\documentclass[10pt,reqno]{amsart}

\usepackage[margin=1in]{geometry}  
\usepackage[utf8]{inputenc}  
\usepackage{fontawesome}
\usepackage{listings}
\usepackage{cancel}
\usepackage{soul}

\usepackage{amsmath,amsthm,amsfonts,amssymb,amscd}
\usepackage{enumitem}

\usepackage{times}

\usepackage{esint,stackrel}
\usepackage{enumitem}
\usepackage{bbm}

\usepackage{float}
\usepackage{emptypage}
\usepackage{xfrac}

\usepackage{stmaryrd}
\usepackage{cases}
\usepackage[mathscr]{euscript}

\usepackage{upgreek}

\usepackage{calligra}

\DeclareMathAlphabet{\mathcalligra}{T1}{calligra}{m}{n}
\DeclareFontShape{T1}{calligra}{m}{n}{<->s*[2.2]callig15}{}

\usepackage{recycle}

\newcommand\ringring[1]{%
  {
   \mathop{\kern0pt #1}\limits^{
     \vbox to-1.85ex{
       \kern-2ex 
       \hbox to 0pt{\hss\normalfont\kern.1em \r{}\kern-.45em \r{}\hss}%
       \vss 
     }
   }
  }
}

\def\p{\partial}
\def\eps{\varepsilon}
\def\les{\lesssim}

\newcommand{\jap}[1]{\langle #1 \rangle} 
\newcommand{\brak}[1]{\langle #1 \rangle} 
\newcommand{\norm}[1]{\left \| #1 \right\|} 
\newcommand{\snorm}[1]{\bigl\| #1 \bigr\|} 
\newcommand{\abs}[1]{\left|#1\right|}
\newcommand{\sabs}[1]{\bigl|#1\bigr|}

\newcommand*{\supp}{\ensuremath{\mathrm{supp\,}}}
\newcommand*{\sgn}{\ensuremath{\mathrm{sgn\,}}}

\newcommand{\TT}{\mathbb T}
\newcommand{\OO}{\mathcal O}

\newcommand{\ZZ}{\mathcal Z}

\renewcommand*{\tilde}{\widetilde}
\renewcommand*{\hat}{\widehat}
\renewcommand*{\bar}{\overline}

\renewcommand{\epsilon}{\varepsilon}

\newcommand{\Bsix}{\langle \mathsf{B}_6 \rangle}

\renewcommand{\theenumi}{(\roman{enumi})}
\renewcommand{\labelenumi}{\theenumi}

\newtheorem{theorem}{Theorem}[section]
\newtheorem{lemma}[theorem]{Lemma}
\newtheorem{proposition}[theorem]{Proposition}
\newtheorem{corollary}[theorem]{Corollary}
\newtheorem{definition}[theorem]{Definition}

\newtheorem{example}[theorem]{Example} 
 
\newtheorem{remark}[theorem]{Remark}

\numberwithin{equation}{section}

\newcommand{\bubu}[1]{#1}

\newcommand{\Cdata}{\mathsf{C_{data}}}
\newcommand{\Cdatatwo}{\mathsf{C}^2_{\mathsf{data}}}
\newcommand{\Csupp}{\mathsf{C_{supp}}}
\newcommand{\Cn}{\mathring{C}}

\def\jump#1{{[\hspace{-2pt}[#1]\hspace{-2pt}]}}

\def\doublecom#1{{( \hspace{-2.5pt}(#1)\hspace{-2.5pt}) }}
\def\comm#1{{( \hspace{-2.5pt}(#1)\hspace{-2.5pt}) }}

\def\n{{\rm n}}

\def\tt{ \mathsf{t} }

\def\U{\mathsf{U}}

\def\Uik{\mathring{\mathsf{U}}^i_k}
\def\Uij{\mathring{\mathsf{U}}^i_j}
\def\Ujk{\mathring{\mathsf{U}}^j_k}
\def\Uii{\mathring{\mathsf{U}}^i_i}

\def\S{\Sigma }
\def\Sk{\mathring\Sigma_k}
\def\Sj{\mathring\Sigma_j}
\def\Si{\mathring\Sigma_i}

\def\Fg{\mathcal{J}} 
\def\Fgss{\Fg^{\! {\frac{3}{2}} } }
\def\Fgtf{\Fg^{\! {\frac{3}{4}} } }
\def\Fgh{\Fg^{\! {\frac{1}{2}} } }

\def\Jg{{J\!_{\scriptscriptstyle g}}}
\def\Jgb{\bar{J}\!_{\scriptscriptstyle g}}

\def\Jgi{\Jg^{\!\!-1} }
\def\Jgt{\Jg^{\!\!2} }

\def\Jgh{\Jg^{\!\! {\frac{1}{2}} } }

\def\Jgtf{\Jg^{\!\! {\frac{3}{4}} } }

\def\mJg{\mathring{\Jg}}
\def\mmJg{\ringring{\Jg}}

\def\JJof{\JJ^{\! {\frac{1}{4}} } }
\def\JJmof{\JJ^{ -\!{\frac{1}{4}} } }
\def\JJhi{\JJ^{-\! {\frac{1}{2}} } }
\def\JJtf{\JJ^{\! {\frac{3}{4}} } }
\def\JJff{\JJ^{\! {\frac{5}{4}} } }
\def\JJss{\JJ^{\! {\frac{3}{2}} } }
\def\JJh{\JJ^{\! {\frac{1}{2}} } }
\def\JJmh{\JJ^{ -\!{\frac{1}{2}} } }

\def\WW{\mathring{\mathsf{W}}}

\def\Wb{{\boldsymbol{ \WW}}}

\def\ZZ{\mathring{\mathsf{Z}}}
\def\Zb{{\boldsymbol{ \ZZ}}}
\def\AA{\mathring{\mathsf{A}}}
\def\Ab{{\boldsymbol{ \AA}}}
\def\Ub{ {\boldsymbol{ \mathring{\U}}}}

\def\Sb{{\boldsymbol{ \mathring{\S}}}}

\def\Wk{\WW_k}
\def\Wi{\WW_i}

\def\Zk{\ZZ_k}
\def\Zi{\ZZ_i}

\def\Ak{\AA_k}
\def\Ai{\AA_i}

\def\nb{\mathsf{D}}
\def\nbs{\tilde{\mathsf{D}}}

\def\nbd{\tilde{\mathfrak{D}}}

\def\cdo{ \! \cdot \!}
\def\cir{ \! \circ \!}

\def\FF{\mathsf{F}}

\def\RR{\mathsf{R}}
\def\II{\mathsf{I}}

\def\Fw{{\FF\!_{\scriptscriptstyle \WW}}}
\def\Fz{{\FF\!_{\scriptscriptstyle \ZZ}}}
\def\Fa{{\FF\!_{\scriptscriptstyle \AA}}}

\def\Rj{{\RR_{\scriptscriptstyle \Jg}}}
\def\Rh{{\RR_{\scriptscriptstyle  \nbs_2 h}}}

\def\Iw{{\II_{\scriptscriptstyle \WW}}}
 
\def\initial{{\mathsf{t_{in}}}}
\def\final{{\mathsf{t_{fin}}}}
\def\medium{{\mathsf{t_{med}}}}

\def\tt{ {\scriptstyle {\mathcal{T}}}} 
\def\ttt{ {\scriptscriptstyle {\mathcal{T}}}} 
\def\nn{ {\scriptstyle {\mathcal{N}}}} 
\def\nnn{ {\scriptscriptstyle {\mathcal{N}}}}

\def\Wbn{{{\boldsymbol{\WW}}_{ \nnn}}}

\def\Wbt{{{\boldsymbol{\WW}}}_{ \ttt}}
\def\Zbn{{\boldsymbol{\ZZ}}_{ \nnn}}

\def\Zbt{{\boldsymbol{\ZZ}}_{ \ttt}}
\def\Abn{{\boldsymbol{\AA}}_{ \nnn}}
\def\Abt{{\boldsymbol{\AA}}_{ \ttt}}
\def\Sbn{{\boldsymbol{\mathring{\S}}}_{\nnn}}
\def\Sbt {{\boldsymbol{\mathring{\S}}}_{\ttt}}

\def\Wbt{{{\boldsymbol{\WW}}}_{ \ttt}}

\def\Zbn{{\boldsymbol{\ZZ}}_{ \nnn}}

\def\Zbt{{\boldsymbol{\ZZ}}_{ \ttt}}

\def\Abn{{\boldsymbol{\AA}}_{ \nnn}}

\def\Abt{{\boldsymbol{\AA}}_{ \ttt}}

\def\Sbn{{\boldsymbol{\mathring{\S}}}_{\nnn}}
\def\Sbt {{\boldsymbol{\mathring{\S}}}_{\ttt}}

\def\Ewn{{{\mathsf{E}}_{\scriptscriptstyle \WW} ^{\nnn}}}

\def\Fwn{{\FF\!_{\scriptscriptstyle \WW} ^{\, \nnn}}}
\def\Fwt{{\FF\!_{\scriptscriptstyle \WW} ^{\, \ttt}}}

\def\Fzn{{\FF\!_{\scriptscriptstyle \ZZ} ^{\, \nnn}}}
\def\Fzt{{\FF\!_{\scriptscriptstyle \ZZ} ^{\, \ttt}}}

\def\Fan{{\FF\!_{\scriptscriptstyle \AA} ^{\, \nnn}}}
\def\Fat{{\FF\!_{\scriptscriptstyle \AA} ^{\, \ttt}}}

\def\Iwn#1{{\Iw}^{\!\!\!\!\nnn,#1}}

\def\Rhn#1{{\Rh}^{\!\!\!\!\!\!\!#1}\, \,}

\def\JJ{ {\mathscr J}}

\def\cir{ \! \circ \!}
\def\tint{\int_0^\s \!\!\!\int \!}

\def\s{{\mathsf{s}}}
\def\Q{\mathsf{Q}}
\def\Qb{ {\bar{\mathsf{Q}} }}
\def\Qc{ {\check{\mathsf{Q}} }}
\def\Qd{ {\hat{\mathsf{Q}} }}
\def\Qr{ {\mathring{\mathsf{Q}} }}

\def\jb{\jmath_{\!{\scriptscriptstyle \beta}}}
 
\def\mPds{ \mathcal{P}^\sharp }
\def\Pdp{{ \mathcal{P}^\sharp_+ }}
\def\Pdm{  { \mathcal{P}^\sharp_-  } }
\def\Pdpm{  { \mathcal{P}^\sharp_\pm } }
\def\Xdp{  { \mathcal{X}_+^\sharp } }
\def\Xdm{  { \mathcal{X}_-^\sharp  } }
\def\Xdpm{  { \mathcal{X}_\pm^\sharp } }
\def\Xds{  { \mathcal{X}_\pm^\sharp(\s) } }
\def\qds{\mathsf{q}}
\def\xstar{ \mathring{x}_1}
\def\xstart{ \mathring{x}_1(x_2)}
\def\xring{ \mathring{x}_1}
\def\xringt{ \mathring{x}_1(x_2)}
\def\tx{ \tilde{x}^*_1}

\def\dl{ {\updelta} }
\def\thd{ {\theta^{\!\dl}} }
\def\Thd{ {\Theta^{\dl}} }
\def\Hdm{ \mathcal{H}^\dl }
\def\Hdmp{ \mathcal{H}^\dl_{+} }
\def\Hdmm{  \mathcal{H}^\dl_{-} }
\def\tHdm{ \mathring{\mathcal{H}}^\dl }
\def\tHdmp{ \mathring{\mathcal{H}}^\dl_{+} }
\def\tHdmm{ \mathring{\mathcal{H}}^\dl_{-} }
\def\sin{{\mathsf{s_{in}}}}
\def\sfin{{\mathsf{s_{fin}}}}

\newcommand{\MGHDB}{\fbox{\rm MGHD}}
 
\usepackage[colorlinks=true, pdfstartview=FitV, linkcolor=blue,citecolor=blue, urlcolor=blue]{hyperref}
 
\title[Geometry of maximal development and shock formation for Euler]{The geometry of maximal development and shock formation for the Euler equations in multiple space dimensions}

\author{Steve  Shkoller}
\address{Department of Mathematics, University of California Davis, Davis, CA 95616.}
\email{\href{shkoller@math.ucdavis.edu}{shkoller@math.ucdavis.edu}}

\author{Vlad Vicol}
\address{Courant Institute of Mathematical Sciences, New York University, New York, NY 10012.}
\email{\href{vicol@cims.nyu.edu}{vicol@cims.nyu.edu}}

\begin{document}

\begin{abstract}
We construct a fundamental piece of the boundary  of the maximal globally hyperbolic development (MGHD) of Cauchy data for the multi-dimensional compressible Euler equations, which is necessary for the local shock development problem.  
For an open set of compressive and generic $H^7$ initial data, we construct unique $H^7$ solutions to the Euler equations in the maximal spacetime 
region below a given time-slice, beyond the time of the first singularity; at any point in this spacetime, the solution can be smoothly and uniquely computed by 
tracing both the fast and slow acoustic characteristic surfaces backward-in-time, until reaching the Cauchy data prescribed along the initial time-slice. 
The future temporal boundary of this spacetime region is a singular hypersurface, containing the union of three sets: first, a co-dimension-$2$ surface
of ``first  singularities'' called the {\em pre-shock}; second, a downstream hypersurface called the {\em singular set} emanating from the pre-shock, 
on which the Euler solution experiences a {\em continuum of gradient catastrophes}; third, an upstream hypersurface consisting of a  
{\em Cauchy horizon} emanating from the pre-shock, which the Euler solution cannot reach. We develop a new geometric framework for the
description of the acoustic characteristic surfaces which is based on the Arbitrary Lagrangian Eulerian (ALE) framework, and combine this with a
new type of differentiated Riemann variables which are linear combinations of gradients of velocity, sound speed,  and the curvature of the fast acoustic characteristic surfaces.  With these new variables, we establish uniform $H^7$ Sobolev bounds for solutions to the Euler equations without derivative loss and with optimal regularity. 
\end{abstract}

\maketitle

\vspace{-0.25in}

\setcounter{tocdepth}{1}
\tableofcontents

\allowdisplaybreaks

\vspace{-0.25in}
\section{Introduction}
We construct a portion of the boundary of the {\it maximal  globally hyperbolic development}  (MGHD) of Cauchy data, during the shock formation process, for solutions of the Euler equations
\begin{subequations} 
\label{euler0}
\begin{align}
\partial_t (\rho u) + \operatorname{div}  (\rho u \otimes u) + \nabla p &= 0 \,,  \\
\partial_t \rho  +  \operatorname{div}  (\rho u)&=0 \,,  \\
\partial_t E +  \operatorname{div}  (u(E+p))&=0 \,,
\end{align}
\end{subequations} 
where
$p = (\gamma -1)  (E- {\frac{1}{2}} \rho |u|^2 )$ is the scalar pressure function, $\gamma>1$ is the adiabatic exponent, 
$u :\mathbb{T}^d  \times \mathbb{R}  \to \mathbb{R}^d$ denotes the $d$-component velocity vector field, 
$\rho: \mathbb{T}^d  \times \mathbb{R}  \to \mathbb{R}  _+$ denotes the strictly positive density function, and
$E: \mathbb{T}^d  \times \mathbb{R}  \to \mathbb{R}$ is the total energy.   
In particular, we develop a new geometric and analytic framework that allows us to obtain uniform Sobolev estimates for the Euler solution, 
evolving past the time of the first gradient singularity uniformly through a spacetime hypersurface of gradient catastrophes, also
known as the {\it singular set}. 
In turn, we are able to give a complete description of a fundamental portion of the boundary of the MGHD
to which (an open set of) compressive Cauchy data can be smoothly and 
uniquely evolved.   

Specifically, we consider a portion of the boundary of the MGHD which contains (a) the set of first singularities, which we call the 
{\it pre-shock set}\footnote{The term {\it pre-shock} has its origins in the work of Fournier \& Frisch~\cite{FoFr1983}, who use the term {\em pr\'echoc}.}, 
(b) the  {\it singular set} which is the downstream hypersurface of gradient catastrophes emanating from the pre-shock set, and (c) the upstream {\it Cauchy horizon} which is the slow acoustic characteristic hypersurface emanating from the pre-shock set.
  Indeed, it is this portion of the boundary of
the MGHD that is essential for the resolution of the shock development problem, which we shall describe below.  In this regard, the results herein
resolve the first step in a two-tier program to establish the existence of unique shock wave solutions for the Euler equations
in multiple space dimensions.

An abbreviated form of our main result, the establishment of a fundamental portion of the MGHD,   can be found in Theorem \ref{thm:mainrough} below, while the detailed statements can be found in Theorems \ref{thm:main:shock}, \ref{thm:main:DS}, and \ref{thm:main:US}.

\subsection{Shock formation and shock development}
\label{sec:two:step}
The system \eqref{euler0} is the quintessential system of nonlinear
hyperbolic conservation laws.   Such systems exhibit shocks;   these are spacetime hypersurfaces of discontinuity
which  emerge in finite time from smooth initial data, and dynamically evolve according to the Rankine-Hugoniot (RH) jump conditions (see, for example,
\cite{Dafermos2016}).   In addition
to the physical variable unknowns in \eqref{euler0} -- velocity, density, and energy -- the location of the shock surface is also an unknown.  A weak solution to the
Euler equations requires that the physical variables satisfy the Euler  equations pointwise on either side of the shock surface, and that the shock
surface propagates with the correct  normal speed.  Moreover, certain physical ``entropy'' conditions must be satisfied to ensure that the solution is
physically meaningful.

While the theory of shock waves and weak solutions to the compressible Euler equations (and more general systems of conservation laws) is fairly 
complete in one space dimension (see \cite{Ri1860}, \cite{Lax1957,Lax1964,Lax1972}, \cite{Glimm1965}, \cite{GlLa1970}, \cite{Li1979}, \cite{Di1983}, 
\cite{Ma1984}, \cite{Lebaud1994}, \cite{Chen1997}, \cite{Yin2004},  \cite{BiBr2005} as well as the fairly complete bibliography in \cite{Dafermos2016}), the problem of 
obtaining and propagating unique shock solutions in more than one space dimension, without symmetry assumptions, remains open.    Detailed shock
formation   under azimuthal symmetry  with the functional description of the solution at the first shock singularity has been extensively studied in
\cite{BuDrShVi2022}, \cite{NeShVi2022} and \cite{NeRiShVi2023}.

The methods of one-dimensional conservation laws have not proven to be easily extendable to multiple space dimensions (see, for example,~\cite{Ra1986}). Moreover, the convex-integration based results originating in~\cite{ChDeKr2015} and refined in~\cite{V-ASS}, have shown that entropy inequalities cannot be used as a uniqueness selection criterion.   As a result, there is yet no general existence theorem, 
describing the  evolution of smooth data towards the creation and  unique propagation of discontinuous shock surfaces for the Euler equations in multiple space dimensions.

Christodoulou \cite{Ch2007,Ch2019}  introduced a novel two-stage
 program for the construction of unique shock wave solutions to the Euler
equations.  Starting from smooth initial data, the first step is called {\em shock formation}, in which  smooth compressive initial data is
evolved up to a cusp-like Eulerian spacetime co-dimension-$2$ hypersurface of ``first singularities'' -- these first singularities are where the gradient of velocity, density, and energy first become infinite.    We term this co-dimension-$2$ hypersurface of ``first singularities''  the {\it pre-shock} set, because
along this set, the Euler solution remains continuous, but forms a $C^{\frac{1}{3}} $ cusp.  The term pre-shock is used, because on this set, the
gradient of the solution has become infinite, but the actual shock {\it discontinuity} is yet to develop.
The second step of the program is called {\it shock development}. Here, one uses the analytical description of the $C^{\frac{1}{3}}$ solution in a 
neighborhood of
the pre-shock as Cauchy data, from which the shock surface of discontinuity instantaneously {\it develops}.   To date, this program remains unresolved in the absence of symmetry assumptions; for the Euler equations, see Christodoulou~\cite{Ch2019} for the so-called restricted shock development problem, Yin \cite{Yin2004}
and  Christodoulou \& Lisibach~\cite{ChLi2016} for shock development in spherical symmetry, and~\cite{BuDrShVi2022} for shock development in
azimuthal symmetry, 
together with the emergence of the {\it weak characteristic discontinuities} conjectured by Landau \& Lifshitz \cite{LaLi1987}.   It is important to
note that Majda's shock stability result \cite{Ma1983, Ma1983b} is neither a {\it shock formation} result nor a {\it shock development} result, but rather
a short-time existence theorem on the propagation of  shock front initial data by the shock speed imposed by the Rankine-Hugoniot jump conditions.
Specifically, Majda assumes the existence of a surface of discontinuity in the data, while the objective of shock development is to dynamically create
this surface of discontinuity from the $C^{\frac{1}{3}} $ cusp-solution at the  pre-shock.
 
To summarize,  the first step of this two-tiered program  necessitates the analysis of a portion of the boundary of the  maximal globally hyperbolic development (MGHD) of smooth compressive Cauchy data. The second step of the program, shock development, then uses the $C^{\frac{1}{3}} $ Euler solution about the pre-shock as Cauchy data, and seeks to evolve a unique shock hypersurface dynamically and instantaneously from the pre-shock, together with a unique weak solution obeying the RH jump conditions across the shock front.  In turn, this requires having unique Euler solutions up to the portion of the boundary of the MGHD containing the pre-shock, the singular set, and the Cauchy horizon. The main objective of this paper is the resolution of the first step of this program, giving the relevant portion of the  boundary of the MGHD.

 \subsection{The evolution of the Euler solution past the time of first blow-up}
With the exception of the recent result of Abbrescia \& Speck \cite{AbSp2022} which we shall describe below,  the analysis of the multi-dimensional shock formation process for the Euler equations considered solutions {\it only up to}  the {\it time of the very first singularity}, the earliest time $t_*$  at  which the solutions' gradient  becomes  infinite (see~\cite{Si1985}, \cite{Ch2007}, \cite{ChMi2014}, \cite{LuSp2018}, 
\cite{LuSp2021}, \cite{BuShVi2022,BuShVi2023a,BuShVi2023b}).
Such an analysis is insufficient to proceed with the problem of shock development.  It is essential to describe the full shock formation process, past the time of the first singularity, and to capture the pre-shock set of ``first singularities'' which successively emerges.   Indeed, it is the description of the solution about the pre-shock
set of  first singularities (see the black curve on both the left and right images in Figure~\ref{fig1} below) that is used as Cauchy data for the development of discontinuous shock waves.
The objective is therefore to create a novel geometric and analytical framework that can provide uniform Sobolev estimates for solutions that are
experiencing successive gradient catastrophes along a hypersurface of spacetime.  Thus, we are not simply trying to  extend the solution past a
single first singularity, but rather  we are evolving the solution through a continuum of gradient blow-ups in appropriately chosen coordinates that allow
for uniform bounds to be maintained.

We  note that while our focus in this work is on the shock-type gradient singularity, solutions to the Euler equations can form finite-time implosions from smooth initial conditions.  Such unstable imploding solutions, which form a finite-time amplitude blow-up, have been proven to exist in~\cite{MeRaRoSz2022a,MeRaRoSz2022b,BuCaGo2023,CaGoShSt2023}.   It was also recently shown in \cite{ChCiShVi2024} that an
axi-symmetric implosion can lead to a finite-time blow-up of the vorticity.

\subsection{Maximal globally hyperbolic development in a box, \MGHDB}
\begin{figure}[htb!]
\centering
\begin{minipage}{.35\linewidth}
  \includegraphics[width=\linewidth]{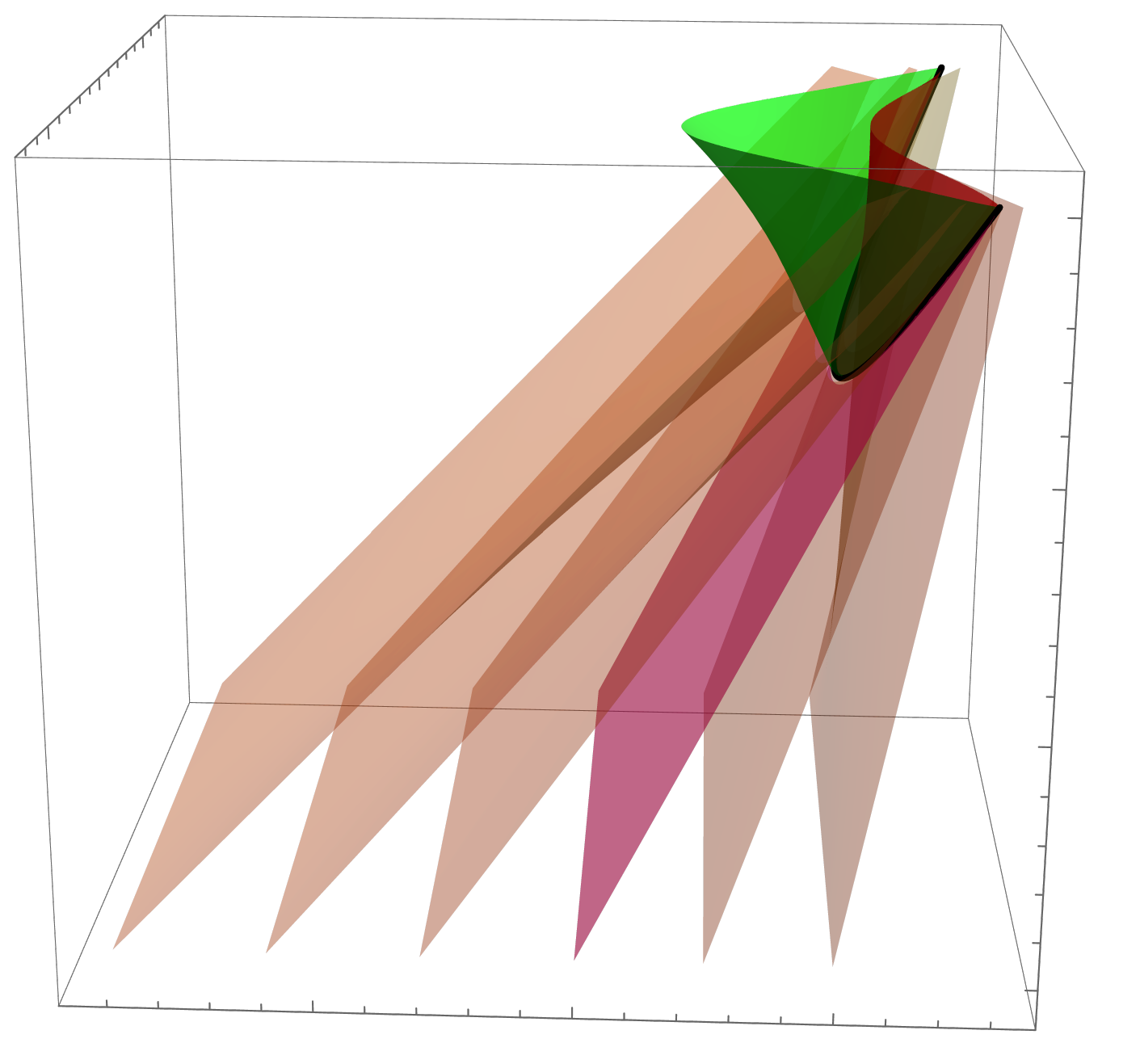}
\end{minipage}
\hspace{.05\linewidth}
\begin{minipage}{.45\linewidth}
  \includegraphics[width=\linewidth]{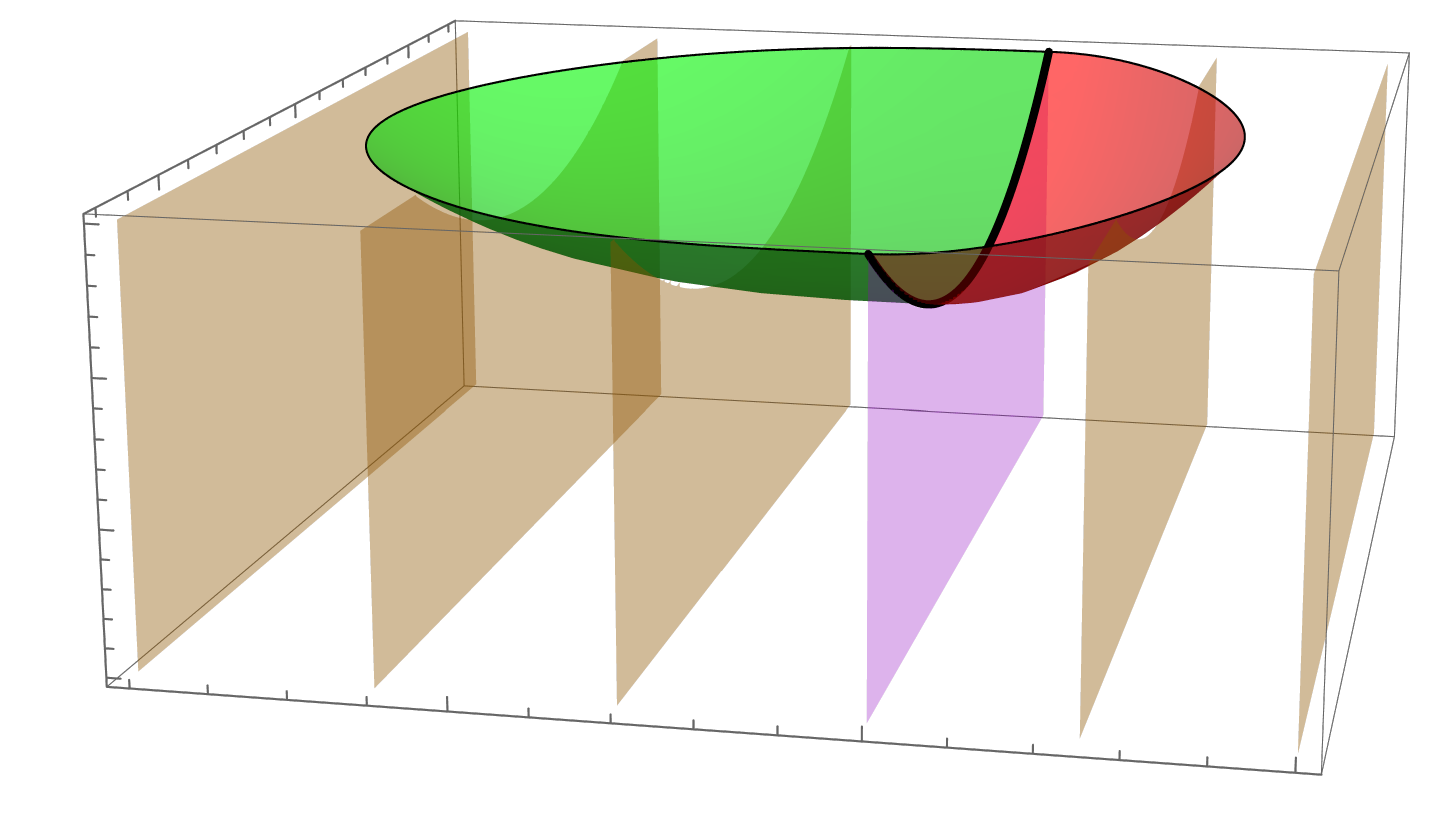}
\end{minipage}
\vspace{-.1 in}
\caption{ \underline{Left}.  \MGHDB\ spacetime in traditional Eulerian coordinates.  The characteristic surfaces are shown to impinge on the  cuspoidal surface of first singularities, which consists of the pre-shock set (shown as the black curve) together with the singular hypersurface
which emanates from the pre-shock in the downstream direction, which consists of a continuum of gradient catastrophes (shown as the red surface).   The slow
acoustic characteristic which emanates from the pre-shock in the upstream direction (shown in green) is a Cauchy horizon which the Euler
solution can never reach. 
\underline{Right}. 
 The  \MGHDB\  spacetime in Arbitrary Lagrangian Eulerian (ALE) coordinates denotes the region below the union of the downstream ALE paraboloid of ``gradient catastrophes'' (in red),  the set of pre-shocks (in black), and the slow
acoustic characteristic surface which emanates from the pre-shock (in green).
In ALE coordinates, each fast acoustic characteristic surface has been flattened.   The hypersurface shown in the magenta color denotes the
fast acoustic characteristic surface that passes through the pre-shock.
}
 \label{fig1}
\end{figure}

In the traditional Cauchy problem in fluid dynamics, initial data is prescribed on a given time-slice 
and the evolution of this data is considered up to some later time-slice. 
For PDEs with finite speed of propagation, localized sets of initial data can be propagated up to space-like submanifolds of spacetime which
do not necessarily coincide with time-slices.  For example, as shown on the left side of Figure~\ref{fig1}, which displays an example of an Eulerian spacetime, the fast acoustic characteristic surfaces (emanating from the initial time-slice) are impinging upon each other; the black curve represents the location in spacetime where the first gradient catastrophes occur, while the union of the red and green surfaces indicate the future temporal boundary of the collection of points that can be smoothly and uniquely reached by the Cauchy data.   The right panel in Figure~\ref{fig1} displays the analogous spacetime set, but in {\it Arbitrary Lagrangian Eulerian} (ALE) coordinates adapted to the fast acoustic characteristic surfaces.  These ALE coordinates, which will be defined in Section \ref{sec:ALE}, provide a smooth geometric framework for our analysis.

There is a notion of
the  {\it maximal globally hyperbolic development} (MGHD)  of a data set, which originated in the study of  general relativity  and can be traced back to the fundamental paper of Choquet-Bruhat \& Geroch \cite{Ch-BrGe1969} (see also the recent discussion in \cite{EpReSb2019}).  For a hyperbolic PDE, the Cauchy data  is prescribed on a spacelike manifold $\mathcal{S}_0$. 
 A {\it development} of this data consists of a spacetime $ \mathcal{M} $, a solution  of
the hyperbolic PDE, together with a diffeomorphism of $\mathcal{S}_0$ onto a spacelike submanifold $ \mathcal{S} $ of  $ \mathcal{M} $
such that the solution restricted to $ \mathcal{S} $ is diffeomorphic to the data on $ \mathcal{S} _0$. In other words, each point on the manifold $ \mathcal{M} $ can be reached by a unique and smooth characteristic emanating from $ \mathcal{S}$.   There exists a  precise notion of a {\it maximal} development. It is known that every initial data set has a hyperbolic development,  and that any two developments of $ \mathcal{S}_0$ are extensions of a common development.  Moreover, for any initial data set on $\mathcal{S}_0$, there exists a development $ \mathcal{M} $ of $ \mathcal{S}_0$ which is an extension of every other development of $\mathcal{S}_0$, and  this development is unique. 

Smooth solutions of the Euler equations generically develop discontinuous shocks.
The construction of such shock-wave solutions to the Euler equations requires more than mere sound-wave analysis in the MGHD: 
a nonlinear free-boundary problem for the motion for the discontinuous shock surface  is to be studied via the Rankine-Hugoniot jump conditions and the second law of thermodynamics.  
Nevertheless, the construction of such physical solutions, even {\em locally-in-time}, necessitates a detailed description of the Euler solution on a {\em local piece} of the boundary of the MGHD; this provides data on the union of the Cauchy horizon and the a priori unknown shock hypersurface, 
 which is itself tangent to the singular set at the pre-shock (see the discussion in Section~\ref{sec:two:step}). With this goal in mind, we introduce a {\em compact-in-time}  version of the maximal globally hyperbolic development (MGHD) of the Cauchy data, which is fundamental to the local-in-time analysis of physical shock-waves, and which we denote as   {\em MGHD in a box}.

The concept of a {\em MGHD  in a box}, written succinctly as \MGHDB, is defined as follows. For initial data prescribed on a time-slice 
$\mathcal{S}_0 = \TT^d \times \{t=\initial\}$ (equivalently, ${\mathbb R}^d\times \{t=\initial\}$), let $t_*$ denote the time at which the first gradient 
singularity occurs. Let $\final>t_*$ be a finite time horizon up to which we wish to evolve the Euler solution, where $\final - t_* = \OO(t_* - \initial)$. 
The \MGHDB\ spacetime is then defined as the intersection of the classical maximal globally hyperbolic development $\mathcal M$ of ${\mathcal S}_0$, with the 
spacetime box $\TT^d \times [\initial,\final]$ (equivalently, the spacetime slab ${\mathbb R}^d\times \{t=\initial\}$). As such, the future temporal 
boundary of the \MGHDB\  captures the portion of the boundary of the classical MGHD which is relevant to the physical shock-wave Euler evolution 
within the compact time interval $[\initial,\final]$.

In traditional Eulerian coordinates, the portion of the boundary of the MGHD for~\eqref{euler0} captured by the \MGHDB\ is the cuspoidal surface in Figure~\ref{fig1} below (see the left panel which displays  the union of the red singular surface, the black pre-shock set, and the green 
Cauchy horizon).  This Eulerian cuspoidal surface acts as a temporal boundary to the spacetime set on which the Cauchy data can be evolved in a 
smooth and unique fashion. As noted above, the black curve, which we will refer to as the  {\it pre-shock}, denotes the set of ``first singularities'' where 
characteristic surfaces first impinge. This {\it pre-shock} set is a collection of spacetime points where the gradient of the solution (velocity, density, 
pressure, etc.) experiences the first  blow-up. The red surface shown in Figure~\ref{fig1} denotes the spacetime hypersurface  on which the fast 
characteristic surfaces impinge (the {\em singular set}), and the green surface shown in Figure~\ref{fig1} denotes the slow acoustic characteristic 
surface emanating from the pre-shock set (the {\em Cauchy horizon}). The \MGHDB\ of $\mathcal{S} _0 = \TT^2 \times \{t=\initial\}$ thus 
consists of the Euler solution in the spacetime set lying ``below'' the  union of the aforementioned cuspoidal surface  in Figure~\ref{fig1} and  the ``lid'' 
at $\{t=\final\}$. For a concise and self-contained introduction to the \MGHDB\ of Cauchy data in gas dynamics, we refer the reader to 
Appendix~\ref{sec:usersguide}, where this concept is discussed in the familiar context of the $1$D Euler equations.

\subsection{Prior results on solutions of the Euler equations past the time of first singularity}
 
\subsubsection{Christodoulou's method \cite{Ch2007}}  
The velocity potential $\varphi$ of an irrotational compressible fluid satisfies the wave
equation $\Box_ g \varphi=0$ (see Synge~\cite{Synge1938}, Unruh~\cite{Unruh1981}, and Visser~\cite{Visser1993});
that is, the irrotational Euler equations are written in terms of the d'Alembertian associated to the Lorentz (acoustic) metric $g$, where
the components of $g^{-1}$ are given by $g^{00}=-1$, $g^{0j}=-u^j$, $g^{i0}= - u^i$, $g^{ij} = c^2\delta^{ij}- u^iu^j$, where $i,j=1,...d$, and the index $0$ denotes the temporal component.  It is quite remarkable that even though the underlying fluid dynamics are Newtonian, non-relativistic, and occur in flat spacetime, the fluctuations (sound waves) are governed by a curved Lorentzian (pseudo-Riemannian) spacetime geometry.   With such a construct, sound plays the role of light, and sound rays follow the {\it null geodesics} of the acoustic metric.   As a result, the analysis of the irrotational compressible Euler equations can be entirely relegated to the study of the second-order nonlinear wave equation $\Box_g \varphi=0$, and thereby placed in the framework of analysis in  General Relativity.
 
 This is indeed the starting point for Christodoulou's method \cite{Ch2007} for the analysis of shock formation in fluids.  The null geodesics and hence the  bicharacteristic\footnote{\label{sec:classical-char-surf} In multiple space dimensions, characteristic surfaces are co-dimension-$1$ submanifolds of spacetime.   
In the classical framework, which can be traced back to Courant \& Friedrichs \cite{CoFr1948},  such manifolds are locally generated by special vector fields, whose directions are termed {\it bicharacteristic directions}.  These vector fields are generators of a cone,  the so-called {\it  bicharacteristic cone}, which is tangent to
the characteristic surfaces. This traditional description of characteristic surfaces is very different from our approach, which we detail in Section \ref{sec:freeboundary}.} directions are obtained as solutions to the {\it eikonal equation}
\begin{equation*}
 g^{\alpha\beta} \p_\alpha q \p_\beta q=0 \,, \qquad \alpha,\beta=0,1, \ldots, d \,.
\end{equation*}
By definition, the derivatives of $q$ then determine the tangent and normal directions to the acoustic characteristic surfaces, as well as the
inverse foliation density $\mu$.  Gradient singularities form as $\mu$ vanishes.    With 
geometric variables defined, the closure of $N$th-order Sobolev-class $\mu$-weighted energy estimates on  second-order wave equations provides
lower-order bounds on the irrotational Euler solutions up to the time of first singularity (for $M$ derivatives, where $M$ is much less than $N$).

Christodoulou's method has been used as the framework for the study of shock  formation (and singularity formation) for a  large class of Lorentzian wave equations;  see, for example,  \cite{Sp2016,Sp2018}, \cite{HoKlSpWo2016}, \cite{MiYu2017}, \cite{LuSp2018,LuSp2021}, \cite{LuYu2023a,LuYu2023b}.

\subsubsection{The recent result of Abbrescia \& Speck \cite{AbSp2022}}
Employing the geometric framework of Christodoulou's method, Abbrescia \& Speck \cite{AbSp2022} have constructed a piece of the boundary of the MGHD for the 3D Euler equations with vorticity and entropy,  which contains (a) the {\it pre-shock} (which they refer to as the {\it crease}), and  (b) a localized portion of the singular set; the analysis in   \cite{AbSp2022} does not  contain {\it any} portion of the boundary of the MGHD containing the Cauchy horizon.  This portion of the boundary of the MGHD (as well as the structure of the MGHD) was discussed in Chapter 15 of
Christodoulou's book \cite{Ch2007}.

While Christodoulou's result in~\cite{Ch2007} obtains uniform bounds on the irrotational Euler solution up to the time of first singularity, 
Abbrescia \& Speck~\cite{AbSp2022} are able to go past this first singular time and to generalize Christodoulou's method to allow for non-trivial
vorticity and non-constant entropy.    To go past the time of first singularity, Abbrescia \& Speck~\cite{AbSp2022} use foliations of spacetime rather
than constant time-slices, thereby enabling the closure of energy estimates up the singular set, downstream of the crease.  Additionally, the energy
method is quantified with a lower-bound for the number of derivatives $N=25$ (the top-order energy scheme) required to close estimates
(see also~\cite{LuSp2018,LuSp2021}).   More
importantly, they are able to couple the transport equations for vorticity and entropy to the second-order wave equations of the irrotational and 
isentropic problem.   Because Christodoulou's original scheme relied crucially on  pure-wave structure, the coupling of vorticity transport is an
important  technical advancement in \cite{AbSp2022} and  \cite{LuSp2018,LuSp2021}.  It should be noted, however, that this transport-wave coupling requires a stringent regularity 
assumption on the initial data for vorticity, requiring one extra derivative of regularity for vorticity relative to the regularity assumed for the velocity.  
This incompatibility in the smoothness of vorticity becomes problematic for the shock development problem.

\subsection{Preparing the Euler equations for shock formation analysis}
\subsubsection{A  symmetric form of the Euler equations}
While the Euler solution stays smooth, the energy equation $\partial_t E +  \operatorname{div}  (u(E+p))=0$
can be replaced by the transport equation for the specific entropy $S$, 
 $$\p_t S + u\cdot \nabla S =0\,,$$ 
 in which case the pressure law can be equivalently written as  
 $$p= {\tfrac{1}{\gamma}} e^S\rho^\gamma \,, \ \ \gamma >1 \,.$$
When the {\it initial entropy function}  is a constant, the entropy remains a constant during the shock formation process and the
 dynamics are termed {\it isentropic}; for isentropic dynamics,  the pressure law is given by $p= {\tfrac{1}{\gamma}} \rho^\gamma$ for $\gamma >1$.
 
It is convenient to rewrite the Euler equations in the symmetric form
\begin{subequations} 
\label{euler1}
\begin{align} 
\p_t u +  u \cdot \nabla u  +   \alpha  \sigma \nabla \sigma & = 0 \,, \label{u1} \\
\p_t\sigma +    u \cdot \nabla \sigma  +  \alpha  \sigma  \operatorname{div} u & = 0  \,,  \label{sigma1}
\end{align} 
\end{subequations} 
where $ \alpha = \tfrac{\gamma-1}{2} >0$ is the adiabatic exponent, and $ \sigma = \tfrac{1}{ \alpha } c$ is the rescaled sound speed.  
The sound speed is $c = \rho ^ \alpha $.

\subsubsection{Acoustic form of the Euler equations}
A key feature in the analysis of the Euler equations, and general systems of nonlinear hyperbolic equations, is the concept of {\it characteristic surfaces} and the geometry describing their evolution (see Courant \& Hilbert \cite{CoHi1962}). Characteristic surfaces in spacetime can be viewed as propagating wave fronts, and for
the Euler equations, these represent either sound waves or vorticity waves.

 We shall focus our presentation on isentropic dynamics in space dimension $d=2$.\footnote{For both simplicity and concision, we present our method of analysis of the \MGHDB\ for the case of $d=2$, in which  there are three independent variables $(x_1,x_2,t)$, each characteristic surface is a $2$-dimensional hypersurface of spacetime which possesses one unit normal vector $n$ and one unit tangent vector $\tau$ at each point. The modifications of our theory to the case that $d=3$, merely require the use of two linearly independent tangent vectors $\tau_1$ and $\tau_2$ to each fast acoustic characteristic hypersurface of spacetime, which turns out to be a fairly trivial generalization.}

Let $(n( \cdot , t),\tau(\cdot , t))$ denote an orthonormal basis for $ \mathbb{T}  ^2$ for each time $t$.
The Euler equations \eqref{euler1} have three {\it distinct} waves speeds
\begin{equation} 
\lambda_1  = u \cdot n - \alpha \sigma \,, \qquad
\lambda_2  = u \cdot n \,, \qquad
\lambda_3  = u \cdot n + \alpha \sigma  \,. \label{wave-speeds}
\end{equation}
Here $\lambda_3$ is the fast acoustic wave speed, and it is along the transport velocity $u+ \alpha \sigma n$ that sound waves steepen
to form shocks.   The normal vector $n$ will be made explicit below as the normal to  dynamically evolving and steepening {\it fast acoustic
characteristic surfaces}.   The wave speed $\lambda_1$ denotes the {\it slow acoustic wave speed}, and plays a prominent role in defining the
future temporal boundary for the \MGHDB.

When spacetime is foliated by the fast acoustic characteristic surfaces,  and if $n$ and $\tau$ denote the normal and tangent vectors, respectively, to the intersection
of these surfaces with each time-slice, then the
propagation of acoustic waves (and shock formation) can be studied by  rewriting \eqref{euler1} as 
\begin{subequations} 
\label{euler-shock}
\begin{align} 
\p_t u +  (u + \alpha \sigma n) \cdot \nabla u  +   \alpha  \sigma (\nabla \sigma - n \cdot \nabla u) & = 0 \,, \label{u-shock1} \\
\p_t\sigma +    (u + \alpha \sigma n)  \cdot \nabla \sigma  +  \alpha  \sigma ( \operatorname{div} u - n \cdot \nabla \sigma)& = 0  \,.  \label{sigma-shock1}
\end{align} 
\end{subequations} 
The system \eqref{euler-shock} can be thought of as the {\it acoustic Euler equations} with fast acoustic  transport velocity $u+ \alpha \sigma n$, along
which the fast characteristic surfaces propagate.

\subsubsection{Classical Riemann variables}
The description and dynamics of these characteristic surfaces greatly simplify in  
one space dimension.  In the 1D case, the characteristic surfaces are curves in two-dimensional spacetime that can propagate in only
two directions; namely, the positive or the  negative spatial direction.  Along these characteristic directions, the solutions to the 1D Euler equations 
possess certain invariant functions which are called {\em Riemann invariants}~\cite{Ri1860}, which (in the case of constant entropy) are exactly
transported along the fast and slow acoustic wave speeds.

In multiple space dimensions, complete invariance is generally not preserved, but Riemann  variables  can nevertheless be defined 
to both capture the dominant sound wave motion and to maintain small deviation of the dominant variable when transported along the fast characteristic surfaces.
The classical Riemann variables are defined as\footnote{For the case that $d=3$, $w$ and $z$ are defined in the identical fashion, while
the tangential velocity is given by the vector $a=(a_1, a_2):=(u\cdot\tau_1, u\cdot \tau_2)$ where $(\tau_1,\tau_2)$ form a basis
for each tangent space to the fast acoustic characteristic hypersurface.}
\begin{equation}
w  = u \cdot n  + \sigma\,, \qquad \qquad z  = u\cdot n -\sigma\,, \qquad \qquad a = u \cdot \tau \,.
\label{wza-intro}  
\end{equation}

\subsubsection{Generic and compressive initial data for shock formation} 
We consider an open set of data $(u_0,\sigma_0) \in H^7$ which satisfy certain generic and compressive properties that are made precise 
in Section \ref{cauchydata} below.
We define the dominant Riemann variable at initial time $t=\initial$ by  $w_0(x) := w(x,\initial) = u_0^1(x,\initial) + \sigma(x,\initial)$, where $n(x,\initial)=e_1$.
We set the maximal negative slope of $w_0$ to occur in
the $x_1$ direction.   
We suppose that for $0< \eps \ll 1$,  the derivative  $\p_1w_0$ takes its minimum value at the origin and is given by $\p_1w_0 (0)= - {\tfrac{1}{\eps}} $.
We further suppose that the initial conditions $z_0(x)=z(x,\initial)$ and $a_0(x)=a(x,\initial)$ are small, and have $\OO(1)$ derivatives.   Such data is called 
non-degenerate, or, generic,  if $\nabla^2 \p_1w_0(0)$ is positive definite.

\begin{remark}[Compressive data for which $\p_n w$ blow-up]\label{rem:pn-w-blowup}
For generic and compressive initial data (described above) which yield sound waves that  steepen in the primary direction of propagation $x_1$, and 
evolve with relatively small changes in the transverse coordinate $x_2$, 
$w$ is the {\it dominant} Riemann variable and $z$ is the {\it subdominant} Riemann variable. 
During the shock formation process, it is the normal derivative $\p_n w$ that blows-up, while $\p_n z$, $\p_n a$, as well as $\p_\tau w$, $\p_\tau z$,
and $\p_\tau a$ all remain bounded.
\end{remark}

\subsection{The geometry and regularity of the fast acoustic characteristic surfaces}\label{sec:freeboundary}
The explicit presence of the unit normal vector $n$ (to the steepening fast characteristic surfaces) in 
the system of equations \eqref{euler-shock}  shows the geometric nature of these equations, when written in the form suitable for shock formation.
Together with $u$ and $\sigma$, the normal vector $n$ is one 
of the fundamental {\it unknowns} in the dynamics of shock formation.   The physical unknowns $u$ and $\sigma$ are directly coupled to the evolution
of the geometry of the fast acoustic characteristic surface, and as such, 
 the dynamics of shock formation can be thought of as a highly nonlinear example of a
{\it free boundary problem in fluid dynamics}.\footnote{Our analysis of characteristic surfaces and their underlying geometry is quite different from the traditional viewpoint
of solving the Eikonal equation to determine the bicharacteristic cone.  Instead, we place the geometry of characteristic surfaces in the context of
free boundary problems in fluid dynamics. See also footnote~\ref{sec:classical-char-surf}.}
Traditional free boundary problems have the location (or shape) of the fluid
boundary as one of the basic unknowns, and the dynamics of this free boundary must be coupled to the evolution of the physical flow fields.   For 
both incompressible \cite{CoSh2007} and compressible \cite{CoSh2012} free boundary problems, the dynamics of the geometry are governed by the 
fluid velocity $u$, corresponding to the wave speed $\lambda_2$ in \eqref{wave-speeds}, while for the shock formation problem it is the location and 
shape of the sound wave which is the geometric unknown, and this wave pattern is carried by the fast wave speed $\lambda_3$.  This creates a 
new level of difficulty for the analysis.

To place this difficulty in perspective, we shall consider the Lagrangian parameterization of the geometry.   For traditional  free boundary problems
in fluid dynamics, one considers the Lagrangian flow $\eta(x,t)$ associated to the fluid velocity; namely, 
\begin{equation}
\p_t\eta(x,t) = u(\eta(x,t),t) \ \text{ for } \ t> \initial \,, 
\qquad 
\eta(x,\initial)=x \,,
 \label{trad-flow}
\end{equation} 
where $\initial$ denotes the initial time for the flow.   A key observation is that \eqref{trad-flow} defines an ordinary differential equation, and Picard 
iteration shows that for (at least Lipschitz) velocity fields $u$, the flow map $\eta$ inherits (at least) the regularity of the velocity field, and can often
be shown to gain regularity such that the associated normal vector $n$ possesses the same regularity as the velocity field.  This is indeed the case
for the classical incompressible and compressible Euler free boundary problem.  We now contrast this scenario with the geometric dynamics of shock formation.

Observe from \eqref{euler-shock} that the transport velocity associated with the fast wave speed $\lambda_3$ is given by
\begin{equation} 
\mathcal{V}_3 = u + \alpha \sigma n =  \lambda_3 n + (u \cdot \tau) \tau \,.  
\end{equation} 
We can foliate spacetime with acoustic characteristic surfaces associated to the ``fast'' wave speed $\lambda_3$.
One way of doing so is by studying the Lagrangian flow map of the transport velocity $\mathcal{V} _3$:
\begin{equation}
\p_t \eta(x_1,x_2,t)  = \mathcal{V} _3(\eta(x_1,x_2,t),t) \,, 
\qquad
\eta( x , \initial) = x  \,.
\label{eta-char-intro}
\end{equation} 
In terms of the standard Cartesian basis, we have that\footnote{We identity vector fields and $1$-forms in Euclidean space; in particular, raised indices for components
$F^i$ of a vector field are obtained as $F^i= \delta^{ij} \tilde F_j$, where $\tilde F_j$ denotes the components of the $1$-form and $\delta^{ij}$ denotes the Kronecker-$\delta$.}
\begin{equation*}
\eta(x_1,x_2,t) = (\eta^1(x_1,x_2,t),\eta^2(x_1,x_2,t))\,,
\end{equation*}
and that
\begin{equation*}
\eta^1(x_1,x_2,\initial) = x_1 \ \ \text{ and } \ \ \eta^2(x_1,x_2,\initial) = x_2 \,.
\end{equation*}
Using the flow map $\eta$ we can  give a geometric description of the fast acoustic characteristics surfaces.

\begin{figure}[htb!]
\centering
  \includegraphics[width=.45\linewidth]{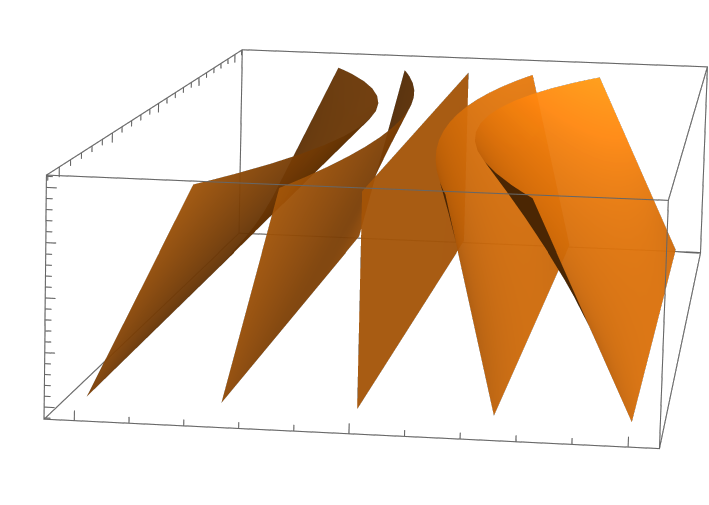}
 \vspace{-.1in}
\caption{For $T>\initial$ which is strictly less than the very first blow-up time, we display the characteristic surfaces $\Gamma_{x_1}(T)$ defined in~\eqref{Cx1-intro} emanating from five different values of $x_1\in\TT$. At $t= \initial$, the curves $\{ \gamma_{x_1}(\initial) \}_{x_1\in\TT}$ are lines which foliate $\TT^2$. The distance between the characteristic surfaces $\Gamma_{x_1}(T)$ is decreasing as $T$ increases, leading to shock formation when this distance vanishes.}
\label{fig:surf-pre}
\end{figure}

At the initial time $t=\initial$, we foliate $\TT^2$ by lines parallel to $e_1=(1,0)$, and denote these lines by
$ \gamma_{x_1}(\initial) = \{x_1\} \times \TT$.   For each $x_1 \in \TT  $ and $t\in[ \initial,T]$, we define the 
characteristic curve (at a fixed time-slice) by
\begin{equation} 
\gamma_{x_1}(t) = \eta( \gamma_{x_1}(\initial), t)  \,,  
\end{equation} 
and the characteristic surfaces up to time $T \ge \initial$, parameterized of $x_1$, by 
\begin{equation} 
\Gamma_{x_1}(T) = \bigcup\nolimits_{t\in[\initial,T]} \gamma_{x_1}(t) \,. \label{Cx1-intro}
\end{equation} 
Figure~\ref{fig:surf-pre}  displays a few such characteristic surfaces $\Gamma_{x_1}$ for five different values of $x_1 \in \TT$.  

The unit tangent vector to $\gamma_{x_1}(t) $ is given by 
\begin{subequations}
\label{eta-n-and-tau}
\begin{equation} 
\tau(\eta(x,t),t):= |\p_2\eta|^{-1}  \p_2\eta= |\p_2\eta|^{-1} \big( \p_2\eta^1,\p_2\eta^2)
\end{equation} 
and by a rotation,  the unit
normal vector $\gamma_{x_1}(t) $ is given by 
\begin{equation} 
n(\eta(x,t),t):=  |\p_2\eta|^{-1}  \p_2\eta^\perp=|\p_2\eta|^{-1} \big(\p_2\eta^2, - \p_2\eta^1) \,.
\end{equation} 
\end{subequations} 
Substituting this identity into the flow equation \eqref{eta-char-intro}, shows that $\eta(x,t)$ is a solution to
\begin{equation}
\label{shock-flow}
\p_t \eta(x_1,x_2,t) = u(\eta(x_1,x_2,t),t) + \alpha \sigma (\eta(x_1,x_2,t),t) \tfrac{ \p_2\eta^\perp(x,t)}{ |\p_2\eta(x,t)| } \,,
\qquad 
\eta( x , \initial) = x  \,.
\end{equation} 
Observe that \eqref{shock-flow} is a PDE for $\eta(x,t)$, while the traditional Lagrangian flow map \eqref{trad-flow} solves an ODE.   As a 
consequence of the structure of the PDE  \eqref{shock-flow},  it is not {\it a priori} clear if $\eta$ maintains the Sobolev regularity of the velocity $u$
and sound speed $\sigma$.  What is clear, however, is that the normal vector $n$
loses one derivative of smoothness relative to the physical variables $u$ and $\sigma$.

This mismatch in regularity between $n$ and $(u,\sigma)$ necessitates studying a differentiated form of the Euler equations, in which the fundamental
unknowns have the same regularity as the normal vector $n$.  Rather than analyzing the system \eqref{euler-shock}, we return to the system
\eqref{euler1}.
We first compute the evolution of the partial (space) derivatives of $u$ and $\sigma$.  We differentiate \eqref{euler1}, and then re-express the
resulting equation using the fast acoustic transport velocity $ \mathcal{V} _3$.  We find
\begin{subequations} 
\label{euler-d00}
\begin{align} 
\p_t u^i,_k + (u\cdot n + \alpha \sigma) u^i,_{kj} n^j + (u \cdot \tau)  u^i,_{kj}\tau^j + \alpha \sigma (  \sigma,_{ki} -  u^i,_{kj}n^j) + \alpha \sigma,_k \sigma,_i + u^j,_k u^i,_j & = 0 \,, \label{uik-eqn}\\
\p_t \sigma,_k +  (u\cdot n + \alpha \sigma) \sigma,_{kj}n^j +(u \cdot \tau) \sigma,_{kj}\tau^j + \alpha \sigma (  u^i,_{ki} -  \sigma,_{kj} n^j)+ \alpha \sigma,_k u^i,_i + u^j,_k \sigma,_j & =0 \,. \label{sk-eqn}
\end{align} 
\end{subequations} 
The system \eqref{euler-d00} constitutes the differentiated (fast) acoustic Euler equations.  It is imperative to study this differentiated form of the
Euler equations in order to avoid  derivative loss in the geometry.   We are using the following derivative notation: for a differentiable function
$f$, we write
\begin{equation*}
f,_k := \p_k f  \ \text{ for } \ k\in \{1,2\} \,.
\end{equation*}
We are also employing the Einstein summation convention in which repeated indices are summed from $1$ to $2$; e.g. 
\begin{equation*}
u^j,_k u^i,_j := \sum\nolimits_{j=1}^2 \p_k u^j \p_j u^i  \,.
\end{equation*}
Therefore, the equation \eqref{uik-eqn} is a matrix equation with indices $(i,k)$, with $i,k \in\{1,2\}$, and the equation \eqref{sk-eqn} is a vector equation with
indices $k$, where $k\in\{1,2\}$.

\subsection{An Arbitrary Lagrangian Eulerian (ALE) parameterization of the fast characteristic surfaces}
\label{sec:intro:ALE}
The differentiated equation set \eqref{euler-d00} will indeed be the foundation for our analysis, but we will not use the Lagrangian flow $\eta$ of
$ \mathcal{V} _3$ to parameterize the  fast acoustic characteristic surfaces.    Instead, we shall use a novel parameterization based
on the Arbitrary Lagrangian Eulerian description of fluid flow which we shall detail in Section \ref{sec:ALE}.   

The use of the natural map $\eta$ provides control of the second fundamental form along fast characteristic surfaces, but due to a mild
degeneracy created by the tangential re-parameterization symmetry, control of the first fundamental form necessitates a cumbersome analysis.
We avoid this problem: 
a tangential re-parameterization of $\eta$ is introduced in the form of the so-called
Arbitrary-Lagrangian-Eulerian (ALE) coordinates.    This tangential re-parameterization, via the ALE maps,  provides a simple identity to control both the curvature and the metric tensors associated to the fast characteristic surfaces, and is reminiscent of DeTurck's simplification \cite{DeTurck1983} of
Hamilton's \cite{Hamilton1982} Ricci flow local existence theorem, in which the infinite-dimensional  kernel of the linearized Ricci flow operator (caused
by the tangential re-parameterization symmetry) lead to derivative loss and Hamilton's application of Nash-Moser iteration.

Our ALE re-parameterization works in the following manner.
Because each curve $\gamma_{x_1}(t)$ is a graph over the set $\{ x_2 \in \mathbb{T} \}$, we introduce the {\it height} function $h(x_1,x_2,t)$ such that
$\gamma_{x_1}(t) = (h(x_1,x_2,t),x_2)$.
The induced metric on $\gamma_{x_1}(t) $  is given by
$g(x_1,x_2,t) = 1 + (h,_2(x_1,x_2,t))^2$
and the unit tangent vectors $\tt$ and normal vectors $\nn$ to the curves $\gamma_{x_1}(t) $ are  
\begin{equation} 
\tt(x_1, x_2,t) = g^{- {\frac{1}{2}} } ( h,_2, 1) 
\qquad \text{ and  } \qquad \nn(x_1, x_2,t) = g^{- {\frac{1}{2}} } (1, -h,_2) \,,
\label{nn-tt-early}
\end{equation} 
 respectively. 
We define the ALE family of maps $\psi$  by $\psi(x_1,x_2,t) = h(x_1,x_2,t) e_1 + x_2 e_2  $, with initial condition
$h(x_1,x_2,\initial)= x_1$. Note that ${\rm det}(\nabla \psi) = h,_1$.
   In order to preserve the shape of the characteristic surfaces $\Gamma_{x_1}$,  the family of diffeomorphisms 
$\psi( \cdot , t)$ must satisfy the constraint that 
$\p_t \psi \cdot \nn  =  (\mathcal{V}_3 \circ \psi) \cdot \nn $.  As we shall explain in Section \ref{sec:ALE}, 
the dynamics of the ALE family of diffeomorphisms $\psi( \cdot ,t)$ are governed by
\begin{equation*} 
\p_t \psi 
= \left( (u \circ \psi) \cdot \nn + \alpha \sigma \circ \psi \right) \nn +
\left( (u \circ \psi) \cdot \nn + \alpha \sigma \circ \psi \right) h,_2 \tt \,.
\end{equation*}
With this definition, the fast acoustic characteristics (see~Figure~\ref{fig:surf-pre}) are parameterized by $(x_2,t) \mapsto (\psi(x_1,x_2,t),t)$.

Shock formation is measured by the  {\it metric-rescaled Jacobian determinant} of the deformation tensor $ \nabla \psi$.   In particular, we define
\begin{equation} 
\Jg(x,t) = g(x,t)^{- {\frac{1}{2}} } \nabla \psi (x,t) = g(x,t)^{- {\frac{1}{2}} }\p_1 h (x,t)   \,. \label{Jg-def-early}
\end{equation} 
As can be seen in Figure \ref{fig1}, the fast acoustic characteristic surfaces impinge on one another due to compression, and  gradient blow-up
occurs exactly when
\begin{equation*}
\Jg(x,t) =0 \,.
\end{equation*}
The pre-shock, represented as the black curve in Figure \ref{fig1}, denotes the spacetime co-dimension-$2$ surface of {\it first singularities}, and is
defined to be the set of points $(x,t)$ such that the following two  conditions simultaneously hold
\begin{equation*}
\Jg(x,t)=0  \qquad \text{ and } \qquad \p_1 \Jg(x,t) =0 \,.
\end{equation*}
This set of pre-shocks indicates the location and time of the first gradient blow-up in the primary direction of wave steepening (the $x_1$-direction), parameterized by the
transverse coordinates $(x_2, ..., x_d)$.  The condition $\p_1\Jg(x,t)=0$ together with the
{\it non-degeneracy condition}
\begin{equation*}
\operatorname{Hess}( \Jg) > 0 \,,
\end{equation*}
indicate that points in the pre-shock set are {\it local minima} for $\Jg(x,t)$.

The red surface in Figure \ref{fig1} corresponds to the level set $\{\Jg(x,t)=0\}$ {\it downstream} of the pre-shock set.   The red surface indicates
the location in spacetime of subsequent gradient catastrophes, and again designates the location of   characteristic-surface-impingement caused by 
the steepening sound waves during compression.  This red surface $\{\Jg(x,t)=0\}$, together with the pre-shock set $\{\Jg(x,t)=0\}\cap \{\p_1\Jg(x,t)=0\}$,
form a portion of the future temporal boundary of the \MGHDB.

The green surface shown in Figure \ref{fig1}, {\it upstream} of the pre-shock, displays the distinguished slow acoustic characteristic surface which
passes through the pre-shock set.   This green surface forms the remaining portion of the future temporal boundary of the \MGHDB.   In effect,   this distinguished slow acoustic characteristic surface (passing through the pre-shock) plays the role of the event horizon.
These notions are discussed in further detail in Section \ref{sec:usersguide} in the simplified setting of the Euler equations in one space dimension.

\subsection{A new set of differentiated Riemann variables that prevent derivative loss}
We can now map the physical variables $(u,\sigma)$,  as well as the classical Riemann variables $(w,z,a)$ defined in \eqref{wza-intro},  into our 
ALE coordinate system.    We define
\begin{alignat*}{2}
U^i& = u^i \circ \psi \,, \qquad  &&\Sigma = \sigma \circ \psi \,,   \\
W&=U\cdo\nn + \Sigma\,, && Z= U\cdo\nn-\Sigma\,, \qquad A=U\cdo\tt \,.
\end{alignat*} 
 Since it is the differentiated form of the Euler equations \eqref{euler-d00} that will be our starting point, the non-differentiated 
variables $(U,\Sigma)$ and $(W,Z,A)$  will play only a secondary role in our analysis.  

 We introduce specially constructed
{\it differentiated Riemann variables} that remove any derivative loss from the resulting analysis.
First, for $i,k=1,2$, we define
\begin{alignat*}{2}
\Uik&=  u^i,_k \circ \psi \,, \qquad && \mathring\Sigma_k =  \sigma,_k \circ \psi \,.
\end{alignat*} 
The variable $\Uik(x,t)$ denotes the $(i,k)$-component of the matrix $ \nabla u(\psi(x,t),t)$, and the variable $\Sk(x,t)$ denotes the $k$-component
of the vector $\nabla \sigma (\psi(x,t),t)$.
Second, we introduce the differentiated Riemann variables as the vector fields (with components)
\begin{equation} 
\Wk = \nn^i\Uik  + \Sk  \,, \qquad \Zk = \nn^i \Uik - \Sk   \,, \qquad \Ak =\tt^i \Uik \,, \qquad k=1,2\,. \label{important}
\end{equation} 
Again, $\Wk$, $\Zk$, $\Ak$ denote the $k$-component of the vectors $\Wb:= \nn^{\mathsf T} \Ub+\Sb$,  $\Zb:= \nn^{\mathsf T} \Ub-\Sb$, and $\Ab=\tt^{\mathsf T} \Ub$, 
respectively.   The vector field $\Wb$ is the {\it dominant} differentiated Riemann variable, the vector field $\Zb$ is the {\it subdominant} differentiated Riemann variable, and the vector field $\Ab$ is the {\it tangential} component of the gradient of the fluid velocity vector.   These three vector fields are 
then further projected unto their normal and tangential components.  We define
\begin{equation*} 
\Wbn=\Wb\cdo\nn \,, 
\qquad 
\Zbn=\Zb\cdo\nn\,, 
\qquad 
\Abn=\Ab\cdo\nn \,, 
\qquad 
\Wbt=\Wb\cdo\tt\,, 
\qquad 
\Zbt=\Zb\cdo\tt \,, 
\qquad 
\Abt=\Ab\cdo\tt \,.
\end{equation*} 
The generic and compressive initial data that we employ are designed to create steepening
sound waves whose dominant direction of propagation is along the $x_1$ coordinate.  By design,  it is the function $\Wbn$ that
encodes the steepening of the sound wave, and which blows-up when the slope of this steepening sound wave becomes infinite.
Meanwhile,  the other five functions $\Zbn,\Abn, \Wbt,\Zbt,\Abt$ remain uniformly bounded throughout the entire shock
formation process.    As noted in
Remark \ref{rem:pn-w-blowup},  it is the Eulerian quantity $\p_n w$ that blows-up during the shock formation process, and  
it is therefore tempting to believe that the normal component of the dominant differentiated Riemann variable $\Wbn$ should be defined as $\p_n w \cir \psi$. This is, in actuality, not the case.   In fact, our new Riemann variables are chosen to satisfy the following important identities:
\begin{subequations} 
\label{good-unknowns}
\begin{align}
\Wbn & =  (\p_n w -a \p_n n \cdot \tau) \cir \psi \,, \ \ 
\Zbn 
=   (\p_n z -a \p_n n \cdot \tau) \cir \psi \,, \ \ 
\Abn 
=    (\p_n a - \tfrac 12 (w+z)  \p_n \tau\cdot n )\cir\psi  \,,
\\
\Wbt & = (\p_\tau w - a  \p_\tau n \cdot \tau )\cir \psi \,, \ \ 
\Zbt 
= (\p_\tau z -  a \p_\tau n \cdot \tau)\cir\psi \,, \ \ 
\Abt 
 = \big(\p_\tau a - \tfrac{1}{2} (w+z) \p_\tau \tau \cdot n\big) \cir \psi
\,.
\end{align}
\end{subequations} 
This {\it unique linear combination} of differentiated classical Riemann variables together with the {\it curvature of the fast acoustic characteristic surfaces}
creates the {\it good variables}  that prevent derivative loss.

 Let us note the importance of using these differentiated Riemann variables for our
analysis.  In addition to preventing derivative loss in energy estimates,  these differentiated Riemann variables are created so that modulo very small errors
\begin{equation} 
\Jg \Wbn (x,t) \sim \p_1 w_0(x) \,,  \label{JgWbn-sim-p1w0}
\end{equation} 
where $w_0(x)$ is the initial condition for the classical dominant Riemann variable defined in \eqref{wza}, and in particular, 
$$
w_0(x) = w(x,\initial) = u^1(x,\initial) + \sigma(x,\initial) \,.
$$
The approximate identity \eqref{JgWbn-sim-p1w0} indicates that $\Jg\Wbn$ is almost {\it frozen into the flow of the fast acoustic characteristic}, i.e., 
it is almost exactly transported by the fast characteristics, and this fact is
fundamental to all of the pointwise bounds that are used in our work.

We shall describe the variables defined in~\eqref{good-unknowns} in greater detail in Section \ref{sec:new:Euler:variables}. Next, we shall explain how~\eqref{JgWbn-sim-p1w0}
is used in our energy methodt.

\subsection{The leading order dynamics}
Our objective is to establish uniform Sobolev-type bounds for the solution of the Euler equations in the \MGHDB.  We work with the following collection of fundamental variables:
\begin{equation*}
(\Jg\Wbn, \Jg\Zbn, \Jg \Abn, \Wbt, \Zbt, \Abt) 
\end{equation*}
whose evolution is coupled to the following  basic geometric variables
\begin{equation*}
(\Jg, h,_2)
\end{equation*}
as well as as the undifferentiated sound speed $\Sigma$.

While the exact evolution equations are given in Section \ref{sec:new:Euler:variables} below, it is instructive at this point to write down the approximate dynamical
 systems in which only leading order terms are displayed.    It is convenient to define the directional derivative operators
 \begin{subequations} 
 \begin{align} 
 \p_\ttt &= g^{- {\frac{1}{2}} } \p_2 \,, \label{p-tt-approx} \\
 \p_\nnn &= \p_1 - \Jg g^{- {\frac{1}{2}} } h,_2 \p_2 \label{p-nn-approx} \,.
 \end{align} 
 \end{subequations} 
We note that  since $h,_2$ is proven to maintain $\OO(\eps)$ size for the duration of the  development of the Cauchy data, to leading order,  
$\p_\nnn \sim \p_1$.  For the same reason,
 $g \sim 1$ so that $\p_\ttt \sim \p_2$.   Thus, with no danger of arriving at an erroneous conclusion, the reader is safe to make these derivative
 replacements  in the discussion that now follows.

 To leading order we have the following system of equations for the normal component variables: 
 \begin{subequations} 
 \label{nn-approx}
 \begin{align} 
&
\tfrac{1}{\Sigma} \p_t(\Jg\Wbn  )
+ \alpha    \p_\ttt (\Jg\Abn) 
- \alpha     \Abn \p_\ttt \Jg
-  \mathsf{P_1}   \p_\ttt\tt\cdot\nn  = \operatorname{l.o.t.}  \,, 
\label{Wbn-approx}
\\
 &
\tfrac{1}{\Sigma} \p_t (\Jg\Zbn)  
- \alpha    \p_\ttt  (\Jg\Abn)
 -2 \alpha \Jgi \p_\nnn (\Jg\Zbn) 
+ \alpha    \Abn  \p_\ttt\Jg 
+  \mathsf{P_1}  \p_\ttt \tt\cdot\nn  
\notag \\
&\qquad \qquad
 -\mathsf{P_2}\Jgi \p_\nnn \tt \cdot \nn  
 +2 \alpha\Zbn \Jgi  \p_\nnn \Jg 
 = \operatorname{l.o.t.} \,, 
 \label{Zbn-approx}  \\
&
\tfrac{1}{\Sigma}  \p_t (\Jg\Abn)  
+ \alpha   \p_\ttt  (\Jg\Sbn)
- \alpha \Jgi  \p_\nnn(\Jg\Abn)
- \alpha    \Sbn \p_\ttt\Jg
+ \alpha \p_\ttt \tt \cdot \nn \Sbt
 \notag \\
&\qquad \qquad
+ \mathsf{P_1} \Jgi  \p_\nnn \tt \cdot\nn 
+ \alpha \Jgi \p_\nnn \Jg \Abn
= \operatorname{l.o.t.}\,,
\label{Abn-approx}
 \end{align} 
\end{subequations} 
where $\mathsf{P_1}=\tfrac{\alpha }{2}  (\Jg\Wbn + \Jg\Zbn - 2\Jg\Abt)$, $\mathsf{P_2}= 2 \alpha ( \Jg\Abn+ \Jg\Zbt)$, 
and $ \operatorname{l.o.t.} $ denotes a polynomial of the fundamental variables.   
We additionally have the following system of evolution equations for the tangential-component variables:
\begin{subequations} 
\label{tt-approx}
\begin{align} 
&
\tfrac{1}{\Sigma} \p_t \Wbt  
+ \alpha\p_\ttt \Abt -   \mathsf{P_3} \p_\ttt \tt \cdot\nn = \operatorname{l.o.t.} 
\,,
  \\
&
\tfrac{1}{\Sigma} \p_t \Zbt 
- \alpha  \p_\ttt \Abt 
- 2 \alpha \Jgi  \p_\nnn \Zbt 
+   \mathsf{P_3} \p_\ttt \tt \cdot\nn
  - \mathsf{P_4} \Jgi \p_\nnn  \tt\cdot\nn
= \operatorname{l.o.t.} \,, 
  \\
 &
\tfrac{1}{\Sigma} \p_t  \Ab_\tt
+ \alpha   \p_\ttt \Sbt
- \alpha \Jgi \p_\nnn\Abt 
-  \alpha   \p_\ttt \Sbn \tt \cdot \nn
+ \mathsf{P_3} \Jgi \p_\nnn \tt\cdot\nn
= \operatorname{l.o.t.} \,,
\end{align} 
\end{subequations} 
where $\mathsf{P_3}=\alpha  (\Omega +\Wbt+\Zbt )$ and $\mathsf{P_4}= 2 \alpha (\Abt-\Zbn)$.
We couple the above systems of equations with the leading-order evolution equations for $\Jg$, $h,_2$, and $\Sigma$:
\begin{subequations} 
\label{rest-approx}
\begin{align} 
\p_t \Jg & =  \tfrac{1+\alpha}{2} \Jg\Wbn   + \operatorname{l.o.t.} \label{Jg-approx} \,, \\
\p_t h,_2 & =  \tfrac{1+\alpha}{2} \Wbt  +  \operatorname{l.o.t.}\,,  \label{p2h-approx} \\
\p_t  \Sigma   & = -\alpha \Sigma  \Zbn +  \operatorname{l.o.t.} \,, \label{Sigma-approx}
\end{align} 
together with the identity
\begin{equation} 
\Sigma,_1 = \tfrac{1}{2} \Jg\Wbn +  \operatorname{l.o.t.}  \,. \label{p1-Sigma-approx}
\end{equation} 
\end{subequations} 

\subsection{An overview of our energy method}
Energy estimates are performed simultaneously for the normal-component variables $(\Jg\Wbn, \Jg\Zbn, \Jg \Abn)$ in \eqref{nn-approx}
and for the tangential-component variables $(\Wbt,\Zbt,\Abt)$ in \eqref{tt-approx}, and then separately for $\Jg$ and  $h,_2$ in
\eqref{rest-approx}.    The actual energy estimates will be done in three different coordinate systems in which the time coordinate $t$ is
transformed in three different ways, but in order to describe the main
difficulties that must be overcome, for pedagogical reasons, our discussion will be in terms of the independent variables $(x_1,x_2,t)$, and we shall 
use the derivative notation $\nb= (\eps \p_t, \eps \p_1, \p_2)$.   Our energy method will be performed at the level of the sixth-order differentiated system.

We begin with the normal-component system  \eqref{nn-approx}.  The rough idea is that we first multiply the equations \eqref{Zbn-approx} and \eqref{Abn-approx} by 
$\Jg$ (thereby eliminating the presence of the inverse power $\Jgi$),  
we then let $D^6$ act on each equation in  \eqref{nn-approx}, and then
test this differentiated equation set with $\Sigma^{-2\beta+1} \big(\Jg \varphi^{2r} \nb^6(\Jg\Wbn),   \varphi^{2r} \nb^6(\Jg\Zbn),2 \varphi^{2r}  \nb^6(\Jg\Abn) \big)$, 
where
$ \varphi $ is a weight function that degenerates to zero at the {\it future temporal boundary}.   The presence of the additional $\Jg$ weight for
$\nb^6(\Jg\Wbn)$ is needed to match the natural weight that will appear for $\nb^6(\Jg\Zbn)$ and $\nb^6(\Jg\Abn)$ (since we multiplied those
equations by $\Jg$ prior to differentiation).
The weight function $\varphi$ will be chosen in three different ways, corresponding to three different spacetime regions that we shall employ for the analysis.   The
values of $\beta$ and  $r$ will be chosen below.

  At this stage, we simply wish to explain why a weight function is necessary for the energy estimates, and why its choice is determined by
the region of spacetime that is being analyzed.   We shall consider the energy estimates term-by-term, and we shall begin with the first term
in \eqref{nn-approx} containing the time-derivative $\p_t$.

 For demonstration purposes only, let us drastically simplify the domain of integration  so that we can explain a few
of the fundamental ideas (in actuality, our spacetimes require certain changes of coordinates and are somewhat more complicated).   Let us suppose that our energy method employs the spacetime $[\initial, \mathsf{t}_{\sf top}] \times  \mathbb{T} ^2 $ with spacetime
integral $ \int_\initial^t \iint_{ \mathbb{T} ^2}$ for $\initial \le t \le  \mathsf{t}_{\sf top}$.

\subsubsection{ Energy estimates for the first term in \eqref{nn-approx}}
The first term  in all three equations in eqref{nn-approx} yields both the 
standard {\it energy norm},  and also a
compression-induced {\it damping norm}.    The standard energy norm is an $L^\infty$-in-time and $L^2$-in-space norm, while the so-called
damping norm is an $L^2$-in-time and $L^2$-in-space norm,  but with a smaller power of the weight function $\varphi$ (and hence better
regularity).   In order to explain this, let us set $\varphi=\Jg$ and let us suppose (again for demonstration purposes only) that
$\Jg(x,  \mathsf{t}_{\sf top})=0$.   With this assumption, the weight function $\Jg$ degenerates to zero at $t= \mathsf{t}_{\sf top}$ and indicates
a gradient catastrophe.   While this spacetime is artificially constructed for demonstration purposes only, it allows us to use the function
$\Jg$ as our example-weight; 
$\Jg$ possesses the   compression-property that all of the actual weight functions $\varphi$ must possess. 

To explain this, we will make crucial use of the following three properties:
\begin{itemize}[leftmargin=16pt] 
\item[(a)] The initial condition $w_0$, for the dominant classical Riemann variable $w$, satisfies 
$- {\tfrac{1}{\eps}} \le \p_1w_0(x) \le - {\tfrac{1}{2\eps}}  $ in a local open neighborhood of $x_1=0$.
\item[(b)] Observe that the evolution equation for $\Jg\Wbn$ in \eqref{Wbn-approx} does not have any normal-derivative terms (whereas the
$\Jg\Zbn$ equation and the $\Jg\Abn$ equation both do). By applying the fundamental theorem of calculus (in time) to equation \eqref{Wbn-approx},
and using that the time-integral of all of the tangential-derivative terms are small relative to $\p_1w_0$, 
we prove that
$$
\Jg\Wbn(x,t)= \p_1 w_0(x) + \OO(\eps) \text{ as } \eps \to 0\,.
$$
\item[(c)] Using (a), (b), and  the approximate identity  \eqref{Jg-approx}, we see that 
\begin{subequations} 
\begin{alignat}{2} 
- \Jg\Wbn & \ge \tfrac{1}{4\eps}  \quad &&\  \text{  in a local open neighborhood of $x_1=0$,} 
\label{compression0}  \\
-\p_t \Jg(x,t) & \ge \tfrac{1}{4\eps} \cdot \tfrac{1+\alpha}{2} \quad &&\  \text{  in a local open neighborhood of $x_1=0$.} 
\label{compression} 
\end{alignat}  
\end{subequations} 
\end{itemize} 
The inequality \eqref{compression}  indicates that the fluid is in {\it compression},\footnote{Recall that $\Jg$ is the metric-scaled Jacobian
determinant of the deformation tensor $ \nabla \psi$ associated to the flow of the fast acoustic transport velocity.  When $\p_t\Jg < 0$, this indicates
that infinitesimal volume elements are shrinking along the fast characteristic flow, which is the natural geometric description of a compressed
gas.} and is responsible for the emergence of the damping norm using
the following argument (with the weight $\varphi$ set to be $\Jg$) for our energy method applied to the first term in all three equations of 
\eqref{nn-approx}:
\begin{align} 
&\int_\initial^t \iint_{ \mathbb{T} ^2} \Sigma^{-2\beta}  \Jg^{\!\!2r+1} \p_t\big(\nb^6 (\Jg\Wbn), \nb^6(\Jg\Zbn), \nb^6(\Jg\Abn)\big)\cdot \big(\nb^6 (\Jg\Wbn), \nb^6(\Jg\Zbn), 2\nb^6(\Jg\Abn)\big)dx dt' \\
& \qquad
=
\int_\initial^t  \frac{d}{2dt} \iint_{ \mathbb{T} ^2} \Sigma^{-2\beta}  \Jg^{\!\!2r+1} \big( \snorm{\nb^6(\Jg\Wbn)}^2 +\snorm{\nb^6(\Jg\Zbn)}^2 
+2\snorm{ \nb^6(\Jg\Abn)}^2 \big) dx dt' 
\notag \\
& \qquad\qquad
+ (r+\tfrac{1}{2} ) \int_\initial^t   \iint_{ \mathbb{T} ^2} \Sigma^{-2\beta}  \Jg^{\!\!2r} \big( -\p_t\Jg) \big( \snorm{\nb^6(\Jg\Wbn)}^2 +\snorm{\nb^6(\Jg\Zbn)}^2 
+2\snorm{ \nb^6(\Jg\Abn)}^2 \big) dx dt'
\notag \\
& \qquad\qquad
+ \beta \int_\initial^t \iint_{ \mathbb{T} ^2} \Sigma^{-2\beta-1} \p_t\Sigma  \Jg \Jg^{\!\!2r} \big( \snorm{\nb^6(\Jg\Wbn)}^2 +\snorm{\nb^6(\Jg\Zbn)}^2 
+ 2\snorm{\nb^6(\Jg\Abn)}^2 \big) dx dt'  \,. \label{marmot-is-on-a-diet1}
\end{align} 
Now, the ALE sound speed $\Sigma$ satisfies $ \tfrac{\kappa_0}{4} \le \Sigma \le \kappa_0$ for a fixed constant $\kappa_0 \ge 1$, and hence
the first integral on the right side of the equality is sign-definite and positive, and  therefore produces the $L^\infty$-in-time and $L^2$-in-space {\it energy norm}.   The second integral on
the right side of the equality  produces the damping norm; in particular, to avoid obfuscation, we shall assume that the compression inequality
\eqref{compression} holds globally in our spacetime set (this is, of course, not the case but it is a technical, rather than fundamental, matter to
contend with).  In this case, the second integral is  is again sign-definite and positive, and bounded from below by
\begin{equation} 
\tfrac{(2r+1)(1+\alpha)}{16\eps}  \int_\initial^t   \iint_{ \mathbb{T} ^2} \Sigma^{-2\beta}  \Jg^{\!\!2r}  \big( \snorm{\nb^6(\Jg\Wbn)}^2 +\snorm{\nb^6(\Jg\Zbn)}^2 
+ 2\snorm{ \nb^6(\Jg\Abn)}^2 \big) dx dt'  \,,
\label{damping-integral-approx}
\end{equation} 
and this is the $L^2$-in-time and $L^2$-in-space {\it damping norm} that arises from compression.    Notice that the pre-factor in front of this integral
is proportional to $ {\tfrac{1}{\eps}} $ and that $\eps>0$ is  a parameter which is taken to be very small.   This means that other integrals that arise from
our energy method that have the same type of integrand, {\it but whose pre-factor is $\OO(1)$ as $\eps \to 0$} can be viewed as
small errors that our damping integral can easily absorb. In particular, the third integral on the right side of \eqref{marmot-is-on-a-diet1} is an example of
such a small error integral.   From \eqref{Sigma-approx}, we have that $\p_t\Sigma $ is approximately equal to $- \alpha \Zbn \Sigma $ 
and, as we will show in Section \ref{sec:pointwise:bootstraps},  $ |\Zbn| \le C_{\Zbn} =\OO(1)$ and $\Jg \le \tfrac{6}{5} $;  therefore, one simply requires that 
$ \tfrac{(2r+1)(1+ \alpha )}{16} $ is larger than 
$ \eps \alpha\beta \tfrac{12 C_{\Zbn}}{5} $ and by choosing $\eps <   \tfrac{5 (2r+1)(1+ \alpha )}{192 \alpha \beta C_{\Zbn}}$ this is indeed achieved.

To summarize, the first term in all three equations in \eqref{nn-approx} have produced two types of regularity via an {\it energy} norm and a {\it damping} norm.  Reverting to the notation
$\varphi$ for the weight function, we record these norms here
as follows:
\begin{subequations} 
\label{norms-approx}
\begin{align} 
\mathcal{E}_{6,\nnn}^2(t)  &= \snorm{  \varphi^{\!\frac 34} \Jgh \nb^6 (\Jg\Wbn,\Jg\Zbn, \Jg\Abn)( \cdot , t)}^2_{L^2_x}  \,, \label{E6n-approx} \\
\mathcal{D}_{6,\nnn}^2(t) 
&=   \int_0^t 
\snorm{\varphi^{\!\frac 14}  \Jgh\nb^6 (\Jg \Wbn, \Jg\Zbn, \Jg\Abn)( \cdot , t')}_{L^2_x}^2 {\rm d}t' \,. \label{D6n-approx}
\end{align} 
\end{subequations} 
The purpose of our energy method is to obtain uniform bounds on the
 norms $\mathcal{E}_{6,\nnn}^2(t)$ and $\mathcal{D}_{6,\nnn}^2(t) $ for $\initial \le t \le  \mathsf{t}_{\sf top}$.   The energy
 estimates for the tangential-component equations \eqref{tt-approx} will produce analogous energy and damping norms 
 $\mathcal{E}_{6,\ttt}^2(t)$ and $\mathcal{D}_{6,\ttt}^2(t) $ which we will define below, and we also obtain uniform bounds for these tangential norms
 for $\initial \le t \le  \mathsf{t}_{\sf top}$.

\subsubsection{Energy estimates for the second term in \eqref{nn-approx}}
Recall that for the purposes of this pedagogical overview, we are setting the weight function $\varphi$ to equal $\Jg$, we multiplying
\eqref{Zbn-approx} and \eqref{Abn-approx} by $\Jg$,  we then let 
 $D^6$ act on each equation in  \eqref{nn-approx}, and test the resulting equations 
with $ \big(\Jg \Jg^{\!\!2r} \nb^6(\Jg\Wbn),  \Jg^{\!\!2r} \nb^6(\Jg\Zbn),2 \Jg^{\!\!2r}  \nb^6(\Jg\Abn) \big)$. We now focus on
the second term in \eqref{nn-approx}.   The second terms in all three equations must be grouped together to form an exact derivative which can then
be moved off of the highest-derivative term by use of integration-by-parts.   Note that for us to be able to group all three terms together to 
form an exact derivative,   it is essential that the weights are identical in all three equations.
Our energy method yields the following combination of highest-order  integrals:
\begin{align*} 
 \alpha   \int_\initial^t   \iint_{ \mathbb{T} ^2}\Jg^{\!\!(2r+1)} \Big(    \p_\ttt \nb^6 (\Jg\Abn) \big( \nb^6(\Jg\Wbn)-  \nb^6(\Jg\Zbn)  \big) 
 + 2 \p_\ttt  \nb^6(\Jg\Sbn) \nb^6 (\Jg\Abn) \Big) {\rm d} x {\rm d} t' \,.
\end{align*} 
By definition, we have that $\nb^6(\Jg\Wbn)-  \nb^6(\Jg\Zbn)  = 2 \nb^6(\Jg\Sbn) $, and hence the above integrand contains the exact derivative
$2 \p_{\ttt} \big(  \nb^6(\Jg\Sbn) \nb^6 (\Jg\Abn)\big)$ and upon integrating-by-parts with respect to $\p_{\ttt}$, we obtain at highest-order, the
resulting integral
\begin{align*} 
 -(2r+1) \alpha   \int_\initial^t   \iint_{ \mathbb{T} ^2} \p_\ttt\Jg \, \Jg^{\!\!2r}
   \big( \nb^6(\Jg\Wbn) - \nb^6(\Jg\Zbn)\big) \,  \nb^6 (\Jg\Abn) {\rm d} x {\rm d} t' \,.
\end{align*} 
We prove in Section \ref{sec:pointwise:bootstraps} that $\p_\ttt\Jg = \OO(1)$ and hence, the above integral is an {\it error integral} which is easily controlled by our damping norm $ \mathcal{D}_{6,\nnn}(t)$  via an
application of the Cauchy-Young inequality.

\subsubsection{Energy estimates for the third term in \eqref{Zbn-approx} and \eqref{Abn-approx} }   
We now focus on the energy estimates for the
 third term in \eqref{Zbn-approx} and \eqref{Abn-approx}.  As can be seen in \eqref{Wbn-approx},   this type of normal derivative term does not exist in 
 the equation for $\Jg\Wbn$, but the presence of such terms in the evolution equations for $\Jg\Zbn$ and $\Jg\Abn$ create the fundamental difficulty
in the analysis of the \MGHDB, and are the terms which are primarily responsible for our use of three different spacetime 
regions with three different weight functions $\varphi$.

Our energy method applied to the third term in \eqref{Zbn-approx} and \eqref{Abn-approx} produces the following integral:
\begin{equation} 
-\alpha \int_\initial^t   \iint_{ \mathbb{T} ^2} \Sigma^{-2\beta+1} \varphi^{2r} \p_\nnn 
\big( \sabs{ \nb^6(\Jg\Zbn)}^2 +  \sabs{ \nb^6(\Jg\Abn)}^2\big) {\rm d}x {\rm d}t' \,.
\label{marmot-is-on-a-diet7}
\end{equation} 
Observe that we have written the weigh function in this integral  as $\varphi^{2r} $. 
While for the purposes of this simplified overview, we  have set $\varphi$  equal to $\Jg$,  for the integral arising from this third term, we 
use the more general $\varphi$ for the weight function, and we will explain the reason for this notational choice below.

Upon integrating-by-parts with respect to $\p_\nnn$ in \eqref{marmot-is-on-a-diet7},  using \eqref{p-nn-approx} and the fact that $h,_2=\OO(\eps)$, 
we have that  \eqref{marmot-is-on-a-diet7},  to leading order,  is given by
\begin{align} 
&
- \alpha (2   \beta-1)  \int_\initial^t   \iint_{ \mathbb{T} ^2} \Sigma^{-2\beta}  \varphi^{2r} \p_1 \Sigma   \, 
\big( \sabs{ \nb^6(\Jg\Zbn)}^2 +  \sabs{ \nb^6(\Jg\Abn)}^2\big) {\rm d}x {\rm d}t 
\notag \\
& \qquad\qquad
+2  \alpha r \int_\initial^t   \iint_{ \mathbb{T} ^2} \Sigma^{-2\beta+1} \varphi^{2r-1} \p_1 \varphi \, 
\big( \sabs{ \nb^6(\Jg\Zbn)}^2 +  \sabs{ \nb^6(\Jg\Abn)}^2\big) {\rm d}x {\rm d}t 
 \,.
\label{marmot-is-on-a-diet2}
\end{align} 
The first integral in \eqref{marmot-is-on-a-diet2} produces another (arbitrarily large) type of damping term thanks to the use of the function
$\Sigma^{-2\beta+1}$ as part of the weighting.   We use the approximate identity \eqref{p1-Sigma-approx} and we again assume (only for this
demonstration) that the compression lower-bound \eqref{compression0} holds globally (rather than locally as stated).  In this case, we see
that the  first integral in \eqref{marmot-is-on-a-diet2} has the positive lower bound
$$
\tfrac{ \alpha  (2\beta-1)}{8 \eps}  \int_\initial^t   \iint_{ \mathbb{T} ^2} \Sigma^{-2\beta-1}  \varphi^{2r}    \, 
\big( \sabs{ \nb^6(\Jg\Zbn)}^2 +  \sabs{ \nb^6(\Jg\Abn)}^2\big) {\rm d}x {\rm d}t 
$$
whose size can be adjusted with the choice of $\beta> \tfrac{1}{2} $.

\subsubsection{Analysis is split into three regions of spacetime} 
 \begin{figure}[htb!]
  \includegraphics[width=.5\linewidth]{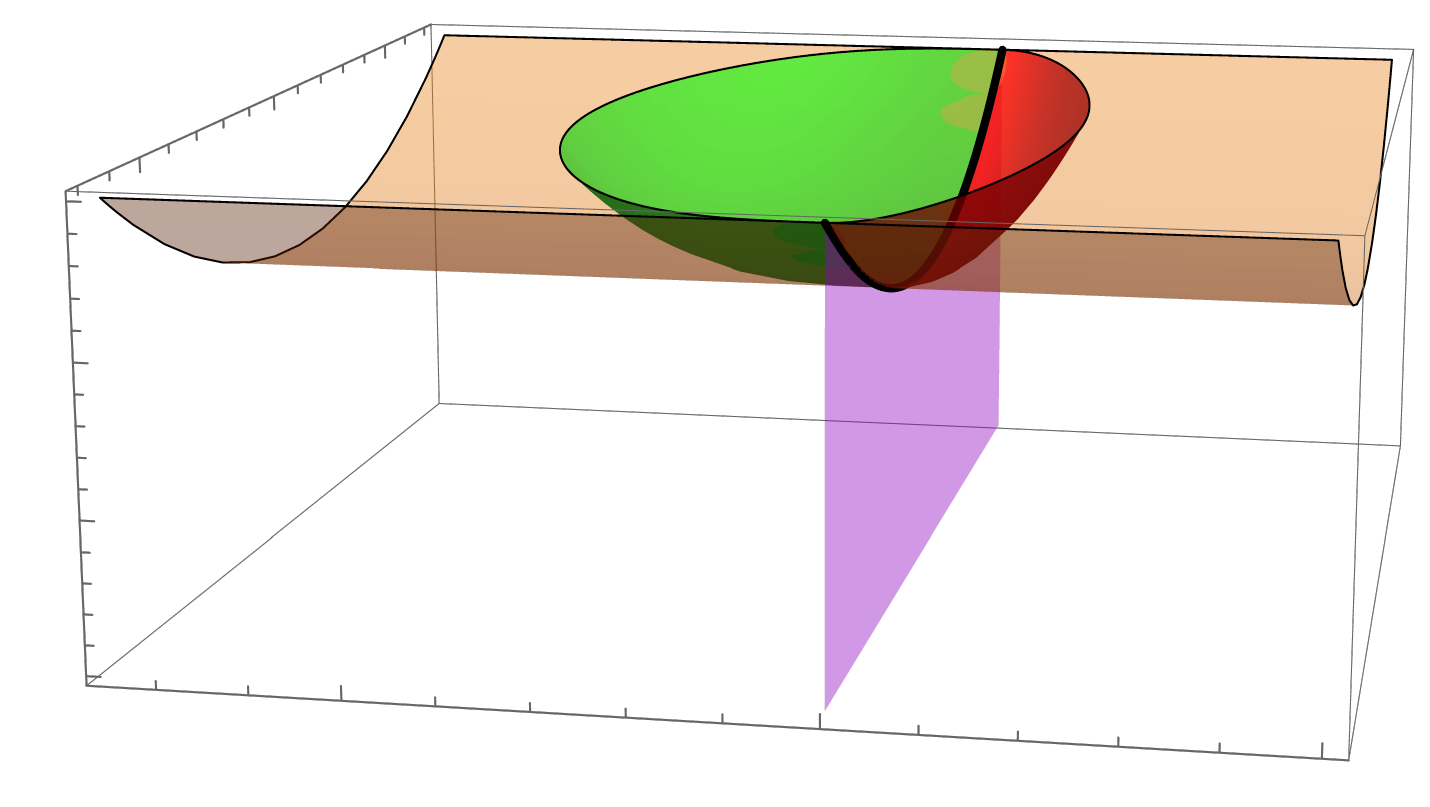}
\vspace{-.1 in}
\caption{Four fundamental hypersurfaces are displayed. The ``nearly vertical'' surface $x_1=x_1^*(x_2,t):= \operatorname{argmin} _{x_1 \in \mathbb{T}  } \Jg(x_1,x_2,t) $ is shown in magenta.  This surface passes through the set of pre-shocks, displayed as the
 black curve.  In red, the downstream surface $\{\Jg(x,t)=0\}$ is displayed, and in green, the upstream slow acoustic characteristic surface that passes
 through the pre-shock set is shown.  In orange, the cylindrical surface $t=t^*(x_2)$ is displayed, where $t^*(x_2)$ denotes the time coordinate along
 the set of pre-shocks.}
\label{fig:2D:Lagrangian:a0}
\end{figure}
The second integral  in \eqref{marmot-is-on-a-diet2} is a highly problematic {\it error integral}. The difficulty in bounding this
error integral is twofold: first,  the damping integral
has $\varphi^{2r}$ as the weight, while the problematic error integral has only $\varphi^{2r-1}$ as a weight and hence, this integral cannot be naively absorbed by our
damping norm; second, the function $\p_1\varphi$ is not a signed function, and therein lies a fundamental difficulty for our analysis.   
Specifically, we consider the hypersurface, shown as the (almost vertical) magenta surface in Figure \ref{fig:2D:Lagrangian:a0}, defined by 
$$
x_1^*(x_2,t) :=  \operatorname{argmin} _{x_1 \in \mathbb{T}  } \Jg(x_1,x_2,t) \,.
$$
We can now define the {\it upstream} and the {\it downstream} regions of spacetime.
The {\it upstream region}  consists of  all triples $(x_1,x_2,t)$ such that $\{ x_1 < x_1^*(x_2,t)\}$ and time $t$ is less
than the times corresponding to the 
green slow characteristic surface in the right panel in Figure \ref{fig:3spacetimes}.  Similarly, the {\it downstream region} consists of 
 all triples $(x_1,x_2,t)$ such that $\{ x_1 > x_1^*(x_2,t)\}$ and time $t$ is less than the times corresponding to the red surface in the
 center panel in 
 Figure \ref{fig:3spacetimes} where $\Jg$ vanishes.   The magenta surface $x_1=x_1^*(x_2,t)$ passes through the pre-shock, denoted by
 the black curve in Figure \ref{fig:2D:Lagrangian:a0} (see also Figure \ref{fig:3spacetimes}).
 
 Clearly, the sign of $\p_1 \varphi$ is of basic importance in estimating the second integral in \eqref{marmot-is-on-a-diet2}. Since we
 have set $\varphi$ to equal $\Jg$ for this overview, we have that $\p_1\varphi=\p_1\Jg$.   As $\p_1\Jg=0$ on
 the surface $x_1= x_1^*(x_2,t)$, it must change sign from  the upstream region  to  the downstream region.  
 In fact, we have that  $\p_1 \Jg < 0$ in the upstream region, and $\p_1\Jg>0$
 in the downstream region.\footnote{This is only true in a local open neighborhood of the surface $x_1= x_1^*(x_2,t)$.  For $x_1 > x_1^*(x_2,t)$ sufficiently large, $\p_1\Jg$ becomes negative, but this technical difficulty can be ignored at this stage of the presentation.}
 Returning to the second integral in \eqref{marmot-is-on-a-diet2}, we see that this error integral acquires a ``good'' sign in the downstream region and
 can be viewed as an additional damping integral with a reduced power in the weight function $\Jg^{\!\!2r-1}$.  Meanwhile, in the upstream region
 in which $\p_1\Jg < 0$, this error integral has the ``bad'' sign and cannot be properly bounded with this choice of weight function.   To be precise,
 {\it we cannot set the weight function $\varphi$  equal to $\Jg$ in the upstream region}; instead, we must devise a weight function $\varphi$ that obeys good properties with respect $\p_1\varphi$ in the upstream region; we will use a weight function which is  essentially  {\it transported} by the slow 
acoustic characteristics.   This type of transport structure produces a cancellation that entirely eliminates the problematic error integral in  
\eqref{marmot-is-on-a-diet2} in the upstream region.

 Clearly, different weights must be used in the upstream and downstream regions of spacetime so as to
 bound the error integral in \eqref{marmot-is-on-a-diet2}.   To that end, we consider a third region of spacetime, whose closure contains 
 the pre-shock set as well as large subsets of both the upstream and downstream regions. See the left panel of Figure \ref{fig:3spacetimes},
 which displays this third region, consisting of the triples $(x_1,x_2,t)$ such
that time $t$ is less than the times corresponding to the set of pre-shocks. 
Because the set of pre-shocks depends only upon the transverse coordinate $x_2$ and time $t$,  in order to bound the error integral in  
\eqref{marmot-is-on-a-diet2} in this third spacetime region, we can set the weight function $\varphi$ to be a function
 that degenerates to zero along this cylindrical surface and which is independent of the $x_1$ coordinate.    This can be done by setting
 $\varphi(x_2,t) = \Jg (x_1^*(x_2,t),x_2,t)$, in which case the problematic error integral in \eqref{marmot-is-on-a-diet2} vanishes.\footnote{The actual weight function  for this region of spacetime uses a modification of  $\Jg$ which is defined in \eqref{def-Jgbar}.}
 
\subsubsection{ Energy estimates for the remaining terms in \eqref{nn-approx}}
It remains for us to explain how our energy method bounds the fourth term in \eqref{Wbn-approx} and the fourth, fifth, sixth, and seventh terms in
\eqref{Zbn-approx} and \eqref{Abn-approx}.   The common feature which these terms share is {\it over-differentiated
geometry}; however, instead of producing derivative-loss,  all of these terms (by design) have exact derivative structure, and some of the 
over-differentiated terms produce new types of damping norms,
as well as {\it anti-damping} norms, the latter requiring the choice of sufficiently large exponent $r$.

Perhaps the most interesting of the {\it over-differentiated} terms is the fourth term in \eqref{Abn-approx}, $- \alpha \Sbn \p_\ttt \Jg$.   Recalling
that the procedure for our energy method requires first multiplying by $\Jg$, then applying $\nb^6$, and then testing with 
$2\Sigma^{-2\beta+1} \varphi^{2r} \nb^6(\Jg\Abn)$, we find that the leading order integral obtained from this term is given
by
\begin{equation} 
- 2\alpha  \int_\initial^t   \iint_{ \mathbb{T} ^2} \Sigma^{-2\beta+1}  \varphi^{2r}  (\Jg \Sbn) \p_\ttt \nb^6 \Jg \,  \nb^6(\Jg\Abn) {\rm d}x{\rm d}t' \,.
\label{marmot-is-on-a-diet3}
\end{equation} 
This integral explains the over-differentiated geometry nomenclature: it appears that there is one too many derivatives in  $\p_\ttt \nb^6\Jg$.  We
can have six derivatives on $\Jg$ but not seven, and as we shall see, there is indeed an exact derivative structure here.   Let us explain why
six derivatives on $\Jg$ is in agreement with the norms \eqref{norms-approx} of our energy method.  

 We return to the approximate dynamics for
$\Jg$ given by \eqref{Jg-approx} and perform a sixth-order energy estimate on this relation.  The resulting differential inequality yields, via the
Gr\"onwall inequality, the bound
\begin{equation} 
\sup_{t' \in [\initial,t]} \eps \snorm{ \varphi^{\! \frac 14} \nbs^6 \Jg (x,t') }_{ L^2_x}^2
+ \int_{\initial}^t \snorm{  \varphi^{-\! \frac 14}\nbs^6 \Jg  (x,t')}_{L^2_{x}}^2 {\rm d}t'
\les  \eps^2 \mathcal{D}_{6,\nnn}^2(t)  \,.  \label{Jg-bound-approx}
\end{equation} 
 Because $\varphi$ is equal go $\Jg$ and because $\Jgi \ge \tfrac{5}{6} $ which is proven in Section \ref{sec:pointwise:bootstraps}, 
we see from \eqref{Jg-bound-approx} that the unweighted function
 $\nb^6 \Jg$ is bounded in the spacetime $L^2$-norm by $\eps  \mathcal{D}_{6,\nnn}^2(t) $, showing not only that six derivatives
of $\Jg$ are bounded, but that the upper bound of the spacetime $L^2$-norm is extremely small.

To estimate the integral \eqref{marmot-is-on-a-diet3}, we integrate-by-parts with respect to $\p_\ttt$.   Using that $\Jg \Sbn = \tfrac{1}{2} (\Jg\Wbn-\Jg\Zbn)$, that $\Jg\Wbn = \OO( - {\tfrac{1}{\eps}} )$, and that $\Jg\Zbn = \OO(1)$, integration-by-parts 
in \eqref{marmot-is-on-a-diet3} produces the highest order integral:
\begin{equation} 
I^{ \operatorname{error} }
=
\int_\initial^t   \iint_{ \mathbb{T} ^2} \Sigma^{-2\beta+1}  \varphi^{2r}  (\Jg \Wbn)  \nb^6 \Jg \,  \alpha  \p_\ttt\nb^6(\Jg\Abn) {\rm d}x{\rm d}t' \,.
\label{marmot-is-on-a-diet4}
\end{equation} 
At this stage, we use the evolution equation \eqref{Wbn-approx} to make the substitution
$$
 \alpha    \p_\ttt \nb^6(\Jg\Abn) 
= 
-\tfrac{1}{\Sigma} \p_t \nb^6(\Jg\Wbn  )
+ \alpha     \Abn \p_\ttt \nb^6 \Jg
+  \mathsf{P_1}   \p_\ttt \nb^6\tt\cdot\nn  
+\operatorname{l.o.t.}  
$$
To leading order,  the integral $I^{ \operatorname{error} }$ in  \eqref{marmot-is-on-a-diet4} is written as
\begin{align} 
I^{ \operatorname{error} }&= I^{ \operatorname{error} }_1+I^{ \operatorname{error} }_2+I^{ \operatorname{error} }_3 \,,
\notag \\
I^{ \operatorname{error} }_1&=
-  \int_\initial^t   \iint_{ \mathbb{T} ^2} \Sigma^{-2\beta}  \varphi^{2r}  (\Jg \Wbn)  \nb^6 \Jg \,   \p_t\nb^6(\Jg\Wbn  )  {\rm d}x{\rm d}t'
\,, \notag \\
I^{ \operatorname{error} }_2&=
\alpha    \int_\initial^t   \iint_{ \mathbb{T} ^2} \Sigma^{-2\beta+1}  \varphi^{2r}  (\Jg \Wbn)  \Abn  \nb^6 \Jg \,   \p_\ttt\nb^6\Jg {\rm d}x{\rm d}t' 
\,, \notag \\
I^{ \operatorname{error} }_3&=
  \int_\initial^t   \iint_{ \mathbb{T} ^2} \Sigma^{-2\beta+1}  \varphi^{2r}  (\Jg \Wbn)   \mathsf{P_1} g^{- {\frac{1}{2}} }  \, \nb^6 \Jg \,   \p_2\nb^6\tt\cdot\nn   {\rm d}x{\rm d}t' \,.
\label{marmot-is-on-a-diet5}
\end{align} 
Let us first explain how we bound the integrals $I^{ \operatorname{error} }_2$ and $I^{ \operatorname{error} }_3$.   The integral
$I^{ \operatorname{error} }_2$ has the obvious exact-derivative structure, which together with an integration-by-parts with respect to $\p_\ttt$
produces the following:
\begin{align*} 
I^{ \operatorname{error} }_2
&=
\tfrac{\alpha }{2}   \int_\initial^t   \iint_{ \mathbb{T} ^2} \Sigma^{-2\beta+1}  \varphi^{2r}  (\Jg \Wbn)  \Abn   \p_\ttt \sabs{\nb^6\Jg}^2 {\rm d}x{\rm d}t'
=
-\tfrac{\alpha }{2}   \int_\initial^t   \iint_{ \mathbb{T} ^2}  \p_\ttt\big( \Sigma^{-2\beta+1}  \varphi^{2r}  (\Jg \Wbn)  \Abn \big)  \sabs{\nb^6\Jg}^2 {\rm d}x{\rm d}t' \,.
\end{align*} 
This integral is now easily bounded using \eqref{Jg-bound-approx} together with our established pointwise bounds for the coefficient function which
are proven to hold in Section \ref{sec:pointwise:bootstraps}.

In the same way that we obtained the bound \eqref{Jg-bound-approx} for $\Jg$, we can obtain the analogous bound for $h,_2$.  We
perform a sixth-order energy estimate on \eqref{p2h-approx} and find that
\begin{equation} 
\sup_{t' \in [\initial,t]} \eps \snorm{ \varphi^{\! \frac 14} \nbs^6 h,_2 (x,t') }_{ L^2_x}^2
+ \int_{\initial}^t \snorm{  \varphi^{-\! \frac 14}\nbs^6 h,_2  (x,t')}_{L^2_{x}}^2 {\rm d}t'
\les \mathsf{K} \eps^2 \mathcal{D}_{6,\ttt}^2(t)  \,, \label{p2h-bound-approx}
\end{equation} 
where $ \mathcal{D}_{6,\ttt}(t)$ is the damping norm that arises from the tangential-component energy estimates for the system 
\eqref{tt-approx} and $\mathsf{K}$ is a constant used in the bounds for those tangential-component energy estimates.   From 
\eqref{nn-tt-early} and \eqref{Jg-def-early}, to leading order, we have that
\begin{equation*}
\nb^6 \Jg = g^{- {\frac{1}{2}} } \nb^6 h,_1 \ \text{ and } \ \p_2\nb^6\tt \cdot \nn = g^{-1} \nbs^6 h,_{22} \,.
\end{equation*}
Therefore, to leading order, the  integral $I^{ \operatorname{error} }_3$ in \eqref{marmot-is-on-a-diet5} can be written as
\begin{equation*} 
I^{ \operatorname{error} }_3
=
\int_\initial^t   \iint_{ \mathbb{T} ^2} \Sigma^{-2\beta+1}  \varphi^{2r}  (\Jg \Wbn)   \mathsf{P_1} g^{- 2 }  \, \nb^6 h,_1 \,   \nb^6 h,_{22} {\rm d}x{\rm d}t' \,.
\end{equation*} 
We integrate-by-parts with respect to $\p_2$, form the exact derivative $\p_1\sabs{ \nbs^6h,_2}^2$, and  integrate-by-parts with respect to $\p_1$.  We
obtain that to leading order,
\begin{equation*} 
I^{ \operatorname{error} }_3
=
\tfrac{1 }{2}   \int_\initial^t   \iint_{ \mathbb{T} ^2}\p_1\big( \Sigma^{-2\beta+1}  \varphi^{2r}  (\Jg \Wbn)   \mathsf{P_1} g^{- 2 } \big)\sabs{  \nb^6 h,_{2}}^2 {\rm d}x{\rm d}t'  \,,
\end{equation*} 
and this integral is easily bounded using \eqref{p2h-bound-approx}.

We now explain how to estimate the  integral $I^{ \operatorname{error} }_1$ in \eqref{marmot-is-on-a-diet5}.
We first integrate-by-parts with respect to $\p_t$, use \eqref{Jg-approx},  and find that to leading order with $\varphi$ set equal to $\Jg$, we have that
\begin{equation*} 
I^{ \operatorname{error} }_1
=
\int_\initial^t   \iint_{ \mathbb{T} ^2} \Sigma^{-2\beta}  \Jg^{\!\!2r}  (\p_t\Jg) \sabs{\nb^6(\Jg\Wbn  )}^2  {\rm d}x{\rm d}t'
\end{equation*} 
Comparing this integral with the second integral on the right side of \eqref{marmot-is-on-a-diet1}, we see that the sign on $\p_t\Jg$ is now positive
while it is negative in  \eqref{marmot-is-on-a-diet1}.   In   \eqref{marmot-is-on-a-diet1}, the second integral on the right side is a damping integral,
while the integral $I^{ \operatorname{error} }_1$ produces an {\it anti-damping} integral, meaning a sign-definite integral but with the wrong sign.
Consequently, this anti-damping integral $I^{ \operatorname{error} }_1$ must be combined with the damping integral on the right side of
\eqref{marmot-is-on-a-diet1}.  The sum of these two integrals yields the same type of integral as in \eqref{damping-integral-approx} for the
function $\nb^6(\Jg\Wbn)$; more precisely, we have the lower bound for this sum of integral given by
\begin{equation} 
\tfrac{(r-\frac{1}{2} )(1+\alpha)}{8\eps}  \int_\initial^t   \iint_{ \mathbb{T} ^2} \Sigma^{-2\beta}  \Jg^{\!\!2r}  \snorm{\nb^6(\Jg\Wbn)}^2 dx dt' \,,
\label{damping-integral-approx2}
\end{equation} 
it is is therefore clear that in order to obtain our damping norm, we must choose $r> \tfrac{1}{2} $.  For our analysis, we set 
$r= \tfrac{3}{4}$.

We have now given an overview of our energy method for the normal-component system \eqref{nn-approx}.   Similar energy estimates
are performed for the tangential-component system, leading us to bound the following norms
\begin{subequations} 
\label{norms-approx-tt}
\begin{align} 
\mathcal{E}_{6,\ttt}^2(t)  &= \snorm{  \varphi^{\!\frac 34} \Jgh \nb^6 (\Wbt,\Zbt, \Abt)( \cdot , t)}^2_{L^2_x}  \,, \label{E6t-approx} \\
\mathcal{D}_{6,\ttt}^2(t) 
&=   \int_0^t 
\snorm{\varphi^{\!\frac 14}  \Jgh\nb^6 ( \Wbt, \Zbt, \Abt)( \cdot , t')}_{L^2_x}^2 {\rm d}t' \,. \label{D6t-approx}
\end{align} 
\end{subequations} 

Defining the {\it total norms}  by
\begin{align*} 
\mathcal{E}_{6}^2(t) &=  \mathcal{E}_{6,\nnn}^2(t)  + (\mathsf{K}\eps)^{-2} \mathcal{E}_{6,\ttt}^2(t) \,, \notag \\
\mathcal{D}_{6}^2(t) &=  \mathcal{D}_{6,\nnn}^2(t)  + (\mathsf{K}\eps)^{-2} \mathcal{D}_{6,\ttt}^2(t) \,,
\end{align*} 
our combined energy estimates prove that $\mathcal{E}_{6}(t)$ and $\mathcal{D}_{6}(t)$ remain uniformly bounded 
for $\initial \le t \le  \mathsf{t}_{\sf top}$.

As we stated above, this overview used a highly simplified spacetime to explain some of the key ideas for our energy method.  The actual
scheme uses three different spacetime regions: (1) the spacetime region bounded from above by the pre-shock cylindrical surface shown
in the left panel of Figure \ref{fig:3spacetimes}, (2) the downstream region bounded from above by the level set $\{\Jg(x,t)=0\}$,  shown in red in
the center panel of Figure \ref{fig:3spacetimes}; and 
(3) the upstream region bounded from above by the distinguished slow acoustic characteristic surface passing through the pre-shock, shown in
green in the right panel of Figure~\ref{fig:3spacetimes} below.


\subsection{A rough statement of the main theorem}
\begin{figure}[htb!]
\centering
\begin{minipage}{.3\linewidth}
  \includegraphics[width=\linewidth]{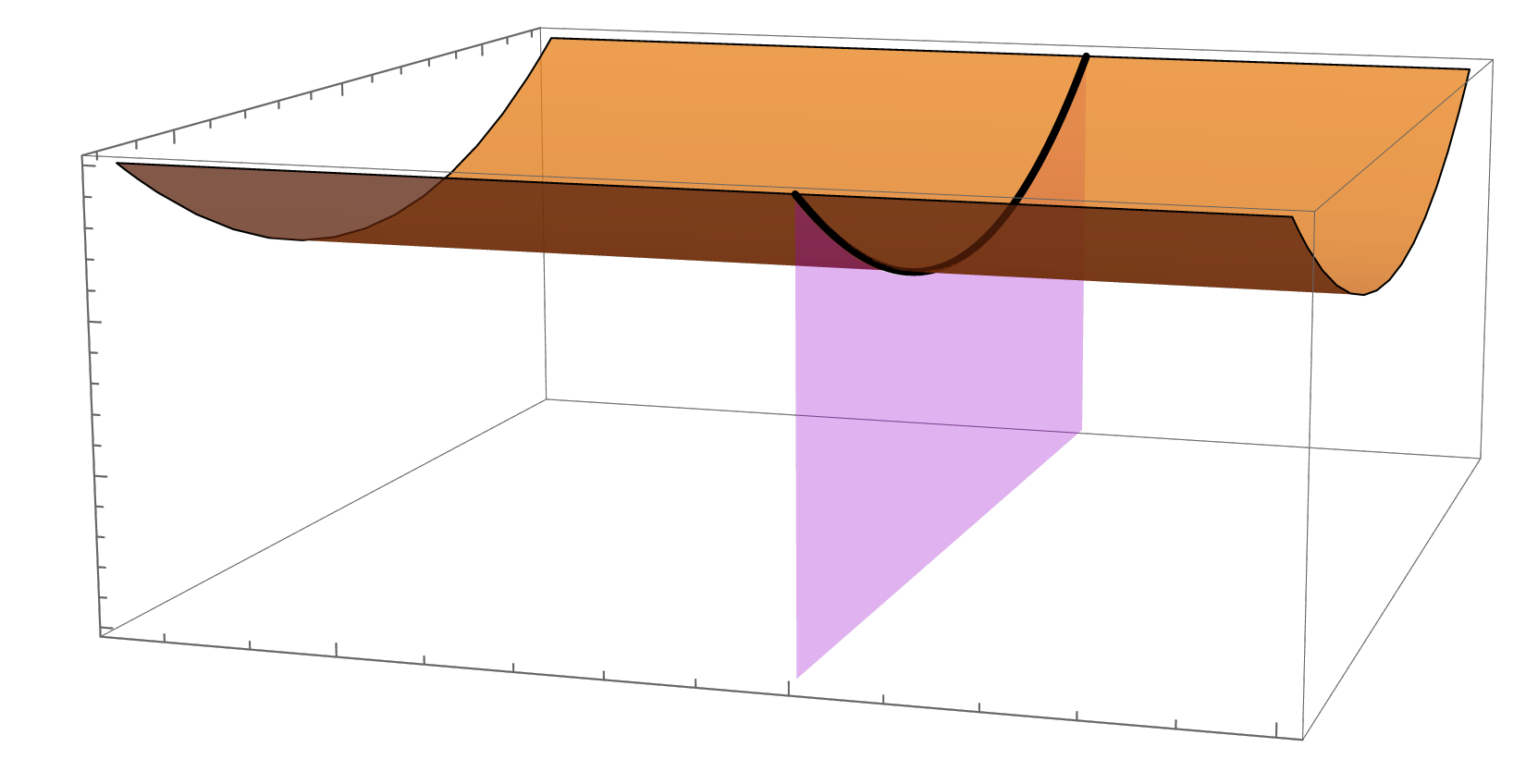}
\end{minipage}
\hspace{.01\linewidth}
\begin{minipage}{.3\linewidth}
  \includegraphics[width=\linewidth]{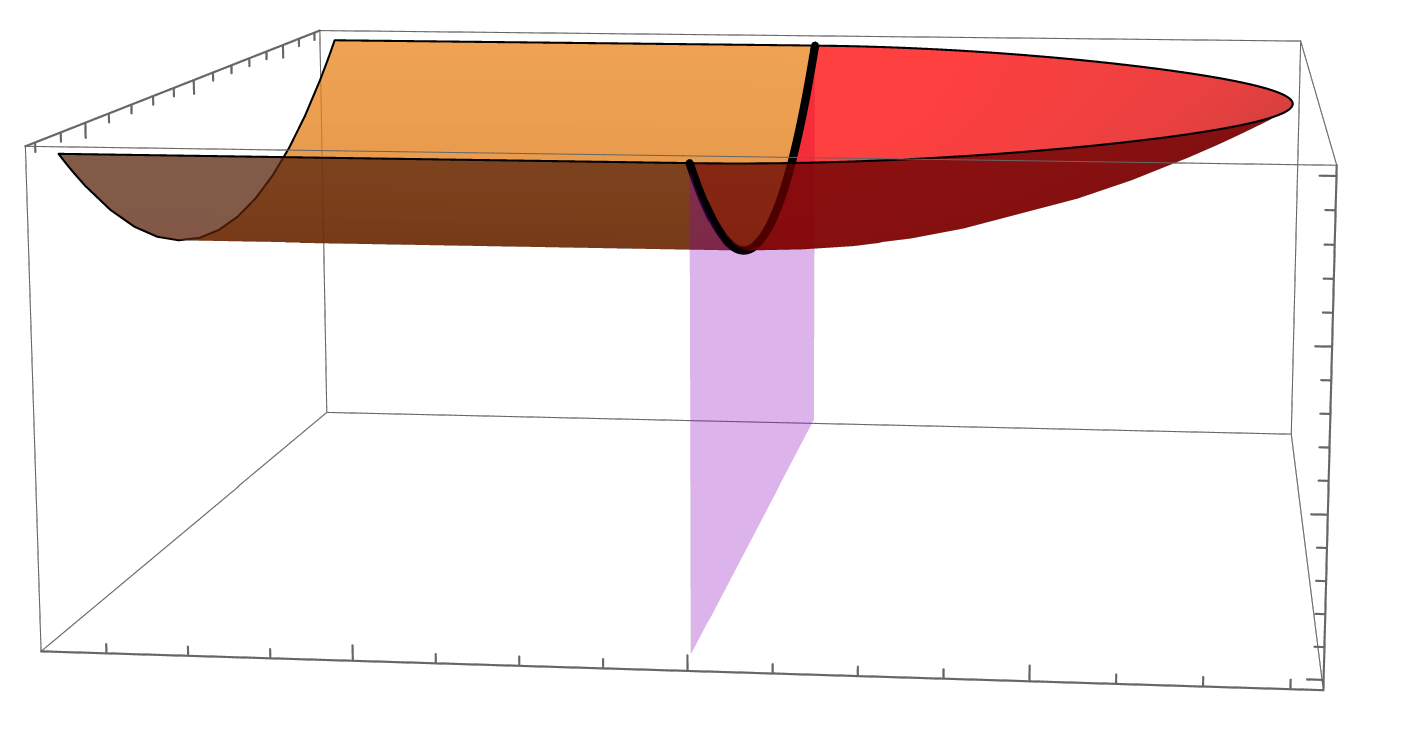}
\end{minipage}
\hspace{.01\linewidth}
\begin{minipage}{.3\linewidth}
  \includegraphics[width=\linewidth]{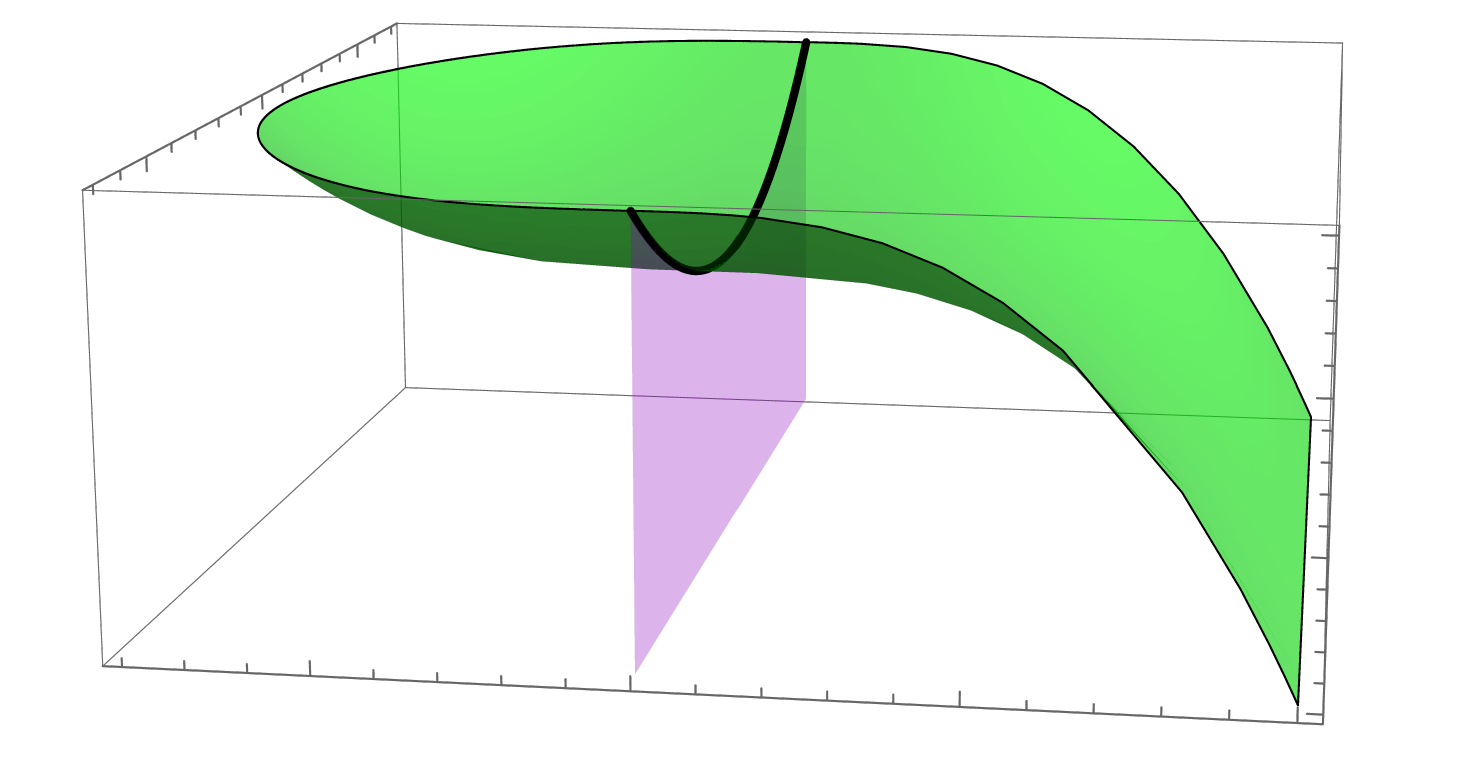}
\end{minipage}
\vspace{-.1 in}
\caption{ \underline{Left}. The spacetime region with the future temporal boundary consisting of the orange pre-shock cylindrical surface passing through the co-dimension-$2$ pre-shock set, shown as the black curve.
\underline{Center}.  The spacetime region used for the analysis of the downstream part of the \MGHDB.  The red surface displays the singular hypersurface consisting of the level set $\{\Jg=0\}$, which emanates from the pre-shock (the black curve)
in the downstream direction. In orange we represent the pre-shock cylindrical surface which also emanates from the pre-shock, but in the upstream direction.
\underline{Right}. The spacetime region used for the analysis of the upstream part of the \MGHDB.  The green surface denotes the distinguished slow acoustic characteristic surface emanating from the pre-shock set (in black), which connects to the initial time slice in the downstream direction.
}
\label{fig:3spacetimes}
\end{figure}

The main results of this paper are 
Theorem~\ref{thm:main:shock}, Theorem~\ref{thm:main:DS}, and Theorem~\ref{thm:main:US}, corresponding to the three panes in Figure~\ref{fig:3spacetimes}. For convenience, we summarize in Theorem~\ref{thm:mainrough} a significantly abbreviated version of our main results.

\begin{theorem}[\bf Maximal globally hyperbolic development in a box] 
\label{thm:mainrough}
Consider the 2D Euler equations~\eqref{euler1} for arbitrary~$\gamma>1$, with $H^7(\TT^2)$-smooth isentropic initial data $(u_0,\sigma_0)$, bounded away from vacuum. Assume that the data is compressive in the $x_1$ direction and generic.\footnote{The precise assumptions on the Cauchy data are given in Section~\ref{cauchydata}, items~\eqref{item:ic:supp}--\eqref{item:ic:w0:d11:positive}. The set of such initial conditions forms an open set in the $H^7$ topology, as discussed in~Remark~\ref{rem:open:set:data}.} For $0<\eps \ll 1$, we assume that the initial dominant Riemann variable $w_0 := u_0^{1} + \sigma_0 = \OO(1)$ is assumed to have a point at which $\p_1 w_0$ attains its global (non-degenerate) minimum, this maximally negative slope in the $x_1$ direction is $\OO(-\frac{1}{\eps})$, while in the $x_2$ direction the slope of $w_0$ is $\OO(1)$. Assume that the initial sub-dominant Riemann variable $z_0 := u_0^1 - \sigma_0$ and the initial tangential velocity $a_0 = u_0^2$ are $ \OO(\eps)$ in amplitude and that their derivatives satisfy $(\p_1 z_0,\p_1 a_0) = \OO(1)$ and $(\p_2 z_0,\p_2 a_0) = \OO(\eps)$.  The initial condition  for the geometry is $(\nn,\tt,\Jg)=(e_1, e_2,1)$. Then, assuming that $\eps$ is sufficiently small, the following hold:
\begin{enumerate}[leftmargin=26pt]

\item There exists a spacetime $\mathcal{M}_{\sf Eulerian}$, the \MGHDB\ of the Cauchy data $(u_0,\sigma_0)$ prescribed at the initial time slice, and a  unique solution $(u,\sigma)$ of~\eqref{euler1} in this spacetime, which propagates the regularity of the initial data; in particular, $(u,\sigma) \in C^0_t H^7_y \cap C^7_t L^2_y$.

\item There exists a family of Arbitrary-Lagrangian-Eulerian (ALE) diffeomorphisms $\psi(\cdot,t)$, indexed by time and defined in Section~\ref{sec:intro:ALE}, which flattens all fast acoustic characteristic surfaces. Under the action of $\psi$, the \MGHDB\ spacetime gets mapped into its ALE counterpart, the spacetime $\mathcal{M}_{\sf ALE}$.

\item The Euler evolution~\eqref{euler1} is equivalent to the evolution of the differentiated ALE Riemann variables $(\Wb,\Zb,\Ab)$, of the geometry $(\nn,\tt,\Jg)$ and of the (rescaled) sound speed $\Sigma$, cf.~\eqref{nn-approx}, \eqref{tt-approx}, and~\eqref{rest-approx}. These new geometric unknowns propagate the regularity of their initial data throughout the spacetime $\mathcal{M}_{\mathsf{ALE}}$; for instance, $\Jg, h,_2, \Jg\Wbn,\Zbn,\Abn,\Wbt,\Zbt,\Abt$ are bounded in $C^0_t H^6_x$, whereas $\psi$ and $\Sigma$ remain bounded in $C^0_t H^7_x$.

\item  The spacetime $\mathcal{M}_{\mathsf{ALE}}$ is the \MGHDB\ of the Cauchy data for the ALE  dynamics~\eqref{nn-approx}, \eqref{tt-approx}, and~\eqref{rest-approx}. The future temporal boundary (``top boundary'') of 
$\mathcal{M}_{\mathsf{ALE}}$ consists of: a co-dimension-$2$ set of ``first gradient singularities'', which is the set of pre-shocks 
$\Xi^* = \{ \Jg =0 \} \cap \{ \Jg,_1 = 0\}$; a co-dimension-$1$ singular 
surface which emerges from the pre-shock set in the downstream region, the set $\{\Jg = 0\}$, which parametrizes a continuum of gradient catastrophes
for the density and the normal velocity; the distinguished slow acoustic characteristic co-dimension-$1$ surface emanating from the pre-shock set
in the upstream direction, which serves as a Cauchy horizon for the ALE Euler evolution. See Figures \ref{fig:2D:Lagrangian:a0} and \ref{fig:3spacetimes} below for a pictorial 
representation of the future temporal boundary of $\mathcal{M}_{\mathsf{ALE}}$.

\end{enumerate}
\end{theorem}

\subsection{Organization of the paper}

In Section \ref{sec:acoustic:geometry}, we introduce the Arbitrary Lagrangian Eulerian (ALE) coordinate system, adapted to the fast acoustic
characteristic surfaces.   
In Section \ref{sec:new:Euler:variables}, we introduce a new type of differentiated Riemann-type variable which is a linear combination of the
gradient of velocity and sound speed and the curvature of the fast acoustic characteristic surfaces (see the identities \eqref{good-unknowns}).  Our analysis will make use of these
new differentiated Riemann variables in the ALE coordinate system.   

In Section \ref{sec:thesetup}, we give a detailed description of  the open set of compressive and generic initial
conditions used for our analysis.    We then state the three main theorems of this work.   Theorem \ref{thm:main:shock} produces the unique
Euler solution up the pre-shock.   Theorem \ref{thm:main:DS} establishes the downstream part of the \MGHDB, and Theorem
\ref{thm:main:US} establishes the upstream part of the \MGHDB.  

In Section \ref{sec:formation:setup}, we explain the geometry of the spacetime region lying below the cylindrical-type surface corresponding
to the set of pre-shocks.   Specifically, we smooth the corner formed at the intersection of the set of pre-shocks with our final time-slice, and create
a new smooth cylindrical-type surface using a modification of the metric-scaled determinant $\Jg$.   We then remap time, $t \mapsto \s$, so as to flatten 
this surface, and redefine all of our variables to now be functions of $(x,\s)$.  We define the  Sobolev norms used for energy estimates and state the
pointwise bootstrap assumptions on low-order derivatives and the Sobolev bootstrap assumptions on high-order derivatives.
Section \ref{sec:first:consequences} provides some of the basic consequences of the pointwise bootstrap assumptions, including 
the fact that $\Jg\Wbn(x,t) \sim \p_1w_0(x)$, which is used on numerous occasions throughout our proof.

Section \ref{sec:geometry:sound:ALE} establishes the sixth-order Sobolev bounds for the geometric quantities $\Jg$, $\p_2h$, $\nn$, $\tt$,  
the sound speed
$\Sigma$, and the tangential reparameterization velocity $V$.
In Section \ref{sec:vorticity}, we prove sixth-order energy estimates for the vorticity.  The resulting vorticity bound allows us to produce
improved $L^2$ bounds for both the fifth-order and sixth-order derivatives of $\Jg\Abn$, $\Jg\Zbn$, and $\Jg\Wbn$.
In Section \ref{sec:pointwise:bootstraps}, we close all of the bootstrap bounds.

Sections \ref{sec:sixth:order:energy-tangential}--\ref{sec:sixth:order:energy} contain the complete set of energy estimates 
used to obtain uniform Sobolev bounds for Euler solutions 
in the spacetime bounded from above by the pre-shock cylindrical surface.   Using the $(x,\s)$ coordinates defined in \eqref{eq:t-to-s-transform:all},
and the spacetime gradient operator $\nbs$ defined in \eqref{nb-s}, sixth-order energy estimates on the tangential-component equation set 
\eqref{energy-WZA-tan-s} are performed in 
Section \ref{sec:sixth:order:energy-tangential}, providing uniform bounds for the tangential norms 
 $\widetilde{\mathcal{E}}_{5,\ttt}(\s)$,  $\widetilde{\mathcal{E}}_{6,\ttt}(\s)$,  $\widetilde{\mathcal{D}}_{5,\ttt}(\s)$, 
and $\widetilde{\mathcal{D}}_{6,\ttt}(\s)$  defined in \eqref{eq:norms:L2:first}.  In Section \ref{sec:sixth:order:energy}, energy estimates for
the normal-component equation set \eqref{energy-WZA-s} are given, with uniform bounds for the norms 
 $\widetilde{\mathcal{E}}_{5,\nnn}(\s)$,  $\widetilde{\mathcal{E}}_{6,\nnn}(\s)$,  $\widetilde{\mathcal{D}}_{5,\nnn}(\s)$, 
and $\widetilde{\mathcal{D}}_{6,\nnn}(\s)$  defined in \eqref{eq:norms:L2:first}.   The closure of these normal-component energy estimates
relies upon improved bounds for the case of six pure time-derivatives, and these improved bounds are proven in Section \ref{sec:pure:time}.

In Section \ref{sec:downstreammaxdev}, we give the construction of the downstream part of the \MGHDB.  
With the modified determinant function $\Jgb$ defined in \eqref{eq:Jgb:identity:0}, we consider the spacetime with
future temporal boundary given by the level-set $\{\Jgb(x,t)=0\}$ in the downstream region, and the pre-shock cylinder in the upstream region.
A new coordinate transformation \eqref{eq:t-to-s-transform:all-P} is introduced which maps $(x,t)\mapsto (x,\s)$ and flattens this future
temporal boundary.   Again, sixth-order energy estimates are closed for the normal-component equations, improved normal-component estimates
are obtained for the case of six pure time derivatives, and 
tangential-component estimates are closed.  The resulting  sixth-derivative uniform bounds establishes the existence of the unique Euler solution for all times
up to the singular boundary where the  metric-scale determinant $\Jg$ vanishes and thus  where gradient blow-up occurs.

The upstream part of the \MGHDB\ is established in  Section \ref{sec:upstreammaxdev}.  A large portion of the upstream spacetime
region is foliated by (what are essentially) slow acoustic characteristic surfaces.    The natural {\it time evolution} of the slow characteristic surfaces is 
somewhat singular when written in our ALE coordinates,  adapted to the fast acoustic characteristic surfaces.  In particular, the temporal
rate of change of each slow characteristic surface is proportional to $\Jgi$ which blows-up at the pre-shock.   As such, we introduce a 
reparameterization of these slow characteristic surfaces but employing the $x_1$ coordinate as the evolutionary independent variable. The resulting
description of the geometry of the slow characteristic surfaces becomes smooth.   We define the weight function $\JJ$ used in upstream energy
method via transport along these slow surfaces of the value of $\Jg$ along the fast characteristic surface passing through the pre-shock.  With
this weight function, we once again close sixth-order energy estimates for the normal-component equations, improved normal-component estimates
are once again obtained for the case of six pure time derivatives, and  tangential-component estimates are again closed. 
Uniform bounds are therefore established in the entire upstream spacetime region, lying below the slow characteristic surface that emanates from
the pre-shock.

Section \ref{sec:optimal:reg} is devoted to establishing the optimal regularity of the velocity, sound speed, and the ALE family of diffeomorphism
$\psi$.   We prove that $H^7$ Sobolev regularity is maintained for the entire development of the  Cauchy data.

In  Appendix \ref{sec:usersguide}, we provide the reader a  self-contained introduction to the notion of \MGHDB\ of Cauchy data for
the Euler equations in the simplified setting of one space dimension.  A complete description  is provided in both the traditional Eulerian setting, as well as the more geometric Lagrangian framework.

Appendix \ref{app:functional} is devoted to the basic functional analysis lemmas in the three different spacetime regions that we employ for our
analysis.  We prove a number of technical lemmas which are spacetime variants of the classical Sobolev and Poincar\'{e} inequalities, the
Gagliardo-Nirenberg inequalities,  Moser inequalities,  and a number of commutator lemmas.

Finally, in 
Appendix \ref{sec:app:transport}, we establish $L^\infty $ bounds for solutions to certain transport equations, by obtaining $p$-independent bounds
for $L^p$ energy estimates and passing to the limit as $p\to \infty $.   Keeping in mind that our ALE coordinate system is adapted to the
fast characteristics associated to the $\lambda_3$ wave speed, the lemmas in this section allow us to obtain pointwise bounds for quantities which are
naturally associated with either the $\lambda_1$ or $\lambda_2$ wave speeds.


\section{Acoustic characteristic surfaces associated to shock formation}
\label{sec:acoustic:geometry}
In this section, we develop the geometry for shock formation.  We shall study the problem on the spacetime $ \TT^2  \times [\initial, T]$, where the initial time $\initial$ shall be made precise below, and $T \in (\initial,\final]$ is arbitrary. The times $\initial$ and $\final$ are defined in~\eqref{tin}, respectively in~\eqref{time-of-existence}.

\subsection{Characteristic surfaces}
\label{sec:ALE}
We denote the $3$-transport velocity by 
\begin{equation} 
\mathcal{V}_3 = u + \alpha \sigma n =  \lambda_3 n + (u \cdot \tau) \tau \,. \label{V3}
\end{equation} 
We shall foliate spacetime with acoustic characteristic surfaces associated to the ``fast'' wave speed $\lambda_3$.
One way of doing so is by studying the Lagrangian flow map of the $3$-transport velocity $\mathcal{V} _3$:
\begin{equation}
 \p_t \eta(x_1,x_2,t)  = \mathcal{V} _3(\eta(x_1,x_2,t),t) \,, 
 \qquad 
\eta( x , \initial) = x  \,.
\label{eta-char}
\end{equation} 
In terms of the standard Cartesian basis, we have that
\begin{equation*}
\eta(x_1,x_2,t) = (\eta^1(x_1,x_2,t),\eta^2(x_1,x_2,t))\,,
\end{equation*}
and that
\begin{equation*}
\eta^1(x_1,x_2,\initial) = x_1 \ \ \text{ and } \ \ \eta^2(x_1,x_2,\initial) = x_2 \,.
\end{equation*}
Using the flow map $\eta$ we can  give a geometric description of the fast acoustic characteristics surfaces.

At initial time $t=\initial$, we foliate $\TT^2$ by lines parallel to $e_2=(0,1)$, and denote these lines by
$ \gamma_{x_1}(\initial) = \{x_1\} \times \TT$.   For each $x_1 \in \TT  $ and $t\in[ \initial,T]$, we define the 
characteristic curve (at a fixed time-slice) by
\begin{equation} 
\gamma_{x_1}(t) = \eta( \gamma_{x_1}(\initial), t)  \,,  \label{Cx1t}
\end{equation} 
and the characteristic surfaces up to time $T \ge \initial$ (which are parameterized of $x_1$) by 
\begin{equation} 
\Gamma_{x_1}(T) = \bigcup\nolimits_{t\in[\initial,T]} \gamma_{x_1}(t) \,. \label{Cx1}
\end{equation} 
Figure~\ref{fig:surf} below displays a few such characteristic surfaces $\Gamma_{x_1}$ for five different values of $x_1 \in \TT$.
Figure~\ref{fig:surf_new} displays the fast-acoustic characteristic surfaces $\Gamma_{x_1}$ {\em jointly} with the slow-acoustic characteristic surfaces emanating from the same initial foliation of $\TT^2$.

\begin{figure}[htb!]
\centering
\includegraphics[width=.6\linewidth]{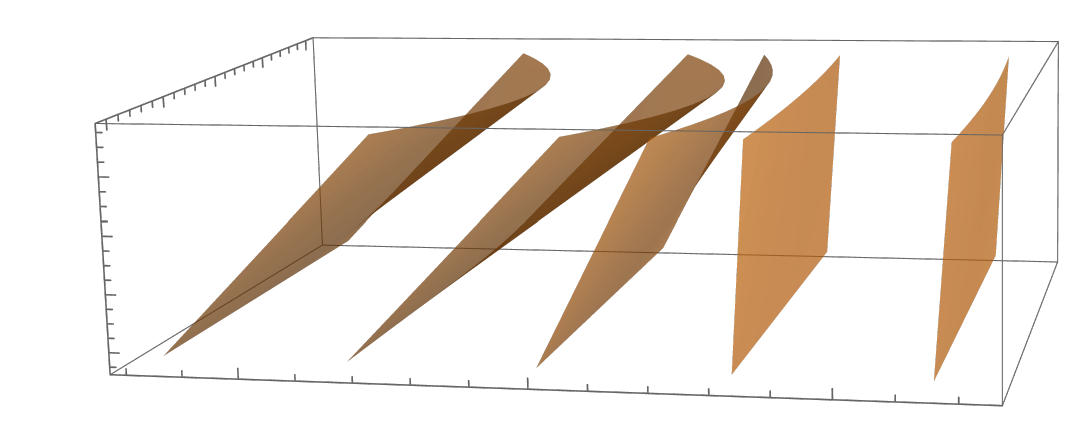}
\vspace{-.1in}
\caption{For $T>\initial$ which is strictly less than the very first blow-up time, we display the characteristic surfaces $\Gamma_{x_1}(T)$ defined in~\eqref{Cx1} emanating from five different values of $x_1\in\TT$. At $t= \initial$, the curves $\{ \gamma_{x_1}(\initial) \}_{x_1\in\TT}$ are lines which foliate $\TT^2$. The distance between the characteristic surfaces $\Gamma_{x_1}(T)$ is decreasing as $T$ increases leading to shock formation when this distance vanishes.}
\label{fig:surf}
\end{figure}

\begin{figure}[htb!]
\centering
\includegraphics[width=.6\linewidth]{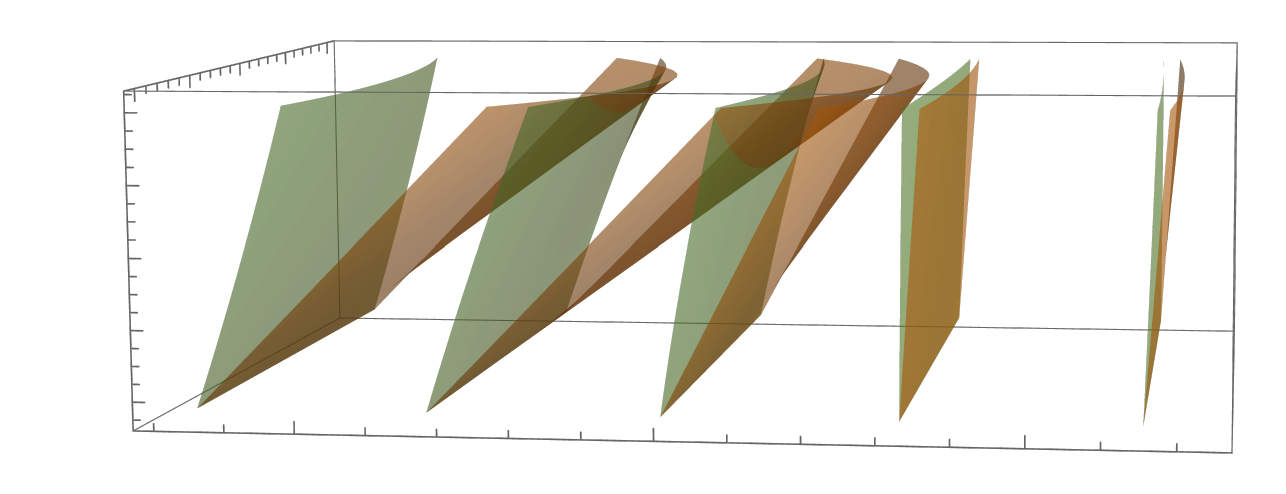}
\vspace{-.1in}
\caption{For $T>\initial$ which is strictly less than the very first blow-up time, and for five different values of $x_1\in\TT$, we display both the fast acoustic characteristic surfaces $\Gamma_{x_1}(T)$ (in orange), and the corresponding slow acoustic characteristic surfaces emanating from the same values of $x_1$ (in olive-green). While the orange fast acoustic characteristic surfaces are close-to-impinging on each other, the olive-green slow acoustic surfaces smoothly foliate spacetime.}
\label{fig:surf_new}
\end{figure}

\subsection{An Arbitrary-Lagrangian-Eulerian (ALE) description of the geometry}
While the Lagrangian flow $\eta$ is a natural  parameterization for the fast acoustic characteristic surfaces $\Gamma_{x_1}$, it is convenient to introduce a tangential re-parameterization in the form of the so-called Arbitrary-Lagrangian-Eulerian (ALE) coordinates.   

Because each curve $\gamma_{x_1}(t)$ is a graph over the set $\TT \ni x_2$, we introduce a {\it height} function $h(x_1,x_2,t)$ such that
\begin{equation*}
\gamma_{x_1}(t) = \{ (h(x_1,x_2,t),x_2) \colon x_2 \in \TT \}, \qquad t \in [\initial ,T] \,.
\end{equation*}
The induced metric on $\gamma_{x_1}(t)$  is given by
\begin{equation} 
g(x_1,x_2,t) = 1 + | h,_2(x_1,x_2,t)|^2 \,, 
\label{g-def}
\end{equation} 
and the unit tangent vectors $\tt$ and normal vectors $\nn$ to the curves $\gamma_{x_1}(t) $ are  then given by
\begin{equation} 
\tt(x_1, x_2,t) = g^{- {\frac{1}{2}} } ( h,_2, 1) \,, \qquad \mbox{and} \qquad \nn(x_1, x_2,t) = g^{- {\frac{1}{2}} } (1, -h,_2) \,. 
\label{tn-lag}
\end{equation}

We define the ALE family of maps $\psi(x_1,x_2,t) = (\psi^1(x_1,x_2,t),\psi^2(x_1,x_2,t))\,,$ 
by\footnote{The tangent and normal vectors to $\gamma_{x_1}(t)$ can be equivalently defined via the map $\psi$.  In particular, we have that
$\tt = g^{- {\frac{1}{2}} } \psi,_2 = g^{- {\frac{1}{2}} } (h,_2,1)$, $ \nn = g^{- {\frac{1}{2}} }  \psi,_2^\perp  = g^{- {\frac{1}{2}} } (1, -h,_2)$, and the induced
metric $g = \psi,_2 \cdot \psi,_2 = 1+ |h,_2|^2$.}
\begin{equation} 
\psi(x_1,x_2,t) = h(x_1,x_2,t) e_1 + x_2 e_2 \,,   \label{psi-def}
\end{equation} 
where 
\begin{equation}
h(x_1,x_2,\initial)= x_1\,.
\label{eq:h:initial:condition}
\end{equation}
In order to preserve the shape of the characteristic surfaces $\Gamma_{x_1}(T)$,  the family of diffeomorphisms 
$\psi( \cdot , t)$ must satisfy the constraint  
\begin{equation} 
\p_t \psi \cdot \nn = (\mathcal{V} _3 \circ \psi) \cdot \nn =  ( u+ \alpha \sigma n) \circ \psi \cdot \nn \,.  
\label{psit-def}
\end{equation} 
Time-differentiating \eqref{psi-def}, we have that
\begin{equation*}
\p_t \psi \cdot \nn = \p_t h (e_1 \cdot \nn) = g^{- {\frac{1}{2}} } \p_t h \,,
\end{equation*}
and from \eqref{psit-def}, we have that
\begin{equation} 
\p_t h = g^ {\frac{1}{2}} \left( (u \circ \psi) \cdot \nn + \alpha \sigma \circ \psi \right) \,.
\label{eq:psi:evo:def}
\end{equation} 
Similarly, 
\begin{equation*}
\p_t \psi \cdot \tt = \p_t h (e_1 \cdot \tt) = g^{- {\frac{1}{2}} } h,_2 \p_t h \,.
\end{equation*}
It follows that
\begin{align} 
\p_t \psi 
&= (\p_t \psi \cdot \nn) \nn + ( \p_t \psi \cdot \tt) \tt  
\notag \\
&= \left( (u \circ \psi) \cdot \nn + \alpha \sigma \circ \psi \right) \nn +
\left( (u \circ \psi) \cdot \nn + \alpha \sigma \circ \psi \right) h,_2 \tt \,. 
\label{pt-psi}
\end{align}

\subsection{The  deformation matrix $\nabla \psi$ and its determinant and inverse}
\begin{figure}[htb!]
\centering
\begin{minipage}{.45\linewidth}
\includegraphics[width=\linewidth]{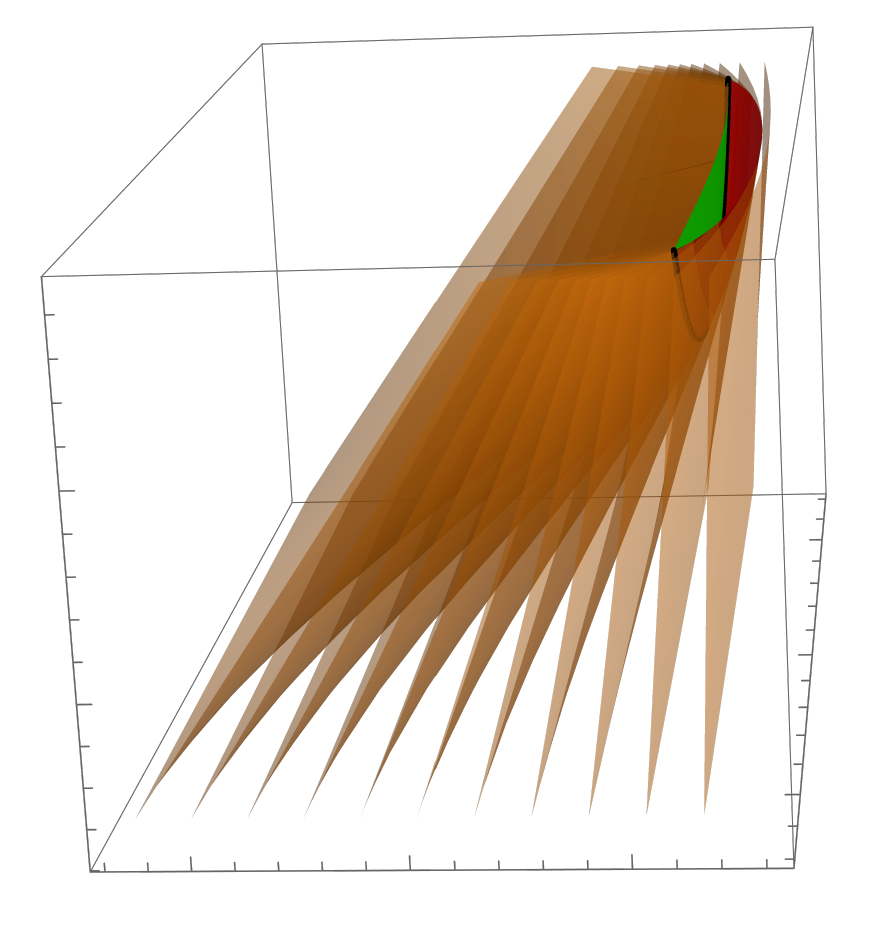}
\end{minipage}
\hspace{.02\linewidth}
\begin{minipage}{.45\linewidth}
\includegraphics[width=\linewidth]{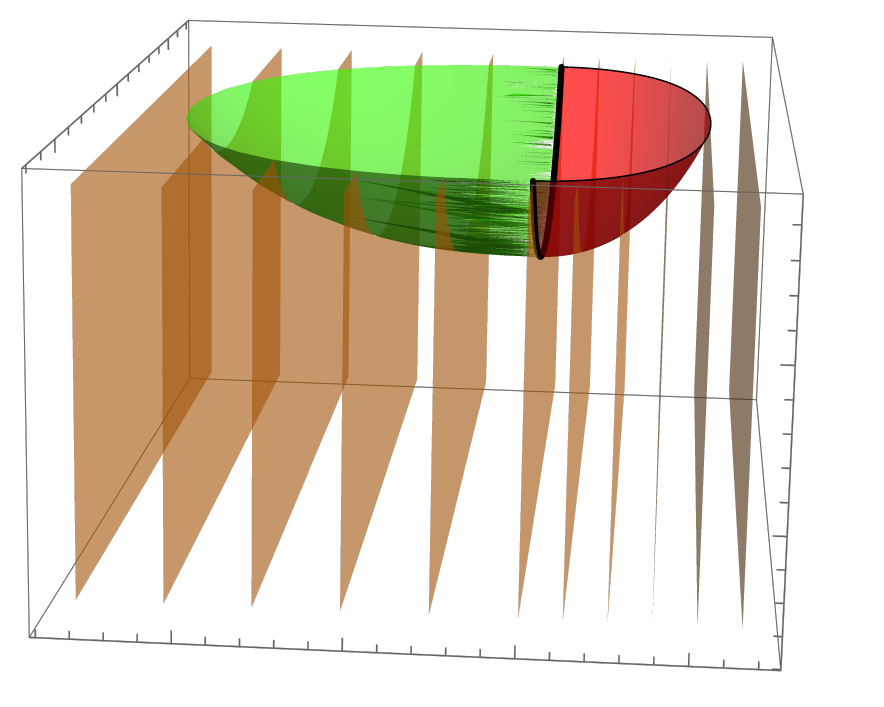}
\end{minipage}
\vspace{-.1 in}
\caption{\underline{Left}: in Eulerian coordinates, we represent the 2D analogue of the 1D image from Figure~\ref{fig:1D:Eulerian} (below). For $T \in [\initial,\final]$ we display in orange the fast acoustic characteristic surfaces $\Gamma_{x_1}(T)$ for various values of $x_1$. The black curve, parametrized by $x_2$ represents the set of pre-shocks in Eulerian coordinates (see Definition~\ref{def:pre-shock} below). The green surface is the slow acoustic characterstic surface emanating from the pre-shock set in the upstream direction. The cuspoidal red surface, which also emanates from the curve of pre-shocks, is the surface of ``first singularities'', the envelope (or ``top'' boundary) of the spacetime region in which the fast acoustic characteristic surfaces remain in one-to-one correspondence with the initial foliation of spacetime. The Eulerian \MGHDB\  of the initial data is the spacetime which lies in the ``temporal past'' of the union of the green and red surfaces (and the time slice $\{t = \final\}$). This spacetime has a boundary which is not smooth; meaning, not differentiable, it has limited H\"older regularity.
\underline{Right}: in ALE coordinates, we represent the 2D analogue of the 1D image from Figure~\ref{fig:1D:Lagrangian} (below). We represent the same surfaces as in the left figure, except that we are composing with the ALE diffeomorphism $\psi$. In these ALE coordinates, the fast acoustic characteristics are flat planes which foliate spacetime, parametrized by $x_1$. The spacetime of \MGHDB\ of the initial data in ALE coordinates is the spacetime which lies ``below'' the union of the green and red surfaces (and the time slice $\{t = \final\}$). In turn, the red surface is the portion of the paraboloid of ``first singularities'' $\{\Jg = 0\}$ which lies downstream of the curve of pre-shocks, whereas the green surface represents the upstream part of the slow acoustic characteristic surface emanating from the pre-shock, composed with the flow $\psi$. The boundary of the ALE spacetime is   $W^{2,\infty}$ smooth.}
\label{fig:2D:Eulerian:a}
\end{figure}

The diffeomorphisms $\psi$ are fundamental in our analysis, since the definition of the paraboloid of first singularities, which describes the downstream \MGHDB, is determined by the vanishing of the Jacobian determinant of $\psi$  (see Figure \ref{fig:2D:Eulerian:a}).

From \eqref{psi-def} we have that
\begin{equation*} 
\nabla \psi = 
\left[
\begin{matrix}
h,_1  & h,_2 \\
0 &  1
\end{matrix}\right] \,,
\end{equation*} 
so that the Jacobian determinant is given by
$$
J = \det \nabla \psi = h,_1 \,.
$$
We introduce the metric-normalized Jacobian determinant as
\begin{align} 
\Jg =  g^ {-\frac{1}{2}} J = g^{- {\frac{1}{2}} } h,_1  \,.   \label{Jg-def}
\end{align} 
The paraboloid of first singularities on the right side of Figure~\ref{fig:2D:Eulerian:a} will be shown to be the level set $\{ \Jg=0\}$.

The cofactor matrix of $ \nabla \psi$ is denoted by
\begin{equation*} 
b= 
\left[
\begin{matrix}
1  & -h,_2 \\
0 &  h,_1
\end{matrix}\right] \,,
\end{equation*} 
and the inverse matrix $B(x,t) = [\nabla \psi(x,t)] ^{-1} $ is defined by
\begin{equation*}
B = J ^{-1}  b \,.
\end{equation*}
The components of $b$ are denoted by $b^i_j$, the upper index for the row, and the lower index for the column.
It is important to observe that
\begin{equation}
b^1_j = {\psi,_2^j}^\perp = g^ {\frac{1}{2}} \nn^j \,, \qquad\mbox{and}\qquad b^2_j =  (0, h,_1) =:  J e_2^j \,,   
\label{b-comp}
\end{equation} 
so that
\begin{equation}
B^1_j = \Jgi \nn^j \,, \qquad \mbox{and}\qquad B^2_j =    \delta^2_j \,.
\label{B-comp}
\end{equation} 
As usual, the columns of $b$ are divergence-free.

\subsection{The relationship between ALE and Eulerian derivatives}
The associated Eulerian unit tangent and normal vectors are given by 
\begin{equation} 
 \tau  = \tt \circ \psi^{-1}  \,, \ \ \ n =  \nn \circ \psi^{-1}  \,.  
 \label{Eulerian-geom}
\end{equation} 
Suppose that $f$ denotes a differentiable Eulerian function, and that let  $F = f \circ \psi$ be the associated ALE function.
Using \eqref{B-comp} together with the chain-rule,  we have that
\begin{equation*}
f,_i \circ \psi = F,_k B^k_i  = \Jgi F,_1 \nn^i + F,_2 e_2^i \,.
\end{equation*}
Since $e_2 \cdot \nn = -g^{- {\frac{1}{2}} } h,_2$ and $e_2 \cdot \tt = g^{- {\frac{1}{2}} } $,  it follows that 
\begin{subequations} 
 \label{dtran}
\begin{align} 
\p_n f \circ \psi  & =  \Jgi F,_1 - g^{- {\frac{1}{2}} } h,_2  F,_2  \,,  \label{dtran-n}\\
\p_\tau f \circ \psi &= g^{-\frac{1}{2}} F,_2   \,.  \label{dtran-tau}
\end{align} 
\end{subequations} 
The identities in~\eqref{dtran} show that the ALE coordinate system characterizes the ``tame'' tangential derivatives of $f$ simply as derivatives with respect to $x_2$ for $f\circ \psi$. The ``singular'' nature of the normal derivatives of $f$ is characterized not as much by the fact that these are derivatives with respect to the $x_1$ or $x_2$ direction; it is the presence of the $\Jgi$ term (which blows up as $\Jg \to 0$)  in front of $F,_1$ which fully characterizes the singular nature of normal derivatives.

\section{A new set of variables for Euler shock formation}
\label{sec:new:Euler:variables}

Having introduced a new system of ALE coordinates, we now introduce a new set of variables that are essential for the analysis of the shock formation process and of the \MGHDB.

\subsection{Euler equations in geometric ALE coordinates}
\subsubsection{The differentiated acoustic Euler equations}
We first compute the evolution of partial (space) derivatives of $u$ and $\sigma$.  We differentiate \eqref{euler-shock} and find that
\begin{subequations} 
\label{euler-d1}
\begin{align} 
\p_t u^i,_k + (u\cdot n + \alpha \sigma) u^i,_{kj} n^j + (u \cdot \tau)  u^i,_{kj}\tau^j + \alpha \sigma (  \sigma,_{ki} -  \p_n u^i,_k) + \alpha \sigma,_k \sigma,_i + u^j,_k u^i,_j & = 0 \,, \\
\p_t \sigma,_k +  (u\cdot n + \alpha \sigma) \sigma,_{kj}n^j +(u \cdot \tau) \sigma,_{kj}\tau^j + \alpha \sigma (  u^i,_{ki} -  \p_n \sigma,_k)+ \alpha \sigma,_k u^i,_i + u^j,_k \sigma,_j & =0 \,.
\end{align} 
\end{subequations} 
The system \eqref{euler-d1} constitutes the differentiated acoustic Euler equations.  It is imperative to study this differentiated form to avoid 
derivative loss in the geometry.

\subsubsection{ALE variables} Using our family of ALE mappings $\psi$ defined in \eqref{psi-def}, we define 
a new set of  ALE variables representing velocity, sound speed, and their gradients along the acoustic characteristic surfaces by
\begin{subequations} 
\label{good-US}
\begin{alignat}{2}
U^i& = u^i \circ \psi \,, \qquad  &&\Sigma = \sigma \circ \psi \,,   \label{vsL} \\
\Uik&=  u^i,_k \circ \psi \,, \qquad && \mathring\Sigma_k =  \sigma,_k \circ \psi \,.  \label{USigma}
\end{alignat} 
Additionally, we define the ALE  wave-speed 
\begin{equation} 
 \Lambda_3 = \lambda_3 \circ \psi  = U \cdot \nn + \alpha \Sigma \,,  \label{lambda3-ale}
\end{equation} 
\end{subequations} 
With respect to the variables \eqref{good-US}, the system \eqref{euler-d1} takes the form
\begin{subequations} 
\label{euler-ALE}
\begin{align} 
&\p_t \Uik + g^{- {\frac{1}{2}} } \left(U\cdo\tt - \Lambda_3 h,_2 \right) \Uik,_2  \notag \\
& \qquad
 +\alpha \Sigma \Jgi ( \Sk,_1 \nn^i - \Uik,_1  ) 
  +\alpha \Sigma ( \Sk,_2 e_2^i + g^{- {\frac{1}{2}} }  \Uik,_2  h,_2)
+ \alpha \Sk \Si + \Uij\Ujk =0 \,,  \label{Euler-ALEa}\\
&\p_t \Sk  + g^{- {\frac{1}{2}} } \left(U\cdo\tt - \Lambda_3 h,_2 \right) \Sk,_2 \notag \\
& \qquad
 + \alpha \Sigma \Jgi (\Uik,_1  \nn^i - \Sk,_1)  + \alpha \Sigma g^{- {\frac{1}{2}} } ( -\Uik,_2 \nn^i h,_2 + \Sk,_2 h,_2 + \Uik,_2 \tt^i)
+ \alpha \Sk\Uii +\Sj\Ujk =0 \,. \label{Euler-ALEb}
\end{align} 
\end{subequations} 

\subsubsection{ALE Riemann variables}
We define the ALE  Riemann variables
\begin{equation} 
W= U \cdot \nn + \Sigma \,, \qquad  Z = U \cdot \nn - \Sigma\,, \qquad A = U \cdot \tt \,,   \label{WZA-ALE}
\end{equation} 
and a new set of Riemann-type variables (which are of fundamental importance to our analysis) by
\begin{equation} 
\label{BigWZA}
\Wk = \nn^i\Uik  + \Sk  \,, \qquad \Zk = \nn^i \Uik - \Sk   \,, \qquad \Ak =\tt^i \Uik \,.
\end{equation} 
We note that $W$, $Z$, and $A$ are scalar functions, while $\Wk$, $\Zk$, and $\Ak$ are the $k^{\text{th}}$ components of vector functions.
The use of the ALE mappings $\psi$ to parameterize the acoustic characteristic surfaces introduces a tangential reparameterization with
the following  transport velocity
\begin{equation} 
V 
= g^{- {\frac{1}{2}} } \bigl( A - \Lambda_3 h,_2 \bigr)  
= g^{- {\frac{1}{2}} } \big( A - h,_2  (\tfrac{1+\alpha}{2} W + \tfrac{1-\alpha}{2} Z  ) \big) \,.   \label{transport-ale}
\end{equation} 
In order to obtain the evolution equations for $\Wk$, $\Zk$, and $\Ak$, we must first compute the dynamics of 
the normal and tangential components of \eqref{Euler-ALEa}.  We find that
\begin{subequations} 
\label{euler-WZA}
\begin{align} 
&
\Jg\Bigl( \nn^i(\p_t + V\p_2)\Uik +(\p_t + V\p_2) \Sk\Bigr)
+ \alpha \Jg  \Sigma g^{- {\frac{1}{2}} }\tt^i  \Uik,_2 
 \notag \\
& \qquad\qquad \qquad\qquad
+  \Jg ( \tfrac{1+ \alpha }{2}  \Wk +  \tfrac{1- \alpha }{2} \Zk ) \Wi \nn^i
+ \Jg \Ak \Wi \tt^i + \tfrac{\alpha }{2} \Jg (\Wk-\Zk) \Ai \tt^i =0 \,, \label{euler-W} \\
&
\Jg\Bigl( \nn^i(\p_t + V\p_2)\Uik -(\p_t + V\p_2) \Sk\Bigr)
- \alpha \Jg  \Sigma g^{- {\frac{1}{2}} } \tt^i \Uik,_2  \notag \\
& \qquad\qquad \qquad\qquad
 -2 \alpha \Sigma \Bigl(\nn^i\Uik,_1  - \Sk,_1\Bigr) + 2 \alpha \Sigma \Jg g^{- {\frac{1}{2}} }h,_2  \Bigl( \nn^i\Uik,_2 - \Sk,_2 \Bigr)
\notag \\
& \qquad\qquad \qquad\qquad
+ \bubu{  \Jg ( \tfrac{1- \alpha }{2} \Wk + \tfrac{1+ \alpha }{2} \Zk ) \Zi \nn^i}
+ \Jg \Ak \Zi \tt^i - \tfrac{\alpha }{2} \Jg (\Wk-\Zk)\Ai\tt^i
=0  \,, \label{euler-Z} \\
&
\Jg \tt^i (\p_t + V\p_2) \Uik 
+ \alpha \Sigma g^{- {\frac{1}{2}} } \Jg  \Sk,_2 
 - \alpha \Sigma \tt^i \Uik,_1 
+ \alpha \Sigma g^{- {\frac{1}{2}} } \Jg  h,_2 \tt^i\Uik,_2   \notag \\
& \qquad\qquad \qquad\qquad
+  \Jg\Bigl(  \tfrac{1}{2} (\Wk + \Zk) \Ai\nn^i + \Ak \Ai \tt^i + \tfrac{\alpha }{4}(\Wk-\Zk)(\Wi-\Zi)\tt^i   \Bigr) 
=0  \,, \label{euler-A} 
\end{align} 
\end{subequations} 
which is coupled with the transport-type equation for $h$ given by
\begin{equation} 
\p_t h  = g^{{\frac{1}{2}} } (U \cdot \nn + \alpha \Sigma) = g^{{\frac{1}{2}} } ( \tfrac{1+\alpha}{2} W + \tfrac{1-\alpha}{2} Z) \,, \label{h-evo}
\end{equation} 
and the evolution equations for $(W,Z,A)$, which are given in~\eqref{WZA-evo} below. 

\subsubsection{Tensor notation}
We use  bold notation for tensors, and write \eqref{BigWZA} as
\begin{equation*}
\Ub = \nabla u \circ \psi \,, \ \ \ \Sb = \nabla \sigma \circ \psi\,, \ \ \ \Wb = \nn   \Ub + \Sb\,, \ \ \Zb = \nn   \Ub - \Sb\,, \ \ \Ab = \tt \Ub \,.
\end{equation*}
By $\nn   \Ub$ we mean the contraction of the vector $\nn$ with the tensor $\Ub$ so that $\nn \Ub$ is a $1$-form, which we identify with a vector field, 
with components  defined by
$$(\nn \Ub)_k = \nn^i \Uik \,.$$
Similarly, the vector $\tt\Ub$ has components $(\tt \Ub)_k = \tt^i \Uik$.

We  next  introduce the following component notation 
$$
\Wbn = \Wb \cdot \nn \,, \ \ \Wbt = \Wb \cdot \tt \,, \ \ \Zbn = \Zb \cdot \nn\,, \ \ \Zbt = \Zb \cdot \tt\,, \ \ 
\Abn = \Ab \cdot \nn\,, \ \ \Abt = \Ab \cdot \tt \,,
$$
and this shall henceforth be used.

\subsubsection{Geometric significance of the variables $\Wb$, $\Zb$, and $\Ab$}
The variables $\Wb$, $\Zb$, and $\Ab$ have been designed to both encode certain components of the characteristic surface curvature tensor and
avoid geometric derivative loss.   The significance of these variables is demonstrated using classical Eulerian Riemann variables. 
Following our definition of the ALE Riemann variables in \eqref{WZA-ALE}, we define the Eulerian Riemann variables by
\begin{equation} 
w= u \cdot n + \sigma \,, \qquad  z = u \cdot n - \sigma\,, \qquad a = u \cdot \tau \,.   
\label{wza}
\end{equation} 
Then, the identities \eqref{dtran}, \eqref{WZA-ALE}, \eqref{BigWZA}, and \eqref{wza} provide the following: 
\begin{subequations} 
\label{component-identities}
\begin{align}
\Wbn 
&=  \Jgi W,_1 - h,_2 g^{-\frac 12}  W,_2
+ g^{-1} A (\Jgi h,_{12} - h,_2 g^{-\frac 12} h,_{22} ) 
=  (\p_n w -a \p_n n \cdot \tau) \cir \psi \,,
\label{comp-id-Wbn}
\\
\Zbn 
&=  \Jgi Z,_1 - h,_2 g^{-\frac 12}  Z,_2
+ g^{-1} A (\Jgi h,_{12} - h,_2 g^{-\frac 12} h,_{22} ) 
=   (\p_n z -a \p_n n \cdot \tau) \cir \psi \,,
\label{comp-id-Zbn}
\\
\Abn 
&=  \Jgi A,_1 - h,_2 g^{-\frac 12}  A,_2
- \tfrac 12 g^{-1} (W+Z) (\Jgi h,_{12} - h,_2 g^{-\frac 12} h,_{22} )
=    (\p_n a - \tfrac 12 (w+z)  \p_n \tau\cdot n )\cir\psi  \,,
\label{comp-id-Abn}
\\
\Wbt 
&= g^{-\frac 12} W,_2 + A g^{-\frac 32} h,_{22}
= (\p_\tau w - a  \p_\tau n \cdot \tau )\cir \psi \,,
\label{comp-id-Wbt}
\\
\Zbt 
&= g^{-\frac 12} Z,_2 + A g^{-\frac 32} h,_{22}
= (\p_\tau z -  a \p_\tau n \cdot \tau)\cir\psi \,,
\label{comp-id-Zbt}
\\
\Abt 
&= g^{-\frac 12} A,_2 - \tfrac 12 (W+Z) g^{-\frac 32} h,_{22} = \big(\p_\tau a - \tfrac{1}{2} (w+z) \p_\tau \tau \cdot n\big) \cir \psi
\,.
\label{comp-id-Abt}
\end{align}
\end{subequations} 

\subsubsection{Dynamics of geometry}
We will use the identities 
\begin{equation} 
D  \nn = - g^{-1} D h,_2 \tt \,, \qquad D  \tt =  g^{-1} D h,_2 \nn \,, \qquad D g= 2 h,_2 D h,_2\,,  \label{d-n-tau}
\end{equation} 
which hold for any differential operator $D$.
From \eqref{BigWZA}, we see that
\begin{equation} 
\nn \Ub + \alpha \Sb = \tfrac{1+ \alpha }{2} \Wb + \tfrac{1- \alpha }{2} \Zb \,.  
\label{d-Lambda3}
\end{equation} 
Differentiating \eqref{h-evo} and using the identity \eqref{d-Lambda3},  we obtain that
\begin{subequations} 
\label{d12-h}
\begin{alignat}{2} 
&(\p_t+V\p_2) h,_1=  h,_1  \bigl(\tfrac{1+ \alpha }{2} \Wb + \tfrac{1- \alpha }{2} \Zb\bigr) \cdot (\nn  + h,_2 \tt)\,, \ && \text{ in } \  \TT^2  \times [\initial, T]
 \,, \label{p1-h} \\
&(\p_t+V\p_2) h,_2=  g  \bigl(\tfrac{1+ \alpha }{2} \Wbt + \tfrac{1- \alpha }{2} \Zbt \bigr)\,,  \ && \text{ in } \  \TT^2  \times [\initial, T] \,, \label{p2-h} \\
& h,_1 =  1 \ \text{ and }\ h,_2 =0\, ,\ && \text{ on } \  \TT^2  \times \{t=\initial\} \,. \label{ics-h,1-h,2}
\end{alignat} 
\end{subequations} 
From \eqref{tn-lag} and \eqref{d12-h}, we then have that
\begin{subequations} 
\label{nn-tt-evo}
\begin{alignat}{2}
&(\p_t+V\p_2)\nn  + \bigl(\tfrac{1+ \alpha }{2} \Wbt + \tfrac{1- \alpha }{2} \Zbt \bigr) \tt  = 0\,, \ && \text{ in } \  \TT^2  \times [\initial, T] \,, 
 \label{nn-evo}\\
&(\p_t+V\p_2)\tt   - \bigl(\tfrac{1+ \alpha }{2} \Wbt + \tfrac{1- \alpha }{2} \Zbt \bigr)  \nn  = 0\,, \ && \text{ in } \  \TT^2  \times [\initial, T] \,,
 \label{tt-evo} \\
 & \nn =  e_1 \ \text{ and }\ \tt =e_2\,,  \ && \text{ on } \  \TT^2  \times \{t=\initial\} \,. \label{ics-nn-tt}
\end{alignat} 
\end{subequations}    

The definition of $\Jg$ in \eqref{Jg-def}, together with the dynamics \eqref{d12-h}, shows that
\begin{subequations} 
\label{Jg-system}
\begin{alignat}{2}
(\p_t+V\p_2)\Jg  &= \Jg (  \tfrac{1+ \alpha }{2} \Wbn + \tfrac{1- \alpha }{2} \Zbn ) \,,  && \text{ in } \  \TT^2  \times [\initial, T] \,, \label{Jg-evo} \\
 \Jg &=  1\,,  \ && \text{ on } \  \TT^2  \times \{t=\initial\} \,. \label{ics-Jg}
\end{alignat} 
\end{subequations}

Similarly, with \eqref{g-def}, we compute that
\begin{subequations} 
\begin{alignat}{2}
(\p_t+V\p_2)g^{{\frac{1}{2}} }   & = g^{{\frac{1}{2}} } \bigl(\tfrac{1+ \alpha }{2} \Wbt + \tfrac{1- \alpha }{2} \Zbt\bigr) h,_2\,, && \text{ in } \  \TT^2  \times [\initial, T]  \,, \label{g12-evo} \\
 g^{\frac{1}{2}}  &=  1\,,  \ && \text{ on } \  \TT^2  \times \{t=\initial\} \,. \label{ics-g12}
\end{alignat} 
\end{subequations}

Thanks to  \eqref{Jg-def} and  \eqref{d-n-tau}, we have that 
\begin{equation} 
\Jg,_2 = g^{- {\frac{1}{2}} } h,_{12}  - g^{- {\frac{3}{2}} } h,_1 h,_2 h,_{22}
= g^{- {\frac{1}{2}} } \big( h,_{12}  -  \Jg g^{-\frac 12}   h,_2 h,_{22}\big)
\,.
\label{p1-tau-and-p2-Jg}
\end{equation} 
and therefore
\begin{equation} 
\tt,_1 \cdot \nn 
= g^{-1} h,_{12}
 = g^{- {\frac{1}{2}} } \Jg,_2  +g^{- {\frac{3}{2}}  }  \Jg h,_2 h,_{22}
 =  - \nn,_1 \cdot \tt \,.
\label{p1-n-tau}
\end{equation}

\subsubsection{Dynamics of the Riemann variables}
With \eqref{good-US} and  \eqref{BigWZA},  the un-differentiated Euler system \eqref{euler-shock} is written in ALE variables as 
\begin{subequations} 
\label{Euler0-ALE}
\begin{alignat}{2}
&(\p_t + V \p_2) U^i 
- \alpha \Sigma \tt^i  \Abn
- \tfrac{\alpha}{2} \Sigma \nn^i (\Wbn +  \Zbn) 
+ \tfrac{\alpha}{2} \Sigma (\Wi - \Zi)
 =0\,,  && \text{ in } \  \TT^2  \times [\initial, T]
\,,   \label{U0-ALE} \\
&(\p_t + V \p_2) \Sigma 
+ \alpha \Sigma (\Zbn+ \Abt ) =0 \,,  && \text{ in } \  \TT^2  \times [\initial, T]
  \,, \label{Sigma0-ALE} \\
& U^i =  u_0^i \ \text{ and } \Sigma =\sigma_0\,,  \ && \text{ on } \  \TT^2  \times \{t=\initial\} \,. \label{ics-U-Sigma}
\end{alignat} 
\end{subequations}

Using the chain-rule and \eqref{B-comp}, we also record that
\begin{subequations}
\label{grad-Sigma}
\begin{align} 
\Sigma,_1 &= \Jg (\Sbn + h,_2 \Sbt ) = \tfrac{1}{2} \Jg(\Wbn -\Zbn)  + \tfrac{1}{2} \Jg h,_2 (\Wbt -\Zbt) \,, \label{p1-Sigma} \\
\Sigma,_2 &=  g^{{\frac{1}{2}} } \Sbt   = \tfrac{1}{2} g^{{\frac{1}{2}} } (\Wbt-\Zbt)\,. \label{p2-Sigma}
\end{align} 
\end{subequations}

\begin{figure}[htb!]
\centering
\begin{minipage}{.35\linewidth}
  \includegraphics[width=\linewidth]{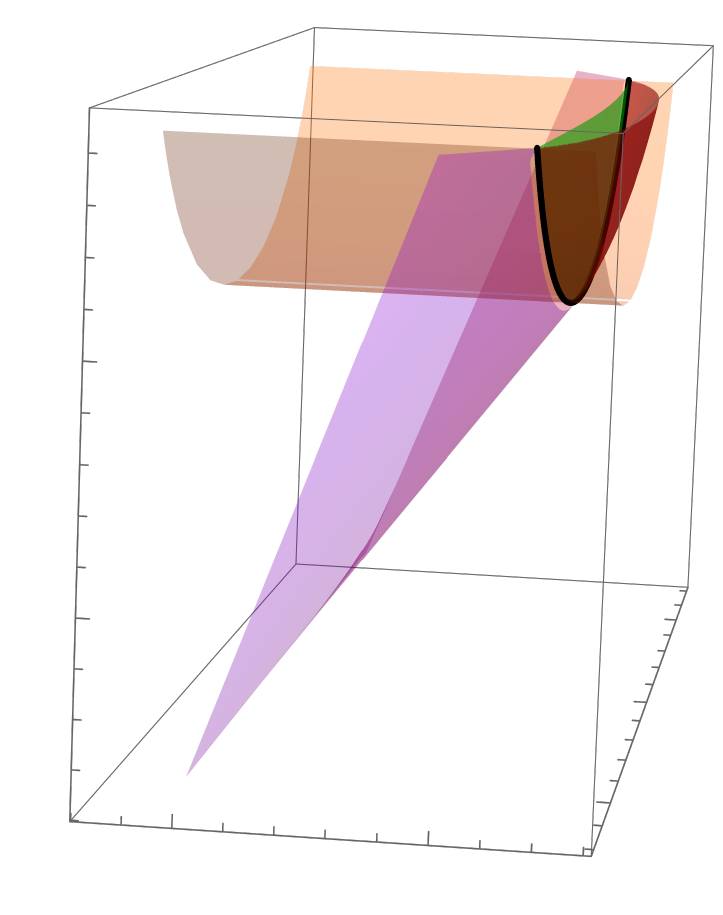}
\end{minipage}
\hspace{.02\linewidth}
\begin{minipage}{.45\linewidth}
  \includegraphics[width=\linewidth]{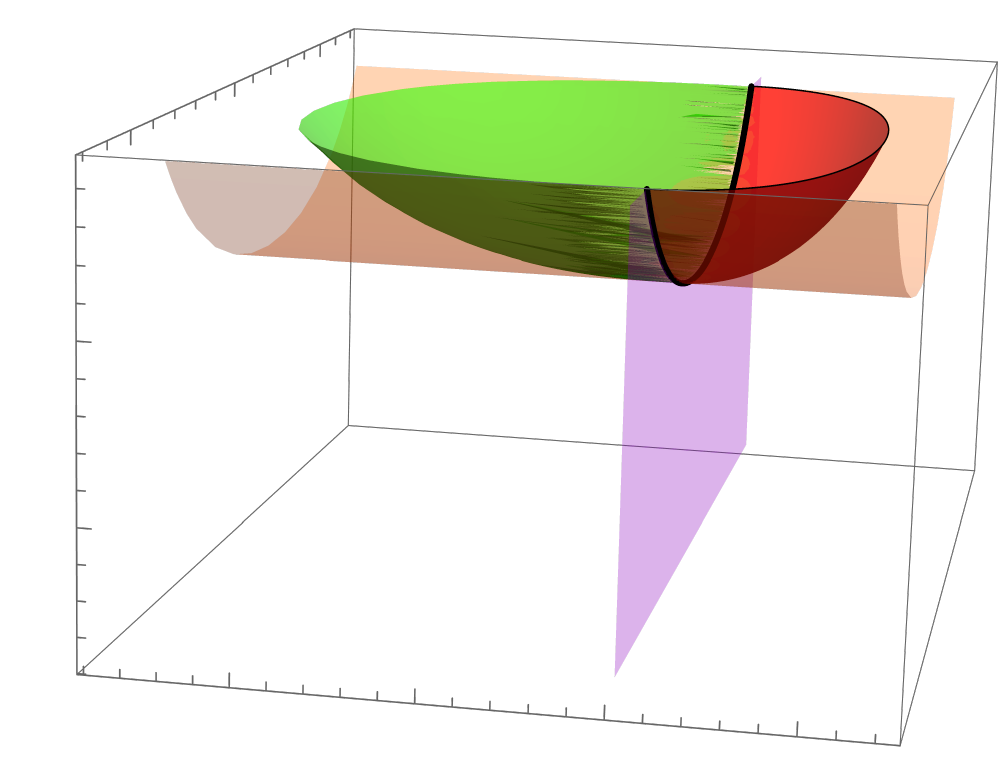}
\end{minipage}
\vspace{-.1 in}
\caption{We revisit Figure~\ref{fig:2D:Eulerian:a}, emphasizing the top boundary of the spacetime of \MGHDB\ of the data, in which the evolution equations~\eqref{nn-tt-evo}, \eqref{Jg-system}, \eqref{Euler0-ALE}, and~\eqref{euler-WZA-real}
will be studied.
On the \underline{left}, we use traditional Eulerian coordinates, while on the 
\underline{right} we use ALE coordinates. In addition to the curve of pre-shocks (in black), the slow acoustic characteristic surface emanating from the pre-shock (in green), and the downstream part of the cuspoidal surface/paraboloid of ``first singularities'' (in red), we display two more surfaces (in both figures). In magenta, we display the fast acoustic characteristic surface which contains the curve of pre-shocks. In orange, we display the cylinder obtained by translating the curve of pre-shocks with respect to $x_1$. The Euler evolution in the spacetime which lies below the orange cylinder is analyzed in Sections~\ref{sec:formation:setup}--\ref{sec:sixth:order:energy}. The Euler evolution on the downstream side of the pre-shock, i.e., in between the orange  and the red surfaces, is analyzed in Section~\ref{sec:downstreammaxdev}.  The Euler evolution on the upstream side of the pre-shock, i.e., in between the orange and the green surfaces, is analyzed in Section~\ref{sec:upstreammaxdev}.}
\label{fig:2D:Lagrangian:a}
\end{figure}

From \eqref{WZA-ALE}, \eqref{Euler0-ALE}, and \eqref{nn-tt-evo}  , we have that
\begin{subequations} 
\label{WZA-evo}
\begin{align}
&(\p_t +V\p_2)W  +  A \bigl(\tfrac{1+ \alpha }{2} \Wbt + \tfrac{1- \alpha }{2} \Zbt\bigr)  +   \alpha \Sigma \Abt
=0\,, 
\label{WZA-evo:a}  \\
&(\p_t +V\p_2)Z+    A \bigl(\tfrac{1+ \alpha }{2} \Wbt + \tfrac{1- \alpha }{2} \Zbt\bigr) -  \alpha \Sigma ( \Abt+
2\Zbn)=0\,, 
\label{WZA-evo:b}  \\
&(\p_t +V\p_2)A   + \tfrac{\alpha}{2} \Sigma (\Wbt - \Zbt - 2 \Abn) - \tfrac 12 (W+Z) \bigl(\tfrac{1+\alpha}{2}\Wbt + \tfrac{1-\alpha}{2} \Zbt\bigr)  =0\,, 
\label{WZA-evo:c}   
\end{align} 
in $\TT^2  \times [\initial, T]$, with initial datum
\begin{equation}
 W =  w_0:= u_0^1+ \alpha \sigma \,,  \qquad 
 Z= z_0:= u_0^1- \alpha \sigma\,, \qquad
 A =a_0:=u_0^2\,,
 \label{ics-WZA}
\end{equation}
\end{subequations} 
on~$\TT^2  \times \{t=\initial\}$.
Similarly, from \eqref{BigWZA}, \eqref{euler-WZA}, and  \eqref{nn-tt-evo}, we deduce
\begin{subequations} 
\label{euler-WZA-real}
\begin{align}
&
\Jg(\p_t +V\p_2)\Wb 
+ \alpha \Jg  \Sigma g^{- {\frac{1}{2}} } \tt \Ub,_2   \notag \\
& \qquad
+  \Jg \Wbn ( \tfrac{1+ \alpha }{2} \Wb + \tfrac{1- \alpha }{2} \Zb )  
+  \Jg \Ab ( \tfrac{3+ \alpha }{2} \Wbt+ \tfrac{1- \alpha }{2} \Zbt ) 
+ \tfrac{\alpha }{2} \Jg (\Wb-\Zb) \Abt =0\,,  
\label{euler-W-real} \\
&
\Jg (\p_t +V\p_2)\Zb
- \alpha \Jg  \Sigma g^{- {\frac{1}{2}} } \tt \Ub,_2   
-2 \alpha \Sigma \bigl(\nn\Ub,_1 - \Sb,_1\bigr) + 2 \alpha \Sigma \Jg g^{- {\frac{1}{2}} }h,_2  \bigl(\nn\Ub,_2 - \Sb,_2 \bigr)
\notag \\
& \qquad
+  \bubu{ \Jg ( \tfrac{1- \alpha }{2} \Wb + \tfrac{1+ \alpha }{2} \Zb ) \Zbn}
+  \Jg \Ab ( \tfrac{1+ \alpha }{2} \Wbt+ \tfrac{3- \alpha }{2} \Zbt ) 
- \tfrac{\alpha }{2} \Jg (\Wb-\Zb)\Abt
=0\,, 
\label{euler-Z-real} \\
&
\Jg(\p_t +V\p_2)\Ab + \alpha \Sigma g^{- {\frac{1}{2}} } \Jg \Sb,_2   - \alpha \Sigma \tt\Ub,_1 
+ \alpha \Sigma g^{- {\frac{1}{2}} } \Jg h,_2 \tt \Ub,_2   \notag \\
& \qquad  
+  \Jg\bigl(  \tfrac{1}{2} (\Wb + \Zb) (\Abn - \tfrac{1+ \alpha }{2} \Wbt - \tfrac{1- \alpha }{2} \Zbt )+ \Ab \Abt + \tfrac{\alpha }{4}(\Wb-\Zb)(\Wbt-\Zbt)    \bigr)  
=0\,,
\label{euler-A-real}
\end{align}
in $\TT^2  \times [\initial, T]$, 
with initial datum
\begin{equation}
 \Wb =  \nabla w_0 := \nabla u_0^1 + \alpha \nabla \sigma_0 \,, 
 \qquad
 \Zb =  \nabla z_0 := \nabla u_0^1 - \alpha \nabla \sigma_0\,, 
 \qquad 
 \Ab = \nabla a_0 := \nabla u_0^2\,,
 \label{ics-bigWZA}
\end{equation} 
\end{subequations} 
on $\TT^2 \times \{t=\initial\}$. The set of equations~\eqref{nn-tt-evo}, \eqref{Jg-system}, \eqref{Euler0-ALE}, and~\eqref{euler-WZA-real} will be studied  in the spacetime region represented in the left panel of Figure~\ref{fig:2D:Lagrangian:a}.

\subsection{Dynamics of $V$}
From \eqref{transport-ale}, \eqref{p2-h}, \eqref{g12-evo}, and \eqref{WZA-evo}, we have that
\begin{equation} 
(\p_t + V \p_2) V = - \alpha \Sigma g^{-\frac 12}  \Bigl(\tfrac{2+ \alpha }{2} \Wbt - \tfrac{\alpha }{2} \Zbt- \Abn - h,_2 \bigl(  \alpha   \Abt -  (1-\alpha)  \Zbn \bigr) \Bigr)    
 \,. \label{V-evo}
\end{equation}

\subsection{Dynamics of normal components $\Wbn$, $\Zbn$, and $\Abn$}
From \eqref{BigWZA}, \eqref{d-n-tau}, \eqref{nn-evo} and \eqref{euler-W-real} we deduce that
\begin{subequations}
\begin{align}
(\p_t + V \p_2) \Wbn
&= - \Wbn \bigl(\tfrac{1+\alpha}{2} \Wbn + \tfrac{1-\alpha}{2} \Zbn + \tfrac{\alpha}{2} \Abt - \tfrac{\alpha}{2} \Sigma g^{-\frac 32} h,_{22}\bigr)  - \alpha \Sigma g^{-\frac 12}  \Abn,_2   + \tfrac{\alpha}{2} \bigl( \Abt + \Sigma g^{-\frac 32} h,_{22} \bigr)  \Zbn 
\notag\\
&\quad
-  \bigl(\tfrac{3+\alpha}{2} \Wbt + \tfrac{1-\alpha}{2} \Zbt\bigr) \Abn
- \bigl(\tfrac{1+\alpha}{2} \Wbt + \tfrac{1-\alpha}{2} \Zbt\bigr) \Wbt  
- \alpha \Sigma g^{-\frac 32} h,_{22} \Abt
\,.
\label{eq:Wb:nn}
\end{align}
Furthermore, using \eqref{Jg-evo} we obtain that
\begin{align}
(\p_t + V \p_2) (\Jg \Wbn)
&= -(\Jg\Wbn) \bigl(\tfrac{\alpha}{2} \Abt - \tfrac{\alpha}{2} \Sigma g^{-\frac 32} h,_{22}\bigr)  
- \alpha \Sigma g^{-\frac 12} \Jg \Abn,_2 + \tfrac{\alpha}{2} \bigl( \Abt + \Sigma g^{-\frac 32} h,_{22} \bigr) \Jg \Zbn 
\notag\\
&\quad
-  \bigl(\tfrac{3+\alpha}{2} \Wbt + \tfrac{1-\alpha}{2} \Zbt\bigr) \Jg\Abn
- \bigl(\tfrac{1+\alpha}{2} \Wbt + \tfrac{1-\alpha}{2} \Zbt\bigr)\Jg \Wbt  
- \alpha \Sigma g^{-\frac 32} h,_{22} \Jg \Abt
\,.
\label{eq:Jg:Wb:nn}
\end{align}
\end{subequations}

From \eqref{BigWZA}, \eqref{d-n-tau}, \eqref{nn-evo} and \eqref{euler-Z-real} we have that
\begin{subequations}
\begin{align}
(\p_t + V \p_2)  \Zbn
&= 
-  \Zbn \bigl(\bubu{\tfrac{1-\alpha}{2} \Wbn + \tfrac{1+\alpha}{2} \Zbn  }
+ \tfrac{\alpha}{2} \Abt + \tfrac{\alpha}{2} \Sigma g^{-\frac 32} h,_{22} \bigr) 
- 2\alpha \Sigma g^{-\frac 12} h,_{2}  {\Zbn,_2}  
+ 2\alpha \Sigma\Jgi    {\Zbn,_1} 
\notag\\
 &\qquad 
 +\alpha \Sigma g^{-\frac 12}    {\Abn,_2}  
+ 2\alpha \Sigma \bubu{g^{-\frac 12} } (\Zbt+ \Abn)\Jgi  \Jg,_2    
+\tfrac{\alpha}{2} \bigl( \Abt - \Sigma g^{-\frac 32}   h,_{22} \bigr) \Wbn 
\notag\\
&\qquad   
-  \bigl(\tfrac{1+\alpha}{2} \Wbt + \tfrac{3-\alpha}{2} \Zbt \bigr) \Abn
- \bigl(\tfrac{1+\alpha}{2} \Wbt + \tfrac{1-\alpha}{2} \Zbt\bigr) \Zbt
+   \alpha  \Sigma g^{-\frac 32}   h,_{22}    \Abt 
\,.
\label{eq:Zb:nn} 
\end{align}
Furthermore, using \eqref{Jg-evo} we obtain that
\begin{align}
(\p_t + V \p_2) (\Jg \Zbn)
&= 
- (\Jg  \Zbn) \bigl(\bubu{- \alpha \Wbn + \alpha \Zbn} + \tfrac{\alpha}{2} \Abt + \tfrac{\alpha}{2} \Sigma g^{-\frac 32} h,_{22} \bigr) 
- 2\alpha \Sigma g^{-\frac 12} h,_{2}   ( \Jg \Zbn),_2
+ 2\alpha \Sigma   \Zbn,_1 
\notag\\
 &\ \
+ \alpha \Sigma g^{-\frac 12}   \Jg  \Abn,_2
+ 2\alpha \Sigma \bubu{g^{-\frac 12}} ( h,_{2} \Zbn +  \Zbt+ \Abn)  \Jg,_2     
+\tfrac{\alpha}{2} \bigl( \Abt - \Sigma g^{-\frac 32}   h,_{22} \bigr) \Jg  \Wbn  
\notag\\
&\ \
-  \bigl(\tfrac{1+\alpha}{2} \Wbt + \tfrac{3-\alpha}{2} \Zbt \bigr) \Jg  \Abn
- \bigl(\tfrac{1+\alpha}{2} \Wbt + \tfrac{1-\alpha}{2} \Zbt\bigr) \Jg  \Zbt
+   \alpha  \Sigma g^{-\frac 32}   h,_{22}    \Jg \Abt 
\,.
\label{eq:Jg:Zb:nn} 
\end{align}
We will later make use of the fact that we can write the $\Zbn$ evolution in terms of a transport velocity based on the $\lambda_1$ wave speed
as follows: 
\begin{align}
&\bigl(\p_t + (V +2\alpha \Sigma g^{-\frac 12} h,_{2}   )  \p_2 -2\alpha \Sigma\Jgi \p_1 \bigr)  \Zbn \notag \\
& \
= 
- \Zbn \bigl(\bubu{ \tfrac{1-\alpha}{2} \Wbn + \tfrac{1+\alpha}{2} \Zbn } + \tfrac{\alpha}{2} \Abt + \tfrac{\alpha}{2} \Sigma g^{-\frac 32} h,_{22} \bigr) 
+ 2\alpha \Sigma (g^{-\frac 12} \Zbt+ \Abn)\Jgi   {\Jg,_2}     
+\alpha \Sigma g^{-\frac 12}   \Abn,_2 
\notag\\
&\ \ 
+\tfrac{\alpha}{2} \bigl( \Abt - \Sigma g^{-\frac 32}   h,_{22} \bigr) \Wbn    
-  \bigl(\tfrac{1+\alpha}{2} \Wbt + \tfrac{3-\alpha}{2} \Zbt \bigr) \Abn
- \bigl(\tfrac{1+\alpha}{2} \Wbt + \tfrac{1-\alpha}{2} \Zbt \bigr) \Zbt
+   \alpha  \Sigma g^{-\frac 32}   h,_{22}    \Abt 
\,.
\label{eq:Zb:nn:alt} 
\end{align}
\end{subequations}

From \eqref{BigWZA}, \eqref{d-n-tau}, \eqref{nn-evo} and \eqref{euler-A-real} we have that
\begin{subequations}
\begin{align}
(\p_t + V \p_2)  \Abn
&=  - \Abn (\tfrac 12 \Wbn + \tfrac 12 \Zbn +  \Abt)
+\alpha \Sigma\Jgi \Abn,_1 
- \alpha \Sigma g^{-\frac 12} h,_2 \Abn,_2 
\notag\\
&\qquad
- \tfrac{\alpha}{2} \Sigma g^{-\frac 12} (\Wbn - \Zbn),_2
+ \tfrac 12 (\Wbn + \Zbn - 2 \Abt) (\tfrac{1+\alpha}{2} \Wbt + \tfrac{1-\alpha}{2} \Zbt )
\notag\\
&\qquad
- \tfrac{\alpha}{4} (\Wbn - \Zbn + 2 \Sigma g^{-\frac 32} h,_{22}) (\Wbt -\Zbt)
- \tfrac{\alpha}{2} \Sigma g^{-\frac 12}\Jgi \Jg,_2 (\Wbn + \Zbn - 2 \Abt)
 \,,
\label{eq:Ab:nn}
\end{align}
and by using \eqref{Jg-evo} we find that
\begin{align}
(\p_t + V \p_2) (\Jg \Abn)
&= - (\Jg \Abn) (- \tfrac{ \alpha}{2} \Wbn + \tfrac{\alpha}{2} \Zbn + \Abt) 
+\alpha  \Sigma \Abn,_1 
- \alpha \Sigma g^{-\frac 12} h,_2  (\Jg \Abn),_2
\notag\\
& 
- \tfrac{\alpha}{2} \Sigma g^{-\frac 12}  (\Jg\Wbn - \Jg\Zbn),_2
+ \tfrac 12 (\Jg \Wbn + \Jg \Zbn - 2 \Jg \Abt) (\tfrac{1+\alpha}{2} \Wbt + \tfrac{1-\alpha}{2} \Zbt)
\notag\\
&
- \tfrac{\alpha}{4} (\Jg \Wbn - \Jg \Zbn + 2 \Sigma g^{-\frac 32} \Jg h,_{22}) (\Wbt -\Zbt)
+  \alpha \Sigma g^{-\frac 12}  \Jg,_2 ( h,_2 \Abn + \Abt - \bubu{\Zbn})
 \,.
\label{eq:Jg:Ab:nn}
\end{align}
 We shall also make use of the $\Abn$ evolution with a transport velocity that encodes the $\lambda_2$ wave speed as
\begin{align}
&\bigl(\p_t + ( V + \alpha \Sigma g^{-\frac 12} h,_2) \p_2 - \alpha \Sigma\Jgi \p_1\bigr) \Abn \notag \\
& \quad
=  - \Abn (\tfrac 12 \Wbn + \tfrac 12 \Zbn +  \Abt)
- \tfrac{\alpha}{2} \Sigma g^{-\frac 12}\Jgi \Jg,_2 (\Wbn + \Zbn - 2 \Abt)
- \tfrac{\alpha}{2} \Sigma g^{-\frac 12} (\Wbn - \Zbn),_2
\notag\\
&\quad
+ \tfrac 12 (\Wbn + \Zbn - 2 \Abt) (\tfrac{1+\alpha}{2} \Wbt + \tfrac{1-\alpha}{2} \Zbt )
- \tfrac{\alpha}{4} (\Wbn - \Zbn + 2 \Sigma g^{-\frac 32} h,_{22}) (\Wbt -\Zbt) \,.
\label{eq:Ab:nn:alt}
\end{align}
\end{subequations}

\subsection{Dynamics of tangential components $\Wbt$, $\Zbt$, and $\Abt$}
Using \eqref{BigWZA}, \eqref{d-n-tau}, \eqref{tt-evo} and \eqref{euler-W-real}, we have that
\begin{align}
(\p_t + V \p_2)  \Wbt
&= - \bigl( \tfrac{3+2\alpha}{2} \Wbt + \tfrac{1-2\alpha}{2} \Zbt \bigr) \Abt
+ \tfrac{\alpha}{2} \Sigma g^{-\frac 32} h,_{22} \bigl(  \Wbt + \Zbt \bigr)
- \alpha \Sigma g^{-\frac 12} {\Ab,_2 \cdot \tt}
\notag\\
&= - \alpha \Sigma g^{-\frac 12} \Abt,_2
 + \tfrac{\alpha}{2} \Sigma g^{-\frac 32}  h,_{22} (\Wbt + \Zbt + 2 \Abn )
 -  \bigl( \tfrac{3+2\alpha}{2} \Wbt + \tfrac{1-2\alpha}{2} \Zbt \bigr) \Abt  \,.
\label{eq:Wb:tt}
\end{align}

In a similar way, from \eqref{BigWZA}, \eqref{d-n-tau}, \eqref{tt-evo} and \eqref{euler-Z-real},  we deduce that
\begin{subequations}
\begin{align}
(\p_t + V \p_2)  \Zbt
&= 
- 2 \alpha \Sigma g^{-\frac 12} h,_2  \Zbt,_2
+ 2\alpha \Sigma\Jgi   \Zbt,_1 
- 2\alpha \Sigma \bubu{g^{-\frac 12}} \bigl(  \Zbn  - \Abt \bigr)\Jg^{-1}   \Jg,_2
+ \alpha \Sigma g^{-\frac 12}  \Abt,_2
\notag\\
&      
- \tfrac{\alpha}{2} \Sigma g^{-\frac 32} h,_{22}(\Wbt + \Zbt + 2 \Abn)
\bubu{+ \alpha \Zbn (\Wbt - \Zbt)}
-   \Abt   \bigl(\tfrac 12 \Wbt +\tfrac 32 \Zbt    \bigr)   \,.
\label{eq:Zb:tt}
\end{align}
Just as we did for $\Zbn$, we can write the $\Zbt$ evolution using the transport velocity based on the $\lambda_1$ wave speed as
\begin{align}
&\bigl(\p_t + (V +2\alpha \Sigma g^{-\frac 12} h,_2) \p_2 - 2 \alpha \Sigma\Jgi \p_1 \bigr)  \Zbt 
\notag\\
&= 2\alpha \Sigma \bubu{g^{-\frac 12}} \Jgi  {\Jg,_2}  \bigl( \Abt  - \Zbn  \bigr)
+ \alpha \Sigma g^{-\frac 12}  \Abt,_2
- \tfrac{\alpha}{2} \Sigma g^{-\frac 32} h,_{22}(\Wbt + \Zbt + 2 \Abn)
\notag\\
&\qquad 
\bubu{+ \alpha \Zbn (\Wbt - \Zbt)}
-   \Abt   \bigl(\tfrac 12 \Wbt +\tfrac 32 \Zbt    \bigr) \,.
\label{eq:Zb:tt:alt}
\end{align}
\end{subequations}

Lastly, from \eqref{BigWZA}, \eqref{d-n-tau}, \eqref{tt-evo} and \eqref{euler-A-real}, we compute that
\begin{subequations}
\begin{align}
(\p_t +V\p_2) \Abt
&= \alpha \Sigma\Jgi \Abt,_1    
- \alpha \Sigma g^{- \frac{1}{2} }  h,_2 \Abt,_2
- \tfrac \alpha 2 \Sigma g^{- \frac{1}{2} }  (\Wbt - \Zbt),_2  
\notag\\
&\  
- \tfrac \alpha 2   \Sigma g^{-\frac 12} (\Wbt + \Zbt + 2 \Abn )\Jgi  \Jg,_2 
+ \tfrac{\alpha}{2} \Sigma g^{- \frac{3}{2} } (\Wbn- \Zbn - h,_2 \Abn) h,_{22}
\notag \\
&\  
 + \tfrac{\alpha}{2} \Abn ( \Wbt - \Zbt ) 
 -  (\Abt)^2 
- \tfrac{\alpha }{4} (\Wbt -\Zbt)^2 
- \tfrac{1}{2}  ( \Wbt+\Zbt) ( \tfrac{1+ \alpha }{2} \Wbt + \tfrac{1- \alpha }{2} \Zbt)
\,,
\label{eq:Ab:tt}
\end{align}
and along the transport velocity using the $\lambda_2$ wave speed, we also have that
\begin{align}
&\bigl(\p_t + (V + \alpha \Sigma g^{-\frac 12} h,_{2}) \p_2 - \alpha \Sigma\Jgi \p_1 \bigr) \Abt
\notag\\
&=
- \tfrac \alpha 2 \Sigma g^{- \frac{1}{2} }  (\Wbt - \Zbt),_2  
- \tfrac \alpha 2   \Sigma g^{-\frac 12} (\Wbt + \Zbt + 2 \Abn )\Jgi  \Jg,_2 
+ \tfrac{\alpha}{2} \Sigma g^{- \frac{3}{2} } (\Wbn- \Zbn - h,_2 \Abn) h,_{22}
\notag \\
&\  
 + \tfrac{\alpha}{2} \Abn ( \Wbt - \Zbt ) 
 -  (\Abt)^2 
- \tfrac{\alpha }{4} (\Wbt -\Zbt)^2 
- \tfrac{1}{2}  ( \Wbt+\Zbt) ( \tfrac{1+ \alpha }{2} \Wbt + \tfrac{1- \alpha }{2} \Zbt)
\,.
\label{eq:Ab:tt:alt}
\end{align}
\end{subequations} 

\subsection{Dynamics of vorticity}
The Eulerian vorticity $\omega := \nabla^\perp \cdot u = \bubu{-} \p_\tau u \cdot n \bubu{+} \p_n u \cdot \tau$ is a solution of 
\begin{equation*} 
\p_t\omega + (u+ \alpha \sigma n) \cdot \nabla \omega + \omega ( \operatorname{div} u - n\cdot \nabla \omega)  =0 \,.
\end{equation*} 
We define the ALE vorticity
\begin{equation}
\Omega  
= \omega \cir\psi 
= \Abn - \tfrac{1}{2} (\Wbt+\Zbt) \,.  
\label{vort-id-good}
\end{equation} 
We next define the Eulerian specific vorticity by  $\upomega=  (\alpha \sigma)^{-\frac{1}{\alpha}}\omega$ and a simple computation verifies that
\begin{equation*}
\p_t \upomega+ u \cdot \nabla \upomega=0 \,.
\end{equation*}
Defining the ALE specific vorticity
\begin{equation} 
 \Upomega=  \upomega\circ\psi = (\alpha \Sigma)^{-\frac{1}{\alpha}}\Omega \,, 
 \label{svort-id-good}
\end{equation} 
an application of the chain-rule, \eqref{B-comp}, and \eqref{BigWZA} shows that
\begin{equation} 
\tfrac{\Jg}{\Sigma}  (\p_t + V \p_2) \Upomega    - \alpha \Upomega,_1 
+ \alpha  g^{- {\frac{1}{2}} } h,_2 \Jg \Upomega,_2  =0  \,. 
\label{vort-t}
\end{equation}

\subsection{Identity for the divergence}
The Eulerian divergence of the fluid velocity is given as ${\rm div}\, u := \p_n u \cdot n + \p_\tau u \cdot \tau$. We show here that the ALE version of the divergence of the fluid velocity, namely $({\rm div}\, u) \cir \psi$, may be written as a linear combination of the differentiated ALE Riemann variables $\Wb, \Zb$, and $\Ab$. More precisely, by appealing to~\eqref{WZA-ALE} and~\eqref{component-identities} we have
\begin{align}
({\rm div}\, u)\cir \psi
&= (\p_n u \cdot n + \p_\tau u \cdot \tau) \cir \psi
\notag\\
&= (\p_n (u \cdot n) - u \cdot \tau \p_n n \cdot \tau + \p_\tau (u \cdot \tau) - u \cdot n \p_\tau \tau \cdot n) \cir \psi
\notag\\
&= (\tfrac 12 \p_n (w+z) - a \p_n n \cdot \tau + \p_\tau a - \tfrac 12 (w+z) \p_\tau \tau \cdot n) \cir \psi
\notag\\
&= \tfrac 12 (\Wbn + \Zbn) + \Abt
\,.
\label{div-id-good}
\end{align}
The above identity is in direct analogy to how the ALE vorticity was written in~\eqref{vort-id-good} as a linear combination of components of $\Wb, \Zb$, and $\Ab$.

\subsection{Identities for $\Jg \Wb, \Jg \Zb, \Jg \Ab$}
It is important to first rewrite the system of equations \eqref{euler-WZA-real} by 
commuting $\Jg$ with  the operator $\p_t+V\p_2$ as follows
\begin{subequations} 
\label{euler-WZA-aug}
\begin{align} 
&
(\p_t +V\p_2)(\Jg\Wb  )
+ \alpha \Sigma g^{- {\frac{1}{2}} }  \Jg  \tt \Ub,_2   =\Fw \,,
 \label{euler-W-aug}  \\
 &
(\p_t +V\p_2)(\Jg\Zb) 
- \alpha \Sigma g^{- {\frac{1}{2}} }  \Jg  \tt\Ub,_2   
 -2 \alpha \Sigma \bigl( \nn\Ub,_1 \ - \Sb,_1\bigr) + 2 \alpha \Sigma \Jg g^{- {\frac{1}{2}} }h,_2  \bigl(\nn\Ub,_2  - \Sb,_2 \bigr)
= \Fz  \,,   \label{euler-Z-aug}  \\
&
(\p_t +V\p_2)(\Jg\Ab)  + \alpha \Sigma  g^{- {\frac{1}{2}} } \Jg \Sb,_2 
- \alpha \Sigma \tt \Ub,_1 
+ \alpha \Sigma  \Jg g^{- {\frac{1}{2}} } h,_2  \tt \Ub,_2  =\Fa  \,, 
 \label{euler-A-aug} \\
 &
 (\p_t+V\p_2) h,_2=  g  \bigl(\tfrac{1+ \alpha }{2} \Wbt + \tfrac{1- \alpha }{2} \Zbt \bigr) \,, \label{euler-h,2-aug} \\
 &
 (\p_t+V\p_2) \Jg=    \bigl(\tfrac{1+ \alpha }{2} \Jg\Wbn + \tfrac{1- \alpha }{2} \Jg\Zbn \bigr) \,, \label{euler-Jg-aug} 
\end{align} 
\end{subequations} 
where
\begin{align*} 
\Fw &= - \tfrac{1- \alpha }{2}\Jg( \Wbn \Zbt - \Wbt \Zbn) \tt 
-  \Jg  ( \tfrac{3+ \alpha }{2} \Wbt+ \tfrac{1- \alpha }{2} \Zbt ) \Ab
- \tfrac{\alpha }{2} \Jg (\Wb-\Zb) \Abt \,, \\
\Fz & =  \bubu{
\alpha \Jg(\Wbn-\Zbn)\Zbn \nn 
- \Jg \big(  \tfrac{1- \alpha }{2}  \Wbt \Zbn  -   \tfrac{1+ \alpha }{2}  \Wbn \Zbt + \alpha \Zbn\Zbt\big) \tt
} \notag \\
&  \qquad
-  \Jg  ( \tfrac{1+ \alpha }{2} \Wbt+ \tfrac{3- \alpha }{2} \Zbt ) \Ab + \tfrac{\alpha }{2} \Jg (\Wb-\Zb)\Abt \,, \\
\Fa & = 
\bubu{
\tfrac{\alpha}{2} \Jg (\Wbn - \Zbn) \Abn \nn 
+ \Jg \big( \Abt (\tfrac{1+\alpha}{2} \Wbn + \tfrac{1-\alpha}{2} \Zbn) - \tfrac 12 \Abn  (\Wbt + \Zbt) \big) \tt
}
\notag \\
& \qquad  
-  \Jg\big(   \Ab \Abt + \tfrac{\alpha }{4}(\Wb-\Zb)(\Wbt-\Zbt)  \big)  
+ \tfrac{1}{2} \Jg ( \Wb+\Zb) ( \tfrac{1+ \alpha }{2} \Wbt+ \tfrac{1- \alpha }{2} \Zbt ) \,.
\end{align*}

In order to close highest-order energy estimates, it is essential to re-weight the equations  \eqref{euler-Z-aug} and \eqref{euler-A-aug}
in a manner specific to the normal and tangential components.  
By  computing the normal components of \eqref{euler-WZA-aug} and the tangential components of \eqref{euler-WZA-real}, we  arrive at the following system of
equations:
\begin{subequations} 
\label{euler-WZA-aug5}
\begin{align} 
&
\tfrac{1}{\Sigma} (\p_t +V\p_2)(\Jg\Wbn  )
+ \alpha    g^{- {\frac{1}{2}} }  (\Jg\Abn),_2 
- \alpha    g^{- {\frac{1}{2}} }  \Abn \Jg,_2
- \tfrac{\alpha}{2} g^{- {\frac{1}{2}} } \tt,_2\cdot\nn(\Jg\Wbn + \Jg\Zbn - 2\Jg\Abt)
   =\Fwn\,,
 \label{euler-Wn}  \\
 &
\tfrac{\Jg}{\Sigma} (\p_t +V\p_2)(\Jg\Zbn)  
- \bubu{  \tfrac{\alpha }{\Sigma} (\Jg \Wbn - \Jg \Zbn)(\Jg\Zbn)}
- \alpha    g^{- {\frac{1}{2}} } \Jg  (\Jg\Abn),_2 
+ \alpha    g^{- {\frac{1}{2}} }  \Abn \Jg \Jg,_2
\notag \\
&\quad
+ \tfrac{\alpha}{2} g^{- {\frac{1}{2}} } \tt,_2\cdot\nn \Jg (\Jg\Wbn + \Jg\Zbn - 2\Jg\Abt)
 -2 \alpha   (\Jg\Zbn),_1  
 -2 \alpha \tt,_1 \cdot \nn \Jg (\Abn+\Zbt) 
 +2 \alpha \Jg,_1 \Zbn
 \notag \\
 &\quad
 +2 \alpha \Jg g^{- {\frac{1}{2}} } h,_2  (\Jg\Zbn),_2
 +2 \alpha \tt,_2 \cdot \nn \Jgt g^{- {\frac{1}{2}} } h,_2  (\Abn+\Zbt) 
 -2 \alpha \Jg,_2 \Jg g^{- {\frac{1}{2}} } h,_2  \Zbn
 = \Fzn \,,   \label{euler-Zn}  \\
&
\tfrac{\Jg}{\Sigma}  (\p_t +V\p_2)(\Jg\Abn)  
\bubu{- \tfrac{\alpha}{2\Sigma} ( \Jg\Wbn - \Jg\Zbn) (\Jg\Abn) } 
+ \alpha   g^{- {\frac{1}{2}} }  \Jg  (\Jg\Sbn),_2 
- \alpha   g^{- {\frac{1}{2}} }  (\Jg\Sbn) \Jg,_2
 \notag \\
&\quad
+ \alpha g^{- {\frac{1}{2}} } \Jgt \tt,_2 \cdot \nn \Sbt
- \alpha (\Jg\Abn),_1 
+ \tfrac{\alpha}{2} \tt,_1 \cdot\nn (\Jg\Wbn +\Jg\Zbn - 2\Jg \Abt)
+ \alpha \Jg,_1 \Abn
 \notag \\
&\quad
+ \alpha \Jg g^{- {\frac{1}{2}} } h,_2 (\Jg\Abn),_2
- \tfrac{\alpha}{2} \tt,_2 \cdot\nn \Jg g^{- {\frac{1}{2}} } h,_2 (\Jg\Wbn +\Jg\Zbn - 2\Jg \Abt)
- \alpha \Jg g^{- {\frac{1}{2}} } h,_2 \Jg,_2 \Abn
 =\Fan \,,
 \label{euler-An} 
\end{align} 
and
\begin{align} 
&
\tfrac{1}{\Sigma} (\p_t +V\p_2)\Wbt  
+ \alpha g^{- {\frac{1}{2}} } \Abt,_2 -  \alpha  g^{- {\frac{1}{2}} } \tt,_2\cdot\nn(\Omega +\Wbt+\Zbt )
 =\Fwt\,,
 \label{euler-Wt}  \\
&
\tfrac{\Jg}{\Sigma} (\p_t +V\p_2)\Zbt 
- \alpha \Jg g^{- {\frac{1}{2}} } \Abt,_2 +   \alpha\Jg  g^{- {\frac{1}{2}} } \tt,_2\cdot\nn(\Omega +\Wbt+\Zbt )
- 2 \alpha \Zbt,_1 - 2 \alpha \tt,_1\cdot\nn(\Abt-\Zbn)
\notag \\
& \quad
 + 2 \alpha  g^{- {\frac{1}{2}} }h,_2 \Jg  \Zbt,_2
 + 2 \alpha  g^{- {\frac{1}{2}} }h,_2 \Jg \tt,_2\cdot\nn (\Abt-\Zbn)
=\Fzt \,, 
\label{euler-Zt}  \\
 &
\tfrac{\Jg}{\Sigma} (\p_t +V\p_2)\Abt 
+ \alpha   g^{- {\frac{1}{2}} } \Jg \Sbt,_2  
-  \alpha   g^{- {\frac{1}{2}} } \Jg \Sbn \tt,_2\cdot \nn
- \alpha \Abt,_1 + \alpha \tt,_1\cdot\nn (\Omega +\Wbt+\Zbt )
\notag \\
& \quad 
+ \alpha   g^{- {\frac{1}{2}} }  h,_2 \Jg\Abt,_2
- \alpha   g^{- {\frac{1}{2}} }  h,_2 \Jg\tt,_2\cdot\nn (\Omega +\Wbt+\Zbt )
 =\Fat  \,,
 \label{euler-At} 
\end{align} 
\end{subequations} 
where
\begin{subequations} 
\label{forcing-nt}
\begin{align} 
\Fwn & = 
-\tfrac{\Jg}{\Sigma} \big(  ( \tfrac{3+ \alpha }{2} \Wbt+ \tfrac{1- \alpha }{2} \Zbt ) \Abn
+ \tfrac{\alpha }{2} (\Wbn-\Zbn) \Abt  + ( \tfrac{1+\alpha}{2} \Wbt + \tfrac{1-\alpha}{2} \Zbt)\Wbt\big)
\,, \\
\Fzn & =   -\tfrac{\Jgt}{\Sigma} \big(
  ( \tfrac{1+ \alpha }{2} \Wbt+ \tfrac{3- \alpha }{2} \Zbt )\Abn
- \tfrac{\alpha }{2}\Abt  (\Wbn-\Zbn)
+ ( \tfrac{1+\alpha}{2} \Wbt + \tfrac{1-\alpha}{2} \Zbt)\Zbt \big)  \,, \\
\Fan & = -\tfrac{\Jgt}{\Sigma} \big(
 \Abn \bubu{\Abt} \bubu{-} \tfrac{1}{4}  \Wbn(\Wbt+\Zbt)
 -\Zbn ( \tfrac{1+2\alpha}{4} \Wbt - \tfrac{1-2\alpha}{4} \Zbt)
 \bubu{+ \Abt (\tfrac{1+\alpha}{2} \Wbt + \tfrac{1-\alpha}{2} \Zbt)}
\big) \,, \\
\Fwt & = -\tfrac{1}{\Sigma}    \Abt ( \tfrac{3+ 2\alpha }{2} \Wbt+ \tfrac{1- 2\alpha }{2} \Zbt ) 
\,, \\
\Fzt & = - \tfrac{\Jg}{\Sigma}   \Abt ( \tfrac{1 }{2} \Wbt+ \tfrac{3 }{2} \Zbt ) + \bubu{\tfrac{\alpha}{\Sigma} \Jg\Zbn (\Wbt - \Zbt)}
 \,, \\
\Fat & =   \tfrac{\Jg}{\Sigma} \big(
 \tfrac{\alpha}{2}  \Omega (\Wbt - \Zbt) 
 + \tfrac 12 (\Wbt + \Zbt) (\tfrac{1+2\alpha}{2} \Wbt + \tfrac {1-2\alpha}{2}\Zbt ) -  \Abt^2   - \tfrac{\alpha }{4}(\Wbt-\Zbt)^2 
\big) 
 \,.
\end{align} 
\end{subequations}


\section{Initial data and main results}
\label{sec:thesetup}

We consider the Euler system posed on the spatial torus $\mathbb{T}^2 = [-\pi,\pi]^2$ with periodic boundary conditions in space.

\subsection{The time interval}
\label{sec:tin:tmed:tfin}
The initial data is prescribed at the initial time 
\begin{equation}
\initial := - \tfrac{2}{1+\alpha} \eps 
\,.
\label{tin}
\end{equation}
The above choice of initial time is made so that the very first time at which is $\Jg$ vanishes, is given by 
\begin{equation}
0\pm \OO(\eps^{2})
\,.
\label{first-time}
\end{equation} 
That is, the very first singularity is expected to arise at a time which is $\OO(\eps)$ past $\initial$. Our goal is to analyze the Euler dynamics for an additional $\OO(\eps)$ amount of time past this very first blow-up time, so that there is ample room (in time) to observe the geometry of the \MGHDB. For this purpose, we define
\begin{equation}
\medium = \tfrac{2}{1+\alpha} \eps \cdot  \tfrac{1}{100}
\,.
\label{tmed}
\end{equation}
We aim to study the Euler evolution on $\TT^2 \times [\initial,\medium]$, i.e., for an $\OO(\eps)$ amount of time past the very first singularity. 

In practice, it is convenient to smoothly cutoff the \MGHDB\ spacetime with respect to time, while leaving the dynamics unaltered on $[\initial,\medium]$. For this purpose we introduce a final time
\begin{equation}
\final:= \tfrac{2}{1+\alpha} \eps \cdot  \tfrac{1}{50}
\,.
\label{time-of-existence}
\end{equation}
We leave the Euler dynamics unaltered on $[\initial,\medium]$, and employ a smooth cutoff procedure in time (see the definition of the function~$\Jgb$ below in~\eqref{eq:Jgb:identity:0}) on the time interval $(\medium,\final]$. We establish bounds which are valid on $\TT^2 \times [\initial,\final]$, but the Euler geometry is only captured on the subset $\TT^2\times[\initial,\medium]$. Throughout the paper, we restrict the time variable to satisfy  $t\leq \final$, so that in particular $0 \leq t-\initial \leq \frac{2\eps}{1+\alpha} \cdot \frac{51}{50}$ for all $t\in [\initial,\final]$. We note that the precise definitions of $0 < \medium < \final$ given above are not relevant; these choices are only made here for convenience, and the only features that matter are that $0 < \medium = \OO(\eps)$ and that $0 < \final - \medium = \OO(\eps)$.

\begin{remark}[\bf Notation: rescaled derivatives]
\label{rem:notation:nb:doublecomm}
Since the evolution occurs on an $\OO(\eps)$ time interval, we expect that each time derivatives to ``cost'' a factor of $\eps^{-1}$. It turns out that the analysis in the paper is performed on  $\OO(\eps) \times \OO(1)$ sub-domain of $\mathbb{T}^2$, and thus we expect $\p_{1}$ derivatives to ``cost'' a factor of $\eps^{-1}$, while $\p_2$ derivatives to not bring in powers of $\eps$. Keeping this in mind, throughout the paper, we will use the {\em rescaled spacetime gradient operator} $\nb$, given  in $(x,t)$ coordinates, by
\begin{equation}
\nb = ( \nb_t, \nb_1, \nb_2) := (\eps \p_t, \eps \p_1, \p_2)
\,.
\label{eq:nb:x:t}
\end{equation}
When $\nb$ acts on a function of space alone, we implicitly identify $\nb$ with $(\nb_1,\nb_2) = (\eps \p_1, \p_2)$.
We will use the following notation when discussing these derivatives:
\begin{itemize}[leftmargin=16pt] 

\item The symbol $\nb^m$ is used to denote any partial derivative $\nb^\gamma$ with $\gamma \in \mathbb{N}_0^3$, where $\nb^\gamma = (\eps \p_t)^{\gamma_0} (\eps \p_1)^{\gamma_1}\p_2^{\gamma_2}$ with $|\gamma|=m$. In particular, throughout this section there is no need to keep track of the specific multi-index $\gamma$, just of the total order $|\gamma|=m$ of the tangential derivative. 

\item Naturally, the symbol $\| \nb^m f\|_{L^2}^2$ denotes $\sum_{|\gamma|=m} \|(\eps \p_t)^{\gamma_0} (\eps \p_1)^{\gamma_1} \p_2^{\gamma_2} f\|_{L^2}^2$. Whenever the aforementioned sum over all pairs $\gamma \in \mathbb{N}_0^3$ with $|\gamma|=m$ is implicit, it will be dropped, so that we do not further clutter the notation.

\item For any scalar function $f$, with the notation for $\nb$ introduced in~\eqref{eq:nb:x:t}, we shall denote commutators as
\begin{equation*}
 \jump{f,\nb^\gamma} g 
 := f \nb^\gamma g -\nb^\gamma (f g) = - \sum\nolimits_{\delta \leq \gamma, |\delta|\leq |\gamma|-1} {\gamma \choose \delta} \nb^{\gamma-\delta} f \; \nb^\delta g \,.
\end{equation*}
Lastly, we shall use the notation 
\begin{equation*}
\doublecom{\nb^\gamma, f , g}
= \nb^{\gamma}(f \, g) - f \nb^\gamma g - g \nb^\gamma f
= \sum\nolimits_{\delta \leq \gamma, 1\leq |\delta|\leq |\gamma|-1} {\gamma \choose \delta} \nb^{\gamma-\delta} f \; \nb^\delta g \,,
\end{equation*}
to denote ``double-commutators''.
\end{itemize}
\end{remark}

\subsection{The initial data}
\label{cauchydata}
It is convenient to state the initial data assumptions in terms of $(w_0,z_0,a_0)$ defined cf.~\eqref{tn-lag},~\eqref{eq:h:initial:condition}, and~\eqref{WZA-ALE} via
\begin{equation}
w_0 = u_0 \cdot e_1 + \sigma_0, \qquad 
z_0 = u_0 \cdot e_1 - \sigma_0, \qquad 
a_0 = u_0 \cdot e_2,
\label{eq:w0:z0:a0:def}
\end{equation}
rather the velocity and rescaled sound speed $(u_0,\sigma_0)$. 
The initial data is taken to satisfy the following properties:

\begin{enumerate}[leftmargin=26pt]

\item \label{item:ic:supp} There exists a constant\footnote{The purpose of $\kappa_0$ is to ensure that the initial data does not have vacuum, see~\eqref{sigma0-bound}. The assumption $\kappa_0 \geq 20$ is made only for convenience.} $\kappa_0 \geq 20$, {\bf which is independent of $\eps$}, such that  
\begin{equation}
\supp(w_0 -\kappa_0) \cup \supp(z_0) \cup \supp(a_0)
\subseteq
\mathcal{X}_{\mathsf{in}}:=
[- 13 \pi \eps, 13 \pi \eps]\times \mathbb{T}
\,.
\label{eq:ic:supp}
\end{equation}
Naturally, $\eps$ is assumed small enough to ensure that $13 \eps \leq 1$, so that $\mathcal{X}_{\mathsf{in}} \subset \TT^2$.

\item \label{item:ic:infinity} 
At the level of no derivatives, we assume that 
\begin{equation*}
\norm{w_0 -\kappa_0}_{L^\infty_x} \leq \tfrac{5}{3}\,,
\qquad 
\norm{z_0}_{L^\infty_x} \leq \eps \kappa_0\,,
\qquad 
\norm{a_0}_{L^\infty_x} \leq \eps \kappa_0 \,.
\end{equation*}
Therefore, the initial rescaled sound speed $\sigma_0$ satisfies
\begin{equation}
\tfrac{1}{3} \kappa_0 \le \sigma _0(x) \le \tfrac{2}{3} \kappa_0  \,,
\qquad \mbox{for all} \qquad
x \in \mathbb{T}^2\,.
\label{sigma0-bound}
\end{equation}  

\item \label{item:ic:reg} Assume that $(w_0,z_0,a_0) \in H^7(\mathbb{T}^2)$, and that there exists a constant $\bar{\mathsf{C}} \geq 1$, {\bf which is independent of $\eps$}, with\footnote{Intuitively, estimate~\eqref{eq:ic:norm} says that (at least up to the seventh derivative) we should think of $w_0(x) \sim \kappa_0 + \mathcal{W}(\frac{x_1}{\eps},x_2)$, $z_0(x) \sim \eps \mathcal{Z}(\frac{x_1}{\eps},x_2)$, and $a_0(x) \sim \eps \mathcal{A}(\frac{x_1}{\eps},x_2)$, for some smooth functions $(\mathcal{W}, \mathcal{Z},\mathcal{A})$ which obey $\OO(1)$ bounds (with respect to $\eps$) for all their derivatives.}
\begin{equation}
\sum_{1 \leq |\gamma| \leq 5}
\bigl(
\snorm{\nb^\gamma w_0}_{L^\infty_x}
+
\eps^{-1} \snorm{\nb^\gamma (z_0,a_0)}_{L^\infty_x}
\bigr)
+
\sum_{6 \leq |\gamma| \leq 7} 
\bigl(
\eps^{-\frac 12} \snorm{\nb^\gamma w_0}_{L^2_x}
+
\eps^{-\frac 32} \snorm{\nb^\gamma (z_0,a_0)}_{L^2_x} 
\bigr)
\leq 
\bar{\mathsf{C}}
\,.
\label{eq:ic:norm}
\end{equation}
Here we have used the notation $\nb_1 = \eps \p_1$, $\nb_2 = \p_2$, and $\nb^\gamma = \nb_1^{\gamma_1} \nb_2^{\gamma_2}$ for $\gamma = (\gamma_1,\gamma_2) \in \mathbb{N}_0^2$.

\item  \label{item:ic:max:w0} 
Recalling the notation $\nb_1 = \eps \p_1$ and $\nb_2 = \p_2$, we assume that  for all $x \in \TT^2$ we have
\begin{equation}
-1 \leq \nb_1 w_0 \leq \tfrac{1}{10}
\,, 
\qquad
|\nb_2 w_0| \leq 1 
\,,
\qquad \mbox{and}\qquad 
|\nb \nb_1 w_0| 
\leq  2
\,.
\label{eq:why:the:fuck:not:0}
\end{equation}

\item  \label{item:ic:W:Taylor:1} 
We assume that the global minimum of $\nb_1 w_0$ is non-degenerate and occurs at $x=0$. By this we mean that $\nb_1 w_0(0) = -1$,  $\nb_2 w_0(0)= 0$,  $\nb \nb_1 w_0(0) = 0$, and $(1-\eps) {\rm Id} \leq \nb^2 \nb_1 w_0(0) \leq (1+\eps) {\rm Id}$. 

\item  \label{item:ic:w0:x2:negative} 
We assume that for each $x_2 \in \mathbb{T}$, the function $x_1 \mapsto \nb_1 w_0(x_1,x_2)$ attains its global minimum at a unique point $x_{1}^{\vee} = x_{1}^{\vee}(x_2)$, such that the   bounds  $\nb_1 w_0(x_{1}^{\vee}(x_2),x_2) \leq - \frac {9}{10}$, and $\nb_1^3 w_0(x_{1}^{\vee}(x_2),x_2) \geq  \frac{9}{10}$ hold.

\item \label{item:ic:w0:x2:special}
We assume\footnote{Assumption~\eqref{item:ic:w0:x2:special} is only used once: in the proof of Lemma~\ref{lem:Q:bnds}; more precisely, in the proof of estimate~\eqref{eq:x1star:x1vee}.} that there exists $\eps_0>0$ such that for all $\eps \in (0,\eps_0]$ and all  $x_1$ such that $|x_1 - x_1^\vee(x_2)|\geq \eps^{\frac 54}$, we have $\nb_1w_0(x_1,x_2) \geq \nb_1w_0(x_1^\vee(x_2),x_2) + \eps^{\frac 34}$.

\item \label{item:ic:w0:d11:positive}
We assume\footnote{Assumption~\eqref{item:ic:w0:d11:positive} is only used in Sections~\ref{sec:downstreammaxdev} and~\ref{sec:upstreammaxdev}, where we consider the downstream \MGHDB, respectively the upstream \MGHDB.} that there exists $\eps_0>0$ such that for all $\eps \in (0,\eps_0]$ the following holds. If $x=(x_1,x_2)$ is such that $x_1 - x_1^\vee(x_2) \geq \eps^{\frac 74}$ and $\nb_1^2 w_0(x) \leq \eps^{\frac 78}$, then $\nb_1 w_0(x) \geq - \frac 13$ and $\nb_1^2 w_0(x) \geq - 1$.  
Symmetrically, if $(x_1,x_2)$ is such that $x_1 - x_1^\vee(x_2) \leq - \eps^{\frac 74}$ and $\nb_1^2 w_0(x_1,x_2) \geq - \eps^{\frac 78}$, then  we  assume that $\nb_1 w_0(x_1,x_2) \geq - \frac 13$ and  $\nb_1^2 w_0(x_1,x_2) \leq 1$. Lastly, we also assume that for $x\in\TT^2$, if $\nb_1w_0 (x) \leq - \frac 13$, then $\nb_1^3 w_0(x) \geq \frac 13$.

\end{enumerate}

\begin{remark}[\bf Size of derivatives of the initial data]
\label{rem:table:derivatives}
Instead of working with~\eqref{eq:ic:norm}, it is convenient to quantify the size of higher order derivatives of the fundamental variables in the analysis. Recall that $(W,Z,A,U,\Sigma)|_{t=\initial} = (w_0,z_0,a_0,u_0,\sigma_0)$, $\Jg|_{t=\initial} = 1$, $h,_1|_{t=\initial} = 1$,  $h,_2|_{t=\initial} = 0$, $g|_{t=\initial} = 1$, $\nn|_{t=\initial} = e_1$, $\tt|_{t=\initial} = e_2$,   that $(\Wb,\Zb,\Ab)|_{t=\initial}$ may be computed from the identities \eqref{component-identities}, while $V|_{t=\initial}$ may be computed from~\eqref{transport-ale}.  Then, using the identities in Section~\ref{sec:new:Euler:variables} (at  $t= \initial$), we may show that \eqref{eq:ic:norm} and assumptions~\eqref{item:ic:supp}, \eqref{item:ic:infinity} imply that there exists a constant $\Cdata = \Cdata(\alpha,\kappa_0,\bar{\mathsf{C}}) \geq 1$, {\bf which is independent of $\eps$}, such that
\begin{align}
& \sum_{1\leq |\gamma|\leq 7} \eps^{-\frac 12} \snorm{\nb^\gamma W(\cdot,\initial)}_{L^2_x} 
+  \eps^{-\frac 32} \snorm{\nb^\gamma (Z,A)(\cdot,\initial)}_{L^2_x}
+ \sum_{|\gamma|\leq 5} \snorm{\nb^\gamma W(\cdot,\initial)}_{L^\infty_x} 
+ \eps^{-1} \snorm{\nb^\gamma (Z,A)(\cdot,\initial)}_{L^\infty_x}
\notag\\
&\quad 
+ \sum_{|\gamma|\leq 6} \eps^{\frac 12} \snorm{\nb^\gamma (\Jg \Wbn)(\cdot,\initial)}_{L^2_x} 
+  \eps^{-\frac 12} \snorm{\nb^\gamma (\Jg \Zbn,\Jg \Abn)(\cdot,\initial)}_{L^2_x}
+  \eps^{-\frac 12} \snorm{\nb^\gamma (\Zbn,\Abn)(\cdot,\initial)}_{L^2_x}
\notag\\
&\quad 
+ \sum_{|\gamma|\leq 4} \eps \snorm{\nb^\gamma (\Jg \Wbn)(\cdot,\initial)}_{L^\infty_x} 
+ \snorm{\nb^\gamma (\Jg \Zbn,\Jg \Abn)(\cdot,\initial)}_{L^\infty_x}
+ \snorm{\nb^\gamma (\Zbn,\Abn)(\cdot,\initial)}_{L^\infty_x}
\notag\\
&\quad 
+ \sum_{|\gamma|\leq 6} \eps^{-\frac 12} \snorm{\nb^\gamma  \Wbt(\cdot,\initial)}_{L^2_x} 
+  \eps^{-\frac 32} \snorm{\nb^\gamma (\Zbt,\Abt)(\cdot,\initial)}_{L^2_x}
\notag\\
&\quad 
+ \sum_{|\gamma|\leq 4}   \snorm{\nb^\gamma  \Wbt(\cdot,\initial)}_{L^\infty_x} 
+ \eps^{-1} \snorm{\nb^\gamma ( \Zbt, \Abt)(\cdot,\initial)}_{L^\infty_x}
\notag\\
&\quad 
+ \sum_{|\gamma|\leq 6} \eps^{-\frac 12} \snorm{\nb^\gamma (\Jg,\Sigma,\Upomega,g) (\cdot,\initial)}_{L^2_x} 
+ \sum_{|\gamma|\leq 4}   \snorm{\nb^\gamma  (\Jg,\Sigma,\Upomega,g) (\cdot,\initial)}_{L^\infty_x} 
\notag\\
&\quad 
+ \sum_{|\gamma|\leq 6} \eps^{-\frac 32} \snorm{\nb^\gamma  (h,h,_2,V)(\cdot,\initial)}_{L^2_x} 
+ \sum_{|\gamma|\leq 4}  \eps^{-1} \snorm{\nb^\gamma (h,h,_2,V) (\cdot,\initial)}_{L^\infty_x} 
\leq \Cdata
\,.
\label{table:derivatives}
\end{align}
Here we use the notation from~\eqref{eq:nb:x:t}, with $\nb^\gamma= \nb_t^{\gamma_0}\nb_1^{\gamma_1}\nb_2^{\gamma_2}$, and $\gamma \in \mathbb{N}_0^3$. We note that the gain of $\eps^{\frac 12}$ that the $L^2_x$ bounds experience over the $L^\infty_x$ bounds are due to the support of size $\OO(\eps)$ in the $x_1$ direction, see assumption~\eqref{item:ic:supp}.
Verifying that the bounds in \eqref{table:derivatives} have a scaling with respect to $\eps$ that is consistent with \eqref{eq:ic:norm}, and also with assumptions~\eqref{item:ic:supp}, \eqref{item:ic:infinity}, is an exercise whose details are omitted here. 
Throughout the paper we shall refer to~\eqref{table:derivatives} instead of~\eqref{eq:ic:norm}, although the former follows from the latter. 
\end{remark}

\begin{example}[\bf The prototypical initial data]
\label{ex:prototypical}
The prototypical example  for $w_0$ is of the type 
\begin{equation}
w_{0,{\rm ex}}(x_1,x_2) = \kappa_0 + \varphi(\tfrac{x_1}{\eps}) \phi(x_2),
\label{eq:w0:ex}
\end{equation}
where $\kappa_0 \geq 20$, and the functions $\varphi \in C_0^\infty(\mathbb{R})$  and $\phi \in C^\infty(\TT)$  have the following properties:
\begin{itemize}[leftmargin=16pt]
\item $\varphi(r) = - r + \frac 16 r^3$ for all $|r| \leq \sqrt{2}$ and $\varphi(r) =0$ for $|r| \geq 13 \pi$. For $\sqrt{2} < |r| < 13 \pi$ we take $\varphi$ such that $-1 \leq {\rm sgn}(r) \varphi(r) \leq \frac 12$, $-\frac 14 \leq \varphi'(r) \leq \frac{1}{11}$, and $- \frac{1}{2} \leq \sgn(r)\varphi''(r) \leq \frac 32$. For $\eps \leq \frac{1}{13}$, we may view $\varphi(\frac{\cdot}{\eps})$ as  $[-\pi,\pi]$-periodic. It is straightforward to construct a function $\varphi$
which satisfies all the above conditions.

\item $\phi(\bar r) = 1 - \frac 12 \bar r^2$ for all $|\bar r| \leq \frac{1}{\sqrt{20}}$ and  $\phi(\bar r) = \frac{19}{20}$ for $\frac{\pi}{2} \leq |\bar r| \leq \pi$.
For  $\frac{1}{\sqrt{20}} < |\bar r| < \frac{\pi}{2}$ we take $\phi$ such that  $- \frac 12 \leq {\rm sgn}(\bar r) \phi'(\bar r) \leq 0$, and $|\phi''(\bar r)| \leq \frac 32$. It is straightforward to  construct a function $\phi$
which satisfies all the above conditions.
\end{itemize}
A plot of a prototypical function $w_{0,{\rm ex}}$ and its derivative $\nb_1w_{0,{\rm ex}}$ are given in Figures~\ref{fig:w0:ex} and~\ref{fig:w0,1:ex} below. 
We then verify that the function defined in \eqref{eq:w0:ex} satisfies the assumptions we imposed on $w_0$:
\begin{itemize}[leftmargin=16pt]

\item \eqref{item:ic:supp} holds because $\varphi(\frac{x_1}{\eps})=0$ when $|x_1|\geq 13\pi \eps$.

\item \eqref{item:ic:infinity} holds since $\kappa_0 \geq 20$ and because $|\varphi| \leq 1$ and $0 \leq \phi \leq 1$.

\item \eqref{item:ic:reg} holds for some $\bar{\mathsf{C}}>0$ because $\varphi$ and $\phi$ are $C^\infty$ smooth and have the correct scaling in $x_1$ and $x_2$.

\item \eqref{item:ic:max:w0} holds because $-1 \leq \varphi' \leq \frac{1}{11}$, $|\phi'| \leq \frac 12$, and $|\varphi''|\leq \frac 32$.

\item \eqref{item:ic:W:Taylor:1} holds because $w_{0,{\rm ex}}(x_1,x_2) = \kappa_0 - \frac{x_1}{\eps}  + \frac {1}{2} \frac{x_1}{\eps} x_2^2 + \frac{1}{6} \frac{x_1^3}{\eps^3}  - \frac{1}{12} \frac{x_1^3}{\eps^3} x_2^2$, for all $|x_1|\leq \sqrt{2} \eps$ and $|x_2|\leq \frac{1}{\sqrt{20}}$. 

\item \eqref{item:ic:w0:x2:negative} holds with $x_{1}^{\vee}(x_2) = 0$ because $\phi\geq 0$ and so the global minimum of the function $\p_1 w_{0,{\rm ex}}(\cdot,x_2) = \frac{1}{\eps} \varphi'(\frac{\cdot}{\eps}) \phi(y_2)$ occurs where $\varphi'$ is most negative, i.e., at $x_1=0$. Moreover, we have $\nb_1 w_{0,{\rm ex}} (x_{1}^{\vee}(x_2),x_2) = \varphi'(0) \phi(y_2) = - \phi(y_2) \leq -  \frac{19}{20} < -\frac{9}{10}$,  and $\nb_1^3 w_{0,{\rm ex}} (x_{1}^{\vee}(x_2),x_2) = \varphi'''(0) \phi(y_2) =  \phi(y_2) \geq    \frac{19}{20} >   \frac{9}{10}$. 

\item \eqref{item:ic:w0:x2:special} holds with $\eps_0 = 6^{-5}$, because for $|x_1 - x_1^\vee(x_2)| = |x_1| \geq \eps^{\frac 54}$, we have   $\nb_1w_0(x_1,x_2) - \nb_1w_0(x_1^\vee(x_2),x_2) = (\varphi'(x_1) - \varphi'(0)) \phi(x_2) \geq \frac{19}{20}  (1+\varphi'(x_1)) =  \frac{19}{20} \min\{ \frac 34, \frac 12 (\frac{x_1}{\eps})^2\} \geq \frac{19}{40} \eps^{\frac 12} > \eps^{\frac 34}$.

\item \eqref{item:ic:w0:d11:positive} holds because 
if $x_1 - x_1^\vee(x_2) = x_1 \geq \eps^{\frac 74}$, it means that $ \frac{x_1}{\eps} \geq \eps^{\frac 34}$, and since $\phi \geq \frac{19}{20}$ we have $\nb_1^2 w_{0,{\rm ex}}(x) \leq \eps^{\frac 78}  \Leftrightarrow \varphi''(\frac{x_1}{\eps} )\leq \frac{20}{19} \eps^{\frac 78} \Rightarrow  \frac{x_1}{\eps} > \sqrt{2}$, since if $\frac{x_1}{\eps} \leq \sqrt{2}$, then $\varphi''(\frac{x_1}{\eps}) = \frac{x_1}{\eps} \geq \eps^{\frac 34}$ which is strictly larger than $\frac{20}{19} \eps^{\frac 78}$ if $\eps\leq \eps_0 \leq (\frac{19}{20})^{8}$. But if $\frac{x_1}{\eps} \geq \sqrt{2}$, then $\varphi'(\frac{x_1}{\eps}) \geq  -\frac14$ and so $\nb_1 w_{0,{\rm ex}}(x) = \varphi'(\frac{x_1}{\eps}) \phi(x_2)  \geq -\frac14 > -\frac 13$. Moreover, $\varphi''(\frac{x_1}{\eps}) \geq - \frac{1}{2}$ and so $\nb_1^2 w_{0,{\rm ex}}(x) = \varphi''(\frac{x_1}{\eps}) \phi(x_2) \geq - \frac{1}{2}$. The symmetric statement for $x_1 - x_1^\vee(x_2) = x_1 \leq - \eps^{\frac 74}$ holds for the same reason. The last condition holds true because if $- \frac 13 \geq \nb_1 w_{0,{\rm ex}}(x) = \varphi'(\frac{x_1}{\eps}) \phi(x_2)$, then $\varphi'(\frac{x_1}{\eps}) \leq -\frac 13 < - \frac 14$, and hence $|\frac{x_1}{\eps}| \leq \sqrt{2}$. But in this region we have that $\varphi'''(\frac{x_1}{\eps}) = 1$, and hence $\nb_1^3 w_{0,{\rm ex}}(x) = \varphi'''(\frac{x_1}{\eps}) \phi(x_2) \geq \phi(x_2) \geq \frac{19}{20} > \frac 13$.
\end{itemize}
\begin{figure}[htb!]
\centering
\begin{minipage}{.45\linewidth}
  \includegraphics[width=\linewidth]{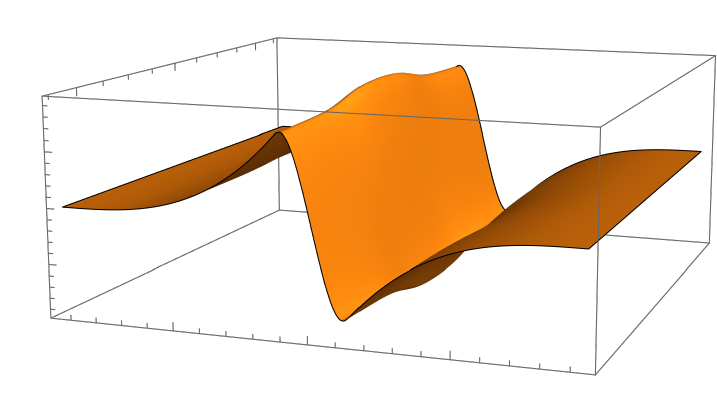}
\end{minipage}
\hspace{.05\linewidth}
\begin{minipage}{.45\linewidth}
  \includegraphics[width=\linewidth]{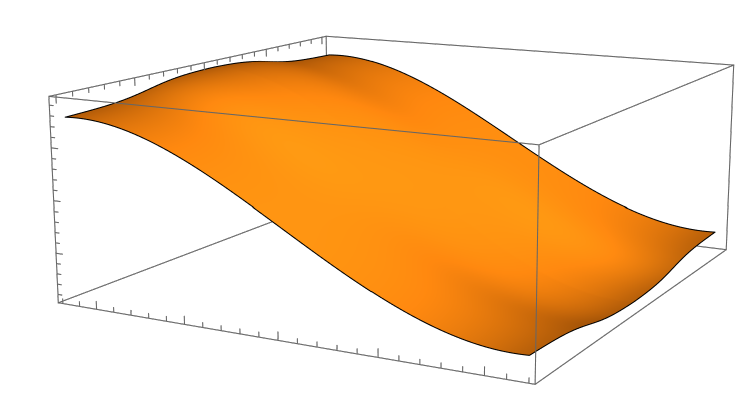}
\end{minipage}
\vspace{-.1 in}
\caption{\underline{Left}: a global view the function $w_{0,{\rm ex}}$, plotted for $|x_1|\leq 10\pi\eps$  and $|x_2| \leq \tfrac 14$, with $\eps=\tfrac{1}{20}$. \underline{Right}: a zoom-in of the function $w_{0,{\rm ex}} $, plotted for $|x_1|\leq \tfrac{\pi\eps}{2}$ and $|x_2| \leq \tfrac 18$, with $\eps=\tfrac{1}{20}$.}
\label{fig:w0:ex}
\end{figure} 

\begin{figure}[htb!]
\centering
\begin{minipage}{.45\linewidth}
  \includegraphics[width=\linewidth]{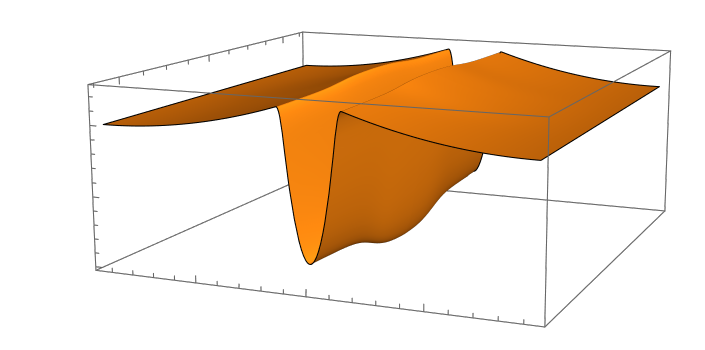}
\end{minipage}
\hspace{.05\linewidth}
\begin{minipage}{.45\linewidth}
  \includegraphics[width=\linewidth]{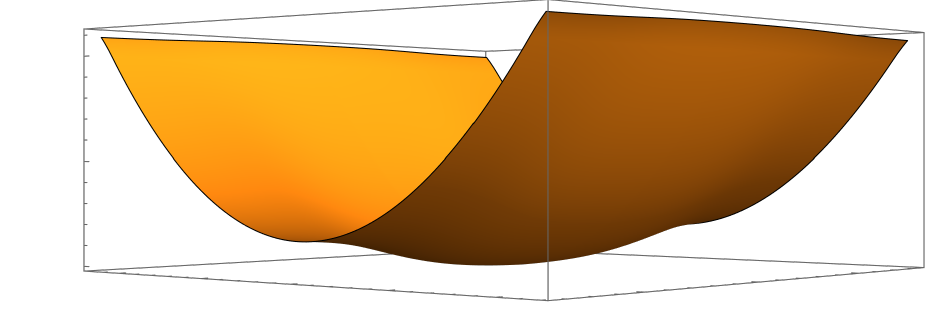}
\end{minipage}
\vspace{-.1 in}
\caption{\underline{Left}: a global view the function $\nb_1 w_{0,{\rm ex}}$, plotted for $|x_1|\leq 10\pi\eps$ and $|x_2| \leq \tfrac 14$, with $\eps=\tfrac{1}{20}$. \underline{Right}: a zoom-in of the function $\nb_1 w_{0,{\rm ex}} $, plotted for $|x_1|\leq \tfrac{\pi\eps}{2}$  and $|x_2| \leq \tfrac 18$, with $\eps=\tfrac{1}{20}$.}
\label{fig:w0,1:ex}
\end{figure} 

Next, we identify prototypical initial data for $a_0$ and $z_0$. Note that these fields only need to satisfy conditions~\eqref{item:ic:infinity} and \eqref{item:ic:reg}. As such, we may for instance take 
\begin{equation}
z_{0,{\rm ex}} (x_1,x_2)  = 0 
\, .
\label{eq:z0:ex}
\end{equation}
We may also take
\begin{subequations}
\label{eq:a0:ex}
\begin{equation}
a_{0,{\rm ex}}(x_1,x_2)  = 0
\,,
\label{eq:a0:ex:a}
\end{equation}
but maybe a more interesting prototypical example of initial data for $a_0$ is one for which the initial velocity $u_0$ is irrotational, which is equivalent to $\frac 12 \partial_2 w_0 = \partial_1 a_{0,{\rm ex}}$, and this is given by
\begin{equation}
a_{0,{\rm ex}} (x_1,x_2)= \tfrac{\eps}{2} \Phi(\tfrac{x_1}{\eps}) \phi'(x_2)
\,,
\label{eq:a0:ex:b}
\end{equation}
\end{subequations}
where $\phi$ is as in \eqref{eq:w0:ex}, and $\Phi$ is the compactly supported primitive of the function $\varphi$ from \eqref{eq:w0:ex}. That is, $\Phi(r) = \int_{-\infty}^r \varphi(r') {\rm d}r'$. By choosing $\varphi$ to be odd, we ensure that $\Phi$ is supported in $|r|\leq 13 \pi$, ensuring condition~\eqref{item:ic:supp}. Moreover, condition~\eqref{item:ic:infinity} holds as long as $\kappa_0 \geq 20$, so that \eqref{eq:a0:ex:b} is a permissible choice of initial data for $a_0$.
\end{example}

\begin{remark}[\bf An open set of initial data]
\label{rem:open:set:data}
We note that the initial data $(w_0,z_0,a_0)$ may be taken in an open set with respect to a suitable topology. The most direct way to see this is to consider a ball of radius $\eps^{N}$ with respect to the $H^7_0(\mathcal{X}_{\rm in})$ topology\footnote{Alternatively, we may take perturbations in the standard $H^7(\TT^2)$ topology which are sufficiently small. Such perturbations can potentially destroy the support assumption~\eqref{eq:ic:supp}. Nonetheless, because these perturbations are infinitesimal, and because the Euler system has finite speed of propagation, a suitable cutting procedure, and the classical local well-posedness theorem for the Euler system reduces the problem to studying a ``cut'' or ``localized'' initial data, which now does satisfy the support assumption~\eqref{eq:ic:supp}. We refer to~\cite[Section 13.7]{BuShVi2023b} for the detailed argument.}, with $N$ sufficiently large, centered at the prototypical functions $(w_{0,{\rm ex}}, z_{0,{\rm ex}}, a_{0,{\rm ex}})$, defined earlier in Example~\ref{ex:prototypical}---see~\eqref{eq:w0:ex}, \eqref{eq:z0:ex}, \eqref{eq:a0:ex}---with $\kappa_0 > 20$, and $\eps \in (0,\eps_0)$, where $\eps_0 = \eps_0(\alpha,\kappa_0)$ is a sufficiently small constant. Indeed, for any function in this ball, conditions~\eqref{item:ic:supp}--\eqref{item:ic:w0:d11:positive} are satisfied, if $N$ is taken to be sufficiently large. To see this, start with conditions~\eqref{item:ic:supp} and~\eqref{item:ic:reg}; these hold automatically for all functions in this ball by the definition of the $H^7_0(\mathcal{X}_{\rm in})$ norm, upon possibly enlarging the value of $\bar{\mathsf{C}}$. Condition~\eqref{item:ic:infinity} holds because for $\kappa_0>20$ the functions $(w_{0,{\rm ex}}, z_{0,{\rm ex}}, a_{0,{\rm ex}})$ satisfy these bounds with strict inequalities, and we have the Sobolev embedding $H^7_0(\mathcal{X}_{\rm in}) \subset L^\infty(\mathcal{X}_{\rm in})$. Similarly, the bounds $\nb_1 w_0 \leq \frac{1}{10}$, $|\nb_2 w_0| \leq 1$, and $|\nb\nb_1 w_0|\leq 2$ appearing in \eqref{item:ic:max:w0} hold if $\eps$ is taken to be sufficiently small, because by construction, the function $w_{0,{\rm ex}}$ satisfies these bounds with strict inequalities. The bound $\nb_1 w_0 (x)\geq -1$ holds in a vicinity of $x=0$ due to assumption~\eqref{item:ic:W:Taylor:1}, while for $x$ away from $0$ it holds because there  $w_{0,{\rm ex}}$ satisfies this bound  with a strict inequality. Similarly, the conditions on the initial data in \eqref{item:ic:w0:x2:negative},  \eqref{item:ic:w0:x2:special}, and~\eqref{item:ic:w0:d11:positive} are all bounds satisfied by $w_{0,{\rm ex}}$ with strict inequalities, and hence small enough smooth perturbations will also satisfy these bounds.  It thus remains to discuss assumption \eqref{item:ic:W:Taylor:1} on the initial data. It is clear that arbitrary small perturbations of $w_{0,{\rm ex}}$ may not anymore attain their global minimum exactly at $x=0$, or this minimum may not anymore equal exactly $-1$, or we may not anymore have the exact equality $\nb_2 w_0 = 0$ or $\nb^2 \nb_1 w_0 = {\rm Id}$ at this global minimum. Nonetheless, as we have previously discussed in~\cite[Section 13.7]{BuShVi2023b}, we may use the Galilean symmetry group and the  scaling invariance of the Euler system to relax the pointwise constraints in~\eqref{item:ic:W:Taylor:1}. For instance, small and smooth perturbations of $w_{0,{\rm ex}}$ will attain their global minimum at a point near $0$, which may then be shifted to be exactly at $0$ using translational invariance. An affine transformation of space may then be used to ensure that $\nb_2 w_0$ and $\nb \nb_1 w_0$ vanish at this point, and scaling may be used to enforce that $\nb_1 w_0$ equals $-1$. Our condition on the Hessian of $\nb_1 w_0$ is already an open condition, so it will be automatically satisfied for small perturbations. This concludes the proof of the fact that all smooth and sufficiently small perturbations of the prototypical functions constructed in Example~\ref{ex:prototypical} satisfy all the assumptions on the initial data: \eqref{item:ic:supp}--\eqref{item:ic:w0:d11:positive}.
\end{remark}

\begin{remark}[\bf Notation: usage of $\les$ and the dependence of of generic constants]
\label{rem:les}
Throughout the paper we shall write $A \les B$ to mean that there exists a constant $\Cn \geq 1$ such that $A \leq \Cn B$, where $\Cn$ is allowed to depend only on $\alpha$, $\kappa_0$, and $\Cdata$, but be independent of $\eps$. Throughout the paper we use $\Cn$ to denote a sufficiently large constant which depends only on $\alpha$, $\kappa_0$, and $\Cdata$, and which may change (increase) from line to line. We emphasize that $\Cn$ is never allowed to depend on $\eps$. Since $\eps$ will be chosen to be sufficiently small with respect to $\alpha,\kappa_0$, and $\Cdata$, we frequently write inequalities of the type $\eps \Cn \leq 1$.
\end{remark}

\subsection{Main results}
\label{subsec:mainresults}
The following three theorems are the main results of this paper. Theorem~\ref{thm:main:shock} concerns the process of shock formation, and is proven in Sections~\ref{sec:formation:setup}--\ref{sec:sixth:order:energy}. Theorem~\ref{thm:main:DS} concerns the spacetime of downstream \MGHDB\ of the initial data, and is proven in Section~\ref{sec:downstreammaxdev}. Theorem~\ref{thm:main:US} concerns the upstream spacetime of \MGHDB\ of the initial data, and is proven in Section~\ref{sec:upstreammaxdev}. See Figure~\ref{fig:3spacetimes} above. Additional optimal bounds for velocity, sound speed, and ALE map are reported in Section~\ref{sec:optimal:reg}  in all cases. 

\begin{theorem}[\bf Shock formation and the set of pre-shocks]
\label{thm:main:shock}
Fix $\alpha = \frac{\gamma-1}{2} > 0$, where $\gamma>1$ is the adiabatic exponent. Let $\kappa_0 \geq 20$ and $\bar{\mathsf{C}} \geq 1$ be two arbitrary constants. Then, there exists a sufficiently small $\eps_0 = \eps_0(\alpha,\kappa_0,\bar{\mathsf{C}}) \in (0,1]$ such that for every $\eps \in (0,\eps_0]$ the following holds. If the initial data $(u_0,\sigma_0)$--or equivalently, $(w_0,z_0,a_0)$ cf.~\eqref{eq:w0:z0:a0:def}--of the Euler equations at time $t=\initial$ (cf.~\eqref{tin}) satisfies assumptions~\eqref{item:ic:supp}--\eqref{item:ic:w0:x2:special} with parameters $(\alpha,\kappa_0,\bar{\mathsf{C}},\eps)$, then there exists a spacetime $\mathcal{P}$ and a time-dependent family of diffeomorphisms $\psi (\cdot,t) \colon \mathcal{P} \cap \{t\} \to \mathbb{R}^2$ such that the following hold:

\begin{enumerate}[leftmargin=26pt]
\renewcommand{\theenumi}{\alph{enumi}}
\renewcommand{\labelenumi}{(\theenumi)}

\item \label{item:shock:a} 
There exists a \underline{unique classical solution} $(u,\sigma)$ of the Cauchy problem for the Euler equations~\eqref{euler1} in the spacetime $\mathcal{P}_{\sf Eulerian} := \{ (\psi(x,t),t) \colon (x,t) \in \mathcal{P}\}$, with data $(u_0,\sigma_0)$. The solution $(u,\sigma$) is as smooth as the initial data, i.e., \underline{it does not lose derivatives}.

\item \label{item:shock:b} 
Each diffeomorphism $\psi(\cdot,t)$ is \underline{invertible} with ${\rm det}(\nabla \psi) > 0$ on $\mathcal{P}$, for every $t$ the map $x\mapsto \psi(x,t) - x$ is $\TT^2$-periodic, and $\psi$ is as smooth as the initial data, i.e., \underline{it does not lose derivatives}.

\item \label{item:shock:c} 
The map $\psi$ defines a smooth ALE coordinate system~\eqref{psi-def} on $\mathcal{P}$, with associated smooth normal \& tangent vectors $\nn$ \& $\tt$ defined via \eqref{tn-lag}, and smooth metric-normalized Jacobian determinant $\Jg \approx {\rm det}(\nabla \psi)$ defined via \eqref{Jg-def}. This \underline{ALE coordinate system flattens every fast acoustic characteristic surface} and allows us to characterize $\mathcal{P} = \{(x,t) \in \TT^2 \times [\initial,\medium] \colon \min_{x_1 \in \TT} \Jg(x_1,x_2,t) > 0\}$, cf.~\eqref{eq:spacetime:smooth}, where $\medium$ is given by \eqref{tmed}.  The spacetime $\mathcal{P}$ describes the Euler evolution for an $\OO(\eps)$ amount of time \underline{past the ``very first'' singularity} and satisfies  $\mathcal{P}  \subset \TT^2 \times [\initial,\medium]$.

\item \label{item:shock:d} 
The ``top'' boundary (future temporal boundary) of $\mathcal{P}$, i.e., $\partial_{\sf top}\mathcal{P} = \{(x,t) \in \TT^2 \times [\initial,\medium] \colon \min_{x_1 \in \TT} \Jg = 0\}$ contains the \underline{set of pre-shocks} $\Xi^*$, which parametrizes a \underline{cascade of first gradient catastrophes}, resulting from the distance between fast acoustic characteristic surfaces collapsing to zero. The set of pre-shocks  is a smooth co-dimension-$2$ subset of spacetime (see Definition~\ref{def:pre-shock}) characterized as the intersection of two co-dimension-$1$ surfaces: $\Xi^* = \{(x,t) \in \TT^2 \times [\initial,\medium] \colon \Jg(x,t) = 0\}  \cap  \{(x,t) \in \TT^2 \times [\initial,\medium] \colon \partial_1 \Jg(x,t) = 0\} \subset \partial_{\sf top}\mathcal{P}$. 

\item \label{item:shock:e}
The ALE coordinate system allows us to define, via \eqref{USigma} and~\eqref{BigWZA}, a new set of smooth multi-dimensional \underline{differentiated geometric Riemann variables} $(\Wb,\Zb,\Ab)$ whose time evolution is given by~\eqref{euler-WZA-real}. On the spacetime $\mathcal{P}$ the Euler equations~\eqref{euler1} for $(u,\sigma)$ are \underline{equivalent} to the evolution of the differentiated geometric Riemann variables, sound speed, and  of the geometry itself, via~\eqref{euler-WZA-real}, \eqref{Sigma0-ALE}, \eqref{nn-tt-evo}, and \eqref{Jg-system}. 

\item \label{item:shock:f}
The unique solution $(\Wb,\Zb,\Ab,\Sigma, \nn,\tt,\Jg)$ of this ALE-Euler system of equations--\eqref{euler-WZA-real}, \eqref{Sigma0-ALE}, \eqref{nn-tt-evo}, \eqref{Jg-system}--on the spacetime $\mathcal{P}$ \underline{maintains uniform $H^6$ Sobolev bounds} throughout the cascade of gradient catastrophes emerging on $\Xi^* \subset \partial_{\sf top}\mathcal{P}$. These Sobolev estimates propagate the regularity of the initial data, \underline{there is no} \underline{derivative loss}. The precise pointwise and energy estimates for $(\Wb,\Zb,\Ab,\Sigma, \nn,\tt,\Jg)$ are found in the bootstraps~\eqref{bootstraps}, the geometry bounds~\eqref{geometry-bounds-new}, the improved estimates~\eqref{eq:Jg:Abn:D5:improve}, \eqref{eq:Jg:Zbn:D5:improve}, and \eqref{eq:madman:2:all}, and in the optimal $H^7$ regularity bounds for $u\circ \psi$, $\sigma \circ \psi$, and $\psi$ reported in \eqref{eq:D7:U:Sigma:h}.

\item \label{item:shock:g}
No gradient singularity occurs at points in the closure of $\mathcal{P}$ which are away from the curve of pre-shocks. That is, for $(x_*,t_*) \in \partial_{\sf top} \mathcal{P} \setminus \Xi^*$ we have~$\lim_{\mathcal{P} \ni (x,t) \to (x_*,t_*)} (|\nabla u|, |\nabla \sigma|)\circ\psi(x,t) < + \infty$. On the other hand, for $(x_*,t_*) \in \Xi^*$  exactly one component of $(\nabla u) \circ \psi$ and one component of $(\nabla \sigma) \circ \psi$ blows up at $(x_*,t_*)$. With $n = \nn\circ \psi^{-1}$ and $\tau = \tt\circ \psi^{-1}$, we have that $\lim_{\mathcal{P} \ni (x,t) \to (x_*,t_*)} (|\tau \cdot \p_n u|, |\p_\tau u|, |{\rm curl}\, u|, |\partial_\tau \sigma|)\circ\psi(x,t) < + \infty$ and also $\lim_{\mathcal{P} \ni (x,t) \to (x_*,t_*)}(n\cdot \p_n u, {\rm div}\, u, \p_n \sigma) \circ \psi(x,t) = - \lim_{\mathcal{P} \ni (x,t) \to (x_*,t_*)} \Jgi(x,t) = - \infty$. That is, \underline{the singularities emerging on $\Xi^*$ are all pre-shocks}, and there are \underline{no other singularities on the closure of $\mathcal{P}$}.

\item \label{item:shock:h} 
With respect to the usual Eulerian variables $(y,t)$, the solution $(u,\sigma)$ \underline{inherits the $H^7$ regularity} from $U = u\circ \psi$, $\Sigma = \sigma \circ \psi$, and the $H^7$ invertible map $\psi$, in the interior of the spacetime $\mathcal{P}_{\sf Eulerian} = \{ (y,t) \colon y = \psi(x,t), (x,t) \in \mathcal{P}\}$. In particular,  $(u,\sigma) \in C^0_t C^5_y \cap C^5_t C^0_y$ is a classical solution of the Cauchy problem for the Euler equations in the interior of $\mathcal{P}_{\sf Eulerian}$. The ``top'' boundary of the spacetime $\mathcal{P}_{\sf Eulerian}$ contains the \underline{Eulerian curve of pre-shocks} defined as $\Xi^*_{\sf Eulerian} := \{(y,t) \colon y = \psi(x,t), (x,t) \in \Xi^*\}$. We have that $|\nabla u|$ and $|\nabla \sigma|$ remain bounded as we approach boundary points away from the curve of pre-shocks. As we approach points on the co-dimension-$2$ set of pre-shocks,  $n \cdot \p_n u$, ${\rm div}\, u$,  and $\p_n \sigma$ diverge towards $-\infty$ at a rate proportional to the spacetime distance to $\Xi^*_{\sf Eulerian}$, while $\tau \cdot \p_n u$, $\p_\tau u$, ${\rm curl}\, u$, and $\p_\tau \sigma$ remain bounded. 
\end{enumerate}
\end{theorem}

\begin{theorem}[\bf Downstream maximal globally hyperbolic development in a box]
\label{thm:main:DS}
Let $0< \eps \leq \eps_0(\alpha,\kappa_0,\bar{\mathsf{C}})$ be as in Theorem~\ref{thm:main:shock}, and assume that the initial data $(w_0,z_0,a_0)$ satisfies the same assumptions as in Theorem~\ref{thm:main:shock}.  If $w_0$ furthermore satisfies \eqref{item:ic:w0:d11:positive}, then there exists a spacetime $\mPds$ and a family of diffeomorphisms $\psi(\cdot,t) \colon \mPds \cap \{t\} \to\mathbb{R}^2$ such that $\mPds \supset \mathcal{P}$, $\psi\bigl|_{\mathcal{P}}$ is the same as the diffeomorphism $\psi$ from Theorem~\ref{thm:main:shock}, and such that the following hold: 

\begin{enumerate}[leftmargin=26pt]
\renewcommand{\theenumi}{\alph{enumi}}
\renewcommand{\labelenumi}{(\theenumi)}

\item \label{item:DS:a}  
There exists a \underline{unique classical solution} $(u,\sigma)$ of the Cauchy problem  for the Euler equations~\eqref{euler1} in the spacetime $\mathcal{P}^\sharp_{\sf Eulerian} := \{ (\psi(x,t),t) \colon (x,t) \in \mathcal{P}^\sharp\}$,  with data $(u_0,\sigma_0)$. The solution $(u,\sigma$) is as smooth as the initial data, i.e., \underline{it does not lose derivatives}, and $(u,\sigma)\bigl|_{\mathcal{P}}$ is the same as the solution $(u,\sigma)$ of Theorem~\ref{thm:main:shock}.

\item \label{item:DS:b}
Each diffeomorphism $\psi(\cdot,)$ is \underline{invertible} with ${\rm det}(\nabla \psi) > 0$ on $\mathcal{P}^\sharp$, for every $t$ the map $x\mapsto \psi(x,t) - x$ is $\TT^2$-periodic, and $\psi$ is as smooth as the initial data, i.e., \underline{it does not lose derivatives}. 
As in Theorem~\ref{thm:main:shock}, the diffeomorphism $\psi$ defines a \underline{smooth ALE coordinate system} on $\mathcal{P}^\sharp$, with associated smooth normal \& tangent vectors $\nn$ \& $\tt$, and smooth metric-normalized Jacobian determinant $\Jg \approx {\rm det}(\nabla \psi)$,  which flattens every fast acoustic characteristic surface.

\item \label{item:DS:c}
There exists a co-dimension-$1$ surface $\Pi$ parametrized as $\Pi = \{ (x_1^*(x_2,t),x_2,t) \colon (\cdot,x_2,t) \in \mathcal{P}\}$ such that the co-dimension-$2$ surface of pre-shocks defined in Theorem~\ref{thm:main:shock} is given by $\Xi^* = \Pi \cap \partial_{\sf top} \mathcal{P}$, and such that $\Pi \subset \{ \Jg,_1 = 0\}$ (cf.~\eqref{eq:x1star:def} and \eqref{eq:x1star:critical}). 
We say that a point $(x,t)$ lies \underline{upstream} of the surface $\Pi$ if $x_1< x_1^*(x_2,t)$, and write this as $(x,t) \in \Pi_-$. 
We say that a point $(x,t)$ lies \underline{downstream} of $\Pi$ if $x_1> x_1^*(x_2,t)$, and we write this as $(x,t) \in \Pi_+$. 

\item \label{item:DS:d}
The spacetime $\mathcal{P}^\sharp$ is characterized as follows. By~\eqref{eq:spacetime:P} and Remark~\ref{rem:oJJ:upstream} we have that $\mathcal{P}^\sharp \cap \Pi_-  = \mathcal{P} \cap \Pi_-$ and $\mathcal{P}^\sharp \cap \Pi = \mathcal{P} \cap \Pi$. In the \underline{downstream region}, the spacetime $\mathcal{P}^\sharp \cap \Pi_+$ \underline{is strictly larger than} $\mathcal{P} \cap \Pi_+$, and by~\eqref{eq:spacetime:P} and Remark~\ref{rem:oJJ:gg:0:maximal} it is given by $\mathcal{P}^\sharp \cap \Pi_+ = \{(x,t) \in \TT^2 \times [\initial,\medium] \colon x_1 > x_1^*(x_2,t), \Jg(x,t) > 0\}$. The new results in the present theorem  (when compared to Theorem~\ref{thm:main:shock}) concern the downstream part of $\mathcal{P}^\sharp$. The Euler evolution
within the spacetime $\mathcal{P}^\sharp \cap \Pi_+$ is the \MGHDB\ of the Cauchy data $(u_0,\sigma_0)$ within $\Pi_+$. In the closure of the upstream region, $\psi$ is the same as in Theorem~\ref{thm:main:shock} and moreover all results from Theorem~\ref{thm:main:shock} apply as is. 

\item \label{item:DS:e} 
The ``top'' boundary (future temporal boundary) of $\mathcal{P}^\sharp$ has global $W^{2,\infty}$ regularity, and is smooth (the level set of an $H^6$ function) on either side of the set of pre-shocks $\Xi^*$, which lies at the intersection of $\partial_{\sf top} \mathcal{P}^\sharp$ with the surface $\Pi$. A  \underline{surface of fast acoustic characteristic singularities} smoothly connects to the set of pre-shocks in the downstream part of the ``top'' boundary of $\mathcal{P}^\sharp$. This is the co-dimension-$1$ surface given explicitly as $\partial_{\sf top}\mathcal{P}^\sharp \cap \Pi_+ = \{(x,t) \in \TT^2 \times [\initial,\medium] \colon \Jg(x,t) = 0, x_1 >x_1^*(x_2,t)\}$. Since $\Jg \approx {\rm det}(\nabla \psi)$, this surface parametrizes gradient catastrophes resulting from impinging fast acoustic characteristic surfaces in $\Pi_+$, i.e., it is the ``envelope'' of the spacetime in which the fast acoustic characteristic surfaces remain in one-to-one correspondence with the initial foliation of spacetime.

\item \label{item:DS:f}
In the spacetime $\mathcal{P}^\sharp$, the smooth evolution~\eqref{euler-WZA-real} of the differentiated geometric Riemann variables $(\Wb,\Zb,\Ab)$, together with the evolution of sound speed and the geoemetry in~\eqref{Sigma0-ALE}, \eqref{nn-tt-evo}, and \eqref{Jg-system}, is \underline{equivalent} to the Euler equations~\eqref{euler1} for $(u,\sigma)$. The unique solution $(\Wb,\Zb,\Ab,\Sigma, \nn,\tt,\Jg)$ of this ALE-Euler system of equations maintains uniform $H^6$ Sobolev bounds on the spacetime $\mathcal{P}^\sharp$. These Sobolev estimates propagate the level regularity of the initial data, i.e., \underline{there is no derivative loss}. The precise pointwise and energy estimates are found in the bootstrap bounds~\eqref{bootstraps-P}, the geometry bounds in Proposition~\ref{prop:geometry-P}, the improved estimates~\eqref{improved-P} and \eqref{eq:madman:2-P}, and in the optimal $H^7$ regularity bounds for $u\circ \psi$, $\sigma \circ \psi$, and $\psi$ reported in \eqref{eq:D7:U:Sigma:h}.

\item \label{item:DS:g}
\underline{Gradient singularities occur at every point} which lies in the \underline{downstream part} of the \underline{``top'' boundary of $\mathcal{P}^\sharp$}. That is, for all $(x_*,t_*) \in \partial_{\sf top} \mathcal{P}^\sharp \cap \Pi_+$ we have that~$\lim_{\mathcal{P}^\sharp \ni (x,t) \to (x_*,t_*)} (|\tau \cdot \p_n u|, |\p_\tau u|, |{\rm curl}\, u|, |\partial_\tau \sigma|)\circ\psi(x,t) < + \infty$ and  $\lim_{\mathcal{P}^\sharp \ni (x,t) \to (x_*,t_*)}(n\cdot \p_n u, {\rm div}\, u, \p_n \sigma) \circ \psi(x,t) = - \lim_{\mathcal{P}^\sharp \ni (x,t) \to (x_*,t_*)} \Jgi(x,t) = - \infty$. The same type of gradient singularities occur on the set of pre-shocks $\partial_{\sf top} \mathcal{P}^\sharp \cap \Pi$. There are no gradient singularities on the upstream part of the ``top'' boundary of $\mathcal{P}^\sharp$, i.e., on $\partial_{\sf top} \mathcal{P}^\sharp \cap \Pi_-$.

\item \label{item:DS:h}
With respect to the Eulerian variables $(y,t)$, the solution $(u,\sigma)$ \underline{inherits the $H^7$ regularity} from $U = u\circ \psi$, $\Sigma = \sigma \circ \psi$, and the $H^7$ invertible map $\psi$, in the interior of the spacetime $\mathcal{P}^\sharp_{\sf Eulerian} = \{ (y,t) \colon y = \psi(x,t), (x,t) \in \mathcal{P}^\sharp\}$. In particular, $(u,\sigma) \in C^0_t C^5_y \cap C^5_t C^0_y$ is a classical solution of the Cauchy problem for the Euler equations in the interior of $\mathcal{P}_{\sf Eulerian}^\sharp$. The ``top'' boundary of the spacetime $\mathcal{P}_{\sf Eulerian}^\sharp$ has global $W^{2,\infty}$ regularity and is smooth on either side of the \underline{Eulerian curve of pre-shocks} $\Xi^*_{\sf Eulerian}$. \underline{Gradient singularities} occur \underline{at all points} which lie on the \underline{downstream part of $\partial_{\sf top} \mathcal{P}_{\sf Eulerian}^\sharp$}. Here,  $n \cdot \p_n u$, ${\rm div}\, u$, and $\p_n \sigma$ diverge towards $-\infty$ at a rate proportional to the spacetime distance to $\partial_{\sf top} \mathcal{P}_{\sf Eulerian}^\sharp$, while $\tau \cdot \p_n u$, $\p_\tau u$, ${\rm curl}\, u$, and $\p_\tau \sigma$ remain  bounded. 
\end{enumerate}
\end{theorem}

\begin{theorem}[\bf Upstream maximal globally hyperbolic development in a box]
\label{thm:main:US}
Fix $\alpha = \frac{\gamma-1}{2} > 0$ for $\gamma>1$. Let $\kappa_0 \geq 20$ be large enough with respect to $\alpha$ to ensure that~\eqref{eq:US:kappa:0:cond:0} holds. Let $\bar{\mathsf{C}} \geq 1$ and $\dl \in (0,\tfrac{1}{2} ]$ be  arbitrary. Then, there exists a sufficiently small $\eps_0 = \eps_0(\alpha,\kappa_0,\bar{\mathsf{C}},\dl) \in (0,1]$ such that for every $\eps \in (0,\eps_0]$ the following holds. If the initial data $(u_0,\sigma_0)$--or equivalently $(w_0,z_0,a_0)$--of the Euler equations at time $t=\initial$ satisfies assumptions~\eqref{item:ic:supp}--\eqref{item:ic:w0:d11:positive} with parameters $(\alpha,\kappa_0,\bar{\mathsf{C}},\eps)$, then there exists a spacetime $\tHdm$ and a time-dependent family of diffeomorphisms $\psi(\cdot,t) \colon \tHdm \cap \{t\} \to \mathbb{R}^2$, such that $\psi\bigl|_{\tHdm \cap \mathcal{P}}$ is the same as the diffeomorphism $\psi$ from Theorem~\ref{thm:main:shock}, and such that the following hold:
\begin{enumerate}[leftmargin=26pt]
\renewcommand{\theenumi}{\alph{enumi}}
\renewcommand{\labelenumi}{(\theenumi)}

\item \label{item:US:a}  
There exists a \underline{unique classical solution $(u,\sigma)$ of the Cauchy problem} for the Euler equations~\eqref{euler1} in the spacetime $\mathring{\mathcal{H}}^\dl_{\mathsf{Eulerian}} 
 := \{ (\psi(x,t),t) \colon (x,t) \in\tHdm\}$,  with data $(u_0,\sigma_0)$. The solution $(u,\sigma$) is as smooth as the initial data, i.e., \underline{it does not lose derivatives}, and $(u,\sigma)\bigl|_{\tHdm \cap \mathcal{P}}$ is the same as the solution $(u,\sigma)$ of Theorem~\ref{thm:main:shock}.

\item \label{item:US:b}  
Each diffeomorphism $\psi$ is \underline{invertible} with ${\rm det}(\nabla \psi) > 0$ on $\tHdm$, for every $t$ the map $x\mapsto \psi(x,t) - x$ is $\TT^2$-periodic, and $\psi$ is as smooth as the initial data, i.e., \underline{it does not lose derivatives}. As in Theorem~\ref{thm:main:shock}, the diffeomorphism $\psi$ defines a \underline{smooth ALE coordinate system} on $\tHdm$, with associated smooth normal \& tangent vectors $\nn$ \& $\tt$, and smooth metric-normalized Jacobian determinant $\Jg \approx {\rm det}(\nabla \psi)$,  which flattens every fast acoustic characteristic surface.

\item \label{item:US:c}  
Upstream of the surface $\Pi$ defined in Theorem~\ref{thm:main:DS}, item~\eqref{item:DS:c}, the spacetime $\Pi_-$ is foliated by slow acoustic characteristic surfaces which emanate from $\Pi = \{x_1 = x_1^*(x_2,t)\}$, at least locally for $x_1 < x_1^*(x_2,t)$. The portion of the \MGHDB\ of the initial data $(u_0,\sigma_0)$, which lies within $\Pi_-$, \underline{has as ``top'' boundary} (future temporal boundary) \underline{the unique slow acoustic characteristic surface emanating} \underline{from the set of pre-shocks} $\Xi^*$. 

\item \label{item:US:d}  
The spacetime $\tHdm$ is characterized as follows. For $\dl \in (0,\frac 12]$ \underline{arbitrary}, we consider (cf.~\eqref{sec:Thd:defn}) $\dl$-approximate slow acoustic characteristic surfaces $\{(x,\Thd(x,t))\}$. Among these surfaces there exists a unique and smooth \underline{distinguished $\dl$-approximate slow acoustic characteristic surface} $\{(x,\bar\Thd(x))\}$, which \underline{emanates from the set of} \underline{pre-shocks} $\Xi^*$. The spacetime $\tHdm$ is then characterized as $\{(x,t) \in \TT^2 \times [\initial,\final) \colon t < \bar{\Thd}(x) \}$. The new results in this theorem concern the upstream part of $\tHdm$, i.e.~$\tHdm \cap \Pi_-$. In the downstream region $\tHdm \cap \Pi_+$, the spacetime considered in 
Section~\ref{sec:upstreammaxdev} is a strict subset of the spacetime $\mathcal{P}$ from Theorem~\ref{thm:main:shock}, and all the results in Theorem~\ref{thm:main:shock} still apply, as is.

\item \label{item:US:e}  
The ``top'' boundary of $\tHdm$ is the surface $\partial_{\mathsf{top}} \tHdm = \{(x,\min\{\bar{\Thd}(x)),\final\} \}$ and \underline{the set of pre-shocks embeds} \underline{this future temporal boundary} as $\partial_{\mathsf{top}} \tHdm \cap \Pi = \Xi^*$. As $\dl \to 0^+$, the surface $\{(x,\bar{\Thd}(x)\}$ converges precisely to the slow acoustic characteristic surface emanating from the curve of pre-shocks, so that in the limit $\dl\to 0^+$ we  recover the entire upstream part of the \MGHDB\ spacetime.

\item \label{item:US:f}  
In the spacetime $\tHdm$ the smooth evolution~\eqref{euler-WZA-real} of the differentiated geometric Riemann variables $(\Wb,\Zb,\Ab)$, together with the evolution of sound speed and the geoemetry in~\eqref{Sigma0-ALE}, \eqref{nn-tt-evo}, and \eqref{Jg-system}, is \underline{equivalent} to the Euler equations~\eqref{euler1} for $(u,\sigma)$. The unique solution $(\Wb,\Zb,\Ab,\Sigma, \nn,\tt,\Jg)$ of this ALE-Euler system of equations maintains uniform $H^6$ Sobolev bounds on the spacetime $\tHdm$. These Sobolev estimates propagate the level regularity of the initial data, i.e., \underline{there is no derivative loss}. The precise pointwise and energy estimates are found in the bootstrap bounds~\eqref{bootstraps-H}, the geometry bounds in Proposition~\ref{prop:geometry-H}, the improved estimates in \eqref{eq:Jg:Abn:D5:improve-H}, \eqref{eq:Jg:Zbn:D5:improve-H}, and \eqref{eq:Jg:Wbn:improve:material-H}, and in the optimal $H^7$ regularity bounds for $u\circ \psi$, $\sigma \circ \psi$, and $\psi$ from~\eqref{eq:D7:U:Sigma:h}.

\item \label{item:US:g}  
No gradient singularity occurs at points in the closure of $\tHdm$ which are away from the curve of pre-shocks. That is, for $(x_*,t_*) \in \partial_{\mathsf{top}} \tHdm \setminus \Xi^*$ we have~$\lim_{\tHdm \ni (x,t) \to (x_*,t_*)} (|\nabla u|, |\nabla \sigma|)\circ\psi(x,t) < + \infty$. A different kind of singular phenomenon occurs in the upstream part of $\partial_{\mathsf{top}}  \tHdm$, in the limit as $\dl \to 0$: the ALE diffeomorphism $\psi$ cannot be extended beyond this Cauchy horizon in a unique and smooth fashion.

\item \label{item:US:h}  
With respect to the Eulerian variables $(y,t)$, the solution $(u,\sigma)$ \underline{inherits the $H^7$ regularity} from $U = u\circ \psi$, $\Sigma = \sigma \circ \psi$, and the $H^7$ invertible map $\psi$, in the interior of the spacetime $\mathring{\mathcal{H}}^\dl_{\mathsf{Eulerian}}  = \{ (y,t) \colon y = \psi(x,t), (x,t) \in \tHdm\}$. In particular, $(u,\sigma) \in C^0_t C^5_y \cap C^5_t C^0_y$ is a classical solution of the Cauchy problem for the Euler equations in the interior of $\mathring{\mathcal{H}}^\dl_{\mathsf{Eulerian}}$. The ``top'' boundary of  $\mathring{\mathcal{H}}^\dl_{\mathsf{Eulerian}}$ is smooth and the only \underline{gradient singularities} occur on the \underline{Eulerian curve of pre-shocks} $\Xi^*_{\sf Eulerian}$ which is embedded in $\partial_{\mathsf{top}} \mathring{\mathcal{H}}^\dl_{\mathsf{Eulerian}}$. 
\end{enumerate}
\end{theorem}

\subsection{The proofs of the main theorems}
\label{sec:proof:roadmap}
The remainder of the paper contains the proofs of Theorems~\ref{thm:main:shock},~\ref{thm:main:DS}, and~\ref{thm:main:US}.
The proofs of these theorems paint a much more detailed picture than what is summarized in the statements above, which only mentions the highlights. Here we provide a roadmap for the structure of these proofs (including the necessary forward references). The precise details are given in subsequent sections.

\subsubsection{The proof of Theorem~\ref{thm:main:shock}}
\label{sec:outline:thm:main:shock}
Assume that the initial data $(w_0,z_0,a_0) = (u_0 \cdot e_1 + \sigma_0, u_0\cdot e_1 - \sigma_0, u_0\cdot e_2)$ satisfies conditions~\eqref{item:ic:supp}--\eqref{item:ic:w0:x2:special} from Section~\ref{cauchydata}, for some parameters $\alpha$, $\kappa_0$, $\bar{\mathsf{C}}$ (independent of $\eps$), and $\eps>0$. As discussed in Remark~\ref{rem:open:set:data}, this constitutes an open set of initial data. These assumptions in particular give $(u_0,\sigma_0) \in H^7(\TT^2)$, and the initial density is bounded away from vacuum. By the classical local well-posedness theory for the isentropic Euler system in Sobolev spaces, we know that there exists a sufficiently small time $T > \initial$ and a unique  classical solution $(u,\sigma) \in C^0([\initial,T];H^7(\TT^2))$ of the Euler equations~\eqref{euler1}, such that all bounds on this solution are inherited from the initial data, up to a deterioration/magnification factor of $1+\eps$ for all norms. 

On $\TT^2\times [\initial,T]$ we may define as in Section~\ref{sec:acoustic:geometry} the Arbitrary-Lagrangian-Eulerian (ALE) coordinates $(\psi(x,t),t) = (h(x_1,x_2,t),x_2,t)$ adapted to the geometry of the fast acoustic characteristics. The family of diffeomorphisms $\psi(\cdot,t)$, which evolves according to \eqref{pt-psi}, induces a normal ($\nn$) and tangent ($\tt$) vector according to \eqref{tn-lag}, and a metric-normalized Jacobian determinant $\Jg$ according to \eqref{Jg-def}. In terms of this ALE geometry, we may then define as in Section~\ref{sec:new:Euler:variables} a new set of differentiated multidimensional Riemann variables $(\Wb,\Zb,\Ab)$, according to \eqref{USigma} and \eqref{BigWZA}. Since on $\TT^2\times [\initial,T]$ we are dealing with $C^2_{x,t}$ functions, and since the map $\psi$ is invertible (${\rm det}(\nabla \psi) \approx \Jg$ is bounded from below by a strictly positive number), by the construction of these differentiated Riemann variables we have that their evolution in~\eqref{euler-WZA-real}, together with the evolution of the rescaled sound speed~\eqref{Sigma0-ALE}, and that of the geometry~\eqref{nn-tt-evo}--\eqref{Jg-system} is in fact equivalent to the original Euler evolution in~\eqref{euler1}.  

The above described short-time analysis may then be extended to a larger spacetime via a classical continuation argument (which relies on local well-posedness of the Euler system and on finite speed of propagation) if we are able to guarantee that in this extended spacetime all the unknowns in the problem retain the regularity of the initial data (in this case, $H^7$ regularity for $U=u\circ\psi, \Sigma = \sigma\circ\psi$ and $\psi$ itself, and $H^6$ regularity for the geometric quantities $\nn,\tt,\Jg$ and for the differentiated Riemann variables $(\Wb,\Zb,\Ab)$), and if we are able to show that in this larger spacetime the family of diffeomorphisms $\psi$ remain invertible (that is, $\Jg>0$). A rigorous implementation of this continuation argument requires quantitative bounds on all unknowns in the problem, which we establish via a series of ``bootstrap inequalities'' for the solutions $(\Wb,\Zb,\Ab,\Sigma,\nn,\tt,\Jg)$ of~\eqref{euler-WZA-real},~\eqref{Sigma0-ALE},~\eqref{nn-tt-evo}--\eqref{Jg-system}. These bootstrap inequalities (see the inequalities in~\eqref{bootstraps} below) consist of pointwise bounds for the fields $(\Wb,\Zb,\Ab,\Sigma,\nn,\tt,\Jg)$ and their derivatives with respect to space and time, and of $L^2$-based energy bounds for derivatives up to order six for $(\Wb,\Zb,\Ab,\nn,\tt,\Jg)$ (the same as the regularity of these fields at the initial time). These bootstrap inequalities are stated in either the ALE spacetime coordinates $(x,t)$, or equivalently using a set of ``flattened'' $(x,\s)$ spacetime coordinates (given by~\eqref{t-to-s-transform} and \eqref{eq:f:tilde:f}), which are more convenient to use for energy estimates (see also Remark~\ref{rem:L2:norms:x:s:x:t:A} below). The spacetime $\mathcal{P}$ mentioned in Theorem~\ref{thm:main:shock} is then defined as a cylinder (meaning, it is invariant under translations in the $x_1$-variable) in which the map $\psi(\cdot,t)$ remains invertible, and which is quantified as $\mathcal{J}(x_2,t) = \min_{x_1\in\TT} \Jg(x_1,x_2,t) > 0$ (see~\eqref{eq:spacetime:smooth}). 

We note that the ``top'' boundary of the spacetime $\{\mathcal{J} > 0\}$ intersects the final time slice mentioned in Theorem~\ref{thm:main:shock}, namely $\{t = \medium\}$, in a Lipschitz (as opposed to $H^6$-smooth) fashion. In order to work with a spacetime which is as smooth as possible, in our case the zero level set of a $H^6$ function, in Section~\ref{sec:tin:tmed:tfin} we have included one more time-slice denoted by $\{t = \final\}$. On $\TT^2 \times [\medium,\final]$, which is a spacetime beyond the scope of Theorem~\ref{thm:main:shock}, we smoothly extend the Euler evolution in an artificial way (by working with the function $\Jgb$ defined in~\eqref{Jg-bar-evo}, instead of the natural Jacobian determinant $\Jg$), in order to ensure a smooth termination of our spacetime before the slice $\{t=\final\}$. While this extension is not seen in the statement of Theorem~\ref{thm:main:shock}, its use is very convenient for the proof, as it for instance ensures the flattening map $(x,t)\mapsto(x,\s)$ given by~\eqref{t-to-s-transform} retains maximal regularity, instead of being merely Lipschitz continuous. It is important to emphasize that this technically useful extension does not alter the Euler dynamics in any way on $\TT^2\times[\initial,\medium]$, and at first reading one should ignore the modifications due to $\Jg \mapsto \Jgb$.

The proof of Theorem~\ref{thm:main:shock} then consists of showing that in this dynamically defined spacetime the bootstrap inequalities may be ``closed'' if $\eps$ is taken to be sufficiently small. By ``closing the bootstrap assumptions'' we mean the standard continuity argument: assuming the inequalities in~\eqref{bootstraps} hold true with a specific constant on $\mathcal{P}$, we use the evolution equations~\eqref{euler-WZA-real},~\eqref{Sigma0-ALE},~\eqref{nn-tt-evo}--\eqref{Jg-system} to prove that the bounds in fact hold true on $\mathcal{P}$ with a strictly smaller constant than what was postulated. This requires a careful fine-tuning of the constants in the bootstrap assumptions, which is detailed in Remark~\ref{rem:order:of:constants} below. It is important here to notice that $\eps$ is the last parameter chosen in the proof, sufficiently small with respect to $\alpha,\kappa_0$, and $\bar{\mathsf{C}}$. 

The closure of the bootstrap inequalities~\eqref{bootstraps} is achieved in Sections~\ref{sec:formation:setup}--\ref{sec:sixth:order:energy} below. This necessitates a careful blend of pointwise bounds and energy estimates, which appeal not just the bootstraps themselves, but also to a number of bounds that are direct consequences of the bootstrap assumptions when combined with the ALE Euler evolution and the functional analytic framework from Appendix~\ref{app:functional} and Appendix~\ref{sec:app:transport}. In particular, the closure of the energy bootstraps requires that we carefully keep track of the behavior of all unknowns in the problem as we reach the ``top'' boundary of the spacetime $\mathcal{P}$. It is here that we encounter the co-dimension-$2$ set of pre-shocks $\Xi^*$ (see Defintion~\ref{def:pre-shock}), on which $\Jg$ vanishes ($\psi$ becomes not invertible). We keep track of the behavior of all unknowns in the vicinity of $\partial_{\sf top} \mathcal{P}$ using carefully chosen weights for the energy norms, in terms of fractional powers of $\mathcal{J}$ and $\Jg$; see the definitions of the energy and damping norms in Subsection~\ref{sec:norms:L2:first}.

Once the bootstraps  are closed on $\mathcal{P}$, we have established optimal  $H^6$ regularity estimates for $(\Wb,\Zb,\Ab,\nn,\tt,\Jg)$, and also the invertibility of the map $\psi$ (guaranteed by $\Jg>0$). Optimal $H^7$ bounds for velocity $U$, sound speed $\Sigma$, and ALE map $\psi$ are reported in Section~\ref{sec:optimal:reg}. By the Sobolev embedding and the (Sobolev) inverse function theorem, this implies the claimed $C^5_{x,t}$ regularity of $u,\sigma,\psi$ and $\psi^{-1}$ in the interior of $\mathcal{P}$, and also the equivalence of the system \eqref{euler-WZA-real},~\eqref{Sigma0-ALE},~\eqref{nn-tt-evo}--\eqref{Jg-system} with the original Euler evolution in~\eqref{euler1}. The claimed properties of the set of pre-shocks are established in Section~\ref{sec:x1star}.  

Lastly, the behavior of gradients of the solution as we approach the ``top'' boundary of the spacetime,  as claimed in~\eqref{item:shock:g} (or equivalently, \eqref{item:shock:h} in Eulerian variables) are now direct consequences of: the identities~\eqref{component-identities}, of the pointwise bootstrap bounds~\eqref{bootstraps}, the properties of $\mathcal{J}$ and $\Jg$ established in Sections~\ref{sec:Wbn:Jg:pointwise} and~\ref{sec:Jg:properties}, and the characterization of $\Xi^*$ in Proposition~\ref{prop:pre-shock}. For example, \eqref{comp-id-Wbn} and~\eqref{comp-id-Zbn} imply that $(n\cdot\p_n u)\circ \psi = \tfrac 12 (\Wbn+\Zbn)$ and $(\p_n \sigma)\circ\psi = \tfrac 12 (\Wbn-\Zbn)$. The bootstrap \eqref{bs-nnZb} gives that $\Zbn$ remains uniformly bounded in $\mathcal{P}$. The bootstrap \eqref{bs-Jg} together with the characterization of the pre-shock in Proposition~\ref{prop:pre-shock} show that as $\mathcal{P} \ni (x,t) \to \Xi^*$, we must have $\Jg(x,t) \to 0^+$ and hence $\Wbn(x,t) \leq - \tfrac{9}{10} \eps^{-1} \Jg(x,t)^{-1} \to -\infty$. This shows that $(n\cdot \p_n u, \p_n \sigma)\circ\psi(x,t) \to -\infty$ as $\mathcal{P} \ni (x,t) \to \Xi^*$.  On the other hand, \eqref{component-identities} shows that the gradients $(\tau \cdot \p_n u)\circ \psi, (\p_\tau u)\circ \psi $ and $(\p_\tau \sigma) \circ \psi$ may be computed solely in terms of $\Abn, \Wbt, \Zbt,\Abt$, whereas the bootstraps \eqref{bs-nnAb}, \eqref{bs-ttWb}, \eqref{bs-ttZb}, and \eqref{bs-ttAb} show that these terms remain uniformly bounded in $\mathcal{P}$. This shows that $(\tau \cdot \p_n u, \p_\tau u, \p_\tau \sigma)\circ\psi(x,t)$ remain bounded as $\mathcal{P} \ni (x,t) \to \Xi^*$. The fact that as $\mathcal{P} \ni (x,t) \to \partial_{\rm top}\mathcal{P} \setminus \Xi^*$ all gradients remain bounded is a consequence of the fact that {\em $\Jg$ does not vanish on $\partial_{\sf top}\mathcal{P} \setminus \Xi^*$}, which is in turn a consequence of the proof of Lemma~\ref{lem:x1*:unique}. The statements concerning ${\rm div}\, u$ and ${\rm curl}\, u$ follow from~\eqref{vort-id-good} and~\eqref{div-id-good}. The asymptotic behavior of gradients in Eulerian variables, as claimed in point~\eqref{item:shock:h}, follows identically; we omit these redundant details.
  
\subsubsection{The proof of Theorem~\ref{thm:main:DS}}
\label{sec:outline:thm:main:DS}

The proof is very similar in both spirit and implementation to the proof of Theorem~\ref{thm:main:shock}, outlined above. It is based on the equivalent formulation of the Euler equations in ALE variables from Sections~\ref{sec:acoustic:geometry} and~\ref{sec:new:Euler:variables}, and a continuation argument which is made quantitative via the propagation of bootstrap inequalities. This close resemblance allows us to  confine the entire analysis to one section, namely Section~\ref{sec:downstreammaxdev}, in which we highlight the details in the analysis which are different from the analysis in Sections~\ref{sec:formation:setup}--\ref{sec:sixth:order:energy}. 

The heart of the proof is to close bootstrap inequalities and to ensure that the map $\psi$ is smooth and invertible in the spacetime considered. The bootstrap inequalities themselves are the same as in Sections~\ref{sec:formation:setup}--\ref{sec:sixth:order:energy} and have been re-stated for convenience in~\eqref{bootstraps-P}.  The principal difference with respect to the analysis in Sections~\ref{sec:formation:setup}--\ref{sec:sixth:order:energy} is that in the downstream region, i.e., for $x_1 > x_1^*(x_2,t)$ (written as $\Pi_+$ in the statement of the theorem), we wish to extend the spacetime $\mathcal{P}$ to a strictly larger spacetime $\mathcal{P}^\sharp$, whose ``top'' boundary should be characterized by $\{ \Jg(x_1,x_2,t) = 0\}$, as opposed to $\{\mathcal{J}(x_2,t) = \min_{x_1} \Jg(x_1,x_2,t) = 0\}$, for times prior to $\medium$. In particular, this means that any parametrization of the downstream part of the ``top'' boundary of the spacetime necessitates $x_1$-dependence. In turn, this $x_1$ dependence enters the weight function $\JJ$  which replaces $\mathcal{J}$ in the downstream region (see definition~\eqref{eq:fake:Jg:def:ds}), and in the definition of the flattening map $\s = \qds(x,t)$ which replaces $\mathfrak{q}$ in the downstream region (see definitions~\eqref{eq:t-to-s-transform:all-P} and~\eqref{eq:f:tilde:f-P}). The closure of the energy bootstraps is then complicated by the appearance of an $\JJ,_1$ term in the energy estimates, and of  the coefficient $\Qb_1 = \eps \JJ,_1$ in the definition of the $\nbs_1$ operator (see~\eqref{QQQ-P} and~\eqref{nb-s-P}). This difficulty is overcome by noting that for $x_1$ which is in the downstream region $\mathcal{P}^\sharp \cap \Pi_+$ and is ``close'' to the co-dimension-$1$ set $\Pi = \{(x_1^*(x_2,t),x_2,t)\}$, we have that $\Jg,_1>0$, while for points in $\mathcal{P}^\sharp \cap \Pi_+$ which are far from $\Pi$, we have that $\Jg$ is bounded from below by a positive constant. A careful design of the weight function $\JJ$ and of the flattening map $\qds$ in the downstream region (see Section~\ref{sec:design:JJ:weight:DS}) then ensures $\JJ,_1$ is related to $\Jg,_1$ and thus has a favorable sign. This information is encoded through the fact that the coefficients $\Qb_1$ and $\Qr_1$ are non-negative (see~\eqref{eq:Qrs1:bbq-DS} and~\eqref{eq:Qr1:bbq-DS}), and hence certain dangerous terms in our energy estimates have a favorable sign. Physically speaking, this desired favorable sign in our energy estimates is a manifestation of the phenomenon of ``compression'', which is natural in the downstream region.

The weight function $\JJ$ and the flattening map $\qds$ are also carefully designed so that the spacetime $\mathcal{P}^\sharp \cap \Pi_+$ captures the full downstream part of the \MGHDB, for times prior to $\medium$ (see Remark~\ref{rem:oJJ:gg:0:maximal}). As in Sections~\ref{sec:formation:setup}--\ref{sec:sixth:order:energy}, we artificially extend  their definitions for times $t\in (\medium,\final]$ in order to ensure the smoothness of the ``top'' boundary of the downstream part of the extended spacetime $\mathcal{P}^\sharp$. As before, this extension of the ALE Euler dynamics to $\TT^2 \times (\medium,\final]$ (past the scope of Theorem~\ref{thm:main:DS}) is done for technical convenience only, and it does not affect the statement of Theorem~\ref{thm:main:DS}. 

The closure of the bootstraps corresponding to the spacetime $\mathcal{P}^\sharp$, is established in Sections~\ref{sec:bootstrap:consequences:DS}--\ref{sec:bootstrap:closed:DS}. Modifications to the argument in Sections~\ref{sec:formation:setup}--\ref{sec:sixth:order:energy} (which already covers the spacetime $\mathcal{P}^\sharp \cap \Pi_-$) are only required for the downstream part $\mathcal{P}^\sharp\cap\Pi_+$. The closure of these bootstraps then implies optimal  $H^6$ regularity estimates for $(\Wb,\Zb,\Ab,\nn,\tt,\Jg)$, the invertibility of the map $\psi$ (guaranteed by $\Jg>0$ in the interior of $\mathcal{P}^\sharp$), and optimal $H^7$ bounds for velocity $(U,\Sigma,\psi)$ (reported in Section~\ref{sec:optimal:reg}). The claimed $C^5_{x,t}$ regularity of $u,\sigma,\psi$ and $\psi^{-1}$ in the interior of $\mathcal{P}^\sharp$, and the equivalence of the system \eqref{euler-WZA-real},~\eqref{Sigma0-ALE},~\eqref{nn-tt-evo}--\eqref{Jg-system} with the original Euler evolution in~\eqref{euler1} directly follows. 

The fact that $\partial_{\sf top}\mathcal{P}^\sharp \cap \Pi_+$ and $\partial_{\sf top}\mathcal{P}^\sharp \cap \Pi_-$ are smooth (the zero level sets of $H^6$ functions) follows by construction. The fact that $\partial_{\sf top}\mathcal{P}^\sharp$ only retains $W^{2,\infty}$ regularity across its intersection with $\Pi$, i.e., at the pre-shock $\Xi^*$, is due to the fact that the second derivative with respect to $x_1$ of the weight function $\JJ$ is equal to $0$ as $x_1 \to x_1^*(x_2,t)^-$ (from the left, the upstream part), while the second $x_1$ derivative of the weight function $\JJ$ is strictly positive (due to \eqref{eq:Jg:11:lower}) as $x_1 \to x_1^*(x_2,t)^+$ (from the right, the downstream part). 

The remaining issue to discuss is the behavior of gradients of the solutions $(u,\sigma)$ discussed in item~\eqref{item:DS:g} (in ALE variables) and item~\eqref{item:DS:h} (in Eulerian variables). The novelty here lies in the statement that the gradient components $(n\cdot \p_n u,\p_n\sigma)\circ\psi(x,t)$ blow up as $(x,t) \in \mathcal{P}^\sharp$ approaches {\em any point} on the downstream part of the ``top'' boundary, $\partial_{\sf top}\mathcal{P}^\sharp \cap \Pi_+  \cap \{ t \leq \medium\}$, not just at points on the pre-shock $\Xi^*$ (as was shown in Theorem~\ref{thm:main:shock}). This fact is in a sense the very definition of downstream \MGHDB: components of $(\nabla u, \nabla \sigma)$ blow up {\em everywhere} on this future temporal boundary of the spacetime. In turn, this blow-up follows by our construction, which implies that this future temporal boundary $\{(x,t)\colon \JJ(x,t) = 0, x_1>x_1^*(x_2,t), t\leq \medium\}$ in fact equals the set $\{(x,t)\colon \Jg(x,t) = 0, x_1>x_1^*(x_2,t), t\leq \medium\}$ (see Remark~\ref{rem:oJJ:gg:0:maximal}), and thus $\Jg$ vanishes identically on this set. As discussed in the last paragraph of Section~\ref{sec:outline:thm:main:shock}, the asymptotic vanishing of $\Jg$ is equivalent to the divergence towards $-\infty$ of $\Wbn$, and thus also of $(n\cdot \p_n u,\p_n\sigma)\circ\psi = (\tfrac 12 (\Wbn+\Zbn), \tfrac 12 (\Wbn-\Zbn) )$. The bootstraps~\eqref{bootstraps-P} also imply uniform bounds for $(\Zbn,\Abn,\Wbt,\Zbt,\Abt)$ on $\mathcal{P}^\sharp$, showing that $(\tau\cdot \p_n u, \p_\tau u, \p_\tau \sigma) \circ \psi$ remain uniformly bounded on $\mathcal{P}^\sharp$. The statements concerning ${\rm div}\, u$ and ${\rm curl}\, u$ follow from~\eqref{vort-id-good} and~\eqref{div-id-good}. The asymptotic behavior of gradients in Eulerian variables, as claimed in point~\eqref{item:DS:h}, follows identically.

\subsubsection{The proof of Theorem~\ref{thm:main:US}}
\label{sec:outline:thm:main:US}
The proof follows the same strategy that was utilized in the proofs of Theorem~\ref{thm:main:shock} and~\ref{thm:main:DS} above: we use the equivalent formulation of the Euler equations in ALE variables from Sections~\ref{sec:acoustic:geometry} and~\ref{sec:new:Euler:variables}, and a continuation argument which is made quantitative via the propagation of bootstrap inequalities. We confine the entire proof to Section~\ref{sec:upstreammaxdev}, where we highlight the details in the analysis which are different from the analysis in Sections~\ref{sec:formation:setup}--\ref{sec:sixth:order:energy}. Nearly all differences arise due to the fact that we need to carefully analyze all slow acoustic characteristic surfaces emanating from the pre-shock, and its vicinity. 

The heart of the proof is to close bootstrap inequalities and to ensure that the map $\psi$ is smooth and invertible in the spacetime $\tHdm$. Both of these require a careful design and analysis of the upstream part of the spacetime, denoted as $\tHdm \cap \Pi_-$ in the statement of the theorem. 
While the bootstrap inequalities themselves are the same as in Sections~\ref{sec:formation:setup}--\ref{sec:sixth:order:energy}, and have been re-stated for convenience in~\eqref{bootstraps-H}, the meaning of the $L^2_x$-based norms present in~\eqref{boots-H} has been adapted to the upstream geometry (cf.~\eqref{norms-H}), and the weight function $\JJ$ present in the definition of the energy (cf.~\eqref{eq:tilde:E5E6-H}) and damping (cf.~\eqref{eq:tilde:D5D6-H}) norms  has undergone a significant transformation (see~\eqref{JJ-def-plus-t} and~\eqref{qps-JJ-tHdmm}) in order to account for the degeneracy in the problem which occurs along the slow acoustic characteristics emanating from the pre-shock. 

The intuition behind the construction of the weight function $\JJ$ in~\eqref{JJ-def-plus-t} and~\eqref{qps-JJ-tHdmm} is as follows. Based on intuition gained from Sections~\ref{sec:formation:setup}--\ref{sec:sixth:order:energy} and Section~\ref{sec:downstreammaxdev}, we need to design  the weight function $\JJ$ such that:
\begin{itemize}[leftmargin=16pt] 
\item $\JJ$ is $H^6$ smooth, the same regularity as $\Jg$,
\item the level set $\{\JJ = 0\}$ perfectly describes the future temporal boundary of the upstream part of the spacetime $\tHdm$, at least in the vicinity of the set of pre-shocks $\Xi^*$, where the gradient singularities are lurking,
\item such that the  action of the $\lambda_i$-transport operators $(\p_t + V \p_2) - (3-i) \alpha \Sigma (\Jgi \p_1 - g^{-\frac 12} h,_2 \p_2)$, gives a sign-definite term when acting on $\JJ$, for all $i \in \{1,2,3\}$.
\end{itemize}
The immediate issue is that as opposed to our earlier analysis we cannot let $\JJ$ equal simply to $\min_{x_1} \Jg$ (cf.~Sections~\ref{sec:formation:setup}--\ref{sec:sixth:order:energy}), or even $\Jg$ itself (cf.~Section~\ref{sec:downstreammaxdev}). This is because upstream of the pre-shock, the level set $\{\Jg =0\}$ describes the future temporal boundary of a spacetime which cannot be accessed by neither $\lambda_1$-characteristic surfaces (suitable for propagation of slow sound waves via $\Zb$) nor $\lambda_2$-characteristic surfaces (suitable for the propagation density waves via $\Sigma$, vorticity waves via $\Omega$, and tangential velocity waves via $\Ab$), which emanate from the initial data $(u_0,\sigma_0)$ at $t=\initial$. 

This begs the question: what is the maximal spacetime upstream of the pre-shock which is accessible by all characteristic surfaces (corresponding to the $\lambda_1$, $\lambda_2$, and $\lambda_3$ transport operators) emanating from the initial data? As discussed in Section~\ref{sec:usersguide}, this is the spacetime whose future temporal boundary is given by the slow acoustic characteristic surface (corresponding to the slowest wave-speed, $\lambda_1$) which emanates from the set or pre-shocks, in the upstream direction. This matches item~\eqref{item:US:c}  in the statement of Theorem~\ref{thm:main:US}. As discussed in~\eqref{1-flow},~\eqref{theta-dynamics-t}, and Figure~\ref{fig:upstream:maxdev} below, this surface would normally be characterized  as a graph over $(x_2,t)$,  by letting $x_1 = \theta(x_2,t)$ for a suitable function $\theta$. The fact that this surface emanates from the pre-shock $\Xi^* = \{ (\mathring{x}_1(x_2),x_2,t^*(x_2))\}$ is then the statement $\theta(x_2, t^*(x_2)) = \mathring{x}_1(x_2)$ for all $x_2 \in \TT$ (cf.~\eqref{theta-constraint}). The immediate issue which emerges is that the evolution equation for the function $\theta$ contains factors of $\Jgi$, which degenerate in as one approaches $\Xi^*$. Our observation is that if we re-parametrize the  slow acoustic characteristic surface emanating from $\Xi^*$ as a graph over $(x_1,x_2$), by letting $t =\Theta(x_1,x_2)$ for the function $\Theta$ such that $\Theta(\theta(x_2,t),x_2) = t$. A consequence of this re-parametrization is that the ``evolution equation'' for $\Theta(x)$ (we view $x_1$ as the evolution direction) now contains only factors of $\Jg$ (cf.~\eqref{p1-Theta-t-old-old}), which merely vanish as one approaches $\Xi^*$. This makes a smooth analysis of $\Theta$ accessible, and with that,  a smooth description of the spacetime of upstream \MGHDB. 

For technical reasons, related to the third bullet in the above-described requirements for $\JJ$, it is convenient to retain a damping term in our energy estimates (see the discussion in Remark~\ref{rem:JJ-formula}). As such, for $\dl>0$, arbitrarily small, we define a $\dl$-approximate slow acoustic characteristic surface passing though the curve of pre-shocks, and replace the $\Theta(x)$ described above by $\bar{\Thd}(x)$, as defined in~\eqref{eq:Thd:def}. Then, the $\dl$-adjusted upstream spacetime of \MGHDB\ of the Cauchy data, $\tHdm$, is characterized as the set $(x,t)$ such that $t< \bar{\Thd(x)}$, matching~\eqref{item:US:d} in the statement of Theorem~\ref{thm:main:US}. For convenience, we also require $t<\final$ in order to cap the time evolution at an $\OO(\eps)$ past the pre-shock. Then, the second bullet described above dictates that $\JJ$ needs to be  designed such that $\JJ(x,\bar{\Thd}(x))=0$ for all $x_1$ in the vicinity of the pre-shock $\Xi^*$ located at $x_1 = \mathring{x}_1(x_2)$ and $t = t^*(x_2)$. In order to ensure that $\JJ$ vanishes on $\partial_{\mathsf{top}} \tHdm$ for times $t< \final$, we thus design $\JJ$ as a $H^6$ smooth function in $\tHdm$, whose {\em zero level-set} is given by $\{ (x,\bar{\Thd}(x)) \}$. This matches the first two bullets in the above list of requirements for $\JJ$. The extension of $\JJ$ from $\partial_{\mathsf{top}} \tHdm$ down into $\tHdm$ is then made to also take into account the third bullet, by ensuring that the $\dl$-modified $\lambda_1$ transport operator $(1-\dl) (\p_t + V \p_2)  - 2\alpha \Sigma(\Jgi \p_1 - g^{-\frac 12} h,_2 \p_2)$ has $\JJ$ in its kernel (see~\eqref{JJ-formula-t-evo}). Additionally, we require that on the plane $\{(\mathring{x}_1(x_2),x_2,t)\}$ which emerges from the pre-shock at earlier times $t< t^*(x_2)$, the weight $\JJ$ precisely matches the function $\Jgb$ (cf.~\eqref{JJ-formula-t-BC}). This ensures that $\JJ$ vanishes not just on the pre-shock, but on the entire surface $(x,\bar{\Thd}(x))$, which is a characteristic surface for $(1-\dl) (\p_t + V \p_2)  - 2\alpha \Sigma(\Jgi \p_1 - g^{-\frac 12} h,_2 \p_2)$. This strategy is implemented by letting $\JJ(x,\Thd(x,t)) = \Jgb(\mathring{x}_1(x_2),x_2,t)$ for $t<t^*(x_2)$ and all $x_2 \in \TT$ (cf.~\eqref{JJ-def-plus-t}), where $\Thd(x,t)$ represents a family of characteristic surfaces for the $\dl$-approximate $\lambda_1$ transport operator $(1-\dl) (\p_t + V \p_2)  - 2\alpha \Sigma(\Jgi \p_1 - g^{-\frac 12} h,_2 \p_2)$, which emanate from the plane $\{(\mathring{x}_1(x_2),x_2,t)\}$ (see~\eqref{eq:Thd:PDE}). In fact, the surfaces $(x,\Thd(x,t))$, for $(x,t)$ as described in~\eqref{eq:Omega:US:+}, smoothly foliate a portion of the spacetime $\tHdm$ which has $\Xi^*$ on its future temporal boundary, labeled as $\tHdmp$ in the analysis, and defined in~\eqref{eq:H:dl:max:+:def}. The fact that $(x,\Thd(x,t))$ smoothly foliates $\tHdmp$ allows us to perform a smooth and sharp analysis upstream of the set of pre-shocks. The spacetime $\tHdm \setminus \tHdmp$ is denoted by $\tHdmm$ in \eqref{eq:H:dl:max:-:def}. Here the analysis is simpler because $\Jg \approx 1$ (see~\eqref{JJ-le-Jg:c}), we are ``far away'' from $\partial_{\mathsf{top}} \tHdm$, and so we just need to ensure that $\JJ$ satisfies the third bullet from the above list of requirements. This is implemented in~\eqref{qps-JJ-tHdmm}.

With the weight $\JJ$ and our spacetime $\tHdm$ defined precisely, the proof turns to the closure of the bootstrap assumptions, establishing the properties stated in item~\eqref{item:US:f} of the Theorem. As before, pointwise bootstraps are closed in $(x,t)\in\tHdm$ coordinates, while energy norms are estimated in flattened $(x,\s) \in \Hdm$ coordinates defined in~\eqref{t-to-s-transform-H}--\eqref{eq:f:tilde:f-H} below. Note that while the parameter $\dl$ does not enter the bootstrap assumptions explicitly, this parameter does affect the dependencies of the bootstrap constants themselves, cf.~Remark~\ref{rem:order:of:constants}. Our analysis shows that the constants appearing in items~\eqref{item:K:sf:choice} (corresponding to $\mathsf{K}$)--\eqref{item:eps:final:choice} (corresponding to $\eps$) of Remark~\ref{rem:order:of:constants} need to be chosen to depend on $\dl$. 

At the level of pointwise estimates, we highlight the lower bounds for $\JJ$ in~\eqref{JJ-and-t}, the fact that $\JJ$ gives a (good) signed contribution when acted upon by the $\{\lambda_i\}_{i=1}^3$ transport operators (this follows from~\eqref{eq:waitin:for:the:bus}), and the fact that $\Jg \geq \tfrac 13 \JJ$ (see~\eqref{JJ-le-Jg}). In particular, this last fact and the fact that $\JJ>0$ in the interior of $\tHdm$ ensures that $\Jg>0$ in the interior of $\tHdm$, implying the invertibility of $\psi$ claimed in item~\eqref{item:US:b} of the theorem.

At the level of energy estimates, several complications arise when compared to the analysis in Sections~\ref{sec:formation:setup}--\ref{sec:sixth:order:energy} and Section~\ref{sec:downstreammaxdev}. The principal new difficulty stems from the fact that  $\Hdm$ has a ``right lateral'' boundary located at $x_1 = \thd(x_2,\s)$ (see the definition in~\eqref{little-theta}). As such, the adjoint operator $\nbs^*$ corresponding to the $L^2$-norms defined in~\eqref{little-theta} contains a number of boundary terms~\eqref{eq:adjoints-H} at $x_1= \thd(x_2,\s)$.  At the top level of the energy estimates,  these boundary terms seem to be out of control; a more careful inspection, which uses fine properties of the spacetime $\Hdm$ and the design of the weight function $\JJ$ shows however that these boundary terms (or, suitable combinations of them) either have a good sign (see e.g.~\eqref{eq:got:my:mojo:workin}), or they vanish altogether since $\JJ$ vanishes (see for instance the proof of Proposition~\eqref{prop:vort:H6-H}). Another difficulty in closing the energy estimates stems from the $\lambda_i$-transport operators, written as $(\Q \p_\s + V \p_2) - (3-i) \alpha \Sigma (\Jgi \p_1 - g^{-\frac 12} \nbs_2 h \nbs_2)$ in $(x,\s)$ coordinates, acting on $\JJ$. Here, by design we obtain helpful signed damping terms, see e.g.~Remark~\ref{rem:JJ-formula} and the lower bound corresponding to $\Zbn$ in~\eqref{nicechickennice} (which is due to the fact that $\dl>0$). We also mention here that the functional analytic framework from Appendix~\ref{app:functional} can be adapted to the $(x,\s)$ coordinates considered in Section~\ref{sec:upstreammaxdev}, as discussed in Section~\ref{rem:app:upstream:flat} below. Similarly, the space-time $L^\infty$ estimates from  Appendix~\ref{sec:app:transport} also hold in the flattened upstream geometry, which changes to the proof of these estimates that are described in~Section~\ref{app:upstream:Lp}.

Concerning points~\eqref{item:US:e} and~\eqref{item:US:g} in the statement of the Theorem, we remark that at points $(x,t)$ such that $x_1 < x_1^*(x_2,t)$ and such that $t > \sup_{\dl\in(0,\frac 12]} \bar{\Thd}(x)$, that is, at points which lie upstream of the pre-shock and above the envelope of the $\dl$-approximate slow acoustic characteristic surfaces $(\cdot,\Thd(\cdot,\cdot))$, an Euler solution cannot be computed in a smooth and unique fashion from the initial data $(u_0,\sigma_0)$ at time $\initial$. This is  because a slow acoustic characteristic surface passing through $(x,t)$ would necessarily have to intersect (backwards in time) the surface $\partial_{\mathsf{top}} \mathcal{P}^\sharp \cap \Pi_+$, the downstream part of the top boundary of the space time $\mathcal{P}^\sharp$ constructed in Theorem~\ref{thm:main:DS}. But according to~Theorem~\ref{thm:main:DS}, item~\eqref{item:DS:g}, at every point on $\partial_{\mathsf{top}} \mathcal{P}^\sharp \cap \Pi_+$ a gradient singularity occurs both in density and in the normal derivative of the normal velocity, precluding a smooth continuation back to the initial data. 

In closing, we mention that the only gradient singularities for the fundamental variables $u$ or $\sigma$ which may be encountered on the closure of the spacetime $\tHdm$ (the closure of $\mathring{\mathcal{H}}^\dl_{\mathsf{Eulerian}}$ in Eulerian variables) occur on the set of pre-shocks $\Xi^*$ (denoted as $\Xi^*_{\sf Eulerian}$ in Eulerian variables), which is embedded in  the future temporal boundary of our spacetime, $\partial_{\mathsf{top}} \tHdm$. This fact was claimed in items~\eqref{item:US:g} and~\eqref{item:US:h} of the statement of the Theorem. Indeed, as was already discussed in proof of Theorem~\ref{thm:main:shock} and the proof of Theorem~\ref{thm:main:DS}, the only potential singularities permitted by the pointwise bootstrap assumptions are in $(n\cdot \p_n u,\p_n\sigma)\circ\psi(x,t)$, because these terms are computed in terms of $\Wbn$, while the bootstraps only  control   $\Jg \Wbn \approx (w_0),_1$. As before, for $(x_*,t_*) \in \Xi^*$ we have that $\lim_{\tHdm \ni (x,t) \to (x_*,t_*)} \Jg(x,t) = 0$, and thus $\Wbn \to -\infty$. However, for any $(x_*,t_*) \in \partial_{\mathsf{top}} \tHdm \setminus \Xi^*$, we either  have $\lim_{\tHdm \ni (x,t) \to (x_*,t_*)} \Jg(x,t) \geq \frac 19$ when ${\rm dist}( (x_*,t_*), \Xi^*) \gtrsim \eps$ due to \eqref{eq:Jg:useful:US:a}, or we have that $\lim_{\tHdm \ni (x,t) \to (x_*,t_*)} \Jg(x,t) \geq \tfrac{(x_{*1} - \mathring{x}_1(x_{*2}))^2}{14 \eps^2} \approx  (\tfrac{1}{\eps} {\rm dist}( (x_*,t_*), \Xi^*) )^2$, when  $0 < {\rm dist}( (x_*,t_*), \Xi^*) \ll \eps$ (due to \eqref{eq:stone:2} and~\eqref{eq:Err:lower:bound:impala}).  As such, $\lim_{\tHdm \ni (x,t) \to (x_*,t_*)} \Jg(x,t) >0$, and so the bounds for $\Jg \Wbn$ do imply a (finite) upper bound for $( |n\cdot \p_n u|, |\p_n\sigma|)\circ\psi(x,t)$, thereby concluding the proof.

\section{Shock formation: spacetime, energy  norms, and bootstrap assumptions}
\label{sec:formation:setup}
 
\subsection{Local well-posedness} 
\label{sec:local:existence:classic}
With the Cauchy data defined in Section~\ref{cauchydata}, the classical local well-posedness of the Euler system gives the existence of a time $ T \in (\initial,\final)$, such that uniform sixth-order energy estimates for solutions $(\Wb,\Zb,\Ab,\Jg,h,_2)$ to 
 \eqref{euler-WZA-aug}  are obtained on the time interval $[\initial,T]$, and the support of each solution at all times $t\in [\initial,T]$ is contained in the set 
\begin{equation}
\mathcal{X}_{\rm fin} 
:= 
\bigl\{ x\in \mathbb{T}^2 \colon {\rm dist}(x, \mathcal{X}_{\rm in}) \leq \Csupp \eps \bigr\}
\,,
\label{eq:supp:fin}
\end{equation}
where the constant $\Csupp>0$  depends only on $\alpha$ and $\kappa_0$ (see~\eqref{eq:Csupp:def} below). That is, solutions to the compressible Euler system have finite speed of propagation, and we are bounding   solutions for an amount of time which is at most $\final-\initial = \frac{2\eps}{1+\alpha} \cdot \frac{51}{50}$;  moreover, the local existence theory implies that 
\begin{equation}
\inf_{(x,t) \in \mathbb{T}^2 \times [\initial,T]} \Jg(x,t) > 0
\label{eq:Jg:pos}
\end{equation} 
which is to say that no collision of characteristics  occurrs on $[\initial,T]$. 

\subsection{A smooth remapping of spacetime}
\label{sec:remapping}

Our initial goal is to extend the ALE Euler solution $(\Wb,\Zb,\Ab,\Jg,h,_2)$ of \eqref{euler-WZA-aug} from the set $\mathcal{X}_{\rm fin}  \times [\initial,T]$ described in Section~\ref{sec:local:existence:classic} to a certain
spacetime $\mathcal{P}  \subset \mathcal{X}_{\rm fin} \times [\initial,\final]$, such  that 
\begin{itemize}[leftmargin=16pt] 
\item in this spacetime $ \mathcal{P} $,  the ALE Euler solution maintains uniform sixth-order Sobolev
bounds;
\item $\Jg >0$ in the interior of $\mathcal{P}$; and
\item  the boundary of $\mathcal{P}$ contains a co-dimension-$2$  surface on which $\Jg$ and $\Jg,_1$ vanish, the so-called {\em curve of pre-shocks}.
\end{itemize} 

\begin{figure}[htb!]
\centering
  \includegraphics[width=0.5\linewidth]{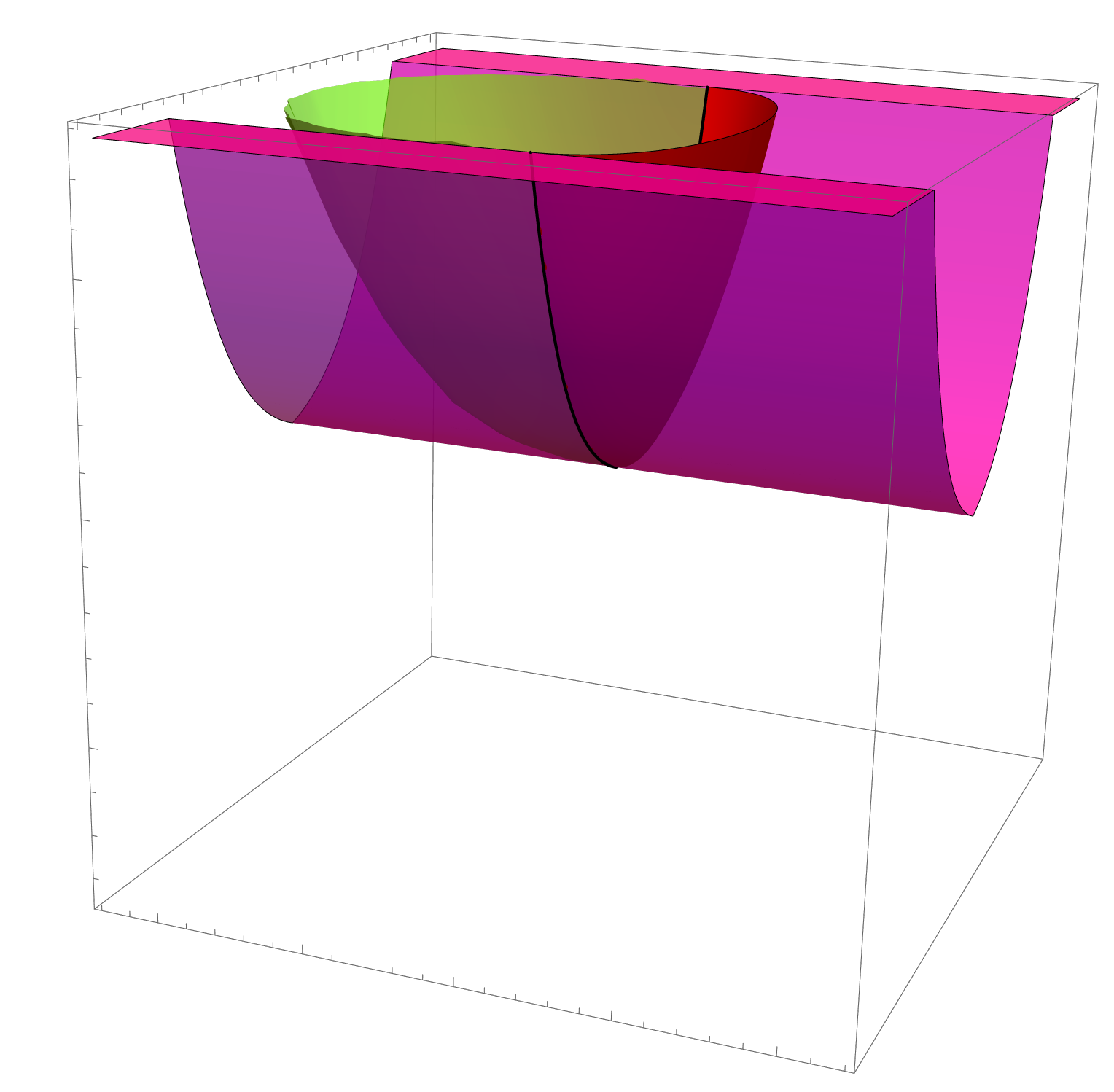}
\vspace{-.1 in}
\caption{Consider the function $\Jg = \Jg(x_1,x_2,t)$ resulting from the evolution of the initial data from Example~\ref{ex:prototypical}, with $(w_0,z_0,a_0) = (w_{0,{\rm ex}},0,0)$, with $\eps=\tfrac{1}{20}$. The bounding box represents the zoomed-in region $|x_1| \leq \tfrac{\pi\eps}{4}$, $|x_2| \leq \tfrac{1}{16}$, and $t\in [-\final,\final]$. In red, we plot the level set $\{(x,t) \in \TT^2 \times [\initial,\final] \colon \Jg(x_1,x_2,t) = 0, x_1>x_1^*(x_2,t)\}$. In magenta, we plot the {\it future} temporal boundary (or ``top'' boundary) of the spacetime $\mathring{ \mathcal{P} }$ from \eqref{eq:spacetime:temp}, which consists of the level  set $\{(x,t) \in \TT^2 \times [\initial,\final] \colon  \min_{x_1} \Jg(x_1,x_2,t) = \Jg(x_1^*(x_2,t),x_2,t)= 0\}$,  together with the flat portion of this temporal boundary $\{(x,\final) \in \TT^2 \times \{\final\} \colon   \min_{x_1} \Jg(x_1,x_2,\final) >0\}$. The {\it curve of pre-shocks}, represented in black, is defined by the intersection of the magenta and the red surfaces, i.e., it is the set $\{(x,t) \colon \Jg(x,t) = 0 = \Jg,_1(x,t) \} $. See also~Definition~\ref{def:pre-shock} below. Lastly, in green we represent the slow acoustic characteristic passing through the curve of pre-shocks.}  
\label{fig:parabolic:spacetime}
\end{figure}

A  priori, it is natural to consider  the  spacetime set 
\begin{equation}
\mathring {\mathcal{P}} := \bigl\{ (x,t) \in \mathbb{T}^2 \times [\initial,\final) \colon \min_{x_1 \in \mathbb{T}} \Jg(x_1,x_2,t) > 0 \bigr\}
 \,.
 \label{eq:spacetime:temp}
\end{equation}
We note that the {\it future} temporal boundary of the spacetime $\mathring{ \mathcal{P} }$ in  \eqref{eq:spacetime:temp} 
is not smooth along the intersection of the parabolic cylinder  $\{ (x_2,t) \colon  \min_{x_1 \in \mathbb{T}} \Jg(x_1,x_2,t)=0   \}$ 
and $\{t= \final\}$ (see the green surface in Figure~\ref{fig:parabolic:spacetime}). 
The lack of smoothness of this {\it future} temporal boundary along this intersection is an artifact of our
choice of the final time $t=\final$;  in particular, any ``final time'' which is $\OO(\eps)$
can be used in place of $\final$.
As such, we introduce a new spacetime, which coincides with $\mathring{ \mathcal{P} }$ for $t\in [\initial,\medium]$, 
but whose future temporal boundary is both smooth and properly contained in the set $\TT^2 \times [\initial,\final]$. 

For this purpose, we introduce a specially constructed modification of $\Jg$, which we denote by $\Jgb$, which has 
the following three properties:
\begin{enumerate}[leftmargin=26pt] 
\item  $\Jgb \equiv \Jg$ for all $t\in [\initial,\medium]$;
\item  $\Jgb$ has the identical regularity as $\Jg$; and
\item  for any $x\in \TT^2$, we have that $\Jgb(x,\cdot)$ vanishes at a time $t_*(x) \leq \final$.
\end{enumerate} 
More precisely, we define $\Jgb$ by modifying \eqref{Jg-system} as follows:
\begin{subequations} 
\label{def-Jgbar}
\begin{alignat}{2}
(\p_t+V\p_2)\Jgb  &=  \Jg  (  \tfrac{1+ \alpha }{2} \Wbn + \tfrac{1- \alpha }{2} \Zbn ) - \mathfrak{J}  \,,  && \text{ in } \  \mathbb{T}^2  \times [\initial, \final] \,, \label{Jg-bar-evo} \\
 \Jgb &=  1\,,  \ && \text{ on } \  \mathbb{T}^2  \times \{t=\initial\} \,, \label{ics-Jg-bar}
\end{alignat} 
\end{subequations} 
where $\mathfrak{J} = \mathfrak{J}(t) \geq 0$ is a smooth time-dependent function (independent of $x$), given by
\begin{equation}
\mathfrak{J}(t) = \tfrac{2(\final-\medium)}{\eps} \mathfrak{C}(\tfrac{t-\medium}{\final-\medium})
\,,
\label{eq:Xi:def}
\end{equation}
where $\mathfrak{C}$ is a $C^5$-smooth cut-off function, with $\mathfrak{C}(r) = 0$ for $r\leq 0$, with $0 < \mathfrak{C}(r) \leq 2$  and $0 < \mathfrak{C}'(r) \leq 4$ for $r \in (0,1]$, with $\int_0^1\mathfrak{C}(r) {\rm d}r = 1$, and with $\|\frac{d^k}{dr^k} \mathfrak{C}\|_{L^\infty(0,1)} \les 1$, where the implicit constant depends only on $k \in \{1,5\}$.

We note that in view of  \eqref{Jg-evo}  and \eqref{Jg-bar-evo} we have $(\p_t + V \p_2) (\Jgb-\Jg) = - \mathfrak{J}$, and thus we arrive at the identity
\begin{equation}
\Jgb  (x,t) = \Jg  (x,t) - {\bf 1}_{t>\medium} \int_{\medium}^t \mathfrak{J} (x_1,\xi(x_1,\xi^{-1}(x,t),t'),t'), {\rm d}t'
\,.
\label{eq:Jgb:identity:1}
\end{equation}
With the choice of $\mathfrak{J}$ in \eqref{eq:Xi:def}, we arrive at the explicit formula
\begin{equation}
\Jgb  (x,t) = \Jg  (x,t) - 2 {\bf 1}_{t>\medium}  \int_0^{\frac{t-\medium}{\final-\medium}} \mathfrak{C}(r) {\rm d}r \,.
\label{eq:Jgb:identity:0}
\end{equation}
This identity, the bootstrap \eqref{bs-Jg-simple} and continuity,  shows that for every $x\in \TT^2$, there exists a time $t_*(x) \in [\initial,\final]$ such that $\Jgb(x,t_*(x)) = 0$. Additionally, \eqref{eq:Jgb:identity:0} shows that  uniformly for $(x,t) \in \mathcal{P}$ we have
\begin{equation}
\p_1 (\Jgb - \Jg) \equiv 0\,,
\qquad
\p_2 (\Jgb -\Jg) \equiv 0\,,
\qquad 
- \tfrac{200(1+\alpha)}{\eps}  {\bf 1}_{t \in [\medium,\final]}\leq  \p_t (\Jgb - \Jg) \leq 0 \,,
\label{eq:Jgb:identity:2}
\end{equation}
and also 
\begin{equation}
\qquad 
|(\eps \p_t)^k (\Jgb - \Jg)| \les  {\bf 1}_{t \in [\medium,\final]}
\label{eq:Jgb:identity:3}
\end{equation}
for all $k \in \{1,\ldots,6\}$, where the implicit constant only depends on $\alpha$ and $k$. Note also that 
\begin{equation}
\Jgb \leq \Jg \,.
\label{Jgb-le-Jg} 
\end{equation}

\begin{figure}[htb!]
\centering
\begin{minipage}{.45\linewidth}
  \includegraphics[width=.93\linewidth]{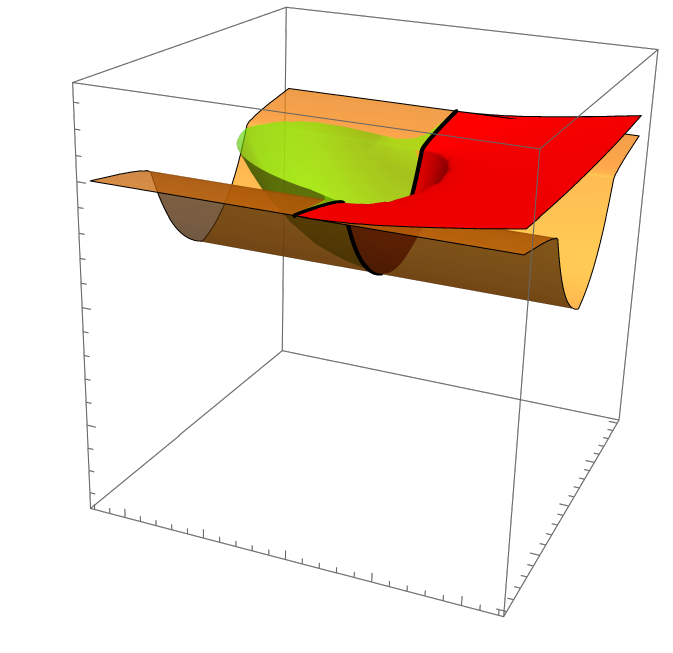}
\end{minipage}
\hspace{.05\linewidth}
\begin{minipage}{.45\linewidth}
  \includegraphics[width=\linewidth]{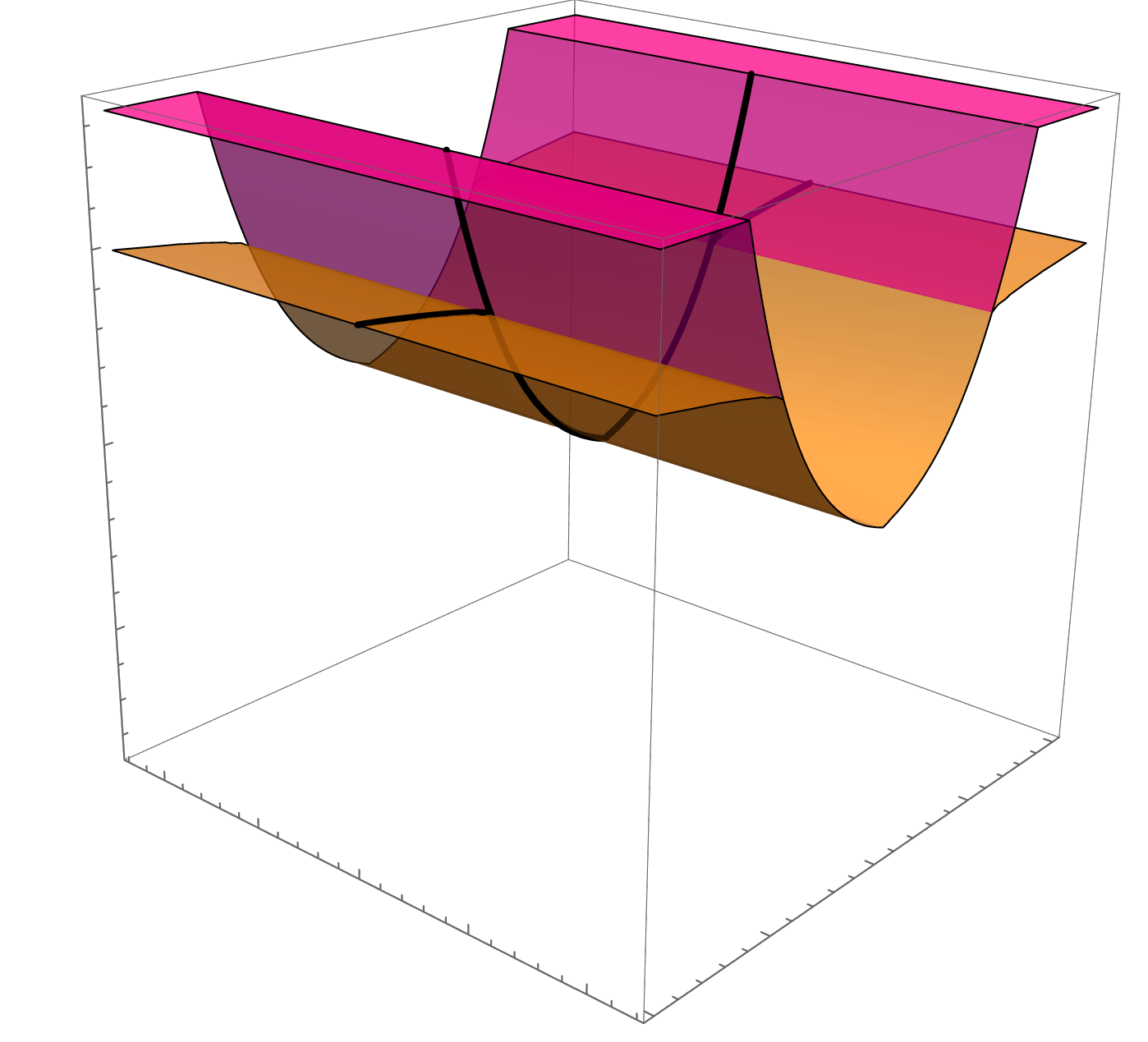}
\end{minipage}
\vspace{-.1 in}
\caption{Consider the function $\Jgb = \Jgb(x_1,x_2,t)$ as defined by~\eqref{eq:Jgb:identity:0}, with $\Jg$ as in Figure~\ref{fig:parabolic:spacetime}. The bounding box represents the zoomed-in region $|x_1| \leq \tfrac{\pi\eps}{4}$, $|x_2| \leq \tfrac{1}{16}$, and $t\in [-\final,\final]$. \underline{Left:} In red, we plot the level set $\{(x,t) \in \TT^2 \times [\initial,\final] \colon \Jgb(x_1,x_2,t) = 0, x_1 > x_1^*(x_2,t)\}$. In orange, we plot the {\it future} temporal boundary (or ``top'' boundary) of the spacetime $\mathring{ \mathcal{P} }$ from \eqref{eq:spacetime:temp}, which consists of the level  set $\{(x,t) \in \TT^2 \times [\initial,\final] \colon \mathcal{J}(x_2,t) = \min_{x_1} \Jgb(x_1,x_2,t)  = \Jgb(x_1^*(x_2,t),x_2,t)= 0\}$. The intersection of the orange and the red surfaces, i.e.,  the set  $\{(x,t) \colon \Jgb(x,t) = 0 = \Jg,_1(x,t) \}$ is represented in black. The slow acoustic characteristic emanating from this black curve is represented in green. \underline{Right:} in order to emphasize the fact that for $t\leq \medium$ the figure on the left side precisely matches Figure~\ref{fig:parabolic:spacetime}, we compare the level set $\{(x,t) \in \TT^2 \times [\initial,\final] \colon \mathcal{J}(x_2,t) = \min_{x_1} \Jgb(x_1,x_2,t) = 0\}$, in orange, with the level set $\{(x,t) \in \TT^2 \times [\initial,\final] \colon (x_2,t) = \min_{x_1} \Jg(x_1,x_2,t) = 0\} \cup  \{(x,\final) \in \TT^2 \times \{\final\} \colon   \min_{x_1} \Jg(x_1,x_2,\final) >0\}$, in magenta. In black we represent both $\{(x,t) \colon \Jg(x,t) = 0 = \Jg,_1(x,t) \}$ and $\{(x,t) \colon \Jgb(x,t) = 0 = \Jg,_1(x,t) \}$. It is clear that the modifications only occur for $t\in (\medium,\final]$, in order to ensure that the top boundary of the spacetime is smooth.
}  
\label{fig:parabolic:spacetime:bar}
\end{figure}

Next, we modify the spacetime $\mathring{ \mathcal{P} }$ of \eqref{eq:spacetime:temp},  and define the spacetime
\begin{equation}
\mathcal{P} := \bigl\{ (x,t) \in \mathbb{T}^2 \times [\initial,\final) \colon \min_{x_1 \in \mathbb{T}} \Jgb(x_1,x_2,t) > 0 \bigr\}
 \,.
 \label{eq:spacetime:smooth}
\end{equation}
A representation of the spacetime $\mathcal{P}$ is given in the left panel of Figure~\ref{fig:parabolic:spacetime:bar}. 
We shall prove in Lemma~\ref{lem:x1*:unique} below that for $(x,t) \in \mathcal{P}$ with $t> \initial$\footnote{Note that since $\Jgb(x,\initial) = 1$, the minimum of $\Jgb$ is not attained at a unique point when $t=\initial$.} the minimum with respect to $x_1$ of $\Jgb$ is attained at a unique point  $x_1^*(x_2,t)$, so that we have
\begin{equation}
 \Jgb(x_1^*(x_2,t),x_2,t) = \min_{x_1 \in \mathbb{T}} \Jgb(x_1,x_2,t)
 \label{eq:x1star:def}
 \,.
\end{equation}
In particular, for $t > \initial$, $x_1^*(x_2,t)$ is a critical point for the function $\Jgb,_1(\cdot,x_2,t)$, i.e.
\begin{equation}
\Jgb,_1(x_1^*(x_2,t),x_2,t) = 0. 
\label{eq:x1star:critical}
\end{equation}
For brevity of notation, throughout Sections~\ref{sec:formation:setup}--\ref{sec:sixth:order:energy} we shall denote
\begin{equation}
\mathcal{J}(x_2,t):= \Jgb(x_1^*(x_2,t),x_2,t).
\label{eq:fake:Jg:def}
\end{equation}
We note at this stage that due to \eqref{Jgb-le-Jg}, we may show that 
\begin{equation}
\mathcal{J}(x_2,t) \leq 1
\,.
\label{eq:Jgb:less:than:1}
\end{equation}
In order to prove \eqref{eq:Jgb:less:than:1}, we refer to \eqref{eq:d1w0:x1:star:temp}--\eqref{eq:d1w0:x1:star} below, which implies the bound $\Jg(x_1^*(x_2,t),x_2,t) \leq 1 + (t-\initial) \frac{1+\alpha}{2} (-\frac{9}{10\eps} + 2 \mathsf{C_{J_t}}) \leq 1 - (t-\initial) \frac{2(1+\alpha)}{5 \eps} \leq 1$, once $\eps$ is chosen sufficiently small.
We note that since $\Jgb - \Jg$ is independent of $x$, {\em the global minimum of $\Jgb(\cdot,x_2,t)$ is attained at the same point where the global minimum of $\Jg(\cdot,x_2,t)$ is attained}, and hence \eqref{eq:x1star:def} may be rewritten as
\begin{equation}
\Jg(x_1^*(x_2,t),x_2,t) = \min_{x_1 \in \mathbb{T}} \Jg(x_1,x_2,t)\,.
\label{eq:Jgb:Jg:identity:3} 
\end{equation}

Note that if $(x_1,x_2,t) \in \mathcal{P}$, then $(  x_1',x_2,t) \in \mathcal{P}$ for all $ x_1'$, and so this spacetime is invariant under shifts in the $x_1$ direction. It is thus convenient to define its projection onto the $(x_2,t)$ coordinates by:
\begin{equation}
\hat{\mathcal{P}} = \bigl\{ (x_2,t) \in \mathbb{T}  \times [\initial,\final) \colon \mathcal{J}(x_2,t) > 0 \bigr\}
 \,.
 \label{eq:spacetime:temp:x2:t}
\end{equation}
In order to perform energy estimates on $\mathcal{P}$, it is convenient to  introduce the transformation  
\begin{subequations}
\label{eq:t-to-s-transform:all}
\begin{align} 
\mathfrak{q} &\colon \hat{\mathcal{P}} \mapsto [0,\eps)\\
\s = \mathfrak{q}(x_2,t)&:= \eps \bigl(1 -  \mathcal{J}(x_2,t) \bigr) \,.
\label{t-to-s-transform}
\end{align} 
\end{subequations}
We have that $\mathfrak{q}(x_2,\initial) = 0$, so that the set  $\{\s = 0\}$ corresponds to the initial time slice $\{t=\initial\}$, which is the {\it past } temporal boundary of the projected spacetime $\hat{\mathcal{P}}$. We also note that the  {\it future}
temporal boundary of  $\hat{\mathcal{P}}$, namely the set 
\begin{equation*}
\partial_{\sf top} \hat{\mathcal{P}}= 
\{(x_2,t) \in \mathbb{T}  \times [\initial,\final] \colon  \mathcal{J}(x_2,t)=0  \} 
\end{equation*}
is mapped under $\mathfrak{q}$ to the set $\{\s = \eps\}$.  

Next, we define a suitable inverse of $\mathfrak{q}$ by 
\begin{subequations} 
\label{s-to-t-transform}
\begin{align} 
\mathfrak{q}^{-1}&: \TT \times [0,\eps)  \to [\initial,\final)  \,, \\
t & = \mathfrak{q}^{-1}(x_2,\s) \,,
\end{align} 
\end{subequations} 
such that $t = \mathfrak{q}^{-1}(x_2,\mathfrak{q}(x_2,t))$ for all $(x_2,t) \in \hat{\mathcal{P}}$, or equivalently, that $\s = {\mathfrak q}(x_2,\mathfrak{q}^{-1}(x_2,\s))$ for all $(x_2,\s) \in \TT \times [0,\eps)$. In \eqref{s-to-t-transform} we are abusing convention: it is the map $(x_2,\s) \mapsto (x_2,t)$ defined from $\TT \times [0,\eps) \to \TT \times [\initial,\final)$ which is the inverse of the map $(x_2,t) \mapsto (x_2,\s) = (x_2,\mathfrak{q}(x_2,t))$. 
The fact that such a map is well-defined is established in Lemma~\ref{lem:q:invertible} below.  

\subsection{Change of coordinates for the remapped spacetime}
\label{sec:remapping:1}
Given any function $f\colon \mathcal{P} \to \mathbb{R}$, where we recall cf.~\eqref{eq:spacetime:smooth} that $\mathcal{P} \subset \TT^2 \times [\initial,\final)$, we define the function $\tilde f \colon \TT^2 \times [0,\eps) \to \mathbb{R}$ by
\begin{equation}
\tilde f(x,\s) := f(x,t), \qquad \mbox{where} \qquad \s = \mathfrak{q}(x_2,t) \,.
\label{eq:f:tilde:f}
\end{equation}
Then, by the chain-rule, \eqref{t-to-s-transform},  and  \eqref{eq:x1star:critical} we obtain  
\begin{subequations}
\label{eq:xt:xs:chain:rule}
\begin{align} 
\p_t f(x,t) &= \Qd(x_2,\s)   \p_\s \tilde f(x,\s) \,, \\
\p_2f(x,t) &=  \bigl(\p_2 - \Qb_2(x_2,\s)  \p_\s\bigr) \tilde f(x,\s)\,,  \\
\p_1f(x,t) &= \p_1 \tilde f (x,\s) \,,
\end{align} 
\end{subequations}
where for compactness of notation we have introduced the functions 
\begin{subequations} 
\label{QQQ}
\begin{align}
\Qd (x_2,\s) &= \p_t \mathfrak{q}(x_2,t) \Big|_{t= \mathfrak{q}^{-1}(x_2,\s)} 
= - \eps (\p_t \Jgb)(x_1^*(x_2,t),x_2,t) \Big|_{t= \mathfrak{q}^{-1}(x_2,\s)}
\label{eq:QQQ:a}
\\
\Qb_2(x_2,\s) &= - \p_2 \mathfrak{q}(x_2,t) \Big|_{t= \mathfrak{q}^{-1}(x_2,\s)} =  \eps  (\p_2 \Jgb)(x_1^*(x_2,t),x_2,t)   \Big|_{t= \mathfrak{q}^{-1}(x_2,\s)}\,.
\label{eq:QQQ:aa}
\end{align}
For later use, it is also convenient to define 
\begin{align}
\Q( x,\s )  
:= \Qd(x_2,\s) - \tilde V(x,\s)  \Qb_2(x_2,\s) 
&= - \eps (\p_t \Jgb + V \p_2 \Jgb)(x_1^*(x_2,t),x_2,t)  \Big|_{t= \mathfrak{q}^{-1}(x_2,\s)}\notag \\
&\qquad -\bigl(V(x_1,x_2,t) - V(x_1^*(x_2,t),x_2,t)\bigr)\Big|_{t= \mathfrak{q}^{-1}(x_2,\s)} \Qb_2(x_2,\s) \,,
\label{eq:QQQ:b}
\end{align} 
and 
\begin{equation} 
\Qc = \p_\s \Q \,, \qquad 
\Qr_\s  = \p_\s \Qd  \,, \qquad 
\Qr_2 =  \p_\s \Qb_2 
\,.
\label{eq:QQQ:c}
\end{equation} 
\end{subequations} 
With the above notation, it follows from \eqref{eq:xt:xs:chain:rule} that the spacetime gradient operator in  $(x,t)$ variables, namely $\nb = (\eps\p_t, \eps \p_1, \p_2)$, becomes the gradient operator $\nbs$ associated with the $(x,\s)$ coordinates, which is defined by
\begin{equation} 
\nbs = (\nbs_\s, \nbs_1, \nbs_2) := \big( \eps \Qd  \p_\s ,  \eps \p_1  ,  \p_2 -\Qb_2 \p_\s \big) \,. 
\label{nb-s}
\end{equation} 
That is, we have that $$\nb f (x,t) = \nbs \tilde f(x,\s).$$ 
Next, we notice that the components of $\nbs$ commute, that is
\begin{equation} 
\jump{\nbs_\s, \nbs_2} =
\jump{\nbs_\s, \nbs_1} =
\jump{\nbs_2, \nbs_1} =
0
\,,
\label{comm-nbs1-nbs2}
\end{equation} 
so that for any $\gamma \in \mathbb{N}_0^3$ we may write unambiguously $\nbs^\gamma = \nbs_\s^{\gamma_0} \nbs_1^{\gamma_1} \nbs_2^{\gamma_2}$, and notice that 
\begin{equation}
\label{nb-to-nbs}
(\nb^\gamma f) (x,t) = (\nbs^\gamma \tilde f)(x,s)
\,.  
\end{equation}
Via the identity $\Q\p_\s +V \p_2= \tfrac{1}{\eps}\nbs_\s + \tilde V\nbs_2$, we note that   material derivatives are mapped into $(x,\s)$ coordinates as
\begin{equation} 
(\p_t  +V \p_2) f(x,t)  = (\Q\p_\s + \tilde V\p_2) \tilde f(x,\s)  = ( \tfrac{1}{\eps}\nbs_\s + \tilde V\nbs_2) \tilde f(x,\s) 
\,. \label{the-time-der}
\end{equation} 
It also follows from \eqref{comm-nbs1-nbs2} and the second equality in \eqref{the-time-der} that 
\begin{equation} 
\jump{ (\Q\p_\s +\tilde V \p_2) , \nbs^k} \tilde f  = \jump{\tilde V, \nbs^k} \nbs_2 \tilde f 
= -\nbs^k \tilde V \, \nbs_2 \tilde f - \doublecom{\nbs^k, \tilde V, \nbs_2 \tilde f} 
\,. \label{good-comm}
\end{equation}
Lastly,  we may identify the adjoint of $\nbs$ with respect to the $L^2$ inner product on $ \TT^2\times[0,\s]$ by 
\begin{subequations}
\label{eq:adjoints}
\begin{align}
\nbs_\s^* &=  - \nbs_\s - \eps  \Qr_\s + \eps \Qd ( \delta_{\s} - \delta_{0} ) \,, \label{adjoint-s}
\\
\nbs_1^* &=  - \nbs_1   \,,  \label{adjoint-1}
\\
\nbs_2^* &=  - \nbs_2 + \Qr_2 -  \Qb_2 \delta_{\s} \,,
 \label{adjoint-2}
\\
(\Q \p_\s + \tilde V \p_2)^* & = - (\Q \p_\s + \tilde V \p_2)  - \Qr_\s + \Q  ( \delta_{\s} - \delta_{0} ) + \tilde V  \Qr_2   - \nbs_2 \tilde V
 \label{adjoint-3}
 \,.
\end{align}
\end{subequations}
Here we have used that $\Jgb(x,\initial)=1$, so that $\Qb_2(x,0) = 0$.

\begin{remark}[\bf Lower bound for $\tilde \Jg$ and definition of ``fake $\Jg$'']
\label{rem:tilde:Jg}
By using \eqref{eq:f:tilde:f} with $f=\Jgb$, and the definition of the map ${\mathfrak q}$, we deduce that  
\begin{equation}
  \tilde \Jgb(x,\s) = \Jgb(x,t) \geq \min_{x_1} \Jgb(x,t) = \mathcal{J}(x_2,t) = \big(1 - \tfrac{\s}{\eps}\bigr)  = \tilde{\mathcal{J}}(x_2,\s)
  \label{eq:fake:Jg}
 \,.
\end{equation}
Throughout the paper, we shall refer to $\mathcal{J}(x_2,t) = \tilde{\mathcal{J}}(x_2,\s)$ as ``fake $\Jg$''.
We shall discuss in Section~\ref{sec:first:consequences} several useful properties of $\tilde{\mathcal{J}}$. Moreover, note that $\tilde{\mathcal{J}}$ does not in fact depend on $x_2$ at all.
\end{remark}

\begin{remark}[\bf Dropping the tildes]
\label{rem:no:tilde}
Rather than working with a new family of variables that depend on the spacetime coordinates $(x,\s)$, namely
$\widetilde \Wb(x,\s) = \Wb(x,t)$, $\widetilde \Zb(x,\s) = \Zb(x,t)$, $\widetilde \Ab(x,\s) = \Ab(x,t)$, $\widetilde \Jg(x,\s) = \Jg(x,t)$, $\widetilde \Jgb(x,\s) = \Jgb(x,t)$, $\tilde{\mathcal{J}}(x,\s) = \mathcal{J}(x,\s)$, 
$\widetilde h(x,\s) = h(x,t)$, $\widetilde g(x,\s) = g(x,t)$, $\widetilde \nn(x,\s) = \nn(x,t)$, and $\widetilde \tt(x,\s) = \tt(x,t)$, for notational simplicity we  drop the  tilde and  abuse notation to
continue using the variables $\Wb,\Zb,\Ab, \Jg, \Jgb, \mathcal{J}, h, g, \nn, \tt$, but now depending on $(x,\s)$ rather than $(x,t)$. This identification is made throughout the rest of the paper and no ambiguity may arise because we shall still use the notation $\nbs$ for the spacetime derivative operator in $(x,\s)$ coordinates. As such, $\nbs f$ means that $f$ is viewed as a function of $(x,\s)$, while $\nb f$ means that $f$ is viewed as a function of $(x,t)$, where $t= \mathfrak{q}^{-1}(x_2,\s)$. 
\end{remark}

At this stage it is convenient to record a few of  the evolution equations transformed into  $(x,\s)$ coordinates. For instance, \eqref{p2-h}, \eqref{Jg-evo}, and \eqref{Jg-bar-evo} imply that (as mentioned in Remark~\ref{rem:no:tilde}, we drop the tildes)
\begin{align} 
(\Q\p_\s+V\p_2) \Jg &= \tfrac{1+ \alpha }{2} \Jg\Wbn + \tfrac{1- \alpha }{2} \Jg\Zbn \,, \label{Jg-evo-s} \\ 
(\Q\p_\s+V\p_2) \Jgb &= \tfrac{1+ \alpha }{2} \Jg\Wbn + \tfrac{1- \alpha }{2} \Jg\Zbn  - \mathfrak{J} \,, \label{Jgb-evo-s} \\
(\Q\p_\s+V\p_2) \nbs_2 h &=  g  \bigl(\tfrac{1+ \alpha }{2} \Wbt + \tfrac{1- \alpha }{2} \Zbt \bigr) \,, \label{p2h-evo-s}
\end{align} 
 from \eqref{ics-U-Sigma} and \eqref{grad-Sigma} we deduce 
\begin{subequations} 
\label{Sigma-s}
\begin{align} 
{\tfrac{1}{\eps}} \nbs_1 \Sigma &=  \tfrac{1}{2} \Jg(\Wbn -\Zbn)  + \tfrac{1}{2} \Jg \nbs_2 h (\Wbt -\Zbt) \,, \label{p1-Sigma-s} \\
\nbs_2 \Sigma &= \tfrac{1}{2} g^{\frac{1}{2}}  (\Wbt-\Zbt)\,,  \label{p2-Sigma-s} \\
(\Q\p_\s+V\p_2) \Sigma  &= -  \alpha \Sigma (\Zbn+ \Abt )  \,,  \label{Sigma0-ALE-s}  \\
(\Q \p_\s +    V \p_2) \Sigma^{-2\beta}   &=   2\alpha\beta \Sigma^{-2\beta}  (\Zbn+ \Abt )  \,,  \label{Sigma0i-ALE-s} 
\end{align} 
\end{subequations}
while transforming \eqref{nn-tt-evo} yields
\begin{subequations} 
\label{nn-tt-evo-s}
\begin{align} 
(\Q\p_\s+V\p_2)\nn  & =-  \bigl(\tfrac{1+ \alpha }{2} \Wbt + \tfrac{1- \alpha }{2} \Zbt \bigr) \tt   \,, \label{nn-evo-s} \\
(\Q\p_\s+V\p_2)\tt  & =   \bigl(\tfrac{1+ \alpha }{2} \Wbt + \tfrac{1- \alpha }{2} \Zbt \bigr) \nn   \,. \label{tt-evo-s}
\end{align} 
\end{subequations} 
The specific vorticity evolution is transformed to the equation
\begin{equation} 
\tfrac{\Jg}{\Sigma}  (\Q\p_\s+V\p_2)  \Upomega    -  \alpha \p_1 \Upomega
+ \alpha \Jg  g^{- {\frac{1}{2}} } \nbs_2h \   \nbs_2 \Upomega  =0  \,. \label{vort-s}
\end{equation} 
Here we have appealed to the abuse of notation mentioned in Remark~\ref{rem:no:tilde}.

\subsection{The $L^2$-based energy norms}
\label{sec:norms:L2:first}
In Sections~\ref{sec:formation:setup}--\ref{sec:sixth:order:energy},  we will make use of the ``energy'' and ``damping'' norms defined  as follows. We use the convention in Remark~\ref{rem:no:tilde}, dropping all tildes for  functions that depend on $(x,\s) \in \TT^2 \times [0,\eps)$. We keep the $\nbs$ notation from \eqref{nb-s} to emphasize this $(x,\s)$ dependence. See also Remark~\ref{rem:L2:norms:x:s:x:t:A} for equivalent norms in term of $(x,t)$ coordinates. 

The energy norms at the sixth derivative level are given by
\begin{subequations}
\label{eq:norms:L2:first}
\begin{align} 
\widetilde{\mathcal{E}}_{6}^2(\s) 
&=  \widetilde{\mathcal{E}}_{6,\nnn}^2(\s)  + (\mathsf{K}\eps)^{-2} \widetilde{\mathcal{E}}_{6,\ttt}^2(\s)
\label{eq:tilde:E6} \\
\widetilde{\mathcal{E}}_{6,\nnn}^2(\s) 
&= \snorm{  \mathcal{J}^{\!\frac 34} \Jgh \nbs^6 (\Jg\Wbn,\Jg\Zbn, \Jg\Abn)( \cdot , \s)}^2_{L^2_x} 
\\
\widetilde{\mathcal{E}}_{6,\ttt}^2(\s)
&=\snorm{ \mathcal{J}^{\!\frac 34} \Jgh  \nbs^6 (\Wbt, \Zbt, \Abt)( \cdot , \s)}^2_{L^2_x} 
\,,
\end{align}
and are defined  at the fifth derivative level by
\begin{align}
\widetilde{\mathcal{E}}_5^2(\s) 
&=  \widetilde{\mathcal{E}}_{5,\nnn}^2(\s) + (\mathsf{K}\eps)^{-2} \widetilde{\mathcal{E}}_{5,\ttt}^2(\s)
\label{eq:tilde:E5} \\
\widetilde{\mathcal{E}}_{5,\nnn}^2(\s)
&= \snorm{ \Jgh \nbs^5 (\Jg\Wbn, \Jg\Zbn, \Jg\Abn)( \cdot , \s)}^2_{L^2_x} 
\\
\widetilde{\mathcal{E}}_{5,\ttt}^2(\s)
&= \snorm{ \Jgh \nbs^5 ( \Wbt, \Zbt, \Abt)( \cdot , \s)}^2_{L^2_x} 
\,,
\end{align}
where  $\mathsf{K} = \mathsf{K}(\alpha) \geq 1$ is a sufficiently large constant chosen   solely in terms of  $\alpha$, see~\eqref{eq:K:choice:1}. In particular, $\mathsf{K}$ is independent of $\eps$. 
The sixth-order damping norms  are given by
\begin{align}
\widetilde{\mathcal{D}}_6^2(\s)
&= \widetilde{\mathcal{D}}_{6,\nnn}^2(\s) + (\mathsf{K} \eps)^{-2} \widetilde{\mathcal{D}}_{6,\ttt}^2(\s) 
\label{eq:tilde:D6}\\
\widetilde{\mathcal{D}}_{6,\nnn}^2(\s) 
&=   \int_0^\s 
\snorm{\mathcal{J}^{\!\frac 14}  \Jgh\nbs^6 (\Jg \Wbn, \Jg\Zbn, \Jg\Abn)( \cdot , \s')}_{L^2_x}^2 {\rm d}\s'
\\
\widetilde{\mathcal{D}}_{6,\ttt}^2(\s) 
&=  \int_0^\s 
\snorm{\mathcal{J}^{\!\frac 14} \Jgh\nbs^6 (\Wbt, \Zbt, \Abt)( \cdot , \s')}_{L^2_x}^2 {\rm d}\s'  
\,,
\end{align}
and the fifth-order damping norms  are
\begin{align}
\widetilde{\mathcal{D}}^2_5(\s)
& = \widetilde{\mathcal{D}}^2_{5,\nnn}(\s) + (\mathsf{K}\eps)^{-2} \widetilde{\mathcal{D}}^2_{5,\ttt}(\s)
\label{eq:tilde:D5} \\
\widetilde{\mathcal{D}}^2_{5,\nnn}(\s)
& =   \int_0^\s   \snorm{ \nbs^5 (\Jg \Wbn, \Jg\Zbn, \Jg\Abn)( \cdot , \s')}^2_{L^2_x} {\rm d}\s'
\\
\widetilde{\mathcal{D}}^2_{5,\ttt}(\s)
&= \int_0^\s    \snorm{ \nbs^5 (\Wbt, \Zbt, \Abt)( \cdot , \s')}^2_{L^2_x}   {\rm d}\s'  
\,,
\end{align}
\end{subequations}
where $\mathsf{K}\geq 1$ is the same constant in \eqref{eq:tilde:E6}, \eqref{eq:tilde:E5}, \eqref{eq:tilde:D6}, and in \eqref{eq:tilde:D5}.

\subsection{Bootstrap assumptions}
The existence of solutions in the spacetime $\mathcal{P}$ (equivalently, on $\TT^2 \times [0,\eps)$ in $(x,\s)$ variables) relies on quantitative bounds on all unknowns in the problem. We establish these quantitative bounds via a series of ``bootstrap inequalities''. Assuming these inequalities hold true with a specific constant on $\mathcal{P}$ (intuitively, this constant is related to the size of various norms of the initial data  multiplied by a constant which only depends on $\alpha$ and $\kappa_0$), we use the equations to prove that the bounds in fact hold true on $\mathcal{P}$ (equivalently, on $\TT^2 \times [0,\eps)$ in $(x,\s)$ variables)  with a constant which is  strictly smaller than what was assumed. A standard continuity argument is then used to justify that the bootstrap inequalities indeed hold true globally on $\mathcal{P}$.

The bootstrap assumptions are as follows. 
There exists a constant $\Csupp>0$, which depends only on $\alpha$ and $\kappa_0$ (see~\eqref{eq:Csupp:def} below), such that 
for all $\s\in [0,\eps)$, we have
\begin{subequations} 
\label{bootstraps}
\begin{equation}
\supp(\Wb,\Zb,\Ab,\nbs \Jg,\nbs \nbs_2 h)(\cdot,\s)
\subset 
\mathcal{X}_{\rm fin} 
:= 
\bigl\{ x\in \mathbb{T}^2 \colon {\rm dist}(x, \mathcal{X}_{\rm in}) \leq \Csupp \eps \bigr\}
\,.
\label{bs-supp}
\end{equation}
Regarding $\Wb$, we assume that pointwise for $(x,t) \in \mathcal{P}$, or equivalently, $(x,\s) \in \mathbb{T}^2 \times [0,\eps)$, we have
\begin{align} 
\Jg  \Wbn  & \ge   - \tfrac{9}{10} \eps^{-1} \ \ \text{ implies that } \ \ \Jg \ge \tfrac{2}{25} \,, 
\label{bs-Jg}  \\
|\Jg  \Wbn| &\le (1 +\eps) \eps^{-1}  \,, 
\label{bs-JgnnWb}\\
|\nbs (\Jg\Wbn)|  &\leq 3 \eps^{-1} \,, 
\label{bs-nnWb,2}\\
|\Wbt|  & \leq  1+\eps  \,,
\label{bs-ttWb} \\
|\nbs \Wbt|  & \leq 2 \Cdata  \,.
\label{bs-D-ttWb}
\end{align}
Regarding $\Zb$ and $\Ab$, we assume that pointwise for $(x,t) \in \mathcal{P}$, or equivalently, $(x,\s) \in \mathbb{T}^2 \times [0,\eps)$, it holds that
\begin{align}
|\Zbn| + |\nbs \Zbn| &\leq  \mathsf{C_{\Zbn}}    \,, 
\label{bs-nnZb}\\
|\Abn| + |\nbs \Abn| &\leq  \mathsf{C_{\Abn}}  \,, 
\label{bs-nnAb}\\
|\Zbt| + |\nbs \Zbt| &\leq  \mathsf{C_{\Zbt}}   \eps\,, 
\label{bs-ttZb}\\
|\Abt| + |\nbs \Abt| &\leq  \mathsf{C_{\Abt}}    \eps\,,
\label{bs-ttAb}
\end{align}
where $\mathsf{C_{\Zbn}}, \mathsf{C_{\Abn}}, \mathsf{C_{\Zbt}}$, and $\mathsf{C_{\Abt}}$ are sufficiently large constants which depend  only on $\alpha, \kappa_0$, and $\Cdata$, and which are fixed throughout the proof (see conditions~\eqref{eq:CZ:cond:1}, \eqref{eq:CZ:cond:2}, \eqref{eq:CZ:cond:3}, \eqref{eq:CZ:cond:4}, \eqref{eq:CZ:cond:5}, \eqref{eq:CZ:cond:6}, \eqref{eq:CZ:cond:7}, and \eqref{eq:CZ:cond:8} below).
Similar pointwise in spacetime bootstraps are assumed on the geometry, ALE-drift, and sound-speed:
\begin{align}
0 &\leq \Jg \leq \tfrac 65
\label{bs-Jg-simple} \\
|\nb  \Jg| &\leq  4(1+\alpha) 
\label{bs-Jg,1} \\
\max\bigr\{ \tfrac{1}{2} |\nb_1 h|  ,  \tfrac{1}{3} |\nb_2 h| , \tfrac{1}{(1+\alpha)\kappa_0 } | \nb_t h|\bigl\}&\le    \eps \,, 
\label{bs-h}\\
\max\{\tfrac{1}{5 (1+\alpha)}\sabs{\nb \nb_1 h} , \tfrac{1}{5 \Cdata}\sabs{\nb \nb_2 h} \}&\leq   \eps \,, 
\label{bs-h,22}\\
|V| + |\nb  V| &\leq  \mathsf{C_V} \eps \,,  
\label{bs-V}\\
\tfrac{\kappa_0}{4} &\leq \Sigma \leq \kappa_0 \,, 
\label{bs-Sigma}
\\
|\nb \Sigma| &\leq 2 \kappa_0 \,,
\label{bs-D-Sigma}
\end{align} 
where $\mathsf{C_V}$ is a sufficiently large constant which depends only on $\alpha,\kappa_0$, and $\Cdata$, which which is fixed throughout the proof  (see~\eqref{eq:C_V:choice}).
Lastly, for the energy bootstrap we assume that there exist constants $\mathsf{B_6}, \mathsf{B_5},  \mathsf{B_J}, \mathsf{B_h} \geq 1$, which only depend on $\alpha$, $\kappa_0$, and $\Cdata$, such that  
\begin{align}
\eps^{\frac 12} \!\!\sup_{\s \in [0,\eps]} \widetilde{\mathcal{E}}_{6} (\s) + \tilde{\mathcal{D}}_{6}(\eps) &\leq \mathsf{B}_6 \,, \label{bootstraps-Dnorm:6} \\
\eps^{\frac 12} \!\!\sup_{\s \in [0,\eps]} \widetilde{\mathcal{E}}_{5} (\s) + \tilde{\mathcal{D}}_{5}(\eps) &\leq \mathsf{B}_5 \,,\label{bootstraps-Dnorm:5}  \\
 \snorm{\nbs^6 \nbs_1 h}_{L^2_{x,\s}([0,\eps)\times \TT^2)} 
 +
 \snorm{\nbs^6 \nbs_2 h}_{L^2_{x,\s}([0,\eps)\times \TT^2)} &\leq \mathsf{B_h}  \eps^2\,, \label{bootstraps-Dnorm:h2} \\
\snorm{\nbs^6 \Jg}_{L^2_{x,\s}([0,\eps)\times \TT^2)} &\leq \mathsf{B_J} \eps.\label{bootstraps-Dnorm:Jg} 
\end{align}
\end{subequations} 
Without loss of generality, we will henceforth assume the ordering $ \mathsf{B}_6 \leq \mathsf{B}_5 \leq \mathsf{B_J} ,\mathsf{B_h}$.

\begin{remark}[\bf Norms with respect to $(x,\s)$ versus $(x,t)$ variables]
\label{rem:L2:norms:x:s:x:t:A}
Using definition~\eqref{eq:f:tilde:f}, in the bootstrap bounds \eqref{bootstraps} we have identified functions $F = F(x,t) \colon \mathcal{P} \to \mathbb{R}$ and their counterparts $\tilde F = \tilde F(x,\s) \colon \TT^2 \times [0,\eps) \to \mathbb{R}$. Additionally, according to Remark~\ref{rem:no:tilde} we have dropped tildes, writing $F$ instead of $\tilde F$, but have kept $\nbs$ instead of $\nb$ to emphasize $(x,\s)$ versus $(x,t)$ dependence. The perspective taken in our proof is that some of the bootstrap assumptions (e.g.~for $\Wbn$ and $\Jg$) are more convenient to close in $(x,t)$ variables, while some others (e.g.~the energy bounds) are more convenient to close in $(x,\s)$ variables. 
It is important however to emphasize that with the exception of the bootstraps for $\sup_{\s} \tilde{\mathcal{E}}_5(\s)$ and $\sup_{\s} \tilde{\mathcal{E}}_6(\s)$, no ambiguity arises from using $(x,\s)$ versus $(x,t)$ variables. Indeed, it is clear that at the level of pointwise bounds,~\eqref{eq:f:tilde:f} and~\eqref{nb-to-nbs} imply that for any function $F(x,t) = \tilde F(x,\s)$ and $k\geq 0$, we have
\begin{subequations}
\begin{equation}
\|\nb^k F\|_{L^\infty_{x,t}(\mathcal{P})} =  \|\nbs^k \tilde F\|_{L^\infty_{x,\s}([0,\eps)\times \TT^2)}\,.
\label{eq:norm:equivalence:a}
\end{equation}
This addresses~\eqref{bs-supp}--\eqref{bs-D-Sigma}. Next, we note that the Jacobian of the map $(x,t) \mapsto (x,\s)$ present in~\eqref{eq:f:tilde:f} is easily seen to equal $|\p_t \mathfrak{q}| = \Qd$, and the bound~\eqref{Qd-lower-upper} below gives global upper and lower bounds for $\Qd$ (which are strictly positive, and depend only on $\alpha$). As such, with the spacetime $\mathcal{P}$ defined in~\eqref{eq:spacetime:smooth}, the change of variable formula gives that for any function $F(x,t) = \tilde F(x,\s)$, any $k\geq 0$, and any weight $\varphi(x,t) = \tilde \varphi(x,\s) \geq 0$, we have
\begin{equation}
C_\alpha^{-1} \|\tilde \varphi \nbs^k \tilde F\|_{L^2_{x,\s}([0,\eps)\times \TT^2)} 
\leq \|\varphi \nb^k F\|_{L^2_{x,t}(\mathcal{P})} 
\leq C_\alpha  \|\tilde \varphi \nbs^k \tilde F\|_{L^2_{x,\s}([0,\eps)\times \TT^2)}\,,
\label{eq:norm:equivalence:b}
\end{equation}
for a constant $C_\alpha\geq 1$ that only depends on $\alpha$.
This addresses the $L^2_{x,\s} = L^2([0,\eps)\times\TT^2)$ norms present in \eqref{bootstraps-Dnorm:6}--\eqref{bootstraps-Dnorm:Jg}. It remains to discuss the $L^\infty_\s L^2_x = L^\infty([0,\eps);L^2(\TT^2))$ norms norms encoded by the bootstraps for $\sup_{\s} \tilde{\mathcal{E}}_5(\s)$ and $\sup_{\s} \tilde{\mathcal{E}}_6(\s)$ in \eqref{bootstraps-Dnorm:6}--\eqref{bootstraps-Dnorm:5}. Here, bounds which correspond to the ``time-slice foliation'' of $\TT^2 \times [0,\eps)$, namely $(\TT^2 \times \{\s\})_{\s\in[0,\eps]}$, do not translate to bounds on a ``time-slice foliation'' of $\mathcal{P}$. Instead, a foliation via level sets of $\mathcal{J}$, namely $\{(x_1,x_2,t) \colon \mathcal{J}(x_2,t) = 1 -\tfrac{\s}{\eps} \} = \{(x_1,x_2,\mathfrak{q}^{-1}(x_2,\s))\}$ for $\s\in[0,\eps]$, must be used. That is, for any function $F(x,t) = \tilde F(x,\s)$, any $k\geq 0$, and any weight $\varphi(x,t) = \tilde \varphi(x,\s) \geq 0$, we have
\begin{equation}
\|(\tilde \varphi \nbs^k \tilde F) (x,\s)\|_{L^2_{x}(\TT^2)} 
= \|(\varphi \nb^k F)(x,\mathfrak{q}^{-1}(x_2,\s))\|_{L^2_{x}(\TT^2)} 
\,,
\label{eq:norm:equivalence:c}
\end{equation}
for any $\s \in [0,\eps)$.
\end{subequations}
The equivalences in~\eqref{eq:norm:equivalence:a} and~\eqref{eq:norm:equivalence:b} will be used throughout the paper. In contrast,~\eqref{eq:norm:equivalence:c} is never used.
\end{remark}

\begin{remark}[\bf The order in which the bootstrap constants are chosen]
\label{rem:order:of:constants}
We will show that the bootstrap assumptions~\eqref{bootstraps} close if the various constants appearing therein, namely 
\begin{equation}
\mathsf{C_{supp}} \,, 
\mathsf{C}_{\Zbn} \,,
\mathsf{C}_{\Abn} \,,
\mathsf{C}_{\Zbt} \,,
\mathsf{C}_{\Abt} \,,
\mathsf{C_V} \,,
\mathsf{K}\,,
\mathsf{B}_6 \,,
\mathsf{B_5} \,,
\mathsf{B_J} \,, 
\mathsf{B_h}\,,
\label{eq:the:bootsrap:constants}
\end{equation}
are chosen suitably. In order to make sure that there is no circular argument, we discuss the precise interdependence of these constants. We first note that:
\begin{enumerate}[leftmargin=26pt]
\item\label{item:alpha:choice} The Euler equation fixes the parameter $\alpha = \frac{\gamma-1}{2} >0$, where $\gamma>1$ is the adiabatic exponent.
\item\label{item:Cdata:choice} The initial data, through assumptions~\eqref{item:ic:supp}--\eqref{item:ic:reg}, fixes two parameters $\kappa_0 \geq 1$ and $\bar{\mathsf{C}} \geq 1$. Moreover, as explained in Remark~\ref{rem:table:derivatives}, the parameters $\alpha, \kappa_0$, and $\bar{\mathsf{C}}$ determine a sufficiently large parameter $\Cdata \geq \kappa_0$. 
\end{enumerate}
It is important to emphasize that $\alpha, \kappa_0$, and $\Cdata$ are independent of $\eps>0$, which will be chosen to be sufficiently small at the end of the proof. Next, the precise order in which the bootstrap constants are chosen is as follows:
\begin{enumerate}[leftmargin=26pt]
\setcounter{enumi}{2} 

\item \label{item:C:supp:choice} $\mathsf{C_{supp}}$ is chosen in~\eqref{eq:Csupp:def} to depend only on $\alpha$ and $\kappa_0$. Subsequent dependence on $\mathsf{C_{supp}}$ is encoded as dependence on $\alpha$ and $\kappa_0$.

\item \label{item:K:sf:choice} $\mathsf{K} \geq 1$ is chosen in~\eqref{eq:K:choice:1} to depend only on $\alpha$. For the downstream \MGHDB\ given in Section~\ref{sec:downstreammaxdev}, $\mathsf{K}$ also needs to depend on $\kappa_0$, cf.~\eqref{eq:K:choice:1-P}.

\item \label{item:B6:choice} $\mathsf{B_6}\geq 1$ is determined by~\eqref{eq:B6:choice:1} and~\eqref{eq:B6:choice:2} to be sufficiently large with respect to only $\alpha$, and $\Cdata$. For the downstream \MGHDB\ in Section~\ref{sec:downstreammaxdev}, $\mathsf{B}_6$ also needs to depend on $\kappa_0$, cf.~\eqref{eq:B6:choice:1-P} and~\eqref{eq:B6:choice:2-P}.

\item $\mathsf{B}_5$ is chosen in~\eqref{eq:B5B6:relation} to depend only on  $\alpha$, $\Cdata$, and $\mathsf{B}_6$. As shown in~\eqref{rem:B5:B6}, the quotient $\mathsf{B}_5 \mathsf{B}_6^{-1}$ is bounded from above by a universal constant, and from below by a constant that only depends on $\alpha$. As such, subsequent bounds of the type $A \leq \mathsf{B}_5$ will be written  as $A \leq \Cn \mathsf{B}_6$ or $A\les \mathsf{B}_6$ (see~Remark~\ref{rem:les}).

\item $\mathsf{B_J}$ is chosen in~\eqref{eq:B:Jg:cond} to depend only on $\alpha,\kappa_0, \Cdata$, and in a linear fashion on $\brak{\mathsf{B_6}}$. Since $\mathsf{B}_6\geq 1$, subsequent bounds of the type $A \leq \mathsf{B_J}$ will be written as $A \leq \Cn \mathsf{B}_6$ or $A \les \mathsf{B}_6$ (see~Remark~\ref{rem:les}).

\item $\mathsf{B_h}$ is chosen in~\eqref{eq:B:h2:cond} to depend only on $\alpha,\kappa_0, \Cdata$, and in a linear fashion on $\brak{\mathsf{K}} \brak{\mathsf{B_6}}$. Since $\mathsf{B}_6\geq 1$ and $\mathsf{K}\geq 1$, subsequent bounds of the type $A \leq \mathsf{B_h}$ will be written as $A \leq \Cn \mathsf{K} \mathsf{B}_6$ or $A \les \mathsf{K}  \mathsf{B}_6$ (see~Remark~\ref{rem:les}).

\item  \label{item:C:Abn:choice}  $\mathsf{C}_{\Abn}$ is determined by~\eqref{eq:CZ:cond:1} and \eqref{eq:CZ:cond:2}, and depends only on $\alpha$ and $\Cdata$. 

\item  \label{item:C:Abt:choice} $\mathsf{C}_{\Abt}$ is determined by~\eqref{eq:CZ:cond:3} and \eqref{eq:CZ:cond:4}, and depends on $\alpha, \kappa_0,\Cdata, \mathsf{C_{\Abn}}$, and $\mathsf{B_6}$. In view of points~\eqref{item:B6:choice} and~\eqref{item:C:Abn:choice} above, dependence on $\mathsf{C}_{\Abt}$ is subsequently encoded as dependence only on $\alpha, \kappa_0$, and $\Cdata$.

\item $\mathsf{C_{V}}$ is determined by~\eqref{eq:C_V:choice} and depends only on $\alpha,\kappa_0,\Cdata,\mathsf{C}_{\Abn}$, and $\mathsf{C}_{\Abt}$.  In view of points~\eqref{item:C:Abn:choice} and~\eqref{item:C:Abt:choice} above, dependence on $\mathsf{C_V}$ is subsequently encoded as dependence only on $\alpha, \kappa_0$, and $\Cdata$.

\item \label{item:C:Zbn:choice} $\mathsf{C}_{\Zbn}$ is determined by~\eqref{eq:CZ:cond:5} and \eqref{eq:CZ:cond:6} and depends only on $\alpha,\kappa_0$ and $\Cdata$.

\item \label{item:C:Zbt:choice} $\mathsf{C}_{\Zbt}$ is determined by~\eqref{eq:CZ:cond:7} and \eqref{eq:CZ:cond:8} and depends only on $\alpha,\kappa_0, \Cdata$, and $\mathsf{C}_{\Zbn}$. In view of point~\eqref{item:C:Zbn:choice} above, dependence on $\mathsf{C}_{\Zbt}$ is subsequently encoded as dependence only on $\alpha, \kappa_0$, and $\Cdata$. 
\end{enumerate}
The last parameter chosen in the proof is:
\begin{enumerate}[leftmargin=26pt]
\setcounter{enumi}{13}
\item  \label{item:eps:final:choice} $\eps >0$, which is taken to be sufficiently small with respect to $\alpha$, $\kappa_0$, and $\Cdata$. That is, $\eps$ is taken to be small enough with respect to the parameters induced by the initial data, cf.~points~\eqref{item:alpha:choice}--\eqref{item:Cdata:choice} above.
\end{enumerate}
We emphasize that in view of points~\eqref{item:C:supp:choice}--\eqref{item:C:Zbt:choice} above, $\eps$ is  sufficiently small enough with respect to any of the constants appearing in the bootstrap assumptions, cf.~\eqref{eq:the:bootsrap:constants}. This fact is used implicitly throughout the paper.
\end{remark}

\section{First consequences of the bootstrap assumptions}
\label{sec:first:consequences}
In this section we collect a few direct consequences of the bootstrap bounds \eqref{bootstraps}, which are then subsequently used throughout the paper. We emphasize that the order in which these consequences are proven is irrelevant, they are all consequences of \eqref{bootstraps}.  As such, we sometimes make forward references to other bounds which are direct consequences of the bootstrap assumptions.

\subsection{Spatial support}
\label{sec:spatial:support}
The goal of this subsection is to prove the bootstrap \eqref{bs-supp}. Recall  that at the initial time $t=\initial$ we have that $(\Wb,\Zb,\Ab)$ and hence $(\Ub,\Sb)$ are compactly supported in $[-13\pi\eps,13\pi\eps] \times \mathbb{T}$; see the set~$\mathcal{X}_{\rm in}$ defined in \eqref{eq:ic:supp}. Then, for a speed $\mathsf{v}$ which is to be determined later in the proof, one may define an expanding set 
\begin{equation*}
\mathcal{X}(\s) := \{ x\in\mathbb{T}^2 \colon |x_1| \leq 26\pi \eps + \mathsf{v} \s  \}.
\end{equation*}
Note that at time $\s=0$ we have that $\mathcal{X}_{\rm in} \subset \mathcal{X}(0)$, giving us a bit of room to operate, at least for some infinitesimally small time. Then, we may use the system \eqref{euler-ALE}, \eqref{Jg-evo}, and \eqref{Sigma0-ALE} to show that there exists a sufficiently large parameter $\mathsf{v}$, depending only on $\alpha$, $\kappa$, such that 
\begin{equation}
\int_{(\mathcal{X}(\s))^\complement}
\tfrac{\Q \Jg}{2\Sigma} \bigl(|\Ub|^2  + |\Sb|^2\bigr) {\rm d}x
= 0 \,,
\label{eq:bs:supp}
\end{equation}
for all $\s \in [0,\eps)$, where we have denoted $|\Ub|^2 = \Uik \Uik$ and $|\Sb|^2 = \Sk \Sk$, with the usual convention of summation over repeated indices. As $\frac{\Q \Jg}{2\Sigma} > 0$ on $\TT^2 \times [0,\eps)$ (cf.~\eqref{eq:spacetime:smooth}, \eqref{Qd-lower-upper}, \eqref{Q-lower-upper}, \eqref{eq:Sigma:sharp})  if we establish the above identity, it means that the solution $(\Ub,\Sb)$ is compactly supported in $\mathcal{X}(\s)$, as  claimed. For the sake of a contradiction, assume that there is a minimal $\bar \s \in (0,\eps)$ such that \eqref{eq:bs:supp} holds on $[0,\bar \s]$, but that   $(\frac{d}{d\s} \int_{(\mathcal{X}(\s))^\complement}
\tfrac{\Q \Jg}{2\Sigma} \bigl(|\Ub|^2  + |\Sb|^2\bigr) {\rm d}x)|_{\s=\bar \s} > 0$. Then, from the chain rule we obtain 
\begin{equation}
\Bigl(\tfrac{d}{d \s} \int_{(\mathcal{X}(\s))^\complement}
\tfrac{\Q\Jg}{2 \Sigma} \bigl(|\Ub|^2  + |\Sb|^2\bigr) {\rm d}x\Bigr) \Bigr|_{\s=\bar \s}
\leq \int_{(\mathcal{X}(\bar \s))^\complement} 
\p_\s \Bigl(\tfrac{\Q\Jg}{2 \Sigma} \bigl(|\Ub|^2  + |\Sb|^2\bigr)\Bigr) {\rm d}x
- \tfrac{\mathsf{v}}{2}
\int_{\partial \mathcal{X}(\bar \s)}
 \tfrac{\Q\Jg}{2 \Sigma} \bigl(|\Ub|^2  + |\Sb|^2\bigr) {\rm d}x
 \,.
 \label{eq:bs:supp:2}
\end{equation}
Next, from \eqref{euler-ALE}, \eqref{Jg-evo}, \eqref{Sigma0-ALE}, \eqref{eq:xt:xs:chain:rule} and \eqref{the-time-der}, we obtain
\begin{align} 
\p_\s  \bigl(\tfrac{\Q \Jg}{2 \Sigma} ( |\Ub|^2 + |\Sb|^2) \bigr)
&=
 \tfrac{1}{\Sigma}( |\Ub|^2 + |\Sb|^2) \bigl( \tfrac{1+\alpha}{2}  \Jg\Wbn + \tfrac{1- \alpha}{2} \Jg\Zbn + \tfrac{\alpha}{2} \Jg \Abt + \tfrac{\Qc + V,_2 }{2} \Jg\bigr)
\notag\\
&\quad
-  \tfrac{\Jg}{\Sigma} \Bigl( (1+\alpha) \Sk \Si  \Uik + \alpha |\Sb|^2 \Uii + \Uij\Ujk \Uik  \Bigr)   
+ \tfrac{\alpha}{2} \Bigl(  |\Ub|^2 + |\Sb|^2 - 2 \nn^i \Sk  \Uik \Bigr),_1 
\notag\\
&\quad
- \Jg   \Bigl(   \tfrac{V  }{2 \Sigma} ( |\Ub|^2 + |\Sb|^2)  
- \alpha    \Sk e_2^i \Uik 
- \tfrac{\alpha}{2}    g^{-\frac 12} h,_2 (|\Ub|^2  + |\Sb|^2)  \Bigr),_2
\notag\\
&\quad  
-   \Bigl( \alpha  \Uik \Sk \tt^i  g^{-\frac 12}   
+  \alpha   \Sk e_2^i \Uik   
+ \tfrac{\alpha}{2}  (|\Ub|^2  + |\Sb|^2)  g^{-1} h,_2 
+ \tfrac{V}{2 \Sigma} (|\Ub|^2  + |\Sb|^2) \Bigr) 
\Jg,_2  
\notag\\
&\quad
- \Jg  \Bigl( \alpha \tt^i g^{-1} h,_2  \Uik \Sk   - \tfrac{\alpha}{2} g^{-\frac 32}  (|\Ub|^2  + |\Sb|^2)   \Bigr)   h,_{22} 
\,. 
\label{eq:supp:energy}
\end{align} 
A remarks is in order. Since for all $s\in [0,\bar \s]$ we have that $(\Ub,\Sb)$, and hence also $(\Wb,\Zb,\Ab)$, are supported in $\mathcal{X}(\s)$, from identity  \eqref{Jg-circ-xi} below, and  the bounds \eqref{eq:xi:nabla:xi}, \eqref{bs-V}, and \eqref{eq:Q:all:bbq}, we obtain that  if  $\mathsf{v} \geq C \eps$, then $\Jg(x,\s)  =1$ and $h,_2(x,\s) = 0$ for all $x \in (\mathcal{X}(\bar \s))^\complement$. Thus, in \eqref{eq:supp:energy} we have that the $\Jg,_2$ and $h,_{22}$ terms vanish identically, at all points in the interior of $(\mathcal{X}(\bar \s))^\complement$. Moreover, we may freely multiply (or divide) the right side of \eqref{eq:supp:energy} by powers of $\Jg$, since we are then multiplying by $1$.  Using  these  observations, we integrate \eqref{eq:supp:energy} over $(\mathcal{X}(\s))^\complement$, integrate by parts the pure derivative terms, and appeal to the  bootstrap bounds~\eqref{bootstraps} to conclude
\begin{align}
&\int_{(\mathcal{X}(\bar \s))^\complement} 
\p_\s \Bigl(\tfrac{\Q \Jg}{2 \Sigma} \bigl(|\Ub|^2  + |\Sb|^2\bigr)\Bigr) {\rm d}x
\notag\\
&\qquad \leq  
\tfrac{\Cn}{\eps} 
\int_{(\mathcal{X}(\bar  \s))^\complement} 
\tfrac{\Q \Jg}{2 \Sigma} \bigl(|\Ub|^2  + |\Sb|^2\bigr) {\rm d}x
+ \left( 32 \alpha (1+\alpha)\kappa_0  + \Cn \eps\right)
\int_{\partial \mathcal{X}(\bar \s)} 
\tfrac{\Q \Jg}{2 \Sigma} \bigl(|\Ub|^2  + |\Sb|^2\bigr) {\rm d}x
\,.
\label{eq:bs:supp:3}
\end{align}
Note that by assumption on the minimality of $\bar s$, the first term on the right side of \eqref{eq:bs:supp:3},  vanishes. As such, if $\mathsf{v}$ is taken to as
$ \mathsf{v} := 65 \alpha (1+\alpha)\kappa_0 $, 
and $\eps$ is chosen sufficiently small, \eqref{eq:bs:supp:2} and \eqref{eq:bs:supp:3}
yield    $ \frac{d}{ds} \int_{(\mathcal{X}(\s))^\complement}
\tfrac{\Q \Jg}{2\Sigma} \bigl(|\Ub|^2  + |\Sb|^2\bigr) {\rm d}x|_{\s=\bar \s} \leq 0$, a contradiction.

We have thus shown that by letting $\mathsf{v} = 65 \alpha (1+\alpha)\kappa_0$, the functions $(\Ub,\Sb,\nbs \Jg,\nbs \nbs,_2 h)(\cdot,\s)$ are compactly supported in the set $\mathcal{X}(\eps)$, for all $\s \in [0,\eps]$. Choosing the constant $\mathsf{C_{supp}}$ in \eqref{eq:supp:fin} to equal $13 \pi + \mathsf{v}$, i.e.,
\begin{equation}
\mathsf{C_{supp}} = 13 \pi + 65 \alpha (1+\alpha)\kappa_0\,,
\label{eq:Csupp:def}
\end{equation}
shows that the bootstrap~\eqref{bs-supp} is closed.

\subsection{Flow of the ALE velocity $V$}
\label{sec:ALE:flow}
Due to the presence of the ALE transport operator $\p_t + V \p_2$, when working in $(x,t)$ coordinates it convenient to consider the flow map:
\begin{subequations} 
\label{xi-flow}
\begin{align} 
\p_t\xi (x_1,x_2,t) & = V(x_1, \xi(x_1,x_2,t),t)\, \qquad \mbox{for} \qquad t \in (\initial,\final) \,, \\
\xi(x_1,x_2,\initial) &= x_2 \,,
\end{align} 
\end{subequations} 
where we recall that $V$ is defined in \eqref{transport-ale}. Given a label $x$, the flow $\xi(x,\cdot)$ is then defined up to the stopping time
\begin{equation*} 
\mathsf{T}_{\xi}(x) = \sup \{ t \in [\initial,\final] \colon (\xi(x,t),t) \in \mathcal{P} \}\,.
\end{equation*} 
Then, for any function $F \colon \mathcal{P} \to \mathbb{R}$, by the composition $F \cir \xi$ we mean
\begin{equation*} 
(F \cir \xi)(x,t) = F(x_1, \xi(x_1,x_2,t),t) \,,
\end{equation*} 
which is well-defined for $t \leq  \mathsf{T}_{\xi}(x)$. Similarly, the composition with the inverse of $\xi$, which is denoted as usual by $\xi^{-1}$ only affects the second space coordinate.

Next, we note that the bootstrap assumptions~\eqref{bootstraps} and the identities 
\begin{subequations}
\label{eq:nabla:V}
\begin{align}
\p_t (\log \xi,_2) 
&= V,_2 \circ \xi \\
\xi,_2 \p_t (\tfrac{\xi,_1}{\xi,_2}) 
&=  V,_1 \circ \xi \\
V,_1 &=
- \alpha \Sigma g^{-\frac 32} h,_{12}
+ g^{-\frac 12} \bigl(\Jg\Abn - h,_2 (\tfrac{1+\alpha}{2} \Jg \Wbn + \tfrac{1-\alpha}{2} \Jg \Zbn ) \bigr)
\notag\\
&\qquad \qquad \qquad\qquad
+ g^{-\frac 12}  h,_2 \Jg \bigl( \Abt  - h,_2 (\tfrac{1+\alpha}{2}  \Wbt + \tfrac{1-\alpha}{2} \Zbt ) \bigr)
\,,
\label{eq:nabla:V:1}
\\
V,_2 &=
- \alpha \Sigma g^{-\frac 32} h,_{22} 
+ \Abt
 - h,_2 \bigl( \tfrac{1+\alpha}{2} \Wbt + \tfrac{1-\alpha}{2} \Zbt \bigr) \,,
\label{eq:nabla:V:2}
\end{align}
\end{subequations}
imply that pointwise in $(x,t) \in \mathcal{P}$  we have
\begin{equation}
\abs{\xi(x,t) - x_2} \les \eps^2
\,, \qquad
 \abs{\xi,_2 - 1} \les \eps^{2}
 \,, \qquad 
 \abs{\p_t \xi} \les \eps
  \,, \qquad 
 \abs{\xi,_1 } \les \eps
 \,.
 \label{eq:xi:nabla:xi}
\end{equation}

Pointwise estimates for $\nabla^2\xi$ may also be obtained upon by noting that \eqref{xi-flow},  \eqref{eq:xi:nabla:xi}, \eqref{eq:Sigma:H6:new}, and \eqref{eq:Sobolev}, imply
\begin{subequations}
\label{eq:xi:hessian:xi}
\begin{align}
\abs{\xi,_{11}} 
 &\les
 \eps^{-1} \snorm{\nbs_1^2 V}_{L^\infty_{x,\s}}
 + \eps  \snorm{\nbs_1 \nbs_2 V }_{L^\infty_{x,\s}} 
 + \eps^{3} \snorm{\nbs_2^2 V }_{L^\infty_{x,\s}} 
 \notag\\
 & \les
 \eps^{-\frac 32} \snorm{\nbs_1 \nbs^3 V(\cdot,0)}_{L^2_x} 
 + \eps^{-2} \snorm{\nbs^4 \nbs_1 V}_{L^2_{x,\s}}
 \les 
 \mathsf{K} \brak{\mathsf{B}_6}
 \\
 \abs{\xi,_{12}} 
 &\les  \snorm{\nbs_1 \nbs_2 V }_{L^\infty_{x,\s}} 
 + \eps^2 \snorm{\nbs_2^2 V }_{L^\infty_{x,\s}} 
  \les \eps  \mathsf{K} \brak{\mathsf{B}_6} 
 \,, \\ 
  \abs{\xi,_{22}} 
&\les \eps   \snorm{\nbs_2^2 V }_{L^\infty_{x,\s}} 
\les \eps^2  \mathsf{K} \brak{\mathsf{B}_6}
\,,
\end{align}
\end{subequations}
pointwise for $(x,t) \in \mathcal{P}$. 
The bounds \eqref{eq:xi:nabla:xi} and \eqref{eq:xi:hessian:xi} imply that for any smooth function $f$ we have
\begin{subequations}
\label{eq:xi:switch}
\begin{align}
\sabs{\nb(f\circ \xi) - (\nb f)\circ \xi} 
&\les \eps^2 |\nb f| \circ \xi 
\\
\sabs{\nb^2(f\circ \xi) - (\nb^2 f)\circ \xi} 
&\les  \eps^2 \mathsf{K} \brak{\mathsf{B}_6} |\nb_2 f|\circ \xi  + \eps^2 |\nb^2 f| \circ \xi 
\end{align} 
\end{subequations}
pointwise for $(x,t) \in \mathcal{P}$. 
For later use, we record the  equations \eqref{eq:nabla:V:1} and \eqref{eq:nabla:V:2} transformed into $(x, \s)$ coordinates:
\begin{subequations} 
\label{D1-and-D2-V-s}
\begin{align} 
\nbs_1 V &=
- \alpha \Sigma g^{-\frac 32}\nbs_1\nbs_2 h 
+ \eps g^{-\frac 12} \bigl(\Jg\Abn - \nbs_2 h (\tfrac{1+\alpha}{2} \Jg \Wbn + \tfrac{1-\alpha}{2} \Jg \Zbn ) \bigr)
\notag \\
&\qquad \qquad
+\eps g^{-\frac 12}  \nbs_2 h\  \Jg \bigl( \Abt  - \nbs_2 h (\tfrac{1+\alpha}{2}  \Wbt + \tfrac{1-\alpha}{2} \Zbt ) \bigr)
\,,  \label{eq:nabla:V:1-s}
\\
\nbs_2V &= - \alpha \Sigma g^{-\frac 32} \nbs_2^2h + \Abt
 - \nbs_2 h \bigl( \tfrac{1+\alpha}{2} \Wbt + \tfrac{1-\alpha}{2} \Zbt \bigr) \,.   \label{eq:nabla:V:2-s}
\end{align} 
\end{subequations} 
 
\subsection{Bounds for $(W,Z,A)$}
\label{sec:WZA:pointwise} 
The goal of this subsection is to establish the pointwise bounds 
\begin{equation}
|W(x,t) - \kappa_0| \leq \tfrac{1}{6} \kappa_0
\,, \qquad
|Z(x,t)| \leq \kappa_0 \eps \bigl(1 + \tfrac{5\alpha }{1+\alpha}  \mathsf{C}_{\Zbn} \bigr)
\,, \qquad 
|A(x,t)| \leq \kappa_0 \eps  \bigl(1 + 3\kappa_0 + \tfrac{3\alpha}{1+\alpha}\mathsf{C}_{\Abn} \bigr)
\,.
\label{eq:WZA:bounds:sup}  
\end{equation}
for all $(x,t) \in \mathcal{P}$.
 
From \eqref{U0-ALE}, \eqref{nn-evo}, and \eqref{WZA-evo:c} we deduce that 
\begin{subequations}
\label{U.n-and-A}
\begin{align}
(\p_t + V \p_2) (U\cdot \nn) - \alpha \Sigma  \Zbn + A \bigl(\tfrac{1+ \alpha }{2} \Wbt + \tfrac{1- \alpha }{2} \Zbt\bigr) &= 0 \,, 
\label{eq:bingo:de:bongo:forever:a}\\
(\p_t +V\p_2)A   + \tfrac{\alpha}{2} \Sigma (\Wbt - \Zbt - 2 \Abn) -  (U\cdot \nn) \bigl(\tfrac{1+\alpha}{2}\Wbt + \tfrac{1-\alpha}{2} \Zbt\bigr)  &=0 \,,
\label{eq:bingo:de:bongo:forever:b}
\end{align}
\end{subequations}
and hence, by also appealing to the bootstrap inequalities \eqref{bootstraps}, we obtain
\begin{align}
\tfrac 12 (\p_t + V \p_2) \bigl((U\cdot \nn)^2 + A^2\bigr)
&= \alpha \Sigma \bigl( (U\cdot\nn) \Zbn - \tfrac 12 A (\Wbt - \Zbt - 2 \Abn)\bigr)
\notag \\
&\leq \alpha \Sigma \bigl((U\cdot \nn)^2 + A^2\bigr)^{\frac 12} \bigl( |\Zbn| +  \tfrac 12 (|\Wbt| + |\Zbt| + 2 |\Abn|)\bigr)
\notag \\
&\leq \Cn \bigl((U\cdot \nn)^2 + A^2\bigr)^{\frac 12} \,.
\end{align}
Upon integrating the above estimate we further obtain 
\begin{equation}
\sabs{ \bigl((U\cdot \nn)^2 + A^2\bigr)^{\frac 12} \circ \xi (x,t) - |u_0(x)|}
\leq \Cn \eps
\label{eq:bingo:de:bongo:forever:aa}
\,.
\end{equation}
On the other hand, appealing to \eqref{item:ic:infinity} we have that 
\begin{equation}
|u_0(x)| = \bigl( \tfrac 14 (w_0(x) + z_0(x))^2 + a_0(x)^2\bigr)^{\frac 12} 
\in [\tfrac{11}{24} - \eps, \tfrac{13}{24} + \eps] \kappa_0
\,.
\label{eq:bingo:de:bongo:forever:bb}
\end{equation}
We deduce in particular the preliminary estimate $\|U\cdot \nn\|_{L^\infty_{x,t}} \leq\frac{13}{24} \kappa_0 + \Cn \eps \leq  \frac{7}{12}\kappa_0$.

The estimate for $A$ now follows by integrating \eqref{eq:bingo:de:bongo:forever:b}, and using the previously established bound for $U \cdot \nn$ along with the bootstrap inequalities~\eqref{bootstraps}, to obtain
\begin{align*}
\sabs{A\circ \xi(x,t)}
&\leq |a_0(x)| + 
\tfrac{3\eps}{1+\alpha} \Bigl( \tfrac{\alpha}{2} \kappa_0 \bigl(1 + \eps  + \eps \mathsf{C_{\Zbt}} + 2 \mathsf{C_{\Abn}}\bigr)   + \tfrac{7}{12} \kappa_0 \bigl( \tfrac{1+\alpha}{2} (1+\eps)   + \tfrac{1-\alpha}{2} \eps \mathsf{C}_{\Zbt}\bigr)  \Bigr)
\notag\\
&
\leq \kappa_0 \eps  + 
\tfrac{2\eps}{1+\alpha} \tfrac{51}{50} \Bigl(\tfrac {7+19\alpha }{24} \kappa_0 +  \alpha \kappa_0 \mathsf{C}_{\Abn}  +\Cn \eps \Bigr)
\leq \kappa_0 \eps  \Bigl(1 + 3 \kappa_0 + \tfrac{3\alpha}{1+\alpha}\mathsf{C}_{\Abn} \Bigr)
\,.
\end{align*}
The $A$ estimate in \eqref{eq:WZA:bounds:sup} follows by composing with $\xi^{-1}$.

The $Z$ estimate is obtained in a similar way by integrating \eqref{WZA-evo:b}, and appealing to the existing pointwise bounds
\begin{equation*}
\sabs{Z\circ \xi(x,t)}
\leq |z_0(x)| + 
\tfrac{2\eps}{1+\alpha} \tfrac{51}{50} \Bigl( \eps \Cn + 2 \alpha \kappa_0 \mathsf{C}_{\Zbn} \Bigr)
\leq \kappa_0 \eps \Bigl(1 + \tfrac{5\alpha }{1+\alpha}  \mathsf{C}_{\Zbn} \Bigr)
\,.
\end{equation*}
The $Z$ estimate in \eqref{eq:WZA:bounds:sup} now follows.

Lastly, the $W$ estimate is obtained by integrating~\eqref{WZA-evo:a}, together with assumption~\eqref{bs-Sigma}, and appealing to the existing pointwise bounds
\begin{equation*}
\sabs{W \circ \xi(x,t) - w_0(x)}
\leq \Cn \eps^2
\Rightarrow
\sabs{W \circ \xi(x,t) - \kappa_0}
\leq \sabs{W \circ \xi(x,t) - w_0(x)} + \sabs{w_0(x)-\kappa_0}
\leq \Cn \eps^2 + \tfrac{\kappa_0}{12}\,.
\end{equation*}
Taking $\eps$ to be sufficiently small and composing with $\xi^{-1}$ yields the $W$ estimate in \eqref{eq:WZA:bounds:sup}.

\subsection{Pointwise bounds for $\nb^k (\Jg \Wbn)$ and $\nb^k   \Jg  $ when $0 \leq k \leq 2$}
\label{sec:Wbn:Jg:pointwise}
For future use, it is convenient to record the following pointwise bounds for the first few derivatives of $\Jg \Wbn$ and $\Jg \Zbn$.
\begin{lemma}
\label{lem:JgWbn:close:to:data}
Assume that the bootstrap bounds \eqref{bootstraps} hold on $\TT^2 \times [0,\eps)$.
If $\eps$ is taken to be sufficiently small with respect to $\alpha, \kappa_0$, and $\Cdata$, then 
\begin{subequations}
\label{eq:broncos:eat:shit:2}
\begin{align}
\sabs{ (\Jg \Wbn)(x,t) -   (w_0),_1(x) }
&\les \eps  
\label{eq:broncos:eat:shit:20}\\
\sabs{\nb   (\Jg \Wbn)(x,t) - \nb  (w_0),_1(x) }
&\les \eps \mathsf{K} \brak{\mathsf{B_6}}
\label{eq:broncos:eat:shit:2a}\\
\sabs{\nb^2  (\Jg \Wbn)(x,t) - \nb^2  (w_0),_1(x) }
&\les \mathsf{K} \brak{\mathsf{B_6}}
\label{eq:broncos:eat:shit:2b}
\end{align}
holds for all $(x,t) \in \mathcal{P}$. 
\end{subequations}
\end{lemma}
\begin{proof}[Proof of Lemma~\ref{lem:JgWbn:close:to:data}]
From \eqref{eq:Jg:Wb:nn}, upon composing with the flow $\xi$ associated to the vector field $\p_t +V \p_2$, we deduce that for each frozen $x$, and for $t\leq \mathsf{T}_{\xi}(x)$, we have
\begin{subequations} 
\label{Jg-Wbn-family}
\begin{equation}
(\Jg  \Wbn) \cir \xi(x,t)
= (w_0),_{1}\!(x)\  \Iwn1(x,t) + \Iwn2(x,t)  \,,  \label{Jg-Wbn-1}
\end{equation}
where we have denoted
\begin{align}
\Iwn1(x,t)&=e^{- \frac{\alpha}{2} \int_{\initial}^t ( \Abt  - g^{-\frac 32} \Sigma h,_{22}) \circ \xi (x,r) {\rm d}r}  \,,  \label{Iwn1}\\
\Iwn2(x,t)&= \int_{\initial}^t \Ewn \cir \xi (x,r) e^{-\frac{\alpha}{2}\int_{r}^t ( \Abt - g^{-\frac 32} \Sigma h,_{22}) \circ \xi (x,r') {\rm d}r'} {\rm d}r \,, \label{Iwn2} \\
\Ewn(x,t)&= - \alpha \Sigma g^{-\frac 12} \Jg \Abn,_2  
+ \tfrac{\alpha}{2} \Sigma g^{-\frac 32} h,_{22} (\Jg \Zbn  -2\Jg \Abt) \notag \\
&\qquad 
+ \tfrac{\alpha}{2}   \Abt   \Jg \Zbn  -  \bigl(\tfrac{3+\alpha}{2} \Wbt + \tfrac{1-\alpha}{2} \Zbt\bigr) \Jg\Abn
- \bigl(\tfrac{1+\alpha}{2} \Wbt + \tfrac{1-\alpha}{2} \Zbt\bigr)\Jg \Wbt   
\label{Ewn}
\,.
\end{align}
\end{subequations}
Using the bootstrap assumptions~\eqref{bootstraps}, we  see that the terms defined in \eqref{Iwn1}--\eqref{Iwn2} satisfy the pointwise estimates
\begin{equation} 
\sabs{1-\Iwn1} \les \eps^{2} \,, \qquad  \sabs{\Iwn2} \les \eps   \,, \label{Iwn1-2-est}
\end{equation} 
from where we deduce that 
\begin{equation}
\sabs{(\Jg  \Wbn) \cir \xi(x,t) - (w_0)_{,1}(x)}
\les  \eps^2 |(w_0)_{,1}(x)| + \eps \les \eps \,.
\end{equation}
Additionally, the implicit function theorem combined with the bounds~\eqref{table:derivatives} and~ \eqref{eq:xi:nabla:xi} give
\begin{equation}
\label{eq:d1:w0:close}
|\nb^k (w_0),_1(x) - \nb^k (w_0),_1\circ \xi^{-1}(x,t)|
\leq \|\nb^k (w_0),_{12}\|_{L^\infty_x} |x_2 - \xi^{-1}(x,t)| 
\les \eps \,,
\end{equation}
for $k\in\{0,1,2\}$, 
from which \eqref{eq:broncos:eat:shit:20} follows.

Establishing  \eqref{eq:broncos:eat:shit:2a} and \eqref{eq:broncos:eat:shit:2b} requires bounds for the derivatives of $\Iwn1$ and $\Iwn2$, which we obtain  as follows.  Using \eqref{eq:xi:nabla:xi}, we have that $|\nb_1 (f\circ\xi)| \les   |\nb_1 f|\circ\xi + \eps^2 |\nb_2 f|\circ \xi \les   |\nb f|\circ \xi$, $|\nb_2 (f\circ\xi)|\les |\nb_2 f|\circ \xi$, and $|\nb_t (f\circ \xi)| \les |\nb_t f|\circ \xi+ \eps^2 |\nb_2 f|\circ \xi \les |\nb f|\circ \xi$. With this in hand, we return to the terms defined in~\eqref{Iwn1}--\eqref{Ewn} and use the bootstrap inequalities~\eqref{bootstraps}, the product and chain rules to estimate
\begin{align*}
\snorm{\nb \, \Iwn1}_{L^\infty_{x,t}}
&\les \eps^2 + \eps \snorm{\nb h,_{22}}_{L^\infty_{x,t}}
\notag\\
\snorm{\nb^2 \, \Iwn1}_{L^\infty_{x,t}}
&\les\eps^3 + \eps^{\frac 12} \snorm{\nb^2 \Abt}_{{L^2_t L^\infty_{x}}}  +  \eps \snorm{\nb^2 h,_{2}}_{L^\infty_{x,t}}  + \eps^2 \snorm{\nb^2 \Sigma}_{L^\infty_{x,t}} + \eps^{\frac 12} \snorm{\nb^3 h,_2}_{{L^2_t L^\infty_x}} + \eps \snorm{\nb^2 h,_{2}}_{L^\infty_{x,t}}^2  
\notag\\
\snorm{\nb\, \Ewn}_{L^\infty_{x,t}}
&\les  1 + \snorm{\nb^2 h,_2 }_{L^\infty_{x,t}} +  \snorm{\nb^2(\Jg\Abn)}_{L^\infty_{x,t}} + \snorm{\nb^2 \Jg }_{L^\infty_{x,t}} 
\notag\\
\snorm{\nb^2\, \Ewn}_{{L^2_t L^\infty_x}}
&\les \eps^{\frac 52} 
+  \snorm{\nb^2 (\Jg \Abn,_2)}_{{L^2_t L^\infty_x}}
+  \snorm{\nb \Abn,_2}_{{L^2_t L^\infty_{x}}}
+ \eps^{\frac 12} \snorm{\nb^2 (\Sigma,\Jg)}_{L^\infty_{x,t}}
+  \snorm{\nb^3 \Jg}_{{L^2_t L^\infty_{x}}} \notag\\
&\qquad + \eps^{\frac 32} \snorm{\nb^2 h,_2}_{L^\infty_{x,t}}
+  \snorm{\nb^2 (\Wbt,\Zbt,\Jg \Abn)}_{{L^2_t L^\infty_{x}}}
+ \eps  \snorm{\nb^2 \Abt}_{{L^2_t L^\infty_{x}}}
+\eps  \snorm{\nb^2 (\Jg \Zbn)}_{{L^2_t L^\infty_{x}}}
\,.
\end{align*}
The $L^\infty_{x,t}$ norms present in the above estimate, which are the same as $L^\infty_{x,\s}$ norms cf.~\eqref{eq:norm:equivalence:a}, are estimated using \eqref{eq:Sobolev}. On the other hand, the $L^2_t L^\infty_x$ norms are not equivalent to $L^2_\s L^\infty_x$ norms, akin to the discussion above~\eqref{eq:norm:equivalence:c}. Nonetheless, the proof of estimate \eqref{eq:Steve:punch} (designed for $L^2_\s L^\infty_x$ norms) can still be applied in this case,\footnote{\label{foot:L2xt} First, at each fixed $t$ one uses fundamental theorem of calculus in the $x_1$ and $x_2$ variables, and the Poincar\'e inequality in the $x_1$ variable to bound the $L^\infty_x( (\{t\}\times \TT^2) \cap \mathcal{P})$ norm of a function $F$ in terms of the $L^2_x((\{t\}\times \TT^2) \cap \mathcal{P})$ norm of $\eps^{-\frac 12} \nb_1 \nb F$ (see~\eqref{eq:L2:FTC:space}). Second, one integrates this expression in $t$, leading to a bound for $\|F\|_{L^2_t L^\infty_x}$ in terms of $\eps^{-\frac 12} \|\nb_1 \nb F\|_{L^2_{x,t}}$. Third, one uses~\eqref{eq:norm:equivalence:b} to transfer an upper bound in terms of $\|\nb_1 \nb F\|_{L^2_{x,t}}$ into an upper bound in terms of $\|\nbs_1 \nbs F\|_{L^2_{x,\s}}$. This argument is used repeatedly throughout this section.} leading to obtain upper bounds in terms of $L^2_{x,\s}$ norms for $5^{th}$ order derivatives. In turn, to bound these we appeal to the bootstraps~\eqref{bootstraps-Dnorm:6}--\eqref{bootstraps-Dnorm:Jg}, the bounds for the geometry and sound speed in~\eqref{geometry-bounds-new}, and to the vorticity estimate~\eqref{eq:vort:H6}.  
For instance, \eqref{eq:Sobolev}, \eqref{D6h2Energy:new}, and the fact that $\mathcal{J} \leq 1$, imply that $\| \nb^2 h,_2\|_{L^\infty_{x,t}} \les \eps^{-\frac 12} \| \nb^4 h,_2(\cdot,\initial)\|_{L^2_x} + \eps^{-1} \|\nb^5 h,_2\|_{L^2_{x,t}} \les \mathsf{K} \eps \brak{\mathsf{B_6}}$. Similar arguments show that $\|\nb^2 \Jg\|_{L^\infty_{x,t}} \leq \brak{\mathsf{B_6}}$, $\|\nb^2 \Sigma\|_{L^\infty_{x,t}} \leq \brak{\mathsf{B_6}}$, $\|\nb^3 h,_2\|_{{L^2_t L^\infty_x}} \les \mathsf{K} \eps^{\frac 32} \brak{\mathsf{B_6}}$, and $\|\nb^3 \Jg\|_{{L^2_t L^\infty_x}} \les  \eps^{\frac 12} \brak{\mathsf{B_6}}$. Similar arguments imply for $(\Wbt,\Zbt,\Abt,\Jg \Zbn,\Jg\Abn)$. For instance, \eqref{eq:Steve:punch} and \eqref{eq:tilde:D5} imply that $\|\nb^2  (\Wbt,\Zbt, \Abt)\|_{{L^2_t L^\infty_x}} \les \eps^{-\frac 12} \|\nb^4 (\Wbt,\Zbt, \Abt)\|_{L^2_{x,t}} \les \eps^{-\frac 12} \tilde{\mathcal{D}}_{5,\ttt} \les \mathsf{K} \eps^{\frac 12} \brak{\mathsf{B_6}}$.  The same argument shows $\|\nb^2  (\Jg \Zbn) \|_{{L^2_t L^\infty_x}} \les \eps^{-\frac 12} \brak{\mathsf{B_6}}$. Some of the bounds for the terms involving $\Abn$ are more delicate and require that we relate $\Abn$ to the vorticity via $\Abn = \Omega + \tfrac 12 \Wbt + \tfrac 12 \Zbt $. These vorticity-improved bounds are obtained in Section~\ref{sec:vorticity} (see~\eqref{eq:Jg:Abn:D5:improve:a} and~\eqref{eq:Jg:Abn:D5:improve:c}). For instance, \eqref{eq:Jg:Abn:D5:improve:c} implies $\|\nb^2 \Abn\|_{{L^2_t L^\infty_x}} \les  \eps^{-\frac 12} \|\nb^4 \Abn\|_{L^2_{x,t}} \les \mathsf{K} \eps^{\frac 12} \brak{\mathsf{B_6}} $. This bound, the bootstraps~\eqref{bootstraps}, and the product rule in turn gives $\|\nb^2(\Jg\Abn)\|_{{L^2_t L^\infty_{x}}} \les \mathsf{K} \eps^{\frac 12} \brak{\mathsf{B_6}}$. Next, by appealing to the improved bounds in Corollary~\ref{cor:Abn:improve}, and to Lemma~\ref{lem:Moser:tangent}, we may show that $\|\nb^2(\Jg \Abn)\|_{L^\infty_{x,t}} \les 1 + \eps^{-1} \|\nb^5(\Jg\Abn)\|_{L^2_{x,t}} \les \brak{\mathsf{B_6}} +  \eps^{-1} \| \nb^5 \Abn \|_{L^2_{x,t}} \les \mathsf{K} \brak{\mathsf{B_6}}$, and that $\|\nb^2(\Jg \Abn,_2)\|_{{L^2_t L^\infty_x}} \les \eps^{-\frac 12} \|\nb^5(\Jg\Abn)\|_{L^2_{x,t}} + \eps^{\frac 12} \brak{\mathsf{B_6}} \les \mathsf{K} \eps^{\frac 12} \brak{\mathsf{B_6}}$. Inserting these bounds into the previously established bounds for $\nb \Iwn1, \nb^2 \Iwn1, \nb\Ewn$, and $\nb^2 \Ewn$, we derive 
\begin{subequations}
\label{eq:broncos:suck}
\begin{align}
\snorm{\nb \, \Iwn1}_{L^\infty_{t,x}}
&\les \mathsf{K}  \eps^2 \brak{\mathsf{B_6}} \\
\snorm{\nb^2 \, \Iwn1}_{L^\infty_{t,x}}
&\les \mathsf{K}  \eps \brak{\mathsf{B_6}}   \\
\snorm{\nb\, \Ewn}_{L^\infty_{t,x}}
&\les \mathsf{K}   \brak{\mathsf{B_6}}   \\
\snorm{\nb^2\, \Ewn}_{{L^2_t L^\infty_x}}
&\les \mathsf{K}  \eps^{\frac 12} \brak{\mathsf{B_6}}
\,.
\end{align}
The bounds on $\Ewn$ and the integrating factor $\Iwn1$ imply also that 
\begin{equation}
\snorm{\nb \, \Iwn2}_{L^\infty_{t,x}}
+ \snorm{\nb^2 \, \Iwn2}_{L^\infty_{t,x}}
\les  \mathsf{K}\eps \brak{\mathsf{B_6}} 
\,. 
\end{equation}
\end{subequations}
With \eqref{Iwn1-2-est} and \eqref{eq:broncos:suck}, we return to \eqref{Jg-Wbn-1} and obtain the pointwise estimates
\begin{subequations}
\label{eq:broncos:eat:shit}
\begin{align}
\sabs{\nb  \bigl( (\Jg \Wbn)\circ \xi(x,t) - (w_0),_1(x)\bigr)}
&\les \mathsf{K}\eps \brak{\mathsf{B_6}} 
\label{eq:broncos:eat:shit:a}
\\
\sabs{\nb^2 \bigl( (\Jg \Wbn)\circ \xi(x,t) - (w_0),_1(x)\bigr)}
&\les \mathsf{K}  \brak{\mathsf{B_6}} 
\label{eq:broncos:eat:shit:b}
\end{align}
\end{subequations}
for all $(x,t) \in \mathcal{P}$.
We note that \eqref{eq:broncos:eat:shit} and \eqref{eq:xi:switch}, together with \eqref{eq:broncos:eat:shit:20} and \eqref{eq:d1:w0:close} imply \eqref{eq:broncos:eat:shit:2a}--\eqref{eq:broncos:eat:shit:2b}. 
\end{proof}

Integrating the bounds obtained in Lemma~\ref{lem:JgWbn:close:to:data}, we obtain pointwise bounds for $\Jg$ as follows:
\begin{corollary}
\label{cor:Jg:initial} 
Assume that the bootstrap bounds \eqref{bootstraps} hold on $\TT^2 \times [0,\eps)$.
If $\eps$ is taken to be sufficiently small with respect to $\alpha, \kappa_0$, and $\Cdata$, then the bounds
\begin{subequations}
\label{eq:Jg:Wbn:horse}
\begin{align}
\sabs{  \Jg(x,t) -  1- (t-\initial) \tfrac{1+\alpha}{2}  (w_0),_1(x) }
&\leq \tfrac{1+\alpha}{2} \mathsf{C_{J_t}} (t-\initial),
\label{bs-Jg-0}\\
\sabs{\nb_i   \Jg  (x,t) - (t-\initial) \tfrac{1+\alpha}{2}  \nb_i  (w_0),_1(x) }
&\les (t-\initial) \eps \mathsf{K}  \brak{\mathsf{B_6}}, \qquad \mbox{for} \qquad i\in\{1,2\},
\label{bs-Jg-1}\\
\sabs{(\p_t + V \p_2) \Jg  (x,t) -  \tfrac{1+\alpha}{2}     (w_0),_1(x) }
&\leq \tfrac{ 5|1-\alpha|}{6} \mathsf{C_{\Zbn}},  
\label{bs-Jg-1a}\\
\sabs{\p_t \Jg  (x,t) -  \tfrac{1+\alpha}{2}  (w_0),_1(x) }
&\leq  \tfrac{1+\alpha}{2} \mathsf{C_{J_t}},  
\label{bs-Jg-1b}\\
\sabs{\nb_i\nb_j   \Jg  (x,t) - (t-\initial) \tfrac{1+\alpha}{2}  \nb_i \nb_j  (w_0),_1(x) }
&\les (t-\initial)   \mathsf{K}  \brak{\mathsf{B_6}},
\qquad \mbox{for} \qquad i,j \in \{1,2\},
\label{bs-Jg-2}\\
\sabs{\nb_i\p_t   \Jg  (x,t) -  \tfrac{1+\alpha}{2}  \nb_i  (w_0),_1(x) }
&\leq   8(1+\alpha) \mathsf{C_{\Zbn}},
\qquad \mbox{for} \qquad i \in \{1,2\},
\label{bs-Jg-2a}\\
\sabs{\p_t \p_t   \Jg  (x,t) }
&\leq  4(1+\alpha)|1-\alpha|  \mathsf{C_{\Zbn}} \eps^{-1},
\label{bs-Jg-2b}
\end{align}
hold  for all $(x,t) \in \mathcal{P}$. Here we have introduced  $\mathsf{C_{J_t}} = \mathsf{C_{J_t}}(\alpha,\kappa_0,\Cdata)>0$, defined by $|1-\alpha| \mathsf{C_{\Zbn}} = \frac{1+\alpha}{2} \mathsf{C_{J_t}}$. 
\end{subequations}
\end{corollary}
\begin{proof}[Proof of Corollary~\ref{cor:Jg:initial}]
Upon composing \eqref{Jg-evo} with the flow $\xi$ and integrating in time, we obtain that for $t \leq \mathsf{T}_{\xi}(x)$, 
\begin{equation}
\Jg \cir \xi (x,t) 
= 1 + \int_{\initial}^t \left( \tfrac{1+\alpha}{2} (\Jg  \Wbn ) \cir \xi + \tfrac{1-\alpha}{2} (\Jg  \Zbn) \cir \xi\right) (x,t') {\rm d}t' 
\,. \label{Jg-circ-xi}
\end{equation}
Using the bootstrap assumptions~\eqref{bootstraps}, and the previously established bound \eqref{eq:broncos:eat:shit:20} , we deduce from the above identity that
\begin{equation}
\sabs{\Jg \cir \xi (x,t) - 1 - (t-\initial) \tfrac{1+\alpha}{2} (w_0),_1(x)}
\leq   (t-\initial) \bigl( \eps \Cn + \tfrac{3|1-\alpha|}{4} \mathsf{C}_{\Zbn} \bigr)  \,.
\end{equation}
Upon composing with $\xi^{-1}(x,t)$ and appealing to \eqref{eq:d1:w0:close} with $k=0$, we deduce \eqref{bs-Jg-0}.

In a similar fashion, we may differentiate \eqref{eq:xi:nabla:xi} with respect to space, appeal to \eqref{Jg-Wbn-1}, and deduce that for $i \in \{1,2\}$ we have
\begin{align}
&\nb_i (\Jg  \cir \xi) (x,t)  - (t-\initial) \tfrac{1+\alpha}{2} \nb_i(w_0),_1(x)\notag\\
&=   \int_{\initial}^t \Bigl( \tfrac{1+\alpha}{2} (w_0),_1(x) \nb_i \Iwn1(x,t') + \tfrac{1+\alpha}{2} \nb_i (w_0),_1(x) (\Iwn1(x,t')-1)  \notag\\
&\qquad \qquad \qquad + \tfrac{1+\alpha}{2} \nb_i \Iwn2(x,t')+ \tfrac{1-\alpha}{2} \nb_i ( (\Jg  \Zbn) \cir \xi)\Bigr) (x,t') {\rm d}t'
\,. 
\label{eq:nb:J:circ:xi}
\end{align}
Using the initial data assumptions, the bounds \eqref{Iwn1-2-est}, \eqref{eq:broncos:suck}, \eqref{eq:xi:switch}, \eqref{bs-nnZb}, \eqref{bs-Jg-simple}, and \eqref{bs-Jg,1}, we deduce
\begin{equation}
\sabs{\nb_i (\Jg  \cir \xi) (x,t)  - (t-\initial) \tfrac{1+\alpha}{2} \nb_i(w_0),_1(x)}
\les  (t-\initial) \eps \mathsf{K}  \brak{\mathsf{B_6}}
\,.
\end{equation}
Next, instead of merely appealing to \eqref{eq:xi:nabla:xi}, we use that similarly to \eqref{eq:xi:nabla:xi} we may show $|\xi,_2(x,t)-1|\les \eps (t-\initial)$ and $|\xi,_1(x,t)|\les (t-\initial)$. The resulting bounds are
\begin{equation}
\sabs{(\nb_i \Jg)  \circ \xi(x,t)  - (t-\initial) \tfrac{1+\alpha}{2} \nb_i(w_0),_1(x)}
\les  (t-\initial) \eps \mathsf{K}  \brak{\mathsf{B_6}}\,,
\qquad \mbox{for} \qquad i \in \{1,2\}\,.
\end{equation}
Combined with \eqref{eq:d1:w0:close} with $k=1$, this bound implies \eqref{bs-Jg-1}.

Differentiating \eqref{Jg-circ-xi} with respect to $t$, recalling the definition of $\xi$ in \eqref{xi-flow}, the identity \eqref{Jg-Wbn-1}, and the bounds \eqref{Iwn1-2-est}, \eqref{bs-nnZb}, and \eqref{bs-Jg-simple}, we obtain
\begin{equation}
\sabs{( (\p_t + V \p_2) \Jg)  \circ \xi(x,t)  -  \tfrac{1+\alpha}{2}  (w_0),_1(x)}
\leq  \Cn \eps + \tfrac{3|1-\alpha|}{4} \mathsf{C_{\Zbn}}
\,.
\end{equation}
The bound \eqref{bs-Jg-1a} now follows upon composing with $\xi^{-1}$ and using \eqref{eq:d1:w0:close} with $k=0$.  The bound \eqref{bs-Jg-1b} follows from \eqref{bs-Jg-1} with $i=2$, \eqref{bs-Jg-1a}, and  \eqref{bs-V}. 

Next, we establish \eqref{bs-Jg-2}. We apply $\nb_j$ to \eqref{eq:nb:J:circ:xi}, appeal to the bounds \eqref{Iwn1-2-est}, \eqref{eq:broncos:suck}, and \eqref{eq:xi:switch}, to obtain
\begin{equation}
\sabs{\nb_j\nb_i (\Jg  \cir \xi) (x,t)  - (t-\initial) \tfrac{1+\alpha}{2} \nb_j\nb_i(w_0),_1(x)}
\les (t-\initial)   \mathsf{K}  \brak{\mathsf{B_6}} + (t-\initial) \|\nb^2(\Jg \Zbn)\|_{L^\infty_{x,t}} \,.
\label{eq:Knicks:1}
\end{equation}
To bound the last term in the above estimate, we  use \eqref{eq:Sobolev} and the improved $\Jg \Zbn$ bound available from \eqref{eq:Jg:Zbn:D5:improve:a}. We obtain 
\begin{equation}
\|\nb^2(\Jg \Zbn)\|_{L^\infty_{x,t}} 
\les \eps^{-\frac 12} \|\nb^4(\Jg \Zbn)(\cdot,\initial)\|_{L^2_{x}} + \eps^{-1} \|\nb^5(\Jg \Zbn)\|_{L^2_{x,\s}}
\les \mathsf{K} \brak{\mathsf{B_6}} \,.
\label{eq:Knicks:2}
\end{equation}
The bounds \eqref{eq:Knicks:1} and \eqref{eq:Knicks:2} thus imply 
\begin{equation}
\sabs{\nb_j\nb_i (\Jg  \cir \xi) (x,t)  - (t-\initial) \tfrac{1+\alpha}{2} \nb_j\nb_i(w_0),_1(x)}
\les  (t-\initial)  \mathsf{K}  \brak{\mathsf{B_6}} \,.
\label{eq:Knicks:3}
\end{equation}
Appealing to \eqref{bs-Jg-1}, \eqref{eq:xi:hessian:xi},   to the bounds $|\xi,_2(x,t)-1|\les \eps (t-\initial)$ and $|\xi,_1(x,t)|\les (t-\initial)$, and to the previously derived estimate $\|\nb^2 \Jg\|_{L^\infty_{x,t}} \leq \brak{\mathsf{B_6}}$,  we   also obtain
\begin{equation}
 \sabs{\nb_j\nb_i (\Jg\circ\xi)(x,t) - (\nb_j\nb_i  \Jg)\cir \xi (x,t)}
\les  (t-\initial)    \mathsf{K}  \brak{\mathsf{B_6}} 
\label{eq:Knicks:4}
\,.
\end{equation} 
Upon combining \eqref{eq:Knicks:3} and \eqref{eq:Knicks:4}, composing with $\xi^{-1}$, and using \eqref{eq:d1:w0:close} with $k=2$, we arrive at the proof of \eqref{bs-Jg-2}.

In order to prove \eqref{bs-Jg-2a} we apply $\nb_i \p_t$ to \eqref{Jg-circ-xi}, use \eqref{eq:broncos:eat:shit:a} and the bootstrap assumptions~\eqref{bootstraps}, and obtain that 
\begin{subequations}
\begin{align}
\sabs{\xi,_2 (\nb_2 (\p_t + V \p_2) \Jg) \circ \xi   - \tfrac{1+\alpha}{2} \nb_2 (w_0),_1}
&\leq \Cn \eps \mathsf{K}  \brak{\mathsf{B_6}} +  7(1+\alpha) \mathsf{C_{\Zbn}}\,, \\
\sabs{\xi,_2 (\nb_1 (\p_t + V \p_2) \Jg) \circ \xi   - \tfrac{1+\alpha}{2} \nb_1 (w_0),_1 - \eps \xi,_1 (\nb_2 (\p_t + V \p_2) \Jg) \circ \xi }
&\leq \Cn \eps \mathsf{K}  \brak{\mathsf{B_6}} +  7(1+\alpha) \mathsf{C_{\Zbn}}  \,.
\end{align}
\end{subequations}
By again appealing to \eqref{eq:xi:nabla:xi} and also to \eqref{bs-V} and to $\|\nb^2 \Jg\|_{L^\infty_{x,t}} \leq \brak{\mathsf{B_6}}$, we derive from the above that
\begin{equation}
\sabs{ (\nb_i  \p_t  \Jg) \circ \xi   - \tfrac{1+\alpha}{2} \nb_i (w_0),_1}
\leq \Cn \eps \mathsf{K}  \brak{\mathsf{B_6}} +  7(1+\alpha) \mathsf{C_{\Zbn}}\,,
\end{equation}
for $i \in \{1,2\}$.
The bound \eqref{bs-Jg-2a} now follows from the above estimate, upon composing with $\xi^{-1}$, and from \eqref{eq:d1:w0:close} with $k=1$.

The last estimate in \eqref{eq:Jg:Wbn:horse}, namely \eqref{bs-Jg-2b}, follows by  differentiating \eqref{Jg-evo} with respect to time, which yields
\begin{align}
 (\p_t\p_t + \tfrac{1}{\eps} \nb_t V \p_2 + V \p_t \p_2)  \Jg  = \tfrac{1+\alpha}{2 }  (\p_t + V \p_2) (\Jg \Wbn) -
 \tfrac{1+\alpha}{2} V \p_2(\Jg \Wbn) + \tfrac{1-\alpha}{2\eps}  \nb_t   (\Jg \Zbn) \,.
 \label{eq:cant:believe:I:need:this:shit}
\end{align}
Using the previously established bound \eqref{bs-Jg-2a},   the bootstrap assumptions~\eqref{bootstraps},   the time differentiated version of \eqref{Jg-Wbn-1}  which gives $(\p_t + V \p_2) (\Jg \Wbn) \circ \xi(x,t) = (w_0),_1(x) \p_t \Iwn1(x,t) + \p_t \Iwn2(x,t)$, and the bounds \eqref{eq:broncos:suck}, we deduce that  
\begin{equation*}
\sabs{\p_t\p_t \Jg} \leq \Cn + \tfrac{|1-\alpha|}{2\eps} \mathsf{C_{\Zbn}} 7(1+\alpha)\,. 
\end{equation*}
The bound  \eqref{eq:Jg:Wbn:horse} now follows, concluding the proof of the Corollary.
\end{proof}

\subsection{Properties of the remapping coefficients}
\label{sec:Q:coeff}

We recall the coefficients $\Q,\Qb_2,\Qd,\Qc,\Qr_\s$, and $\Qr_2$ introduced in \eqref{QQQ}. These coefficients are bounded  as follows:

\begin{lemma}
\label{lem:Q:bnds}
Assume that the bootstrap bounds \eqref{bootstraps} hold on $\TT^2 \times [0,\eps)$.
If $\eps$ is taken to be sufficiently small with respect to $\alpha, \kappa_0$, and $\Cdata$, then the   various functions appearing in  in \eqref{QQQ} bounds are bounded as
\begin{subequations}
\label{eq:Q:all:bbq}
\begin{align} 
\tfrac{17 (1+\alpha) }{40}      &\le \Qd(x_2,\s)  \le  401  (1+\alpha)  \,,  
\label{Qd-lower-upper} \\ 
\sabs{\Q - \Qd} &\leq 3 \mathsf{C_V} \eps \s  \leq 3 \mathsf{C_V} \eps^2 \,,  
\label{Q-lower-upper} \\
\sabs{\Qb_2} &\leq 3 \s \leq 3 \eps\,, 
\label{eq:Qb:bbq} \\
\sabs{\Qr_\s} &\leq   \tfrac{2 \cdot 250^2 \Q}{\eps}   \,, 
\label{eq:Qrs:bbq} \\
 \sabs{\Qr_2 }  &\leq  \Cn  \,, 
\label{eq:Qr2:bbq} 
\\
\sabs{\Qc - \Qr_\s} &\leq \Cn \eps
\,, 
\label{eq:Qc:bbq} 
\end{align} 
hold uniformly for all $(x,\s) \in \TT^2 \times [0,\eps)$.  
In particular, \eqref{Qd-lower-upper} and \eqref{Q-lower-upper} imply that 
\begin{equation}
\tfrac{2(1+\alpha)}{5} \leq \Q \leq  402 (1+\alpha)  \, .
\label{eq:Q:bbq}
\end{equation}
\end{subequations}
\end{lemma}
\begin{proof}[Proof of Lemma~\ref{lem:Q:bnds}]
From~\eqref{eq:QQQ:a}, \eqref{bs-Jg-1b}, and \eqref{eq:Jgb:identity:0} it follows that in order to obtain a bound on $\Qd$, we first need a bound on $(w_0),_1(x_1^*(x_2,t),x_2)$. 
To obtain such an estimate, we recall that by \eqref{eq:Jgb:Jg:identity:3}  $x_1^*(x_2,t)$ is the point at which the global minimum of $\Jg(\cdot,x_2,t)$ is attained, and hence by \eqref{bs-Jg-0}  we have
\begin{align}
1 + (t-\initial) \tfrac{1+\alpha}{2} \bigl( (w_0),_1(x_1^*(x_2,t),x_2) - \mathsf{C_{J_t}} \bigr)
&\leq \Jg(x_1^*(x_2,t),x_2,t)
\notag\\
&\leq \Jg(x_1^\vee(x_2),x_2,t)
\notag\\
&\leq 
1 + (t-\initial) \tfrac{1+\alpha}{2} \bigl( (w_0),_1(x_1^\vee(x_2),x_2) + \mathsf{C_{J_t}} \bigr) \,,
 \label{eq:d1w0:x1:star:temp}
\end{align}
where we recall cf.~assumption \eqref{item:ic:w0:x2:negative} that $x_1^\vee(x_2)$ is the point at which $(w_0),_1(\cdot,x_2)$ attains a minimum. Therefore,
\begin{equation}
 -\tfrac{1}{\eps}
 \leq
 (w_0),_1(x_1^\vee(x_2),x_2)
 \leq 
 (w_0),_1(x_1^*(x_2,t),x_2)
 \leq 
 (w_0),_1(x_1^\vee(x_2),x_2) + 2 \mathsf{C_{J_t}} 
 \leq 
 - \tfrac{9}{10\eps} +  2 \mathsf{C_{J_t}} 
 \label{eq:d1w0:x1:star}
 \,.
\end{equation}
Inserting the above obtained bounds for  $(w_0),_1(x_1^*(x_2,t),x_2)$ into  \eqref{bs-Jg-1b} yields
\begin{align}
\tfrac{1+\alpha}{2} \bigl(-\tfrac{1}{\eps}  - \mathsf{C_{J_t}} \bigr) 
&\leq
\tfrac{1+\alpha}{2} \bigl((w_0),_1(x_1^*(x_2,t),x_2) - \mathsf{C_{J_t}} \bigr)
\notag\\
&\leq 
(\p_t \Jg)(x_1^*(x_2,t),x_2,t) 
\notag\\
&\leq 
\tfrac{1+\alpha}{2} \bigl((w_0),_1(x_1^*(x_2,t),x_2) + \mathsf{C_{J_t}} \bigr)
\notag\\
&\leq 
\tfrac{1+\alpha}{2} \bigl(-\tfrac{9}{10\eps}  + \mathsf{C_{J_t}} \bigr) 
\,.
 \label{eq:dtJg:x1:star}
\end{align}
Next, we use the bounds \eqref{eq:d1w0:x1:star} and \eqref{eq:dtJg:x1:star} in order to prove \eqref{Qd-lower-upper}.

Differentiating \eqref{t-to-s-transform} with respect to $t$ and appealing to \eqref{eq:Jgb:identity:0}, we have 
\begin{align}
\Qd(x_2,\mathfrak{q}(x_2,t)) 
=
(\p_t \mathfrak{q})(x_2,t) 
&= 
-\eps  (\p_t \Jgb)(x_1^*(x_2,t), x_2,t) \notag\\
&=
-\eps  (\p_t \Jg)(x_1^*(x_2,t), x_2,t) 
+200(1+\alpha) \mathfrak{C}\bigl(\tfrac{t-\medium}{\final-\medium}\bigr)
\,.
\end{align}
We observe that \eqref{eq:dtJg:x1:star} and the definition of $\mathfrak{C}$ implies
\begin{equation}
\tfrac{1+\alpha}{2} (\tfrac{9}{10} - \eps \mathsf{C_{J_t}}) 
\leq
\Qd(x_2,\s) 
\leq 
\tfrac{1+\alpha}{2} (1 + \eps \mathsf{C_{J_t}}) + 400 (1+\alpha) {\bf 1}_{t\in (\medium,\final]}
\,.
\label{eq:Qcal:bbq:temp:3}
\end{equation}
The bound \eqref{Qd-lower-upper}  now immediately follows.

Due to \eqref{eq:QQQ:b} we have that  
\begin{equation*}
|\Q - \Qd|(x,\s) = \eps |V(x,\s)| |\p_2 \Jgb|(x_1^*(x_2,t),x_2,t)|_{t=\mathfrak{q}^{-1}(x_2,\s)} = \eps |V(x,\s)| |\p_2 \Jg|(x_1^*(x_2,t),x_2,t)|_{t=\mathfrak{q}^{-1}(x_2,\s)}
\end{equation*}
and so by appealing to \eqref{bs-V} and  \eqref{bs-Jg-1} with $i=2$, we obtain
\begin{align}
|\Q - \Qd| 
&\leq 
\eps^2 \mathsf{C_V} (\mathfrak{q}^{-1}(x_2,\s)-\initial) \tfrac{1+\alpha}{2} \bigl( |(w_0),_{12}(x_1^*,x_2)| + \Cn \mathsf{K} \eps \jap{\mathsf{B_6}}\bigr) 
\notag\\
&\leq \eps (1+\alpha) \mathsf{C_V} (\mathfrak{q}^{-1}(x_2,\s)-\initial)
\label{eq:Qcal:bbq:temp:1}
\,.
\end{align}
The bound \eqref{Q-lower-upper} will follow if we are able to establish $\mathfrak{q}^{-1}(x_2,\s) - \initial \les \s$. Using that $\mathfrak{q}(x_2,\initial) = 0$, we write $\s = \mathfrak{q}(x_2,t) =\mathfrak{q}(x_2,t)- \mathfrak{q}(x_2,\initial) = \int_{\initial}^t (\p_t \mathfrak{q})(x_2,t') {\rm d}t'$.  
Integrating the lower bound in \eqref{eq:Qcal:bbq:temp:3} with respect to time, 
\begin{equation}
 \s 
 = \int_{\initial}^t (\p_t \mathfrak{q})(x_2,t') {\rm d}t'
 \geq \tfrac{2(1+\alpha)}{5} (t-\initial)
 \,.
 \label{eq:s:t:lower}
\end{equation}
Note that a reverse inequality  also holds:
\begin{equation}
 \s 
 = \int_{\initial}^t (\p_t \mathfrak{q})(x_2,t') {\rm d}t'
\leq 401 ( 1+\alpha) (t-\initial)
 \,.
 \label{eq:s:t:upper}
\end{equation}
From \eqref{eq:s:t:lower} and \eqref{eq:s:t:upper} we obtain that 
\begin{equation}
\tfrac{2(1+\alpha)}{5} \leq \tfrac{\s}{\mathfrak{q}^{-1}(x_2,\s) -\initial} \leq  401( 1+\alpha) 
\,.
\label{eq:Qcal:bbq:temp:2}
\end{equation}
Combining \eqref{eq:Qcal:bbq:temp:1} and \eqref{eq:Qcal:bbq:temp:2}, we immediately arrive at \eqref{Q-lower-upper}.

Next, we recall the definition of $\Qb_2(x_2,\s)$ in \eqref{eq:QQQ:aa}. Using \eqref{eq:why:the:fuck:not:0}, \eqref{bs-Jg-1}, \eqref{eq:Jgb:identity:2}, and \eqref{eq:Qcal:bbq:temp:2} we deduce that 
\begin{equation}
|\Qb_2(x_2,\s)|\leq   (t-\initial) \tfrac{1+\alpha}{2} \bigl( |\nb_1 \nb_2 (w_0)(x_1^*,x_2)| + \Cn \eps^2 \mathsf{K}  \jap{\mathsf{B_6}}\bigr) 
\leq (\tfrac{102}{50} + \Cn \eps^2\mathsf{K}  \brak{\mathsf{B_6}}) \s
\label{Qb2-bound} 
\end{equation}  
The bound \eqref{eq:Qb:bbq}  follows.

In order to prove \eqref{eq:Qrs:bbq}--\eqref{eq:Qc:bbq}, we first note that implicit differentiation of the relation $\Jgb,_1(x_1^*(x_2,t),x_2,t) = \Jg,_1(x_1^*(x_2,t),x_2,t) = 0$, with respect to $x_2$ and with respect to $t$ gives the relations 
\begin{equation}
\p_2 x_1^*(x_2,t) = - \bigl(\tfrac{\p_1 \p_2 \Jg }{\p_1 \p_1 \Jg }\bigr)(x_1^*(x_2,t),x_2,t)
\,, \qquad \mbox{and} \qquad 
\p_t x_1^*(x_2,t) = - \bigl(\tfrac{\p_t \p_1 \Jg }{\p_1 \p_1 \Jg }\bigr)(x_1^*(x_2,t),x_2,t)
\,.
\label{p2-and-pt-x1star}
\end{equation}
In turn, cf.~\eqref{eq:Jgb:identity:2} and~\eqref{QQQ} this implies 
\begin{subequations}
\label{eq:Qring:2s:hat}
\begin{align}
 \Qr_\s (x_2,\s)&= \tfrac{-\eps}{\Qd(x_2,\s)} 
\Bigl(  \p_t \p_t  \Jgb  - \tfrac{ \p_1 \p_t \Jg \cdot \p_1\p_t\Jg}{\p_1 \p_1 \Jg}    \Bigr)(x_1^*(x_2,t),x_2,t) \Bigr|_{t = \mathfrak{q}^{-1}(x_2,\s)}
 \, \\
 \Qr_2 (x_2,\s)&= \tfrac{\eps}{\Qd(x_2,\s)}
 \Bigl(  \p_t \p_2 \Jg - \tfrac{ \p_1 \p_2 \Jg 
 \cdot \p_1 \p_t \Jg }{\p_1 \p_1 \Jg}   \Bigr)(x_1^*(x_2,t),x_2,t) \Bigr|_{t = \mathfrak{q}^{-1}(x_2,\s)}
 \, \\
 \Qc(x_2,\s)&= \Qr_\s (x_2,\s) - \Qr_2 (x_2,\s) V(x,\s)
 -  \Qb_2(x_2,\s)   \p_\s V  (x,\s) \,.
\end{align}
\end{subequations}
It thus becomes apparent that we require a positive lower bound for $\Jg,_{11}(x_1^*(x_2,t),x_2,t)$; we note that this is the only place in the argument where assumption~\eqref{item:ic:w0:x2:special} on the initial data enters. We revisit the second and third inequalities in \eqref{eq:d1w0:x1:star}, which show that $0 \leq \nb_1 w_0(x_1^*(x_2,t),x_2) -  \nb_1 w_0(x_1^\vee(x_2,t),x_2) \leq 2 \eps \mathsf{C_{J_t}}$. Since $2 \eps \mathsf{C_{J_t}} < \eps^{\frac 34}$ for $\eps$ sufficiently small, by assumption~\eqref{item:ic:w0:x2:special} this implies that $|x_1^*(x_2,t)- x_1^\vee(x_2)| < \eps^{\frac 54}$. We can however obtain an improved bound for this difference. We recall that $\Jg,_1(x_1^*(x_2,t),x_2,t) = 0$ (see~\eqref{eq:x1star:critical} and \eqref{eq:Jgb:identity:2}), while \eqref{bs-Jg-1}, \eqref{bootstraps-Dnorm:6}, and \eqref{bootstraps-Dnorm:5} imply that $|\Jg,_1(x_1^*(x_2,t),x_2,t) - (t-\initial) (w_0),_{11}(x_1^*(x_2,t),x_2)|\leq \Cn (t-\initial) \mathsf{K}  \jap{\mathsf{B}_6}$. Together, these bounds yield
\begin{equation}
\sabs{(w_0),_{11}(x_1^*(x_2,t),x_2)} \leq \Cn \mathsf{K} \jap{\mathsf{B}_6}.
\label{eq:d1d1w0}
\end{equation}
On the other hand, using assumption~\eqref{item:ic:w0:x2:negative} and~\eqref{table:derivatives} we may perform a second order Taylor expansion in $x_1$ (at fixed $x_2$) around $x_1^{\vee} = x_1^{\vee}(x_2)$, to obtain
\begin{equation}
\underbrace{(w_0),_{11}(x_1^*,x_2)}_{|\cdot|\leq \Cn \jap{\mathsf{B}_6}} = \underbrace{(w_0),_{11}(x_1^\vee,x_2)}_{=0} + (x_1^*-x_1^\vee) \underbrace{(w_0),_{111}(x_1^\vee,x_2)}_{\geq \frac{9}{10\eps^3}} + \tfrac 12 (x_1^*-x_1^\vee)^2 \underbrace{(w_0),_{1111}(x_1^\sharp,x_2)}_{|\cdot|\leq \frac{\Cdata}{\eps^4}}
\end{equation}
for some $x_1^\sharp$ that lies in between $x_1^*$ and $x_1^\vee$. Moreover, using that $|x_1^*(x_2,t)- x_1^\vee(x_2)| < \eps^{\frac 54}$ we have that $|(x_1^*-x_1^\vee) (w_0),_{1111}(x_1^\sharp,x_2)| \leq \Cdata \eps^{-3 + \frac 14}$ and thus
\begin{equation}
 |x_1^*(x_2,t)- x_1^\vee(x_2)|
\leq \Cn \mathsf{K} \jap{\mathsf{B_6}} \bigl( \tfrac{9}{10\eps^3} - \tfrac{\Cdata \eps^{\frac 14}}{2 \eps^{3}}\bigr)^{-1}  
 \leq  2 \Cn \mathsf{K} \jap{\mathsf{B_6}}  \eps^3 \,.
 \label{eq:x1star:x1vee}
\end{equation}
Notice the improvement of $\OO(\eps^{\frac 54}) \mapsto \OO(\eps^3)$ that the above bound gives over assumption \eqref{item:ic:w0:x2:special}. Using \eqref{eq:x1star:x1vee} we return to \eqref{bs-Jg-2} with $i=j=1$, use the mean value theorem in $x_1$ and assumption~\eqref{item:ic:w0:x2:negative}, and deduce that 
\begin{align}
\nb_{1}\nb_1 \Jg(x_1^*,x_2,t) 
&\geq (t-\initial)\tfrac{1+\alpha}{2\eps} \nb_1^3 (w_0) (x_1^*,x_2) - \Cn (t-\initial) \mathsf{K}  \jap{\mathsf{B_6}}
\notag\\
&\geq (t-\initial)\tfrac{1+\alpha}{2\eps} \nb_1^3 (w_0) (x_1^\vee,x_2) 
- (t-\initial)\tfrac{1+\alpha}{2\eps} \tfrac{|x_1^*-x_1^\vee|}{\eps} \|\nb_1^4 w_0\|_{L^\infty}
- \Cn (t-\initial)  \mathsf{K}  \jap{\mathsf{B_6}}
\notag\\
&\geq (t-\initial)\tfrac{9(1+\alpha)}{20\eps} 
- (t-\initial)\tfrac{1+\alpha}{2\eps} \tfrac{2 \Cn \mathsf{K}  \jap{\mathsf{B_6}}  \eps^3}{\eps} \Cdata
- \Cn (t-\initial) \mathsf{K}  \jap{\mathsf{B_6}}
\notag\\
&\geq (t-\initial)\tfrac{2(1+\alpha)}{5\eps}\,,
\label{eq:Jg:11:lower}
\end{align}
upon taking $\eps$ sufficiently small. We note that the above lower bound vanishes as $t \to \initial$. We may obtain a matching bound for $\nb_1 \p_t \Jg(x_1^*,x_2,t)$ as follows. By Taylor's theorem in time, and using that $\Jg(\cdot,\initial) \equiv 1$, we obtain that
\begin{equation}
0 = \nb_1 \Jg(x,\initial) = \nb_1 \Jg(x,t) + (\initial-t) \nb_1 \p_t \Jg(x,t) +  \int_{\initial}^t \nb_1 \p_t^2 \Jg(x,t') (t' - \initial) {\rm d}t'  
\,.
\end{equation}
Evaluating the above expression at $x= (x_1^*(x_2,t),x_2)$, and using that $\Jg,_1(x_1^*(x_2,t),x_2,t) = 0$, we obtain that 
\begin{equation}
 \nb_1 \p_t \Jg(x_1^*(x_2,t),x_2,t) = \tfrac{1}{t-\initial} \int_{\initial}^t \nb_1 \p_t^2 \Jg(x_1^*(x_2,t),x_2,t') (t' - \initial) {\rm d}t' 
\end{equation}
and therefore
\begin{equation}
\sabs{ \nb_1 \p_t \Jg(x_1^*(x_2,t),x_2,t) } \leq  \tfrac{t-\initial}{2} \| \nb_1 \p_t^2 \Jg \|_{L^\infty_{x,t}}
\,.
\label{eq:I:cant:believe:I:need:this:shit:2}
\end{equation}
In order to bound the right side of the above identity, we appeal to  Lemma~\ref{lem:anisotropic:sobolev} and the available bounds at the $6^{th}$ derivative level which are given by the bootstrap~\eqref{bootstraps-Dnorm:Jg} and the initial data assumption. More precisely, we have $\|\nb^3 \Jg\|_{L^\infty_{x,\s}} \les \eps^{-\frac 12} \|\nb^5 \Jg(\cdot,\initial)\|_{L^2_x} + \eps^{-1} \|\nb^6 \Jg\|_{L^2_{x,\s}} \les \brak{\mathsf{B_J}}$. This  estimate gives a sub-optimal bound since  upon noting that $\nb_1 \p_t^2 = \eps^{-2} \nb_1  \nb_t^2$, we infer $\| \nb_1 \p_t^2 \Jg \|_{L^\infty_{x,t}} \les \eps^{-2} \brak{\mathsf{B_J}}$. To obtain the correct bound, which is sharper by a full power of $\eps$, we note that \eqref{eq:cant:believe:I:need:this:shit} implies
\begin{align}
\nb_1 \p_t^2  \Jg
&=
- \tfrac{1}{\eps} \nb_1(\nb_t V \p_2 \Jg)
- \nb_1(V \p_t \p_2\Jg)
+ \tfrac{1+\alpha}{2 }  \nb_1(\p_t + V \p_2) (\Jg \Wbn) 
\notag\\
&\qquad 
- \tfrac{1+\alpha}{2} \nb_1(V \p_2(\Jg \Wbn) )
+ \tfrac{1-\alpha}{2\eps} \nb_1 \nb_t   (\Jg \Zbn) 
\,.
\label{eq:I:cant:believe:I:need:this:shit:3}
\end{align}
Next, recall that the $\nb_1 \p_t$ differentiated version of \eqref{Jg-Wbn-1}  gives $\nb_1 (\p_t + V \p_2) (\Jg \Wbn) \circ \xi(x,t) = -  \eps \xi,_1 \p_2 (\p_t + V \p_2) (\Jg \Wbn) \circ \xi(x,t) + \eps (w_0),_{11}(x) \p_t \Iwn1(x,t) + (w_0),_1(x) \nb_1 \p_t \Iwn1(x,t) + \nb_1 \p_t \Iwn2(x,t)$. Therefore, using \eqref{eq:xi:nabla:xi}, \eqref{bootstraps-Dnorm:5}, the Sobolev embedding bound $\|\nb^2 (\Jg \Wbn)\|_{L^\infty_{x,\s}} \les \eps^{-\frac 12} \|\nb^4 (\Jg \Wbn)(\cdot,\initial)\|_{L^2_x} + \eps^{-1} \| \nb^5(\Jg \Wbn)\|_{L^2_{x,\s}} \les \eps^{-1} \brak{\mathsf{B_5}}$, and \eqref{eq:broncos:suck}, we deduce 
\begin{equation}
\| \nb_1(\p_t + V \p_2) (\Jg \Wbn) \|_{L^\infty_{x,t}} \les \eps^{-1} \brak{\mathsf{B_6}}\,.
\end{equation}
Using the above estimate, we return to \eqref{eq:I:cant:believe:I:need:this:shit:3}, and bound the remaining term by brute force using the bootstraps~\eqref{bootstraps} and the improved $\Zbn$ bounds from \eqref{eq:Jg:Zbn:D5:improve:a}, to obtain that
\begin{equation}
\| \nb_1 \p_t^2 \Jg \|_{L^\infty_{x,t}} \les \eps^{-1} \brak{\mathsf{B_6}} 
\,.
\end{equation}
Then, the above estimate and \eqref{eq:I:cant:believe:I:need:this:shit:2} imply that 
\begin{equation}
\sabs{ \nb_1 \p_t \Jg(x_1^*(x_2,t),x_2,t) }  \les  \tfrac{t-\initial}{\eps } \brak{\mathsf{B_6}}
\label{eq:Jg:1t:upper}
\end{equation}
holds pointwise for $(x_2,t) \in \hat{\mathcal{P}}$.

Using \eqref{eq:Jg:11:lower} and \eqref{eq:Jg:1t:upper}, we return to bound three terms in \eqref{eq:Qring:2s:hat}. By also appealing to \eqref{eq:Jgb:identity:0} we rewrite $\p_t \p_t \Jgb = \p_t \p_t (\Jgb - \Jg) + \p_t \p_t \Jg$,  the fact that $|\mathfrak{C}'| \leq 4$ and $\final-\medium = \frac{\eps}{50(1+\alpha)}$, and to the bounds \eqref{eq:Jg:Wbn:horse}, \eqref{Qd-lower-upper}, \eqref{Q-lower-upper}, \eqref{eq:Qb:bbq}, \eqref{eq:s:t:upper}, and \eqref{eq:d1d1w0}, we deduce
\begin{subequations}
\begin{align}
\sabs{\Qr_\s(x_2,\s)} &\leq \tfrac{\eps }{\Qd(x_2,\s)} \bigl( \bigl( \tfrac{8}{(\final-\medium)^2} + \tfrac{\Cn}{\eps}\bigl) + \tfrac{\Cn \s}{\eps}   \brak{\mathsf{B_6}}^2  \bigr)
\leq \tfrac{2 \cdot 5^2 \cdot 50^2 \Q}{\eps}
 \, \\
 \sabs{\Qr_2 (x_2,\s)}&\les\tfrac{\eps }{\Qd(x_2,\s)} \bigl(\tfrac{1}{\eps} + \tfrac{\s}{\eps}   \brak{\mathsf{B_6}}   \bigr) \les 1
 \, \\
 \sabs{\Qc(x_2,\s) - \Qr_\s (x_2,\s) }&\les \eps + \s  \les \eps 
 \,.
\end{align}
\end{subequations}
The above bounds establish \eqref{eq:Qrs:bbq}--\eqref{eq:Qc:bbq}, thereby concluding the proof of the lemma.
\end{proof}

\subsection{Properties of $\Jgb$, $\mathfrak{q}$, and the definition of the curve of pre-shocks}
\label{sec:x1star}

First, we show that the minimum of $\Jg(\cdot,x_2,t)$ is attained at a unique point as soon as $t>\initial$, justifying the definition \eqref{eq:x1star:def}.

\begin{lemma}[\bf The point $x_1^*(x_2,t)$ is uniquely defined]
\label{lem:x1*:unique}
Assume that the bootstraps~\eqref{bootstraps} hold, and assume that $\eps$ is sufficiently small with respect to $\alpha,\kappa_0$, and $\Cdata$. Then, 
for all $(x_2,t) \in \hat{\mathcal{P}}$ with $t>\initial$, there exists a unique $x_1^*(x_2,t) $ at which the minimum of $\Jg(\cdot,x_2,t)$ is attained.
\end{lemma}
\begin{proof}[Proof of Lemma~\ref{lem:x1*:unique}]
Recall cf.~\eqref{eq:Jgb:Jg:identity:3} and the discussion in the paragraph above that equation that $x_1^*$ can be equivalently defined as the location of the global minimum of $\Jgb(\cdot,x_2,t)$ or $\Jg(\cdot,x_2,t)$. 
Fix $(x_2,t) \in \hat{\mathcal{P}}$, with $t>\initial$. By \eqref{eq:ic:supp}, \eqref{bs-supp}, \eqref{eq:xi:nabla:xi}, and \eqref{Jg-circ-xi}, we have that the continuous map $x_1 \mapsto \Jg(x_1,x_2,t) - 1$ is supported in $\{|x_1|\leq (9 + \Csupp) \eps\} \subset \TT$. As such the minimum of this function is attained at at least one point. For any such point $x_1^*$, we have that \eqref{eq:d1w0:x1:star} holds, and thus the argument which as lead to the bound \eqref{eq:x1star:x1vee} holds true, yielding the estimate $|x_1^* - x_1^\vee| \les \eps^3$. Now assume that $x_1^*$ was not unique, rendering the existence of two such points $x_{1,a}^*$ and $x_{1,b}^*$. Then we must have $\Jg,_1(x_{1,a}^*,x_2,t) = \Jg,_1(x_{1,b}^*,x_2,t) = 0$, and by the mean value theorem there must exist $x_1^\sharp$ which lies in between $x_{1,a}^*$ and $x_{1,b}^*$, such that $\Jg,_{11}(x_{1}^\sharp,x_2,t)=0$. But note that  \eqref{eq:x1star:x1vee} implies $|x_1^\sharp- x_1^\vee| \les \eps^3$. Therefore, we may repeat the bounds in \eqref{eq:Jg:11:lower}, with $x_1^*$ replaced by $x_1^\sharp$, and deduce that $ \Jg,_{11}(x_{1}^\sharp,x_2,t) \geq (t-\initial) \tfrac{2(1+\alpha)}{5\eps^3} > 0$. This is a contradiction,  concluding the proof.
\end{proof}

\begin{lemma}[\bf The map $\mathfrak{q}$ is invertible]
\label{lem:q:invertible}
Assume that the bootstraps~\eqref{bootstraps} hold, and assume that $\eps$ is sufficiently small with respect to $\alpha,\kappa_0$, and $\Cdata$. Then, the map $\mathfrak{q}$ defined by \eqref{eq:t-to-s-transform:all} is invertible, with inverse $\mathfrak{q}^{-1}$ defined by \eqref{s-to-t-transform}.
\end{lemma}
\begin{proof}[Proof of Lemma~\ref{lem:q:invertible}]
Fix $x_2 \in \mathbb{T}$ and $\s \in [0,\eps)$. We need to show that the equation $\s - \mathfrak{q}(x_2,t) = 0$ has a unique solution $t\in [\initial,\final)$. The uniqueness part is easy: if two solutions $\initial \leq t_a < t_b< \final$ would exist, then we'd have $\mathfrak{q}(x_2,t_a) = \mathfrak{q}(x_2,t_b)$, and so by the mean value theorem  $\p_t \mathfrak{q}(x_2,t^\sharp)=0$ for some $t_a < t^\sharp < t_b$. But we have already shown earlier, see~\eqref{eq:Qcal:bbq:temp:3}, that $\p_t \mathfrak{q}(x_2,t^\sharp) \geq \frac{2(1+\alpha)}{5}  >0$, a contradiction. When $\s=0$, then $t=\initial$ clearly solves the desired equation since $\mathfrak{q}(x_2,\initial) = 0$. When $\s\in(0,\eps)$, the existence of $t$ follows by the intermediate value theorem, the continuity of $\mathfrak{q}$ with respect to time, and the bounds \eqref{bs-Jg-simple} and \eqref{eq:s:t:upper}. Indeed, \eqref{eq:s:t:upper} shows that $\mathfrak{q}(x_2,t) \leq 401(1+\alpha)  (t-\initial)$, so we can find $t_a>\initial$ such that $\s- \mathfrak{q}(x_2,t_a) >0$. For the other bound, we recall cf.~\eqref{eq:Jgb:identity:0} and the bootstrap \eqref{bs-Jg-simple} that for every $x\in \TT^2$ there exists $t_*(x)<\final$ such that $\Jgb(x,t_*(x))=0$. This implies that for every $x_2$ there exists $t_b\in [\initial,\final)$ with $\mathcal{J}(x_2,t_b) = 0$. Hence, $\s-\mathfrak{q}(x_2,t_b) = \s - \eps < 0$. Thus, by the intermediate value theorem there must exist a root $t \in (t_a,t_b)$ of $\s - \mathfrak{q}(x_2,t) = 0$ .
\end{proof}

\begin{definition}[\bf The curve of pre-shocks]
\label{def:pre-shock}
For all $x_2 \in \TT$ we define, extending the definition of $\mathfrak{q}^{-1}$ by continuity, 
\begin{equation}
\xstart := \lim_{\s \to \eps^{-}} x_1^*\bigl(x_2,\mathfrak{q}^{-1}(x_2,\s)\bigr) 
= x_1^*\bigl(x_2,\mathfrak{q}^{-1}(x_2,\eps)\bigr). 
\label{x1-star}
\end{equation}
he parametrized curve
\begin{equation*}
\Xi^* := \Bigl\{ \bigl(\xstart,x_2, \mathfrak{q}^{-1}(x_2,\eps) \bigr) \colon x_2\in\TT\Bigr\}
\end{equation*}
is called the {\em curve of pre-shocks}.    We define $t^*(x_2):= \mathfrak{q}^{-1}(x_2,\eps) $ so that 
$\Xi^* := \Bigl\{ \bigl(\xstart,x_2,  t^*(x_2) \bigr) \colon x_2\in\TT\Bigr\}$.
\end{definition}

\begin{proposition}[\bf Equivalent characterization of the curve of pre-shocks]
\label{prop:pre-shock}
The curve of pre-shocks $\Xi^*$ is precisely the intersection of the two-dimensional surfaces $\{\Jgb =0 \}$ and $\{\Jg,_1 = 0\}$.
\end{proposition}

\begin{proof}[Proof of Proposition~\ref{prop:pre-shock}]
We first establish the inclusion $\Xi^* \subseteq \{\Jgb =0 \} \cap \{\Jg,_1 = 0\}$. 
The fact that $\Jg,_1$ vanishes on $\Xi^*$ follows from the definition of $x_1^*$ (see~\eqref{eq:x1star:critical}), which  gives that $\Jg,_1(x_1^*(x_2,t),x_2,t))|_{t= \mathfrak{q}^{-1}(x_2,\s)} = 0$, and thus this equality also holds as $\s\to \eps$, by continuity. The fact that $\Jgb$ vanishes on $\Xi^*$ is a consequence of the definition of the map $\mathfrak{q}^{-1}$ (see \eqref{t-to-s-transform} and \eqref{s-to-t-transform}). Indeed, as $\s\to \eps$, we have that $\mathcal{J}(x_2,\mathfrak{q}^{-1}(x_2,\s)) \to 0$, which means via~\eqref{eq:fake:Jg:def} that $\Jgb(x_1^*(x_2,\mathfrak{q}^{-1}(x_2,\s)),x_2,\mathfrak{q}^{-1}(x_2,\s)) \to 0$ as $\s\to \eps$. By continuity of $\Jgb$ in the $x_1$ and $t$ entries, it follows that $\Jgb \equiv 0$ on $\Xi^*$.

The proof is completed once we establish the reverse inclusion, namely $\{\Jgb =0 \} \cap \{\Jg,_1 = 0\} \subseteq \Xi^*$. Let $(x_1,x_2,t) \in \{\Jgb =0 \} \cap \{\Jg,_1 = 0\}$. We need to show two things: $t = \mathfrak{q}^{-1}(x_2,\eps)$ and $x_1 = x_1^* (x_2,\mathfrak{q}^{-1}(x_2,\eps) )$. Since $\Jgb,_1 = \Jg,_1$, we have that $\Jgb(x_1,x_2,t) = \Jgb,_1(x_1,x_2,t) = 0$, and therefore the map $x_1 \mapsto \Jgb(x_1,x_2,t)$, with $(x_2,t)$ frozen, has a global minimum at $x_1$ (indeed, $\Jgb$ cannot attain strictly negative values in the closure of the spacetime considered here). By \eqref{eq:x1star:def} and the uniqueness statement established in Lemma~\ref{lem:x1*:unique}, it  follows that $x_1 = x_1^*(x_2,t)$. Moreover, by the definition~\eqref{eq:fake:Jg:def} it follows that $0 = \Jgb(x_1,x_2,t) = \mathcal{J}(x_2,t)$, which gives in light of \eqref{t-to-s-transform} that $t = \mathfrak{q}^{-1}(x_2,\eps)$. This concludes the proof.
\end{proof}

\subsection{Damping properties of $\Jg$ and $\mathcal{J}$}
\label{sec:Jg:properties}
In this section we record the properties of $\Jg$  and $\mathcal{J}$ that are most important to the analysis, especially to the energy estimates in subsequent sections.

\begin{lemma}[\bf Damping]
\label{lem:damping}
Assume that the bootstraps~\eqref{bootstraps} hold, and assume that $\eps$ is sufficiently small with respect to $\alpha,\kappa_0$, and $\Cdata$. Then, 
for all $(x,t) \in \mathcal{P}$, we have that
\begin{equation}
(\Jg \Wbn)(x,t) \leq - \tfrac{9}{10\eps}  +  \tfrac{13}{\eps} \Jg(x,t)\,.
 \label{eq:signed:Jg}
\end{equation}
\end{lemma}
\begin{proof}[Proof of Lemma~\ref{lem:damping}]
If $(x,t) \in \mathcal{P}$ is such that $\Jg \Wbn(x,t) \leq - \frac{9}{10\eps}$, then \eqref{eq:signed:Jg} holds automatically due to \eqref{bs-Jg-simple}. If on the other hand $\Jg \Wbn(x,t) \geq - \frac{9}{10\eps}$, then by \eqref{bs-Jg} we have $\Jg(x,t) \geq \frac{2}{25}$. In this case, \eqref{eq:broncos:eat:shit:20} and \eqref{eq:why:the:fuck:not:0} give 
\begin{equation*}
\tfrac{9}{10\eps} + \Jg \Wbn(x,t) \leq  \tfrac{9}{10\eps} + (w_0),_1(x) + \Cn \eps \leq  \tfrac{9}{10\eps} + \tfrac{1}{10\eps}  + \Cn \eps 
\leq \tfrac{1 + \Cn \eps^2}{\eps} \Jg  \Jg^{-1} \leq \tfrac{1 + \Cn \eps^2}{\eps} \Jg \tfrac{25}{2}
\,,
\end{equation*} 
thereby proving \eqref{eq:signed:Jg}.
\end{proof}

For the purpose of closing energy estimates, we will make use of the following crucial lemma:
\begin{lemma}[\bf Damping and anti-damping]
\label{lem:damping:anti:damping}
Assume that the bootstraps~\eqref{bootstraps} hold, and assume that $\eps$ is sufficiently small with respect to $\alpha,\kappa_0$, and $\Cdata$.
For all $(x,\s) \in \mathbb{T}^2 \times [0,\eps)$, we have
\begin{equation}
- \Jg (\Q \p_\s + V\p_2) \mathcal{J}^{\frac 32}   
+ 
 \mathcal{J}^{\frac 32}(\Q \p_\s + V\p_2)   \Jg
\geq \tfrac{1+\alpha}{8\eps} \mathcal{J}^\frac 12     \Jg 
\,.
\label{eq:fakeJg:LB}
\end{equation}
\end{lemma}
\begin{proof}[Proof of Lemma~\ref{lem:damping:anti:damping}]
We recall that by \eqref{eq:fake:Jg} we have 
\begin{equation}
(\Q \p_\s + V\p_2) \mathcal{J} = - \tfrac{\Q}{\eps} 
\,.
\label{evo-fakeJg}
\end{equation}
As such, using~\eqref{Jg-evo-s} we rewrite the left side of \eqref{eq:fakeJg:LB} as
\begin{equation}
- \Jg (\Q \p_\s + V\p_2) \mathcal{J}^{\frac 32}   
+ 
 \mathcal{J}^{\frac 32}(\Q \p_\s + V\p_2)   \Jg
=
\tfrac 32 \mathcal{J}^{\frac 12} \Jg  \tfrac{\Q}{\eps}
+ 
\mathcal{J}^{\frac 32} \bigl( \tfrac{1+\alpha}{2} \Jg \Wbn + \tfrac{1-\alpha}{2} \Jg \Zbn \bigr)
\,.
\label{eq:fakeJg:LB:x}
\end{equation}
Using the bootstraps~\eqref{bootstraps}, the coefficient bounds~\eqref{eq:Q:all:bbq}, and the fact that $0 \leq \mathcal{J} \leq \Jg$, we then bound from below the right side of \eqref{eq:fakeJg:LB:x} as
\begin{align}
\tfrac 32 \mathcal{J}^{\frac 12} \Jg  \tfrac{\Q}{\eps}
+ 
\mathcal{J}^{\frac 32} \bigl( \tfrac{1+\alpha}{2} \Jg \Wbn + \tfrac{1-\alpha}{2} \Jg \Zbn \bigr)
&=
\tfrac 32 \mathcal{J}^{\frac 12} \Jg  \tfrac{\Qd}{\eps}
+
\tfrac 32 \mathcal{J}^{\frac 12} \Jg  \tfrac{\Q - \Qd}{\eps}
+ 
\mathcal{J}^{\frac 32} \bigl( \tfrac{1+\alpha}{2} \Jg \Wbn + \tfrac{1-\alpha}{2} \Jg \Zbn \bigr)
\notag\\
&\geq 
\tfrac 32 \mathcal{J}^{\frac 12} \Jg  \tfrac{17(1+\alpha)}{40 \eps}
- 
\tfrac{9 \mathsf{C_V}}{2} \eps  \mathcal{J}^{\frac 12} \Jg  
- 
\mathcal{J}^{\frac 32} \bigl( \tfrac{1+\alpha}{2 \eps} + \Cn \bigr)
\notag\\
&\geq 
\mathcal{J}^{\frac 12} \Jg  \tfrac{51(1+\alpha)}{80 \eps}
- 
\mathcal{J}^{\frac 12}\Jg \bigl( \tfrac{1+\alpha}{2 \eps} + \Cn \bigr)
\notag\\
&\geq 
\mathcal{J}^{\frac 12} \Jg  \bigl( \tfrac{11(1+\alpha)}{80 \eps} - \Cn \bigr)
\notag\\
&\geq 
\mathcal{J}^{\frac 12} \Jg   \tfrac{1+\alpha}{8 \eps} 
\,.
\label{eq:fakeJg:LB:xx}
\end{align}
Combining~\eqref{eq:fakeJg:LB:x} and \eqref{eq:fakeJg:LB:xx} proves \eqref{eq:fakeJg:LB}.
\end{proof}

\subsection{Closure of the \eqref{bootstraps-Dnorm:5} bootstrap}
\label{sec:D5:bootstrap}
In this section we show that \eqref{bootstraps-Dnorm:5} follows from \eqref{bootstraps-Dnorm:6}.

We first consider a general function $F$ such that $\mathcal{J}^r \p_\s F \in L^2_{x,\s}$ for some $r \in \mathbb{R}$. Then, by the fundamental theorem of calculus in time we have that for $r' \in \mathbb{R}$ to be determined, pointwise in $\s$ it holds that 
\begin{equation}
\mathcal{J}^{r'}(\s) \| F(\cdot,\s)\|_{L^2_x}
\leq 
\mathcal{J}^{r'}(\s) \| F(\cdot,0)\|_{L^2_x}
+ \mathcal{J}^{r'}(\s) 
\int_0^{\s} \mathcal{J}^{r}(\s') \|\p_\s F(\cdot,\s')\|_{L^2_x} \mathcal{J}^{-r}(\s') {\rm d}\s'
\,,
\label{eq:killing:me:softly:1}
\end{equation}
and therefore we have
\begin{equation}
 \| \mathcal{J}^{r'} F \|_{L^2_{x,\s}} 
 \leq \| F(\cdot,0)\|_{L^2_x} \| \mathcal{J}^{r'}   \|_{L^2_\s} 
 +  \| \mathcal{J}^{r} \p_\s F \|_{L^2_{\s,x}} 
 \left( \int_0^\eps \mathcal{J}^{2r'}(\s)  \int_0^{\s} \mathcal{J}^{-2r}(\s') {\rm d}\s' {\rm d}\s\right)^{\frac 12}
 \,.
 \label{eq:killing:me:softly:2}
\end{equation}
At this point, we recall from the definition \eqref{t-to-s-transform} that we have the identity $\mathcal{J}(\s) = 1 - \frac{\s}{\eps}$. Therefore,  an explicit computation shows that \eqref{eq:killing:me:softly:2} yields
\begin{equation}
\mbox{if } r > \tfrac 12 \mbox{ and } r'  > r -1, \mbox{ then:} \qquad 
\| \mathcal{J}^{r'} F \|_{L^2_{x,\s}} 
 \leq  \tfrac{\eps^{\frac 12}}{\sqrt{2r'+1}} \| F(\cdot,0)\|_{L^2_x} 
 + \tfrac{\eps}{\sqrt{2(2r-1)(1+r'-r)}}   \| \mathcal{J}^{r} \p_\s F \|_{L^2_{x,\s}} 
 \,,
\label{eq:r:r':L2:time}
\end{equation}
and 
\begin{equation}
\mbox{if } r < \tfrac 12 \mbox{ and } r'  > - \tfrac 12, \mbox{ then:} \qquad 
\| \mathcal{J}^{r'} F \|_{L^2_{x,\s}} 
 \leq  \tfrac{\eps^{\frac 12}}{\sqrt{2r'+1}} \| F(\cdot,0)\|_{L^2_x} 
 + \tfrac{\eps}{\sqrt{(2r'+1)(1-2r)}}   \| \mathcal{J}^{r} \p_\s F \|_{L^2_{x,\s}} 
 \,.
\label{eq:r:r':L2:time:new}
\end{equation}
Returning to \eqref{eq:killing:me:softly:1}, we also note that upon taking $r'=0$, we have  
\begin{equation}
\mbox{if } r < \tfrac 12, \mbox{ then:} \qquad  
\| F \|_{L^\infty_\s L^2_{x}} 
 \leq  \| F(\cdot,0)\|_{L^2_x} 
 +  \tfrac{\eps^{\frac 12}}{\sqrt{1-2r}}   \| \mathcal{J}^{r} \p_\s F \|_{L^2_{x,\s}} 
 \,.
\label{eq:r:Linfty:time}
\end{equation}
Lastly we note that if only a bound on $\mathcal{J}^r \p_\s F \in L^2_{x,\s}$ is available with $r>\frac 12$, then \eqref{eq:r:Linfty:time} is not available. In this situation, we require knowledge of $\mathcal{J}^{r-\frac 12} \Jgh \p_\s F \in L^2_{x,\s}$ and conclude $\Jgh F \in L^\infty_\s L^2_x$. The argument is as follows. Using \eqref{Jg-evo-s}, \eqref{bootstraps}, \eqref{eq:Q:all:bbq}, and \eqref{eq:signed:Jg}, we conclude
\begin{align}
 \tfrac{d}{d\s} \| \Jgh F \|_{L^2_x}^2 
 &= \int \tfrac{1}{\Q} |F|^2 \Q \p_\s \Jg + 2 \int \Jg F \p_\s F 
 \notag\\
 &= \int \tfrac{1}{\Q} |F|^2 \bigl( - V \p_2\Jg + \tfrac{1+\alpha}{2} \Jg \Wbn + \tfrac{1-\alpha}{2} \Jg \Zbn \bigr) + 2 \int \Jg F \p_\s F
 \notag\\
 &\leq \int \tfrac{1}{\Q} |F|^2 \bigl(5 \eps (1+\alpha) \mathsf{C_V} - \tfrac{9(1+\alpha)}{20 \eps} + \tfrac{13(1+\alpha)}{2\eps} \Jg + \tfrac{3(1-\alpha)}{4} \mathsf{C}_{\Zbn} \bigr) + 2 \|\Jgh F\|_{L^2_x} \|\Jgh \p_\s F\|_{L^2_x}
 \notag\\
 &\leq \tfrac{18}{\eps} \| \Jgh F \|_{L^2_x}^2  + 2 \|\Jgh F\|_{L^2_x} \|\Jgh \p_\s F\|_{L^2_x}
 \,,
 \label{eq:r:Linfty:time:1}
\end{align}
upon taking $\eps$ to be sufficiently small. Integrating the above inequality in time, and recalling that $\mathcal{J}(\s) = 1- \frac{\s}{\eps}$, we are lead to conclude that 
\begin{align}
\mbox{if } \tfrac 12 < r < 1, \mbox{ then:} \qquad 
\sup_{\s \in [0,\eps]}  \| \Jgh F (\cdot,\s) \|_{L^2_x} 
&\leq 
e^9 \| \Jgh F (\cdot,0) \|_{L^2_x}  
+ e^9 \int_0^\eps \|\Jgh \p_\s  F(\cdot,\s)\|_{L^2_x} {\rm d}\s
\notag\\
&=
e^9 \| \Jgh F (\cdot,0) \|_{L^2_x}  
+ e^9 \int_0^\eps \|\mathcal{J}^{r-\frac 12} \Jgh \p_\s  F(\cdot,\s)\|_{L^2_x} \mathcal{J}^{-r+\frac 12} (\s) {\rm d}\s
\notag\\
&\leq 
e^9 \| \Jgh F (\cdot,0) \|_{L^2_x} 
+ \tfrac{\eps^{\frac 12} e^9}{\sqrt{2(1-r)}} \|\mathcal{J}^{r-\frac 12} \Jgh \p_\s  F(\cdot,\s)\|_{L^2_{x,\s}} 
\,.
\label{eq:r:Linfty:time:2}
\end{align}

Having established the bounds \eqref{eq:r:r':L2:time}, \eqref{eq:r:Linfty:time}, and \eqref{eq:r:Linfty:time:2}, we show that \eqref{bootstraps-Dnorm:5} follows from \eqref{bootstraps-Dnorm:6}, assuming $\mathsf{B}_5$ is sufficiently large with respect to $\mathsf{B}_6$. Indeed, from \eqref{eq:r:r':L2:time} with $r'=0$ and $r=\frac 34$,   using \eqref{table:derivatives}, 
recalling that $\p_\s = \frac{1}{\eps \Qd} \nbs_\s$, and that $\Qd^{-1}$ is bounded according to \eqref{Qd-lower-upper}, we deduce 
\begin{align*}
\widetilde{\mathcal{D}}_{5,\nnn} (\eps)
&\leq   \Cdata
+ \tfrac{5}{1+\alpha}  \widetilde{\mathcal{D}}_{6,\nnn}(\eps)\,,
\\
\widetilde{\mathcal{D}}_{5,\ttt} (\eps)
&\leq  \eps \Cdata
+ \tfrac{5}{1+\alpha}  \widetilde{\mathcal{D}}_{6,\ttt}(\eps)
\,.
\end{align*}
From the above bound and the definitions \eqref{eq:tilde:D6} and \eqref{eq:tilde:D5}, it follows that 
\begin{equation}
\widetilde{\mathcal{D}}^2_5(\eps)  
\leq 2 \Cdatatwo + \tfrac{50}{(1+\alpha)^2}  \widetilde{\mathcal{D}}_{6,\nnn}^2(\eps)
+ (\mathsf{K} \eps)^{-2} \bigl(2 \eps^2 \Cdatatwo + \tfrac{50}{(1+\alpha)^2} \widetilde{\mathcal{D}}_{6,\ttt}^2(\eps)\bigr)
\leq 4 \Cdatatwo  + \tfrac{50}{(1+\alpha)^2} \widetilde{\mathcal{D}}^2_6(\eps) 
\,,
\label{eq:D5:D6:dominate}
\end{equation}
since $\mathsf{K}\geq 1$.
Similarly, by appealing to \eqref{eq:r:Linfty:time:2} with $r= \frac 34$, using \eqref{table:derivatives} and the fact that $\Jg(\cdot,0) = 1$, we deduce that 
\begin{align*}
\sup_{\s\in [0,\eps]} \widetilde{\mathcal{E}}_{5,\nnn}(\s)
&\leq e^9 \Cdata \eps^{-\frac 12}  + (\tfrac{2}{\eps})^{\frac 12} e^9 \tfrac{5}{1+\alpha} \widetilde{\mathcal{D}}_{6,\nnn}(\eps)
\\
\sup_{\s\in [0,\eps]} \widetilde{\mathcal{E}}_{5,\ttt}(\s)
&\leq e^9 \Cdata \eps^{\frac 12}  + (\tfrac{2}{\eps})^{\frac 12} e^9 \tfrac{5}{1+\alpha} \widetilde{\mathcal{D}}_{6,\ttt}(\eps)
\end{align*}
and therefore, upon recalling \eqref{eq:tilde:E6} and \eqref{eq:tilde:E5}, we obtain
\begin{align}
\eps \sup_{\s\in [0,\eps]} \widetilde{\mathcal{E}}_{5}^2(\s)
&\leq 2 e^{18} \Cdatatwo + 4 e^{18}\tfrac{5^2}{(1+\alpha)^2} \widetilde{\mathcal{D}}_{6,\nnn}^2(\eps)
+
(\mathsf{K}\eps)^{-2} \bigl( 2 e^{18} \Cdatatwo \eps^2 + 4   e^{18}\tfrac{5^2}{(1+\alpha)^2} \widetilde{\mathcal{D}}_{6,\ttt}^2(\eps)\big)
\notag\\
&\leq 4 e^{18} \Cdatatwo+ 4   e^{18}\tfrac{5^2}{(1+\alpha)^2}\widetilde{\mathcal{D}}^2_6(\eps) 
\,.
\label{eq:E5:D6:dominate}
\end{align}
From \eqref{eq:D5:D6:dominate} and \eqref{eq:E5:D6:dominate} we thus obtain
\begin{equation}
\eps^{\frac 12} \sup_{\s\in [0,\eps]} \widetilde{\mathcal{E}}_{5}(\s) 
+ \widetilde{\mathcal{D}}_5(\eps)  
\leq 2 (1+ e^9) \Cdata + \tfrac{10}{1+\alpha} (1+ e^9) \widetilde{\mathcal{D}}_6(\eps)
\end{equation}
and so the bootstrap \eqref{bootstraps-Dnorm:6} implies that \eqref{bootstraps-Dnorm:5} holds (with a strict inequality) as soon as 
\begin{equation}
  2 (1+ e^9) \Cdata + \tfrac{10}{1+\alpha} (1+ e^9) \mathsf{B}_6 =: \mathsf{B}_5
  \,.
  \label{eq:B5B6:relation}
\end{equation}

\begin{remark}[\bf $\mathsf{B}_5$ and $\mathsf{B}_6$ are proportional]
\label{rem:B5:B6}
 Note that since we will choose $\mathsf{B}_6 \geq \Cdata$ (see~\eqref{eq:B6:choice:1}), the relation \eqref{eq:B5B6:relation} implies   
\begin{equation}
\tfrac{10}{1+\alpha} (1+ e^9) \mathsf{B_6} \leq  \mathsf{B}_5 \leq 12 (1+ e^9) \mathsf{B}_6
\,.
\end{equation}
As such, any upper bound of the type $A \les \mathsf{B}_5$ may be written as $A \les \mathsf{B}_6$, upon changing the implicit constant.
\end{remark}

\section{Bounds for the geometry, sound speed, and the tangential reparameterization velocity}
\label{sec:geometry:sound:ALE}
The purpose of this section is to establish the following bounds, which are then subsequently used throughout the paper. Additionally, we close the bootstrap assumptions for the $\nbs^6$ level bounds on the geometry,~\eqref{bootstraps-Dnorm:Jg}--\eqref{bootstraps-Dnorm:h2}.
\begin{proposition}[Bounds for the geometry, sound speed, and ALE velocity]
\label{prop:geometry}
Assume   the bootstrap assumptions~\eqref{bootstraps} hold, and that $\eps$ is taken to be sufficiently small to ensure $\eps^{\frac 12} \bigl( \brak{\mathsf{B_J}} +  \brak{\mathsf{B_h}} +     \brak{\mathsf{B_6}}\bigr)  \leq 1$.  
Then, assuming $\eps$ is sufficiently small with respect to $\alpha,\kappa_0,$ and $\Cdata$, we have
\begin{subequations}
\label{geometry-bounds-new}
\begin{align} 
\eps^{\frac 12} \snorm{ \mathcal{J}^{\! \frac 14} \nbs^6 \Jg }_{L^\infty_\s L^2_x}
+ \snorm{  \mathcal{J}^{-\! \frac 14}\nbs^6 \Jg }_{L^2_{x,\s}} 
&\les  \eps  \brak{\mathsf{B}_6}   \,,   \label{D6JgEnergy:new}
\\
 \snorm{ \nbs^5 \Jg }_{L^\infty_\s L^2_x} 
 &\les  
 \eps^{\frac 12}  \brak{\mathsf{B}_6}   \,,   \label{D5JgEnergy} 
 \\
\eps^{\frac 12} \| \mathcal{J}^{\! \frac 14} \nbs^{6} \nbs_2h  \|^2_{L^\infty_\s L^2_x}
+ \|\mathcal{J}^{-\! \frac 14} \nbs^{6}  \nbs_2h \|_{L^2_{x,\s}}^2
&\les 
 \mathsf{K}  \eps^2 \brak{\mathsf{B}_6}
\,,
\label{D6h2Energy:new}
\\
\eps^{\frac 12} \| \mathcal{J}^{\frac 14} \nbs^6 \nbs_1 h\|_{L^\infty_{\s} L^2_x}
+ \|\mathcal{J}^{-\frac 14}  \nbs^{6}  \nbs_1 h \|_{L^2_{x,\s}}
&\les  \eps^2  \brak{\mathsf{B}_6} 
\,,
\label{D6h1Energy:new}
\\
 \snorm{ \nbs^5 \nbs_2 h(\cdot,\s) }_{L^\infty_\s L^2_x} 
 &\les  
  \mathsf{K}  \eps^{\frac 32}  \brak{\mathsf{B}_6}
 \,,   \label{D5h2Energy} 
 \\
 \eps^{\frac 12} \|\mathcal{J}^{\!\frac 14}  \nbs^{6} g  \|_{L^\infty_{\s} L^2_x}
+ \| \mathcal{J}^{-\!\frac 14} \nbs^{6}  g \|_{L^2_{x,\s}}
&\les \mathsf{K}^2 \eps^3 \brak{\mathsf{B_6}}^2
\,, \label{D6gbound} 
\\
 {\textstyle\sum}_{3\leq |\gamma|\leq 6} 
 \|\mathcal{J}^{-\frac 14} \bigl(\nbs^{|\gamma|}  \nn + g^{-1}  \tt \nbs^{|\gamma|} \nbs_2 h\bigr) \|_{{L^2_{x,\s}} }
+
\|\mathcal{J}^{-\frac 14} \bigl( \nbs^{|\gamma|}  \tt - g^{-1}  \nn \nbs^{|\gamma|} \nbs_2 h \bigr) \|_{{L^2_{x,\s}} }
&\les  \mathsf{K}  \eps^3 \brak{\mathsf{B}_6}  \,,
\label{D6n-bound:a:new}
\\
\| \mathcal{J}^{-\frac 14}  \nbs^6  \nn  \|_{{L^2_{x,\s}} }
+
\|\mathcal{J}^{-\frac 14}  \nbs^6 \tt \|_{{L^2_{x,\s}} }
&\les  \mathsf{K}  \eps^2 \brak{\mathsf{B}_6}   \,,
\label{D6n-bound:b:new}
\\
\|\mathcal{J}^{\frac 14} \nbs^6  \nn  \|_{L^\infty_\s L^2_x}
+
\|\mathcal{J}^{\frac 14}  \nbs^6 \tt \|_{L^\infty_\s L^2_x}
&\les  \mathsf{K}  \eps^{\frac 32} \brak{\mathsf{B}_6}   \,,
\label{D6n-bound:b:new:bdd}
\\
\snorm{\nbs^6 \Sigma^{\pm1}}_{L^2_{x,\s} } 
&\les 
\eps \brak{\mathsf{B}_6} \,,
\label{eq:Sigma:H6:new}
\\
\snorm{\Jgh \nbs^6 \Sigma^{\pm1}}_{L^\infty_\s L^2_{x} } 
&\les 
\eps^{\frac 12} \brak{\mathsf{B}_6} \,,
\label{eq:Sigma:H6:new:bdd}
\\
\snorm{ \nbs^6 V}_{L^2_{x,\s}} 
&\les  \mathsf{K}  \eps^2   \brak{\mathsf{B}_6} \,,
\label{eq:V:H6:new}
\\
\snorm{\Jgh \nbs^6 V}_{L^\infty_\s L^2_{x}} 
&\les  \mathsf{K}  \eps^{\frac 32} \brak{\mathsf{B}_6} \,,
\label{eq:V:H6:new:bdd}
\end{align} 
where the implicit constants in all the above inequalities depend only on $\alpha$, $\kappa_0$, and $\Cdata$. 
\end{subequations}
\end{proposition}
\begin{proof}[Proof of Proposition~\ref{prop:geometry}]
The proof consists in combining the bounds contained in Lemmas~\ref{lem:Sigma:V},~\ref{lem:D5-Jg},~\ref{lem:D5-h2},~\ref{lem:D6:V:Sigma:bdd:time}, in Remarks~~\ref{rem:B5:B6},~\ref{rem:D5:Jg:energy},~\ref{rem:D5:h2:energy}, \ref{rem:D6:g},~\ref{rem:D6:D1:h},~\ref{rem:D5:Sigma:V}, and 
Corollary~\ref{cor:D5n-bound}, which are all proven below. 
\end{proof}

One immediate consequence of the above proposition (we recall that the implicit constants therein only depend on $\alpha,\kappa_0,$ and $\Cdata$), and of the bound \eqref{eq:Jgb:less:than:1} is that the bootstraps~\eqref{bootstraps-Dnorm:h2}--\eqref{bootstraps-Dnorm:Jg} are closed.  
\begin{corollary}
\label{cor:Bj:Bh}
Assume that $\mathsf{B_J}$ and $\mathsf{B_h}$ are sufficiently large with respect to $\mathsf{B_6}$ and $\mathsf{K}$. More precisely, define
\begin{subequations}
\begin{align}
6 C_{\eqref{D6JgEnergy:new}}^{\frac 12} \brak{\mathsf{B_6}}  
&=:   \mathsf{B_J}
\label{eq:B:Jg:cond} \\
6 \bigl( C_{\eqref{D6h2Energy:new}}^{\frac 12}  \mathsf{K}  + C_{\eqref{D6h1Energy:new}}^{\frac 12} \bigr) \brak{\mathsf{B_6}}  
&=:   \mathsf{B_h}
\label{eq:B:h2:cond}
\,.
\end{align}
\end{subequations}
Then the bounds~\eqref{bootstraps-Dnorm:h2}--\eqref{bootstraps-Dnorm:Jg} hold with constants $\frac 12 \mathsf{B_h}$ and $\frac 12  \mathsf{B_J}$ respectively, thereby closing these bootstraps.
\end{corollary}

The following inequalities will be used  several times throughout the proof, and  provide bounds for the norms of $V$ and $\Sigma^{\pm 1}$ in terms 
of norms of $\nbs_2 h$ and $\Jg$, and of  $\brak{\mathsf{B_6}}$.
\begin{lemma}
\label{lem:Sigma:V}
Under the same assumptions as Proposition~\ref{prop:geometry},
the rescaled sound speed $\Sigma$ and the tangential transport velocity $V$ satisfy
\begin{subequations}
\label{D6-Sigma-V}
\begin{align} 
 \snorm{\nbs^6 \Sigma}_{L^2_{x,\s}} 
 & \les \eps \brak{\widetilde{\mathcal{D}}_{5}}
+\eps \bigl(\eps +   \snorm{\nbs^6 \nbs_2 h}_{L^2_{x,\s}}  + \eps \snorm{\nbs^5 \Jg}_{L^2_{x,\s}} \bigr)  \bigl( \brak{\mathsf{B}_6}
 +   \snorm{\nbs^6 \nbs_2 h}_{L^2_{x,\s}}   
 + \eps \snorm{\nbs^5 \Jg}_{L^2_{x,\s}} \bigr) \,,  
 \label{eq:Sigma:H5} \\
  \snorm{\nbs^6 V}_{L^2_{x,\s}}
 & \les \eps^2 \brak{\mathsf{B}_6}
 +  \snorm{\nbs^6 \nbs_2 h}_{L^2_{x,\s}}   
 + \eps \snorm{\nbs^5 \Jg}_{L^2_{x,\s}}  \,.
 \label{eq:V:H6}
 \end{align} 
\end{subequations}
\end{lemma}
\begin{proof}[Proof of Lemma~\ref{lem:Sigma:V}]
We prove \eqref{eq:Sigma:H5} only for $\Sigma$. Comparing \eqref{Sigma0i-ALE-s} with $\beta = \frac 12$ and \eqref{Sigma0-ALE-s}, and appealing to the bootstrap \eqref{bs-Sigma}, it is clear that the bounds for $\Sigma^{-1}$ are proven in exactly the same way.

In the case that $\nbs^6=\nbs_2\nbs^5$, we apply $\nbs^5$ to \eqref{p2-Sigma-s}, and apply \eqref{table:derivatives}, \eqref{bootstraps}, and \eqref{eq:Lynch:1}, to obtain
\begin{equation} 
\snorm{\nbs^5 \nbs_2 \Sigma}_{L^2_{x,\s}} 
\le
\tfrac{1}{2} \snorm{  \nbs^5\big( g^{\frac{1}{2}}  (\Wbt-\Zbt) \big)}_{L^2_{x,\s}}  
 \les \widetilde{\mathcal{D}}_{5,\ttt} + \eps  \snorm{\nbs^5 \nbs_2 h}_{L^2_{x,\s}}   + \eps
\les \eps \widetilde{\mathcal{D}}_{5}
+ \eps  \snorm{\nbs^5 \nbs_2 h}_{L^2_{x,\s}}   
\,.
\label{D2D5Sigma}
\end{equation} 
Similarly, in the case that $\nbs^6 = \nbs_1\nbs^5$, we let $\nbs^5$ act upon to \eqref{p1-Sigma-s} and  
and apply the inequalities \eqref{table:derivatives}, \eqref{bootstraps},   and \eqref{eq:Lynch:1},  to obtain that
\begin{align}
\snorm{\nbs^5 \nbs_1 \Sigma}_{L^2_{x,\s}} 
&\les 
\eps \snorm{\nbs^5 (\Jg \Wbn,\Jg \Zbn)}_{L^2_{x,\s}} 
+ \eps \snorm{ \nbs^5(\Jg \nbs_2 h (\Wbt -\Zbt) )}_{L^2_{x,\s}} 
\notag\\
& \les
\eps \widetilde{\mathcal{D}}^2_{5,\nnn}
+
\eps^2 \widetilde{\mathcal{D}}^2_{5,\ttt}
 + \eps  \snorm{\nbs^5 \nbs_2 h}_{L^2_{x,\s}}   
 + \eps^2 \snorm{\nbs^5 \Jg}_{L^2_{x,\s}} 
 + \eps^2  
 \notag\\
 &\les \eps \widetilde{\mathcal{D}}_{5}
 + \eps  \snorm{\nbs^5 \nbs_2 h}_{L^2_{x,\s}}   
 + \eps^2 \snorm{\nbs^5 \Jg}_{L^2_{x,\s}} 
 \label{D1D5Sigma}
\end{align} 
Finally, in the case that $\nbs^6 = \nbs_\s^6$, 
from \eqref{the-time-der}, \eqref{Sigma0-ALE-s}, we have the identity
\begin{align*}
\nbs_\s^6  \Sigma  
& = - \eps \nbs_2 \Sigma \nbs_\s^5 V  - \eps V \nbs_\s^5  \nbs_2 \Sigma 
 -  \alpha \eps (\Zbn+\Abt) \nbs^5\Sigma - \alpha  \eps  \Sigma  (\nbs_\s^5\Zbn + \nbs_\s^5\Abt) 
 \notag \\
 & \qquad
 - \eps \doublecom{ \nbs_\s^5  ,V, \nbs_2 \Sigma } -  \alpha \eps \doublecom{\nbs_\s^5, \Zbn+\Abt ,\Sigma } \,.
\end{align*}
Using  \eqref{table:derivatives}, \eqref{eq:tilde:D5}, \eqref{bootstraps}, \eqref{D2D5Sigma}, \eqref{D1D5Sigma}, 
 \eqref{eq:Jg:Zbn:D5:improve:a},   \eqref{eq:x1:Poincare}, \eqref{eq:Sobolev}, and \eqref{eq:Lynch:2} we find that
\begin{align}
\snorm{\nbs_\s^6 \Sigma}_{L^2_{x,\s}}
& \les \eps \snorm{\nbs_\s^5 V}_{L^2_{x,\s}}
+\eps^2 \snorm{\nbs^5 \nbs_2 \Sigma}_{L^2_{x,\s}} 
+\eps \snorm{\nbs^5 \nbs_1 \Sigma}_{L^2_{x,\s}} 
+ \eps \snorm{\nbs_\s^5 (\Zbn, \Abt)}_{L^2_{x,\s}}
\notag \\
&
\quad + \eps \bigl( \eps + \| \nbs_\s^4 V\|_{L^2_{x,\s}}\bigr) 
\bigl( \eps^{-\frac 12}  \| \nbs_\s \nbs_2 \Sigma(\cdot,0)\|_{L^2_{x}}
+ \eps^{-1} \| \nbs^2 \nbs_1 \nbs_s \nbs_2 \Sigma\|_{L^2_{x,\s}} \bigr)
+\eps \bigl( 1 + \| \nbs_\s^4 (\Zbn,\Abt)\|_{L^2_{x,\s}}  \bigr)
\notag\\
&\les 
\eps \snorm{\nbs^5 \nbs_1 V}_{L^2_{x,\s}}
+\eps^2 \snorm{\nbs^5 \nbs_2 \Sigma}_{L^2_{x,\s}} 
+\eps \snorm{\nbs^5 \nbs_1 \Sigma}_{L^2_{x,\s}} 
+ \| \nbs_\s^4 V\|_{L^2_{x,\s}}\| \nbs^4 \nbs_1 \Sigma\|_{L^2_{x,\s}}
+ \eps^2 \brak{\mathsf{B_6}} + \eps
\notag\\
&\les 
\eps
+ \eps \snorm{\nbs^5 \nbs_1 V}_{L^2_{x,\s}}
+ \eps
\bigl( \eps + \snorm{\nbs^5 \nbs_1 V}_{L^2_{x,\s}} \bigr)
\bigl( \brak{\mathsf{B_6}}
 +   \snorm{\nbs^5 \nbs_2 h}_{L^2_{x,\s}}   
 + \eps \snorm{\nbs^5 \Jg}_{L^2_{x,\s}} \bigr)
 \,.
  \label{DsD5Sigma}
\end{align} 
In light of \eqref{D2D5Sigma}, \eqref{D1D5Sigma}, and \eqref{DsD5Sigma}, the $\nbs^6 \Sigma$ bound is completed, once a bound for $\nbs^5 \nbs_1 V$ is available.

The bound for $\nbs^6 V$ is obtained in a similar fashion. If $\nbs^6 =\nbs_2\nbs^5$, then using the identity \eqref{eq:nabla:V:2-s} together with
\eqref{table:derivatives}, \eqref{bootstraps}, \eqref{D1D5Sigma}, \eqref{eq:x1:Poincare}, and \eqref{eq:Lynch:1} shows that
\begin{align}
\snorm{\nbs^5 \nbs_2 V}_{L^2_{x,\s}}
&\les 
\snorm{\nbs^5 (\Sigma g^{-\frac 32}  \nbs_2^2 h)}_{L^2_{x,\s}}
+ \snorm{\nbs^5 \Abt}_{L^2_{x,\s}}
+ \snorm{\nbs^5\big( \nbs_2h ( \tfrac{1+\alpha}{2} \Wbt + \tfrac{1-\alpha}{2} \Zbt)\big)}_{L^2_{x,\s}}
\notag \\
&\les \eps \snorm{\nbs^5 \Sigma}_{L^2_{x,\s}}
+ \snorm{\nbs^6 \nbs_2 h}_{L^2_{x,\s}}
+ \eps^2 
+ \eps \widetilde{\mathcal{D}}_{5,\ttt}
\notag\\
&\les  \eps^2 \brak{\widetilde{\mathcal{D}}_{5}}
 +  \snorm{\nbs^6 \nbs_2 h}_{L^2_{x,\s}}   
 + \eps^3 \snorm{\nbs^5 \Jg}_{L^2_{x,\s}} 
\,.
\label{D2D5V}
\end{align}
Next, if $\nbs^6 =\nbs_1\nbs^5$ we apply $\nbs^5$ to  \eqref{eq:nabla:V:1-s}, and use \eqref{table:derivatives}, \eqref{bootstraps}, \eqref{eq:Jg:Abn:D5:improve:c}, \eqref{eq:Jg:Zbn:D5:improve}, \eqref{eq:x1:Poincare}, and \eqref{eq:Lynch:1}
 to obtain
\begin{align} 
\snorm{\nbs_1 \nbs^5 V}_{L^2_{x,\s}} 
&\les \snorm{\nbs^5 (\Sigma g^{-\frac 32}\nbs_1\nbs_2 h )}_{L^2_{x,\s}}
+\eps \snorm{\nbs^5 (g^{-\frac 12} \Jg \Abn)}_{L^2_{x,\s}}
+\eps \snorm{\nbs^5 (g^{-\frac 12}\nbs_2 h ( \tfrac{1+\alpha}{2}\Jg \Wbn + \tfrac{1-\alpha}{2} \Jg \Zbn) )}_{L^2_{x,\s}}
\notag\\
&\quad 
+\eps \snorm{\nbs^5 (g^{-\frac 12}\nbs_2 h \Jg \Abt) }_{L^2_{x,\s}}
+\eps \snorm{\nbs^5 (g^{-\frac 12}(\nbs_2 h)^2 \Jg ( \tfrac{1+\alpha}{2} \Wbt + \tfrac{1-\alpha}{2}  \Zbt) )}_{L^2_{x,\s}}
\notag\\
&\les \eps \snorm{\nbs^5\Sigma}_{L^2_{x,\s}} 
+  \snorm{\nbs^6\nbs_2 h}_{L^2_{x,\s}}
+\eps \snorm{\nbs^5 \Jg }_{L^2_{x,\s}}
+ \eps^2
\notag\\
&\quad
+\eps \snorm{\nbs^5 \Abn}_{L^2_{x,\s}}
+ \eps^2 \snorm{\nbs^5 (\Jg \Wbn,\Jg \Zbn)}_{L^2_{x,\s}}
+ \eps^2 \snorm{\nbs^5 (\Wbt,\Zbt,\Abt)}_{L^2_{x,\s}}
\notag\\
&\les
\eps^2 \brak{\mathsf{B_6}}
+  \snorm{\nbs^6\nbs_2 h}_{L^2_{x,\s}}
+\eps \snorm{\nbs^5 \Jg }_{L^2_{x,\s}}
\label{D1D5V}
\end{align} 
Lastly, if $\nbs^6 =\nbs_\s^6$, we
transform \eqref{V-evo} into $(x,\s)$  coordinates, apply $\nbs_\s^5$ and find that
\begin{equation*} 
\snorm{\nbs_s^6 V}_{L^2_{x,\s}}
\les \eps \snorm{\nbs_\s^5 (V\nbs_2 V)}_{L^2_{x,\s}}
+\eps  \snorm{\nbs_\s^5  \Bigl( \Sigma g^{-\frac 12}  \Bigl(\tfrac{2+ \alpha }{2} \Wbt - \tfrac{\alpha }{2} \Zbt- \Abn - \nbs_2 h  \bigl(  \alpha   \Abt -  (1-\alpha)  \Zbn \bigr) \Bigr)   \Bigr)}_{L^2_{x,\s}} 
 \,.
\end{equation*} 
By appealing to \eqref{table:derivatives}, \eqref{bootstraps}, 
\eqref{D1D5Sigma}, \eqref{D2D5V}, \eqref{D1D5V}, 
\eqref{eq:Jg:Abn:D5:improve:c}, \eqref{eq:Jg:Zbn:D5:improve}, \eqref{eq:x1:Poincare}, and \eqref{eq:Lynch:1},
we obtain
\begin{align} 
\snorm{\nbs_s^6 V}_{L^2_{x,\s}}
&\les \eps^2 \snorm{\nbs_\s^5 (\nbs_1 V,\nbs_2 V)}_{L^2_{x,\s}}
+\eps  \snorm{\nbs_\s^5  \nbs_1 \Sigma }_{L^2_{x,\s}} 
+\eps^2  \snorm{\nbs_\s^5 \nbs_2 h }_{L^2_{x,\s}}
+ \eps^2 
\notag\\
&\quad 
+\eps  \snorm{\nbs_\s^5 ( \Wbt,\Zbt,\Abt, \Abn ,\Zbn )  }_{L^2_{x,\s}} 
+\eps \snorm{\nbs_\s^5 \nbs_2 h }_{L^2_{x,\s}}
\notag\\
&\les 
\eps^2 \brak{\mathsf{B_6}}
+ \eps^2 \snorm{\nbs^6\nbs_2 h}_{L^2_{x,\s}}
+\eps^3 \snorm{\nbs^5 \Jg }_{L^2_{x,\s}} 
\,.
 \label{DsD5V}
\end{align} 
In the last inequality we have used \eqref{the-time-der} and \eqref{p2h-evo-s} to bound
\begin{align*}
\snorm{\nbs_\s^5 \nbs_2 h }_{L^2_{x,\s}} 
&\les  \eps \snorm{\nbs_\s^4 \bigl(g (\tfrac{1+\alpha}{2} \Wbt + \tfrac{1-\alpha}{2} \Zbt)\bigr)}_{L^2_{x,\s}}
+ \eps \snorm{\nbs_\s^4 \bigl(V \nbs_2^2 h\bigr)}_{L^2_{x,\s}}
\notag\\
&\les \eps^2 \snorm{\nbs_\s^4 \nbs_2 h }_{L^2_{x,\s}} + \eps^2 \brak{\mathsf{B_6}} + \eps^2 \snorm{\nbs_\s^4 \nbs_2^2 h }_{L^2_{x,\s}} 
+ \eps^2 \snorm{\nbs_\s^4 V}_{L^2_{x,\s}}
\,.
\end{align*}

To conclude the proof of the lemma, we first sum \eqref{D2D5V}, \eqref{D1D5V},  \eqref{DsD5V}  to obtain
\begin{equation*} 
\snorm{ \nbs^6 V}_{L^2_{x,\s}}
\les
\eps^2 \brak{\mathsf{B_6}}
 +  \snorm{\nbs^6 \nbs_2 h}_{L^2_{x,\s}}   
 + \eps \snorm{\nbs^5 \Jg}_{L^2_{x,\s}} 
 \, ,
\end{equation*} 
which proves \eqref{eq:V:H6}. Finally, inserting \eqref{D1D5V} in the sum of 
\eqref{D2D5Sigma}, \eqref{D1D5Sigma}, and \eqref{DsD5Sigma}, we obtain that
 \begin{align*} 
 \snorm{\nbs^6 \Sigma}_{L^2_{x,\s}} 
 &  \les 
 \eps \brak{\widetilde{\mathcal{D}}_{5}}
+ \eps  \snorm{\nbs^5 \nbs_2 h}_{L^2_{x,\s}}   
 + \eps^2 \snorm{\nbs^5 \Jg}_{L^2_{x,\s}} 
 \notag\\
 &\quad 
+ \eps \bigl(\eps^2 \brak{\mathsf{B_6}}
+  \snorm{\nbs^6\nbs_2 h}_{L^2_{x,\s}}
+\eps \snorm{\nbs^5 \Jg }_{L^2_{x,\s}}\bigr)
\notag\\
&\quad 
+ \eps
\bigl( \eps + \eps^2 \brak{\mathsf{B_6}}
+  \snorm{\nbs^6\nbs_2 h}_{L^2_{x,\s}}
+\eps \snorm{\nbs^5 \Jg }_{L^2_{x,\s}} \bigr)
\bigl( \brak{\mathsf{B_6}}
 +   \snorm{\nbs^5 \nbs_2 h}_{L^2_{x,\s}}   
 + \eps \snorm{\nbs^5 \Jg}_{L^2_{x,\s}} \bigr)
 \notag\\
  &  \les 
 \eps \brak{\widetilde{\mathcal{D}}_{5}}
+\eps \bigl(\eps +   \snorm{\nbs^6 \nbs_2 h}_{L^2_{x,\s}}  + \eps \snorm{\nbs^5 \Jg}_{L^2_{x,\s}} \bigr)  \bigl( \brak{\mathsf{B_6}}
 +   \snorm{\nbs^6 \nbs_2 h}_{L^2_{x,\s}}   
 + \eps \snorm{\nbs^5 \Jg}_{L^2_{x,\s}} \bigr)
 \end{align*} 
which proves \eqref{eq:Sigma:H5}. Here we have used that $\eps \brak{\mathsf{B}_6} \leq 1$.
\end{proof}

Next, we turn to the energy estimate for $\Jg$, aiming to prove \eqref{D6JgEnergy:new}.
\begin{lemma}
\label{lem:D5-Jg}
Under the same assumptions as Proposition~\ref{prop:geometry}, we have that  
\begin{equation} 
\sup_{\s\in[0,\eps]} \snorm{\mathcal{J}^{\!\frac 14}  \nbs^6 \Jg(\cdot,\s) }^2_{L^2_x}
+ \tfrac{1}{\eps} \int_0^\eps \snorm{\mathcal{J}^{-\!\frac 14} \nbs^6 \Jg(\cdot,\s) }^2_{L^2_x} {\rm d}\s
 \les  
 \eps \brak{\mathsf{B}_6}^2 
 \,,   \label{D6JgEnergy}
\end{equation} 
where the implicit constant depends only on $\alpha$ and $\Cdata$.
\end{lemma} 
\begin{proof}[Proof of Lemma \ref{lem:D5-Jg}]
We let $\nbs^6$ act on \eqref{Jg-evo-s} and write this as
\begin{equation} 
(\Q\p_\s +V\p_2)(\nbs^6\Jg) = \tfrac{1+\alpha}{2} \nbs^6(\Jg \Wbn)+ \tfrac{1-\alpha}{2}\nbs^6 (\Jg\Zbn) + \Rj\,, 
\label{D5-Jg-s}
\end{equation} 
where $\Rj  = -\nbs^6 V \, \nbs_2 \Jg - \doublecom{\nbs^6, V, \nbs_2 \Jg}$
thanks to \eqref{good-comm}.

We compute the $L^2$-inner product of \eqref{D5-Jg-s} with $\mathcal{J}^{\!\frac 12} \nbs^6\Jg$ to obtain that
\begin{equation}
\tfrac{1}{2}\! \int \! \mathcal{J}^{\!\frac 12}  (\Q\p_\s +V\p_2) |\nbs^6\Jg|^2 
= \tfrac{1+\alpha}{2} \! \int \! \mathcal{J}^{\!\frac 12} \nbs^6(\Jg \Wbn)\nbs^6\Jg
+ \tfrac{1-\alpha}{2} \! \int \! \mathcal{J}^{\!\frac 12}   \nbs^6(\Jg\Zbn)\nbs^6\Jg 
+ \int \mathcal{J}^{\!\frac 12}   \Rj \ \nbs^6\Jg \,.
\label{eq:D6-Jg-new-label}
\end{equation} 
We next commute $\mathcal{J}^{\!\frac 12} $ around $(\Q\p_\s +V\p_2)$.
For this purpose we note that identity \eqref{evo-fakeJg} shows that for any function $f = f(x,\s)$ and any $r\in \mathbb{R}$ we have
\begin{align}
 \mathcal{J}^r  (\Q\p_\s+V\p_2)  f 
 &=\p_\s \bigl( \mathcal{J}^r \Q f\bigr) + \p_2 \bigl(\mathcal{J}^r V f\bigr)
 +  r f \mathcal{J}^{r-1} \tfrac{\Q}{\eps}
 - f \mathcal{J}^{r} \bigl(\Qc + \p_2 V\bigr)
 \,.
 \label{eq:diff:by:parts}
\end{align}
Using \eqref{eq:diff:by:parts} with $r=\frac 12$ and $f = \frac 12  |\nbs^6\Jg|^2$, we rewrite \eqref{eq:D6-Jg-new-label} as
\begin{align*}
&\tfrac{1}{2 d\s} \int \mathcal{J}^{\!\frac 12}  \Q  |\nbs^6\Jg|^2 
+ \tfrac 14 \int \mathcal{J}^{-\frac 12}  \tfrac{\Q}{\eps}  |\nbs^6\Jg|^2 
- \tfrac 12  \int \mathcal{J}^{\!\frac 12}   \bigl(\Qc + \p_2 V \bigr)  |\nbs^6\Jg|^2 
\notag\\
&= \tfrac{1+\alpha}{2}\int \mathcal{J}^{\!\frac 12} \nbs^6(\Jg \Wbn) \nbs^6\Jg
+ \tfrac{1-\alpha}{2}\int \mathcal{J}^{\!\frac 12}   \nbs^6(\Jg\Zbn)\nbs^6\Jg 
+ \int \mathcal{J}^{\!\frac 12}   \Rj \ \nbs^6\Jg \,.
\end{align*} 
Then, using the bootstraps~\eqref{bootstraps}, and the bounds \eqref{eq:Q:all:bbq} we obtain that
\begin{align} 
&\tfrac{d}{2d\s}  \int \Q \mathcal{J}^{\!\frac 12}  |\nbs^6\Jg|^2   
+ \tfrac{1+ \alpha }{10\eps} (1- \Cn \eps) \| \mathcal{J}^{-\!\frac 14}  \nbs^6\Jg  \|_{L^2}^2 
- \tfrac{2\cdot 250^2}{\eps} \int \Q \mathcal{J}^{\!\frac 12}  |\nbs^6\Jg|^2   
\notag\\
&\leq
\| \mathcal{J}^{-\!\frac 14}  \nbs^6\Jg  \|_{L^2} \Bigl(
\tfrac{1+\alpha}{2} \|\mathcal{J}^{\!\frac 34}  \nbs^6(\Jg \Wbn)\|_{L^2}
+ \tfrac{|1-\alpha|}{2} \|\mathcal{J}^{\!\frac 34} \nbs^6(\Jg\Zbn)\|_{L^2} 
+\| \mathcal{J}^{\!\frac 34}   \Rj \|_{L^2} \Bigr)
\notag\\
&\leq
\tfrac{1+\alpha}{20 \eps} \| \mathcal{J}^{-\!\frac 14}  \nbs^6\Jg  \|_{L^2}^2
+4(1+\alpha) \eps \|\mathcal{J}^{\!\frac 34}  \nbs^6(\Jg \Wbn,\Jg \Zbn)\|_{L^2}^2 
+ \tfrac{16\eps}{1+\alpha} \| \mathcal{J}^{\!\frac 34}   \Rj \|_{L^2}^2
\,.
\label{D6JgEnergy:temp1}
\end{align} 
In order to bound the commutator term appearing on the right side of the above estimate, we appeal to  the bootstraps~\eqref{bootstraps}, the bound $ \mathcal{J} \les 1$, and to Lemmas~\ref{lem:Sigma:V},~\ref{lem:anisotropic:sobolev} and~\ref{lem:comm:tangent},
to conclude that 
\begin{align}
\snorm{\mathcal{J}^{\!\frac 34} \Rj}_{L^2_{x,\s}}
&\leq
\|  \nbs^6 V \|_{L^2_{x,\s}} \| \nbs_2 \Jg \|_{L^\infty_{x,\s}} 
+ \|  \doublecom{\nbs^6, V, \nbs_2 \Jg} \|_{L^2_{x,\s}}
\notag\\
&\les
\bigl(\eps^2 \brak{\mathsf{B_6}}
 +  \snorm{\nbs^6 \nbs_2 h}_{L^2_{x,\s}}   
 + \eps \snorm{\nbs^5 \Jg}_{L^2_{x,\s}} 
 \bigr) \bigl( 1 + \|\nbs^2 \Jg\|_{L^\infty_{x,\s}}\bigr)
+ \eps \|  \nbs^6 \Jg \|_{L^2_{x,\s}}
+ \eps^2  \bigl( 1 + \|\nbs^2 \Jg\|_{L^\infty_{x,\s}}\bigr)
\notag\\
&\les \eps^2 \brak{\mathsf{B_J}} \bigl( \brak{\mathsf{B_6}} + \brak{\mathsf{B_J}} + \brak{\mathsf{B_h}}\bigr)
\,.
\label{D6JgEnergy:temp2}
\end{align}
Using Gr\"onwall's inequality in time for $\s \in [0,\eps]$ in the bound \eqref{D6JgEnergy:temp1}, inserting the commutator bound \eqref{D6JgEnergy:temp2} with the assumption that $\eps^{\frac 12}  ( \brak{\mathsf{B_J}} +   \brak{\mathsf{B_h}} +     \brak{\mathsf{B_6}}) \leq 1$, 
appealing to the upper and lower bound on $\Q$ which arises from \eqref{eq:Q:all:bbq}, we deduce
\begin{equation} 
\sup_{\s\in[0,\eps]}  \|\mathcal{J}^{\!\frac 14} \nbs^6\Jg(\cdot , \s)\|^2_{L^2} 
+ \tfrac{1 }{\eps} \!\! \int_0^\eps  \|\mathcal{J}^{-\!\frac 14} \nbs^6\Jg (\cdot , \s)\|^2_{L^2} {\rm d}\s
\les  
\| \nbs^6\Jg(\cdot , 0)\|^2_{L^2_x}
+ \eps \brak{\mathsf{B}_6}^2
 \,, \label{D6JgEnergy:temp3}
\end{equation}
where the implicit constant depends only on $\alpha$.
The bound \eqref{table:derivatives}, which gives $\|\nbs^6\Jg(\cdot , 0)\|^2_{L^2_x} \les \eps$,  and the bootstraps \eqref{bootstraps-Dnorm:6}--\eqref{bootstraps-Dnorm:5} conclude the proof of \eqref{D6JgEnergy}.
\end{proof} 

\begin{remark}
We shall frequently make use of the fact that \eqref{eq:Sobolev}, \eqref{table:derivatives}, \eqref{D6JgEnergy}, and the bound $1 \leq \mathcal{J}^{-\frac 14}$ imply
\begin{equation} 
\snorm{ \nbs^3\Jg}_{L^\infty _{x,\s}}\les \brak{\mathsf{B_6}} \,. \label{D2-Jg-Linfty}
\end{equation}  
\end{remark}

\begin{remark}
\label{rem:D5:Jg:energy}
In the proof of Lemma~\ref{lem:D5-Jg}, we have tested the equation \eqref{D5-Jg-s} with $\mathcal{J}^{\frac 12} \nbs^6 \Jg$. The presence of the $\mathcal{J}^{\frac 12}$ is what allowed us to obtain the damping term on the left side of \eqref{D6JgEnergy:temp1}, by commuting $\Q \p_\s + V \p_2$ past $\mathcal{J}^{\frac 12}$. In addition, this factor was necessary in order to bound $\nbs^6 (\Jg \Wbn, \Jg \Zbn)$. Other than this, the $\mathcal{J}^{\frac 12}$ factor did not play any role in the bounds; indeed, already in the first line of  \eqref{D6JgEnergy:temp2} we have discarded the extra factor of $\mathcal{J}^{\frac 34}$. With this in mind, we may return to \eqref{Jg-evo-s}, act on it with $\nbs^5$, and this time test it with $\nbs^5 \Jg$. By repeating the same bounds as in the  proof of Lemma~\ref{lem:D5-Jg}, in analogy to \eqref{D6JgEnergy:temp1} and \eqref{D6JgEnergy:temp2},  we may thus establish the bound
\begin{equation*} 
\tfrac{d}{2d\s}  \int  \Q    |\nbs^5 \Jg |^2 
\leq \Cn \eps^{-1}  \int  \Q   |\nbs^5 \Jg |^2
+ \Cn \eps \snorm{ \nbs^5(\Jg  \Wbn, \Jg \Zbn)}_{L^2}^2
+ \Cn \eps  \snorm{ \Rj }_{L^2}^2
\,,
\end{equation*} 
and therefore 
\begin{equation} 
\sup_{\s\in[0,\eps]} e^{- \frac{\Cn \s}{\eps}} \|  \nbs^{5} \Jg (\cdot, \s) \|_{L^2}^2 
\les   \eps  \brak{\mathsf{B_6}}^2
\,.
\label{D5JgEnergy:sup}
\end{equation} 
This concludes the proof of \eqref{D5JgEnergy}.
\end{remark}

\begin{lemma}\label{lem:D5-h2}
Under the same assumptions as Proposition~\ref{prop:geometry},   we have that  
\begin{equation} 
\sup_{\s\in[0,\eps]} \|\mathcal{J}^{\!\frac 14}  \nbs^{6} \nbs_2h(\cdot, \s) \|^2_{L^2}
+ \tfrac{1}{\eps} \int_0^\eps  \| \mathcal{J}^{-\!\frac 14} \nbs^{6}  \nbs_2h \|^2_{L^2}
\les \mathsf{K}^2 \eps^3 \brak{\mathsf{B_6}}^2
\,,
\label{D6h2Energy}
\end{equation} 
where the implicit constant  depends only on $\alpha$ and $\Cdata$.
\end{lemma} 
\begin{proof}[Proof of Lemma~\ref{lem:D5-h2}]
The proof is similar to that of Lemma~\ref{lem:D5-Jg}. We let $\nbs^6$ act on \eqref{p2h-evo-s} and write this as
\begin{align} 
&(\Q\p_\s +V\p_2)\nbs^6 \nbs_2h \notag\\
&= \bigl((1+\alpha)\Wbt + (1- \alpha) \Zbt \bigr) \nbs_2h \nbs^6\nbs_2h
+g \bigl( \tfrac{1+\alpha}{2}  \nbs^6 \Wbt+ \tfrac{1-\alpha}{2} \nbs^6 \Zbt\bigr)
 + {\textstyle\sum}_{k=1}^3 \Rhn k
 \,, \label{D5-h2-s}
\end{align} 
where the remainder terms are given by
\begin{align*} 
\Rhn1 & = - \jump{\nbs^6,V} \nbs_2^2 h \,,  \ \ 
\Rhn2  =   \comm{\nbs^6, g, \tfrac{1+ \alpha }{2} \Wbt + \tfrac{1- \alpha }{2} \Zbt } \,,  \ \
\Rhn3  =   \bigl((1+\alpha)\Wbt + (1- \alpha) \Zbt \bigr) \jump{\nbs^5, \nbs_2h} \nbs  \nbs_2h \,.
\end{align*} 
We compute the $L^2$-inner product of \eqref{D5-h2-s} with $\mathcal{J}^{\!\frac 12}  \nbs^6 \nbs_2h$, and use \eqref{eq:diff:by:parts} with $r = \frac 12$ and $f= \frac 12 |\nbs^6 \nbs_2 h|^2$ to obtain
\begin{align*} 
&\tfrac{d}{2d\s}  \int  \Q \mathcal{J}^{\!\frac 12}  |\nbs^6\nbs_2h|^2 
+ \tfrac{1}{4} \int  \mathcal{J}^{-\!\frac 12} \tfrac{\Q}{\eps}  |\nbs^6\nbs_2h|^2 
\notag\\
&\quad 
=\int  \mathcal{J}^{\frac 12}  \bigl( \tfrac 12  (\Qc + \p_2 V) +  ((1+\alpha)\Wbt + (1- \alpha) \Zbt  ) \nbs_2h \bigr) |\nbs^6\nbs_2h|^2 
\notag \\
& \qquad 
 + \int \mathcal{J}^{\!\frac 12} g \bigl(  \tfrac{1+ \alpha }{2} \nbs^6 \Wbt + \tfrac{1- \alpha }{2} \nbs^6 \Zbt \bigr) 
 \nbs^6 \nbs_2 h 
+ {\textstyle\sum}_{k=1}^3\int\mathcal{J}^{\!\frac 12}   \Rhn k \ \nbs^6\nbs_2h \,.
\end{align*} 
Thanks to \eqref{bootstraps}, \eqref{eq:Q:all:bbq}, and the Cauchy-Young inequality, it follows that 
\begin{align} 
&\tfrac{d}{2d\s}  \int  \Q \mathcal{J}^{\!\frac 12}  |\nbs^6\nbs_2h|^2 
+ \tfrac{1+\alpha}{10 \eps} (1- \Cn \eps) \snorm{\mathcal{J}^{-\!\frac 14}   \nbs^6\nbs_2h}_{L^2}^2 
- \tfrac{2\cdot 250^2}{\eps} \int \Q \mathcal{J}^{\!\frac 12}  |\nbs^6\nbs_2h|^2   
\notag\\
&\quad 
\leq \tfrac{1+\alpha}{2} (1+ \Cn \eps^2) \snorm{\mathcal{J}^{-\!\frac 14}   \nbs^6\nbs_2h}_{L^2}
\snorm{\mathcal{J}^{\frac 34}   ( \nbs^6 \Wbt, \nbs^6 \Zbt)}_{L^2}
+ \snorm{\mathcal{J}^{-\!\frac 14}   \nbs^6\nbs_2h}_{L^2} 
{\textstyle\sum}_{k=1}^3 \snorm{\mathcal{J}^{\!\frac 34}   \Rhn k }_{L^2}
\notag\\
&\quad 
\leq \tfrac{1+\alpha}{20 \eps} \snorm{\mathcal{J}^{-\!\frac 14}   \nbs^6\nbs_2h}_{L^2}^2
+ 4 (1+\alpha) \eps \snorm{\mathcal{J}^{\frac 34}   ( \nbs^6 \Wbt, \nbs^6 \Zbt)}_{L^2}^2
+ \tfrac{15  \eps}{1+\alpha} {\textstyle\sum}_{k=1}^3 \snorm{ \Rhn k }_{L^2}^2
\,.
\label{D6h2Energy:temp:1}
\end{align} 
In anticipation of integrating \eqref{D6h2Energy:temp:1} in time (via Gr\"onwall), we bound the three commutators appearing on the right side as follows. 
First, using \eqref{table:derivatives}, \eqref{bootstraps}, \eqref{eq:V:H6}, \eqref{eq:Sobolev}, and Lemma~\ref{lem:time:interpolation}, Lemma~\ref{lem:comm:tangent}, we have
that
\begin{equation*}
\snorm{  \Rhn1}_{L^2_{x,\s}} 
\leq 
\snorm{\nbs_2^2 h}_{L^\infty_{x,\s}} \snorm{  \nbs^6 V}_{L^2_{x,\s}} 
+
\snorm{\doublecom{\nbs^6,V,\nbs_2^2 h}}_{L^2_{x,\s}}
\les  \eps^2 \brak{\mathsf{B_h}} \brak{\mathsf{B_6}}
\,.
\end{equation*}
Similarly, since $\eps^{\frac 12} \brak{\mathsf{B_h}} \leq 1$, we have
\begin{equation*}
\snorm{\Rhn2}_{L^2_{x,\s}}  
\leq  \tfrac{1+ \alpha }{2} \snorm{\comm{\nbs^6, g,  \Wbt}}_{L^2_{x,\s}} 
+ \tfrac{1- \alpha }{2}  \snorm{\comm{\nbs^6, g, \Zbt }}_{L^2_{x,\s}}
\les
\eps^3 \brak{\mathsf{B_h}} + \eps^3 \mathsf{K} \mathsf{B_5}
\les \eps^2 
\,,
\end{equation*}
and also
\begin{align*}
\snorm{\Rhn3}_{L^2_{x,\s}}
&\leq 
 \bigl((1+\alpha) \snorm{\Wbt}_{L^\infty_{x,\s}} + (1- \alpha) \snorm{\Zbt}_{L^\infty_{x,\s}} \bigr) \bigl( \snorm{\doublecom{\nbs^5, \nbs_2 h, \nbs  \nbs_2 h}}_{L^2_{x,\s}} + \snorm{\nbs_2h}_{L^\infty_{x,\s}} \snorm{\nbs^6  \nbs_2h}_{L^2_{x,\s}} \bigr)
 \notag\\
&\les \eps^3 \brak{\mathsf{B_h}} \les \eps^2 \,.
\end{align*}
Inserting the bounds obtained in the previous three displays into the time-integrated form of \eqref{D6h2Energy:temp:1}, and appealing to the upper and lower bound on $\Q$ which arises from \eqref{eq:Q:all:bbq}, we obtain
\begin{equation} 
\sup_{\s\in[0,\eps]} \| \mathcal{J}^{\frac 14} \nbs^6\nbs_2h(\cdot, \s)\|^2_{L^2} 
+ \tfrac{1 }{\eps} \|\mathcal{J}^{ -\frac 14} \nbs^6\nbs_2 h\|_{L^2_{x,\s}}^2
\les  \| \mathcal{J}^{\frac 14} \nbs^6\nbs_2h(\cdot, 0)\|^2_{L^2} 
+ \mathsf{K}^2 \eps^3  \brak{\mathsf{B_6}}^2  
+ \eps^5 \brak{\mathsf{B_h}}^2 \brak{\mathsf{B_6}}^2 \,,
\label{D6h2Energy:temp:2}
\end{equation}
where the implicit constant depends only on $\alpha$.
The fact that  $\eps \brak{\mathsf{B_h}}^2 + \eps \brak{\mathsf{B_6}}^2\leq 1$, and the bound on the initial data which arises from \eqref{table:derivatives}, concludes the proof of \eqref{D6h2Energy}.
\end{proof}

\begin{remark}
We shall also make use of the fact that from \eqref{eq:Sobolev}, combined with \eqref{table:derivatives}, the bound $1 \leq \mathcal{J}^{-\frac 14}$, and \eqref{D6h2Energy}, we have
\begin{align} 
\snorm{ \nbs^3 \nbs_2 h}_{L^\infty _{x,\s}}\les \mathsf{K} \eps \brak{\mathsf{B_6}} \,. 
\label{eq:h_2:D2-bound}
\end{align}  
\end{remark}

\begin{remark}
\label{rem:D5:h2:energy}
In the proof of Lemma~\ref{lem:D5-h2}, we have tested the equation \eqref{D5-h2-s} with $\mathcal{J}^{\frac 12} \nbs^6 \nbs_2 h$. The presence of the $\mathcal{J}^{\frac 12}$ is what allowed us to obtain the damping term on the left side of \eqref{D6h2Energy}, by commuting $\Q \p_\s + V \p_2$ past $\mathcal{J}^{\frac 12}$. In addition, this factor was necessary in order to bound $\nbs^6 (\Wbt, \Zbt)$. Other than this, the $\mathcal{J}^{\frac 12}$ factor did not play any role in the bounds; indeed, already in~\eqref{D6h2Energy:temp:1} we have discarded the extra factor of $\mathcal{J}^{\frac 34}$ from the commutator terms on the right side. Keeping this in mind, we return to \eqref{p2h-evo-s}, act on it with $\nbs^5$, and test the resulting equation with $\nbs^5 \nbs_2 h$. By repeating the same arguments as in the  proof of Lemma~\ref{lem:D5-h2}, in analogy to \eqref{D6h2Energy:temp:1}--\eqref{D6h2Energy:temp:2}   we may   establish  
\begin{equation*} 
\tfrac{d}{2d\s}  \int  \Q    |\nbs^5\nbs_2h|^2 
\leq \Cn \eps^{-1}  \int  \Q   |\nbs^5\nbs_2h|^2
+ \Cn \eps \snorm{ ( \nbs^5 \Wbt, \nbs^5 \Zbt)}_{L^2}^2
+ \Cn \eps {\textstyle\sum}_{k=1}^3 \snorm{\Rhn k }_{L^2}^2
\,,
\end{equation*} 
and therefore 
\begin{equation} 
\sup_{\s\in[0,\eps]} e^{- \frac{\Cn \s}{\eps}} \|  \nbs^{5} \nbs_2h(\cdot, \s) \|_{L^2}^2 
\les  \mathsf{K}^2 \eps^3 \brak{\mathsf{B_6}}^2
\,.
\label{D5h2Energy:sup}
\end{equation} 
This concludes the proof of \eqref{D5h2Energy}.
\end{remark}

\begin{remark}[\bf Estimates for $g$]
\label{rem:D6:g}
The $\nbs_2 h$ estimates established in Lemma~\ref{lem:D5-h2} have as a direct consequence estimates for $g = 1 + (\nbs_2 h)^2$. More precisely, the identity $\nbs^6 g = 2 \nbs^5( \nbs_2 h \nbs_2^2 h) = 2 \nbs_2 h \nbs_2^7 h + 2 \jump{\nbs^5, \nbs_2 h} \nbs_2^2 h$ combined with the bounds \eqref{bs-h}, \eqref{bs-h,22}, \eqref{D6h2Energy}, \eqref{eq:h_2:D2-bound}, \eqref{D5h2Energy:sup}, \eqref{eq:Lynch:2}, \eqref{eq:Lynch:3}, and \eqref{eq:Lynch:3:bdd} (with $a=a'=b=0$) imply that
\begin{equation*} 
\eps^{\frac 12} \sup_{\s\in[0,\eps]} \|\mathcal{J}^{\!\frac 14}  \nbs^{6} g(\cdot, \s) \|_{L^2}
+ \| \mathcal{J}^{-\!\frac 14} \nbs^{6}  g \|_{L^2_{x,\s}}
\les \mathsf{K}^2 \eps^3 \brak{\mathsf{B_6}}^2
\,,
\end{equation*} 
which proves~\eqref{D6gbound}.
Since $1\leq g \leq 1 + \Cn \eps^2$, by appealing to the chain rule it is clear that any rational power of $g$ appearing in the proof (e.g.~$g^{\frac 32}, g^{\frac 12}, g^{-\frac 12}, g^{-\frac 32}$) satisfies the same bound as $g$.
\end{remark}

\begin{remark}[\bf Estimates for $\nbs_1 h$]
\label{rem:D6:D1:h}
We note that the bound \eqref{D6h1Energy:new} is a direct consequence of the identity~$\nbs_1 h = \eps g^{\frac 12} \Jg$, of the bounds~\eqref{bs-Jg-simple}--\eqref{bs-h,22}, \eqref{D6JgEnergy:new}--\eqref{D6gbound} and of the product Moser-type Lemmas~\ref{lem:Moser:tangent} and~\ref{lem:comm:bdd}.
\end{remark}

\begin{remark}
\label{rem:D5:Sigma:V}
We note that with  the bound $1\les \mathcal{J}^{-\frac 14}$, inserting estimates~\eqref{D6JgEnergy} and~\eqref{D6h2Energy} into \eqref{D6-Sigma-V}, gives the proof of \eqref{eq:Sigma:H6:new} and \eqref{eq:V:H6:new}. Indeed, \eqref{eq:Sigma:H5} becomes
\begin{equation*} 
 \snorm{\nbs^6 \Sigma}_{L^2_{x,\s}} 
 \les \eps \brak{\widetilde{\mathcal{D}}_{5}}
+\eps \bigl(\eps + \mathsf{K} \eps^2 \brak{\mathsf{B_6}} + \eps^2 \brak{\mathsf{B_6}} \bigr)  \bigl( \brak{\mathsf{B_6}}
 + \mathsf{K} \eps^2 \brak{\mathsf{B_6}} 
 +  \eps^2 \brak{\mathsf{B_6}} \bigr) 
 \les \eps \brak{\widetilde{\mathcal{D}}_{5}}  + \eps^2 \brak{\mathsf{B_6}}
 \les \eps \brak{\mathsf{B_6}}
 \,,
 \end{equation*}
 since $\eps \brak{\mathsf{B_6}} \leq 1$, and recalling Remark~\ref{rem:B5:B6}.
Similarly, \eqref{eq:V:H6} becomes  
\begin{equation*}
\snorm{\nbs^6 V}_{L^2_{x,\s}}
\les \eps^2 \brak{\mathsf{B_6}}
+ \mathsf{K} \eps^2 \brak{\mathsf{B_6}}  
+  \eps^2 \brak{\mathsf{B_6}} 
\les \mathsf{K} \eps^2 \brak{\mathsf{B_6}}  \,.
\end{equation*} 
\end{remark}

\begin{corollary}\label{cor:D5n-bound}
Under the assumptions of Proposition~\ref{prop:geometry}, we have that for $3\leq |\gamma|\leq 6$,
\begin{subequations}
\label{D5n-bound}
\begin{align} 
\|\mathcal{J}^{-\frac 14} \bigl( \nbs^{|\gamma|}  \nn + g^{-1}  \tt \nbs^{|\gamma|} \nbs_2 h \bigr) \|_{{L^2_{x,\s}} }
+
\|\mathcal{J}^{-\frac 14} \bigl( \nbs^{|\gamma|}  \tt - g^{-1}  \nn \nbs^{|\gamma|} \nbs_2 h \bigr) \|_{{L^2_{x,\s}} }
&\les \mathsf{K} \eps^3 \brak{\mathsf{B_6}} 
\label{D5n-bound:a}\\
\|\mathcal{J}^{-\frac 14}  \nbs^{|\gamma|}  \nn  \|_{{L^2_{x,\s}} }
+
\|\mathcal{J}^{-\frac 14} \nbs^{|\gamma|}  \tt \|_{{L^2_{x,\s}} }
&\les \mathsf{K}\eps^2 \brak{\mathsf{B_6}}   
\label{D5n-bound:b}
\,.
\end{align} 
\end{subequations}
\end{corollary} 
\begin{proof}[Proof of Corollary~\ref{cor:D5n-bound}]
The bounds for $\nn$ and $\tt$ are symmetric, so we only prove the estimates for $\nn$.
From \eqref{d-n-tau}, we have that $\nbs \nn = - g^{-1} \nbs \nbs_2 h \, \tt$. 
Therefore, we have 
\begin{equation}
\nbs^{|\gamma|}  \nn  + g^{-1}  \tt \nbs^{|\gamma|} \nbs_2 h = - \jump{\nbs^{|\gamma|-1},g^{-1} \tt} \nbs \nbs_2 h
.
\label{eq:chain:rule:dont:stop}
\end{equation}
In order to bound the commutator term on the right side of \eqref{eq:chain:rule:dont:stop}, we appeal to 
inequality \eqref{eq:r:r':L2:time:new} with $r'=-\frac 14> - \frac 12$, $r=0<\frac 14$, and $F = \jump{\nbs^{|\gamma|-1},g^{-1} \tt} \nbs \nbs_2 h$, together with the chain and product rules, and the $\nbs^k$ bounds for $\nbs_2 h(\cdot,0)$ contained in \eqref{table:derivatives}, to deduce
\begin{align}
\snorm{\mathcal{J}^{-\frac 14} \jump{\nbs^{|\gamma|-1},g^{-1} \tt} \nbs \nbs_2 h}_{L^2_{x,\s}}
&\les  
\eps^{\frac 12} \snorm{\jump{\nbs^{|\gamma|-1},g^{-1} \tt} \nbs \nbs_2 h(\cdot,0)}_{L^2_{x}}
+ \eps \snorm{\p_\s \jump{\nbs^{|\gamma|-1},g^{-1} \tt} \nbs \nbs_2 h}_{L^2_{x,\s}}
\notag\\
&\les \eps^3 + \snorm{\nbs_\s \jump{\nbs^{|\gamma|-1},g^{-1} \tt} \nbs \nbs_2 h}_{L^2_{x,\s}}.
\end{align}
Using the product rule and \eqref{d-n-tau}, we may further rewrite
\begin{equation*}
\nbs_\s \jump{\nbs^{|\gamma|-1},g^{-1} \tt} \nbs \nbs_2 h
=  \jump{\nbs_\s\nbs^{|\gamma|-1},g^{-1} \tt} \nbs \nbs_2 h 
+ \bigl(2   \nbs_2 h   \tt - \nn \bigr) g^{-2} \nbs_\s \nbs_2 h  \nbs^{|\gamma|} \nbs_2 h
\end{equation*}
Due to \eqref{eq:Jgb:less:than:1}, \eqref{bs-h},  \eqref{bs-h,22}, and \eqref{D6h2Energy:new}, the second term in the above estimate may be bounded as 
\begin{equation*}
\snorm{\bigl(2   \nbs_2 h   \tt - \nn \bigr) g^{-2} \nbs_\s \nbs_2 h  \nbs^{|\gamma|} \nbs_2 h}_{L^2_{x,\s}}
\les \eps  \snorm{\nbs^{|\gamma|} \nbs_2 h}_{L^2_{x,\s}} 
\les \eps^3 \mathsf{K} \brak{\mathsf{B}_6}\,.
\end{equation*}
Similarly, by  by appealing to the chain rule, the bootstrap inequalities~\eqref{bootstraps}, the bounds~\eqref{D6h2Energy} and~\eqref{eq:h_2:D2-bound}, and Lemma~\ref{lem:comm:tangent}, that
\begin{equation*}
\|\jump{\nbs_\s\nbs^{|\gamma|-1},g^{-1} \tt} \nbs \nbs_2 h \|_{{L^2_{x,\s}} }
\leq \snorm{\nbs \nbs_2 h}_{L^\infty_{x,\s}} \snorm{\nbs^{|\gamma|} (g^{-1} \tt)}_{L^2_{x,\s}}
+  \snorm{\doublecom{\nbs^{|\gamma|},g^{-1} \tt} \nbs \nbs_2 h}_{L^2_{x,\s}}
\les \eps^3 \mathsf{K} \brak{\mathsf{B}_6}
\,.
\end{equation*}
Combining the above five displays concludes the proof of \eqref{D5n-bound:a}. The bound \eqref{D5n-bound:b} follows from \eqref{bs-h}, \eqref{D6h2Energy}, and \eqref{D5n-bound:a}.
\end{proof}

It remains to prove estimates~\eqref{eq:Sigma:H6:new:bdd} and~\eqref{eq:V:H6:new:bdd}, which is achieved next. 

\begin{lemma}
\label{lem:D6:V:Sigma:bdd:time} 
Under the assumptions of Proposition~\ref{prop:geometry}, we have that
\begin{subequations}
\begin{align}
\snorm{\Jgh \nbs^6 \Sigma^{\pm1}}_{L^\infty_\s L^2_{x} } 
&\les 
\mathsf{K}  \eps^{\frac 12} \brak{\mathsf{B}_6},
\\
\snorm{\Jgh \nbs^6 V}_{L^\infty_s L^2_{x}} 
&\les \mathsf{K} \eps^{\frac 32} \brak{\mathsf{B}_6},
\\
\|\mathcal{J}^{\frac 14} \nbs^6  \nn  \|_{L^\infty_\s L^2_x}
+
\|\mathcal{J}^{\frac 14}  \nbs^6 \tt \|_{L^\infty_\s L^2_x}
&\les \mathsf{K} \eps^{\frac 32} \brak{\mathsf{B}_6} .
\end{align} 
\end{subequations}
\end{lemma}
\begin{proof}[Proof of Lemma~\ref{lem:D6:V:Sigma:bdd:time}]
In view of Remark~\ref{rem:shit:is:similar}, the proof is very similar to that of Lemma~\ref{lem:Sigma:V}, and so we only present here the differences.  As such, by repeating the argument used to establish \eqref{D2D5Sigma}, we deduce from \eqref{eq:Lynch:1:bdd} that 
\begin{equation*} 
\snorm{\Jgh \nbs^5 \nbs_2 \Sigma}_{L^\infty_\s L^2_{x}} 
\les \snorm{\Jgh  \nbs^5 (\Wbt,\Zbt) }_{L^\infty_\s L^2_{x}} 
+ \eps \snorm{\Jgh  \nbs^5 \nbs_2 h  }_{L^\infty_\s L^2_{x}} 
+ \eps^{\frac 12} \brak{\mathsf{B}_6} 
\les \mathsf{K} \eps^{\frac 12} \brak{\mathsf{B}_6} 
\,.
\end{equation*} 
Similarly to \eqref{D1D5Sigma}, with \eqref{eq:Lynch:1:bdd} we have
\begin{align} 
\snorm{\Jgh \nbs^5 \nbs_1 \Sigma}_{L^\infty_\s L^2_{x}} 
&\les 
\eps \snorm{\Jgh  \nbs^5 (\Jg \Wbn,\Zbn) }_{L^\infty_\s L^2_{x}} 
+ \eps^2 \snorm{\Jgh  \nbs^5  \Jg  }_{L^\infty_\s L^2_{x}} 
\notag\\
&\qquad 
+ \eps \snorm{\Jgh  \nbs^5   \nbs_2 h    }_{L^\infty_\s L^2_{x}} 
+ \eps^2 \snorm{\Jgh   (\Wbt ,\Zbt)  }_{L^\infty_\s L^2_{x}} 
+ \eps^{\frac 12} \brak{\mathsf{B}_6} 
\les \eps^{\frac 12} \brak{\mathsf{B}_6} 
\,,
\end{align}
while in analogy to \eqref{DsD5Sigma} from \eqref{eq:Lynch:1:bdd} and \eqref{eq:sup:in:time:L2} we have
\begin{align} 
\snorm{\Jgh \nbs_{\s}^6 \Sigma}_{L^\infty_\s L^2_{x}} 
&\les 
\eps \snorm{ \nbs_{\s}^5  V}_{L^\infty_\s L^2_{x}} 
+ \eps^2 \snorm{\Jgh \nbs_{\s}^5  \nbs_2 \Sigma}_{L^\infty_\s L^2_{x}} 
+ \eps \snorm{ \nbs_{\s}^5    \Sigma}_{L^\infty_\s L^2_{x}} 
\notag\\
&\qquad
+ \eps \snorm{\Jgh \nbs_{\s}^5 \Zbn}_{L^\infty_\s L^2_{x}} 
+ \eps \snorm{\Jgh \nbs_{\s}^5 \Abt}_{L^\infty_\s L^2_{x}} 
+ \eps^{\frac 12} \brak{\mathsf{B}_6}
\notag\\
&\les 
\eps^{\frac 12} \snorm{  \nbs_{\s}^6  V}_{L^2_{x,\s}} 
+ \eps^2 \snorm{\Jgh \nbs_{\s}^5  \nbs_2 \Sigma}_{L^\infty_\s L^2_{x}} 
+ \eps^{\frac 12} \snorm{  \nbs_{\s}^6    \Sigma}_{L^2_{x,\s}} 
\notag\\
&\qquad
+ \eps \snorm{\Jgh \nbs_{\s}^5 \Zbn}_{L^\infty_\s L^2_{x}} 
+ \eps \snorm{\Jgh \nbs_{\s}^5 \Abt}_{L^\infty_\s L^2_{x}} 
+ \eps^{\frac 12} \brak{\mathsf{B}_6}
\les 
 \eps^{\frac 12} \brak{\mathsf{B}_6}
 \,.
\end{align}
In the last inequality we have additionally appealed to \eqref{eq:Sigma:H6:new}, \eqref{eq:V:H6:new}, 
the improved $\Zbn$ estimate in \eqref{eq:Jg:Zbn:D5:improve:a}, and to the previously established bound for 
$\snorm{\Jgh \nbs_{\s}^5  \nbs_2 \Sigma}_{L^\infty_\s L^2_{x}}$. 

The bound for $\snorm{\Jgh \nbs^6 V}_{L^\infty_\s L^2_{x}}$ is obtained similarly to \eqref{D2D5V} when $\nbs^6 = \nbs^5 \nbs_2$, in analogy to \eqref{D1D5V} when $\nbs^6 = \nbs^5\nbs_1$, and in parallel to \eqref{DsD5V} when $\nbs^6 = \nbs_\s^6$. To avoid redundancy we omit these details.

The bound for $\snorm{\mathcal{J}^{\frac 14} \nbs^6 \nn}_{L^\infty_\s L^2_x}$ follows from the identity~\eqref{eq:chain:rule:dont:stop}, the previously established bounds \eqref{D6h2Energy:new} and~\eqref{D5h2Energy}, the bound~\eqref{eq:h_2:D2-bound}, and the commutator bound~\eqref{eq:Lynch:3:bdd} with $m=5$ and $a=b=0$:
\begin{align}
\snorm{\mathcal{J}^{\frac 14} \nbs^6 \nn}_{L^\infty_\s L^2_x}
&\leq  
\snorm{\mathcal{J}^{\frac 14} \nbs^6 \nbs_2 h}_{L^\infty_\s L^2_x}
+ 
\Cn \eps \snorm{\nbs^5(g^{-1} \tt)}_{L^\infty_\s L^2_x} 
+ 
\Cn \eps^{-\frac 12} \bigl(\eps \snorm{\nbs^5(g^{-1} \tt)}_{L^2_{x,\s}} + \snorm{\nbs^6 \nbs_2 h}_{L^2_{x,\s}} + \eps^2 \bigr)
\notag\\
&\leq \Cn \eps^{\frac 32} \mathsf{K} \brak{\mathsf{B}_6}\,.
\end{align}
This concludes the proof of the lemma.
\end{proof}

\section{Vorticity energy estimates and the resulting improved estimates}
\label{sec:vorticity}

\subsection{Bounds for the vorticity}

The goal of this subsection is to establish the following bound. 
\begin{proposition}[$H^6$ estimates for the vorticity]
\label{prop:vort:H6}
Let $\Omega$ be the ALE vorticity, defined in \eqref{vort-id-good}, and $\Upomega$ be the ALE specific vorticity given by \eqref{svort-id-good}.
Assume that the bootstrap assumptions~\eqref{bootstraps} hold, and that $\eps$ is taken to be sufficiently small to ensure $\eps^{\frac 12} \brak{\mathsf{B_J}} + \eps^{\frac 12} \brak{\mathsf{B_h}} +   \eps^{\frac 12} \brak{\mathsf{B_6}} \leq 1$.  Then, assuming $\eps$ is sufficiently small with respect to $\alpha,\kappa_0$, and $\Cdata$, we have the bound
\begin{equation}
\sup_{\s\in[0,\eps]} \snorm{\Jgh  \nbs^6\Upomega (\cdot,\s)}_{L^2_x}^2 
+
\tfrac{1}{ \eps}  \int_0^{\eps}\snorm{ \nbs^6\Upomega (\cdot,\s)}_{L^2_x}^2 {\rm d}\s
\les \eps \brak{\mathsf{B_6}}^2  \,,
\label{eq:svort:H6}
\end{equation}
where the implicit constant depends only on $\alpha,\kappa_0$, and $\Cdata$.
Additionally, we have that 
\begin{equation}
 \sup_{\s\in[0,\eps]} \snorm{\Jgh  \nbs^6\Omega (\cdot,\s)}_{L^2_x}^2  
+
\tfrac{1}{ \eps}  \int_0^{\eps}\snorm{\nbs^6\Omega (\cdot,\s)}_{L^2_x}^2 {\rm d}\s
\les \eps \brak{\mathsf{B_6}}^2 \,,
\label{eq:vort:H6}
\end{equation}
where the implicit constant depends only on $\alpha,\kappa_0$, and $\Cdata$.
Moreover, we have that 
\begin{subequations}
\label{eq:vorticity:pointwise}
\begin{align}
\snorm{\Omega}_{L^\infty_{x,\s}} 
&\leq  2^{3+\frac{2}{\alpha}} e^{18}  \Cdata
\,,
\label{eq:vorticity:pointwise:a}
\\
\snorm{\nbs \Omega}_{L^\infty_{x,\s}} 
&\leq 2 (4 e^{18})^{\frac{20 \cdot 23 (1+\alpha)}{\alpha}} \Cdata
\,.
\label{eq:vorticity:pointwise:b}
\end{align}
\end{subequations}
\end{proposition}

\begin{proof}[Proof of Proposition~\ref{prop:vort:H6}]
Recall from \eqref{svort-id-good} that  the vorticity is obtained from the specific vorticity by $\Omega = (\alpha \Sigma)^{\frac 1\alpha} \Upomega$. 
In light of the already established bound \eqref{eq:Sigma:H6:new}  and the product and chain rules, the bound \eqref{eq:vort:H6} follows from \eqref{eq:svort:H6}. As such, we shall only establish the later.

Applying the operator $\nbs^6$ to \eqref{vort-s}, we have that
\begin{equation} 
\tfrac{\Jg}{\Sigma} (\Q \p_\s + V \p_2)  \nbs^6  \Upomega -  \alpha    \p_1 \nbs^6  \Upomega + \alpha \Jg   g^{-\frac 12} \nbs_2 h  \nbs_2 \nbs^6  \Upomega = \RR_\Upomega
\,, \label{D6-vort-s}
\end{equation} 
where the remainder term $\RR_\Upomega$ is given by appealing to \eqref{the-time-der} and \eqref{good-comm} as 
\begin{align} 
\RR_\Upomega
&=
- \jump{\nbs^6,\tfrac{\Jg}{\Sigma}} \bigl(\tfrac{1}{\eps} \nbs_\s  + V \nbs_2\bigr) \Upomega
- \tfrac{\Jg}{\Sigma} \jump{\nbs^6,V} \nbs_2 \Upomega 
- \alpha \jump{\nbs^6, \Jg g^{- {\frac{1}{2}} }\nbs_2 h } \nbs_2 \Upomega 
\notag\\
&=
- \tfrac{1}{\eps}  \jump{\nbs^6,\tfrac{\Jg}{\Sigma}} \nbs_s  \Upomega
- \jump{\nbs^6,\tfrac{\Jg}{\Sigma}}  \bigl( V \nbs_2\Upomega\bigr) 
- \tfrac{\Jg}{\Sigma} \jump{\nbs^6,V} \nbs_2 \Upomega 
-  \alpha \jump{\nbs^6, \Jg g^{- {\frac{1}{2}} }\nbs_2 h } \nbs_2 \Upomega 
\notag\\
&=: \RR_\Upomega^{(1)} + \RR_\Upomega^{(2)} + \RR_\Upomega^{(3)} + \RR_\Upomega^{(4)}
  \,. \label{D4-vort-s-remainder}
\end{align} 
For $\beta>0$ to be chosen below, we compute the spacetime $L^2$ inner-product of \eqref{D6-vort-s} with $\Sigma^{-2\beta+1} \nbs^6 \Upomega$, 
use the identities \eqref{adjoint-1}, \eqref{adjoint-2}, and \eqref{adjoint-3}, to obtain that
\begin{align}
&  \snorm{\tfrac{(\Jg \Q)^{\frac 12}}{ \Sigma^{\beta}}  \nbs^6\Upomega (\cdot,\s)}_{L^2_x}^2 
-  \snorm{\tfrac{(\Jg \Q)^{\frac 12}}{ \Sigma^{\beta}}  \nbs^6\Upomega  (\cdot,0)}_{L^2_x}^2 
-   \tint \sabs{\nbs^6\Upomega}^2 (\Q \p_\s + V \p_2) \tfrac{\Jg}{\Sigma^{2\beta}}  
- \alpha (2\beta-1) \tint \sabs{\nbs^6\Upomega}^2 \tfrac{\Sigma,_1}{\Sigma^{2\beta}} 
\notag\\
&=   \tint \sabs{\nbs^6\Upomega}^2 \tfrac{\Jg}{\Sigma^{2\beta}} \bigl(\Qr_\s -  V \Qr_2 + \nbs_2 V  \bigr) 
+ \alpha  \tint \sabs{\nbs^6\Upomega}^2 \bigl(\nbs_2 - \Qr_2\bigr) \tfrac{\Jg  g^{-\frac 12} \nbs_2 h}{\Sigma^{2\beta-1}}
\notag\\
&\qquad
+ \alpha  \int \sabs{\nbs^6\Upomega}^2 \Qb_2 \tfrac{\Jg  g^{-\frac 12} \nbs_2 h}{\Sigma^{2\beta-1}} \Big|_{\s}
+  \tint\tfrac{2}{\Sigma^{2\beta-1}} \RR_\Upomega \nbs^6 \Upomega 
\,.
\label{eq:vorticity:energy:0}
\end{align}
By appealing to \eqref{Jg-evo-s}, \eqref{p1-Sigma-s}, \eqref{Sigma0-ALE-s}, the above becomes
\begin{align}
&  \snorm{\tfrac{(\Jg \Q)^{\frac 12}}{ \Sigma^{\beta}}  \nbs^6\Upomega (\cdot,\s)}_{L^2_x}^2 
-  \snorm{\tfrac{(\Jg \Q)^{\frac 12}}{ \Sigma^{\beta}}  \nbs^6\Upomega  (\cdot,0)}_{L^2_x}^2 
+ \tint \sabs{\nbs^6\Upomega}^2 \tfrac{1}{\Sigma^{2\beta}}  \mathsf{G}_{\Upomega} 
\notag\\
&=    \alpha  \int \tfrac{\Jg}{\Sigma^{2\beta}} \sabs{\nbs^6\Upomega}^2 \Qb_2  \Sigma  g^{-\frac 12} \nbs_2 h  \Big|_{\s}
+  \tint\tfrac{2}{\Sigma^{2\beta-1}} \RR_\Upomega \nbs^6 \Upomega 
\,.
\end{align}
where 
\begin{align}
\mathsf{G}_{\Upomega}
&:=  
-  \Bigl( \bigl( \tfrac{1+\alpha}{2} \Jg \Wbn + \tfrac{1-\alpha}{2} \Jg \Zbn\bigr) + 2 \alpha \beta  \Jg   \bigl( \Zbn + \Abt\bigr)\Bigr)
- \alpha (\beta-\tfrac 12) \bigl(  \Jg \Wbn -  \Jg \Zbn + \Jg \nbs_2 h (\Wbt - \Zbt) \bigr)
\notag\\
&\qquad 
- \Jg \bigl(\Qr_\s -  V \Qr_2 + \nbs_2 V  \bigr) 
+ (2\beta-1) \alpha \Jg  g^{-\frac 12} \nbs_2 h \nbs_2 \Sigma 
+ \alpha  \Qr_2 \Sigma \Jg  g^{-\frac 12} \nbs_2 h 
\label{eq:vorticity:energy:00}
\end{align}
Using the bounds \eqref{eq:Q:all:bbq}, the damping inequality \eqref{eq:signed:Jg}, and the bootstrap 
assumptions~\eqref{bootstraps}, we conclude that  
\begin{subequations}
\label{justbecause0}
\begin{equation}
\mathsf{G}_{\Upomega} 
\geq - (\alpha \beta + \tfrac 12)  \Jg \Wbn 
- \tfrac{2 \cdot 250^2}{\eps} \Q \Jg 
- \Cn \brak{\beta}
\geq  (\alpha \beta + \tfrac 12) \bigl( \tfrac{9}{10\eps} - \tfrac{13}{\eps} \Jg \bigr) 
- \tfrac{2 \cdot 250^2}{\eps} \Q \Jg  - \Cn \brak{\beta}
\,,
\end{equation}
and that
\begin{equation}
\sabs{ \Qb_2  \Sigma  g^{-\frac 12} \nbs_2 h } \leq \Cn \eps^2 
\,.
\end{equation}
\end{subequations}
From~\eqref{eq:vorticity:energy:0}--\eqref{justbecause0} we conclude that 
\begin{align}
&(1 - \Cn \eps^2)  \snorm{\tfrac{(\Jg \Q)^{\frac 12}}{ \Sigma^{\beta}}  \nbs^6\Upomega (\cdot,\s)}_{L^2_x}^2 
-  \snorm{\tfrac{(\Jg \Q)^{\frac 12}}{ \Sigma^{\beta}}  \nbs^6\Upomega  (\cdot,0)}_{L^2_x}^2 
+ \tfrac{9}{10\eps} \bigl( \alpha \beta + \tfrac 12 - \Cn \eps \brak{\beta} \bigr) 
\int_0^{\s} \snorm{\tfrac{1}{\Sigma^\beta} \nbs^6\Upomega (\cdot,\s')}_{L^2_x}^2 {\rm d}\s'   
\notag\\
&\leq \tfrac{33 (\alpha \beta + \frac 12) + 2 \cdot 250^2 (1+\alpha)}{\eps (1+\alpha)}  
\int_0^{\s} \snorm{\tfrac{(\Jg \Q)^{\frac 12}}{ \Sigma^{\beta}}  \nbs^6\Upomega (\cdot,\s')}_{L^2_x}^2  {\rm d}\s'
+ \tfrac{2}{\kappa_0} \int_0^\s  \snorm{\tfrac{1}{\Sigma^\beta} \nbs^6\Upomega (\cdot,\s')}_{L^2_x} \snorm{\tfrac{1}{ \Sigma^{\beta}}   \RR_\Upomega  (\cdot,\s')}_{L^2_x}  {\rm d}\s'  
\,,
\label{eq:vorticity:energy:1}
\end{align}
where, as usual, we have that $\Cn$ is a function only of $\alpha,\kappa_0$, and $\Cdata$, but is independent of $\beta$.

Next, we estimate $\snorm{\tfrac{1}{\Sigma^\beta} \RR_\Upomega}_{L^2_{x,\s}}$. The most delicate term is $\RR_\Upomega^{(1)}$. From H\"older's inequality and the bootstrap assumption~\eqref{bs-Sigma},  
\begin{equation*}
 \snorm{\tfrac{1}{\Sigma^\beta} \RR_\Upomega^{(1)}}_{L^2_{x,\s}} 
\les
\tfrac{1}{\eps} \snorm{\nbs \tfrac{\Jg}{\Sigma}}_{L^\infty_{x,\s}} \snorm{\tfrac{1}{\Sigma^\beta}  \nbs^6 \Upomega}_{L^2_{x,\s}}
+ \tfrac{1}{\eps}
(4 \kappa_0^{-1})^\beta \sum_{i=0}^{4} \snorm{\nbs^{5-i} \nbs\tfrac{\Jg}{\Sigma}}_{L^{\frac{10}{5-i}}_{x,\s}} \snorm{\nbs^{i}  \nbs_\s \Upomega}_{L^{\frac{10}{i}}_{x,\s}} 
\,,
\end{equation*}
where the implicit constant depends only on $\mathsf{C_{supp}}$, hence on $\alpha$ and $\kappa_0$. Using the initial data assumptions~\eqref{table:derivatives}, the bootstrap assumptions~\eqref{bootstraps}, and the bounds in  Lemmas~\ref{lem:time:interpolation}, we arrive at
\begin{align}
 \snorm{\tfrac{1}{\Sigma^\beta} \RR_\Upomega^{(1)}}_{L^2_{x,\s}} 
&\leq
\tfrac{\Cn}{\eps}  \snorm{\tfrac{1}{\Sigma^\beta}  \nbs^6 \Upomega}_{L^2_{x,\s}} 
+ \tfrac{\Cn}{\eps} (4 \kappa_0^{-1})^\beta \sum_{i=0}^{4} 
\bigl(\snorm{\nbs^{6}  \tfrac{\Jg}{\Sigma}}_{L^{2}_{x,\s}}^{\frac{5-i}{5}}  + \eps^{\frac{5-i}{5}}    \bigr)
\bigl((\kappa_0^\beta \snorm{\tfrac{1}{\Sigma^\beta} \nbs^{6}  \Upomega}_{L^{2}_{x,\s}})^{\frac{i}{5}}  +  \eps^{\frac{i}{5}}  \bigr)
\notag\\
&\leq\tfrac{\Cn}{\eps} \Bigl(
 \snorm{\tfrac{1}{\Sigma^\beta}  \nbs^6 \Upomega}_{L^2_{x,\s}} 
+ (\tfrac{4^5}{\kappa_0})^\beta   \snorm{\nbs^{6} \tfrac{\Jg}{\Sigma}}_{L^{2}_{x,\s}}
+ \eps      (\tfrac{4^5}{\kappa_0})^\beta \Bigr)
 \,,
 \label{eq:vorticity:reminder:1}
\end{align}
where as usual,  $\Cn = \Cn(\alpha,\kappa_0,\Cdata)$ is independent of $\beta$. 
For the remainder terms $\{\RR_\Upomega^{(i)}\}_{i=2,3,4}$ we can afford rougher upper bounds since no inverse powers of $\epsilon$ are present. As such, using the initial data assumptions~\eqref{table:derivatives}, the bootstrap assumptions~\eqref{bootstraps}, the Poincar\'e inequality~\eqref{eq:x1:Poincare}, the bounds in  Lemmas~\ref{lem:time:interpolation} and~\ref{lem:comm:tangent}, and the definition of $\Upomega$, we may thus estimate that
\begin{align}
&{\textstyle \sum}_{i=2}^{4} \snorm{\tfrac{1}{\Sigma^\beta} \RR_\Upomega^{(i)}}_{L^2_{x,\s}}
\leq 
(4 \kappa_0^{-1})^\beta {\textstyle \sum}_{i=2}^{4} \snorm{  \RR_\Upomega^{(i)}}_{L^2_{x,\s}}
\notag\\
&\leq \Cn (4 \kappa_0^{-1})^\beta
\Bigl(
 \eps \kappa_0^\beta \snorm{\tfrac{1}{\Sigma^\beta} \nbs^6 \Upomega}_{L^2_{x,\s}}
+ \eps  \snorm{\nbs^6 \tfrac{\Jg}{\Sigma} }_{L^2_{x,\s}}
+ \snorm{\nbs^6  V }_{L^2_{x,\s}}
+ \eps \snorm{\nbs^6 \Jg}_{L^2_{x,\s}}
+  \snorm{\nbs^6 \nbs_2 h}_{L^2_{x,\s}} + \eps\Bigr)
 \label{eq:vorticity:reminder:2}
\,.
\end{align}
Note that Lemmas~\ref{lem:Moser:tangent}, and the Poincar\'e inequality~\eqref{eq:x1:Poincare} and the bootstrap bounds give that 
\begin{equation}
\snorm{\nbs^6 \tfrac{\Jg}{\Sigma} }_{L^2_{x,\s}}
\les \snorm{\nbs^6 \Jg}_{L^2_{x,\s}} + \snorm{ \nbs^6 \Sigma^{-1}}_{L^2_{x,\s}} + \eps
 \label{eq:vorticity:reminder:3}
\,.
\end{equation}
By further appealing to Proposition~\ref{prop:geometry} and the bound $1 \leq \mathcal{J}^{-\frac 14}$, we deduce from \eqref{eq:vorticity:reminder:1}--\eqref{eq:vorticity:reminder:3} that 
\begin{equation}
\snorm{\tfrac{1}{\Sigma^\beta}  \RR_\Upomega}_{L^2_{x,\s}}
\leq 
\tfrac{\Cn}{\eps} (1 + \eps^2 4^\beta) \snorm{\tfrac{1}{\Sigma^{\beta}}  \nbs^6 \Upomega}_{L^2_{x,\s}}  
+ \Cn (4^5 \kappa_0^{-1})^\beta \brak{\mathsf{B_6}}
\,,
\label{eq:vorticity:remainder:4}
\end{equation}
where $\Cn = \Cn(\alpha,\kappa_0,\Cdata) \geq 1$ is a constant independent of $\beta$. 

Next, by combining \eqref{eq:vorticity:energy:1} and \eqref{eq:vorticity:remainder:4}, with the Cauchy-Young inequality we arrive at
\begin{align}
&(1 - \Cn \eps^2)  \snorm{\tfrac{(\Jg \Q)^{\frac 12}}{ \Sigma^{\beta}}  \nbs^6\Upomega (\cdot,\s)}_{L^2_x}^2 
-  \snorm{\tfrac{(\Jg \Q)^{\frac 12}}{ \Sigma^{\beta}}  \nbs^6\Upomega  (\cdot,0)}_{L^2_x}^2 
\notag\\
&\qquad 
+ \tfrac{9}{10\eps} \bigl( \alpha \beta + \tfrac 12 - \Cn \eps \brak{\beta} -  \Cn  (1 + \eps^2 4^\beta) \bigr) 
\int_0^{\s} \snorm{\tfrac{1}{\Sigma^\beta} \nbs^6\Upomega (\cdot,\s')}_{L^2_x}^2 {\rm d}\s'   
\notag\\
&\leq \tfrac{33 (\alpha \beta + \frac 12)+ 2 \cdot 250^2 (1+\alpha)}{\eps (1+\alpha)}  
\int_0^{\s} \snorm{\tfrac{(\Jg \Q)^{\frac 12}}{ \Sigma^{\beta}}  \nbs^6\Upomega (\cdot,\s')}_{L^2_x}^2  {\rm d}\s'
+ \Cn \eps (4^5 \kappa_0^{-1})^{2\beta} \brak{\mathsf{B_6}}^2      
\,,
\label{eq:vorticity:energy:2}
\end{align}
where $\Cn = \Cn(\alpha,\kappa_0,\Cdata) \geq 1$ is a constant independent of $\beta$. 
 Since $\Cn_{\eqref{eq:vorticity:energy:2}}$ is independent of $\beta$, we may first choose  $\beta$  to be sufficiently large (in terms of $\alpha,\kappa_0,\Cdata$), and then $\eps$ to be sufficiently small (in terms of $\alpha,\kappa_0,\Cdata$), to ensure that
\begin{equation}
\label{eq:vorticity:energy:beta:def}
\alpha \beta   -  \Cn_{\eqref{eq:vorticity:energy:2}} \geq 0,
\qquad
\mbox{and}
\qquad 
\tfrac 14   - \eps \Cn_{\eqref{eq:vorticity:energy:2}}  \brak{\beta} - \eps^2  \Cn_{\eqref{eq:vorticity:energy:2}}    4^\beta
\geq 0
\,.
\end{equation}
This makes $\beta = \beta(\alpha,\kappa_0,\Cdata)$.  With this choice of $\beta$, \eqref{eq:vorticity:energy:2} implies
\begin{align}
&\tfrac 12  \snorm{\tfrac{(\Jg \Q)^{\frac 12}}{ \Sigma^{\beta}}  \nbs^6\Upomega (\cdot,\s)}_{L^2_x}^2 
+ \tfrac{1}{8\eps}  
\int_0^{\s} \snorm{\tfrac{1}{\Sigma^\beta} \nbs^6\Upomega (\cdot,\s')}_{L^2_x}^2 {\rm d}\s'   
\notag\\
&\leq    \snorm{\tfrac{(\Jg \Q)^{\frac 12}}{ \Sigma^{\beta}}  \nbs^6\Upomega  (\cdot,0)}_{L^2_x}^2 
+ \tfrac{\Cn}{\eps}  
\int_0^{\s} \snorm{\tfrac{(\Jg \Q)^{\frac 12}}{ \Sigma^{\beta}}  \nbs^6\Upomega (\cdot,\s')}_{L^2_x}^2  {\rm d}\s'
+ \Cn  \eps (4^5 \kappa_0^{-1})^{2\beta} \brak{\mathsf{B_6}}^2      
\,,
\label{eq:vorticity:energy:3}
\end{align}
Using a standard Gr\"onwall argument,  and using the initial data bound provided by \eqref{table:derivatives}, we obtain from \eqref{eq:vorticity:energy:3} that 
\begin{equation*}
\sup_{\s\in[0,\eps]} \snorm{\tfrac{(\Jg \Q)^{\frac 12}}{\Sigma^{\beta}}   \nbs^6\Upomega (\cdot,\s)}_{L^2_x}^2 
+
 \tfrac{1}{ \eps}  \int_0^{\eps}\snorm{\tfrac{1}{\Sigma^{\beta}}   \nbs^6\Upomega (\cdot,\s)}_{L^2_x}^2 {\rm d}\s
\leq \Cn \eps  \brak{\mathsf{B_6}}^2
\,.
\end{equation*}
where $\Cn = \Cn(\alpha,\kappa_0,\Cdata)>0$ is a constant (the   $\beta$ dependence is included in the dependence on $ \alpha,\kappa_0,\Cdata$). The proof of \eqref{eq:svort:H6} is concluded upon multiplying the above estimate by $\kappa_0^{2\beta}$ and appealing to \eqref{bs-Sigma} and \eqref{Qd-lower-upper}--\eqref{Q-lower-upper}.

The first inequality in \eqref{eq:vorticity:pointwise} follows from \eqref{table:derivatives} since we may apply Proposition~\ref{prop:transverse:bounds} to the evolution \eqref{vort-t} of $\Upomega$, to deduce that $\|\Upomega\|_{L^\infty_{x,\s}} \leq 4 e^{18} \Cdata$.  The bootstrap \eqref{bs-Sigma} allows us to convert this into $\|\Omega\|_{L^\infty_{x,\s}} \leq 4 \cdot 4^{\frac{1}{\alpha}} \cdot e^{16} \Cdata$. 

The second inequality in \eqref{eq:vorticity:pointwise} follows in a similar manner, but we need to differentiate \eqref{vort-s} before applying Proposition~\ref{prop:transverse:bounds}. We have that 
\begin{equation}
 \Jg \Qd\p_\s (\nbs \Upomega) + \Jg\big( V+ \alpha\Sigma  g^{- {\frac{1}{2}} } \nbs_2 h \big) \nbs_2 (\nbs \Upomega)   - \tfrac{\alpha}{\eps}\Sigma \nbs_1 (\nbs\Upomega)
  = {\rm Error}  \,, \label{nbsUpomega}
\end{equation}
where
\begin{equation}
{\rm Error} = -  \Sigma \nbs (\tfrac{ \Jg}{\Sigma})  (\tfrac{1}{\eps} \nbs_\s   \Upomega  + V \nbs_2 \Upomega) 
- \Bigl(  \Jg \nbs  V  + \alpha \Sigma \nbs(\Jg g^{-\frac 12} \nbs_2h) \Bigr) \nbs_2 \Upomega  \,.  \label{nbsUpomega-Error}
\end{equation}
Notice that \eqref{nbsUpomega} is in the form of the equation  \eqref{eq:abstract:2:xs} with $\tilde f = \nbs\Upomega$ and with
$\tilde m \tilde f + \tilde q ={\rm Error}$.
Using the bootstrap assumptions~\eqref{bootstraps}, we have that  
\begin{equation}
\|\tilde m\|_{L^\infty_{x,\s}}
\leq \tfrac{1}{\eps} \snorm{\Sigma \nb(\tfrac{\Jg}{\Sigma})}_{L^\infty_{x,t}}
\le \tfrac{1}{\eps}  \bigl( 5(1+\alpha) + 18 \bigr) \leq \tfrac{23(1+\alpha)}{\eps} \,, \quad  \text{ and } \quad
\|\tilde q\|_{L^\infty_{x,t}} \le \Cn \eps \| \nb_2 \Upomega\|_{L^\infty_{x,\s}} 
\,.
\end{equation}
In turn, this means that the parameter $\beta$ appearing in \eqref{eq:abstract:transport:2} may be chosen to be  $\beta = \max\{\frac{20 \cdot 23 (1+\alpha)}{\alpha},1\} = \frac{20 \cdot 23 (1+\alpha)}{\alpha}$, which is a constant that depends only on  $\alpha$. Then, \eqref{eq:abstract:transport:2} leads to the bound
\begin{equation*}
 \snorm{\nb\Upomega}_{L^\infty_{x,\s}} 
\leq (4 e^{18})^\beta \snorm{\nb\Upomega(\cdot,0)}_{L^\infty_x} 
+\Cn\eps   \tfrac{20 \eps}{\alpha \beta} (4 e^{18})^\beta  \|\nb_2 \Upomega\|_{L^\infty_{x,t}}
\,.
\end{equation*}
Taking $\eps$ to be sufficiently small to absorb the second term on the right side  into the left side, concludes the proof.
\end{proof}

\subsection{Improved estimates for $\Abn$}
The boundedness of the damping norms  $\tilde{\mathcal{D}}_6$  and $\tilde{\mathcal{D}}_5$ assumed in the bootstraps \eqref{bootstraps-Dnorm:6} and \eqref{bootstraps-Dnorm:5} implies that the $L^2_{x,\s}$ norm of $\mathcal{J}^{\frac 14} \Jgh \nbs^6(\Jg \Abn)$ and $\nbs^5(\Jg \Abn)$ are bounded by $\mathsf{B_6}$ and respectively $\mathsf{B_5}$. The goal of this subsection is to show that these bounds for $\Abn$ are greatly improved (one fewer powers of $\Jg$, and one gain of $\eps$) as a consequence of the  vorticity estimate obtained earlier in Proposition~\ref{prop:vort:H6}.  
\begin{corollary}
\label{cor:Abn:improve}
Under the standing bootstrap assumptions, we have that 
\begin{subequations}
\label{eq:Jg:Abn:D5:improve}
\begin{align}
\snorm{\Jgh  \nbs^5 \Abn}_{L^\infty_\s L^2_{x}}
&\les \eps^{\frac 12} \mathsf{K} \brak{\mathsf{B_6}}
\label{eq:Jg:Abn:D5:improve:a}
\,, \\
\snorm{{\mathcal J}^{\! \frac 34} \Jgh  \nbs^6 \Abn}_{L^\infty_\s L^2_{x}}
&\les \eps^{\frac 12} \mathsf{K} \brak{\mathsf{B_6}}
\label{eq:Jg:Abn:D6:improve:b}
\,, \\
\snorm{\nbs^5 \Abn}_{L^2_{x,\s}} 
&\les 
\eps \mathsf{K} \brak{\mathsf{B_6}}
\label{eq:Jg:Abn:D5:improve:c}
\,, \\
\snorm{{\mathcal J}^{\! \frac 14} \Jgh \nbs^6  \Abn}_{L^2_{x,\s}} 
&\les 
\eps \mathsf{K} \brak{\mathsf{B_6}}
\label{eq:Jg:Abn:D6:improve:d}
\,.
\end{align}
\end{subequations}
where the implicit constant depends only on $\alpha,\kappa_0$, and $\Cdata$.
\end{corollary}
\begin{proof}[Proof of Corollary~\ref{cor:Abn:improve}]
We first prove \eqref{eq:Jg:Abn:D6:improve:b} and  \eqref{eq:Jg:Abn:D6:improve:d}. 
From \eqref{vort-id-good}, \eqref{eq:tilde:E6}, \eqref{eq:tilde:D6},   \eqref{bootstraps-Dnorm:6}, \eqref{bootstraps-Dnorm:5}, \eqref{eq:vort:H6} and the bounds $\mathcal{J}\leq 1$ and $\Jg \leq \frac 32$,  we deduce that  
\begin{align*}
\snorm{{\mathcal J}^{\! \frac 34} \Jgh  \nbs^6 \Abn}_{L^\infty_\s L^2_{x}}
& \les
\snorm{ \Jgh  \nbs^6 \Omega}_{L^\infty_\s L^2_{x}}
+
\snorm{{\mathcal J}^{\! \frac 34} \Jgh  \nbs^6 (  \Wbt,  \Zbt)}_{L^\infty_\s L^2_{x}}
\les \eps^\frac 12  \brak{\mathsf{B_6}} +  \mathsf{K} \eps \brak{\mathsf{B_6}} \eps^{-\frac 12} \,, \\
\snorm{{\mathcal J}^{\! \frac 14} \Jgh \nbs^6  \Abn }_{L^2_{x,\s}} 
& \les \snorm{ \nbs^6     \Omega }_{L^2_{x,\s}}  + \snorm{{\mathcal J}^{\! \frac 14} \Jgh \nbs^6    (  \Wbt,  \Zbt)}_{L^2_{x,\s}} 
\les \eps   \brak{\mathsf{B_6}}  + \mathsf{K} \eps \brak{\mathsf{B_6}}\,,
\end{align*}
thereby proving  \eqref{eq:Jg:Abn:D6:improve:b} and  \eqref{eq:Jg:Abn:D6:improve:d}, since $\mathcal{J} \leq \Jgb \leq \Jg$.

Next, we note that  \eqref{eq:Jg:Abn:D5:improve:c} immediately follows  from \eqref{eq:Jg:Abn:D6:improve:d} and  \eqref{eq:r:r':L2:time}. Indeed, we may apply \eqref{eq:r:r':L2:time} with $r'=0$ and $r=\frac 34$,  and recall that $\p_\s = \frac{1}{\eps \Qd} \nbs_\s$  so that with \eqref{table:derivatives}, \eqref{Qd-lower-upper}, \eqref{eq:Jg:Abn:D6:improve:d}, and $F = \nbs^5 \Abn$, we obtain
\begin{equation*}
 \| \nbs^5 \Abn \|_{L^2_{x,\s}} 
 \leq \eps^{\frac 12}  \| \nbs^5 \Abn(\cdot,0)\|_{L^2_x} 
 +2   \snorm{\Qd^{-1}}_{L^\infty_{x,\s}}  \| \mathcal{J}^{\!\frac 34} \nbs_\s \nbs^5 \Abn \|_{L^2_{x,\s}} 
 \les \eps \mathsf{K} \brak{\mathsf{B}_6}
 \,,
\end{equation*}
thereby establishing \eqref{eq:Jg:Abn:D5:improve:c}. Similarly,~\eqref{eq:Jg:Abn:D5:improve:a} follows from \eqref{eq:Jg:Abn:D6:improve:d} and  \eqref{eq:r:Linfty:time:2}. Applying  \eqref{eq:r:Linfty:time:2} with $r=\frac 34$ and $F = \nbs^5 \Abn$ we arrive at
\begin{equation*}
 \sup_{\s \in [0,\eps]} \|\Jgh \nbs^5 \Abn (\cdot,\s)\|_{L^2_{x}} 
 \les   \| \nbs^5 \Abn(\cdot,0)\|_{L^2_x} 
 + \eps^{-\frac 12} \| \mathcal{J}^{\!\frac 14} \Jgh \nbs_\s \nbs^5 \Abn \|_{L^2_{x,\s}} 
 \les \eps^{\frac 12} \mathsf{K} \brak{\mathsf{B}_6}
 \,,
\end{equation*}
concluding the proof of the lemma.
\end{proof}

\subsection{Improved estimate for $\Zbn$}
Since $\Zbn$ does not appear in the definition of the vorticity, we take a different approach to improve the bounds for the derivatives of 
$\Zbn$ in $L^2_{x,\s}$.   By using an argument similar to that in the proof of  Proposition~\ref{prop:vort:H6}, and by appealing to a key step in
the bound \eqref{eq:Jg:Abn:D5:improve}, we are able to show that when $\eps$ is sufficiently small (in terms of $\alpha,\kappa_0,\Cdata$), we have the following: 
\begin{lemma}
\label{lem:Zbn:improve}
Under the  assumptions of Proposition~\ref{prop:vort:H6}, we have that  for any \bubu{$\bar \beta \geq 0$}, and $\bar a \in [0,\frac 12]$, 
\begin{subequations}
\label{eq:Jg:Zbn:D5:improve}
\begin{align}
\snorm{ \Jg^{\!\frac 12}  \nbs^5 \Zbn }_{L^\infty _\s L^2_x}^2
+ \tfrac{1}{\eps} \int_0^\eps
\snorm{\nbs^5  \Zbn (\s)}_{L^2_{x}}^2 {\rm d}\s
&\les \eps \brak{\mathsf{B}_6}^2 \,, 
\label{eq:Jg:Zbn:D5:improve:a}\\
\snorm{\nbs^4 \Zbn }_{L^\infty _\s L^2_x}
&\les \eps^{\frac{1}{2}}   \brak{\mathsf{B}_6} \,, 
\label{eq:Jg:Zbn:D3:improve} \\
\snorm{\Sigma^{-\bar \beta} \mathcal{J}^{\frac 34 - \bar a} \Jg^{\! \bar a} \nbs_1 \nbs^5 \Zbn }_{L^2_{x,\s}} 
&\le \tfrac{1}{2\alpha} \| \Sigma^{-1-\bar \beta} \mathcal{J}^{\frac 34 - \bar a} \Jg^{\! \bar a} \nbs^5 \nbs_\s (\Jg \Zbn) \|_{L^2_{x,\s}}
+ \Cn \eps  (\tfrac{4}{\kappa_0})^{\bar \beta}    \mathsf{K} \brak{\mathsf{B_6}}\,,
\label{eq:Jg:Zbn:p1D5:improve}  \\
\snorm{\Sigma^{-\bar \beta} \mathcal{J}^{\frac 34 - \bar a} \Jg^{\! \bar a}  \nbs_2 \nbs^5 \Zbn }_{L^2_{x,\s}}  
&\le \tfrac{1}{\alpha  \eps} \|\Sigma^{-1 - \bar \beta } \mathcal{J}^{\frac 34 - \bar a} \Jg^{\! \bar a} \nbs^5 \nbs_\s \Zbt \|_{L^2_{x,\s}} + \Cn \eps (\tfrac{4}{\kappa_0})^{\bar \beta} \mathsf{K} \brak{\mathsf{B_6}}  
\label{eq:Jg:Zbn:p2D5:improve:new} \,,
\end{align}
where the implicit constant and the constant $\Cn$ depend only on $\alpha,\kappa_0$, and $\Cdata$. 
It is also convenient to record the estimates 
\begin{align}
\snorm{\mathcal{J}^{\frac 34} \Jgh \nbs_1 \nbs^5 \Zbn }_{L^\infty_\s L^2_{x}} 
&\les \eps^{-\frac 12} \brak{\mathsf{B_6}}\,,
\label{eq:Jg:Zbn:p1D6:sup:improve:new}  \\
\snorm{\mathcal{J}^{\frac 34} \Jgh \nbs_2 \nbs^5 \Zbn }_{L^\infty_\s L^2_{x}}  
&\les  \eps^{-\frac 12}   \brak{\mathsf{B_6}} \,,
\label{eq:Jg:Zbn:p2D6:sup:improve:new}
\end{align}
\end{subequations}
where the implicit constant and the constant depends only on $\alpha,\kappa_0$, and $\Cdata$. 
\end{lemma}
\begin{proof}[Proof of Lemma~\ref{lem:Zbn:improve}]
We first prove \eqref{eq:Jg:Zbn:D5:improve:a}.  Letting $\nbs^5$ act on \bubu{the $(x,\s)$-variable form of \eqref{eq:Zb:nn} which has been pre-multiplied by $\frac{\Jg}{\Sigma}$},  and find
that
\begin{equation}
\tfrac{\Jg}{\Sigma}  (\Q \p_\s + V \p_2) (\nbs^5 \Zbn) - 2\alpha  (\nbs^5 \Zbn),_1 + 2 \alpha   \Jg g^{-\frac 12} \nbs_2 h \nbs_2 (\nbs^5 \Zbn)  
= 
- \tfrac{1}{\Sigma} (\nbs^5 \Zbn) \bubu{(\tfrac{1-\alpha}{2} \Jg \Wbn + \tfrac{1+\alpha}{2} \Jg \Zbn)}
+ \mathsf{Err_Z} \,,
\label{eq:Zbn:improv:1} 
\end{equation}
where the error term $\mathsf{Err_Z}$ is given by 
\begin{align}
\mathsf{Err_Z}
&= 
-  \jump{\nbs^5, \tfrac{\Jg}{\Sigma}} (\tfrac{1}{\eps} \nbs_\s + V \nbs_2) \Zbn
-  \tfrac{\Jg}{\Sigma}  \jump{\nbs^5, V} \nbs_2 \Zbn
- 2\alpha    \jump{\nbs^5, \Jg g^{-\frac 12} \nbs_2 h} \nbs_2 \Zbn
\notag\\
&
-   \jump{\nbs^5, \tfrac{1}{\Sigma} \bubu{(\tfrac{1-\alpha}{2} \Jg \Wbn + \tfrac{1+\alpha}{2} \Jg \Zbn)}} \Zbn
- \tfrac{\alpha}{2 } \nbs^5 \bigl( (\Jg \Zbn) \bubu{(\tfrac{1}{\Sigma} \Abt+ g^{-\frac 32} \nbs_2^2 h )} \bigr)
\notag\\
&  
+ \tfrac{\alpha}{2 } \nbs^5 \bigl( (\Jg \Wbn) \bubu{(\tfrac{1}{\Sigma} \Abt - g^{-\frac 32} \nbs_2^2 h )}\bigr)
+ \alpha  \nbs^5 \bigl(   g^{-\frac 12} (\Jg \Abn),_2 \bigr)
+2 \alpha  \nbs^5 \bigl( \bubu{g^{-\frac 12}  ( \Zbt + \tfrac 12 \Abn)} \nbs_2 \Jg \bigr)
\notag\\
&  
-  \nbs^5 \bigl(\bubu{\tfrac{1}{\Sigma}} (\Jg \Abn) (\tfrac{1+\alpha}{2} \Wbt + \tfrac{3-\alpha}{2} \Zbt) \bigr)
-  \nbs^5 \bigl(\bubu{\tfrac{1}{\Sigma}} (\Jg \Zbt) (\tfrac{1+\alpha}{2} \Wbt + \tfrac{1-\alpha}{2} \Zbt) \bigr)
+  \alpha  \nbs^5 \bigl( \bubu{(\Jg \Abt)}   g^{-\frac 32} \nbs_2^2 h  \bigr)
\notag\\
&=: {\textstyle \sum}_{i=1}^{11}\mathsf{Err_Z}^{(i)}
\label{eq:Zbn:improv:2}
\,.
\end{align}
The first three terms appearing in \eqref{eq:Zbn:improv:2} should be compared to the remainder terms in \eqref{D4-vort-s-remainder}.
Next, we test \eqref{eq:Zbn:improv:1}  with $\Sigma^{-2\beta+1} \nbs^5 \Zbn$, for $\beta>0$ to be determined, and similarly to \eqref{eq:vorticity:energy:0} we arrive at
\begin{align}
&  \snorm{\tfrac{(\Jg \Q)^{\frac 12}}{ \Sigma^{\beta}}  \nbs^5\Zbn (\cdot,\s)}_{L^2_x}^2 
-  \snorm{\tfrac{(\Jg \Q)^{\frac 12}}{ \Sigma^{\beta}}  \nbs^5\Zbn  (\cdot,0)}_{L^2_x}^2 
-   \tint \sabs{\nbs^5\Zbn}^2 (\Q \p_s + V \p_2) \tfrac{\Jg}{\Sigma^{2\beta}}  
- 2 \alpha (2\beta-1) \tint \sabs{\nbs^5\Zbn}^2 \tfrac{\Sigma,_1}{\Sigma^{2\beta}} 
\notag\\
&=   \tint \sabs{\nbs^5\Zbn}^2 \tfrac{\Jg}{\Sigma^{2\beta}} \bigl(\Qr_\s -  V \Qr_2 + \nbs_2 V  \bigr) 
+ 2 \alpha  \tint \sabs{\nbs^5\Zbn}^2 \bigl(\nbs_2 - \Qr_2\bigr) \tfrac{\Jg  g^{-\frac 12} \nbs_2 h}{\Sigma^{2\beta-1}}
\notag\\
&\quad
+ 2 \alpha  \int \sabs{\nbs^5 \Zbn}^2 \Qb_2 \tfrac{\Jg  g^{-\frac 12} \nbs_2 h}{\Sigma^{2\beta-1}} \Big|_{\s}
-   \tint\tfrac{2}{\Sigma^{2\beta}} |\nbs^5 \Zbn|^2 \bubu{(\tfrac{1-\alpha}{2} \Jg \Wbn + \tfrac{1+\alpha}{2} \Jg \Zbn)}
+  \tint\tfrac{2}{\Sigma^{2\beta-1}} \nbs^5 \Zbn \mathsf{Err_Z}  
\,.
\end{align}
By appealing to  \eqref{p1-Sigma-s} and \eqref{Sigma0-ALE-s} the above becomes
\begin{align}
&  \snorm{\tfrac{(\Jg \Q)^{\frac 12}}{ \Sigma^{\beta}}  \nbs^5\Zbn (\cdot,\s)}_{L^2_x}^2 
-  \snorm{\tfrac{(\Jg \Q)^{\frac 12}}{ \Sigma^{\beta}}  \nbs^5\Zbn  (\cdot,0)}_{L^2_x}^2 
+  \tint \sabs{\nbs^5\Zbn}^2   \tfrac{1}{\Sigma^{2\beta}}  \mathsf{G}
\notag\\
&=   
 2 \alpha  \int \sabs{\nbs^5 \Zbn}^2 \Qb_2 \tfrac{\Jg  g^{-\frac 12} \nbs_2 h}{\Sigma^{2\beta-1}} \Big|_{\s}
+  \tint\tfrac{2}{\Sigma^{2\beta-1}} \nbs^5 \Zbn \mathsf{Err_Z}  
\,,
\end{align}
where
\begin{align}
\mathsf{G}
&:=  
 \Bigl( \bigl( \bubu{\tfrac{1-3\alpha}{2}} \Jg \Wbn + \bubu{\tfrac{1+3\alpha}{2}} \Jg \Zbn\bigr) - 2 \alpha \beta  \Jg   \bigl( \Zbn + \Abt\bigr)\Bigr)
- 2\alpha (\beta-\tfrac 12) \bigl(  \Jg \Wbn -  \Jg \Zbn + \Jg \nbs_2 h (\Wbt - \Zbt) \bigr)
\notag\\
&\qquad 
- \Jg \bigl(\Qr_\s -  V \Qr_2 + \nbs_2 V  \bigr) 
+  2 \alpha (2\beta-1) \Jg  g^{-\frac 12} \nbs_2 h \nbs_2 \Sigma 
+ 2 \alpha  \Qr_2 \Sigma \Jg  g^{-\frac 12} \nbs_2 h 
\end{align}
By using the pointwise bootstrap inequalities in \eqref{bootstraps} and~\eqref{eq:Q:all:bbq}, using the bound \eqref{eq:signed:Jg}, upon taking 
\begin{equation} 
\beta > \tfrac{2+3\alpha}{2\alpha}\,, \label{betachoice1}
\end{equation} 
and $\eps$ sufficiently small with respect to $\alpha,\kappa_0$, and $\Cdata$ (but not $\beta$), we may derive the  bounds
\begin{subequations} 
\label{justbecause1}
\begin{align}
&\mathsf{G} 
\geq - (2 \alpha \beta - \alpha - \bubu{\tfrac{1-3\alpha}{2}})  \Jg \Wbn - \tfrac{2 \cdot 250^2}{\eps} \Q \Jg - \Cn \brak{\beta}
\geq  (\alpha \beta + \tfrac 12) \bigl( \tfrac{9}{10\eps} - \tfrac{13}{\eps} \Jg \bigr) - \tfrac{2 \cdot 250^2}{\eps} \Q \Jg  - \Cn \brak{\beta}
\,,\\
&\sabs{ \Qb_2  \Sigma  g^{-\frac 12} \nbs_2 h } \leq \Cn \eps^2 
\,.
\end{align}
\end{subequations}
We note that the bounds in \eqref{justbecause1} are identical to those in \eqref{justbecause0}, and hence, as in \eqref{eq:vorticity:energy:1} we deduce
\begin{align}
&(1 - \Cn \eps^2)  \snorm{\tfrac{(\Jg \Q)^{\frac 12}}{ \Sigma^{\beta}}  \nbs^5\Zbn (\cdot,\s)}_{L^2_x}^2 
-  \snorm{\tfrac{(\Jg \Q)^{\frac 12}}{ \Sigma^{\beta}}  \nbs^5\Zbn  (\cdot,0)}_{L^2_x}^2 
+ \tfrac{9}{10\eps} \bigl( \alpha \beta + \tfrac 12 - \Cn \eps \brak{\beta} \bigr) 
\int_0^{\s} \snorm{\tfrac{1}{\Sigma^\beta} \nbs^5\Zbn (\cdot,\s')}_{L^2_x}^2 {\rm d}\s'   
\notag\\
&\leq \tfrac{33 (\alpha \beta + \frac 12) + 2 \cdot 250^2(1+\alpha)}{\eps (1+\alpha)}  
\int_0^{\s} \snorm{\tfrac{(\Jg \Q)^{\frac 12}}{ \Sigma^{\beta}}  \nbs^5\Zbn (\cdot,\s')}_{L^2_x}^2  {\rm d}\s'
+ \tfrac{2}{\kappa_0} \int_0^\s  \snorm{\tfrac{1}{\Sigma^\beta} \nbs^5\Zbn (\cdot,\s')}_{L^2_x} \snorm{\tfrac{1}{ \Sigma^{\beta}}  \mathsf{Err_Z}  (\cdot,\s')}_{L^2_x}  {\rm d}\s'  
\,,
\label{eq:Zbn:improv:3}
\end{align}
where $\Cn$ is a function only of $\alpha,\kappa_0$, and $\Cdata$, but is independent of $\beta$.
It thus remains to bound  the error term $\mathsf{Err_Z}$ appearing in \eqref{eq:Zbn:improv:2}, in $L^2_{x,\s}$. We note that the  first three terms in the definition of 
$\mathsf{Err_Z}$ are in direct correspondence with the three terms in the remainder $\RR_\Upomega$ from \eqref{D4-vort-s-remainder}; 
as such, these terms are estimated in a nearly identical fashion, leading to a bound analogous to \eqref{eq:vorticity:remainder:4}, namely
\begin{equation}
{\textstyle \sum}_{i=1}^{3} \snorm{\tfrac{1}{\Sigma^\beta}\mathsf{Err_Z}^{(i)}}_{L^2_{x,\s}} 
\leq \tfrac{\Cn}{\eps} (1 + \eps^2 4^\beta) \snorm{\tfrac{1}{\Sigma^{\beta}}  \nbs^5 \Zbn}_{L^2_{x,\s}}  
+ \Cn (4^5 \kappa_0^{-1})^\beta \brak{\mathsf{B_6}} \,.
\label{eq:Zbn:improv:4}
\end{equation}
The remaining terms  on the right side of \eqref{eq:Zbn:improv:2}, namely $\{\mathsf{Err_Z}^{(i)}\}_{i=4}^{11}$, are bounded similarly, using the Moser-type inequality in Lemma~\ref{lem:Moser:tangent}, the 
Gagliardo-Nirenberg inequality in Lemma~\ref{lem:time:interpolation}, the bootstrap inequalities~\eqref{bootstraps}, the  bounds in 
Proposition~\ref{prop:geometry}, and the above established estimates~\eqref{eq:Jg:Abn:D5:improve}. 
We shall detail the estimates for the two most difficult terms: $\mathsf{Err_Z}^{(4)}$ and $\mathsf{Err_Z}^{(7)}$.  For the fourth error term, we have that  
\begin{align}
\snorm{\tfrac{1}{\Sigma^\beta}\mathsf{Err_Z}^{(4)}}_{L^2_{x,\s}} 
&\leq (4 \kappa_0^{-1})^{\beta}
\snorm{\jump{\nbs^5, \tfrac{1}{\Sigma} (\bubu{\tfrac{1-\alpha}{2}} \Jg \Wbn + \bubu{\tfrac{1+\alpha}{2}} \Jg \Zbn) } \Zbn}_{L^2_{x,\s}} 
\notag\\
&  \leq \Cn (4 \kappa_0^{-1})^{\beta}\sum_{i=0}^4 \snorm{\nbs^{5-i} (\tfrac{1}{\Sigma} (\bubu{\tfrac{1-\alpha}{2}} \Jg \Wbn + \bubu{\tfrac{1+\alpha}{2}} \Jg \Zbn))}_{L^{\frac{10}{5-i}}_{x,\s}} \snorm{\nbs^i \Zbn}_{L^{\frac{10}{i}}_{x,\s}}
\notag\\
&\leq \Cn (4\kappa_0^{-1})^\beta \sum_{i=0}^4
\bigl(\snorm{\nbs^{5} (\tfrac{1}{\Sigma}(\Jg\Wbn,\Jg\Zbn))}_{L^{2}_{x,\s}}^{\frac{5-i}{5}}  \eps^{-\frac{i}{5}}+ \eps^{-\frac{i}{5}}    \bigr)
\bigl((\kappa_0^\beta \snorm{\tfrac{1}{\Sigma^\beta} \nbs^{5}  \Zbn}_{L^{2}_{x,\s}})^{\frac{i}{5}}  +  \eps^{\frac{i}{5}}  \bigr)
\notag\\
&   \leq \tfrac{\Cn}{\eps}  \snorm{\tfrac{1}{\Sigma^\beta}  \nbs^5 \Zbn }_{L^2_{x,\s}}
+ (\tfrac{4^5}{\kappa_0})^\beta \brak{\mathsf{B}_6}\,,
\label{eq:Zbn:improv:4b}
\end{align}
and for the seventh error term, we obtain that
\begin{equation}
\snorm{\tfrac{1}{\Sigma^\beta}\mathsf{Err_Z}^{(7)}}_{L^2_{x,\s}}
\leq \snorm{\tfrac{1}{\Sigma^{\beta}} \nbs^5 (  g^{-\frac 12} \nbs_2 (\Jg \Abn)) }_{L^2_{x,\s}} 
 \leq 
\Cn \eps (4\kappa_0^{-1})^\beta  \brak{\mathsf{B_6}} 
 \,.
\label{eq:Zbn:improv:4c}
\end{equation}
A straightforward application of \eqref{bootstraps} and \eqref{eq:Lynch:2}  shows that all of the remaining terms on the right side of  \eqref{eq:Zbn:improv:2} satisfy even better bounds, 
 which combined with  \eqref{eq:Zbn:improv:4}--\eqref{eq:Zbn:improv:4c} leads to 
\begin{equation}
\snorm{\Sigma^{-\beta} \mathsf{Err_Z}}_{L^2_{x,\s}}
\leq  \tfrac{\Cn}{\eps} (1 + \eps^2 4^\beta) \snorm{\tfrac{1}{\Sigma^{\beta}}  \nbs^5 \Zbn}_{L^2_{x,\s}}  
+ \Cn (4^5 \kappa_0^{-1})^\beta \brak{\mathsf{B_6}} 
\label{eq:Zbn:improv:4d}
\end{equation}
where $\Cn = \Cn(\alpha,\kappa_0,\Cdata)$ is a suitably large constant, which is independent of $\beta$.
Inserting the bound \eqref{eq:Zbn:improv:4d} into \eqref{eq:Zbn:improv:3}, similarly to \eqref{eq:vorticity:energy:2} we are lead to 
\begin{align}
&(1 - \Cn \eps^2)  \snorm{\tfrac{(\Jg \Q)^{\frac 12}}{ \Sigma^{\beta}}  \nbs^5\Zbn (\cdot,\s)}_{L^2_x}^2 
-  \snorm{\tfrac{(\Jg \Q)^{\frac 12}}{ \Sigma^{\beta}}  \nbs^5\Zbn  (\cdot,0)}_{L^2_x}^2 
\notag\\
&\qquad 
+ \tfrac{9}{10\eps} \bigl( \alpha \beta + \tfrac 12 - \Cn \eps \brak{\beta} -  \Cn  (1 + \eps^2 4^\beta) \bigr) 
\int_0^{\s} \snorm{\tfrac{1}{\Sigma^\beta} \nbs^5\Zbn (\cdot,\s')}_{L^2_x}^2 {\rm d}\s'   
\notag\\
&\leq \tfrac{33 (\alpha \beta + \frac 12)+2 \cdot 250^2(1+\alpha)}{\eps (1+\alpha)}  
\int_0^{\s} \snorm{\tfrac{(\Jg \Q)^{\frac 12}}{ \Sigma^{\beta}}  \nbs^5\Zbn (\cdot,\s')}_{L^2_x}^2 {\rm d}\s'
+ \Cn \eps (4^5 \kappa_0^{-1})^{2\beta} \brak{\mathsf{B_6}}^2      
\,,
\label{eq:Zbn:improv:5}
\end{align}
where $\Cn = \Cn(\alpha,\kappa_0,\Cdata) \geq 1$ is a constant independent of $\beta$. 
Since $\Cn_{\eqref{eq:Zbn:improv:5}}$ is independent of $\beta$, we may  choose first $\beta$ to be sufficiently large (in terms of $\alpha,\kappa_0,\Cdata$), and then $\eps$ to be sufficiently small (in terms of $\alpha,\kappa_0,\Cdata$), to ensure that
\begin{equation}
\alpha \beta   -  \Cn_{\eqref{eq:Zbn:improv:5}} \geq 0,
\qquad
\mbox{and}
\qquad 
\tfrac 14   - \eps \Cn_{\eqref{eq:Zbn:improv:5}}  \brak{\beta} - \eps^2  \Cn_{\eqref{eq:Zbn:improv:5}}    4^\beta
\geq 0
\,.
\label{betachoice2}
\end{equation}
The choice \eqref{betachoice1} and \eqref{betachoice2} makes $\beta = \beta(\alpha,\kappa_0,\Cdata)$.  With this choice of $\beta$, \eqref{eq:Zbn:improv:5} implies
\begin{align}
&\tfrac 12  \snorm{\tfrac{(\Jg \Q)^{\frac 12}}{ \Sigma^{\beta}}   \nbs^5\Zbn (\cdot,\s)}_{L^2_x}^2 
+ \tfrac{1}{8\eps}  
\int_0^{\s} \snorm{\tfrac{1}{\Sigma^\beta}  \nbs^5\Zbn (\cdot,\s')}_{L^2_x}^2 {\rm d}\s'   
\notag\\
&\leq    \snorm{\tfrac{(\Jg \Q)^{\frac 12}}{ \Sigma^{\beta}}   \nbs^5\Zbn (\cdot,0)}_{L^2_x}^2 
+ \tfrac{\Cn}{\eps}  
\int_0^{\s} \snorm{\tfrac{(\Jg \Q)^{\frac 12}}{ \Sigma^{\beta}}   \nbs^5\Zbn (\cdot,\s')}_{L^2_x}^2  {\rm d}\s'
+ \Cn  \eps (4^5 \kappa_0^{-1})^{2\beta} \brak{\mathsf{B_6}}^2      
\,,
\label{eq:Zbn:improv:6}
\end{align}
Using a standard Gr\"onwall argument,  and using the initial data bound provided by \eqref{table:derivatives}, we obtain from \eqref{eq:Zbn:improv:6} that 
\begin{equation}
\sup_{\s\in[0,\eps]} \snorm{\tfrac{(\Jg \Q)^{\frac 12}}{\Sigma^{\beta}}    \nbs^5\Zbn(\cdot,\s)}_{L^2_x}^2 
+
 \tfrac{1}{ \eps}  \int_0^{\eps}\snorm{\tfrac{1}{\Sigma^{\beta}}    \nbs^5\Zbn (\cdot,\s)}_{L^2_x}^2 {\rm d}\s
\leq \Cn \eps  \brak{\mathsf{B_6}}^2
\,.
\end{equation}
where $\Cn = \Cn(\alpha,\kappa_0,\Cdata)$   (the   $\beta$ dependence is included in the dependence on $ \alpha,\kappa_0,\Cdata$). The proof of \eqref{eq:Jg:Zbn:D5:improve:a} is concluded upon multiplying the above estimate by $\kappa_0^{2\beta}$ and appealing to \eqref{bs-Sigma} and \eqref{Qd-lower-upper}--\eqref{Q-lower-upper}.

The bound \eqref{eq:Jg:Zbn:D3:improve} follows directly from  \eqref{eq:Jg:Zbn:D5:improve:a} by the fundamental theorem of calculus; in 
particular, using \eqref{eq:sup:in:time:L2}, we obtain that
\begin{equation}
\| \nbs^4 \Zbn\|_{L^\infty_\s L^2_x} 
\les 
\| \nbs^4 \Zbn (\cdot,0)\|_{L^2_x} 
+ 
\eps^{-\frac 12} \| \nbs^5 \Zbn\|_{L^2_{x,\s}}
\les 
\eps^{\frac 12} \Cdata + \eps^{-\frac 12} \cdot \eps \brak{\mathsf{B_6}}
\,,
\end{equation}
which proves \eqref{eq:Jg:Zbn:D3:improve}.

Next, we turn to the improved estimates for $\nbs^6 \Zbn$ stated in \eqref{eq:Jg:Zbn:p1D5:improve}--\eqref{eq:Jg:Zbn:p2D5:improve:new}. We emphasize that the parameter $\bar \beta\geq 0$ appearing in these estimates is arbitrary, it is not the same as the $\beta$ appearing in \eqref{eq:Zbn:improv:6}. Note that we are only able at this point to obtain such an improved estimate if either $\nbs^6 = \nbs^5 \nbs_1$, or $\nbs^6 = \nbs^5 \nbs_2$, i.e., if at least one space derivative is present. To see this, we return to \eqref{Sigma0-ALE}, which we first differentiate in space, and then convert to $(x,\s)$-variables, to obtain
\begin{equation}
\nbs_i \Zbn
= -\nbs_i \Abt - \tfrac{1}{\alpha \Sigma} (\tfrac{1}{\eps} \nbs_\s + V \nbs_2)  \nbs_i \Sigma 
- \tfrac{1}{\Sigma}  (   \Zbn +   \Abt  ) \nbs_i \Sigma
- \tfrac{1}{\alpha \Sigma}  \nbs_i V \nbs_2 \Sigma
\,, \qquad \mbox{for} \qquad i \in \{1,2\} \,.
\end{equation}
Upon applying $\nbs^5$ to the above identity, and recalling  from \eqref{grad-Sigma} that 
\begin{equation}
\nbs_1 \Sigma =  \tfrac{\eps}{2} \Jg(\Wbn -\Zbn)  + \tfrac{\eps}{2} \Jg \nbs_2 h (\Wbt -\Zbt)
\,, \qquad \mbox{and} \qquad
\nbs_2 \Sigma = \tfrac{1}{2} g^{{\frac{1}{2}} } (\Wbt-\Zbt)
\,,
\end{equation}
we deduce  
\begin{subequations}
\label{eq:Zn:nb:Sigma:identity}
\begin{align}
\nbs^5 \nbs_1 \Zbn
&= - \nbs^5 \nbs_1 \Abt 
- \tfrac{1}{2\alpha  }  \Sigma^{-1} \nbs^5 \nbs_\s   \bigl( \Jg \Wbn - \Jg \Zbn + \Jg \nbs_2 h(\Wbt-\Zbt) \bigr)
- \tfrac{1}{2\alpha} V \Sigma^{-1}  \nbs^5  \nbs_1 \bigl(g^{{\frac{1}{2}} } (\Wbt-\Zbt)\bigr)
\notag\\
&\qquad 
- \tfrac{1}{2\alpha  } \jump{\nbs^5,  \Sigma^{-1}} \nbs_\s     \bigl( \Jg \Wbn - \Jg \Zbn + \Jg \nbs_2 h(\Wbt-\Zbt) \bigr)
- \tfrac{1}{\alpha} \jump{\nbs^5, V \Sigma^{-1}}\nbs_1 \nbs_2 \Sigma
\notag\\
&\qquad 
- \nbs^5 \bigl(\Sigma^{-1} (   \Zbn +   \Abt  )  \nbs_1 \Sigma\bigr)
- \tfrac{1}{\alpha}  \nbs^5 \bigl(\Sigma^{-1} \nbs_1 V \nbs_2 \Sigma\bigr)
\label{eq:Zn:nb1:Sigma:identity}
\\
\nbs^5 \nbs_2 \Zbn
&= - \nbs^5\nbs_2 \Abt 
- \tfrac{1}{2\alpha \eps} \Sigma^{-1} \nbs^5 \nbs_\s  \bigl( g^{{\frac{1}{2}} } (\Wbt-\Zbt) \bigr)
- \tfrac{1}{2\alpha}   V\Sigma^{-1} \nbs^5 \nbs_2 \bigl(g^{{\frac{1}{2}} } (\Wbt-\Zbt) \bigr)
\notag\\
&\qquad 
- \tfrac{1}{2 \alpha \eps} \jump{\nbs^5,\Sigma^{-1}} \nbs_\s \bigl( g^{{\frac{1}{2}} } (\Wbt-\Zbt) \bigr)
- \tfrac{1}{\alpha}  \jump{\nbs^5, V\Sigma^{-1}}  \nbs_2^2 \Sigma 
\notag\\
&\qquad
- \nbs^5 \bigl( \tfrac{1}{\Sigma}  (   \Zbn +   \Abt  + \tfrac{1}{\alpha} \nbs_2 V) \nbs_2 \Sigma\bigr)
\label{eq:Zn:nb2:Sigma:identity}
\end{align}
\end{subequations}
Taking into account the bootstraps  \eqref{bootstraps},  the previously established bounds \eqref{geometry-bounds-new}, \eqref{eq:Jg:Abn:D5:improve},  \eqref{eq:Jg:Zbn:D5:improve:a}, and the Moser-type bound \eqref{eq:Lynch:1}, 
we deduce from \eqref{eq:Zn:nb1:Sigma:identity} that
 \begin{align}
\snorm{\Sigma^{-\bar \beta} \mathcal{J}^{\frac 34 - \bar a} \Jg^{\! \bar a} \nbs^5 \nbs_1 \Zbn}_{L^2_{x,\s}}
&\leq 
\tfrac{1}{2\alpha} \| \Sigma^{-1-\bar \beta} \mathcal{J}^{\frac 34 - \bar a} \Jg^{\! \bar a} \nbs^5 \nbs_\s (\Jg \Zbn) \|_{L^2_{x,\s}}
\notag\\
&\qquad
 +   (4\kappa_0^{-1})^{\bar \beta} (1+ \Cn \eps) \eps \mathsf{K} \mathsf{B_6}
  + \Cn (4\kappa_0^{-1})^{\bar \beta}  \eps^2 \mathsf{K} \brak{\mathsf{B_6}}  
  + \Cn  (4\kappa_0^{-1})^{\bar \beta} \eps \brak{\mathsf{B_6}}  
 \notag\\
&\qquad 
+ \Cn   (4\kappa_0^{-1})^{\bar \beta} \eps  \| \mathcal{J}^{\frac 34 - \bar a} \Jg^{\! \bar a} \nbs^5 (\Q\p_\s + V\p_2)( \Jg \Wbn) \|_{L^2_{x,\s}}
\notag\\
&\qquad 
+ \Cn  (4\kappa_0^{-1})^{\bar \beta} \eps^2 \brak{\mathsf{B}_6} \|(\Q\p_\s + V\p_2) ( \Jg \Wbn)\|_{L^\infty_{x,\s}}
\notag\\
&\qquad 
+ \Cn  (4\kappa_0^{-1})^{\bar \beta} \eps \|  \nbs^4(\Q\p_\s + V\p_2) (\Jg \Wbn)\|_{L^2_{x,\s}}
\,.
\label{eq:Zn:nb1:Sigma:bnd:1}
\end{align}
In the above estimate we have chosen not to leave the terms involving $\tfrac{1}{\eps} \nbs_\s (\Jg \Wbn)$ as is, since they would give sub-optimal bounds, and instead to write them in terms of $(\Q \p_\s + V\p_2) (\Jg \Wbn)$, which satisfies better bounds. Indeed, from \eqref{eq:Jg:Wb:nn} written in $(x,\s)$-variables, combined with \eqref{bootstraps}, \eqref{geometry-bounds-new}, \eqref{eq:Jg:Abn:D5:improve},  \eqref{eq:Jg:Zbn:D5:improve:a}, and \eqref{eq:Lynch:1}, 
we obtain that 
\begin{subequations}
\label{eq:Jg:Wbn:improve:material}
\begin{align}
 \|(\Q\p_\s + V\p_2) ( \Jg \Wbn)\|_{L^\infty_{x,\s}} 
 &\les 1
 \label{eq:Jg:Wbn:improve:material:a}
 \\
\|  \nbs^4(\Q\p_\s + V\p_2) (\Jg \Wbn)\|_{L^2_{x,\s}}
&\les  \mathsf{K} \brak{\mathsf{B_6}}
\label{eq:Jg:Wbn:improve:material:b}
\\
\| \mathcal{J}^{\frac 14} \nbs^5 (\Q\p_\s + V\p_2)( \Jg \Wbn) \|_{L^2_{x,\s}}
&\les \mathsf{K} \brak{\mathsf{B_6}}
\label{eq:Jg:Wbn:improve:material:c}
\end{align}
\end{subequations}
holds, under the assumptions of Proposition~\ref{prop:vort:H6}. Inserting the above bounds into \eqref{eq:Zn:nb1:Sigma:bnd:1}, and using that $\eps$ is sufficiently small with respect to $\alpha,\kappa_0$ and $\Cdata$, it follows that for any $\bar a \in [0,\frac 12]$
 \begin{equation}
\snorm{\Sigma^{-\bar \beta} \mathcal{J}^{\frac 34 - \bar a} \Jg^{\! \bar a} \nbs^5 \nbs_1 \Zbn}_{L^2_{x,\s}}
\leq 
\tfrac{1}{2\alpha} \| \Sigma^{-1-\bar \beta} \mathcal{J}^{\frac 34 - \bar a} \Jg^{\! \bar a} \nbs^5 \nbs_\s (\Jg \Zbn) \|_{L^2_{x,\s}}
+ \Cn   (4\kappa_0^{-1})^{\bar \beta} \eps   \mathsf{K} \brak{\mathsf{B_6}}
\,,
\label{eq:Zn:nb1:Sigma:bnd:2}
\end{equation}
where the implicit constant is independent of $\bar \beta$ and $\bar a$.

Similarly to \eqref{eq:Zn:nb1:Sigma:bnd:1}, we may deduce from \eqref{eq:Zn:nb2:Sigma:identity} that 
\begin{align}
 \snorm{\Sigma^{-\bar \beta} \mathcal{J}^{\frac 34 - \bar a} \Jg^{\! \bar a} \nbs^5 \nbs_2  \Zbn}_{L^2_{x,\s}}
&\leq 
\Cn  \eps  (4\kappa_0^{-1})^{\bar \beta} \mathsf{K} \brak{\mathsf{B_6}}
+ \tfrac{1}{\alpha }  \|\Sigma^{-1-\bar \beta} \mathcal{J}^{\frac 34 - \bar a} \Jg^{\! \bar a} \nbs^5 (\Q \p_\s + V\p_2) \bigl( g^{{\frac{1}{2}} } (\Wbt-\Zbt) \bigr) \|_{L^2_{x,\s}}
\notag\\
&\qquad + \Cn   \eps  (4\kappa_0^{-1})^{\bar \beta} \brak{\mathsf{B_6}}  \| (\Q \p_\s + V\p_2) \bigl( g^{{\frac{1}{2}} } (\Wbt-\Zbt) \bigr) \|_{L^\infty_{x,\s}}\,.
\label{eq:Zn:nb2:Sigma:bnd:1}
\end{align}
Here we have again chosen to rewrite some of the $\frac{1}{\eps} \nbs_\s$ terms in terms of $(\Q \p_s + V\p_2)$, as otherwise we are obtaining sub-optimal bounds. By using  \eqref{g12-evo} and \eqref{eq:Wb:tt},  written in terms of $(x,\s)$ variables, and by appealing to \eqref{bootstraps}, \eqref{geometry-bounds-new}, \eqref{eq:Jg:Abn:D5:improve},  \eqref{eq:Jg:Zbn:D5:improve:a}, and \eqref{eq:Lynch:1}, we may show that 
\begin{subequations}
\begin{align}
 &\|\Sigma^{-1-\bar \beta} \mathcal{J}^{\frac 34 - \bar a} \Jg^{\! \bar a}  \nbs^5 (\Q \p_\s + V\p_2) \bigl( g^{{\frac{1}{2}} } (\Wbt-\Zbt) \bigr) \|_{L^2_{x,\s}}
 \notag\\
 &\qquad \leq \Cn \eps (4\kappa_0^{-1})^{\bar \beta}  \mathsf{K} \brak{\mathsf{B}_6} 
 + \tfrac{1}{\eps} \|\Sigma^{-1-\bar \beta} \mathcal{J}^{\frac 34 - \bar a} \Jg^{\! \bar a}   \nbs^5 \nbs_\s \Zbt \|_{L^2_{x,\s}}
\end{align}
and
\begin{equation}
 \| (\Q \p_\s + V\p_2) \bigl( g^{{\frac{1}{2}} } (\Wbt-\Zbt) \bigr) \|_{L^\infty_{x,\s}}
 \les \eps   + \| (\Q \p_\s + V\p_2)  \Zbt  \|_{L^\infty_{x,\s}} \les 1
 \,.
\end{equation} 
\end{subequations}
The bound \eqref{eq:Jg:Zbn:p2D5:improve:new} now follows by combining the above two estimates and~\eqref{eq:Zn:nb2:Sigma:bnd:1}. 

It remains to prove the bounds~\eqref{eq:Jg:Zbn:p1D6:sup:improve:new} and~\eqref{eq:Jg:Zbn:p2D6:sup:improve:new}.
Notice that here we do not aim for a sharp pre-factor in front of the leading order term, as done previously for \eqref{eq:Jg:Zbn:p1D5:improve} and \eqref{eq:Jg:Zbn:p2D5:improve:new}. 
First, we revisit \eqref{eq:Zn:nb1:Sigma:identity}. Taking into account the bootstraps  \eqref{bootstraps},  the previously established bounds \eqref{geometry-bounds-new}, \eqref{eq:Jg:Abn:D5:improve},  \eqref{eq:Jg:Zbn:D5:improve:a}, and the bounds in Lemma~\ref{lem:comm:bdd}, we deduce
\begin{align}
\snorm{\mathcal{J}^{\frac 34} \Jgh \nbs_1 \nbs^5 \Zbn }_{L^\infty_\s L^2_{x}} 
&
\les \eps^{\frac 12} \mathsf{K} \mathsf{B}_6 
+ \eps^{-\frac 12} \mathsf{B}_6 + 
\snorm{\mathcal{J}^{\frac 34} \Jgh  \nbs^6 (\Jg \nbs_2 h (\Wbt-\Zbt) )}_{L^\infty_\s L^2_{x}} 
\notag\\
&\quad 
+ \eps \snorm{\mathcal{J}^{\frac 34} \Jgh  \nbs^6 (g^{\frac 12} (\Wbt-\Zbt) )}_{L^\infty_\s L^2_{x}} 
+ \eps^{-\frac 12}  \brak{\mathsf{B}_6} 
\notag\\
&\quad 
+ \snorm{\mathcal{J}^{\frac 34} \Jgh  \jump{\nbs^5,V\Sigma^{-1}} \nbs_1\nbs_2 \Sigma}_{L^\infty_\s L^2_{x}}  
+ \eps^{\frac 32}  \brak{\mathsf{B}_6}  
+ \eps^{\frac 12}  \brak{\mathsf{B}_6} 
\,.
\label{eq:Live:Earls:Court:75a}
\end{align}
For the three terms on the right side of the \eqref{eq:Live:Earls:Court:75a} which still involve norms of products and commutators, the results in Lemma~\ref{lem:comm:bdd} do not apply directly; instead, the desired bound follows from the argument which was used to prove Lemma~\ref{lem:comm:bdd}. For example, the available bounds, Sobolev interpolation, and the fundamental theorem of calculus in time, imply 
\begin{subequations}
\label{eq:Live:Earls:Court:75}
\begin{align}
\snorm{\mathcal{J}^{\frac 34} \Jgh  \nbs^6 (g^{\frac 12} (\Wbt-\Zbt) )}_{L^\infty_\s L^2_{x}}  
&\les \mathsf{K} \eps^{\frac 12}\brak{\mathsf{B}_6} 
+ \displaystyle{\sum}_{k=2}^{4} \snorm{\nbs^{6-k} (g^{\frac 12}) \nbs^{k}(\Wbt,\Zbt) }_{L^\infty_\s L^2_{x}}
\notag\\
&\les \mathsf{K} \eps^{\frac 12}\brak{\mathsf{B}_6} +  \eps^{\frac 32} 
+ \eps^{-\frac 12} \displaystyle{\sum}_{k=2}^{4} \snorm{\nbs_\s \nbs^{6-k} (g^{\frac 12}) \nbs^{k}(\Wbt,\Zbt) }_{L^2_{x,\s}} 
\notag\\
&\qquad 
+ \eps^{-\frac 12} \displaystyle{\sum}_{k=2}^{4}
\snorm{\nbs^{6-k} (g^{\frac 12}) \nbs_\s \nbs^{k}(\Wbt,\Zbt) }_{L^2_{x,\s}} 
\notag\\
&\les \mathsf{K} \eps^{\frac 12}\brak{\mathsf{B}_6} \,.
\label{eq:Live:Earls:Court:75b}
\end{align}
Similarly to \eqref{eq:Live:Earls:Court:75b} one may show that 
\begin{align}
 \snorm{\mathcal{J}^{\frac 34} \Jgh  \nbs^6 (\Jg \nbs_2 h (\Wbt-\Zbt) )}_{L^\infty_\s L^2_{x}}
&\les \mathsf{K} \eps^{\frac 12}\brak{\mathsf{B}_6} \,,
\label{eq:Live:Earls:Court:75c}\\
\snorm{\mathcal{J}^{\frac 34} \Jgh  \jump{\nbs^5,V\Sigma^{-1}} \nbs_1\nbs_2 \Sigma}_{L^\infty_\s L^2_{x}}  
&\les \mathsf{K} \eps^{\frac 32}\brak{\mathsf{B}_6} 
\label{eq:Live:Earls:Court:75d}
\,.
\end{align}
\end{subequations}
Inserting \eqref{eq:Live:Earls:Court:75} into the bound~\eqref{eq:Live:Earls:Court:75a}
leads to
\begin{equation*}
\snorm{\mathcal{J}^{\frac 34} \Jgh \nbs_1 \nbs^5 \Zbn }_{L^\infty_\s L^2_{x}} 
\les \eps^{-\frac 12} \brak{\mathsf{B_6}}
\,, 
\end{equation*}
thereby proving~\eqref{eq:Jg:Zbn:p1D6:sup:improve:new}.
Estimate~\eqref{eq:Jg:Zbn:p2D6:sup:improve:new} is established in a nearly identical manner, by bounding the right side of \eqref{eq:Zn:nb2:Sigma:identity}, instead of \eqref{eq:Zn:nb1:Sigma:identity}. We omit these redundant details.
\end{proof}

\subsection{Improved estimates for $\Jg \Wbn$}
We next obtain a few improved estimates for  $\nbs^5(\Jg \Wbn)$.
\begin{lemma}
\label{lem:duck:a:fck}
Under the  assumptions of Proposition~\ref{prop:vort:H6}, in addition to the improved bounds for $(\Q \p_\s + V \p_2) (\Jg \Wbn)$ from \eqref{eq:Jg:Wbn:improve:material},  we have the estimates
\begin{subequations}
\label{Wbn:improved}
\begin{align}
\snorm{\nbs^5 (\Jg \Wbn)}_{L^\infty_\s L^2_x}
&\leq 2 \Cdata \eps^{-\frac 12}
\,,
\label{eq:D5:JgWbn}
\\
\snorm{{\mathcal J}^{\! \frac 14} \Jgh \nbs_\s \nbs^5 (\Jg\Wbn )}_{L^2_{x,\s}} 
&\leq \Cn \eps \mathsf{K} \brak{\mathsf{B}_6}
\,.   
\label{eq:Ds:D4:Jg:Wbn:new}
\end{align} 
\end{subequations}
\end{lemma}

\begin{proof}[Proof of Lemma~\ref{lem:duck:a:fck}]
In order to prove \eqref{eq:D5:JgWbn}, we apply $\nb^5$ to \eqref{eq:Jg:Wb:nn} and then transform to $(x,\s)$ variables, to obtain
\begin{align}
&(\Q  \p_\s + V \p_2) \nbs^5 (\Jg \Wbn)
+\alpha \Sigma g^{-\frac 12} \Jg \nbs^5  \nbs_2 \Abn 
\notag\\
&= - \jump{\nbs^5,V} \nbs_2 (\Jg \Wbn)
- \nbs^5  \bigl( \bigl(\tfrac{\alpha}{2} \Abt - \tfrac{\alpha}{2} \Sigma g^{-\frac 32} \nbs_2^2 h \bigr)  (\Jg\Wbn) \bigr) - \mathsf{Rem}
\label{eq:duck:a:fck:1}
\end{align}
where the remainder term is given by 
\begin{align}
\mathsf{Rem}
&=   
\alpha \jump{\nbs^5, \Sigma g^{-\frac 12} \Jg }   \nbs_2 \Abn 
- \tfrac{\alpha}{2}  \nbs^5 \bigl( \bigl( \Abt + \Sigma g^{-\frac 32} \nbs_2^2 h \bigr) \Jg \Zbn \bigr) 
\notag\\
&\quad 
+ \nbs^5 \bigl( \bigl(\tfrac{3+\alpha}{2} \Wbt + \tfrac{1-\alpha}{2} \Zbt\bigr) \Jg\Abn \bigr)
+ \nbs^5 \bigl(\bigl(\tfrac{1+\alpha}{2} \Wbt + \tfrac{1-\alpha}{2} \Zbt\bigr)\Jg \Wbt   \bigr)
+ \alpha \nbs^5 \bigl( \Sigma g^{-\frac 32} \nbs_2^2 h  \Jg \Abt \bigr)
\,.
\label{eq:duck:a:fck:2}
\end{align}
Thus, establishing \eqref{eq:D5:JgWbn} amounts to bounding the forcing terms in \eqref{eq:duck:a:fck:1} and the remainder from \eqref{eq:duck:a:fck:2} in $L^2_{x,\s}$. For this purpose, we note that \eqref{bootstraps}, \eqref{eq:broncos:eat:shit:2b}, \eqref{geometry-bounds-new},  \eqref{eq:Jg:Abn:D5:improve}, \eqref{eq:x1:Poincare}, \eqref{eq:Lynch:1}, and \eqref{eq:Lynch:3}, imply
\begin{equation}
\snorm{\jump{\nbs^5,V} \nbs_2 (\Jg \Wbn)}_{L^2_{x,\s}}
+ \snorm{\nbs^5  \bigl( \bigl(\tfrac{\alpha}{2} \Abt - \tfrac{\alpha}{2} \Sigma g^{-\frac 32} \nbs_2^2 h \bigr)  (\Jg\Wbn) \bigr)}_{L^2_{x,\s}}
\les \mathsf{K}  \brak{\mathsf{B}_6}
\label{eq:duck:a:fck:2a}
\end{equation}
and
\begin{equation}
\snorm{\mathsf{Rem}}_{L^2_{x,\s}}
\les \mathsf{K} \eps \brak{\mathsf{B}_6}
\label{eq:duck:a:fck:2b}
\,.
\end{equation}

Next we return to \eqref{eq:duck:a:fck:1}, which we test with $\nbs^5 (\Jg\Wbn)$. By appealing to \eqref{adjoint-2}, \eqref{eq:Q:all:bbq}, \eqref{eq:duck:a:fck:2a}, and \eqref{eq:duck:a:fck:2b}, we obtain
\begin{align}
\snorm{\Q^{\frac 12}  \nbs^5 (\Jg\Wbn)(\cdot,\s)}_{L^2_x}^2
&\leq \snorm{\Q^{\frac 12}  \nbs^5 (\Jg\Wbn)(\cdot,0)}_{L^2_x}^2 
-2 \alpha \int_0^\s \!\!\! \int \Sigma g^{-\frac 12} \Jg \nbs^5 \nbs_2 \Abn \nbs^5 (\Jg\Wbn) 
\notag\\
&\qquad 
+ \Cn \int_0^\s  \snorm{\Q^{\frac 12}  \nbs^5 (\Jg\Wbn)(\cdot,\s')}_{L^2_x}^2 {\rm d}\s'
+ \Cn  \mathsf{K}^2 \brak{\mathsf{B}_6}^2
\,.
\label{eq:duck:a:fck:3}
\end{align}
The only tricky term is the second term on the right side of \eqref{eq:duck:a:fck:3}. For this term, we recall the bootstrap~\eqref{bootstraps-Dnorm:5}, Remark~\ref{rem:B5:B6}, and the improved estimate \eqref{eq:Jg:Abn:D6:improve:d}, to bound
\begin{align}
\left| \int_0^\s \!\!\! \int \Sigma g^{-\frac 12} \Jg \nbs^5 \nbs_2 \Abn \nbs^5 (\Jg\Wbn) \right| 
&\leq \|\Sigma g^{-\frac 12}\|_{L^\infty_{x,\s}} \snorm{\Jgh  \nbs^5 (\Jg\Wbn)}_{L^\infty_\s L^2_x}
\|\mathcal{J}^{\frac 14} \Jgh \nbs^5 \nbs_2 \Abn\|_{L^2_{x,\s}} \snorm{\mathcal{J}^{-\frac 14}}_{L^2_{\s}}
\notag\\
&\leq \Cn \eps^{\frac 12} \mathsf{K} \brak{\mathsf{B}_6}^2 \left( \int_0^\eps \mathcal{J}(\s)^{-\frac 12} {\rm d}\s\right)^{\frac 12}
\notag\\
&\leq \Cn \eps  \mathsf{K} \brak{\mathsf{B}_6}^2
\,,
\label{last-pain-in-the-ass}
\end{align} 
where in the last inequality we have used the fact that $\s\leq \eps$ and that $\mathcal{J}(\s) = 1 - \frac{\s}{\eps}$ is a function independent of $x$. Using this bound and taking a supremum in time in \eqref{eq:duck:a:fck:3} for $\s\in[0,\eps]$, we deduce that 
\begin{equation*}
\sup_{\s \in [0,\eps]} \snorm{\Q^{\frac 12}  \nbs^5 (\Jg\Wbn)(\cdot,\s)}_{L^2_x}^2
\leq  \snorm{\Q^{\frac 12}  \nbs^5 (\Jg\Wbn)(\cdot,0)}_{L^2_x}^2 
+ \Cn \eps \sup_{\s \in [0,\eps]} \snorm{\Q^{\frac 12}  \nbs^5 (\Jg\Wbn)(\cdot,\s)}_{L^2_x}^2
+ \Cn  \mathsf{K}^2 \brak{\mathsf{B}_6}^2
\,.
\end{equation*}
The proof of \eqref{eq:D5:JgWbn} now follows upon absorbing the second term on the right side into the left side, by appealing to \eqref{table:derivatives}, to the upper and lower bounds for $\Q$ at time $\s=0$ which follow from  \eqref{eq:Qcal:bbq:temp:3}--\eqref{eq:Qcal:bbq:temp:1}, and by taking $\eps$ to be sufficiently small with respect to $\mathsf{K}$ and $\mathsf{B}_6$.

In order to prove \eqref{eq:Ds:D4:Jg:Wbn:new}, we appeal to \eqref{eq:Jg:Wbn:improve:material:c}, to the identity $(\Q \p_\s + V \p_2) = \frac{1}{\eps} \nbs_\s + V \nbs_2$, and to the bound $\Jg \leq \frac 32$, to conclude
\begin{equation}
\tfrac{1}{\eps}  \| \mathcal{J}^{\frac 14}\Jgh \nbs^5 \nbs_\s  ( \Jg \Wbn) \|_{L^2_{x,\s}}
\les
\eps \| \mathcal{J}^{\frac 14} \Jgh  \nbs^5 \nbs_2 ( \Jg \Wbn) \|_{L^2_{x,\s}}
+ \|  \jump{\nbs^5, V}\nbs_2 ( \Jg \Wbn) \bigr)\|_{L^2_{x,\s}}
+ \mathsf{K} \brak{\mathsf{B_6}} 
\label{eq:cying:out:at:night}
\,.
\end{equation}
Note that \eqref{bootstraps-Dnorm:6} gives $\| \mathcal{J}^{\frac 14} \Jgh  \nbs^5 \nbs_2 ( \Jg \Wbn) \|_{L^2_{x,\s}} \leq \mathsf{B}_6$. Moreover, from \eqref{eq:Lynch:3}, \eqref{bootstraps}, \eqref{eq:V:H6:new}, and \eqref{eq:D5:JgWbn} we obtain $\|   \jump{\nbs^5, V}\nbs_2 ( \Jg \Wbn) \bigr)\|_{L^2_{x,\s}} \les \eps \mathsf{K} \brak{\mathsf{B}_6}$. With these   bounds, multiplying both sides of  \eqref{eq:cying:out:at:night} by $\eps$, gives \eqref{eq:Ds:D4:Jg:Wbn:new}.
\end{proof}


\section{Closing the pointwise bootstrap inequalities}
\label{sec:pointwise:bootstraps}

The goal of this section is to close the pointwise bootstrap inequalities \eqref{bs-Jg}--\eqref{bs-D-Sigma}. We first claim that:
 
\begin{proposition}
\label{cor:Jg:positive}
Assume that the bootstraps~\eqref{bootstraps} hold, and assume that $\eps$ is sufficiently small with respect to $\alpha,\kappa_0$, and $\Cdata$. Then, for all $(x,t) \in \mathcal{P}$ we have that 
\begin{subequations}
\begin{align}
 \Jg  \Wbn  & \geq   - \tfrac{9}{10} \eps^{-1} \ \ \text{ implies that } \ \ \Jg \geq \tfrac{81}{1000}  \,,\\
 \sabs{\Jg \Wbn} &\leq (1+ \tfrac{\eps}{2}) \eps^{-1} \,,\\
 \sabs{\nb(\Jg \Wbn)}  &\leq \tfrac 52 \eps^{-1} \,,\\
 \Jg &\leq \tfrac{21}{20}\,,\\
 \sabs{\nb \Jg} &\leq 4+ \alpha \,.
\end{align}
\end{subequations}
In particular,  the bootstraps \eqref{bs-Jg}, \eqref{bs-JgnnWb}, \eqref{bs-nnWb,2}, \eqref{bs-Jg-simple}, and \eqref{bs-Jg,1} are closed.
\end{proposition}
\begin{proof}[Proof of Proposition~\ref{cor:Jg:positive}]
The bound \eqref{eq:broncos:eat:shit:20} implies that 
\begin{equation*}
|\Jg \Wbn(x,t)| \leq |(w_0),_1(x)| + \Cn \eps \leq \eps^{-1} (1 + \Cn \eps^2) \leq \eps^{-1} (1+\tfrac{\eps}{2}).
\end{equation*}
The same estimate also shows that if $\Jg \Wbn(x,t) \geq - \frac{9}{10} \eps^{-1}$, then we must have $(w_0),_1(x) \geq  - \frac{9}{10} \eps^{-1} - \Cn \eps$, and therefore \eqref{bs-Jg-0} and \eqref{time-of-existence} imply
\begin{align*}
\Jg(x,t) 
&\geq 1 + (t-\initial) \tfrac{1+\alpha}{2} \bigl(  (w_0),_1(x)  - \mathsf{C_{J_t}} \bigr)
\notag\\
&\geq 1 + (t-\initial) \tfrac{1+\alpha}{2}  \bigl( - \tfrac{9}{10 \eps}   -2 \mathsf{C_{J_t}} \bigr)
\notag\\
& \geq 1 - \tfrac{2\eps}{1+\alpha} \tfrac{51}{50} \tfrac{1+\alpha}{2\eps} \bigl( \tfrac{9}{10}  + 2 \eps \mathsf{C_{J_t}} \bigr)
\geq 1 - \tfrac{459}{500}   - 3 \eps \mathsf{C_{J_t}} = \tfrac{41}{500} - 3 \mathsf{C_{J_t}} \geq \tfrac{81}{1000} \,,
\end{align*}
assuming that $\eps$ is sufficiently small. In a similar fashion, the upper bound for $\Jg$ is also obtained from  \eqref{bs-Jg-0}, combined with the upper bound in \eqref{eq:why:the:fuck:not:0}
\begin{equation*}
\Jg(x,t) 
\leq 1 + (t-\initial)  \tfrac{1+\alpha}{2}  \bigl( (w_0),_1(x) + \mathsf{C_{J_t}}  \bigr)
\leq 1 + \tfrac{2\eps}{1+\alpha} \tfrac{51}{50}  \tfrac{1+\alpha}{2 \eps} \bigl(\tfrac{1}{10}   + \eps \mathsf{C_{J_t}} \bigr)
\leq 1 +  \tfrac{51}{500} + 2 \eps \mathsf{C_{J_t}} \leq \tfrac{21}{20}\,,
\end{equation*}
assuming that $\eps$ is sufficiently small. The bound for $\nb (\Jg \Wbn)$ follows from \eqref{eq:broncos:eat:shit:2a} and \eqref{eq:why:the:fuck:not:0}, which together give 
\begin{equation*}
|\nb(\Jg \Wbn)(x,t)| \leq \tfrac{1}{\eps} |\nb \nb_1 w_0(x)| + \Cn \eps \mathsf{K} \brak{\mathsf{B}_6}
\leq \tfrac{2}{\eps}+ \Cn \eps \mathsf{K} \brak{\mathsf{B}_6}
\leq \tfrac{5}{2\eps}
\end{equation*}
assuming that $\eps$ is sufficiently small.  It remains to estimate $\nb \Jg$. When $\nb = \nb_1$ or $\nb = \nb_2$, from \eqref{bs-Jg-1} and \eqref{eq:why:the:fuck:not:0} we deduce
\begin{equation}
|\nb_i \Jg(x,t)| 
\leq (t-\initial) \tfrac{1+\alpha}{2\eps} |\nb \nb_1 w_0(x)| + \Cn (t-\initial) \eps \mathsf{K} \brak{\mathsf{B}_6}
\leq \tfrac{2\eps}{1+\alpha} \tfrac{51}{50} \tfrac{1+\alpha}{\eps}  + \Cn \eps^2 \mathsf{K} \brak{\mathsf{B}_6}
\leq 3,
\label{eq:Jg,_1:sharp}
\end{equation}
assuming that $\eps$ is sufficiently small. When $\nb = \eps \p_t$, from \eqref{bs-Jg-1b} we deduce 
\begin{equation*}
|\eps \p_t \Jg(x,t)|  \leq \tfrac{1+\alpha}{2} |\nb_1w_0(x)| + \Cn \eps  \leq \tfrac{1+\alpha}{2} + \Cn \eps \leq 1+\alpha.
\end{equation*}
 The claimed bound for $\nb \Jg$ thus holds in both cases, concluding the proof of the Proposition.
\end{proof}

\subsection{The $\Wbt$ bootstraps}
Next we turn to the bounds for $\Wbt$. From \eqref{eq:Wb:tt} we obtain that
\begin{align}
\Wbt\circ\xi(x,t)
&= 
(w_0),_2(x)
e^{\int_{\initial}^t (\frac{\alpha}{2} \Sigma g^{-\frac 32} h,_{22} - \frac{3+2\alpha}{2} \Abt)\circ \xi(x,r) {\rm d}r} 
\notag\\
&+
\int_{\initial}^t
\Bigl( -\alpha \Sigma g^{-\frac 12} \Abt,_2 + \tfrac{\alpha}{2} \Sigma g^{-\frac 32} h,_{22} (\Zbt + 2 \Abn) - \tfrac{1-2\alpha}{2} \Zbt \Abt\Bigr) \circ \xi(x,r) 
\notag\\
&\qquad \qquad \qquad \times 
e^{\int_{r}^t (\frac{\alpha}{2} \Sigma g^{-\frac 32} h,_{22} - \frac{3+2\alpha}{2} \Abt)\circ \xi(x,r') {\rm d}r'} {\rm d}r
\end{align}
By appealing to the bootstrap inequalities~\eqref{bootstraps} we deduce from the above formula that
\begin{equation}
\sabs{\Wbt \circ \xi(x,t) - (w_0),_2(x)}
\les \eps^2 \sabs{(w_0),_2(x)} + \eps^2
\,.
\label{eq:itchabod:0}
\end{equation}
Upon composing with $\xi^{-1}(x,t)$ and using that \eqref{eq:xi:nabla:xi} implies $|(w_0),_2(\xi^{-1}(x,t)) - (w_0),_2(x)|\les \eps^2 \|(w_0),_{22}\|_{L^\infty_x} \les \eps^2$, we arrive at 
\begin{equation}
\sabs{\Wbt(x,t) - (w_0),_2(x)} \leq \Cn \eps^2
\label{eq:itchabod:1}
\end{equation}
where the implicit constant depends solely on $\alpha,\kappa_0$, and $\Cdata$. In view of assumption \eqref{item:ic:max:w0} on the initial data, which gives $\|(w_0),_2\|_{L^\infty} \leq 1$, the bound \eqref{eq:itchabod:1} closes the bootstrap \eqref{bs-ttWb}, upon taking $\eps$ to be sufficiently small.

In order to obtain a  bound for $\nb \Wbt$, we differentiate \eqref{eq:Wb:tt} to obtain
\begin{align}
(\p_t + V \p_2) \nb \Wbt
&= - \nb V \nb_2 \Wbt  + \bigl(\tfrac{\alpha}{2} \Sigma g^{-\frac 32}  h,_{22} -  \tfrac{3+2\alpha}{2}\Abt  \bigr) \nb \Wbt  
\notag\\
&\quad 
- \alpha \nb(\Sigma g^{-\frac 12}) \nb_2 \Abt - \alpha \Sigma g^{-\frac 12} \nb\nb_2 \Abt
+ \bigl(\tfrac{\alpha}{2} \nb(\Sigma g^{-\frac 32})  h,_{22} + \tfrac{\alpha}{2}  \Sigma g^{-\frac 32}  \nb h,_{22} -  \tfrac{3+2\alpha}{2} \nb \Abt  \bigr) \Wbt  
\notag\\
&\quad 
+ \tfrac{\alpha}{2} \nb \bigl(\Sigma g^{-\frac 32}   (\Zbt + 2 \Abn )\bigr) h,_{22}
+ \tfrac{\alpha}{2} \Sigma g^{-\frac 32} \nb h,_{22} (\Zbt + 2 \Abn )
-  \tfrac{1-2\alpha}{2} \nb(\Zbt \Abt)  \,.
\label{eq:itchabod:2}
\end{align}
In order to derive an estimate in the spirit of \eqref{eq:itchabod:0}, we need to estimate in $L^\infty_{x,t}$, the last two lines of the above identity. The most difficult terms are $\nb \nb_2 \Abt$ and $\nb h,_{22}$, for which we do not have pointwise bootstrap inequalities readily available. Instead, we use \eqref{eq:h_2:D2-bound} to bound $\|\nb h,_{22}\|_{L^\infty_{x,t}} \les \mathsf{K} \eps \brak{\mathsf{B}_6}$, and we appeal to the Sobolev embedding \eqref{eq:Sobolev}, which gives $\|\nb \nb_2 \Abt\|_{L^\infty_{x,t}} \les \mathsf{K} \brak{\mathsf{B}_6}$. The remaining terms on the second and third line on the right side of \eqref{eq:itchabod:1} are then bounded using bootstrap inequalities~\eqref{bootstraps}, in $L^\infty_{x,t}$, by $\Cn \eps \brak{\mathsf{B_6}}$. Similarly to \eqref{eq:itchabod:0}, we then obtain
\begin{equation}
\sabs{(\nb \Wbt) \circ \xi(x,t) - \nb (w_0),_2(x)}
\les \eps^2 \sabs{\nb (w_0),_2(x)} + \mathsf{K} \eps \brak{\mathsf{B}_6}
\end{equation}
which in turn implies upon composing with $\xi^{-1}$ that
\begin{equation}
\sabs{(\nb \Wbt) (x,t) - \nb (w_0),_2(x)} \les   \mathsf{K} \eps \brak{\mathsf{B}_6}
\,.
\label{eq:itchabod:3}
\end{equation}
By combining  \eqref{eq:itchabod:3} with \eqref{table:derivatives}, 
we thus obtain
\begin{equation}
|\nb \Wbt(x,t)| 
\leq   \|\nb (w_0),_2\|_{L^\infty_x} + \Cn \mathsf{K} \eps \brak{\mathsf{B}_6}
\leq \Cdata + \Cn \mathsf{K} \eps \brak{\mathsf{B}_6}
\end{equation}
where as usual we have $\Cn = \Cn(\alpha,\kappa_0,\Cdata)$.  Taking $\eps$ to be sufficiently small with respect to $\alpha,\kappa_0,\Cdata$, and $\mathsf{B_6}$, closes the bootstrap inequality \eqref{bs-D-ttWb}.

\subsection{The $\Sigma$ bootstraps}
We seek sharp estimates for the quantities 
\begin{equation*}
\Sigma - \sigma_0\,, 
\qquad  \mbox{and} \qquad 
\tfrac{\nb \Sigma}{\Sigma} - \tfrac{\nb \sigma_0}{\sigma_0} \,\,. 
\end{equation*}
For this purpose, we appeal to \eqref{eq:xi:nabla:xi}, \eqref{eq:xi:hessian:xi}, and the $\Sigma$ evolution in \eqref{Sigma0-ALE}. First, by using the bootstraps~\eqref{bs-nnZb} and \eqref{bs-ttAb}, we deduce that 
\begin{equation*}
\sabs{\Sigma \circ \xi (x,t) - \sigma_0(x)} 
= \sigma_0(x) \Bigl|e^{-\alpha \int_{\initial}^t \Zbn \circ \xi + \Abt \circ \xi {\rm d}r} - 1 \Bigr|
\leq \Cn \eps \sigma_0(x) e^{\Cn \eps}\,.
\end{equation*}
Since the mean value theorem gives 
$|\Sigma \circ \xi (x,t) - \Sigma(x,t)| \leq |\xi(x,t) - x_2| \| \Sigma,_2\|_{L^\infty_{x,t}}$, by \eqref{sigma0-bound}, \eqref{eq:xi:nabla:xi}, and \eqref{bs-D-Sigma}, we deduce that
\label{eq:the:sharp:Sigma}
\begin{equation}
\abs{\tfrac{\Sigma(x,t)}{\sigma_0(x)} - 1} 
\leq \Cn \eps  e^{\Cn \eps} + \tfrac{\Cn \eps^2 }{\sigma_0(x)}
\leq  \Cn \eps
\label{eq:Sigma:sharp}
\end{equation}
upon taking $\eps$ sufficiently small with respect to $\alpha,\kappa_0$, and $\Cdata$. The initial data assumption \eqref{sigma0-bound} closes the bootstrap \eqref{bs-Sigma} if we choose $\eps$ to be small enough to ensure $\Cn \eps \leq \frac{1}{12} \kappa_0$.

Next, differentiating \eqref{Sigma0-ALE} we deduce  
\begin{equation*}
(\p_t + V \p_2) \tfrac{\nb \Sigma}{\Sigma} = - \nb V \tfrac{\Sigma,_2}{\Sigma} - \alpha \nb \Zbn - \alpha \nb \Abt 
\,,
\end{equation*}
and thus by appealing to \eqref{bs-nnZb}, \eqref{bs-nnAb}, and \eqref{bs-V}, we obtain that 
\begin{align*}
\sabs{\tfrac{\nb \Sigma}{\Sigma} \circ \xi (x,t) - \tfrac{\nb \sigma_0}{\sigma_0}(x)} 
&\leq \tfrac{|\nb \sigma_0(x)|}{\sigma_0(x)} \Bigl|e^{ \int_{\initial}^t |\nb V| \circ \xi {\rm d}r} - 1 \Bigr| + \alpha e^{ \int_{\initial}^t |\nb V| \circ \xi {\rm d}r} \int_{\initial}^t |\nb \Zbn\circ \xi| + |\nb \Abt \circ \xi| {\rm d}r
\notag\\
&\leq \Cn \eps^2   \tfrac{|\nb \sigma_0(x)|}{\sigma_0(x)} + \Cn \eps  
\,.
\end{align*}
Using the bound 
$\| \nb_2 \tfrac{\nb \Sigma}{\Sigma} \|_{L^\infty_{x,t}} \les 1+ \|\nb^2 \Sigma\|_{L^\infty_{x,t}} 
\les 1 + \eps^{-1} \|\nb^4 \nb_1 \Sigma\|_{L^2_{x,t}} \les \brak{\mathsf{B}_6}$, and the estimate \eqref{eq:xi:nabla:xi}, we deduce from the above bound that 
\begin{equation*}
 \sabs{\tfrac{\nb \Sigma}{\Sigma}  (x,t) - \tfrac{\nb \sigma_0}{\sigma_0}(x)} 
 \leq \Cn \eps^2     \tfrac{|\nb \sigma_0(x)|}{\sigma_0(x)} +  \Cn \eps   + \Cn \eps^2   \brak{\mathsf{B}_6}
 \leq \Cn \eps^2     \tfrac{|\nb \sigma_0(x)|}{\sigma_0(x)} +  \Cn \eps .
\end{equation*}
Note that from \eqref{item:ic:infinity}, \eqref{item:ic:max:w0},   \eqref{eq:why:the:fuck:not:0}, and \eqref{table:derivatives}, upon taking $\eps$ to be sufficiently small with respect to $\kappa_0$ and $\Cdata$, we have that 
\begin{equation}
\sabs{\tfrac{\nb_1 \sigma_0}{\sigma_0}(x)}
\leq \tfrac{1}{\kappa_0}
\leq 1 \,,
\qquad\mbox{and}\qquad 
\sabs{\tfrac{\nb_2 \sigma_0}{\sigma_0}(x)}
\leq \tfrac{1}{\kappa_0}
\leq 1\,,
\label{Dsigma0-bound} 
\end{equation}
where we have used that $\kappa_0 \geq 1$. 
By also appealing to \eqref{eq:Sigma:sharp}, and by taking $\eps$ to be sufficiently small with respect to $\alpha,\kappa_0,\Cdata$, and $\brak{\mathsf{B_6}}$, we deduce   
\begin{equation}
\sabs{\nb \Sigma(x,t)}
\leq |\nb \sigma_0(x)| (1 + \eps \Cn) + \eps \Cn   
\leq \tfrac 23 + \Cn \eps
\leq 1 \leq \kappa_0 \,.
\end{equation}
This closes the bootstrap \eqref{bs-D-Sigma}.

\subsection{The $h$ bootstraps}
We first note that the bootstrap \eqref{bs-h} implies
\begin{equation*}
\sabs{g(x,t)-1} \leq (1+\alpha)^2\kappa_0^2 \eps^2
\,.
\end{equation*}
In order to close the bootstrap \eqref{bs-h} we first recall from \eqref{h-evo} that $\nb_t h = \eps \p_t h = \eps g^{-\frac 12} (\frac{1+\alpha}{2}W+\frac{1-\alpha}{2}Z)$.
Therefore, with \eqref{eq:WZA:bounds:sup} we arrive at 
\begin{equation*}
|\nb_t h| \leq \eps (1- \Cn \eps^2)^{-\frac 12} \bigl(\tfrac{3(1+\alpha)}{4} \kappa_0 + \Cn \eps \bigr)
< \tfrac{5(1+\alpha)}{6} \eps \kappa_0 
\,,
\end{equation*}
which closes the $\nb_t h$ part of the bootstrap \eqref{bs-h}. Moreover, since $h,_{1} = g^{\frac 12} \Jg$, from \eqref{bs-Jg-simple} we deduce 
\begin{equation*}
|\nb_1 h| \leq \tfrac 32 \eps (1 + \Cn \eps^2)^{\frac 12} 
\,,
\end{equation*}
which closes the $\nb_1 h$ part of the bootstrap \eqref{bs-h} upon taking $\eps$ to be sufficiently small.
On the other hand, from \eqref{p2-h}, \eqref{bs-ttWb}, and \eqref{bs-ttZb}, we obtain that 
\begin{equation*}
\sabs{h,_{2} \circ \xi (x,t)}
\leq \tfrac{4\eps}{1+\alpha} (1 + \Cn \eps^2)^{\frac 12}  \Bigl(\tfrac{1+\alpha}{2} \cdot (1+\eps) + \tfrac{1-\alpha}{2} \cdot  \mathsf{C}_{\Zbt} \eps  \Bigr)
\leq 2  \eps   + \Cn \eps^2
\,.
\end{equation*}
which improves the $\nbs_2 h$ part of \eqref{bs-h}.

Next, we estimate $\nb h,_2$. Upon applying $\nb$ to \eqref{p2-h}, we derive 
\begin{equation*}
(\p_t + V \p_2) \nb h,_2 + \nb V h,_{22} = \tfrac{1+\alpha}{2} g \nb \Wbt  + \tfrac{1-\alpha}{2}  g \nb \Zbt
+ 
2 (\tfrac{1+\alpha}{2}  \Wbt  + \tfrac{1-\alpha}{2}   \Zbt) h,_2 \nb h,_2
\,.
\end{equation*}
By appealing to the bootstraps \eqref{bs-ttWb}, \eqref{bs-D-ttWb}, \eqref{bs-nnZb}, \eqref{bs-V}, 
we deduce from the above identity that
\begin{equation*}
\sabs{(\nb h,_2) \circ \xi(x,t)}
\leq \tfrac{4\eps}{1+\alpha} (1+ \Cn \eps^2) \bigl( \tfrac{1+\alpha}{2} \cdot 2 \Cdata + \Cn \eps\bigr)  e^{\Cn \eps}
\leq 4 \Cdata \eps + \Cn \eps^2 \,.
\end{equation*}
Taking $\eps$ to be sufficiently small, and composing with $\xi^{-1}$ then  improves the $\nb \n_2 h$ part of the bootstrap \eqref{bs-h,22}. In order to improve the $\nb \nb_1 h$ part of \eqref{bs-h,22}, we just note that 
\begin{equation*}
|\nb \nb_1 h|
\leq \eps g^{\frac 12} |\nb \Jg| + \eps \Jg  g^{-\frac 12} |h,_2 \nb h,_2|
\leq \eps (1 + \Cn \eps^2) 4(1+\alpha) + \Cn \eps^3 
\leq 4(1+\alpha) \eps + \Cn \eps^3\,.
\end{equation*}

\subsection{The $V$  bootstraps}
From \eqref{transport-ale},  \eqref{bs-h}, and \eqref{eq:WZA:bounds:sup}  we obtain that 
\begin{align}
\sabs{V} 
&\leq (1 + \Cn \eps^2)^{-\frac 12} 
\Bigl( \kappa_0 \eps  \bigl(1 + 6 \kappa_0 + \tfrac{4\alpha}{1+\alpha}\mathsf{C}_{\Abn} \bigr)  + 3  \eps \bigl(\tfrac{1+\alpha}{2} \tfrac{3}{2} \kappa_0 + \Cn \eps \bigr)\Bigr)
\notag\\
&\leq  \kappa_0 \eps \bigl(1 + 10 (1+\alpha) \kappa_0 + \tfrac{4\alpha}{1+\alpha}\mathsf{C}_{\Abn} \bigr)  \,,
\label{eq:seus:1}
\end{align}
upon taking $\eps$ to be sufficiently small.
Next, the bootstrap inequalities \eqref{bootstraps} and the identities \eqref{eq:nabla:V:1}--\eqref{eq:nabla:V:2}, yield
\begin{subequations}
\label{eq:seus:2}
\begin{align}
|\nb_1 V|
&\leq  5 \Cdata \alpha \eps \kappa_0 (1 - \Cn \eps^2)^{-\frac32}
+ \eps (1 - \Cn \eps^2)^{-\frac 12} \bigl(\tfrac 32 \mathsf{C}_{\Abn} + 5 \kappa_0 \eps \tfrac{1+\alpha}{\eps} \bigr)
+ \Cn \eps^2    \notag\\
&\leq \eps \kappa_0 \bigl(6 \Cdata  + 5(1+\alpha )  + \tfrac 32 \mathsf{C}_{\Abn} \bigr)
\,, 
\label{eq:seus:2a}
\\
|\nb_2 V|
&\leq 5 \Cdata \alpha \eps \kappa_0 (1 - \Cn \eps^2)^{-\frac32}
+ \mathsf{C}_{\Abt} \eps + 10 \kappa_0^2 \eps
\notag\\
&\leq \eps \kappa_0 \bigl( 6 \Cdata \alpha   + 10 \kappa_0  + \mathsf{C}_{\Abt} \bigr)
\,,
\label{eq:seus:2b}
\end{align}
where we have used that $\kappa_0\geq 1$.
It remains to estimate $\nb_t V = \eps \p_t V$. For this purpose, we appeal to \eqref{V-evo}, which combined with the bootstrap inequalities \eqref{bootstraps} yields (for a constanat $C_1 = C_1(\alpha,\kappa_0,\Cdata)>0$)
\begin{align}
\sabs{\nb_t V}
&\leq \Cn \eps^3 + \alpha \kappa_0 \eps (1 - \Cn \eps^2)^{-\frac 12}
\Bigl((2+\alpha)\kappa_0 +\Cn \eps + \mathsf{C}_{\Abn}  \Bigr)
\notag\\
&\leq   \alpha  \kappa_0  \eps   \bigl(2(1+\alpha)\kappa_0 + \mathsf{C}_{\Abn} \bigr)
\label{eq:seus:2c}
\end{align}
\end{subequations}
by choosing $\eps$ to be sufficiently small, in terms of $\alpha,\kappa_0$, and $\Cdata$. 
By combining \eqref{eq:seus:1} and \eqref{eq:seus:2}, and choosing the constant $\mathsf{C_V}$ to satisfy
\begin{equation}
\mathsf{C_V} \geq 2 \Bigl(
1 + 27 (1+\alpha) \kappa_0  
+ 6 (1+\alpha) \Cdata   + \tfrac{13}{2} \mathsf{C}_{\Abn}  + \mathsf{C}_{\Abt}
\Bigr)
\,,
\label{eq:C_V:choice}
\end{equation}
we have thus improved the bootstrap \eqref{bs-V}. We note that since $\mathsf{C}_{\Abn}$ and 
$\mathsf{C}_{\Abt}$ only depend on $\alpha,\kappa_0,$ and $\Cdata$, and therefore so does $\mathsf{C_V}$.

\subsection{The $\Ab$ bootstraps}
The estimate for $\Abn$ is rather direct to obtain, since we've already estimated the vorticity $\Omega = \Abn - \frac 12 (\Wbt + \Zbt)$ in \eqref{eq:vorticity:pointwise}.
Indeed, by using \eqref{eq:vorticity:pointwise}, \eqref{bs-ttWb},  and \eqref{bs-ttZb},   we obtain that
\begin{equation}
\sabs{\Abn} 
\leq C_{\eqref{eq:vorticity:pointwise:a}} \Cdata + \kappa_0 + \Cn \eps
\leq \bigl(C_{\eqref{eq:vorticity:pointwise:a}} + 1 \bigr) \Cdata + \Cn \eps
\end{equation}
where we have appealed to the bound $\Cdata \geq \kappa_0$, and recall that $C_{\eqref{eq:vorticity:pointwise:a}}  = C_{\eqref{eq:vorticity:pointwise:a}} (\alpha) = 2^{3+\frac{2}{\alpha}} e^{18}$. As such, if $\eps$ is taken to be sufficiently small, and if 
\begin{equation}
 \mathsf{C_{\Abn}} \geq 4 \bigl(C_{\eqref{eq:vorticity:pointwise:a}} + 1 \bigr) \Cdata + 1
 \label{eq:CZ:cond:1}
\end{equation}
then 
\begin{equation}
\sabs{\Abn} 
\leq \tfrac 14 \mathsf{C_{\Abn}}     
\,.
\label{eq:wooket:Abn}
\end{equation}
Similarly, by using \eqref{eq:vorticity:pointwise},   \eqref{bs-D-ttWb}, and \eqref{bs-ttZb},   we obtain that 
\begin{equation}
\sabs{\nbs \Abn} 
\leq C_{\eqref{eq:vorticity:pointwise:b}} \Cdata + \Cdata + \Cn \eps
\leq \bigl(C_{\eqref{eq:vorticity:pointwise:b}} + 1  \bigr) \Cdata +\Cn \eps  
\,,
\end{equation}
where $C_{\eqref{eq:vorticity:pointwise:b}}  = C_{\eqref{eq:vorticity:pointwise:b}} (\alpha) = 2 (4 e^{20})^{\frac{20 \cdot 23 (1+\alpha)}{\alpha}}$. 
Thus, if we also ensure that 
\begin{equation}
 \mathsf{C_{\Abn}} \geq 4 (C_{\eqref{eq:vorticity:pointwise:b}} + 1 \bigr) \Cdata + 1
 \label{eq:CZ:cond:2}
\end{equation}
then 
\begin{equation}
\sabs{\nbs \Abn} 
\leq \tfrac 14 \mathsf{C_{\Abn}}   \,.
\label{eq:wooket:D:Abn}
\end{equation}
These estimates thus close the bootstrap assumption~\eqref{bs-nnAb}. Note that both $C_{\eqref{eq:vorticity:pointwise:a}}$ and $C_{\eqref{eq:vorticity:pointwise:b}}$ are explicit functions of the parameter $\alpha$ alone, and thus the smallest number $\mathsf{C}_{\Abn}$ which satisfies  both \eqref{eq:CZ:cond:1} and \eqref{eq:CZ:cond:2}, depends on $\alpha$ and $\Cdata$ alone. This parameter is henceforth fixed.

The estimates for $\Abt$ and its derivative does not follow from the vorticity estimate, but it is established in a similar manner, by appealing to Proposition~\ref{prop:transverse:bounds} and  Corollary~\ref{cor:transverse:bounds}. We start from \eqref{eq:Ab:tt:alt}, which we multiply by $\Jg$:
\begin{equation}
(\Jg \p_t + \Jg (V + \alpha \Sigma g^{-\frac 12} h,_2) \p_2 - \alpha \Sigma \p_1) \Abt
= q_{\Abt} 
\label{eq:pucket:Abt}
\end{equation}
where
\begin{align}
q_{\Abt} 
&:= 
- \tfrac \alpha 2 \Sigma g^{- \frac{1}{2} }  \Jg (\Wbt - \Zbt),_2  
- \tfrac \alpha 2   \Sigma g^{-\frac 12} (\Wbt + \Zbt + 2 \Abn )  \Jg,_2 
+ \tfrac{\alpha}{2} \Sigma g^{- \frac{3}{2} } (\Jg \Wbn- \Jg \Zbn - h,_2 \Jg \Abn) h,_{22}
\notag \\
&\qquad 
 + \tfrac{\alpha}{2} \Jg \Abn ( \Wbt - \Zbt ) 
 - \Jg (\Abt)^2 
- \tfrac{\alpha }{4} \Jg (\Wbt -\Zbt)^2 
- \tfrac{1}{2}  \Jg ( \Wbt+\Zbt) ( \tfrac{1+ \alpha }{2} \Wbt + \tfrac{1- \alpha }{2} \Zbt)
\,.
\end{align}
Using the bootstrap inequalities~\eqref{bootstraps}, we have that 
\begin{align}
\snorm{q_{\Abt}}_{L^\infty_{x,t}}  
&\leq 
\tfrac{6\alpha \kappa_0  }{5} \Cdata 
+ 2\alpha (1+\alpha) \kappa_0 (1+ 2 \mathsf{C}_{\Abn})
+ \tfrac{5 \alpha \kappa_0}{2} \Cdata 
+ \tfrac{3\alpha}{5} \mathsf{C}_{\Abn}
+ \tfrac{3(1+2\alpha)}{10} 
+ \Cn \eps
\notag\\
&\leq  4\alpha \kappa_0    \Cdata 
+ 2\alpha (1+\alpha) \kappa_0 (1+ 2 \mathsf{C}_{\Abn})
+ \tfrac{3\alpha}{5} \mathsf{C}_{\Abn}
+ \tfrac{2 (1+2\alpha)}{5} 
=:C_{\eqref{eq:constant:Abt}},
\label{eq:constant:Abt}
\end{align}
upon taking $\eps$ to be sufficiently small. Here we recall that $\mathsf{C}_{\Abn}$ was already chosen (see~ \eqref{eq:CZ:cond:1} and \eqref{eq:CZ:cond:2}) to depend solely on $\alpha$ and $\Cdata$, and so $C_{\eqref{eq:constant:Abt}}$ in fact just dependens on $\alpha$, $\kappa_0$, and $\Cdata$. Then, since \eqref{eq:pucket:Abt} takes precisely the form of \eqref{eq:abstract:0}, with $f = \Abt$, we deduce from the above estimate, from \eqref{table:derivatives}, and appealing to \eqref{eq:abstract:transport}, that 
\begin{equation}
\snorm{\Abt}_{L^\infty_{x,\s}} \leq 4 e^{18} \snorm{\Abt(\cdot,0)}_{L^\infty_{x}}+ \tfrac{20 e^{18}}{\alpha}  C_{\eqref{eq:constant:Abt}} \eps 
\leq \eps \bigl(4 e^{18} \Cdata  + \tfrac{20 e^{18}}{\alpha}  C_{\eqref{eq:constant:Abt}}  \bigr) =:  \eps C_{\eqref{eq:wooket:Abt}} \,,
\label{eq:wooket:Abt}
\end{equation}
for some computable constant $C_{\eqref{eq:wooket:Abt}}$ which only depends on $\alpha,\kappa_0$, and $\Cdata$. 
Thus, if we ensure that 
\begin{equation}
 \mathsf{C_{\Abt}} \geq 4  C_{\eqref{eq:wooket:Abt}}   \,,
 \label{eq:CZ:cond:3}
\end{equation}
where we note that $\mathsf{C}_{\Abn}$ has already been chosen at this stage by \eqref{eq:CZ:cond:1} and \eqref{eq:CZ:cond:2}, in terms of $\alpha, \kappa_0,\Cdata$, 
then 
\begin{equation}
\sabs{\Abt} 
\leq \tfrac 14 \eps \mathsf{C_{\Abt}}   \,.
\label{eq:wooket:Abt:2}
\end{equation}

Lastly, we turn to the $\nb \Abt$ bound. We differentiate \eqref{eq:pucket:Abt} to obtain
\begin{align}
&(\Jg \p_t + \Jg (V + \alpha \Sigma g^{-\frac 12} h,_2) \p_2 - \alpha \Sigma \p_1) \nb \Abt
\notag\\
&= - \underbrace{ \nb \Jg  \p_t  \Abt  - \alpha \nb\Sigma \p_1  \Abt}_{=: m_{\nb \Abt} \nb \Abt}
- \underbrace{ \nb \Jg   (V + \alpha \Sigma g^{-\frac 12} h,_2) \p_2 \Abt 
- \Jg \bigl(\nb V +\alpha \nb (\Sigma g^{-\frac 12} h,_2) \bigr) \p_2 \Abt 
+ \nb q_{\Abt}}_{=:q_{\nb \Abt}}
\,
\label{eq:pucket:D:Abt}
\end{align}
which thus takes the form of equation \eqref{eq:abstract:2} for $f= \nb \Abt$, with suitably defined functions $m = m_{\nb \Abt}$ and $q= q_{\nb \Abt}$.
First, we note that 
\begin{equation}
\snorm{m_{\nb \Abt}}_{L^\infty_{x,t}} \leq \tfrac{1}{\eps} \snorm{ \nb \Jg }_{L^\infty_{x,t}} + \tfrac{\alpha}{\eps}  \snorm{\nb \Sigma}_{L^\infty_{x,t}}
\leq \tfrac{4(1+\alpha)+2\alpha \kappa_0}{\eps} \,.
\end{equation}
As such, the parameter $\beta$ in Corollary~\ref{cor:transverse:bounds} may be taken to equal 
$\beta = \frac{20}{\alpha} (4(1+\alpha)+2\alpha \kappa_0)$, so that $\beta = \beta(\alpha,\kappa_0)$.
In order to apply \eqref{eq:abstract:transport:2}, we then use the bootstraps~\eqref{bootstraps}, the Sobolev embedding~\eqref{eq:Sobolev}, the geometry bounds~\eqref{geometry-bounds-new}, and the initial data bounds \eqref{table:derivatives}, to estimate
\begin{align}
\snorm{q_{\nb \Abt}}_{L^\infty_{x,t}}
&\leq  \Cn  \bigl(1 + \mathsf{C}_{\Abn}  +  \eps \mathsf{C_{\Abt}} + \eps \mathsf{C}_{\Zbt} + \eps^2 \mathsf{C}_{\Abt}^2 \bigr)
\notag\\
&\qquad + \Cn \bigl( \snorm{\nb^2 (\Jg \Wbt,\Jg \Zbt)}_{L^\infty_{x,t}} + (1 + \mathsf{C}_{\Abn} ) \snorm{\nb^2  \Jg}_{L^\infty_{x,t}} \bigr) + \tfrac{C_1}{\eps} \snorm{\nb h,_{22}}_{L^\infty_{x,t}}
\notag\\
&\leq C_{\eqref{eq:constant:D:Abt}} \brak{\mathsf{C}_{\Abn}}  \brak{\mathsf{B_6}}
\,,
\label{eq:constant:D:Abt}
\end{align}
where $C_{\eqref{eq:constant:D:Abt}}$ is a computable constant that depends only on $\alpha,\kappa_0$, and $\Cdata$.
Using the above two estimates, we apply Corollary~\ref{cor:transverse:bounds} to the evolution \eqref{eq:pucket:D:Abt}, estimate the initial data via~\eqref{table:derivatives} and deduce that 
\begin{equation}
\snorm{\nb \Abt}_{L^\infty_{x,\s}}
\leq  (4 e^{18})^{\frac{20}{\alpha}(4(1+\alpha)+2\alpha\kappa_0)}  \Bigl( \eps \Cdata 
+  \eps \tfrac{  C_{\eqref{eq:constant:D:Abt}} }{ 4(1+\alpha)+2\alpha\kappa_0 } \brak{\mathsf{C}_{\Abn}}  \brak{\mathsf{B_6}} \Bigr)
\,.
\end{equation}
Taking $\eps$ to be sufficiently small, in terms of $\alpha,\kappa_0$ and $\Cdata$, we deduce that there exists a computable constant $C_{\eqref{eq:wooket:D:Abt}}> 0$, with the same dependences, such that 
\begin{equation}
\snorm{\nb \Abt}_{L^\infty_{x,\s}} \leq \eps C_{\eqref{eq:wooket:D:Abt}}  \brak{\mathsf{C}_{\Abn}} \brak{\mathsf{B_6}}\,.
\label{eq:wooket:D:Abt} 
\end{equation}
Thus, if we ensure that 
\begin{equation}
 \mathsf{C_{\Abt}} \geq 4 C_{\eqref{eq:wooket:D:Abt}}  \brak{\mathsf{C}_{\Abn}} \brak{\mathsf{B_6}} \,,
 \label{eq:CZ:cond:4}
\end{equation}
then 
\begin{equation}
\sabs{\nb \Abt} 
\leq \tfrac 14 \eps \mathsf{C_{\Abt}}   \,.
\label{eq:wooket:D:Abt:2}
\end{equation}
We note that \eqref{eq:wooket:Abt:2} and \eqref{eq:wooket:D:Abt:2} close the bootstrap \eqref{bs-ttAb}. Moreover, we emphasize that the conditions \eqref{eq:CZ:cond:3} and \eqref{eq:CZ:cond:4} on $\mathsf{C}_{\Abt}$ mandate that this constant is chosen to be sufficiently large with respect to $\alpha,\kappa_0,\Cdata$, but also with respect to $\mathsf{C}_{\Abn}$ (which has been chosen already in terms of $\alpha ,\Cdata$; cf.~\eqref{eq:CZ:cond:1} and \eqref{eq:CZ:cond:2}) and $\mathsf{B_6}$ (which will be chosen to be dependent only on $\alpha,\kappa_0,\Cdata$; cf.~\eqref{eq:B6:choice:1}, and~\eqref{eq:B6:choice:2}). 

\subsection{The $\Zb$ bootstraps}
From \eqref{eq:Zb:nn:alt} we deduce that $\Zbn$ solves the forced transport equation
\begin{equation}
\bigl(\Jg \p_t + \Jg (V +2\alpha \Sigma g^{-\frac 12} h,_{2}   )  \p_2 -2\alpha \Sigma \p_1 \bigr)  \Zbn   =  m_{\Zbn} \Zbn + q_{\Zbn}\,,
\label{eq:Zb:nn:alt:alt}
\end{equation}
where we have denoted
\begin{equation*}
m_{\Zbn} := - \bubu{\tfrac{1-\alpha}{2}} \Jg \Wbn - \bubu{\tfrac{1+\alpha}{2}} \Jg \Zbn
- \tfrac{\alpha}{2} \Abt - \tfrac{\alpha}{2} \Sigma g^{-\frac 32} h,_{22}
\end{equation*}
and 
\begin{align}
q_{\Zbn}
&:=
 2\alpha \Sigma \bubu{g^{-\frac 12}} ( \Zbt+ \Abn)  {\Jg,_2}     
+\alpha \Sigma g^{-\frac 12}\Jg   \Abn,_2
+\tfrac{\alpha}{2} \bigl( \Abt - \Sigma g^{-\frac 32}   h,_{22} \bigr) \Jg \Wbn     
\\
&
-  \bigl(\tfrac{1+\alpha}{2} \Wbt + \tfrac{3-\alpha}{2} \Zbt \bigr) \Jg \Abn
- \bigl(\tfrac{1+\alpha}{2} \Wbt + \tfrac{1-\alpha}{2} \Zbt \bigr) \Jg \Zbt
+   \alpha  \Sigma g^{-\frac 32}   h,_{22}  \Jg \Abt 
.
\label{eq:umbrella:1}
\end{align}
Since \eqref{eq:Zb:nn:alt:alt} takes the form of \eqref{eq:abstract:2}, with $f= \Zbn$ and $\alpha$ replaced by $2\alpha$, the desired upper bound for $\Zbn$ will follow from \eqref{eq:abstract:transport:2}. We note that since \eqref{bootstraps} implies $\| m_{\Zbn}\|_{L^\infty_{x,t}} \leq (1+\alpha) \eps^{-1}$ (upon letting $\eps$ be sufficiently small), we may choose the exponent $\beta$ appearing in \eqref{eq:abstract:transport:2} to equal $\beta = \beta(\alpha) = \frac{20(1+\alpha)}{\alpha}$. Moreover, \eqref{table:derivatives} implies that $\|\Zbn(\cdot,\initial)\|_{L^\infty_x} \leq \Cdata$. Thus, it remains to estimate the space-time $L^\infty$ norm of $q_{\Zbn}$. From the bootstrap inequalities~\eqref{bootstraps} we deduce   
\begin{align}
\|q_{\Zbn}\|_{L^\infty_{x,t}}
&\leq 8\alpha(1+\alpha) \kappa_0 \mathsf{C}_{\Abn}
+ \alpha \kappa_0 \tfrac 65 \mathsf{C}_{\Abn}
+ \tfrac{\alpha}{2} (\mathsf{C}_{\Abt} + 5  \kappa_0 \Cdata)
+ \tfrac{1+\alpha}{2}   \tfrac 65 \mathsf{C}_{\Abn}
+ \Cn \eps 
\notag\\
&\leq (1+\alpha)(1+9\alpha) \kappa_0 \mathsf{C}_{\Abn}
+ \tfrac{\alpha}{2} (\mathsf{C}_{\Abt} + 5  \kappa_0 \Cdata)
=: C_{\eqref{eq:constant:Zbn}}
\,,
\label{eq:constant:Zbn}
\end{align}
where $ C_{\eqref{eq:constant:Zbn}}>0$ is explicitly computable solely in terms of $\alpha,\kappa_0$, and $\Cdata$.  With the above estimate, we deduce from \eqref{eq:abstract:transport:2} the pointwise estimate 
\begin{equation}
\sabs{\Zbn}
\leq 
(4 e^{18})^{\frac{20(1+\alpha)}{\alpha}} 
\Bigl( \snorm{\Zbn(\cdot,\initial)}_{L^\infty_x}
+ 
 \tfrac{\eps}{ 1+\alpha }  C_{\eqref{eq:constant:Zbn}}
\Bigr)
\leq 2 \Cdata (4 e^{18})^{\frac{20(1+\alpha)}{\alpha}} 
\label{eq:umbrella:3}
\end{equation}
upon taking $\eps$ to be sufficiently small.
This justifies the first condition we need to impose on $\mathsf{C}_{\Zbn}$, namely 
\begin{equation}
 \mathsf{C}_{\Zbn} \geq 8 \Cdata (4 e^{18})^{\frac{20(1+\alpha)}{\alpha}} 
 \,,
 \label{eq:CZ:cond:5}
\end{equation}
which in turn implies the pointwise estimate
\begin{equation}
\sabs{\Zbn} \leq  \tfrac 14 \mathsf{C}_{\Zbn}
\,.
\label{eq:wooket:Zbn}
\end{equation}

The estimate for $\nb \Zbn$ is obtained by differentiating \eqref{eq:Zb:nn:alt:alt}, which leads to 
\begin{align}
&\bigl(\Jg \p_t + \Jg (V +2\alpha \Sigma g^{-\frac 12} h,_{2}   )  \p_2 -2\alpha \Sigma \p_1 \bigr)  \nb \Zbn   
\notag\\
&\qquad = \underbrace{m_{\Zbn} \nb \Zbn -  \nb \Jg \p_t  \Zbn -\nb \bigl(  \Jg (V +2\alpha \Sigma g^{-\frac 12} h,_{2}   )\bigr)  \p_2  \Zbn + 2\alpha \nb \Sigma \p_1 \Zbn }_{:= m_{\nb \Zbn} \nb \Zbn} + q_{\nb \Zbn}  \,,
\label{eq:D:Zb:nn:alt:alt}
\end{align}
where 
\begin{equation}
q_{\nb \Zbn} =  \nb m_{\Zbn} \Zbn + \nb q_{\Zbn}
\,.
\label{eq:umbrella:2}
\end{equation}
As before, it follows that $\nb \Zbn$ solves a forced transport equation of the type \eqref{eq:abstract:2}. We start by using \eqref{bootstraps} to estimate the stretching factor
\begin{equation}
\|m_{\nb \Zbn}\|_{L^\infty_{x,t}}
\leq \tfrac{1+\alpha}{\eps} + \tfrac{5(1+\alpha)}{\eps} + \Cn \eps + \tfrac{6 \alpha \kappa_0}{\eps} 
\leq 
 \tfrac{6(1+\alpha)(1+\kappa_0)}{\eps} 
\end{equation}
by taking $\eps$ to be sufficiently small.
As such, we may take the constant $\beta$ appearing in \eqref{eq:abstract:transport:2} to equal $\beta = \tfrac{120}{\alpha} (1+\alpha)(1+\kappa_0)$. In order to estimate $q_{\nb \Zbn}$, as defined by \eqref{eq:umbrella:1} and \eqref{eq:umbrella:2},  we use the bootstrap inequalities~\eqref{bootstraps}, the previously derived bound \eqref{eq:umbrella:1}, the Sobolev estimate~\eqref{eq:Sobolev}, the initial data bounds~\eqref{table:derivatives},   the bounds for the geometry~\eqref{geometry-bounds-new}, and the improved estimate~\eqref{eq:Jg:Abn:D5:improve}, to deduce 
\begin{equation}
\|q_{\nb \Zbn}\|_{L^\infty_{x,t}}
\leq 2 \Cdata (4 e^{18})^{\frac{20(1+\alpha)}{\alpha}} \bigl(\tfrac{1+\alpha}{2}(1+\alpha) \Cdata \eps^{-1} + \Cn \bigr)  + \Cn \brak{\mathsf{B_6}} .
\end{equation}
By taking $\eps$ to be sufficiently small with respect to $\alpha,\kappa_0,\Cdata$ and $\mathsf{B_6}$, we deduce from the above two estimates and 
\eqref{eq:abstract:transport:2} that 
\begin{equation}
\sabs{\nb \Zbn}
\leq 
(4 e^{18})^{\beta} \Cdata
+ 
 \tfrac{(1+\alpha)}{6(1+\kappa_0)}  
    (4 e^{18})^{\frac{20(1+\alpha)}{\alpha} + \beta}   \Cdatatwo 
\end{equation}
where we recall that $ \beta = \beta(\alpha,\kappa_0) =  \frac{120}{\alpha}  (1+\alpha)(1+\kappa_0)$.
Thus, by further imposing that
\begin{equation}
 \mathsf{C}_{\Zbn} \geq 4 \Bigl(
 (4 e^{18})^{\beta} \Cdata
+ 
 \tfrac{1(1+\alpha)}{20(1+\kappa_0)}  
    (4 e^{18})^{\frac{20(1+\alpha)}{\alpha} + \beta}   \Cdatatwo  \Bigr)
 \,,
 \label{eq:CZ:cond:6}
\end{equation}
we obtain
\begin{equation}
\sabs{\nb \Zbn} \leq  \tfrac 14 \mathsf{C}_{\Zbn}
\end{equation}
which together with \eqref{eq:wooket:Zbn} closes the bootstrap inequality \eqref{bs-nnZb}. 
We note that the constant $\mathsf{C}_{\Zbn}$ depends only on $\alpha,\kappa_0$, and $\Cdata$, through \eqref{eq:CZ:cond:5} and \eqref{eq:CZ:cond:6}.

It thus remains to estimate $\Zbt$ and $\nb \Zbt$. For this purpose, \eqref{eq:Zb:tt:alt} yields 
\begin{equation}
\bigl(\Jg \p_t + \Jg (V +2\alpha \Sigma g^{-\frac 12} h,_2) \p_2 - 2 \alpha \Sigma  \p_1 \bigr)  \Zbt 
= - \bigl(\underbrace{\bubu{\alpha \Jg \Zbn  +}  \tfrac{\alpha}{2} \Sigma g^{-\frac 32} h,_{22} \Jg + \tfrac 32 \Jg \Abt}_{=: m_{\Zbt}} \bigr) \Zbt + q_{\Zbt}
\,,
\label{eq:pucket:Zbt}
\end{equation}
where we have denoted 
\begin{equation}
\bubu{q_{\Zbt}} =
2\alpha \Sigma   {\Jg,_2} \bubu{g^{-\frac 12}} \bigl( \Abt  -  \Zbn  \bigr)
+ \alpha \Sigma g^{-\frac 12} \Jg  \Abt,_2
- \tfrac{\alpha}{2} \Sigma g^{-\frac 32} h,_{22} \Jg (\Wbt  + 2 \Abn)
\bubu{+ \alpha \Jg \Zbn \Wbt }
- \tfrac 12  \Jg \Abt    \Wbt
\,.
\end{equation}
Thus, $\Zbt$ solves a forced transport equation of the type \eqref{eq:abstract:2}, with $\alpha$ replaced by $2\alpha$.
From \eqref{bootstraps} we note that 
\begin{equation}
\| m_{\Zbt} \|_{L^\infty_{x,t}} \leq 
\bubu{\tfrac{6\alpha}{5} \mathsf{C}_{\Zbn} } + \Cn \eps
\leq \bubu{2\alpha \mathsf{C}_{\Zbn}},
\end{equation}
and thus by taking $\eps$ to be sufficiently small we may ensure that $\beta =1$ in \eqref{eq:abstract:transport:2}. Next, we estimate
\begin{equation}
\|  q_{\Zbt}\|_{L^\infty_{x,t}} \leq 
\bubu{12} \alpha (1+\alpha) \kappa_0 (\mathsf{C}_{\Zbn} + \Cn \eps) + \Cn \eps
\leq 
\bubu{12} (1+\alpha)^2 \kappa_0  \mathsf{C}_{\Zbn} 
\end{equation}
upon taking $\eps$ to be sufficiently small in terms of  $\alpha,\kappa_0,\Cdata$.
From the above bound, \eqref{table:derivatives}, and \eqref{eq:abstract:transport:2}, we deduce
\begin{equation}
\sabs{\Zbt} 
\leq 4 e^{18} \eps \Cdata + \eps \tfrac{20}{\alpha} 4 e^{18} \cdot \bubu{12} (1+\alpha)^2 \kappa_0  \mathsf{C}_{\Zbn} 
\,.
\end{equation}
As such, if we ensure   
\begin{equation}
\mathsf{C}_{\Zbt} 
\geq 16 e^{18} \bigl(  \Cdata +   \tfrac{\bubu{240}  (1+\alpha)^2 \kappa_0}{\alpha}     \mathsf{C}_{\Zbn} \bigr)
\,,
\label{eq:CZ:cond:7}
\end{equation}
then we have the uniform estimate
\begin{equation}
 \sabs{\Zbt} \leq  \tfrac 14 \eps \mathsf{C}_{\Zbt} 
 \,.
 \label{eq:wooket:Zbt}
\end{equation}
The estimate for $\nb \Zbt$ is obtained by differentiating \eqref{eq:pucket:Zbt}, which yields
\begin{align}
& \bigl(\Jg \p_t + \Jg (V +2\alpha \Sigma g^{-\frac 12} h,_2) \p_2 - 2 \alpha \Sigma  \p_1 \bigr)  \nb \Zbt 
\notag\\
&=\underbrace{ - m_{\Zbt} \nb \Zbt - \nb \Jg \p_t \Zbt - \nb \bigl(\Jg(V + 2\alpha \Sigma g^{-\frac 12} h,_2)\bigr) \nb_2 \Zbt + 2\alpha \nb \Sigma \bubu{\p_1} \Zbt}_{=: m_{\nb \Zbt} \nb \Zbt} + q_{\nb \Zbt}
\end{align}
where we have denoted 
\begin{equation}
q_{\nb \Zbt} = \Zbt \nb m_{\Zbt} + \nb q_{\Zbt}  
\,.
\end{equation}
Using \eqref{bootstraps} we first bound 
\begin{equation}
\| m_{\nb \Zbt}\|_{L^\infty_{x,t}} 
\leq  \tfrac{5(1+\alpha)}{\eps} + \tfrac{6\alpha\kappa_0}{\eps}
+ \Cn    
\leq \tfrac{6 (1+\alpha + \alpha \kappa_0)}{\eps}
\end{equation}
so that the constant $\beta$ in \eqref{eq:abstract:transport:2} may be taken to equal 
$\beta = \frac{120}{\alpha} (1+\alpha + \alpha \kappa_0)$. 
Next, by using  the bootstraps~\eqref{bootstraps},  the Sobolev estimate~\eqref{eq:Sobolev}, the initial data bounds~\eqref{table:derivatives},   the bounds for the geometry~\eqref{geometry-bounds-new},  we obtain
\begin{equation}
\| q_{\nb \Zbt}\|_{L^\infty_{x,t}}
\leq  \Cn \brak{\mathsf{B}_6}   \eps + \Cn \mathsf{C}_{\Zbn} 
\leq C_{\eqref{eq:wooket:final}} \brak{\mathsf{C}_{\Zbn} }
\label{eq:wooket:final}
\end{equation}
for some computable constant $C_{\eqref{eq:wooket:final}}$, which depends only on $\alpha,\kappa_0$, and $\Cdata$.

From the above bound, \eqref{table:derivatives}, and \eqref{eq:abstract:transport:2}, we deduce the pointwise estimate
\begin{equation}
\sabs{\nb \Zbt} 
\leq (4 e^{18} \eps)^{\beta} \eps \Cdata + \eps \tfrac{1}{6 (1+\alpha+\alpha \beta)} (4 e^{18})^\beta C_{\eqref{eq:wooket:final}} \brak{\mathsf{C}_{\Zbn} }
\,,
\end{equation}
where we recall that $\beta = \frac{120}{\alpha} (1+\alpha + \alpha \kappa_0)$.
Thus, choosing $\mathsf{C}_{\Zbt}$ large enough to ensure
\begin{equation}
\mathsf{C}_{\Zbt} 
\geq 4 (4 e^{18} \eps)^{\beta} \bigl(\Cdata + C_{\eqref{eq:wooket:final}} \brak{\mathsf{C}_{\Zbn} }\bigr)
\,,
\label{eq:CZ:cond:8}
\end{equation}
we obtain the uniform bound
\begin{equation}
 \sabs{\nb \Zbt} \leq  \tfrac 14 \eps \mathsf{C}_{\Zbt} 
\,,
\end{equation}
which together with \eqref{eq:wooket:Zbt} closes the bootstrap \eqref{bs-ttAb}. 
We note that since $C_{\eqref{eq:wooket:final}}$ and $ \brak{\mathsf{C}_{\Zbn} }$ only depend on $\alpha,\kappa_0$, and $\Cdata$, so does $\mathsf{C}_{\Zbt}$, via \eqref{eq:CZ:cond:7} and \eqref{eq:CZ:cond:8}.


\section{The sixth order energy estimates for the tangential components}
\label{sec:sixth:order:energy-tangential}

From \eqref{euler-Wt}, \eqref{euler-Zt}, and \eqref{euler-At}, we obtain the sixth-order differentiated $(\Wbt,\Zbt, \Abt)$ equations 
\begin{subequations} 
\label{energy-WZA-tan-s}
\begin{align} 
&
\tfrac{1}{\Sigma}(\Q\p_\s +V\p_2) \nbs^6\Wbt  
+ \alpha    g^{- {\frac{1}{2}} }\nbs_2\nbs^6 \Abt
- \alpha    g^{- {\frac{1}{2}} } \nbs_2\nbs^6 \tt\cdot\nn (\Omega+\Wbt+\Zbt)
 = \nbs^6\Fwt  + \mathcal{R}_\Wb^\tt + \mathcal{C}_\Wb^\tt \,,   \label{energy-Wt-s}
 \\
 &
\tfrac{\Jg}{\Sigma} (\Q\p_\s +V\p_2) \nbs^6\Zbt
-  \alpha  \Jg  g^{- {\frac{1}{2}} } \nbs_2\nbs^6 \Abt
+ \alpha   \Jg  g^{- {\frac{1}{2}} } \nbs_2\nbs^6 \tt\cdot\nn (\Omega+\Wbt+\Zbt)
\notag \\
& \qquad 
-2 \alpha \nbs^6 \Zbt,_1 
- 2 \alpha (\Abt -\Zbn) \nbs^6 \tt,_1\cdot\nn 
+ 2 \alpha  \Jg g^{- {\frac{1}{2}} } \nbs_2 h \, \nbs_2\nbs^6\Zbt 
+ 2 \alpha  \Jg g^{- {\frac{1}{2}} } \nbs_2 h\,  (\Abt-\Zbn) \nbs_2\nbs^6\tt\cdot \nn 
   \notag \\
& \qquad\qquad\qquad\qquad
 = \nbs^6 \Fzt   + \mathcal{R}_\Zb^\tt + \mathcal{C}_\Zb^\tt\ \,,     \label{energy-Zt-s} \\
&
\tfrac{\Jg}{\Sigma}  (\Q\p_\s +V\p_2) \nbs^6\Abt 
+ \alpha \Jg   g^{- {\frac{1}{2}} } \nbs_2 \nbs^6  \Sbt
-  \alpha    g^{- {\frac{1}{2}} } (\Jg\Sbn) \nbs_2\nbs^6\tt\cdot \nn
 - \alpha  \nbs^6 \Abt,_1 
 + \alpha  (\Omega +\Wbt+\Zbt) \nbs^6\tt,_1\cdot\nn 
\notag \\
& \qquad 
+ \alpha\Jg g^{- {\frac{1}{2}} } \nbs_2h\,  \nbs_2\nbs^6\Abt
-  \alpha   \Jg g^{- {\frac{1}{2}} } \nbs_2 h\, (\Omega+\Wbt+\Zbt)  \nbs_2\nbs^6\tt\cdot \nn
 = \nbs^6 \Fat  + \mathcal{R}_\Ab^\tt + \mathcal{C}_\Ab^\tt  \,,   \label{energy-At-s}
 \end{align} 
\end{subequations} 
where
\begin{subequations} 
\label{Cw-Rw-comm-tan}
\begin{align} 
\mathcal{R} _\Wb^\tt & 
= \nbs^6 \Sigma ^{-1} (\Q\p_\s+V\p_2) \Wbt 
+ \Sigma ^{-1} \nbs^6V  \nbs_2\Wbt  
+ \alpha \nbs^6 g^{- {\frac{1}{2}} }  \nbs_2\Abt
- \alpha \nbs_2\tt^k \nbs^6 \big(\nn_k g^{- {\frac{1}{2}} } (\Omega +\Wbt+\Zbt)\big)
\,, \\
\mathcal{C}_\Wb^\tt & 
= \Sigma^{-1} \doublecom{\nbs^6, V, \nbs_2\Wbt}
+ \doublecom{\nbs^6, \Sigma^{-1}  , (\Q\p_s+V\p_2)  \Wbt}  
+ \alpha \doublecom{\nbs^6, g^{- {\frac{1}{2}} } , \nbs_2 \Abt}
\notag \\
& \qquad 
- \alpha \doublecom{ \nbs^6, \nbs_2\tt^k,  \nn_k g^{- {\frac{1}{2}} } (\Omega +\Wbt+\Zbt)}
 \,, 
\end{align} 
\end{subequations} 
and
\begin{subequations} 
\label{Cz-Rz-comm-tan}
\begin{align} 
\mathcal{R} _\Zb^\tt & 
= \nbs^6(\tfrac{\Jg}{\Sigma}) (\Q\p_\s+V\p_2) \Zbt 
+ \tfrac{1}{\Sigma} \Jg\nbs^6 V  \nbs_2\Zbt 
- \alpha \nbs^6 (g^{- {\frac{1}{2}} } \Jg)  \nbs_2\Abt
+ \alpha  \nbs_2\tt^k \nbs^6 \big(\nn_k \Jg g^{- {\frac{1}{2}} } (\Omega +\Wbt+\Zbt)\big) 
\notag \\
& \qquad
-2 \alpha \nbs^6\big( (\Abt-\Zbn) \nn_k\big)\tt^k,_1
+ 2 \alpha \nbs^6  \Jg g^{- {\frac{1}{2}} } \nbs_2 h  \nbs_2\Zbt 
+ 2 \alpha \nbs^6(  \Jg g^{- {\frac{1}{2}} } \nbs_2 h\,  (\Abt-\Zbn) \nn_k) \nbs_2\tt^k
\,, \\
 \mathcal{C}_\Zb^\tt & 
= \doublecom{\nbs^6, \tfrac{\Jg}{\Sigma} , (\Q\p_\s+V\p_2)  \Zbt}
+ \tfrac{\Jg}{\Sigma} \doublecom{\nbs^6, V, \nbs_2\Zbt}
- \alpha \doublecom{\nbs^6, g^{- {\frac{1}{2}} } \Jg, \nbs_2 \Abt}
\notag \\
&\qquad
+\alpha \doublecom{ \nbs_5, \nbs_2\tt^k, \nn_k \Jg g^{- {\frac{1}{2}} } (\Omega +\Wbt+\Zbt)}
-2 \alpha \doublecom{ \nbs^6 ,  (\Abt-\Zbn) \nn_k , \tt^k,_1}
  \notag \\
&\qquad
+ 2 \alpha \doublecom{ \nbs^6,   \Jg g^{- {\frac{1}{2}} } \nbs_2 h,  \nbs_2\Zbt }
+ 2 \alpha \doublecom{ \nbs^6, \Jg g^{- {\frac{1}{2}} } \nbs_2 h\,  (\Abt-\Zbn) \nn_k,  \nbs_2\tt^k}
 \,, 
\end{align} 
\end{subequations} 
and
\begin{subequations} 
\label{Ca-Ra-comm-tan}
\begin{align} 
\mathcal{R} _\Ab^\tt & 
=\tfrac{1}{\Sigma} \Jg\nbs^6 V  \nbs_2\Abt
+ \nbs^6(\tfrac{\Jg}{\Sigma}) (\Q\p_\s+V\p_2) \Abt
+ \alpha \nbs^6(\Jg g^{- {\frac{1}{2}} } ) \nbs_2\Sbt
- \alpha \nbs_2\tt^k \nbs^6(\nn_k  g^{- {\frac{1}{2}} } \Jg\Sbn) 
\notag \\
&\qquad
 + \alpha   \tt^k,_1 \nbs^6\big( \nn_k(\Omega +\Wbt+\Zbt)\big)
\notag \\
& \qquad
+ \alpha  \Abt \nbs^6(\Jg g^{- {\frac{1}{2}} } \nbs_2h\,  \nbs_2)
-  \alpha  \nbs_2\tt^k \nbs^6\big( \nn_k  \Jg g^{- {\frac{1}{2}} } \nbs_2 h\, (\Omega+\Wbt+\Zbt)  \big)
\,, \\
 \mathcal{C}_\Ab^\tt & 
=  \tfrac{\Jg}{\Sigma} \tt^k \doublecom{\nbs^6, V, \nbs_2 \Ak}
+ \doublecom{\nbs^6, \tfrac{\Jg}{\Sigma} \tt^k , (\Q\p_\s+V\p_2) \Ak)}  
- \alpha \doublecom{\nbs^6,\Jg g^{- {\frac{1}{2}} } , \nbs_2 \Sbt } 
\notag \\
&\qquad
- \alpha \doublecom{ \nbs^6,   g^{- {\frac{1}{2}} } \Jg\Sbn, \nbs_2\tt\cdot\nn}
+ \alpha \doublecom{\nbs^6 , \Jg g^{- {\frac{1}{2}} } \nbs_2h\,  \nbs_2, \Abt}
\notag\\
&\qquad 
-  \alpha \doublecom{ \nbs^6, \nn_k  \Jg g^{- {\frac{1}{2}} } \nbs_2 h\, (\Omega+\Wbt+\Zbt)  , \nbs_2\tt^k }
 \,.
\end{align} 
\end{subequations}

\subsection{The $\nbs^6$ tangential energy identity}
In order to obtain the fundamental $\nbs^6$ energy identities, we shall use $(x,\s)$ coordinates where $x\in \mathbb{R} ^2$ and $\s \in [0,\eps]$. We
compute the  following spacetime $L^2$ inner-product:
\begin{equation} 
\tint \jb \mathcal{J}^{\frac 32} \Big(  \underbrace{\eqref{energy-Wt-s} \ \Jg  \nbs^6\Wbt } _{ I^{\WW_\tau}} 
+  \underbrace{\eqref{energy-Zt-s} \ \nbs^6\Zbt} _{ I^{\ZZ_\tau}} 
+  \underbrace{\eqref{energy-At-s} \ 2\nbs^6\Abt} _{ I^{\AA_\tau}} \Big) {\rm d} x {\rm d}\s' =0 \,, 
\label{D6-L2-tan}
\end{equation} 
where
\begin{equation*}
\jb = \Sigma^{-2\beta+1}  
\end{equation*} 
and $\beta>0$ is a constant which will be chosen to be sufficiently large, only with respect to $\alpha$, in \eqref{eq:tangential:bounds:beta} below. Throughout this analysis, for notational convenience we will mostly omit the spacetime Lebesgue measure ${\rm d} x {\rm d}\s'$ from these integrals. 

The goal of this section is to show that identity \eqref{D6-L2-tan}, together with a good choice of $\beta= \beta(\alpha)$, and a choice of $\eps$ sufficiently small with respect to $\alpha,\kappa_0$ and $\Cdata$, implies a differential inequality of the type
\begin{align}
&\snorm{\tfrac{\mathcal{J}^{\frac 34} (\Jg \Q)^{\frac 12}}{\Sigma^{\beta}} \nbs^6(\Wbt,\Zbt,\Abt)(\cdot,\s)}_{L^2_x}^2
+ \tfrac{1}{\eps}   \int_0^{\s}  
\snorm{\tfrac{\mathcal{J}^{\frac 14} \Jg^{\frac 12}}{\Sigma^{\beta}} \nbs^6 (\Wbt,\Zbt,\Abt)(\cdot,\s')}_{L^2_x}^2 
 {\rm d} \s'
\notag\\
&\qquad \leq C(\alpha) \eps (\tfrac{4}{\kappa_0})^{2\beta}\mathsf{B}_6^2  
+ \tfrac{C(\alpha)}{\eps} \int_0^{\s}
\snorm{\tfrac{\mathcal{J}^{\frac 34}(\Jg \Q)^{\frac 12}}{\Sigma^{\beta}} \nbs^6 (\Wbt,\Zbt,\Abt)(\cdot,\s')}_{L^2_x}^2
 {\rm d} \s'
 \,,
\label{D6-L2-tan:wish}
\end{align}
where $C(\alpha) > 0$ is a constant that only depends on $\alpha$. The true inequality we establish, see~\eqref{eq:hate:11} below, is a bit more complicated because it turns out we need to augment the tangential energy estimate with the energy term $\tfrac{1}{\eps^2} \snorm{ \tfrac{\Q \mathcal{J}^{\frac 14}}{\Sigma^{\beta_\alpha}}\nn \cdot  \nbs^6 \tt (\cdot,\s)}_{L^2_{x}}^2$. Leaving this complication aside, if we were to establish \eqref{D6-L2-tan:wish}, then a standard Gr\"onwall inequality in $\s \in [0,\eps]$ shows that 
\begin{equation*}
\sup_{\s\in[0,\eps]} \snorm{\tfrac{\mathcal{J}^{\frac 34} (\Jg \Q)^{\frac 12}}{\Sigma^{\beta}} \nbs^6(\Wbt,\Zbt,\Abt)(\cdot,\s)}_{L^2_x}^2
+ \tfrac{1}{\eps}   \int_0^{\s}  
\snorm{\tfrac{\mathcal{J}^{\frac 14} \Jg^{\frac 12}}{\Sigma^{\beta}} \nbs^6 (\Wbt,\Zbt,\Abt)(\cdot,\s')}_{L^2_x}^2 {\rm d}\s'
\leq C'(\alpha) \eps (\tfrac{4}{\kappa_0})^{2\beta}\mathsf{B}_6^2 \,,
\end{equation*}
for another constant $C'(\alpha)>0$ which only depends on $\alpha$.
Upon multiplying the above estimate by $\kappa^{2\beta}$, noting that cf.~\eqref{bs-Sigma} we have $1 \leq \kappa_0^{2\beta} \Sigma^{-2\beta}$, and recalling that $\beta = \beta(\alpha)$ we deduce that the tangential part of the bootstrap \eqref{bootstraps-Dnorm:6} may be closed as soon as $K$ is taken to be sufficiently large with respect to $\alpha$. This argument is made precise in~\eqref{sec:D6:tau:final} below.

\subsection{The integral  $I^{\WW_\tau}$} 
 We additively decompose the integral  $I^{\WW_\tau}$ as 
 \begin{subequations} 
 \label{Integral-Wbt}
\begin{align}
I^{\WW_\tau}&= I^{\WW_\tau}_1+I^{\WW_\tau}_2+I^{\WW_\tau}_3 
+I^{\WW_\tau}_4   
\notag \,, \\
 I^{\WW_n}_1 &=
\tint \tfrac{1}{\Sigma^{2\beta}} \mathcal{J}^{\frac 32} \Jg (\Q\p_\s +V\p_2)\nbs^6\Wbt \ \nbs^6\Wbt
\,, \label{I1-Wbt} \\
 I^{\WW_\tau}_2 &=
\alpha \tint  \jb \mathcal{J}^{\frac 32} g^{- {\frac{1}{2}} } \Jg \nbs_2\nbs^6 \Abt \ 
  \nbs^6\Wbt
  \,,   \label{I2-Wbt} \\
I^{\WW_\tau}_3 &=
-\alpha \tint \jb \mathcal{J}^{\frac 32}  g^{-\frac 12} \Jg (\Omega+\Wbt+\Zbt) \nbs_2\nbs^6\tt \cdot \nn 
 \ \nbs^6\Wbt
  \,,   \label{I3-Wbt} \\
I^{\WW_\tau}_4 &=
-\tint  \jb   \mathcal{J}^{\frac 32} \Jg \big(\nbs^6\Fwt + \mathcal{R}_\Wb^\tau + \mathcal{C}_\Wb^\tau  \big)\  \nbs^6\Wbt
 \,, \label{I4-Wbt} 
 \end{align} 
\end{subequations}

\subsection{The integral  $I^{\ZZ_\tau}$} 
 We additively decompose the integral  $I^{\ZZ_\tau}$ as 
 \begin{subequations} 
 \label{Integral-Zbt}
 \begin{align}
  I^{\ZZ_\tau} & = I^{\ZZ_\tau}_1 + I^{\ZZ_\tau}_2 + I^{\ZZ_\tau}_3+ I^{\ZZ_\tau}_4+ I^{\ZZ_\tau}_5+ I^{\ZZ_\tau}_6+ I^{\ZZ_\tau}_7+ I^{\ZZ_\tau}_8
  \,,  \notag \\
 I^{\ZZ_\tau}_1& =
 \tint \tfrac{1}{\Sigma^{2\beta}}  \mathcal{J}^{\frac 32} \Jg(\Q\p_\s +V\p_2)\nbs^6\Zbt \   \nbs^6\Zbt
  \,, \label{I1-Zbt}\\
 I^{\ZZ_\tau}_2 &=
-\alpha \tint \jb  \mathcal{J}^{\frac 32}   \Jg g^{- {\frac{1}{2}} } \nbs_2\nbs^6 \Abt \ 
  \nbs^6\Zbt
  \,,   \label{I2-Zbt} \\
I^{\ZZ_\tau}_3 &=
\alpha \tint \jb \mathcal{J}^{\frac 32} \Jg  g^{- {\frac{1}{2}} } (\Omega+\Wbt+\Zbt) \nbs_2\nbs^6\tt \cdot \nn 
 \ \nbs^6\Zbt
  \,,   \label{I3-Zbt} \\
 I^{\ZZ_\tau}_4& =
- 2 \alpha  \tint \jb \mathcal{J}^{\frac 32} \nbs^6 \Zbt,_1 \   \nbs^6\Zbt
\,,  \label{I4-Zbt}\\
I^{\ZZ_\tau}_5& =
- 2 \alpha  \tint \jb \mathcal{J}^{\frac 32}(\Abt-\Zbn) \nbs^6\tt,_1\cdot\nn   \   \nbs^6\Zbt
\,,  \label{I5-Zbt}\\
 I^{\ZZ_\tau}_6& =
 2 \alpha \tint \jb \mathcal{J}^{\frac 32} \Jg g^{- {\frac{1}{2}} } \nbs_2h\,  \nbs_2\nbs^6\Zbt  \  \nbs^6\Zbt
\,,  \label{I6-Zbt} \\
I^{\ZZ_\tau}_7& =
 2 \alpha \tint \jb \mathcal{J}^{\frac 32}\Jg g^{- {\frac{1}{2}} } \nbs_2h\,  (\Abt-\Zbn) \nbs_2\nbs^6\tt\cdot \nn \  \nbs^6\Zbt
\,,  \label{I7-Zbt} \\
 I^{\ZZ_\tau}_8 & =
 - \tint \jb \mathcal{J}^{\frac 32} \big(  \nbs^6 \Fzt + \mathcal{R}^\tau_{\Zb} + \mathcal{C}^\tau_{\Zb}\big)
\ \nbs^6\Zbt
\,.  \label{I8-Zbt}
\end{align} 
\end{subequations} 
 
 \subsection{The integral  $I^{\AA_\tau}$}  
We additively decompose the integral  $I^{\AA_\tau}$ as
\begin{subequations} 
\label{Integral-Abt}
\begin{align}
 I^{\AA_\tau} & = I^{\AA_\tau}_1 + I^{\AA_\tau}_2 + I^{\AA_\tau}_3+ I^{\AA_\tau}_4+ I^{\AA_\tau}_5+ I^{\AA_\tau}_6+ I^{\AA_\tau}_7+ I^{\AA_\tau}_8
  \,, \notag \\
 I^{\AA_\tau}_1& =
 \tint \tfrac{2}{\Sigma^{2\beta}} \mathcal{J}^{\frac 32}  \Jg(\Q\p_\s +V\p_2)\nbs^6\Abt \   \nbs^6\Abt
  \,, \label{I1-Abt}\\
 I^{\AA_\tau}_2 &=
2\alpha \tint  \jb g^{- {\frac{1}{2}} } \mathcal{J}^{\frac 32} \Jg \nbs_2\nbs^6 \Sbt \ 
  \nbs^6\Abt
  \,,   \label{I2-Abt} \\
I^{\AA_\tau}_3 &=
-2\alpha \tint  \jb g^{- {\frac{1}{2}} }\mathcal{J}^{\frac 32} ( \Jg\Sbn) \nbs_2\nbs^6\tt \cdot \nn 
 \ \nbs^6\Abt    \,,   \label{I3-Abt} \\
 I^{\AA_\tau}_4& =
- 2 \alpha  \tint \jb \mathcal{J}^{\frac 32}   \Abt,_1 \   \nbs^6\Abt
\,,  \label{I4-Abt}\\
I^{\AA_\tau}_5& =
2 \alpha  \tint  \jb \mathcal{J}^{\frac 32} (\Omega+\Wbt+\Zbt) \nbs^6\tt,_1\cdot\nn   \   \nbs^6\Abt
\,,  \label{I5-Abt}\\
 I^{\AA_\tau}_6& =
 2 \alpha \tint  \jb \mathcal{J}^{\frac 32} \Jg g^{- {\frac{1}{2}} } \nbs_2h\,  \nbs_2 \nbs^6\Abt  \  \nbs^6\Abt
\,,  \label{I6-Abt} \\
I^{\AA_\tau}_7& =
 -2 \alpha \tint  \jb \mathcal{J}^{\frac 32} \Jg g^{- {\frac{1}{2}} } \nbs_2h\,  (\Omega+\Wbt+\Zbt) \nbs_2\nbs^6\tt\cdot \nn \  \nbs^6\Abt
\,,  \label{I7-Abt} \\
  I^{\AA_\tau}_8 & =
 - \tint  \jb \mathcal{J}^{\frac 32}\big(  \nbs^6 \Fat + \mathcal{R}^\tau_{\Ab} + \mathcal{C}^\tau_{\Ab}\big)
\ \nbs^6\Abt
\,.  \label{I8-Abt}
\end{align} 
\end{subequations} 

\subsection{The exact derivative terms}
For the terms involving a time derivative, we note that summing \eqref{I1-Wbt}, \eqref{I1-Zbt}, and \eqref{I1-Abt}, integrating by parts and appealing to \eqref{p2h-evo-s}, \eqref{adjoint-3}, \eqref{Jg-evo-s}, and \eqref{Sigma0i-ALE-s}, we obtain
\begin{align}
I^{\WW_\tau}_1 + I^{\ZZ_\tau}_1  + I^{\AA_\tau}_1 
&=\tint \tfrac{ \mathcal{J}^{\frac 32} \Jg}{2\Sigma^{2\beta}}  (\Q\p_\s +V\p_2) \Bigl( (\nbs^6\Wbt)^2 + (\nbs^6\Zbt)^2 + 2 (\nbs^6\Abt)^2 \Bigr)
\notag\\
&= \tfrac 12 \snorm{\tfrac{\mathcal{J}^{\frac 34}(\Jg \Q)^{\frac 12}}{\Sigma^\beta} \nbs^6 \Wbt(\cdot,\s)}_{L^2_x}^2
+ \tfrac 12 \snorm{\mathcal{J}^{\frac 34}\tfrac{(\Jg \Q)^{\frac 12}}{\Sigma^\beta} \nbs^6 \Zbt(\cdot,\s)}_{L^2_x}^2
+   \snorm{\tfrac{\mathcal{J}^{\frac 34}(\Jg \Q)^{\frac 12}}{\Sigma^\beta} \nbs^6 \Abt(\cdot,\s)}_{L^2_x}^2
\notag\\
&\qquad 
- \tfrac 12 \snorm{\tfrac{\mathcal{J}^{\frac 34}(\Jg \Q)^{\frac 12}}{\Sigma^\beta} \nbs^6 \Wbt(\cdot,0)}_{L^2_x}^2
- \tfrac 12 \snorm{\tfrac{\mathcal{J}^{\frac 34}(\Jg \Q)^{\frac 12}}{\Sigma^\beta} \nbs^6 \Zbt(\cdot,0)}_{L^2_x}^2
-   \snorm{\tfrac{\mathcal{J}^{\frac 34}(\Jg \Q)^{\frac 12}}{\Sigma^\beta} \nbs^6 \Abt(\cdot,0)}_{L^2_x}^2
\notag\\
&\qquad
+ \tint \tfrac{1}{\Sigma^{2\beta}} \mathsf{G}_0  \Bigl( (\nbs^6\Wbt)^2 + (\nbs^6\Zbt)^2 + 2 (\nbs^6\Abt)^2 \Bigr)
\label{eq:heavy:fuel:1}
\end{align}
where we have defined  
\begin{equation*}
 \mathsf{G}_0 
 :=  - \tfrac 12 (\Q \p_\s + V\p_2) \bigl( \mathcal{J}^{\frac 32} \Jg \bigr) 
 + \tfrac 12 \bigl(V \Qr_2 - \Qr_\s - \nbs_2 V - 2 \alpha \beta (\Zbn+\Abt) \bigr) \mathcal{J}^{\frac 32} \Jg
 \,.
\end{equation*}
At this stage, we record the pointwise bound
\begin{equation}
\mathsf{G}_0 
\geq 
\underbrace{- \tfrac 12 (\Q \p_\s + V\p_2) \bigr( \mathcal{J}^{\frac 32} \Jg \bigl)}_{=: \mathsf{G_{good}}}
- \Cn \brak{\beta}   \mathcal{J}^{\frac 12} \Jg  
- \tfrac{250^2}{\eps} \Q  \mathcal{J}^{\frac 32} \Jg  
  \label{eq:G0:tau:lower}
  \,,
\end{equation}
which follows from \eqref{eq:Jgb:less:than:1}, \eqref{bootstraps}, and \eqref{eq:Q:all:bbq}.

For the terms involving a $\p_1$ derivative of the fundamental variables, we add \eqref{I4-Zbt} and \eqref{I4-Abt}, and integrate by parts with respect to $\p_1$ (here, recall that $\p_1 \mathcal{J} = 0$) to arrive at
\begin{equation}
I^{\ZZ_\tau}_4  + I^{\AA_\tau}_4 
= - \alpha \tint \jb \mathcal{J}^{\frac 32} \p_1 \bigl( (\nbs^6 \Zbt)^2 + (\nbs^6 \Abt)^2 \bigr)
=   \tint \tfrac{1}{\Sigma^{2\beta}} \mathsf{G}_1 \bigl( (\nbs^6 \Zbt)^2 + (\nbs^6 \Abt)^2 \bigr)
\label{eq:heavy:fuel:2}
\end{equation}
where 
\begin{equation}
\mathsf{G}_1 
:=  - \alpha (2\beta-1) \mathcal{J}^{\frac 32} \Sigma,_1 
\geq
\alpha (\beta-\tfrac 12)  \bigl(\tfrac{9}{10\eps} - \tfrac{33}{(1+\alpha) \eps}  \Jg \Q -  \Cn   \bigr)  \mathcal{J}^{\frac 32}
\label{eq:heavy:fuel:G1:lower}
\end{equation}
as soon as $\beta\geq \frac 12$, 
in light of \eqref{p1-Sigma}, \eqref{bootstraps}, \eqref{eq:Q:all:bbq}, and \eqref{eq:signed:Jg}.

For the terms involving a $\p_2$ derivative of the fundamental variables, we add \eqref{I6-Zbt} and \eqref{I6-Abt}, and integrate by parts using~\eqref{adjoint-2} to arrive at
\begin{align}
I^{\ZZ_\tau}_6  + I^{\AA_\tau}_6
&=  \alpha \tint \jb \mathcal{J}^{\frac 32} \Jg g^{- {\frac{1}{2}} } \nbs_2h\,  \nbs_2 \bigl( (\nbs^6\Zbt )^2 +  (\nbs^6\Abt )^2 \bigr)
\notag\\
&= -   \alpha \int \Qb_2 \jb \mathcal{J}^{\frac 32} \Jg g^{- {\frac{1}{2}} } \nbs_2h\,   \bigl( (\nbs^6\Zbt )^2 +  (\nbs^6\Abt )^2 \bigr)\Bigr|_{\s} 
+ \tint \tfrac{1}{\Sigma^{2\beta}} \mathsf{G}_2 \bigl( (\nbs^6\Zbt )^2 +  (\nbs^6\Abt )^2 \bigr) \,  
\label{eq:heavy:fuel:3}
\end{align}
where 
\begin{equation}
\mathsf{G}_2 
=  \alpha\Sigma^{2\beta}    (\Qr_2 - \nbs_2) \bigl( \jb \mathcal{J}^{\frac 32} \Jg g^{- {\frac{1}{2}} } \nbs_2h \bigr)
\geq -  \Cn  \brak{\beta} \eps \mathcal{J}^{\frac 12} \Jg 
\label{eq:heavy:fuel:G2:lower}
\end{equation}
and we have appealed to  \eqref{bootstraps}, \eqref{eq:Q:all:bbq}, \eqref{eq:signed:Jg}, and the bound $\mathcal{J} \leq \Jg$. At this stage we also note that the bootstrap inequalities imply
\begin{align}
&-   \alpha \int \Qb_2 \jb\mathcal{J}^{\frac 32}\Jg g^{- {\frac{1}{2}} } \nbs_2h\,   \bigl( (\nbs^6\Zbt )^2 +  (\nbs^6\Abt )^2 \bigr)\Bigr|_{\s} 
\notag\\ 
&\qquad
\geq - \Cn \eps^2 \Bigl(
\snorm{\tfrac{\mathcal{J}^{\frac 34} (\Jg \Q)^{\frac 12}}{\Sigma^\beta} \nbs^6 \Zbt(\cdot,\s)}_{L^2_x}^2
+   \snorm{\tfrac{\mathcal{J}^{\frac 34}  (\Jg \Q)^{\frac 12}}{\Sigma^\beta} \nbs^6 \Abt(\cdot,\s)}_{L^2_x}^2 \Bigr)
\,.
\label{eq:heavy:fuel:G2:lower:a}
\end{align}

Lastly, we note that there are three terms with seven derivatives landing on the fundamental variables; these terms combine to yield an exact derivative, which we then integrate by parts. Adding \eqref{I2-Wbt}, \eqref{I2-Zbt}, \eqref{I2-Abt}, and recalling that $\Sbt = \tfrac 12 \Wbt - \tfrac 12 \Zbt$, using~\eqref{adjoint-2} we have that 
\begin{align}
&I^{\WW_\tau}_2 + I^{\ZZ_\tau}_2 + I^{\AA_\tau}_2
= \alpha \tint \jb \mathcal{J}^{\frac 32} g^{- {\frac{1}{2}} } \Jg \nbs_2 \Bigl( \nbs^6 \Abt \; \bigl( \nbs^6\Wbt - \nbs^6 \Zbt\bigr) \Bigr)
\notag\\
&=  \alpha \tint (\Qr_2 - \nbs_2) \bigl( \jb \mathcal{J}^{\frac 32} g^{- {\frac{1}{2}} } \Jg \bigr)  \;  \nbs^6 \Abt \; \bigl( \nbs^6\Wbt - \nbs^6 \Zbt\bigr) 
\notag\\
&\qquad  
- \alpha   \int \Qb_2 \jb g^{- {\frac{1}{2}} }  \mathcal{J}^{\frac 32}  \Jg \; \nbs^6 \Abt \; \bigl( \nbs^6\Wbt - \nbs^6 \Zbt\bigr) \Bigr|_{\s} 
\,.
\label{eq:heavy:fuel:4}
\end{align}
Note that the available pointwise bootstrap bounds imply 
\begin{equation}
\sabs{\alpha \Sigma^{2\beta} (\Qr_2 - \nbs_2)  \bigl( \jb \mathcal{J}^{\frac 32} g^{- {\frac{1}{2}} } \Jg \bigr)}
\leq  \Cn \brak{\beta} \mathcal{J}^{\frac 12} \Jg \,,
\qquad \mbox{and} \qquad 
\sabs{\alpha \Qb_2 \Sigma g^{- {\frac{1}{2}} }\Q^{-1} } 
\leq \Cn \eps
\,.
\label{eq:heavy:fuel:G3:lower}
\end{equation}

Summarizing the identities~\eqref{eq:heavy:fuel:1}, \eqref{eq:heavy:fuel:2}, \eqref{eq:heavy:fuel:3}, \eqref{eq:heavy:fuel:4} and bounds~\eqref{eq:G0:tau:lower}, \eqref{eq:heavy:fuel:G1:lower}, \eqref{eq:heavy:fuel:G2:lower}, \eqref{eq:heavy:fuel:G2:lower:a}, \eqref{eq:heavy:fuel:G3:lower} upon taking $\eps$ to be sufficiently small in terms of $\alpha,\kappa_0$ and $\Cdata$, and taking $\beta \geq 1 $, gives
\begin{align}
& I^{\WW_\tau}_1 + I^{\ZZ_\tau}_1  + I^{\AA_\tau}_1 
 + I^{\WW_\tau}_2 + I^{\ZZ_\tau}_2 + I^{\AA_\tau}_2
 + I^{\ZZ_\tau}_4  + I^{\AA_\tau}_4 
 + I^{\ZZ_\tau}_6  + I^{\AA_\tau}_6 
 \notag\\
&\quad\geq \bigl( \tfrac 12 - \Cn \eps\bigr) \Bigl( \snorm{\tfrac{\mathcal{J}^{\frac 34} (\Jg \Q)^{\frac 12}}{\Sigma^\beta} \nbs^6 \Wbt(\cdot,\s)}_{L^2_x}^2
+  \snorm{\tfrac{\mathcal{J}^{\frac 34} (\Jg \Q)^{\frac 12}}{\Sigma^\beta} \nbs^6 \Zbt(\cdot,\s)}_{L^2_x}^2
+ 2  \snorm{\tfrac{\mathcal{J}^{\frac 34} (\Jg \Q)^{\frac 12}}{\Sigma^\beta} \nbs^6 \Abt(\cdot,\s)}_{L^2_x}^2 \Bigr)
\notag\\
&\qquad 
- \tfrac 12 \Bigl( \snorm{ \tfrac{\mathcal{J}^{\frac 34}(\Jg \Q)^{\frac 12}}{\Sigma^\beta} \nbs^6 \Wbt(\cdot,0)}_{L^2_x}^2
+ \snorm{\tfrac{\mathcal{J}^{\frac 34} (\Jg \Q)^{\frac 12}}{\Sigma^\beta} \nbs^6 \Zbt(\cdot,0)}_{L^2_x}^2
+2  \snorm{\tfrac{\mathcal{J}^{\frac 34} (\Jg \Q)^{\frac 12}}{\Sigma^\beta} \nbs^6 \Abt(\cdot,0)}_{L^2_x}^2
\Bigr)
\notag\\
&\qquad
+ \tint \tfrac{1}{\Sigma^{2\beta}} \Bigl( \mathsf{G_{good}} - \Cn \beta\mathcal{J}^{\frac 12} \Jg\Bigr)
\Bigl(  |\nbs^6\Wbt|^2 
+ |\nbs^6\Zbt|^2
+ 2 |\nbs^6\Abt|^2 \Bigr)
\notag\\
&\qquad
- \tfrac{250^2}{\eps} \int_0^{\s} 
\Bigl(\snorm{\tfrac{\mathcal{J}^{\frac 34}(\Jg \Q)^{\frac 12}}{\Sigma^\beta} \nbs^6 \Wbt(\cdot,\s')}_{L^2_x}^2
+  \snorm{\tfrac{\mathcal{J}^{\frac 34}(\Jg \Q)^{\frac 12}}{\Sigma^\beta} \nbs^6 \Zbt(\cdot,\s')}_{L^2_x}^2
+    \snorm{\tfrac{\mathcal{J}^{\frac 34}(\Jg \Q)^{\frac 12}}{\Sigma^\beta} \nbs^6 \Abt(\cdot,\s')}_{L^2_x}^2 \Bigr)
{\rm d} \s'
\notag\\
&\qquad 
+   \tfrac{4\alpha( \beta - \frac 12)}{5\eps}  
\int_0^{\s} 
\Bigl( \snorm{\tfrac{ \mathcal{J}^{\frac 34} }{\Sigma^\beta} \nbs^6\Zbt (\cdot,\s')}_{L^2_x}^2
+  \snorm{\tfrac{ \mathcal{J}^{\frac 34} }{\Sigma^\beta} \nbs^6\Abt (\cdot,\s')}_{L^2_x}^2 \Bigr)
{\rm d} \s'
\notag\\
&\qquad 
- \tfrac{33 \alpha  (\beta -\frac 12) }{(1+\alpha)\eps}  \int_0^{\s} \Bigl(  \snorm{\tfrac{\mathcal{J}^{\frac 34}(\Jg \Q)^{\frac 12}}{\Sigma^\beta} \nbs^6 \Zbt(\cdot,\s')}_{L^2_x}^2
+    \snorm{\tfrac{\mathcal{J}^{\frac 34}(\Jg \Q)^{\frac 12}}{\Sigma^\beta} \nbs^6 \Abt(\cdot,\s')}_{L^2_x}^2 \Bigr)
{\rm d} \s' 
\,,
\label{eq:I:tau:1246}
\end{align}
where as usual $\Cn = \Cn (\alpha,\kappa_0,\Cdata)$ is a positive computable constant. Note   that $\Cn$ is independent of $\beta$.

\subsection{The terms with over-differentiated geometry}
Next, we consider the terms in \eqref{Integral-Wbt}--\eqref{Integral-Abt} which contain seven derivatives on the tangent vector $\tt$ (or equivalently, the normal vector $\nn$).

\subsubsection{The sum $ I^{\ZZ_\tau}_5+I^{\ZZ_\tau}_7$}
Integrating-by-parts in \eqref{I5-Zbt} and \eqref{I7-Zbt} we find that
\begin{align} 
 I^{\ZZ_\tau}_{5} +  I^{\ZZ_\tau}_{7}
 &= 2 \alpha  \tint \jb \mathcal{J}^{\frac 32} (\Abt-\Zbn) \nbs^6\tt\cdot\nn   \   \Big( \nbs^6\Zbt,_1- \Jg g^{- {\frac{1}{2}} } \nbs_2h\, \nbs_2  \nbs^6\Zbt\Big) 
 \notag \\
 &  
\qquad + 2 \alpha  \tint \Big( \p_1 \big( \jb\mathcal{J}^{\frac 32}  (\Abt-\Zbn) \nn_k\big) -\nbs_2\big( \jb \mathcal{J}^{\frac 32}  \Jg g^{- {\frac{1}{2}} } \nbs_2h\,  (\Abt-\Zbn) \nn_k\big)   \Big) \nbs^6\tt^k \   \nbs^6\Zbt 
\notag \\
&  
\qquad + 2 \alpha \tint \Qr_2 \jb \mathcal{J}^{\frac 32}  \Jg g^{- {\frac{1}{2}} } \nbs_2h\,  (\Abt-\Zbn) \nbs^6\tt\cdot \nn \  \nbs^6\Zbt
\notag\\
&
\qquad - 2\alpha  \int \Qb_2 \jb \mathcal{J}^{\frac 32} \Jg g^{- {\frac{1}{2}} } \nbs_2h\,  (\Abt-\Zbn) \nbs^6\tt\cdot \nn \  \nbs^6\Zbt\Big|_\s
\notag \\
& =:  I^{\ZZ_\tau}_{5+7,a}+I^{\ZZ_\tau}_{5+7,b}+I^{\ZZ_\tau}_{5+7,c}+I^{\ZZ_\tau}_{5+7,d} \,.
\label{eq:hate:1}
\end{align} 
At this stage we note that the first term in \eqref{eq:hate:1}, namely $I^{\ZZ_\tau}_{5+7,a}$, contains over-differentiated terms, i.e.~seven derivatives on $\Zbt$, but that the remaining three terms are under control. Indeed, from \eqref{bootstraps}, \eqref{eq:Q:all:bbq}, \eqref{D6n-bound:b:new}, and the fact that $\mathcal{J}\leq \Jg$,  we have
\begin{subequations}
\label{eq:I:Zt:5+7:bcd}
\begin{align}
\sabs{I^{\ZZ_\tau}_{5+7,b}}
&\les   \eps^{-1} \beta (\tfrac{4}{\kappa_0})^{\beta} \|\nbs^6\tt\|_{L^2_{x,\s}} 
\snorm{\tfrac{\mathcal{J}^{\frac 14} \Jg^{\frac 12}}{\Sigma^\beta} \nbs^6 \Zbt}_{L^2_{x,\s}}
\les   \eps \beta  (\tfrac{4}{\kappa_0})^{\beta} \mathsf{K} \brak{\mathsf{B}_6}
\snorm{\tfrac{\mathcal{J}^{\frac 14} \Jg^{\frac 12}}{\Sigma^\beta} \nbs^6 \Zbt}_{L^2_{x,\s}}
\\
\sabs{I^{\ZZ_\tau}_{5+7,c}}
&\les   \eps (\tfrac{4}{\kappa_0})^{\beta} \|\nbs^6\tt\|_{L^2_{x,\s}} 
\snorm{\tfrac{\mathcal{J}^{\frac 14} \Jg^{\frac 12}}{\Sigma^\beta} \nbs^6 \Zbt}_{L^2_{x,\s}}
\les \eps^3  (\tfrac{4}{\kappa_0})^{\beta} \mathsf{K} \brak{\mathsf{B}_6}
\snorm{\tfrac{\mathcal{J}^{\frac 14} \Jg^{\frac 12}}{\Sigma^\beta} \nbs^6 \Zbt}_{L^2_{x,\s}}
\\
\sabs{I^{\ZZ_\tau}_{5+7,d}}
&\les \eps^2  (\tfrac{4}{\kappa_0})^\beta \snorm{\mathcal{J}^{\frac 14} \nbs^6\tt}_{L^\infty_\s L^2_x} \snorm{\tfrac{\mathcal{J}^{\frac 34} (\Jg \Q)^{\frac 12}}{\Sigma^\beta} \nbs^6 \Zbt(\cdot,\s)}_{L^2_x}
\les \eps^{\frac 72} (\tfrac{4}{\kappa_0})^\beta \mathsf{K} \brak{\mathsf{B}_6} \snorm{\tfrac{\mathcal{J}^{\frac 34} (\Jg \Q)^{\frac 12}}{\Sigma^\beta} \nbs^6 \Zbt(\cdot,\s)}_{L^2_x}
\,.
\end{align}
\end{subequations}

In order to handle the first term on the right side of \eqref{eq:hate:1}, we use equation \eqref{energy-Zt-s} and \eqref{energy-Wt-s} to rewrite
\begin{align*} 
& 2 \alpha \big( \nbs^6 \Zbt,_1  -   \Jg g^{- {\frac{1}{2}} } \nbs_2 h \, \nbs_2\nbs^6\Zbt \big)
\notag\\
&\qquad=
-  \alpha  \Jg  g^{- {\frac{1}{2}} } \nbs_2\nbs^6 \Abt
+ \alpha   \Jg  g^{- {\frac{1}{2}} } \nbs_2\nbs^6 \tt\cdot\nn (\Omega+\Wbt+\Zbt)
+\tfrac{\Jg}{\Sigma} (\Q\p_\s +V\p_2) \nbs^6\Zbt
\notag \\
& \qquad \qquad
- 2 \alpha (\Abt -\Zbn) \nbs^6 \tt,_1\cdot\nn 
+ 2 \alpha  \Jg g^{- {\frac{1}{2}} } \nbs_2 h\,  (\Abt-\Zbn) \nbs_2\nbs^6\tt\cdot \nn 
-\big(\nbs^6 \Fzt   + \mathcal{R}_\Zb^\tt + \mathcal{C}_\Zb^\tt\big) 
\notag\\
&\qquad=
\tfrac{\Jg}{\Sigma} (\Q\p_\s +V\p_2) \nbs^6\Wbt 
+\tfrac{\Jg}{\Sigma} (\Q\p_\s +V\p_2) \nbs^6\Zbt
\notag \\
& \qquad \qquad
- 2 \alpha (\Abt -\Zbn) \nbs^6 \tt,_1\cdot\nn 
+ 2 \alpha  \Jg g^{- {\frac{1}{2}} } \nbs_2 h\,  (\Abt-\Zbn) \nbs_2\nbs^6\tt\cdot \nn 
\notag\\
&\qquad \qquad
- \big( \nbs^6\Fwt  + \mathcal{R}_\Wb^\tt + \mathcal{C}_\Wb^\tt\big)
-\big(\nbs^6 \Fzt   + \mathcal{R}_\Zb^\tt + \mathcal{C}_\Zb^\tt\big) 
\,.
\end{align*} 
Substitution of this identity into the integral $I^{\ZZ_\tau}_{5+7,a}$, gives us the further decomposition
\begin{subequations}
\label{eq:hate:2}
\begin{align} 
I^{\ZZ_\tau}_{5+7,a} & = I^{\ZZ_\tau}_{5+7,a,i}+I^{\ZZ_\tau}_{5+7,a,ii}+I^{\ZZ_\tau}_{5+7,a,iii}\,,
\notag \\
I^{\ZZ_\tau}_{5+7,a,i} 
 &=\tint \tfrac{1}{\Sigma^{2\beta}} \mathcal{J}^{\frac 32} \Jg (\Abt-\Zbn) \nbs^6\tt\cdot\nn   \  (\Q\p_\s+V\p_2) \bigl(\nbs^6\Wbt + \nbs^6 \Zbt\bigr)\,,
\label{eq:hate:2a} \\
I^{\ZZ_\tau}_{5+7,a,ii }
&= - 2\alpha  \tint \jb \mathcal{J}^{\frac 32} (\Abt-\Zbn)^2 \nbs^6\tt\cdot\nn   \  \big(\nbs^6\tt,_1\cdot\nn 
-   \Jg g^{- {\frac{1}{2}} } \nbs_2 h\,  \nbs_2\nbs^6\tt\cdot \nn \big)   \,,
\label{eq:hate:2b} \\
I^{\ZZ_\tau}_{5+7,a,iii} 
&= - \tint \jb \mathcal{J}^{\frac 32} (\Abt-\Zbn) \nbs^6\tt\cdot\nn   \ \big( \nbs^6\Fwt  + \mathcal{R}_\Wb^\tt + \mathcal{C}_\Wb^\tt + \nbs^6 \Fzt   + \mathcal{R}_\Zb^\tt + \mathcal{C}_\Zb^\tt \big) \,.
\label{eq:hate:2c}
\end{align} 
\end{subequations}
As we shall see next, the first two terms terms in \eqref{eq:hate:2} have a good structure (in spite of them having terms with seven derivatives), while the forcing terms present in \eqref{eq:hate:2c} may be estimated directly (see~Subsection~\ref{sec:D6:tau:forcing:commutator} below). The term $I^{\ZZ_\tau}_{5+7,a,ii }$ defined in \eqref{eq:hate:2b} contains perfect derivatives which may be integrated by parts:
\begin{align}
I^{\ZZ_\tau}_{5+7,a,ii }
&= -  \alpha  \tint \jb \mathcal{J}^{\frac 32} (\Abt-\Zbn)^2    \big(\p_1 (\nbs^6\tt \cdot\nn)^2 
-   \Jg g^{- {\frac{1}{2}} } \nbs_2 h\,  \nbs_2(\nbs^6\tt\cdot \nn)^2 \big) 
\notag\\
&\qquad + 2\alpha  \tint \jb \mathcal{J}^{\frac 32} (\Abt-\Zbn)^2 \nbs^6\tt\cdot\nn   \  \big(\nbs^6\tt \cdot \tt \nn,_1 \cdot\tt
-   \Jg g^{- {\frac{1}{2}} } \nbs_2 h\,  \nbs^6\tt\cdot \tt \nbs_2\nn \cdot \tt\big) 
\notag\\
&= \alpha  \tint \Bigl( \p_1\bigl( \jb \mathcal{J}^{\frac 32} (\Abt-\Zbn)^2 \bigr) + (\nbs_2 - \Qr_2) \bigl( \jb \mathcal{J}^{\frac 32} (\Abt-\Zbn)^2 \Jg g^{- {\frac{1}{2}} } \nbs_2 h \bigr)\Bigr)   (\nbs^6\tt \cdot\nn)^2 
\notag\\
&\qquad 
+ \alpha  \int \Qb_2 \jb \mathcal{J}^{\frac 32} (\Abt-\Zbn)^2      \Jg g^{- {\frac{1}{2}} } \nbs_2 h\,   (\nbs^6\tt\cdot \nn)^2 \Bigl|_{\s}
\notag\\
&\qquad + 2\alpha  \tint \jb \mathcal{J}^{\frac 32} (\Abt-\Zbn)^2 \nbs^6\tt\cdot\nn   \  \big(\nbs^6\tt \cdot \tt \nn,_1 \cdot\tt
-   \Jg g^{- {\frac{1}{2}} } \nbs_2 h\,  \nbs^6\tt\cdot \tt \nbs_2\nn \cdot \tt\big) 
\,.
\end{align}
By appealing to  \eqref{bootstraps}, \eqref{eq:Q:all:bbq}, \eqref{D6n-bound:b:new}, and the bound $\mathcal{J} \leq 1$, we conclude
\begin{equation}
\sabs{I^{\ZZ_\tau}_{5+7,a,ii}} 
\les 
\eps^3 \beta (\tfrac{4}{\kappa_0})^{2\beta} \mathsf{K}^2  \brak{\mathsf{B}_6}^2
\label{eq:hate:2b:bound}
\, .
\end{equation}
The term $I^{\ZZ_\tau}_{5+7,a,i}$ defined in \eqref{eq:hate:2a} requires that we integrate by parts the $(\Q \p_\s + V \p_2)$ operator, and then appeal to the $\nbs^6$-differentiated variant of \eqref{tt-evo-s}. Indeed, from \eqref{adjoint-3}, and \eqref{tt-evo-s}, we may rewrite
\begin{align}
I^{\ZZ_\tau}_{5+7,a,i}
&= 
- \tint \bigl(\nbs^6\Wbt + \nbs^6 \Zbt\bigr) \tfrac{1}{\Sigma^{2\beta}} \mathcal{J}^{\frac 32} \Jg (\Abt-\Zbn) 
  \bigl(\tfrac{1+\alpha}{2} \nbs^6 \Wbt + \tfrac{1-\alpha}{2}  \nbs^6\Zbt\bigr) 
 \notag\\
&\qquad 
- \tint \bigl(\nbs^6\Wbt + \nbs^6 \Zbt\bigr) \tfrac{1}{\Sigma^{2\beta}} \mathcal{J}^{\frac 32} \Jg (\Abt-\Zbn) 
  \bigl(\tfrac{1+\alpha}{2} \Wbt + \tfrac{1-\alpha}{2} \Zbt\bigr) \nn \cdot \nbs^6  \nn  
\notag\\
&\qquad 
- \tint \bigl(\nbs^6\Wbt + \nbs^6 \Zbt\bigr) \tfrac{1}{\Sigma^{2\beta}} \mathcal{J}^{\frac 32} \Jg (\Abt-\Zbn) \nn_k  \doublecom{\nbs^6 , \bigl(\tfrac{1+\alpha}{2} \Wbt + \tfrac{1-\alpha}{2} \Zbt\bigr) , \nn_k}
\notag\\
&\qquad - \tint \bigl(\nbs^6\Wbt + \nbs^6 \Zbt\bigr) \tfrac{1}{\Sigma^{2\beta}} \mathcal{J}^{\frac 32} \Jg (\Abt-\Zbn) \nn_k \jump{\nbs^6,V} \nbs_2 \tt_k 
\notag\\
&\qquad 
+ \tint \bigl(\nbs^6\Wbt + \nbs^6 \Zbt\bigr) \bigl(V \Qr_2 - \nbs_2 V - \Qr_\s \bigr) \Bigl( \tfrac{1}{\Sigma^{2\beta}} \mathcal{J}^{\frac 32} \Jg (\Abt-\Zbn) \nn \cdot \nbs^6\tt   \Bigr)   
\notag\\
&\qquad 
-  \tint \bigl(\nbs^6\Wbt + \nbs^6 \Zbt\bigr) \nbs^6\tt_k (\Q\p_\s+V\p_2) \Bigl( \tfrac{1}{\Sigma^{2\beta}} \mathcal{J}^{\frac 32} \Jg (\Abt-\Zbn) \nn_k   \Bigr)   
\notag\\
&\qquad 
+ \int \bigl(\nbs^6\Wbt + \nbs^6 \Zbt\bigr) \Q \Bigl( \tfrac{1}{\Sigma^{2\beta}} \mathcal{J}^{\frac 32} \Jg (\Abt-\Zbn) \nn \cdot \nbs^6\tt  \Bigr)  \Bigr|_{\s}
\notag\\
&\qquad 
-  \int \bigl(\nbs^6\Wbt + \nbs^6 \Zbt\bigr) \Q \Bigl( \tfrac{1}{\Sigma^{2\beta}} \mathcal{J}^{\frac 32} \Jg (\Abt-\Zbn) \nn \cdot \nbs^6\tt   \Bigr)  \Bigr|_{0} 
\,.
\label{eq:hate:2a:rewrite}
\end{align}
We have isolated on the right side of \eqref{eq:hate:2a:rewrite} the first term as the most dangerous term, with the remaining seven terms being lower order. Using the available bootstrap bounds and improved estimates, we deduce from \eqref{eq:hate:2a:rewrite} that
\begin{align}
\sabs{I^{\ZZ_\tau}_{5+7,a,i}}
&\les  \int_0^{\s} \Bigl(  \snorm{\tfrac{\mathcal{J}^{\frac 34}(\Jg \Q)^{\frac 12}}{\Sigma^\beta} \nbs^6 \Wbt(\cdot,\s')}_{L^2_x}^2
+    \snorm{\tfrac{\mathcal{J}^{\frac 34}(\Jg \Q)^{\frac 12}}{\Sigma^\beta} \nbs^6 \Zbt(\cdot,\s')}_{L^2_x}^2 \Bigr)
{\rm d} \s'
\notag\\
&\qquad 
+ \eps^2 (\tfrac{4}{\kappa_0})^{2\beta} \mathsf{K}^2  \brak{\mathsf{B}_6}^2
+ \eps \mathsf{K} \brak{\mathsf{B}_6} (\tfrac{4}{\kappa_0})^\beta \Bigl(\snorm{\tfrac{\mathcal{J}^{\frac 14} \Jg^{\frac 12}}{\Sigma^\beta} \nbs^6 \Wbt}_{L^2_{x,\s}} + \snorm{\tfrac{\mathcal{J}^{\frac 14} \Jg^{\frac 12}}{\Sigma^\beta} \nbs^6 \Zbt}_{L^2_{x,\s}} \Bigr)
\notag\\
&\qquad 
+ \eps^{\frac 32} (\tfrac{4}{\kappa_0})^\beta \mathsf{K} \brak{\mathsf{B}_6} \Bigl(\snorm{\tfrac{\mathcal{J}^{\frac 34} (\Jg \Q)^{\frac 12}}{\Sigma^\beta} \nbs^6 \Wbt(\cdot,\s)}_{L^2_x} + \snorm{\tfrac{\mathcal{J}^{\frac 34} (\Jg \Q)^{\frac 12}}{\Sigma^\beta} \nbs^6 \Zbt(\cdot,\s)}_{L^2_x}\Bigr)
\notag\\
&\qquad 
+ \eps^{\frac 32} (\tfrac{4}{\kappa_0})^\beta \mathsf{K} \brak{\mathsf{B}_6} \Bigl(\snorm{\tfrac{\mathcal{J}^{\frac 34} (\Jg \Q)^{\frac 12}}{\Sigma^\beta} \nbs^6 \Wbt(\cdot,0)}_{L^2_x} + \snorm{\tfrac{\mathcal{J}^{\frac 34} (\Jg \Q)^{\frac 12}}{\Sigma^\beta} \nbs^6 \Zbt(\cdot,0)}_{L^2_x}\Bigr)
\,.
\label{eq:hate:2a:bound}
\end{align}
Summarizing the decompositions~\eqref{eq:hate:1} and~\eqref{eq:hate:2},   appealing to the bounds~\eqref{eq:I:Zt:5+7:bcd}, \eqref{eq:hate:2b:bound}, \eqref{eq:hate:2a:bound}, the Cauchy-Schwartz inequality, and taking $\eps$ to be sufficiently small (at this stage, only in terms of $\alpha,\kappa_0,\Cdata$), for $\beta \geq 1$ we arrive at
\begin{align}
\sabs{I^{\ZZ_\tau}_{5} +  I^{\ZZ_\tau}_{7}}
&
\leq
\sabs{I^{\ZZ_\tau}_{5+7,a,iii}}
+ 
\Cn \eps^2 \beta^2 (\tfrac{4}{\kappa_0})^{2\beta} \mathsf{K}^2  \brak{\mathsf{B}_6}^2
\notag\\
&\qquad 
+ \Cn \eps \Bigl(\snorm{\tfrac{\mathcal{J}^{\frac 34} (\Jg \Q)^{\frac 12}}{\Sigma^\beta} \nbs^6 \Wbt(\cdot,\s)}_{L^2_x}^2 + \snorm{\tfrac{\mathcal{J}^{\frac 34} (\Jg \Q)^{\frac 12}}{\Sigma^\beta} \nbs^6 \Zbt(\cdot,\s)}_{L^2_x}^2\Bigr)
\notag\\
& \qquad 
+ \Cn \eps  \Bigl(\snorm{\tfrac{\mathcal{J}^{\frac 34} (\Jg \Q)^{\frac 12}}{\Sigma^\beta} \nbs^6 \Wbt(\cdot,0)}_{L^2_x}^2 + \snorm{\tfrac{\mathcal{J}^{\frac 34} (\Jg \Q)^{\frac 12}}{\Sigma^\beta} \nbs^6 \Zbt(\cdot,0)}_{L^2_x}^2\Bigr) 
\notag\\
&\qquad 
+\Cn \int_0^{\s} \Bigl(  \snorm{\tfrac{\mathcal{J}^{\frac 14} \Jg^{\frac 12}}{\Sigma^\beta} \nbs^6 \Wbt(\cdot,\s')}_{L^2_x}^2
+    \snorm{\tfrac{\mathcal{J}^{\frac 14} \Jg^{\frac 12}}{\Sigma^\beta} \nbs^6 \Zbt(\cdot,\s')}_{L^2_x}^2 \Bigr)
{\rm d} \s'
\,,
\label{eq:I:Zt:5+7:final}
\end{align}
where $I^{\ZZ_\tau}_{5+7,a,iii}$ is defined in \eqref{eq:hate:2c}, and will be shown to be dominated by the other terms on the right side of \eqref{eq:I:Zt:5+7:final}.

\subsubsection{The sum $ I^{\WW_\tau}_3 +  I^{\ZZ_\tau}_3+ I^{\AA_\tau}_5+ I^{\AA_\tau}_7$}
Integrating-by-parts in~\eqref{I5-Abt} and~\eqref{I7-Abt}, in a similar fashion to \eqref{eq:hate:1}, we find that
\begin{align} 
I^{\AA_\tau}_5+I^{\AA_\tau}_7
 &= - 2\alpha \tint \jb \mathcal{J}^{\frac 32} (\Omega+\Wbt+\Zbt) \nbs^6\tt\cdot\nn     \Big( \nbs^6\Abt,_1- \Jg g^{- {\frac{1}{2}} } \nbs_2h\, \nbs_2  \nbs^6\Abt\Big) 
 \notag \\
 &  \quad 
- 2 \alpha  \tint \Big( \p_1 \big( \jb \mathcal{J}^{\frac 32} (\Omega+\Wbt+\Zbt) \nn_k\big) -\nbs_2\big( \jb\mathcal{J}^{\frac 32} \Jg g^{- {\frac{1}{2}} } \nbs_2h\,  (\Omega+\Wbt+\Zbt) \nn_k\big)   \Big) \nbs^6\tt^k \   \nbs^6\Abt 
\notag \\
&  \quad
- 2 \alpha\!\! \tint \!\Qr_2 \jb \mathcal{J}^{\frac 32}  \Jg g^{- {\frac{1}{2}} } \nbs_2h\,  (\Omega\!+\!\Wbt\!+\!\Zbt) \nbs^6\tt\cdo \nn \  \nbs^6\Abt 
\notag\\
&\quad
+2\alpha \!\! \int \!\! \Qb_2 \jb  \mathcal{J}^{\frac 32}  \Jg g^{- {\frac{1}{2}} } \nbs_2h\,  (\Omega\!+\!\Wbt\!+\!\Zbt) \nbs^6\tt\cdot\nn \  \nbs^6\Abt\Big|_\s
\notag \\
& =:  I^{\AA_\tau}_{5+7,a}+I^{\AA_\tau}_{5+7,b}+I^{\AA_\tau}_{5+7,c}+I^{\AA_\tau}_{5+7,d} \,.
\label{eq:hate:3}
\end{align} 
Note that the first term in \eqref{eq:hate:2}, namely $I^{\AA_\tau}_{5+7,a}$, contains  terms with seven derivatives on $\Abt$, but that the remaining three terms are under control. Indeed, from \eqref{bootstraps}, \eqref{eq:Q:all:bbq}, \eqref{D6n-bound:b:new}, \eqref{eq:vorticity:pointwise} and the fact that $\mathcal{J}\leq \Jg$, in analogy to~\eqref{eq:I:Zt:5+7:bcd}  we have 
\begin{subequations}
\label{eq:I:At:5+7:bcd}
\begin{align}
\sabs{I^{\AA_\tau}_{5+7,b}}
&  
\les   \eps \beta  (\tfrac{4}{\kappa_0})^{\beta} \mathsf{K} \brak{\mathsf{B}_6}
\snorm{\tfrac{\mathcal{J}^{\frac 14} \Jg^{\frac 12}}{\Sigma^\beta} \nbs^6 \Abt}_{L^2_{x,\s}}
\\
\sabs{I^{\AA_\tau}_{5+7,c}}
& 
\les \eps^3  (\tfrac{4}{\kappa_0})^{\beta} \mathsf{K} \brak{\mathsf{B}_6}
\snorm{\tfrac{\mathcal{J}^{\frac 14} \Jg^{\frac 12}}{\Sigma^\beta} \nbs^6 \Abt}_{L^2_{x,\s}}
\\
\sabs{I^{\AA_\tau}_{5+7,d}}
& 
\les \eps^{\frac 72} (\tfrac{4}{\kappa_0})^\beta \mathsf{K} \brak{\mathsf{B}_6} \snorm{\tfrac{\mathcal{J}^{\frac 34} (\Jg \Q)^{\frac 12}}{\Sigma^\beta} \nbs^6 \Abt(\cdot,\s)}_{L^2_x}
\,.
\end{align}
\end{subequations}
In order to handle the first term on the right side of \eqref{eq:hate:3}, using equation \eqref{energy-At-s}, we see that
\begin{align*} 
&
-2 \alpha \big( \nbs^6 \Abt,_1  -   \Jg g^{- {\frac{1}{2}} } \nbs_2 h \, \nbs_2\nbs^6\Abt \big)
\notag\\
&=
- 2 \alpha  \Jg  g^{- {\frac{1}{2}} } \nbs_2\nbs^6 \Sbt
+2 \alpha    g^{- {\frac{1}{2}} } (\Jg\Sbn)\nbs_2\nbs^6 \tt\cdot\nn 
-2\tfrac{\Jg}{\Sigma} (\Q\p_\s +V\p_2) \nbs^6\Abt
\notag \\
& \qquad 
- 2 \alpha (\Omega+\Wbt+\Zbt) \nbs^6 \tt,_1\cdot\nn 
+2 \alpha  \Jg g^{- {\frac{1}{2}} } \nbs_2 h\,  (\Omega+\Wbt+\Zbt) \nbs_2\nbs^6\tt\cdot \nn 
+\big(\nbs^6 \Fat   + \mathcal{R}_\Ab^\tt + \mathcal{C}_\Ab^\tt  \big)\,.
\end{align*} 
Substitution of this identity into the integral defining $I^{\AA_\tau}_{5+7,a}$ in \eqref{eq:hate:3}, gives us the further decomposition
\begin{subequations}
\label{eq:hate:4}
\begin{align} 
I^{\AA_\tau}_{5+7,a} & = I^{\AA_\tau}_{5+7,a,i}+I^{\AA_\tau}_{5+7,a,ii}+I^{\AA_\tau}_{5+7,a,iii}+I^{\AA_\tau}_{5+7,a,iv}+I^{\AA_\tau}_{5+7,a,v}
 \,,
\notag \\
I^{\AA_\tau}_{5+7,a,i} 
&= -2 \alpha  \tint \jb  \mathcal{J}^{\frac 32} \Jg g^{- {\frac{1}{2}} } (\Omega+\Wbt+\Zbt) \nbs^6\tt\cdot\nn   \   \nbs_2\nbs^6\Sbt \,,
\label{eq:hate:4a}
 \\
I^{\AA_\tau}_{5+7,a,ii} 
&=   \alpha  \tint \jb  \mathcal{J}^{\frac 32} g^{- {\frac{1}{2}} } (\Jg\Wbn - \Jg \Zbn) (\Omega+\Wbt+\Zbt)\nbs^6\tt\cdot\nn   \ \nbs_2\nbs^6 \tt\cdot\nn    \,,
\label{eq:hate:4b} \\
I^{\AA_\tau}_{5+7,a,iii} 
 &=-2\tint \tfrac{ 1}{\Sigma^{2\beta}}  \mathcal{J}^{\frac 32} \Jg (\Omega+\Wbt+\Zbt) \nbs^6\tt\cdot\nn   \  (\Q\p_s+V\p_2) \nbs^6\Abt \,,
\label{eq:hate:4c} \\
I^{\AA_\tau}_{5+7,a,iv} 
&= - 2\alpha  \tint \jb  \mathcal{J}^{\frac 32}  g^{- {\frac{1}{2}} } (\Omega+\Wbt+\Zbt)^2 \nbs^6\tt\cdot\nn   \  \big(\nbs^6\tt,_1\cdot\nn 
-   \Jg g^{- {\frac{1}{2}} } \nbs_2 h\,  \nbs_2\nbs^6\tt\cdot \nn \big)   \,,
\label{eq:hate:4d} \\
I^{\AA_\tau}_{5+7,a,v} 
&=  \tint \jb  \mathcal{J}^{\frac 32} (\Omega+\Wbt+\Zbt) \nbs^6\tt\cdot\nn   \ \big(\nbs^6 \Fat   + \mathcal{R}_\Ab^\tt + \mathcal{C}_\Ab^\tt \big) \,.
\label{eq:hate:4e}
\end{align} 
\end{subequations}
The first term in the above expression, $I^{\AA_\tau}_{5+7,a,i}$, requires careful analysis since it involves seven derivatives on $\Sbt$ (hence $\Wbt$ and $\Zbt$). We may however estimate the remaining terms  by drawing an analogy with \eqref{eq:hate:2}. Save for the difference between $g^{-\frac 12}(\Omega + \Wbt + \Zbt)$ and $\Abt - \Zbn$ (these terms satisfy the same bounds in $L^\infty_{x,\s}$ even with one $\nbs$ derivative acting on them), we have that $I^{\AA_\tau}_{5+7,a,iii}$ (see \eqref{eq:hate:4c}) is nearly identical to $I^{\ZZ_\tau}_{5+7,a,i}$ (see \eqref{eq:hate:2a}),  $I^{\AA_\tau}_{5+7,a,iv}$ (see \eqref{eq:hate:4d}) is nearly identical to $I^{\ZZ_\tau}_{5+7,a,ii}$ (see \eqref{eq:hate:2b}), and $I^{\AA_\tau}_{5+7,a,ii}$ (see~\eqref{eq:hate:4b}) is very similar to the $\nbs_2$-part of $I^{\ZZ_\tau}_{5+7,a,ii}$ (see \eqref{eq:hate:2b}). To avoid redundancy, we do not repeat the steps which have lead to the bounds~\eqref{eq:hate:2a} and~\eqref{eq:hate:2b}; instead, we argue similarly and deduce that 
\begin{subequations}
\label{eq:I:At:5+7:a:ii:iv}
\begin{align}
\sabs{I^{\AA_\tau}_{5+7,a,ii}} + \sabs{I^{\AA_\tau}_{5+7,a,iv}} 
&\les 
\eps^3 \beta (\tfrac{4}{\kappa_0})^{2\beta} \mathsf{K}^2  \brak{\mathsf{B}_6}^2
, \\
\sabs{I^{\AA_\tau}_{5+7,a,iii}} 
&\les  \int_0^{\s}  \snorm{\tfrac{\mathcal{J}^{\frac 34}(\Jg \Q)^{\frac 12}}{\Sigma^\beta} \nbs^6 \Abt(\cdot,\s')}_{L^2_x}^2 {\rm d} \s'
+ \eps \mathsf{K} \brak{\mathsf{B}_6} (\tfrac{4}{\kappa_0})^\beta  
\snorm{\tfrac{\mathcal{J}^{\frac 14} \Jg^{\frac 12}}{\Sigma^\beta} \nbs^6 \Abt}_{L^2_{x,\s}} \notag\\
&\qquad 
+ \eps^{\frac 32} (\tfrac{4}{\kappa_0})^\beta \mathsf{K} \brak{\mathsf{B}_6} \snorm{\tfrac{\mathcal{J}^{\frac 34} (\Jg \Q)^{\frac 12}}{\Sigma^\beta} \nbs^6 \Abt(\cdot,\s)}_{L^2_x}  
\notag\\
&\qquad 
+ \eps^{\frac 32} (\tfrac{4}{\kappa_0})^\beta \mathsf{K} \brak{\mathsf{B}_6} \snorm{\tfrac{\mathcal{J}^{\frac 34} (\Jg \Q)^{\frac 12}}{\Sigma^\beta} \nbs^6 \Abt(\cdot,0)}_{L^2_x}  
+ \eps^4   (\tfrac{4}{\kappa_0})^{2\beta} \mathsf{K}^2  \brak{\mathsf{B}_6}^2
.
\end{align} 
\end{subequations}
The forcing and commutator terms   in $I^{\AA_\tau}_{5+7,a,v}$ will be shown in~Subsection~\ref{sec:D6:tau:forcing:commutator} below to satisfy similar bounds. 

This leaves us with the most delicate term, namely $I^{\AA_\tau}_{5+7,a,i}$ (see~\eqref{eq:hate:4a}). It turns out that while a good bound for this term is not available, this term completely cancels with a piece of the sum $ I^{\WW_\tau}_3 +  I^{\ZZ_\tau}_3$. To see this, we next consider the sum $ I^{\WW_\tau}_3 +  I^{\ZZ_\tau}_3$, given by \eqref{I3-Wbt} and \eqref{I3-Zbt}.  Using \eqref{adjoint-2}, we have that
\begin{subequations}
\label{eq:I:Wt:Zt:3}
\begin{align} 
 I^{\WW_\tau}_3 +  I^{\ZZ_\tau}_3 
  &= -2\alpha \tint \jb  \mathcal{J}^{\frac 32} \Jg g^{-\frac 12} (\Omega+\Wbt+\Zbt) \nbs_2\nbs^6\tt \cdot \nn 
 \ \nbs^6\Sbt 
 \notag \\
  &=: I^{\WW_\tau+\ZZ_\tau}_{3,a} + I^{\WW_\tau+\ZZ_\tau}_{3,b} + I^{\WW_\tau+\ZZ_\tau}_{3,c} + I^{\WW_\tau+\ZZ_\tau}_{3,d} 
\end{align}
where the four terms arising from integrating by parts the $\nbs_2$ operator are defined as
\begin{align}
I^{\WW_\tau+\ZZ_\tau}_{3,a}
&= 2\alpha \tint \jb \mathcal{J}^{\frac 32} \Jg g^{-\frac 12} (\Omega+\Wbt+\Zbt) \nbs^6\tt \cdot \nn  \ \nbs_2\nbs^6\Sbt
\\
I^{\WW_\tau+\ZZ_\tau}_{3,b}
&= 2\alpha \tint \nbs_2\big( \jb \mathcal{J}^{\frac 32} \Jg g^{-\frac 12} (\Omega+\Wbt+\Zbt) \nn_k\big) \nbs^6\tt^k  \ \nbs^6\Sbt
\\
I^{\WW_\tau+\ZZ_\tau}_{3,c} 
&= -  2\alpha \tint \Qr_2 \jb  \mathcal{J}^{\frac 32} \Jg g^{-\frac 12} (\Omega+\Wbt+\Zbt) \nbs^6\tt \cdot \nn  \ \nbs^6\Sbt
\\
I^{\WW_\tau+\ZZ_\tau}_{3,d}
&=  2\alpha   \int \Qb_2 \jb \mathcal{J}^{\frac 32} \Jg g^{-\frac 12}  (\Omega+\Wbt+\Zbt) \nbs^6\tt \cdot \nn  \ \nbs^6\Sbt\Big|_\s \,.
\end{align} 
\end{subequations}
We notice that the term $I^{\WW_\tau+\ZZ_\tau}_{3,a}$ precisely cancels the term $I^{\AA_\tau}_{5+7,a,i}$, as expected:
\begin{equation*}
I^{\AA_\tau}_{5+7,a,i}  + I^{\WW_\tau+\ZZ_\tau}_{3,a} = 0
\,.
\end{equation*}
The remaining terms in \eqref{eq:I:Wt:Zt:3} are bounded using  \eqref{bootstraps}, \eqref{eq:Q:all:bbq}, \eqref{D6n-bound:b:new}, \eqref{D6n-bound:b:new:bdd}, \eqref{eq:vorticity:pointwise}, and the bound $\mathcal{J} \leq \Jg$, as 
\begin{subequations}
\label{eq:I:Wt:Zt:3:bcd}
\begin{align}
\sabs{I^{\WW_\tau+\ZZ_\tau}_{3,b}} 
+ 
\sabs{I^{\WW_\tau+\ZZ_\tau}_{3,c}}
& 
\les \eps^2  (\tfrac{4}{\kappa_0})^{\beta} \mathsf{K} \brak{\mathsf{B}_6}
\Bigl(
\snorm{\tfrac{\mathcal{J}^{\frac 14} \Jg^{\frac 12}}{\Sigma^\beta} \nbs^6 \Wbt}_{L^2_{x,\s}}
+
\snorm{\tfrac{\mathcal{J}^{\frac 14} \Jg^{\frac 12}}{\Sigma^\beta} \nbs^6 \Zbt}_{L^2_{x,\s}}
\Bigr)
\\
\sabs{I^{\WW_\tau+\ZZ_\tau}_{3,d}}
& 
\les \eps^{\frac 52} (\tfrac{4}{\kappa_0})^\beta \mathsf{K} \brak{\mathsf{B}_6} \Bigl(
\snorm{\tfrac{\mathcal{J}^{\frac 34} (\Jg \Q)^{\frac 12}}{\Sigma^\beta} \nbs^6 \Wbt(\cdot,\s)}_{L^2_x}
+ 
\snorm{\tfrac{\mathcal{J}^{\frac 34} (\Jg \Q)^{\frac 12}}{\Sigma^\beta} \nbs^6 \Zbt(\cdot,\s)}_{L^2_x}
\Bigr)
\,.
\end{align}
\end{subequations}

Summarizing the decompositions~\eqref{eq:hate:3}, \eqref{eq:hate:4},~\eqref{eq:I:Wt:Zt:3}, appealing to the bounds~\eqref{eq:I:At:5+7:bcd}, \eqref{eq:I:At:5+7:a:ii:iv}, \eqref{eq:I:Wt:Zt:3:bcd}, the Cauchy-Schwartz inequality, and taking $\eps$ to be sufficiently small (only in terms of $\alpha,\kappa_0,\Cdata$), for $\beta \geq 1$ we arrive at 
\begin{align} 
&\sabs{I^{\WW_\tau}_3 +  I^{\ZZ_\tau}_3+ I^{\AA_\tau}_5+I^{\AA_\tau}_7} 
\notag\\
&\leq
\sabs{I^{\AA_\tau}_{5+7,a,v} }
+ \Cn \eps^2 \beta^2 (\tfrac{4}{\kappa_0})^{2\beta} \mathsf{K}^2  \brak{\mathsf{B}_6}^2
+ \Cn \eps \snorm{\tfrac{\mathcal{J}^{\frac 34} (\Jg \Q)^{\frac 12}}{\Sigma^\beta} \nbs^6 \Abt(\cdot,0)}_{L^2_x}^2  
\notag\\&\qquad
+ \Cn \int_0^{\s} 
\Bigl(
 \snorm{\tfrac{\mathcal{J}^{\frac 14} \Jg^{\frac 12}}{\Sigma^\beta} \nbs^6 \Wbt(\cdot,\s')}_{L^2_x}^2 {\rm d} \s'
+
 \snorm{\tfrac{\mathcal{J}^{\frac 14} \Jg^{\frac 12}}{\Sigma^\beta} \nbs^6 \Zbt(\cdot,\s')}_{L^2_x}^2 {\rm d} \s'
+
 \snorm{\tfrac{\mathcal{J}^{\frac 14} \Jg^{\frac 12}}{\Sigma^\beta} \nbs^6 \Abt(\cdot,\s')}_{L^2_x}^2 
 \Bigr) {\rm d} \s'
\notag\\&\qquad
+ \Cn \eps  
\Bigl(
\snorm{\tfrac{\mathcal{J}^{\frac 34} (\Jg \Q)^{\frac 12}}{\Sigma^\beta} \nbs^6 \Wbt(\cdot,\s)}_{L^2_x}
+ 
\snorm{\tfrac{\mathcal{J}^{\frac 34} (\Jg \Q)^{\frac 12}}{\Sigma^\beta} \nbs^6 \Zbt(\cdot,\s)}_{L^2_x}
+
\snorm{\tfrac{\mathcal{J}^{\frac 34} (\Jg \Q)^{\frac 12}}{\Sigma^\beta} \nbs^6 \Abt(\cdot,\s)}_{L^2_x}^2 
\Bigr)
 \,.
\label{eq:hate:5}
\end{align} 

\subsubsection{The integral $I^{\AA_\tau}_3$ }
\label{sec:nasty:tangential}
The only term with over-differentiated geometry that we have not yet considered is $I^{\AA_\tau}_3$, as defined in~\eqref{I3-Abt}. When compared to the $I^{\WW_\tau}_3$ and $I^{\ZZ_\tau}_3$ terms which we have just bounded in \eqref{eq:hate:5}, we expect $I^{\AA_\tau}_3$ to be worst behaved because it contains a copy of $\Jg \Sbn$, instead of $\Jg \Sbt$. We start by integrating by parts the $\nbs_2$ derivative from \eqref{I3-Abt}, and by using  \eqref{adjoint-2}, we arrive at
\begin{align} 
I^{\AA_\tau}_3 &=
-2\alpha \tint \jb \mathcal{J}^{\frac 32}  g^{- {\frac{1}{2}} }( \Jg\Sbn) \nbs_2\nbs^6\tt \cdot \nn  \ \nbs^6\Abt
\notag\\
&=
 \alpha \tint \jb \mathcal{J}^{\frac 32} g^{- {\frac{1}{2}} }  ( \Jg\Wbn - \Jg \Zbn) \nbs^6\tt \cdot \nn  \  \nbs_2\nbs^6\Abt 
\notag\\&\qquad
+ \tfrac{3\alpha}{2} \tint \jb \mathcal{J}^{\frac 12} g^{- {\frac{1}{2}} }(\Jg\Wbn - \Jg \Zbn) (\nbs_2 \mathcal{J})\nbs^6\tt \cdot \nn  \ \nbs^6\Abt
\notag\\&\qquad
+  \alpha \tint \mathcal{J}^{\frac 32} \nbs_2\big( \jb   g^{- {\frac{1}{2}} }(\Jg\Wbn - \Jg \Zbn) \nn_k \big)\nbs^6\tt^k  \ \nbs^6\Abt
\notag \\ & \qquad 
- \alpha \tint \Qr_2 \jb \mathcal{J}^{\frac 32}  g^{- {\frac{1}{2}} }(\Jg\Wbn - \Jg \Zbn) \nbs^6\tt \cdot \nn  \ \nbs^6\Abt
\notag\\&\qquad
+ \alpha  \int \Qb_2 \jb \mathcal{J}^{\frac 32}  g^{- {\frac{1}{2}} }(\Jg\Wbn - \Jg \Zbn) \nbs^6\tt \cdot \nn  \ \nbs^6\Abt\Big|_\s 
\notag \\
&=: I^{\AA_\tau}_{3,a} + I^{\AA_\tau}_{3,b} + I^{\AA_\tau}_{3,c}+ I^{\AA_\tau}_{3,d}+ I^{\AA_\tau}_{3,e}
\,.
\label{eq:hate:6}
\end{align} 
In analogy to \eqref{eq:hate:1}--\eqref{eq:I:Zt:5+7:bcd}, \eqref{eq:hate:3}--\eqref{eq:I:At:5+7:bcd}, and \eqref{eq:I:Wt:Zt:3}--\eqref{eq:I:Wt:Zt:3:bcd}, the last thee terms on the right side of \eqref{eq:hate:6} may be bounded directly (recall that $\mathcal{J}\leq \Jg$) as
\begin{subequations}
\label{eq:I:At:3:bcde}
\begin{align}
\sabs{I^{\AA_\tau}_{3,c}} + \sabs{I^{\AA_\tau}_{3,d}}
& 
\les \eps \beta (\tfrac{4}{\kappa_0})^{\beta} \mathsf{K} \brak{\mathsf{B}_6}
\snorm{\tfrac{\mathcal{J}^{\frac 14} \Jg^{\frac 12}}{\Sigma^\beta} \nbs^6 \Abt}_{L^2_{x,\s}}
\\
\sabs{I^{\AA_\tau}_{3,e}}
& 
\les \eps^{\frac 32} (\tfrac{4}{\kappa_0})^\beta \mathsf{K} \brak{\mathsf{B}_6} \snorm{\tfrac{\mathcal{J}^{\frac 34} (\Jg \Q)^{\frac 12}}{\Sigma^\beta} \nbs^6 \Abt(\cdot,\s)}_{L^2_x}
\,.
\end{align}
The second term on the right side of \eqref{eq:hate:6} requires special care because we only have $\mathcal{J}^{\frac 12}$  to be paired with $\nbs^6 \Abt$, but the former requires $\mathcal{J}^{\frac 34}$ in order to be bounded suitably in $L^2_{x,\s}$.  The additional power of $\mathcal{J}^{\frac 14}$ that we are missing arises from the fact that the $\nbs^6 \tt$ bound in \eqref{D6h2Energy:new} in fact carries $\mathcal{J}^{-\frac 14}$. As such, we obtain
\begin{equation}
\sabs{I^{\AA_\tau}_{3,b}}  
\leq (\tfrac{4}{\kappa_0})^{\beta} \|\Jg \Wbn - \Jg \Zbn\|_{L^\infty_{x,\s}}  \|\mathcal{J}^{-\frac 14} \nbs^6 \tt\|_{L^2_{x,\s}} \snorm{\tfrac{\mathcal{J}^{\frac 34}}{\Sigma^\beta} \nbs^6 \Abt}_{L^2_{x,\s}} 
\les   \eps    (\tfrac{4}{\kappa_0})^{\beta} \mathsf{K} \brak{\mathsf{B}_6}
\snorm{\tfrac{\mathcal{J}^{\frac 34}}{\Sigma^\beta} \nbs^6 \Abt}_{L^2_{x,\s}}
\,.
\end{equation}
\end{subequations}
The only term that we are left to estimate in the decomposition \eqref{eq:hate:6} is $I^{\AA_\tau}_{3,a}$. For this term, we use the equation \eqref{energy-Wt-s} to further rewrite the term over-differentiated term $\nbs_2 \nbs^6 \Abt$. To be precise, we may derive the identity 
\begin{subequations}
\label{eq:hate:7}
\begin{align} 
I^{\AA_\tau}_{3,a} & = I^{\AA_\tau}_{3,a,i} + I^{\AA_\tau}_{3,a,ii} + I^{\AA_\tau}_{3,a,iii} \,,
\notag \\
 I^{\AA_\tau}_{3,a,i} &= 
 \alpha \tint \jb \mathcal{J}^{\frac 32} g^{- {\frac{1}{2}} }  (\Jg\Wbn - \Jg \Zbn)  (\Omega+\Wbt+\Zbt) \nbs^6\tt \cdot \nn  \   \nbs_2\nbs^6 \tt\cdot\nn \,,
\label{eq:hate:7a} \\
I^{\AA_\tau}_{3,a,ii}  &= 
-  \tint \tfrac{1}{\Sigma^{2\beta}}  \mathcal{J}^{\frac 32}  (\Jg\Wbn - \Jg \Zbn) \nbs^6\tt \cdot \nn  \ (\Q\p_\s+V\p_2) \nbs^6\Wbt \,,
\label{eq:hate:7b} \\
  I^{\AA_\tau}_{3,a,iii} &= 
    \tint \jb\mathcal{J}^{\frac 32}   (\Jg\Wbn - \Jg \Zbn) \nbs^6\tt \cdot \nn  \  \big(\nbs^6\Fwt  + \mathcal{R}_\Wb^\tt + \mathcal{C}_\Wb^\tt\big)\,.
  \label{eq:hate:7c}
\end{align} 
\end{subequations}
The bound for the $I^{\AA_\tau}_{3,a,i}$ is similar to bounds which were obtained earlier. We rewrite 
\begin{equation*}
\nbs^6\tt \cdot \nn \nbs_2\nbs^6 \tt\cdot\nn = \tfrac 12 \nbs_2 (\nbs^6\tt \cdot \nn)^2 - \tt \cdot \nbs_2 \nn \nbs^6 \tt \cdot \nn (\nbs^6 \tt - g^{-1} \nn \nbs^6 \nbs_2 h) \cdot \tt
\end{equation*}
and for the first term integrate by parts the $\nbs_2$ derivative. Using the bootstrap inequalities~\eqref{bootstraps}, the estimates~\eqref{eq:Q:all:bbq}, and the bounds~\eqref{geometry-bounds-new} for the geometry, we deduce
\begin{equation}
\sabs{ I^{\AA_\tau}_{3,a,i} }
\les \eps^3 \beta (\tfrac{4}{\kappa_0})^{2\beta} \mathsf{K}^2 \brak{\mathsf{B}_6}^2
\,.
\label{eq:hate:7a:bound}
\end{equation}
In order to obtain a bound for the $I^{\AA_\tau}_{3,a,ii}$ term in \eqref{eq:hate:7}, similarly to \eqref{eq:hate:2a} and \eqref{eq:hate:2a:rewrite}, we need to integrate by parts the $(\Q\p_\s + V \p_2)$ off from $\nbs^6 \Wbt$, leading to the identity 
\begin{align}
I^{\AA_\tau}_{3,a,ii}
&= 
\tint \nbs^6\Wbt \tfrac{1}{\Sigma^{2\beta}}  \mathcal{J}^{\frac 32}  (\Jg\Wbn - \Jg \Zbn)
  \bigl(\tfrac{1+\alpha}{2} \nbs^6 \Wbt + \tfrac{1-\alpha}{2}  \nbs^6\Zbt\bigr)
 \notag\\
 &\qquad 
+  \tint \nbs^6\Wbt \nbs^6\tt_k (\Q\p_\s+V\p_2) \Bigl( \tfrac{1}{\Sigma^{2\beta}}  \mathcal{J}^{\frac 32}  (\Jg\Wbn - \Jg \Zbn) \nn_k   \Bigr)   
\notag\\
&\qquad 
- \int \nbs^6\Wbt \Q \Bigl( \tfrac{1}{\Sigma^{2\beta}}  \mathcal{J}^{\frac 32}  (\Jg\Wbn - \Jg \Zbn) \nn \cdot \nbs^6\tt   \Bigr)  \Bigr|_{\s}
\notag\\
&\qquad 
+ \tint \nbs^6\Wbt \tfrac{1}{\Sigma^{2\beta}}  \mathcal{J}^{\frac 32}  (\Jg\Wbn - \Jg \Zbn)
  \bigl(\tfrac{1+\alpha}{2} \Wbt + \tfrac{1-\alpha}{2} \Zbt\bigr) \nn \cdot \bigl(\nbs^6  \nn + g^{-1} \tt \nbs^6 \nbs_2 h \bigr) 
\notag\\
&\qquad 
+ \tint \nbs^6\Wbt \tfrac{1}{\Sigma^{2\beta}}  \mathcal{J}^{\frac 32}  (\Jg\Wbn - \Jg \Zbn) \nn_k  \doublecom{\nbs^6 , \bigl(\tfrac{1+\alpha}{2} \Wbt + \tfrac{1-\alpha}{2} \Zbt\bigr) , \nn_k}
\notag\\
&\qquad 
+ \tint \nbs^6\Wbt \tfrac{1}{\Sigma^{2\beta}}  \mathcal{J}^{\frac 32}  (\Jg\Wbn - \Jg \Zbn) \nn_k \jump{\nbs^6,V} \nbs_2 \tt_k 
\notag\\
&\qquad 
- \tint \nbs^6\Wbt \bigl(V \Qr_2 - \nbs_2 V -   \Qr_\s  \bigr) \Bigl( \tfrac{1}{\Sigma^{2\beta}}  \mathcal{J}^{\frac 32}  (\Jg\Wbn - \Jg \Zbn) \nn \cdot \nbs^6\tt   \Bigr)   
\notag\\
&\qquad 
+  \int \nbs^6\Wbt \Q \Bigl( \tfrac{1}{\Sigma^{2\beta}}  \mathcal{J}^{\frac 32}  (\Jg\Wbn - \Jg \Zbn) \nn \cdot \nbs^6\tt   \Bigr)  \Bigr|_{0} 
\,.
\label{eq:hate:7b:rewrite}
\end{align}
At this stage, we denote the first three lines on the right side of \eqref{eq:hate:7b:rewrite} by $\mathsf{M}_1$, $\mathsf{M}_2$, and $\mathsf{M}_3$, and note  the existing bounds and the initial data assumption~\eqref{table:derivatives} imply that 
\begin{align}
\sabs{I^{\AA_\tau}_{3,a,ii} -  \mathsf{M}_1 - \mathsf{M}_2 - \mathsf{M}_3}
&\leq \Cn \eps  (\tfrac{4}{\kappa_0})^\beta  \mathsf{K} \brak{\mathsf{B}_6} \snorm{\tfrac{\mathcal{J}^{\frac 14} \Jg^{\frac 12}}{\Sigma^\beta} \nbs^6 \Wbt}_{L^2_{x,\s}}
+
 6(1+\alpha) \eps    (\tfrac{3}{\kappa_0})^{2\beta}    \Cdatatwo
\notag\\&\qquad
+\tfrac{500^2}{\eps^2 \sqrt{1+\alpha}} 
\snorm{\tfrac{\mathcal{J}^{\frac 34} (\Q \Jg)^{\frac 12}}{\Sigma^\beta} \nbs^6 \Wbt}_{L^2_{x,\s}} 
\snorm{\tfrac{\mathcal{J}^{\frac 14} \Q}{\Sigma^\beta} \nn \cdot \nbs^6 \tt}_{L^2_{x,\s}} 
 \label{eq:hate:7b:bound:initial}
 \,.
\end{align}
It thus remains to estimate the three bad terms in \eqref{eq:hate:7b:rewrite}, namely $\mathsf{M}_1,\mathsf{M}_2, \mathsf{M}_3$. Concerning the second and third one, we have
\begin{align}
&\left| \mathsf{M}_2
 - \tint \nbs^6\Wbt \nn \cdot \nbs^6\tt \tfrac{1}{\Sigma^{2\beta}} \Jg\Wbn  (\Q\p_\s+V\p_2)\mathcal{J}^{\frac 32} \right|
\les \eps  (\tfrac{4}{\kappa_0})^\beta  \mathsf{K} \brak{\mathsf{B}_6} \snorm{\tfrac{\mathcal{J}^{\frac 14} \Jg^{\frac 12}}{\Sigma^\beta} \nbs^6 \Wbt}_{L^2_{x,\s}}
\label{eq:WTF:*}
\\
&\left| \mathsf{M}_3 + \int \nbs^6\Wbt  \nn \cdot \nbs^6\tt   \tfrac{1}{\Sigma^{2\beta}}     \Jg\Wbn    \Q \mathcal{J}^{\frac 32}   \Bigr|_{\s} \right|
\les  \eps^{\frac 32}  (\tfrac{4}{\kappa_0})^\beta \mathsf{K} \brak{\mathsf{B}_6}  \snorm{\tfrac{\mathcal{J}^{\frac 34} (\Jg \Q)^{\frac 12}}{\Sigma^\beta} \nbs^6 \Wbt(\cdot,\s)}_{L^2_x} 
\label{eq:WTF:**}
\,,
\end{align}
and the terms on the RHS of the above have an acceptable size. 
At this point we note that by the definition of $\mathcal{J}$, 
\begin{equation*}
(\Q\p_\s+V\p_2)\mathcal{J}^{\frac 32}
=
- \tfrac{3}{2\eps}  \Q  \mathcal{J}^{\frac 12}
\,.
\end{equation*}
Thus, we can write the left side of \eqref{eq:WTF:*} as $| \mathsf{M}_2 - \tilde{\mathsf{M}}_2|$ and the left side of \eqref{eq:WTF:**} as $| \mathsf{M}_3 - \tilde{\mathsf{M}}_3|$, where we define 
\begin{align}
 \tilde{\mathsf{M}}_2
 &= - \tfrac{3}{2\eps} \tint \nbs^6\Wbt \nn \cdot \nbs^6\tt \tfrac{1}{\Sigma^{2\beta}} \Jg\Wbn  \Q  \mathcal{J}^{\frac 12}
  \label{eq:Really:Bad:2:def}
 \,, \\
 \tilde{\mathsf{M}}_3
 &=  - \int \nbs^6\Wbt  \nn \cdot \nbs^6\tt   \tfrac{1}{\Sigma^{2\beta}}     \Jg\Wbn    \Q \mathcal{J}^{\frac 32}   \Bigr|_{\s}  \, .
 \label{eq:Really:Bad:3:def}
\end{align}
Next, we try to manipulate the $\tilde{\mathsf{M}}_2$ term. 
We rewrite 
\begin{align}
 \tilde{\mathsf{M}}_2
 &= - \tfrac{3}{2\eps} \tint \nbs^6\Wbt \nn \cdot \nbs^6\tt \tfrac{1}{\Sigma^{2\beta}} \bigl(\Jg\Wbn - \tfrac{13}{\eps} \Jg \bigr)  \Q  \mathcal{J}^{\frac 12} 
 - \tfrac{39}{2 \eps^2} \tint \nbs^6\Wbt \nn \cdot \nbs^6\tt \tfrac{1}{\Sigma^{2\beta}}     \Q  \mathcal{J}^{\frac 12} \Jg 
 \notag\\
 &=:  \tilde{\mathsf{M}}_2' +  \tilde{\mathsf{M}}_2''.
 \label{eq:really:bad:2}
\end{align}
The second term in \eqref{eq:really:bad:2} is estimated using Cauchy-Schwartz, \eqref{eq:Q:bbq}, and \eqref{bs-Jg-simple} as
\begin{equation}
\sabs{\tilde{\mathsf{M}}_2''}
\leq \tfrac{34}{\eps^2 \sqrt{1+\alpha}} \int_0^{\s}  \snorm{\tfrac{\mathcal{J}^{\frac 34} (\Jg \Q)^{\frac 12}}{\Sigma^\beta} \nbs^6 \Wbt(\cdot,\s')}_{L^2_x}  \snorm{\tfrac{\mathcal{J}^{-\frac 14}  \Q}{\Sigma^\beta} \nn \cdot \nbs^6 \tt(\cdot,\s')}_{L^2_x} {\rm d} \s'
\,. \label{eq:really:bad:2:double:prime}
\end{equation}
Next, we consider the term $\tilde{\mathsf{M}}_2'$ in \eqref{eq:really:bad:2}. We apply $\nbs^6$ to \eqref{tt-evo-s},  and obtain 
\begin{equation}
\nn \cdot (\Q\p_\s+V\p_2)\nbs^6 \tt  
=  \nn^k   \nbs^6 \Bigl( \bigl(\tfrac{1+ \alpha }{2} \Wbt + \tfrac{1- \alpha }{2} \Zbt \bigr) \nn^k \Bigr)  - \nn^k \jump{ \nbs^6,  V } \nbs_2 \tt^k \,,
\end{equation}
which can be further manipulated to read
\begin{align}
(\Q\p_\s+V\p_2)\bigl( \nn \cdot  \nbs^6 \tt  \bigr)
& =  \tfrac{1+ \alpha }{2}  \nbs^6 \Wbt  + \tfrac{1- \alpha }{2}  \nbs^6 \Zbt  
+  \bigl(\tfrac{1+ \alpha }{2} \Wbt  + \tfrac{1- \alpha }{2} \Zbt \bigr) 
\bigl( \nn \cdot   \nbs^6 \nn   - \tt \cdot \nbs^6 \tt \bigr)
\notag\\
&\qquad 
+  \nn^k \doublecom{  \nbs^6, \bigl(\tfrac{1+ \alpha }{2} \Wbt  + \tfrac{1- \alpha }{2} \Zbt \bigr), \nn^k}
- \nn^k \jump{ \nbs^6,  V } \nbs_2 \tt^k \,.
 \label{eq:hate:10-old}
\end{align}
We deduce from \eqref{eq:hate:10-old} that 
\begin{equation*}
\snorm{ \tfrac{1+ \alpha }{2}  \nbs^6 \Wbt + \tfrac{1- \alpha }{2}  \nbs^6 \Zbt  - (\Q\p_\s+V\p_2)\bigl( \nn \cdot  \nbs^6 \tt  \bigr) }_{L^2_{x,\s}}
 \les 
\mathsf{K} \eps^2 \brak{\mathsf{B}_6}
\,,
\end{equation*}
and hence,
\begin{equation}
\left|
\tilde{\mathsf{M}}_2'
+ \tfrac{3}{(1+\alpha)\eps} \tint \bigl((\Q\p_\s+V\p_2)\bigl( \nn \cdot  \nbs^6 \tt  \bigr) - \tfrac{1- \alpha }{2}  \nbs^6 \Zbt \bigr) \nn \cdot \nbs^6\tt \tfrac{1}{\Sigma^{2\beta}} \bigl(\Jg\Wbn - \tfrac{13}{\eps} \Jg \bigr) \Q  \mathcal{J}^{\frac 12} 
\right| 
\les  \mathsf{K}^2 \eps^3 \brak{\mathsf{B}_6}^2 
.
\label{eq:hate:11-very-old}
\end{equation}
Note that we also have the estimate
\begin{equation}
\left|
\tfrac{3}{(1+\alpha)\eps} \tint  \tfrac{1- \alpha }{2}  \nbs^6 \Zbt  \nn \cdot \nbs^6\tt \tfrac{1}{\Sigma^{2\beta}} \bigl(\Jg\Wbn - \tfrac{13}{\eps} \Jg \bigr)  \Q  \mathcal{J}^{\frac 12} 
\right| 
\leq   
\tfrac{25}{\eps^2}  
\snorm{\tfrac{\mathcal{J}^{-\frac 14} \Q }{\Sigma^\beta} \nbs^6 \tt }_{L^2_{x,\s}} 
\snorm{\tfrac{\mathcal{J}^{\frac 34} }{\Sigma^\beta} \nbs^6 \Zbt }_{L^2_{x,\s}} 
\,.
\label{eq:hate:12-very-old}
\end{equation}
From \eqref{eq:hate:11-very-old} and \eqref{eq:hate:12-very-old}, we deduce that 
\begin{align}
&\left|
\tilde{\mathsf{M}}_2'
+ \tfrac{3}{2(1+\alpha)\eps} \tint (\Q\p_\s+V\p_2) ( \nn \cdot  \nbs^6 \tt)^2    \tfrac{1}{\Sigma^{2\beta}} \bigl(\Jg\Wbn - \tfrac{13}{\eps} \Jg \bigr)  \Q  \mathcal{J}^{\frac 12} 
\right| 
\notag\\
&\qquad 
\leq
\tfrac{25}{\eps^2}  
\snorm{\tfrac{\mathcal{J}^{-\frac 14} \Q }{\Sigma^\beta} \nbs^6 \tt }_{L^2_{x,\s}} 
\snorm{\tfrac{\mathcal{J}^{\frac 34} }{\Sigma^\beta} \nbs^6 \Zbt }_{L^2_{x,\s}} 
+ 
\Cn  \mathsf{K}^2 \eps^3 \brak{\mathsf{B}_6}^2 
\,.
\label{eq:hate:13}
\end{align}
Thus, we are left to consider the following term (which we rewrite using \eqref{adjoint-3})
\begin{align*}
\tilde{\mathsf{M}}_2'''
&:=-  \tfrac{3}{2(1+\alpha)\eps} \tint (\Q\p_\s+V\p_2) ( \nn \cdot  \nbs^6 \tt)^2    \tfrac{1}{\Sigma^{2\beta}} \bigl(\Jg\Wbn - \tfrac{13}{\eps} \Jg \bigr)  \Q  \mathcal{J}^{\frac 12}  
\notag\\
&=  \tfrac{3}{2(1+\alpha)\eps} \tint  ( \nn \cdot  \nbs^6 \tt)^2  \tfrac{1}{\Sigma^{2\beta}} \bigl(\Jg\Wbn - \tfrac{13}{\eps} \Jg \bigr)  \Q  (\Q\p_\s+V\p_2) \mathcal{J}^{\frac 12}  
\notag\\&\qquad 
+ \tfrac{3}{2(1+\alpha)\eps} \tint  ( \nn \cdot  \nbs^6 \tt)^2   \tfrac{1}{\Sigma^{2\beta}} \Q  \mathcal{J}^{\frac 12}  (\Q\p_\s+V\p_2)  \bigl(\Jg\Wbn - \tfrac{13}{\eps}\Jg \bigr)    
\notag\\&\qquad 
+ \tfrac{3}{2(1+\alpha)\eps} \tint  ( \nn \cdot  \nbs^6 \tt)^2 \bigl(\Jg\Wbn - \tfrac{13}{\eps} \Jg \bigr)  \mathcal{J}^{\frac 12} (\Q\p_\s+V\p_2) \Bigl(   \tfrac{1}{\Sigma^{2\beta}}  \Q   \Bigr)
\notag\\&\qquad 
-  \tfrac{3}{2(1+\alpha)\eps} \tint (V \Qr_2 - \nbs_2 V - \Qr_\s) ( \nn \cdot  \nbs^6 \tt)^2    \tfrac{1}{\Sigma^{2\beta}} \bigl(\Jg\Wbn - \tfrac{13}{\eps} \Jg \bigr)  \Q  \mathcal{J}^{\frac 12}  
\notag\\&\qquad 
-  \tfrac{3}{2(1+\alpha)\eps} \int \Q ( \nn \cdot  \nbs^6 \tt)^2    \tfrac{1}{\Sigma^{2\beta}} \bigl(\Jg\Wbn - \tfrac{13}{\eps} \Jg \bigr)  \Q  \mathcal{J}^{\frac 12}  \Bigr|_{\s}
\notag\\
&\qquad 
+  \tfrac{3}{2(1+\alpha)\eps} \tint  \Q ( \nn \cdot  \nbs^6 \tt)^2    \tfrac{1}{\Sigma^{2\beta}} \bigl(\Jg\Wbn - \tfrac{13}{\eps} \Jg \bigr)  \Q  \mathcal{J}^{\frac 12}  \Bigr|_{0}
\notag\\
&
=: \tilde{\mathsf{M}}_{2,a}''' + \tilde{\mathsf{M}}_{2,b}''' + \tilde{\mathsf{M}}_{2,c}''' + \tilde{\mathsf{M}}_{2,d}''' + \tilde{\mathsf{M}}_{2,e}''' + \tilde{\mathsf{M}}_{2,f}'''
\,.
\end{align*}
At this stage, we note that since $(\Q \p_\s + V\p_2 ) \mathcal{J} = - \tfrac{\Q}{\eps}$, and by also appealing to \eqref{Jg-evo-s} we have
\begin{subequations}
\label{eq:really:bad:2''':rewrite}
\begin{align}
 \tilde{\mathsf{M}}_{2,a}'''
 &=\tfrac{3}{4(1+\alpha)\eps^2} \tint  ( \nn \cdot  \nbs^6 \tt)^2  \tfrac{1}{\Sigma^{2\beta}} \bigl(- \Jg\Wbn + \tfrac{13}{\eps} \Jg \bigr)  \Q^2  \mathcal{J}^{-\frac 12}
 \\
 \tilde{\mathsf{M}}_{2,b}'''
 &= \tfrac{3}{2(1+\alpha)\eps} \tint  ( \nn \cdot  \nbs^6 \tt)^2   \tfrac{1}{\Sigma^{2\beta}} \Q  \mathcal{J}^{\frac 12}  (\Q\p_\s+V\p_2)   (\Jg\Wbn)
\notag\\
&\qquad 
-  \tfrac{39}{ 2(1+\alpha)\eps^2} \tint  ( \nn \cdot  \nbs^6 \tt)^2   \tfrac{1}{\Sigma^{2\beta}} \Q  \mathcal{J}^{\frac 12}  \bigl(\tfrac{1+\alpha}{2} \Jg \Wbn + \tfrac{1-\alpha}{2} \Jg \Zbn \bigr)  
 \,.     
\end{align}
\end{subequations}
Appealing to \eqref{eq:why:the:fuck:not:0},  \eqref{Sigma0i-ALE-s}, \eqref{eq:broncos:eat:shit:20}, the identity $(\Q \p_\s + V \p_2) \Q =  \Q \Qc  + V \p_2 \Q$, the bootstrap bounds~\eqref{bootstraps}, the initial data bounds~\eqref{table:derivatives}, the bounds \eqref{eq:Q:all:bbq} for the $\Q$-related coefficients,   the bounds \eqref{geometry-bounds-new} on the geometry, and the improved estimate for $(\Q\p_\s + V \p_2)\Jg \Wbn$ in~\eqref{eq:Jg:Wbn:improve:material:a}, we have
\begin{subequations}
\label{eq:how:the:fuck?}
\begin{align}
 \tilde{\mathsf{M}}_{2,b}'''
 &\geq -  \Cn \eps^2 (\tfrac{4}{\kappa_0})^{2\beta} \mathsf{K}^2 \brak{\mathsf{B}_6}^2 
 -  \tfrac{25}{(1+\alpha)\eps^3} \int_0^\s \snorm{ \tfrac{\Q \mathcal{J}^{\frac 14}}{\Sigma^\beta}\nn \cdot  \nbs^6 \tt (\cdot,\s')}_{L^2_{x}}^2 {\rm d} \s'
 \label{eq:how:the:fuck?:b}
 \,, \\ 
 \sabs{\tilde{\mathsf{M}}_{2,c}'''}  +  \sabs{\tilde{\mathsf{M}}_{2,d}'''}
 &\leq \Cn \eps^2 \beta (\tfrac{4}{\kappa_0})^{2\beta} \mathsf{K}^2 \brak{\mathsf{B}_6}^2
 +  \tfrac{100\cdot 250^2}{(1+\alpha)\eps^3} \int_0^\s \snorm{ \tfrac{\Q \mathcal{J}^{\frac 14}}{\Sigma^\beta}\nn \cdot  \nbs^6 \tt (\cdot,\s')}_{L^2_{x}}^2 {\rm d} \s'
 \label{eq:how:the:fuck?:cd}
 \,, \\
 \sabs{\tilde{\mathsf{M}}_{2,f}'''}
 &\leq  25 (1+\alpha)\eps (\tfrac{3}{\kappa_0})^{2\beta} \Cdatatwo 
 \label{eq:how:the:fuck?:f}
 \,.
\end{align}
In the last inequality we have used that $\|\Q(\cdot,0)\|_{L^\infty} \leq 1+\alpha$.
On the other hand, \eqref{eq:signed:Jg} yields
\begin{align}
 \tilde{\mathsf{M}}_{2,a}'''
 &\geq \tfrac{27}{40(1+\alpha)\eps^3} 
 \int_0^{\s} \snorm{ \tfrac{\Q \mathcal{J}^{-\frac 14}}{\Sigma^\beta}\nn \cdot  \nbs^6 \tt (\cdot,\s')}_{L^2_{x}}^2 {\rm d}\s'
 \label{eq:how:the:fuck?:a}
 \,,
 \\
 \tilde{\mathsf{M}}_{2,e}'''
 &\geq
  \tfrac{27}{20(1+\alpha)\eps^2} \snorm{ \tfrac{\Q \mathcal{J}^{\frac 14}}{\Sigma^\beta}\nn \cdot  \nbs^6 \tt (\cdot,\s)}_{L^2_{x}}^2
  \label{eq:how:the:fuck?:e}
  \,.
\end{align}
\end{subequations}
We note here that the second term on the right side of \eqref{eq:how:the:fuck?:b} is the correct ``Gr\"onwall term'' that corresponds to the energy given by \eqref{eq:how:the:fuck?:e}.
Combining the estimates in~\eqref{eq:really:bad:2}--\eqref{eq:how:the:fuck?}, we deduce that  
\begin{align}
 \tilde{\mathsf{M}}_2 
 &\geq  
 - \tfrac{34^2}{\eps} \int_0^{\s}  \snorm{\tfrac{\mathcal{J}^{\frac 34} (\Jg \Q)^{\frac 12}}{\Sigma^\beta} \nbs^6 \Wbt(\cdot,\s')}_{L^2_x}^2  {\rm d} \s'
- \tfrac{25^2(1+\alpha)}{\eps}  
\int_0^\s  \snorm{\tfrac{\mathcal{J}^{\frac 34} }{\Sigma^\beta} \nbs^6 \Zbt (\cdot,\s')}_{L^2_{x}}^2 {\rm d} \s' 
\notag\\
&\qquad
 - 25 (1+\alpha)\eps (\tfrac{3}{\kappa_0})^{2\beta} \Cdatatwo 
-  \Cn \eps^2 (\tfrac{4}{\kappa_0})^{2\beta} \mathsf{K}^2 \brak{\mathsf{B}_6}^2 
 \notag\\
 &\qquad 
 + \tfrac{7}{40(1+\alpha)\eps^3} \int_0^{\s} \snorm{ \tfrac{\Q \mathcal{J}^{-\frac 14}}{\Sigma^\beta}\nn \cdot  \nbs^6 \tt (\cdot,\s')}_{L^2_{x}}^2 {\rm d} \s'
 \notag\\
 &\qquad
 + \tfrac{27}{20(1+\alpha)\eps^2} \snorm{ \tfrac{\Q \mathcal{J}^{\frac 14}}{\Sigma^\beta}\nn \cdot  \nbs^6 \tt (\cdot,\s)}_{L^2_{x}}^2
  -  \tfrac{25 + 100\cdot 250^2}{(1+\alpha)\eps^3} \int_0^\s \snorm{ \tfrac{\Q \mathcal{J}^{\frac 14}}{\Sigma^\beta}\nn \cdot  \nbs^6 \tt (\cdot,\s')}_{L^2_{x}}^2 {\rm d} \s' \,.
 \label{eq:Really:Bad:2:final}
\end{align}
The first and last term in \eqref{eq:Really:Bad:2:final} are Gr\"onwall terms, while the second term is a damping term, which will be handled by taking $\beta$ to be sufficiently large in terms of $\alpha$.

Next, we return to the $\tilde{\mathsf{M}}_3$ term defined in \eqref{eq:Really:Bad:3:def}. We estimate this term simply using \eqref{eq:Q:bbq}, the Cauchy-Schwartz inequality, and the bound $\mathcal{J}\leq \Jg$, as to obtain
\begin{align}
 \tilde{\mathsf{M}}_3
 &\geq - \tfrac{\sqrt{5}}{\eps\sqrt{2(1+\alpha)}}  \snorm{\tfrac{\mathcal{J}^{\frac 34} (\Jg \Q)^{\frac 12}}{\Sigma^\beta} \nbs^6 \Wbt(\cdot,\s)}_{L^2_x} \snorm{ \tfrac{\Q \mathcal{J}^{\frac 14}}{\Sigma^\beta}\nn \cdot  \nbs^6 \tt (\cdot,\s)}_{L^2_{x}} 
 \notag\\
 &\geq
  - \tfrac{25}{52} \snorm{\tfrac{\mathcal{J}^{\frac 34} (\Jg \Q)^{\frac 12}}{\Sigma^\beta} \nbs^6 \Wbt(\cdot,\s)}_{L^2_x}^2
   - \tfrac{13}{10(1+\alpha)\eps^2} \snorm{ \tfrac{\Q \mathcal{J}^{\frac 14}}{\Sigma^\beta}\nn \cdot  \nbs^6 \tt (\cdot,\s)}_{L^2_{x}}^2
    \label{eq:Really:Bad:3:final}
   \,.
\end{align} 
We emphasize at this stage that $\frac{25}{52}  < \frac 12$, so that the above bound is compatible with \eqref{eq:I:tau:1246} and \eqref{eq:Really:Bad:2:final}. Combining \eqref{eq:WTF:*}, \eqref{eq:WTF:**}, \eqref{eq:Really:Bad:2:def}, \eqref{eq:Really:Bad:3:def}, \eqref{eq:Really:Bad:2:final}, and \eqref{eq:Really:Bad:3:final}, we obtain
\begin{align}
 \mathsf{M}_2  + \mathsf{M}_3
 &\geq  
   - (\tfrac{25}{52} + \Cn \eps) \snorm{\tfrac{\mathcal{J}^{\frac 34} (\Jg \Q)^{\frac 12}}{\Sigma^\beta} \nbs^6 \Wbt(\cdot,\s)}_{L^2_x}^2
 - \Cn \int_0^{\s}  \snorm{\tfrac{\mathcal{J}^{\frac 14} \Jg^{\frac 12}}{\Sigma^\beta} \nbs^6 \Wbt(\cdot,\s')}_{L^2_x}^2  {\rm d} \s'
 \notag\\
 &\qquad 
  - \tfrac{34^2}{\eps} \int_0^{\s}  \snorm{\tfrac{\mathcal{J}^{\frac 34} (\Jg \Q)^{\frac 12}}{\Sigma^\beta} \nbs^6 \Wbt(\cdot,\s')}_{L^2_x}^2  {\rm d} \s'
- \tfrac{25^2(1+\alpha)}{\eps}  
\int_0^\s  \snorm{\tfrac{\mathcal{J}^{\frac 34} }{\Sigma^\beta} \nbs^6 \Zbt (\cdot,\s')}_{L^2_{x}}^2 {\rm d} \s' 
\notag\\
&\qquad
 - 25 (1+\alpha)\eps (\tfrac{3}{\kappa_0})^{2\beta} \Cdatatwo 
-  \Cn \eps^2 (\tfrac{4}{\kappa_0})^{2\beta} \mathsf{K}^2 \brak{\mathsf{B}_6}^2 
 \notag\\
 &\qquad 
 + \tfrac{7}{40(1+\alpha)\eps^3} \int_0^{\s} \snorm{ \tfrac{\Q \mathcal{J}^{-\frac 14}}{\Sigma^\beta}\nn \cdot  \nbs^6 \tt (\cdot,\s')}_{L^2_{x}}^2 {\rm d} \s'
 \notag\\
 &\qquad
 + \tfrac{1}{20(1+\alpha)\eps^2} \snorm{ \tfrac{\Q \mathcal{J}^{\frac 14}}{\Sigma^\beta}\nn \cdot  \nbs^6 \tt (\cdot,\s)}_{L^2_{x}}^2
  -  \tfrac{25   + 100\cdot 250^2}{(1+\alpha)\eps^3} \int_0^\s \snorm{ \tfrac{\Q \mathcal{J}^{\frac 14}}{\Sigma^\beta}\nn \cdot  \nbs^6 \tt (\cdot,\s')}_{L^2_{x}}^2 {\rm d} \s' \,.
 \label{eq:Bad:2:Bad:3:final}
\end{align} 
In view of \eqref{eq:hate:7b:bound:initial}, it is left to obtain a good lower bound for the term $\mathsf{M}_1$ defined in \eqref{eq:hate:7b:rewrite}.
We have that 
\begin{align}
\mathsf{M}_1
&\geq 
 \tint  \tfrac{1}{\Sigma^{2\beta}} \mathsf{G_{bad}} \bigl(\nbs^6\Wbt\bigr)^2  
 - \tfrac{1+\alpha }{\eps} \int_0^\s   \snorm{\tfrac{\mathcal{J}^{\frac 14} \Jg^{\frac 12}}{\Sigma^\beta} \nbs^6 \Wbt(\cdot,\s')}_{L^2_x}  \snorm{\tfrac{\mathcal{J}^{\frac 34}  }{\Sigma^\beta} \nbs^6 \Zbt(\cdot,\s')}_{L^2_x} {\rm d} \s'
\notag\\
 &\qquad - \Cn \int_0^{\s} 
\Bigl(
 \snorm{\tfrac{\mathcal{J}^{\frac 14} \Jg^{\frac 12}}{\Sigma^\beta} \nbs^6 \Wbt(\cdot,\s')}_{L^2_x}^2 {\rm d} \s'
+
 \snorm{\tfrac{\mathcal{J}^{\frac 14} \Jg^{\frac 12}}{\Sigma^\beta} \nbs^6 \Zbt(\cdot,\s')}_{L^2_x}^2 {\rm d} \s'
 \Bigr) {\rm d} \s'
 \label{eq:Bad:1:lower:bound}
\end{align}
where we have introduced
\begin{equation}
\label{eq:G:bad:tangential}
\mathsf{G_{bad}}
= \mathcal{J}^{\frac 32}  (\Q \p_\s + V \p_2) \Jg
\, .
\end{equation}
We emphasize that the first term on the right side of \eqref{eq:Bad:1:lower:bound} is to be combined in \eqref{eq:hate:9} with the third term on the right side of \eqref{eq:I:tau:1246}, which contains the damping factor $\mathsf{G_{good}}$ defined in \eqref{eq:G0:tau:lower}. 
In particular, we emphasize that \eqref{eq:fakeJg:LB} gives
\begin{align}
 \mathsf{G_{good}} + \mathsf{G_{bad}}
 &= - \tfrac 12 (\Q \p_s + V\p_2) \bigr( \mathcal{J}^{\frac 32} \Jg \bigl) + \mathcal{J}^{\frac 32}  (\Q \p_\s + V \p_2) \Jg
 \notag\\
 & = \tfrac 12 \Bigl(-  \Jg (\Q \p_s + V\p_2)  \mathcal{J}^{\frac 32}   + \mathcal{J}^{\frac 32}  (\Q \p_\s + V \p_2) \Jg\Bigr)
 \geq \tfrac{1+\alpha}{16\eps} \mathcal{J}^{\frac 12} \Jg
 \,.
 \label{eq:G:good:G:bad}
\end{align}

To summarize the bounds in this subsection, we combine \eqref{eq:hate:7}, \eqref{eq:hate:7a:bound}, \eqref{eq:hate:7b:bound:initial}, \eqref{eq:Bad:2:Bad:3:final}, and \eqref{eq:Bad:1:lower:bound}, and deduce that  
\begin{align}
I^{\AA_\tau}_{3,a}
&\geq  I^{\AA_\tau}_{3,a,iii}
- (\tfrac{25}{52} + \Cn \eps) \snorm{\tfrac{\mathcal{J}^{\frac 34} (\Jg \Q)^{\frac 12}}{\Sigma^\beta} \nbs^6 \Wbt(\cdot,\s)}_{L^2_x}^2
+ \tint  \tfrac{1}{\Sigma^{2\beta}} \mathsf{G_{bad}} \bigl(\nbs^6\Wbt\bigr)^2  
 \notag\\
 &\qquad 
 - \tfrac{1+\alpha }{\eps} \int_0^\s   \snorm{\tfrac{\mathcal{J}^{\frac 14} \Jg^{\frac 12}}{\Sigma^\beta} \nbs^6 \Wbt(\cdot,\s')}_{L^2_x}  \snorm{\tfrac{\mathcal{J}^{\frac 34}  }{\Sigma^\beta} \nbs^6 \Zbt(\cdot,\s')}_{L^2_x} {\rm d} \s'
\notag\\
 &\qquad - \Cn \int_0^{\s} 
\Bigl(
 \snorm{\tfrac{\mathcal{J}^{\frac 14} \Jg^{\frac 12}}{\Sigma^\beta} \nbs^6 \Wbt(\cdot,\s')}_{L^2_x}^2 {\rm d} \s'
+
 \snorm{\tfrac{\mathcal{J}^{\frac 14} \Jg^{\frac 12}}{\Sigma^\beta} \nbs^6 \Zbt(\cdot,\s')}_{L^2_x}^2 {\rm d} \s'
 \Bigr) {\rm d} \s'
 \notag\\
 &\qquad 
  - \tfrac{34^2+500^2}{\eps} \int_0^{\s}  \snorm{\tfrac{\mathcal{J}^{\frac 34} (\Jg \Q)^{\frac 12}}{\Sigma^\beta} \nbs^6 \Wbt(\cdot,\s')}_{L^2_x}^2  {\rm d} \s'
- \tfrac{25^2(1+\alpha)}{\eps}  
\int_0^\s  \snorm{\tfrac{\mathcal{J}^{\frac 34} }{\Sigma^\beta} \nbs^6 \Zbt (\cdot,\s')}_{L^2_{x}}^2 {\rm d} \s' 
\notag\\
&\qquad
 - 31 (1+\alpha)\eps (\tfrac{3}{\kappa_0})^{2\beta} \Cdatatwo 
-  \Cn \eps^2 (\tfrac{4}{\kappa_0})^{2\beta} \mathsf{K}^2 \brak{\mathsf{B}_6}^2 
 \notag\\
 &\qquad 
 + \tfrac{7}{40(1+\alpha)\eps^3} \int_0^{\s} \snorm{ \tfrac{\Q \mathcal{J}^{-\frac 14}}{\Sigma^\beta}\nn \cdot  \nbs^6 \tt (\cdot,\s')}_{L^2_{x}}^2 {\rm d} \s'
 \notag\\
 &\qquad
 + \tfrac{1}{20(1+\alpha)\eps^2} \snorm{ \tfrac{\Q \mathcal{J}^{\frac 14}}{\Sigma^\beta}\nn \cdot  \nbs^6 \tt (\cdot,\s)}_{L^2_{x}}^2
  -  \tfrac{25+ 100\cdot 250^2+500^2}{(1+\alpha)\eps^3} \int_0^\s \snorm{ \tfrac{\Q \mathcal{J}^{\frac 14}}{\Sigma^\beta}\nn \cdot  \nbs^6 \tt (\cdot,\s')}_{L^2_{x}}^2 {\rm d} \s' \,.
  \label{eq:hate:7b:bound}
\end{align}
The forcing and commutator terms appearing in $I^{\AA_\tau}_{3,a,iii}$ will be shown in~Subsection~\ref{sec:D6:tau:forcing:commutator} to satisfy similar bounds. 

\subsubsection{Combining all the terms with over-differentiated geometry}
Here we summarize the identities for the terms which contain over-differentiated geometry. 
By combining \eqref{eq:I:Zt:5+7:final}, \eqref{eq:hate:5}, and \eqref{eq:hate:7b:bound}, we arrive at the identity
\begin{align}
&I^{\WW_\tau}_3 +  I^{\ZZ_\tau}_3 +  I^{\AA_\tau}_3 + I^{\AA_\tau}_5+I^{\AA_\tau}_7  +  I^{\ZZ_\tau}_{5} +  I^{\ZZ_\tau}_{7} 
\notag\\
&\geq  
- \sabs{I^{\ZZ_\tau}_{5+7,a,iii}}  - \sabs{I^{\AA_\tau}_{5+7,a,v} } -  \sabs{ I^{\AA_\tau}_{3,a,iii}}
 - 31 (1+\alpha)\eps (\tfrac{3}{\kappa_0})^{2\beta} \Cdatatwo 
-  \Cn \eps^2 \beta^2 (\tfrac{4}{\kappa_0})^{2\beta} \mathsf{K}^2  \brak{\mathsf{B}_6}^2
\notag\\
&\qquad 
- (\tfrac{25}{52} + \Cn \eps) 
\snorm{\tfrac{\mathcal{J}^{\frac 34} (\Jg \Q)^{\frac 12}}{\Sigma^\beta} \nbs^6 \Wbt(\cdot,\s)}_{L^2_x}^2 
- \Cn \eps 
\Bigl( 
\snorm{\tfrac{\mathcal{J}^{\frac 34} (\Jg \Q)^{\frac 12}}{\Sigma^\beta} \nbs^6 \Zbt(\cdot,\s)}_{L^2_x}^2
+ 
\snorm{\tfrac{\mathcal{J}^{\frac 34} (\Jg \Q)^{\frac 12}}{\Sigma^\beta} \nbs^6 \Abt(\cdot,\s)}_{L^2_x}^2 
\Bigr)
\notag\\
& \qquad
- \Cn \eps  \Bigl(\snorm{\tfrac{\mathcal{J}^{\frac 34} (\Jg \Q)^{\frac 12}}{\Sigma^\beta} \nbs^6 \Wbt(\cdot,0)}_{L^2_x}^2 + \snorm{\tfrac{\mathcal{J}^{\frac 34} (\Jg \Q)^{\frac 12}}{\Sigma^\beta} \nbs^6 \Zbt(\cdot,0)}_{L^2_x}^2 + \snorm{\tfrac{\mathcal{J}^{\frac 34} (\Jg \Q)^{\frac 12}}{\Sigma^\beta} \nbs^6 \Abt(\cdot,0)}_{L^2_x}^2  
\Bigr) 
\notag\\
&\qquad 
- \Cn \int_0^{\s} 
\Bigl(
 \snorm{\tfrac{\mathcal{J}^{\frac 14} \Jg^{\frac 12}}{\Sigma^\beta} \nbs^6 \Wbt(\cdot,\s')}_{L^2_x}^2 {\rm d} \s'
+
 \snorm{\tfrac{\mathcal{J}^{\frac 14} \Jg^{\frac 12}}{\Sigma^\beta} \nbs^6 \Zbt(\cdot,\s')}_{L^2_x}^2 {\rm d} \s'
+
 \snorm{\tfrac{\mathcal{J}^{\frac 14} \Jg^{\frac 12}}{\Sigma^\beta} \nbs^6 \Abt(\cdot,\s')}_{L^2_x}^2 
 \Bigr) {\rm d} \s'
\notag\\
&\qquad
+ \tint  \tfrac{1}{\Sigma^{2\beta}} \mathsf{G_{bad}} \bigl(\nbs^6\Wbt\bigr)^2  
 - \tfrac{1+\alpha }{\eps} \int_0^\s   \snorm{\tfrac{\mathcal{J}^{\frac 14} \Jg^{\frac 12}}{\Sigma^\beta} \nbs^6 \Wbt(\cdot,\s')}_{L^2_x}  \snorm{\tfrac{\mathcal{J}^{\frac 34}  }{\Sigma^\beta} \nbs^6 \Zbt(\cdot,\s')}_{L^2_x} {\rm d} \s'
\notag\\
 &\qquad 
  - \tfrac{34^2+500^2}{\eps} \int_0^{\s}  \snorm{\tfrac{\mathcal{J}^{\frac 34} (\Jg \Q)^{\frac 12}}{\Sigma^\beta} \nbs^6 \Wbt(\cdot,\s')}_{L^2_x}^2  {\rm d} \s'
- \tfrac{25^2(1+\alpha)}{\eps}  
\int_0^\s  \snorm{\tfrac{\mathcal{J}^{\frac 34} }{\Sigma^\beta} \nbs^6 \Zbt (\cdot,\s')}_{L^2_{x}}^2 {\rm d} \s' 
\notag\\
&\qquad 
 + \tfrac{7}{40(1+\alpha)\eps^3} \int_0^{\s} \snorm{ \tfrac{\Q \mathcal{J}^{-\frac 14}}{\Sigma^\beta}\nn \cdot  \nbs^6 \tt (\cdot,\s')}_{L^2_{x}}^2 {\rm d} \s'
 \notag\\
 &\qquad
 + \tfrac{1}{20(1+\alpha)\eps^2} \snorm{ \tfrac{\Q \mathcal{J}^{\frac 14}}{\Sigma^\beta}\nn \cdot  \nbs^6 \tt (\cdot,\s)}_{L^2_{x}}^2
  -  \tfrac{25+ 100\cdot 250^2+500^2}{(1+\alpha)\eps^3} \int_0^\s \snorm{ \tfrac{\Q \mathcal{J}^{\frac 14}}{\Sigma^\beta}\nn \cdot  \nbs^6 \tt (\cdot,\s')}_{L^2_{x}}^2 {\rm d} \s' 
\,.
\label{eq:hate:8}
\end{align}

\subsection{The forcing and commutator terms}
\label{sec:D6:tau:forcing:commutator}
The only terms left to bound are: $I^{\WW_\tau}_4$ (cf.~\eqref{I4-Wbt}), $I^{\ZZ_\tau}_8$ (cf.~\eqref{I8-Zbt}), $I^{\AA_\tau}_8$ (cf.~\eqref{I8-Abt}), $I^{\ZZ_\tau}_{5+7,a,iii}$ (cf.~\eqref{eq:hate:2c}), $I^{\AA_\tau}_{5+7,a,v}$ (cf.~\eqref{eq:hate:4e}), and $I^{\AA_\tau}_{3,a,iii}$ (cf.~\eqref{eq:hate:7c}). 
These terms do not contain factors with derivative loss, or with a deficient power of $\Jg$ or $\mathcal{J}$; as such, these terms are bounded directly using Cauchy-Schwartz and the available estimates for the fundamental unknowns and for the geometry.

For instance, from \eqref{I4-Wbt}, \eqref{bs-Sigma}, and \eqref{eq:Q:bbq} we have that 
\begin{subequations}
\label{eq:I:W:tt:4:all}
\begin{align}
\sabs{I^{\WW_\tau}_4}
&\leq\tfrac{5\kappa_0}{2(1+\alpha)} \int_0^\s  \snorm{\tfrac{\mathcal{J}^{\frac 34} (\Jg \Q)^{\frac 12}}{\Sigma^\beta} \nbs^6 \Wbt(\cdot,\s')}_{L^2_x}  
\notag\\
&\qquad   \times \Bigl( \snorm{\tfrac{\mathcal{J}^{\frac 34} (\Jg \Q)^{\frac 12}}{\Sigma^\beta} \nbs^6\Fwt (\cdot,\s')}_{L^2_x}  + \snorm{\tfrac{\mathcal{J}^{\frac 34} (\Jg \Q)^{\frac 12}}{\Sigma^\beta} \mathcal{R}_\Wb^\tau (\cdot,\s')}_{L^2_x}  + \snorm{\tfrac{\mathcal{J}^{\frac 34} (\Jg \Q)^{\frac 12}}{\Sigma^\beta} \mathcal{C}_\Wb^\tau(\cdot,\s')}_{L^2_x}  \Bigr){\rm d} \s'
\,.
\label{eq:I:W:tt:4:a}
\end{align}
By the definition of $\Fwt$ in \eqref{forcing-nt}, and appealing to the $\mathcal{J}$ estimate~\eqref{eq:Jgb:less:than:1}, the bootstrap inequalities~\eqref{bootstraps}, the coefficient bounds~\eqref{eq:Q:all:bbq}, the bounds for the geometry~\eqref{geometry-bounds-new}, and the double-commutator bound in \eqref{eq:Lynch:2}, we deduce
\begin{equation}
\snorm{\tfrac{\mathcal{J}^{\frac 34} (\Jg \Q)^{\frac 12}}{\Sigma^\beta} \nbs^6\Fwt}_{L^2_{x,\s}}
\les  \snorm{\tfrac{\mathcal{J}^{\frac 34} (\Jg \Q)^{\frac 12}}{\Sigma^\beta} \nbs^6(\Wbt,\Zbt,\Abt)}_{L^2_{x,\s}}
+   \eps (\tfrac{4}{\kappa_0})^\beta  \mathsf{K} \brak{\mathsf{B}_6}
\les 
 \eps (\tfrac{4}{\kappa_0})^\beta  \mathsf{K} \brak{\mathsf{B}_6}
\,.
\label{eq:I:W:tt:4:b}
\end{equation}
Similarly, from the definition of $ \mathcal{R}_\Wb^\tau$ in \eqref{Cw-Rw-comm-tan}, and by additionally appealing to \eqref{eq:vort:H6} and \eqref{eq:Wb:tt}, which gives the estimate $\|(\Q \p_\s + V \p_2) \Wbt\|_{L^\infty_{x,\s}} \les \eps$, we deduce 
\begin{equation}
\snorm{\tfrac{\mathcal{J}^{\frac 34} (\Jg \Q)^{\frac 12}}{\Sigma^\beta}\mathcal{R}_\Wb^\tau}_{L^2_{x,\s}}
\les   \eps^2 (\tfrac{4}{\kappa_0})^\beta  \mathsf{K} \brak{\mathsf{B}_6}
+   \eps \Bigl(\snorm{\tfrac{\mathcal{J}^{\frac 34} (\Jg \Q)^{\frac 12}}{\Sigma^\beta} \nbs^6\Wbt}_{L^2_{x,\s}}
+ \snorm{\tfrac{\mathcal{J}^{\frac 34} (\Jg \Q)^{\frac 12}}{\Sigma^\beta} \nbs^6\Zbt}_{L^2_{x,\s}}\Bigr)
\les   \eps^2 (\tfrac{4}{\kappa_0})^\beta  \mathsf{K} \brak{\mathsf{B}_6}
\label{eq:I:W:tt:4:c}
\,.
\end{equation}
Lastly, using the definition of $\mathcal{C}_\Wb^\tau$ in \eqref{Cw-Rw-comm-tan}, identity~\eqref{eq:Wb:tt},   the aforementioned  bounds,  Lemma~\ref{lem:anisotropic:sobolev}, the Leibniz rule and Lemma~\ref{lem:time:interpolation},  we may also obtain 
\begin{align}
\snorm{\tfrac{\mathcal{J}^{\frac 34} (\Jg \Q)^{\frac 12}}{\Sigma^\beta}\mathcal{C}_\Wb^\tau}_{L^2_{x,\s}}
&\les  \eps^2 (\tfrac{4}{\kappa_0})^\beta  \mathsf{K}^2 \brak{\mathsf{B}_6}^2
+  \eps \snorm{\tfrac{\mathcal{J}^{\frac 34} (\Jg \Q)^{\frac 12}}{\Sigma^\beta} \nbs^6\Wbt}_{L^2_{x,\s}}
+ \snorm{\tfrac{\mathcal{J}^{\frac 34} (\Jg \Q)^{\frac 12}}{\Sigma^\beta} \nbs^6\Abt}_{L^2_{x,\s}}
\les   \eps  (\tfrac{4}{\kappa_0})^\beta  \mathsf{K} \brak{\mathsf{B}_6}
\,.
\label{eq:I:W:tt:4:d}
\end{align}
\end{subequations}
Combining the bounds in \eqref{eq:I:W:tt:4:all}, we thus deduce
\begin{align}
\label{eq:I:W:tt:4:final}
\sabs{I^{\WW_\tau}_4} 
&\leq
 \Cn \eps^2 (\tfrac{4}{\kappa_0})^{2\beta}  \mathsf{K}^2 \brak{\mathsf{B}_6}^2
\,.
\end{align}

Next, we consider the bounds for $I^{\ZZ_\tau}_8$ and $I^{\AA_\tau}_8$, which are nearly identical. From \eqref{bs-Sigma}, \eqref{I8-Zbt}, and \eqref{I8-Abt}, we obtain
\begin{subequations}
\label{eq:I:ZA:tt:8:all}
\begin{align}
&\sabs{I^{\ZZ_\tau}_8}
+
\sabs{I^{\AA_\tau}_8}
\notag\\
&\leq \int_0^\s \snorm{\tfrac{\mathcal{J}^{\frac 34}}{\Sigma^\beta} \nbs^6 \Zbt(\cdot,\s')}_{L^2_x}
\Bigl(\kappa_0  \snorm{\tfrac{\mathcal{J}^{\frac 34}}{\Sigma^\beta} \nbs^6 \Fzt(\cdot,\s')}_{L^2_x}
+ \snorm{\tfrac{\mathcal{J}^{\frac 34}}{\Sigma^{\beta-1}} \mathcal{R}^\tau_{\Zb}(\cdot,\s')}_{L^2_x}
+ \snorm{\tfrac{\mathcal{J}^{\frac 34}}{\Sigma^{\beta-1}} \mathcal{C}^\tau_{\Zb}(\cdot,\s')}_{L^2_x}
\Bigr) {\rm d} \s'
\notag\\
&\qquad +  \int_0^\s \snorm{\tfrac{\mathcal{J}^{\frac 34}}{\Sigma^\beta} \nbs^6 \Abt(\cdot,\s')}_{L^2_x}
\Bigl(\kappa_0  \snorm{\tfrac{\mathcal{J}^{\frac 34}}{\Sigma^\beta} \nbs^6 \Fat(\cdot,\s')}_{L^2_x}
+ \snorm{\tfrac{\mathcal{J}^{\frac 34}}{\Sigma^{\beta-1}} \mathcal{R}^\tau_{\Ab}(\cdot,\s')}_{L^2_x}
 + \snorm{\tfrac{\mathcal{J}^{\frac 34}}{\Sigma^{\beta-1}} \mathcal{C}^\tau_{\Ab}(\cdot,\s')}_{L^2_x}
\Bigr) {\rm d} \s'
\,.
\label{eq:I:ZA:tt:8:a}
\end{align}
From the definitions of $\Fzt$ and $\Fat$ in \eqref{forcing-nt}, the bootstrap inequalities~\eqref{bootstraps},  the estimate $\mathcal{J}\leq \Jg$, the bounds for the geometry~\eqref{geometry-bounds-new}, the vorticity estimates~\eqref{eq:vort:H6}--\eqref{eq:vorticity:pointwise},  the double-commutator in \eqref{eq:Lynch:2}, \bubu{and the improved estimate for $\nbs^5 \Zbn$ in \eqref{eq:Jg:Zbn:D5:improve:a}}, we deduce
\begin{equation}
\snorm{\tfrac{\mathcal{J}^{\frac 34}}{\Sigma^\beta} \nbs^6 \Fzt}_{L^2_{x,\s}}
+ 
\snorm{\tfrac{\mathcal{J}^{\frac 34}}{\Sigma^\beta} \nbs^6 \Fat}_{L^2_{x,\s}}
\leq  \bubu{\tfrac{\alpha}{\kappa_0} (\tfrac{4}{\kappa_0})^\beta \mathsf{B}_6} + \bubu{\Cn} \eps(\tfrac{4}{\kappa_0})^\beta  \mathsf{K} \brak{\mathsf{B}_6}
\,.
\label{eq:I:ZA:tt:8:b}
\end{equation}
Similarly, the definitions of $\mathcal{R}^\tau_{\Zb}$ and $\mathcal{R}^\tau_{\Ab}$ in \eqref{Cz-Rz-comm-tan} and \eqref{Ca-Ra-comm-tan} furthermore allow us to estimate  
\begin{align}
\snorm{\tfrac{\mathcal{J}^{\frac 34}}{\Sigma^{\beta-1}} \mathcal{R}^\tau_{\Zb}}_{L^2_{x,\s}}
&\leq \Cn \eps  (\tfrac{4}{\kappa_0})^\beta  \mathsf{K} \brak{\mathsf{B}_6} 
+2\alpha \| \nn \cdot \tt,_1 - \Jg g^{-\frac 12} \nbs_2h \nn \cdot \nbs_2 \tt\|_{L^\infty_{x,\s}} 
\snorm{\tfrac{\mathcal{J}^{\frac 34} }{\Sigma^{\beta-1}} \nbs^6 \Zbn}_{L^2_{x,\s}}
\label{eq:I:ZA:tt:8:c:0} \\
\snorm{\tfrac{\mathcal{J}^{\frac 34}}{\Sigma^{\beta-1}} \mathcal{R}^\tau_{\Ab}}_{L^2_{x,\s}}
&\leq \Cn \eps  (\tfrac{4}{\kappa_0})^\beta  \mathsf{K} \brak{\mathsf{B}_6} 
\,.
\label{eq:I:ZA:tt:8:c:00}
\end{align}
The term highlighted on the second line of the right side of \eqref{eq:I:ZA:tt:8:c:0} is not necessarily small, and so it must be treated with care. First, according to \eqref{p1-n-tau} and \eqref{bootstraps} we have that 
\begin{equation*}
\| \nn \cdot \tt,_1 - \Jg g^{-\frac 12} \nbs_2h \nn \cdot \nbs_2 \tt\|_{L^\infty_{x,\s}}  
\leq \| g^{-\frac 12} \nbs_2 \Jg \|_{L^\infty_{x,\s}}  
 + \Cn \eps
 \leq 4(1+\alpha) + \Cn \eps \leq 5(1+\alpha)
 \,.
\end{equation*}
Then, we note that if $\nbs^6$ contains a single copy of $\nbs_1$, or a single copy of $\nbs_2$, then by \eqref{eq:Jg:Zbn:p1D5:improve}, \eqref{eq:Jg:Zbn:p2D5:improve:new}, and the inequality $\mathcal{J} \leq \Jg$, we have that 
\begin{align}
 \snorm{\tfrac{\mathcal{J}^{\frac 34} }{\Sigma^{\beta-1}} \nbs_1 \nbs^5 \Zbn}_{L^2_{x,\s}}
&\leq \tfrac{1}{2 \alpha }  (\tfrac{4}{\kappa_0})^{\beta} \snorm{ \mathcal{J}^{\frac 14} \Jgh \nbs^6 (\Jg \Zbn)}_{L^2_{x,\s}}
+ \Cn \eps (\tfrac{4}{\kappa_0})^\beta \mathsf{K} \brak{\mathsf{B}_6}
\label{eq:I:ZA:tt:8:c:00000i}
\\
 \snorm{\tfrac{\mathcal{J}^{\frac 34} }{\Sigma^{\beta-1}} \nbs_2 \nbs^5 \Zbn}_{L^2_{x,\s}}
&\leq  \tfrac{1}{\alpha \eps}  \snorm{ \tfrac{\mathcal{J}^{\frac 34} }{\Sigma^\beta} \nbs^6  \Zbt}_{L^2_{x,\s}}
+ \Cn \eps (\tfrac{4}{\kappa_0})^\beta \mathsf{K} \brak{\mathsf{B}_6}
\label{eq:I:ZA:tt:8:c:00000ii}
\,.
\end{align}
On the other hand, if $\nbs^6 = \nbs_\s^6$, then by estimate \eqref{eq:madman:2} below we obtain that
\begin{equation}
 2\alpha \snorm{\tfrac{\mathcal{J}^{\frac 34} }{\Sigma^{\beta-1}} \nbs_\s^6 \Zbn}_{L^2_{x,\s}}
\leq \Cn \eps (\tfrac{4}{\kappa_0})^\beta \mathsf{K} \brak{\mathsf{B}_6}
\,.
\label{eq:I:ZA:tt:8:c:00000}
\end{equation}
By combining \eqref{eq:I:ZA:tt:8:c:00000i}--\eqref{eq:I:ZA:tt:8:c:00000} with \eqref{eq:I:ZA:tt:8:c:0}, the bootstrap assumption~\eqref{bootstraps-Dnorm:6}, and definition \eqref{eq:tilde:D6}, we obtain
\begin{equation}
\snorm{\tfrac{\mathcal{J}^{\frac 34}}{\Sigma^{\beta-1}} \mathcal{R}^\tau_{\Zb}}_{L^2_{x,\s}}
\leq \Cn \eps  (\tfrac{4}{\kappa_0})^\beta  \mathsf{K} \brak{\mathsf{B}_6} 
+ \tfrac{5(1+\alpha)}{\alpha \eps} \snorm{ \tfrac{\mathcal{J}^{\frac 34} }{\Sigma^\beta} \nbs^6  \Zbt}_{L^2_{x,\s}}
+6(1+\alpha) (\tfrac{4}{\kappa_0})^\beta \mathsf{B}_6 
\,.
\label{eq:I:ZA:tt:8:c}
\end{equation}
We emphasize here that the last term on the right side of \eqref{eq:I:ZA:tt:8:c} does not contain a factor of $\mathsf{K}$. This fact is crucial for the proof: we can absorb this term not because it has a helpful power of $\eps$, but because it is missing a damaging factor of $\mathsf{K}$. 
Lastly, using the definitions of $\mathcal{C}_\Zb^\tt$ and $\mathcal{C}_\Ab^\tt$ in \eqref{Cz-Rz-comm-tan} and \eqref{Ca-Ra-comm-tan}, and by appealing to all available bounds, we obtain
\begin{align}
\snorm{\tfrac{\mathcal{J}^{\frac 34}}{\Sigma^{\beta-1}} \mathcal{C}_\Zb^\tt}_{L^2_{x,\s}}
& \leq \Cn  \eps  (\tfrac{4}{\kappa_0})^\beta  \mathsf{K} \brak{\mathsf{B}_6}  
+\tfrac{1}{\eps} \snorm{\tfrac{\mathcal{J}^{\frac 34}}{\Sigma^{\beta-1}}  \doublecom{\nbs^6, \tfrac{\Jg}{\Sigma} , \nbs_{\s} \Zbt} }_{L^2_{x,\s}}
\notag\\
& \leq  
\tfrac{32(1+\alpha)}{\eps}  \snorm{\tfrac{\mathcal{J}^{\frac 34}}{\Sigma^{\beta}}   \nbs^6 \Zbt  }_{L^2_{x,\s}}
+ \Cn  \eps^{\frac 14}  (\tfrac{4}{\kappa_0})^\beta  \mathsf{K} \brak{\mathsf{B}_6} 
\label{eq:I:ZA:tt:8:d} \\
\snorm{\tfrac{\mathcal{J}^{\frac 34}}{\Sigma^{\beta-1}} \mathcal{C}_\Ab^\tt}_{L^2_{x,\s}}
& \leq \Cn  \eps  (\tfrac{4}{\kappa_0})^\beta  \mathsf{K} \brak{\mathsf{B}_6} 
+\tfrac{1}{\eps} \snorm{\tfrac{\mathcal{J}^{\frac 34}}{\Sigma^{\beta-1}}  \doublecom{\nbs^6, \tfrac{\Jg}{\Sigma} \tt^k , \nbs_\s \Ak)}   }_{L^2_{x,\s}}
\notag\\
&   \leq 
 \tfrac{32(1+\alpha)}{\eps}  \snorm{\tfrac{\mathcal{J}^{\frac 34}}{\Sigma^{\beta}}   \nbs^6 \Abt  }_{L^2_{x,\s}}
+ \Cn  \eps^{\frac 14}  (\tfrac{4}{\kappa_0})^\beta  \mathsf{K} \brak{\mathsf{B}_6} 
\,.
\label{eq:I:ZA:tt:8:dd}
\end{align}
\end{subequations}
In the last inequality in both \eqref{eq:I:ZA:tt:8:d} and \eqref{eq:I:ZA:tt:8:dd} we have appealed to \eqref{eq:poop:stinks:4} with $m=6$, and note that the most dangerous term comes from $i=2$; this is the cause of the $\eps^{\frac 14}$ term, instead of the usual $\eps$. Upon combining all the bounds we have obtained in \eqref{eq:I:ZA:tt:8:all}, we deduce that 
\begin{align}
\sabs{I^{\ZZ_\tau}_8}
+
\sabs{I^{\AA_\tau}_8}
&\leq 
\tfrac{\bubu{35}(1+\alpha)^2}{\alpha \eps} 
\int_0^\s \snorm{\tfrac{\mathcal{J}^{\frac 34}}{\Sigma^\beta} \nbs^6 \Zbt(\cdot,\s')}_{L^2_x}^2 {\rm d}\s'
+ \tfrac{\bubu{35}(1+\alpha)}{\eps} 
\int_0^\s \snorm{\tfrac{\mathcal{J}^{\frac 34}}{\Sigma^\beta} \nbs^6 \Abt(\cdot,\s')}_{L^2_x}^2 {\rm d}\s'
\notag\\
&\qquad +
\bubu{10} \alpha \eps (\tfrac{4}{\kappa_0})^{2\beta} \mathsf{B}_6^2
 +   \Cn  \eps^{\frac 32}  (\tfrac{4}{\kappa_0})^{2\beta}  \mathsf{K}^2 \brak{\mathsf{B}_6}^2
\,.
\label{eq:I:ZA:tt:8:final}
\end{align}
 
At last, the terms $I^{\ZZ_\tau}_{5+7,a,iii}$ (defined in~\eqref{eq:hate:2c}), $I^{\AA_\tau}_{5+7,a,v}$ (defined in~\eqref{eq:hate:4e}), and $I^{\AA_\tau}_{3,a,iii}$ (defined in~\eqref{eq:hate:7c}), may be estimated using the forcing and commutator estimates that we have just obtained in \eqref{eq:I:W:tt:4:b}, \eqref{eq:I:W:tt:4:c}, and \eqref{eq:I:W:tt:4:d} for $\Wbt$, \eqref{eq:I:ZA:tt:8:b}, \eqref{eq:I:ZA:tt:8:c}, and \eqref{eq:I:ZA:tt:8:d} for $\Zbt$, and \eqref{eq:I:ZA:tt:8:b}, \eqref{eq:I:ZA:tt:8:c:00}, and \eqref{eq:I:ZA:tt:8:dd} for $\Abt$. We record the bounds
\begin{subequations}
\label{eq:money:dont:jiggle}
\begin{align}
\sabs{I^{\ZZ_\tau}_{5+7,a,iii}}
&\leq 
\Cn \int_0^{\s} 
\snorm{ \tfrac{\mathcal{J}^{\frac 14}}{\Sigma^\beta}\nn \cdot  \nbs^6 \tt (\cdot,\s')}_{L^2_{x}}
\snorm{ \tfrac{\mathcal{J}^{\frac 54}}{\Sigma^\beta} \big( \nbs^6\Fwt  + \mathcal{R}_\Wb^\tt + \mathcal{C}_\Wb^\tt + \nbs^6 \Fzt   + \mathcal{R}_\Zb^\tt + \mathcal{C}_\Zb^\tt \big) (\cdot,\s')}_{L^2_{x}} {\rm d}\s'
\notag\\
&\leq \tfrac{1}{(1+\alpha)\eps^3} \int_0^{\s} \snorm{ \tfrac{\Q \mathcal{J}^{\frac 14}}{\Sigma^\beta}\nn \cdot  \nbs^6 \tt (\cdot,\s')}_{L^2_{x}}^2 {\rm d} \s' 
+ \Cn \eps^3 (\tfrac{4}{\kappa_0})^{2\beta} \mathsf{K}^2 \brak{\mathsf{B}_6}^2
\\
\sabs{I^{\AA_\tau}_{5+7,a,v}}
&\leq \Cn \int_0^{\s} 
\snorm{ \tfrac{\mathcal{J}^{\frac 14}}{\Sigma^\beta}\nn \cdot  \nbs^6 \tt (\cdot,\s')}_{L^2_{x}}
\snorm{ \tfrac{\mathcal{J}^{\frac 54}}{\Sigma^\beta} \big( \nbs^6 \Fat   + \mathcal{R}_\Ab^\tt + \mathcal{C}_\Ab^\tt   \big) (\cdot,\s')}_{L^2_{x}} {\rm d}\s'
\notag\\
&\leq 
\tfrac{1}{(1+\alpha)\eps^3} \int_0^{\s} \snorm{ \tfrac{\Q \mathcal{J}^{\frac 14}}{\Sigma^\beta}\nn \cdot  \nbs^6 \tt (\cdot,\s')}_{L^2_{x}}^2 {\rm d} \s' 
+ \Cn \eps^3 (\tfrac{4}{\kappa_0})^{2\beta} \mathsf{K}^2 \brak{\mathsf{B}_6}^2
\\
\sabs{I^{\AA_\tau}_{3,a,iii}}
&\leq \tfrac{8}{\eps \kappa_0}
\int_0^{\s} 
\snorm{ \tfrac{\mathcal{J}^{\frac 14}}{\Sigma^\beta}\nn \cdot  \nbs^6 \tt (\cdot,\s')}_{L^2_{x}}
\snorm{ \tfrac{\mathcal{J}^{\frac 54}}{\Sigma^\beta} \big(\nbs^6\Fwt  + \mathcal{R}_\Wb^\tt + \mathcal{C}_\Wb^\tt \big) (\cdot,\s')}_{L^2_{x}} {\rm d}\s'
\notag\\
&\leq 
\tfrac{1}{(1+\alpha)\eps^3} \int_0^{\s} \snorm{ \tfrac{\Q \mathcal{J}^{\frac 14}}{\Sigma^\beta}\nn \cdot  \nbs^6 \tt (\cdot,\s')}_{L^2_{x}}^2 {\rm d} \s' 
+ \Cn \eps^3 (\tfrac{4}{\kappa_0})^{2\beta} \mathsf{K}^2 \brak{\mathsf{B}_6}^2
\,.
\end{align}
\end{subequations}

Combining the estimates \eqref{eq:I:W:tt:4:final}, \eqref{eq:I:ZA:tt:8:final}, and \eqref{eq:money:dont:jiggle}, we may summarize the forcing and commutator bounds as
\begin{align}
&\sabs{I^{\WW_\tau}_4}
+ 
\sabs{I^{\ZZ_\tau}_8}
+
\sabs{I^{\AA_\tau}_8}
+
\sabs{I^{\ZZ_\tau}_{5+7,a,iii}}
+
\sabs{I^{\AA_\tau}_{5+7,a,v}}
+
\sabs{I^{\AA_\tau}_{3,a,iii}}
\notag\\
&\leq \tfrac{\bubu{35}(1+\alpha)^2}{\alpha \eps} 
\int_0^\s \snorm{\tfrac{\mathcal{J}^{\frac 34}}{\Sigma^\beta} \nbs^6 \Zbt(\cdot,\s')}_{L^2_x}^2 {\rm d}\s'
+ \tfrac{\bubu{35}(1+\alpha)}{\eps} 
\int_0^\s \snorm{\tfrac{\mathcal{J}^{\frac 34}}{\Sigma^\beta} \nbs^6 \Abt(\cdot,\s')}_{L^2_x}^2 {\rm d}\s'
\notag\\
&\qquad 
+ \tfrac{3}{(1+\alpha)\eps^3} \int_0^{\s} \snorm{ \tfrac{\Q \mathcal{J}^{\frac 14}}{\Sigma^\beta}\nn \cdot  \nbs^6 \tt (\cdot,\s')}_{L^2_{x}}^2 {\rm d} \s' 
+ \bubu{10} \alpha \eps (\tfrac{4}{\kappa_0})^{2\beta} \mathsf{B}_6^2
+ \Cn  \eps^{\frac 32}  (\tfrac{4}{\kappa_0})^{2\beta}  \mathsf{K}^2 \brak{\mathsf{B}_6}^2
\,.
\label{eq:hate:9}
\end{align}

\subsection{Conclusion of the six derivative tangential energy bounds}
\label{sec:D6:tau:final}
It remains to merge the identity \eqref{D6-L2-tan}, the lower bounds from \eqref{eq:I:tau:1246} and \eqref{eq:hate:8}, the estimate~\eqref{eq:G:good:G:bad}, the fact that $\mathsf{G_{good}} \geq - \frac{1+\alpha}{3\eps} \mathcal{J}^{\frac 32}$, and  the upper bounds from \eqref{eq:hate:8} and \eqref{eq:hate:9}, we obtain
\begin{align}
0
&\geq  ( \tfrac{1}{52} - \Cn \eps )  \snorm{\tfrac{\mathcal{J}^{\frac 34} (\Jg \Q)^{\frac 12}}{\Sigma^\beta} \nbs^6 \Wbt(\cdot,\s)}_{L^2_x}^2
+
( \tfrac 12 - \Cn \eps )
\Bigl(  \snorm{\tfrac{\mathcal{J}^{\frac 34} (\Jg \Q)^{\frac 12}}{\Sigma^\beta} \nbs^6 \Zbt(\cdot,\s)}_{L^2_x}^2
+ 2  \snorm{\tfrac{\mathcal{J}^{\frac 34} (\Jg \Q)^{\frac 12}}{\Sigma^\beta} \nbs^6 \Abt(\cdot,\s)}_{L^2_x}^2 \Bigr)
\notag\\
&\qquad 
- (\tfrac 12  + \Cn \eps) \Bigl( \snorm{ \tfrac{\mathcal{J}^{\frac 34}(\Jg \Q)^{\frac 12}}{\Sigma^\beta} \nbs^6 \Wbt(\cdot,0)}_{L^2_x}^2
+ \snorm{\tfrac{\mathcal{J}^{\frac 34} (\Jg \Q)^{\frac 12}}{\Sigma^\beta} \nbs^6 \Zbt(\cdot,0)}_{L^2_x}^2
+2  \snorm{\tfrac{\mathcal{J}^{\frac 34} (\Jg \Q)^{\frac 12}}{\Sigma^\beta} \nbs^6 \Abt(\cdot,0)}_{L^2_x}^2
\Bigr)
\notag\\
&\qquad
+ (\tfrac{1+\alpha}{16\eps} - \Cn \eps \beta) \int_0^{\s}  
\snorm{\tfrac{\mathcal{J}^{\frac 14} \Jg^{\frac 12}}{\Sigma^\beta} \nbs^6 \Wbt(\cdot,\s')}_{L^2_x}^2 
 {\rm d} \s'
\notag\\
&\qquad
+ (\tfrac{1+\alpha}{16\eps} - \Cn \eps \beta) \int_0^{\s}  
\Bigl(
\snorm{\tfrac{\mathcal{J}^{\frac 14} \Jg^{\frac 12}}{\Sigma^\beta} \nbs^6 \Zbt(\cdot,\s')}_{L^2_x}^2 
+
\snorm{\tfrac{\mathcal{J}^{\frac 14} \Jg^{\frac 12}}{\Sigma^\beta} \nbs^6 \Abt(\cdot,\s')}_{L^2_x}^2 
\Bigr) {\rm d} \s'
\notag\\
&\qquad 
+  \bigl( \tfrac{4\alpha( \beta - \frac 12)}{5\eps}   -  \tfrac{4(1+\alpha)}{\eps} - \tfrac{25^2(1+\alpha)}{\eps}  - \tfrac{\bubu{35}(1+\alpha)^2}{\alpha \eps}  \bigr)
\int_0^{\s} 
\Bigl( \snorm{\tfrac{ \mathcal{J}^{\frac 34} }{\Sigma^\beta} \nbs^6\Zbt (\cdot,\s')}_{L^2_x}^2
+  \snorm{\tfrac{ \mathcal{J}^{\frac 34} }{\Sigma^\beta} \nbs^6\Abt (\cdot,\s')}_{L^2_x}^2 \Bigr)
{\rm d} \s'
\notag\\ 
&\qquad 
- \tfrac{250^2}{\eps} \int_0^{\s} 
\Bigl(\snorm{\tfrac{\mathcal{J}^{\frac 34}(\Jg \Q)^{\frac 12}}{\Sigma^\beta} \nbs^6 \Wbt(\cdot,\s')}_{L^2_x}^2
+  \snorm{\tfrac{\mathcal{J}^{\frac 34}(\Jg \Q)^{\frac 12}}{\Sigma^\beta} \nbs^6 \Zbt(\cdot,\s')}_{L^2_x}^2
+    \snorm{\tfrac{\mathcal{J}^{\frac 34}(\Jg \Q)^{\frac 12}}{\Sigma^\beta} \nbs^6 \Abt(\cdot,\s')}_{L^2_x}^2 \Bigr)
{\rm d} \s'
\notag\\
&\qquad 
- \tfrac{33 \alpha  (\beta -\frac 12) }{(1+\alpha)\eps}  \int_0^{\s} \Bigl(  \snorm{\tfrac{\mathcal{J}^{\frac 34}(\Jg \Q)^{\frac 12}}{\Sigma^\beta} \nbs^6 \Zbt(\cdot,\s')}_{L^2_x}^2
+    \snorm{\tfrac{\mathcal{J}^{\frac 34}(\Jg \Q)^{\frac 12}}{\Sigma^\beta} \nbs^6 \Abt(\cdot,\s')}_{L^2_x}^2 \Bigr)
{\rm d} \s'
\notag\\
&\qquad 
- 31 (1+\alpha)\eps (\tfrac{3}{\kappa_0})^{2\beta} \Cdatatwo 
- \bubu{10} \alpha \eps (\tfrac{4}{\kappa_0})^{2\beta}\mathsf{B}_6^2
-  \Cn \eps^{\frac 32} \beta^2 (\tfrac{4}{\kappa_0})^{2\beta} \mathsf{K}^2  \brak{\mathsf{B}_6}^2
\notag\\
&\qquad 
- \tfrac{34^2+500^2}{\eps} \int_0^{\s}  \snorm{\tfrac{\mathcal{J}^{\frac 34} (\Jg \Q)^{\frac 12}}{\Sigma^\beta} \nbs^6 \Wbt(\cdot,\s')}_{L^2_x}^2  {\rm d} \s'
\notag\\
&\qquad 
+ \tfrac{7}{40(1+\alpha)\eps^3} \int_0^{\s} \snorm{ \tfrac{\Q \mathcal{J}^{-\frac 14}}{\Sigma^\beta}\nn \cdot  \nbs^6 \tt (\cdot,\s')}_{L^2_{x}}^2 {\rm d} \s'
\notag\\
&\qquad
+ \tfrac{1}{20(1+\alpha)\eps^2} \snorm{ \tfrac{\Q \mathcal{J}^{\frac 14}}{\Sigma^\beta}\nn \cdot  \nbs^6 \tt (\cdot,\s)}_{L^2_{x}}^2
-  \tfrac{28+ 100\cdot 250^2+500^2}{(1+\alpha)\eps^3} \int_0^\s \snorm{ \tfrac{\Q \mathcal{J}^{\frac 14}}{\Sigma^\beta}\nn \cdot  \nbs^6 \tt (\cdot,\s')}_{L^2_{x}}^2 {\rm d} \s' 
\,,
\label{eq:hate:10}
\end{align}
where $\Cn = \Cn (\alpha,\kappa_0,\Cdata)$ is independent of $\beta$ (and $\eps$, as always).

At this stage, we choose $\beta = \beta(\alpha)$ to be sufficiently large to ensure that the damping for $\nbs^6(\Zbt,\Abt)$ is strong enough, i.e., so that 
\begin{equation*}
 \tfrac{4\alpha( \beta - \frac 12)}{5\eps}   -  \tfrac{4(1+\alpha)}{\eps} - \tfrac{25^2(1+\alpha)}{\eps}  - \tfrac{\bubu{35}(1+\alpha)^2}{\alpha \eps} \geq  0 \,.
\end{equation*}
More precisely, we choose $\beta$ to ensure equality in the above inequality; namely, as 
\begin{equation}
\beta_\alpha  := \tfrac 12 + \tfrac{5(1+\alpha)}{4\alpha} \Bigl(  4 + 25^2 + \tfrac{\bubu{35}(1+\alpha)}{\alpha}  \Bigr) \,.
\label{eq:tangential:bounds:beta}
\end{equation}
With this choice of $\beta = \beta_\alpha$, we return to \eqref{eq:hate:10}, choose $\eps$ to be sufficiently small in terms of $\alpha,\kappa_0,\Cdata$, and use \eqref{table:derivatives} to bound the initial data term. After re-arranging we deduce that
\begin{align}
&\tfrac{1}{53}    \snorm{\tfrac{\mathcal{J}^{\frac 34} (\Jg \Q)^{\frac 12}}{\Sigma^{\beta_\alpha}} \nbs^6(\Wbt,\Zbt,\Abt)(\cdot,\s)}_{L^2_x}^2
+ \tfrac{1}{20(1+\alpha)\eps^2} \snorm{ \tfrac{\Q \mathcal{J}^{\frac 14}}{\Sigma^{\beta_\alpha}}\nn \cdot  \nbs^6 \tt (\cdot,\s)}_{L^2_{x}}^2
\notag\\
&\quad
+ \tfrac{1+\alpha}{18\eps}   \int_0^{\s}  
\snorm{\tfrac{\mathcal{J}^{\frac 14} \Jg^{\frac 12}}{\Sigma^{\beta_\alpha}} \nbs^6 (\Wbt,\Zbt,\Abt)(\cdot,\s')}_{L^2_x}^2 
 {\rm d} \s'
+ \tfrac{7}{40(1+\alpha)\eps^3} \int_0^{\s} \snorm{ \tfrac{\Q \mathcal{J}^{-\frac 14}}{\Sigma^{\beta_\alpha}}\nn \cdot  \nbs^6 \tt (\cdot,\s')}_{L^2_{x}}^2 {\rm d} \s'
\notag\\
&\leq  \tfrac{54}{53}  \snorm{ \tfrac{\mathcal{J}^{\frac 34}(\Jg \Q)^{\frac 12}}{\Sigma^{\beta_\alpha}} \nbs^6 (\Wbt,\Zbt,\Abt)(\cdot,0)}_{L^2_x}^2
\notag\\
&\quad 
+ 31 (1+\alpha)\eps (\tfrac{3}{\kappa_0})^{2\beta_\alpha} \Cdatatwo 
+  \bubu{10} \alpha \eps (\tfrac{4}{\kappa_0})^{2\beta_\alpha}\mathsf{B}_6^2
+ \Cn \eps^{\frac 32} \beta_\alpha^2 (\tfrac{4}{\kappa_0})^{2\beta_\alpha} \mathsf{K}^2  \brak{\mathsf{B}_6}^2
\notag\\
&\quad 
+ \tfrac{C_1(\alpha)}{\eps}     \int_0^{\s}
\snorm{\tfrac{\mathcal{J}^{\frac 34}(\Jg \Q)^{\frac 12}}{\Sigma^{\beta_\alpha}} \nbs^6 (\Wbt,\Zbt,\Abt)(\cdot,\s')}_{L^2_x}^2
 {\rm d} \s'
 + \tfrac{C_2(\alpha)}{\eps^3} \int_0^\s \snorm{ \tfrac{\Q \mathcal{J}^{\frac 14}}{\Sigma^{\beta_\alpha}}\nn \cdot  \nbs^6 \tt (\cdot,\s')}_{L^2_{x}}^2 {\rm d} \s' 
\notag\\
&\leq 33 (1+\alpha) \eps (\tfrac{3}{\kappa_0})^{2\beta_\alpha} \Cdatatwo
+  \bubu{10} \alpha \eps (\tfrac{4}{\kappa_0})^{2\beta_\alpha}\mathsf{B}_6^2
+ \Cn \eps^{\frac 32} \beta_\alpha^2 (\tfrac{4}{\kappa_0})^{2\beta_\alpha} \mathsf{K}^2  \brak{\mathsf{B}_6}^2
\notag\\
&\quad 
+ \tfrac{C_1(\alpha)}{\eps}    \int_0^{\s}
\snorm{\tfrac{\mathcal{J}^{\frac 34}(\Jg \Q)^{\frac 12}}{\Sigma^{\beta_\alpha}} \nbs^6 (\Wbt,\Zbt,\Abt)(\cdot,\s')}_{L^2_x}^2
 {\rm d} \s'
 + \tfrac{C_2(\alpha)}{\eps^3} \int_0^\s \snorm{ \tfrac{\Q \mathcal{J}^{\frac 14}}{\Sigma^{\beta_\alpha}}\nn \cdot  \nbs^6 \tt (\cdot,\s')}_{L^2_{x}}^2 {\rm d} \s' \,,
\label{eq:hate:11}
\end{align}
where $C_1(\alpha) = \frac{33 \alpha  (\beta_\alpha -\frac 12)}{1+\alpha} + (34^2 + 2 \cdot 500^2) $ and $C_2(\alpha) =\tfrac{28+ 100\cdot 250^2+500^2}{(1+\alpha)}$ are two constants that depend only on $\alpha$. 
By inspecting the first line on the left side and the last line on the right side of \eqref{eq:hate:11}, we observe that we may apply Gr\"onwall's inequality for $\s \in [0,\eps]$. More precisely, there exists a constant 
\begin{equation}
\hat{\mathsf{c}}_\alpha >0
\label{eq:hate:11:a}
\end{equation}
which only depends on $\alpha$, and may be computed explicitly from \eqref{eq:Q:bbq}, \eqref{eq:tangential:bounds:beta}, and \eqref{eq:hate:11}, such that 
\begin{align}
& 
\sup_{\s \in [0,\eps]}
\snorm{\tfrac{\mathcal{J}^{\frac 34}  \Jg^{\frac 12}}{\Sigma^{\beta_\alpha}} \nbs^6(\Wbt,\Zbt,\Abt)(\cdot,\s)}_{L^2_x}^2 
 + \tfrac{1}{\eps}   \int_0^{\eps}  
\snorm{\tfrac{\mathcal{J}^{\frac 14} \Jg^{\frac 12}}{\Sigma^{\beta_\alpha}} \nbs^6 (\Wbt,\Zbt,\Abt)(\cdot,\s)}_{L^2_x}^2   {\rm d} \s
 \notag\\
 &\qquad 
+ \tfrac{1}{\eps^2} \sup_{\s \in [0,\eps]} 
 \snorm{ \tfrac{\mathcal{J}^{\frac 14}}{\Sigma^{\beta_\alpha}}\nn \cdot  \nbs^6 \tt (\cdot,\s)}_{L^2_{x}}^2
+ \tfrac{1}{\eps^3} \int_0^{\eps} \snorm{ \tfrac{\mathcal{J}^{-\frac 14}}{\Sigma^{\beta_\alpha}}\nn \cdot  \nbs^6 \tt (\cdot,\s)}_{L^2_{x}}^2 {\rm d} \s
\notag\\
 &  
 \leq \hat{\mathsf{c}}_\alpha \eps (\tfrac{4}{\kappa_0})^{2\beta_\alpha} 
 \Big( \Cdatatwo + \mathsf{B}_6^2 + \Cn \eps^{\frac 12}  \mathsf{K}^2  \brak{\mathsf{B}_6}^2\Bigr)
 \,.
 \label{eq:hate:12}
\end{align}
At last, we multiply the above estimate by $\kappa_0^{2 \beta_\alpha}$, appeal to \eqref{bs-Sigma}, drop the energy and damping terms for $\nn \cdot \nbs^6 \tt$ (since these were bounded already in Proposition~\ref{prop:geometry}), use the inequality $\mathcal{J} \leq \Jg$,  and recall the definitions of $\widetilde{\mathcal{E}}_{6,\ttt}^2(\s)$ and $\widetilde{\mathcal{D}}_{6,\ttt}^2(\s) $ to deduce that 
\begin{align}
\eps \sup_{\s \in [0,\eps]} \widetilde{\mathcal{E}}_{6,\ttt}^2(\s)
+\widetilde{\mathcal{D}}_{6,\ttt}^2(\eps) 
& 
\leq \hat{\mathsf{c}}_\alpha \eps^2   4^{2\beta_\alpha} 
 \Big( \Cdatatwo + \mathsf{B}_6^2 + \Cn \eps^{\frac 12}  \mathsf{K}^2  \brak{\mathsf{B}_6}^2\Bigr)
 \notag\\
&
\leq (\eps \mathsf{K})^2 \mathsf{B}_6^2 \cdot \hat{\mathsf{c}}_\alpha   4^{2\beta_\alpha} 
 \Big(\tfrac{ \Cdatatwo + \mathsf{B}_6^2}{\mathsf{K}^2 \mathsf{B}_6^2} + \Cn \eps^{\frac 12}   \tfrac{\brak{\mathsf{B}_6}^2}{\mathsf{B}_6^2}\Bigr) 
 \,.
 \label{eq:hate:13:aa}
\end{align}
Upon defining
\begin{equation}
\mathsf{K} := 8 \max\{1, \hat{\mathsf{c}}_\alpha^{\frac 12}  4^{\beta_\alpha} \} \,,
\label{eq:K:choice:1}
\end{equation}
where $\beta_\alpha$ is as defined in~\eqref{eq:tangential:bounds:beta}, 
and ensuring that 
\begin{equation}
\mathsf{B}_6 \geq \max\{1, \Cdata \} \,,
\label{eq:B6:choice:1}
\end{equation}
and $\eps$ sufficiently small in terms of $\alpha,\kappa_0,\Cdata$, we deduce from \eqref{eq:hate:13:aa} that 
\begin{equation}
\eps \sup_{\s \in [0,\eps]} \widetilde{\mathcal{E}}_{6,\ttt}^2(\s)
+\widetilde{\mathcal{D}}_{6,\ttt}^2(\eps) 
\leq \tfrac{1}{8} (\eps \mathsf{K})^2 \mathsf{B}_6^2 
 \,,
 \label{eq:hate:13:a}
\end{equation}
which closes the ``tangential part'' of the remaining bootstrap \eqref{bootstraps-Dnorm:6}.


\section{Improved normal-component estimates for six pure time derivatives}
\label{sec:pure:time}
We now consider a new set of sixth-order energy estimates for $\Zbn$ and $\Abn$, which involve only time, or ``material'' derivatives. The bounds in this section, estimate \eqref{eq:madman:2} to be more precise, were used in Section~\ref{sec:sixth:order:energy-tangential} to bound the $2\alpha ( \nn \cdot \tt,_1 - \Jg g^{-\frac 12} \nbs_2h \nn \cdot \nbs_2 \tt) \nbs^6 \Zbn$ contribution to $\mathcal{R}^\tau_{\Zb}$ (see~\eqref{eq:I:ZA:tt:8:c:00000}), and are used in  Section~\ref{sec:sixth:order:energy} to bound the $\tfrac{2\alpha }{\eps} \nbs_1\Jg \nbs^6\Zbn$ contribution to $\mathcal{R}_\Zb^\nn$ (see~\eqref{10hard3}). Except for these remainder terms, the bounds in this section are not used anywhere else in the argument.

We define the operator $\eps$-rescaled ALE transport operator in $(x,\s)$ coordinates by
\begin{equation*}
\nbd =   \eps (\Q\p_\s+V\p_2) =  \nbs_\s + \eps V\nbs_2  \,.
\end{equation*} 
With the above notation, the goal of this section is to establish:
\begin{proposition}
\label{prop:madman}
Under the standing bootstrap assumptions~\eqref{bootstraps}, assuming that $\eps$ is sufficiently small with respect to $\alpha,\kappa_0$, and $\Cdata$, we have
\begin{equation} 
\sup_{\s\in[0,\eps]} \snorm{ \Jgh \mathcal{J}^{\! \frac 34} \nbd^6 (\Zbn, \Abn)( \cdot , \s)}_{L^2_x}^2
+
\tfrac{1}{\eps} \int_0^\eps   \snorm{\mathcal{J}^{\! \frac 14} \Jgh \nbd^6 (\Zbn, \Abn)( \cdot ,\s)}^2_{L^2_x}  {\rm d} \s
\les \eps \mathsf{K}^2 \brak{\mathsf{B}_6}^2 \,,
\label{eq:madman}
\end{equation}
where the implicit constant in \eqref{eq:madman} only depends on $\alpha,\kappa_0$, and $\Cdata$.
\end{proposition}
 
Before turning to the proof of the above estimate, we record a corollary which is used throughout our proof. 
\begin{corollary}
\label{cor:beyond:madman}
Under the standing bootstrap assumptions~\eqref{bootstraps}, assuming that $\eps$ is sufficiently small with respect to $\alpha,\kappa_0$, and $\Cdata$, we have that 
\begin{subequations}
\label{eq:madman:2:all}
\begin{align} 
\snorm{\mathcal{J}^{\! \frac 34} \Jgh \nbs_{\s}^6  \Zbn }_{L^\infty_{\s} L^2_{x}} 
&\les
\eps^{\frac 12} \mathsf{K} \brak{\mathsf{B}_6} \,,
\label{eq:madman:2:75}
\\
\snorm{\mathcal{J}^{\! \frac 34} \nbs_{\s}^6  \Zbn }_{L^2_{x,\s}} 
\leq
\snorm{\mathcal{J}^{\! \frac 14} \Jgh \nbs_{\s}^6  \Zbn }_{L^2_{x,\s}} 
&\les \eps \mathsf{K} \brak{\mathsf{B}_6} 
\,,
\label{eq:madman:2}
\end{align}
\end{subequations}
where the implicit constant depends only on $\alpha,\kappa_0$, and $\Cdata$.
\end{corollary}
\begin{proof}[Proof of Corollary~\ref{cor:beyond:madman}]
We note that for $k\geq 2$ and sufficiently smooth functions $f$, we have
\begin{align} 
\nbd^k f 
- \nbs_\s^k f 
&= \sum_{i=0}^{k-1} {k \choose i} (\eps V)^{k-i} \nbs_\s^i \nbs_2^{k-i}  f  
\notag\\
&+ 
\eps \sum_{i=0}^{k-2} \sum_{n=0}^{k-i-1} \sum_{j=0}^i  \nbs_{\s}^{i-j}  \nbs_2^{j+1} f \cdot 
 \eps^{j+n}
  \sum_{|\alpha| = k-i-1-n , |\beta| = n} c_{k,i,n,j,\alpha,\beta} \prod_{\ell=1}^{j+1+n} \nbs_\s^{\alpha_\ell} \nbs_2^{\beta_\ell}  V 
\,,
\label{eq:expand:nbs:k}
\end{align}
for   suitable combinatorial coefficients $c_{k,i,n,j,\alpha,\beta}\geq 0$. 
Identity~\eqref{eq:expand:nbs:k} with $k=6$ shows that $\nbd^6 \Zbn - \nbs_{\s}^6 \Zbn$ consists of a sum of terms with at most five derivatives on $\nbs_2 \Zbn$  times a power of $\eps$ which is at least equal to one. These terms are already bounded by~\eqref{eq:Jg:Zbn:D5:improve:a}, by \eqref{eq:Jg:Zbn:p2D5:improve:new} (with $\bar \beta = 0$ and $\bar a = \frac 12$), which for instance gives $\snorm{\mathcal{J}^{\frac 14} \Jgh \nbs_2 \nbs^5 \Zbn }_{L^2_{x,\s}}  
\le \tfrac{4 \mathsf{K}}{\alpha \kappa_0 } \mathsf{B_6}$, and by \eqref{eq:Jg:Zbn:p2D6:sup:improve:new}. Using the triangle inequality, the bounds for $V$ established in~\eqref{eq:V:H6:new}--\eqref{eq:V:H6:new:bdd}, and the interpolation bounds in~Lemma~\ref{lem:time:interpolation} to treat the non-endpoint cases, it is then clear that \eqref{eq:madman} implies both \eqref{eq:madman:2:75} and \eqref{eq:madman:2}.
\end{proof}

\begin{corollary}
\label{cor:useful:junk}
Expansion~\eqref{eq:expand:nbs:k} implies that we have the bounds
\begin{subequations}
\label{eq:very:useful:junk}
\begin{align}
\snorm{\nbd^k f}_{L^2_{x,\s}} 
&\leq 2 \snorm{\nbs^k f}_{L^2_{x,\s}} 
+ \Cn \eps^3 \mathsf{K} \brak{\mathsf{B}_6}  \mathfrak{B}_f
\,,
\quad  0 \leq k \leq 6
\,,
\label{eq:useful:junk}
\\
\snorm{\mathcal{J}^a \Jg^{\! b} \nbd^5 f}_{L^2_{x,\s}} 
&\leq 2 \snorm{\mathcal{J}^a \Jg^{\! b} \nbs^5 f}_{L^2_{x,\s}} 
+ \Cn \eps^3 \mathsf{K} \brak{\mathsf{B}_6} 
\bigl(\snorm{\nbs^3 f(\cdot,0)}_{L^\infty_x} + \eps^{-1} \snorm{\nbs^4 f}_{L^2_{x,\s}} \bigr)
\,,
\quad  a,b\geq 0
\,,
\label{eq:useful:junk:L2L2}
\\
\snorm{\mathcal{J}^a \Jg^{\!b} \nbd^5 f}_{L^\infty_\s L^2_x} 
&\leq 2 \snorm{\mathcal{J}^a \Jg^{\!b}  \nbs^5 f}_{L^\infty_\s L^2_x} 
+ \Cn \eps^{\frac 52} \mathsf{K} \brak{\mathsf{B}_6} 
\bigl(\snorm{\nbs^3 f(\cdot,0)}_{L^\infty_x} + \eps^{-1} \snorm{\nbs^4 f}_{L^2_{x,\s}} \bigr)
\,,
\quad  a,b\geq 0
\,,
\label{eq:useful:junk:LinftyL2}
\end{align}
\end{subequations}
where we use the notation in \eqref{eq:pink:shirt}. 
\end{corollary}
\begin{proof}[Proof of Corollary~\ref{cor:useful:junk}]
For simplicity, we only give the proof of \eqref{eq:useful:junk} when $k=6$, the most difficult case. Due to \eqref{bs-V}, the first line of \eqref{eq:expand:nbs:k} contributes   at most $(1+ \Cn \eps^2) \|\nbs^k f\|_{L^2_{x,\s}}$ to the upper bound. In order to deal with the second line in \eqref{eq:expand:nbs:k}, note that the $V$ estimates in \eqref{bs-V} and \eqref{eq:V:H6:new} imply that $\|\nbs^j V\|_{L^\infty_{x,\s}} \les\eps \mathsf{K} \brak{\mathsf{B}_6}$ for all $0\leq j \leq 3$. As such, by also appealing to the Poincar\'e-type inequality \eqref{eq:x1:Poincare}, we see that the second line of \eqref{eq:expand:nbs:k} only has a nontrivial contribution when $i=0$ or $i=1$. For these special cases, we use the interpolation inequality  \eqref{se3:time}, to obtain
\begin{equation*}
\snorm{\nbd^6 f}_{L^2_{x,\s}} \leq (1 + \Cn \eps^2\mathsf{K} \brak{\mathsf{B}_6})  \snorm{\nbs^6 f}_{L^2_{x,\s}} + \Cn \eps \sum_{i=0}^{1} (\|\nbs^6 f\|_{L^2_{x,\s}}^{\frac{i+1}{6}} {\mathfrak{B}_f}^{\frac{5-i}{6}} + \eps^{\frac{i+1}{6}} \mathfrak{B}_f) (\eps \mathsf{K}\brak{\mathsf{B}_6})^{\frac{5-i}{6}} \eps.
\end{equation*}
The bound \eqref{eq:useful:junk} with $k=6$ now follows from the $\eps$-Young inequality. 

The bounds \eqref{eq:useful:junk:L2L2} and \eqref{eq:useful:junk:LinftyL2} follow since the $L^\infty_{x,\s}$ bound on $\nbs^3 V$ implies 
\begin{align*}
\snorm{\mathcal{J}^a \Jg^{\! b} \nbd^5 f}_{L^2_{x,\s}} 
&\leq (1 + \Cn \eps^2\mathsf{K} \brak{\mathsf{B}_6})  \snorm{\mathcal{J}^a \Jg^{\! b} \nbs^5 f}_{L^2_{x,\s}} + \Cn \eps \|   \nbs_2 f \nbs^5 V \|_{L^2_{x,\s}}
,
\\
\snorm{\mathcal{J}^a \Jg^{\! b} \nbd^5 f}_{L^\infty_\s L^2_x} 
&\leq (1 + \Cn \eps^2\mathsf{K} \brak{\mathsf{B}_6})  \snorm{\mathcal{J}^a \Jg^{\! b} \nbs^5 f}_{L^\infty_\s L^2_x} + \Cn \eps \|  \nbs_2 f  \nbs^5 V \|_{L^\infty_\s L^2_x}
.
\end{align*} 
We conclude the proof of  \eqref{eq:useful:junk:L2L2}--\eqref{eq:useful:junk:LinftyL2} by appealing to \eqref{eq:Sobolev}, \eqref{eq:V:H6:new}, and to \eqref{eq:r:Linfty:time} with $r=0$.
\end{proof}

The remainder of this section is dedicated to the proof of Proposition~\ref{prop:madman}. 
The proof of \eqref{eq:madman} is based on an energy estimate for $\Abn$ and $\Zbn$ (see~\eqref{D_sD^5-L2}), in which the $\Jg \Wbn$ evolution is used passively as a ``constitutive relation''. The proof consists of several steps, which are then finalized in Section~\ref{sec:proof:eq:madman}. 

In preparation, we write equations \eqref{eq:Jg:Wb:nn}, \bubu{\eqref{eq:Zb:nn}, and \eqref{eq:Ab:nn}} in $(x,\s)$ coordinates as
\begin{subequations} 
\begin{align}
&(\Q\p_\s+V\p_2) (\Jg \Wbn) + \alpha \Sigma g^{-\frac 12} \Jg \nbs_2  \Abn 
= \mathcal{G}^\nn_\Wb 
\,, \label{eq:Jg:Wb:nn-s} \\
&\tfrac{\Jg}{\Sigma} (\Q\p_\s+V\p_2)  \Zbn 
-\alpha g^{-\frac 12}   \Jg  \nbs_2  \Abn
- 2\alpha    \Zbn,_1 
+ 2\alpha  g^{-\frac 12} \nbs_2 h   \Jg \nbs_2   \Zbn 
 = \mathcal{G}^\nn_\Zb  
\,, \label{eq:Jg:Zb:nn-s}  \\
&\tfrac{\Jg}{\Sigma} (\Q\p_\s+V\p_2)   \Abn  
- \tfrac{\alpha}{2}  g^{-\frac 12} \Jg \nbs_2  \Zbn 
-\alpha  \Abn,_1 
+ \alpha   g^{-\frac 12} \nbs_2 h \Jg \nbs_2   \Abn 
+ \tfrac{\alpha}{2} g^{-\frac 12}  \nbs_2 (\Jg\Wbn) 
= \mathcal{G}^\nn_\Ab 
 \,,
\label{eq:Jg:Ab:nn-s}
\end{align}
\end{subequations} 
where
\begin{subequations} 
\begin{align} 
\mathcal{G}^\nn_\Wb 
& := 
\tfrac{\alpha}{2} \Sigma g^{-\frac 32}  \bigl(   \Jg\Wbn + \Jg \Zbn - 2\Jg\Abt \bigr)  \nbs_2^2 h
- \tfrac{\alpha}{2} \bigl( \Jg\Wbn - \Jg \Zbn  \bigr) \Abt 
\notag\\
&\qquad
-  \bigl(\tfrac{3+\alpha}{2} \Wbt + \tfrac{1-\alpha}{2} \Zbt\bigr) \Jg\Abn
- \bigl(\tfrac{1+\alpha}{2} \Wbt + \tfrac{1-\alpha}{2} \Zbt\bigr)\Jg \Wbt  
 \,, \label{Gwn} \\
\mathcal{G}^\nn_\Zb 
& := 
- \tfrac{\alpha}{2} g^{-\frac 32}  \bigl(\Jg \Wbn +  \Jg \Zbn - 2 \Jg \Abt\bigr) \nbs_2^2 h
+
 2\alpha \bubu{g^{-\frac 12}}  ( \Zbt+ \Abn)  \nbs_2 \Jg     
-   \tfrac{1}{\Sigma} \bigl( \bubu{\tfrac{1-\alpha}{2}} \Jg \Wbn + \bubu{\tfrac{1+\alpha}{2}} \Jg \Zbn \bigr)  \Zbn
\notag\\
&\qquad  
+
\tfrac{1}{\Sigma} \Big(
\tfrac{\alpha}{2}   \Abt \bigl(   \Jg  \Wbn   - \Jg \Zbn \bigr)
-  \bigl(\tfrac{1+\alpha}{2} \Wbt + \tfrac{3-\alpha}{2} \Zbt \bigr) \Jg  \Abn
- \bigl(\tfrac{1+\alpha}{2} \Wbt + \tfrac{1-\alpha}{2} \Zbt\bigr) \Jg  \Zbt
 \Big) 
\,, \label{Gzn} \\
\mathcal{G}^\nn_\Ab 
&:= 
-\tfrac{\alpha}{2} g^{-\frac 32} \Jg (\Wbt - \Zbt) \nbs_2^2 h
+  \alpha   g^{-\frac 12}   ( \tfrac 12 \Zbn +  \Abt ) \nbs_2 \Jg
- \tfrac{1}{2 \Sigma}    (  \Jg\Wbn + \tfrac{1}{2}\Jg \Zbn + 2 \Jg\Abt)   \Abn
\notag\\
&\qquad +
\tfrac{1}{2\Sigma} \Big(
  (\Jg \Wbn + \Jg \Zbn - 2 \Jg \Abt) (\tfrac{1+\alpha}{2} \Wbt + \tfrac{1-\alpha}{2} \Zbt)
- \tfrac{\alpha}{2} (\Jg \Wbn - \Jg \Zbn) (\Wbt -\Zbt)
\Big) \,. \label{Gan}
\end{align} \
\end{subequations} 
Next, we let $\nbd^6$
act upon equations 
\eqref{eq:Jg:Zb:nn-s} and \eqref{eq:Jg:Zb:nn-s} to obtain that
\begin{subequations} 
\begin{align} 
&\tfrac{\Jg}{\Sigma} (\Q\p_s+V\p_2) \nbd^6 \Zbn 
-\alpha g^{-\frac 12}   \Jg   \nbs_2\nbd^6\Abn
- 2\alpha  \nbd^6\Zbn,_1  
+ 2\alpha \Jg g^{-\frac 12} \nbs_2 h   \nbs_2 \nbd^6    \Zbn 
\notag \\
& \qquad\qquad
=  \nbd^6\mathcal{G}^\nn_\Zb  +  \mathfrak{C}^\nn_\Zb \,,  
\label{D_sD^5:Jg:Zb:nn-s}  \\
&\tfrac{\Jg}{\Sigma} (\Q\p_s+V\p_2)  \nbd^6   \Abn  
- \tfrac{\alpha}{2} g^{-\frac 12}  \Jg \nbs_2 \nbd^6 \Zbn 
-\alpha  \nbd^6\Abn,_1 
+ \alpha  g^{-\frac 12} \nbs_2 h \Jg \nbs_2   \nbd^6 \Abn  
+ \tfrac{\alpha}{2} g^{-\frac 12}  \nbs_2 \nbd^6(\Jg\Wbn) 
\notag \\
& \qquad\qquad
= \nbd^6\mathcal{G}^\nn_\Ab    +  \mathfrak{C}^\nn_\Ab   
\,,  \label{D_sD^5:Jg:Ab:nn-s} 
\end{align} 
\end{subequations} 
where the commutator terms are given by
\begin{subequations} 
\label{D_sD^5:Jg:Zn:An-commutators} 
\begin{align} 
\mathfrak{C}^\nn_\Zb
&:=- \tfrac{1}{\eps} \jump{\nbd^6,\Jg \Sigma^{-1}} \nbd   \Zbn 
+ \alpha \jump{\nbd^6, g^{-\frac 12} \Jg} \nbs_2 \Abn
-2\alpha \jump{\nbd^6, g^{-\frac 12} \nbs_2 h \Jg} \nbs_2  \Zbn \,, \\
\mathfrak{C}^\nn_\Ab 
&:= - \tfrac{1}{\eps} \jump{\nbd^6,\Jg \Sigma^{-1}} \nbd   \Abn  
+ \tfrac{\alpha}{2} \jump{\nbd^6, g^{-\frac 12} \Jg} \nbs_2 \Zbn 
- \alpha \jump{\nbd^6, g^{-\frac 12} \nbs_2 h \Jg} \nbs_2  \Abn\,.
\end{align} 
\end{subequations} 
Our goal is to compute the  spacetime $L^2$ inner-product:
\begin{equation} 
\tint \Sigma^{-2\beta+1} \mathcal{J}^{\! \frac 32} \Big(
 \underbrace{\eqref{D_sD^5:Jg:Zb:nn-s} \ \nbd^5 \Zbn} _{ \mathfrak{I}^{\ZZ_n}}
+ \underbrace{\eqref{D_sD^5:Jg:Ab:nn-s} \ 2\nbd^5\Abn } _{ \mathfrak{I}^{\AA_n}}
 \Big) {\rm d}x {\rm d}r=0 \,, \label{D_sD^5-L2}
\end{equation} 
where $\beta = \beta(\alpha) >0$ is a constant whose value will be made precise below, see~\eqref{eq:Zbn:Abn:6:time:beta:choice}. For convenience, we will abbreviate $\Sigma^{-2\beta+1} =: \jb$.

\subsection{The integral  $\mathfrak{I} ^{\ZZ_n}$} 
 We additively decompose the integral  $\mathfrak{I} ^{\ZZ_n}$ as 
\begin{subequations} 
 \label{IntZbn}
 \begin{align}
  \mathfrak{I}^{\ZZ_n} & = \mathfrak{I} ^{\ZZ_n}_1 + \mathfrak{I} ^{\ZZ_n}_2 + \mathfrak{I} ^{\ZZ_n}_3+ \mathfrak{I} ^{\ZZ_n}_4+ \mathfrak{I} ^{\ZZ_n}_5
  \,,  \notag \\
 \mathfrak{I} ^{\ZZ_n}_1& =
 \tint \tfrac{\mathcal{J}^{\! \frac 32} \Jg}{2 \Sigma^{2\beta}}  (\Q\p_\s +V\p_2)  \bigl( (\nbd^6 \Zbn  )^2\bigr)
  \,, \label{I1-Zbn-imp}\\
 \mathfrak{I} ^{\ZZ_n}_2 &=
-\alpha \tint \jb \mathcal{J}^{\! \frac 32} g^{- {\frac{1}{2}} } \Jg \nbs_2 \nbd^6\Abn \   \nbd^6\Zbn
  \,,   \label{I2-Zbn-imp} \\
 \mathfrak{I} ^{\ZZ_n}_3 &=
 - 2 \alpha  \tint \jb\mathcal{J}^{\! \frac 32} \nbd^6 \Zbn,_1 \     \nbd^6\Zbn
  \,,   \label{I3-Zbn-imp} \\
 \mathfrak{I} ^{\ZZ_n}_4& =
 \alpha  \tint \jb \mathcal{J}^{\! \frac 32} \Jg g^{-\frac 12} \nbs_2 h    \nbs_2  \bigl( ( \nbd^6 \Zbn)^2 \bigr)
\,,  \label{I4-Zbn-imp}\\
 \mathfrak{I} ^{\ZZ_n}_5& =
-   \tint \jb \mathcal{J}^{\! \frac 32} \big( \nbd^6 \mathcal{G}^\nn_\Zb   +  \mathfrak{C}^\nn_\Zb\big)  \     \nbd^6\Zbn
\,.  \label{I5-Zbn-imp}
 \end{align} 
\end{subequations} 

\subsection{The integral  $\mathfrak{I} ^{\AA_n}$} 
 We additively decompose the integral  $\mathfrak{I} ^{\AA_n}$ as 
 \begin{subequations} 
 \label{IntAbn}
 \begin{align}
  \mathfrak{I}^{\AA_n} & = \mathfrak{I} ^{\AA_n}_1 + \mathfrak{I} ^{\AA_n}_2 + \mathfrak{I} ^{\AA_n}_3+ \mathfrak{I} ^{\AA_n}_4+ \mathfrak{I} ^{\AA_n}_5 + \mathfrak{I} ^{\AA_n}_6
  \,,  \notag \\
 \mathfrak{I} ^{\AA_n}_1& =
 \tint \tfrac{\Jg \mathcal{J}^{\! \frac 32}}{\Sigma^{2\beta}}   (\Q\p_\s +V\p_2) \bigl( (\nbd^6 \Abn)^2 \bigr) 
  \,, \label{I1-Abn-imp}\\
 \mathfrak{I} ^{\AA_n}_2 &=
- \alpha \tint \jb \mathcal{J}^{\! \frac 32} \Jg g^{-\frac 12}  \nbs_2\nbd^6 \Zbn \ 
 \nbd^6\Abn
  \,,   \label{I2-Abn-imp} \\
 \mathfrak{I} ^{\AA_n}_3 &=
 - 2 \alpha  \tint \jb\mathcal{J}^{\! \frac 32} \nbd^6 \Abn,_1 \   \nbd^6\Abn
  \,,   \label{I3-Abn-imp} \\
 \mathfrak{I} ^{\AA_n}_4& =
 \alpha  \tint \jb\mathcal{J}^{\! \frac 32}  \Jg  g^{-\frac 12} \nbs_2 h    \nbs_2 \bigl( (\nbd^6 \Abn)^2 \bigr)
\,,  \label{I4-Abn-imp}\\
 \mathfrak{I} ^{\AA_n}_5& =
-   \tint \jb\mathcal{J}^{\! \frac 32} \big( \nbd^6\mathcal{G}^\nn_\Ab   +  \mathfrak{C}^\nn_\Ab\big)  \   \nbd^6\Abn
\,,  \label{I5-Abn-imp} \\
 \mathfrak{I} ^{\AA_n}_6& =
 \alpha   \tint \jb\mathcal{J}^{\! \frac 32}  g^{-\frac 12} \nbs_2 \nbd^6(\Jg\Wbn) \   \nbd^6\Abn
\,.  \label{I6-Abn-imp}
 \end{align} 
\end{subequations} 

\subsection{The exact derivative terms} 
We note that using \eqref{eq:adjoints} we may integrate by parts certain exact derivative terms present in \eqref{D_sD^5-L2}, so that analogously to \eqref{eq:heavy:fuel:1}, \eqref{eq:heavy:fuel:G1:lower}, and \eqref{eq:heavy:fuel:G2:lower}, we have 
\begin{align}
  &
  \mathfrak{I} ^{\ZZ_n}_1 +  \mathfrak{I} ^{\AA_n}_1
  + 
  \mathfrak{I} ^{\ZZ_n}_2 +  \mathfrak{I} ^{\AA_n}_2
  +
  \mathfrak{I} ^{\ZZ_n}_3 +  \mathfrak{I} ^{\AA_n}_3
  +
  \mathfrak{I} ^{\ZZ_n}_4 +  \mathfrak{I} ^{\AA_n}_4
  \notag\\
  &=
  \tint \tfrac{\mathcal{J}^{\frac 32}\Jg}{\Sigma^{2\beta}}   
  (\Q\p_\s +V\p_2) \bigl(\tfrac 12 (\nbd^6 \Zbn)^2 + (\nbd^6 \Abn)^2 \bigr)
 - \alpha \tint \jb \mathcal{J}^{\frac 32} \Jg g^{-\frac 12}  
  \nbs_2\bigl( \nbd^6 \Zbn \,  \nbd^6\Abn\bigr)
 \notag\\
 &\qquad 
 -   \alpha  \tint \jb \mathcal{J}^{\frac 32} \p_1 \bigl( (\nbd^6 \Zbn)^2 + (\nbd^6 \Abn)^2 \bigr)
 +   \alpha  \tint \jb \mathcal{J}^{\frac 32} \Jg g^{-\frac 12} \nbs_2 h \nbs_2 \bigl( (\nbd^6 \Zbn)^2 + (\nbd^6 \Abn)^2 \bigr)
 \notag\\
 &\qquad 
  -  2 \alpha  \tint \jb \mathcal{J}^{\frac 32} \bigl(\nbd^6 \Zbn  \jump{\nbd^6,\p_1}\Zbn  + \nbd^6 \Abn  \jump{\nbd^6,\p_1}\Abn \bigr)
  \notag\\
   &= \tfrac 12 \snorm{ \tfrac{\mathcal{J}^{\frac 34} (\Q \Jg)^{\frac 12}}{\Sigma^\beta} \nbd^6 \Zbn(\cdot,\s)}_{L^2_x}^2 
   + \snorm{ \tfrac{\mathcal{J}^{\frac 34}  (\Q \Jg)^{\frac 12}}{\Sigma^\beta} \nbd^6 \Abn(\cdot,\s)}_{L^2_x}^2 
   \notag\\
   &\qquad 
   - \tfrac 12 \snorm{ \tfrac{\mathcal{J}^{\frac 34}  (\Q  \Jg)^{\frac 12}}{\Sigma^\beta} \nbd^6 \Zbn(\cdot,0)}_{L^2_x}^2 
   - \snorm{ \tfrac{\mathcal{J}^{\frac 34} (\Q  \Jg)^{\frac 12}}{\Sigma^\beta} \nbd^6 \Abn(\cdot,0)}_{L^2_x}^2 
   \notag\\
   &\qquad 
   + \tint  \tfrac{1}{\Sigma^{2\beta}} \bigl( \mathsf{G}^{\ZZ_n} (\nbd^6 \Zbn)^2 + \mathsf{G}^{\AA_n} (\nbd^6 \Abn)^2 \bigr)
  -  \alpha \tint \nbd^6 \Zbn \,  \nbd^6\Abn (\Qr_2 - \nbs_2) \bigl( \jb \mathcal{J}^{\frac 32} \Jg g^{-\frac 12}  \bigr)
   \notag\\
   &\qquad
  + \alpha \int \Qb_2 \jb \mathcal{J}^{\frac 32} \Jg g^{-\frac 12}    \nbd^6 \Zbn  \nbd^6\Abn \Bigr|_{\s}
-   \alpha  \int  \Qb_2 \jb \mathcal{J}^{\frac 32} \Jg g^{-\frac 12} \nbs_2 h  \bigl( (\nbd^6 \Zbn)^2 + (\nbd^6 \Abn)^2 \bigr)  \Bigr|_{\s}
   \notag\\
   &\qquad 
     -  2 \alpha  \tint \jb \mathcal{J}^{\frac 32} \bigl(\nbd^6 \Zbn  \jump{\nbd^6,\p_1}\Zbn  + \nbd^6 \Abn  \jump{\nbd^6,\p_1}\Abn \bigr)
   \label{eq:D6:An:Zn:1}
 \end{align}
where we introduced the coefficients
\begin{align}
\mathsf{G}^{\ZZ_n}
&=  
- \alpha (2\beta -1) \mathcal{J}^{\frac 32} \Sigma,_1 
- \tfrac12  (\Q\p_\s +V\p_2) \bigl( \mathcal{J}^{\frac 32} \Jg\bigr) 
\notag\\
&\qquad 
+ \alpha \Sigma^{2\beta} (\Qr_2 - \nbs_2) \bigl(   \jb \mathcal{J}^{\frac 32} \Jg g^{-\frac 12} \nbs_2 h  \bigr) 
+ \tfrac12 \bigl( V \Qr_2 - \Qr_\s - \nbs_2 V - 2 \alpha \beta (\Zbn + \Abt) \bigr) \mathcal{J}^{\frac 32} \Jg
\\
\mathsf{G}^{\AA_n}
&= 
- \alpha (2\beta -1) \mathcal{J}^{\frac 32} \Sigma,_1 
-  (\Q\p_\s +V\p_2) \bigl( \mathcal{J}^{\frac 32} \Jg\bigr) 
\notag\\
&\qquad 
+ \alpha \Sigma^{2\beta} (\Qr_2 - \nbs_2) \bigl(   \jb \mathcal{J}^{\frac 32} \Jg g^{-\frac 12} \nbs_2 h  \bigr) 
+ \bigl(V \Qr_2 - \Qr_\s - \nbs_2 V - 2 \alpha \beta (\Zbn + \Abt) \bigr) \mathcal{J}^{\frac 32} \Jg
\,.
\end{align}
We note at this stage that for $\beta \geq 1$, analogously to \eqref{eq:heavy:fuel:1}, \eqref{eq:heavy:fuel:G1:lower}, and \eqref{eq:heavy:fuel:G2:lower}, using the lower bound in~\eqref{eq:Q:bbq} and choosing $\eps$ to be sufficiently small in terms of $\alpha,\kappa_0$, and $\Cdata$,  the following lower bounds hold:
\begin{subequations}
\label{eq:D6:An:Zn:1:useful}
\begin{align}
 \mathsf{G}^{\ZZ_n}
 &\geq \bigl( \alpha (\beta-\tfrac 12) + \tfrac{1+\alpha}{4} \bigr) \bigl(\tfrac{4}{5\eps} - \tfrac{33}{(1+\alpha) \eps}  \Jg \Q  \bigr)  \mathcal{J}^{\frac 32} 
 - \tfrac{250^2}{\eps} \Q \Jg \mathcal{J}^{\frac 32}
 + \tfrac{3}{4\eps}\tfrac{2(1+\alpha)}{5} \mathcal{J}^{\frac 12} \Jg  
 \,,
 \\
 \mathsf{G}^{\AA_n}
 &\geq \bigl( \alpha (\beta-\tfrac 12) + \tfrac{1+\alpha}{2} \bigr) \bigl(\tfrac{4}{5\eps} - \tfrac{33}{(1+\alpha) \eps}  \Jg \Q    \bigr)  \mathcal{J}^{\frac 32} 
  - \tfrac{250^2}{\eps} \Q \Jg \mathcal{J}^{\frac 32}
 + \tfrac{3}{2\eps}\tfrac{2(1+\alpha)}{5}  \mathcal{J}^{\frac 12} \Jg   
 \,.
\end{align}
\end{subequations}
The sixth, seventh, and eight terms on the right side of \eqref{eq:D6:An:Zn:1} are bounded by appealing to \eqref{eq:heavy:fuel:G3:lower}. In order to estimate the ninth (and last) term on the right side of \eqref{eq:D6:An:Zn:1}, we observe that expression~\eqref{eq:expand:nbs:k} and the commutator identity \eqref{comm-nbs1-nbs2} also implies
\begin{align*} 
\jump{\nbd^6,\nbs_1} f 
&= - \eps \sum_{i=0}^5 {6 \choose i} (6-i) \nbs_1 V (\eps V)^{5-i} \nbs_2 (\nbs_s^i \nbs_2^{5-i}  f) \notag\\
& \qquad - \eps  \sum_{i=0}^{4} \sum_{n=0}^{5-i} \sum_{j=0}^i \nbs_{\s}^{i-j}  \nbs_2^{j+1} f \cdot \eps^{j+n}   \sum_{|\alpha| = 5-i -n , |\beta| = n} c_{k,i,n,j,\alpha,\beta} \nbs_1 \prod_{\ell=1}^{j+1+n} \nbs_\s^{\alpha_\ell} \nbs_2^{\beta_\ell}  V 
\,,
\end{align*}
Note in particular that for all terms in the above commutator, we have at least a $\nbs_2$ present on $f$. Since in the sixth term on the right side of \eqref{eq:D6:An:Zn:1} $f$ is either $\Abn$ or $\Zbn$, we may appeal to the $\nbs^6 \Abn$ and $\nbs_2 \nbs^5 \Zbn$ bounds in \eqref{eq:Jg:Abn:D6:improve:d} and respectively \eqref{eq:Jg:Zbn:p2D5:improve:new}, to the pointwise bounds on $(\nbs_2 \Abn,\nbs_2 \Zbn)$ in \eqref{bs-nnZb}, to the $V$ estimates in \eqref{bs-V} and \eqref{eq:V:H6:new}, and to the interpolation inequality 
\eqref{se3:time}, to deduce 
\begin{equation}
\left| 2 \alpha  \tint \jb \mathcal{J}^{\frac 32} \bigl(\nbd^6 \Zbn  \jump{\nbd^6,\p_1}\Zbn  + \nbd^6 \Abn  \jump{\nbd^6,\p_1}\Abn \bigr) \right|
\les \eps (\tfrac{4}{\kappa_0})^\beta  \snorm{ \tfrac{\mathcal{J}^{\frac 34} }{\Sigma^{\beta}}  (\nbd^6 \Zbn,\nbd^6 \Abn)}_{L^2_{x,\s}} \mathsf{K} \brak{\mathsf{B}_6}
\,.
\label{eq:D5:An:Zn:2:useful}
\end{equation}

By combining \eqref{eq:D6:An:Zn:1} with \eqref{eq:D6:An:Zn:1:useful}, \eqref{eq:heavy:fuel:G3:lower}, \eqref{eq:D5:An:Zn:2:useful}, and taking $\eps$ to be sufficiently small in terms of $\alpha,\kappa_0,\Cdata$ (not on $\beta \geq 1$),  we arrive at
\begin{align}
  &
  \mathfrak{I} ^{\ZZ_n}_1 +  \mathfrak{I} ^{\AA_n}_1
  + 
  \mathfrak{I} ^{\ZZ_n}_2 +  \mathfrak{I} ^{\AA_n}_2
  +
  \mathfrak{I} ^{\ZZ_n}_3 +  \mathfrak{I} ^{\AA_n}_3
  +
  \mathfrak{I} ^{\ZZ_n}_4 +  \mathfrak{I} ^{\AA_n}_4
  \notag\\
  &\geq
  \bigl(\tfrac 12 - \Cn \eps \bigr) \snorm{ \tfrac{\mathcal{J}^{\frac 34} (\Q \Jg)^{\frac 12}}{\Sigma^\beta} ( \nbd^6 \Zbn, \nbd^6 \Abn)(\cdot,\s)}_{L^2_x}^2 
  - \snorm{ \tfrac{\mathcal{J}^{\frac 34} (\Q \Jg)^{\frac 12}}{\Sigma^\beta}  ( \nbd^6 \Zbn, \nbd^6 \Abn)(\cdot,0)}_{L^2_x}^2 
   \notag\\
   &\qquad 
   +\tfrac{4\alpha(\beta-\frac 12)}{5\eps} \int_0^{\s} \snorm{ \tfrac{\mathcal{J}^{\frac 34}}{\Sigma^{\beta}}  (\nbd^6 \Zbn,\nbd^6 \Abn) (\cdot,\s')}_{L^2_x}^2 {\rm d}\s' 
   \notag\\
   &\qquad
   - \tfrac{33 (\alpha\beta+\frac 12) + 250^2(1+\alpha)}{(1+\alpha)\eps} \int_0^{\s} \snorm{ \tfrac{\mathcal{J}^{\frac 34}(\Q \Jg)^{\frac 12}}{\Sigma^{\beta}}  (\nbd^6 \Zbn,\nbd^6 \Abn) (\cdot,\s')}_{L^2_x}^2 {\rm d}\s' 
      \notag\\
   &\qquad    
   - \Cn \eps^3 (\tfrac{4}{\kappa_0})^{2\beta} \mathsf{K}^2 \brak{\mathsf{B}_6}^2
  + \bigl( \tfrac{3(1+\alpha)}{5\eps}  -  \Cn \beta  \bigr) \tint  \tfrac{\mathcal{J}^{\frac 12} \Jg}{\Sigma^{2\beta}} \bigl( \tfrac 12 (\nbd^6 \Zbn)^2 +  (\nbd^6 \Abn)^2 \bigr)
   \label{eq:D5:An:Zn:2}
   \,,
 \end{align}
where $\Cn = \Cn(\alpha,\kappa_0,\Cdata)>0$ is a constant.

\subsection{The integral  $\mathfrak{I} ^{\AA_n}_6$}

In order to estimate the $\mathfrak{I} ^{\AA_n}_{6}$ term defined in \eqref{I6-Abn-imp}, which we note 
contains the seventh order derivative $\nbs_2 \nbd^6 (\Jg \Wbn) = \eps \nbs_2 \nbd^5 (\Q \p_\s + V \p_2)(\Jg\Wbn)$, we apply $\nbs_2\nbd^5$ to \eqref{eq:Jg:Wb:nn-s} and find that
\begin{align*} 
\nbs_2\nbd^6(\Jg\Wbn)
& = -  \alpha\eps \Sigma g^{-\frac 12} \Jg \nbs_2^2 \nbd^5  \Abn 
- \alpha \eps \Sigma g^{-\frac 12} \Jg \nbs_2 \jump{\nbd^5 ,\nbs_2}\Abn
\notag \\
& \qquad
-  \alpha \eps \nbs_2\nbd^5( \Sigma g^{-\frac 12} \Jg) \nbs_2  \Abn 
-  \alpha \eps \doublecom{\nbs_2\nbd^5, \Sigma g^{-\frac 12} \Jg,  \nbs_2  \Abn   }
+ \eps \nbs_2\nbd^5 \mathcal{G}^\nn_\Wb  \,.
\end{align*} 
Substitution of this identity into the integral $\mathfrak{I} ^{\AA_n}_{6}$  in \eqref{I6-Abn-imp} shows that 
\begin{align} 
\mathfrak{I} ^{\AA_n}_{6} 
&= 
-\alpha^2 \eps \tint \Sigma g^{-1} \jb\mathcal{J}^{\! \frac 32} \Jg    \nbs_2^2 \nbd^5  \Abn   \    \nbd^6\Abn 
-\alpha^2 \eps \tint \jb\Sigma g^{-1}  \mathcal{J}^{\! \frac 32} \Jg    \nbs_2 \jump{\nbd^5 ,\nbs_2}\Abn   \    \nbd^6\Abn 
 \notag \\
& \quad
-\alpha^2 \eps \tint  \jb g^{-\frac 12} \mathcal{J}^{\!\frac 32} \nbs_2\nbd^5( \Sigma g^{-\frac 12} \Jg) \nbs_2  \Abn  \   \nbd^6\Abn
-\alpha^2 \eps  \tint\jb g^{-\frac 12} \mathcal{J}^{\!\frac 32}   \doublecom{\nbs_2\nbd^5 , \Sigma g^{-\frac 12} \Jg,  \nbs_2  \Abn   }  \   \nbd^6\Abn
\notag \\
& \quad
+\alpha \eps \tint \jb g^{-\frac 12} \mathcal{J}^{\!\frac 32} \nbs_2\nbd^5 \mathcal{G}^\nn_\Wb    \   \nbd^6\Abn
\notag \\
&=: \mathfrak{I} ^{\AA_n}_{6,a} + \mathfrak{I} ^{\AA_n}_{6,b} + \mathfrak{I} ^{\AA_n}_{6,c}  + \mathfrak{I} ^{\AA_n}_{6,d}  + \mathfrak{I} ^{\AA_n}_{6,e} \,.
\label{eq:I:An:6:decomposition}
\end{align} 
The key terms are $\mathfrak{I} ^{\AA_n}_{6,a}$ and $\mathfrak{I} ^{\AA_n}_{6,e}$, because these terms involve objects with seven derivatives: $\nbs_2^2 \nbd^5 \Abn$, and respectively $\nbs_2 \nbd^5 \nbs_2 (\nbs_2 h)$. We deal with these over-differentiated terms first.

Using \eqref{eq:adjoints} and the fact that $\jump{\nbs_2 ,\nbd} = \eps \nbs_2 V \nbs_2$, we rewrite
\begin{align}
\mathfrak{I} ^{\AA_n}_{6,a}
&= 
\tfrac{\alpha^2}{2}  \eps^2 \tint  \Sigma \jb g^{-1} \mathcal{J}^{\frac 32} \Jg   (\Q \p_\s + V \p_2) \bigl( (\nbs_2 \nbd^5 \Abn)^2 \bigr)
+ \alpha^2  \eps^2 \tint  \Sigma \jb g^{-1} \mathcal{J}^{\frac 32}  \Jg  \nbs_2 V   (\nbs_2 \nbd^5 \Abn)^2
\notag\\
&\qquad
- \alpha^2  \eps  \tint \nbs_2 \nbd^5 \Abn \nbd^6 \Abn  (\Qr_2 - \nbs_2)  \bigl(\Sigma \jb g^{-1}  \mathcal{J}^{\frac 32}  \Jg   \bigr)  
+ \alpha^2 \eps^2 \int \Qb_2 \Sigma \jb g^{-1}  \mathcal{J}^{\frac 32} \Jg \nbs_2 \nbd^5 \Abn \nbd^6 \Abn \Bigr|_{\s}
\notag\\
&= 
\tfrac{\alpha^2}{2}  \eps^2 \snorm{\tfrac{\Sigma g^{-\frac 12}\mathcal{J}^{\frac 34} (\Q  \Jg)^{\frac 12}}{\Sigma^{\beta}}     \nbs_2 \nbd^5 \Abn(\cdot,\s)}_{L^2_x}^2 
- 
\tfrac{\alpha^2}{2}  \eps^2 \snorm{\tfrac{\Sigma g^{-\frac 12}\mathcal{J}^{\frac 34} (\Q  \Jg)^{\frac 12}}{\Sigma^{\beta}}     \nbs_2 \nbd^5 \Abn(\cdot,0)}_{L^2_x}^2 
\notag\\
&\qquad 
+ \tfrac{\alpha^2}{2}  \eps^2 \tint  \tfrac{1}{\Sigma^{2\beta}} \mathsf{G}_*   (\nbs_2 \nbd^5 \Abn)^2  
- \alpha^2  \eps  \tint \nbs_2 \nbd^5 \Abn \nbd^6 \Abn  (\Qr_2 - \nbs_2)  \bigl(\Sigma \jb g^{-1}  \mathcal{J}^{\frac 32}  \Jg   \bigr)  
\notag\\
&\qquad
+ \alpha^2 \eps^2 \int \Qb_2 \Sigma \jb g^{-1}  \mathcal{J}^{\frac 32} \Jg \nbs_2 \nbd^5 \Abn \nbd^6 \Abn \Bigr|_{\s}
\label{eq:mathfrak:G:emerges:first}
\end{align}
where  
\begin{align}
\mathsf{G}_*
&:=
\Sigma^{2}    g^{-1} \mathcal{J}^{\frac 32} \Jg \bigl(\nbs_2 V + V \Qr_2 -  \Qr_\s  \bigr)  
- \Sigma^{2\beta} (\Q \p_\s + V \p_2)\bigl( \Sigma \jb g^{-1} \mathcal{J}^{\frac 32} \Jg \bigr) 
\notag\\
&\geq - \Cn \bigl(\tfrac 1{\eps} + \brak{\beta}\bigr) \mathcal{J}^{\frac 32} \Jg   
+ \tfrac 32 \Sigma^{2} g^{-1}  \mathcal{J}^{\frac 12} \Jg \tfrac{\Q}{\eps}
+ \tfrac{1+\alpha}{2} \Sigma^{2}  g^{-1} \mathcal{J}^{\frac 32} \bigl( \tfrac{9}{10\eps} - \tfrac{13}{\eps} \Jg\bigr)
\notag\\
&\geq 
- \Cn \bigl(\tfrac 1{\eps} + \brak{\beta}\bigr) \mathcal{J}^{\frac 32} \Jg   
+ \tfrac 32 \Sigma^{2} g^{-1}  \mathcal{J}^{\frac 12} \Jg \tfrac{\Q}{\eps}
+ \tfrac{1+\alpha}{2} \Sigma^{2}  g^{-1} \mathcal{J}^{\frac 32}  \tfrac{9}{10\eps}  
\notag\\
&\geq 
-    \Cn \bigl(\tfrac 1{\eps} + \brak{\beta}\bigr) \mathcal{J}^{\frac 32} \Jg    
\,.
\label{eq:mathfrak:G:lower:bound}
\end{align}
In the second-to-last inequality above we have appealed to the bootstraps~\eqref{bootstraps}, the bounds \eqref{eq:Q:all:bbq},  and \eqref{eq:signed:Jg}, have assumed that $\beta \geq 1$, and have taken $\eps$ to be sufficiently small with respect to $\alpha,\kappa_0,\Cdata$ (but not with respect to $\beta$). As such, by also appealing to \eqref{eq:useful:junk:L2L2} with $f = \nbs_2 \Abn$ and $a = \frac 34, b = 0$, and with \eqref{eq:Jg:Abn:D5:improve:c}--\eqref{eq:Jg:Abn:D6:improve:d}, we have
\begin{align}
\tfrac{\alpha^2}{2}  \eps^2 \tint  \tfrac{1}{\Sigma^{2\beta}} \mathsf{G}_*   (\nbs_2 \nbd^5 \Abn)^2 
&\geq -   \eps^2 \Cn \bigl(\tfrac 1{\eps} + \brak{\beta}\bigr) (\tfrac{4}{\kappa_0})^{2\beta} \| \mathcal{J}^{\frac 34}  \nbd^5 ( \nbs_2 \Abn)\|_{L^2_{x,\s}}^2
\notag\\
&\geq - \Cn \eps^3 (1 + \Cn \eps \brak{\beta}) (\tfrac{4}{\kappa_0})^{2\beta} \mathsf{K}^2 \brak{\mathsf{B}_6}^2
\,.
 \label{eq:mathfrak:G:lower:bound:no:fuck}
\end{align}

In order to bound the fourth term on the right side of \eqref{eq:mathfrak:G:emerges:first}, we use \eqref{eq:useful:junk:L2L2} with $f = \nbs_2 \Abn$ and $a = \frac 34, b = 0$, in conjunction with \eqref{eq:Jg:Abn:D5:improve:c}--\eqref{eq:Jg:Abn:D6:improve:d}, to deduce 
\begin{align}
\alpha^2  \eps \left| \tint \nbs_2 \nbd^5 \Abn \nbd^6 \Abn  (\Qr_2 - \nbs_2)  \bigl(\Sigma \jb g^{-1}  \mathcal{J}^{\frac 32}  \Jg   \bigr) \right|
& \leq  \Cn \eps (\tfrac{4}{\kappa_0})^\beta \!\! \int_0^{\s} \snorm{ \tfrac{\mathcal{J}^{\frac 34}}{\Sigma^{\beta}}   \nbd^6 \Abn (\cdot,\s')}_{L^2_x} (\eps \mathsf{K} \brak{\mathsf{B}_6} + \eps^3 \mathsf{K}^2 \brak{\mathsf{B}_6}^2 ){\rm d}\s'
\notag\\
&\leq \Cn \int_0^{\s} \snorm{ \tfrac{\mathcal{J}^{\frac 34}}{\Sigma^{\beta}}   \nbd^6 \Abn (\cdot,\s')}_{L^2_x}^2 {\rm d}\s'
+ \Cn  \eps^3  (\tfrac{4}{\kappa_0})^{2\beta}  \mathsf{K}^2 \brak{\mathsf{B}_6}^2.
\label{eq:mathfrak:G:lower:bound:fuck}
\end{align}
Similarly,  to bound the fifth term on the right side of \eqref{eq:mathfrak:G:emerges:first}, we use \eqref{eq:useful:junk:LinftyL2} with $f = \nbs_2 \Abn$ and $a = \frac 34, b = \frac 12$, in conjunction with \eqref{eq:Jg:Abn:D6:improve:b}--\eqref{eq:Jg:Abn:D5:improve:c}, to deduce 
\begin{align}
 \alpha^2 \eps^2 \left| \int \Qb_2 \Sigma \jb g^{-1}  \mathcal{J}^{\frac 32} \Jg \nbs_2 \nbd^5 \Abn \nbd^6 \Abn \Bigr|_{\s} \right|
 &\leq \Cn \eps^3 (\tfrac{4}{\kappa_0})^\beta \snorm{ \tfrac{\mathcal{J}^{\frac 34} (\Q \Jg)^{\frac 12}}{\Sigma^\beta}  \nbd^6 \Abn(\cdot,\s)}_{L^2_x} (\eps^{\frac 12} \mathsf{K} \brak{\mathsf{B}_6} + \eps^{\frac 52} \mathsf{K}^2 \brak{\mathsf{B}_6}^2) 
 \notag\\
  &\leq \Cn \eps^2  \snorm{ \tfrac{\mathcal{J}^{\frac 34} (\Q \Jg)^{\frac 12}}{\Sigma^\beta}  \nbd^6 \Abn(\cdot,\s)}_{L^2_x}^2 + \Cn  \eps^5 (\tfrac{4}{\kappa_0})^{2\beta} \mathsf{K}^2 \brak{\mathsf{B}_6}^2
  \,.
  \label{eq:mathfrak:G:lower:bound:fuck:fuck}
\end{align}

Thus, by combining \eqref{eq:mathfrak:G:emerges:first}, \eqref{eq:mathfrak:G:lower:bound:no:fuck}, \eqref{eq:mathfrak:G:lower:bound:fuck}, and \eqref{eq:mathfrak:G:lower:bound:fuck:fuck},  we deduce 
\begin{align}
\mathfrak{I} ^{\AA_n}_{6,a}
&\geq 
\tfrac{\alpha^2}{2}  \eps^2 \snorm{\tfrac{\Sigma g^{-\frac 12}\mathcal{J}^{\frac 34} (\Q  \Jg)^{\frac 12}}{\Sigma^{\beta}}     \nbs_2 \nbd^5 \Abn(\cdot,\s)}_{L^2_x}^2 
- 
\tfrac{\alpha^2}{2}  \eps^2 \snorm{\tfrac{\Sigma g^{-\frac 12}\mathcal{J}^{\frac 34} (\Q  \Jg)^{\frac 12}}{\Sigma^{\beta}}     \nbs_2 \nbd^5 \Abn(\cdot,0)}_{L^2_x}^2 
\notag\\
&\qquad 
- 
\Cn \eps^2  \snorm{ \tfrac{\mathcal{J}^{\frac 34} (\Q \Jg)^{\frac 12}}{\Sigma^\beta}  \nbd^6 \Abn(\cdot,\s)}_{L^2_x}^2
-
\Cn \int_0^{\s} \snorm{ \tfrac{\mathcal{J}^{\frac 34}}{\Sigma^{\beta}}   \nbd^6 \Abn (\cdot,\s')}_{L^2_x}^2 {\rm d}\s'
- \Cn \eps^3 (1 + \Cn \eps \brak{\beta}) (\tfrac{4}{\kappa_0})^{2\beta} \mathsf{K}^2 \brak{\mathsf{B}_6}^2
\notag\\
&
\geq 
- 
\Cn \eps^2  \snorm{ \tfrac{\mathcal{J}^{\frac 34} (\Q \Jg)^{\frac 12}}{\Sigma^\beta}  \nbd^6 \Abn(\cdot,\s)}_{L^2_x}^2
-
\Cn  \int_0^{\s} \snorm{ \tfrac{\mathcal{J}^{\frac 34}}{\Sigma^{\beta}}   \nbd^6 \Abn (\cdot,\s')}_{L^2_x}^2 {\rm d}\s'
\notag\\
&\qquad 
- \Cn \eps^3 (1 + \Cn \eps \brak{\beta}) (\tfrac{4}{\kappa_0})^{2\beta} \mathsf{K}^2 \brak{\mathsf{B}_6}^2
\,.
\label{eq:mathfrak:G:has:gone}
\end{align}
In the second inequality in~\eqref{eq:mathfrak:G:has:gone} we have appealed to \eqref{table:derivatives}, the bootstraps~\eqref{bootstraps}, and to \eqref{eq:useful:junk:LinftyL2} with $f=\nbs_2 \Abn$, $a=b=0$.

Next, we turn to the other delicate term in \eqref{eq:I:An:6:decomposition}, namely $\mathfrak{I} ^{\AA_n}_{6,e}$. Recalling \eqref{Gwn}, \eqref{eq:Jg:Wb:nn},  and   \eqref{p2h-evo-s}, we may write 
\begin{align} 
&\nbs_2 \nbd^5  \mathcal{G}^\nn_\Wb 
\notag\\
&=
\eps \nbs_2 \nbd^4(\Q \p_\s +   V\p_2) \mathcal{G}^\nn_\Wb  
\notag\\
&= 
\eps \tfrac{\alpha}{2} \Sigma g^{-\frac 32}  \bigl(   \Jg\Wbn + \Jg \Zbn - 2\Jg\Abt \bigr)\nbs_2 \nbd^4 \nbs_2 \bigl( g ( \tfrac{1+\alpha}{2} \Wbt + \tfrac{1-\alpha}{2} \Zbt ) \bigr)
\notag\\
&  
- \eps \tfrac{\alpha}{2} \bigl( \nbs_2^2 h \Sigma g^{-\frac 32}  - \Abt\bigr) \, \nbs_2 \nbd^4  
\Bigl((\Jg\Wbn) \bigl(\tfrac{\alpha}{2} \Abt  - \tfrac{\alpha}{2} \Sigma g^{-\frac 32} \nbs_2^2 h \bigr)  
+ \alpha \Sigma g^{-\frac 12} \Jg \nbs_2 \Abn
- \tfrac{\alpha}{2} \bigl( \Abt + \Sigma g^{-\frac 32} \nbs_2^2 h  \bigr) \Jg \Zbn 
\notag\\
&  \qquad  \qquad\qquad\qquad\qquad\qquad\quad 
+  \bigl(\tfrac{3+\alpha}{2} \Wbt + \tfrac{1-\alpha}{2} \Zbt\bigr) \Jg\Abn
+ \bigl(\tfrac{1+\alpha}{2} \Wbt + \tfrac{1-\alpha}{2} \Zbt\bigr)\Jg \Wbt  
+ \alpha \Sigma g^{-\frac 32} \nbs_2^2 h  \Jg \Abt \Bigr)
\notag\\
&  
- \eps  \tfrac{\alpha}{2} \Sigma g^{-\frac 32}  \bigl(   \Jg\Wbn + \Jg \Zbn - 2\Jg\Abt \bigr)\nbs_2 \nbd^4 (\nbs_2 V  \nbs_2 h )
-   \tfrac{\alpha}{2} \bigl( \nbs_2^2 h \Sigma g^{-\frac 32}  - \Abt\bigr) \, \nbs_2 \nbd^5 (   \Jg \Zbn   ) 
 \notag\\
 &  
 + \alpha  \nbs_2^2 h \, \nbs_2 \nbd^5 \Bigl( \Sigma g^{-\frac 32}  \bigl( \Jg \Zbn -  \Jg\Abt \bigr)  \Bigr)
 + \tfrac{\alpha}{2} \doublecom{\nbs_2 \nbd^5, \Sigma g^{-\frac 32}  \bigl(   \Jg\Wbn + \Jg \Zbn - 2\Jg\Abt \bigr), \nbs_2^2 h}
 \notag\\
 &  - \nbs_2 \nbd^5 \Bigl(    \bigl(\tfrac{3+\alpha}{2} \Wbt + \tfrac{1-\alpha}{2} \Zbt\bigr) \Jg\Abn
+ \bigl(\tfrac{1+\alpha}{2} \Wbt + \tfrac{1-\alpha}{2} \Zbt\bigr)\Jg \Wbt  \Bigr) 
 \notag\\
 &  
+ \tfrac{\alpha}{2} \nbs_2^2 h \jump{ \nbs_2 \nbd^5, \Sigma g^{-\frac 32}}  \bigl(   \Jg\Wbn - \Jg \Zbn  \bigr) 
- \tfrac{\alpha}{2} \jump{\nbs_2 \nbd^5 , \Abt} \bigl( \Jg\Wbn - \Jg \Zbn  \bigr)
\,.
\label{eq:WTF:is:this}
\end{align} 
The precise form of the above expression is not relevant. Its important features are: the term which previously contained too many derivatives, namely $\nbs_2 \nbd^5 \nbs_2 (\nbs_2 h)$, has been rewritten in terms of factors with at most six derivatives on them. In the above expression, six derivatives land on either the geometry (these terms are bounded due to Proposition~\ref{prop:geometry}), or on $\Zbn$ and $\Abn$ (these terms we are currently writing evolution equations for, or they contain at least on $\nbs_1$ or $\nbs_2$, in which case they are bounded via \eqref{eq:very:useful:junk}, \eqref{eq:Jg:Abn:D5:improve}, \eqref{eq:Jg:Zbn:D5:improve}), or on the tangential components $(\Wbt,\Zbt,\Abt)$ (which require one-less power of $\Jg$, and are bounded using \eqref{bootstraps-Dnorm:6} and \eqref{bootstraps-Dnorm:5}). 
Using \eqref{eq:WTF:is:this}, and by appealing to the bounds \eqref{bootstraps}, \eqref{geometry-bounds-new}, \eqref{eq:h_2:D2-bound}, \eqref{eq:Jg:Abn:D5:improve},  \eqref{eq:Jg:Zbn:D5:improve}, \eqref{eq:D5:JgWbn}, \eqref{eq:very:useful:junk}, \eqref{eq:Lynch:1}, \eqref{eq:Lynch:2}, we may thus bound the $\mathfrak{I} ^{\AA_n}_{6,e}$ integral  in \eqref{eq:I:An:6:decomposition} as 
\begin{subequations}
\begin{align}
\sabs{\mathfrak{I} ^{\AA_n}_{6,e}}
&\leq 2\alpha\eps \kappa_0 (\tfrac{4}{\kappa_0})^\beta 
\int_0^{\s} 
\snorm{ \mathcal{J}^{\frac 34}  \nbs_2 \nbd^5  \mathcal{G}^\nn_\Wb (\cdot,\s')}_{L^2_{x}} 
\snorm{ \tfrac{\mathcal{J}^{\frac 34}}{\Sigma^{\beta}}   \nbd^6 \Abn (\cdot,\s')}_{L^2_x} 
{\rm d}\s'
\notag\\
&\leq  \Cn \eps (\tfrac{4}{\kappa_0})^\beta \mathsf{K} \brak{\mathsf{B}_6} 
\int_0^{\s} \snorm{ \tfrac{\mathcal{J}^{\frac 34}}{\Sigma^{\beta}}   \nbd^6 \Abn (\cdot,\s')}_{L^2_x} 
{\rm d}\s' 
\notag\\
&\leq \Cn  \eps^3 (\tfrac{4}{\kappa_0})^{2\beta} \mathsf{K}^2 \brak{\mathsf{B}_6}^2 
+ \tfrac{\alpha}{5\eps}  \int_0^{\s} \snorm{ \tfrac{\mathcal{J}^{\frac 34}}{\Sigma^{\beta}}   \nbd^6 \Abn (\cdot,\s')}_{L^2_x}^2 
{\rm d}\s' 
\, .
\end{align}
Here we single out the term $\tfrac{\alpha}{2}  \Jg\Wbn \mathcal{J}^{\frac 34} \nbs_2 \nbd^5 \Abt  $ arising in $\tfrac{\alpha}{2} \mathcal{J}^{\frac 34} \jump{\nbs_2 \nbd^5, \Abt} \bigl( \Jg\Wbn - \Jg \Zbn  \bigr)$, as the only term responsible for the factor of $\mathsf{K} \brak{\mathsf{B}_6}$ in the bound for $\nbs_2 \nbd^5  \mathcal{G}^\nn_\Wb $; one may verify that all other terms in \eqref{eq:WTF:is:this} contribute at most a factor of $\eps \mathsf{K} \brak{\mathsf{B}_6}$. 

In a similar fashion, we may bound the remaining (easy) terms $\mathfrak{I} ^{\AA_n}_{6,b}, \mathfrak{I} ^{\AA_n}_{6,c}, \mathfrak{I} ^{\AA_n}_{6,d}$  in \eqref{eq:I:An:6:decomposition} as follows:
\begin{align}
\sabs{\mathfrak{I} ^{\AA_n}_{6,b}} 
+
\sabs{\mathfrak{I} ^{\AA_n}_{6,c}} 
+
\sabs{\mathfrak{I} ^{\AA_n}_{6,d}} 
&\leq \Cn \eps^2 (\tfrac{4}{\kappa_0})^\beta \mathsf{K} \brak{\mathsf{B}_6} 
\int_0^{\s} \snorm{ \tfrac{\mathcal{J}^{\frac 34}}{\Sigma^{\beta}}   \nbd^6 \Abn (\cdot,\s')}_{L^2_x} 
{\rm d}\s' 
\notag\\
&\leq  \Cn  \eps^3 (\tfrac{4}{\kappa_0})^{2\beta} \mathsf{K}^2 \brak{\mathsf{B}_6}^2 
+ \Cn \eps \int_0^{\s} \snorm{ \tfrac{\mathcal{J}^{\frac 34}}{\Sigma^{\beta}}   \nbd^6 \Abn (\cdot,\s')}_{L^2_x}^2 
{\rm d}\s'  
\,,
\end{align}
\label{eq:junk:is:king}
\end{subequations}
where the implicit constant depends only on $\alpha,\kappa_0$, and $\Cdata$. 

In summary, from the \eqref{eq:I:An:6:decomposition} and the bounds \eqref{eq:mathfrak:G:has:gone} and \eqref{eq:junk:is:king} we deduce the lower bound
\begin{align}
\mathfrak{I} ^{\AA_n}_{6}
&\geq  
- 
\Cn \eps^2  \snorm{ \tfrac{\mathcal{J}^{\frac 34} (\Q \Jg)^{\frac 12}}{\Sigma^\beta}  \nbd^6 \Abn(\cdot,\s)}_{L^2_x}^2
-
\tfrac{2\alpha}{5\eps}  \int_0^{\s} \snorm{ \tfrac{\mathcal{J}^{\frac 34}}{\Sigma^{\beta}}   \nbd^6 \Abn (\cdot,\s')}_{L^2_x}^2 {\rm d}\s'
\notag\\
&\qquad 
- \Cn \eps^3 (1 +  \eps \brak{\beta}) (\tfrac{4}{\kappa_0})^{2\beta} \mathsf{K}^2 \brak{\mathsf{B}_6}^2
\,.
\label{eq:mathfrak:G:has:gone:2}
\end{align}

\subsection{The forcing and commutator terms}
\label{sec:An:Zn:improve:forcing:commutators}

Returning to \eqref{D_sD^5-L2} and the decompositions~\eqref{IntZbn}--\eqref{IntAbn}, it remains to estimate the forcing and commutator terms $ \mathfrak{I} ^{\ZZ_n}_5$ and $ \mathfrak{I} ^{\AA_n}_5$, defined in \eqref{I5-Zbn-imp} and respectively \eqref{I5-Abn-imp}. 

For the contributions from the commutator terms defined in \eqref{D_sD^5:Jg:Zn:An-commutators}, by appealing to Lemmas~~\ref{lem:anisotropic:sobolev}, \ref{lem:time:interpolation} and~\ref{lem:comm:tangent}, to the boostraps~\eqref{bootstraps}, to the bounds on the geometry in Proposition~\ref{prop:geometry}, the improved $\Abn$ estimates~\eqref{eq:Jg:Abn:D5:improve}, to the improved $\Zbn$ estimates~\eqref{eq:Jg:Zbn:D5:improve}, and to the bounds \eqref{eq:very:useful:junk}, we obtain 
\begin{subequations}
\label{eq:useful:junk:8}
\begin{align}
&\snorm{\tfrac{\mathcal{J}^{\frac 34}}{\Sigma^{\beta-1}} \mathfrak{C}^\nn_\Zb }_{L^2_{x,\s}} 
\notag\\
&\leq \tfrac{6}{\eps} \|\Sigma \nbd (\Jg \Sigma^{-1}) \|_{L^\infty_{x,\s}} \snorm{\tfrac{\mathcal{J}^{\frac 34}}{\Sigma^{\beta}} \nbd^6 \Zbn}_{L^2_{x,\s}}  
+ \tfrac{\Cn}{\eps} (\tfrac{4}{\kappa_0})^\beta \displaystyle{\sum}_{k=0}^4 \snorm{ \nbd^{6-k}(\Jg \Sigma^{-1}) \nbd^{k+1}\Zbn}_{L^2_{x,\s}} 
+ \Cn \eps (\tfrac{4}{\kappa_0})^\beta \mathsf{K} \brak{\mathsf{B}_6} 
\notag\\
&\leq \tfrac{6 (11+\alpha)}{\eps} \snorm{\tfrac{\mathcal{J}^{\frac 34}}{\Sigma^{\beta}} \nbd^6 \Zbn}_{L^2_{x,\s}} 
+  \tfrac{\Cn}{\eps} (\tfrac{4}{\kappa_0})^\beta \displaystyle{\sum}_{k=0}^4 \snorm{\nbd^{4-k}  \nbd^2(\Jg \Sigma^{-1})}_{L^{\frac{8}{4-k}}_{x,\s}} \snorm{\nbd^{k} \nbd \Zbn}_{L^{\frac{8}{k}}_{x,\s}} 
+ \Cn \eps (\tfrac{4}{\kappa_0})^\beta \mathsf{K} \brak{\mathsf{B}_6}    
\notag\\
&\leq \tfrac{6 (11+\alpha)}{\eps} \snorm{\tfrac{\mathcal{J}^{\frac 34}}{\Sigma^{\beta}} \nbd^6 \Zbn}_{L^2_{x,\s}} 
+ \Cn    \mathsf{K} (\tfrac{4}{\kappa_0})^\beta  \brak{\mathsf{B}_6} 
\,, 
\end{align}
and in a similar fashion
\begin{equation}
\snorm{\tfrac{\mathcal{J}^{\frac 34}}{\Sigma^{\beta-1}} \mathfrak{C}^\nn_\Ab }_{L^2_{x,\s}}
\leq 
\tfrac{6 (11+\alpha)}{\eps} \snorm{\tfrac{\mathcal{J}^{\frac 34}}{\Sigma^{\beta}} \nbd^6 \Abn}_{L^2_{x,\s}} 
+ \Cn   \mathsf{K}   (\tfrac{4}{\kappa_0})^\beta  \brak{\mathsf{B}_6}
\,. 
\end{equation}
\end{subequations}

For the contributions from $\nbd^6$ acting on the forcing terms defined in \eqref{Gzn}--\eqref{Gan}, we need to be quite careful, especially the contributions from $\nbd^6 (\Jg \Wbn)$, $\nbd^6 \nbs_2^2 h$, and $\nbd^6 \nbs_2 \Jg$. For these terms, we use that identities 
\eqref{eq:Jg:Wb:nn-s} and \eqref{Gwn}, along with the bootstraps~\eqref{bootstraps}, the bounds on the geometry in Proposition~\ref{prop:geometry}, the improved $\Abn$ estimates~\eqref{eq:Jg:Abn:D5:improve}, to the improved $\Zbn$ estimates~\eqref{eq:Jg:Zbn:D5:improve}, to the bounds~\eqref{eq:very:useful:junk}, and the Moser inequality~\eqref{eq:Lynch:1}, imply
\begin{subequations}
\label{eq:junk:loggers:1}
\begin{align}
 \snorm{\mathcal{J}^{\frac 34} \nbd^6 (\Jg \Wbn) }_{L^2_{x,\s}}
 &\leq 
 \eps \snorm{\mathcal{J}^{\frac 34}  \nbd^5 (\Q\p_\s + V \p_2) (\Jg \Wbn) }_{L^2_{x,\s}} 
 \notag\\
 &\leq 
 \alpha \eps \|\Sigma g^{-\frac 12} \Jg\|_{L^\infty_{x,\s}} \snorm{\mathcal{J}^{\frac 34}   \nbd^5 \nbs_2 \Abn}_{L^2_{x,\s}} 
 +  \alpha \eps \snorm{ \jump{\nbd^5, \Sigma g^{-\frac 12} \Jg} \nbs_2 \Abn}_{L^2_{x,\s}} 
 + \eps \snorm{\nbd^5 \mathcal{G}^\nn_\Wb}_{L^2_{x,\s}}
 \notag\\
 &\les \eps \mathsf{K} \brak{\mathsf{B}_6} \,.
\end{align}
By additionally appealing to the $\nbs_2$ differentiated version of \eqref{p2h-evo-s}, which gives $(\Q\p_\s+V\p_2) \nbs_2^2 h  = - \nbs_2 V \nbs_2^2 h +  \nbs_2 h \nbs_2^2 h  \bigl(\tfrac{1+ \alpha }{2} \Wbt + \tfrac{1- \alpha }{2} \Zbt \bigr) +  g  \bigl(\tfrac{1+ \alpha }{2} \nbs_2 \Wbt + \tfrac{1- \alpha }{2} \nbs_2 \Zbt \bigr)$, we may similarly obtain
\begin{equation}
 \snorm{\mathcal{J}^{\frac 34} \nbd^6 \nbs_2^2 h }_{L^2_{x,\s}}  
 =
 \eps
 \snorm{\mathcal{J}^{\frac 34} \nbd^5 (\Q \p_\s + V \p_2) \nbs_2^2 h }_{L^2_{x,\s}}  
 \les \eps^2 \mathsf{K} \brak{\mathsf{B}_6}
 \,,
\end{equation}
and by using the $\nbs_2$ differentiated version of  \eqref{Jg-evo-s}, which gives $(\Q\p_\s+V\p_2) \nbs_2 \Jg  = - \nbs_2 V \nbs_2 \Jg +  \tfrac{1+ \alpha }{2} \nbs_2(\Jg\Wbn) + \tfrac{1- \alpha }{2} \nbs_2(\Jg\Zbn)$,  we have
\begin{equation}
 \snorm{\mathcal{J}^{\frac 34} \nbd^6 \nbs_2 \Jg }_{L^2_{x,\s}}  
 = \eps 
  \snorm{\mathcal{J}^{\frac 34} \nbd^5 (\Q \p_\s + V \p_2)\nbs_2 \Jg }_{L^2_{x,\s}} 
 \les \eps  \mathsf{K} \brak{\mathsf{B}_6}
 \,.
\end{equation}
\end{subequations}
Using these above estimates in \eqref{eq:junk:loggers:1}, the bound \eqref{eq:very:useful:junk}, the product and commutator lemmas bounds in~\ref{eq:Lynch:1} and~\ref{eq:Lynch:2}, we may derive from \eqref{Gzn}--\eqref{Gan} the bounds
\begin{subequations}
\label{eq:useful:junk:9}
\begin{align}
\snorm{\tfrac{\mathcal{J}^{\frac 34}}{\Sigma^{\beta-1}}\nbd^6\mathcal{G}^\nn_\Zb}_{L^2_{x,\s}}
&\leq \tfrac{1+\alpha}{\eps}   \snorm{\tfrac{\mathcal{J}^{\frac 34}}{\Sigma^{\beta}} \nbd^6 \Zbn}_{L^2_{x,\s}} 
+  \Cn   \mathsf{K}   (\tfrac{4}{\kappa_0})^\beta  \brak{\mathsf{B}_6} \,, \\
\snorm{\tfrac{\mathcal{J}^{\frac 34}}{\Sigma^{\beta-1}} \nbd^6\mathcal{G}^\nn_\Ab}_{L^2_{x,\s}}
&\leq \tfrac{1}{\eps}   \snorm{\tfrac{\mathcal{J}^{\frac 34}}{\Sigma^{\beta}} \nbd^6 \Abn}_{L^2_{x,\s}} 
+  \Cn   \mathsf{K}   (\tfrac{4}{\kappa_0})^\beta  \brak{\mathsf{B}_6}\,.
\end{align}
\end{subequations}
Inserting the bounds \eqref{eq:useful:junk:8} and \eqref{eq:useful:junk:9}, into definitions \eqref{I5-Zbn-imp} and \eqref{I5-Abn-imp}, we arrive at 
\begin{equation}
\sabs{\mathfrak{I} ^{\ZZ_n}_5} + \sabs{\mathfrak{I} ^{\AA_n}_5}
\leq  
\tfrac{100(1+\alpha)}{\eps}   \int_0^{\s} \snorm{ \tfrac{\mathcal{J}^{\frac 34}}{\Sigma^{\beta}}  (\nbd^6 \Zbn, \nbd^6 \Abn) (\cdot,\s')}_{L^2_x}^2 
{\rm d}\s' 
+\Cn   \eps    \mathsf{K}^2   (\tfrac{4}{\kappa_0})^{2\beta}  \brak{\mathsf{B}_6}^2 \,,
\label{eq:useful:junk:7}
\end{equation}
with $\Cn = \Cn(\alpha,\kappa_0,\Cdata)>0$.

\subsection{Conclusion of the proof of Proposition~\ref{prop:madman}}
\label{sec:proof:eq:madman}
At this stage we gather the identity \eqref{D_sD^5-L2} and the lower bounds \eqref{eq:D5:An:Zn:2}, \eqref{eq:mathfrak:G:has:gone:2}, and \eqref{eq:useful:junk:7}, to derive that 
\begin{align}
0
&\geq
\bigl(\tfrac 12 - \Cn \eps \bigr) \snorm{ \tfrac{\mathcal{J}^{\frac 34} (\Q \Jg)^{\frac 12}}{\Sigma^\beta} ( \nbd^6 \Zbn, \nbd^6 \Abn)(\cdot,\s)}_{L^2_x}^2 
- \snorm{ \tfrac{\mathcal{J}^{\frac 34} (\Q \Jg)^{\frac 12}}{\Sigma^\beta}  ( \nbd^6 \Zbn, \nbd^6 \Abn)(\cdot,0)}_{L^2_x}^2 
\notag\\
&\qquad 
+\Bigl( \tfrac{4\alpha(\beta-1)}{5\eps} -   \tfrac{100(1+\alpha)}{\eps}  \Bigr)
\int_0^{\s} \snorm{ \tfrac{\mathcal{J}^{\frac 34}}{\Sigma^{\beta}}  (\nbd^6 \Zbn,\nbd^6 \Abn) (\cdot,\s')}_{L^2_x}^2 {\rm d}\s' 
\notag\\
&\qquad 
- \tfrac{33 (\alpha\beta+\frac 12)+250^2(1+\alpha)}{(1+\alpha)\eps} \int_0^{\s} \snorm{ \tfrac{\mathcal{J}^{\frac 34}(\Q \Jg)^{\frac 12}}{\Sigma^{\beta}}  (\nbd^6 \Zbn,\nbd^6 \Abn) (\cdot,\s')}_{L^2_x}^2 {\rm d}\s' 
\notag\\
&\qquad    
- \Cn \eps  (1 +   \eps^3 \brak{\beta}) (\tfrac{4}{\kappa_0})^{2\beta} \mathsf{K}^2 \brak{\mathsf{B}_6}^2
+\bigl( \tfrac{3(1+\alpha)}{5\eps}  -  \Cn \beta  \bigr)  \tint  \tfrac{\mathcal{J}^{\frac 12} \Jg}{\Sigma^{2\beta}} \bigl( \tfrac 12 (\nbd^6 \Zbn)^2 +  (\nbd^6 \Abn)^2 \bigr)
\,.
\label{eq:Zbn:Abn:6:time:final:bound:0}
\end{align}
In view of the above estimate, we first choose  $\beta$ such that 
\begin{equation}
 \tfrac{4\alpha(\beta-1)}{5 } =  100(1+\alpha) 
 \qquad 
 \Leftrightarrow
 \qquad 
 \beta = \beta_\alpha := \tfrac{500(1+\alpha)}{4\alpha} + 1 
\,,
 \label{eq:Zbn:Abn:6:time:beta:choice}
\end{equation}
and then, for this fixed $\beta_\alpha$, we choose $\eps$ to be sufficiently small, in terms of $\alpha,\kappa_0,\Cdata$, such that 
\begin{equation*}
\tfrac{2(1+\alpha)}{5\eps} \geq \Cn_{\eqref{eq:Zbn:Abn:6:time:final:bound:0}} \beta_\alpha,
\qquad\mbox{and}\qquad
1+ \Cn_{\eqref{eq:Zbn:Abn:6:time:final:bound:0}} \eps^3 \brak{\beta_\alpha} \leq 2,
\qquad\mbox{and}\qquad
\Cn_{\eqref{eq:Zbn:Abn:6:time:final:bound:0}} \eps \leq \tfrac 14.
\end{equation*}
With this choice of $\beta = \beta_\alpha$ and $\eps$ sufficiently small, we return to \eqref{eq:Zbn:Abn:6:time:final:bound:0}, appeal to \eqref{table:derivatives} and \eqref{eq:expand:nbs:k} in order to bound the initial data term (recall that at $\s=0$ we have $\mathcal{J} \equiv \Jg \equiv 1$), and deduce (since $\mathsf{B}_6 \geq \Cdata$ by~\eqref{eq:B6:choice:1})
\begin{align}
&\tfrac 14 \snorm{ \tfrac{\mathcal{J}^{\frac 34} (\Q \Jg)^{\frac 12}}{\Sigma^\beta} ( \nbd^6 \Zbn, \nbd^6 \Abn)(\cdot,\s)}_{L^2_x}^2 
+ \tfrac{ 1+\alpha}{5\eps} \tint  \tfrac{\mathcal{J}^{\frac 12} \Jg}{\Sigma^{2\beta}} \bigl( \tfrac 12 (\nbd^6 \Zbn)^2 +  (\nbd^6 \Abn)^2 \bigr) 
\notag\\
&\leq \Cn \eps  (\tfrac{4}{\kappa_0})^{2\beta_\alpha} \mathsf{K}^2 \brak{\mathsf{B}_6}^2
+ \tfrac{33 (\alpha\beta_\alpha+\frac 12)+250^2(1+\alpha)}{(1+\alpha)\eps} \int_0^{\s} \snorm{ \tfrac{\mathcal{J}^{\frac 34}(\Q \Jg)^{\frac 12}}{\Sigma^{\beta}}  (\nbd^6 \Zbn,\nbd^6 \Abn) (\cdot,\s')}_{L^2_x}^2 {\rm d}\s' 
\,,
\label{eq:useful:junk:6}
\end{align}
for some constant $\Cn$ which depends only on $\alpha,\kappa_0$, and $\Cdata$.
To conclude the proof of  \eqref{eq:madman}, we apply the Gr\"onwall inequality to \eqref{eq:useful:junk:6} for $\s \in [0,\eps]$, then multiply the resulting bound by $\kappa_0^\beta$, and appeal to \eqref{bs-Sigma} and \eqref{eq:Q:bbq}.


\section{The sixth order energy estimates for the normal components}
\label{sec:sixth:order:energy}
 
\subsection{The   sixth-order differentiated equations for $(\Jg\Wbn, \Jg \Zbn, \Jg\Abn)$}

From~\eqref{euler-Wn},~\eqref{euler-Zn},~\eqref{euler-An} we write the sixth-order differentiated $(\Jg\Wbn, \Jg\Zbn, \Jg \Abn)$ equations in $(x,\s)$ coordinates as
\begin{subequations} 
\label{energy-WZA-s}
\begin{align} 
&
\tfrac{1}{\Sigma} (\Q\p_\s +V\p_2) \nbs^6(\Jg\Wbn  )
+ \alpha    g^{- {\frac{1}{2}} } \nbs_2 \nbs^6 (\Jg\Abn)
- \alpha    g^{- {\frac{1}{2}} }  \Abn \nbs_2\nbs^6 \Jg
 \notag \\
& \qquad
- \tfrac{\alpha}{2} g^{- {\frac{1}{2}} } (\Jg\Wbn + \Jg\Zbn - 2\Jg\Abt)\nbs_2 \nbs^6 \tt\cdo\nn
=  \nbs^6\Fwn + \mathcal{R}_\Wb^\nn + \mathcal{C}_\Wb^\nn
\,,  \label{energy-Wn-s}
 \\
  &
\tfrac{1}{\Sigma} \Jg(\Q\p_\s +V\p_2) \nbs^6(\Jg\Zbn) 
-\bubu{  \tfrac{\alpha}{\Sigma} ( \Jg \Wbn - \Jg \Zbn)\nbs^6(\Jg\Zbn)}
- \alpha \Jg   g^{- {\frac{1}{2}} } \nbs_2 \nbs^6 (\Jg\Abn)
+ \alpha \Jg   g^{- {\frac{1}{2}} }  \Abn \nbs_2\nbs^6 \Jg
\notag \\
& \qquad
+ \tfrac{\alpha}{2} \Jg g^{- {\frac{1}{2}} }(\Jg\Wbn + \Jg\Zbn - 2\Jg\Abt) \nbs_2 \nbs^6 \tt\cdo\nn
- \tfrac{2 \alpha}{\eps} \nbs_1 \nbs^6 (\Jg\Zbn)
 - \tfrac{2 \alpha}{\eps} \Jg (\Abn+\Zbt)  \nbs_1 \nbs^6 \tt \cdot \nn
\notag \\
& \qquad
 + \tfrac{2 \alpha}{\eps} \Zbn  \nbs_1 \nbs^6 \Jg
 + 2 \alpha \Jg g^{- {\frac{1}{2}} } \nbs_2 h  \nbs_2\nbs^6 (\Jg\Zbn)
 + 2 \alpha  \nbs_2\nbs^6\tt \cdot \nn \Jgt g^{- {\frac{1}{2}} } \nbs_2 h  (\Abn+\Zbt) 
 \notag \\
& \qquad
 - 2 \alpha \Jg g^{- {\frac{1}{2}} } \nbs_2h  \Zbn \nbs_2\nbs^6\Jg 
 = \nbs^6 \Fzn  + \mathcal{R}_\Zb^\nn + \mathcal{C}_\Zb^\nn\,,     \label{energy-Zn-s}  \\
 &
\tfrac{1}{\Sigma} \Jg(\Q\p_\s +V\p_2) \nbs^6(\Jg\Abn) 
\bubu{ - \tfrac{\alpha }{2\Sigma} (  \Jg\Wbn - \Jg \Zbn)\nbs^6(\Jg\Abn)}
+ \alpha   g^{- {\frac{1}{2}} }  \Jg  \nbs_2\nbs^6(\Jg\Sbn)
 \notag \\
& \qquad
- \alpha   g^{- {\frac{1}{2}} }  (\Jg\Sbn) \nbs_2\nbs^6 \Jg
+ \alpha g^{- {\frac{1}{2}} } \Jgt  \Sbt \nbs_2\nbs^6 \tt \cdot \nn
+ \tfrac{\alpha}{2\eps}  \nbs_1\nbs^6\tt \cdot\nn (\Jg\Wbn +\Jg\Zbn - 2\Jg \Abt)
 \notag \\
& \qquad
- \tfrac{\alpha}{\eps}\nbs_1\nbs^6 (\Jg\Abn)
+ \tfrac{\alpha}{\eps} \nbs_1\nbs^6 \Jg \Abn
+ \alpha \Jg g^{- {\frac{1}{2}} } \nbs_2 h \nbs_2\nbs^6 (\Jg\Abn)
- \alpha \Jg g^{- {\frac{1}{2}} } \nbs_2 h \Abn \nbs_2\nbs^6\Jg 
\notag \\
& \qquad 
- \tfrac{\alpha}{2} \nbs_2 \nbs^6 \tt \cdot\nn \Jg g^{- {\frac{1}{2}} } \nbs_2 h (\Jg\Wbn +\Jg\Zbn - 2\Jg \Abt)
 = \nbs^6 \Fan  + \mathcal{R}_\Ab^\nn + \mathcal{C}_\Ab^\nn 
  \,,   \label{energy-An-s}  
  \end{align} 
\end{subequations} 
where the forcing terms $\Fwn, \Fzn, \Fan$ 
are defined by~\eqref{forcing-nt}, while the remainder and commutator terms are given by 
\begin{subequations} 
\label{Cw-Rw-comm}
\begin{align} 
\mathcal{R} _\Wb^\nn & 
=\Sigma ^{-1} \nbs^6V \nbs_2 (\Jg\Wbn) 
+ \nbs^6 \Sigma ^{-1} (\Q\p_\s+V\p_2) (\Jg\Wbn) 
+ \alpha \nbs^6 g^{- {\frac{1}{2}} } \nbs_2 (\Jg\Abn)
\notag \\
& \qquad
- \alpha \nbs^6( g^{- {\frac{1}{2}} } \Abn) \nbs_2 \Jg
- \tfrac{\alpha}{2}  \nbs^6\big( g^{- {\frac{1}{2}} } (\Jg\Wbn + \Jg\Zbn - 2\Jg\Abt) \nn_i \big)\nbs_2  \tt^i
\,, \\
\mathcal{C}_\Wb^\nn 
&= \doublecom{\nbs^6, \Sigma^{-1}  , (\Q\p_\s+V\p_2)  (\Jg\Wbn)}
+ \Sigma^{-1} \doublecom{\nbs^6, V, \nbs_2(\Jg\Wbn)}
+ \alpha \doublecom{ \nbs^6 , g^{- {\frac{1}{2}} },  \nbs_2 (\Jg\Abn)}
\notag \\
& \qquad
- \alpha  \doublecom{\nbs^6,  g^{- {\frac{1}{2}} } \Abn,  \nbs_2 \Jg}
- \tfrac{\alpha}{2} \doublecom{  \nbs^6,  g^{- {\frac{1}{2}} } (\Jg\Wbn + \Jg\Zbn - 2\Jg\Abt) \nn_i , \nbs_2  \tt^i}
 \,, 
\label{com-Cwn}
\end{align} 
\end{subequations} 
for \eqref{energy-Wn-s}, and
\begin{subequations} 
\label{Cz-Rz-comm}
\begin{align} 
\mathcal{R} _\Zb^\nn & 
= \Sigma ^{-1}\Jg\nbs^6V \nbs_2(\Jg\Zbn) 
+ \nbs^6( \Sigma ^{-1}\Jg) (\Q\p_\s+V\p_2) (\Jg\Zbn) 
-\bubu{  \nbs^6\big( \tfrac{\alpha }{\Sigma} (  \Jg\Wbn - \Jg\Zbn)\big) \Jg\Zbn}
\notag \\
& \qquad
- \alpha \nbs^6 (\Jg g^{- {\frac{1}{2}} } )\nbs_2 (\Jg\Abn)
+ \alpha \nbs^6(\Jg  g^{- {\frac{1}{2}} } \Abn) \nbs_2 \Jg
+ \tfrac{\alpha}{2}  \nbs^6\big( \Jg g^{- {\frac{1}{2}} } (\Jg\Wbn + \Jg\Zbn - 2\Jg\Abt) \nn_i \big)\nbs_2  \tt^i
 \notag \\
& \qquad
 - \tfrac{2 \alpha}{\eps} \nbs^6\big( \Jg (\Abn+\Zbt) \nn_i \big)  \nbs_1  \tt^i
 +2 \alpha  \nbs^6\big( \Jg g^{- {\frac{1}{2}} } \nbs_2 h \big) \nbs_2(\Jg\Zbn)
 +2 \alpha \nbs^6\big(   \Jgt g^{- {\frac{1}{2}} } \nbs_2 h  (\Abn+\Zbt)\nn_i\big)  \nbs_2\tt^i
 \notag \\
& \qquad
 -2 \alpha  \nbs^6\big( \Jg g^{- {\frac{1}{2}} } \nbs_2h  \Zbn \big) \nbs_2\Jg 
 - \tfrac{2 \alpha}{\eps} \nbs_1  \Jg \, \nbs^6 \Zbn  
\,, \label{Zbn-remainder} \\
\mathcal{C}_\Zb^\nn & 
= \doublecom{\nbs^6, \Sigma^{-1} \Jg , (\Q\p_\s+V\p_2)  (\Jg\Zbn)}
+ \Sigma^{-1}\Jg \doublecom{\nbs^6, V,  \nbs_2 (\Jg\Zbn)}
- \bubu{ \doublecom{ \nbs^6, \tfrac{\alpha }{\Sigma} (  \Jg\Wbn - \Jg\Zbn), \Jg\Zbn}  }
\notag\\
&\qquad
+ \alpha \doublecom{ \nbs^6 , \Jg g^{- {\frac{1}{2}} },  \nbs_2 (\Jg\Abn)}
- \alpha  \doublecom{\nbs^6, \Jg  g^{- {\frac{1}{2}} } \Abn,  \nbs_2 \Jg}
 \notag \\
&\qquad
- \tfrac{\alpha}{2} \doublecom{  \nbs^6, \Jg g^{- {\frac{1}{2}} } (\Jg\Wbn + \Jg\Zbn - 2\Jg\Abt) \nn_i , \nbs_2  \tt^i}
 \notag \\
&\qquad
 - \tfrac{2 \alpha}{\eps} \doublecom{ \nbs^6,  \Jg (\Abn+\Zbt) \nn_i ,  \nbs_1  \tt^i}
 + \tfrac{2 \alpha}{\eps} \doublecom{ \nbs^6 , \Zbn,   \nbs_1  \Jg}
 +2 \alpha \doublecom{  \nbs^6,  \Jg g^{- {\frac{1}{2}} } \nbs_2 h , \nbs_2(\Jg\Zbn)}
  \notag \\
& \qquad
 +2 \alpha\doublecom{ \nbs^6,    \Jgt g^{- {\frac{1}{2}} } \nbs_2 h  (\Abn+\Zbt)\nn_i ,   \nbs_2\tt^i}
 -2 \alpha  \doublecom{\nbs^6,  \Jg g^{- {\frac{1}{2}} } \nbs_2h  \Zbn ,  \nbs_2\Jg }
 \,,
 \end{align} 
\end{subequations} 
for \eqref{energy-Zn-s}, and
\begin{subequations} 
\label{Ca-Ra-comm}
\begin{align} 
\mathcal{R} _\Ab^\nn & 
= \Sigma ^{-1}\Jg\nbs^6V \nbs_2(\Jg \Abn) 
+ \nbs^6( \Sigma ^{-1}\Jg)(\Q\p_\s+V\p_2) (\Jg \Abn)
- \bubu{ \nbs^6\big( \tfrac{\alpha }{2\Sigma} (  \Jg\Wbn - \Jg\Zbn)\big) \Jg \Abn }
 \notag \\
&\qquad
+ \alpha \nbs^6\big(  \Jg g^{- {\frac{1}{2}} }  \big)  \nbs_2 (\Jg\Sbn)
- \alpha  \nbs^6\big( g^{- {\frac{1}{2}} }  (\Jg\Sbn) \big)  \nbs_2 \Jg
+ \alpha \nbs^6\big(   \Jgt g^{- {\frac{1}{2}} } \Sbt \nn_i\big) \nbs_2\nbs^6 \tt^i
\notag \\
&\qquad
+ \tfrac{\alpha}{2\eps} \nbs^6\big(  \Jg\Wbn +\Jg\Zbn - 2\Jg \Abt) \nn_i \big)  \nbs_1 \tt^i
+ \alpha \nbs^6\big( \Jg g^{- {\frac{1}{2}} } \nbs_2 h \nbs_2 \big) \nbs_2 (\Jg\Abn)
+ \alpha \nbs^6\big( \Jg g^{- {\frac{1}{2}} } \nbs_2 h\big)  \Abn \nbs_2\Jg 
 \notag \\
& \qquad
- \tfrac{\alpha}{2} \nbs^6\big( \Jg g^{- {\frac{1}{2}} } \nbs_2 h (\Jg\Wbn +\Jg\Zbn - 2\Jg \Abt) \nn_i \big)   \nbs_2  \tt^i
- \tfrac{\alpha}{\eps}\nbs_1 \Jg \,  \nbs^6\Abn 
\,, \label{Abn-remainder} \\
\mathcal{C}_\Ab^\nn & 
= \doublecom{\nbs^6, \Sigma^{-1} \Jg , (\Q\p_\s+V\p_2)  (\Jg \Abn)}
+ \Sigma^{-1}\Jg \doublecom{\nbs^6, V, \nbs_2(\Jg \Abn)}
- \bubu{ \doublecom{ \nbs^6,   \tfrac{\alpha}{2\Sigma}( \Jg\Wbn - \Jg\Zbn) , \Jg \Abn} }
 \notag \\
&\qquad
+ \alpha \doublecom{ \nbs^6,   \Jg g^{- {\frac{1}{2}} }  ,  \nbs_2 (\Jg\Sbn) }
- \alpha  \doublecom{ \nbs^6, g^{- {\frac{1}{2}} }  (\Jg\Sbn) ,   \nbs_2 \Jg }
+ \alpha  \doublecom{ \nbs^6,    \Jgt g^{- {\frac{1}{2}} } \Sbt \nn_i ,  \nbs_2\nbs^6 \tt^i}
\notag \\
&\qquad
+ \tfrac{\alpha}{2\eps} \doublecom{\nbs^6,  (  \Jg\Wbn +\Jg\Zbn - 2\Jg \Abt) \nn_i ,  \nbs_1 \tt^i}
- \tfrac{\alpha}{\eps}\doublecom{ \nbs^6, \Abn,  \nbs_1 \Jg }
+ \alpha \doublecom{ \nbs^6 ,  \Jg g^{- {\frac{1}{2}} } \nbs_2 h \nbs_2 ,  \nbs_2 (\Jg\Abn) }
\notag \\
& \qquad
- \tfrac{\alpha}{2} \doublecom{ \nbs^6 ,  \Jg g^{- {\frac{1}{2}} } \nbs_2 h (\Jg\Wbn +\Jg\Zbn - 2\Jg \Abt) \nn_i ,   \nbs_2  \tt^i }
+ \alpha \doublecom{ \nbs^6 ,  \Jg g^{- {\frac{1}{2}} } \nbs_2 h ,  \Abn \nbs_2\Jg  }
 \,,
\end{align} 
\end{subequations} 
for \eqref{energy-An-s}. 
 
\subsection{The $\nbs^6$ normal energy identity}
We compute the  following spacetime $L^2$ inner-product:
\begin{equation} 
\tint \jb \Fgss \Big( \underbrace{\eqref{energy-Wn-s} \ \Jg  \nbs^6(\Jg\Wbn)} _{ I^{\WW_n}}
+ \underbrace{\eqref{energy-Zn-s} \  \nbs^6(\Jg\Zbn)} _{ I^{\ZZ_n}}
+ \underbrace{2\eqref{energy-An-s} \ \nbs^6(  \Jg \Abn) } _{ I^{\AA_n}}
 \Big)  {\rm d} x {\rm d}\s'=0 \,, \label{D6-L2-N}
\end{equation} 
where
\begin{equation*}
\jb = \Sigma^{-2\beta+1} 
\end{equation*} 
and $\beta>0$ is a constant which will be chosen to be sufficiently large, only with respect to $\alpha$, in \eqref{eq:normal:bounds:beta}  below. For notational convenience we will mostly omit the spacetime Lebesgue measure ${\rm d} x {\rm d}\s'$ from these integrals.  We proceed to the  analysis of  each of the three integrals in \eqref{D6-L2-N}.   

\subsection{The integral  $I^{\WW_n}$} 
 We additively decompose the integral  $I^{\WW_n}$ as 
 \begin{subequations} 
 \label{Integral-Wbn}
\begin{align}
I^{\WW_n}&= I^{\WW_n}_1+I^{\WW_n}_3+I^{\WW_n}_4 +I^{\WW_n}_5+I^{\WW_n}_6
\notag \,, \\
 I^{\WW_n}_1 &=
\tint \tfrac{1}{\Sigma^{2\beta}} \Fgss \Jg (\Q\p_\s +V\p_2)\nbs^6( \Jg\Wbn  ) \  \nbs^6(\Jg\Wbn)
\,, \label{I1-Wbn} \\
I^{\WW_n}_3 &=
 \alpha \tint \jb \Fgss \Jg g^{- {\frac{1}{2}} }  \nbs_2\nbs^6(\Jg\Abn)  \  \nbs^6(\Jg\Wbn)
\,, \label{I3-Wbn} \\
 I^{\WW_n}_4 &=
-\alpha \tint \jb  \Fgss \Jg g^{- {\frac{1}{2}} }  \Abn \nbs_2\nbs^6 \Jg \  \nbs^6(\Jg\Wbn)
  \,,   \label{I4-Wbn} \\
 I^{\WW_n}_5 &=
-\tfrac{\alpha }{2}  \tint \jb  \Fgss  \Jg g^{- {\frac{1}{2}} } (\Jg\Wbn + \Jg\Zbn - 2\Jg\Abt)\nbs_2 \nbs^6 \tt\cdo\nn \ \nbs^6(\Jg\Wbn)
 \,, \label{I5-Wbn} \\
 I^{\WW_n}_6 &=
-\tfrac{\alpha }{2}  \tint \jb \Fgss \Jg \big(  \nbs^6 \Fwn + \mathcal{R}^\nn_{\Wb} + \mathcal{C}^\nn_{\Wb}\big) \ \nbs^6(\Jg\Wbn)
 \,. \label{I6-Wbn}
\end{align} 
\end{subequations} 
Note that the integral $I^{\WW_n}_2$ is missing since there is no pure damping term in the $\nbs^6(\Jg \Wbn)$ evolution~\eqref{energy-Wn-s}.

\subsection{The integral  $I^{\ZZ_n}$} 
In analogy to~\eqref{Integral-Zbn}, we additively decompose the integral  $I^{\ZZ_n}$ as 
 \begin{subequations} 
 \label{Integral-Zbn}
 \begin{align}
  I^{\ZZ_n} & = I^{\ZZ_n}_1 + I^{\ZZ_n}_2 + I^{\ZZ_n}_3+ I^{\ZZ_n}_4+ I^{\ZZ_n}_5+ I^{\ZZ_n}_6+ I^{\ZZ_n}_7+ I^{\ZZ_n}_8
  + I^{\ZZ_n}_9+ I^{\ZZ_n}_{10}
  \,,  \notag \\
 I^{\ZZ_n}_1& =
 \tint \tfrac{1}{\Sigma^{2\beta}} \Fgss  \Jg(\Q\p_\s +V\p_2)\nbs^6( \Jg\Zbn) \  \nbs^6(\Jg\Zbn)
  \,, \label{I1-Zbn}\\
 I^{\ZZ_n}_2& =
-\bubu{  \tint \tfrac{\alpha }{\Sigma^{2\beta}} ( \Jg\Wbn -\Jg\Zbn) } \Fgss\sabs{\nbs^6(\Jg\Zbn)  }^2
  \,, \label{I2-Zbn} \\
I^{\ZZ_n}_3 &=
- \alpha \tint \jb  g^{- {\frac{1}{2}} }\Fgss \Jg  \nbs_2\nbs^6(\Jg\Abn)  \  \nbs^6(\Jg\Zbn)
\,, \label{I3-Zbn} \\
 I^{\ZZ_n}_4 &=
\alpha \tint \jb\Fgss \Jg g^{- {\frac{1}{2}} }  \Abn \nbs_2\nbs^6 \Jg \  \nbs^6(\Jg\Zbn)
  \,,   \label{I4-Zbn} \\
 I^{\ZZ_n}_5 &=
\tfrac{\alpha }{2}  \tint \jb \Fgss \Jg g^{- {\frac{1}{2}} } (\Jg\Wbn + \Jg\Zbn - 2\Jg\Abt)\nbs_2 \nbs^6 \tt\cdo\nn \ \nbs^6(\Jg\Zbn)
 \,, \label{I5-Zbn} 
\\
 I^{\ZZ_n}_6& =
-{\tfrac{2 \alpha }{\eps}} \tint  \jb \Fgss   \nbs_1 \nbs^6 (\Jg\Zbn) \ \nbs^6(\Jg\Zbn)
\,,  \label{I6-Zbn}\\
 I^{\ZZ_n}_7& =
-{\tfrac{2 \alpha }{\eps}} \tint  \jb  \Fgss  \Jg(\Abn+\Zbt) \Big( \nbs_1\nbs^6\tt \cdo\nn - \eps \Jg g^{- {\frac{1}{2}} } \nbs_2h\, \nbs_2\nbs^6\tt \cdo\nn   \Big)   \ \nbs^6(\Jg\Zbn)
\,,  \label{I7-Zbn}\\
 I^{\ZZ_n}_8& =
{\tfrac{2 \alpha }{\eps}} \tint  \jb \Fgss  \Zbn  \Big( \nbs_1\nbs^6\Jg -  \eps \Jg g^{- {\frac{1}{2}} } \nbs_2 h\,  \nbs_2\nbs^6\Jg \Big) \ \nbs^6(\Jg\Zbn)
\,,  \label{I8-Zbn}\\
 I^{\ZZ_n}_9& =
2 \alpha \tint  \jb  \Fgss \Jg g^{- {\frac{1}{2}} } \nbs_2h  \nbs_2 \nbs^6 (\Jg\Zbn) \ \nbs^6(\Jg\Zbn)
\,,  \label{I9-Zbn}\\
 I^{\ZZ_n}_{10}& =
 - \tint  \jb \Fgss \big(  \nbs^6 \Fzn + \mathcal{R}^\nn_{\Zb} + \mathcal{C}^\nn_{\Zb}\big) \  \nbs^6(\Jg\Zbn)
\,.  \label{I10-Zbn}
\end{align} 
\end{subequations} 

 \subsection{The integral  $I^{\AA_n}$}  
We additively decompose the integral  $I^{\AA_n}$ as
\begin{subequations} 
\label{Integral-Abn}
 \begin{align}
  I^{\AA_n} & = I^{\AA_n}_1 + I^{\AA_n}_2 + I^{\AA_n}_3+ I^{\AA_n}_4+ I^{\AA_n}_5+ I^{\AA_n}_6+ I^{\AA_n}_7+ I^{\AA_n}_8
  + I^{\AA_n}_9+ I^{\AA_n}_{10}
  \,,  \notag \\
 I^{\AA_n}_1& =
2 \tint \tfrac{1}{\Sigma^{2\beta}}  \Fgss \Jg(\Q\p_\s +V\p_2)\nbs^6( \Jg\Abn) \  \nbs^6(\Jg\Abn)
  \,, \label{I1-Abn}\\
 I^{\AA_n}_2& =
- \bubu{ \tint \tfrac{\alpha }{\Sigma^{2\beta}}   \big( \Jg\Wbn - \Jg\Zbn\big) } \Fgss\sabs{\nbs^6(\Jg\Abn)}^2
  \,, \label{I2-Abn} \\
I^{\AA_n}_3 &=
2 \alpha \tint \jb  g^{- {\frac{1}{2}} } \Fgss \Jg  \nbs_2\nbs^6(\Jg\Sbn)  \  \nbs^6(\Jg\Abn)
\,, \label{I3-Abn} \\
 I^{\AA_n}_4 &=
-2\alpha \tint \jb  g^{- {\frac{1}{2}} } \Fgss (\Jg \Sbn) \nbs_2\nbs^6 \Jg \  \nbs^6(\Jg\Abn)
  \,,   \label{I4-Abn} \\
 I^{\AA_n}_5 &=
2 \alpha \tint \jb \Fgss \Jg g^{- {\frac{1}{2}} }  \Sbt \nbs_2 \nbs^6 \tt\cdo\nn \ \nbs^6(\Jg\Abn)
 \,, \label{I5-Abn} 
\\
 I^{\AA_n}_6& =
-{\tfrac{2 \alpha }{\eps}} \tint  \jb   \Fgss \nbs_1 \nbs^6 (\Jg\Abn) \ \nbs^6(\Jg\Abn)
\,,  \label{I6-Abn}\\
 I^{\AA_n}_7& =
{\tfrac{\alpha }{\eps}} \tint  \jb \Fgss   (\Jg\Wbn+\Jg\Zbn - 2\Jg\Abt) \Big( \nbs_1\nbs^6\tt \cdo\nn - \eps \Jg g^{- {\frac{1}{2}} } \nbs_2h\, \nbs_2\nbs^6\tt \cdo\nn   \Big)  \ \nbs^6(\Jg\Abn)
\,,  \label{I7-Abn}\\
 I^{\AA_n}_8& =
{\tfrac{2 \alpha }{\eps}} \tint  \jb  \Fgss  \Abn \Big( \nbs_1\nbs^6\Jg -  \eps \Jg g^{- {\frac{1}{2}} } \nbs_2 h\,  \nbs_2\nbs^6\Jg \Big)  \ \nbs^6(\Jg\Abn)
\,,  \label{I8-Abn}\\
 I^{\AA_n}_9& =
2 \alpha \tint  \jb  \Fgss \Jg g^{- {\frac{1}{2}} } \nbs_2h  \nbs_2 \nbs^6 (\Jg\Abn) \ \nbs^6(\Jg\Abn)
\,,  \label{I9-Abn}\\
 I^{\AA_n}_{10}& =
 - \tint  \jb \Fgss \big(  \nbs^6 \Fzn + \mathcal{R}^\nn_{\Zb} + \mathcal{C}^\nn_{\Zb}\big) \  \nbs^6(\Jg\Abn)
\,.  \label{I10-Abn}
 \end{align} 
 \end{subequations} 
 
 \subsection{The exact derivative terms}
For the terms involving a time derivative, we note that summing \eqref{I1-Wbn}, \eqref{I1-Zbn}, and \eqref{I1-Abn}, integrating by parts and appealing to \eqref{p2h-evo-s}, \eqref{adjoint-3}, \eqref{Jg-evo-s}, and \eqref{Sigma0i-ALE-s}, we obtain
\begin{align}
I^{\WW_n}_1& + I^{\ZZ_n}_1  + I^{\AA_n}_1 \notag \\
&=\tint \tfrac{ \mathcal{J}^{\frac 32} \Jg}{2\Sigma^{2\beta}}  (\Q\p_\s +V\p_2) \Bigl( \sabs{\nbs^6(\Jg\Wbn)}^2 + \sabs{\nbs^6(\Jg\Zbn)}^2 + 2 \sabs{\nbs^6(\Jg\Abn)}^2 \Bigr)
\notag\\
&= \tfrac 12 \snorm{\tfrac{\mathcal{J}^{\frac 34}(\Jg \Q)^{\frac 12}}{\Sigma^\beta} \nbs^6 (\Jg\Wbn)\cdot,\s)}_{L^2_x}^2
+ \tfrac 12 \snorm{\mathcal{J}^{\frac 34}\tfrac{(\Jg \Q)^{\frac 12}}{\Sigma^\beta} \nbs^6(\Jg \Zbn)(\cdot,\s)}_{L^2_x}^2
+   \snorm{\tfrac{\mathcal{J}^{\frac 34}(\Jg \Q)^{\frac 12}}{\Sigma^\beta} \nbs^6(\Jg \Abn)(\cdot,\s)}_{L^2_x}^2
\notag\\
&\quad 
- \tfrac 12 \snorm{\tfrac{\mathcal{J}^{\frac 34}(\Jg \Q)^{\frac 12}}{\Sigma^\beta} \nbs^6(\Jg \Wbn)(\cdot,0)}_{L^2_x}^2
- \tfrac 12 \snorm{\tfrac{\mathcal{J}^{\frac 34}(\Jg \Q)^{\frac 12}}{\Sigma^\beta} \nbs^6(\Jg \Zbn)(\cdot,0)}_{L^2_x}^2
-   \snorm{\tfrac{\mathcal{J}^{\frac 34}(\Jg \Q)^{\frac 12}}{\Sigma^\beta} \nbs^6 (\Jg\Abn)(\cdot,0)}_{L^2_x}^2
\notag\\
&\quad
+ \tint \tfrac{1}{\Sigma^{2\beta}} \mathsf{G}_0  \Bigl( \sabs{\nbs^6(\Jg\Wbn)}^2 + \sabs{\nbs^6(\Jg\Zbn)}^2 + 2 \sabs{\nbs^6(\Jg\Abn)}^2 \Bigr)
\label{eq:heavy:fuel:1n}
\end{align}
where we have defined  
\begin{equation*}
 \mathsf{G}_0 
 :=  - \tfrac 12 (\Q \p_\s + V\p_2) \bigr( \mathcal{J}^{\frac 32} \Jg \bigl) 
 + \tfrac 12 (V \Qr_2 - \Qr_\s  - \nbs_2 V - 2 \alpha \beta (\Zbn+\Abt) ) \mathcal{J}^{\frac 32} \Jg
 \,.
\end{equation*}
At this stage, we record the pointwise bound
\begin{equation}
\mathsf{G}_0 
\geq 
\underbrace{- \tfrac 12 (\Q \p_\s + V\p_2) \bigr( \mathcal{J}^{\frac 32} \Jg \bigl)}_{=: \mathsf{G_{good}}}
- \tfrac{250^2}{\eps} \Q \Jg  \mathcal{J}^{\frac 32}
- \Cn \brak{\beta}   \mathcal{J}^{\frac 12} \Jg  
\label{eq:G0:n:lower}
  \,,
\end{equation}
which follows from \eqref{eq:Jgb:less:than:1}, \eqref{bootstraps}, and \eqref{eq:Q:all:bbq}.

For the terms involving a $\p_1$ derivative of the fundamental variables, we add \eqref{I6-Zbn} and \eqref{I6-Abn}, and integrate by parts with respect to $\p_1$ (here, recall that $\p_1 \mathcal{J} = 0$) to arrive at
\begin{align}
I^{\ZZ_n}_6  + I^{\AA_n}_6
&= - \alpha \tint \jb \mathcal{J}^{\frac 32} \p_1 \bigl( \sabs{\nbs^6(\Jg \Zbn)}^2 + \sabs{\nbs^6(\Jg\Abn)}^2 \bigr)
\notag\\
&=   \tint \tfrac{1}{\Sigma^{2\beta}} \mathsf{G}_1 \bigl( \sabs{\nbs^6(\Jg \Zbn)}^2 + \sabs{\nbs^6(\Jg\Abn)}^2 \bigr)
\label{eq:heavy:fuel:2n}
\end{align}
where 
\begin{equation}
\mathsf{G}_1 
:=  - \alpha (2\beta-1) \mathcal{J}^{\frac 32} \Sigma,_1 
\geq
\alpha (\beta-\tfrac 12)  \bigl(\tfrac{1}{4\eps} - \tfrac{16}{(1+\alpha) \eps}  \Jg \Q -  \Cn   \bigr)  \mathcal{J}^{\frac 32}
\label{eq:heavy:fuel:G1:lowern}
\end{equation}
as soon as $\beta\geq \frac 12$, 
in light of \eqref{p1-Sigma}, \eqref{bootstraps}, \eqref{eq:Q:all:bbq}, and \eqref{eq:signed:Jg}.

For the terms involving a $\p_2$ derivative of the fundamental variables, we add \eqref{I9-Zbn} and \eqref{I9-Abn}, and integrate by parts using~\eqref{adjoint-2} to arrive at
\begin{align}
&I^{\ZZ_n}_9  + I^{\AA_n}_9
\notag \\
&=  \alpha \tint \jb \mathcal{J}^{\frac 32} \Jg g^{- {\frac{1}{2}} } \nbs_2h\,  \nbs_2 
\bigl( \sabs{\nbs^6(\Jg\Zbn) }^2 +  \sabs{\nbs^6(\Jg\Zbn) }^2 \bigr)
\notag\\
&= -   \alpha \int \Qb_2 \jb \mathcal{J}^{\frac 32} \Jg g^{- {\frac{1}{2}} } \nbs_2h\,   
\bigl( \sabs{\nbs^6(\Jg\Zbn) }^2 +  \sabs{\nbs^6(\Jg\Zbn) }^2 \bigr)\Bigr|_{\s} 
+ \tint \tfrac{1}{\Sigma^{2\beta}} \mathsf{G}_2 
\bigl( \sabs{\nbs^6(\Jg\Zbn) }^2 +  \sabs{\nbs^6(\Jg\Zbn) }^2 \bigr)\,  
\label{eq:heavy:fuel:3n}
\end{align}
where 
\begin{equation}
\mathsf{G}_2 
=  \alpha\Sigma^{2\beta}    (\Qr_2 - \nbs_2) \bigl( \jb \mathcal{J}^{\frac 32} \Jg g^{- {\frac{1}{2}} } \nbs_2h \bigr)
\geq -  \Cn  \brak{\beta} \eps \mathcal{J}^{\frac 12} \Jg 
\label{eq:heavy:fuel:G2:lowern}
\end{equation}
and we have appealed to  \eqref{bootstraps}, \eqref{eq:Q:all:bbq}, \eqref{eq:signed:Jg}, and the bound $\mathcal{J} \leq \Jg$. At this stage we also note that the bootstrap inequalities imply
\begin{align}
&-   \alpha \int \Qb_2 \jb\mathcal{J}^{\frac 32}\Jg g^{- {\frac{1}{2}} } \nbs_2h\,   
\bigl( \sabs{\nbs^6(\Jg\Zbn) }^2 +  \sabs{\nbs^6(\Jg\Zbn) }^2 \bigr)\Bigr|_{\s} 
\notag\\ 
&\qquad
\geq - \Cn \eps^2 \Bigl(
 \snorm{\tfrac{\mathcal{J}^{\frac 34} (\Jg \Q)^{\frac 12}}{\Sigma^\beta} \nbs^6 (\Jg\Zbn)(\cdot,\s)}_{L^2_x}^2
+   \snorm{\tfrac{\mathcal{J}^{\frac 34}  (\Jg \Q)^{\frac 12}}{\Sigma^\beta} \nbs^6(\Jg \Abn)(\cdot,\s)}_{L^2_x}^2 \Bigr)
\,.
\label{eq:heavy:fuel:G2:lower:an}
\end{align}

Finally, we have the pure damping terms $I^{\ZZ_n}_2$ and $I^{\AA_n}_2$.  Summing  \eqref{I2-Zbn} and \eqref{I2-Abn}
yields
\begin{equation}
I^{\ZZ_n}_2  + I^{\AA_n}_2
=  \tint \tfrac{1}{\Sigma^{2\beta}}  \mathsf{G_3}
\bigl( \sabs{\nbs^6(\Jg\Zbn) }^2 +  \sabs{\nbs^6(\Jg\Abn) }^2 \bigr) \,,
\label{eq:heavy:fuel:4n}
\end{equation} 
and by using \eqref{bootstraps}, \eqref{eq:Q:bbq}, and  \eqref{eq:signed:Jg}, we choosing $\eps $ sufficiently small, we have that
\begin{equation}
\bubu{ 
\mathsf{G_3} 
:= - \alpha ( \Jg \Wbn -  \Jg \Zbn) \Fgss
\ge   \tfrac{9 \alpha }{10\eps}\Fgss  -  \tfrac{33}{\eps} \Q \Fgss \Jg   \,. }
\label{eq:heavy:fuel:G3:lower:an}
\end{equation}

There are three terms with seven derivatives acting on the fundamental variables, and once again, these terms combine 
to produce an exact derivative, which we then integrate by parts. 
Adding \eqref{I3-Wbn}, \eqref{I3-Zbn}, \eqref{I3-Abn}, employing the identity $\Sbn = \tfrac 12 \Wbn - \tfrac 12 \Zbn$, 
using~\eqref{adjoint-2} we have that 
\begin{align} 
I^{\WW_n}_3+I^{\ZZ_n}_3+I^{\AA_n}_3
&= 2\alpha \tint \jb g^{- {\frac{1}{2}} } \Fgss \Jg \nbs_2 \Big( \nbs^6(\Jg\Abn)  \  \nbs^6(\Jg\Sbn)\Big)
\notag \\
&= \alpha \tint (\Qr_2- \nbs_2)( \jb g^{- {\frac{1}{2}} } \Fgss \Jg)   \nbs^6(\Jg\Abn)  \ 
\big(  \nbs^6(\Jg\Wbn)- \nbs^6(\Jg\Zbn) \big)
\notag \\
& \qquad \qquad
-\alpha  \int \jb \Qb_2 g^{- {\frac{1}{2}} }\Fgss \Jg \nbs^6(\Jg\Abn)  \  \big(  \nbs^6(\Jg\Wbn)- \nbs^6(\Jg\Zbn) \big)\Big|_\s \,.
\label{eq:heavy:fuel:5n}
\end{align} 
These two integrals are bounded in the usual way by appealing to \eqref{eq:heavy:fuel:G3:lower}.

 Summarizing,  the identities~\eqref{eq:heavy:fuel:1n}, \eqref{eq:heavy:fuel:2n}, \eqref{eq:heavy:fuel:3n},
 \eqref{eq:heavy:fuel:4n},   \eqref{eq:heavy:fuel:5n}  
 and bounds~ \eqref{eq:heavy:fuel:G3:lower}, \eqref{eq:G0:n:lower},
  \eqref{eq:heavy:fuel:G1:lowern}, 
 \eqref{eq:heavy:fuel:G2:lower:an},
 \eqref{eq:heavy:fuel:G2:lowern}, 
 \eqref{eq:heavy:fuel:G3:lower:an}, upon taking 
 $\eps$ to be sufficiently small in terms of $\alpha,\kappa_0$ and $\Cdata$, and taking $\beta \geq 1 $, gives
\begin{align}
& I^{\WW_n}_1 + I^{\ZZ_n}_1  + I^{\AA_n}_1 + I^{\ZZ_n}_2 + I^{\AA_n}_2
 + I^{\WW_n}_3 + I^{\ZZ_n}_3 + I^{\AA_n}_3
 + I^{\ZZ_n}_6  + I^{\AA_n}_6 
 + I^{\ZZ_n}_9  + I^{\AA_n}_9 
 \notag\\
&\geq \bigl( \tfrac 12 - \Cn \eps\bigr) 
\Bigl( \snorm{\tfrac{\mathcal{J}^{\frac 34} (\Jg \Q)^{\frac 12}}{\Sigma^\beta} \nbs^6(\Jg \Wbn)(\cdot,\s)}_{L^2_x}^2
+  \snorm{\tfrac{\mathcal{J}^{\frac 34} (\Jg \Q)^{\frac 12}}{\Sigma^\beta} \nbs^6(\Jg \Zbn)(\cdot,\s)}_{L^2_x}^2
+ 2  \snorm{\tfrac{\mathcal{J}^{\frac 34} (\Jg \Q)^{\frac 12}}{\Sigma^\beta} \nbs^6( \Jg\Abn)(\cdot,\s)}_{L^2_x}^2 \Bigr)
\notag\\
&\qquad 
- \tfrac 12 \Bigl( \snorm{ \tfrac{\mathcal{J}^{\frac 34}(\Jg \Q)^{\frac 12}}{\Sigma^\beta} \nbs^6 (\Jg\Wbn)(\cdot,0)}_{L^2_x}^2
+ \snorm{\tfrac{\mathcal{J}^{\frac 34} (\Jg \Q)^{\frac 12}}{\Sigma^\beta} \nbs^6 (\Jg\Zbn)(\cdot,0)}_{L^2_x}^2
+2  \snorm{\tfrac{\mathcal{J}^{\frac 34} (\Jg \Q)^{\frac 12}}{\Sigma^\beta} \nbs^6(\Jg \Abn)(\cdot,0)}_{L^2_x}^2
\Bigr)
\notag\\
&\qquad
+ \tint \tfrac{1}{\Sigma^{2\beta}} \Bigl( \mathsf{G_{good}} - \Cn \beta\mathcal{J}^{\frac 12} \Jg\Bigr)
\Bigl(  |\nbs^6(\Jg\Wbn)|^2 + |\nbs^6(\Jg\Zbn)|^2 + 2 |\nbs^6(\Jg\Abn)|^2 \Bigr)
\notag\\
&\qquad 
+ \Bigl(  \tfrac{\alpha( \beta - \frac 12)}{8\eps}   + \bubu{   \tfrac{9\alpha }{10\eps}  } \Bigr)
\int_0^{\s} 
\Bigl( \snorm{\tfrac{ \mathcal{J}^{\frac 34} }{\Sigma^\beta} \nbs^6(\Jg\Zbn) (\cdot,\s')}_{L^2_x}^2
+  \snorm{\tfrac{ \mathcal{J}^{\frac 34} }{\Sigma^\beta} \nbs^6(\Jg\Abn) (\cdot,\s')}_{L^2_x}^2 \Bigr)
{\rm d} \s'
\notag\\
&\qquad 
- \Bigl( \tfrac{16 \alpha  (\beta -\frac 12) }{(1+\alpha)\eps}  + \tfrac{\bubu{ 33 } + 2 \cdot 250^2}{\eps} \Bigr) 
\int_0^{\s} 
\Bigl(  \snorm{\tfrac{\mathcal{J}^{\frac 34}(\Jg \Q)^{\frac 12}}{\Sigma^\beta} \nbs^6(\Jg \Zbn)(\cdot,\s')}_{L^2_x}^2
+    \snorm{\tfrac{\mathcal{J}^{\frac 34}(\Jg \Q)^{\frac 12}}{\Sigma^\beta} \nbs^6 (\Jg\Abn)(\cdot,\s')}_{L^2_x}^2 \Bigr)
{\rm d} \s'
\notag\\
&\qquad
- \tfrac{250^2}{\eps}
\int_0^{\s} 
 \snorm{\tfrac{\mathcal{J}^{\frac 34}(\Jg \Q)^{\frac 12}}{\Sigma^\beta} \nbs^6(\Jg \Wbn)(\cdot,\s')}_{L^2_x}^2 {\rm d} \s'
\,,
\label{eq:I:n:12369}
\end{align}
where as usual $\Cn = \Cn (\alpha,\kappa_0,\Cdata)$ is a positive computable constant. Note   that $\Cn$ is independent of $\beta$.

\subsection{The terms with over-differentiated geometry}
\label{sec:normal:overdiff}
Next, we consider the terms in \eqref{Integral-Wbn}--\eqref{Integral-Abn} which contain seven derivatives on either $\tt$ or $\Jg$.
Specifically, we shall study the integrals 
$ I^{\WW_n}_4$, 
$ I^{\WW_n}_5$,
$ I^{\ZZ_n}_4$,  
$ I^{\ZZ_n}_5$, 
$ I^{\ZZ_n}_7$, 
$ I^{\ZZ_n}_8$, and 
$ I^{\AA_n}_4$, 
$ I^{\AA_n}_5$, 
$ I^{\AA_n}_7$, 
$ I^{\AA_n}_8$.  
Various grouping of these integrals shall combine to produce important cancellations and thereby yield the  desired bounds.

\subsubsection{Important geometric identities}

\begin{lemma} \label{lem:tau-Jg-D2}
For a function $f(x,\s)$, we have that
\begin{equation} 
\left| \tint f \, \nn \cdo \nbs^6\tt  \ \nbs_2 \nbs^6 \Jg \right| 
+
\left| \tint f \, \nn \cdo \nbs_2 \nbs^6\tt  \  \nbs^6 \Jg \right| 
\les  \eps^{3} \mathsf{K}^2 \brak{\mathsf{B}_6}^2 
\Bigl(\snorm{\mathcal{J}^{\frac 12}  \nbs f}_{L^\infty_{x,\s}}  
+\eps \snorm{\mathcal{J}^{-\frac 12} f}_{L^\infty_{x,\s}}  
   \Bigr) \,.
\label{eq:tau-Jg-D2}
\end{equation} 
\end{lemma} 
\begin{proof} [Proof of Lemma \ref{lem:tau-Jg-D2}]
Using the identities \eqref{d-n-tau} and \eqref{p1-n-tau}, we have that
\begin{align}
\nbs^6 \nbs_2 \Jg
&= 
\tfrac{1}{\eps} g^{\frac 12} \nbs_1(\nn \cdo \nbs^6 \tt) 
-  \Jg \nbs_2 h \nbs_2(\nn \cdo \nbs^6 \tt) 
\notag\\
&\qquad
- \tfrac{1}{\eps} g^{\frac 12} \tt \cdo \nbs_1  \nn\; \tt \cdo \nbs^6 \tt
+ \tfrac{1}{\eps} \nn \cdo \nbs^6(g^{\frac 12} \nn) \;  \nn \cdo \nbs_1 \tt
+ \tfrac{1}{\eps} \doublecom{\nbs^6, g^{\frac 12} \nn_i, \nbs_1 \tt_i}
\notag\\
&\qquad 
+ \Jg \nbs_2 h \tt \cdo \nbs_2 \nn \; \tt \cdo \nbs^6 \tt
- \nn \cdo \nbs^6(\Jg \nbs_2 h \nn) \; \nn \cdo \nbs_2 \tt 
- \doublecom{\nbs^6, \Jg \nbs_2 h \nn_i, \nbs_2 \tt_i}
\label{big-chicken1}
\end{align}
Hence, integrating-by-parts using \eqref{eq:adjoints} and \eqref{big-chicken1}, we find that
\begin{align}
\tint f \, \nn \cdo \nbs^6\tt  \, \nbs_2 \nbs^6 \Jg
&= 
- \tfrac{1}{2\eps} \tint  \nbs_1 (f \, g^{\frac 12})   \, (\nn \cdo \nbs^6 \tt)^2 
- \tfrac 12 \tint (\Qr_2 - \nbs_2) (f \, \Jg \nbs_2 h) \, (\nn \cdo \nbs^6 \tt)^2
\notag\\
&\qquad 
+ \tfrac 12 \int \Qb_2 f \, \Jg \nbs_2 h   (\nn \cdo \nbs^6 \tt)^2 \Bigr|_{\s}
- \tfrac{1}{\eps} \tint f \, \nn \cdo \nbs^6\tt  \,g^{\frac 12} \tt \cdo \nbs_1  \nn\; \tt \cdo \nbs^6 \tt
\notag\\
&\qquad
+ \tfrac{1}{\eps} \tint f \, \nn \cdo \nbs^6\tt  \,\nn \cdo \nbs^6(g^{\frac 12} \nn) \;  \nn \cdo \nbs_1 \tt
+ \tfrac{1}{\eps} \tint f \, \nn \cdo \nbs^6\tt  \,\doublecom{\nbs^6, g^{\frac 12} \nn_i, \nbs_1 \tt_i}
\notag\\
&\qquad 
+ \tint f \, \nn \cdo \nbs^6\tt  \,\Jg \nbs_2 h \tt \cdo \nbs_2 \nn \; \tt \cdo \nbs^6 \tt
- \tint f \, \nn \cdo \nbs^6\tt  \,\nn \cdo \nbs^6(\Jg \nbs_2 h \nn) \; \nn \cdo \nbs_2 \tt 
\notag\\
&\qquad
- \tint f \, \nn \cdo \nbs^6\tt  \, \doublecom{\nbs^6, \Jg \nbs_2 h \nn_i, \nbs_2 \tt_i}  
\,.
\label{big-chicken2}
\end{align}
To conclude the bound for the first term in~\eqref{eq:tau-Jg-D2}, we appeal to the bootstraps~\eqref{bootstraps}, the estimates~\eqref{eq:Q:all:bbq}, to the bounds~\eqref{D6n-bound:a:new}--\eqref{D6n-bound:b:new:bdd}, and to Lemmas~\ref{lem:time:interpolation}--\ref{lem:comm:tangent}, to obtain
\begin{align*}
\left| \tint f \, \nn \cdo \nbs^6\tt  \, \nbs_2 \nbs^6 \Jg\right|
&\les 
\snorm{\mathcal{J}^{\frac 12}  \nbs_1  f}_{L^\infty_{x,\s}}  
\mathsf{K}^2 \eps^3 \brak{\mathsf{B}_6}^2
+\bigl( \snorm{\mathcal{J}^{\frac 12}   f }_{L^\infty_{x,\s}}  
+ \snorm{\mathcal{J}^{\frac 12}  \nbs_2 f}_{L^\infty_{x,\s}} 
\bigr)
\mathsf{K}^2 \eps^5 \brak{\mathsf{B}_6}^2
\notag\\
&\qquad 
+  \snorm{\mathcal{J}^{-\frac 12}  f}_{L^\infty_{x,\s}}   
\mathsf{K}^2 \eps^5 \brak{\mathsf{B}_6}^2
+  \snorm{ f}_{L^\infty_{x,\s}} 
\mathsf{K}^2 \eps^4 \brak{\mathsf{B}_6}^2
\,.
\end{align*}
Since $\mathcal{J}\leq 1$, this concludes the estimate for the first term in~\eqref{eq:tau-Jg-D2}. The bound for the second term follows similarly since \eqref{adjoint-2} gives
\begin{align}
 \tint f \, \nn \cdo \nbs_2 \nbs^6\tt  \  \nbs^6 \Jg
 &= - \tint f \, \nn \cdo  \nbs^6\tt  \ \nbs_2 \nbs^6 \Jg
 - \tint \bigl(\nbs_2 f \, \nn \cdo  \nbs^6\tt + f \tt \cdot \nbs_2 \nn \tt \cdo \nbs^6 \tt \bigr)  \nbs^6 \Jg
 \notag\\
 &\qquad 
 + \tint \Qr_2 f \, \nn \cdo  \nbs^6\tt  \  \nbs^6 \Jg
 - \int \Qb_2 f \, \nn \cdo  \nbs^6\tt  \  \nbs^6 \Jg \Bigr|_{\s}
 \,.
\label{big-chicken-worm2}
\end{align}
Using the bootstrap inequalities~\eqref{bootstraps}, the estimates~\eqref{eq:Q:all:bbq}, and the bounds for the geometry~\eqref{geometry-bounds-new}, we deduce
\begin{align*}
 \left|\tint f \, \nn \cdo \nbs_2 \nbs^6\tt  \  \nbs^6 \Jg\right|
 &\leq \left|\tint f \, \nn \cdo  \nbs^6\tt  \ \nbs_2 \nbs^6 \Jg\right|
 + \Cn \snorm{\mathcal{J}^{\frac 12}  \nbs_2 f}_{L^\infty_{x,\s}}  \eps^3 \mathsf{K}   \brak{\mathsf{B}_6}^2  
 \notag\\
 &\qquad 
 + \Cn \snorm{\mathcal{J}^{\frac 12}  f}_{L^\infty_{x,\s}}  \eps^3 \mathsf{K}   \brak{\mathsf{B}_6}^2  
 + \Cn \snorm{\mathcal{J}^{-\frac 12}  f}_{L^\infty_{x,\s}}  \eps^3 \mathsf{K}   \brak{\mathsf{B}_6}^2 
 \,.
\end{align*}
This concludes the proof of the lemma.
\end{proof} 

\begin{lemma} \label{lem:tau-Jg-D1}
For a function $f(x,\s)$, we have that
\begin{equation} 
\left|\tint f \, \nn \cdo \nbs_1 \nbs^6\tt   \  \nbs^6 \Jg\right|
+
\left|\tint f \, \nn \cdo \nbs^6\tt  \ \nbs_1 \nbs^6 \Jg\right|
\les
\eps^3 \mathsf{K} \brak{\mathsf{B}_6}^2
\Bigl( \|\mathcal{J}^{\frac 12} \nbs f\|_{L^\infty_{x,\s}} +  \|\mathcal{J}^{-\frac 12} f\|_{L^\infty_{x,\s}} \Bigr) 
\,. 
\label{eq:tau-Jg-D1}
\end{equation} 
\end{lemma} 
\begin{proof} [Proof of Lemma \ref{lem:tau-Jg-D1}]
From \eqref{Jg-def}  and \eqref{d-n-tau}  we have that
\begin{equation}
\nbs_1 \tt 
= \nn g^{-1} \nbs_1 \nbs_2 h 
= \eps \nn g^{-1} \nbs_2 (g^{\frac 12} \Jg ) 
= \eps \nn g^{-\frac 12} \nbs_2  \Jg 
+
\eps g^{-\frac 12} \nbs_2 h \Jg  \nbs_2 \tt  
\,.
\end{equation}
Applying $\nbs^6$ to this identity and then taking a dot product with $\nn$, we derive
\begin{align}
\nn \cdo \nbs_1 \nbs^6 \tt 
&= 
\eps   g^{-\frac 12} \nbs_2  \nbs^6 \Jg 
+
\eps \nn_i \nbs^6(\nn_i g^{-\frac 12}) \nbs_2  \Jg 
+
\eps \nn_i  \doublecom{\nbs^6, \nn_i g^{-\frac 12}, \nbs_2 \Jg}
\notag\\
&\quad
+
\eps g^{-\frac 12} \nbs_2 h \Jg  \nn \cdo \nbs_2 \nbs^6  \tt  
+
\eps \nbs^6  \bigl(g^{-\frac 12} \nbs_2 h \Jg\bigr)  \nn \cdo \nbs_2 \tt  
+
\eps \nn_i \doublecom{\nbs^6, g^{-\frac 12} \nbs_2 h \Jg,\nbs_2  \tt_i}  
\,.
\label{eq:tau-Jg-D1-chicken-worm}
\end{align}
Hence, integrating-by-parts in the first term arising from the above expression (using \eqref{eq:adjoints}), we find that
\begin{align}
\tfrac{1}{\eps} \tint f \, \nn \cdo \nbs_1 \nbs^6\tt   \  \nbs^6 \Jg
&=
\tfrac 12 \tint (\Qr_2 - \nbs_2) \bigl( f \,  g^{-\frac 12}\bigr)  ( \nbs^6 \Jg)^2 
-
\tfrac 12 \int \Qb_2 f \,  g^{-\frac 12}  ( \nbs^6 \Jg)^2 \Bigr|_{\s}
\notag\\
&\quad 
+
\tint f \, \nn_i \nbs^6(\nn_i g^{-\frac 12}) \nbs_2  \Jg   \  \nbs^6 \Jg 
+
\tint f \,   \nn_i  \doublecom{\nbs^6, \nn_i g^{-\frac 12}, \nbs_2 \Jg} \  \nbs^6 \Jg 
\notag\\
&\quad 
+
\tint \bigl( f \, g^{-\frac 12} \nbs_2 h \Jg\bigr)  \nn \cdo \nbs_2 \nbs^6  \tt   \  \nbs^6 \Jg 
+
\tint f \,  \nbs^6  \bigl(g^{-\frac 12} \nbs_2 h \Jg\bigr)  \nn \cdo \nbs_2 \tt  \  \nbs^6 \Jg 
\notag\\
&\quad 
+
\tint f \,   \nn_i \doublecom{\nbs^6, g^{-\frac 12} \nbs_2 h \Jg,\nbs_2  \tt_i}  \  \nbs^6 \Jg 
\,.
\end{align}
We note that the fifth term on the right side of the above expression is precisely an integral of the type bounded earlier in~\eqref{eq:tau-Jg-D2}, except that $f$ is replaced by  $ f g^{-\frac 12} \nbs_2 h \Jg$.

Using the bootstraps~\eqref{bootstraps}, the estimates~\eqref{eq:Q:all:bbq},   the bounds~\eqref{geometry-bounds-new}, the previously established estimate~\eqref{eq:tau-Jg-D2}, and the Lemmas~\ref{lem:time:interpolation}--\ref{lem:comm:tangent}, we arrive at
\begin{align}
\tfrac{1}{\eps} \left| \tint f \, \nn \cdo \nbs_1 \nbs^6\tt   \  \nbs^6 \Jg\right|
&\les   
\|\mathcal{J}^{\frac 12} \nbs_2 f\|_{L^\infty_{x,\s}}   \eps^2 \brak{\mathsf{B}_6}^2
+
\|\mathcal{J}^{-\frac 12} f\|_{L^\infty_{x,\s}}  \eps^2 \brak{\mathsf{B}_6}^2
\notag\\
&\quad 
+
\Bigl(\snorm{\mathcal{J}^{\frac 12}  \nbs  f  }_{L^\infty_{x,\s}}  
+\eps \snorm{\mathcal{J}^{-\frac 12} f  }_{L^\infty_{x,\s}}  \Bigr)
\eps^{4} \mathsf{K}^2 \brak{\mathsf{B}_6}^2 
\,.
\end{align}
The bound for the first term on the left side of~\eqref{eq:tau-Jg-D1}, now follows from the above estimate and $\mathcal{J}\leq 1$. In order to bound the second term on the left side of~\eqref{eq:tau-Jg-D1}, we integrate by parts with respect to $\nbs_1$ (using \eqref{adjoint-1}), 
so that 
\begin{equation*}
\tint f \, \nn \cdo \nbs^6\tt  \ \nbs_1 \nbs^6 \Jg
= - \tint f \, \nn \cdo \nbs_1 \nbs^6\tt  \  \nbs^6 \Jg 
- \tint \nbs_1(f \, \nn_i)  \, \nbs^6\tt_i  \  \nbs^6 \Jg
\,.
\end{equation*}
The first term on the right side of the above identity was already previously bounded, while the second term is bounded from above by 
$
\| \mathcal{J}^{\frac 12} \nbs_1 f\|_{L^\infty}  \eps^3 \mathsf{K} \brak{\mathsf{B}_6}^2
+ 
\| \mathcal{J}^{\frac 12} f\|_{L^\infty}  \eps^5 \mathsf{K} \brak{\mathsf{B}_6}^2
$.
This concludes the proof of the lemma.
\end{proof} 

\subsubsection{Bounds for the forcing, remainder, and commutator functions}
\label{subsub:FRC}
We first bound the forcing, remainder, and commutator functions associated to the equation \eqref{energy-Wn-s}.
By the definition of $\Fwn$ in \eqref{forcing-nt}, and appealing to the $\mathcal{J}$ estimate~\eqref{eq:Jgb:less:than:1}, the bootstrap inequalities~\eqref{bootstraps}, the coefficient bounds~\eqref{eq:Q:all:bbq}, the bounds for the geometry~\eqref{geometry-bounds-new}, and the double-commutator bound in \eqref{eq:Lynch:2}, we deduce that
\begin{subequations}
\label{eq:I:W:nn:6:all}
\begin{equation}
\snorm{\tfrac{\mathcal{J}^{\frac 34} (\Jg \Q)^{\frac 12}}{\Sigma^\beta} \nbs^6\Fwn}_{L^2_{x,\s}}
\les  \snorm{\tfrac{\mathcal{J}^{\frac 34} (\Jg \Q)^{\frac 12}}{\Sigma^\beta} \nbs^6(\Jg\Wbn, \Jg\Zbn, \Jg\Abn)}_{L^2_{x,\s}}
+   \eps (\tfrac{4}{\kappa_0})^\beta   \brak{\mathsf{B}_6}
\les  (\tfrac{4}{\kappa_0})^\beta  \brak{\mathsf{B}_6}
\,.
\label{eq:I:W:nn:6:b}
\end{equation}
Similarly, from the definition of $ \mathcal{R}_\Wb^\nn$ in \eqref{Cw-Rw-comm}, and by additionally appealing to \eqref{eq:Jg:Wbn:improve:material}, we find that
\begin{equation}
\snorm{\tfrac{\mathcal{J}^{\frac 34} (\Jg \Q)^{\frac 12}}{\Sigma^\beta}\mathcal{R}_\Wb^\nn}_{L^2_{x,\s}}
\les   \snorm{\tfrac{\mathcal{J}^{\frac 34} (\Jg \Q)^{\frac 12}}{\Sigma^\beta} \nbs^6(\Jg\Wbn,\Jg \Zbn)}_{L^2_{x,\s}}
+ (\tfrac{4}{\kappa_0})^\beta  \brak{\mathsf{B}_6}
\les  (\tfrac{4}{\kappa_0})^\beta  \brak{\mathsf{B}_6}
\label{eq:I:W:nn:6:c}
\,.
\end{equation}
Lastly, using the definition of $\mathcal{C}_\Wb^\nn$ in \eqref{Cw-Rw-comm}, identity~\eqref{eq:Jg:Wb:nn},   the aforementioned  bounds,  Lemma~\ref{lem:anisotropic:sobolev}, the Leibniz rule and Lemma~\ref{lem:time:interpolation},  we may also obtain 
\begin{align}
\snorm{\tfrac{\mathcal{J}^{\frac 34} (\Jg \Q)^{\frac 12}}{\Sigma^\beta}\mathcal{C}_\Wb^\nn}_{L^2_{x,\s}}
\les    (\tfrac{4}{\kappa_0})^\beta  \brak{\mathsf{B}_6}
\,.
\label{eq:I:W:nn:6:d}
\end{align}
\end{subequations}
To give a more detailed explanation of how we arrived at \eqref{eq:I:W:nn:6:d}, we examine a typical commutator term, namely: $\doublecom{\nbs^6, \Sigma^{-1}  , (\Q\p_\s+V\p_2)  (\Jg\Wbn)}$.  Using \eqref{euler-Wn}, for $k=1,2,3,4,5$ we have
\begin{align*} 
\nbs^k&(\Q\p_\s+V\p_2) (\Jg\Wbn  )
\notag \\
&
=\nbs^k\big(- \alpha    g^{- {\frac{1}{2}} }\Sigma  \nbs_2(\Jg\Abn) 
+ \alpha    g^{- {\frac{1}{2}} } \Sigma \Abn \nbs_2 \Jg
+ \tfrac{\alpha}{2} g^{- {\frac{1}{2}} } \Sigma(\Jg\Wbn + \Jg\Zbn - 2\Jg\Abt)\nbs_2 \tt\cdot\nn
 +\Sigma\Fwn\Big)\,.
\end{align*} 
Using  \eqref{bootstraps}, \eqref{geometry-bounds-new}, \eqref{eq:Jg:Abn:D5:improve},  \eqref{eq:Jg:Zbn:D5:improve}, \eqref{eq:D5:JgWbn},    \eqref{D2-Jg-Linfty}, \eqref{se3:time},  and \eqref{eq:Lynch:1}, we obtain the bounds
\begin{align}
&\snorm{\tfrac{\mathcal{J}^{\frac 34}}{\Sigma^{\beta-1}} \doublecom{\nbs^6, \Sigma^{-1}  , (\Q\p_\s+V\p_2)  (\Jg\Wbn)} }_{L^2_{x,\s}} 
\notag\\
&\le (\tfrac{4}{\kappa_0})^\beta  \| \Sigma\nbs (\Jg \Sigma^{-1}) \|_{L^\infty_{x,\s}} \snorm{\nbs^5(\Jg \Wbn)}_{L^2_{x,\s}}  
+ \Cn (\tfrac{4}{\kappa_0})^\beta \displaystyle{\sum}_{k=1}^3 \snorm{ \nbs^{4-k} \nbs(\Jg \Sigma^{-1}) \nbs^{k}\nbs(\Jg\Wbn)}_{L^2_{x,\s}} 
+ \Cn \eps (\tfrac{4}{\kappa_0})^\beta \brak{\mathsf{B}_6} 
\notag\\
&\le \Cn (\tfrac{4}{\kappa_0})^\beta \snorm{\nbs^5(\Jg \Wbn)}_{L^2_{x,\s}} 
+ \Cn(\tfrac{4}{\kappa_0})^\beta \displaystyle{\sum}_{k=1}^3 \snorm{\nbs^{4-k}  \nbs(\Jg \Sigma^{-1})}_{L^{\frac{8}{4-k}}_{x,\s}} \snorm{\nbs^{k} \nbs(\Jg \Wbn)}_{L^{\frac{8}{k}}_{x,\s}} 
+ \Cn \eps (\tfrac{4}{\kappa_0})^\beta \brak{\mathsf{B}_6}    
\notag\\
&\le 
 \Cn   (\tfrac{4}{\kappa_0})^\beta  \brak{\mathsf{B}_6} 
\,.  \label{I-Wbn6-comm1-bound}
\end{align}
The bound~\eqref{I-Wbn6-comm1-bound} is typical of how the most difficult terms in \eqref{eq:I:W:nn:6:d} are estimated.

We next consider the bounds for the forcing, remainder, and commutator functions associated to the equation \eqref{energy-Zn-s}.
A straightforward application of \eqref{bootstraps} and \eqref{eq:Lynch:1} shows that
\begin{subequations}
\label{eq:I:Z:nn:10:all}
\begin{equation}
\snorm{\tfrac{\mathcal{J}^{\frac 34} }{\Sigma^{\beta-1}} \nbs^6\Fzn}_{L^2_{x,\s}}
\les  (\tfrac{4}{\kappa_0})^\beta  \mathsf{K} \brak{\mathsf{B}_6}
\,.
\label{eq:I:Z:nn:10:b}
\end{equation}
Using  the definition of $ \mathcal{R}_\Zb^\nn$ in \eqref{Cz-Rz-comm}, we shall prove that
\begin{equation}
\snorm{\tfrac{\mathcal{J}^{\frac 34} }{\Sigma^{\beta-1}}\mathcal{R}_\Zb^\nn}_{L^2_{x,\s}}
\leq \Cn    (\tfrac{4}{\kappa_0})^\beta  \mathsf{K} \brak{\mathsf{B}_6}
+   \tfrac{4 (1+ \alpha )}{\eps} \snorm{ \tfrac{ \mathcal{J}^{\frac 34}}{\Sigma^{\beta}}   \nbs^6 (\Jg \Zbn)}_{L^2_{x,\s}}
\label{eq:I:Z:nn:10:c}
\,.
\end{equation}
Lastly, using the definition of $\mathcal{C}_\Zb^\nn$ in \eqref{Cz-Rz-comm}, identity~\eqref{eq:Jg:Wb:nn},   the aforementioned  bounds,  Lemma~\ref{lem:anisotropic:sobolev}, the Leibniz rule and Lemma~\ref{lem:time:interpolation},  we may also obtain that
\begin{align}
\snorm{\tfrac{\mathcal{J}^{\frac 34} }{\Sigma^{\beta-1}}\mathcal{C}_\Zb^\nn}_{L^2_{x,\s}}
\les    (\tfrac{4}{\kappa_0})^\beta  \mathsf{K} \brak{\mathsf{B}_6}
\,.
\label{eq:I:Z:nn:10:d}
\end{align}
\end{subequations}

To prove \eqref{eq:I:Z:nn:10:c}, we  write $ \mathcal{R}_\Zb^\nn = ( \mathcal{R}_\Zb^\nn- \tfrac{2\alpha }{\eps} \nbs_1\Jg \nbs^6\Zbn) + \tfrac{2\alpha }{\eps} \nbs_1\Jg \nbs^6\Zbn$.  From the bounds \eqref{bootstraps}, \eqref{eq:Q:all:bbq}, and \eqref{eq:Lynch:1}, we obtain that
\begin{equation} 
\snorm{\tfrac{\mathcal{J}^{\frac 34} }{\Sigma^{\beta-1}} \big( \mathcal{R}_\Zb^\nn- \tfrac{2\alpha }{\eps} \nbs_1\Jg \nbs^6\Zbn\big) }_{L^2_{x,\s}}
\les  (\tfrac{4}{\kappa_0})^\beta  \mathsf{K} \brak{\mathsf{B}_6} \,.
\label{eq:I:Z:nn:10:dd}
\end{equation} 
In order to estimate $ \tfrac{2\alpha }{\eps}\snorm{\tfrac{\mathcal{J}^{\frac 34} }{\Sigma^{\beta-1}} \nbs_1\Jg \nbs^6\Zbn }_{L^2_{x,\s}}$,
we consider the three cases that (a) $\nbs^6= \nbs_1 \nbs^5$, (b) $\nbs^6= \nbs_2 \nbs^5$, and (c) $\nbs^6= \nbs_\s^6$.
For case (a), 
from the bound \eqref{eq:Jg:Zbn:p1D5:improve}, we have that 
\begin{equation} 
\snorm{\tfrac{1}{\Sigma^{\beta-1}}  \Fgtf \nbs_1 \nbs^5 \Zbn }_{L^2_{x,\s}} 
\le \tfrac{1}{2\alpha} \snorm{ \tfrac{1}{\Sigma^{\beta}} \mathcal{J}^{\frac 34}  \nbs^6 (\Jg \Zbn)}_{L^2_{x,\s}}
+ \Cn \eps (\tfrac{4}{\kappa_0})^\beta   \mathsf{K} \brak{\mathsf{B_6}}\,, \label{Shnoblinsky-maneuver1}
\end{equation} 
Then, using \eqref{bs-Jg,1} and \eqref{Shnoblinsky-maneuver1}, we have that
\begin{align} 
 \tfrac{2\alpha }{\eps}\snorm{\tfrac{\mathcal{J}^{\frac 34} }{\Sigma^{\beta-1}} \nbs_1\Jg \nbs_1\nbs^5\Zbn }_{L^2_{x,\s}}
& \le  \tfrac{2 \alpha }{\eps}   \snorm{ \nbs_1 \Jg }_{L^\infty_{x,\s}}  
\snorm{ \tfrac{1}{\Sigma^{\beta-1}} \mathcal{J}^{\frac 34} \nbs_1 \nbs^5\Zbn }_{L^2_{x,\s}}
\notag \\
& \le
\tfrac{4 (1+ \alpha )}{\eps} \snorm{ \tfrac{1}{\Sigma^{\beta}}\mathcal{J}^{\frac 34} \nbs^6 (\Jg \Zbn)}_{L^2_{x,\s}}
+ \Cn \eps (\tfrac{4}{\kappa_0})^\beta    \mathsf{K} \brak{\mathsf{B_6}}
 \,.
 \label{10hard1}
\end{align} 
For case (b),  it follows from \eqref{eq:Jg:Zbn:p2D5:improve:new} that
\begin{equation} 
\snorm{\tfrac{1}{\Sigma^{\beta-1}}  \mathcal{J}^{\frac 34}  \nbs_2 \nbs^5 \Zbn }_{L^2_{x,\s}}  
\le \tfrac{1}{\alpha  \eps} \|\tfrac{1}{\Sigma^{\beta}} \mathcal{J}^{\frac 34}   \nbs^6\Zbt \|_{L^2_{x,\s}} 
+ \Cn \eps (\tfrac{4}{\kappa_0})^\beta \mathsf{K} \brak{\mathsf{B_6}} 
\le \Cn  (\tfrac{4}{\kappa_0})^\beta \mathsf{K} \brak{\mathsf{B_6}}   \,,
\label{Shnoblinsky-maneuver2}
\end{equation}
and hence using \eqref{Shnoblinsky-maneuver2}, we find that
\begin{equation} 
\tfrac{2\alpha }{\eps}\snorm{\tfrac{\mathcal{J}^{\frac 34} }{\Sigma^{\beta-1}} \nbs_1\Jg \nbs_2 \nbs^5\Zbn }_{L^2_{x,\s}}
\le \Cn  (\tfrac{4}{\kappa_0})^\beta \mathsf{K}\brak{\mathsf{B_6}} \,.
\label{10hard2}
\end{equation} 
Finally for case (c), from \eqref{eq:madman:2}, we also have that
\begin{equation} 
\tfrac{2\alpha }{\eps}\snorm{\tfrac{\mathcal{J}^{\frac 34} }{\Sigma^{\beta-1}}\nbs_1 \Jg \nbs_\s^6\Zbn }_{L^2_{x,\s}}
\le \Cn (\tfrac{4}{\kappa_0})^\beta \mathsf{K}\brak{\mathsf{B_6}} \,.
\label{10hard3}
\end{equation} 
Combining \eqref{eq:I:Z:nn:10:dd}, \eqref{10hard1}, \eqref{10hard2}, and  \eqref{10hard3} proves the inequality \eqref{eq:I:Z:nn:10:c}.
A similar decomposition of $ \mathcal{C}_\Zb^\nn$ yields \eqref{eq:I:Z:nn:10:d}.

The bounds for the forcing, remainder, and commutator functions for the $\Abn$-equation \eqref{energy-An-s} are obtained in the identical 
fashion as for the $\Zbn$-equation, and we have that 
\begin{subequations}
\label{eq:I:A:nn:10:all}
\begin{align}
\snorm{\tfrac{\mathcal{J}^{\frac 34} }{\Sigma^{\beta-1}} \nbs^6\Fan}_{L^2_{x,\s}}
&\leq \Cn  (\tfrac{4}{\kappa_0})^\beta  \mathsf{K} \brak{\mathsf{B}_6}
\,,  \label{eq:I:A:nn:10:b} \\
\snorm{\tfrac{\mathcal{J}^{\frac 34} }{\Sigma^{\beta-1}}\mathcal{R}_\Ab^\nn}_{L^2_{x,\s}}
&\leq \Cn    (\tfrac{4}{\kappa_0})^\beta  \mathsf{K} \brak{\mathsf{B}_6}
+   \tfrac{4 (1+ \alpha )}{\eps} \snorm{ \tfrac{1}{\Sigma^{\beta}} \mathcal{J}^{\frac 34}  \nbs^6 (\Jg \Abn)}_{L^2_{x,\s}} \,,
\label{eq:I:A:nn:10:c} \\
\snorm{\tfrac{\mathcal{J}^{\frac 34} }{\Sigma^{\beta-1}}\mathcal{C}_\Ab^\nn}_{L^2_{x,\s}}
&\leq \Cn    (\tfrac{4}{\kappa_0})^\beta  \mathsf{K} \brak{\mathsf{B}_6}
\,.
\label{eq:I:A:nn:10:d}
\end{align}
\end{subequations}

\subsubsection{The sum $I^{\WW_n}_5+I^{\ZZ_n}_5+I^{\AA_n}_7$}
\label{sec:IW5+IZ5+IA7}
We first note that by using $\Sbn= \tfrac{1}{2} (\Wbn-\Zbn)$, we have that
\begin{equation} 
I^{\WW_n}_5+I^{\ZZ_n}_5 =
-\alpha  \tint \jb \Fgss  \Jg g^{- {\frac{1}{2}} } (\Jg\Wbn + \Jg\Zbn - 2\Jg\Abt)\nbs_2 \nbs^6 \tt\cdo\nn \ \nbs^6(\Jg\Sbn) \,.
\label{IWn5+IZn5}
\end{equation} 
For the integral $I^{\AA_n}_7$ in \eqref{I7-Abn}, we use integration-by-parts; in particular, with \eqref{adjoint-1} and \eqref{adjoint-2}, to obtain that
\begin{align} 
I^{\AA_n}_7& = I^{\AA_n}_{7,a} + I^{\AA_n}_{7,b} + I^{\AA_n}_{7,c} +  I^{\AA_n}_{7,d }  \,, 
\label{eq:IAAn:7:decompose}\\
I^{\AA_n}_{7,a} &=
-{\tfrac{\alpha }{\eps}} \tint  \jb  \Fgss (\Jg\Wbn+\Jg\Zbn - 2\Jg\Abt) \nbs^6\tt \cdo\nn  
\ \Big( \nbs_1\nbs^6(\Jg\Abn) -  \eps\Jg g^{- {\frac{1}{2}} } \nbs_2h\, \nbs_2\nbs^6(\Jg\Abn) \Big) \,, 
\notag \\
I^{\AA_n}_{7,b} &=
-{\tfrac{\alpha }{\eps}} \tint  \Big(\nbs_1 - \eps \Jg g^{- {\frac{1}{2}} }  \nbs_2 h \nbs_2 \Big)\Big( \jb  \Fgss (\Jg\Wbn+\Jg\Zbn - 2\Jg\Abt) \nn_i \Big) \nbs^6\tt_i  \ \nbs^6(\Jg\Abn) \,, 
\notag \\
I^{\AA_n}_{7,c} &=
- \alpha \tint  \jb \Fgss  (\Jg\Wbn+\Jg\Zbn - 2\Jg\Abt) \bigl(\Qr_2 - \nbs_2\bigr)\bigl(\Jg g^{- {\frac{1}{2}} } \nbs_2 h\bigr) \, \nbs^6\tt \cdo\nn  \ \nbs^6(\Jg\Abn) \,, 
\notag \\
I^{\AA_n}_{7,d} &=
\alpha  \int   \Jg g^{- {\frac{1}{2}} } \nbs_2 h \, \Qb_2 \  \jb\Fgss   (\Jg\Wbn+\Jg\Zbn - 2\Jg\Abt) \nbs^6\tt \cdo\nn  \ \nbs^6(\Jg\Abn)\Big|_\s 
\notag\,.
\end{align} 

First, let us observe that by employing the bounds \eqref{bootstraps}, \eqref{eq:Q:all:bbq},  and  \eqref{geometry-bounds-new},  we have that
\begin{equation} 
\sabs{I^{\AA_n}_{7,b} } +\sabs{I^{\AA_n}_{7,c} } + \sabs{I^{\AA_n}_{7,d} } 
\les  (\tfrac{4}{\kappa_0})^{2\beta}  \mathsf{K} \Bsix^2 \,.
\label{I7-An-bcd}
\end{equation}

To study the integral $I^{\AA_n}_{7,a}$, we use equation \eqref{energy-An-s} to substitute for 
$\tfrac{\alpha}{\eps} \Big( \nbs_1\nbs^6(\Jg\Abn) -  \eps\Jg g^{- {\frac{1}{2}} } \nbs_2h\, \nbs_2\nbs^6(\Jg\Abn) \Big)$ 
and we find  that
\begin{subequations} 
\label{IAbn-7a-i-viii}
\begin{align} 
I^{\AA_n}_{7,a}&= I^{\AA_n}_{7,a,i}+ I^{\AA_n}_{7,a,ii} + I^{\AA_n}_{7,a,iii} + I^{\AA_n}_{7,a,iv} + I^{\AA_n}_{7,a,v}
+ I^{\AA_n}_{7,a,vi} + I^{\AA_n}_{7,a,vii} + I^{\AA_n}_{7,a,viii}  \,,
\notag \\
I^{\AA_n}_{7,a,i} &= 
 - \alpha  \tint  \jb \Fgss \Jg g^{- {\frac{1}{2}} }   (\Jg\Wbn+\Jg\Zbn - 2\Jg\Abt) \nbs^6\tt \cdo\nn  \     \nbs_2\nbs^6(\Jg\Sbn) \,,
\label{IAbn-7a-i} \\
I^{\AA_n}_{7,a,ii} &= - \tint  \tfrac{1}{\Sigma^{2\beta}} \Fgss   (\Jg\Wbn+\Jg\Zbn - 2\Jg\Abt) \nbs^6\tt \cdo\nn  \  \Jg(\Q\p_\s +V\p_2) \nbs^6(\Jg\Abn)  \,,
\label{IAbn-7a-ii}  \\
I^{\AA_n}_{7,a,iii} &=  \tint  \tfrac{1}{\Sigma^{2\beta}}   \Fgss (\Jg\Wbn+\Jg\Zbn - 2\Jg\Abt) \nbs^6\tt \cdo\nn  
\ \bubu{ ( \tfrac{\alpha}{2} \Jg\Wbn - \tfrac{\alpha}{2} \Jg \Zbn)}\nbs^6(\Jg\Abn) \,,
\label{IAbn-7a-iii}  \\
I^{\AA_n}_{7,a,iv} &= 
 \alpha  \tint  \jb \Fgss  (\Jg\Wbn+\Jg\Zbn - 2\Jg\Abt) \nbs^6\tt \cdo\nn  \    g^{- {\frac{1}{2}} }  (\Jg\Sbn) \nbs_2\nbs^6 \Jg \,,
\label{IAbn-7a-iv}  \\
I^{\AA_n}_{7,a,v} &= 
-\alpha  \tint  \jb  \Fgss (\Jg\Wbn+\Jg\Zbn - 2\Jg\Abt) \nbs^6\tt \cdo\nn  \    g^{- {\frac{1}{2}} } \Jgh  \Sbt \nbs_2\nbs^6 \tt \cdot \nn
\,,
\label{IAbn-7a-v}  \\
I^{\AA_n}_{7,a,vi} &= 
-{\tfrac{\alpha }{2\eps}} \tint  \jb \Fgss  (\Jg\Wbn+\Jg\Zbn - 2\Jg\Abt)^2 \nbs^6\tt \cdo\nn 
 \ \Jg \Big( \nbs_1\nbs^6\tt \cdot\nn - \eps \Jg g^{- {\frac{1}{2}} } \nbs_2h\, \nbs_2\nbs^6\tt \cdot\nn  \Big) \,,
\label{IAbn-7a-vi} \\
I^{\AA_n}_{7,a,vii} &= 
-{\tfrac{\alpha }{\eps}} \tint  \jb \Fgss  \Abn (\Jg\Wbn+\Jg\Zbn - 2\Jg\Abt) \nbs^6\tt \cdo\nn  
\ \Big(  \nbs_1\nbs^6 \Jg - \eps \Jg g^{- {\frac{1}{2}} } \nbs_2h\, \nbs_2\nbs^6 \Jg \Big)  \,,
\label{IAbn-7a-vii}\\
I^{\AA_n}_{7,a,viii} &= 
 \tint  \jb \Fgss  (\Jg\Wbn+\Jg\Zbn - 2\Jg\Abt) \nbs^6\tt \cdo\nn  \  \big( \nbs^6 \Fan  + \mathcal{R}_\Ab^\nn + \mathcal{C}_\Ab^\nn\big) \,.
 \label{IAbn-7a-viii}
\end{align} 
\end{subequations} 
Again, using \eqref{bootstraps} and  \eqref{geometry-bounds-new}, we see that
\begin{equation}
\sabs{I^{\AA_n}_{7,a,iii} } \les  (\tfrac{4}{\kappa_0})^{2\beta}  \mathsf{K} \Bsix^2  \,.
\label{small-chicken1}
\end{equation}
By additionally using \eqref{eq:I:A:nn:10:all}, 
we also have that
\begin{equation}
\sabs{I^{\AA_n}_{7,a,viii} } \les   (\tfrac{4}{\kappa_0})^{2\beta}  \mathsf{K} \Bsix^2  \,.
\label{small-chicken2}
\end{equation}
Next, using exact-derivative structure and integration by parts, and appealing to Lemmas \ref{lem:tau-Jg-D2} and \ref{lem:tau-Jg-D1}, we
have that
\begin{equation}
\sabs{I^{\AA_n}_{7,a,iv} } +\sabs{I^{\AA_n}_{7,a,v} } +\sabs{I^{\AA_n}_{7,a,vi} } +\sabs{I^{\AA_n}_{7,a,vii} } 
\les  (\tfrac{4}{\kappa_0})^{2\beta}  \mathsf{K}^2 \Bsix^2  \,.
\label{small-chicken3}
\end{equation} 
Let us next explain the procedure for bounding the integral $I^{\AA_n}_{7,a,ii}$.   We integrate by parts with respect to the operator
$ (\Q\p_\s+V\p_2) $ and use \eqref{adjoint-3} with \eqref{good-comm}, \eqref{Jg-evo-s},  and \eqref{tt-evo-s} to  find that
\begin{align} 
I^{\AA_n}_{7,a,ii}&= I^{\AA_n}_{7,a,ii_1} + I^{\AA_n}_{7,a,ii_2} + I^{\AA_n}_{7,a,ii_3} + I^{\AA_n}_{7,a,ii_4} +  I^{\AA_n}_{7,a,ii_5} +   I^{\AA_n}_{7,a,ii_6} + I^{\AA_n}_{7,a,ii_7} \,,
\label{eq:I:AAn:7:ii:decomp} \\
I^{\AA_n}_{7,a,ii_1} 
& = \tint  \tfrac{1}{\Sigma^{2\beta}} \Fgss  \Jg (\Jg\Wbn+\Jg\Zbn - 2\Jg\Abt) \  \nbs^6 \Big((\tfrac{1+\alpha}{2} \Wbt
+ \tfrac{1-\alpha}{2} \Zbt)\nn \Big)\cdo\nn  \   \nbs^6(\Jg\Abn)\,,
\notag \\
I^{\AA_n}_{7,a,ii_2} 
& = \tint  \tfrac{1}{\Sigma^{2\beta}} \Fgss  \big(\tfrac{1+\alpha}{2} \Jg\Wbn + \tfrac{1-\alpha}{2} \Jg\Zbn \big) \ 
(\Jg\Wbn+\Jg\Zbn - 2\Jg\Abt) \  \nbs^6\tt \cdo\nn  \   \nbs^6(\Jg\Abn)\,,
\notag \\
I^{\AA_n}_{7,a,ii_3} 
& = \tint  \Jg  (\Q\p_\s+V\p_2) \big( \tfrac{1}{\Sigma^{2\beta}} \Fgss   \ (\Jg\Wbn+\Jg\Zbn - 2\Jg\Abt) \nn_i\big) \  \nbs^6\tt_i  \   \nbs^6(\Jg\Abn) \,,
\notag \\
I^{\AA_n}_{7,a,ii_4} 
& = \tint  \Jg  (\nbs_2 V +  \Qr_\s  - V \Qr_2)  \tfrac{1}{\Sigma^{2\beta}} \Fgss   \ (\Jg\Wbn+\Jg\Zbn - 2\Jg\Abt) \  \nbs^6\tt \cdo\nn  \   \nbs^6(\Jg\Abn) \,,
\notag \\
I^{\AA_n}_{7,a,ii_5} 
& = -\tint  \tfrac{1}{\Sigma^{2\beta}} \Fgss  \Jg (\Jg\Wbn+\Jg\Zbn - 2\Jg\Abt) \ \nbs^6(\Jg\Abn) \  \big( \nbs^6V \nbs_2\tt
+\doublecom{\nbs^6,V, \nbs_2 \tt} \big)\cdo\nn   
\,,
\notag \\
I^{\AA_n}_{7,a,ii_6} 
& = - \int   \tfrac{\Q}{\Sigma^{2\beta}}  \Jg  \Fgss   \ (\Jg\Wbn+\Jg\Zbn - 2\Jg\Abt) \  \nbs^6\tt \cdo\nn  \   \nbs^6(\Jg\Abn) \Big|_\s 
\notag \\
I^{\AA_n}_{7,a,ii_7} 
& =   \int   \tfrac{\Q}{\Sigma^{2\beta}}  \Jg  \Fgss   \ (\Jg\Wbn+\Jg\Zbn - 2\Jg\Abt) \  \nbs^6\tt \cdo\nn  \   \nbs^6(\Jg\Abn) \Big|_0
\notag
\,.
\end{align} 
By using \eqref{table:derivatives}, \eqref{bootstraps}, \eqref{eq:Q:all:bbq},  and \eqref{geometry-bounds-new} together with
\eqref{eq:Lynch:1} and \eqref{eq:Lynch:2}, all of the seven  integrals above are directly estimated and we obtain that
\begin{equation} 
\sabs{ I^{\AA_n}_{7,a,ii} } \les (\tfrac{4}{\kappa_0})^{2\beta}  \mathsf{K}^2 \Bsix^2 \,.
\label{small-chicken5}
\end{equation}

As we have shown  in \eqref{small-chicken1}--\eqref{small-chicken5}, 
the integrals $I^{\AA_n}_{7,a,ii}$ through $I^{\AA_n}_{7,a,viii}$, defined in \eqref{IAbn-7a-ii}--\eqref{IAbn-7a-viii}, are bounded rather directly. In contrast,  the integral
$I^{\AA_n}_{7,a,i}$ defined in~\eqref{IAbn-7a-i} cannot be bounded, but rather requires a cancellation to occur.   In particular, employing integration-by-parts for 
this integral, we will arrive at a cancellation with the sum $I^{\WW_n}_5+I^{\ZZ_n}_5$ given in \eqref{IWn5+IZn5}.   
By once again using \eqref{adjoint-2}, we have that
\begin{align} 
I^{\AA_n}_{7,a,i} &= I^{\AA_n}_{7,a,i_1} + I^{\AA_n}_{7,a,i_2} + I^{\AA_n}_{7,a,i_3}+ I^{\AA_n}_{7,a,i_4} \,,
\notag \\
I^{\AA_n}_{7,a,i_1} & = 
\alpha  \tint  \jb \Fgss \Jg g^{- {\frac{1}{2}} }   (\Jg\Wbn+\Jg\Zbn - 2\Jg\Abt) \nbs_2\nbs^6\tt \cdo\nn  \     \nbs^6(\Jg\Sbn) \,,
\notag \\
I^{\AA_n}_{7,a,i_2} & = 
\alpha  \tint \nbs_2\Big(  \jb\Fgss  \Jg g^{- {\frac{1}{2}} }   (\Jg\Wbn+\Jg\Zbn - 2\Jg\Abt)\nn_i \Big)  \nbs^6\tt_i \     \nbs^6(\Jg\Sbn) \,,
\notag \\
I^{\AA_n}_{7,a,i_3} & = 
- \alpha  \tint  \Qr_2  \jb\Fgss  \Jg g^{- {\frac{1}{2}} }   (\Jg\Wbn+\Jg\Zbn - 2\Jg\Ab) \nbs^6\tt \cdo\nn  \     \nbs^6(\Jg\Sbn) \,,
\notag \\
I^{\AA_n}_{7,a,i_4} & = 
 \alpha   \int  \Qb_2 \jb \Fgss \Jg g^{- {\frac{1}{2}} }   (\Jg\Wbn+\Jg\Zbn - 2\Jg\Abt) \nbs^6\tt \cdo\nn  \     \nbs^6(\Jg\Sbn)\Big|_\s \,.
 \label{IAn7,a,i}
\end{align} 
The first term in the above decomposition precisely equals $-(I^{\WW_n}_5+I^{\ZZ_n}_5)$. 
Therefore, we sum  \eqref{IWn5+IZn5} with \eqref{IAn7,a,i} and find that
\begin{equation} 
I^{\WW_n}_5+I^{\ZZ_n}_5+I^{\AA_n}_{7,a,i} = I^{\AA_n}_{7,a,i_2}+I^{\AA_n}_{7,a,i_3}+I^{\AA_n}_{7,a,i_4} \,.
\label{IWn5+IZn5+IAn7cancel}
\end{equation} 
Using that $\Sbn = \tfrac{1}{2} (\Wbn-\Zbn)$, the bounds \eqref{bootstraps}, \eqref{eq:Q:all:bbq},  and \eqref{geometry-bounds-new} show that
\begin{equation} 
\sabs{I^{\AA_n}_{7,a,i_2}}+ \sabs{ I^{\AA_n}_{7,a,i_3}} + \sabs{ I^{\AA_n}_{7,a,i_4}}
\les \eps  (\tfrac{4}{\kappa_0})^{2\beta}  \mathsf{K}^2 \Bsix^2 \,.
\label{small-chicken4}
\end{equation} 
Combining \eqref{I7-An-bcd}, \eqref{small-chicken1}--\eqref{small-chicken5},
\eqref{IWn5+IZn5+IAn7cancel}, and~\eqref{small-chicken4}, we arrive at
\begin{equation} 
\sabs{I^{\WW_n}_5+I^{\ZZ_n}_5+I^{\AA_n}_{7} }\les (\tfrac{4}{\kappa_0})^{2\beta}  \mathsf{K}^2 \Bsix^2  \,.
\label{IWn5+IZn5+IAn7}
\end{equation} 

\subsubsection{The sum $I^{\WW_n}_4+I^{\ZZ_n}_4+I^{\AA_n}_8$}
\label{sec:IW4+IZ4+IA8}
Summing \eqref{I4-Wbn} and \eqref{I4-Zbn}, we 
 define $I^{{\scriptscriptstyle \WW_n+\ZZ_n}}_4  := I^{\WW_n}_4  +I^{\ZZ_n}_4 $, and 
again using \eqref{adjoint-2}, we obtain that
\begin{align*} 
I^{{\scriptscriptstyle \WW_n+\ZZ_n}}_4 
&=  
I^{{\scriptscriptstyle \WW_n+\ZZ_n}}_{4,a}+I^{{\scriptscriptstyle \WW_n+\ZZ_n}}_{4,b}
+I^{{\scriptscriptstyle \WW_n+\ZZ_n}}_{4,c} \,, 
 \notag \\
I^{{\scriptscriptstyle \WW_n+\ZZ_n}}_{4,a} 
&=
2\alpha \tint \jb  \Fgss \Abn \nbs^6 \Jg \  g^{- {\frac{1}{2}} }\Jg \nbs_2 \nbs^6(\Jg\Sbn) \,,
\notag \\
I^{{\scriptscriptstyle \WW_n+\ZZ_n}}_{4,b}
&=
2\alpha \tint  (\nbs_2-\Qr_2)(\jb\Fgss \Jg g^{- {\frac{1}{2}} }  \Abn) \nbs^6 \Jg \   \nbs^6(\Jg\Sbn) \,,
\notag \\
I^{{\scriptscriptstyle \WW_n+\ZZ_n}}_{4,c}
&=
2\alpha \int  \Qb_2 \jb \Fgss \Jg g^{- {\frac{1}{2}} }  \Abn \nbs^6 \Jg \   \nbs^6(\Jg\Sbn)\Big|_\s \,.
\end{align*} 
The bounds \eqref{bootstraps}, \eqref{eq:Q:all:bbq},  and \eqref{geometry-bounds-new} show that
\begin{equation}
 \sabs{I^{{\scriptscriptstyle \WW_n+\ZZ_n}}_{4,b}}
 +
 \sabs{I^{{\scriptscriptstyle \WW_n+\ZZ_n}}_{4,c}}
 \les 
 \eps
 (\tfrac{4}{\kappa_0})^{2\beta} \Bsix^2
 \,.
 \label{eq:small:duck:1}
\end{equation}
In order to estimate the integral $I^{{\scriptscriptstyle \WW_n+\ZZ_n}}_{4,a}$, which contains a term with seven derivatives, we shall again make use of the equation
\eqref{energy-An-s} and substitute for $\alpha g^{-\frac 12} \Jg \nbs_2 \nbs^6(\Jg\Sbn)$ to obtain that
\begin{subequations} 
\label{I4-Wbn+Zbn-a}
\begin{align} 
I^{{\scriptscriptstyle \WW_n+\ZZ_n}}_{4,a}
& = 
I^{{\scriptscriptstyle \WW_n+\ZZ_n}}_{4,a,i} 
+ I^{{\scriptscriptstyle \WW_n+\ZZ_n}}_{4,a,ii}
+ I^{{\scriptscriptstyle \WW_n+\ZZ_n}}_{4,a,iii} 
+ I^{{\scriptscriptstyle \WW_n+\ZZ_n}}_{4,a,iv} 
+ I^{{\scriptscriptstyle \WW_n+\ZZ_n}}_{4,a,v}
+ I^{{\scriptscriptstyle \WW_n+\ZZ_n}}_{4,a,vi} 
+ I^{{\scriptscriptstyle \WW_n+\ZZ_n}}_{4,a,vii} + I^{{\scriptscriptstyle \WW_n+\ZZ_n}}_{4,a,viii}\,,
\notag \\
I^{{\scriptscriptstyle \WW_n+\ZZ_n}}_{4,a,i} &= 
-2 \tint \tfrac{1}{\Sigma^{2\beta}} \Fgss \Jg \Abn \nbs^6 \Jg \  (\Q\p_\s +V\p_2) \nbs^6(\Jg\Abn)   \,,
 \\
I^{{\scriptscriptstyle \WW_n+\ZZ_n}}_{4,a,ii} &= 
2 \tint \tfrac{1}{\Sigma^{2\beta}} \Fgss  \Abn \nbs^6 \Jg \   \bubu{ ( \tfrac{\alpha}{2} \Jg\Wbn - \tfrac{\alpha}{2} \Jg \Zbn)} \nbs^6(\Jg\Abn)  \,,
  \\
I^{{\scriptscriptstyle \WW_n+\ZZ_n}}_{4,a,iii} &= 
\alpha  \tint \jb  \Fgss g^{- {\frac{1}{2}} }  (\Jg\Sbn) \Abn  \  \nbs_2 \Bigl( \bigl( \nbs^6 \Jg \bigr)^2 \Bigr) \,,
  \\
I^{{\scriptscriptstyle \WW_n+\ZZ_n}}_{4,a,iv} &= 
-2 \alpha  \tint \jb \Fgss g^{- {\frac{1}{2}} }  \Jgt   \Sbt \Abn \nbs^6 \Jg \   \nbs_2\nbs^6 \tt \cdot \nn\,,
  \\
I^{{\scriptscriptstyle \WW_n+\ZZ_n}}_{4,a,v} &= 
{\tfrac{2\alpha }{\eps}}   \tint \jb \Fgss  \Abn \nbs^6 \Jg \  \Big(  \nbs_1\nbs^6 (\Jg\Abn) - \eps \Jg g^{- {\frac{1}{2}} } \nbs_2 h\,  \nbs_2\nbs^6 (\Jg\Abn)\Big)  \,,
\label{I4-Wbn+Zbn-a:v}  \\
I^{{\scriptscriptstyle \WW_n+\ZZ_n}}_{4,a,vi} &= 
-{\tfrac{\alpha }{\eps}}   \tint \jb \Fgss  \Abn (\Jg\Wbn +\Jg\Zbn - 2\Jg \Abt)  \nbs^6 \Jg
\  \Big(  \nbs_1\nbs^6\tt \cdot\nn -  \eps\Jg g^{- {\frac{1}{2}} } \nbs_2h\, \nbs_2\nbs^6\tt \cdot\nn \Big) \,,
  \\
I^{{\scriptscriptstyle \WW_n+\ZZ_n}}_{4,a,vii} &= 
-{\tfrac{ \alpha }{\eps}}   \tint \jb \Fgss \Abn^2   \Big( \nbs_1  - \eps \Jg g^{- {\frac{1}{2}} } \nbs_2 h\, \nbs_2 \Big) \Bigl( \bigl(\nbs^6 \Jg \bigr)^2 \Bigr) \,,
 \\
I^{{\scriptscriptstyle \WW_n+\ZZ_n}}_{4,a,viii} &= 
2 \tint \jb \Fgss \Abn \nbs^6 \Jg \  ( \nbs^6 \Fan  + \mathcal{R}_\Ab^\nn + \mathcal{C}_\Ab^\nn) \,.
\end{align} 
\end{subequations} 
The eight terms present in \eqref{I4-Wbn+Zbn-a} may all be bounded directly, except for $I^{{\scriptscriptstyle \WW_n+\ZZ_n}}_{4,a,v}$, which must be cancelled with a contribution from $ I^{\AA_n}_8$, see~\eqref{eq:ducky:von:ducken:1} below. Indeed, from \eqref{bootstraps}, \eqref{eq:Q:all:bbq}, \eqref{geometry-bounds-new}, \eqref{eq:I:A:nn:10:all} we deduce
\begin{equation}
\sabs{I^{{\scriptscriptstyle \WW_n+\ZZ_n}}_{4,a,ii}}
+
\sabs{I^{{\scriptscriptstyle \WW_n+\ZZ_n}}_{4,a,viii}}
\les  (\tfrac{4}{\kappa_0})^{2\beta} \Bsix^2
  \label{eq:small:duck:2}
  \,.
\end{equation}
By also appealing to Lemmas~\ref{lem:tau-Jg-D2} and~\ref{lem:tau-Jg-D1} we have
\begin{equation}
\sabs{I^{{\scriptscriptstyle \WW_n+\ZZ_n}}_{4,a,iv}}
+
\sabs{I^{{\scriptscriptstyle \WW_n+\ZZ_n}}_{4,a,vi}}
\les
\eps (\tfrac{4}{\kappa_0})^{2\beta} \mathsf{K}^2 \brak{\mathsf{B}_6}^2  
\,.
\label{eq:small:duck:3}
\end{equation}
By furthermore integrating-by-parts the exact derivative structure via~\eqref{adjoint-1} and~\eqref{adjoint-2}, we deduce
\begin{equation}
\sabs{I^{{\scriptscriptstyle \WW_n+\ZZ_n}}_{4,a,iii}}
+
\sabs{I^{{\scriptscriptstyle \WW_n+\ZZ_n}}_{4,a,vii}}
\les 
\eps (\tfrac{4}{\kappa_0})^{2\beta}  \brak{\mathsf{B}_6}^2
\,.
  \label{eq:small:duck:4}
\end{equation}
Lastly, for the first term in \eqref{I4-Wbn+Zbn-a}, we integrate-by-parts using \eqref{adjoint-3} and appeal to \eqref{good-comm} and \eqref{Jg-evo-s}, obtaining
\begin{equation}
  \label{eq:small:duck:5}
\sabs{I^{{\scriptscriptstyle \WW_n+\ZZ_n}}_{4,a,i} }
\les   (\tfrac{4}{\kappa_0})^{2\beta}  \brak{\mathsf{B}_6}^2
\,.
\end{equation}
This concludes the bounds for the straightforward terms in \eqref{I4-Wbn+Zbn-a}.

We next use \eqref{adjoint-1} and  \eqref{adjoint-2} to write  $I^{\AA_n}_8$ as
\begin{subequations}
 \label{I4-Wbn+Zbn-a:minus}
\begin{align} 
 I^{\AA_n}_8& =  I^{\AA_n}_{8,a} + I^{\AA_n}_{8,b} + I^{\AA_n}_{8,c} + I^{\AA_n}_{8,d} \,,
 \notag \\
 I^{\AA_n}_{8,a}&= 
-{\tfrac{2 \alpha }{\eps}} \tint  \jb \Fgss  \Abn \nbs^6\Jg \ \Big( \nbs_1\nbs^6(\Jg\Abn) -  \eps \Jg g^{- {\frac{1}{2}} } \nbs_2 h\,  \nbs_2\nbs^6(\Jg\Abn) \Big)  \,,
 \label{I4-Wbn+Zbn-a:minus:v} \\
I^{\AA_n}_{8,b}&= 
-{\tfrac{2 \alpha }{\eps}} \tint  \Big( \nbs_1(\jb\Fgss\Abn)- \eps \nbs_2\big( \Fgss \Jg g^{- {\frac{1}{2}} } \nbs_2 h\, \nbs_2(\jb\Abn) \big)  \Big)
 \nbs^6\Jg \ \nbs^6(\Jg\Abn)  \,,
  \\
I^{\AA_n}_{8,c}&= 
-2 \alpha   \tint  \jb \Fgss \Abn   \Jg g^{- {\frac{1}{2}} } \nbs_2 h\,  \Qr_2 
 \nbs^6\Jg \ \nbs^6(\Jg\Abn)  \,,
  \\
I^{\AA_n}_{8,d}&= 
2 \alpha   \int  \jb \Fgss\Abn    \Jg g^{- {\frac{1}{2}} } \nbs_2 h\,  \Qb_2 
 \nbs^6\Jg \ \nbs^6(\Jg\Abn) \Big|_\s \,.
\end{align} 
\end{subequations}
For the last three terms above, by appealing to \eqref{bootstraps}, \eqref{eq:Q:all:bbq}, and \eqref{geometry-bounds-new}, we deduce
\begin{equation}
  \label{eq:small:duck:6}
\sabs{I^{\AA_n}_{8,b}}
+
\sabs{I^{\AA_n}_{8,c}}
+
\sabs{I^{\AA_n}_{8,d}}
\les   (\tfrac{4}{\kappa_0})^{2\beta}  \brak{\mathsf{B}_6}^2
\,,
\end{equation}
while for the first term we note the cancelation structure obtained by adding \eqref{I4-Wbn+Zbn-a:v} and \eqref{I4-Wbn+Zbn-a:minus:v}:
\begin{equation}
\label{eq:ducky:von:ducken:1} 
I^{{\scriptscriptstyle \WW_n+\ZZ_n}}_{4,a,v} +I^{\AA_n}_{8,a}=0 \,.
\end{equation}
From~\eqref{I4-Wbn+Zbn-a}, \eqref{I4-Wbn+Zbn-a:minus}, and~\eqref{eq:ducky:von:ducken:1} we thus obtain
\begin{equation*} 
I^{{\scriptscriptstyle \WW_n+\ZZ_n}}_{4,a} +  I^{\AA_n}_8
= I^{{\scriptscriptstyle \WW_n+\ZZ_n}}_{4,a,i} + I^{{\scriptscriptstyle \WW_n+\ZZ_n}}_{4,a,ii} + I^{{\scriptscriptstyle \WW_n+\ZZ_n}}_{4,a,iii}
+ I^{{\scriptscriptstyle \WW_n+\ZZ_n}}_{4,a,iv} + I^{{\scriptscriptstyle \WW_n+\ZZ_n}}_{4,a,vi} + I^{{\scriptscriptstyle \WW_n+\ZZ_n}}_{4,a,vii}
+ I^{{\scriptscriptstyle \WW_n+\ZZ_n}}_{4,a,viii}
+I^{\AA_n}_{8,b} +I^{\AA_n}_{8,c} +I^{\AA_n}_{8,d}
\,,
\end{equation*} 
which may be combined with the bounds \eqref{eq:small:duck:1}, \eqref{eq:small:duck:2}--\eqref{eq:small:duck:5}, and \eqref{eq:small:duck:6} to deduce
\begin{equation}
\sabs{I^{\WW_n}_4+I^{\ZZ_n}_4+I^{\AA_n}_8} 
\les  
(\tfrac{4}{\kappa_0})^{2\beta}  \brak{\mathsf{B}_6}^2
\,.
\label{eq:small:duck:0}
\end{equation}

\subsubsection{The integral $I^{\AA_n}_4$}
\label{sec:IA4}
For the integral $ I^{\AA_n}_4$ defined in~\eqref{I4-Abn} we first integrate-by-parts the $\nbs_2$ derivative.  Using \eqref{adjoint-2}, we have that
\begin{align*} 
 I^{\AA_n}_{4}&=  I^{\AA_n}_{4,a} +  I^{\AA_n}_{4,b} +  I^{\AA_n}_{4,c}  \,,
 \notag \\
I^{\AA_n}_{4,a} &= 
2\alpha \tint \jb  \Fgss ( \Jg \Sbn) \nbs^6 \Jg \  g^{- {\frac{1}{2}} } \nbs^6\nbs_2(\Jg\Abn) \,,
\notag \\
I^{\AA_n}_{4,b} &= 
2\alpha \tint (\nbs_2 - \Qr_2)( \jb \Fgss g^{- {\frac{1}{2}} }  \Jg \Sbn ) \nbs^6 \Jg \   \nbs^6(\Jg\Abn) \,,
\notag \\
I^{\AA_n}_{4,c} &= 
2\alpha   \int  \Qb_2  \jb  g^{- {\frac{1}{2}} }  \Fgss \Jg \Sbn  \nbs^6 \Jg \   \nbs^6(\Jg\Abn) \Big|_\s \,.
\end{align*} 
We bound  $I^{\AA_n}_{4,b}$ and $I^{\AA_n}_{4,c}$ in a straightforward manner by using \eqref{bootstraps}, \eqref{eq:Q:all:bbq},   \eqref{geometry-bounds-new}, and the bound $\mathcal{J} \leq \Jg$, to obtain
\begin{align}
\sabs{I^{\AA_n}_{4,b}} + \sabs{I^{\AA_n}_{4,c}}
\les 
(\tfrac{4}{\kappa_0})^{2\beta}  \brak{\mathsf{B}_6}^2
\,.
\label{eq:small:duck:7} 
\end{align}
For the integral $I^{\AA_n}_{4,a} $, we use equation \eqref{energy-Wn-s} to substitute for  $\alpha g^{- {\frac{1}{2}} } \nbs^6\nbs_2(\Jg\Abn)$.   This leads to an
additive decomposition for $I^{\AA_n}_{4,a} $ which we write as
\begin{subequations} 
\label{I4-Abn-a}
\begin{align} 
I^{\AA_n}_{4,a} & = I^{\AA_n}_{4,a,i} + I^{\AA_n}_{4,a,ii} + I^{\AA_n}_{4,a,iii} + I^{\AA_n}_{4,a,iv} \,,
\notag \\
 I^{\AA_n}_{4,a,i}&= - 2\tint \tfrac{1}{\Sigma^{2\beta}}  \Fgss  ( \Jg \Sbn) \nbs^6 \Jg \ (\Q\p_\s +V\p_2) \nbs^6(\Jg\Wbn  ) \,,
\label{I4-Abn-a,i} \\
 I^{\AA_n}_{4,a,ii}&= 
 \alpha \tint \jb\Fgss    (\Jg \Sbn)  \  g^{- {\frac{1}{2}} } \Abn \nbs_2 \Bigl( \bigl(\nbs^6 \Jg\bigr)^2 \Bigr) \,,
\label{I4-Abn-a,ii} \\
 I^{\AA_n}_{4,a,iii}&= 
\alpha \tint \jb  \Fgss  (\Jg \Sbn) \  g^{- {\frac{1}{2}} } (\Jg\Wbn + \Jg\Zbn - 2\Jg\Abt)\nbs^6 \Jg  \ \nbs_2 \nbs^6 \tt\cdo\nn \,,
\label{I4-Abn-a,iii} \\
 I^{\AA_n}_{4,a,iv}&= 
2 \tint \jb  \Fgss ( \Jg \Sbn) \nbs^6 \Jg \   ( \nbs^6\Fwn + \mathcal{R}_\Wb^\nn + \mathcal{C}_\Wb^\nn) \,.
\label{I4-Abn-a,iv}
\end{align} 
\end{subequations} 
The last three terms in \eqref{I4-Abn-a} may be estimated directly. Indeed, from \eqref{bootstraps}, \eqref{eq:Q:all:bbq}, \eqref{geometry-bounds-new}, \eqref{eq:I:W:nn:6:all}, and the inequality $\mathcal{J} \leq \Jg$ we deduce
\begin{subequations}
\label{eq:small:duck:8}
\begin{equation}
\sabs{I^{\AA_n}_{4,a,iv}}
\les    (\tfrac{4}{\kappa_0})^{2\beta}  \brak{\mathsf{B}_6}^2
\,,
\end{equation}
by additionally integrating by parts the $\nbs_2$ term via \eqref{adjoint-2} we obtain
\begin{equation}
\sabs{I^{\AA_n}_{4,a,ii}}
\les \eps   (\tfrac{4}{\kappa_0})^{2\beta}  \brak{\mathsf{B}_6}^2
\,,
\end{equation}
and by also appealing to Lemma~\ref{lem:tau-Jg-D2}  we have
\begin{equation}
\sabs{I^{\AA_n}_{4,a,iii}}
\les \eps  \mathsf{K} (\tfrac{4}{\kappa_0})^{2\beta}  \brak{\mathsf{B}_6}^2
\,.
\end{equation}
\end{subequations}

We now focus on the integral $ I^{\AA_n}_{4,a,i}$ defined in \eqref{I4-Abn-a,i} which produces an anti-damping term.  By once again using integration-by-parts via \eqref{adjoint-3}, we see that, 
\begin{subequations}
\label{I4-Abn-a-i}
\begin{align} 
 I^{\AA_n}_{4,a,i} &= J^{\AA_n}_{1} + J^{\AA_n}_{2}+ J^{\AA_n}_{3}+ J^{\AA_n}_{4}+ J^{\AA_n}_{5}+ J^{\AA_n}_{6}+ J^{\AA_n}_{7}
 + J^{\AA_n}_{8} \,
\notag \\
J^{\AA_n}_{1} &  = 2 \tint \tfrac{1}{\Sigma^{2\beta}}  \Fgss   (\Jg\Sbn) (\Q\p_\s+V\p_2) \nbs^6\Jg  \  \nbs^6(\Jg\Wbn)  \,,
  \\
J^{\AA_n}_{2} &  =2 \tint \tfrac{1}{\Sigma^{2\beta}}  \Fgss  (\Q\p_\s+V\p_2) (\Jg\Sbn)  \nbs^6\Jg  \  \nbs^6(\Jg\Wbn)   \,,
  \\
J^{\AA_n}_{3} &  =  \tint \tfrac{1}{\Sigma^{2\beta}} (\Q\p_\s+V\p_2) \Fgss  \  (\Jg\Wbn)  \nbs^6\Jg  \  \nbs^6(\Jg\Wbn)  \,,
  \\
J^{\AA_n}_{4} &  =  2\tint      (\Q\p_\s +V\p_2)\big(  \tfrac{1}{\Sigma^{2\beta}}\big) \ \Fgss (\Jg \Sbn) \   \nbs^6 \Jg \  \nbs^6(\Jg\Wbn  ) 
  \\
J^{\AA_n}_{5} & = 2 \tint \tfrac{1}{\Sigma^{2\beta}}  \Fgss   (\Jg \Sbn)( \Qr_\s  -V \Qr_2 + \nbs_2 V) \nbs^6 \Jg \ \nbs^6(\Jg\Wbn  ) \,,
  \\
J^{\AA_n}_{6} & =
-  \int \tfrac{\Q}{\Sigma^{2\beta}} \Fgss    \Jg \Wbn \nbs^6 \Jg \  \nbs^6(\Jg\Wbn  ) \Big|_\s
 \,,
  \\
J^{\AA_n}_{7} & =
\int \tfrac{\Q}{\Sigma^{2\beta}} \Fgss    \Jg \Zbn \nbs^6 \Jg \  \nbs^6(\Jg\Wbn  ) \Big|_\s
  \\  
J^{\AA_n}_{8} & =
 2 \int \tfrac{\Q}{\Sigma^{2\beta}} \Fgss    \Jg \Sbn \nbs^6 \Jg \  \nbs^6(\Jg\Wbn  ) \Big|_0 \,,
  \\    
J^{\AA_n}_{9} &  = - \tint \tfrac{1}{\Sigma^{2\beta}} (\Q\p_\s+V\p_2) \Fgss  \  (\Jg\Zbn)  \nbs^6\Jg  \  \nbs^6(\Jg\Wbn)  \,.
\end{align} 
\end{subequations}
Most terms in the decomposition~\eqref{I4-Abn-a-i} are estimated handily using \eqref{bootstraps}, \eqref{eq:Q:all:bbq}, \eqref{geometry-bounds-new}, \eqref{eq:Jg:Wbn:improve:material:a} as 
\begin{align}
&\sabs{J^{\AA_n}_{2}}
+ 
\sabs{J^{\AA_n}_{4}}
+
\sabs{J^{\AA_n}_{5}}
+
\sabs{J^{\AA_n}_{7}}
+
\sabs{J^{\AA_n}_{9}}
\notag\\
&\qquad
\les 
(\tfrac{4}{\kappa_0})^{2\beta}  \brak{\mathsf{B}_6}^2
+ 
\tfrac{500^2}{\eps^2 \sqrt{1+\alpha}} 
\int_0^{\s} 
\snorm{\tfrac{\mathcal{J}^{\frac 34} (\Jg \Q)^{\frac 12}}{\Sigma^\beta} \nbs^6 (\Jg \Wbn)(\cdot,\s')}_{L^2_x} 
\snorm{\tfrac{\Q \mathcal{J}^{\frac 14}}{\Sigma^\beta} \nbs^6 \Jg(\cdot,\s')}_{L^2_x} 
{\rm d}\s'
\,,
\label{eq:small:duck:8a} 
\end{align}
and by also appealing to \eqref{sigma0-bound}, \eqref{table:derivatives} and \eqref{eq:Qcal:bbq:temp:3} with $\s=0$, we have
\begin{equation}
\sabs{J^{\AA_n}_{8}}
\leq \bigl( \tfrac{1}{\eps} +\Cn \bigr) \tfrac{1+\alpha}{2} (1+\Cn \eps)
(\tfrac{3}{\kappa_0})^{2\beta}  \Cdatatwo
\leq \tfrac{1+\alpha}{\eps} (\tfrac{3}{\kappa_0})^{2\beta}  \Cdatatwo
\,.
\label{eq:small:duck:8b} 
\end{equation}

The remaining terms, namely $J^{\AA_n}_{1}$, $J^{\AA_n}_{3}$, and $J^{\AA_n}_{6}$, need to be handled carefully.
The integral $J^{\AA_n}_{1}$ produces an anti-damping term that must be combined with the last integral on the right
side of \eqref{eq:heavy:fuel:1n}.   Using \eqref{good-comm}, \eqref{Jg-evo-s},  a short computation shows that
\begin{align} 
J^{\AA_n}_{1} & = J^{\AA_n}_{1,a} + J^{\AA_n}_{1,b}+ J^{\AA_n}_{1,c}  \,, 
\label{eq:small:duck:big:duck:bear} \\
J^{\AA_n}_{1,a} &  = \tint \tfrac{1}{\Sigma^{2\beta}}  \Fgss   (\Q\p_\s+V\p_2) \Jg\    \sabs{  \nbs^6(\Jg\Wbn)  }^2
\notag \\
J^{\AA_n}_{1,b} &  = \tfrac{1-\alpha}{2}
 \tint \tfrac{1}{\Sigma^{2\beta}}  \Fgss    (\Jg\Wbn) \nbs^6(\Jg\Zbn)  
\nbs^6(\Jg\Wbn) \,,
\notag \\
J^{\AA_n}_{1,c} &  = 
- \tint \tfrac{1}{\Sigma^{2\beta}}  \Fgss   (\Jg \Zbn) \bigl( \nbs^6(\Jg\Wbn) 
\bigr)^2 \,,
\notag \\
J^{\AA_n}_{1,d} &  = 
- 2 \tint \tfrac{1}{\Sigma^{2\beta}}   \Fgss  ( \Jg \Sbn)  \big(\nbs^6 V \nbs_2\Jg + \doublecom{\nbs^6,V, \nbs_2\Jg} \big)   \  \nbs^6(\Jg\Wbn  ) 
\notag \,.
\end{align} 
Using \eqref{bootstraps}, \eqref{eq:Q:all:bbq}, \eqref{geometry-bounds-new}, and Lemma~\ref{lem:comm:tangent} we deduce
\begin{subequations}
\label{eq:small:duck:big:duck}
\begin{equation}
\sabs{J^{\AA_n}_{1,c}} + \sabs{J^{\AA_n}_{1,d}}
\les (\tfrac{4}{\kappa_0})^{2\beta}  \brak{\mathsf{B}_6}^2 
+ \mathsf{K} \eps (\tfrac{4}{\kappa_0})^{2\beta}  \brak{\mathsf{B}_6}^2
\les (\tfrac{4}{\kappa_0})^{2\beta}  \brak{\mathsf{B}_6}^2
\label{eq:small:duck:9} 
\,.
\end{equation}
Moreover, by Cauchy-Schwartz and \eqref{bs-JgnnWb} we deduce 
\begin{equation}
\sabs{J^{\AA_n}_{1,b}}
\leq \tfrac{1+\alpha}{\eps} 
\int_0^{\s}  \snorm{\tfrac{\mathcal{J}^{\frac 34}}{\Sigma^\beta} \nbs^6(\Jg\Zbn) (\cdot,\s')}_{L^2_x}
\snorm{\tfrac{\mathcal{J}^{\frac 34}}{\Sigma^\beta} \nbs^6(\Jg\Wbn)(\cdot,\s')}_{L^2_x} {\rm d} \s'
\,,
\label{eq:small:duck:10}  
\end{equation}
whereas 
\begin{equation}
J^{\AA_n}_{1,a}
= \tint \tfrac{1}{\Sigma^{2\beta}}   \mathsf{G_{bad}}    \sabs{  \nbs^6(\Jg\Wbn)  }^2
\label{eq:small:duck:11}  
\end{equation}
where the anti-damping term $\mathsf{G_{bad}} = \mathcal{J}^{\frac 32} (\Q\p_\s + V \p_2) \Jg$ is precisely the same as in the case of the tangential derivatives (see~\eqref{eq:G:bad:tangential}). This last term will be combined with the damping term $\mathsf{G_{good}}$ present on the right side of \eqref{eq:I:n:12369}. For this purpose, precisely in the same way as for the tangential derivatives (see~\eqref{eq:G:good:G:bad}), we record the inequality
\begin{equation}
 \mathsf{G_{good}} + \mathsf{G_{bad}}
 \geq \tfrac{1+\alpha}{16\eps} \mathcal{J}^{\frac 12} \Jg
 \,,
\label{eq:small:duck:11:a}
\end{equation}
\end{subequations}
which follows directly from \eqref{eq:fakeJg:LB}.

Among the nine terms in \eqref{I4-Abn-a-i} it thus remains to consider the integral $J^{\AA_n}_{3} $ (which produces a new type of damping and energy norm which are displayed, respectively, in the integrals  $J^{\AA_n}_{3,a}$ and $J^{\AA_n}_{3,b}$ below), and the boundary term $J^{\AA_n}_{6} $ (which is to be dealt with by a careful application of the $\eps$-Young inequality). It is important to make an analogy between $J^{\AA_n}_{3}$ and $J^{\AA_n}_{6} $, respectively the terms $\tilde{\mathsf{M}}_2$ an $\tilde{\mathsf{M}}_3$ which have previously appeared in the tangential energy estimates, cf.~\eqref{eq:Really:Bad:2:def} and~\eqref{eq:Really:Bad:3:def}. Our treatment of the terms $J^{\AA_n}_{3}$ and $J^{\AA_n}_{6} $ closely resembles the analysis in Subsection~\ref{sec:nasty:tangential}.

Concerning $J^{\AA_n}_{3,a}$, we note that by the definition of $\mathcal{J}$ in \eqref{t-to-s-transform}, and upon using~\eqref{Jg-evo-s} and \eqref{good-comm} to rewrite
\begin{align*}
\nbs^6 (\Jg \Wbn) 
&= \tfrac{2}{1+\alpha} \bigl(\nbs^6 (\Q \p_\s + V\p_2) \Jg - \tfrac{1-\alpha}{2} \nbs^6 (\Jg \Zbn) \bigr)
\\
&=  \tfrac{2}{1+\alpha}(\Q \p_\s + V\p_2)  \nbs^6  \Jg + \tfrac{2}{1+\alpha} \bigl(\nbs^6 V \nbs_2 \Jg + \doublecom{\nbs^6,V,\nbs_2 \Jg} \bigr)- \tfrac{1-\alpha}{1+\alpha} \nbs^6 (\Jg \Zbn)  
\,,
\end{align*}
and finally, integrating-by-parts the material derivative using~\eqref{adjoint-3}, we arrive at
\begin{align} 
J^{\AA_n}_{3} 
&=- \tfrac{3}{2\eps} \tint \tfrac{\Q}{\Sigma^{2\beta}} \mathcal{J}^{\frac 12}   (\Jg\Wbn)  \nbs^6\Jg  \  \nbs^6(\Jg\Wbn) 
\notag\\
&= J^{\AA_n}_{3,a} + J^{\AA_n}_{3,b}+ J^{\AA_n}_{3,c} + J^{\AA_n}_{3,d} + J^{\AA_n}_{3,e}+ J^{\AA_n}_{3,f}+ J^{\AA_n}_{3,g}+ J^{\AA_n}_{3,h}+ J^{\AA_n}_{3,i}\,, 
\label{eq:J3:An:decompose}\\
J^{\AA_n}_{3,a}
&=  \tfrac{3}{2(1+ \alpha )}  {\tfrac{1}{\eps}} \int \tfrac{\Q^2}{\Sigma^{2\beta}} \mathcal{J}^{\frac 12} (-\Jg\Wbn + \tfrac{13}{\eps} \Jg )\sabs{  \nbs^6\Jg}^2\Big|_\s \,,
\notag \\
J^{\AA_n}_{3,b}
&=  \tfrac{3}{4(1+ \alpha )}   \tfrac{1}{\eps^2}  \tint \tfrac{\Q^2}{\Sigma^{2\beta}} \mathcal{J}^{-\frac 12} (-\Jg\Wbn + \tfrac{13}{\eps} \Jg) \sabs{\nbs^6\Jg}^2 \,,
\notag \\
J^{\AA_n}_{3,c}
&=  \tfrac{3}{2(1+ \alpha )}  {\tfrac{1}{\eps}} \tint \tfrac{\Q}{\Sigma^{2\beta}} \Fgh (\Q\p_\s+V\p_2) (\Jg\Wbn - \tfrac{13}{\eps} \Jg) 
\ \sabs{ \nbs^6\Jg}^2 \,,
\notag \\
J^{\AA_n}_{3,d}
&=  \tfrac{3}{2(1+ \alpha )}  {\tfrac{1}{\eps}} \tint  (\Q\p_\s+V\p_2) \big( \tfrac{\Q}{\Sigma^{2\beta}}\big) \Fgh (\Jg\Wbn -   \tfrac{13}{\eps} \Jg) 
\ \sabs{ \nbs^6\Jg}^2 \,,
\notag \\
J^{\AA_n}_{3,e}
&=  \tfrac{3}{2(1+ \alpha )}  {\tfrac{1}{\eps}} \tint    \tfrac{\Q}{\Sigma^{2\beta}} ( \Qr_\s  + \nbs_2 V - V \Qr_2 ) \Fgh (\Jg\Wbn -   \tfrac{13}{\eps} \Jg) 
\ \sabs{ \nbs^6\Jg}^2 \,,
\notag \\
 J^{\AA_n}_{3,f} 
&=  \tfrac{3}{2(1+ \alpha )}  {\tfrac{1}{\eps}} \int \tfrac{\Q^2}{\Sigma^{2\beta}} \mathcal{J}^{\frac 12} (\Jg\Wbn -  \tfrac{13}{\eps} \Jg)\sabs{ \nbs^6\Jg}^2\Big|_0
\,,
\notag\\
J^{\AA_n}_{3,g}
&= - \tfrac{39}{2\eps^2} \tint \tfrac{\Q}{\Sigma^{2\beta}} \mathcal{J}^{\frac 12}   \Jg  \nbs^6\Jg  \  \nbs^6(\Jg\Wbn) 
\,,
\notag\\
J^{\AA_n}_{3,h}
&= - \tfrac{3}{(1+\alpha) \eps} \tint \tfrac{\Q}{\Sigma^{2\beta}} \mathcal{J}^{\frac 12}   (\Jg\Wbn -  \tfrac{13}{\eps} \Jg)  \nbs^6\Jg  \  \bigl(\nbs^6 V \nbs_2 \Jg + \doublecom{\nbs^6,V,\nbs_2 \Jg} \bigr) \,,
\notag \\
J^{\AA_n}_{3,i}
&=  \tfrac{3(1-\alpha)}{2(1+\alpha)\eps} \tint \tfrac{\Q}{\Sigma^{2\beta}} \mathcal{J}^{\frac 12}   (\Jg\Wbn  -  \tfrac{13}{\eps} \Jg )  \nbs^6\Jg  \  \nbs^6(\Jg\Zbn)
\notag
\,.
\end{align} 
As in the tangential case, the point of the add-and-subtract of the term $ \frac{13}{\eps} \Jg$ is to allow for the applicability of \eqref{eq:signed:Jg}, which yields (in direct analogy to \eqref{eq:how:the:fuck?:a} and~\eqref{eq:how:the:fuck?:e}) 
\begin{subequations}
\label{eq:fuck:yeah:0}
\begin{align}
J^{\AA_n}_{3,a}
&\geq 
\tfrac{27}{20(1+\alpha)} \tfrac{1}{\eps^2} \snorm{\tfrac{\Q \mathcal{J}^{\frac 14}}{\Sigma^\beta} \nbs^6 \Jg(\cdot,\s)}_{L^2_x}^2 
\,,
\label{eq:fuck:yeah:1}
\\
J^{\AA_n}_{3,b}
&\geq 
\tfrac{27}{40(1+\alpha)} \tfrac{1}{\eps^3} \int_0^{\s} \snorm{\tfrac{\Q \mathcal{J}^{-\frac 14}}{\Sigma^\beta} \nbs^6 \Jg(\cdot,\s')}_{L^2_x}^2 {\rm d} \s'
\,.
\label{eq:fuck:yeah:2}
\end{align}
For the remaining terms in the additive decomposition of $J^{\AA_n}_{3}$, using the bounds \eqref{table:derivatives}, \eqref{bootstraps}, \eqref{eq:Q:all:bbq}, \eqref{eq:signed:Jg}, \eqref{geometry-bounds-new},  \eqref{eq:Jg:Wbn:improve:material:a}, and the identities $(\Q \p_\s + V \p_2) \Q = \Q  \Qc  + V \p_2 \Q$ and \eqref{Jg-evo-s}, to obtain the bounds 
\begin{align}
J^{\AA_n}_{3,c}
&\geq \tfrac{3}{2(1+ \alpha )}  {\tfrac{1}{\eps}} \tint \tfrac{\Q}{\Sigma^{2\beta}} \Fgh (\Q\p_\s+V\p_2) (\Jg\Wbn)
\ \sabs{ \nbs^6\Jg}^2
- \tfrac{39(1-\alpha)}{4(1+ \alpha )}  {\tfrac{1}{\eps^2}} \tint \tfrac{\Q}{\Sigma^{2\beta}} \Fgh ( \Jg \Zbn)  \sabs{ \nbs^6\Jg}^2
\notag\\
&\qquad
+ \tfrac{39 \cdot 9}{40}  {\tfrac{1}{\eps^3}} \tint \tfrac{\Q}{\Sigma^{2\beta}} \Fgh    \sabs{ \nbs^6\Jg}^2
- \tfrac{3 \cdot 13^2}{4}  {\tfrac{1}{\eps^3}} \tint \tfrac{\Q}{\Sigma^{2\beta}} \Fgh  \Jg \sabs{ \nbs^6\Jg}^2
\notag\\
&\geq
\bigl(8 - \Cn \eps \bigr)  {\tfrac{1}{\eps^3}}  \tint \tfrac{\Q}{\Sigma^{2\beta}} \Fgh    \sabs{ \nbs^6\Jg}^2
-  \tfrac{39^2}{4 (1+\alpha)}  {\tfrac{1}{\eps^3}} \int_0^{\s}   \snorm{\tfrac{\Q \mathcal{J}^{\frac 14}}{\Sigma^\beta} \nbs^6 \Jg(\cdot,\s')}_{L^2_x}^2 {\rm d} \s'
\label{eq:fuck:yeah:3}
\\
\sabs{J^{\AA_n}_{3,d}}
+
\sabs{J^{\AA_n}_{3,e}}
&\leq \tfrac{40 \cdot 250^2 + \eps \Cn \brak{\beta}}{\eps^3} \tint \tfrac{\Q}{\Sigma^{2\beta}} \Fgh    \sabs{ \nbs^6\Jg}^2
\label{eq:fuck:yeah:4}
\,,
\\
\sabs{J^{\AA_n}_{3,f}}
&\leq \tfrac{3}{2(1+\alpha)\eps} \tfrac{(1+\alpha)^2}{3} \cdot \tfrac{14}{\eps} (\tfrac{3}{\kappa_0})^{2\beta} \eps \Cdatatwo
\leq
\tfrac{ 7(1+\alpha) }{\eps}  (\tfrac{3}{\kappa_0})^{2\beta}   \Cdatatwo
\label{eq:fuck:yeah:5}
\,,
\\
\sabs{J^{\AA_n}_{3,g}}
&\leq \tfrac{39}{2\eps^2} \|\Jgh \Q^{-\frac 12}\|_{L^\infty_{x,\s}} \int_0^{\s}
\snorm{\tfrac{\Q \mathcal{J}^{-\frac 14}}{\Sigma^\beta} \nbs^6 \Jg(\cdot,\s')}_{L^2_x} 
\snorm{\tfrac{\mathcal{J}^{\frac 34} (\Jg\Q)^{\frac 12}}{\Sigma^{\beta}} \nbs^6 (\Jg \Wbn)(\cdot,\s')}_{L^2_x} {\rm d}\s' 
\notag\\
&\leq 
\tfrac{1}{4(1+\alpha)} \tfrac{1}{\eps^3} \int_0^{\s} \snorm{\tfrac{\Q \mathcal{J}^{-\frac 14}}{\Sigma^\beta} \nbs^6 \Jg(\cdot,\s')}_{L^2_x}^2 {\rm d} \s'
+ 
\tfrac{3 \cdot 39^2}{4} \tfrac{1}{\eps}
\int_0^{\s} 
\snorm{\tfrac{\mathcal{J}^{\frac 34} (\Jg\Q)^{\frac 12}}{\Sigma^{\beta}} \nbs^6 (\Jg \Wbn)(\cdot,\s')}_{L^2_x}^2 {\rm d}\s' 
\label{eq:fuck:yeah:6}
\,,
\\
\sabs{J^{\AA_n}_{3,h}}
&\les \eps \mathsf{K} (\tfrac{4}{\kappa_0})^{2\beta} \brak{\mathsf{B}_6}^2
\label{eq:fuck:yeah:7}
\,,
\\
\sabs{J^{\AA_n}_{3,i}}
&\leq  \tfrac{ 25 }{\eps^2} \int_0^{\s}  
\snorm{\tfrac{\mathcal{J}^{\frac 34}}{\Sigma^\beta} \nbs^6(\Jg\Zbn) (\cdot,\s')}_{L^2_x}
\snorm{\tfrac{\Q \mathcal{J}^{-\frac 14}}{\Sigma^\beta} \nbs^6 \Jg (\cdot,\s')}_{L^2_x}
{\rm d} \s'
\notag\\
&\leq \tfrac{1}{4(1+\alpha)} \tfrac{1}{\eps^3} 
\int_0^{\s} 
\snorm{\tfrac{\Q \mathcal{J}^{-\frac 14}}{\Sigma^\beta} \nbs^6 \Jg(\cdot,\s')}_{L^2_x}^2 
{\rm d} \s'
+ 
\tfrac{ 25^2 (1+\alpha)}{\eps} 
\int_0^{\s} 
\snorm{\tfrac{\mathcal{J}^{\frac 34}}{\Sigma^\beta} \nbs^6(\Jg\Zbn) (\cdot,\s')}_{L^2_x}^2 
{\rm d} \s'
\,.
\label{eq:fuck:yeah:8}
\end{align}
\end{subequations}
The bounds in~\eqref{eq:fuck:yeah:0}  complete our estimates for the nine terms in the decomposition of $J^{\AA_n}_{3}$. 

It remains to bound the term~$J^{\AA_n}_{6}$ in~\eqref{I4-Abn-a-i}. This term requires a special application of the $\eps$-Young inequality, akin to \eqref{eq:Really:Bad:3:final} for the tangential estimates. More precisely, by \eqref{bs-JgnnWb}, the lower bound in \eqref{Qd-lower-upper}, \eqref{Q-lower-upper} and the bound $\mathcal{J}\leq \Jg$, we have
\begin{align}
\sabs{J^{\AA_n}_{6}}
&\leq \tfrac{\sqrt{5}}{\eps \sqrt{2(1+\alpha)}}
\snorm{\tfrac{\Q \mathcal{J}^{\frac 14}}{\Sigma^\beta} \nbs^6 \Jg(\cdot,\s)}_{L^2_x}
\snorm{\tfrac{\mathcal{J}^{\frac 34} (\Jg \Q)^{\frac 12}}{\Sigma^\beta} \nbs^6(\Jg\Wbn)(\cdot,\s)}_{L^2_x} 
\notag\\
&\leq 
\tfrac{25}{52}
\snorm{\tfrac{\mathcal{J}^{\frac 34} (\Jg \Q)^{\frac 12}}{\Sigma^\beta} \nbs^6(\Jg\Wbn)(\cdot,\s)}^2
+
\tfrac{13}{10(1+\alpha)\eps^2}
\snorm{\tfrac{\Q \mathcal{J}^{\frac 14}}{\Sigma^\beta} \nbs^6 \Jg(\cdot,\s)}_{L^2_x}^2
\,.
\label{eq:fuck:yeah:9}
\end{align}
Again, we recall that the fact $\tfrac{25}{52} < \frac 12$ allows us to close the energy estimate.

We summarize the bounds in this subsection, namely~\eqref{eq:small:duck:7}, \eqref{eq:small:duck:8}, \eqref{eq:small:duck:8a}, \eqref{eq:small:duck:8b}, \eqref{eq:small:duck:big:duck}, \eqref{eq:fuck:yeah:0}, and~\eqref{eq:fuck:yeah:9}, together with $\mathcal{J}\leq \Jg$ and $\eps$-Young, to obtain that
\begin{align}
\label{eq:fuck:yeah:end}
 I^{\AA_n}_{4}
&\geq  
- \Cn (\tfrac{4}{\kappa_0})^{2\beta}  \brak{\mathsf{B}_6}^2
- \tfrac{8(1+\alpha)}{\eps}   (\tfrac{3}{\kappa_0})^{2\beta}  \Cdatatwo
\notag\\
&\qquad 
+ \tint \tfrac{1}{\Sigma^{2\beta}} \bigl(\tfrac{1+\alpha}{24\eps} \mathcal{J}^{\frac 12} \Jg - \mathsf{G_{good}} \bigr)   \sabs{  \nbs^6(\Jg\Wbn)  }^2
- \tfrac{(12 +25^2 )(1+\alpha)}{\eps} 
\int_0^{\s} \snorm{\tfrac{\mathcal{J}^{\frac 34}}{\Sigma^\beta} \nbs^6(\Jg\Zbn) (\cdot,\s')}_{L^2_x}^2 {\rm d} \s'
\notag\\
&\qquad
+ \tfrac{1}{20(1+\alpha)} \tfrac{1}{\eps^2} \snorm{\tfrac{\Q \mathcal{J}^{\frac 14}}{\Sigma^\beta} \nbs^6 \Jg(\cdot,\s)}_{L^2_x}^2 
-  \tfrac{20^2 + 500^2 +100 \cdot 250^2 + \eps \Cn \brak{\beta}}{(1+\alpha)}  {\tfrac{1}{\eps^3}} \int_0^{\s}   \snorm{\tfrac{\Q \mathcal{J}^{\frac 14}}{\Sigma^\beta} \nbs^6 \Jg(\cdot,\s')}_{L^2_x}^2 {\rm d} \s'
\notag\\
&\qquad
+ \tfrac{7}{40(1+\alpha)} \tfrac{1}{\eps^3} \int_0^{\s} \snorm{\tfrac{\Q \mathcal{J}^{-\frac 14}}{\Sigma^\beta} \nbs^6 \Jg(\cdot,\s')}_{L^2_x}^2 {\rm d} \s'
+ \bigl(8 - \Cn \brak{\beta} \eps \bigr)  {\tfrac{1}{\eps^3}}  \tint \tfrac{\Q}{\Sigma^{2\beta}} \Fgh    \sabs{ \nbs^6\Jg}^2
\notag\\
&\qquad
- \tfrac{25}{52}
\snorm{\tfrac{\mathcal{J}^{\frac 34} (\Jg \Q)^{\frac 12}}{\Sigma^\beta} \nbs^6(\Jg\Wbn)(\cdot,\s)}_{L^2_x}^2
- 
\tfrac{ 39^2+ 500^2 }{\eps}  
\int_0^{\s} 
\snorm{\tfrac{\mathcal{J}^{\frac 34} (\Jg\Q)^{\frac 12}}{\Sigma^{\beta}} \nbs^6 (\Jg \Wbn)(\cdot,\s')}_{L^2_x}^2 {\rm d}\s' 
\,.
\end{align}

\subsubsection{The integral $I^{\ZZ_n}_7$} 
\label{sec:I:Zn:7:integral}
We next study the integral $I^{\ZZ_n}_7$ in \eqref{I7-Zbn}.
Using \eqref{eq:adjoints}, we have that
\begin{align} 
 I^{\ZZ_n}_7 &  =   I^{\ZZ_n}_{7,a} +  I^{\ZZ_n}_{7,b}+  I^{\ZZ_n}_{7,c}+  I^{\ZZ_n}_{7,d} \,,
 \label{eq:I:Zn:7:decompose} \\
 I^{\ZZ_n}_{7,a} 
 & = {\tfrac{2 \alpha }{\eps}} \tint  \jb  \Fgss  \Jg(\Abn+\Zbt) \nbs^6\tt\cdo\nn 
 \big( \nbs_1 - \eps \Jg g^{- {\frac{1}{2}} } \nbs_2h \, \nbs_2  \big)\nbs^6(\Jg\Zbn) \,,
\notag \\
 I^{\ZZ_n}_{7,b} 
 & = {\tfrac{2 \alpha }{\eps}} \tint      \big( \nbs_1 - \eps \Jg g^{- {\frac{1}{2}} } \nbs_2h \, \nbs_2  \big)
 \big(  \jb \Fgss  \Jg (\Abn+ \Zbt)\nn_i\big) \  \nbs^6\tt_i
 \nbs^6(\Jg\Zbn) \,,
 \notag \\
 I^{\ZZ_n}_{7,c} 
 & = 2 \alpha \tint  \jb \Fgss \Jg(\Abn + \Zbt)    (\Qr_2- \nbs_2)( \Jg g^{- {\frac{1}{2}} } \nbs_2h) 
\  \nbs^6\tt\cdo\nn  \, \nbs^6(\Jg\Zbn) \,,
 \notag \\
 I^{\ZZ_n}_{7,d} 
 & = -2 \alpha \int  \jb \Qb_2 \Fgss \Jg(\Abn + \Zbt)     \Jg g^{- {\frac{1}{2}} } \nbs_2h
\  \nbs^6\tt\cdo\nn  \, \nbs^6(\Jg\Zbn)\Big|_\s \,.
\notag
\end{align} 
With this additive decomposition, only the integral $ I^{\ZZ_n}_{7,a} $ remains over-differentiated  while the integrals $I^{\ZZ_n}_{7,b}$, $I^{\ZZ_n}_{7,c}$, and $I^{\ZZ_n}_{7,d}$ can be directly estimated using \eqref{bootstraps}, \eqref{eq:Q:all:bbq}, and \eqref{geometry-bounds-new}, as
\begin{equation}
\sabs{I^{\ZZ_n}_{7,b} } 
+
\sabs{I^{\ZZ_n}_{7,c} } 
+
\sabs{I^{\ZZ_n}_{7,d} } 
\les
\mathsf{K} \eps  (\tfrac{4}{\kappa_0})^{2\beta} \brak{\mathsf{B}_6}^2
\label{eq:duck:a:duck:1}
\,.
\end{equation}
For the interesting integral, namely $ I^{\ZZ_n}_{7,a} $, we use equation \eqref{energy-Zn-s} to replace ${\tfrac{2 \alpha }{\eps}} \big( \nbs_1 - \eps \Jg g^{- {\frac{1}{2}} } \nbs_2h \, \nbs_2  \big) \nbs^6 (\Jg \Zbn)$.  This replacement yields integrals possessing either ``exact-derivative'' structure or a derivative reduction which arises from terms containing the operator $ (\Q\p_\s+V\p_2)$.   Using \eqref{energy-Zn-s}, we proceed with this substitution, and find that
\begin{subequations} 
\begin{align} 
I^{\ZZ_n}_{7,a} & = {\textstyle\sum_{i=1}^{8}} L^{\ZZ_n}_{i}  \,, 
\label{I-Zbn-7a-decomp} \\
L^{\ZZ_n}_{1} &= \tint \tfrac{1}{\Sigma^{2\beta}} \Fgss\Jg(\Abn+\Zbt) \nbs^6\tt\cdo\nn \ \Jg (\Q\p_\s+V\p_2) \nbs^6(\Jg\Zbn)\,,
 \label{L-Zbn-1} \\
L^{\ZZ_n}_{2} &= -\tint \tfrac{1}{\Sigma^{2\beta}} \Fgss\Jg(\Abn+\Zbt) \nbs^6\tt\cdo\nn \
\bubu{  \alpha ( \Jg\Wbn - \Jg\Zbn)} \nbs^6(\Jg\Zbn) \,,
 \label{L-Zbn-2} \\
L^{\ZZ_n}_{3} &= -\alpha \tint \jb \Fgss\Jg(\Abn+\Zbt) \nbs^6\tt\cdo\nn \
 \Jg g^{- {\frac{1}{2}} } \nbs_2 \nbs^6(\Jg\Abn) \,,
 \label{L-Zbn-3}\\
L^{\ZZ_n}_{4} &=  \alpha  \tint \jb \Fgss\Jg(\Abn+\Zbt) \nbs^6\tt\cdo\nn \
 \Jg g^{- {\frac{1}{2}} } \Abn \nbs_2 \nbs^6\Jg \,,
 \label{L-Zbn-4}\\
L^{\ZZ_n}_{5} &= \tfrac{\alpha}{2} \tint \jb \Fgss\Jg(\Abn+\Zbt) \nbs^6\tt\cdo\nn \
 \Jg g^{- {\frac{1}{2}} }  (\Jg\Wbn+\Jg\Zbn - 2 \Jg\Abt) \nbs_2 \nbs^6\tt\cdo\nn \,,
 \label{L-Zbn-5}\\
L^{\ZZ_n}_{6} &= -\tfrac{2\alpha}{\eps} \tint \jb \Fgss \big(\Jg(\Abn+\Zbt)\big)^2 \nbs^6\tt\cdo\nn \
 \big( \nbs_1 - \eps g^{- {\frac{1}{2}} } \Jg \nbs_2h \nbs_2\big)(\nbs^6\tt)\cdo\nn \,,
 \label{L-Zbn-6} \\
L^{\ZZ_n}_{7} &= \tfrac{2\alpha}{\eps} \tint \jb \Fgss\Jg(\Abn+\Zbt) \Zbn \nbs^6\tt\cdo\nn \
 \big( \nbs_1 - \eps g^{- {\frac{1}{2}} } \Jg \nbs_2h \nbs_2\big)\nbs^6\Jg\,,
 \label{L-Zbn-7}\\
L^{\ZZ_n}_{8} &= - \tint \jb \Fgss\Jg(\Abn+\Zbt)  \nbs^6\tt\cdo\nn \
 \big(  \nbs^6 \Fzn  + \mathcal{R}_\Zb^\nn + \mathcal{C}_\Zb^\nn \big)\,. 
 \label{L-Zbn-8}
\end{align} 
\end{subequations} 
A few of the terms present in the additive decomposition~\eqref{I-Zbn-7a-decomp} may be bounded directly by using \eqref{bootstraps}, \eqref{eq:Q:all:bbq}, \eqref{geometry-bounds-new}, and~\eqref{eq:I:Z:nn:10:all}:
\begin{subequations}
\label{eq:duck:a:duck:2} 
\begin{equation}
\sabs{L^{\ZZ_n}_{2}}
+
\sabs{L^{\ZZ_n}_{8}}
\les
\mathsf{K} \eps  (\tfrac{4}{\kappa_0})^{2\beta} \brak{\mathsf{B}_6}^2
\,.
\end{equation}
By additionally exploring exact derivative structure, which requires integration by parts of either $\nbs_1$ via~\eqref{adjoint-1} or $\nbs_2$ via~\eqref{adjoint-2}, we may also bound
\begin{equation}
\sabs{L^{\ZZ_n}_{5}} 
+
\sabs{L^{\ZZ_n}_{6}}
\les
\mathsf{K} \eps^3  (\tfrac{4}{\kappa_0})^{2\beta} \brak{\mathsf{B}_6}^2
\,.
\end{equation}
Using the geometric lemmas in Lemma~\ref{lem:tau-Jg-D2} and~\ref{lem:tau-Jg-D1},
we may directly establish that
\begin{equation}
\sabs{L^{\ZZ_n}_{4}} 
+
\sabs{L^{\ZZ_n}_{7}}
\les
\mathsf{K}^2 \eps^2  (\tfrac{4}{\kappa_0})^{2\beta} \brak{\mathsf{B}_6}^2
\,.
\end{equation}
\end{subequations}
The remaining two terms, $L^{\ZZ_n}_{1}$ and $L^{\ZZ_n}_{3}$ require special care. For the $L^{\ZZ_n}_{1}$ integral, we integrate by parts the $(\Q \p_s + V\p_2)$ derivative using~\eqref{adjoint-3}, and use \eqref{eq:hate:10-old} to rewrite the over-differentiated term $(\Q \p_s + V\p_2)(\nn \cdot \nbs^6 \tt)$. We obtain 
\begin{align*}
L^{\ZZ_n}_{1} 
&= L^{\ZZ_n}_{1,a} + L^{\ZZ_n}_{1,b} + L^{\ZZ_n}_{1,c} + L^{\ZZ_n}_{1,d}
\\
L^{\ZZ_n}_{1,a}
&= - \tint \tfrac{1}{\Sigma^{2\beta}} \Fgss\Jgt (\Abn+\Zbt) \nbs^6(\Jg\Zbn)   
\biggl(  \tfrac{1+ \alpha }{2}  \nbs^6 \Wbt  + \tfrac{1- \alpha }{2}  \nbs^6 \Zbt  
\notag\\
&\qquad 
+  \bigl(\tfrac{1+ \alpha }{2} \Wbt  + \tfrac{1- \alpha }{2} \Zbt \bigr) 
\bigl( \nn \cdot   \nbs^6 \nn   - \tt \cdot \nbs^6 \tt \bigr)
+  \nn^k \doublecom{  \nbs^6, \bigl(\tfrac{1+ \alpha }{2} \Wbt  + \tfrac{1- \alpha }{2} \Zbt \bigr), \nn^k}
- \nn^k \jump{ \nbs^6,  V } \nbs_2 \tt^k \biggr),
\\
L^{\ZZ_n}_{1,b}
&= - \tint (\Q\p_\s+V\p_2)  \Bigl(\tfrac{1}{\Sigma^{2\beta}} \Fgss\Jgt (\Abn+\Zbt)\Bigr) \nbs^6\tt\cdo\nn \,   \nbs^6(\Jg\Zbn) ,
\\
L^{\ZZ_n}_{1,c}
&= \tint \bigl(V \Qr_2 - \nbs_2 V -  \Qr_\s \bigr) \tfrac{1}{\Sigma^{2\beta}} \Fgss\Jgt (\Abn+\Zbt) \nbs^6\tt\cdo\nn \, \nbs^6(\Jg\Zbn) ,
\\
L^{\ZZ_n}_{1,d}
&= \int \tfrac{\Q}{\Sigma^{2\beta}} \Fgss\Jgt (\Abn+\Zbt) \nbs^6\tt\cdo\nn  \, \nbs^6(\Jg\Zbn)\Bigr|_{0}^{\s} 
\,.
\end{align*}
From the above decomposition it is clear that no over-differentiated terms are present in the $L^{\ZZ_n}_{1}$ integral, and as such we may use \eqref{bootstraps}, \eqref{eq:Q:all:bbq}, and \eqref{geometry-bounds-new}, to deduce
\begin{equation}
\label{eq:duck:a:duck:3}
\sabs{L^{\ZZ_n}_{1}}
\les \mathsf{K} \eps   (\tfrac{4}{\kappa_0})^{2\beta} \brak{\mathsf{B}_6}^2
\,.
\end{equation}
The last integral we need to consider from \eqref{I-Zbn-7a-decomp} is $L^{\ZZ_n}_{3}$; this term requires further refinement in order to contend with the over-differentiation.  To that end,
we now use equation \eqref{energy-Wn-s} to substitute for the term $ \alpha g^{- {\frac{1}{2}} } \nbs_2 \nbs^6(\Jg\Abn) $, and 
we obtain the further decomposition of this integral as follows: 
\begin{subequations} 
\begin{align} 
L^{\ZZ_n}_{3} & = L^{\ZZ_n}_{3,a}+L^{\ZZ_n}_{3,b}+ L^{\ZZ_n}_{3,c}+ L^{\ZZ_n}_{3,d} \,, 
\notag \\
L^{\ZZ_n}_{3,a} 
&= \tint \tfrac{1}{\Sigma^{2\beta}} \Fgss \Jg(\Abn +\Zbt)  \nbs^6\tt\cdo\nn \Jg (\Q\p_\s+V\p_2) \nbs^6(\Jg\Wbn) \,,
 \label{L-Zbn-3a} \\
L^{\ZZ_n}_{3,b} 
&=-\alpha  \tint \jb\Fgss \Jg(\Abn +\Zbt)  \nbs^6\tt\cdo\nn  g^{- {\frac{1}{2}} } \Jg \Abn \nbs_2 \nbs^6\Jg \,,
 \label{L-Zbn-3b}\\
L^{\ZZ_n}_{3,c} 
&=- \tfrac{\alpha}{2} \tint \jb\Fgss \Jg(\Abn +\Zbt)  \nbs^6\tt\cdo\nn  g^{- {\frac{1}{2}} } \Jg
(\Jg\Wbn +\Jg\Zbn -2\Jg\Abt) \nbs_2\nbs^6\tt\cdo \nn \,, 
 \label{L-Zbn-3c}\\
L^{\ZZ_n}_{3,d} &= - \tint \jb \Fgss\Jg(\Abn+\Zbt)  \nbs^6\tt\cdo\nn \
\Jg \big(  \nbs^6 \Fwn  + \mathcal{R}_\Wb^\nn + \mathcal{C}_\Wb^\nn \big)\,.
  \label{L-Zbn-3d}
\end{align} 
\end{subequations} 
The term $L^{\ZZ_n}_{3,a}$ in \eqref{L-Zbn-3a} is nearly identical to the term $L^{\ZZ_n}_{1}$ defined in~\eqref{L-Zbn-1}: the only difference is that $\nbs^6(\Jg \Zbn)$ is now replaced by $\nbs^6(\Jg \Wbn)$. As such, the bound for $L^{\ZZ_n}_{3,a}$ is the same as for $L^{\ZZ_n}_{1}$, namely~\eqref{eq:duck:a:duck:3}. For the term $L^{\ZZ_n}_{3,b}$ in~\eqref{L-Zbn-3b} we appeal to \eqref{eq:tau-Jg-D2}, in addition to \eqref{bootstraps}, \eqref{eq:Q:all:bbq}, and \eqref{geometry-bounds-new}, while for the  $L^{\ZZ_n}_{3,c}$ term in~\eqref{L-Zbn-3c} we note the exact derivative structure involving $\nbs_2 ( (\nn \cdot \nbs^6 \tt)^2 )$, which is dealt with as usual by integrating-by-parts the $\nbs_2$ term. Lastly, for $L^{\ZZ_n}_{3,d}$ term in~\eqref{L-Zbn-3d} we additionally  appeal to \eqref{eq:I:W:nn:6:all}. Putting this together leads to 
\begin{equation}
\label{eq:duck:a:duck:4}
\sabs{L^{\ZZ_n}_{3}}
\les \mathsf{K} \eps   (\tfrac{4}{\kappa_0})^{2\beta} \brak{\mathsf{B}_6}^2
\,.
\end{equation}
Combining the bounds in this subsection, namely~\eqref{eq:duck:a:duck:1}, \eqref{eq:duck:a:duck:2}, \eqref{eq:duck:a:duck:3}, and~\eqref{eq:duck:a:duck:4}, 
we obtain
\begin{equation}
\label{eq:duck:a:duck:all}
\sabs{I^{\ZZ_n}_{7}}
\les \mathsf{K} \eps   (\tfrac{4}{\kappa_0})^{2\beta} \brak{\mathsf{B}_6}^2
\,.
\end{equation}

\subsubsection{The integral $I^{\AA_n}_5$}
\label{sec:IA5}
Using \eqref{adjoint-2} we may rewrite the integral defined in~\eqref{I5-Abn} as 
\begin{align*} 
 I^{\AA_n}_{5} &= I^{\AA_n}_{5,a} + I^{\AA_n}_{5,b} + I^{\AA_n}_{5,c}  + I^{\AA_n}_{5,d} \,, 
 \notag \\
 I^{\AA_n}_{5,a} &=
-2 \alpha \tint \jb \Fgss \Jg  \Sbt  \nbs^6 \tt\cdo\nn \  g^{- {\frac{1}{2}} }  \nbs_2 \nbs^6(\Jg\Abn) \,,
 \notag \\
 I^{\AA_n}_{5,b} &=
- 2 \alpha \tint  \nbs_2( \jb \Fgss\Jg g^{- {\frac{1}{2}} }  \Sbt \nn_i)  \nbs^6 \tt_i \ \nbs^6(\Jg\Abn) \,,
 \notag \\
 I^{\AA_n}_{5,c} &=
2 \alpha \tint   \jb \Qr_2 \Fgss\Jg g^{- {\frac{1}{2}} }  \Sbt   \nbs^6 \tt\cdo\nn \ \nbs^6(\Jg\Abn) \,,
 \notag \\
 I^{\AA_n}_{5,d} &=
-2 \alpha \int   \jb \Qb_2 \Fgss\Jg g^{- {\frac{1}{2}} }  \Sbt   \nbs^6 \tt\cdo\nn \ \nbs^6(\Jg\Abn)\Big|_\s \,.
\end{align*} 
The integrals $ I^{\AA_n}_{5,b}$, $ I^{\AA_n}_{5,c}$, and $ I^{\AA_n}_{5,d}$ can be estimated directly by appealing to \eqref{bootstraps}, \eqref{eq:Q:all:bbq}, and \eqref{geometry-bounds-new} as
\begin{equation*}
\sabs{I^{\AA_n}_{5,b}} + \sabs{I^{\AA_n}_{5,c}} + \sabs{I^{\AA_n}_{5,d}}
\les
\mathsf{K} \eps^2   (\tfrac{4}{\kappa_0})^{2\beta} \brak{\mathsf{B}_6}^2
\,.
\end{equation*}
On the other hand, the 
integral $ I^{\AA_n}_{5,a}$ contains an over-differentiated term. We notice that $I^{\AA_n}_{5,a}$ is nearly identical to the integral $L^{\ZZ_n}_{3} $ in \eqref{L-Zbn-3}: $\Sbt$ is replacing $\Jg (\Abn + \Zbt)$. These terms however satisfy the same bounds, even at the first derivative level, in view of \eqref{bootstraps}. As such, we may bound $I^{\AA_n}_{5,b}$ in the identical fashion as the integral $L^{\ZZ_n}_{3} $ (cf.~\eqref{eq:duck:a:duck:4}) and obtain
\begin{equation*}
\sabs{I^{\AA_n}_{5,a}} 
\les 
\mathsf{K} \eps  (\tfrac{4}{\kappa_0})^{2\beta} \brak{\mathsf{B}_6}^2
\,.
\end{equation*}
Putting together the previous two estimates we thus obtain
\begin{equation} 
\label{eq:insane:in:the:membrane}
\sabs{ I^{\AA_n}_{5}} 
\les 
\mathsf{K} \eps  (\tfrac{4}{\kappa_0})^{2\beta} \brak{\mathsf{B}_6}^2
\,.
\end{equation}

\subsubsection{The integral $I^{\ZZ_n}_8$}  
\label{sec:IZ8}
We next study the integral $I^{\ZZ_n}_8$, defined in \eqref{I8-Zbn}. By comparing the integral $I^{\ZZ_n}_8$ with $I^{\ZZ_n}_7$ (as defined in~\eqref{I7-Zbn}), we notice two differences. First, the lower order term $- (\Jg \Abn  + \Jg \Zbt)$ in $I^{\ZZ_n}_7$, is replaced by the term $\Zbn$ in $I^{\ZZ_n}_8$. These terms, and their $\nbs$ derivatives, satisfy the same upper bounds in $L^\infty_{x,\s}$, in light of \eqref{bootstraps}. Second, the differential operator $\nbs_1 - \eps \Jg g^{-\frac 12} \nbs_2 h \nbs_2$ acting on $\nbs^6 \tt \cdot \nn$ in $I^{\ZZ_n}_7$, is now acting on $\nbs^6 \Jg$ in $I^{\ZZ_n}_8$. These terms satisfy nearly the same bounds in view of~\eqref{D6JgEnergy:new}, \eqref{D6n-bound:b:new}, and~\eqref{D6n-bound:b:new:bdd}, with only one difference: the bounds for the sixth order derivatives of $\Jg$ are worse by a power of $(\mathsf{K} \eps)^{-1}$ than the sixth order derivatives of $\tt$. As such, it is clear that by repeating the decompositions and the bounds in Subsection~\ref{sec:I:Zn:7:integral}, we arrive at a bound for $I^{\ZZ_n}_8$ which is worse than that for $I^{\ZZ_n}_7$ (see~\eqref{eq:duck:a:duck:all}) by a power of $(\mathsf{K} \eps)^{-1}$, namely 
\begin{equation} 
\label{eq:duck:a:duck:all:again}
\sabs{I^{\ZZ_n}_8} 
\le (\tfrac{4}{\kappa_0})^{2\beta} \brak{\mathsf{B}_6}^2
\,.
\end{equation} 
To avoid redundancy, we do not repeat the argument which proves~\eqref{eq:duck:a:duck:all:again}.

\subsection{The forcing and commutator terms}\label{sec:11:forcing:comm}
We now turn to the bounds for the remaining integrals $ I^{\WW_n}_6$, $ I^{\ZZ_n}_{10}$, $ I^{\AA_n}_{10}$ in 
\eqref{I6-Wbn}, \eqref{I10-Zbn}, and \eqref{I10-Abn}, respectively.
In order to estimate these integrals, we use  the definitions  for the forcing functions in \eqref{forcing-nt} together with
the definitions of the so-called remainder and commutator functions in \eqref{Cw-Rw-comm}, \eqref{Cz-Rz-comm}, and  \eqref{Ca-Ra-comm}.  
Bounds for these quantities were obtained earlier in~\eqref{eq:I:W:nn:6:all},~\eqref{eq:I:Z:nn:10:all}, and~\eqref{eq:I:A:nn:10:all}. Using Cauchy-Schwartz, the bound $\mathcal{J} \leq \Jg$, and~\eqref{bootstraps-Dnorm:6}, we deduce from the aforementioned bounds that
\begin{subequations}
\begin{align}
\sabs{I^{\WW_n}_6}
&\leq \Cn (\tfrac{4}{\kappa_0})^{2\beta}  \brak{\mathsf{B}_6}^2 
\label{eq:I:W:nn:6:final}
\,.
\end{align}
Similarly, using \eqref{eq:I:Z:nn:10:all} we also have that
\begin{align}
\sabs{I^{\ZZ_n}_{10}}
&\le  \int_0^\s  \snorm{\tfrac{ \mathcal{J}^{\frac 34} }{\Sigma^{\beta-1}} \nbs^6(\Jg \Zbn)(\cdot,\s')}_{L^2_x}  
 \Bigl( \snorm{\tfrac{\mathcal{J}^{\frac 34} }{\Sigma^{\beta-1}} \nbs^6\Fzn (\cdot,\s')}_{L^2_x}  
+ \snorm{\tfrac{\mathcal{J}^{\frac 34} }{\Sigma^{\beta-1}} \mathcal{R}_\Zb^\nn (\cdot,\s')}_{L^2_x}  
+ \snorm{\tfrac{\mathcal{J}^{\frac 34} }{\Sigma^{\beta-1}} \mathcal{C}_\Zb^\nn(\cdot,\s')}_{L^2_x}  \Bigr){\rm d} \s'
\notag\\
&\le \tfrac{4 (1+ \alpha )}{\eps} \int_0^\s \snorm{ \tfrac{\mathcal{J}^{\frac 34}}{\Sigma^{\beta}}   \nbs^6 (\Jg \Zbn)(\cdot,\s')}_{L^2_{x}}^2 {\rm d}\s'
+ \Cn   (\tfrac{4}{\kappa_0})^{2\beta}    \mathsf{K} \brak{\mathsf{B_6}}^2
   \,.  
\label{I10-Zn}
\end{align} 
In the same way, using \eqref{eq:I:A:nn:10:all} we  obtain the bound
\begin{equation} 
\sabs{ I^{\AA_n}_{10} } \le 
\tfrac{4 (1+ \alpha )}{\eps} \int_0^\s \snorm{ \tfrac{\mathcal{J}^{\frac 34}}{\Sigma^{\beta}}   \nbs^6 (\Jg \Abn)(\cdot,\s')}_{L^2_{x}}^2 {\rm d}\s'
+ \Cn   (\tfrac{4}{\kappa_0})^{2\beta}    \mathsf{K} \brak{\mathsf{B_6}}^2
   \,.  \label{I10-An}
\end{equation} 
\end{subequations}
Summarizing the bounds \eqref{eq:I:W:nn:6:final}, \eqref{I10-Zn}, and \eqref{I10-An}  we thus deduce that
\begin{equation} 
\sabs{ I^{\WW_n}_6}+ \sabs{ I^{\ZZ_n}_{10}} + \sabs{ I^{\AA_n}_{10}}
  \le
\tfrac{4 (1+ \alpha )}{\eps} \int_0^\s \snorm{ \tfrac{\mathcal{J}^{\frac 34}}{\Sigma^{\beta}}   \nbs^6 (\Jg\Zbn,\Jg \Abn)(\cdot,\s')}_{L^2_{x}}^2 {\rm d}\s'
+ \Cn   (\tfrac{4}{\kappa_0})^{2\beta}  \mathsf{K} \brak{\mathsf{B_6}}^2\,.
\label{eq:I:forcing:normal}
\end{equation}

\subsection{Conclusion of the six derivative normal energy bounds}
\label{sec:D6:n:final}
We return to the energy identity~\eqref{D6-L2-N}, with the decompositions~\eqref{Integral-Wbn}, \eqref{Integral-Zbn}, and~\eqref{Integral-Abn}. We collect the lower bounds~\eqref{eq:I:n:12369}, \eqref{eq:fuck:yeah:end}, the estimate $\mathsf{G_{good}} \geq - \frac{1+\alpha}{3\eps} \mathcal{J}^{\frac 32}$,   the upper bounds~\eqref{IWn5+IZn5+IAn7}, \eqref{eq:small:duck:6}, \eqref{eq:duck:a:duck:all}, \eqref{eq:insane:in:the:membrane}, \eqref{eq:duck:a:duck:all:again}, \eqref{eq:I:forcing:normal}, and the initial data assumption~\eqref{table:derivatives}, to obtain
\begin{align}
0
&\geq 
\bigl( \tfrac{1}{52} - \Cn \eps\bigr) 
\snorm{\tfrac{\mathcal{J}^{\frac 34} (\Jg \Q)^{\frac 12}}{\Sigma^\beta} \nbs^6(\Jg \Wbn,\Jg \Zbn,\Jg \Abn)(\cdot,\s)}_{L^2_x}^2
- \tfrac{9(1+\alpha)}{\eps}   (\tfrac{3}{\kappa_0})^{2\beta}  \Cdatatwo
- \Cn (\tfrac{4}{\kappa_0})^{2\beta}  \mathsf{K}^2 \Bsix^2
\notag\\
&\qquad
+\Bigl(\tfrac{1+\alpha}{24} - \Cn \eps \beta\Bigr) \tfrac{1}{\eps} 
\int_0^{\s} \snorm{\tfrac{\mathcal{J}^{\frac 14} \Jgh}{\Sigma^\beta}   \nbs^6(\Jg\Wbn,\Jg\Zbn,\Jg\Abn)(\cdot,\s')}_{L^2_x}^2 {\rm d}\s'
\notag\\
&\qquad 
+ \Bigl(  \tfrac{\alpha( \beta - \frac 12)}{8 }   +  \bubu{ \tfrac{9\alpha }{10 }}  -  (16 + 25^2) (1+ \alpha )  \Bigr) \tfrac{1}{\eps}
\int_0^{\s}  \snorm{\tfrac{ \mathcal{J}^{\frac 34} }{\Sigma^\beta} \nbs^6(\Jg\Zbn,\Jg\Abn) (\cdot,\s')}_{L^2_x}^2
{\rm d} \s'
\notag\\
&\qquad 
- \Bigl( \tfrac{16 \alpha  (\beta -\frac 12) }{(1+\alpha)}  +  \bubu{ 33 } + 2 \cdot 250^2 +   39^2 + 500^2\Bigr) \tfrac{1}{\eps} \int_0^{\s} 
\snorm{\tfrac{\mathcal{J}^{\frac 34}(\Jg \Q)^{\frac 12}}{\Sigma^\beta} \nbs^6(\Jg\Wbn,\Jg \Zbn,\Jg\Abn)(\cdot,\s')}_{L^2_x}^2
{\rm d} \s'
\notag\\
&\qquad 
+ \tfrac{1}{20(1+\alpha)} \tfrac{1}{\eps^2} \snorm{\tfrac{\Q \mathcal{J}^{\frac 14}}{\Sigma^\beta} \nbs^6 \Jg(\cdot,\s)}_{L^2_x}^2 
-  \tfrac{20^2 + 500^2 + 100 \cdot 250^2}{1+\alpha}  {\tfrac{1}{\eps^3}} \int_0^{\s}   \snorm{\tfrac{\Q \mathcal{J}^{\frac 14}}{\Sigma^\beta} \nbs^6 \Jg(\cdot,\s')}_{L^2_x}^2 {\rm d} \s'
\notag\\
&\qquad
+ \tfrac{7}{40(1+\alpha)} \tfrac{1}{\eps^3} \int_0^{\s} \snorm{\tfrac{\Q \mathcal{J}^{-\frac 14}}{\Sigma^\beta} \nbs^6 \Jg(\cdot,\s')}_{L^2_x}^2 {\rm d} \s'
+ \bigl(8 - \Cn \brak{\beta} \eps \bigr)  {\tfrac{1}{\eps^3}}  \tint \tfrac{\Q}{\Sigma^{2\beta}} \Fgh    \sabs{ \nbs^6\Jg}^2
\,,
\label{eq:normal:conclusion:1}
\end{align}
where $\Cn = \Cn(\alpha,\kappa_0,\Cdata)$ is independent of $\beta$ (and $\eps$, as always).

At this stage, we choose $\beta = \beta(\alpha)$ to be sufficiently large to ensure that the damping for $\nbs^6(\Jg\Zbn,\Jg\Abn)$ is strong enough, i.e., so that 
\begin{equation*}
\tfrac{\alpha( \beta - \frac 12)}{8 }   +  \bubu{ \tfrac{9\alpha }{10 }  } -  (16+ 25^2 ) (1+ \alpha ) \geq 0 \,.
\end{equation*}
More precisely, we choose $\beta$ to ensure equality in the above inequality; namely, we have that 
\begin{equation}
\bubu{ \beta_\alpha  := \tfrac{8(1+\alpha)}{\alpha} \bigl( 16+ 25^2     \bigr)  - \tfrac{67}{10} }\,.
\label{eq:normal:bounds:beta}
\end{equation}
With this choice of $\beta = \beta_\alpha$, we return to \eqref{eq:normal:conclusion:1}, and choose $\eps$ to be sufficiently small in terms of $\alpha,\kappa_0,\Cdata$. After re-arranging we deduce that
\begin{align}
& \tfrac{1}{53}  \snorm{\tfrac{\mathcal{J}^{\frac 34} (\Jg \Q)^{\frac 12}}{\Sigma^{\beta_\alpha}} \nbs^6(\Jg \Wbn,\Jg \Zbn,\Jg \Abn)(\cdot,\s)}_{L^2_x}^2
 + \tfrac{1}{20(1+\alpha)} \tfrac{1}{\eps^2} \snorm{\tfrac{\Q \mathcal{J}^{\frac 14}}{\Sigma^{\beta_\alpha}} \nbs^6 \Jg(\cdot,\s)}_{L^2_x}^2 
 \notag\\
&\qquad
+ \tfrac{1+\alpha}{48 \eps} 
\int_0^{\s} \snorm{\tfrac{\mathcal{J}^{\frac 14} \Jg^{\frac12}}{\Sigma^{\beta_\alpha}}   \nbs^6(\Jg\Wbn,\Jg\Zbn,\Jg\Abn)(\cdot,\s')}_{L^2_x}^2 {\rm d}\s'
+ \tfrac{7}{40(1+\alpha)} \tfrac{1}{\eps^3} \int_0^{\s} \snorm{\tfrac{\Q \mathcal{J}^{-\frac 14}}{\Sigma^{\beta_\alpha}} \nbs^6 \Jg(\cdot,\s')}_{L^2_x}^2 {\rm d} \s'
\notag\\
&\leq
\tfrac{9(1+\alpha)}{\eps}  (\tfrac{3}{\kappa_0})^{2\beta_\alpha}  \Cdatatwo
+
\Cn (\tfrac{4}{\kappa_0})^{2\beta_\alpha}  \mathsf{K}^2 \Bsix^2
\notag\\
 &\qquad
+ \tfrac{C}{\eps} \int_0^{\s} 
\snorm{\tfrac{\mathcal{J}^{\frac 34}(\Jg \Q)^{\frac 12}}{\Sigma^{\beta_\alpha}} \nbs^6(\Jg\Wbn,\Jg \Zbn,\Jg\Abn)(\cdot,\s')}_{L^2_x}^2
{\rm d} \s'
+  \tfrac{C}{\eps^3} \int_0^{\s}   \snorm{\tfrac{\Q \mathcal{J}^{\frac 14}}{\Sigma^{\beta_\alpha}} \nbs^6 \Jg(\cdot,\s')}_{L^2_x}^2 {\rm d} \s'
\,,
\label{eq:normal:conclusion:2}
\end{align}
where $C$ is a universal constant (in particular, independent of $\alpha,\kappa_0,\Cdata$), and $\Cn$ is as usual.

By inspecting the first line on the left side and the last line on the right side of \eqref{eq:normal:conclusion:2}, we observe that we may apply Gr\"onwall's inequality for $\s \in [0,\eps]$. More precisely, there exists a constant 
\begin{equation}
\check{\mathsf{c}}_\alpha >0
\end{equation}
which only depends on $\alpha$, and may be computed explicitly from \eqref{eq:Q:bbq}, \eqref{eq:normal:bounds:beta}, and \eqref{eq:normal:conclusion:2}, such that 
\begin{align}
& 
\sup_{\s \in [0,\eps]}
\snorm{\tfrac{\mathcal{J}^{\frac 34}  \Jg^{\frac 12}}{\Sigma^{\beta_\alpha}} \nbs^6(\Jg\Wbn,\Jg\Zbn,\Jg\Abn)(\cdot,\s)}_{L^2_x}^2 
 + \tfrac{1}{\eps}   \int_0^{\eps}  
\snorm{\tfrac{\mathcal{J}^{\frac 14} \Jg^{\frac 12}}{\Sigma^{\beta_\alpha}} \nbs^6 (\Jg\Wbn,\Jg\Zbn,\Jg\Abn)(\cdot,\s)}_{L^2_x}^2   {\rm d} \s
 \notag\\
 &\qquad 
+ \tfrac{1}{\eps^2} \sup_{\s \in [0,\eps]} 
 \snorm{ \tfrac{\mathcal{J}^{\frac 14}}{\Sigma^{\beta_\alpha}}  \nbs^6 \Jg (\cdot,\s)}_{L^2_{x}}^2
+ \tfrac{1}{\eps^3} \int_0^{\eps} \snorm{ \tfrac{\mathcal{J}^{-\frac 14}}{\Sigma^{\beta_\alpha}} \nbs^6 \Jg (\cdot,\s)}_{L^2_{x}}^2 {\rm d} \s
\notag\\
 &  
 \leq \check{\mathsf{c}}_\alpha \tfrac{1}{\eps} (\tfrac{4}{\kappa_0})^{2\beta_\alpha} 
 \Big( \Cdatatwo  + \Cn \eps \mathsf{K}^2  \brak{\mathsf{B}_6}^2\Bigr)
 \,.
 \label{eq:normal:conclusion:3}
\end{align}
 
At last, we multiply the above estimate by $\kappa_0^{2 \beta_\alpha}$, appeal to \eqref{bs-Sigma}, drop the energy and damping terms for $ \nbs^6 \Jg$ (since these were bounded already in Proposition~\ref{prop:geometry}),  and recall the definitions of $\widetilde{\mathcal{E}}_{6,\nnn}^2(\s)$ and $\widetilde{\mathcal{D}}_{6,\nnn}^2(\s) $ to deduce that 
\begin{align}
\eps \sup_{\s \in [0,\eps]} \widetilde{\mathcal{E}}_{6,\nnn}^2(\s)
+\widetilde{\mathcal{D}}_{6,\nnn}^2(\eps) 
&  \leq \check{\mathsf{c}}_\alpha  4^{2\beta_\alpha} 
 \Big( \Cdatatwo  + \Cn \eps \mathsf{K}^2  \brak{\mathsf{B}_6}^2\Bigr)
 \notag\\
 &  \leq  \mathsf{B}_6^2 \check{\mathsf{c}}_\alpha  4^{2\beta_\alpha} 
 \Big(\tfrac{\Cdatatwo}{ \mathsf{B}_6^2}  + \Cn \eps \mathsf{K}^2  \tfrac{\brak{\mathsf{B}_6}^2} {\mathsf{B}_6^2}\Bigr)
 \,.
 \label{eq:normal:conclusion:4}
\end{align}
Since $ \mathsf{B}_6 \geq 1$ (cf.~\eqref{eq:B6:choice:1}) and since $\mathsf{K}$ was already defined solely in terms of $\alpha$ (cf.~\eqref{eq:K:choice:1}),  upon ensuring that 
\begin{align}
\mathsf{B}_6 &\geq 
4 \check{\mathsf{c}}_\alpha^{\frac 12} 
4^{\beta_\alpha}
\Cdata\,,
\label{eq:B6:choice:2}
\end{align}
where $\beta_\alpha$ is as defined in~\eqref{eq:normal:bounds:beta}, 
and taking $\eps$ sufficiently small in terms of $\alpha,\kappa_0,\Cdata$, we deduce from \eqref{eq:normal:conclusion:4} that 
\begin{equation}
\eps \sup_{\s \in [0,\eps]} \widetilde{\mathcal{E}}_{6,\nnn}^2(\s)
+\widetilde{\mathcal{D}}_{6,\nnn}^2(\eps) 
\leq \tfrac{1}{8}   \mathsf{B}_6^2 
 \,,
\label{eq:normal:conclusion:5}
\end{equation}
which closes the ``normal part'' of the remaining bootstrap \eqref{bootstraps-Dnorm:6}.

Indeed, we note in closing that adding \eqref{eq:hate:13:a} and \eqref{eq:normal:conclusion:5}
gives
\begin{equation}
\eps \sup_{\s \in [0,\eps]} \widetilde{\mathcal{E}}_{6}(\s)
+\widetilde{\mathcal{D}}_{6}(\eps) 
\leq \tfrac{1}{2}   \mathsf{B}_6 
\,, \label{eq:normfinal1}
\end{equation}
which closes the bootstrap \eqref{bootstraps-Dnorm:6}.

\section{Downstream maximal globally hyperbolic development in a box}
\label{sec:downstreammaxdev}

In this section we give the proof of Theorem~\ref{thm:main:DS}.
We continue to use the notation in Section~\ref{sec:remapping}.
Consider the level-set $\{ \Jgb,_1 =0\}$. This is a hypersurface parameterized by the graph $(x_2,t)\mapsto (x_1^*(x_2,t),x_2,t)$; see the
definition of $x_1^*(x_2,t)$ given in \eqref{eq:x1star:def}. 
In this section we consider the Euler evolution in the spacetime geometry which is on the {\em downstream side} of the hypersurface $\{\Jgb,_1 =0\}$ (namely, for $x_1 >x_1^*(x_2,t)$), and which is bounded from above by (a small modification of) the parabolic hypersurface $\{\Jgb =0\}$. On the {\em upstream side} (namely, for
$x_1 < x_1^*(x_2,t)$), the spacetime we consider is  bounded from above by the cylindrical hypersurface $\{\mathcal{J} = 0\} = \{ \Jgb(x_1^*(x_2,t),x_2,t)=0\}$, which is the same as the top boundary considered in Sections~\ref{sec:formation:setup}--\ref{sec:sixth:order:energy} (see Figures~\ref{fig:downstream} and~\ref{fig:parabolic-cylinder-flat} below).

\begin{figure}[htb!]
\centering
  \includegraphics[width=0.4\linewidth]{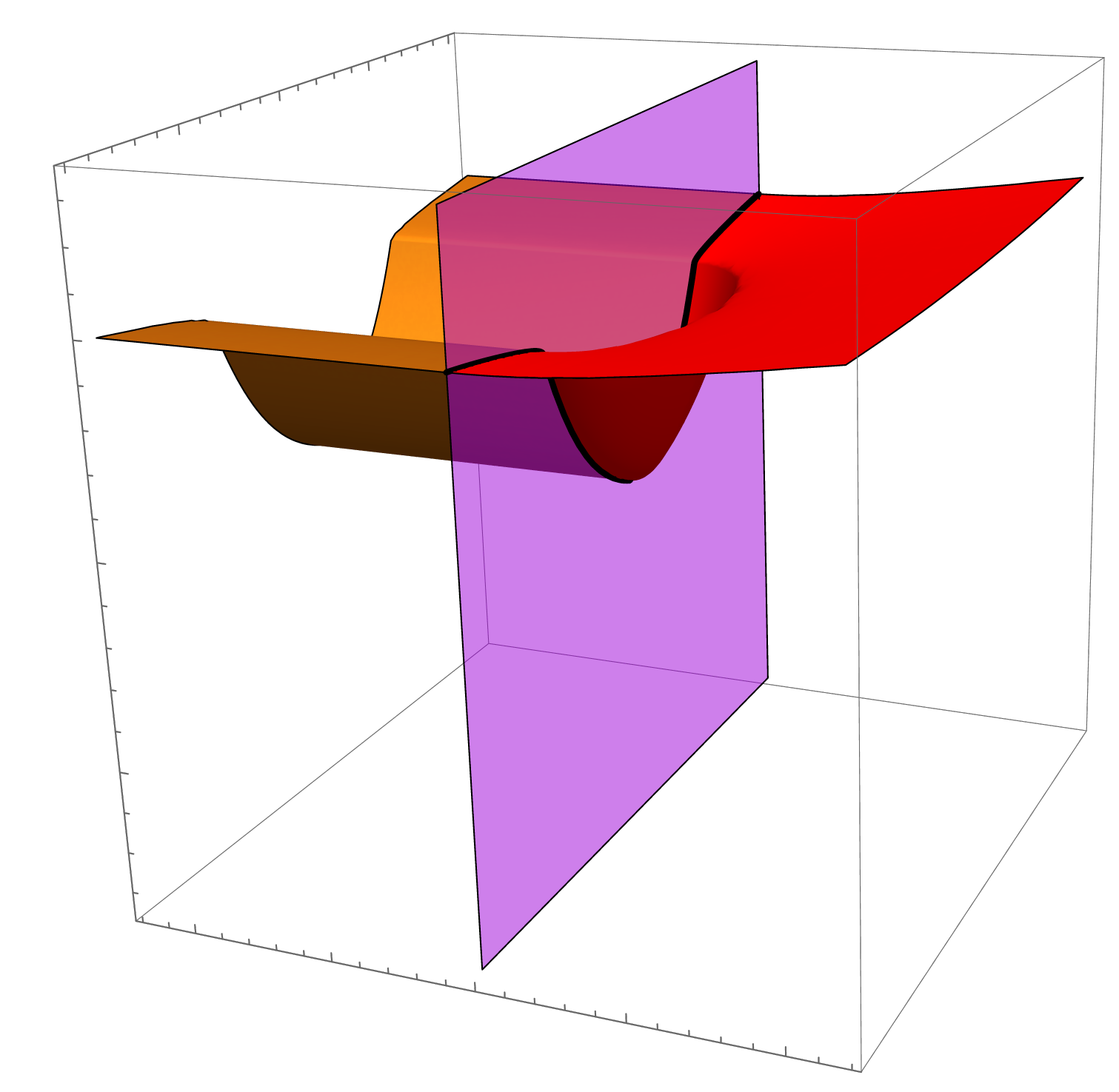}
\vspace{-.1 in}
\caption{Consider the function $\Jgb = \Jgb(x_1,x_2,t)$, as defined by~\eqref{eq:Jgb:identity:0}, with $\Jg$ as in Figure~\ref{fig:parabolic:spacetime}. The bounding box represents the zoomed-in region $|x_1| \leq \tfrac{\pi\eps}{4}$, $|x_2| \leq \tfrac{1}{16}$, and $t\in [-\final,\final]$. In magenta, we plot the surface $\{x_1 = x_1^*(x_2,t)\}$, which separates the downstream side (to the right) from the upstream side (to the left). In orange, we plot the surface $\{(x,t) \in \TT^2 \times [\initial,\final] \colon x_1< x_1^*(x_2,t), \mathcal{J}(x_2,t) = \Jgb(x_1^*(x_2,t),x_2,t)=0 \}$, which represents the future temporal boundary (``top'' boundary) of the spacetime considered in this section on the upstream side. In red, we plot the surface $\{(x,t) \in \TT^2 \times [\initial,\final] \colon x_1 > x_1^*(x_2,t),  \Jgb(x_1,x_2,t)=0 \}$, which represents the future temporal boundary of the spacetime considered in this section on the downstream side.}  
\label{fig:downstream}
\end{figure}

\subsection{Flattening the top of the spacetime}
\label{sec:design:JJ:weight:DS}
For $(x,t) \in \TT^2 \times (\initial,\final]$, we define 
\begin{align*} 
\mJg(x_1,x_2,t)
:= 
\begin{cases}
\Jgb(x_1,x_2,t), &x_1^*(x_2,t) \le x_1 \le x_1^\sharp(x_2)\,, \\
\Jgb(x_1^\sharp(x_2),x_2,t)  , &x_1 \ge x_1^\sharp(x_2)\,,
\end{cases}
\end{align*} 
where $x_1^\sharp(x_2)$ is defined in terms of the function $w_0(x)$ as
\begin{equation} 
x_1^{\sharp}(x_2) = \bigl\{ x_1 \colon x_1 > x_1^\vee(x_2), \nb_1w_0(x_1,x_2) = - \tfrac{17}{20} \bigr\} 
\label{eq:x1:sharp:def}
\,.
\end{equation} 
We also define $\mJg(x,\initial):= 1$ for all $x\in \TT^2$. 

\begin{proposition}[\bf $x_1^\sharp$ is well-defined and $\mJg \leq \Jgb$]
\label{prop:mJg:well:def}
The function $\TT \ni x_2 \mapsto x_1^\sharp(x_2)$ given by \eqref{eq:x1:sharp:def} is well-defined and differentiable, and satisfies 
\begin{equation}
\label{eq:x1:sharp:props}
x_1^{\sharp}(x_2) - x_1^*(x_2,t) \geq \tfrac{\eps}{41}\,,
\qquad
|\p_2 x_1^\sharp(x_2)| \leq 240 \eps \,,
\end{equation}
pointwise in $(x_2,t)$. Moreover we have that 
\begin{equation}
 \mJg \leq \Jgb
 \label{eq:mJg:props}
 \,.
\end{equation}
\end{proposition}
\begin{proof}[Proof of Proposition~\ref{prop:mJg:well:def}]
In order to check that $\mJg$ is well-defined, we need to verify $x_1^\sharp(x_2) > x_1^*(x_2,t)$, pointwise in $(x_2,t)$. We recall the assumptions on $w_0(x)$ given in Section~\ref{cauchydata}. To show that $x_1^{\sharp}(x_2) - x_1^*(x_2,t) \geq \frac{\eps}{41} $, uniformly in $t$, we first recall from \eqref{eq:x1star:x1vee} that $|x_1^*(x_2,t)-x_1^\vee(x_2)|\leq \Cn \mathsf{K} \brak{\mathsf{B}_6}\eps^3$, so that it is sufficient to show $x_1^\sharp > x_1^\vee + \frac{\eps}{40}$.  Using \eqref{eq:why:the:fuck:not:0}, \eqref{item:ic:w0:x2:negative}, and \eqref{eq:x1:sharp:def}, we have that $\frac{1}{20} \leq |\nb_1 w_0(x_1^\sharp(x_2),x_2) - \nb_1 w_0(x_1^\vee(x_2),x_2) | = |x_1^\sharp(x_2) - x_1^\vee(x_2)| \|\p_1 \nb_1  w_0\|_{L^\infty_x} \leq \frac{2}{\eps} |x_1^\sharp(x_2) - x_1^\vee(x_2)| $. Since by assumption $x_1^\sharp(x_2) > x_1^\vee(x_2)$ (see~\eqref{eq:x1:sharp:def}), we deduce $x_1^\sharp > x_1^\vee + \frac{\eps}{40}$.

Next, we show that $x_1^\sharp(x_2)$ is well-defined. Note that from~\eqref{item:ic:supp},~\eqref{item:ic:w0:x2:negative}, the continuity of $\nb_1 w_0$, and the intermediate value theorem, for every $x_2 \in \TT$ there exists at least one $x_1^\vee(x_2) < x_1^{\sharp}(x_2) < 13 \pi \eps$ satisfying~$\nb_1 w_0(x_1,x_2) = - \frac{17}{20}$. To establish uniqueness, we first define $x_1^\sharp(x_2)$ as the smallest value $x_1>x_1^\vee(x_2)$ such that $\nb_1w_0(x) = -\frac{17}{40}$, and then show that for $x_1 \in (x_1^\sharp(x_2),13 \pi \eps]$, we must have $\nb_1w_0(x) > -\frac{17}{40}$. Since $x_1^{\sharp}(x_2) - x_1^\vee(x_2) \geq \frac{\eps}{40} \geq \eps^{\frac 74}$, from assumption \eqref{item:ic:w0:d11:positive} we know that for all $x_1 \in (x_1^\sharp(x_2),13 \pi \eps]$ with $\nb_1 w_0(x) < - \frac 13$, we must have $\nb_1^2 w_0(x)\geq \eps^{\frac 78}$, showing that as $x_1$ increases, $\nb_1 w_0(x)$ strictly increases from the value $-\frac{17}{20}$ when $x_1 = x_1^\sharp(x_2)$, until it reaches the value $-\frac 13$ at some point $x_1= x_1^\top(x_2)$. Additionally, \eqref{item:ic:w0:d11:positive} implies that for $x_1 > x_1^\top(x_2)$ we have that $\nb_1 w_0(x) \geq -\frac 13 $: this is because if $\nb_1 w_0(x)$ wanted to dip below the value $-\frac 13$, then it would need to decrease as a function of $x_1$, but \eqref{item:ic:w0:d11:positive} implies that $\nb_1 w_0(x)$ can only decrease in $x_1$ if $\nb_1 w_0(x) \geq - \frac 13$. Therefore, $\nb_1 w_0(x_1,x_2) > - \frac{17}{20}$ for all $x_1 > x_1^\sharp(x_2)$, giving uniqueness.

Next, we establish the second bound in~\eqref{eq:x1:sharp:props}. 
By implicitly differentiating the relation $\nb_1 w_0(x_1^\sharp(x_2),x_2) = - \frac{17}{20}$, we obtain that $\p_2 x_1^\sharp(x_2) = - \eps \frac{\nb_1 \nb_2 w_0}{\nb_1\nb_1 w_0}(x_1^\sharp(x_2),x_2)$. The numerator of the above fraction is bounded from above by $2$, in absolute value (due to~\eqref{eq:why:the:fuck:not:0}). The denominator is estimated by noting that via the mean value theorem, $\nb_1^2 w_0(x_1^\sharp(x_2),x_2) = \nb_1^2 w_0(x_1^\vee(x_2),x_2) + \frac{x_1^\sharp(x_2)- x_1^\vee(x_2)}{\eps} \nb_1^3 w_0(x_1^\prime,x_2) = \frac{x_1^\sharp(x_2) - x_1^\vee(x_2)}{\eps} \nb_1^3 w_0(x_1^\prime,x_2) \geq \frac{1}{40} \cdot \frac{1}{3}$, for some $x_1^\prime \in (x_1^\vee(x_2),x_1^\sharp(x_2))$. Here we have used the previously established estimate $x_1^\sharp(x_2) - x_1^\vee \geq \frac{\eps}{40}$, the fact that by definition $\nb_1^2(x_1^\vee(x_2) ,x_2)= 0$, and the fact that at all points $x_1^\prime$ in between $x_1^\vee$ and $x_1^\sharp$ we have that $\nb_1 w_0(x_1^\prime,x_2) \leq - \frac 13$ (in fact $\leq - \frac{17}{20}$), and thus by assumption~\eqref{item:ic:w0:d11:positive} we know that $\nb_1^3 w_0(x_1^\prime,x_2) \geq \frac 13$. This shows that $\nb_1^2 w_0(x_1^\sharp(x_2),x_2) \geq \frac{1}{120}$, and therefore $|\p_2 x_1^\sharp(x_2)| \leq 240 \eps$, as claimed.

Lastly, we establish~\eqref{eq:mJg:props}. That is, we need to show that $0 \leq \Jgb(x_1,x_2,t) - \Jgb(x_1^\sharp(x_2),x_2,t) = \Jg(x_1,x_2,t) - \Jg(x_1^\sharp(x_2),x_2,t)$ whenever $x_1 \ge x_1^\sharp(x_2)$; here, in the second equality we have used \eqref{eq:Jgb:identity:0}. In turn, this fact  follows from assumption \eqref{item:ic:w0:d11:positive} and Corollary~\ref{cor:Jg:initial}, as follows. From the mean value theorem in $x_1$, the bound \eqref{bs-Jg-1} with $i=1$, and the considerations in the second paragraph of this proof, we have that $\Jg(x_1,x_2,t) - \Jg(x_1^\sharp(x_2),x_2,t)$ cannot vanish for $x_1 \in (x_1^\sharp(x_2),x_1^\top(x_2)]$. But if $x_1 > x_1^\top(x_2)$ then necessarily $(w_0),_1(x_1,x_2) \geq -\frac{1}{3\eps}$ and since $(w_0),_1(x_1^\sharp(x_2),x_2) = -\frac{17}{20 \eps}$, we have that $\Jg(x_1,x_2,t) - \Jg(x_1^\sharp(x_2),x_2,t)$ cannot vanish due to \eqref{bs-Jg-0}.
\end{proof}

The only issue is that while the map $x_1 \mapsto \mJg(x_1,x_2,t)$ is Lipschitz continuous, it is not $C^1$ smooth. As such, we need to consider a variant of $\mJg$, denoted by $\mmJg$, which is $H^6$ smooth in space and time in the set $\{(x,t) \colon t \in (\initial,\final), x_2 \in \TT, x_1> x_1^*(x_2,t)\}$. It is convenient to introduce 
\begin{equation}
x_{1,+}^\sharp(x_2) = x_1^\sharp(x_2) + \tfrac{\eps}{1000}\,,
\qquad\mbox{and}\qquad 
x_{1,-}^\sharp(x_2) = x_1^\sharp(x_2) - \tfrac{\eps}{1000}\,,
\end{equation}
where $x_1^\sharp(x_2)$ is as in \eqref{eq:x1:sharp:def}.
Then, we let $\mmJg(x,\initial) =1$
and for $(x,t) \in \TT^2\times (\initial,\final]$ with $x_1 \geq x_1^*(x_2,t)$ we define
\begin{subequations}
\begin{align} 
\label{eq:mmJg:def}
\mmJg(x_1,x_2,t)
:= 
\begin{cases}
\mJg(x_1,x_2,t) = \Jgb(x_1,x_2,t) &x_1^*(x_2,t) \le x_1 \le x_{1,-}^\sharp(x_2)\,, \\
\mbox{a smooth connection satisfying } \eqref{eq:mmJg:def:middle}, &x_{1,-}^\sharp(x_2) < x_1 <x_{1,+}^\sharp(x_2) \,, \\
\mJg(x_1,x_2,t) = \Jgb(x_1^\sharp(x_2),x_2,t), &x_{1,+}^\sharp(x_2) \leq x_1 \,,
\end{cases}
\end{align} 
where the middle branch in the definition of $\mmJg$ is taken to satisfy  the properties  
\begin{equation}
0 \leq \p_1 \p_t \mmJg(x,t) \leq 2 (1+\alpha) \eps^{-2}\,,
\qquad 
\sabs{\p_2 \p_t \mmJg(x,t)} \leq 250 (1+\alpha) \eps^{-1} \,,
\qquad 
\mmJg (x,t) \leq \mJg(x,t)\,,
\label{eq:mmJg:def:middle} 
\end{equation}
for $x_{1,-}^\sharp(x_2) < x_1 <x_{1,+}^\sharp(x_2) \Leftrightarrow |x_1 - x_1^\sharp(x_2)| < \frac{\eps}{1000}$.
\end{subequations}
Note, in particular, that Proposition~\ref{prop:mJg:well:def} and the last condition in \eqref{eq:mmJg:def:middle} imply that $\mmJg \leq \Jgb$.

\begin{remark}[\bf $\mmJg$ is well-defined and $\mmJg \leq \mJg$]
\label{rem:mmJg:well:def}
In order to see that a smooth (space and time) connection as claimed in the middle branch of the definition \eqref{eq:mmJg:def} is even possible, we first discuss the condition $0 \leq \p_1 \p_t \mmJg(x_1,x_2,t)$, which is equivalent to the fact that for every fixed $(x_2,t)$, the map $x_1 \mapsto \p_t \mmJg(x_1,x_2,t)$ is monotone increasing. For this purpose, we first 
note that since the third branch of \eqref{eq:mmJg:def} is independent of $x_1$, we have that $\p_1 \p_t \mmJg(x_{1,+}^\sharp(x_2)^+,x_2,t) = 0$. Second, we
check that $\p_1 \p_t \mmJg(x_{1,-}^\sharp(x_2)^-,x_2,t) > 0$. To do so, we observe that $x_{1,-}^\sharp(x_2)  \geq x_1^\vee(x_2) + \frac{\eps}{40} - \frac{\eps}{1000} \geq x_1^\vee(x_2)  + \eps^{\frac 74}$, and hence we are in the range of $x_1$ for which assumption~\eqref{item:ic:w0:d11:positive} applies. Moreover, for any $|x_1 - x_1^\sharp(x_2)| \leq \tfrac{\eps}{1000}$, by \eqref{eq:x1:sharp:def} and \eqref{eq:why:the:fuck:not:0} we know that $|(w_0),_1(x_1,x_2) + \frac{17}{20\eps}| \leq \frac{\eps}{1000} \cdot \frac{2}{\eps^2} = \frac{1}{500 \eps}$, and hence $(w_0),_1(x_1,x_2) \leq - \frac{17}{20\eps} + \frac{1}{500 \eps} \leq - \frac{33}{40\eps} < -\frac{1}{3\eps} $. Hence, assumption~\eqref{item:ic:w0:d11:positive} implies that $(w_0),_{11}(x_1,x_2) > \eps^{-\frac 98}$. Combining this information with \eqref{bs-Jg-2a}, we deduce $\p_1 \p_t \mmJg(x_{1,-}^\sharp(x_2)^-,x_2,t)  \geq \frac{1+\alpha}{2}  (w_0),_{11}(x_{1,-}^\sharp(x_2),x_2) - \Cn \eps^{-1} \geq \frac{1+\alpha}{2} \eps^{-\frac 98}  - \Cn \eps^{-1} \geq \frac{1+\alpha}{4} \eps^{-\frac 98}$. Third, we verify that $\p_t \mmJg(x_{1,-}^\sharp(x_2) ,x_2,t) < \p_t \mmJg (x_{1,+}^\sharp(x_2),x_2,t)$. In light of~\eqref{eq:Jgb:identity:0} and of~\eqref{bs-Jg-1b}, this inequality would follow from $(w_0),_{1}(x_{1,-}^\sharp(x_2),x_2) < (w_0),_1(x_1^\sharp(x_2),x_2) - 2 \mathsf{C_{J_t}}$. 
But by the mean value theorem and the previously obtained $(w_0),_{11}(x_1,x_2) > \eps^{-\frac 98}$ we have $(w_0),_1(x_1^\sharp(x_2),x_2) - (w_0),_1(x_{1,-}^\sharp(x_2) ,x_2) > \frac{\eps}{1000} \cdot \eps^{-\frac 98} = \frac{1}{1000} \eps^{-\frac 18} > 2 \mathsf{C_{J_t}}$, assuming $\eps$ is sufficiently small. These three points show that the first condition in \eqref{eq:mmJg:def:middle} (the map $x_1 \mapsto \p_t \mmJg(x_1,x_2,t)$ is monotone increasing) is permissible. Next, we discuss the second condition in \eqref{eq:mmJg:def:middle}, namely $\p_1 \p_t \mmJg(x_1,x_2,t) \leq 2(1+\alpha) \eps^{-2}$. This condition is easy to ensure because $\p_1 \p_t \mmJg(x_{1,+}^\sharp(x_2)^+,x_2,t) = 0$ and by \eqref{bs-Jg-2a} and \eqref{eq:why:the:fuck:not:0} we have $|\p_1 \p_t \mmJg(x_{1,-}^\sharp(x_2)^-,x_2,t)|  \leq \frac{1+\alpha}{2\eps^2}   \|\nb_1^2 w_0\|_{L^\infty} + \Cn \eps^{-1} \leq (1+\alpha) \eps^{-2} + \Cn \eps^{-1} <  2(1+\alpha) \eps^{-2}$. Next, we discuss the third condition in \eqref{eq:mmJg:def:middle}, namely $|\p_2 \p_t \mmJg(x_1,x_2,t)| \leq 2(1+\alpha) \eps^{-2}$. Again, using \eqref{bs-Jg-2a} and \eqref{eq:why:the:fuck:not:0} we have $|\p_2 \p_t \mmJg(x_{1,-}^\sharp(x_2)^-,x_2,t)|  \leq \frac{1+\alpha}{2\eps}   \|\nb_1 \nb_2 w_0\|_{L^\infty} + \Cn \leq (1+\alpha) \eps^{-1} + \Cn < 2(1+\alpha) \eps^{-1}$. On the other hand, by the chain rule and~\eqref{eq:x1:sharp:props} we have that $|\p_2 \p_t \mmJg(x_{1,+}^\sharp(x_2)^+,x_2,t)| \leq |\p_2 \p_t \Jg(x_1^\sharp(x_2),x_2,t)| + |\p_1 \p_t \Jg(x_1^\sharp(x_2),x_2,t)| \cdot |\p_2 x_1^\sharp(x_2)| \leq ((1+\alpha)\eps^{-1} + \Cn)(1+ \eps^{-1} |\p_2 x_1^\sharp(x_2)|) \leq 241 ((1+\alpha)\eps^{-1} + \Cn)  \leq 242 (1+\alpha) \eps^{-1}$. Lastly, we discuss the fourth requirement in \eqref{eq:mmJg:def:middle}, namely the fact that $\mmJg\leq \mJg$. This requirement is already satisfied for $x_1^\sharp(x_2) \leq x_1 \leq x_{1,+}^\sharp(x_2) $ since $\mmJg(x_{1,+}^\sharp(x_2),x_2) = \mJg(x_{1,+}^\sharp(x_2),x_2)$ and $\p_t \p_1 (\mmJg - \mJg) \geq 0$. Using this fact, and that by construction we have $\mmJg(x_{1,-}^\sharp(x_2),x_2) = \mJg(x_{1,-}^\sharp(x_2),x_2)$ it is straightforward to verify that the bound $\mmJg\leq \mJg$ may also be ensured to hold for $x_{1,-}^\sharp(x_2)   \leq x_1 \leq x_1^\sharp(x_2) $. 
\end{remark}

Then,  with $\mmJg$ defined as in~\eqref{eq:mmJg:def}, for $(x,t) \in \TT^2 \times (\initial,\final]$, we define
\begin{align}
\JJ(x_1,x_2,t)
:= 
\begin{cases}
\Jgb(x_1^*(x_2,t),x_2,t), &x_1 \leq x_1^*(x_2,t) \,, \\
\mmJg(x_1,x_2,t)  , &x_1 > x_1^*(x_2,t) \,,
\end{cases}
\label{eq:fake:Jg:def:ds}
\end{align}
and we let $\JJ(x,\initial)=1$ for all $x\in \TT^2$. It is important to observe that $\JJ$ has limited Lipschitz regularity across the hypersurface $\{ x_1 = x_1^*(x_2,t) \}$, but that it is $H^6$ smooth on either side of this hypersurface. We also note that in the upstream region $x_1 \leq x_1^*(x_2,t)$, we have that $\JJ(x,t) = \mathcal{J}(x_2,t)$, as defined in \eqref{eq:fake:Jg:def}. 

\begin{figure}[htb!]
\centering
\begin{minipage}{.4\linewidth}
  \includegraphics[width=.93\linewidth]{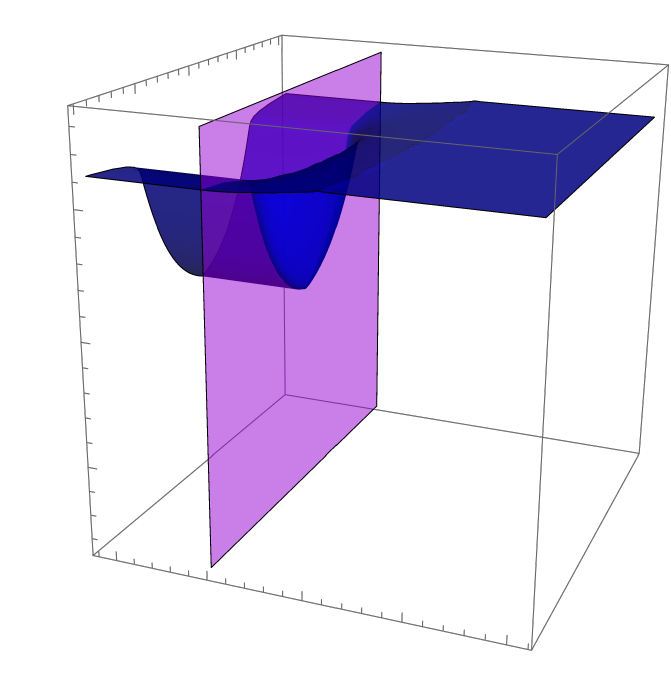}
\end{minipage}
\hspace{.05\linewidth}
\begin{minipage}{.4\linewidth}
  \includegraphics[width=.93\linewidth]{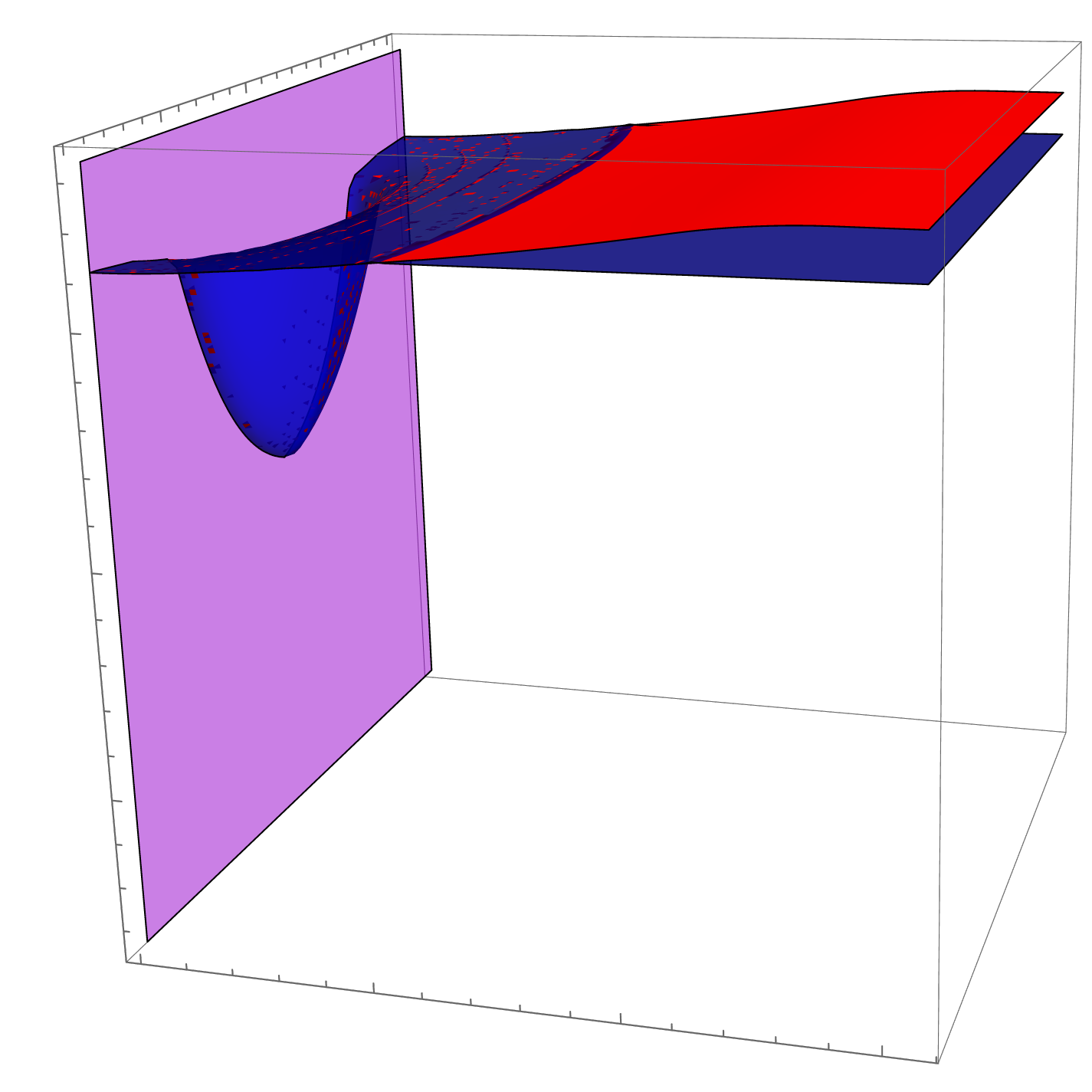}
\end{minipage}
\vspace{-.1 in}
\caption{Consider the functions $\JJ(x_1,x_2,t)$, as defined by~\eqref{eq:fake:Jg:def:ds},  $\Jgb(x_1,x_2,t)$, as defined by~\eqref{eq:Jgb:identity:0}, and $\Jg(x_1,x_2,t)$ as in Figure~\ref{fig:parabolic:spacetime}. The bounding box represents the zoomed in region $-\tfrac{\pi \eps}{4} \leq x_1  \leq \pi\eps $, $|x_2| \leq \tfrac{1}{16}$, and $t\in [-\final,\final]$. 
\underline{Left}: in blue we plot the level-set $\{(x,t) \in \TT^2\times[\initial,\final] \colon \JJ(x_1,x_2,t)=0\}$ which is the  future  temporal boundary (``top'' boundary) of the spacetime $\mPds$ defined in \eqref{eq:spacetime:P} below. This spacetime is the one analyzed in this section. In magenta we plot the surface $\{x_1 = x_1^*(x_2,t)\}$, which separates the upstream from the downstream side.
\underline{Right}: we focus on the downstream side. In blue, we plot the level-set $\{(x,t) \in \TT^2\times[\initial,\final] \colon x_1>x_1^*(x_2,t), \JJ(x_1,x_2,t)=0\}$. As in Figure~\ref{fig:downstream}, in magenta we plot the surface $\{x_1 = x_1^*(x_2,t)\}$, while in red we plot the surface $\{(x,t) \in \TT^2 \times [\initial,\final] \colon x_1 > x_1^*(x_2,t),  \Jgb(x_1,x_2,t)=0 \}$. We emphasize that the blue and red surfaces, which correspond to $\{\JJ=0\}$ and $\{\Jgb =0\}$, respectively, coincide for times $t\leq \medium$. Only for $t\in(\medium,\final]$ do the sets $\{\JJ=0\}$ and $\{\Jgb =0\}$ differ, with the ordering given by $\JJ \leq \Jg$.}  
\label{fig:parabolic-cylinder-flat}
\end{figure}

\begin{remark}[\bf Properties of $\JJ$ on the upstream side are the same as those of $\mathcal{J}$]
\label{rem:oJJ:upstream}
We note that \eqref{eq:fake:Jg:def} implies the equality $\JJ(x,t) = \Jgb(x_1^*(x_2,t),x_2,t) = \mathcal{J}(x_2,t)$,  for all $x_1 < x_1^*(x_2,t)$. As such, all the properties of $\mathcal{J}$ which were established in Sections~\ref{sec:formation:setup} and~\ref{sec:first:consequences} directly carry over to properties of $\JJ$ in the spacetime $\{(x,t) \in \TT^2 \times [\initial,\final] \colon \JJ(x,t) > 0, x_1 < x_1^*(x_2,t) \}$. To avoid redundancy, these estimates will not be re-proven in this section. 
\end{remark}
 
\begin{remark}[\bf The spacetime $\{\JJ>0\}$ terminates before $\final$]
\label{rem:downstream:terminates}
We note that for consistency with the rest of the paper, the spacetime $\{\JJ > 0\}$ is designed to terminate before $\final$. To see this, note that in the upstream region this fact was already established (see Remark~\ref{rem:oJJ:upstream}), while in the downstream region we have
$\JJ \leq \Jgb$, since $\mmJg \leq \Jgb$. By the proof of Lemma~\ref{lem:q:invertible}, for all $x\in \TT^2$ there exists $t_*(x) < \final$ with $\Jgb(x,t_*(x)) = 0$. Hence, for all $x\in \TT^2$ there exists a time $t^\sharp(x)\leq t_*(x) < \final$ such that $\JJ(x,t^\sharp(x)) = 0$. 
\end{remark}

\begin{remark}[\bf The spacetime $\{\JJ >0\}$ captures the downstream \MGHDB\ prior to $\medium$]
\label{rem:oJJ:gg:0:maximal}  
The spacetime $\{\JJ > 0, x_1> x_1^*(x_2,t), t\leq \medium \}$ coincides with the spacetime $\{\Jg > 0, x_1> x_1^*(x_2,t),t\leq \medium \}$, so that by studying $\{\JJ >0, x_1> x_1^*(x_2,t), t\leq \medium\}$, we are indeed analyzing the full development on the downstream side, for times $t\in [\initial,\medium]$. To prove this fact, since $\JJ \leq \Jgb = \Jg$ in the aforementioned spacetime, we only need to show that if $t \leq \medium$ and $x_1> x_1^*(x_2,t)$, then $\JJ(x,t) = 0 \Rightarrow \Jg(x,t) = 0$. To see this, recall that for $t\leq \medium$ we have that $\Jg(x,t) = \Jgb(x,t)$. If we additionally impose $x_1 \leq x_{1,-}^\sharp(x_2)$, then by \eqref{eq:mmJg:def} and \eqref{eq:fake:Jg:def:ds} we have that $\JJ(x,t) = \mmJg(x,t) = \Jgb(x,t) = \Jg(x,t)$, so that $\JJ$ vanishes if and only if $\Jg$ vanishes. On the other hand, if $x_1 > x_{1,+}^\sharp(x_2)$ and $t<\medium$, then  $\JJ$ cannot vanish at $(x,t)$. Indeed, by \eqref{bs-Jg-0} and the definition of $x_1^\sharp(x_2)$ in~\eqref{eq:x1:sharp:def}, we have that $\JJ(x,t) = \Jgb(x_1^\sharp(x_2),x_2,t) = \Jg(x_1^\sharp(x_2),x_2,t) \geq 1 + (t-\initial) \frac{1+\alpha}{2} ( (w_0),_1(x_1^\sharp(x_2),x_2) - \mathsf{C_{J_t}}) 
\geq 1 - (\medium-\initial) \frac{1+\alpha}{2} (\frac{17}{20 \eps} +\mathsf{C_{J_t}}) 
\geq 1 - \frac{51}{50} ( \frac{17}{20} + \eps \mathsf{C_{J_t}}) \geq \frac 18$ if $\eps$ is sufficiently small. By the mean value theorem and the bound \eqref{eq:Jg,_1:sharp},  this also implies that for $x_{1,-}^\sharp(x_2) \leq x_1 \leq x_{1,+}^\sharp(x_2)$ and $t<\medium$, we have $\JJ(x,t) = \Jg(x,t) \geq \Jg(x_1^\sharp(x_2),x_2,t) - \frac{\eps}{1000} \|\Jg,_1\|_{L^\infty} \geq \frac 18 - \frac{3}{1000} \geq \frac 19$. As such, $\JJ$ also cannot vanish in this region, completing the proof of $\JJ = 0 \Rightarrow \Jg = 0$ when $t\leq \medium$.
\end{remark}

For downstream \MGHDB, we use the spacetime set given by
\begin{align}
\mPds := \bigl\{ (x,t) \in \mathbb{T}^2 \times [\initial,\final) \colon \JJ (x_1,x_2,t) > 0   \bigr\}
 \,,
 \label{eq:spacetime:P}
\end{align}
We define our spacetime ``flattening'' transformation by
\begin{subequations}
\label{eq:t-to-s-transform:all-P}
\begin{align} 
\qds &\colon \mPds \mapsto   [0,\eps) \\
\s = \qds(x_1, x_2,t)&:= \eps \bigl(1 - \JJ(x_1,x_2,t) \bigr) \,.
\label{t-to-s-transform-P}
\end{align} 
\end{subequations}
We have that $\qds(x_1,x_2,\initial) = 0$, so that the set  $\{\s = 0\}$ corresponds to the initial time slice $\{t=\initial\}$, or, the  {\it past} temporal boundary of $\mPds$; meanwhile, the {\it future} temporal boundary of  
$\mPds$ is flattened to the set $\{\s = \eps\}$.  
We note that in the upstream region $x_1 \leq x_1^*(x_2,t)$, the map $\qds = \qds(x,t)$ defined in \eqref{t-to-s-transform-P} is precisely equal to the map $\mathfrak{q}= \mathfrak{q}(x_2,t)$ defined in \eqref{t-to-s-transform}. Only in the downstream region $x_1 > x_1^*(x_2,t)$ do $\qds(x,t)$ and $\mathfrak{q}(x_2,t)$ differ.

The  inverse of $\qds$ is defined by
\begin{subequations} 
\label{s-to-t-transform-P}
\begin{align} 
\qds^{-1}&: \TT^2 \times [0,\eps)  \to [\initial,\final)  \,, \\
t & = \qds^{-1}(x_1,x_2,\s) \,,
\end{align} 
\end{subequations} 
such that $t = \qds^{-1}(x_1,x_2,\qds(x_1,x_2,t))$ for all $(x_1,x_2,t) \in \mPds$, or equivalently, that 
$\s = \qds(x_1,x_2,\qds^{-1}(x_1,x_2,\s))$ for all $(x_1,x_2,\s) \in \TT^2 \times [0,\eps)$. 
In \eqref{s-to-t-transform-P}, 
we are once again abusing convention: it is the map $(x_1,x_2,\s) \mapsto (x_1,x_2,t)$ defined from 
$\TT^2 \times [0,\eps) \to \TT^2 \times [\initial,\final)$ which is the inverse of the map 
$(x_1,x_2,t) \mapsto (x_1,x_2,\s) = (x_1, x_2,\qds(x_1,x_2,t))$. 
The fact that such a map is well-defined is established in Lemma~\ref{lem:q:invertible-DS} below.  

\subsection{Change of coordinates for the remapped spacetime}
\label{sec:remapping:2}
Given any function $f\colon \mPds \to \mathbb{R}$, we define the function $\tilde f \colon \TT^2 \times [0,\eps) \to \mathbb{R}$ by
\begin{equation}
\tilde f(x,\s) := f(x,t), \qquad \mbox{where} \qquad \s = \qds(x,t)   \,.
\label{eq:f:tilde:f-P}
\end{equation}
Then, by the chain-rule and \eqref{t-to-s-transform-P},  we obtain  that
\begin{subequations}
\label{eq:xt:xs:chain:rule-P}
\begin{align} 
\p_t f(x,t) &= \Qd(x,\s)   \p_\s \tilde f(x,\s) \,, \label{ddt-dds-P} \\
\p_if(x,t) &=  \bigl(\p_i - \Qb_i(x,\s)  \p_\s\bigr) \tilde f(x,\s)\,,  \qquad i=1,2\,,
\end{align} 
\end{subequations}
where for compactness of notation we have introduced the functions 
\begin{subequations} 
\label{QQQ-P}
\begin{align}
\Qd (x,\s) &
= - \eps (\p_t \JJ)(x,t) \Big|_{t= \qds^{-1}(x,\s)}
\label{eq:QQQ:a-P}
\\
\Qb_i(x,\s) &=   \eps  (\p_i \JJ)(x,t)   \Big|_{t= \qds^{-1}(x,\s)} \ \text{ for } \ i=1,2 \,.
\label{eq:QQQ:aa-P}
\end{align}
Note that $\Qb_1=0$ for $x_1 \leq x_1^*(x_2, t)|_{t= \qds^{-1}(x,\s)}$.
We also  define 
\begin{equation}
\Q(x,\s )  := \Qd(x,\s) - \tilde V(x,\s)  \Qb_2(x,\s) 
= - \eps (\p_t \JJ + V \p_2 \JJ)(x,t)  \Big|_{t= \qds^{-1}(x,\s)}\,,
\label{eq:QQQ:b-P}
\end{equation} 
and
\begin{equation} 
\Qc = \p_\s \Q \,, \qquad 
\Qr_\s  = \p_\s \Qd  \,, \qquad 
\Qr_1= \p_\s\Qb_1 \,,  \qquad
\Qr_2= \p_\s\Qb_2 
\,.
\label{eq:QQQ:c-P}
\end{equation} 
\end{subequations} 
With the above notation, it follows from \eqref{eq:xt:xs:chain:rule} that the spacetime gradient operator in  $(x,t)$ variables, namely $\nb = (\eps\p_t, \eps \p_1, \p_2)$, becomes the gradient operator $\nbs$ associated with the $(x,\s)$ coordinates, which is defined by
\begin{align} 
\nbs = (\nbs_\s, \nbs_1, \nbs_2) := \big( \eps \Qd  \p_\s ,  \eps( \p_1 -\Qb_1\p_\s) ,  \p_2 -\Qb_2 \p_\s \big) \,. 
\label{nb-s-P}
\end{align} 

Therefore,  $\nb f (x,t) = \nbs \tilde f(x,\s)$, and the components of $\nbs$ commute:
\begin{align} 
\jump{\nbs_\s, \nbs_2} =
\jump{\nbs_\s, \nbs_1} =
\jump{\nbs_2, \nbs_1} =
0
\,.
\label{comm-nbs1-nbs2-P}
\end{align} 
For any $\gamma \in \mathbb{N}_0^3$,  $\nbs^\gamma = \nbs_\s^{\gamma_0} \nbs_1^{\gamma_1} \nbs_2^{\gamma_2}$, and 
\begin{equation}
\label{nb-to-nbs-P}
(\nb^\gamma f) (x,t) = (\nbs^\gamma \tilde f)(x,s)
\,.  
\end{equation}
From the identity $\Q\p_\s +V \p_2= \tfrac{1}{\eps}\nbs_s + V\nbs_2$, we note that   material derivatives are mapped into $(x,\s)$ coordinates as
\begin{align} 
(\p_t  +V \p_2) f(x,t)  = (\Q\p_\s + \tilde V\p_2) \tilde f(x,\s)  = ( \tfrac{1}{\eps}\nbs_\s + \tilde V\nbs_2) \tilde f(x,\s) 
\,. \label{the-time-der-P}
\end{align} 
It also follows from \eqref{comm-nbs1-nbs2-P} and the second equality in \eqref{the-time-der-P} that 
\begin{align} 
\jump{ (\Q\p_\s +\tilde V \p_2) , \nbs^k} \tilde f  = \jump{\tilde V, \nbs^k} \nbs_2 \tilde f 
= -\nbs^k \tilde V \, \nbs_2 \tilde f - \doublecom{\nbs^k, \tilde V, \nbs_2 \tilde f} 
\,. \label{good-comm-P}
\end{align}

With the notation in~\eqref{eq:f:tilde:f-P} and~\eqref{nb-s-P}, the definition~\eqref{t-to-s-transform-P} implies the following identities for $\tilde \JJ$ and $\nbs \tilde \JJ$:
\begin{equation}
\label{eq:D:JJ} 
\tilde \JJ(x,\s) = 1 - \tfrac{\s}{\eps}\,,
\qquad
\nbs_\s \tilde \JJ  = - \Qd\,,
\qquad
\nbs_1 \tilde \JJ = \Qb_1 \,,
\qquad 
\nbs_2 \tilde \JJ  =  \Qb_2 \,.
\end{equation}
The identities in \eqref{eq:D:JJ} will be used implicitly throughout the section.

\subsection{Adjoint formulas}
\label{sec:adjoint-P}
We define 
\begin{equation}
\tx(x_2,\s) = x_1^*(x_2,t)|_{t= \qds^{-1}(x,\s)}
\,.
\label{eq:tx:x2:s}
\end{equation}  
Note that by construction we have $\Qb_1(x,\s)=0$ for $x_1 \leq x_1^*(x_2, t)|_{t= \qds^{-1}(x,\s)} = \tx(x_2,\s)$, and so it follows via~\eqref{eq:D:JJ} that
\begin{equation} 
 \nbs_1\tilde{\JJ} (\tx(x_2,\s), x_2, \s) = \Qb_1(\tx(x_2,\s),x_2,\s) = 0 \,.
 \label{Qb1-xstar}
\end{equation} 
Using the notation in~\eqref{eq:tx:x2:s}, we define
\begin{equation} 
\Gamma(\s) =  \bigcup_{\s' \in [0,\s]} \big( \tx(x_2,\s') ,x_2,\s'\big) \label{interface}
\end{equation} 
and hence \eqref{Qb1-xstar} implies $\nbs_1\tilde{\JJ}=0$ on $\Gamma(\eps)$.

For energy estimates we shall work on the union of the following two sets:
\begin{subequations} 
\label{twosets}
\begin{align}
\Pdp(\s) &:= \{ (x,\s') \in  \mathbb{T} ^2 \times [0,\s] \colon x_1 > \tx(x_2,\s') \} \\
\Pdm(\s) &:= \{ (x,\s') \in  \mathbb{T} ^2 \times [0,\s] \colon x_1 < \tx(x_2,\s') \} 
 \,, \\
 \Pdpm(\s) &= \Pdp(\s) \cup \Pdm(\s) \,.
\end{align}
\end{subequations} 
With \eqref{interface} and \eqref{twosets}, we see that $\Pdpm(\eps) \cup \Gamma(\eps) =   \mathbb{T} ^2 \times [0,\eps]$
is the entire spacetime set, but energy estimates will be set on $\Pdm(\s) \cup \Pdp(\s)$, thus avoiding differentiation across $\Gamma(\s)$.  We similarly consider the associated time-slices
\begin{subequations}
 \label{twosets-slices}
\begin{align}
\Xdp(\s) &:= \{ x \in  \mathbb{T} ^2 \colon x_1 > \tx(x_2,\s) \} \,,
\\
\Xdm(\s) &:= \{ x \in  \mathbb{T} ^2  \colon x_1 < \tx(x_2,\s) \} 
 \,, \\
\Xdpm(\s) &= \Xdp(\s) \cup \Xdm(\s) \,.
\end{align}
\end{subequations}
Using the notation in~\eqref{twosets} and~\eqref{twosets-slices}, throughout this section we shall use the following notation for norms:
for $0 \le \s \le \eps$,  we  denote 
\begin{subequations}
\label{eq:norm:notation:ds}
\begin{align}
\| F(\cdot ,\s)  \|_{L^2} = \| F(\cdot ,\s)  \|_{L^2_x} &:= \| F(\cdot ,\s) \|_{L^2(\Xdpm(\s))}
\,,\\
\| F  \|_{L^\infty_\s L^2_x} &:= \sup_{\s \in [0,\eps]} \| F(\cdot ,\s) \|_{L^2(\Xdpm(\s))}
\,,\\
\| F  \|_{L^2_{x,\s}} &:= \| F \|_{L^2(\Pdpm(\s))} \,,
\\
\| F  \|_{L^\infty_{x,\s}} &:= \| F \|_{L^\infty(\Pdpm(\s))} \,.
\end{align}  
\end{subequations}
In view of \eqref{eq:Qd:bbq-DS}, the same comments as in Remark~\ref{rem:L2:norms:x:s:x:t:A} (regarding equivalence of $L^2_{x,\s}$ and $L^2_{x,t}$, respectively $L^\infty_{x,\s}$ and $L^\infty_{x,t}$) also apply to the norms in~\eqref{eq:norm:notation:ds}.

Next, we compute adjoints of the differential operator $\nbs$ (as defined by \eqref{nb-s-P}), with respect to the $L^2$ inner product on $\Pdpm(\s)$. We claim that
\begin{subequations}
\label{eq:adjoints-P}
\begin{align}
\nbs_\s^* &=  - \nbs_\s - \eps  \Qr_\s + \eps \Qd ( \delta_{\s} - \delta_{0} ) \,, \label{adjoint-s-P}
\\
\nbs_1^* &=  - \nbs_1   + \eps \Qr_1 - \eps \Qb_1 \delta_{\s}   \,,  \label{adjoint-1-P}
\\
\nbs_2^* &=  - \nbs_2 + \Qr_2 -  \Qb_2 \delta_{\s} \,,
 \label{adjoint-2-P}
\\
(\Q \p_\s + \tilde V \p_2)^* & = - (\Q  \p_\s + \tilde V \p_2)  - \Qr_\s + \Q  ( \delta_{\s} - \delta_{0} ) + \tilde V  \Qr_2   
- \nbs_2 \tilde V
 \label{adjoint-3-P}
 \,.
\end{align}
\end{subequations}
We note that from the definition \eqref{QQQ-P}, we have
$\Qb_1(x,\s) = \Qr_1(x ,\s) =0$ for $ x_1 < \tx(x_2,\s)$ and so the identities \eqref{eq:adjoints-P} precisely match the previously established formulas \eqref{eq:adjoints} in the upstream region.
The identities in~\eqref{eq:adjoints-P} follow from \eqref{nb-s-P}, 
and from \eqref{nbs1-IBP}, \eqref{nbss-IBP}, and \eqref{nbs2-IBP2}, which we establish below. 

We first consider $\nbs_1^*$. From \eqref{nb-s-P}, we have that
\begin{align}
\nbs_1
:= 
\begin{cases}
 \eps \p_1 \,,  &x_1 < \tx(x_2,\s) \,, \\
 \eps( \p_1 -\Qb_1\p_\s)\,&x_1 > \tx(x_2,\s)  \,.
\end{cases}
\label{nbs1-P}
\end{align}
By also appealing to \eqref{Qb1-xstar}, it follows that\footnote{Our computation in \eqref{nbs1-IBP} makes use of the equality
$(FG)(\tx(x_2,\s')^+,x_2,\s')=(FG)(\tx(x_2,\s')^-,x_2,\s')$, which holds whenever the trace of $FG$ on $\Gamma(\s)$ from the domains 
$\Pdm(\s)$  and $\Pdp(\s)$ exist and agree.   Note that $\Gamma(s)$ is contained in the spacetime $\hat{ \mathcal{P} }$ defined in
\eqref{eq:spacetime:temp:x2:t}, and that we have established the existence of a unique $H^7$-class Euler solution in $\hat{ \mathcal{P} }$ with uniform bounds
in the norms \eqref{eq:norms:L2:first} for the $H^7$-class initial data specified in Section \ref{cauchydata}.  We observe that by increasing
the regularity of our initial data to $H^k$ for $k \ge 8$, then our argument would yield a unique  $H^k$-class Euler solution in $\hat{ \mathcal{P} }$ with uniform bounds in the $(k-1)$-order version of the sixth-order norms in \eqref{eq:norms:L2:first}.  In particular, in applying the trace to the functions
$F$ and $G$ in \eqref{nbs1-IBP}, we can assume that we have sufficiently regularity for the cancellation of the two-sided trace to hold.
\label{footnote-trace}
}
\begin{align} 
\int_{\Pdpm(\s)} \nbs_1 F \, G  \, {\rm d} x {\rm d}\s'
&=  \eps  \int_0^\s \int_{-\pi}^\pi\int_{\tx(x_2,\s')}^\pi  ( \p_1 -\Qb_1\p_\s)F \, G  \, {\rm d} x_1 {\rm d}x_2 {\rm d}\s'
+  \eps  \int_0^\s \int_{-\pi}^\pi\int_{-\pi}^{\tx(x_2,\s')}   \p_1F \, G \, {\rm d} x_1 {\rm d}x_2 {\rm d}\s'\notag \\
&= \eps  \int_0^\s \int_{-\pi}^\pi  (F\,G)(y_1,x_2,\s')\Big|_{y_1=\tx(x_2,\s')^+} ^{y_1=\tx(x_2,\s')^-} {\rm d}x_2 {\rm d}\s' \notag \\
& \qquad \qquad
  -\eps \int_{\Pdpm(\s)}  F \, \p_1G \, {\rm d} x_1 {\rm d}x_2 {\rm d}\s'
 - \eps  \int_0^\s \int_{-\pi}^\pi\int_{\tx(x_2,\s')}^\pi  \Qb_1 G\, \p_\s F \, {\rm d} x_1 {\rm d}x_2 {\rm d}\s' \notag \\
 & =- \int_{\Pdpm(\s)}  F \, \nbs_1G \, {\rm d} x_1 {\rm d}x_2 {\rm d}\s'
 +  \eps  \int_{\Pdp(\s)} \Qr_1 \,  G\,  F \, {\rm d} x_1 {\rm d}x_2 {\rm d}\s'
 -  \eps  \int_{\Xdp(\s)}  \Qb_1 G\,  F\, {\rm d} x_1 {\rm d}x_2  
  \,,  \label{nbs1-IBP}
\end{align} 
where we have used \eqref{eq:QQQ:c-P} to write $\Qr_1 = \p_\s \Qb_1$ and \eqref{eq:QQQ:aa-P} to conclude that
$\Qb_1(x,0) =0$. This proves~\eqref{adjoint-1-P}.

Second, we consider $\nbs_\s^*$. 
Employing the fundamental theorem of calculus, we obtain the identities
\begin{align} 
\frac{d}{d\s}  \int_{-\pi}^\pi\int_{\tx(x_2,\s')}^\pi  F(x,\s) \mathsf{G}(x,\s) {\rm d}x_1 {\rm d}x_2
&= \int_{-\pi}^\pi\int_{\tx(x_2,\s)}^\pi  \p_\s \big( F(x,\s) \mathsf{G}(x,\s) \big){\rm d}x_1 {\rm d}x_2
\notag \\
& \ \
- \int_{-\pi}^\pi  \p_\s \tx(x_2,\s)  F(\tx(x_2,\s),x_2,\s) \mathsf{G}(\tx(x_2,\s),x_2,\s) {\rm d}x_2 \,,
\label{ftoc1}
\end{align} 
and 
\begin{align} 
\frac{d}{d\s}  \int_{-\pi}^\pi\int_{-\pi}^{\tx(x_2,\s')} F(x,\s) \mathsf{G}(x,\s) {\rm d}x_1 {\rm d}x_2
&= \int_{-\pi}^\pi\int_{-\pi}^{\tx(x_2,\s)}  \p_\s \big( F(x,\s) \mathsf{G}(x,\s) \big){\rm d}x_1 {\rm d}x_2
\notag \\
& \ \
+\int_{-\pi}^\pi  \p_\s \tx(x_2,\s)  F(\tx(x_2,\s),x_2,\s) \mathsf{G}(\tx(x_2,\s),x_2,\s) {\rm d}x_2 \,.
\label{ftoc2}
\end{align} 
The $\p_\s \tx(x_2,\s) $ present in \eqref{ftoc1}--\eqref{ftoc2} may be computed from \eqref{ddt-dds-P}, \eqref{eq:QQQ:a-P}, and \eqref{p2-and-pt-x1star}, as
\begin{equation} 
\p_\s \tx(x_2,\s) 
= \tfrac{1 }{\Qd(x,\s)}  \p_t x_1^*(x_2,t)|_{t= \qds^{-1}(x,\s)}
= - \tfrac{1}{\Qd(x,\s)} \bigl(  \tfrac{ \p_t \p_1 \Jg }{\p_1\p_1 \Jg}\bigr) (x_1^*(x_2,t),t) \big|_{t= \qds^{-1}(x,\s)}  \,.
\label{dds-x1star}
\end{equation} 
From~\eqref{eq:Jg:11:lower}, \eqref{eq:Jg:1t:upper}, and \eqref{eq:Qd:bbq-DS} below, we may deduce that $|\p_\s \tx| \les \eps$, so that this term is finite. This allows us to add together \eqref{ftoc1} and \eqref{ftoc2}, which shows that the two integrals evaluated along $\Gamma(\s)$ cancel each other, and upon integrating in $\s$ deduce that
\begin{align*}
\int_{\Xdpm(\s')} F(x,\s') \mathsf{G}(x,\s') {\rm d}x_1 {\rm d}x_2\Big|_{\s'=0}^{\s'=\s}
&= \int_{\Pdpm(\s)} \p_\s \big( F(x,\s') \mathsf{G}(x,\s') \big){\rm d}x_1 {\rm d}x_2{\rm d}\s'
\end{align*} 
By setting $\mathsf{G} = \eps\Qd G$, we have that
\begin{align}
 \int_{\Pdpm(\s)} \nbs_\s F \, G {\rm d}x_1 {\rm d}x_2{\rm d}\s'
&=
\eps \int_{ \Xds }  \Qd \, F\, G   {\rm d}x_1 {\rm d}x_2\Big|_{0}^{\s}
-  \int_{\Pdpm(\s)}   F\, \nbs_\s G   {\rm d}x_1 {\rm d}x_2 {\rm d}\s'
- \eps \int_{\Pdpm(\s)}     \Qr_\s \,  F G   {\rm d}x_1 {\rm d}x_2 {\rm d}\s' \,, \label{nbss-IBP}
\end{align} 
which proves~\eqref{adjoint-s-P}.

Third, we consider $\nbs_2^*$. By repeating the computations leading to \eqref{ftoc1} and \eqref{ftoc2} but with $\p_\s$ replaced by $\p_2$, we find that
\begin{align} 
\int_{\Xdpm(\s)} \p_2 \big( F(x,\s') \mathsf{G}(x,\s') \big) {\rm d}x_1 {\rm d}x_2 {\rm d}\s' =0 \,.
\label{eq:vomit:cascade:1}
 \end{align} 
Letting $\mathsf{G} = \Qb_2 G$, we find that
 \begin{align} 
\int_{\Pdpm(\s)} \nbs_2 F \, G 
& =  - \int_{\Pdpm(\s)}  F \, \nbs_2G 
 +   \int_{\Pdp(\s)} \Qr_2   G\,  F 
 -  \int_{-\pi}^\pi\int_{\tx(x_2,\s)}^\pi (\Qb_2 G\,  F)(x,\s) {\rm d}x_1 {\rm d}x_2\,.  \label{nbs2-IBP2}
\end{align} 
and thus~\eqref{adjoint-2-P} follows. Identity~\eqref{adjoint-3-P} is a consequence of~\eqref{adjoint-s-P}--\eqref{adjoint-2-P}, since $\Q\p_\s+ \tilde V \p_2 = \frac{1}{\eps} \nbs_\s + V \nbs_2$.

Finally, we note that by virtue of the cancellation  of the ``boundary'' integrals evaluated along $\Gamma(\s)$,  the sum of the identities \eqref{ftoc1} and \eqref{ftoc2} shows that
\begin{align} 
\frac{d}{d\s}  \int_{\Xdpm(\s)} F(x,\s) {\rm d}x =  \int_{\Xdpm(\s)} \p_\s F(x,\s) {\rm d}x  \,. \label{ds-exchange}
\end{align}

\begin{remark}[\bf Dropping the tildes]
\label{rem:no:tilde-P}
For notational convenience, we shall once again drop the tildes from all the variables in $(x,\s)$ coordinates. Dropping tildes on the fundamental Euler variables and on the geometric variables is done in direct analogy with~Remark~\ref{rem:no:tilde}. Notably, we shall denote $\tilde \Jg, \tilde \JJ, \widetilde{\Jgb}$ simply as $\Jg, \JJ,  \Jgb$. This identification is made throughout the rest of theis section and no ambiguity may arise because we shall still use the notation $\nbs$ for the spacetime derivative operator in $(x,\s)$ coordinates. As such, $\nbs f$ means that $f$ is viewed as a function of $(x,\s)$, while $\nb f$ means that $f$ is viewed as a function of $(x,t)$, where $t= \qds^{-1}(x,\s)$.
\end{remark}

\subsection{The $L^2$-based energy norms}
\label{sec:norms:L2:P}
For downstream \MGHDB, with the notation in~\eqref{eq:norm:notation:ds} we  define the energy norms by
\begin{subequations}
\label{eq:tilde:E5E6} 
\begin{alignat}{2}
\widetilde{\mathcal{E}}_{6}^2(\s) 
&=  \widetilde{\mathcal{E}}_{6,\nnn}^2(\s)  + (\mathsf{K}\eps)^{-2} \widetilde{\mathcal{E}}_{6,\ttt}^2(\s) \,,
\ \ \ 
&&\widetilde{\mathcal{E}}_5^2(\s) 
=  \widetilde{\mathcal{E}}_{5,\nnn}^2(\s) + (\mathsf{K}\eps)^{-2} \widetilde{\mathcal{E}}_{5,\ttt}^2(\s)
\label{eq:tilde:E5E6:N+T} 
\\
\widetilde{\mathcal{E}}_{6,\nnn}^2(\s) 
&= \snorm{   \JJtf \Jgh\nbs^6 (\Jg\Wbn,\Jg\Zbn, \Jg\Abn)(\cdot,\s)}^2_{L^2}  \,,
\ \ \ 
&&\widetilde{\mathcal{E}}_{5,\nnn}^2(\s)
= \snorm{ \Jgh \nbs^5 (\Jg\Wbn, \Jg\Zbn, \Jg\Abn)(\cdot,\s)}^2_{L^2}
 \\
\widetilde{\mathcal{E}}_{6,\ttt}^2(\s)
&=\snorm{ \JJtf \Jgh \nbs^6 (\Wbt, \Zbt, \Abt)(\cdot,\s)}^2_{L^2}  \,, 
\ \ \
&&\widetilde{\mathcal{E}}_{5,\ttt}^2(\s)
= \snorm{ \Jgh \nbs^5 ( \Wbt, \Zbt, \Abt)(\cdot,\s)}^2_{L^2} \,, 
\end{alignat}
\end{subequations}
and the damping norms by
\begin{subequations}
\label{eq:tilde:D5D6}
\begin{alignat}{2}
\widetilde{\mathcal{D}}_6^2(\s)
&= \widetilde{\mathcal{D}}_{6,\nnn}^2(\s) + (\mathsf{K} \eps)^{-2} \widetilde{\mathcal{D}}_{6,\ttt}^2(\s)  \,,
\ \ \
&&\widetilde{\mathcal{D}}^2_5(\s)
 = \widetilde{\mathcal{D}}^2_{5,\nnn}(\s) + (\mathsf{K}\eps)^{-2} \widetilde{\mathcal{D}}^2_{5,\ttt}(\s) 
 \label{eq:tilde:D5D6:N+T} 
 \\
\widetilde{\mathcal{D}}_{6,\nnn}^2(\s) 
&=   
\int_0^\s\snorm{ \JJof \Jgh \nbs^6 (\Jg \Wbn, \Jg\Zbn, \Jg\Abn)}^2_{L^2} {\rm d}\s'  \,, 
\ \ \ 
&&\widetilde{\mathcal{D}}^2_{5,\nnn}(\s)
 =    \int_0^\s  \snorm{ \nbs^5 (\Jg \Wbn, \Jg\Zbn, \Jg\Abn)}^2_{L^2} {\rm d}\s' \,, \\
\widetilde{\mathcal{D}}_{6,\ttt}^2(\s) 
&=  \int_0^\s \snorm{ \JJof \Jgh\nbs^6 (\Wbt, \Zbt, \Abt)}^2_{L^2} {\rm d}\s'  \,,
\ \ \
&&\widetilde{\mathcal{D}}^2_{5,\ttt}(\s)
=   \int_0^\s \snorm{ \nbs^5 (\Wbt, \Zbt, \Abt)}^2_{L^2}{\rm d}\s'  \,.
\end{alignat}
\end{subequations}
where once again $\mathsf{K} \geq 1$ is a sufficiently large constant, independent of $\eps$, chosen at the end of the 
proof, solely in terms of  $\alpha$ and $\kappa_0$ (see~\eqref{eq:K:choice:1-P} below).

\subsection{Bootstrap assumptions} 
We continue to use the same bootstrap assumptions as in~\eqref{bootstraps}, but instead of assuming that these bootstraps hold in the spacetime $\mathcal{P}$ (cf.~\eqref{eq:spacetime:smooth}), we now assume that these bootstraps hold in the $\mPds$ spacetime (cf.~\eqref{eq:spacetime:P}). As such, in this section all pointwise bootstraps are assumed to hold for $(x,t) \in \mPds$, or equivalently, for all  $(x,\s) \in   \TT^2\times[0,\eps)$ via the flattening map $\qds$ (cf.~\eqref{eq:t-to-s-transform:all-P}), and for the energy and damping norms defined earlier in~\eqref{eq:tilde:E5E6} and~\eqref{eq:tilde:D5D6}. 
We continue to follow the convention in Remark~\ref{rem:no:tilde-P} and drop the tildes from all fundamental Euler variables, all the geometric variables, and on the flattening coefficients.  

To be more precise, the working bootstrap assumptions in this section are that 
\begin{subequations}
\label{bootstraps-P}
\begin{align}
&( {\Wb}, {\Zb}, {\Ab}, {\Jg}, {h}, {V}, {\Sigma}) 
  \mbox{ satisfy the pointwise  bootstraps }   
\eqref{bs-supp}\!-\!\eqref{bs-D-Sigma} \mbox{ in } \mPds \Leftrightarrow \TT^2 \times [0,\eps)\label{boots-PP}
\,,\\  
&\widetilde{\mathcal{E}}_{6} ,
\widetilde{\mathcal{D}}_{6} ,
\widetilde{\mathcal{E}}_{5} ,
\widetilde{\mathcal{D}}_{5} ,
\snorm{\nbs^6 \nbs_1 \tilde h}_{L^2_{x,\s}},
\snorm{\nbs^6 \nbs_2 \tilde h}_{L^2_{x,\s}}, 
\snorm{\nbs^6 \tilde \Jg}_{L^2_{x,\s}} 
  \mbox{ satisfy the energy  bootstraps }   
\eqref{bootstraps-Dnorm:6}\!-\!\eqref{bootstraps-Dnorm:Jg}. \label{boots-P}
\end{align}
\end{subequations}
Here $( {\Wb}, {\Zb}, {\Ab}, {\Jg}, {h}, {V}, {\Sigma})$ are defined according to the 
flattening~\eqref{eq:f:tilde:f-P}, and the energy and damping norms are defined in~\eqref{eq:tilde:E5E6} and~\eqref{eq:tilde:D5D6}, 
respectively. Since the bootstraps~\eqref{bootstraps-P} in this section  are the same as the bootstraps~\eqref{bootstraps} used in Sections~\ref{sec:formation:setup}--\ref{sec:sixth:order:energy}, save for the different weights in the $L^2$ norms (see~\eqref{eq:tilde:E5E6} and~\eqref{eq:tilde:D5D6}), we shall sometimes (more frequently for the pointwise bootstraps) make reference to \eqref{bootstraps} instead of \eqref{bootstraps-P}. 

As in Sections~\ref{sec:formation:setup}--\ref{sec:sixth:order:energy}, the burden of the proof in the current section is to close the bootstrap assumptions~\eqref{bootstraps-P}, i.e., to show that these bounds hold with $<$ symbols instead of $\leq$ symbols. To avoid redundancy, we do not repeat the arguments of how bootstraps are closed when the proof is either identical to that in given earlier in Sections~\ref{sec:formation:setup}--\ref{sec:sixth:order:energy}, or if it requires infinitesimal and straightforward adjustments. Instead, we focus on the proofs of those bootstraps which are different because of either the $x_1$-dependence of $\qds$  manifested through the fact that $\nbs_1^* \neq - \nbs_1$ (see~\eqref{adjoint-1-P}), or because of the fact that the weight $\JJ$ used in the energy estimates  satisfies  $\nbs_1 \JJ \neq 0$ (see~\eqref{eq:D:JJ}). The remainder of this section is dedicated to closing the bootstrap assumptions~\eqref{bootstraps-P}.

In the process of closing the bootstrap estimates~\eqref{bootstraps-P} we make use of the functional analytic framework in the flattened domain developed in Appendix~\ref{app:functional}, with the modifications described in Section~\ref{rem:app:downstream:flat}. We shall also utilize the pointwise estimates for objects that naturally flow along the $1$-and $2$-characteristics, as developed in Appendix~\ref{sec:app:transport}, keeping in mind Section~\ref{app:downstream:Lp}, where the modifications due to the $\qds$ flattening map are discussed.

\subsection{Consequences of the bootstraps and updated estimates for downstream development}
\label{sec:bootstrap:consequences:DS}
Many of the bounds established in Sections~\ref{sec:first:consequences}--\ref{sec:pointwise:bootstraps} are direct consequences of the bootstrap assumptions, the functional analytic setup in Appendix~\ref{app:functional}, and of the $L^\infty$ estimates from Appendix~\ref{sec:app:transport}. We emphasize that many of these arguments apply {\bf as is} in the geometry of the downstream \MGHDB, except that we refer to the bootstraps~\eqref{bootstraps-P}, to Section~\ref{rem:app:downstream:flat}, and to~Section~\ref{app:downstream:Lp}. Examples include: the bounds for the diameter of the support in Section~\ref{sec:spatial:support}, the bounds for the ALE flow in Section~\ref{sec:ALE:flow}, the pointwise bounds for $(W,Z,A)$ in Section~\ref{sec:WZA:pointwise}, the pointwise bounds for $\nb^k (\Jg \Wbn)$ and $\nb^k   \Jg  $ when $0 \leq k \leq 2$ from Section~\ref{sec:Wbn:Jg:pointwise}, the properties of $x_1^*(x_2,t)$ and of the curve of pre-shocks in Section~\ref{sec:x1star}, the damping properties of $\Jg$ and $\mathcal{J}$ from Section~\ref{sec:Jg:properties}, 
and the closure of the bootstrap for the fifth order derivatives (cf.~\eqref{bootstraps-Dnorm:5}) discussed in Section~\ref{sec:D5:bootstrap} (here we only use that $\JJ = 1- \frac{\s}{\eps}$). These arguments are not repeated here. In fact, because these bounds are the same, throughout this section we abuse notation and make reference to equation numbers from Sections~\ref{sec:formation:setup}--\ref{sec:first:consequences}.

The only estimates from Section~\ref{sec:first:consequences} which require a substantial modification specific to downstream \MGHDB\ are the bounds for the remapping coefficients $(\Qd,\Q,\Qb_i,\Qc,\Qr_i)$ defined in \eqref{QQQ-P}. The analysis in Section~\ref{sec:Q:coeff}, or more precisely, the bounds from~Lemma~\ref{lem:Q:bnds}, are to be replaced with the estimates obtained in Lemma~\ref{lem:Q:bnds-DS} below. 
 
\begin{lemma}
\label{lem:Q:bnds-DS}
Assume that the bootstrap bounds \eqref{bootstraps-P} hold on $\mPds$.
If $\eps$ is taken to be sufficiently small with respect to $\alpha, \kappa_0$, and $\Cdata$, then the   functions
$(\Qd,\Q,\Qb_i,\Qc,\Qr_i)$ defined in \eqref{QQQ-P} satisfy the bounds 
\begin{subequations}
\label{eq:Q:all:bbq-DS}
\begin{align} 
\tfrac{17(1+\alpha)}{40} - \Cn \eps   \mathsf{C_{J_t}}
&\leq \Qd \leq 401(1+\alpha)\,,
\label{eq:Qd:bbq-DS}
\\
\tfrac{2(1+\alpha)}{5} &\leq \Q \leq  402 (1+\alpha)  \, ,
\label{eq:Q:bbq-DS} \\
0 &\leq \Qb_1 \leq 5 {\bf 1}_{\Pdp(\s)} \,, 
\label{eq:Qrs1:bbq-DS} \\
0 &\leq \Qr_1 \leq  5 \eps^{-1} {\bf 1}_{\Pdp(\s)}\,,
\label{eq:Qr1:bbq-DS} \\
\sabs{\Qb_2} & \leq 1500 \eps\,, 
\label{eq:Qrs2:bbq-DS} \\
|\Qr_2|  &\leq  1500  \,,
\label{eq:Qr2:bbq-DS} \\
|\Qr_\s|, |\Qc| &\leq  2 \cdot 250^2 \Q \eps^{-1} + \Cn \,.
\label{eq:Qrs:bbq-DS}
\end{align} 
\end{subequations} 
hold uniformly for all  $(x,\s) \in \Pdpm(\eps)$.  
\end{lemma}

We note that the lower bound in \eqref{eq:Qd:bbq-DS} matches that in \eqref{Qd-lower-upper}, up to an inconsequential $\OO(\eps)$ term. All other bounds precisely match those obtained previously in~\eqref{eq:Q:all:bbq}, including the all-important lower bound for $\Q$ in~\eqref{eq:Q:bbq-DS}. The bounds~\eqref{eq:Qrs1:bbq-DS} and~\eqref{eq:Qr1:bbq-DS} are new; their key feature is the $\geq 0$ lower bounds for both $\Qb_1$ and $\Qr_1$.

\begin{proof}[Proof of Lemma~\ref{lem:Q:bnds-DS}]  
For $x_1\leq x_1^*(x_2,\s)$, the functions $(\Qd,\Qb_1, \Qb_2, \Q, \Qc, \Qr_\s, \Qr_1,\Qr_2)$ defined in \eqref{QQQ-P} exactly equal the functions $(\Qd, 0, \Qb_2, \Q, \Qc, \Qr_\s, 0 ,\Qr_2)$ defined in~\eqref{QQQ}, because $\JJ = \mathcal{J}$ (recall the paragraph below~\eqref{eq:fake:Jg:def:ds}). As such, for these values of $x_1$, all the bounds established earlier in \eqref{eq:Q:all:bbq} continue to hold. Thus for $x_1\leq x_1^*(x_2,\s)$ no proof is required. Next, we consider the region $x_1> x_1^*(x_2,\s)$.
 
We start with the bound for $\Qd$ which is defined in \eqref{eq:QQQ:a-P}. In light of~\eqref{eq:mmJg:def} and~\eqref{eq:fake:Jg:def:ds}, we first consider values of $x_1$ such that $x_1^*(x_2,t) \le x_1 \le x_{1,-}^\sharp(x_2)$. Here, via~\eqref{eq:Jgb:identity:1}, the fact that $\mathfrak{C}\geq 0$, using \eqref{bs-Jg-1b}, and the fact that the map $x_1 \mapsto \nb_1 w_0(x_1,x_2)$ is monotone increasing (at fixed $x_2$) at least until it reaches the value of $-\frac 13$ (see the paragraph below \eqref{eq:x1:sharp:def}), we have that 
\begin{align}
\Qd(x,\s)
= - \eps (\p_t \JJ)(x,t) \big|_{t= \qds^{-1}(x,\s)} 
&
\geq - \eps (\p_t \Jg)(x,t)  \big|_{t= \qds^{-1}(x,\s)} 
\notag\\
&\geq - \tfrac{1+\alpha}{2} \big(\nb_1 w_0(x) + \eps \mathsf{C_{J_t}} \big)
\notag\\
&\geq - \tfrac{1+\alpha}{2} \big(\nb_1 w_0(x_{1}^\sharp(x_2),x_2) + \eps \mathsf{C_{J_t}} \big)
\notag\\
&\geq  \tfrac{1+\alpha}{2} \big( \tfrac{17}{20} - \eps \mathsf{C_{J_t}}  \big) 
\geq  \tfrac{17(1+\alpha)}{40} - \eps \tfrac{1+\alpha}{2} \mathsf{C_{J_t}}
 \,.
 \label{eq:Qd:lower:sec:12:a}
\end{align}
In identical fashion, for $x_1 \geq x_{1,+}^\sharp(x_2)$, by using the third branch in~\eqref{eq:mmJg:def} we have that 
\begin{equation}
\Qd(x,\s)
= - \eps (\p_t \JJ)(x,t) \big|_{t= \qds^{-1}(x,\s)} 
\geq - \tfrac{1+\alpha}{2} \big(\nb_1 w_0(x_1^\sharp(x_2),x_2) + \eps \mathsf{C_{J_t}} \big)
\geq  \tfrac{17(1+\alpha)}{40} - \eps \tfrac{1+\alpha}{2} \mathsf{C_{J_t}}
 \,.
 \label{eq:Qd:lower:sec:12:b}
\end{equation}
It thus remains to consider values of $x_1$ such that $x_{1,-}^\sharp(x_2) < x_1 < x_{1,+}^\sharp(x_2)$, which represents the middle branch of \eqref{eq:mmJg:def}. In this region, by construction we have that $\p_t \JJ = \p_t \mmJg$ is monotone increasing in $x_1$ (see~\eqref{eq:mmJg:def:middle}), and hence $\Qd$ is monotone decreasing in $x_1$. Its minimum value is thus attained when $x_1 = x_{1,+}^\sharp(x_2)$, a point at which the bound \eqref{eq:Qd:lower:sec:12:b} was previously established. 
This concludes the proof of the lower bound in~\eqref{eq:Qd:bbq-DS}. 
The upper bound in \eqref{eq:Qd:bbq-DS} is obtained by computing $-\eps  \p_t \mmJg$, which by its definition in~\eqref{eq:mmJg:def} and by~\eqref{eq:mmJg:def:middle} attains its maximum in the region $x_1^*(x_2,t) \leq  x_1 \leq x_{1,-}^\sharp(x_2)$. In this region, we conclude via the last bound in \eqref{eq:Jgb:identity:2} and the fact that  $(-\eps \p_t \Jg) \leq \tfrac{1+\alpha}{2} (1 + \eps \mathsf{C_{J_t}})$, via \eqref{bs-Jg-1b}.

The bound~\eqref{eq:Q:bbq-DS} follows from \eqref{eq:Qd:bbq-DS} and~\eqref{eq:Qrs2:bbq-DS}, since $|\Q - \Qd| \leq |V| |\Qb_2| \leq \Cn \eps^2$. 

Next, we turn to the bounds~\eqref{eq:Qrs1:bbq-DS} and~\eqref{eq:Qr1:bbq-DS} which concern $1$-derivatives of the geometric flattening coefficients,~$\Qb_1$ and $\Qr_1$. The definitions \eqref{eq:QQQ:aa-P} and \eqref{eq:fake:Jg:def:ds} yield 
\begin{subequations}
\label{eq:Qb1:DS}
\begin{equation}
\Qb_1(x,\s) = \eps (\p_1 \JJ)(x,t) = \eps (\p_1 \mmJg)(x,t)\,, \qquad \mbox{where} \qquad  t= \qds^{-1}(x,\s)\,,
\end{equation}
in the region of interest, $\{x_1 > x_1^*(x_2,t)\}$. We first note that the definition of $\mmJg$ in \eqref{eq:mmJg:def} yields 
\begin{equation}
(\p_1 \mmJg)(x,t) = 0\,, \qquad\mbox{for}\qquad x_1 > x_{1,+}^\sharp(x_2)\, .
\end{equation}
For the middle branch in \eqref{eq:mmJg:def}, namely $|x_1 - x_1^\sharp(x_2)| < \frac{\eps}{1000}$, we use the first two inequalities in~\eqref{eq:mmJg:def:middle}, which may be integrated in time (since $x_1^\sharp$ does not depend on time) to show that  
\begin{align}
\p_1 \mmJg(x,t) &= \underbrace{\p_1 \mmJg(x,\initial)}_{=0} + \int_{\initial}^t (\p_t \p_1 \mmJg)(x,t') {\rm d} t' \notag\\
&=  \int_{\initial}^t (\p_t \p_1 \mmJg)(x,t') {\rm d} t' \in [0, (\final-\initial) \cdot 2 (1+\alpha) \eps^{-2}] \subseteq [0, 5 \eps^{-1}].
\end{align}
It remains to consider the $\p_1$-derivative of the first branch in \eqref{eq:mmJg:def}, which concerns $x_1^* < x_1 < x_{1,-}^\sharp(x_2)$. For these values of $x_1$ we have that $\p_1 \mmJg = \p_1 \Jg$, while~\eqref{eq:why:the:fuck:not:0} and \eqref{bs-Jg-1}, imply 
\begin{equation}
\|\p_1 \Jg\|_{L^\infty_{x,t}} \leq (\final - \initial) \tfrac{1+\alpha}{2\eps^2} (\|\nb_1^2 w_0\|_{L^\infty_x} + \Cn \eps^2) \leq \tfrac{51}{25 \eps} + \Cn \eps \leq \tfrac{3}{\eps}
\,.
\end{equation}
This gives the desired upper bound. For the lower bound, as already noted in Proposition~\ref{prop:mJg:well:def}, for $x_1 \in [x_1^*(x_2,t),x_1^\sharp(x_2)]$ we have that $\nb_1 w_0(x) \leq - \frac{17}{20} < -\frac 13$, and therefore $\nb_1^2 w_0(x) \geq \eps^{\frac 78}$ (due to assumption~\eqref{item:ic:w0:d11:positive}). With this information,~\eqref{eq:why:the:fuck:not:0} and \eqref{bs-Jg-1}, now imply 
\begin{equation}
\p_1 \Jg (x,t) \geq (t - \initial) \tfrac{1+\alpha}{2\eps^2} ( \nb_1^2 w_0(x) -  \Cn \eps^2) 
\geq (t - \initial) \tfrac{1+\alpha}{2\eps^2} (\eps^{\frac 78} -  \Cn \eps^2) 
\geq (t - \initial) \tfrac{1+\alpha}{4} \eps^{-\frac 98} \geq 0
\,.
\end{equation}
\end{subequations}
Collecting the bounds in~\eqref{eq:Qb1:DS}, we have thus established~\eqref{eq:Qrs1:bbq-DS}.

Next, we consider the bound~\eqref{eq:Qr1:bbq-DS}. From \eqref{QQQ-P} and the chain rule we have that 
\begin{equation*}
\Qr_1(x,\s) =  \tfrac{\eps}{\Qd(x,\s)} (\p_t \p_1 \JJ)(x,\qds^{-1}(x,\s)) \,.
\end{equation*}
Inspecting the proofs of the bounds in~\eqref{eq:Qb1:DS}, and upon referring to~\eqref{bs-Jg-2a} instead of \eqref{bs-Jg-1}, we may deduce that $0 \leq \p_t \p_1 \JJ (x,t)\leq 2(1+\alpha)\eps^{-2}$, 
for all $x_1> x_1^*(x_2,t)$. Using the lower bound for $\Qd$ obtained in \eqref{eq:Qd:bbq-DS}, the proof of \eqref{eq:Qr1:bbq-DS} now follows. 

Next, we turn to the proof of~\eqref{eq:Qrs2:bbq-DS} for 
\begin{equation*}
\Qb_2(x,\s) = \tfrac{\eps}{\Qd(x,\s)} (\p_2 \JJ)(x,\qds^{-1}(x,\s)) =  \tfrac{\eps}{\Qd(x,\s)} (\p_2 \mmJg)(x,\qds^{-1}(x,\s))\,,
\end{equation*}
for $x_1 > x_1^*(x_2,t)$. For the first branch in \eqref{eq:mmJg:def}, we use ~\eqref{eq:why:the:fuck:not:0} and \eqref{bs-Jg-1}, to deduce $|(\p_2 \mmJg)(x,t)| = |(\nb_2 \Jg)(x,t)| \leq   (\final - \initial) \frac{1+\alpha}{2\eps} (\|\nb_1 \nb_2 w_0\|_{L^\infty_x} + \Cn \eps^2) \leq \frac{51}{25}   + \Cn \eps^2  \leq 3 $. For the middle branch in \eqref{eq:mmJg:def} we appeal to the second bound in \eqref{eq:mmJg:def:middle}, which may be integrated in time to yield $|(\p_2 \mmJg)(x,t)| \leq 250(1+\alpha)\eps^{-1} (\final-\initial) \leq 510$. For the last branch in \eqref{eq:mmJg:def}, we note that for all $x_1 \geq x_{1,+}^\sharp(x_2)$ by the argument in Remark~\ref{rem:mmJg:well:def} we have $|\p_2 \p_t \mmJg(x_1,x_2,t)| = |\p_2 \p_t \mmJg(x_{1,+}^\sharp(x_2),x_2,t)| \leq |\p_t \p_2 \Jg(x_1^\sharp(x_2),x_2,t)| + |\p_2 x_1^\sharp(x_2)| | \p_t \p_1 \Jg(x_1^\sharp(x_2),x_2,t)| \leq 242 (1+\alpha) \eps^{-1}$. Integrating this bound in time, similarly implies $|(\p_2 \mmJg)(x,t)| \leq 242 (1+\alpha) \eps^{-1} (\final-\initial) \leq 510$.   Using the lower bound for $\Qd$ obtained in \eqref{eq:Qd:bbq-DS}, the proof of \eqref{eq:Qrs2:bbq-DS} now follows. 

Inspecting the arguments in the previous paragraph, and refferring to the bound~\eqref{bs-Jg-2a} instead of~\eqref{bs-Jg-1}, we may deduce that $|\p_2 \p_t \mmJg(x,t)| \leq 250(1+\alpha)\eps^{-1}$ for all $x_1>x_1^*(x_2,t)$. Since $\Qr_2(x,\s) =  \tfrac{\eps}{\Qd(x,\s)} (\p_t \p_2 \JJ)(x,\qds^{-1}(x,\s))$, using the lower bound for $\Qd$ obtained in \eqref{eq:Qd:bbq-DS}, the proof of \eqref{eq:Qr2:bbq-DS} now follows.

Lastly, the estimates in~\eqref{eq:Qrs:bbq-DS} are obtained by repeating the proofs of \eqref{eq:Qrs:bbq},  and \eqref{eq:Qc:bbq}, except that these proofs simplify in downstream region. For instance, when compared to the $\Qr_\s$  expression  from \eqref{eq:Qring:2s:hat}, in the downstream region we have from \eqref{QQQ-P} and the chain rule that 
\begin{equation*}
\Qr_\s(x,\s) =  \tfrac{-\eps}{\Qd(x,\s)} (\p_t \p_t \JJ)(x,\qds^{-1}(x,\s)) \, .
\end{equation*}
When compared to  \eqref{eq:Qring:2s:hat}, the terms which involve an $\nb_1^2 \Jg$ denominator are absent. Given that we have already bounded $\Qd$ from below in \eqref{eq:Qd:bbq-DS}, 
repeating the same arguments as in the proofs of \eqref{eq:Qrs:bbq} we obtain the bound  for $\Qr_\s$ claimed in~\eqref{eq:Qr2:bbq-DS} and \eqref{eq:Qrs:bbq-DS}. The bound for $\Qc$ claimed in~\eqref{eq:Qrs:bbq-DS} follows from the $\Qr_\s$, $\Qr_2$, and $\Qb_2$ estimates, since $|\Qc - \Qr_\s| \leq |V| |\Qr_2| +|\p_\s V| |\Qb_2| \les \eps$. 
\end{proof} 

For the sake of completeness, we also record here an analogue of Lemma~\ref{lem:q:invertible}, showing that the map $\qds$ is invertible.

\begin{lemma}[\bf The map $\qds$ is invertible]
\label{lem:q:invertible-DS}
Assume that the bootstraps~\eqref{bootstraps-P} hold and that $\eps$ is sufficiently small with respect to $\alpha,\kappa_0$, and $\Cdata$. Then, the map $\qds$ defined by \eqref{eq:t-to-s-transform:all-P} is invertible, with inverse $\qds^{-1}$ defined by \eqref{s-to-t-transform-P}.
\end{lemma}
\begin{proof}[Proof of Lemma~\ref{lem:q:invertible-DS}] The proof is nearly identical to that of Lemma~\ref{lem:q:invertible}. The fact that the equation $ \s = \qds(x,t)= \eps  (1 - \JJ(x,t)  )$ has at most one solution $t\in [\initial,\final)$, for every fixed $\s \in [0,\eps)$ and $x\in \TT^2$, follows from the fact that $\p_t \qds = - \eps \p_t \JJ = \Qd \geq  \tfrac{2(1+\alpha)}{5} > 0$. The fact that the equation $ \s = \qds(x,t)= \eps  (1 - \JJ(x,t)  )$ has at least one solution $t\in (\initial,\final)$ for $\s \in (0,\eps)$, again follows from the intermediate value theorem because $\JJ \leq \Jgb$  in the downstream region (see the first paragraph in Remark~\ref{rem:oJJ:gg:0:maximal}). The case $\s=0$ is trivial, yielding $t=\initial$.
\end{proof}

\subsection{Identities in the downstream coordinate system}
With respect to the coordinates $(x,\s)$ given by \eqref{eq:t-to-s-transform:all-P}, with the transformation~\eqref{eq:f:tilde:f-P}, and upon dropping the tildes (see Remark~\ref{rem:no:tilde-P}), we have the following fundamental identities, which are translations of the identities in Section~\ref{sec:new:Euler:variables}  into $(x,\s)$ coordinates (see also~\eqref{Jg-evo-s}--\eqref{vort-s}):
\begin{subequations} 
 \label{fundamental-P}
\begin{align} 
(\Q\p_\s+V\p_2) \Jg &= \tfrac{1+ \alpha }{2} \Jg\Wbn + \tfrac{1- \alpha }{2} \Jg\Zbn \,, \label{Jg-evo-s-P} \\ 
(\Q\p_\s+V\p_2)  \JJ&= - \tfrac{\Q}{\eps} \,, \label{Jgb-evo-s-P} \\
(\Q\p_\s+V\p_2) \nbs_2 h &=  g  \bigl(\tfrac{1+ \alpha }{2} \Wbt + \tfrac{1- \alpha }{2} \Zbt \bigr) \,, \label{p2h-evo-s-P} \\
{\tfrac{1}{\eps}} \nbs_1 \Sigma &=  \tfrac{1}{2} \Jg(\Wbn -\Zbn)  + \tfrac{1}{2} \Jg \nbs_2 h (\Wbt -\Zbt) \,, \label{p1-Sigma-s-P} \\
\nbs_2 \Sigma &= \tfrac{1}{2} g^{\frac{1}{2}}  (\Wbt-\Zbt)\,,  \label{p2-Sigma-s-P} \\
(\Q\p_\s+V\p_2) \Sigma  &= -  \alpha \Sigma (\Zbn+ \Abt )  \,,  \label{Sigma0-ALE-s-P}  \\
(\Q \p_\s +    V \p_2) \Sigma^{-2\beta}   &=   2\alpha\beta \Sigma^{-2\beta}  (\Zbn+ \Abt )  \,,  \label{Sigma0i-ALE-s-P} \\
(\Q\p_\s+V\p_2)\nn  & =-  \bigl(\tfrac{1+ \alpha }{2} \Wbt + \tfrac{1- \alpha }{2} \Zbt \bigr) \tt   \,, \label{nn-evo-s-P} \\
(\Q\p_\s+V\p_2)\tt  & =   \bigl(\tfrac{1+ \alpha }{2} \Wbt + \tfrac{1- \alpha }{2} \Zbt \bigr) \nn   \,, \label{tt-evo-s-P} \\
\tfrac{\Jg}{\Sigma}  (\Q\p_\s+V\p_2)  \Upomega & =      {\tfrac{\alpha }{\eps}}  \nbs_1 \Upomega
- \alpha \Jg  g^{- {\frac{1}{2}} } \nbs_2h \   \nbs_2 \Upomega   \,. \label{vort-s-P}
\end{align} 
\end{subequations} 
By using \eqref{eq:Jgb:identity:2}, we also note that $\nbs_1 (\Jg-\Jgb) = \nbs_2 (\Jg-\Jgb)=0$.

\subsection{Bounds for the geometry, sound speed, and ALE velocity for downstream development}
  
We record here the modifications to the bounds obtained earlier in Section~\ref{sec:geometry:sound:ALE}, due to the downstream geometry. It turns out that the only modification is due to the change of the weight in the energy estimates, namely $\mathcal{J} \mapsto \JJ$: all the bounds which do not involve this weight function remain identical to those in Section~\ref{sec:geometry:sound:ALE} (and for those bounds we abuse notation and make reference to equation numbers from Section~\ref{sec:geometry:sound:ALE}), while all the bounds which do involve this weight function need to be modified by exchanging $\mathcal{J}$ with $\JJ$. For instance, the bounds in Proposition~\ref{prop:geometry} now become the bounds given in Proposition~\ref{prop:geometry-P} below. The corollaries and remarks which follow this proposition (in particular, the closure of the \eqref{bootstraps-Dnorm:h2} and~\eqref{bootstraps-Dnorm:Jg} bootstraps in Corollary~\ref{cor:Bj:Bh}) remain the same as in Section~\ref{sec:geometry:sound:ALE}, and to avoid redundancy we do not repeat those arguments here.
 
\begin{proposition}[Bounds for the geometry, sound speed, and ALE velocity]
\label{prop:geometry-P}
Assume that  the bootstrap assumptions \eqref{bootstraps-P} hold, and that $\eps$ is taken to be sufficiently small to ensure $\eps^{\frac 12}  ( \brak{\mathsf{B_J}} +  \brak{\mathsf{B_h}} +     \brak{\mathsf{B_6}} )  \leq 1$.  
Then, assuming $\eps$ is sufficiently small with respect to $\alpha,\kappa_0,$ and $\Cdata$, 
the functions $(\Jg, \nbs_1 h, \nbs_2 h, \Sigma, V)$ satisfy the bounds \eqref{D5JgEnergy}, \eqref{D6h1Energy:new},  
\eqref{D5h2Energy},  \eqref{eq:Sigma:H6:new}--\eqref{eq:Sigma:H6:new:bdd}, and \eqref{eq:V:H6:new}--\eqref{eq:V:H6:new:bdd} respectively. Additionally

\begin{subequations}
\label{geometry-bounds-new-P}
\begin{align} 
\sup_{\s \in [0,\eps]}  \snorm{ \JJof \nbs^6 \Jg(\cdot,\s) }^2_{L^\infty_\s L^2_{x}} 
+ \tfrac{1}{\eps}  \snorm{  \JJmof \nbs^6 \Jg }_{L^2_{x,\s}}^2 
&\les  \eps  \brak{\mathsf{B}_6}^2    \,,   \label{D6JgEnergy:new-P}
\\
\sup_{\s \in [0,\eps]} \| \JJof \nbs^{6} \nbs_2h(\cdot, \s) \|^2_{L^\infty_\s L^2_{x}} 
+ \tfrac{1}{\eps}  \|\JJmof \nbs^{6}  \nbs_2h \|_{L^2_{x,\s}}^2
&\les 
\mathsf{K}^2 \eps^3 \brak{\mathsf{B}_6}^2
\label{D6h2Energy:new-P}
\\
 \sup_{\s \in [0,\eps]}\sum_{ |\gamma|=3}^6 
 \|\JJmof\bigl(\nbs^{|\gamma|}  \nn + \tfrac{1}{g}   \tt \nbs^{|\gamma|} \nbs_2 h\bigr) \|_{L^2_{x,\s}}
+
\|\JJmof \bigl( \nbs^{|\gamma|}  \tt - \tfrac{1}{g}  \nn \nbs^{|\gamma|} \nbs_2 h \bigr) \|_{L^2_{x,\s}}
&\les \mathsf{K} \eps^3 \brak{\mathsf{B}_6}  \,,
\label{D6n-bound:a:new-P}
\\
 \sup_{\s \in [0,\eps]}\| \JJmof  \nbs^6  \nn  \|_{L^2_{x,\s}}
+
\|\JJmof  \nbs^6 \tt \|_{L^2_{x,\s}}
&\les \mathsf{K} \eps^2 \brak{\mathsf{B}_6}   \,,
\label{D6n-bound:b:new-P}
\\
 \sup_{\s \in [0,\eps]} \big( \|\JJof \nbs^6  \nn  \|_{L^\infty_\s L^2_{x}} 
+
\|\JJof  \nbs^6 \tt \|_{L^\infty_\s L^2_{x}} \big)
&\les \mathsf{K} \eps^{\frac 32} \brak{\mathsf{B}_6}   \,,
\label{D6n-bound:b:new:bdd-P}
\end{align} 
where the implicit constants in all the above inequalities depend only on $\alpha$, $\kappa_0$, and $\Cdata$. 
\end{subequations}
\end{proposition}
\begin{proof}[Proof of Proposition~\ref{prop:geometry-P}]
We explain the downstream modifications required for the proof of the inequality  \eqref{D6JgEnergy:new-P}.   Just as the equation
\eqref{D5-Jg-s} was obtained,  letting $\nbs^6$ act on \eqref{Jg-evo-s-P} in the set $\Pdpm(\eps)$, we have that
\begin{align} 
(\Q\p_\s +V\p_2)(\nbs^6\Jg) = \tfrac{1+\alpha}{2} \nbs^6(\Jg \Wbn)+ \tfrac{1-\alpha}{2}\nbs^6 (\Jg\Zbn) + \Rj\,, \label{D5-Jg-s-P}
\end{align} 
where $\Rj  = -\nbs^6 V \, \nbs_2 \Jg - \doublecom{\nbs^6, V, \nbs_2 \Jg}$.  For each $\s \in (0,\eps)$, 
we compute the $L^2(\Xdpm(\s))$-inner product of \eqref{D5-Jg-s-P} with $\JJh \nbs^6\Jg$ to obtain that 
\begin{align}
&\tfrac{1}{2} \int_{\Xdpm(\s)} \JJh (\Q\p_\s +V\p_2) |\nbs^6\Jg|^2 
 \notag \\
& \qquad \qquad
= \tfrac{1+\alpha}{2}\int_{\Xdpm(\s)} \JJh \nbs^6(\Jg \Wbn)\nbs^6\Jg
+ \tfrac{1-\alpha}{2}\int_{\Xdpm(\s)} \JJh  \nbs^6(\Jg\Zbn)\nbs^6\Jg 
+ \int_{\Xdpm(\s)} \JJh   \Rj \ \nbs^6\Jg \,.
\label{eq:vomit:cascade:4}
\end{align} 
Next, we commute the operator $ (\Q\p_\s+V\p_2) $ past $\JJh$ separately in the regions $\Xdm(\s)$ and $\Xdp(\s)$; here, recall the definition of $\JJ$ in  \eqref{eq:fake:Jg:def:ds}. 
In analogy to \eqref{eq:diff:by:parts}, in both $\Xdm(\s)$ and $\Xdp(\s)$ we have that 
for any function $f = f(x,\s) \geq 0 $ and any $r\in \mathbb{R}$ we have the following identity and subsequent estimate (by appealing to \eqref{Jgb-evo-s-P})
\begin{align}
 \JJ^r  (\Q\p_\s+V\p_2)  f 
 &=
 \p_\s \bigl( \JJ^r \Q f\bigr) + \p_2 \bigl(\JJ^r V f\bigr)
 -  r f \JJ^{r-1} (\Q \p_\s + V \p_2) \JJ
 - f \JJ^{r} \bigl(\Qc + \p_2 V\bigr)
 \notag\\
 &= 
 \p_\s \bigl( \JJ^r \Q f\bigr) + \p_2 \bigl(\JJ^r V f\bigr)
 +r  \tfrac{\Q}{\eps}   f \JJ^{r-1}   
 - f \JJ^{r} \bigl(\Qc + \p_2 V  \bigr)
 \,,
 \label{eq:diff:by:parts-P}
\end{align}
Integrating the above expression over $\Xdpm(\s)$,  appealing to \eqref{ds-exchange} for the $\p_\s$ term, to \eqref{eq:vomit:cascade:1} for the $\p_2$ term, appealing to the $\Q$ and $\Qc$ bounds in \eqref{eq:Q:all:bbq-DS},  to the $V,_2 = \OO(\eps)$ bootstrap in \eqref{bootstraps-P}, and upon taking $\eps$ to be sufficiently small, we arrive at the bound 
\begin{align}
\int_{\Xdpm(\s)}  \JJ^r  (\Q\p_\s+V\p_2)  f 
\geq 
\tfrac{d}{d\s} \int_{\Xdpm(\s)}  \JJ^r  \Q f 
+ 
\tfrac{2 r (1+\alpha)}{5\eps} \int_{\Xdpm(\s)}  \JJ^{r-1} f 
-
\tfrac{2 \cdot 250^2 + \Cn \eps}{\eps}\int_{\Xdpm(\s)}  \JJ^{r} \Q f 
\,.
\label{eq:vomit:cascade:3} 
\end{align}
With $r=\frac 12$ and $f= \frac 12 |\nbs^6 \Jg|^2$, we use \eqref{eq:vomit:cascade:3} to lower bound the left side of \eqref{eq:vomit:cascade:4}, resulting in  
\begin{align*}
&\tfrac{1}{2 d\s} \int_{\Xdpm(\s)} \JJh  \Q  |\nbs^6\Jg|^2 
+ \tfrac{1+\alpha}{10\eps} \int_{\Xdm(\s)}  \JJ^{-\frac 12}   |\nbs^6\Jg|^2 
- \tfrac{1+250^2}{\eps}  \int_{\Xdpm(\s)}\JJh  \Q  |\nbs^6\Jg|^2 \\
&\leq \tfrac{1+\alpha}{2}\int_{\Xdpm(\s)} \JJh \nbs^6(\Jg \Wbn)\nbs^6\Jg
+ \tfrac{1-\alpha}{2}\int_{\Xdpm(\s)} \JJh  \nbs^6(\Jg\Zbn)\nbs^6\Jg 
+ \int_{\Xdpm(\s)} \JJh   \Rj \ \nbs^6\Jg \,.
\end{align*} 
The integrals on the right side of the above inequality are analyzed in the identical manner as in the proof of Lemma~\ref{lem:D5-Jg} (see~\eqref{D6JgEnergy:temp1}--\eqref{D6JgEnergy:temp3}). By using Gr\"onwall's inequality on the interval $[0,\eps]$ to handle with the third term on the right side of the above inequality, and using the lower bound in \eqref{eq:Q:bbq-DS}, we arrive at
\begin{align*} 
\sup_{\s\in[0,\eps]}  \|\JJof \nbs^6\Jg(\cdot , \s)\|^2_{L^2} 
+ \tfrac{1 }{\eps}  \|\JJmof  \nbs^6\Jg \|^2_{L^2_{x,\s}} 
\les  
\| \nbs^6\Jg(\cdot , 0)\|^2_{L^2_x}
+ \eps \brak{\mathsf{B}_6}^2  \,.
\end{align*} 
Using \eqref{table:derivatives},  we arrive at  \eqref{D6JgEnergy:new-P}.

The downstream modifications (when compared to the proof of Proposition~\ref{prop:geometry}) required for proof of the inequalities \eqref{D6h2Energy:new-P}--\eqref{D6n-bound:b:new:bdd-P}  are identical (going through~\eqref{eq:vomit:cascade:3} instead of \eqref{eq:diff:by:parts}) and these details will be omitted. To avoid redundancy, we also omit the proofs of the unweighted bounds for $(\Jg, \nbs_1 h, \nbs_2 h, \Sigma, V)$, as these are established exactly as in the proof of Proposition~\ref{prop:geometry}.
\end{proof}

\subsection{Estimates in the downstream geometry which are improved due to vorticity bounds}
In the downstream geometry given by the transformation \eqref{eq:t-to-s-transform:all-P} we still obtain an improved estimates for the vorticity, and as a consequence we have that: 
\begin{align} 
(\Jg\Wbn,\Zbn,\Abn) \text{ satisfy the improved bounds  
\eqref{eq:Jg:Abn:D5:improve},   \eqref{eq:Jg:Zbn:D5:improve}, \eqref{eq:Jg:Wbn:improve:material} and \eqref{Wbn:improved} } \,, \text{ with } \mathcal{J} \text{ replaced by } \JJ \,.
\label{improved-P}
\end{align} 
As shown in Section~\ref{sec:vorticity}, these improved bounds are a consequence of the bootstrap assumptions (these are causing the need to exchange any power of $\mathcal{J}$ for the same power of $\JJ$, which in light of \eqref{eq:fake:Jg} and \eqref{eq:D:JJ} are both equal to the same power of $(1-\frac{\s}{\eps})$) and an energetic bound for the vorticity, with optimal weights (independent of $\JJ$). This vorticity bound was obtained previously in Proposition~\ref{prop:vort:H6}; it remains unchanged in the downstream geometry, though the proof needs to be slightly altered due to the presence of a $\p_1$ derivative (see~\eqref{vort-s}) in the specific vorticity transport (see Section~\ref{app:downstream:Lp}). We record these vorticity estimates next:

\begin{proposition}[$H^6$ estimates for the vorticity]
\label{prop:vort:H6-DS}
Assume that the bootstrap assumptions~\eqref{bootstraps-P} hold and that $\eps$ is taken to be sufficiently small with respect to $\alpha,\kappa_0$ and $\Cdata$.  
The ALE vorticity $\Omega$ satisfies the $H^6$ energy bound \eqref{eq:vort:H6}, the $L^\infty$ bound \eqref{eq:vorticity:pointwise:a} and the $W^{1,\infty}$ bound \eqref{eq:vorticity:pointwise:b}. 
\end{proposition}
\begin{proof}[Proof of Proposition~\ref{prop:vort:H6-DS}]
The claimed vorticity estimates are established identically to the proof of Proposition~\ref{prop:vort:H6}, and are a consequence of the following bound for the specific vorticity
\begin{align}
\sup_{\s\in[0,\eps]} \snorm{\Jgh  \nbs^6\Upomega (\cdot,\s)}_{L^2_x}^2 
+
\tfrac{1}{ \eps}  \int_0^{\eps}\snorm{ \nbs^6\Upomega (\cdot,\s)}_{L^2_x}^2 {\rm d}\s
\les \eps \brak{\mathsf{B_6}}^2  \,.
\label{eq:svort:H6-DS}
\end{align}
To avoid redundancy we only sketch here the modifications (when compared to the proof of~\eqref{eq:svort:H6}) needed to establish \eqref{eq:svort:H6-DS}, which arise to the $\nbs_1$-derivative. Applying $\nbs^6$ to \eqref{vort-s-P}, as in~\eqref{D6-vort-s} we arrive at 
\begin{align} 
&\tfrac{\Jg}{\Sigma} (\Q \p_\s + V \p_2)  \nbs^6  \Upomega -  \tfrac{\alpha}{\eps} \nbs_1 \nbs^6  \Upomega + \alpha \Jg   g^{-\frac 12} \nbs_2 h  \nbs_2 \nbs^6  \Upomega = \RR_\Upomega
\,, \label{D6-vort-s-DS}
\end{align} 
where the remainder term $\RR_\Upomega$ is defined exactly as in~\eqref{D4-vort-s-remainder}. We test~\eqref{D6-vort-s-DS} with $\Sigma^{-2\beta+1} \nbs^6 \Upomega$, with $\beta>0$ to be chosen appropriately (in terms of $\alpha,\kappa_0,\Cdata$), and integrate over $\Pdpm(\s)$. Note that here we do not multiply by powers of $\JJ$ (as opposed to~\eqref{eq:vomit:cascade:4} earlier), and so appealing to \eqref{eq:vomit:cascade:3} is not necessary. Additionally, note that the term containing a $\p_1$ derivative, namely $-\frac{\alpha}{2} \Sigma^{-2\beta+1} \nbs_1 |\nbs^6 \Upomega|^2$ does not contain a $\Jg$ or $\JJ$ term.  Instead, we appeal directly to the adjoint formulas for $(\Q\p_\s + V\p_2)^*$, $\nbs_1^*$, and $\nbs_2^*$ from~\eqref{eq:adjoints-P}, and to the bounds~\eqref{eq:Q:all:bbq-DS} to deduce in analogy with~\eqref{eq:vorticity:energy:0}--\eqref{eq:vorticity:energy:00} that 
\begin{align}
&  \snorm{\tfrac{(\Jg \Q)^{\frac 12}}{ \Sigma^{\beta}}  \nbs^6\Upomega (\cdot,\s)}_{L^2_x}^2 
-  \snorm{\tfrac{(\Jg \Q)^{\frac 12}}{ \Sigma^{\beta}}  \nbs^6\Upomega  (\cdot,0)}_{L^2_x}^2 
+ \tint \sabs{\nbs^6\Upomega}^2 \tfrac{1}{\Sigma^{2\beta}}  \mathsf{G}_{\Upomega} 
+ \alpha  \int \underbrace{\tfrac{1}{\Sigma^{2\beta}} \sabs{\nbs^6\Upomega}^2 \Sigma  \Qb_1}_{\geq 0 \mbox{ due to } \eqref{eq:Qrs1:bbq-DS}}  \Big|_{\s}
\notag\\
&=    \alpha  \int \tfrac{\Jg}{\Sigma^{2\beta}} \sabs{\nbs^6\Upomega}^2\Sigma \Qb_2    g^{-\frac 12} \nbs_2 h  \Big|_{\s}
+  \tint\tfrac{2}{\Sigma^{2\beta-1}} \RR_\Upomega \nbs^6 \Upomega 
\,.
\label{eq:vorticity:DS:1}
\end{align}
where 
\begin{align}
\mathsf{G}_{\Upomega}
&:=  
-  \Bigl( \bigl( \tfrac{1+\alpha}{2} \Jg \Wbn + \tfrac{1-\alpha}{2} \Jg \Zbn\bigr) + 2 \alpha \beta  \Jg   \bigl( \Zbn + \Abt\bigr)\Bigr)
- \alpha (\beta-\tfrac 12) \bigl(  \Jg \Wbn -  \Jg \Zbn + \Jg \nbs_2 h (\Wbt - \Zbt) \bigr)
\notag\\
&\qquad 
- \Jg \bigl(\Qr_\s -  V \Qr_2 + \nbs_2 V  \bigr) 
-  \alpha \Sigma \Qr_1 
+ (2\beta-1) \alpha \Jg  g^{-\frac 12} \nbs_2 h \nbs_2 \Sigma 
+ \alpha  \Qr_2 \Sigma \Jg  g^{-\frac 12} \nbs_2 h 
\,.
\label{eq:vorticity:DS:2}
\end{align}
The remainder term $\RR_\Upomega$ present on the right side of \eqref{eq:vorticity:DS:2} satisfies the same $\Sigma^{-\beta}$-weighted $L^2_{x,\s}$ bound as in \eqref{eq:vorticity:remainder:4}.  
When compared to~\eqref{eq:vorticity:energy:0}--\eqref{eq:vorticity:energy:00}, the new contributions are those due to $\Qb_1$ for the temporal boundary term in~\eqref{eq:vorticity:DS:1}, and $\Qr_1$ in the definition of $\mathsf{G}_{\Upomega}$~\eqref{eq:vorticity:DS:2}.  According to \eqref{eq:Qrs1:bbq-DS} we have that $\Qb_1 \geq 0$, and thus the additional temporal term in~\eqref{eq:vorticity:DS:1} has a favorable sign, allowing us to ignore it. Also, according to~\eqref{eq:Qr1:bbq-DS} we have that $0\leq \Qr_1 \leq 5 \eps^{-1}$, so that similarly to \eqref{justbecause0} we may deduce 
\begin{align*}
\mathsf{G}_{\Upomega} 
\geq  (\alpha \beta + \tfrac 12) \bigl( \tfrac{9}{10\eps} - \tfrac{13}{\eps} \Jg \bigr) 
- \tfrac{2 \cdot 250^2}{\eps} \Q \Jg  - \Cn \brak{\beta}
- \tfrac{5 \alpha \kappa_0}{\eps} \,.
\end{align*}
Due to the last term appearing in the above lower bound on $\mathsf{G}_{\Upomega}$, the size of $\beta$ (cf.~\eqref{eq:vorticity:energy:beta:def}) needs to be increased by an additive factor of $6 \kappa_0$. Thus, for $\beta = \beta(\alpha,\kappa_0,\Cdata)$ chosen suitably, and for $\eps$ sufficiently small (in terms of $\alpha, \kappa_0,\Cdata$), in analogy to \eqref{eq:vorticity:energy:3}, we may deduce from \eqref{eq:vorticity:DS:1} that 
\begin{align*}
&\tfrac 12  \snorm{\tfrac{(\Jg \Q)^{\frac 12}}{ \Sigma^{\beta}}  \nbs^6\Upomega (\cdot,\s)}_{L^2_x}^2 
+ \tfrac{1}{8\eps}  
\int_0^{\s} \snorm{\tfrac{1}{\Sigma^\beta} \nbs^6\Upomega (\cdot,\s')}_{L^2_x}^2 {\rm d}\s'   
\notag\\
&\leq    \snorm{\tfrac{(\Jg \Q)^{\frac 12}}{ \Sigma^{\beta}}  \nbs^6\Upomega  (\cdot,0)}_{L^2_x}^2 
+ \tfrac{\Cn}{\eps}  
\int_0^{\s} \snorm{\tfrac{(\Jg \Q)^{\frac 12}}{ \Sigma^{\beta}}  \nbs^6\Upomega (\cdot,\s')}_{L^2_x}^2  {\rm d}\s'
+ \Cn  \eps (4^5 \kappa_0^{-1})^{2\beta} \brak{\mathsf{B_6}}^2      
\,.
\end{align*}
Using Gr\"onwall for $\s \in [0,\eps]$, multiplying the resulting estimate by $\kappa_0^{2\beta}$, and appealing to \eqref{bs-Sigma} and \eqref{eq:Q:bbq-DS}, concludes the proof of \eqref{eq:svort:H6-DS}.
\end{proof}

\subsection{Closing the pointwise bootstraps in the downstream geometry}
 
The pointwise bounds which were previously established in Section~\ref{sec:pointwise:bootstraps}, only relied on the evolution equations in $(x,\s)$ coordinates (cf.~\eqref{Jg-evo-s}--\eqref{nn-tt-evo-s}), the bootstrap assumptions, the functional analytic setup in Appendix~\ref{app:functional}, and of the $L^\infty$ estimates from Appendix~\ref{sec:app:transport}. 
These arguments apply {\bf as is} in the geometry of the downstream development, except that we refer to the evolution identities~\eqref{fundamental-P}, the bootstraps~\eqref{bootstraps-P}, to  the functional analytic bounds from Section~\ref{rem:app:downstream:flat}, and to the $L^\infty$ estimates from Section~\ref{app:downstream:Lp}. We omit these redundant details concerning the closure of the pointwise bootstraps.

\subsection{Downstream energy estimates}
At this stage in the proof, it only remains to close the bootstrap~\eqref{boots-P} (see~\eqref{bootstraps-Dnorm:6}) for the sixth order energy $\widetilde{\mathcal{E}}_6$ and damping $\widetilde{\mathcal{D}}_6$ norms, defined earlier in~\eqref{eq:tilde:E5E6:N+T} and~\eqref{eq:tilde:D5D6:N+T}. 

Previously, this was achieved by separately establishing a bound for the tangential parts of the energy $\widetilde{\mathcal{E}}_{6,\ttt}$ and damping  $\widetilde{\mathcal{D}}_{6,\ttt}$ in~Section~\ref{sec:sixth:order:energy-tangential}, and the normal parts of the energy $\widetilde{\mathcal{E}}_{6,\nnn}$ and damping  $\widetilde{\mathcal{D}}_{6,\nnn}$ in~Section~\ref{sec:sixth:order:energy}. In turn, these estimates required that we established improved energy bounds for six ``pure time derivatives'' in Section~\ref{sec:pure:time}. 

For the downstream \MGHDB, we follow the same exact strategy. 
As before, the tangential bounds from~Section~\ref{sec:sixth:order:energy-tangential} and normal energy estimates from Section~\ref{sec:sixth:order:energy} run in parallel, the only difference being that the fundamental variables are un-weighted for the tangential part (i.e.~$(\Wbt,\Zbt,\Abt)$) and are $\Jg$-weighted for the normal part (i.e.~$(\Jg \Wbn,\Jg \Zbn, \Jg \Abn)$). The special estimates for six pure time derivatives from Section~\ref{sec:pure:time} are used in the same way, to treat the remainders $\mathcal{R} _\Zb^\tt$ and $\mathcal{R} _\Zb^\nn$, in the $\Zbt$, and respectively $\Jg \Zbn$ equations. 

Since the tangential and normal energy estimates run in parallel (similarities and differences may be seen by comparing Sections~\ref{sec:sixth:order:energy-tangential} and~\ref{sec:sixth:order:energy}), we do not repeat both of these two sets of energy estimates for the downstream geometry. Indeed, the downstream modifications to the tangential component energy estimates are identical to the modifications made to the normal component energy estimates. For the sake of brevity, we chose to only give details for the downstream modifications to the normal energy estimates  (see Section~\ref{eq:downstream:normal} below). 

\subsubsection{Sixth order tangential energy estimates}
For the tangential energy estimates, at this point we simply record that by repeating the arguments from Section~\ref{sec:sixth:order:energy-tangential}, with the modifications outlined in Section~\ref{eq:downstream:normal} below (see the argument leading to~\eqref{eq:normal:conclusion:1-P}--\eqref{eq:normal:conclusion:3-P}), 
similarly to \eqref{eq:hate:11:a}--\eqref{eq:hate:12}, we obtain that there exists a constant 
\begin{equation*}
 \hat{\mathsf{c}}_{\alpha,\kappa_0} > 0 \,,
\end{equation*}
which depends only on $\alpha$ and also on $\kappa_0$, and may be computed explicitly, such that 
\begin{align}
& 
\sup_{\s \in [0,\eps]}
\snorm{ \JJ^{\frac 34}  \Jg^{\!\frac 12} \nbs^6(\Wbt,\Zbt,\Abt)(\cdot,\s)}_{L^2_x}^2 
 + \tfrac{1}{\eps}   \int_0^{\eps}  
\snorm{\JJ^{\frac 14} \Jg^{\!\frac 12} \nbs^6 (\Wbt,\Zbt,\Abt)(\cdot,\s)}_{L^2_x}^2   {\rm d} \s
 \notag\\
 &\qquad 
+ \tfrac{1}{\eps^2} \sup_{\s \in [0,\eps]} 
 \snorm{ \JJ^{\frac 14} \nn \cdot  \nbs^6 \tt (\cdot,\s)}_{L^2_{x}}^2
+ \tfrac{1}{\eps^3} \int_0^{\eps} \snorm{ \JJ^{-\frac 14} \nn \cdot  \nbs^6 \tt (\cdot,\s)}_{L^2_{x}}^2 {\rm d} \s
\notag\\
 &  
 \leq \hat{\mathsf{c}}_{\alpha,\kappa_0} \eps  \Big( \Cdatatwo + \mathsf{B}_6^2 + \Cn \eps^{\frac 12}  \mathsf{K}^2  \brak{\mathsf{B}_6}^2\Bigr)
 \,.
 \label{eq:hate:12-P}
\end{align}
Then, as in~\eqref{eq:hate:13:aa}--\eqref{eq:hate:13:a}, upon ensuring that 
\begin{align}
\mathsf{B}_6 &\geq \max\{1, \Cdata \} \,,
\label{eq:B6:choice:1-P}
\end{align}
and upon defining
\begin{align}
\mathsf{K} &:= 8 \max\{1, \hat{\mathsf{c}}_{\alpha,\kappa_0}^{\frac 12} \} \,,
\label{eq:K:choice:1-P}
\end{align}
by letting $\eps$ be sufficiently small in terms of $\alpha,\kappa_0$, and $\Cdata$, we deduce from~\eqref{eq:hate:12-P} that 
\begin{align}
\eps \sup_{\s \in [0,\eps]} \widetilde{\mathcal{E}}_{6,\ttt}^2(\s)
+\widetilde{\mathcal{D}}_{6,\ttt}^2(\eps) 
\leq \tfrac{1}{8} (\eps \mathsf{K})^2 \mathsf{B}_6^2 
 \,.
 \label{eq:hate:13-P}
\end{align}
This bound is the same as~\eqref{eq:hate:13:a}. It closes the ``tangential part'' of the remaining bootstrap~\eqref{boots-P} for $\widetilde{\mathcal{E}}_{6}$ and $\widetilde{\mathcal{D}}_{6}$.

\subsubsection{Sixth order pure-time energy estimates}
For the energy estimates concerning pure time derivatives,  we record that by repeating the arguments from Section~\ref{sec:pure:time}, with the modifications outlined in Section~\ref{eq:downstream:normal} below, the same bound as given in \eqref{eq:madman:2} holds, namely 
\begin{equation} 
\eps^{\frac 12} 
\snorm{\JJ^{\! \frac 34}\Jgh \nbs_{\s}^6  \Zbn }_{L^\infty_x L^2_{x}}
+ 
\snorm{\JJ^{\! \frac 34} \nbs_{\s}^6  \Zbn }_{L^2_{x,\s}} 
\leq
\eps^{\frac 12} 
\snorm{\JJ^{\! \frac 34}\Jgh \nbs_{\s}^6  \Zbn }_{L^\infty_x L^2_{x}}
+
\snorm{\JJ^{\! \frac 14} \Jgh \nbs_{\s}^6  \Zbn }_{L^2_{x,\s}} 
\les \eps \mathsf{K} \brak{\mathsf{B}_6} 
\,.
\label{eq:madman:2-P}
\end{equation}


\subsection{Downstream energy estimates for normal components}
\label{eq:downstream:normal}
It thus remains to outline the modifications to the normal energy estimates in Section~\ref{sec:sixth:order:energy}, required for the downstream geometry. 
We continue to use the equation set \eqref{energy-WZA-s} for $(\Jg\Wbn, \Jg\Zbn, \Jg \Abn)$ in $(x,\s)$ coordinates, with the operator $\nbs^6$  used only in the set $\Pdpm(\eps)$, so that differentiation does not take place across the surface 
$\Gamma(\eps)$. According to definitions~\eqref{eq:tilde:E5E6} and~\eqref{eq:tilde:D5D6}, the energy identity \eqref{D6-L2-N} is replaced with the downstream energy identity
\begin{equation} 
\tint \jb \JJss \Big( \underbrace{\eqref{energy-Wn-s} \ \Jg  \nbs^6(\Jg\Wbn)} _{ I^{\WW_n}}
+ \underbrace{\eqref{energy-Zn-s} \  \nbs^6(\Jg\Zbn)} _{ I^{\ZZ_n}}
+ \underbrace{2\, \eqref{energy-An-s} \ \nbs^6(  \Jg \Abn) } _{ I^{\AA_n}}
 \Big)  {\rm d}x {\rm d} \s'=0 \,, \label{D5-L2-P}
\end{equation} 
where once again $\jb = \Sigma^{-2\beta+1}$ and $\beta = \beta(\alpha,\kappa_0)>0$ is a sufficiently large constant chosen in the proof (see~\eqref{eq:normal:bounds:beta-P}). Here and throughout the remainder of the section we use the integral notation
\begin{equation}
\int \mbox{ to denote } \int_{\Xdpm(\s)} {\rm d} x \,,
\qquad \mbox{and} \qquad 
\int_0^{\s} \!\!\! \int \mbox{ to denote } \int_0^{\s} \!\!\! \int_{\Xdpm(\s')} {\rm d}x {\rm d}\s' =\int_{\Pdpm(\s)} {\rm d}x {\rm d}\s'\,.
\label{eq:lets:confuse:the:reader}
\end{equation}
With this notation, we frequently appeal to~\eqref{eq:vomit:cascade:1},~\eqref{ds-exchange}, and to the adjoint formulae~\eqref{eq:adjoints-P}.

\subsubsection{The additive decompositions of integrals  $I^{\WW_n}$, $I^{\ZZ_n}$, and $I^{\AA_n}$} 
In analogy to~\eqref{Integral-Wbn}, we additively decompose the integral  $I^{\WW_n}$ as 
 \begin{subequations} 
 \label{Integral-Wbn-P}
\begin{align}
I^{\WW_n}&= I^{\WW_n}_1+I^{\WW_n}_3+I^{\WW_n}_4 +I^{\WW_n}_5+I^{\WW_n}_6
\notag \,, \\
 I^{\WW_n}_1 &=
\tint \tfrac{1}{\Sigma^{2\beta}} \JJss \Jg (\Q\p_\s +V\p_2)\nbs^6( \Jg\Wbn  ) \  \nbs^6(\Jg\Wbn)
\,, \label{I1-Wbn-P} \\
I^{\WW_n}_3 &=
 \alpha \tint \jb \JJss \Jg g^{- {\frac{1}{2}} }  \nbs_2\nbs^6(\Jg\Abn)  \  \nbs^6(\Jg\Wbn)
\,, \label{I3-Wbn-P} \\
 I^{\WW_n}_4 &=
-\alpha \tint \jb  \JJss \Jg g^{- {\frac{1}{2}} }  \Abn \nbs_2\nbs^6 \Jg \  \nbs^6(\Jg\Wbn)
  \,,   \label{I4-Wbn-P} \\
 I^{\WW_n}_5 &=
-\tfrac{\alpha }{2}  \tint \jb  \JJss  \Jg g^{- {\frac{1}{2}} } (\Jg\Wbn + \Jg\Zbn - 2\Jg\Abt)\nbs_2 \nbs^6 \tt\cdo\nn \ \nbs^6(\Jg\Wbn)
 \,, \label{I5-Wbn-P} \\
 I^{\WW_n}_6 &=
-\tfrac{\alpha }{2}  \tint \jb \JJss \Jg \big(  \nbs^6 \Fwn + \mathcal{R}^\nn_{\Wb} + \mathcal{C}^\nn_{\Wb}\big) \ \nbs^6(\Jg\Wbn)
 \,, \label{I6-Wbn-P}
\end{align} 
\end{subequations} 
where we have used the notation in~\eqref{eq:lets:confuse:the:reader}.
Next, in analogy to~\eqref{Integral-Zbn}, we additively decompose the integral  $I^{\ZZ_n}$ as 
 \begin{subequations} 
 \label{Integral-Zbn-P}
 \begin{align}
  I^{\ZZ_n} & = I^{\ZZ_n}_1 + I^{\ZZ_n}_2 + I^{\ZZ_n}_3+ I^{\ZZ_n}_4+ I^{\ZZ_n}_5+ I^{\ZZ_n}_6+ I^{\ZZ_n}_7+ I^{\ZZ_n}_8
  + I^{\ZZ_n}_9+ I^{\ZZ_n}_{10}
  \,,  \notag \\
 I^{\ZZ_n}_1& =
 \tint \tfrac{1}{\Sigma^{2\beta}} \JJss  \Jg(\Q\p_\s +V\p_2)\nbs^6( \Jg\Zbn) \  \nbs^6(\Jg\Zbn)
  \,, \label{I1-Zbn-P}\\
 I^{\ZZ_n}_2& =
- \tint \tfrac{1}{\Sigma^{2\beta}}  \bubu{  \alpha  ( \Jg\Wbn - \Jg\Zbn)} \JJss\sabs{\nbs^6(\Jg\Zbn)  }^2
  \,, \label{I2-Zbn-P} \\
I^{\ZZ_n}_3 &=
- \alpha \tint \jb  g^{- {\frac{1}{2}} }\JJss \Jg  \nbs_2\nbs^6(\Jg\Abn)  \  \nbs^6(\Jg\Zbn)
\,, \label{I3-Zbn-P} \\
 I^{\ZZ_n}_4 &=
\alpha \tint \jb\JJss \Jg g^{- {\frac{1}{2}} }  \Abn \nbs_2\nbs^6 \Jg \  \nbs^6(\Jg\Zbn)
  \,,   \label{I4-Zbn-P} \\
 I^{\ZZ_n}_5 &=
\tfrac{\alpha }{2}  \tint \jb \JJss \Jg g^{- {\frac{1}{2}} } (\Jg\Wbn + \Jg\Zbn - 2\Jg\Abt)\nbs_2 \nbs^6 \tt\cdo\nn \ \nbs^6(\Jg\Zbn)
 \,, \label{I5-Zbn-P} 
\\
 I^{\ZZ_n}_6& =
-{\tfrac{2 \alpha }{\eps}} \tint  \jb \JJss   \nbs_1 \nbs^6 (\Jg\Zbn) \ \nbs^6(\Jg\Zbn)
\,,  \label{I6-Zbn-P}\\
 I^{\ZZ_n}_7& =
-{\tfrac{2 \alpha }{\eps}} \tint  \jb  \JJss  \Jg(\Abn+\Zbt) \Big( \nbs_1\nbs^6\tt \cdo\nn - \eps \Jg g^{- {\frac{1}{2}} } \nbs_2h\, \nbs_2\nbs^6\tt \cdo\nn   \Big)   \ \nbs^6(\Jg\Zbn)
\,,  \label{I7-Zbn-P}\\
 I^{\ZZ_n}_8& =
{\tfrac{2 \alpha }{\eps}} \tint  \jb \JJss  \Zbn  \Big( \nbs_1\nbs^6\Jg -  \eps \Jg g^{- {\frac{1}{2}} } \nbs_2 h\,  \nbs_2\nbs^6\Jg \Big) \ \nbs^6(\Jg\Zbn)
\,,  \label{I8-Zbn-P}\\
 I^{\ZZ_n}_9& =
2 \alpha \tint  \jb  \JJss \Jg g^{- {\frac{1}{2}} } \nbs_2h  \nbs_2 \nbs^6 (\Jg\Zbn) \ \nbs^6(\Jg\Zbn)
\,,  \label{I9-Zbn-P}\\
 I^{\ZZ_n}_{10}& =
 - \tint  \jb \JJss \big(  \nbs^6 \Fzn + \mathcal{R}^\nn_{\Zb} + \mathcal{C}^\nn_{\Zb}\big) \  \nbs^6(\Jg\Zbn)
\,,  \label{I10-Zbn-P}
\end{align} 
\end{subequations} 
where we have used the notation in~\eqref{eq:lets:confuse:the:reader}.
Lastly, in analogy to~\eqref{Integral-Abn}, we additively decompose the integral  $I^{\AA_n}$ as
\begin{subequations} 
\label{Integral-Abn-P}
 \begin{align}
  I^{\AA_n} & = I^{\AA_n}_1 + I^{\AA_n}_2 + I^{\AA_n}_3+ I^{\AA_n}_4+ I^{\AA_n}_5+ I^{\AA_n}_6+ I^{\AA_n}_7+ I^{\AA_n}_8
  + I^{\AA_n}_9+ I^{\AA_n}_{10}
  \,,  \notag \\
 I^{\AA_n}_1& =
2 \tint \tfrac{1}{\Sigma^{2\beta}}  \JJss \Jg(\Q\p_\s +V\p_2)\nbs^6( \Jg\Abn) \  \nbs^6(\Jg\Abn)
  \,, \label{I1-Abn-P}\\
 I^{\AA_n}_2& =
- \tint \tfrac{1}{\Sigma^{2\beta}} \bubu{  \alpha  \big(\Jg\Wbn - \Jg\Zbn\big)} \JJss\sabs{\nbs^6(\Jg\Abn)}^2
  \,, \label{I2-Abn-P} \\
I^{\AA_n}_3 &=
2 \alpha \tint \jb  g^{- {\frac{1}{2}} } \JJss \Jg  \nbs_2\nbs^6(\Jg\Sbn)  \  \nbs^6(\Jg\Abn)
\,, \label{I3-Abn-P} \\
 I^{\AA_n}_4 &=
-2\alpha \tint \jb  g^{- {\frac{1}{2}} } \JJss (\Jg \Sbn) \nbs_2\nbs^6 \Jg \  \nbs^6(\Jg\Abn)
  \,,   \label{I4-Abn-P} \\
 I^{\AA_n}_5 &=
2 \alpha \tint \jb \JJss \Jg g^{- {\frac{1}{2}} }  \Sbt \nbs_2 \nbs^6 \tt\cdo\nn \ \nbs^6(\Jg\Abn)
 \,, \label{I5-Abn-P} 
\\
 I^{\AA_n}_6& =
-{\tfrac{2 \alpha }{\eps}} \tint  \jb   \JJss \nbs_1 \nbs^6 (\Jg\Abn) \ \nbs^6(\Jg\Abn)
\,,  \label{I6-Abn-P}\\
 I^{\AA_n}_7& =
{\tfrac{\alpha }{\eps}} \tint  \jb \JJss   (\Jg\Wbn+\Jg\Zbn - 2\Jg\Abt) \Big( \nbs_1\nbs^6\tt \cdo\nn - \eps \Jg g^{- {\frac{1}{2}} } \nbs_2h\, \nbs_2\nbs^6\tt \cdo\nn   \Big)  \ \nbs^6(\Jg\Abn)
\,,  \label{I7-Abn-P}\\
 I^{\AA_n}_8& =
{\tfrac{2 \alpha }{\eps}} \tint  \jb  \JJss  \Abn \Big( \nbs_1\nbs^6\Jg -  \eps \Jg g^{- {\frac{1}{2}} } \nbs_2 h\,  \nbs_2\nbs^6\Jg \Big)  \ \nbs^6(\Jg\Abn)
\,,  \label{I8-Abn-P}\\
 I^{\AA_n}_9& =
2 \alpha \tint  \jb  \JJss \Jg g^{- {\frac{1}{2}} } \nbs_2h  \nbs_2 \nbs^6 (\Jg\Abn) \ \nbs^6(\Jg\Abn)
\,,  \label{I9-Abn-P}\\
 I^{\AA_n}_{10}& =
 - \tint  \jb \JJss \big(  \nbs^6 \Fzn + \mathcal{R}^\nn_{\Zb} + \mathcal{C}^\nn_{\Zb}\big) \  \nbs^6(\Jg\Abn)
\,,  \label{I10-Abn-P}
\end{align} 
\end{subequations} 
where we have used the notation in~\eqref{eq:lets:confuse:the:reader}.

The majority of the $26$ integrals listed in~\eqref{Integral-Wbn-P}, \eqref{Integral-Zbn-P}, and \eqref{Integral-Abn-P}, are estimated identically as in Section \ref{sec:sixth:order:energy}. We will explain the modifications that are required for those integrals in which we must integrate-by-parts with respect to $\nbs_1$, both because $\nbs_1^*$ is contains two extra terms (see~\eqref{adjoint-1-P}) which require the new bounds \eqref{eq:Qrs1:bbq-DS}--\eqref{eq:Qr1:bbq-DS}, and because $\nbs_1 \JJ \neq 0$, which requires appealing to~\eqref{eq:D:JJ}.

\subsubsection{Downstream modifications to the exact derivative terms}
The identity~\eqref{eq:heavy:fuel:1n} for $I^{\WW_n}_1  + I^{\ZZ_n}_1 + I^{\AA_n}_1$ remains the same, except that the damping term $\mathsf{G}_0$ carries the natural modification $\mathcal{J} \mapsto \JJ$, i.e., \eqref{eq:G0:n:lower} becomes
\begin{equation}
\mathsf{G}_0 
\geq
\underbrace{- \tfrac 12 (\Q \p_\s + V\p_2) \bigr(\JJ^{\frac 32} \Jg \bigl)}_{=: \mathsf{G_{good}}}
- \tfrac{250^2}{\eps} \Q \Jg  \JJ^{\frac 32}
- \Cn \brak{\beta}  \JJ^{\frac12} \Jg  
\,.
\label{eq:G0:n:lower-P}
\end{equation}
The first real modifications are to the term \eqref{eq:heavy:fuel:2n} concerning $I^{\ZZ_n}_6  + I^{\AA_n}_6$. Here $\p_1$ needs to be replaced by $\frac{1}{\eps} \nbs_1$, and using \eqref{adjoint-1-P}, we have that 
\begin{align}
I^{\ZZ_n}_6  + I^{\AA_n}_6 
&= - {\tfrac{\alpha }{\eps}}  \tint \jb \JJss \nbs_1 \sabs{\nbs^6(\Jg \Zbn,\Jg\Abn )}^2
\notag \\
&
=   \tint \tfrac{1}{\Sigma^{2\beta}} \mathsf{G}_1  \sabs{\nbs^6(\Jg \Zbn,\Jg\Abn)}^2  
+ \alpha
\underbrace{\int\jb \Qb_1 \JJss  \sabs{\nbs^6(\Jg \Zbn,\Jg\Abn )}^2\Big|_\s}_{\geq 0 \mbox{ due to } \eqref{eq:Qrs1:bbq-DS}}
\,,
\label{eq:heavy:fuel:2n-P}
\end{align}
where 
\begin{align*}
\mathsf{G}_1 
:=  - \alpha (2\beta-1) \JJss  \tfrac{1}{\eps} \nbs_1\Sigma 
- \alpha \Sigma \Qr_1 \JJss
+ \tfrac{3\alpha }{2\eps} \Sigma \JJh \nbs_1\JJ 
 \,.
\end{align*} 
Taking $\beta\ge \frac 12$,  and appealing to \eqref{eq:D:JJ},  \eqref{bootstraps-P},  \eqref{eq:Qrs1:bbq-DS}, \eqref{eq:Qr1:bbq-DS}, and \eqref{p1-Sigma-s-P},   we then have that 
\begin{align} 
\mathsf{G}_1 
\geq
\alpha (\beta-\tfrac 12)  \bigl(\tfrac{1}{4\eps} - \tfrac{16}{(1+\alpha) \eps}  \Jg \Q -  \Cn   \bigr)  \JJss 
-  \tfrac{ 5 \alpha \kappa_0 }{\eps} \JJss \,.
\label{eq:heavy:fuel:G1:lowern-P}
\end{align}
The last two terms appearing on the right side of \eqref{eq:heavy:fuel:G1:lowern-P} are the new terms caused by the downstream geometry. The first term will be absorbed by choosing $\beta$ sufficiently large, while the second term is obviously a Gr\"onwall term. 
Combining \eqref{eq:heavy:fuel:2n-P} and \eqref{eq:heavy:fuel:G1:lowern-P}, we obtain that
\begin{align}
I^{\ZZ_n}_6  + I^{\AA_n}_6 
& \ge
\tfrac{\alpha( \beta - \frac 12) -  40 \alpha \kappa_0 }{8 \eps}   
\int_0^{\s} 
\snorm{\tfrac{\JJ^{\frac 34} }{\Sigma^\beta} \nbs^6(\Jg\Zbn,\Jg\Abn) (\cdot,\s')}_{L^2_x}^2
{\rm d} \s'
\notag\\
&
\qquad 
- 
 \tfrac{16 \alpha( \beta - \frac 12)}{(1+\alpha)\eps}   
\int_0^{\s} 
\snorm{\tfrac{\JJ^{\frac 34}(\Jg\Q)^{\frac 12} }{\Sigma^\beta} \nbs^6(\Jg\Zbn,\Jg\Abn) (\cdot,\s')}_{L^2_x}^2
{\rm d} \s' \,.
 \label{I6-Zbn+Abn-lowerbound-P}
\end{align} 

The identity \eqref{eq:heavy:fuel:3n} for $I^{\ZZ_n}_9  + I^{\AA_n}_9$ remains the same. The lower bound~\eqref{eq:heavy:fuel:G2:lowern} for $\mathsf{G}_2$ and \eqref{eq:heavy:fuel:G2:lower:an} for the temporal boundary term, remain unchanged (except for $\mathcal{J}^{\! \frac 12} \mapsto \JJ^{\frac 12}$).
The identity~\eqref{eq:heavy:fuel:4n} for $I^{\ZZ_n}_2  + I^{\AA_n}_2$ remains the same, and the damping coefficient $\mathsf{G}_3$ satisfies the same lower bound as in \eqref{eq:heavy:fuel:G3:lower:an}, except for the usual modification $\mathcal{J}^{\frac 32} \mapsto \JJ^{\frac 32}$. At last, the identity~\eqref{eq:heavy:fuel:5n} for $I^{\WW_n}_3  + I^{\ZZ_n}_3 + I^{\AA_n}_3$ remains the same, except that the  weight function is $\JJ^{\frac 32}$ instead of $\mathcal{J}^{\frac 32}$. With the definition of $\mathsf{G_{good}}$ in \eqref{eq:G0:n:lower-P}, and with the updated lower bound~\eqref{I6-Zbn+Abn-lowerbound-P}, the estimate~\eqref{eq:I:n:12369} summarizing the contribution of all exact derivative terms becomes
\begin{align}
& I^{\WW_n}_1 + I^{\ZZ_n}_1  + I^{\AA_n}_1 + I^{\ZZ_n}_2 + I^{\AA_n}_2
 + I^{\WW_n}_3 + I^{\ZZ_n}_3 + I^{\AA_n}_3
 + I^{\ZZ_n}_6  + I^{\AA_n}_6 
 + I^{\ZZ_n}_9  + I^{\AA_n}_9 
 \notag\\
&\geq \bigl( \tfrac 12 - \Cn \eps\bigr) 
 \snorm{\tfrac{\JJ^{\frac 34} (\Jg \Q)^{\frac 12}}{\Sigma^\beta} \nbs^6(\Jg \Wbn,\Jg \Zbn,\Jg \Abn)(\cdot,\s)}_{L^2_x}^2 
-  \snorm{ \tfrac{\JJ^{\frac 34}(\Jg \Q)^{\frac 12}}{\Sigma^\beta} \nbs^6 (\Jg \Wbn,\Jg \Zbn,\Jg \Abn)(\cdot,0)}_{L^2_x}^2 
\notag\\
&\qquad
+ \tint \tfrac{1}{\Sigma^{2\beta}} \bigl( \mathsf{G_{good}} - \Cn \beta\JJ^{\frac 12} \Jg\bigr)
\sabs{\nbs^6(\Jg \Wbn,\Jg \Zbn,\Jg \Abn)}^2  
\notag\\
&\qquad 
+ \Bigl(  \tfrac{\alpha( \beta - \frac 12) -  40 \alpha\kappa_0 }{8\eps}   + \bubu{  \tfrac{9 \alpha }{10\eps} } \Bigr)
\int_0^{\s}  \snorm{\tfrac{ \JJ^{\frac 34} }{\Sigma^\beta} \nbs^6(\Jg\Zbn,\Jg\Abn) (\cdot,\s')}_{L^2_x}^2 
{\rm d} \s'
\notag\\
&\qquad 
- \Bigl( \tfrac{16 \alpha  (\beta -\frac 12) }{(1+\alpha)\eps}  + \tfrac{\bubu{ 33} + 2 \cdot 250^2}{\eps} \Bigr) 
\int_0^{\s} \snorm{\tfrac{\JJ^{\frac 34}(\Jg \Q)^{\frac 12}}{\Sigma^\beta} \nbs^6(\Jg \Zbn,\Jg\Abn)(\cdot,\s')}_{L^2_x}^2 {\rm d} \s'
\notag\\
&\qquad
- \tfrac{250^2}{\eps}
\int_0^{\s} 
 \snorm{\tfrac{\JJ^{\frac 34}(\Jg \Q)^{\frac 12}}{\Sigma^\beta} \nbs^6(\Jg \Wbn)(\cdot,\s')}_{L^2_x}^2 {\rm d} \s'
\,,
\label{eq:I:n:12369-P}
\end{align}
where  $\Cn = \Cn (\alpha,\kappa_0,\Cdata)$ is a positive constant independent of $\beta$. 

\subsubsection{Downstream update to geometric lemmas}
In order to deal with terms that contain over-differentiated geometry, we next generalize Lemmas~\ref{lem:tau-Jg-D2} and \ref{lem:tau-Jg-D1} to the geometry of $\Pdpm(\s)$.

\begin{lemma} \label{lem:tau-Jg-D2-P}
For a function $f(x,\s)$, we have that
\begin{equation} 
\left| \tint f \, \nn \cdo \nbs^6\tt  \ \nbs_2 \nbs^6 \Jg \right| 
+
\left| \tint f \, \nn \cdo \nbs_2 \nbs^6\tt  \  \nbs^6 \Jg \right| 
\les  \eps^{3} \mathsf{K}^2 \brak{\mathsf{B}_6}^2 
\Bigl(\snorm{\JJh \nbs f}_{L^\infty_{x,\s}}  
+   \snorm{\JJmh f}_{L^\infty_{x,\s}}  
   \Bigr) \,.
\label{eq:tau-Jg-D2-P}
\end{equation} 
\end{lemma} 
We note that when compared to \eqref{eq:tau-Jg-D2}, the bound~\eqref{eq:tau-Jg-D2-P} is missing a factor of $\eps$ next to $\|\JJmh f\|_{L^\infty_{x,\s}}$. This helpful factor of $\eps$ was however never used in any application of \eqref{lem:tau-Jg-D2}, so that we may use~\eqref{eq:tau-Jg-D2-P} as is in all energy estimates.
\begin{proof}[Proof of Lemma \ref{lem:tau-Jg-D2-P}]
We  explain the modification of the proof of Lemma \ref{lem:tau-Jg-D2} that will lead to the inequality \eqref{eq:tau-Jg-D2-P}.  In
particular, we explain the modifications to \eqref{big-chicken2} which occur when replacing the formula \eqref{adjoint-1} with the new
identity \eqref{adjoint-1-P}.   The identity  \eqref{big-chicken2} is replaced with
\begin{align}
\tint f \, \nn \cdo \nbs^6\tt  \, \nbs_2 \nbs^6 \Jg 
&= 
- \tfrac{1}{2\eps} \tint  \nbs_1 (f \, g^{\frac 12})   \, (\nn \cdo \nbs^6 \tt)^2 
- \tfrac 12 \tint \Big( (\Qr_2 - \nbs_2) (f \, \Jg \nbs_2 h) -  \Qr_1  f \, g^{\frac 12}\Big) \, (\nn \cdo \nbs^6 \tt)^2
\notag\\
&\qquad 
+ \tfrac 12 \int\big(  \Qb_2 f \, \Jg \nbs_2 h  - \Qb_1  f \, g^{\frac 12} \big)  (\nn \cdo \nbs^6 \tt)^2 \Bigr|_{\s}
- \tfrac{1}{\eps} \tint f \, \nn \cdo \nbs^6\tt  \,g^{\frac 12} \tt \cdo \nbs_1  \nn\; \tt \cdo \nbs^6 \tt
\notag\\
&\qquad
+ \tfrac{1}{\eps} \tint f \, \nn \cdo \nbs^6\tt  \,\nn \cdo \nbs^6(g^{\frac 12} \nn) \;  \nn \cdo \nbs_1 \tt
+ \tfrac{1}{\eps} \tint f \, \nn \cdo \nbs^6\tt  \,\doublecom{\nbs^6, g^{\frac 12} \nn_i, \nbs_1 \tt_i}
\notag\\
&\qquad 
+ \tint f \, \nn \cdo \nbs^6\tt  \,\Jg \nbs_2 h \tt \cdo \nbs_2 \nn \; \tt \cdo \nbs^6 \tt
- \tint f \, \nn \cdo \nbs^6\tt  \,\nn \cdo \nbs^6(\Jg \nbs_2 h \nn) \; \nn \cdo \nbs_2 \tt 
\notag\\
&\qquad
- \tint f \, \nn \cdo \nbs^6\tt  \, \doublecom{\nbs^6, \Jg \nbs_2 h \nn_i, \nbs_2 \tt_i}  \,.
\label{big-chicken2-P}
\end{align}
Using the bounds established in Proposition~\ref{prop:geometry-P}, the new terms in the above expression (when compared to~\eqref{big-chicken2}), due to $\Qr_1$ and $\Qb_1$, are estimated from above by $5 \| \JJ^{-\frac 12} f \|_{L^\infty_{x,\s}} (\eps^{-1} \|\JJ^{\frac 12}\|_{L^\infty_{x,\s}}^2 \|\JJ^{-\frac 14} \nn\cdot\nbs^6 \tt\|_{L^2_{x,\s}}^2 +  \|\JJ^{\frac 14} \nn\cdot\nbs^6 \tt\|_{L^\infty_{\s}L^2_{x}}^2) 
\les \| \JJ^{-\frac 12} f \|_{L^\infty_{x,\s}} \eps^3 \mathsf{K}^2 \brak{\mathsf{B}_6}^2$. This estimate is clearly consistent with \eqref{eq:tau-Jg-D2-P}. Using Lemma~\ref{lem:Q:bnds-DS} and the identical procedure that we used in the proof of Lemma~\ref{lem:tau-Jg-D2} (except that any reference to the bounds in  Proposition~\ref{prop:geometry} need to be changed to the bounds  in Proposition~\ref{prop:geometry-P}), all other terms present on the right side of \eqref{big-chicken2-P} may be shown to be bounded from above by the right side of \eqref{eq:tau-Jg-D2-P}. The bound for the second term on the left side of \eqref{eq:tau-Jg-D2-P} is obtained by using $\nbs_2^*$ to convert $\nbs_2\nbs^6\tt \nbs^6 \Jg$ into $-\nbs^6 \tt \nbs_2 \nbs^6\Jg$, as in \eqref{big-chicken-worm2}. Since the formula for $\nbs_2^*$ in \eqref{adjoint-2-P} is the same as the one used in \eqref{big-chicken-worm2}, the bounds are identical to each other.
\end{proof}

\begin{lemma} \label{lem:tau-Jg-D1-P}
For a function $f(x,\s)$, we have that
\begin{equation} 
\left|\tint f \, \nn \cdo \nbs_1 \nbs^6\tt   \  \nbs^6 \Jg\right|
+
\left|\tint f \, \nn \cdo \nbs^6\tt  \ \nbs_1 \nbs^6 \Jg\right|
\les
\eps^3 \mathsf{K} \brak{\mathsf{B}_6}^2
\Bigl( \|\JJh \nbs f\|_{L^\infty_{x,\s}} +  \|\JJmh f\|_{L^\infty_{x,\s}} \Bigr) 
\,. 
\label{eq:tau-Jg-D1-P}
\end{equation} 
\end{lemma} 
\begin{proof}[Proof of Lemma \ref{lem:tau-Jg-D1-P}]
The bound for the integral $\int_0^\s \!   \int f \, \nn \cdo \nbs_1 \nbs^6\tt   \  \nbs^6 \Jg$ is unchanged from the proof of Lemma
\ref{lem:tau-Jg-D1} and uses~\eqref{eq:tau-Jg-D1-chicken-worm}.  To bound the integral $\int_0^\s \!   \int f \, \nn \cdot \nbs^6\tt  \ \nbs_1 \nbs^6 \Jg$, we use 
the identity \eqref{adjoint-1-P}, and obtain
\begin{align*} 
\tint f \, \nn \cdo \nbs^6\tt   \  \nbs_1 \nbs^6 \Jg  
&= -\tint f \, \nn \cdo  \nbs_1 \nbs^6\tt   \  \nbs^6 \Jg
- \tint \big( \nbs_1(f \, \nn_i) - \eps\Qr_1 f \, \nn_i\big) \, \nbs^6\tt _i  \  \nbs^6 \Jg 
\notag \\
& \qquad\qquad
-\eps \int\Qb_1 f \, \nn\cdo \nbs^6\tt  \  \nbs^6 \Jg\Big|_\s \,.
\end{align*} 
We have established that the first integral on the right side has the bound \eqref{eq:tau-Jg-D1-P}.  The remaining integrals
are bounded using
\eqref{bootstraps-P}, \eqref{geometry-bounds-new-P}, and Lemma \ref{lem:Q:bnds-DS}, and the
desired inequality follows.
\end{proof} 

\subsubsection{Downstream modifications to bounds for the forcing, remainder, and commutator functions}
\label{subsubsec:chicken0}
The bounds obtained earlier in Section~\ref{subsub:FRC} hold without change. Note that in these terms there is no integration by parts being used. One only uses the bootstraps~\eqref{bootstraps-P}, the updated bounds for the geometry from Proposition~\ref{prop:geometry-P}, the improved bounds for normal components~\eqref{improved-P}, especially for $\Zbn$, the estimate~\eqref{eq:madman:2-P} for pure time derivatives of $\Zbn$, and the commutator and product estimates from Appendix~\ref{app:functional}. The bounds~\eqref{eq:I:W:nn:6:all}, \eqref{eq:I:Z:nn:10:all}, and \eqref{eq:I:A:nn:10:all} hold as is, except that the weight $\mathcal{J}$ needs to be updated with the new weight $\JJ$.

\subsubsection{Downstream modifications to the terms with over-differentiated geometry}
\label{subsubsec:chicken1}

We discuss the modifications to the following sets of terms, all of which contain over-differentiated geometry: $I^{\WW_n}_5+I^{\ZZ_n}_5+I^{\AA_n}_7$ (see Section~\ref{sec:IW5+IZ5+IA7}), $I^{\WW_n}_4+I^{\ZZ_n}_4+I^{\AA_n}_8$ (see Section~\ref{sec:IW4+IZ4+IA8}), $I^{\AA_n}_4$ (see Section~\ref{sec:IA4}), $I^{\ZZ_n}_7$ (see Section~\ref{sec:I:Zn:7:integral}), $I^{\AA_n}_5$ (see Section~\ref{sec:IA5}), and $I^{\ZZ_n}_8$ (see Section~\ref{sec:IZ8}). Several of these terms are bounded in precisely the same way as in Section~\ref{sec:normal:overdiff}, either because they do not contain the operator $\frac{1}{\eps} \nbs_1$, or because they are dealt with by directly appealing to Lemmas~\ref{lem:tau-Jg-D2-P} or~\ref{lem:tau-Jg-D1-P}, which have already been updated to the downstream geometry. For the convenience of the reader, we go through these terms one by one.

We offer more details for the analysis of the combination $I^{\WW_n}_5+I^{\ZZ_n}_5+I^{\AA_n}_7$, which was previously analyzed in Section~\ref{sec:IW5+IZ5+IA7}. We view these modifications as a template for the necessary changes to all other terms discussed in the above paragraph. First, consider the term $I^{\AA_n}_7$, as defined in~\eqref{I7-Abn-P}. Taking into account~\eqref{eq:adjoints-P}, the decomposition~\eqref{eq:IAAn:7:decompose} becomes 
\begin{align*} 
I^{\AA_n}_7& = I^{\AA_n}_{7,a} + I^{\AA_n}_{7,b} + I^{\AA_n}_{7,c} +  I^{\AA_n}_{7,d }  \,, 
\\
I^{\AA_n}_{7,a} &=
-{\tfrac{\alpha }{\eps}} \tint  \jb  \JJ^{\! \frac 32} (\Jg\Wbn+\Jg\Zbn - 2\Jg\Abt) \nbs^6\tt \cdo\nn  
\ \Big( \nbs_1\nbs^6(\Jg\Abn) -  \eps\Jg g^{- {\frac{1}{2}} } \nbs_2h\, \nbs_2\nbs^6(\Jg\Abn) \Big) \,, 
\notag \\
I^{\AA_n}_{7,b} &=
-{\tfrac{\alpha }{\eps}} \tint  \Big(\nbs_1 - \eps \Jg g^{- {\frac{1}{2}} }  \nbs_2 h \nbs_2 \Big)\Big( \jb  \JJ^{\! \frac 32} (\Jg\Wbn+\Jg\Zbn - 2\Jg\Abt) \nn_i \Big) \nbs^6\tt_i  \ \nbs^6(\Jg\Abn) \,, 
\notag \\
I^{\AA_n}_{7,c} &=
\alpha \tint  \Bigl(\Qr_1 -  \Jg g^{- {\frac{1}{2}} } \nbs_2 h \, \Qr_2 \Bigr) \jb \JJ^{\! \frac 32}  (\Jg\Wbn+\Jg\Zbn - 2\Jg\Abt) \nbs^6\tt \cdo\nn  \ \nbs^6(\Jg\Abn) \,, 
\notag \\
I^{\AA_n}_{7,d} &=
- \alpha  \int \Bigl(\Qb_1 - \Jg g^{- {\frac{1}{2}} } \nbs_2 h \, \Qb_2 \Bigr)  \jb\JJ^{\! \frac 32}   (\Jg\Wbn+\Jg\Zbn - 2\Jg\Abt) \nbs^6\tt \cdo\nn  \ \nbs^6(\Jg\Abn)\Big|_\s 
\notag\,.
\end{align*} 
The only change  to $I^{\AA_n}_{7,b}$ is that the differential operator $\nbs_1$ may act on the weight $\JJ^{\frac 32}$. This term is handled using \eqref{eq:D:JJ}, which gives $\nbs_1 \JJ = \Qb_1$, and the upper bound in~\eqref{eq:Qrs1:bbq-DS}.
The changes to $I^{\AA_n}_{7,c}$ and $I^{\AA_n}_{7,d}$ are the emergence of the $\Qr_1$ and respectively $\Qb_1$ terms. Note however that these terms are only nonzero in $\Pdp(\s)$ (see~\eqref{eq:Qrs1:bbq-DS} and~\eqref{eq:Qr1:bbq-DS}), and we have that 
\begin{equation*}
\JJ \leq \Jgb \leq \Jg\,,
\end{equation*}
due to Remark~\ref{rem:downstream:terminates}, and to~\eqref{Jgb-le-Jg}--\eqref{eq:fake:Jg:def}. Using the bootstraps~\eqref{bootstraps-P}, the bounds in Lemma~\ref{lem:Q:bnds-DS} for the remapping coefficients, and the bounds for the geometry from Proposition~\ref{prop:geometry-P}, we deduce 
\begin{equation*}
|I^{\AA_n}_{7,b}| + 
|I^{\AA_n}_{7,c}| + 
|I^{\AA_n}_{7,d}| 
\les  (\tfrac{4}{\kappa_0})^{2\beta}  \mathsf{K} \Bsix^2 
\,,
\end{equation*} 
which precisely matches~\eqref{I7-An-bcd}. 
The analysis of $I^{\AA_n}_{7,a}$ requires the decomposition in~\eqref{IAbn-7a-i-viii} (with $\mathcal{J}$ replaced by $\JJ$). Among the eight terms appearing in the decomposition~\eqref{IAbn-7a-i-viii}, the last five terms require no modifications beyond those already given by Lemmas~\ref{lem:tau-Jg-D2-P} and~\ref{lem:tau-Jg-D1-P}, so that the bounds~\eqref{small-chicken1}--\eqref{small-chicken3} remain the same. To see this, consider for instance a term which involves the $\nbs_1$ operator, such as $I^{\AA_n}_{7,a,vi}$. For this term we form an exact derivative and integrate-by-parts with respect to  $\nbs_1$ using \eqref{adjoint-1-P} to obtain that
\begin{align*} 
& - {\tfrac{\alpha }{2\eps}} \tint  \jb \JJ^{\!\frac 32}  (\Jg\Wbn+\Jg\Zbn - 2\Jg\Abt)^2 \nbs^6\tt \cdo\nn 
 \ \Jg   \nbs_1\nbs^6\tt \cdot\nn  
 \notag\\ 
 &= 
{\tfrac{\alpha }{4\eps}} \tint ( \nbs_1-\eps \Qr_1)\Big(\jb \JJ^{\!\frac 32}  (\Jg\Wbn+\Jg\Zbn - 2\Jg\Abt)^2 \nn^i\nn^j \Big)\nbs^6\tt_i \nbs^6\tt_j 
\notag \\
&\qquad
+{\tfrac{\alpha }{4}} \int \Qb_1 \jb \JJ^{\!\frac 32}  (\Jg\Wbn+\Jg\Zbn - 2\Jg\Abt)^2 \, \sabs{\nbs^6\tt\cdo\nn}^2\Big|_\s \,.
\end{align*} 
From \eqref{small-chicken3},  \eqref{bootstraps-P}, \eqref{geometry-bounds-new-P}, and  Lemma \ref{lem:Q:bnds-DS},  we see that the above term is bounded by $(\tfrac{4}{\kappa_0})^{2\beta}  \mathsf{K}^2 \Bsix^2$, consistent with \eqref{small-chicken3}. For the second term  in the decomposition~\eqref{IAbn-7a-i-viii}, namely $I^{\AA_n}_{7,a,ii}$, the only downstream modification stems from the term $I^{\AA_n}_{7,a,ii_3}$, in which $(\Q\p_\s + V\p_2)$ acts on $\JJ^{\frac 32}$ instead of $\mathcal{J}^{\frac 32}$; this term  satisfies however the same bound (up to a universal constant), so that \eqref{small-chicken5} also remains the same. This leads to 
\begin{equation*}
|I^{\AA_n}_{7,a,ii}| 
+
|I^{\AA_n}_{7,a,iii}| +|I^{\AA_n}_{7,a,iv}| +|I^{\AA_n}_{7,a,v}| +|I^{\AA_n}_{7,a,vi}| +|I^{\AA_n}_{7,a,vii}| + |I^{\AA_n}_{7,a,viii}| 
\les  (\tfrac{4}{\kappa_0})^{2\beta}  \mathsf{K} \Bsix^2 
\,.
\end{equation*} 
The last term to be considered in the decomposition~\eqref{IAbn-7a-i-viii} is $I^{\AA_n}_{7,a,i}$. We use $\nbs_2^*$ to decompose this term into four pieces, as in~\eqref{IAn7,a,i}. The first of these terms precisely cancels the sum $I^{\WW_n}_5+I^{\ZZ_n}_5$ which together with the bound \eqref{small-chicken4} for the remaining three terms implies
\begin{equation*}
|I^{\WW_n}_5+I^{\ZZ_n}_5 + I^{\AA_n}_{7,a,i}| 
\les \eps (\tfrac{4}{\kappa_0})^{2\beta}  \mathsf{K} \Bsix^2 
\,.
\end{equation*}
Adding the estimates in the above three displayed inequalities leads to a bound for $I^{\WW_n}_5+I^{\ZZ_n}_5 + I^{\AA_n}_{7}$ which precisely matches~\eqref{IWn5+IZn5+IAn7}.

Next, we consider the combination~$I^{\WW_n}_4+I^{\ZZ_n}_4 + I^{\AA_n}_{8}$. A close inspection of the analysis in Section~\ref{sec:IW4+IZ4+IA8}, shows that the all-important cancellation~\eqref{eq:ducky:von:ducken:1} remains the same, and that upon appealing to the updated Lemmas~\ref{lem:tau-Jg-D2-P} and~\ref{lem:tau-Jg-D1-P} all other bounds remain unchanged, leading to a bound for $|I^{\WW_n}_4+I^{\ZZ_n}_4 + I^{\AA_n}_{8}|$  which exactly matches~\eqref{eq:small:duck:0}.

Next, we consider the term $I^{\AA_n}_{4}$, which was clearly the most challenging one to handle in Section~\ref{sec:normal:overdiff}. The estimates~\eqref{eq:small:duck:7} and~\eqref{eq:small:duck:8}  remain unchanged, so that we only need to estimate the term $I^{\AA_n}_{4,a,i}$ which was defined in~\eqref{I4-Abn-a,i}. Exactly as in~\eqref{I4-Abn-a-i}, we decompose $I^{\AA_n}_{4,a,i}$ into nine parts. Among these nine parts, six of them are bounded in precisely the same way, via~\eqref{eq:small:duck:8a} and~\eqref{eq:small:duck:8b}. It thus remains to carefully analyze the $J^{\AA_n}_{1}$, $J^{\AA_n}_{3}$, and $J^{\AA_n}_{6}$ terms appearing in the decomposition~\eqref{I4-Abn-a-i}, which we recall here for convenience
\begin{subequations}
\label{I4-Abn-a-i-DS}
\begin{align} 
J^{\AA_n}_{1} &  =  \tint \tfrac{1}{\Sigma^{2\beta}}  \JJ^{\!\frac 32}   (\Jg\Wbn - \Jg \Zbn) (\Q\p_\s+V\p_2) \nbs^6\Jg  \  \nbs^6(\Jg\Wbn)  \,,
  \\
J^{\AA_n}_{3} &  = - \tfrac{3}{2\eps} \tint \tfrac{\Q}{\Sigma^{2\beta}}   \JJ^{\!\frac 12}  \  (\Jg\Wbn)  \nbs^6\Jg  \  \nbs^6(\Jg\Wbn)  \,,
  \\
J^{\AA_n}_{6} & =
-  \int \tfrac{\Q}{\Sigma^{2\beta}} \JJ^{\!\frac 32}    \Jg \Wbn \nbs^6 \Jg \  \nbs^6(\Jg\Wbn  ) \Big|_\s
  \,.
\end{align}
\end{subequations}
For $J^{\AA_n}_{1}$, we first commute $\Q\p_\s+V\p_2$ past $\nbs^6$ acting on $\Jg$ using  \eqref{good-comm-P}. The commutator is lower order, while the principal term gives via \eqref{Jg-evo-s-P} a contribution of the type $(\frac{1+\alpha}{2} \nbs^6 (\Jg \Wbn) + \frac{1-\alpha}{2} \nbs^6 (\Jg \Zbn) )$. The next trick is to rewrite $\Jg \Sbn = \frac 12 (\Jg \Wbn) - \frac 12 (\Jg \Zbn) = \frac{1}{1+\alpha} (\Q \p_\s + V \p_2)\Jg - \frac{1}{1+\alpha} (\Jg \Zbn)$. This allows us to rewrite $J^{\AA_n}_{1}$ exactly as in \eqref{eq:small:duck:big:duck:bear}, and appeal to \eqref{eq:small:duck:9}--\eqref{eq:small:duck:10} to estimate
\begin{align}
 J^{\AA_n}_{1} 
 &\geq \tint \tfrac{1}{\Sigma^{2\beta}} 
 \underbrace{\JJ^{\!\frac 32} (\Q\p_\s+V\p_2) \Jg}_{=:\mathsf{G_{bad}}}    \sabs{  \nbs^6(\Jg\Wbn)  }^2
 \notag\\
 &\qquad
 - \Cn (\tfrac{4}{\kappa_0})^{2\beta}  \brak{\mathsf{B}_6}^2 
 - \tfrac{1+\alpha}{\eps} 
\int_0^{\s}  \snorm{\tfrac{\mathcal{J}^{\frac 34}}{\Sigma^\beta} \nbs^6(\Jg\Zbn) (\cdot,\s')}_{L^2_x}
\snorm{\tfrac{\mathcal{J}^{\frac 34}}{\Sigma^\beta} \nbs^6(\Jg\Wbn)(\cdot,\s')}_{L^2_x} {\rm d} \s'
\,.
\label{chicken-choop33}
\end{align}
The last two terms on the right side of the above display are standard, while the first one, containing the ``anti-damping term'' $\mathsf{G_{bad}}$, needs to be combined with the ``damping term'' containing $\mathsf{G_{good}}$, which is already present in \eqref{eq:I:n:12369-P}. Recalling the definition of $\mathsf{G_{good}}$ from \eqref{eq:G0:n:lower-P}, and using the definition of $\mathsf{G_{bad}}$ given above, we next claim that 
\begin{align}
 \mathsf{G_{good}} + \mathsf{G_{bad}}
 &= 
 - \tfrac 12 (\Q \p_\s + V\p_2)  (\JJ^{\frac 32} \Jg)
 + \JJ^{\!\frac 32} (\Q\p_\s+V\p_2) \Jg
 \notag\\
 &= \tfrac 12  \Bigl( \JJ^{\! \frac 32}(\Q \p_\s + V\p_2)  \Jg 
  -   \Jg (\Q \p_\s + V\p_2)  \JJ^{\! \frac 32} \Bigr)
 \geq 
  \tfrac{1+ \alpha }{16\eps} \JJ^{\!\frac 12} \Jg  
 \,,
\label{eq:small:duck:11:a-P}
\end{align}
which serves as the replacement of the lower bound \eqref{eq:small:duck:11:a} in our downstream analysis. To see that \eqref{eq:small:duck:11:a-P} holds, we consider separately the regions $\Pdm(\eps)$ and $\Pdp(\eps)$. For $(x,\s) \in \Pdm(\eps)$, by construction (see~\eqref{eq:fake:Jg:def:ds}) we have that  $\JJ = \mathcal{J}$, and so we can apply~\eqref{eq:fakeJg:LB}, to obtain the lower bound $\tfrac{1+\alpha}{16\eps} \mathcal{J}^{\frac 12} \Jg = \tfrac{1+\alpha}{16 \eps} \JJ^{\!\frac 12} \Jg$, which is  identical to \eqref{eq:small:duck:11:a}, and is consistent with \eqref{eq:small:duck:11:a-P}. For $(x,\s) \in \Pdp(\eps)$, we cannot directly appeal to \eqref{eq:fakeJg:LB}, and instead need to revisit the proof of this bound, which we detail as follows.
Using the bounds $0 < \JJ \leq \Jg$, the identities~\eqref{QQQ-P}, \eqref{eq:D:JJ}, \eqref{Jg-evo-s-P}, the bootstrap assumptions~\eqref{boots-PP}, and the estimates~\eqref{eq:broncos:eat:shit:20} and~\eqref{eq:Qd:bbq-DS}, we deduce
\begin{align}
&\tfrac 12 \JJ^{\!\frac 12} \Bigl( \JJ (\Q \p_\s + V\p_2)  \Jg 
- \tfrac 32  \Jg (\Q \p_\s + V\p_2)  \JJ \Bigr)
\notag\\
&\qquad =  
\tfrac 12 \JJ^{\!\frac 12} \Bigl( 
\JJ \bigl( \tfrac{1+\alpha}{2} \Jg \Wbn + \tfrac{1-\alpha}{2} \Jg \Zbn \bigr)
+ \tfrac 32  \Jg \tfrac{\Q}{\eps}
\Bigr)
\notag \\
&\qquad \geq
\tfrac 12 \JJ^{\!\frac 12} \Bigl(
- \JJ  \bigl( \tfrac{1+\alpha}{2\eps} + \Cn \bigr)
+ \tfrac 32  \Jg \bigl(\tfrac{\Qd}{\eps} - \Cn \bigr)
\Bigr)
\notag \\
&\qquad \geq
\tfrac 12 \JJ^{\!\frac 12} \Bigl(
- \Jg  \bigl( \tfrac{1+\alpha}{2\eps} + \Cn \bigr)
+ \tfrac 32  \Jg \bigl(\tfrac{17(1+\alpha)}{40\eps} - \Cn \bigr)
\Bigr)
= \tfrac 12 \JJ^{\!\frac 12} \Bigl(
\Jg \tfrac{11(1+\alpha)}{80 \eps} - \Cn \Jg
\Bigr)
\geq \JJ^{\!\frac 12}\Jg \tfrac{1+\alpha}{16\eps}
\,.
\label{chicken-choop3}
\end{align}
The above bound is identical to \eqref{eq:fakeJg:LB}, and concludes the proof of~\eqref{eq:small:duck:11:a-P}. Indeed, combining~\eqref{chicken-choop33}--\eqref{chicken-choop3}, we arrive at
\begin{align}
 J^{\AA_n}_{1} 
 &\geq \tint \tfrac{1}{\Sigma^{2\beta}} 
\bigl( \tfrac{1+\alpha}{24 \eps} \JJ^{\!\frac 12} \Jg - \mathsf{G_{good}}\bigr)  \sabs{  \nbs^6(\Jg\Wbn)  }^2
 - \Cn (\tfrac{4}{\kappa_0})^{2\beta}  \brak{\mathsf{B}_6}^2 
 - \tfrac{12 (1+\alpha)}{\eps} 
\int_0^{\s}  \snorm{\tfrac{\mathcal{J}^{\frac 34}}{\Sigma^\beta} \nbs^6(\Jg\Zbn) (\cdot,\s')}_{L^2_x}^2 {\rm d} \s'
\,.
\label{chicken-choop333}
\end{align}
The above bound is identical to the one obtained earlier in Section~\ref{sec:IA4}. 
We return now to the terms $J^{\AA_n}_{3}$ and $J^{\AA_n}_{6}$ given in \eqref{I4-Abn-a-i-DS}.
For $J^{\AA_n}_{3}$ the trick is to again rewrite $\nbs^6 (\Jg \Wbn)  =  \tfrac{2}{1+\alpha}(\Q \p_\s + V\p_2)  \nbs^6  \Jg + \tfrac{2}{1+\alpha} \bigl(\nbs^6 V \nbs_2 \Jg + \doublecom{\nbs^6,V,\nbs_2 \Jg} \bigr)- \tfrac{1-\alpha}{1+\alpha} \nbs^6 (\Jg \Zbn)$, leading to a decomposition identical (except for changing  $\mathcal{J}$ to $\JJ$) to~\eqref{eq:J3:An:decompose}, because the operators $\p_1$ or $\nbs_1$ are not involved here. Since in the downstream development we still have $(\Q\p_\s+V \p_2)\JJ = -\frac{\Q}{\eps}$, and since $\Q$ satisfies the bound \eqref{eq:Q:bbq-DS}, which is identical to~\eqref{eq:Q:bbq}, the bounds~\eqref{eq:fuck:yeah:0} for the nine terms appearing in the decomposition \eqref{eq:J3:An:decompose} of $J^{\AA_n}_{3}$ remain {\bf as is}. This means that we may use precisely the same Cauchy-Young inequality for the $J^{\AA_n}_{6}$ term as in \eqref{eq:fuck:yeah:9}. Together with~\eqref{chicken-choop333}, we may summarize the lower bound for $I^{\AA_n}_{4}$ as
\begin{align}
\label{eq:fuck:yeah:end-P}
 I^{\AA_n}_{4}
&\geq  
- \Cn (\tfrac{4}{\kappa_0})^{2\beta}  \brak{\mathsf{B}_6}^2
- \tfrac{1+\alpha}{\eps} \bigl(1 +\tfrac{132}{(1+\alpha)^4}  \bigr) (\tfrac{3}{\kappa_0})^{2\beta}  \Cdatatwo
\notag\\
&\qquad 
+ \tint \tfrac{1}{\Sigma^{2\beta}} \bigl(\tfrac{1+\alpha}{24\eps} \JJ^{\!\frac 12} \Jg - \mathsf{G_{good}} \bigr)   \sabs{  \nbs^6(\Jg\Wbn)  }^2
- \tfrac{(12+ 25^2)(1+\alpha)}{\eps} 
\int_0^{\s} \snorm{\tfrac{\JJ^{\frac 34}}{\Sigma^\beta} \nbs^6(\Jg\Zbn) (\cdot,\s')}_{L^2_x}^2 {\rm d} \s'
\notag\\
&\qquad
+ \tfrac{1}{20(1+\alpha)} \tfrac{1}{\eps^2} \snorm{\tfrac{\Q \JJ^{\frac 14}}{\Sigma^\beta} \nbs^6 \Jg(\cdot,\s)}_{L^2_x}^2 
-  \tfrac{20^2 + 500^2 +100 \cdot 250^2 + \eps \Cn \brak{\beta}}{(1+\alpha)}  {\tfrac{1}{\eps^3}} \int_0^{\s}   \snorm{\tfrac{\Q \JJ^{\frac 14}}{\Sigma^\beta} \nbs^6 \Jg(\cdot,\s')}_{L^2_x}^2 {\rm d} \s'
\notag\\
&\qquad
+ \tfrac{7}{40(1+\alpha)} \tfrac{1}{\eps^3} \int_0^{\s} \snorm{\tfrac{\Q \JJ^{-\frac 14}}{\Sigma^\beta} \nbs^6 \Jg(\cdot,\s')}_{L^2_x}^2 {\rm d} \s'
+ \bigl(8 - \Cn \brak{\beta} \eps \bigr)  {\tfrac{1}{\eps^3}}  \tint \tfrac{\Q}{\Sigma^{2\beta}} \JJ^{\!\frac 12}   \sabs{ \nbs^6\Jg}^2
\notag\\
&\qquad
- \tfrac{25}{52}
\snorm{\tfrac{\JJ^{\frac 34} (\Jg \Q)^{\frac 12}}{\Sigma^\beta} \nbs^6(\Jg\Wbn)(\cdot,\s)}_{L^2_x}^2
- 
\tfrac{ 39^2+ 500^2 }{\eps}  
\int_0^{\s} 
\snorm{\tfrac{\JJ^{\frac 34} (\Jg\Q)^{\frac 12}}{\Sigma^{\beta}} \nbs^6 (\Jg \Wbn)(\cdot,\s')}_{L^2_x}^2 {\rm d}\s' 
\,,
\end{align}
a bound which is identical to the one obtained earlier in \eqref{eq:fuck:yeah:end}.
 
Next, we consider the term $I^{\ZZ_n}_{7}$ defined in~\eqref{I7-Zbn-P}. Previously, this term was handled in Section~\ref{sec:I:Zn:7:integral}. In this term we need to appeal to the new formula for $\nbs_1^*$ (see~\eqref{adjoint-1-P}). Accordingly, the decomposition~\eqref{eq:I:Zn:7:decompose} is such that $I^{\ZZ_n}_{7,a}$  remains  identical, $I^{\ZZ_n}_{7,b}$ remains the same except that $\nbs_1$ may act on $\JJ^{\frac 32}$, while $I^{\ZZ_n}_{7,c}$ and $I^{\ZZ_n}_{7,d}$ become
\begin{align*}
I^{\ZZ_n}_{7,c} 
&= \mbox{the old }I^{\ZZ_n}_{7,c} \mbox{ from }\eqref{eq:I:Zn:7:decompose} - 2 \alpha \tint   \Qr_1 \jb \JJ^{\!\frac 32} \Jg(\Abn + \Zbt)
\  \nbs^6\tt\cdo\nn  \, \nbs^6(\Jg\Zbn) 
\,, \\
I^{\ZZ_n}_{7,d} 
&= \mbox{the old }I^{\ZZ_n}_{7,d} \mbox{ from }\eqref{eq:I:Zn:7:decompose} + 2 \alpha \int  \Qb_1 \jb  \JJ^{\!\frac 32} \Jg(\Abn + \Zbt) \  \nbs^6\tt\cdo\nn  \, \nbs^6(\Jg\Zbn)\Big|_\s
\,.
\end{align*}
By appealing to~\eqref{eq:D:JJ}, \eqref{bootstraps-P}, \eqref{eq:Qrs1:bbq-DS}, \eqref{eq:Qr1:bbq-DS}, and~\eqref{geometry-bounds-new-P}, we see that the bound \eqref{eq:duck:a:duck:1} for $|I^{\ZZ_n}_{7,b}  + I^{\ZZ_n}_{7,c} + I^{\ZZ_n}_{7,d} |$ remains the same. It thus remains to consider the term $I^{\ZZ_n}_{7,d}$, which is decomposed in eight parts according to~\eqref{I-Zbn-7a-decomp}. By once again employing the updated inequalities \eqref{eq:tau-Jg-D2-P} and \eqref{eq:tau-Jg-D1-P}, using exact-derivative structure, the updated adjoint identities~\eqref{eq:adjoints-P},   the updated coefficient~\eqref{eq:Q:all:bbq-DS} and geometry~\eqref{geometry-bounds-new-P} bounds, we find that following step by step the procedure outlined in~\eqref{eq:duck:a:duck:2}--\eqref{eq:duck:a:duck:4}, we arrive at the bound $|I^{\ZZ_n}_{7,a}| \les  \mathsf{K} \eps   (\tfrac{4}{\kappa_0})^{2\beta} \brak{\mathsf{B}_6}^2$. Putting this together with the bound for $|I^{\ZZ_n}_{7,b}  + I^{\ZZ_n}_{7,c} + I^{\ZZ_n}_{7,d} |$ discussed earlier, we arrive at the same conclusion as~\eqref{eq:duck:a:duck:all}, namely that
\begin{align} 
\sabs{ I^{\ZZ_n}_{7}} \les  \mathsf{K} \eps   (\tfrac{4}{\kappa_0})^{2\beta} \brak{\mathsf{B}_6}^2  \,.
 \label{fatchickenfat5-DS}
\end{align} 
 
The only terms with over-differentiated geometry which remains to be discussed are $I^{\ZZ_n}_{8}$ (see~\eqref{I8-Zbn-P}), and $I^{\AA_n}_{5}$ (see~\eqref{I5-Abn-P}). A close inspection of the analysis of these terms Sections~\ref{sec:IZ8} and~\ref{sec:IA5}, respectively, reveals that the bounds~\eqref{eq:duck:a:duck:all:again} and \eqref{eq:insane:in:the:membrane} remain unchanged. 

\subsubsection{Downstream modifications to the forcing and commutator terms}
\label{subsubsec:chicken2}
A close inspection of the analysis of the forcing, remainder, and commutator terms from Section~\ref{sec:11:forcing:comm}, leading to the bounds for the integrals  $ I^{\WW_n}_6$, $ I^{\ZZ_n}_{10}$, and $ I^{\AA_n}_{10}$), shows that no modification is required, and that the bounds~\eqref{eq:I:W:nn:6:final},~\eqref{I10-Zn}, and~\eqref{I10-An} hold as is (with the weight $\mathcal{J}$ being replaced by $\JJ$).

\subsubsection{Conclusion of the downstream normal component estimates}
We collect the bounds for the integrals that required downstream modifications.  Combining 
\eqref{D5-L2-P}, with the downstream modified bound~\eqref{eq:I:n:12369-P}, 
with the unmodified bounds discussed in Sections~\ref{subsubsec:chicken0}, \ref{subsubsec:chicken1}, \ref{subsubsec:chicken2}, exactly as in \eqref{eq:normal:conclusion:1} we conclude that 
\begin{align}
0
&\geq 
\bigl( \tfrac{1}{52} - \Cn \eps\bigr) 
\snorm{\tfrac{\JJ^{\frac 34} (\Jg \Q)^{\frac 12}}{\Sigma^\beta} \nbs^6(\Jg \Wbn,\Jg \Zbn,\Jg \Abn)(\cdot,\s)}_{L^2_x}^2
- \tfrac{1+\alpha}{\eps} \bigl(2 +\tfrac{132}{(1+\alpha)^4}  \bigr) (\tfrac{3}{\kappa_0})^{2\beta}  \Cdatatwo
- \Cn (\tfrac{4}{\kappa_0})^{2\beta}  \mathsf{K}^2 \Bsix^2
\notag\\
&\qquad
+\Bigl(\tfrac{1+\alpha}{24} - \Cn \eps \beta\Bigr) \tfrac{1}{\eps} 
\int_0^{\s} \snorm{\tfrac{\JJ^{\frac 14} \Jgh}{\Sigma^\beta}   \nbs^6(\Jg\Wbn,\Jg\Zbn,\Jg\Abn)(\cdot,\s')}_{L^2_x}^2 {\rm d}\s'
\notag\\
&\qquad 
+ \Bigl(  \tfrac{\alpha( \beta - \frac 12) - 40\alpha \kappa_0}{8}   + \bubu{  \tfrac{9\alpha }{10 }} -  (16+ 25^2 ) (1+ \alpha )  \Bigr) \tfrac{1}{\eps}
\int_0^{\s}  \snorm{\tfrac{\JJ^{\frac 34} }{\Sigma^\beta} \nbs^6(\Jg\Zbn,\Jg\Abn) (\cdot,\s')}_{L^2_x}^2
{\rm d} \s'
\notag\\
&\qquad 
- \Bigl( \tfrac{16 \alpha  (\beta -\frac 12) }{(1+\alpha)}  +  \bubu{ 33} + 2 \cdot 250^2 +   39^2 + 500^2\Bigr) \tfrac{1}{\eps} \int_0^{\s} 
\snorm{\tfrac{\JJ^{\frac 34}(\Jg \Q)^{\frac 12}}{\Sigma^\beta} \nbs^6(\Jg\Wbn,\Jg \Zbn,\Jg\Abn)(\cdot,\s')}_{L^2_x}^2
{\rm d} \s'
\notag\\
&\qquad 
+ \tfrac{1}{20(1+\alpha)} \tfrac{1}{\eps^2} \snorm{\tfrac{\Q \JJ^{\frac 14}}{\Sigma^\beta} \nbs^6 \Jg(\cdot,\s)}_{L^2_x}^2 
-  \tfrac{20^2 + 500^2 + 100 \cdot 250^2}{1+\alpha}  {\tfrac{1}{\eps^3}} \int_0^{\s}   \snorm{\tfrac{\Q \JJ^{\frac 14}}{\Sigma^\beta} \nbs^6 \Jg(\cdot,\s')}_{L^2_x}^2 {\rm d} \s'
\notag\\
&\qquad
+ \tfrac{7}{40(1+\alpha)} \tfrac{1}{\eps^3} \int_0^{\s} \snorm{\tfrac{\Q \JJ^{-\frac 14}}{\Sigma^\beta} \nbs^6 \Jg(\cdot,\s')}_{L^2_x}^2 {\rm d} \s'
+ \bigl(8 - \Cn \brak{\beta} \eps \bigr)  {\tfrac{1}{\eps^3}}  \tint \tfrac{\Q}{\Sigma^{2\beta}} \JJ^{\!\frac 12}  \sabs{ \nbs^6\Jg}^2
\,,
\label{eq:normal:conclusion:1-P}
\end{align}
where $\Cn = \Cn(\alpha,\kappa_0,\Cdata)$ is independent of $\beta$  and $\eps$.
The bound is nearly identical to~\eqref{eq:normal:conclusion:1}. Besides a factor of $12(1+\alpha)$ which has now become $17(1+\alpha)$, the only other modification comes in the coefficient for the damping term on $\nbs^6(\Jg \Zbn,\Jg \Abn)$~(see the third line on the right side of~\eqref{eq:normal:conclusion:1-P}): instead of $\alpha(\beta-\frac 12)$, the modified bound contains $\alpha (\beta-\frac 12) - 40 \alpha \kappa_0$. 
These modifications have as a consequence the following choice for the parameter $\beta$:
\begin{align}
\bubu{ \beta_{\alpha,\kappa_0}  := 40 \kappa_0  + \tfrac{8(1+\alpha)}{\alpha} \bigl( 16 +  25^2    \bigr)- \tfrac{67 }{10} } \,.
\label{eq:normal:bounds:beta-P}
\end{align}
When compared to~\eqref{eq:normal:bounds:beta}, we notice the linear factor $40\kappa_0$ in the definition of $\beta = \beta_{\alpha,\kappa_0}$. With the choice~\eqref{eq:normal:bounds:beta-P}, we may now return to~\eqref{eq:normal:conclusion:1-P}, choose $\eps$ to be sufficiently small in terms of $\alpha,\kappa_0,\Cdata$, so that in analogy with \eqref{eq:normal:conclusion:2} we deduce
\begin{align}
& \tfrac{1}{53}  \snorm{\tfrac{\JJ^{\frac 34} (\Jg \Q)^{\frac 12}}{\Sigma^{\beta_{\alpha,\kappa_0}}} \nbs^6(\Jg \Wbn,\Jg \Zbn,\Jg \Abn)(\cdot,\s)}_{L^2_x}^2
 + \tfrac{1}{20(1+\alpha)} \tfrac{1}{\eps^2} \snorm{\tfrac{\Q \mathcal{J}^{\frac 14}}{\Sigma^{\beta_{\alpha,\kappa_0}}} \nbs^6 \Jg(\cdot,\s)}_{L^2_x}^2 
 \notag\\
&\qquad
+ \tfrac{1+\alpha}{48 \eps} 
\int_0^{\s} \snorm{\tfrac{\JJ^{\frac 14} \Jg^{\frac12}}{\Sigma^{\beta_{\alpha,\kappa_0}}}   \nbs^6(\Jg\Wbn,\Jg\Zbn,\Jg\Abn)(\cdot,\s')}_{L^2_x}^2 {\rm d}\s'
+ \tfrac{7}{40(1+\alpha)} \tfrac{1}{\eps^3} \int_0^{\s} \snorm{\tfrac{\Q \JJ^{-\frac 14}}{\Sigma^{\beta_{\alpha,\kappa_0}}} \nbs^6 \Jg(\cdot,\s')}_{L^2_x}^2 {\rm d} \s'
\notag\\
&\leq
\tfrac{1+\alpha}{\eps} \bigl(2 +\tfrac{132}{(1+\alpha)^4}  \bigr) (\tfrac{3}{\kappa_0})^{2\beta_{\alpha,\kappa_0}}  \Cdatatwo
+
\Cn 
(\tfrac{4}{\kappa_0})^{2\beta_{\alpha,\kappa_0}}  
\mathsf{K}^2 \Bsix^2
\notag\\
 &\qquad
+ \tfrac{C \brak{\alpha \kappa_0}}{\eps} \int_0^{\s} 
\snorm{\tfrac{\JJ^{\frac 34}(\Jg \Q)^{\frac 12}}{\Sigma^{\beta_{\alpha,\kappa_0}}} \nbs^6(\Jg\Wbn,\Jg \Zbn,\Jg\Abn)(\cdot,\s')}_{L^2_x}^2
{\rm d} \s'
+  \tfrac{C}{\eps^3} \int_0^{\s}   \snorm{\tfrac{\Q \JJ^{\frac 14}}{\Sigma^{\beta_{\alpha,\kappa_0}}} \nbs^6 \Jg(\cdot,\s')}_{L^2_x}^2 {\rm d} \s'
\,,
\label{eq:normal:conclusion:2-P}
\end{align}
where $C$ is a universal constant (in particular, independent of $\alpha,\kappa_0,\Cdata$), and $\Cn = \Cn(\alpha,\kappa_0,\Cdata)$ is as usual.
Comparing to \eqref{eq:normal:conclusion:2}, the main modification (besides $\beta_\alpha \mapsto \beta_{\alpha,\kappa_0}$) arising in \eqref{eq:normal:conclusion:2-P} is the dependence on $\brak{\alpha\kappa_0}$ of the Gr\"onwall constant appearing in the third term on the right side of \eqref{eq:normal:conclusion:2-P} (in~\eqref{eq:normal:conclusion:2} this constant was universal). Nonetheless, we may apply Gr\"onwall's inequality for $\s\in [0,\eps]$ to \eqref{eq:normal:conclusion:2-P} and deduce that there exists an explicitly computable constant 
\begin{equation}
 \check{\mathsf{c}}_{\alpha,\kappa_0} > 0 
\end{equation}
which only depends on $\alpha$ and $\kappa_0$, such that 
after multiplying by $\kappa_0^{\beta_{\alpha,\kappa_0}}$ and using that $\frac{\kappa_0}{4} \leq \Sigma \leq \kappa_0$, and using that $\mathsf{K} = \mathsf{K}(\alpha,\kappa_0)$ was already fixed by the tangential energy estimates (see~\eqref{eq:K:choice:1-P}), as in~\eqref{eq:normal:conclusion:3} we obtain
\begin{align}
&
\sup_{\s\in[0,\eps]} \snorm{\JJ^{\frac 34} \Jg^{\! \frac 12} \nbs^6(\Jg \Wbn,\Jg \Zbn,\Jg \Abn)(\cdot,\s)}_{L^2_x}^2
+ \tfrac{1}{\eps} 
\int_0^{\eps} \snorm{ \JJ^{\frac 14} \Jg^{\! \frac12} \nbs^6(\Jg\Wbn,\Jg\Zbn,\Jg\Abn)(\cdot,\s)}_{L^2_x}^2 {\rm d}\s
\notag\\
&\qquad 
+ 
\sup_{\s\in[0,\eps]}   \tfrac{1}{\eps^2} \snorm{\mathcal{J}^{\frac 14} \nbs^6 \Jg(\cdot,\s)}_{L^2_x}^2 
+\tfrac{1}{\eps^3} \int_0^{\eps} \snorm{\JJ^{-\frac 14} \nbs^6 \Jg(\cdot,\s)}_{L^2_x}^2 {\rm d} \s
\notag\\
&\leq \tfrac{1}{\eps} \check{\mathsf{c}}_{\alpha,\kappa_0} \Bigl(
 \Cdatatwo
+
\Cn \eps  \Bsix^2
\Bigr)
\,.
\label{eq:normal:conclusion:3-P}
\end{align}
Dropping the energy and damping terms for $\nbs^6 \Jg$ (since these were bounded already in Proposition~\ref{prop:geometry-P}), and recalling the definitions of $\widetilde{\mathcal{E}}_{6,\nnn}^2(\s)$ and $\widetilde{\mathcal{D}}_{6,\nnn}^2(\s)$ (in~\eqref{eq:tilde:E5E6}--\eqref{eq:tilde:D5D6}), as in~\eqref{eq:normal:conclusion:4}--\eqref{eq:normal:conclusion:5} we deduce that
\begin{equation}
\eps \sup_{\s \in [0,\eps]} \widetilde{\mathcal{E}}_{6,\nnn}^2(\s)
+\widetilde{\mathcal{D}}_{6,\nnn}^2(\eps) 
\leq \check{\mathsf{c}}_{\alpha,\kappa_0} 
 \Big( \Cdatatwo  + \Cn \eps  \brak{\mathsf{B}_6}^2\Bigr)
\leq 2 \check{\mathsf{c}}_{\alpha,\kappa_0} \Cdatatwo
\leq \tfrac{1}{8}   \mathsf{B}_6^2 
 \,,
 \label{eq:normal:conclusion:5-P}
\end{equation}
once $\eps$ is taken sufficiently small in terms of $\alpha,\kappa_0$, and $\Cdata$, and $\mathsf{B}_6$ is chosen sufficiently large in terms of $\alpha,\kappa_0$, and $\Cdata$ to ensure that $\mathsf{B}_6 \geq \max\{1,\Cdata\}$ (see~\eqref{eq:B6:choice:1-P}) and 
\begin{equation}
\mathsf{B}_6 \geq 4 \check{\mathsf{c}}_{\alpha,\kappa_0}^{\frac 12} \Cdata
\,.
\label{eq:B6:choice:2-P}
\end{equation}
The choice \eqref{eq:B6:choice:2-P} is the downstream-modified version of \eqref{eq:B6:choice:2}, and this closes the proof of ``normal part'' of the remaining bootstrap \eqref{boots-P} for $\mathcal{E}_{6}$ and $\mathcal{D}_{6}$. 

\subsection{Closing the bootstrap for the sixth order energy}
\label{sec:bootstrap:closed:DS}
Combining~\eqref{eq:normal:conclusion:5-P}  with \eqref{eq:hate:13-P}
we arrive at the same inequality as obtained in
\eqref{eq:normfinal1}
\begin{align}
\eps \sup_{\s \in [0,\eps]} \widetilde{\mathcal{E}}_{6}(\s)
+\widetilde{\mathcal{D}}_{6}(\eps) 
\leq \tfrac{1}{2}   \mathsf{B}_6 
\,,  \label{eq:normfinal2}
\end{align}
which closes the bootstrap \eqref{boots-P}  (cf.~\eqref{bootstraps-Dnorm:6}) in the downstream coordinate system  \eqref{eq:t-to-s-transform:all-P}.


\def\tThd{ { {\tilde \Theta}^{\dl}} }
\def\tstar{{t^*(x_2)}}
\def\tht{ { \theta(x_2,t)} }
\def\ths{ {\thd(x_2,\s)} }
\def\thsd{ { \thd(x_2,\s)} }
\def\Thsd{ { \Thd(x_1,x_2)} }
\def\dint{\iint^\thd\!\!}
\def\tint{\int_{\!\sin}^{\sfin} \!\!\! \!\iint^\thd \!\! }
\def\tints{\int_{\!\sin}^{\s} \!\!\iint^\thd \!\! }
\def\Qb{ {\bar{\mathsf{Q}}_2 }}
\def\omd{ (1\!-\! \dl)}
\def\aa{ {\mathsf a}_\alpha }

\def\tin{{\mathsf{t_{in}}}}
\def\tfin{{\mathsf{t_{fin}}}}
\def\ftime{ \mathsf{f} }
\def\thsin{ \thd_{\!\!\!\sin} }
\def\JJr{ \JJ^r}
\def\JJrr{ \JJ^{2r}}
\def\QQ{\mathcal{Q}}

\def\b#1{{\mathsf{b}_{#1} }}
\newenvironment{frcseries}{\fontfamily{frc}\selectfont}{}
\newcommand{\mathfrc}[1]{\text{\textfrc{#1}}}

\section{Upstream maximal globally hyperbolic development in a box}
\label{sec:upstreammaxdev}

\subsection{The slow acoustic characteristic surface} 
Upstream of the pre-shock, the \MGHDB of the Cauchy data is limited by the unique slow acoustic characteristic surface passing
through the co-dimension $2$ surface of pre-shocks $\Xi^*$.   With respect to our fast-characteristic-geometry, the slow acoustic characteristic flow map  $\Upupsilon=(\Upupsilon_1, \Upupsilon_2)$ evolves according to
\begin{align} 
\p_t \Upupsilon_1(x,t) & = -2 \alpha ( \Sigma \Jgi) ( \Upupsilon(x,t),t)) \,, \quad
\p_t \Upupsilon_2(x,t)  = ( V + 2\alpha \Sigma g^{-\frac 12} h,_2) (\Upupsilon(x,t),t) \,.
\label{1-flow}
\end{align} 
As can be seen from equations \eqref{eq:Zb:nn:alt} and \eqref{eq:Zb:tt:alt}, the vector
$ (-2 \alpha  \Sigma \Jgi\,, V + 2\alpha \Sigma g^{-\frac 12} h,_2 )$ is the slow acoustic transport velocity associated to
the wave speed $\lambda_1$ in \eqref{wave-speeds}, but written in the frame of the fast acoustic geometry.
\begin{figure}[htp]
\centering
\includegraphics[width=.4\textwidth]{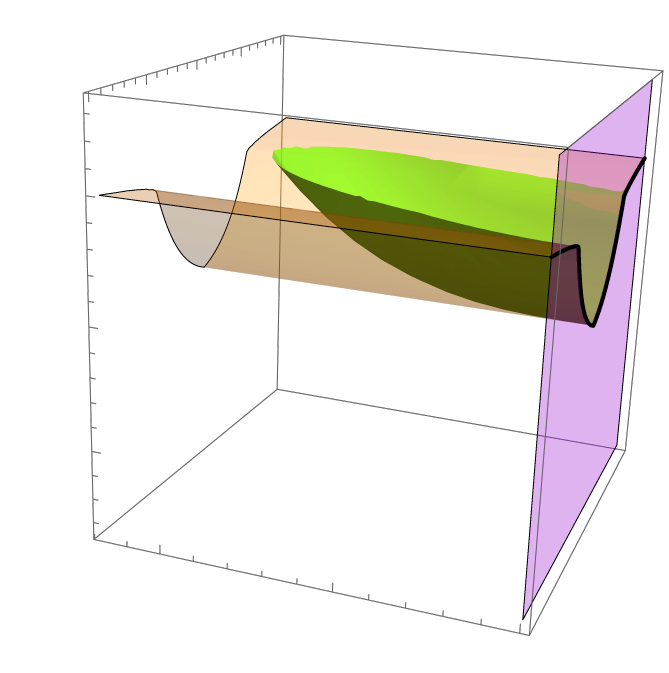};
\vspace{-.1 in}
\caption{The bounding box represents the zoomed-in region $-\tfrac{\pi \eps}{4} \leq x_1  \leq \pi\eps $, $|x_2| \leq \tfrac{1}{16}$, and $t\in [-\final,\final]$. In orange we have plotted the level set $\{\mathcal{J}(x_2,t)= \min_{x_1} \Jgb(x_1,x_2,t) =0 \}$, in magenta the surface $\{x_1 = \xringt\}$, in black the set of pre-shocks $\Xi^*$, and in green we have represented the ``upstream part'' of the slow acoustic characteristic surface passing
through $\Xi^*$, namely the set $\{(\Upupsilon(x,t),t) \colon \Upupsilon_1(x,t) \leq \mathring{x}_1(\Upupsilon_2(x,t)), t\leq \final\}$. We emphasize that the distinguished slow characteristic surface passing through $\Xi^*$ lies ``above'' the orange surface in the upstream side, so that the \MGHDB\ spacetime  is indeed an extension of the spacetime $\{\mathcal{J}>0\}$. }
\label{fig:upstream:maxdev}
\end{figure}

Recalling Definition \ref{def:pre-shock},  the set of pre-shocks is given by
\begin{equation*}
\Xi^* := \Bigl\{ \bigl(\xringt ,x_2,  t^*(x_2) \bigr) \colon x_2\in\TT\Bigr\} \,.
\end{equation*}
We shall at first be concerned with the specific slow characteristic surface that {\it passes though the pre-shock} $\Xi^*$.   This surface
consists of the union of trajectories of $\Upupsilon$ with starting position along $\Xi^*$.  The variable $t$ in \eqref{1-flow} denotes
the flow time of the dynamical system, with initial time $t=\tstar$ corresponding to the time of intersection with $\Xi^*$.  As such, the image of the
pre-shock $\Xi^*$ by the flow of the slow acoustic characteristic is given by
the green surface in Figure~\ref{fig:upstream:maxdev}, namely $\Upupsilon(\Xi^*, t)$ for $\tstar \leq t \leq \final$. Here, as in the previous sections, we use the notation
\begin{subequations}
\begin{equation}
\initial := - \tfrac{2}{1+\alpha} \eps 
\,,
\label{tin-H}
\end{equation}
for the initial time, and
\begin{equation}
\final:= \tfrac{2}{1+\alpha} \eps \cdot  \tfrac{1}{50}
\,,
\label{tfin-H}
\end{equation}
\end{subequations}
for the final time.

\subsection{Parameterizing the slow acoustic characteristic through the pre-shock as a graph}
\label{sec:upstream:reparametrize}
It is important in our analysis to parameterize this {\it distinguished}  slow-acoustic 
characteristic surface passing through $\Xi^*$
as the graph $x_1=\theta(x_2,t)$ over the $(x_2,t)$-plane. 
The dynamics of $\theta(x_2,t)$ are determined from the dynamics of the flow map $\Upupsilon$ as
\begin{equation*} 
\p_t \theta = \big(\p_t \Upupsilon \circ  \Upupsilon^{-1} \circ \theta \big) \cdot \tilde N \,, 
\end{equation*} 
where the normal vector $\tilde N$ to the surface $(\theta(x_2,t),x_2)$ is given by $\tilde N=(1,- \theta,_2)$.   As such, we compute that
\begin{equation} 
\p_t \theta
= -2 \alpha \big(\Sigma \Jgi \big) \cir \theta  - \big(V+ 2\alpha \Sigma   g^{- {\frac{1}{2}} } h,_2\big) \cir \theta \ \theta,_2 \,. 
\label{theta-dynamics-t}
\end{equation} 
The surface $x_1=\theta(x_2,t)$ is a graph-type reparamaterization
of the distinguished slow acoustic characteristic surface passing through the pre-shock, and hence must 
verify the constraint
\begin{align} 
\theta(x_2, \tstar) = \xringt \,.  \label{theta-constraint}
\end{align} 
The graph $x_1=\theta(x_2,t)$ will play the role of the ``right spatial boundary'' for our spacetime.

The characteristic surface given by the graph $x_1=\theta(x_2,t)$ can be alternatively
reparameterized as the graph $t=\Theta(x_1,x_2)$, where for each $x_2$ fixed, $\Theta$ is the inverse of $\theta$, i.e.
\begin{equation} 
\Theta(\theta(x_2,t),x_2)=t \,. \label{BigTheta-t}
\end{equation} 
By differentiating the identity \eqref{BigTheta-t} with  $\p_t$, applying the chain-rule, and using that $x_1=\theta(x_2,t)$, we have that
$$
\p_1\Theta(x_1,x_2) \p_t\theta(x_2, t) =1 \,.
$$
Substituting the dynamics \eqref{theta-dynamics-t} into this relation, we obtain that
\begin{align} 
\Big(2 \alpha \big(\Sigma \Jgi\big) (\theta(x_2,t),x_2,t) 
+ \big(V+ 2 \alpha \Sigma g^{- {\frac{1}{2}} } \p_2 h\big)(\theta(x_2,t),x_2,t)\p_2\theta(x_2,t)\Big) \p_1\Theta(x_1,x_2)=-1 \,.
\label{Elvis-in-the-house-1-t}
\end{align} 
Next, we differentiating the identity \eqref{BigTheta-t} with  $\p_2$ and again use that $x_1=\theta(x_2,t)$; this yields
the identity
\begin{align} 
\p_1\Theta(x_1,x_2)\p_2\theta(x_2,\s) = -\p_2\Theta(x_1,x_2) \,.
\label{Elvis-in-the-house-2-t}
\end{align} 
 Substitution of \eqref{Elvis-in-the-house-2-t} into
\eqref{Elvis-in-the-house-1-t} together with the fact that $t=\Theta(x_1,x_2)$ then shows that
\begin{align} 
\p_1\Theta(x) &= - \tfrac{\Jg}{2\alpha \Sigma} (x_1,x_2,\Theta(x))
+   \big(\tfrac{\Jg V}{2\alpha \Sigma}  + g^{- {\frac{1}{2}} } \p_2h \Jg \big)(x,\Theta(x))    \p_2\Theta(x) \,.
\label{p1-Theta-t-old-old}
\end{align}   
From \eqref{theta-constraint}, the reparameterization $\Theta$ verifies  the boundary condition
\begin{equation*} 
\Theta(\xringt,x_2) = \tstar \,.
\end{equation*} 

There is an important observation  to make about the ``evolution equation'' \eqref{p1-Theta-t-old-old} for the slow acoustic characteristic 
surface $\Theta$ passing through the pre-shock.   Note that with the inversion of $\theta$, the parameterization $\Theta$ defines
the evolution of the slow acoustic characteristic surface via $\Jg$ rather than $\Jgi$ (as was the case in \eqref{theta-dynamics-t}). Avoiding the
degeneracy of $\Jgi$ at the pre-shock maintains our smooth analysis.   This distinguished slow characteristic surface denotes the {\it future 
temporal boundary} of the spacetime for the upstream \MGHDB.  For technical reasons, it is convenient to
use an arbitrarily small perturbation of this characteristic surface, and we shall explain this approximation in what follows.

\subsection{A foliation of spacetime by a family of approximate $1$-characteristic surfaces}  
\label{sec:Thd:defn}
Fix an arbitrary 
\begin{equation}
\dl \in (0, \tfrac{1}{2})\,.
\label{def-dl} 
\end{equation}
We define a family of approximate $1$-characteristic surfaces $\Thd(x_1,x_2,t)$ (they would not be ``approximate'' if $\dl=0$)
  as follows. For $\initial\leq t \leq \tstar$, we define $\Thd$ as the solution of the Cauchy problem
\begin{subequations} 
\label{eq:Thd:PDE}
\begin{align} 
& \p_1\Thd(x,t) = - \tfrac{\omd \Jg}{2\alpha \Sigma} (x,\Thd(x,t))
\notag\\
&\qquad \qquad \qquad \times
\Bigl(
1 
- \bigl(V + \tfrac{2\alpha}{\omd} \Sigma g^{-\frac 12} h,_2 \bigr)(x,\Thd(x,t)) 
\bigl( \p_2 \Thd(x,t) - \tfrac{\p_2 \mathcal{B}}{\p_t \mathcal{B}}(x_2,t) \p_t \Thd(x,t)\bigr)
\Bigr) \,,
\label{p1-Theta-t} \\
&\Thd(\xstart,x_2,t)=  t \,, \qquad  \text{ for each } t \in [\initial,\tstar] \,,
 \label{Theta-BC-t}
\end{align}  
where
\begin{equation}
\mathcal{B}(x_2,t) = \Jgb(\xstart,x_2,t) \,.  
\label{Theta-nastyB-t}
\end{equation}
\end{subequations} 
The solution $\Thd$ of~\eqref{eq:Thd:PDE} is defined (see Section~\ref{sec:Thd:props} below) on the domain
\begin{subequations}
\begin{equation}
\mathring{\Omega}_{\mathsf{US},+} := \left\{ (x,t) \colon x_2 \in \TT\,, \initial \leq t < \tstar\,, \mathfrak{X}_1^-(x_2,t) \leq x_1 \leq \mathfrak{X}_1^+(x_2,t) \right\} \,,
\label{eq:Omega:US:+}
\end{equation}
where the ``stopping-times'' $\mathfrak{X}_1^-$ and $\mathfrak{X}_1^+$ are defined by
\begin{align}
 \mathfrak{X}_1^+ &= \mathfrak{X}_1^+(x_2,t) = \max\left\{ x_1 \in \TT \colon \Thd(x_1,x_2,t) \geq \initial \right\}\,,
 \label{eq:X:plus:stopping:time}
 \\
  \mathfrak{X}_1^- &= \mathfrak{X}_1^-(x_2,t) = \min\left\{ x_1 \in \TT \colon \Thd(x_1,x_2,t) \leq \final \right\}\,.
   \label{eq:X:minus:stopping:time}
\end{align}
\end{subequations}
These stopping times are well defined by continuity of the function $\Thd$ and the compactness of the constraints. Note that $\mathfrak{X}_1^-(x_2,t) \leq \xringt \leq \mathfrak{X}_1^+(x_2,t)$ in light of \eqref{Theta-BC-t}. 

\begin{remark}[\bf Spatial support]
We note that due to the bootstrap~\eqref{bs-supp} present in \eqref{boots-HH}, throughout this section are only interested in points $x\in \mathcal{X}_{\rm fin} 
=  \{ x\in \mathbb{T}^2 \colon {\rm dist}(x, \mathcal{X}_{\rm in}) \leq \Csupp \eps\}$, which in view of \eqref{eq:ic:supp} and \eqref{eq:Csupp:def}, amounts to 
\begin{equation}
|x_1 - \xringt|\leq 2(13\pi + 65\alpha(1+\alpha)\kappa_0)\eps
\,.
\label{eq:driftin:3} 
\end{equation}
This is because for  $x \not \in \mathcal{X}_{\rm fin}$, by~\eqref{eq:ic:supp} we have $\Jg = 1$, $\Sigma = \tfrac 12 \kappa_0$, and $(\Wb,\Zb,\Ab,\nb\Jg,h,_2) \equiv 0$, and so there is no analysis required here (all functions are in fact constants there). Throughout this section we shall implicitly assume that \eqref{eq:driftin:3} holds.
\end{remark}

The distinguished  surface passing through the pre-shock  (which corresponds to $t=\tstar$) is parametrized as $\{ (x_1,x_2,\bar\Thd(x_1,x_2))\}$, where 
\begin{equation}
\bar\Thd(x_1,x_2) := \Thd(x_1,x_2,\tstar) \,.
\label{eq:Thd:def}
\end{equation}
That is, for this distinguished surface, the explicit dependence on time is dropped. Here $x_2 \in \TT$ and $\mathfrak{X}_1^-(x_2,\tstar) \leq x_1 \leq \mathfrak{X}_1^+(x_2,\tstar)$.

Throughout this section we work on the $\dl$-adjusted upstream spacetime 
\begin{equation}
\tHdm := \bigl\{ (x, t) \in  \mathcal{X}_{\rm fin} \times [\initial , \final) \colon  
 \initial \leq t <    \bar \Thd(x)   \bigr\} \,.
 \label{eq:spacetime-Theta-t}
\end{equation}
The surface
$\{t =  \min\{\bar\Thd(x),\final\} \}$ defines the ``top'' temporal boundary of $\tHdm$.  
For times $\initial \le t \le \tstar$, we
have a well-defined foliation (see Section~\ref{sec:Thd:props} for details) of the spacetime subset of $\Hdm$ given by
\begin{subequations}
\label{eq:H:dl:max:split}
\begin{equation} 
\tHdmp:= \{ (x,t) \in \Hdm \colon \Theta^\dl(x,\initial) < t <   \bar \Thd(x)   \}  
= \{ (x,\Thd(x,t)) \colon (x,t) \in \mathring{\Omega}_{\mathsf{US},+} \} \,,
\label{eq:H:dl:max:+:def}
\end{equation} 
where $\mathring{\Omega}_{\mathsf{US},+}$ is defined in~\eqref{eq:Omega:US:+}.
We define the complimentary set by 
\begin{equation} 
\tHdmm := \{ (x,t) \in \Hdm \colon \initial \leq t < \Theta^\dl(x,\initial)\}  \,.
\label{eq:H:dl:max:-:def}
\end{equation} 
so that 
\begin{equation}
\tHdm = \tHdmp \cup \Theta^\dl(x,\initial)\cup \tHdmm 
\,.
\label{eq:H:dl:max:split:def}
\end{equation}
\end{subequations}
The decomposition~\eqref{eq:H:dl:max:split:def} is represented in Figure~\ref{fig:upstream:Hdm} below.

\begin{figure}[htp]
\centering
\includegraphics[width=3.5 in]{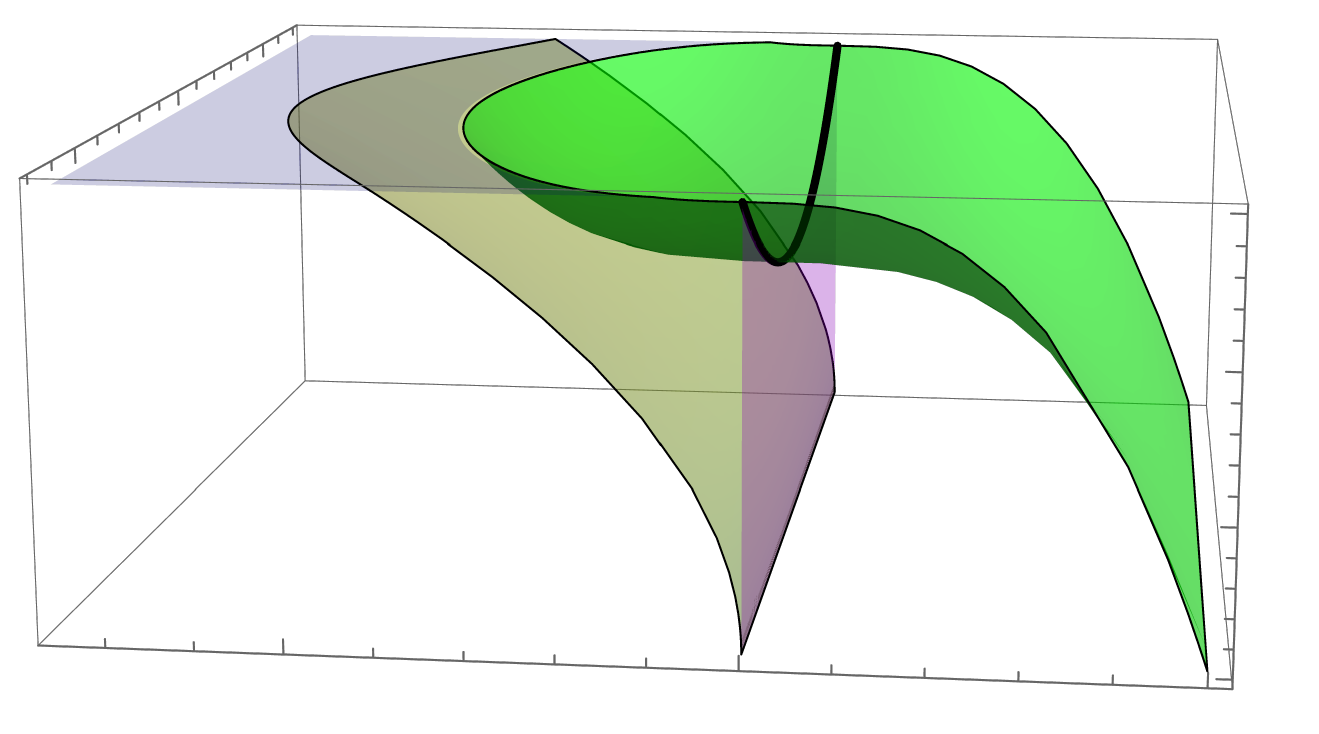};
\vspace{-.1 in}
\caption{We revisit Figure~\ref{fig:upstream:maxdev}, with the aim of describing the spacetimes $\tHdm$, $\tHdmp$, and $\tHdmm$. To retain a point of comparison with Figure~\ref{fig:upstream:maxdev}, we have kept in magenta we the surface $\{x_1 = \xringt, t \leq \tstar\}$, and have represented in black the set of pre-shocks $\Xi^*$. In green we have represented the portion of the surface $\{t = \bar{\Thd}(x) = \Thd(x,\tstar) \}$ which lies below the final times slice $\{t = \final\}$. In light-blue we have represented the portion of the final time slice which lies to the left of the green surface, i.e.~the set $\{(x,\final) \colon \bar{\Thd}(x)>\final\}$. In olive-green we have represented the portion of the surface $\{t = \Thd(x,\initial)\}$ which lies below the final times slice $\{t = \final\}$. 
The ``top'' temporal boundary of the spacetime $\tHdm$, as defined in~\eqref{eq:spacetime-Theta-t}, is thus the union of the green surface $\{t = \bar{\Thd}(x)\}$ and of the light-blue lid. The spacetime $\tHdmp$, as defined in~\eqref{eq:H:dl:max:+:def}, is the region  to the left of the green surface $\{t = \bar{\Thd}(x)\}$, to the right of the olive-green surface $\{t =  \Thd(x,\initial)\}$, and below the lid at $\{t=\final\}$. Finally, the spacetime $\tHdmm$, as defined by~\eqref{eq:H:dl:max:-:def} is the region to the left of the olive-green surface $\{t =  \Thd(x,\initial)\}$, and below the lid at $\{t=\final\}$.}
\label{fig:upstream:Hdm}
\end{figure} 

\subsection{Bounds for derivatives of $\tstar$, $\xringt$, and $\mathcal{B}$} 
Before analyzing the properties of the solution $\Thd$ of \eqref{eq:Thd:PDE}, we need to estimate various derivatives of the functions $\tstar$ and $\xringt$ appearing in the boundary condition~\eqref{Theta-BC-t}, and of the function $\mathcal{B}$ defined in~\eqref{Theta-nastyB-t}.

We recall from Definition~\ref{def:pre-shock} that $ {\Jgb}(\xringt ,x_2,  \tstar)=0$. By employing the chain-rule, the fact that also ${\Jg},_1(\xringt ,x_2,  \tstar)=0$, and the identity~\eqref{eq:Jgb:identity:2}, we deduce that 
\begin{equation} 
\p_2\tstar
=
\tfrac{\Jg,_1(\xringt ,x_2,  \tstar ) \p_2 \xringt +  \Jg,_2(\xringt ,x_2,  \tstar )}{ \p_t {\Jgb}(\xringt ,x_2,  \tstar )}
=
\tfrac{\Jg,_2(\xringt ,x_2,  \tstar )}{ \p_t {\Jgb}(\xringt ,x_2,  \tstar )}
\,.
\label{p2-tstar-def}
\end{equation} 
While $|\Jg,_2| = |\nb_2 \Jg| \leq 4(1+\alpha)$ follows from the bootstrap assumptions (see~\eqref{bs-Jg,1}), the denominator appearing in the above identity was previously estimated \eqref{eq:Qcal:bbq:temp:3} (specialized at $t= \tstar$ so that $\xringt = x_1^*(x_2,\tstar)$), 
resulting in $\tfrac{2(1+\alpha)}{5} \leq - \eps \p_t \Jgb(\xringt,x_2,\tstar) \leq (1+\alpha) (1+ 400 {\bf 1}_{\tstar\in(\medium,\final]})$. From these bounds and~\eqref{p2-tstar-def} we deduce that
\begin{equation} 
\sabs{\p_2 \tstar} \leq 10 \eps \,. \label{p2-tstar}
\end{equation} 
Similarly, from \eqref{p2-and-pt-x1star} evaluated at $t=\tstar$, we have that
\begin{align} 
\p_2 \xringt = - \bigl(\tfrac{\p_1 \p_2  \Jg}{\p_1 \p_1 \Jg}\bigr)(\xringt,x_2,\tstar) \,, \label{p2-xring}
\end{align} 
and using \eqref{eq:Jg:Wbn:horse}, \eqref{eq:Jg:11:lower}, and the fact that $- \Cn \eps \leq \tstar \leq \final $,  we find  that
\begin{equation} 
\sabs{\p_2\xstart}
\leq 
\tfrac 54 \eps \bigr( |\nb_1^2 \nb_2 w_0(\xstart,x_2)| + \eps \Cn \mathsf{K}\brak{\mathsf{B}_6}\bigl)
\leq 2\eps \|\nb^3 w_0\|_{L^\infty} 
\leq 2\eps \Cdata
\,.   \label{p2-xstar}
\end{equation} 

Next, we turn to the second order derivatives of $\tstar$ and $\xringt$.
Differentiating \eqref{p2-tstar-def} once more, we find that
\begin{align*} 
\p_2^2 \tstar =\Big( \tfrac{(\p_2 \xring \p_1 + \p_2 + \p_2 t^* \p_t) \Jg,_2}{\p_t{\Jgb}}  - \tfrac{{\Jg},_2 (\p_2 \xring \p_1 + \p_2 + \p_2 t^* \p_t)\p_t \Jgb}{(\p_t{\Jgb})^2} 
 \Big)(\xring(x_2), x_2, t^*(x_2)) \,.
\end{align*} 
Using the assumptions on $w_0$ given in Section~\ref{cauchydata} together with the bounds \eqref{eq:Jg:Wbn:horse},  \eqref{D2-Jg-Linfty}, \eqref{p2-tstar}, and~\eqref{p2-xstar}, in a similar manner to the proof of  Lemma~\ref{lem:Q:bnds} we deduce
\begin{equation} 
\sabs{\p_2^2 \tstar} \les \eps \,. \label{p22-tstar}
\end{equation} 
Analogously, differentiating \eqref{p2-xring}, we obtain that
\begin{align*} 
\p_2^2\xring (x_2)= -\Big( \tfrac{(\p_2 \xring \p_1 + \p_2 + \p_2 t^* \p_t) \Jg,_{12}}{{\Jg},_{11}} 
- \tfrac{{\Jg},_{12}(\p_2 \xring \p_1 + \p_2 + \p_2 t^* \p_t)\Jg,_{11}}{({\Jg},_{11})^2}  \Big)(\xring(x_2), x_2, t^*(x_2)) \,.
\end{align*} 
With the assumptions on $w_0$ given in Section~\ref{cauchydata}, the bounds \eqref{eq:Jg:Wbn:horse}, \eqref{D2-Jg-Linfty}, \eqref{p2-tstar}, and~\eqref{p2-xstar}, we find that
\begin{equation} 
\sabs{\p_2^2\xstart} \les \eps \langle \mathsf{B}_6\rangle \,.   \label{p22-xstar}
\end{equation} 

Lastly, we turn to the third order derivatives of $\tstar$ and $\xringt$.
Differentiating \eqref{p2-tstar-def} two times with respect to $x_2$, we have that
\begin{align*} 
\p_2^3 \tstar  
&=\Big( \tfrac{(\p_2 \xring \p_1 + \p_2 + \p_2 t^* \p_t)^2 \Jg,_2}{\p_t{\Jgb}}  - \tfrac{2(\p_2 \xring \p_1 + \p_2 + \p_2 t^* \p_t){\Jg},_2 (\p_2 \xring \p_1 + \p_2 + \p_2 t^* \p_t)\p_t \Jgb}{(\p_t{\Jgb})^2} 
\notag\\
&\qquad \qquad 
- \tfrac{ {\Jg},_2 (\p_2 \xring \p_1 + \p_2 + \p_2 t^* \p_t)^2\p_t \Jgb}{(\p_t{\Jgb})^2} 
+ \tfrac{2 {\Jg},_2 ((\p_2 \xring \p_1 + \p_2 + \p_2 t^* \p_t) \p_t \Jgb)^2}{(\p_t{\Jgb})^3} 
 \Big)(\xring(x_2), x_2, t^*(x_2))  \,.
\end{align*} 
With the assumptions on $w_0$ given in Section~\ref{cauchydata}, the bounds \eqref{eq:Jg:Wbn:horse}, \eqref{D2-Jg-Linfty}, and \eqref{p2-xstar}--\eqref{p22-xstar}  yield \begin{equation} 
\sabs{\p_2^3 \tstar} \les    \eps \langle \mathsf{B}_6\rangle  \,. \label{p222-tstar}
\end{equation} 
Finally, differentiating \eqref{p2-xring} two times with respect to $x_2$, we arrive at
\begin{align*} 
\p_2^3 \xringt  
&=-\Big( \tfrac{(\p_2 \xring \p_1 + \p_2 + \p_2 t^* \p_t)^2 \Jg,_{12}}{\Jg,_{11}}  
- \tfrac{2(\p_2 \xring \p_1 + \p_2 + \p_2 t^* \p_t){\Jg},_{12} (\p_2 \xring \p_1 + \p_2 + \p_2 t^* \p_t)\Jg,_{11}}{(\Jg,_{11})^2} 
\notag\\
&\qquad \qquad 
- \tfrac{ {\Jg},_{12} (\p_2 \xring \p_1 + \p_2 + \p_2 t^* \p_t)^2 \Jg,_{11}}{(\Jg,_{11})^2} 
+ \tfrac{2 {\Jg},_{12} ((\p_2 \xring \p_1 + \p_2 + \p_2 t^* \p_t) \Jg,_{11})^2}{(\Jg,_{11})^3} 
 \Big)(\xring(x_2), x_2, t^*(x_2))  \,.
\end{align*} 
The noticeable difference we encounter in the above identity is the appearance of third order derivatives of $\Jg,_1$ along the pre-shock. Pointwise bounds for $\nb^4 \Jg$ are not covered by the $L^\infty$ bound in \eqref{D2-Jg-Linfty}. Instead, we   note that the bootstrap~\eqref{boots-H} postulates $\Jg \in H^6(\tHdm)$, and thus the Sobolev embedding in $2+1$ (space$+$time) dimensions gives that $\Jg \in L^\infty(\tHdm)$. The Sobolev embedding is however not necessarily sharp in terms of the scaling with respect to $\eps$. Indeed, the classical Gagliardo-Nirenberg inequality for a function $f$  which is $H^2$ smooth on a space-time domain $\Omega \subset \mathcal{X}_{\rm fin} \times [\initial,\final] \subset \mathbb{R}^3$: is $\|f\|_{L^\infty(\Omega)} \les \|f\|_{L^2(\Omega)}^{\frac 14} \|\nabla^2 f\|_{L^2(\Omega)}^{\frac 34}+ |\Omega|^{-\frac 12} \|f\|_{L^2(\Omega)}$, where the implicit constant is universal. In terms of the differential operators $\nb$, this implies via the Poincar\'e inequality~\eqref{eq:x1:Poincare} that $\|f\|_{L^\infty_{x,t}} \les \eps^{-\frac 32} \|\nb^2 f\|_{L^2_{x,t}}$. Applying this bound with $f =\nb^4 \Jg$, and appealing to the $\nb^6 \Jg$ bootstrap~\eqref{boots-H}, we deduce $\|\nb^4\Jg \|_{L^\infty_{x,t}} \les \mathsf{B_J} \eps^{-\frac 12}$. Since this term is $\OO(\eps^{-\frac 12})$ instead of $\OO(1)$, the bound for $\p_2 \xringt$ loses a factor of $\eps^{\frac 12}$, resulting in 
\begin{equation} 
\sabs{\p_2^3 \xringt} \les    \eps^{\frac 12} \langle \mathsf{B}_6\rangle  \,. \label{p222-xring}
\end{equation} 

\subsection{Solvability of \eqref{eq:Thd:PDE} and properties of $\Thd$ and its derivatives}
\label{sec:Thd:props}
Solving for $\Thd$ amounts to a standard application of the method of characteristics.
Letting 
\begin{subequations}
\label{eq:frak:MNF}
\begin{align}
\mathfrak{M}(x,t) &:= - \tfrac{(1-\dl)\Jg(x,t)}{2\alpha \Sigma(x,t)} \bigl( V + \tfrac{2\alpha}{1-\dl} \Sigma g^{-\frac 12} h,_2\bigr)(x,t) \,,
\label{eq:frak:MNF:a} \\
\mathfrak{N}(x,t) &:= \tfrac{(1-\dl)\Jg(x,t)}{2\alpha \Sigma(x,t)} \bigl( V + \tfrac{2\alpha}{1-\dl} \Sigma g^{-\frac 12} h,_2\bigr) (x,t) \tfrac{\p_2 \mathcal{B}}{\p_t\mathcal{B}}(x_2,t) \,,
\label{eq:frak:MNF:b}\\
\mathfrak{F}(x,t) &:=  - \tfrac{(1-\dl)\Jg(x,t)}{2\alpha \Sigma(x,t)}\,,
\label{eq:frak:MNF:c}
\end{align}
\end{subequations}
we may write \eqref{p1-Theta-t} as 
\begin{equation}
\label{p1-Theta-t-alt} 
\p_1 \Thd(x,t) + \mathfrak{M}(x,\Thd(x,t)) \p_2 \Thd(x,t) +  \mathfrak{N}(x,\Thd(x,t)) \p_t \Thd(x,t) = \mathfrak{F}(x, \Thd(x,t) )
\,.
\end{equation}
This is a semilinear first order PDE with smooth coefficients. The regularity of $\mathfrak{M},\mathfrak{N}$, and $\mathfrak{F}$ may be seen as follows. First, in analogy to the bound \eqref{p2-tstar} we deduce from~\eqref{eq:Jg:Wbn:horse}, \eqref{eq:Qcal:bbq:temp:3}, \eqref{eq:x1star:x1vee}, and~\eqref{p2-xstar} that the term $\mathcal{B}$ defined in~\eqref{Theta-nastyB-t} satisfies
\begin{subequations}
\label{eq:MNF:bounds}
\begin{equation}
\sabs{\tfrac{\p_2 \mathcal{B}}{\p_t \mathcal{B}}(x_2,t) } 
\leq \tfrac{|\Jg,_1(\xringt,x_2,t)| \cdot |\p_2 \xringt| + |\Jg,_2(\xringt,x_2,t)|}{-\p_t \Jgb(\xringt,x_2,t)}
\leq (\Cn \eps^2 + 4(1+\alpha)) \tfrac{5\eps}{2(1+\alpha)} \leq 11 \eps \,,
\label{eq:MNF:bounds:a}
\end{equation}
for all $t\in [\initial,\tstar)$. Using the above estimate and the pointwise bootstrap assumptions~\eqref{bs-Jg-simple}--\eqref{bs-D-Sigma} present in \eqref{boots-HH}, we deduce that
\begin{equation}
\sabs{\mathfrak{F}(x,t)} \leq \tfrac{2}{\alpha\kappa_0} \Jg(x,t)\,,
\quad 
\sabs{\mathfrak{M}(x,t)} \leq   (16 + \tfrac{2\mathsf{C_V}}{\alpha\kappa_0}) \eps \Jg(x,t)  \,,
\quad 
\sabs{\mathfrak{N}(x,t)} \leq   (176 + \tfrac{22\mathsf{C_V}}{\alpha\kappa_0}) \eps^2 \Jg(x,t) \,,
\label{eq:MNF:bounds:b}
\end{equation}
which become uniform bounds since $\Jg(x,t) \leq \tfrac 65$.
In order to bound derivatives of $\mathfrak{M},\mathfrak{N}$, and $\mathfrak{F}$, we use that from~\eqref{eq:Jg:Wbn:horse}, \eqref{eq:Qcal:bbq:temp:3}, \eqref{eq:x1star:x1vee}, \eqref{p2-xstar}, and~\eqref{p22-xstar} we have
\begin{equation}
\sabs{(\nb_t,\nb_2) \tfrac{\p_2 \mathcal{B}}{\p_t \mathcal{B}}(x_2,t) } 
\les \eps \brak{\mathsf{B}_6}\,,
\label{eq:MNF:bounds:c}
\end{equation}
where we recall that $(\nb_t,\nb_2) = (\eps \p_t, \p_2)$.
Combining this estimate with  the pointwise bootstrap assumptions~\eqref{bs-Jg-simple}--\eqref{bs-D-Sigma} present in \eqref{boots-HH}, we deduce
\begin{equation}
\sabs{\nb \mathfrak{F}(x,t)} \les 1\,,
\quad 
\sabs{\nb \mathfrak{M}(x,t)} \les    \eps \,,
\quad 
\sabs{\nb \mathfrak{N}(x,t)} \les \eps^2 \brak{\mathsf{B}_6}\,,
\label{eq:MNF:bounds:d}
\end{equation}
where we recall that $\nb = (\eps \p_t, \eps \p_1, \p_2)$. Lastly, by appealing to \eqref{eq:Jg:Wbn:horse}, \eqref{p2-xstar}, \eqref{p22-xstar}, \eqref{p222-xring}, 
and the estimate, $|\Jg,_1(\xringt,x_2,t)| \leq |x_1^*(x_2,t) - \xringt| \cdot \|\Jg,_{11}\|_{L^\infty_{x,t}} \les \eps \mathsf{K} \brak{\mathsf{B}_6}$, 
similarly to~\eqref{eq:MNF:bounds:a} and \eqref{eq:MNF:bounds:c} we deduce
\begin{equation}
\sabs{(\nb_t^2,\nb_t \nb_2,\nb_2^2) \tfrac{\p_2 \mathcal{B}}{\p_t \mathcal{B}}(x_2,t) } 
\les \eps \brak{\mathsf{B}_6}\,.
\label{eq:MNF:bounds:2}
\end{equation}
Combined with the bootstraps~\eqref{boots-HH}--\eqref{boots-H}, and the  anisotropic Sobolev estimate in~\eqref{eq:Sobolev}, similarly to \eqref{eq:MNF:bounds:d} we obtain
\begin{equation}
\sabs{\nb^2 \mathfrak{F}(x,t)} \les \brak{\mathsf{B}_6}\,,
\quad 
\sabs{\nb^2 \mathfrak{M}(x,t)} \les    \eps \brak{\mathsf{B}_6}\,,
\quad 
\sabs{\nb^2 \mathfrak{N}(x,t)} \les \eps^2 \brak{\mathsf{B}_6}\,.
\label{eq:MNF:bounds:f}
\end{equation}
\end{subequations}

With the bounds in~\eqref{eq:MNF:bounds}, we turn to solving~\eqref{p1-Theta-t-alt}. Treating $x_1$ as  time  and $(x_2,t)$ as parameters, we introduce the flows $(\zeta_2(x_1,x_2,t),\zeta_t(x_1,x_2,t))$ and the function $\Thd\circ\zeta(x_1,x_2,t) := \Thd(x_1,\zeta_2(x_1,x_2,t),\zeta_t(x_1,x_2,t))$. The flows $(\zeta_2, \zeta_t) $ are the solutions of the characteristic ODEs
\begin{subequations}
\label{eq:characteristics}
\begin{align}
\p_1 \zeta_2(x_1,x_2,t) 
&= \mathfrak{M}\bigl(x_1,\zeta_2(x_1,x_2,t),\Thd\cir\zeta(x_1,x_2,t) \bigr),
\qquad \zeta_2(\xringt,x_2,t) = x_2\,,
\label{eq:characteristics:a}
\\
\p_1 \zeta_t(x_1,x_2,t) 
&= \mathfrak{N}\bigl(x_1,\zeta_2(x_1,x_2,t),\Thd\cir\zeta(x_1,x_2,t) \bigr),
\qquad \zeta_t(\xringt,x_2,t) = t\,,
\label{eq:characteristics:b}
\end{align}
while \eqref{p1-Theta-t-alt}  and \eqref{Theta-BC-t} may be rewritten as
\begin{align}
&\p_1 \bigl(\Thd\cir\zeta(x_1,x_2,t)\bigr)
= \mathfrak{F}\bigl(x_1,\zeta_2(x_1,x_2,t),\Thd\cir\zeta(x_1,x_2,t)\bigr)\,,
\label{eq:characteristics:c}
\\
&\Thd\cir\zeta(\xringt,x_2,t)
= \Thd(\xringt,x_2,t) = t\,.
\label{eq:characteristics:d}
\end{align}
\end{subequations}
Note that the boundary condition at $\{\xringt,x_2,t\}$ is non-characteristic. Moreover, the fields $(\mathfrak{F},\mathfrak{M},\mathfrak{N})$ defined in \eqref{eq:frak:MNF} are uniformly $C^1$ in both space and time in our spacetime (see~\eqref{eq:MNF:bounds}). This ensures unique and smooth solvability of the system \eqref{eq:characteristics}: first by solving the two-dimensional system of coupled ODEs for $\zeta_2(x_1,x_2,t)$ and $\Thd(x_1,\zeta_2(x_1,x_2,t),\eta_t(x_1,x_2,t))$ obtained from \eqref{eq:characteristics:a} and \eqref{eq:characteristics:c}--\eqref{eq:characteristics:d}, with $(x_2,t)$ as parameters, and then afterwards integrating the ODE for $\zeta^t$ in \eqref{eq:characteristics:b}. The global solvability of the characteristic ODEs in the interval $x_1\in[\mathfrak{X}_1^-(x_2,t),\mathfrak{X}_1^+(x_2,t)]$ is a consequence of the $C^1_{x,t}$ regularity of $(\mathfrak{F},\mathfrak{M},\mathfrak{N})$ and the fact that the boundary data at $x_1 = \xringt$ is smooth and non-characteristic. Moreover, using~\eqref{eq:MNF:bounds:b} and~\eqref{eq:MNF:bounds:d}, we have that the map $(x_2,t) \mapsto (\zeta_2(\cdot,x_2,t),\zeta_t(\cdot,x_2,t))$ is invertible and  the bounds
\begin{equation}
|\zeta_2(x_1,x_2,t) -x_2|\leq \Cn \eps^2 \,,
\qquad  
|\zeta_t(x_1,x_2,t) -t|\leq \Cn \eps^3 \,,
\end{equation} 
hold for each $x_1\in[\mathfrak{X}_1^-(x_2,t),\mathfrak{X}_1^+(x_2,t)]$. 
Next, we turn to bounding the derivatives of $\Thd$. We establish the following
\begin{lemma}[\bf Bounds for the derivatives of $\Thd$]
\label{lem:Thd:derivs}
Let $\kappa_0$ be sufficiently large with respect to $\alpha$ to ensure that  \eqref{eq:US:kappa:0:cond:0} holds.
Assume that the bootstrap assumptions~\eqref{bootstraps-H} hold in $\tHdm$ and that $\eps$ is taken sufficiently small with respect to $\alpha,\kappa_0$, and $\Cdata$. Then, for all $(x,t) \in \mathring{\Omega}_{\mathsf{US},+}$, the spacetime defined in~\eqref{eq:Omega:US:+}, we have 
\begin{subequations}
\label{eq:Thd:derivs}
\begin{align} 
-  \tfrac{4 \Jg (x,\Thd(x,t))}{\alpha \kappa_0}
&\leq \p_1\Thd(x,t) \leq -  \tfrac{\Jg (x,\Thd(x,t))}{4 \alpha \kappa_0} < 0\,, 
\label{p1-Theta-sign-t} \\
\sabs{\p_2\Thd(x,t)} 
&\le 5 \cdot 10^3 (1+\alpha)^2   \eps\,, 
\label{boot-p2Theta-t} \\
|\p_t\Thd(x,t)-1|
& \le 3 \cdot 10^{-5} \,,
\label{boot-psTheta-t} \\
\sabs{\p_{22}\Thd(x,t)} 
&\le \b{22}\, \eps \brak{\mathsf{B}_6} \,, 
\label{boot-p22Theta-t} \\
\sabs{\p_{2t}\Thd(x,t)} 
&\le \b{2\s}  \,, 
\label{boot-p2sTheta-t} \\
\sabs{\p_{tt}\Thd(x,t)} 
&\le \b{\s\s}  \,.
\label{boot-pssTheta-t}  
\end{align}
The constants $\b{22}, \b{2\s}, \b{\s\s}$ appearing in~\eqref{boot-p22Theta-t}--\eqref{boot-pssTheta-t} only depend on $\alpha$, $\kappa_0$, and $\Cdata$, and are defined  in~\eqref{eq:b22:constraint}, \eqref{eq:b2s:constraint}, and~\eqref{eq:bss:constraint}.
\end{subequations}
\end{lemma}

Before giving the proof of Lemma~\ref{lem:Thd:derivs}, we record a few immediate consequences.
First, we note that from \eqref{p1-Theta-sign-t}, \eqref{eq:US:kappa:0:cond:0}, and \eqref{bs-Jg-simple} we may deduce
\begin{subequations}
\label{eq:jesus:just:left:chicago}
\begin{equation} 
-\tfrac{1}{10^5(1+\alpha)} \leq 
-  \tfrac{5 }{\alpha \kappa_0}
\leq \p_1\Thd(x,t)  < 0\,, 
\label{p1-Theta-sign-t-new}
\end{equation}
for all $(x,t) \in \mathring{\Omega}_{\mathsf{US},+}$.
Second, we note that by differentiating \eqref{p1-Theta-t-alt} with respect to $x_2$ or $t$, and appealing to the bootstrap inequalities in~\eqref{boots-HH} bounds in~\eqref{eq:MNF:bounds:b}, \eqref{eq:MNF:bounds:d}, and also to identity~\eqref{eq:MNF:Rams:1} and bound~\eqref{eq:US:kappa:0:cond:0} below, we deduce
\begin{align}
\sabs{\p_{12} \Thd(x,t)}
&\leq |(\p_2 \mathfrak{F})\cir \Thd(x,t)| + |\p_2 \Thd \cdot (\p_t \mathfrak{F})\cir \Thd(x,t)| + \Cn \eps\brak{\mathsf{B}_6} 
\notag\\
&\leq \tfrac{32(1+\alpha)}{\alpha\kappa_0} + \tfrac{5\cdot10^3(1+\alpha)^3}{\alpha\kappa_0} + \Cn \eps\brak{\mathsf{B}_6} 
\leq \tfrac{(1+\alpha)^2}{50} 
\label{eq:jesus:just:left:chicago:a}
\,,
\\
\sabs{\p_{1t} \Thd(x,t)}
&\leq |\p_t \Thd \cdot (\p_t \mathfrak{F})\cir \Thd(x,t)| + \Cn \eps\brak{\mathsf{B}_6} 
\notag\\
&\leq  (1 + 3\cdot 10^{-5}) \tfrac{ (1+\alpha)}{\alpha \eps \kappa_0}+ \Cn \eps\brak{\mathsf{B}_6}
\leq \tfrac{1}{10^5 \eps}
\label{eq:jesus:just:left:chicago:b}
\,.
\end{align}
By using these estimates, we may also differentiate \eqref{p1-Theta-t-alt} with respect to $x_1$, and similarly deduce
\begin{align}
\sabs{\p_{11} \Thd(x,t)}
&\leq |(\p_1 \mathfrak{F})\cir \Thd(x,t)| + |\p_1 \Thd \cdot (\p_t \mathfrak{F})\cir \Thd(x,t)| + \Cn \eps\brak{\mathsf{B}_6} 
\notag\\
&\leq  \tfrac{30}{\alpha \eps \kappa_0} + 
\tfrac{5(1+\alpha)}{\eps (\alpha \kappa_0)^2}+ \Cn \eps\brak{\mathsf{B}_6}
\leq \tfrac{1}{10^4 \eps (1+\alpha)}
\label{eq:jesus:just:left:chicago:c}
\,. 
\end{align}
for all $(x,t) \in \mathring{\Omega}_{\mathsf{US},+}$.
\end{subequations}

\begin{proof}[Proof of Lemma~\ref{lem:Thd:derivs}]
First, we prove \eqref{p1-Theta-sign-t}. Recall that $\p_1\Thd$ is computed from 
\eqref{p1-Theta-t-alt}. By appealing to the bounds in~\eqref{eq:MNF:bounds:b}, along with the bounds~\eqref{boot-p2Theta-t}--\eqref{boot-psTheta-t}, we obtain
\begin{equation*}
\sabs{\p_1 \Thd(x,t) + \tfrac{(1-\dl)\Jg}{2\alpha \Sigma}(x,\Thd(x,t))}
\leq \Cn \eps^2 \Jg(x,\Thd(x,t))\,.
\end{equation*}
Combining this bound with the $\Sigma$ bootstrap in~\eqref{bs-Sigma} and taking $\eps$ to be sufficiently small, proves \eqref{p1-Theta-sign-t}.

Next, we prove~\eqref{boot-p2Theta-t}--\eqref{boot-psTheta-t}. We establish these bounds via a bootstrap/continuity argument with respect to $x_1$, starting at $\xringt$. At $x_1 =\xringt$ we have we have $|\p_2 \Thd(\xringt,x_2,t)|  = |\Jg,_1(\xringt,x_2,t)| \cdot |\p_2 \xringt| \leq \Cn \eps^2 \mathsf{K} \brak{\mathsf{B}_6}$ in light of~\eqref{p2-xstar}, and also $\p_t \Thd(\xringt,x_2,t) = 1$. Thus, at $x_1 = \xringt$, the bounds \eqref{boot-p2Theta-t}--\eqref{boot-psTheta-t} hold with a strict inequality. We then continue to propagate these bounds for $x_1$ away from $\xringt$, and prove that they still hold with a strict inequality, yielding the global bound  in $\mathring{\Omega}_{\mathsf{US},+}$.

We  differentiate~\eqref{p1-Theta-t-alt} with respect to the $t$ and $x_2$, and deduce that 
\begin{align}
&\bigl(\p_1 + \mathfrak{M}\cir\Thd \p_2 +\mathfrak{N} \cir \Thd \p_t\bigl)
\begin{pmatrix}
\p_t \Thd \\
\p_2\Thd 
\end{pmatrix} 
- 
(\p_t \mathfrak{F}) \cir \Thd
\begin{pmatrix}
\p_t \Thd \\
\p_2\Thd 
\end{pmatrix} 
\notag\\
&
= \begin{pmatrix}
- (\p_t \mathfrak{M})\cir \Thd \p_2 \Thd \p_t \Thd 
- (\p_t \mathfrak{N})\cir \Thd (\p_t \Thd)^2
\\
- (\p_t \mathfrak{M})\cir\Thd (\p_2 \Thd)^2 
- (\p_2 \mathfrak{M})\cir \Thd \p_2 \Thd 
- (\p_t \mathfrak{N})\cir \Thd \p_2 \Thd \p_t \Thd
- (\p_2 \mathfrak{N})\cir \Thd \p_t \Thd
+ (\p_2 \mathfrak{F})\cir \Thd 
\end{pmatrix} 
\label{eq:MNF:Rams:0}
.
\end{align}
Time differentiating~\eqref{eq:frak:MNF:c}, and by appealing to \eqref{Jg-evo}, \eqref{Sigma0-ALE}, we obtain \eqref{p1-Sigma}, 
\begin{align}
\p_t \mathfrak{F} 
&= - \tfrac{(1-\dl)(1+\alpha)}{2\alpha} \p_1 (\log \Sigma)
\notag\\
&\quad
- \tfrac{(1-\dl)}{2\alpha\Sigma} 
\bigl((1-\alpha) \Jg \Zbn - \alpha \Jg \Abt - \tfrac{1+\alpha}{2} \Jg h,_2 (\Wbt-\Zbt) - V \Jg,_2 - V \Jg \Sigma^{-1} \Sigma,_2 \bigr)\,,
\label{eq:Ophelia:1}
\end{align}
so that the bootstrap assumptions imply 
\begin{equation}
 \sabs{\p_t \mathfrak{F}+\tfrac{(1-\dl)(1+\alpha)}{2\alpha} \p_1 (\log \Sigma)} \leq \Cn \,.
 \label{eq:MNF:Rams:1}
\end{equation}
Using the characteristic flow $(\zeta_2,\zeta_t)$ introduced in \eqref{eq:characteristics}, we additionally note the identity
\begin{align}
\bigl(\p_1(\log \Sigma)\bigr)(x_1,\zeta_2(x_1,x_2,t),\Thd\cir\zeta(x_1,x_2,t)) 
&= 
\p_1 \bigl( \log \Sigma (x_1,\zeta_2(x_1,x_2,t),\Thd\cir\zeta(x_1,x_2,t)) \bigr)   
\notag\\
&\qquad 
- 
\bigl(\p_2(\log \Sigma) \cdot \mathfrak{M} \bigr)(x_1,\zeta_2(x_1,x_2,t),\Thd\cir\zeta(x_1,x_2,t)) 
\notag\\
&\qquad 
- 
\bigl(\p_t(\log \Sigma) \cdot \mathfrak{F} \bigr)(x_1,\zeta_2(x_1,x_2,t),\Thd\cir\zeta(x_1,x_2,t)) 
\,.
 \label{eq:MNF:Rams:2}
\end{align}
Thus, composing \eqref{eq:MNF:Rams:0} with $(\zeta_2,\zeta_t)$, using the bound~\eqref{eq:MNF:Rams:1}, identity~\eqref{eq:MNF:Rams:2}, the bounds~\eqref{eq:MNF:bounds:b}, \eqref{eq:MNF:bounds:d}, the bootstrap assumptions relating to $\Sigma$ in~\eqref{bs-Sigma}--\eqref{bs-D-Sigma}, and the bounds~\eqref{boot-p2Theta-t}--\eqref{boot-psTheta-t} in a bootstrap fashion, we deduce that
\begin{subequations}
\label{eq:MNF:Rams:3}
\begin{align}
&\bigl| \p_1 \bigl( (\p_t \Thd) \cir \zeta (x_1,x_2,t) \bigr)  
\notag\\
&\quad
+ \tfrac{(1-\dl)(1+\alpha)}{2\alpha} \p_1 \bigl( \log \Sigma (x_1,\zeta_2(x_1,x_2,t),\Thd\cir\zeta(x_1,x_2,t)) \bigr)   \cdot (\p_t \Thd) \cir \zeta (x_1,x_2,t) \bigr| \leq \Cn 
\,,
\label{eq:MNF:Rams:3a}
\\
&\bigl| \p_1 \bigl( (\p_2 \Thd) \cir \zeta (x_1,x_2,t) \bigr)  
\notag\\
&\quad
+ \tfrac{(1-\dl)(1+\alpha)}{2\alpha} \p_1 \bigl( \log \Sigma (x_1,\zeta_2(x_1,x_2,t),\Thd\cir\zeta(x_1,x_2,t)) \bigr)   \cdot (\p_2 \Thd) \cir \zeta (x_1,x_2,t)\bigr| \leq \Cn \eps + \| \p_2 \mathfrak{F}\|_{L^\infty_{x,t}}
\label{eq:MNF:Rams:3b}
\,.
\end{align}
\end{subequations}
The estimates in~\eqref{eq:MNF:Rams:3} are set up so that we introduce the integrating factor 
\begin{align}
&\mathcal{I}(x_1,x_2,t):= \exp\Bigl(- \tfrac{(1-\dl)(1+\alpha)}{2\alpha} \int_{\xringt}^{x_1}\p_1 \bigl( \log \Sigma (x_1^\prime,\zeta_2(x_1^\prime,x_2,t),\Thd\cir\zeta(x_1^\prime,x_2,t)) \bigr) {\rm d} x_1^\prime\Bigr)
\notag\\
&=\exp\Bigl(- \tfrac{(1-\dl)(1+\alpha)}{2\alpha} \log \bigl(\tfrac{\Sigma (x_1,\zeta_2(x_1,x_2,t),\Thd\circ\zeta(x_1,x_2,t))}{\Sigma (\xringt,x_2,t)} \bigr)\Bigr)
= \Bigl(\tfrac{\Sigma (\xringt,x_2,t)}{\Sigma (x_1,\zeta_2(x_1,x_2,t),\Thd\circ\zeta(x_1,x_2,t))} \Bigr)^{\tfrac{(1-\dl)(1+\alpha)}{2\alpha}}
\,.
\label{eq:MNF:Rams:4}
\end{align}
We note that since the pointwise bootstraps in~\eqref{bootstraps} and identity~\eqref{Sigma0-ALE} imply $|\p_t \Sigma| \les 1$, and thus via assumption~\eqref{item:ic:infinity} we arrive at $|\Sigma(x,t) - \tfrac 12 w_0(x)| \leq |\Sigma(x,t) - \sigma_0(x)| + \tfrac 12 \eps \kappa_0 \leq \Cn \eps$. Applying  assumption~\eqref{item:ic:infinity} once more, together with the assumption that $\kappa_0\geq 20$, we deduce 
\begin{equation}
\sabs{\tfrac{\Sigma (\xringt,x_2,t)}{\Sigma (x_1,\zeta_2(x_1,x_2,t),\Thd\circ\zeta(x_1,x_2,t))}
 - 1}
\leq 
\tfrac{2 \|w_0 -\kappa_0\|_{L^\infty_x} +\Cn \eps}{\kappa_0 - \|w_0 - \kappa_0\|_{L^\infty_x} - \Cn \eps}
\leq \tfrac{4}{\kappa_0 - 2} \leq \tfrac{5}{\kappa_0}\,.
 \label{eq:MNF:Rams:5}
\end{equation}
Using the fact that $| (1+r)^\beta - 1 |\leq 2 |r|\beta$ for $\beta>0$ and $|r|\leq \tfrac{1}{10}$, such that $\beta |r| \leq \tfrac{1}{10}$, if we take $\kappa_0$ sufficiently large with respect to $\alpha$ to ensure 
\begin{equation}
 \kappa_0 \geq \tfrac{5 \cdot 10^5(1+\alpha)}{\alpha}
 \label{eq:US:kappa:0:cond:0}
 \,,
\end{equation}
then we deduce from \eqref{eq:MNF:Rams:4} and \eqref{eq:MNF:Rams:5} that
\begin{equation}
\sabs{\mathcal{I}(x_1,x_2,t) - 1}
\leq \tfrac{5(1-\dl)(1+\alpha)}{\alpha \kappa_0} \leq 10^{-5}\,,
 \label{eq:MNF:Rams:6}
\end{equation}
uniformly for $\dl \in (0,1)$.

With the bound~\eqref{eq:MNF:Rams:6}, we return to \eqref{eq:MNF:Rams:3a} and estimate for $\p_t \Thd$. Integrating  \eqref{eq:MNF:Rams:3a} with respect to $x_1$, and using the boundary condition $(\p_t \Thd)(\xringt,x_2,t) = 1$, we deduce that
\begin{equation*}
\sabs{ (\p_t \Thd) \cir \zeta (x_1,x_2,t) - 1\cdot \mathcal{I}(x_1,x_2,t) }
\les |x_1 - \xringt| \les \eps \,,
\end{equation*}
which gives via~\eqref{eq:MNF:Rams:6} and upon composing with $\zeta^{-1}$ that
\begin{equation*}
\sabs{ \p_t \Thd -1 }
\leq 10^{-5} + \Cn \eps 
\leq 2 \cdot 10^{-5}\,,
\end{equation*}
upon taking $\eps$ to be sufficiently small. This proves~\eqref{boot-psTheta-t}.

Integrating  \eqref{eq:MNF:Rams:3b} with respect to $x_1$, using that the boundary condition satisfies $|(\p_2 \Thd)(\xringt,x_2,t)| \leq \Cn \eps^2 \mathsf{K} \brak{\mathsf{B}_6}$, appealing to the bound~\eqref{eq:MNF:Rams:6} for the integrating factor, and also using the estimate 
$\| \p_2 \mathfrak{F}\|_{L^\infty_{x,t}} \leq \tfrac{2}{\alpha\kappa_0} \|\Jg,_2\|_{L^\infty_{x,t}} + \tfrac{8}{\alpha\kappa_0^2}\|\Jg \Sigma,_2\|_{L^\infty_{x,t}}\leq \tfrac{30(1+\alpha)}{\alpha\kappa_0}$ (which is a consequence of \eqref{bs-Jg-simple}, \eqref{bs-Jg,1}, \eqref{bs-Sigma}, \eqref{bs-D-Sigma}), we deduce 
\begin{align*}
\sabs{ (\p_2 \Thd) \cir \zeta (x_1,x_2,t)}
&\leq  \mathcal{I}(x_1,x_2,t) |(\p_2 \Thd)(\xringt,x_2,t)| + \bigl(\tfrac{30(1+\alpha)}{\alpha\kappa_0} + \Cn \eps\bigr) \cdot \tfrac{1+10^{-5}}{1-10^{-5}} \cdot|x_1 - \xringt| \notag\\
&\leq  \Cn \eps^2 \mathsf{K} \brak{\mathsf{B}_6} + \tfrac{31(1+\alpha)}{\alpha\kappa_0} \cdot|x_1 - \xringt| \,.
\end{align*}
Appealing to the support assumption~\eqref{eq:driftin:3}, to the fact that $\kappa_0$ is taken sufficiently large with respect to $\alpha$ cf.~\eqref{eq:US:kappa:0:cond:0}, taking $\eps$ to be sufficiently small, and composing with $\zeta^{-1}$ we obtain
\begin{equation*}
\sabs{ \p_2 \Thd }
\leq 4050 (1+\alpha)^2  \eps\,.
\end{equation*}
This proves~\eqref{boot-p2Theta-t}.

It thus remains to establish~\eqref{boot-p22Theta-t}--\eqref{boot-pssTheta-t}. As with~\eqref{boot-p2Theta-t}--\eqref{boot-psTheta-t}, we prove these estimates via a bootstrap/continuity argument originating at $x_1 = \xringt$. First, we need to obtain good bounds for the boundary conditions at $x_1 = \xringt$, upon differentiating \eqref{Theta-BC-t} twice, we deduce 
\begin{subequations}
\label{eq:Stafford}
\begin{align}
\p_{tt} \Thd(\xringt,x_2,t) &= 0\,, \\
\p_{2t} \Thd(\xringt,x_2,t) &= - \p_{1t}\Thd(\xringt,x_2,t) \cdot \p_2 \xstart \,, \\
\p_{22} \Thd(\xringt,x_2,t) &= - \p_{11} \Thd (\xringt,x_2,t) (\p_2 \xstart)^2 - 2 \p_{12} \Thd (\xringt,x_2,t) \p_2 \xstart \notag\\
&\qquad \qquad - \p_1 \Thd(\xringt,x_2,t) \cdot \p_{22} \xstart  \,.
\end{align}
In order to compute $\p_{12} \Thd$ and $\p_{11} \Thd$ at $(\xringt,x_2,t)$, we differentiate \eqref{p1-Theta-t-alt} with respect to $x_2,t$, and $x_1$, appeal to the bootstraps, the bounds~\eqref{p2-xstar},~\eqref{p22-xstar}, \eqref{eq:MNF:bounds:b}, \eqref{eq:MNF:bounds:d}, \eqref{eq:MNF:Rams:1}, and the fact that $\p_t \Thd(\xringt,x_2,t) = 1$ and $|\p_2 \Thd(\xringt,x_2,t)|\leq \Cn \eps^2 \mathsf{K} \brak{\mathsf{B}_6}$, to deduce
\begin{align}
|\p_{1t} \Thd(\xringt,x_2,t)| &\leq \Cn \eps^2 |\p_{12}\Thd (\xringt,x_2,t)|  + \tfrac{2(1+\alpha)}{\alpha \eps \kappa_0} \,,\\
|\p_{11} \Thd(\xringt,x_2,t)|&\leq \Cn \eps |\p_{12}\Thd (\xringt,x_2,t)|   + \tfrac{10(1+\alpha)}{\eps(\alpha \kappa_0)}  + \tfrac{8(1+\alpha)}{ \eps (\alpha\kappa_0)^2}\,,\\
|\p_{12} \Thd(\xringt,x_2,t)|&\leq \Cn \eps  |\p_{22}\Thd (\xringt,x_2,t)|  + \tfrac{28(1+\alpha)}{\alpha \kappa_0}\,.
\end{align}
\end{subequations}
Combining the bounds in~\eqref{eq:Stafford} with~\eqref{p2-xstar}, and~\eqref{p22-xstar}, we obtain that 
\begin{subequations}
\begin{align}
\p_{tt} \Thd(\xringt,x_2,t) &= 0\,, 
\label{eq:Stafford:sucks:a}\\
|\p_{2t} \Thd(\xringt,x_2,t)| &\leq  \tfrac{5(1+\alpha)}{\alpha \kappa_0}    \Cdata   \,, 
\label{eq:Stafford:sucks:b}\\
|\p_{22} \Thd(\xringt,x_2,t)| &\leq \tfrac{200(1+\alpha)}{\alpha \kappa_0}  \Cdata^2 \eps  + \tfrac{5}{\alpha\kappa_0} \cdot \Cn_{\eqref{eq:Stafford:sucks:c}} \eps \brak{\mathsf{B}_6}\,.
\label{eq:Stafford:sucks:c}
\end{align} 
\end{subequations}
Here the constant $\Cn_{\eqref{eq:Stafford:sucks:c}}$ only depends on the implicit constant from \eqref{p22-xstar}, and thus only depends on $\alpha,\kappa_0$, and $\Cdata$.
In particular, recalling that $\Cdata \leq \mathsf{B}_6$, if we let the constant $\b{22}$ appearing in~\eqref{boot-p22Theta-t} satisfy 
\begin{equation}
 \tfrac{200(1+\alpha)}{\alpha \kappa_0}  \Cdata     + \tfrac{5}{\alpha\kappa_0} \cdot \Cn_{\eqref{eq:Stafford:sucks:c}}  \leq \tfrac 12 \b{22}
 \label{eq:b22:constraint:1}
 \,,
\end{equation}
and we let the constant $\b{2\s}$ appearing in~\eqref{boot-p2sTheta-t} satisfy 
\begin{equation}
\tfrac{5(1+\alpha)}{\alpha \kappa_0}    \Cdata   \leq \tfrac 12 \b{2\s}
 \,,
  \label{eq:b2s:constraint:1}
\end{equation}
we are ensured that \eqref{boot-p22Theta-t}--\eqref{boot-pssTheta-t} hold at $x_1=\xringt$, with strict inequalities. Note that no constraint on $\b{\s\s}$ is imposed at this stage. We next show via a bootstrap / continuity argument that these bounds still hold, with strict inequalities, globally in $\mathring{\Omega}_{\mathsf{US},+}$.

We first establish~\eqref{boot-pssTheta-t}. Differentiating the first component of \eqref{eq:MNF:Rams:0} with respect to the $t$ variable, we obtain
\begin{align}
&\bigl(\p_1 + \mathfrak{M}\cir\Thd \p_2 +\mathfrak{N} \cir \Thd \p_t\bigl) (\p_{tt} \Thd)
- (\p_t \mathfrak{F}) \cir \Thd \p_{tt} \Thd
\notag\\
&=
-  (\p_t \mathfrak{M})\cir\Thd \p_t \Thd \p_{2t} \Thd 
- (\p_t \mathfrak{N}) \cir \Thd \p_t \Thd \p_{tt} \Thd
+ (\p_{tt} \mathfrak{F}) \cir \Thd (\p_{t} \Thd)^2
\notag\\
&\qquad 
- (\p_t \mathfrak{M})\cir \Thd \bigl( \p_{2t} \Thd \p_t \Thd  + \p_2 \Thd \p_{tt} \Thd \bigr)
- (\p_{tt} \mathfrak{M})\cir \Thd \p_2 \Thd (\p_t \Thd)^2
\notag\\
&\qquad 
- 2  (\p_t \mathfrak{N})\cir \Thd   \p_t \Thd \p_{tt} \Thd
-  (\p_{tt} \mathfrak{N})\cir \Thd (\p_t \Thd)^3
\,.
\label{eq:p:tt:Thd}
\end{align}
By appealing to \eqref{eq:MNF:bounds:d}, \eqref{eq:MNF:bounds:f}, and \eqref{eq:Thd:derivs}
we deduce that the right side of \eqref{eq:p:tt:Thd} may be bounded as 
\begin{equation}
\sabs{ \mathsf{RHS}_{\eqref{eq:p:tt:Thd}} } 
\leq (1+ 3 \cdot 10^{-5})^2 |(\p_{tt}\mathfrak{F} )\cir \Thd| + \Cn \brak{\mathsf{B}_6}\,,
\label{eq:p:tt:Thd:a}
\end{equation}
where we have used that the constants $\b{22}, \b{2\s}, \b{\s\s}$  only depend on $\alpha$, $\kappa_0$, and $\Cdata$, and thus may be absorbed into $\Cn$. In analogy to \eqref{eq:MNF:Rams:1}, we may use the bootstraps,~\eqref{eq:Ophelia:1} and~\eqref{Sigma0-ALE} to show that 
\begin{equation*}
 \sabs{\p_{tt} \mathfrak{F}
-\tfrac{(1-\dl)(1+\alpha)}{2} \p_{1} \tfrac{\Zbn}{\Sigma}
+ \tfrac{(1-\dl)(1-\alpha)}{2\alpha} \p_t \tfrac{\Jg \Zbn}{\Sigma}  } \leq \Cn \,,
\end{equation*}
and therefore
\begin{equation}
 \sabs{\p_{tt} \mathfrak{F}} \leq \tfrac{64(1+\alpha)}{\eps \kappa_0}\mathsf{C} _{\Zbn}
 \,.
 \label{eq:p:tt:mathfrak:F:improve}
\end{equation}
From the above estimate and~\eqref{eq:p:tt:Thd:a}, upon taking $\eps$ to be sufficiently small, we obtain
\begin{equation}
\sabs{ \mathsf{RHS}_{\eqref{eq:p:tt:Thd}} } 
\leq \tfrac{65(1+\alpha)}{\eps \kappa_0}\mathsf{C} _{\Zbn}\,.
\label{eq:p:tt:Thd:b}
\end{equation}
With the above estimate available, we compose \eqref{eq:p:tt:Thd} with the flow $(\zeta_2,\zeta_t)$, integrate in $x_1$ starting at $\xringt$, use the boundary condition \eqref{eq:Stafford:sucks:a}, the integrating factor $\mathcal{I}$ from~\eqref{eq:MNF:Rams:4} which satisfies~\eqref{eq:MNF:Rams:6}, and the bound~\eqref{eq:driftin:3}, to deduce 
\begin{align*}
\sabs{(\p_{tt} \Thd)\cir \zeta(x_1,x_2,t)}
&\leq \mathcal{I}(x_1,x_2,t) \cdot |(\p_{tt} \Thd)(\xringt,x_2,t)| 
+\sabs{ \mathsf{RHS}_{\eqref{eq:p:tt:Thd}}\cir \zeta(x_1,x_2,t)}  \cdot \tfrac{1+10^{-5}}{1-10^{-5}} \cdot|x_1 - \xringt| 
\notag\\
&\leq  \tfrac{65(1+\alpha)}{\eps \kappa_0}\mathsf{C} _{\Zbn} \cdot \tfrac{1+10^{-5}}{1-10^{-5}} \cdot 2(13\pi + 65\alpha(1+\alpha)\kappa_0)\eps
\notag\\
&\leq  \tfrac{66(1+\alpha)}{\kappa_0}\mathsf{C} _{\Zbn}  \cdot  (26\pi + 130\alpha(1+\alpha)\kappa_0) 
\,.
\end{align*}
Since the right side in the above bound depends only on $\alpha,\kappa_0$, and $\Cdata$, upon composing with the inverse flow of $\zeta$, and upon defining
\begin{equation}
 \tfrac{66(1+\alpha)}{\kappa_0}\mathsf{C} _{\Zbn}  \cdot  (26\pi + 130\alpha(1+\alpha)\kappa_0)  
 =: \tfrac 12 \b{\s\s}\,,
 \label{eq:bss:constraint}
\end{equation}
we have completed the proof of~\eqref{boot-pssTheta-t}.

Next, we establish~\eqref{boot-p2sTheta-t}. The proof is similar to \eqref{boot-pssTheta-t}, except that the boundary condition at $x_1= \xringt$ is now satisfying~\eqref{eq:Stafford:sucks:b} instead of \eqref{eq:Stafford:sucks:a}, and the evolution equation~\eqref{eq:p:tt:Thd} is now replaced by
\begin{align}
&\bigl(\p_1 + \mathfrak{M}\cir\Thd \p_2 +\mathfrak{N} \cir \Thd \p_t\bigl) (\p_{2t} \Thd)
- (\p_t \mathfrak{F}) \cir \Thd \p_{2t} \Thd
\notag\\
&=
- \bigl( (\p_2 \mathfrak{M})\cir\Thd + (\p_t \mathfrak{M})\cir \Thd \p_2 \Thd\bigr)
  \p_{2t} \Thd
- \bigl( (\p_2 \mathfrak{N} \cir \Thd + (\p_t \mathfrak{N}) \cir \Thd \p_2 \Thd \bigr)
\p_{tt} \Thd
\notag\\
&\qquad 
+ \bigl((\p_{2t}\mathfrak{F})\cir \Thd + (\p_{tt} \mathfrak{F}) \cir \Thd \p_2 \Thd\bigr)  
\p_{t} \Thd 
\notag\\
&\qquad 
- (\p_t \mathfrak{M})\cir \Thd \bigl( \p_{22} \Thd \p_t \Thd  + \p_2 \Thd \p_{2t} \Thd \bigr)
- \bigl( (\p_{2t}\mathfrak{M})\cir \Thd + (\p_{tt} \mathfrak{M})\cir \Thd \p_2 \Thd \bigr) \p_2 \Thd (\p_t \Thd)^2
\notag\\
&\qquad 
- 2  (\p_t \mathfrak{N})\cir \Thd   \p_t \Thd \p_{2t} \Thd
- \bigl((\p_{2t} \mathfrak{N})\cir \Thd + (\p_{tt} \mathfrak{N})\cir \Thd \p_2 \Thd \bigr)(\p_t \Thd)^2
\,,
\label{eq:p:2t:Thd}
\end{align}
which is obtained by differentiating the first component of~\eqref{eq:MNF:Rams:0} with respect to $x_2$.
By appealing to \eqref{eq:MNF:bounds:d}, \eqref{eq:MNF:bounds:f}, \eqref{eq:Thd:derivs} and~\eqref{eq:p:tt:mathfrak:F:improve} 
we deduce that the right side of \eqref{eq:p:2t:Thd} may be bounded as 
\begin{equation}
\sabs{ \mathsf{RHS}_{\eqref{eq:p:2t:Thd}} } 
\leq  (1+ 3 \cdot 10^{-5})\bigl( |(\p_{2t}\mathfrak{F} )\cir \Thd| + \tfrac{64(1+\alpha)}{\eps \kappa_0}\mathsf{C} _{\Zbn} \cdot 5 \cdot 10^3(1+\alpha)^2 \eps\bigr)+ \Cn \eps \brak{\mathsf{B}_6}\,.
\label{eq:p:2t:Thd:a}
\end{equation}
On the other hand, by differentiating \eqref{eq:Ophelia:1} with respect to $x_2$ and appealing to the bootstrap assumptions, we may show that 
\begin{equation*}
 \sabs{\p_{2t} \mathfrak{F}
+\tfrac{(1-\dl)(1+\alpha)}{4\alpha} \p_{2} \tfrac{\Jg\Wbn-\Jg\Zbn}{\Sigma}
+ \tfrac{(1-\dl)(1-\alpha)}{2\alpha} \p_2 \tfrac{\Jg \Zbn}{\Sigma}  } \leq \Cn \eps \,,
\end{equation*}
and therefore
\begin{equation}
 \sabs{\p_{2t} \mathfrak{F}} \leq \tfrac{11(1+\alpha)}{ \alpha \eps \kappa_0}  + \Cn 
 \leq \tfrac{12(1+\alpha)}{ \alpha \eps \kappa_0}  
 \,.
 \label{eq:p:2t:mathfrak:F:improve}
\end{equation}
From the above estimate and \eqref{eq:p:2t:Thd:a}, upon taking $\eps$ to be sufficiently small, we obtain
\begin{equation}
\sabs{ \mathsf{RHS}_{\eqref{eq:p:2t:Thd}} } 
\leq  \tfrac{13(1+\alpha)}{ \alpha \eps \kappa_0} \,.
\label{eq:p:2t:Thd:b}
\end{equation}
With the above estimate available, we compose \eqref{eq:p:2t:Thd} with the flow $(\zeta_2,\zeta_t)$, integrate in $x_1$ starting at $\xringt$, use the boundary condition \eqref{eq:Stafford:sucks:b}, the integrating factor $\mathcal{I}$ from~\eqref{eq:MNF:Rams:4} which satisfies~\eqref{eq:MNF:Rams:6}, and the bound~\eqref{eq:driftin:3}, to deduce 
\begin{align*}
\sabs{(\p_{2t} \Thd)\cir \zeta(x_1,x_2,t)}
&\leq \mathcal{I}(x_1,x_2,t) \cdot |(\p_{2t} \Thd)(\xringt,x_2,t)| 
+\sabs{ \mathsf{RHS}_{\eqref{eq:p:2t:Thd}}\cir \zeta(x_1,x_2,t)}  \cdot \tfrac{1+10^{-5}}{1-10^{-5}} \cdot|x_1 - \xringt| 
\notag\\
&\leq\tfrac{(1+10^{-5})5(1+\alpha)}{\alpha \kappa_0}    \Cdata  +  \tfrac{13(1+\alpha)}{ \alpha \eps \kappa_0} \cdot \tfrac{1+10^{-5}}{1-10^{-5}} \cdot 2(13\pi + 65\alpha(1+\alpha)\kappa_0)\eps
\notag\\
&\leq  \tfrac{6(1+\alpha)}{\alpha\kappa_0}\Cdata +  \tfrac{14(1+\alpha)}{ \alpha   \kappa_0}  \cdot  (26\pi + 130\alpha(1+\alpha)\kappa_0) 
\,.
\end{align*}
Since the right side in the above bound depends only on $\alpha,\kappa_0$, and $\Cdata$, upon composing with the inverse flow of $\zeta$, and upon defining
\begin{equation}
\tfrac{6(1+\alpha)}{\alpha\kappa_0}\Cdata +  \tfrac{14(1+\alpha)}{ \alpha   \kappa_0}  \cdot  (26\pi + 130\alpha(1+\alpha)\kappa_0) 
=: \tfrac 12 \b{2\s}\,,
 \label{eq:b2s:constraint}
\end{equation}
which automatically implies~\eqref{eq:b2s:constraint:1},  we have thus completed the proof of~\eqref{boot-p2sTheta-t}.

At last, we establish~\eqref{boot-p22Theta-t}. The proof is nearly identical to that of~\eqref{boot-p2sTheta-t} and~\eqref{boot-pssTheta-t}. In analogy to~\eqref{eq:p:2t:Thd}, by differentiating with respect to $x_2$ the second component of~\eqref{eq:MNF:Rams:0}, we deduce
\begin{align}
&\bigl(\p_1 + \mathfrak{M}\cir\Thd \p_2 +\mathfrak{N} \cir \Thd \p_t\bigl) (\p_{22} \Thd)
- (\p_t \mathfrak{F}) \cir \Thd \p_{22} \Thd
\notag\\
&=
- \bigl( (\p_2 \mathfrak{M})\cir\Thd + (\p_t \mathfrak{M})\cir \Thd \p_2 \Thd\bigr)
  \p_{22} \Thd
- \bigl( (\p_2 \mathfrak{N} \cir \Thd + (\p_t \mathfrak{N}) \cir \Thd \p_2 \Thd \bigr)
\p_{2t} \Thd
\notag\\
&\qquad 
+ \bigl((\p_{2t}\mathfrak{F})\cir \Thd + (\p_{tt} \mathfrak{F}) \cir \Thd \p_2 \Thd\bigr)  
\p_{2} \Thd 
+ \bigl((\p_{22}\mathfrak{F})\cir \Thd + (\p_{2t} \mathfrak{F}) \cir \Thd \p_2 \Thd\bigr)  
\notag\\
&\qquad 
- 2 (\p_t \mathfrak{M})\cir \Thd  \p_{22} \Thd \p_2 \Thd 
- \bigl( (\p_{2t}\mathfrak{M})\cir \Thd + (\p_{tt} \mathfrak{M})\cir \Thd \p_2 \Thd \bigr)  (\p_2 \Thd)^2
\notag\\
&\qquad 
-  (\p_2 \mathfrak{M})\cir \Thd  \p_{22} \Thd 
- \bigl( (\p_{22}\mathfrak{M})\cir \Thd + (\p_{2t}\mathfrak{M})\cir \Thd \p_2 \Thd\bigr) \p_2 \Thd
\notag\\
&\qquad 
-  (\p_t \mathfrak{N})\cir \Thd  \bigl( \p_{22} \Thd \p_{t} \Thd +  \p_2 \Thd \p_{2t} \Thd \bigr)
- \bigl((\p_{2t} \mathfrak{N})\cir \Thd + (\p_{tt} \mathfrak{N})\cir \Thd \p_2 \Thd \bigr) \p_2\Thd \p_t \Thd  
\notag\\
&\qquad 
- (\p_2 \mathfrak{N})\cir \Thd \p_{2t} \Thd
- \bigl((\p_{22} \mathfrak{N})\cir \Thd + (\p_{2t} \mathfrak{N})\cir \Thd \p_2 \Thd \bigr)  \p_t \Thd  
\,.
\label{eq:p:22:Thd}
\end{align}
As before, by appealing to \eqref{eq:MNF:bounds:d}, \eqref{eq:MNF:bounds:f}, \eqref{eq:Thd:derivs}, \eqref{eq:p:tt:mathfrak:F:improve}, and~\eqref{eq:p:2t:mathfrak:F:improve} we may show that 
\begin{equation}
 \sabs{ \mathsf{RHS}_{\eqref{eq:p:22:Thd}} } 
\leq   \bigl( |(\p_{22}\mathfrak{F} )\cir \Thd| + 2 \cdot \tfrac{12(1+\alpha)}{ \alpha \eps \kappa_0}\cdot 5\cdot10^3(1+\alpha)^2 \eps   \bigr)+ \Cn \eps \brak{\mathsf{B}_6}\,.
\label{eq:p:22:Thd:a}
\end{equation}
Using the bootstraps, the definition~\eqref{eq:frak:MNF:c}, and the bounds~\eqref{bs-Jg-2} and~\eqref{p2-Sigma} we may show that
\begin{equation*}
 \sabs{\p_{22} \mathfrak{F}} \leq \tfrac{32}{\alpha\kappa_0}  \bigl(2   \Cdata   +   28(1+\alpha)  \bigr)
 \,.
\end{equation*}
The above estimate and~\eqref{eq:p:22:Thd:a} imply that 
\begin{equation}
 \sabs{ \mathsf{RHS}_{\eqref{eq:p:22:Thd}} } 
\leq \tfrac{32}{\alpha\kappa_0}  \bigl(2   \Cdata   +   28(1+\alpha)  \bigr) +  \tfrac{12 }{ \alpha   \kappa_0}\cdot  10^4(1+\alpha)^3    \,.
\label{eq:p:22:Thd:b}
\end{equation}
With the above estimate available, we compose \eqref{eq:p:22:Thd} with the flow $(\zeta_2,\zeta_t)$, integrate in $x_1$ starting at $\xringt$, use the boundary condition \eqref{eq:Stafford:sucks:c}, the integrating factor $\mathcal{I}$ from~\eqref{eq:MNF:Rams:4} which satisfies~\eqref{eq:MNF:Rams:6}, and the bound~\eqref{eq:driftin:3}, to deduce 
\begin{align}
\sabs{(\p_{22} \Thd)\cir \zeta(x_1,x_2,t)}
&\leq \mathcal{I}(x_1,x_2,t) \cdot |(\p_{22} \Thd)(\xringt,x_2,t)| 
+\sabs{ \mathsf{RHS}_{\eqref{eq:p:2t:Thd}}\cir \zeta(x_1,x_2,t)}  \cdot \tfrac{1+10^{-5}}{1-10^{-5}} \cdot|x_1 - \xringt| 
\notag\\
&\leq \eps \Bigl( \tfrac{200(1+\alpha)}{\alpha \kappa_0}  \Cdata^2 + \tfrac{5}{\alpha\kappa_0} \cdot \Cn_{\eqref{eq:Stafford:sucks:c}}   \brak{\mathsf{B}_6}    
\notag\\
&\qquad\qquad
+ \bigl( \tfrac{32}{\alpha\kappa_0}  \bigl(2   \Cdata   +   28(1+\alpha)  \bigr) +  \tfrac{12 }{ \alpha   \kappa_0}\cdot  10^4(1+\alpha)^3 \bigr) (26\pi + 130\alpha(1+\alpha)\kappa_0) 
\Bigr)\notag\\
&\leq \eps \brak{\mathsf{B}_6} \Bigl( \tfrac{200(1+\alpha)}{\alpha \kappa_0}  \Cdata  + \tfrac{5}{\alpha\kappa_0} \cdot \Cn_{\eqref{eq:Stafford:sucks:c}}     
\notag\\
&\qquad\qquad
+ \bigl( \tfrac{32}{\alpha\kappa_0}  \bigl(2    +   28(1+\alpha)  \bigr) +  \tfrac{12 }{ \alpha   \kappa_0}\cdot  10^4(1+\alpha)^3 \bigr) (26\pi + 130\alpha(1+\alpha)\kappa_0) 
\Bigr)\notag\\
&=: \eps \brak{\mathsf{B}_6} \cdot \tfrac 12 \b{22} \label{eq:b22:constraint}
\,.
\end{align}
The above defined $\b{22}$ depends only on $\alpha,\kappa_0$, and $\Cdata$, and that this choice automatically implies~\eqref{eq:b22:constraint:1}. Upon composing with the inverse flow of $\zeta$, \eqref{eq:b22:constraint} completes the proof of~\eqref{boot-p22Theta-t}, and thus of the Lemma.
\end{proof}

\subsection{The upstream weight function $\JJ$ in $\tHdmp$} 

In analogy with Sections~\ref{sec:formation:setup}--\ref{sec:downstreammaxdev}, we introduce a weight function  (denoted by $\JJ$) that will be used in the upstream energy estimates, and which is a suitable extension of $\Jgb$ away from the pre-shock. According to the decomposition~\eqref{eq:H:dl:max:split:def} of the upstream spacetime $\tHdm$, we  separately define the weight $\JJ$ in $\tHdmp$ (see~\eqref{JJ-def-plus-t} below) and $\tHdmm$ (see~\eqref{qps-JJ-tHdmm} below), ensuring the continuity of certain derivatives across the surface $\Thd(x,\initial)$. 

In this subsection we define the weight function $\JJ$ on $\tHdmp$, which we recall is foliated by the surfaces $(x,\Thd(x,t))$, according to \eqref{eq:H:dl:max:+:def}. In light of this foliation, we may define the upstream weight function $\JJ$ as
\begin{align}
\JJ \big(x_1,x_2, \Thd(x_1,x_2,t)\big) = \mathcal{B}(x_2,t) = \Jgb(\xringt,x_2,t)   \, ,
\label{JJ-def-plus-t}
\end{align}
for all  $(x,t) \in \mathring{\Omega}_{\mathsf{US},+}$, 
where we have used the notation in~\eqref{Theta-nastyB-t} for $\mathcal{B}$.
In order to simplify our exposition, we shall sometimes use the notation $\JJ(x,\Thd)$ to mean $\JJ  (x_1,x_2,\Thd(x_1,x_2,t) ) $.  

The weight $\JJ$ was defined in~\eqref{JJ-def-plus-t}  in order to ensure that is satisfies a PDE, which is a $\dl$-modification of the $1$-characteristic transport PDE. To see this, we differentiate \eqref{JJ-def-plus-t} and obtain that in $\tHdmp$ we have
\begin{subequations} 
\label{p-JJ-Theta}
\begin{align} 
\p_1\big( \JJ (x,\Thd) \big) &=  \p_1\JJ(x,\Thd)+  \p_t \JJ (x,\Thd)\p_1\Thd =0 \,,  \label{p1-JJ-Theta}\\
\p_2 \big( \JJ (x,\Thd) \big) &=  \p_2\JJ(x,\Thd)+  \p_t \JJ (x,\Thd)\p_2\Thd
= \p_2 \mathcal{B}(x_2,t) \,, \label{p2-JJ-Theta} \\
\p_t \big( \JJ (x,\Thd) \big) &=   \p_t \JJ (x,\Thd)  \p_t\Thd = \p_t \mathcal{B}(x_2,t) \,.
 \label{ps-JJ-Theta} 
\end{align} 
\end{subequations} 
The identities in~\eqref{p-JJ-Theta} are substituted into the definition of $\Thd$ in \eqref{p1-Theta-t} to yield
\begin{subequations}
\begin{equation}
\omd (\p_t+V\p_2) \JJ 
-(2 \alpha \Sigma \Jgi)   \p_1\JJ 
+ \big(2\alpha\Sigma g^{- {\frac{1}{2}} } h,_2 \big) \p_2\JJ   
=
0
\,.\label{JJ-formula-t-evo}
\end{equation}
The boundary condition associated to the $\JJ$ evolution \eqref{JJ-formula-t-evo} is deduced from \eqref{JJ-def-plus-t} and \eqref{Theta-BC-t}, leading to 
\begin{equation}
\JJ(\xringt,x_2,t) = \mathcal{B}(x_2,t) =  \Jgb(\xstart,x_2,t)
 \label{JJ-formula-t-BC}
\end{equation}
\end{subequations}
for $t \in[\initial, \tstar)$.

In the energy estimates, the form of the $\JJ$ evolution \eqref{JJ-formula-t-evo} which is most frequently used is
\begin{equation}
 2 \alpha \Sigma  \p_1\JJ 
- 2\alpha\Sigma\Jg g^{- {\frac{1}{2}} } h,_2  \p_2\JJ   
-  \Jg (\p_t+V\p_2) \JJ 
= 
- \dl \Jg (\p_t+V\p_2) \JJ 
\,.
\label{JJ-formula-t}
\end{equation}
We will show (cf.~\eqref{qps-JJ-bound-t} below) that  $(\p_t + V\p_2)\JJ < 0$, and so the condition $\dl>0$ makes the term on the right side of \eqref{JJ-formula-t} strictly positive. In turn, this  induces a strictly positive damping term in certain energy norms; see Remark~\ref{rem:JJ-formula} below. 

\subsection{The upstream weight function $\JJ$ in $\tHdmm$}  
Let us now define the upstream weight function $\JJ$ in the spacetime region $\Hdmm$.   We set
\begin{subequations}
\label{qps-JJ-tHdmm}
\begin{equation} 
(\p_t+V\p_2) \JJ(x,t) = ((\p_t+V\p_2) \JJ)(x,\Thd(x,\initial)) \qquad \text{ for all } \qquad (x,t) \in \tHdmm \,,
\label{qps-JJ-tHdmm:a}
\end{equation} 
with the boundary condition set at the ``top'' boundary of the spacetime $\tHdmm$
\begin{equation}
\JJ(x,t) = 1\qquad \text{ for all } \qquad t  =  \Thd(x,\initial) 
\label{qps-JJ-tHdmm:b}
 \,.
\end{equation}
\end{subequations}
We note from the start that the definition of $\JJ$ in \eqref{qps-JJ-tHdmm} ensures that $\JJ$ is continuous across $\Thd(x,\initial)$ (because $\mathcal{B}(x_2,\initial) = \Jgb(\xringt,x_2,\initial) = 1$), and moreover $(\p_t + V\p_2) \JJ$ is continuous across $\Thd(x,\initial)$ (because of~\eqref{qps-JJ-tHdmm:a}).

It is convenient to solve the boundary value problem \eqref{qps-JJ-tHdmm}.
First, we compute using \eqref{p2-JJ-Theta}, \eqref{ps-JJ-Theta},   the fact that $\Jgb(x,\initial) = 1$, and the fact that $(\p_t \Jgb)(x,\initial) = \tfrac{1+\alpha}{2} (w_0),_1(x) + \tfrac{1-\alpha}{2} (z_0),_1(x)$, that for $t \in [\initial,\tstar]$, we 
\begin{align} 
\bigl( (\p_t+V\p_2) \JJ\bigr) (x,\Thd(x,\initial)) 
&= \p_t \mathcal{B}(x_2,\initial) \tfrac{1 -  V(x,\Thd(x,\initial))  \p_2 \Thd(x,\initial)}{\p_t\Thd(x,\initial)} + V(x,\Thd(x,\initial)) \p_2 \mathcal{B}(x_2,\initial)
\notag\\
&= \bigl(\tfrac{1+\alpha}{2} (w_0),_1 + \tfrac{1-\alpha}{2} (z_0),_1 \bigr)(\xringt,x_2) \tfrac{1 -  V(x,\Thd(x,\initial))  \p_2 \Thd(x,\initial)}{\p_t\Thd(x,\initial)} 
\notag\\
&=: \mathfrak{f}(x_1,x_2)
  \,.
\label{qps-JJ-H}
\end{align} 
Note importantly that the $x_1$ dependence of $\mathfrak{f}$ enters only through the argument of $\Thd$ and its derivatives.

Returning to the evolution equation for $\JJ$ in \eqref{qps-JJ-tHdmm}, we introduce a slight modification of the flow  $\xi$ defined in \eqref{xi-flow} (the modification is that the new flow is the identity at a given time $t$ instead of $\initial$), as follows: for $\initial\leq t \leq t^{\prime}$, and for $x\in \TT^2$ such that $(x,t) \in \tHdmm$, we let 
\begin{equation} 
\label{xi-flow-new}
\p_{t^{\prime}}\xi_{t}(x_1,x_2,t^{\prime}) = V(x_1, \xi_{t}(x_1,x_2,t^{\prime}),t^{\prime})\,, 
\qquad
\xi_{t}(x_1,x_2,t) = x_2 \,.
\end{equation} 
It is clear that $\xi_{t}$ satisfies the bounds~\eqref{eq:xi:nabla:xi} and~\eqref{eq:xi:hessian:xi}.
This flow is defined for times $t^{\prime}$ less than the stopping time 
\begin{subequations}
\label{xi-flow-new-stopping}
\begin{equation}
\mathsf{T}_\xi(x,t)
= \sup \bigl\{ t^{\prime} \in [t,\final) \colon (x_1,\xi_{t}(x_1,x_2,t^{\prime}),t^{\prime}) \in \tHdmm \bigr\}\,.
\label{xi-flow-new-stopping:a}
\end{equation}
Since the ``top'' temporal boundary of $\tHdmm$ is the surface $\Thd(x,\initial)$, we may alternatively characterize the stopping time $\mathsf{T}_\xi(x,t)$ as the implicit solution of
\begin{equation}
\mathsf{T}_\xi(x_1,x_2,t) = \Thd\Bigl(x_1,\xi_{t}\bigl(x_1,x_2,\mathsf{T}_\xi(x_1,x_2,t)\bigr),\initial\Bigl)
\,,
\qquad
\mathsf{T}_\xi(x_1,x_2,t) \in [t,\final)\,.
\label{xi-flow-new-stopping:b}
\end{equation}
\end{subequations}
Note in particular that $\mathsf{T}_\xi(x,\Thd(x,\initial)) = \Thd(x,\initial)$.

Using this notation, we solve for $\JJ$ in \eqref{qps-JJ-tHdmm}. For $(x_1,x_2,t)\in\tHdmm$ fixed, we compose \eqref{qps-JJ-tHdmm:a} with the flow $\xi_{t}$, and using the notation in \eqref{qps-JJ-H} deduce that
\begin{equation*}
\p_{t^\prime} \bigl( \JJ(x_1, \xi_{t}(x_1,x_2,t^\prime),t^\prime) \bigr) = \mathfrak{f}(x_1,\xi_{t}(x_1,x_2,t^\prime))
\,.
\end{equation*}
Integrating the above equation from $t^\prime=t$ and until $t^\prime = \mathsf{T}_{\xi}(x,t)$, and using the boundary condition~\eqref{qps-JJ-tHdmm:b}, we get
\begin{equation}
\JJ(x_1, x_2,t)
=
1- 
\int_{t}^{\mathsf{T}_{\xi}(x_1,x_2,t)}
\mathfrak{f}(x_1,\xi_{t}(x_1,x_2,t^\prime))
{\rm d} t^\prime
\,.
\label{JJ-def-minus-t}
\end{equation}
Identity~\eqref{JJ-def-minus-t} gives the definition of the upstream weight function $\JJ$ in the spacetime $\tHdmm$, where $\xi_t$ is defined as in~\eqref{xi-flow-new}, $\mathfrak{f}$ is given by~\eqref{qps-JJ-H}, and $T_\xi$ is given by \eqref{xi-flow-new-stopping}.

\subsection{Properties of the weight function $\JJ$}
\label{sec:JJ:US:properties}
The upstream weight function $\JJ$ is now fully defined, according to~\eqref{JJ-def-plus-t} in $\tHdmp$, and~\eqref{JJ-def-minus-t} in $\tHdmm$. We collect a number of useful properties of this weight, which will be used throughout the remainder of this section. 

\subsubsection{Lower bounds for $\JJ$}
We claim that for all $(x,t) \in \tHdm$, we have 
\begin{equation}
\JJ(x,t) \geq 
\begin{cases}
\bigl( \bar\Thd(x) - t\bigr) \tfrac{1+\alpha}{2\eps}  \cdot \tfrac{89}{100}\,, &(x,t) \in \tHdmp\,,\\
1\,, &(x,t) \in \tHdmm\,.
\end{cases}
\label{JJ-and-t}
\end{equation}
We prove \eqref{JJ-and-t} separately in the spacetime $\tHdmp$ and $\tHdmm$.  In $\tHdmp$ it is sometimes convenient to use~\eqref{JJ-and-t-new}.

According to~\eqref{eq:H:dl:max:+:def}, for any $(x,t) \in \tHdmp$ there exists $t^\prime = t^\prime(x,t)  \in [\initial,\tstar)$ such that $t = \Thd(x,t^\prime)$. The definition~\eqref{JJ-def-plus-t} then gives 
\begin{equation*}
\JJ(x,t) = \JJ(x,\Thd(x,t^\prime)) = \mathcal{B}(x_2,t^\prime) \,. 
\end{equation*}
On the other hand, from assumptions~\eqref{item:ic:max:w0} and~\eqref{item:ic:w0:x2:negative}, and bounds \eqref{eq:Jgb:identity:2}, \eqref{bs-Jg-1b}, \eqref{eq:x1star:x1vee}, and Definition~\ref{def:pre-shock}, we have
\begin{align}
\JJ(x,\Thd(x,t^\prime))
 &= \mathcal{B}(x_2,t^\prime) 
 = \Jgb(\xringt,x_2,t^\prime)
 = \Jgb(\xringt,x_2,\tstar) - \int_{t^\prime}^{\tstar} \p_t \Jgb(\xringt,x_2,t^{\prime\prime}) {\rm d} t^{\prime\prime}\notag\\
 &\geq - \int_{t^\prime}^{\tstar} \p_t \Jg(\xringt,x_2,t^{\prime\prime}) {\rm d} t^{\prime\prime}
\geq -  (\tstar-t^\prime )\tfrac{1+\alpha}{2} \bigl( (w_0),_1(\xringt,x_2) + \mathsf{C_{J_t}}\bigr)
 \notag\\
 &\geq  (\tstar-t^\prime ) \tfrac{1+\alpha}{2} \bigl(\tfrac{9}{10 \eps} -2 \mathsf{C_{J_t}}\bigr) 
\geq (\tstar-t^\prime ) \tfrac{1+\alpha}{2\eps} \cdot \tfrac{895}{1000}
 \label{JJ-and-t-new}
 \,,
\end{align}
whenever $\initial \leq t^\prime < \tstar$.
It thus remains to appeal to the intermediate value theorem in time, along with the bound~\eqref{boot-psTheta-t}, and obtain that 
\begin{equation*}
 \bar\Thd(x) - t
 =\Thd(x,\tstar) - \Thd(x,t^\prime)
 = (\tstar-t^\prime) \underbrace{\p_t \Thd(x,t^{\prime\prime})}_{\in  [1-3\cdot10^{-5},1+3\cdot10^{-5}]}\,.
\end{equation*}
From the two estimates above, we obtain the bound~\eqref{JJ-and-t}, in the spacetime $\tHdmp$.

For $(x,t) \in \tHdmm$, we appeal to the formula~\eqref{JJ-def-minus-t}. We note that by the definition~\eqref{qps-JJ-H}, the bounds~\eqref{boot-p2Theta-t}--\eqref{boot-psTheta-t}, the fact that $V = \OO(\eps)$, the assumptions on $w_0$ and $z_0$ in Section~\ref{cauchydata}, and the bound~\eqref{eq:x1star:x1vee} the function $\mathfrak{f}$ satisfies the bound
\begin{equation}
\tfrac{1+\alpha}{2\eps} \cdot  \tfrac{101}{100}
\geq 
\bigl(\tfrac{1+\alpha}{2} \cdot \tfrac{1}{ \eps} + \Cn\bigr) \cdot \tfrac{1+ \Cn \eps^2}{1 - 3 \cdot 10^{-5}} \geq -\mathfrak{f}(x_1,x_2) 
\geq \bigl(\tfrac{1+\alpha}{2} \cdot \tfrac{9}{10 \eps} - \Cn\bigr) \cdot \tfrac{1- \Cn \eps^2}{1 + 3 \cdot 10^{-5}}
\geq \tfrac{1+\alpha}{2\eps} \cdot  \tfrac{89}{100}
\,.
\label{eq:mathfrak:f:bound:0}
\end{equation}
Inserting the lower bound~\eqref{eq:mathfrak:f:bound:0} into~\eqref{JJ-def-minus-t}, results in
\begin{equation}
1 + (\mathsf{T}_\xi(x,t) - t)  \tfrac{1+\alpha}{2\eps} \cdot  \tfrac{101}{100} 
  \geq \JJ(x,t) \geq 1 + (\mathsf{T}_\xi(x,t) - t)  \tfrac{1+\alpha}{2\eps} \cdot  \tfrac{89}{100} 
\,.
\label{eq:aubergine:1}
\end{equation}
Since by definition (recall~\eqref{xi-flow-new-stopping:a}) we always have $\mathsf{T}_\xi(x,t)\geq t$, the above estimate proves \eqref{JJ-and-t} in $\tHdmm$.

\subsubsection{Comparison of $\JJ$ and $\Jg$}
Next, we compare the weight function $\JJ$ defined by~\eqref{JJ-def-plus-t} and~\eqref{JJ-def-minus-t} to $\Jg$ itself.

\begin{lemma}[\bf $\JJ$ and $\Jg$] \label{lem:JJ-le-Jg}
Assume that $\kappa_0$ is taken sufficiently large with respect to $\alpha>0$ to ensure that~\eqref{eq:US:kappa:0:cond:0} holds. Assume that the pointwise bootstraps~\eqref{boots-HH} hold in $\tHdm$, and that $\eps$ is taken to be sufficiently small, with respect to $\alpha, \kappa_0$, and $\Cdata$. Then, for $(x,t) \in \tHdm$ we have 
\begin{subequations}
\label{JJ-le-Jg}
\begin{align} 
&\JJ(x,t)\le   \tfrac{101}{100}  \Jg(x,t)\,, \qquad \mbox{whenever} \qquad |x_1| \leq 13 \pi \eps \,,
\label{JJ-le-Jg:a}
\\
&\JJ(x,t)\le \tfrac{21}{10} \Jg(x,t) 
\quad \mbox{and} \quad
|\Jg \Wbn(x,t)|\le \Cn \eps\,, 
\qquad \mbox{whenever} \qquad |x_1| > 13 \pi \eps \,.
\label{JJ-le-Jg:b}
\end{align} 
Moreover, we have
\begin{equation}
0 \le \JJ(x,t) - 1 \leq \tfrac{1}{10^3} {\bf 1}_{|x_1|<13\pi\eps}  + \tfrac{53}{50} {\bf 1}_{|x_1|>13\pi\eps} 
\quad \mbox{and} \quad 
\sabs{\Jg(x,t)-1}  \leq \tfrac{5}{10^{4}}
 \,,
 \label{JJ-le-Jg:c}
\end{equation}
\end{subequations}
for $(x,t) \in \tHdmm$. Since $\JJ > 0$ in $\tHdm$, the bounds in \eqref{JJ-le-Jg} also show $\Jg>0$ in $\tHdm$.
\end{lemma} 

\begin{proof}[Proof of Lemma \ref{lem:JJ-le-Jg}]
Let us first note that  the condition $|x_1| > 13\pi \eps$ implies, via~\eqref{eq:ic:supp}, that $(w_0),_1(x)=0$. As such, the bound~\eqref{eq:broncos:eat:shit:20} immediately implies that for such values of $x$ we have $|\Jg\Wbn| \leq \Cn \eps$. This proves the second bound in \eqref{JJ-le-Jg:b}. It thus only remains to prove the estimates that involve $\JJ$ and $\Jg$.

{\em The proof of~\eqref{JJ-le-Jg:c}.} 
We recall that by its definition in~\eqref{eq:H:dl:max:-:def}, we have $\tHdmm = \{ (x,t) \in \tHdm \colon \initial \leq t < \Thd(x,\initial)\}$. We first prove the $\Jg$ estimate in~\eqref{JJ-le-Jg:c}. Using~\eqref{bs-Jg-0}, \eqref{p1-Theta-sign-t-new}, and the fact that $\initial = \Thd(\xringt,x_2,\initial)$, we have that for all $(x,t) \in \tHdmm$, 
\begin{align*}
\abs{\Jg(x,t) - 1}
&\leq (t - \initial) \tfrac{1+\alpha}{2\eps}\bigl( |\nb_1 w_0(x)| + \eps \mathsf{C_{J_t}} \bigr)
\notag\\
&\leq (\Thd(x_1,x_2,\initial) - \Thd(\xringt,x_2,\initial) \tfrac{1+\alpha}{2\eps} \bigl( |\nb_1 w_0(x)| + \eps \mathsf{C_{J_t}} \bigr)
\notag\\
&\leq  \abs{x_1 - \xringt}  \cdot \tfrac{1}{10^5(1+\alpha)} \cdot \tfrac{1+\alpha}{2\eps} \bigl( {\bf 1}_{|x_1|\leq 13 \pi \eps} + \eps \mathsf{C_{J_t}} \bigr)
\notag\\
&\leq {\bf 1}_{|x_1|\leq 13\pi\eps} \tfrac{13\pi}{10^5} +\Cn \eps {\bf 1}_{|x_1|> 13\pi\eps} 
\leq 5 \cdot 10^{-4}
\,,
\end{align*}
upon taking $\eps$ sufficiently small. This establishes the bounds for $\Jg$ claimed in~\eqref{JJ-le-Jg:c}.

In order to estimate $\JJ(x,t)$ in the spacetime $\Hdmm$ we appeal to \eqref{eq:aubergine:1}, which thus necessitates an upper bound for the non-negative quantity $\mathsf{T}_\xi(x,t) - t$, where $\mathsf{T}_\xi$ is as defined by~\eqref{xi-flow-new-stopping:a}. For this purpose, we note that by construction (see the line below \eqref{xi-flow-new-stopping:b}), we have that 
\begin{equation}
\bigl(\mathsf{T}_{\xi}(x,t)-t\bigr) \bigl|_{t=\Thd(x,\initial)} = 0 \,.
\label{eq:aubergine:2}
\end{equation}
Moreover, implicitly differentiating \eqref{xi-flow-new-stopping:b} with respect to $t$ shows that 
\begin{equation*}
 \p_t \mathsf{T}_{\xi}(x,t)  
= \tfrac{\p_2 \Thd(x_1,\xi_t(x,\mathsf{T}_\xi(x,t)),\initial) (\p_t \xi_t)(x,\mathsf{T}_\xi(x,t))}{1 - V(x_1,\xi_t(x,\mathsf{T}_\xi(x,t)),\mathsf{T}_\xi(x,t)) \p_2 \Thd(x_1,\xi_t(x,\mathsf{T}_\xi(x,t)),\initial)}
\end{equation*}
while differentiation of \eqref{xi-flow-new} with respect to $t$ shows that 
\begin{equation*}
(\p_t \xi)(x,t^\prime) = - V (x,t) \exp\Bigl(\int_{t}^{t^\prime} \p_2 V(x_1,\xi_t(x,t^{\prime\prime}),t^{\prime\prime}) {\rm d} t^{\prime\prime} \Bigr)
\,,
\end{equation*}
Combining the two identities above with \eqref{boot-p2Theta-t} and the fact that $V,\nb V = \OO(\eps)$ by the bootstraps, we deduce that 
\begin{equation}
\sabs{1 + \p_t \bigl(\mathsf{T}_\xi(x,t) -t \bigr)} = \sabs{\p_t \mathsf{T}_{\xi}(x,t)} \les \eps^2
\,.
\label{eq:aubergine:3}
\end{equation}
From \eqref{eq:aubergine:2} and \eqref{eq:aubergine:3} we deduce in turn that 
\begin{align}
(1- \Cn \eps^2) (\Thd(x,\initial) - t) \leq \mathsf{T}_\xi(x,t) - t \leq (1+ \Cn \eps^2) (\Thd(x,\initial) - t)
\,,
\label{eq:aubergine:4}
\end{align}
for all $(x,t) \in \tHdm$. At this stage we consider two cases, $|x_1|\leq 13 \pi \eps$, and $|x_1|>13\pi \eps$. In the first case, when $|x_1| \leq 13\pi \eps$, we use that~\eqref{p1-Theta-sign-t-new} implies 
\begin{align}
{\bf 1}_{|x_1|\leq 13\pi \eps} (\mathsf{T}_{\xi}(x,t) - t)
&\leq 
{\bf 1}_{|x_1|\leq 13\pi \eps} |\Thd(x,\initial) - \initial| (1 + \Cn \eps^2)
\notag\\
&\leq {\bf 1}_{|x_1|\leq 13\pi \eps}  |\Thd(x_1,x_2,\initial)-\Thd(\xringt,x_2,\initial)|
(1 + \Cn \eps^2)
\notag\\
&\leq {\bf 1}_{|x_1|\leq 13\pi \eps}  |x_1-\xringt| \cdot \tfrac{1}{10^5(1+\alpha)}(1 + \Cn \eps^2)
\leq \tfrac{27\pi\eps}{10^5(1+\alpha)} \leq \tfrac{\eps}{10^3(1+\alpha)}
\label{eq:aubergine:5}
\,.
\end{align}
In the second case, when  $|x_1|>13\pi \eps$, we recal the definition of the domain $\mathring{\Omega}_{\mathsf{US},+}$ on which $\Thd$ is defined, cf.~\eqref{eq:Omega:US:+}, to note that $\Thd(x,\initial) \leq \final$, and therefore 
\begin{equation}
 {\bf 1}_{|x_1|\leq 13\pi \eps} (\mathsf{T}_\xi(x,t) - t )
\leq (1+ \Cn \eps^2) (\Thd(x,\initial) - \initial)
\leq (1+ \Cn \eps^2) (\final-\initial) \leq \tfrac{2\eps}{1+\alpha} \cdot \tfrac{52}{50}\,.
 \label{eq:aubergine:6}
\end{equation}
Combining the estimates \eqref{eq:aubergine:5}--\eqref{eq:aubergine:6} with the upper bound in \eqref{eq:aubergine:1}, we deduce
\begin{align*}
\JJ(x,t) \leq 1 + \bigl(\mathsf{T}_{\xi}(x,t)-t) \tfrac{1+\alpha}{2\eps} \cdot \tfrac{101}{100}
&\leq 1 +  \tfrac{2\eps}{1+\alpha} \bigl(  {\bf 1}_{|x_1|<13\pi\eps} \tfrac{1}{2 \cdot 10^3} + {\bf 1}_{|x_1|>13\pi\eps} \tfrac{52}{50} \bigr) \tfrac{1+\alpha}{2\eps} \cdot \tfrac{101}{100}
\notag\\
&\leq 1 + {\bf 1}_{|x_1|<13\pi\eps} \tfrac{1}{10^3} + {\bf 1}_{|x_1|>13\pi\eps} \tfrac{53}{50} \,.
\end{align*}
This proves the upper bound for $\JJ$ in~\eqref{JJ-le-Jg:c}.

{\em The proof of~\eqref{JJ-le-Jg:a}--\eqref{JJ-le-Jg:b}.} 
We note that due to~\eqref{JJ-le-Jg:c}, for $(x,t) \in \tHdm$ (by continuity, on the closure of this spacetime), the bounds~\eqref{JJ-le-Jg:a}--\eqref{JJ-le-Jg:b} are already known to hold. This is because for $|x_1|\leq 13\pi \eps$ we have $(1+ 10^{-3}) \cdot (1 + 5 \cdot 10^{-4}) \leq 1.01$ while for $|x_1|>13 \pi \eps$ we have $(1+ \tfrac{53}{50}) \cdot (1 + 5 \cdot 10^{-4}) \leq 2.1$. Moreover, since we have already established the $\Jg \Wbn$ bound stated in~\eqref{JJ-le-Jg:b}, it only remains to prove the pointwise upper bound for $\JJ \Jgi$ stated in \eqref{JJ-le-Jg:a}--\eqref{JJ-le-Jg:b}, on the spacetime $\tHdmp$. 

The proof consists of two parts. The first one is to obtain rough lower bounds of $\Jg(x,t)$ away from $t=\tstar$ and $x_1=\xringt$. The second consists in comparing $\JJ$ with $\Jg$ using various decompositions of the spacetime $\tHdmp$.

The arguments in the first part of the proof, are reminiscent of the argument in the proof of Proposition~\ref{prop:mJg:well:def}, except that instead of the downstream geometry, we consider the upstream one. In analogy to \eqref{eq:x1:sharp:def} we define 
\begin{equation} 
x_1^{\flat}(x_2) = \bigl\{ x_1 \colon x_1 < x_1^\vee(x_2), \nb_1w_0(x_1,x_2) = - \tfrac{17}{20} \bigr\} 
\label{eq:x1:flat:def}
\,.
\end{equation} 
The fact that the function $\TT \ni x_2 \mapsto x_1^\flat(x_2)$ given by \eqref{eq:x1:flat:def} is well-defined and differentiable requires a proof. This proof is nearly identical to the one given in the proof of Proposition~\ref{prop:mJg:well:def} with signs changed; the changes are as follows. The existence of at least one value $x_1^\flat(x_2) \in (-13\pi \eps, x_1^\vee(x_2) )$ satisfying $\nb_1w_0(x_1^\flat(x_2),x_2) = - \tfrac{17}{20}$ follows from the intermediate value theorem, because $-\frac{17}{20} \in (-\frac{9}{10},0)$. Then, any such possible value $x_1^\flat(x_2)$ must satisfy $x_1^\flat(x_2) \leq x_1^\vee(x_2) - \frac{1}{40} \eps$. This is because by the intermediate value theorem, \eqref{eq:why:the:fuck:not:0}, and \eqref{item:ic:w0:x2:negative}, we have $\frac{1}{20} \leq |\nb_1 w_0(x_1^\flat(x_2),x_2) - \nb_1 w_0(x_1^\vee(x_2),x_2) | = |x_1^\flat(x_2) - x_1^\vee(x_2)| \|\p_1 \nb_1  w_0\|_{L^\infty_x} \leq \frac{2}{\eps} |x_1^\flat(x_2) - x_1^\vee(x_2)| $. In particular, with \eqref{eq:x1star:x1vee} we also have $x_1^\flat(x_2) \leq x_1^*(x_2,t) - \frac{\eps}{41}$. We may then define $x_1^\flat(x_2)$ as the largest value of $x_1$ for which $\nb_1w_0(x_1,x_2) = - \tfrac{17}{20}$ and prove that for all $x_1 < x_1^\flat(x_2)$ we must have $\nb_1w_0(x_1,x_2) > - \tfrac{17}{20}$, yielding also the uniqueness of $x_1^\flat(x_2)$. To see this, note that when $x_1< x_1^\flat(x_2)$, then $x_1 -x_1^\vee(x_2) < - \frac{1}{40} \eps < - \eps^{\frac 74}$. Hence, by assumption~\eqref{item:ic:w0:d11:positive} on the initial data 
we know that for all $x_1 \in [-13 \pi \eps, x_1^\flat(x_2))$ with $\nb_1 w_0(x) < - \frac 13$, we must have $\nb_1^2 w_0(x)\leq - \eps^{\frac 78} < 0$, showing that as $x_1$ decreases, $\nb_1 w_0(x)$ strictly increases from the value $-\frac{17}{20}$ when $x_1 = x_1^\flat(x_2)$, until it reaches the value $-\frac 13$ at some point $x_1= x_1^\bot(x_2)$. Additionally, \eqref{item:ic:w0:d11:positive} implies that for $x_1 < x_1^\bot(x_2)$ we have that $\nb_1 w_0(x) \geq -\frac 13 $: this is because if $\nb_1 w_0(x)$ wanted to dip below the value $-\frac 13$, then it would need to increase as a function of $x_1$, but \eqref{item:ic:w0:d11:positive} implies that $\nb_1 w_0(x)$ can only increase  in $x_1$ if $\nb_1 w_0(x) \geq - \frac 13$.
Thus, we have shown that \eqref{eq:x1:flat:def} gives a well-defined object, and that 
\begin{subequations}
\label{eq:x1:flat:prop}
\begin{align}
&x_1^\flat(x_2) \leq x_1^\vee(x_2) - \tfrac{1}{40} \eps\,,
\qquad 
x_1^\flat(x_2) \leq \xringt - \tfrac{1}{41} \eps\,,
\label{eq:x1:flat:prop:a}
\\
&\nb_1w_0(x)> - \tfrac{17}{20} \mbox{ for all } x_1 < x_1^\flat(x_2)\,,
\qquad
-1 < \nb_1w_0(x) < - \tfrac{17}{20} \mbox{ for all } x_1^\flat(x_2)<x_1<x_1^\vee(x_2)\,.
\label{eq:x1:flat:prop:b}
\end{align}
In particular, \eqref{eq:x1:flat:prop:b}, the inequality $-\tfrac{17}{20} < -\frac 13$, and assumption \eqref{item:ic:w0:d11:positive} imply that 
\begin{align}
\nb_1^3 w_0(x) \geq  \tfrac{1}{3} \mbox{ for all } x_1^\flat(x_2)\leq x_1 \leq x_1^\vee(x_2)\,.
\label{eq:x1:flat:prop:c} 
\end{align}
Using \eqref{eq:x1:flat:prop:c} we may also obtain a lower bound for the value of $x_1^\flat(x_2)$. Since $\nb_1^2 w_0(x_1^\vee(x_2),x_2)=0$, we may write $\nb_1 w_0(x_1^\flat(x_2),x_2) - \nb_1w_0(x_1^\vee(x_2),x_2) = \tfrac 12 (\tfrac{x_1^\flat(x_2) - x_1^\vee(x_2)}{\eps})^2 \nb_1^3 w_0(x_1^\prime,x_2)$ for some $x_1^\prime \in (x_1^\flat(x_2),x_1^\vee(x_2))$. Therefore, we deduce $|x_1^\flat(x_2) - x_1^\vee(x_2)|\leq \sqrt{6} \eps  |\nb_1 w_0(x_1^\flat(x_2),x_2) - \nb_1w_0(x_1^\vee(x_2),x_2)|^{\frac 12} \leq \sqrt{6} \eps \sqrt{3/20} \leq \tfrac{19}{20}\eps$, so
\begin{align}
x_1^\flat(x_2) \geq x_1^\vee(x_2) - \tfrac{19}{20} \eps\,,
\qquad 
x_1^\flat(x_2) \geq \xringt - \eps\,.
\label{eq:x1:flat:prop:d} 
\end{align}
\end{subequations}

We also recall from Section~\ref{sec:design:JJ:weight:DS} that there exists a unique $x_1^\sharp(x_2) > x_1^\vee(x_2)$ such that $\nb_1 w_0(x_1^\sharp(x_2),x_2) = - \tfrac{17}{20}$; cf.~\eqref{eq:x1:sharp:def}. Moreover, in analogy to~\eqref{eq:x1:flat:prop} it is shown in Section~\ref{sec:design:JJ:weight:DS} (just like in the above argument), that 
\begin{subequations}
\label{eq:x1:sharp:prop}
\begin{align}
&x_1^\sharp(x_2) \geq x_1^\vee(x_2) + \tfrac{1}{40} \eps\,,
\qquad 
x_1^\sharp(x_2) \geq \xringt + \tfrac{1}{41} \eps\,,
\label{eq:x1:sharp:prop:a}
\\
&\nb_1w_0(x) > - \tfrac{17}{20} \mbox{ for all } x_1 > x_1^\sharp(x_2)\,,
\qquad
-1 < \nb_1w_0(x) < - \tfrac{17}{20} \mbox{ for all } x_1^\vee(x_2)<x_1<x_1^\sharp(x_2)\,,
\label{eq:x1:sharp:prop:b}
\\
&\nb_1^3 w_0(x) \geq  \tfrac{1}{3} \mbox{ for all } x_1^\vee(x_2)\leq x_1 \leq x_1^\sharp(x_2)\,,
\label{eq:x1:sharp:prop:c} 
\\
&x_1^\sharp(x_2) \leq x_1^\vee(x_2) + \tfrac{19}{20} \eps\,,
\qquad 
x_1^\sharp(x_2) \leq \xringt + \eps\,.
\label{eq:x1:sharp:prop:d} 
\end{align}
\end{subequations}

The bounds obtained in \eqref{eq:x1:flat:prop}  for $w_0$ have two immediate consequences regarding the behavior of $\Jg$, when combined with Corollary~\ref{cor:Jg:initial}.  The first bound in~\eqref{eq:x1:flat:prop:b} together with the support assumption~\eqref{eq:ic:supp} for $(w_0),_1$, and the bound \eqref{bs-Jg-0},  show that whenever $(x,t) \in \tHdm$ is such that $x_1\not\in(x_1^\flat(x_2),x_1^\sharp(x_2))$, 
\begin{subequations}
\label{eq:Jg:useful:US}
\begin{align}
\Jg(x,t)
&\geq 1 + (t - \initial) \tfrac{1+\alpha}{2\eps}\bigl(\nb_1 w_0(x) - \eps \mathsf{C_{J_t}} \bigr)
\notag\\
&\geq 1 - \tfrac{2\eps}{1+\alpha} \cdot \tfrac{51}{50}\cdot \tfrac{1+\alpha}{2\eps}\bigl(\tfrac{17}{20} {\bf 1}_{|x_1|\leq 13 \pi \eps} +  \eps \mathsf{C_{J_t}} \bigr)
\notag\\
&\geq 
\tfrac{1}{9} \bigl({\bf 1}_{-13\pi\eps\leq x_1 \leq x_1^\flat(x_2)} + {\bf 1}_{x_1^\sharp(x_2) \leq x_1 \leq 13 \pi \eps}\bigr) +
(1-\Cn \eps) {\bf 1}_{13\pi \eps < |x_1| \leq \pi} \,.
\label{eq:Jg:useful:US:a}
\end{align}
Similarly, assumption~\eqref{eq:why:the:fuck:not:0},  and the fact that $t^*(x_2) \leq \final$ (which holds by the construction of $\Jgb$, see the proof of Lemma~\ref{lem:q:invertible}),  gives that whenever $(x,t) \in \tHdm$ is such that $\initial \leq t\leq \tfrac 12 (\final+\initial)$, we have 
\begin{equation}
\Jg(x,t)
\geq 1 + (t - \initial) \tfrac{1+\alpha}{2\eps}\bigl(\nb_1 w_0(x) - \eps \mathsf{C_{J_t}} \bigr)
\geq 1 - \tfrac{2\eps}{1+\alpha} \cdot \tfrac{51}{100}\cdot \tfrac{1+\alpha}{2\eps}\bigl(1 +  \eps \mathsf{C_{J_t}} \bigr)
 \geq 
\tfrac 49\,.
\label{eq:Jg:useful:US:b}
\end{equation}
It thus remains to consider points $(x,t) \in \tHdm$ such that $x_1\in (x_1^\flat(x_2),x_1^\sharp(x_2))$ and $\tfrac 12 (\final+\initial) < t \leq \tstar$. In this region, the bounds~\eqref{eq:x1:flat:prop:c}, \eqref{bs-Jg-2} with $i=j=1$,  \eqref{eq:x1star:x1vee}, \eqref{x1-star}, and \eqref{eq:ic:norm},  give
\begin{align}
\nb_1^2 \Jg(x,t)
&\geq (t-\initial) \tfrac{1+\alpha}{2\eps} \bigl(\nb_1^3 w_0(x)- \eps \Cn \mathsf{K} \brak{\mathsf{B_6}} \bigr)
\notag\\
&\geq \tfrac{2\eps}{1+\alpha} \cdot \tfrac{51}{100} \cdot \tfrac{1+\alpha}{2\eps} \bigl(\tfrac 13 - \eps \Cn \mathsf{K} \brak{\mathsf{B_6}} - \eps^2 \Cn \mathsf{K} \brak{\mathsf{B_6}} \bigr)
 \geq \tfrac{1}{6}  \,.
 \label{eq:Jg:useful:US:c}
\end{align}
\end{subequations}

With \eqref{eq:Jg:useful:US} in hand, we turn to the second part of the proof and establish \eqref{JJ-le-Jg:a}--\eqref{JJ-le-Jg:b} for $(x,t) \in \tHdmp$.
 
We note that by definition, this is the set $\{ (x,t) \in \tHdm \colon \Thd(x,\initial) < t < \Thd(x,\tstar)\}$, which may also be written as $ \{ (x,t) \in \Hdm \colon \exists t^\prime \in [\initial,\tstar), \mbox{ with } t = \Thd(x,t^\prime)\}$. We split this spacetime into three different regions, as follows:

{\em (i)\;The case $\initial \leq t \leq \tfrac 12 (\initial+\final)$.}  
According to \eqref{JJ-def-plus-t} and \eqref{Jgb-le-Jg}, for all $t^\prime \in [\initial,\tstar)$ and all $x\in \TT^2$  such that $(x,\Thd(x,t^\prime)) \in \tHdm$,  
\begin{equation*}
\JJ(x,t) = \JJ(x,\Thd(x,t^\prime)) = \Jgb(\xringt,x_2,t^\prime)  \,.
\end{equation*}
On the other hand, using \eqref{eq:Jgb:identity:2}, \eqref{eq:x1star:def}, \eqref{eq:x1star:x1vee}, \eqref{Jgb-le-Jg}, the bootstrap~\eqref{bs-Jg,1} present in~\eqref{boots-HH}, and the mean value theorem, we have
\begin{align*}
\Jgb(\xringt,x_2,t^\prime) 
&\leq 
\Jgb(x_1^*(x_2,t^\prime),x_2,t^\prime)
+ |\xringt-x_1^*(x_2,t^\prime)| \cdot \|\Jg,_1\|_{L^\infty_{x,g}}
\notag\\
&\leq \Jgb(x_1,x_2,t^\prime) + \Cn \mathsf{K} \brak{\mathsf{B}_6} \eps^3 \cdot 4(1+\alpha) \eps^{-1}
\,.
\end{align*}
Combining the two bounds above results in 
\begin{equation}
\label{eq:driftin:1}
 \JJ(x,t)\leq  \Jg(x,t^\prime) + \Cn \mathsf{K} \brak{\mathsf{B}_6} \eps^2\,, \qquad \mbox{where} \qquad t = \Thd(x,t^\prime),
\end{equation}
for all $t^\prime \in [\initial,\tstar)$ and all $x\in \TT^2$  such that $(x,\Thd(x,t^\prime)) \in \tHdm$.

Next, we bound $\Jg(x,t^\prime)$ from above in terms of $\Jg(x,t)=\Jg(x, \Thd(x,t^\prime))$, via the fundamental theorem of calculus in time and \eqref{bs-Jg-1b} as
\begin{equation*}
\Jg(x,t^\prime) 
= \Jg(x,t) - \int_{t^\prime}^{t} \p_{t} \Jg(x,t^{\prime\prime}) {\rm d}t^{\prime\prime}
\leq \Jg(x,t)  + \tfrac{1+\alpha}{2} \bigl( | (w_0),_1(x)| + \mathsf{C_{J_t}} \bigr) |t - t^\prime|
\,.
\end{equation*}
It is thus clear that we need a bound for $t -t^\prime = \Thd(x_1,x_2,t^\prime) -t^\prime$. Recall cf.~\eqref{Theta-BC-t} that $\Theta^\dl(\xringt,x_2,t^\prime) = t^\prime$. The mean value theorem in $x_1$, the bound~\eqref{p1-Theta-sign-t-new}, and the bootstrap~\eqref{bs-Jg-simple} present in~\eqref{boots-HH} give
\begin{align*}
|t - t^\prime| = |\Thd(x_1,x_2,t^\prime)-t^\prime|
&=
|\Thd(x_1,x_2,t^\prime)-\Theta^\dl(\xringt,x_2,t^\prime)|
\notag\\
&=|x_1 - \xringt| \, |\p_1 \Theta^{\dl}(x_1^\prime,x_2,t^\prime)|
\leq |x_1 - \xringt| \tfrac{5}{\alpha \kappa_0}
\,.
\end{align*}
Combining the above two estimates  we thus arrive at the bound
\begin{equation}
\Jg(x,t^\prime) \leq \Jg(x,t) + \tfrac{5(1+\alpha)}{2\alpha \kappa_0}  |x_1 - \xringt| \bigl( | (w_0),_1(x)| + \mathsf{C_{J_t}} \bigr)  
 \label{eq:driftin:2} 
 \,,
\end{equation}
valid for all $t^\prime \in [\initial,\tstar)$ and all $x\in \TT^2$  such that $(x,t) = (x,\Thd(x,t^\prime)) \in \tHdm$.

At this stage, combining \eqref{eq:ic:supp}, \eqref{eq:why:the:fuck:not:0}, \eqref{eq:driftin:1}, \eqref{eq:driftin:2}, and~\eqref{eq:driftin:3}, we obtain
\begin{align}
\JJ(x,t) 
&\leq 
\Jg(x,t) + \tfrac{5(1+\alpha)}{2\alpha \kappa_0}  |x_1 - \xringt| \bigl( | (w_0),_1(x)| + \mathsf{C_{J_t}} \bigr)   + \Cn \mathsf{K} \brak{\mathsf{B}_6} \eps^2
\notag\\
&\leq 
\Jg(x,t) 
+ {\bf 1}_{|x_1|\leq 13\pi\eps} \tfrac{5(1+\alpha)}{2\alpha \kappa_0} \cdot 26\pi\eps  \cdot\bigl( \eps^{-1} + \mathsf{C_{J_t}} \bigr) 
\notag\\
&\qquad
+ {\bf 1}_{|x_1|> 13\pi\eps} \tfrac{5(1+\alpha)}{2\alpha \kappa_0} \cdot 2(13\pi + 65\alpha(1+\alpha)\kappa_0)\eps \cdot 
\bigl(0 + \mathsf{C_{J_t}} \bigr)   
+ \Cn \mathsf{K} \brak{\mathsf{B}_6} \eps^2
\notag\\
&\leq 
\Jg(x,t) 
+ {\bf 1}_{|x_1|\leq 13\pi\eps} \tfrac{210(1+\alpha)}{\alpha \kappa_0} +\Cn \eps {\bf 1}_{|x_1|> 13\pi\eps}  
\,,
\label{eq:driftin:driftin} 
\end{align}
upon taking $\eps$ sufficiently small, and 
for all $t^\prime \in [\initial,\tstar)$ and all $x\in \mathcal{X}_{\rm fin} $  such that $(x,\Thd(x,t^\prime)) \in \tHdm$.

Note that up to this point, the condition  $\initial \leq t \leq \tfrac 12 (\initial+\final)$ was not used. We use this condition  now, in order to bound $\Jg(x,t)$ from below via~\eqref{eq:Jg:useful:US:b} by $\tfrac 49$. Combined with~\eqref{eq:driftin:driftin}, this results in the bound
\begin{equation}
\JJ(x,t) 
\leq 
\Jg(x,t) \bigl( 1   
+ {\bf 1}_{|x_1|\leq 13\pi\eps} \tfrac{9}{4} \cdot \tfrac{210(1+\alpha)}{\alpha \kappa_0} +\Cn \eps {\bf 1}_{|x_1|> 13\pi\eps}  \bigr)
\label{eq:driftin:4}
\,.
\end{equation}
Thus, if we ensure that $\kappa_0$ is taken sufficiently large with respect to $\alpha$ (only), in order to ensure that \eqref{eq:US:kappa:0:cond:0} holds, 
then $ \tfrac{9}{4} \cdot \tfrac{210(1+\alpha)}{\alpha \kappa_0} \leq \tfrac{1}{500} $ and 
upon taking $\eps$ sufficiently small with respect to $\alpha,\kappa_0$, and $\Cdata$, we deduce $\JJ(x,t) \leq \tfrac{101}{100} \Jg(x,t)$, which is consistent with both \eqref{JJ-le-Jg:a} and \eqref{JJ-le-Jg:b}.

{\em (ii) \; The case $\tfrac 12 (\initial + \final) \leq t \leq \min\{ \Thd(x,\tstar),\final \}$, $x_1 \not \in (x_1^\flat(x_2),x_1^\sharp(x_2))$.}  
We still refer to the previously obtained bound~\eqref{eq:driftin:driftin}, except that in this region we use a different lower bound for $\Jg$. This lower bound is coming from~\eqref{eq:Jg:useful:US:a}, which gives $\Jg(x,t) \geq \tfrac{1}{9}$ whenever $x_1 \not \in (x_1^\flat(x_2),x_1^\sharp(x_2))$. Similarly to~\eqref{eq:driftin:4}, we thus obtain
\begin{equation}
\JJ(x,t) 
\leq 
\Jg(x,t) \bigl( 1   
+ {\bf 1}_{|x_1|\leq 13\pi\eps} 9 \cdot \tfrac{210(1+\alpha)}{\alpha \kappa_0} +\Cn \eps {\bf 1}_{|x_1|> 13\pi\eps}  \bigr)
\label{eq:driftin:5}
\,.
\end{equation}
Thus, if $\kappa_0$ is taken sufficiently large with respect to $\alpha$ to ensure that \eqref{eq:US:kappa:0:cond:0} holds, then   $9 \cdot \tfrac{210(1+\alpha)}{\alpha \kappa_0} \leq \tfrac{1}{200} $ and 
upon taking $\eps$ sufficiently small with respect to $\alpha,\kappa_0$, and $\Cdata$, we deduce $\JJ(x,t) \leq \tfrac{101}{100} \Jg(x,t)$, which is consistent with both \eqref{JJ-le-Jg:a} and \eqref{JJ-le-Jg:b}. 
 
{\em (iii) \; The case $\tfrac 12 (\initial + \final) \leq t \leq \min\{ \Thd(x,\tstar),\final \}$, $x_1 \in (x_1^\flat(x_2),x_1^\sharp(x_2))$.}  
This region requires a different analysis, which is not based on~\eqref{eq:driftin:driftin}; this is because we cannot obtain a uniform strictly positive lower bound for $\Jg$ in this spacetime. Instead, both $\JJ$ and $\Jg$ may degenerate towards $0$, and we need to carefully track this degeneracy. To track this degeneracy, we employ a Taylor series expansion in $x_1$, which in light of~\eqref{JJ-def-plus-t} and the fact that by construction we have $\Thd(\xringt,x_2,t^\prime)= t^\prime$, gives 
\begin{align}
&\tfrac{100}{101} \JJ(x_1,x_2,t)
=
\tfrac{100}{101} \JJ(x_1,x_2,\Thd(x_1,x_2,t^\prime))
= 
\tfrac{100}{101} \Jgb(\xringt,x_2,t^\prime)
\notag\\
&\leq 
\tfrac{100}{101} \Jg(\xringt,x_2,t^\prime)
= 
\Jg(\xringt,x_2,t^\prime) -  \tfrac{1}{101} \Jg(\xringt,x_2,t^\prime)
\notag\\
&=
\Jg(\xringt,x_2,\Thd(\xringt,x_2,t^\prime))
 - \tfrac{1}{101} \Jg(\xringt,x_2,t^\prime)
\notag\\
&= 
\Jg(x_1,x_2,\Thd(x_1,x_2,t^\prime))
\notag\\
&\ \ \  
- \tfrac{1}{101} \Jg(\xringt,x_2,t^\prime)
- (x_1-\xringt) 
\bigl( 
\Jg,_1
+ 
\p_{t} \Jg  \p_1 \Thd  
\bigr)(\xringt,x_2,t^\prime)
\notag\\
&\ \ \  
-\tfrac 12 (x_1 - \xringt)^2 
\Bigl(
\Jg,_{11}(x_1^\prime,x_2,\Thd(x_1^\prime,x_2,t^\prime))
+ 2 \p_{t} \Jg,_1(x_1^\prime,x_2,\Thd(x_1^\prime,x_2,t^\prime))  \p_1 \Theta^{\dl}(x_1^\prime,x_2,t^\prime)
\notag\\
&
\qquad \qquad \qquad \qquad \qquad \quad   
+ \p_{t}^2 \Jg(x_1^\prime,x_2,\Thd(x_1^\prime,x_2,t^\prime))  (\p_1 \Theta^{\dl})^2(x_1^\prime,x_2,t^\prime)
\notag\\
&
\qquad \qquad \qquad \qquad \qquad \quad 
+\p_{t}  \Jg(x_1^\prime,x_2,\Thd(x_1^\prime,x_2,t^\prime))\p_1^2 \Theta^{\dl}(x_1^\prime,x_2,t^\prime)
\Bigr)
\notag\\
&=:
\Jg(x_1,x_2,t) - \mathsf{Err}
\label{eq:stone:2}
\end{align}
for some $x_1^\prime$ which lies in between $x_1$ and $\xringt$, and where the error term $\mathsf{Err}$ is defined by the last three terms (containing the minus signs) in the earlier equality. Our claim is that for the range of $x_1$ and $t = \Thd(x,t^\prime)$ considered here, we have that $\mathsf{Err} \geq 0$, which then would result in the estimate $\JJ(x,t) \leq \tfrac{101}{100} \Jg(x,t)$, thereby completing the proof of~\eqref{JJ-le-Jg:a} and~\eqref{JJ-le-Jg:b} in the spacetime $\tHdmp$.
  
It thus remains to show that $\mathsf{Err} \geq 0$. To see this, we first analyze the coefficient of the term which is linear in $x_1 - \xringt$, namely $(\Jg,_1+ \p_t \Jg  \p_1 \Theta^{\dl}  )(\xringt,x_2,t^\prime)$.
Due to Proposition~\ref{prop:pre-shock} we have that 
$\Jg,_1(\xringt,x_2,t^\prime) = - \int_{t^\prime}^{t^*(x_2)} \p_t \Jg,_1(\xringt,x_2,t^{\prime\prime}) {\rm d} t^{\prime\prime}$. By also appealing to \eqref{bs-Jg-2a} and \eqref{eq:d1d1w0} (at time $t^*(x_2)$) we deduce that for $t^\prime \leq \tstar$, we have 
\begin{subequations}
\begin{equation}
|\Jg,_1(\xringt,x_2,t^\prime)|  \leq (t^*(x_2)-t^\prime)   \tfrac{9(1+\alpha)}{\eps} \mathsf{C}_{\Zbn} 
\,.
\label{eq:stone:3a}
\end{equation} 
On the other hand, $\Jg(\xringt,x_2,t^\prime) = - \int_{t^\prime}^{t^*(x_2)} \p_t \Jg(\xringt,x_2,t^{\prime\prime}) {\rm d} t^{\prime\prime}$, and by also appealing to \eqref{bs-Jg-1b} and~\eqref{eq:x1star:x1vee} (at time $t^*(x_2)$) we deduce that for $t^\prime \leq \tstar$, we have
\begin{align}
\Jg(\xringt,x_2,t^\prime)  
&\geq   (t^*(x_2)-t^\prime) \tfrac{1+\alpha}{2} (-(w_0),_1(\xringt,x_2) - \mathsf{C_{J_t}}) \notag\\
&\geq  (t^*(x_2)-t^\prime) \tfrac{1+\alpha}{2} (-(w_0),_1(x_1^\vee(x_2),x_2) - 2 \mathsf{C_{J_t}}) \geq  (t^*(x_2)-t^\prime) \tfrac{2(1+\alpha)}{5\eps}.
\label{eq:stone:3b}
\end{align}
The inequalities \eqref{eq:stone:3a} and \eqref{eq:stone:3b} combined give 
\begin{equation}
\label{eq:stone:3c}
|\Jg,_1(\xringt,x_2,t^\prime)|\leq 23 \mathsf{C}_{\Zbn} \Jg(\xringt,x_2,t^\prime)\,,
\end{equation} for $t^\prime\leq \tstar$. 
\end{subequations}
In order to estimate the $\p_{t} \Jg  \p_1 \Theta^{\dl} $ term, we use the bootstrap~\eqref{bs-Jg,1} contained in~\eqref{boots-HH}, together with  \eqref{p1-Theta-sign-t}, to deduce 
\begin{equation}
|\p_{t} \Jg(\xringt,x_2,t^\prime)  \p_1 \Theta^{\dl}(\xringt,x_2,t^\prime)| 
\leq  \tfrac{4(1+\alpha)}{\eps} \tfrac{4}{\alpha \kappa_0} \Jg(\xringt,x_2,t^\prime) 
= \tfrac{16(1+\alpha)}{\alpha \kappa_0 \eps}\Jg(\xringt,x_2,t^\prime) \,.
\label{eq:stone:4}
\end{equation}
Combining \eqref{eq:stone:3c} and \eqref{eq:stone:4} with $x_1  \in (x_1^\flat(x_2),x_1^\sharp(x_2))$, \eqref{eq:x1:flat:prop:d}, and \eqref{eq:x1:sharp:prop:d}, we obtain
\begin{align}
|x_1 - \xringt| \sabs{(\Jg,_1+ \p_{t} \Jg  \p_1 \Theta^{\dl}  )(\xringt,x_2,t^\prime)}
&\leq  |x_1 - \xringt| \tfrac{17(1+\alpha)}{\alpha \kappa_0 \eps}\Jg(\xringt,x_2,t^\prime)
\notag\\
&\leq   \tfrac{34(1+\alpha)}{\alpha \kappa_0}\Jg(\xringt,x_2,t^\prime)
\leq 10^{-4} \Jg(\xringt,x_2,t^\prime)
\,,\label{eq:stone:4a}
\end{align}
where in the last inequality we have used that $\kappa_0$ satisfies~\eqref{eq:US:kappa:0:cond:0}. Since $10^{-4} < \tfrac{1}{101}$, the above estimate turns out to be sufficiently strong.
 
Second, we analyze the coefficient of the term in \eqref{qinv-H} which is quadratic in $x_1 - \xringt$. For this purpose, we note that since $x_1,\xringt \in (x_1^\flat(x_2),x_1^\sharp(x_2))$, we also have $x_1^\prime \in (x_1^\flat(x_2),x_1^\sharp(x_2))$. Moreover, since we care about $\tfrac 12 (\initial + \final) \leq t \leq \min\{ \Thd(x,\tstar),\final \}$, we  may apply  \eqref{eq:Jg:useful:US:c} to deduce 
\begin{subequations}
\label{eq:stone:5}
\begin{equation}
\label{eq:stone:5a}
\Jg,_{11}(x_1^\prime,x_2,\Theta^{\dl}(x_1^\prime,x_2,t^\prime)) \geq \tfrac{1}{6 \eps^2} \,.
\end{equation} 
Next, we appeal to~\eqref{bs-Jg-2a}, \eqref{eq:why:the:fuck:not:0}, \eqref{p1-Theta-sign-t-new}, and~\eqref{eq:US:kappa:0:cond:0}, to deduce  
\begin{equation}
2 \sabs{ \p_{t} \Jg,_{1}(x_1^\prime,x_2,\Thd(x_1^\prime,x_2,t^\prime)) \p_1 \Theta^{\dl}(x_1^\prime,x_2,t^\prime)}
\leq 2  \cdot \tfrac{1+\alpha}{2\eps^2} \bigl(|\nb_1^2 w_0(x_1^\prime,x_2)| + \Cn \eps\bigr) \cdot \tfrac{5}{\alpha \kappa_0}
\leq \tfrac{11(1+\alpha)}{\alpha \kappa_0 \eps^2} 
\leq \tfrac{10^{-4}}{\eps^2}  \,.
\label{eq:stone:5b}
\end{equation} 
In a similar fashion, by appealing to \eqref{bs-Jg-2b}, we deduce
\begin{equation}
\sabs{\p_t^2 \Jg(x_1^\prime,x_2,\Thd(x_1^\prime,x_2,t^\prime))  (\p_1 \Thd)^2(x_1^\prime,x_2,t^\prime)}
\leq
 \tfrac{4(1+\alpha)|1-\alpha|  \mathsf{C_{\Zbn}}}{\eps} \cdot \bigl(\tfrac{5}{\alpha \kappa_0}\bigr)^2
\leq  \tfrac{10^{-4}}{\eps^2} 
\,,
\label{eq:stone:5c}
\end{equation}
upon taking $\eps$ to be sufficiently small. 
Lastly, upon applying~\eqref{eq:jesus:just:left:chicago:c},  and taking into account the bootstrap assumptions,~\eqref{eq:Q:all:bbq-H},  we deduce
\begin{align}
\sabs{\p_{t}  \Jg(x_1^\prime,x_2,\Thd(x_1^\prime,x_2,t^\prime))\p_1^2 \Thd(x_1^\prime,x_2,t^\prime)}
&\leq     \tfrac{ 1+\alpha }{\eps} \cdot \tfrac{10^{-4}}{\eps (1+\alpha)}
= \tfrac{10^{-4}}{\eps^2}
\label{eq:stone:5d} 
\,.
\end{align}
\end{subequations}
Thus, putting together the bounds established in~\eqref{eq:stone:5}, we obtain that the coefficient of $\tfrac 12 (x_1-\xringt)^2$ appearing in the definition of $\mathsf{Err}$ in \eqref{eq:stone:2}, satisfies
\begin{align}
&\Jg,_{11}(x_1^\prime,x_2,\Thd(x_1^\prime,x_2,t^\prime))
+ 2 \p_{t} \Jg,_1(x_1^\prime,x_2,\Thd(x_1^\prime,x_2,t^\prime))  \p_1\Thd(x_1^\prime,x_2,t^\prime)
\notag\\
&
\qquad
+ \p_{t}^2 \Jg(x_1^\prime,x_2,\Thd(x_1^\prime,x_2,t^\prime))  (\p_1\Thd)^2(x_1^\prime,x_2,t^\prime)
+\p_{t}  \Jg(x_1^\prime,x_2,\Thd(x_1^\prime,x_2,t^\prime))\p_1^2 \Thd(x_1^\prime,x_2,t^\prime)
\notag\\
&\geq \tfrac{1}{6 \eps^2} - \tfrac{3 \cdot 10^{-4}}{\eps^2} 
\geq \tfrac{1}{7 \eps^2}\,.
\label{eq:stone:6}
\end{align}
To summarize, we deduce from~\eqref{eq:stone:4a} and \eqref{eq:stone:6} that the term $\mathsf{Err}$ appearing in~\eqref{eq:stone:2} satisfies
\begin{equation}
 \mathsf{Err}
 \geq \tfrac{1}{101} \Jg(\xringt,x_2,t^\prime) - \tfrac{1}{10^4} \Jg(\xringt,x_2,t^\prime) + \tfrac 12 (x_1-\xringt)^2 \cdot \tfrac{1}{7 \eps^2}
\,,
\label{eq:Err:lower:bound:impala}
\end{equation}
and hence this term is non-negative. In light of~\eqref{eq:stone:2}, this completes the proof of~\eqref{JJ-le-Jg:a} and~\eqref{JJ-le-Jg:b} in the spacetime $\Hdmp \cup \Theta^\dl(x,\sin)$.
\end{proof}

\subsubsection{Damping and anti-damping properties of $\JJ$} One of the crucial properties of the weight $\JJ$ is that its time derivative is strictly negative. This fact is recorded next.

\begin{lemma}
\label{lem:JJdamping}
Assume that $\kappa_0$ satisfies~\eqref{eq:US:kappa:0:cond:0}, that the bootstraps~\eqref{bootstraps-H} hold, and  that  $\eps$ is  then taken sufficiently small with respect to $\alpha,\kappa_0$, and $\Cdata$. Then
for all $(x,t) \in \tHdm$, we have
\begin{subequations}
\label{eq:waitin:for:the:bus}
\begin{equation}
- (\p_t + V \p_2) \JJ(x,t) \geq \tfrac{1+\alpha}{2 \eps} \cdot \tfrac{899}{1000}
\,.
 \label{qps-JJ-bound-t}
\end{equation}
We also have the first derivative upper bounds 
\begin{align}
-\tfrac{402(1+\alpha)}{\eps} \leq  \p_t \JJ(x,t)  &\leq - \tfrac{2(1+\alpha)}{5\eps}\,,
\label{mod-ps-JJ-t}
\\
\sabs{\p_2 \JJ(x,t)} &\leq   3 \cdot 10^6(1+\alpha)^3 + 9 \Cdata + 3 \b{2\s}\,,
\label{fat-marmot2b-t}
\\
\sabs{\p_1 \JJ(x,t)} &\leq \tfrac{1}{10^3 \eps}\,,
\label{fat-marmot2-t}
\\
\sabs{(\p_t + V \p_2)^2 \JJ(x,t)}  
&\leq \tfrac{202^2(1+\alpha)^2}{\eps^2}\,,
\label{qps2JJ-bound-t}
\end{align}
\end{subequations}
for all $(x,t) \in \tHdmp \cup \tHdmm$. Here $\b{2\s} = \b{2\s}(\alpha,\kappa_0,\Cdata)>0$ is the constant from \eqref{boot-p2sTheta-t}.
\end{lemma}
\begin{proof}[Proof of Lemma~\ref{lem:JJdamping}]

In order to prove~\eqref{qps-JJ-bound-t}, we recall that $(\p_t + V\p_2) \JJ$ is continuous across $\Thd(x,\initial)$, and that from the definitions~\eqref{qps-JJ-tHdmm:a},~\eqref{qps-JJ-H} and the lower bound in~\eqref{eq:mathfrak:f:bound:0}, we have
\begin{equation*}
- (\p_t+V\p_2) \JJ(x,t)
= - \mathfrak{f}(x) \geq \tfrac{1+\alpha}{2\eps} \cdot \bigl(\tfrac{9}{10 (1 + 3 \cdot 10^{-5})}  - \Cn \eps\bigr)\,,
\end{equation*}
for all $(x,t) \in \tHdmm$. This bound is consistent with \eqref{qps-JJ-bound-t} upon taking $\eps$ to be sufficiently small, since $\tfrac{9}{10 + 4\cdot 10^{-4}} > \tfrac{899}{1000}$. For the $(x,t) \in \tHdmp$, we write $t=\Thd(x,t^\prime)$ for some $t^\prime\in(\initial,\tstar)$, and use identities~\eqref{p-JJ-Theta} to write
\begin{equation*}
(\p_t+V\p_2) \JJ(x,t)
=
(\p_t+V\p_2) \JJ(x,\Thd(x,t^\prime) 
= 
\p_t \mathcal{B}(x_2,t^\prime) \tfrac{1 -  V(x,\Thd(x,t^\prime))  \p_2 \Thd(x,t^\prime)}{\p_t\Thd(x,t^\prime)} + V(x,\Thd(x,t^\prime)) \p_2 \mathcal{B}(x_2,t^\prime)
\,.
\end{equation*}
Appealing to the bounds established earlier in~\eqref{eq:Qcal:bbq:temp:3}, \eqref{eq:MNF:bounds:a}, and~\eqref{eq:Thd:derivs}, for each $(x,t) \in \tHdmp$ we obtain
\begin{align*}
- (\p_t+V\p_2) \JJ(x,t) 
&\geq -\p_t \Jgb(\xringt,x_2,t^\prime) \bigl(\tfrac{1-\Cn \eps^2}{1 + 3 \cdot 10^{-5}} - \Cn \eps^2 \bigr)
\notag\\
&\geq \tfrac{1+\alpha}{2\eps} \bigl( \tfrac{9}{10} - \eps \Cn\bigr)\bigl(\tfrac{1-\Cn \eps^2}{1 + 3 \cdot 10^{-5}} - \Cn \eps^2 \bigr)
\geq 
 \tfrac{1+\alpha}{2\eps} \bigl( \tfrac{9}{10} \cdot{1}{1+ 3 \cdot 10^{-5}} - \Cn \eps\bigr).
\end{align*}
This bound is consistent with \eqref{qps-JJ-bound-t} upon taking $\eps$ to be sufficiently small, since $\tfrac{9}{10 + 4\cdot 10^{-4}} > \tfrac{899}{1000}$. This completes the proof of~\eqref{qps-JJ-bound-t}.

We next turn to the bounds~\eqref{mod-ps-JJ-t}--\eqref{fat-marmot2-t}. In the spacetime $\tHdmp$ upon writing $(x,t) = (x,\Thd(x,t^\prime))$ and appealing to \eqref{p-JJ-Theta}, \eqref{eq:Thd:derivs}, \eqref{p1-Theta-sign-t-new}, and~\eqref{eq:Qcal:bbq:temp:3}, we have 
\begin{subequations}
\label{eq:la:grange:1}
\begin{align}
\p_t \JJ(x,t) 
&= \p_t \JJ(x,\Thd(x,t^\prime)) 
= \tfrac{\p_t \mathcal{B}(x_2,t)}{\p_t \Thd(x,t^\prime)}  \in \bigl[-\tfrac{402(1+\alpha)}{\eps}, - \tfrac{2(1+\alpha)}{5\eps}\bigr] \,,
\\
\sabs{\p_1 \JJ(x,t)} 
&= \sabs{\p_1 \JJ(x,\Thd(x,t^\prime))} 
= \sabs{\p_1 \Thd(x,t^\prime) \p_t \JJ(x,\Thd(x,t^\prime))}
\leq \tfrac{1}{10^5 (1+\alpha)} \cdot \tfrac{402(1+\alpha)}{\eps} 
\leq \tfrac{1}{10^3 \eps}  \,,
\\
\sabs{\p_2 \JJ(x,t)} 
&= \sabs{\p_2 \JJ(x,\Thd(x,t^\prime))} 
\leq \sabs{\p_2 \mathcal{B}(x_2,t^\prime)} + \sabs{\p_2 \Thd(x,t^\prime) \p_t \JJ(x,\Thd(x,t^\prime))}\notag\\
&\qquad \qquad \qquad \qquad\quad
\leq 5(1+\alpha) + 5 \cdot 10^3 (1+\alpha)^2 \eps \cdot \tfrac{402(1+\alpha)}{\eps} 
\leq 3 \cdot 10^6 (1+\alpha)^3\,.
\end{align}
\end{subequations}
On the other hand, for $(x,t)\in \tHdmm$, we may directly differentiate~\eqref{JJ-def-minus-t}, which in turn requires derivative estimates for the stopping time $\mathsf{T}_\xi$, the flow $\xi_t$, and the function $\mathfrak{f}$. Estimates for the flow $\xi_t$ and its first order derivatives were previously obtained in~\eqref{eq:xi:nabla:xi}. Concerning the space derivatives of $\mathsf{T}_\xi$, implicitly differentiating~\eqref{xi-flow-new-stopping:b} and appealing to \eqref{eq:Thd:derivs}--\eqref{p1-Theta-sign-t-new}, 
we deduce
\begin{subequations}
\label{eq:rough:boy:1}
\begin{align}
\sabs{\p_1 \mathsf{T}_\xi(x,t)}
&\leq \tfrac{1}{1-\Cn \eps^2} \bigl(\tfrac{1}{10^5(1+\alpha)} + \Cn \eps^2\bigr) \leq \tfrac{2}{10^5(1+\alpha)} \,, \\
\sabs{\p_2 \mathsf{T}_\xi(x,t)}
&\leq \tfrac{1}{1-\Cn \eps^2}\cdot 5 \cdot 10^3 (1+\alpha)^2   \eps \cdot(1+ \Cn \eps^2)
\leq 6 \cdot 10^3 (1+\alpha)^2   \eps
\,,
\end{align}
\end{subequations}
for all $(x,t) \in \tHdmm$. 
For the function~$\mathfrak{f}$, the pointwise estimate~\eqref{eq:mathfrak:f:bound:0} is supplemented with the uniform derivative bounds on $\TT^2$ 
\begin{subequations}
\label{eq:rough:boy:2}
\begin{align}
\sabs{\p_1 \mathfrak{f}(x)} 
&\leq \bigl( \tfrac{1+\alpha}{2\eps} + \Cn\bigr) \cdot \tfrac{|\p_{1t} \Thd(x,\initial) + \Cn \eps}{(1-3\cdot 10^{-5})^2}
\leq \tfrac{1+\alpha}{10^5 \eps^2}
\,,\\
\sabs{\p_2 \mathfrak{f}(x)} 
&\leq \bigl( \tfrac{1+\alpha}{2\eps} (4 \Cdata + 2) + \Cn\bigr) \cdot \tfrac{1+\Cn \eps^2}{1-3 \cdot10^{-5}}
+  \bigl( \tfrac{1+\alpha}{2\eps}  + \Cn\bigr) \cdot \tfrac{\b{2\s} + \Cn \brak{\mathsf{B}_6}\eps^2}{(1-3 \cdot10^{-5})^2}  \leq \tfrac{(1+\alpha)}{\eps} ( 4\Cdata + 2  + \b{2\s} )
\,,
\end{align}
\end{subequations}
which in turn follow from the pointwise bootstrap assumptions~\eqref{boots-HH}, the assumptions on the initial data in Section~\ref{cauchydata}, and the bounds \eqref{p2-xstar}, \eqref{eq:Thd:derivs}--\eqref{eq:jesus:just:left:chicago}. Finally, upon differentiating~\eqref{JJ-def-minus-t}, and appealing to the bounds~\eqref{eq:mathfrak:f:bound:0},  \eqref{eq:aubergine:5}--\eqref{eq:aubergine:6}, \eqref{eq:rough:boy:1}--\eqref{eq:rough:boy:2}, and~\eqref{eq:xi:nabla:xi}, we obtain
\begin{subequations}
\label{eq:la:grange:2}
\begin{align}
|\p_1 \JJ(x,t)|
&\leq \sabs{\p_1 \mathsf{T}_\xi(x,t) \mathfrak{f}(x_1,\xi_t(x,\mathsf{T}_\xi(x,t)))}
+\sabs{\mathsf{T}_\xi(x,t)-t} \bigl(\|\p_1 \mathfrak{f}\|_{L^\infty_x} + \|\p_1 \xi_t\|_{L^\infty_{x,t}}  \|\p_2 \mathfrak{f}\|_{L^\infty_x}\bigr)
\notag\\
&\leq \tfrac{2}{10^5 \eps} + \tfrac{2\eps}{1+\alpha} \cdot \tfrac{52}{50} \cdot \bigl( \tfrac{(1+\alpha)}{10^5 \eps^2} + \Cn \bigr)
\leq  \tfrac{1}{10^3 \eps} 
\,,
\\
|\p_2 \JJ(x,t)|
&\leq \sabs{\p_2 \mathsf{T}_\xi(x,t) \mathfrak{f}(x_1,\xi_t(x,\mathsf{T}_\xi(x,t)))}
+\sabs{\mathsf{T}_\xi(x,t)-t}  \|\p_2 \xi_t\|_{L^\infty_{x,t}}  \|\p_2 \mathfrak{f}\|_{L^\infty_x} 
\notag\\
&\leq 6 \cdot10^3 (1+\alpha)^2 \eps 
\cdot\tfrac{1+\alpha}{2\eps} \cdot \tfrac{101}{100} +  \tfrac{2\eps}{1+\alpha} \cdot \tfrac{52}{50} \cdot \bigl(1 + \Cn \eps \bigr) \cdot \tfrac{(1+\alpha)}{\eps} ( 4\Cdata + 2  + \b{2\s} )
\notag\\
&\leq 4 \cdot10^3 (1+\alpha)^3  +   \tfrac{21}{10}  ( 4\Cdata + 2  + \b{2\s} ) \,,
\end{align}
for all $(x,t) \in \tHdmm$. The $\p_2 \JJ$ estimate above, together with \eqref{qps-JJ-tHdmm:a} and \eqref{eq:mathfrak:f:bound:0} shows that 
\begin{equation}
- \tfrac{1+\alpha}{\eps}  
\leq - \tfrac{1+\alpha}{2\eps} \cdot \tfrac{101}{100} - \Cn \eps \leq \p_t \JJ(x,t) \leq -\tfrac{1+\alpha}{2\eps} \cdot \tfrac{89}{100} + \Cn \eps \leq -\tfrac{2(1+\alpha)}{5\eps} 
\end{equation}
for all $(x,t) \in \tHdmm$
\end{subequations}
Together, the estimates \eqref{eq:la:grange:1} and \eqref{eq:la:grange:2} complete the proof of the bounds~\eqref{mod-ps-JJ-t}--\eqref{fat-marmot2-t}, for all $(x,t) \in \tHdmp \cup \tHdmm$.

In order to complete the proof of the Lemma, we need to establish~\eqref{qps2JJ-bound-t}. For $(x,t) \in \tHdmm$, we apply $(\p_t +V\p_2)$ to \eqref{qps-JJ-tHdmm:a}, appeal to the fact that $\mathfrak{f}$ is independent of $t$, that $|V|\leq \Cn \eps$, and that $\p_t \mathfrak{f}$ satisfies~\eqref{eq:rough:boy:2}, to deduce that 
\begin{equation*}
\sabs{(\p_t +V\p_2)^2 \JJ(x,t)} \leq \Cn\,,
\end{equation*}
for all $(x,t) \in \tHdmm$. This estimate is better than what is required by~\eqref{qps2JJ-bound-t}. Lastly, for $(x,t) = (x,\Thd(x,t^\prime)) \in \tHdmp$ we simply decompose 
\begin{equation*}
(\p_t + V \p_2)^2 \JJ = \p_{tt} \JJ + 2 V \p_{2t} \JJ + V^2 \p_{22} \JJ + \p_t V \p_2 \JJ
\,.
\end{equation*}
Then, upon differentiating \eqref{p2-JJ-Theta} and \eqref{ps-JJ-Theta}, we derive
\begin{align*}
\p_{tt} \JJ(x,\Thd(x,t^\prime)) 
&= \tfrac{\p_{tt} \mathcal{B}(x_2,t^\prime) - \p_{tt} \Thd(x,t^\prime) \p_t \JJ(x,\Thd(x,t^\prime))}{(\p_t \Thd(x,t^\prime))^2} 
\,,\\
\p_{2t} \JJ(x,\Thd(x,t^\prime)) 
&= \tfrac{\p_{2t} \mathcal{B}(x_2,t^\prime) - \p_{2t} \Thd(x,t^\prime) \p_t \JJ(x,\Thd(x,t^\prime)) - \p_2 \Thd(x,t^\prime) \p_{tt} \JJ(x,\Thd(x,t^\prime))}{\p_t \Thd(x,t^\prime)} 
\,,\\
\p_{22} \JJ(x,\Thd(x,t^\prime)) 
&=  \p_{22} \mathcal{B}(x_2,t^\prime) - \p_{22} \Thd(x,t^\prime) \p_t \JJ(x,\Thd(x,t^\prime)) 
\notag\\
&\qquad 
- 2 \p_2 \Thd(x,t^\prime) \p_{2t}\JJ(x,\Thd(x,t^\prime)) - (\p_2 \Thd(x,t^\prime))^2 \p_{tt}\JJ(x,\Thd(x,t^\prime)) 
\,.
\end{align*}
Using the estimates~\eqref{eq:Jgb:identity:2}, \eqref{eq:Jgb:identity:3}, \eqref{eq:Qcal:bbq:temp:3}, \eqref{bs-Jg-2b}, \eqref{eq:MNF:bounds:a}, \eqref{eq:MNF:bounds:c}, and \eqref{eq:Thd:derivs}, we deduce 
\begin{align*}
\sabs{\p_{tt} \JJ(x,\Thd(x,t^\prime))}
&= \bigl( \tfrac{200(1+\alpha)}{\eps (1 - 3 \cdot 10^{-5})} \bigr)^2 + \tfrac{\Cn}{\eps}  \leq \tfrac{201^2(1+\alpha)^2}{\eps^2}
\,,\\
\sabs{\p_{2t} \JJ(x,\Thd(x,t^\prime))}
&\leq  \tfrac{\Cn}{\eps}  
\,,\\
\sabs{\p_{22} \JJ(x,\Thd(x,t^\prime))}
&\leq \Cn
\,.
\end{align*}
Combining the above estimates, we obtain that 
\begin{equation*}
\sabs{(\p_t + V \p_2)^2 \JJ (x,t)}
\leq \tfrac{201^2(1+\alpha)^2}{\eps^2} + \Cn \,,
\end{equation*}
for all $(x,t) = (x,\Thd(x,t^\prime)) \in \tHdmp$. This completes the proof of\eqref{qps2JJ-bound-t}, and thus of the lemma.
\end{proof}

We conclude this subsection by providing an upstream analogue of Lemma~\ref{lem:damping:anti:damping}.

\begin{lemma}[\bf Upstream damping and anti-damping]
\label{lem:damping:anti:damping-H-t}
Assume that $\kappa_0$ satisfies~\eqref{eq:US:kappa:0:cond:0}, that the bootstraps~\eqref{bootstraps-H} hold, and  that  $\eps$ is  then taken sufficiently small with respect to $\alpha,\kappa_0$, and $\Cdata$.
Then, we have that
\begin{equation}
-\tfrac 32 \JJh \Jg  (\p_t +V\p_2) \JJ+   \JJss(\p_t + V\p_2)   \Jg 
\ge \tfrac{1+\alpha}{6 \eps} \JJh\Jg     
\,,
\label{eq:fakeJg:LB-H-t}
\end{equation}
pointwise for all $(x,t) \in \tHdm$.
\end{lemma}
\begin{proof}[Proof of Lemma~\ref{lem:damping:anti:damping-H-t}]
Using~\eqref{qps-JJ-bound-t}, \eqref{JJ-le-Jg}, and~\eqref{bs-Jg-1a}, we have that 
\begin{align*}
 -\tfrac 32 \JJh \Jg  (\p_t +V\p_2) \JJ+   \JJss(\p_t + V\p_2)   \Jg 
 &\geq 
 \tfrac 32 \JJh \Jg  \cdot \tfrac{1+\alpha}{2\eps} \cdot \tfrac{899}{1000}
 - \JJss \bigl(\tfrac{1+\alpha}{2\eps} {\bf 1}_{|x_1|\leq 13\pi\eps}+ \Cn \bigr)
 \notag\\
 &\geq 
 \JJh \Jg  \cdot \tfrac{1+\alpha}{2\eps} \bigl( \tfrac 32 \cdot \tfrac{899}{1000} - \tfrac{101}{100} - \Cn \eps\bigr)
 \notag\\
 &\geq 
 \JJh \Jg  \cdot \tfrac{1+\alpha}{2\eps} \cdot \tfrac 13
 \,,
\end{align*}
where  $\eps$ was taken sufficiently small to obtain the last inequality.  This proves
the claimed lower-bound \eqref{eq:fakeJg:LB-H-t}.
\end{proof}

\def\xstar{x_1^*}
\def\xstart{x_1^*(x_2,\s)}

\subsection{Pre-shock flattening and the upstream spacetime set}  
It is convenient to  flatten the ``set  of times'' $ \tstar$ at which the pre-shock occur.
We define the transformation
\begin{align} 
\s = \mathfrak{q}(x_2,t)
&:= \frac{2\eps}{1+ \alpha } \left(  \frac{ t- \tstar }{\tstar -\initial}\right)
= \initial \left( \frac{  \tstar-t }{\tstar -\initial}\right)\,,
\label{t-to-s-transform-H}
\end{align} 
for all $(x,t) \in \tHdm$. Sometimes it is more convenient to express \eqref{t-to-s-transform-H} as 
\begin{equation} 
t= \mathfrak{q}^{-1}(x_2,\s) = \tstar -(\tstar-\initial) \tfrac{\s}{\initial}
\,.
\label{qinv-H}
\end{equation} 
From \eqref{t-to-s-transform-H}, we see that
the initial $\s$-time is given by
\begin{equation} 
\sin := \mathfrak{q}(x_2,\initial) = \initial =  - \tfrac{2\eps}{1+ \alpha } \,.
\label{sin}
\end{equation} 
By design, the set of pre-shock times $\{t = \tstar\}$ gets mapped to the time slice $\{\s=0\}$, which is to say that in  $(x,\s)$ coordinates, the  pre-shock set is given by 
\begin{equation} 
\Xi^* := \Bigl\{ \bigl(\xringt ,x_2,  0 \bigr) \colon x_2\in\TT\Bigr\} \,. \label{Xi*-s}
\end{equation} 

Next, we note that  for all $(x,t) \in \tHdm$ we have 
\begin{equation*}
\mathfrak{q}(x_2,t) \leq \tfrac{2\eps}{1+ \alpha }\cdot  \tfrac{\final- \tstar }{\tstar -\initial} 
\leq \tfrac{2\eps}{1+ \alpha }\cdot  \tfrac{\tfrac{1}{50} \cdot \tfrac{2\eps}{1+\alpha} + \Cn \eps^2}{\tfrac{2\eps}{1+\alpha} - \Cn \eps^2} \leq \tfrac{2\eps}{1+ \alpha }\cdot \tfrac{1+\Cn \eps}{50} \,.
\end{equation*}
As such, the  $\s$-time variable can never exceed 
\begin{equation} 
\sfin :=\tfrac{2\eps}{1+ \alpha }\cdot  \max_{x_2\in\TT}\bigl(  \tfrac{\final- \tstar }{\tstar -\initial}  \bigr) \leq  \tfrac{2\eps}{1+ \alpha }\cdot \tfrac{1+\Cn \eps}{50} \,. \label{sfin}
\end{equation}

Taking into account~\eqref{eq:spacetime-Theta-t}, for upstream development in $(x,\s)$ variables, we use the spacetime set $\Hdm$ defined by 
\begin{align}
\Hdm := \bigl\{ (x, \s) \in \mathcal{X}_{\rm fin} \times [\sin , \sfin) \colon  \s = \mathfrak{q}(x_2,t)\,, (x,t) \in \tHdm \bigr\} \,.
 \label{eq:spacetime:H0}
\end{align}
Moreover, we similarly define
\begin{subequations}
\begin{align}
\Hdmp := \bigl\{ (x, \s) \in \mathcal{X}_{\rm fin} \times [\sin , \sfin) \colon  \s = \mathfrak{q}(x_2,t)\,, (x,t) \in \tHdmp \bigr\} \,,
 \label{eq:spacetime:H0:plus}
 \\
 \Hdmm := \bigl\{ (x, \s) \in \mathcal{X}_{\rm fin} \times [\sin , \sfin) \colon  \s = \mathfrak{q}(x_2,t)\,, (x,t) \in \tHdmm \bigr\} \,.
 \label{eq:spacetime:H0:minus}
\end{align}
\end{subequations}

Given any function $f\colon\tHdm \to \mathbb{R}$, we define the
corresponding function $\tilde f \colon \Hdm \to \mathbb{R}$ by
\begin{equation}
\tilde f(x,\s) := f(x,t), \qquad \mbox{where} \qquad \s = \mathfrak{q}(x_2,t) \,.
\label{eq:f:tilde:f-H}
\end{equation}
From~\eqref{eq:f:tilde:f-H}, and the chain-rule, \eqref{t-to-s-transform-H},  and  \eqref{eq:x1star:critical} we obtain  
\begin{subequations}
\label{eq:xt:xs:chain:rule-H}
\begin{align} 
\p_t f(x,t) &= \Qd(x_2)   \p_\s \tilde f(x,\s) =: \tfrac{1}{\eps} \nbs_\s \tilde f(x,\s)\,, \label{ptf-psf} \\
\p_2f(x,t) &=  \bigl(\p_2 - \Qb(x_2,\s)  \p_\s\bigr) \tilde f(x,\s) =: \nbs_2 \tilde f(x,\s) \,,  \\
\p_1f(x,t) &= \p_1 \tilde f (x,\s) =: \tfrac{1}{\eps} \nbs_1 \tilde f(x,\s)\,,
\end{align} 
\end{subequations}
where for compactness of notation we have introduced the functions 
\begin{subequations} 
\label{QQQ-H}
\begin{alignat}{2}
\Qd(x_2)  &= \p_t \mathfrak{q}(x_2,t) &&=  \tfrac{2\eps}{1+ \alpha } \tfrac{1}{\tstar-\initial}  =  \tfrac{-\initial}{\tstar-\initial}  \,,
\label{eq:QQQ:a-H}
\\
\Qb(x_2,\s) &= -\p_2 \mathfrak{q}(x_2,t) &&=   \tfrac{2\eps}{1+ \alpha } \tfrac{\p_2 \tstar}{\tstar-\initial} \big( 1 - \tfrac{\tstar -t}{\tstar-\initial} \big)  =  \p_2 \tstar \Qd(x_2) \big( 1 -\tfrac{\s}{\sin}\bigr)
 \,.
\label{eq:QQQ:aa-H}
\end{alignat}
For later use, it is also convenient to define 
\begin{align}
\Q( x,\s )  
:= \Qd(x_2) - \tilde V(x,\s)  \Qb(x_2,\s) 
&=  \tfrac{-\initial}{\tstar-\initial} \Big(1 - \tilde V(x,\s) \p_2 \tstar \big(1- \tfrac{\s}{\sin}  \big) \Big) 
\,,
\label{eq:QQQ:b-H}
\end{align} 
and 
\begin{align} 
\Qc = \p_\s \Q = - \p_\s\tilde V  \Qb + \tilde V \Qd \tfrac{\p_2 t^*(x_2)}{\sin}  \,, \qquad 
\Qr_\s  = \p_\s \Qd =0 \,, \qquad 
\Qr_2 =  \p_\s \Qb= - \Qd \tfrac{\p_2 t^*(x_2)}{\sin} 
\,.
\label{eq:QQQ:c-H}
\end{align} 
\end{subequations} 
With the above notation, it follows from \eqref{eq:xt:xs:chain:rule-H} that the spacetime gradient operator in  $(x,t)$ variables, namely $\nb = (\eps\p_t, \eps \p_1, \p_2)$, becomes the gradient operator $\nbs$ associated with the $(x,\s)$ coordinates, which is defined by
\begin{subequations}
\begin{equation} 
\nbs = (\nbs_\s, \nbs_1, \nbs_2) := \big( \eps \Qd \p_\s ,  \eps \p_1  ,  \p_2 -\Qb \p_\s \big) \,. 
\label{nb-s-H}
\end{equation} 
Additionally, the ALE transport operator $(\p_t + V \p_2)$ is written in $(x,\s)$ coordinates as 
\begin{equation} 
\Q \p_\s + \tilde V \p_2 = \Qd \p_\s + \tilde V \nbs_2 = \tfrac{1}{\eps} \nbs_\s + \tilde V \nbs_2 \,. 
\label{transport-s-H}
\end{equation} 
\end{subequations}

\subsection{The  approximate $1$-characteristic surfaces and the upstream weight function in $(x,\s)$ coordinates}  
Using the transformation \eqref{t-to-s-transform-H}, since $\JJ \colon \tHdm \to \mathbb{R}$ we may define upstream weight function in $(x,\s)$ coordinates as
\begin{equation} 
\tilde \JJ(x,\s) = \JJ(x,t) \,,
\label{eq:JJ:x:s:def}
\end{equation} 
for all $(x,\s) \in \Hdm$. Moreover, according to~\eqref{JJ-formula-t}, \eqref{qps-JJ-tHdmm:a}, and~\eqref{eq:xt:xs:chain:rule-H}, we have that $\tilde \JJ$ satisfies the equations
\begin{subequations}
\begin{align}
2 \alpha \Sigma  \p_1 \tilde \JJ 
- 2\alpha \tilde \Sigma \tilde \Jg  \tilde g^{- {\frac{1}{2}} } \nbs_2 \tilde h\,   \nbs_2 \tilde \JJ   
- \tilde \Jg (\Q\p_\s+\tilde V\p_2) \tilde \JJ 
&= 
- \dl \tilde \Jg (\Q \p_\s+\tilde V\p_2) \tilde \JJ 
\,,
\quad \mbox{in} \quad 
\Hdmp\,,
\label{JJ-formula}
\\
(\Q \p_\s+\tilde V\p_2) \tilde \JJ 
&= \mathfrak{f} 
\,,
\quad \mbox{in} \quad 
\Hdmm\,,
\end{align} 
\end{subequations}
where $\mathfrak{f}$ is as defined in \eqref{qps-JJ-H} is the right side of \eqref{qps-JJ-tHdmm:a}.

It is also convenient to obtain a pointwise expression for $\JJ$ in $(x,\s) \in \Hdmp$ coordinates that mimics \eqref{JJ-def-plus-t}. For this purpose, we need to introduce an $(x,\s)$ variant of the approximate $1$-characteristic surfaces $\Thd(x,t)$ which were defined earlier in~\eqref{eq:Thd:PDE}. It is important to note however that the domain of $\Thd$ is not $\tHdmp$, and instead the spacetime $\mathring{\Omega}_{\mathsf{US},+}$ defined in \eqref{eq:Omega:US:+}. This is because $\Thd(x,t)$ acts like a ``time'' variable itself, and so it should not be transformed according to~\eqref{eq:f:tilde:f-H}. Instead, the correct way to define an $(x,\s)$-variable analogue of $\Thd(x,t)$ is similar to~\eqref{t-to-s-transform-H}, and so is given by
\begin{equation}
\tThd(x,\s)= \mathfrak{q}(x_2,\Thd(x,t)) = \initial \cdot \tfrac{\tstar - \Thd(x,t)}{\tstar - \initial}
\,,
\qquad 
\mbox{where}
\qquad
\s = \mathfrak{q}(x_2,t)
\,.
\label{eq:Thd:x:s:def}
\end{equation}
In particular, the boundary condition~\eqref{Theta-BC-t} and definition~\eqref{t-to-s-transform-H} shows that
\begin{equation}
\tThd(\xringt,x_2,\s) =    \initial \cdot \tfrac{\tstar - \Thd(\xringt,x_2,t)}{\tstar - \initial}
 =   \initial \cdot \tfrac{\tstar -t}{\tstar - \initial} = \s\,,
 \label{eq:Thd:x:s:def:BC}
\end{equation}
for all $\sin \leq \s < 0$.
According to~\eqref{eq:Omega:US:+} and \eqref{eq:Thd:x:s:def}, the  natural domain of $\tThd$ is the spacetime
\begin{align}
\Omega_{\mathsf{US},+} 
&:= \left\{ (x,\s) \colon (x,\mathfrak{q}^{-1}(x_2,\s)) \in \mathring{\Omega}_{\mathsf{US},+} \right\}
\notag\\
&= \left\{(x,\s)  \colon x_2 \in \TT\,, \sin \leq \s < 0\,, \mathfrak{X}_1^-(x_2,\mathfrak{q}^{-1}(x_2,\s)) \leq x_1 \leq \mathfrak{X}_1^+(x_2,\mathfrak{q}^{-1}(x_2,\s)) \right\} \,,
\label{eq:tilde:Omega:DS:+}
\end{align}
Note importantly that $\tThd(x,\s)$ is only defined for $\s \in [\sin,0)$, and that the derivatives of $\tThd$ do not transform according to~\eqref{eq:xt:xs:chain:rule-H}, and instead according to
\begin{subequations}
\label{eq:xt:xs:chain:rule-Thd}
\begin{align}
 \p_\s  \tThd(x,\s) &= \p_t \Thd(x,t) \,,
\label{eq:xt:xs:chain:rule-Thd:a}\\
\p_2 \tThd(x,\s) &= \Qd(x_2) \Bigl( \p_2  \Thd(x,t) + \p_2 \tstar \tfrac{t-\initial}{\tstar - \initial} \bigl(\p_t \Thd(x,t)-1\bigr) -   \p_2 \tstar \tfrac{\Thd(x,t) - t}{\tstar-\initial} \Bigr)\,,
\label{eq:xt:xs:chain:rule-Thd:b}\\
\p_1 \tThd(x,\s) &=\Qd(x_2) \p_1 \Thd(x,t)
\label{eq:xt:xs:chain:rule-Thd:c}
\,,
\end{align}
\end{subequations}
Above, identity~\eqref{eq:xt:xs:chain:rule-Thd:a} is a consequence of the fact that $\Qd =\Qd(x_2)$ is independent of time (cf.~\eqref{eq:QQQ:a-H}), while in~\eqref{eq:xt:xs:chain:rule-Thd:b}--\eqref{eq:xt:xs:chain:rule-Thd:c} we have appealed also to the computation in~\eqref{eq:QQQ:aa-H}.

One of the consequences of definition~\eqref{eq:Thd:x:s:def}, when combined with~\eqref{JJ-def-plus-t} and~\eqref{eq:JJ:x:s:def}, is that for all $(x,\s) \in\Omega_{\mathsf{US},+} $ we have 
\begin{subequations}
\begin{align}
\tilde \JJ(x_1,x_2,\tThd(x_1,x_2,\s)) 
&= \tilde{\mathcal{B}}(x_2,\s)\,,
\label{JJ-def-plus}
\\
\tilde{\mathcal{B}}(x_2,\s) 
&= \mathcal{B}(x_2,t)\bigr|_{t = \mathfrak{q}^{-1}(x_2,\s)}
= \Jgb(\xringt,x_2,t)\bigr|_{t = \mathfrak{q}^{-1}(x_2,\s)}
\,.
\end{align}
\end{subequations}
In particular, for every $(x,\s) \in \Hdmp$, we may find $\s^\prime \in [\sin,0)$ such that $\s = \tThd(x,\s^\prime)$, and thus \eqref{JJ-def-plus} yields
\begin{equation*}
\tilde \JJ(x,\s) = \tilde \JJ(x,\tThd(x,\s^\prime)) = \tilde{\mathcal{B}}(x_2,\s^\prime)
\,. 
\end{equation*}

With the above notation, the distinguished surface passing through the pre-shock is parametrized as the graph  
\begin{equation*}
\bigl\{(x_1,x_2,\bar{\tilde{\Theta}^{\dl}}(x_1,x_2)) \colon x_2 \in \TT, x_1 \in [ \mathfrak{X}_1^-(x_2,\tstar) ,\mathfrak{X}_1^+(x_2,\tstar)  ]\bigr\}\,,
\end{equation*}
where
\begin{subequations}
\label{eq:bar:tilde:Theta:delta:def}
\begin{equation}
\bar{\tilde{\Theta}^{\dl}}(x_1,x_2):= \tThd(x_1,x_2 ,0) = \mathfrak{q}(x_2,\bar{\Thd}(x_1,x_2))\,,
\quad 
x_1\in [ \mathfrak{X}_1^-(x_2,\tstar) ,\mathfrak{X}_1^+(x_2,\tstar)  ]
\,,
\end{equation}
and $\mathfrak{X}_1^\pm$ were defined in~\eqref{eq:X:plus:stopping:time}--\eqref{eq:X:minus:stopping:time}. Abusing notation, we extend $\bar{\tilde{\Theta}^{\dl}}$ to be continuous and constant for $x_1 \not \in [ \mathfrak{X}_1^-(x_2,\tstar) ,\mathfrak{X}_1^+(x_2,\tstar)]$; that is, we define
\begin{align}
\bar{\tilde{\Theta}^{\dl}}(x_1,x_2) 
&:= \mathfrak{q}(x_2,\final) \leq \sfin \,,
\,\quad 
x_1 < \mathfrak{X}_1^-(x_2,\tstar) \,,
\\
\bar{\tilde{\Theta}^{\dl}}(x_1,x_2)
&:= \sin
\,,\quad
x_1 > \mathfrak{X}_1^+(x_2,\tstar) \,.
\label{eq:bar:tilde:Theta:delta:def:c}
\end{align}
\end{subequations}

As before, we shall need to represent the surface $\{(x,\bar{\tilde{\Theta}^{\dl}}(x)) \}$ defined via \eqref{eq:bar:tilde:Theta:delta:def} as a graph over the $(x_2,\s)$ plane, i.e.~as
$$
\bigl\{ (\tilde\theta^\dl(x_2,\s),x_2,\s) \colon x_2 \in \TT, \s \in [\sin,\sfin) \bigr\}\,,
$$
where for each $x_2$ fixed, $x_1 = \tilde\theta^\dl(x_2,\s)$ denotes the inverse of $\s = \bar{\tilde{\Theta}^{\dl}}(\cdot,x_2) = \mathfrak{q}(x_2,\bar\Thd(\cdot,x_2))$.  In particular, \begin{equation} 
\bar{\tilde{\Theta}^{\dl}}(\tilde\theta^\dl(x_2,\s),x_2)
=\s \qquad \text{ and } \qquad  
\tilde\theta^\dl(x_2,\bar{\tilde{\Theta}^{\dl}}(x_1,x_2))
=x_1 \,.
\label{little-theta}
\end{equation} 
It follows  from \eqref{JJ-def-plus} and the fact that $\tilde{\mathcal{B}}(x_2,0) = \mathcal{B}(x_2,\tstar) = \Jgb(\xringt,x_2,\tstar) = 0$, that 
\begin{align} 
\tilde\JJ (  \tilde\theta^\dl(x_2,\s), x_2, \s)=
\tilde\JJ \big(x_1,x_2,\bar{\tilde{\Theta}^{\dl}}(x_1,x_2)\big)=0 \,.
\label{JJ-vanish}
\end{align}

\begin{remark}[\bf Dropping the tildes]
\label{rem:US:drop:tilde}
Just as before, see Remarks~\ref{rem:no:tilde} and~\ref{rem:no:tilde-P}, we shall henceforth drop the use of the tilde-notation for all variables defined as functions of the flattened coordinates $(x,\s)$. Notably, besides dropping tildes on the fundamental Euler variables and the geometric variables, we shall denote $\tilde \Jg$, $\tilde \JJ$, and $\tThd$ simply as $\Jg$, $\JJ$, and $\Thd$, keeping the arguments as $(x,\s)$. We shall keep referring to the $\nbs$ operator defined in~\eqref{nb-s-H} as the rescaled spacetime differential operator in $(x,\s)$ coordinates. This identification is made throughout the rest of this section. 
\end{remark}

\subsection{The decomposition of the upstream spacetime in $(x,\s)$ coordinates} 
Recall that in $(x,t)$ coordinates, the upstream spacetime $\tHdm$ was decomposed according to \eqref{eq:H:dl:max:split:def}. In $(x,\s)$ coordinates, this decomposition carries over as follows. Using the notation introduced in~\eqref{eq:bar:tilde:Theta:delta:def}, the upstream spacetime $\Hdm$ defined in \eqref{eq:spacetime:H0} can be equivalently described as
\begin{subequations}
\begin{equation}
\Hdm = \left\{ (x, \s) \in \mathcal{X}_{\rm fin} \times [\sin , \sfin ) \colon  \s <  \bar{\Thd}(x_1,x_2)  \right\} \,.
 \label{eq:spacetime-Theta}
\end{equation}
In analogy to~\eqref{eq:H:dl:max:+:def}, we have a well-defined foliation of the spacetime subset $\Hdmp \subset \Hdm$, as defined in~\eqref{eq:spacetime:H0:plus}, which is given by 
\begin{equation} 
\Hdmp = \left\{ (x,\s) \in \Hdm \colon \Thd(x,\sin) < \s < \Thd(x,0) \right\} 
= \left\{ (x,\Thd(x,\s)) \colon (x,\s) \in \Omega_{\mathsf{US},+} \right\} \,.
\label{eq:spacetime-Theta:aa}
\end{equation} 
The  complimentary spacetime $\Hdmm$ defined in~\eqref{eq:spacetime:H0:minus} may be characterized similarly to~\eqref{eq:H:dl:max:-:def} as 
\begin{equation} 
\Hdmm = \left\{ (x,\s) \in \Hdm \colon \sin \leq \s < \Thd(x,\sin)\right\}  \,.
\label{eq:spacetime-Theta:bb}
\end{equation} 
Moreover, in analogy to~\eqref{eq:H:dl:max:split:def} we have that 
\begin{equation}
\Hdm = \Hdmp \cup \Thd(x,\sin)\cup \Hdmm 
\label{eq:spacetime-Theta:cc}
\,.
\end{equation}
\end{subequations}

While the characterization~\eqref{eq:spacetime-Theta} of $\Hdm$ is most useful for pointwise estimates in $(x,\s)$ coordinates, for energy estimates it is convenient to foliate $\Hdm$ by $\s$-time-slices. This $\s$-slice foliation was precisely the purpose for introducing the function $\theta^\dl(x_2,\s)$ (recall, we have dropped the tilde notation) in~\eqref{little-theta}. By construction, we have 
\begin{equation}
\Hdm = \bigcup_{\s \in [\sin,\sfin]}  \left\{(x_1,x_2,\s) \colon x_2 \in \TT, -\pi \leq x_1 < \theta^{\dl}(x_2,\s)  \right\} \,.
\label{eq:spacetime-Theta-new}
\end{equation}
In fact, the constraint that $x\in \mathcal{X}_{\rm fin}$ present in \eqref{eq:spacetime:H0} gives a more precise lower bound for $x_1$. In particular, according to \eqref{eq:driftin:3} we have $x_1 \geq \xringt - 2 (13 \pi + 65 \alpha(1+\alpha) \kappa_0) \eps$. Since in all our energy estimates the integrands vanish identically for $x\not \in \mathcal{X}_{\rm fin}$, there is no harm done by considering the $\s$-time-slices from~\eqref{eq:spacetime-Theta-new}.

\subsection{Notation for integrals and norms} 
Recalling the definition of $\thd$ given in \eqref{little-theta},
For $\s \in[\sin,\sfin] $ fixed, we use the following double-integral notation for integrals over the $\s$-time-slices of spacetime $\Hdm$, defined in~\eqref{eq:spacetime-Theta-new}:
\begin{equation*}
\dint F 
:= \dint F(x_1,x_2,\s) {\rm d} x_1{\rm d} x_2 
:= \int_{x_2=-\pi}^\pi \int_{x_1=-\pi}^{\thsd} F(x_1,x_2,\s) {\rm d} x_1 {\rm d} x_2 \,,
\end{equation*} 
and the triple integral notation to denote 
\begin{equation*}
\int_{\Hdm}F 
:= \tints F := \tints F(x_1,x_2,\s^\prime) {\rm d} x_1 {\rm d} x_2{\rm d} \s^\prime
:=  \int_{\s^\prime=\sin}^{\s} \int_{x_2=-\pi}^\pi   \int_{x_1=-\pi}^{\thd(x_2,\s^\prime)}F(x_1,x_2,\s^\prime){\rm d} x_1{\rm d} x_2{\rm d} \s^\prime \,.
\end{equation*} 
We then have the associated space and spacetime $L^2$-norms, respectively, defined as
\begin{subequations} 
\label{norms-H}
\begin{align}
\| F(\cdot ,\s)  \|_{L^2}^2 = \| F(\cdot ,\s)  \|_{L^2_x}^2 &:=  \dint F^2
\label{norms-H:a}
\,,\\
\| F  \|_{L^2_{x,\s}}^2 := \| F  \|_{L^2_{x,\s}(\Hdm)}^2 &:= \tint F^2 \,.
\label{norms-H:b}
\end{align}  
\end{subequations} 

Let us note that for integrable functions $F(x_1,x_2,\s)$ over $\Hdm$, the order of integration can be exchanged using the notation in~\eqref{eq:bar:tilde:Theta:delta:def}, so that
\begin{align} 
 \int_{\s=\sin}^{\sfin} \int_{x_2=-\pi}^\pi   \int_{x_1=-\pi}^{\thsd}F(x_1,x_2,\s){\rm d} x_1{\rm d} x_2{\rm d} \s
 = \int_{x_1=-\pi}^\pi  \int_{x_2=-\pi}^\pi \int_{\s=\sin}^{\bar\Thd(x_1,x_2)} F(x_1,x_2,\s) {\rm d} \s {\rm d} x_2  {\rm d} x_1 \,.
 \label{fubini-Thd}
\end{align} 
 
Let us also consider the surface $\Thd(x,\sin)$ (cf.~definition~\eqref{eq:Thd:x:s:def}) passing through $(\xringt,x_2,\sin)$.   While  $\s=\Thd(x,\sin)$ is a  graph over 
$(x_1,x_2)$, it is convenient to view this same surface as a graph over the  $(x_2,\s)$-plane.  In analogy to~\eqref{little-theta}, we let $x_1 = \thsin(x_2,\s)$ denote this surface, where $\thsin$ is defined as the inverse map
\begin{equation} 
\Thd(\thsin(x_2,\s),x_2,\sin)=\s \,,
\qquad\mbox{and}\qquad
\thsin(x_2, \Thd(x_1,x_2,\sin)) = x_1\,. \label{thsin-def}
\end{equation} 
Then, similarly to~\eqref{fubini-Thd}, we have that
\begin{align} 
 \int_{\s=\sin}^{\sfin} \int_{x_2=-\pi}^\pi   \int_{x_1=-\pi}^{\thsin(x_2,\s)}F(x_1,x_2,\s){\rm d}x_1 {\rm d} x_2 {\rm d}\s
 = \int_{x_1=-\pi}^\pi  \int_{x_2=-\pi}^\pi \int_{\s=\sin}^{\Thd(x_1,x_2,\sin)} F(x_1,x_2,\s) {\rm d}\s {\rm d}x_2  {\rm d}x_1 \,,
 \label{fubini-thsin}
\end{align} 
where in analogy to \eqref{eq:bar:tilde:Theta:delta:def} we have abused notation and by continuity defined
\begin{align*}
 \Thd(x_1,x_2,\sin) &:= \mathfrak{q}(x_2,\final) \leq \sfin\,, \qquad x_1 < \mathfrak{X}_1^-(x_2,\initial)\,, \\
  \Thd(x_1,x_2,\sin) &:= \sin\,, \qquad x_1 > \mathfrak{X}_1^+(x_2,\initial)\,,
\end{align*}
where the stopping times $\mathfrak{X}_1^{\pm}$ are defined in~\eqref{eq:X:plus:stopping:time}--\eqref{eq:X:minus:stopping:time}.

For $\s\in[\sin,\sfin]$ fixed, we will also make use of the following integral and norm notation:
\begin{subequations} 
\label{int-Hdmp-Hdmm}
\begin{align} 
\int_{\Hdmp}F &  :=   \int_{\s^\prime=\sin}^{\s} \int_{x_2=-\pi}^\pi   \int_{x_1=\thsin(x_2,\s^\prime)}^{\thd(x_2,\s^\prime)}F(x_1,x_2,\s^\prime){\rm d}x_1 {\rm d}x_2 {\rm d}\s^\prime \,, 
\label{int-Hdmp}\\
\int_{\Hdmm}F &  :=   \int_{\s'=\sin}^{\s} \int_{x_2=-\pi}^\pi   \int_{x_1=-\pi}^{\thsin(x_2,\s^\prime)}F(x_1,x_2,\s^\prime){\rm d}x_1 {\rm d}x_2 {\rm d}\s^\prime \,,
\label{int-Hdmm}
\end{align} 
\end{subequations} 
and in analogy to~\eqref{norms-H:b}, 
\begin{subequations} 
\label{norms-H+-}
\begin{align}
\| F  \|_{L^2_{x,\s,+}}^2 &:= \| F  \|_{L^2_{x,\s}(\Hdmp)}^2 :=  \int_{\Hdmp}F^2 \,, \\
\| F  \|_{L^2_{x,\s,-}}^2 &:= \| F  \|_{L^2_{x,\s}(\Hdmm)}^2 := \int_{\Hdmm}F^2 \,.
\end{align}  
\end{subequations} 

\subsection{$L^2$ adjoint of $\nbs$  in the upstream spacetime $\Hdm$} 
In computing adjoints (with respect to the $L^2$ inner product on $\Hdm$), we make use of the following calculus identities:
\begin{subequations} 
\label{calc-dsd2}
\begin{equation} 
\int_{-\pi}^\thsd \p_\s f(x_1,x_2,\s){\rm d}x_1 
= \frac{\p}{\p\s} \int_{-\pi}^\thsd  f(x_1,x_2,\s) {\rm d}x_1 
- \p_\s\thsd \,  f(\thsd,x_2, \s) \,. \label{calc-ds}
\end{equation} 
Similarly, we have that
\begin{align} 
\int_{-\pi}^\thsd \p_2 f(x_1,x_2,\s){\rm d} x_1 &= \frac{\p}{\p x_2} \int_{-\pi}^\thsd  f(x_1,x_2,\s){\rm d}x_1 
- \p_2\thsd \,  f(\thsd,x_2, \s)
 \,, \label{calc-d2}
\end{align} 
\end{subequations} 
where $\thsd$ is as defined in~\eqref{little-theta}.
With \eqref{calc-dsd2}, the adjoint of $\nbs$ with respect to the $L^2$ inner product  on $\Hdm \cap (\mathcal{X}_{\rm fin} \times [\sin,\s])$, with $\sin < \s <\sfin$, is given by
\begin{subequations}
\label{eq:adjoints-H}
\begin{align}
\nbs_\s^* &=  - \nbs_\s  + \eps \Qd ( \delta_{\s^\prime = \s} - \delta_{\s^\prime = \sin} ) - \nbs_\s\theta^\dl \delta_{x_1=\theta^{\!\dl}\!(x_2,\s^\prime)}  \,, \label{adjoint-s-H}
\\
\nbs_1^* &=  - \nbs_1 +\eps \delta_{x_1=\theta^{\!\dl}\!(x_2,\s^\prime)}   \,,  \label{adjoint-1-H}
\\
\nbs_2^* &=  - \nbs_2 + \Qr_2 
-  \Qb \delta_{\s^\prime = \s}  - \nbs_2\theta^{\!\dl} \delta_{x_1=\theta^{\!\dl}\!(x_2,\s^\prime)} 
\,,
 \label{adjoint-2-H}
\\
(\Q \p_\s + \ V \p_2)^* & = - (\Q \p_\s +  V \p_2)  - (\Qc + V,_2) + \Q  ( \delta_{\s^\prime = \s} - \delta_{\s^\prime = \sin} )  
- (\Q\p_\s \theta^\dl +V\p_2 \theta^\dl )\delta_{x_1=\theta^{\!\dl}\!(x_2,\s^\prime)}
 \,.
  \label{adjoint-3-H}
\end{align}
\end{subequations}
Here we have used that $\Hdm \subset \mathcal{X}_{\rm fin} \times[\sin,\sfin)$, so that only one boundary term emerges when integrating by parts with respect to $x_1$\footnote{See the footnote \ref{footnote-trace}.}, at $x_1 = \thd(x_2,\s)$.
We have also used the fact that \eqref{eq:QQQ:aa-H} implies $\Qb(x_2,\sin) = 0$.

\subsection{The $L^2$-based energy norms} 
\label{sec:norms:L2:H}
Using the weight function $\JJ$ defined by~\eqref{JJ-def-plus-t}, \eqref{JJ-def-minus-t}, and~\eqref{eq:JJ:x:s:def}, and with the notation introduced in~\eqref{norms-H:a} for the $L^2$ norm, for upstream \MGHDB\ we work with the energy norms defined by  
\begin{subequations}
\label{eq:tilde:E5E6-H} 
\begin{alignat}{2}
\widetilde{\mathcal{E}}_{6}^2(\s) 
&=  \widetilde{\mathcal{E}}_{6,\nnn}^2(\s)  + (\mathsf{K}\eps)^{-2} \widetilde{\mathcal{E}}_{6,\ttt}^2(\s) \,,
\ \ \ 
&&\widetilde{\mathcal{E}}_5^2(\s) 
=  \widetilde{\mathcal{E}}_{5,\nnn}^2(\s) + (\mathsf{K}\eps)^{-2} \widetilde{\mathcal{E}}_{5,\ttt}^2(\s)
\label{eq:tilde:E5E6:N+T-H} 
\\
\widetilde{\mathcal{E}}_{6,\nnn}^2(\s) 
&= \snorm{   \JJtf \Jgh \nbs^6 (\Jg\Wbn,\Jg\Zbn, \Jg\Abn)(\cdot,\s)}^2_{L^2}  \,,
\ \ \ 
&&\widetilde{\mathcal{E}}_{5,\nnn}^2(\s)
= \snorm{ \Jgh \nbs^5 (\Jg\Wbn, \Jg\Zbn, \Jg\Abn)(\cdot,\s)}^2_{L^2}
 \\
\widetilde{\mathcal{E}}_{6,\ttt}^2(\s)
&=\snorm{ \JJtf \Jgh \nbs^6 (\Wbt, \Zbt, \Abt)(\cdot,\s)}^2_{L^2}  \,, 
\ \ \
&&\widetilde{\mathcal{E}}_{5,\ttt}^2(\s)
= \snorm{ \Jgh \nbs^5 ( \Wbt, \Zbt, \Abt)(\cdot,\s)}^2_{L^2} \,, 
\end{alignat}
\end{subequations}
and the damping norms by
\begin{subequations}
\label{eq:tilde:D5D6-H}
\begin{alignat}{2}
\widetilde{\mathcal{D}}_6^2(\s)
&= \widetilde{\mathcal{D}}_{6,\nnn}^2(\s) + (\mathsf{K} \eps)^{-2} \widetilde{\mathcal{D}}_{6,\ttt}^2(\s)  \,,
\ \ \
&&\widetilde{\mathcal{D}}^2_5(\s)
 = \widetilde{\mathcal{D}}^2_{5,\nnn}(\s) + (\mathsf{K}\eps)^{-2} \widetilde{\mathcal{D}}^2_{5,\ttt}(\s) 
 \,,
 \label{eq:tilde:D5D6:N+T-H}  \\
\widetilde{\mathcal{D}}_{6,\nnn}^2(\s) 
&=   
\int_\sin^\s\snorm{ \JJof\Jgh \nbs^6 (\Jg \Wbn, \Jg\Zbn, \Jg\Abn)}^2_{L^2} {\rm d}\s'  
\ \ \ 
&&\widetilde{\mathcal{D}}^2_{5,\nnn}(\s)
 =    \int_\sin^\s  \snorm{ \nbs^5 (\Jg \Wbn, \Jg\Zbn, \Jg\Abn)}^2_{L^2} {\rm d}\s' 
  \,, \notag \\
 &\quad\quad +   
\int_\sin^\s\snorm{ \JJtf \nbs^6 (\Jg \Wbn, \Jg\Zbn, \Jg\Abn)}^2_{L^2} {\rm d}\s'  \,, 
\ \ \ 
&& 
  \\
\widetilde{\mathcal{D}}_{6,\ttt}^2(\s) 
&=  \int_\sin^\s \snorm{ \JJof \Jgh \nbs^6 (\Wbt, \Zbt, \Abt)}^2_{L^2} {\rm d}\s'  
\ \ \
&&\widetilde{\mathcal{D}}^2_{5,\ttt}(\s)
=   \int_\sin^\s\snorm{ \nbs^5 (\Wbt,  \Zbt, \Abt)}^2_{L^2}{\rm d}\s'  \,.
\notag \\
&\quad \quad
+   
\int_\sin^\s\snorm{ \JJtf \nbs^6 ( \Wbt, \Zbt,\Abt)}^2_{L^2} {\rm d}\s'  \,, 
\ \ \ 
\end{alignat}
\end{subequations}
Once again $\mathsf{K} \geq 1$ is a sufficiently large constant, independent of $\eps$,  chosen at the end of the 
proof, solely in terms of  $\alpha$ and $\kappa_0$    (see~\eqref{eq:K:choice:1-H} below).

\subsection{Upstream bootstrap assumptions} 
We continue to use the same bootstrap assumptions as in~\eqref{bootstraps}, but instead of assuming that these bootstraps hold 
for $(x,t)$ in the spacetime $\mathcal{P}$ (cf.~\eqref{eq:spacetime:smooth}), we now assume that these bootstraps hold for $(x,t) \in \tHdm $ (cf.~definition~\eqref{eq:spacetime-Theta-t}). As such, in this section all pointwise bootstraps are assumed to hold for $(x,t) \in \tHdm $, or equivalently, for all  $(x,\s) \in  \Hdm$ via the flattening map \eqref{t-to-s-transform-H}, and for the energy and damping norms defined earlier in~\eqref{eq:tilde:E5E6-H} and~\eqref{eq:tilde:D5D6-H}.

To be more precise, the working bootstrap assumptions in this section are that 
\begin{subequations}
\label{bootstraps-H}
\begin{align}
&( {\Wb}, {\Zb}, {\Ab}, {\Jg}, {h}, {V}, {\Sigma}) 
  \mbox{ satisfy the pointwise  bootstraps }   
\eqref{bs-supp}\!-\!\eqref{bs-D-Sigma} \mbox{ in } \tHdm \label{boots-HH}
\,,\\  
&\widetilde{\mathcal{E}}_{6} ,
\widetilde{\mathcal{D}}_{6} ,
\widetilde{\mathcal{E}}_{5} ,
\widetilde{\mathcal{D}}_{5} ,
\snorm{\nbs^6 \nbs_1   h}_{L^2_{x,\s}},
\snorm{\nbs^6 \nbs_2   h}_{L^2_{x,\s}}, 
\snorm{\nbs^6   \Jg}_{L^2_{x,\s}} 
  \mbox{ satisfy the energy  bootstraps }   
\eqref{bootstraps-Dnorm:6}\!-\!\eqref{bootstraps-Dnorm:Jg}  \label{boots-H} \,.
\end{align}
\end{subequations}
Here $( {\Wb}, {\Zb}, {\Ab}, {\Jg}, {h}, {V}, {\Sigma})$ are defined according to the 
flattening~\eqref{eq:f:tilde:f-H} (dropping the tildes as discussed in~Remark~\ref{rem:US:drop:tilde}), while the energy and damping norms are defined in~\eqref{eq:tilde:E5E6-H} and~\eqref{eq:tilde:D5D6-H}, 
respectively. Since the bootstraps~\eqref{bootstraps-H} in this section  are the same as the bootstraps~\eqref{bootstraps} used in Sections~\ref{sec:formation:setup}--\ref{sec:sixth:order:energy}, save for the different weights in the $L^2$ norms (with $\JJ$ replacing $\mathcal{J}$ in~\eqref{eq:tilde:E5E6-H} and~\eqref{eq:tilde:D5D6-H}), we shall sometimes (more frequently for the pointwise bootstraps) make reference to \eqref{bootstraps} instead of \eqref{bootstraps-H}. 

As in Sections~\ref{sec:formation:setup}--\ref{sec:sixth:order:energy}, the burden of the proof in the current section is to close the bootstrap assumptions~\eqref{bootstraps-H}, i.e., to show that these bounds hold with $<$ symbols instead of $\leq$ symbols. To avoid redundancy, we do not repeat the arguments of how bootstraps are closed when the proof is either identical to that in given earlier in Sections~\ref{sec:formation:setup}--\ref{sec:sixth:order:energy}, or if it requires infinitesimal and straightforward adjustments. Instead, we focus on the proofs of those bootstraps which are different due to the upstream spacetime and weights.  In particular, the upstream weight function $\JJ(x_1,x_2,\s)$ that appears in our energy
estimates is defined using a family of surfaces $\Thd(x,\s)$ which closely approximate the slow acoustic characteristic surfaces.   A key feature
of this weight function is the identity \eqref{JJ-formula} which is fundamental to the upstream geometry.
The remainder of this section is dedicated to closing the bootstrap assumptions~\eqref{bootstraps-H}.

\subsection{Identities in the downstream coordinate system}
With respect to the coordinates $(x,\s)$ given by \eqref{t-to-s-transform-H}, with the transformation~\eqref{eq:f:tilde:f-H}, and upon dropping the 
tildes (see Remark~\ref{rem:US:drop:tilde}), we have the following fundamental identities, which are translations of the identities in Section~\ref{sec:new:Euler:variables}  into $(x,\s)$ coordinates (see also~\eqref{Jg-evo-s}--\eqref{vort-s}):
\begin{subequations} 
 \label{fundamental-H}
\begin{align} 
(\Q\p_\s+V\p_2) \Jg &= \tfrac{1+ \alpha }{2} \Jg\Wbn + \tfrac{1- \alpha }{2} \Jg\Zbn \,, \label{Jg-evo-s-P-US} \\ 
(\Q\p_\s+V\p_2) \nbs_2 h &=  g  \bigl(\tfrac{1+ \alpha }{2} \Wbt + \tfrac{1- \alpha }{2} \Zbt \bigr) \,, \label{p2h-evo-s-P-US} \\
{\tfrac{1}{\eps}} \nbs_1 \Sigma &=  \tfrac{1}{2} \Jg(\Wbn -\Zbn)  + \tfrac{1}{2} \Jg \nbs_2 h (\Wbt -\Zbt) \,, \label{p1-Sigma-s-P-US} \\
\nbs_2 \Sigma &= \tfrac{1}{2} g^{\frac{1}{2}}  (\Wbt-\Zbt)\,,  \label{p2-Sigma-s-P-US} \\
(\Q\p_\s+V\p_2) \Sigma  &= -  \alpha \Sigma (\Zbn+ \Abt )  \,,  \label{Sigma0-ALE-s-P-US}  \\
(\Q \p_\s +    V \p_2) \Sigma^{-2\beta}   &=   2\alpha\beta \Sigma^{-2\beta}  (\Zbn+ \Abt )  \,,  \label{Sigma0i-ALE-s-P-US} \\
(\Q\p_\s+V\p_2)\nn  & =-  \bigl(\tfrac{1+ \alpha }{2} \Wbt + \tfrac{1- \alpha }{2} \Zbt \bigr) \tt   \,, \label{nn-evo-s-P-US} \\
(\Q\p_\s+V\p_2)\tt  & =   \bigl(\tfrac{1+ \alpha }{2} \Wbt + \tfrac{1- \alpha }{2} \Zbt \bigr) \nn   \,, \label{tt-evo-s-P-US} \\
\tfrac{\Jg}{\Sigma}  (\Q\p_\s+V\p_2)  \Upomega & =      {\tfrac{\alpha }{\eps}}  \nbs_1 \Upomega
- \alpha \Jg  g^{- {\frac{1}{2}} } \nbs_2h \   \nbs_2 \Upomega   \,. \label{vort-s-P-US}
\end{align} 
\end{subequations} 
By using \eqref{eq:Jgb:identity:2}, we also note that $\nbs_1 (\Jg-\Jgb) = \nbs_2 (\Jg-\Jgb)=0$.  

\subsection{Bounds for $\Q$, $\Qd$, $\Qb$, $\Qr_2$, and $\Qc$}
We list a few properties of the coefficients defined in~\eqref{QQQ-H} in conjunction with the flattening map $ \s = \mathfrak{q}(x_2,t) $.
\begin{lemma}
\label{lem:Q:bnds-H}
Assume that the bootstrap bounds \eqref{bootstraps-H} hold on $\Hdm$.
If $\eps$ is taken to be sufficiently small with respect to $\alpha, \kappa_0$, and $\Cdata$, then  
\begin{subequations}
\label{eq:Q:all:bbq-H}
\begin{align} 
\tfrac{50}{51}      &\le \Qd(x_2)  \le  \tfrac{1001}{1000}   \,,  
\label{Qd-lower-upper-H} \\ 
\sabs{\Q(x,\s) - \Qd(x_2)} &\leq \Cn \eps^2 \,,  
\label{Q-lower-upper-H} \\
\sabs{\Qb(x_2,\s)} &\leq 11 \eps\,, 
\label{eq:Qb:bbq-H} \\
 \sabs{\Qr_2(x_2)}  &\leq 6(1+\alpha) \,, 
\label{eq:Qr2:bbq-H} 
\\
\sabs{\Qc(x,\s)} &\leq \Cn \eps
\,, 
\label{eq:Qc:bbq-H} 
\\
\sabs{\p_2\Qd(x_2) } &\leq 7(1+\alpha)
\,, 
\label{p2-Qd-H} 
\\
\sabs{\p_1\Q(x,\s) } &\leq \Cn \eps
\,, 
\label{p1-Q-H} 
\\\sabs{\p_2\Q(x,\s) } &\leq \Cn
\,, 
\label{p2-Q-H}
\end{align} 
hold uniformly for all $(x,\s) \in \Hdm$.  
\end{subequations}
\end{lemma}
\begin{proof}[Proof of Lemma \ref{lem:Q:bnds-H}]
Using the definition \eqref{eq:QQQ:a-H} together with the fact that 
$ -\Cn \eps^2 - \initial  <
\tstar - \initial \le \final-\initial
$, which follows from the analysis in Section \ref{sec:x1star}, 
we have that
\begin{equation} 
\tfrac{50}{51}
= \tfrac{-\initial}{\final-\initial}   
\le \Qd \le
\tfrac{-\initial}{-\Cn\eps^2-\initial}  
\le 1+ \Cn \eps^2  \,,
\label{tstar-minus-tin}
\end{equation} 
so by choosing $\eps$ sufficiently small,  \eqref{Qd-lower-upper-H} holds.
From the definition \eqref{eq:QQQ:aa-H} and the bounds \eqref{p2-tstar} and  \eqref{tstar-minus-tin}, we deduce
\begin{equation*}
|\Qb(x_2,\s)| \leq 10 \eps \cdot (1 + \Cn \eps^2) \cdot  \sabs{1 - \tfrac{\s}{\sin}} \leq 10 \eps \cdot (1 + \Cn \eps^2) \cdot (\tfrac{51}{50} + \Cn \eps)\,,
\end{equation*}
and so the inequality
\eqref{eq:Qb:bbq-H} immediately follows.   Now, since $\Q = \Qd - V \Qb$, from \eqref{eq:Qb:bbq-H}, we obtain the bound 
\eqref{Q-lower-upper-H}.
Using the definition \eqref{eq:QQQ:c-H}, we see that $\Qc = - \p_\s V \, \Qb + V \Qd \tfrac{\p_2\tstar}{\sin} $.  The bounds
\eqref{bs-V}, \eqref{p2-tstar}, \eqref{tstar-minus-tin}, and \eqref{eq:Qb:bbq-H} show that 
\eqref{eq:Qc:bbq-H} holds.
Next, we have that $\Qr_2 = - \Qd \tfrac{\p_2 \tstar}{\sin} $ so that  with \eqref{p2-tstar} and  \eqref{tstar-minus-tin}, we see that \eqref{eq:Qr2:bbq-H} holds. Since $\p_2\Qd = - \Qd \Qr_2 $, we also have that \eqref{p2-Qd-H} holds.
Finally, from \eqref{eq:QQQ:b-H}, $\p_1\Q =  - {\tfrac{1}{\eps}} \nbs_1 V \,  \Qb$.  Hence, \eqref{bs-V} and \eqref{eq:Qb:bbq-H}
give \eqref{p1-Q-H}.  Similarly,  from \eqref{eq:QQQ:b-H}, $\p_2\Q = \p_2 \Qd  - \p_2 V    \Qb - V \p_2 \Qb$. From \eqref{bs-V}, \eqref{p2-tstar}, \eqref{p22-tstar}, \eqref{eq:Qb:bbq-H}, 
 \eqref{eq:Qr2:bbq-H}, and \eqref{p2-Qd-H},  \eqref{p2-Q-H} follows.  
\end{proof}

\subsection{Bounds for $\Thd$ and $\JJ$ in $(x,\s)$ coordinates} 
In Sections~\ref{sec:Thd:props} and~\ref{sec:JJ:US:properties} we have obtained very precise bounds for $\Thd(x,t)$ (and its derivatives), and for $\JJ(x,t)$ (and its derivatives), in the spacetime $\tHdm$. In this subsection we revisit a few of these estimates and re-state them (using~\eqref{eq:Thd:x:s:def} and \eqref{eq:JJ:x:s:def}), for $\Thd(x,\s)$ and $\JJ(x,\s)$ in the spacetime $\Hdm$.

\subsubsection{Bounds for the derivatives of $\Thd$}
From Lemma~\ref{lem:Thd:derivs} and the identities in~\eqref{eq:xt:xs:chain:rule-Thd}, we immediately derive the first derivative bounds  
\begin{subequations}
\begin{align}
\tfrac{999}{1001} &\leq \p_\s \Thd(x,\s) \leq \tfrac{1001}{999}
\,,
\label{ps-Theta}
\\
-\tfrac{2}{10^5 (1+\alpha)} &\leq \p_1 \Thd(x,\s) <0 
\,,
\label{p1-Theta-sign}
\\
|\p_2 \Thd(x,\s)| &\leq 6 \cdot 10^3 (1+\alpha)^2 \eps
\,,
\label{boot-p2Theta}
\end{align} 
\end{subequations}
for all $(x,\s) \in \Omega_{\mathsf{US},+}$, the spacetime defined in~\eqref{eq:tilde:Omega:DS:+}. Indeed, the bound~\eqref{ps-Theta} follows from~\eqref{boot-psTheta-t} and \eqref{eq:xt:xs:chain:rule-Thd:a}, the estimate~\eqref{p1-Theta-sign} is a consequence of~\eqref{p1-Theta-sign-t-new}, \eqref{eq:xt:xs:chain:rule-Thd:c}, and~\eqref{Qd-lower-upper-H}, while~\eqref{boot-p2Theta} follows from~\eqref{boot-p2Theta-t}, \eqref{boot-psTheta-t}, \eqref{eq:xt:xs:chain:rule-Thd:b}, and \eqref{p2-tstar}.

Bounds for the second order derivatives of $\Thd(x,\s)$ also follow directly from~\eqref{boot-p22Theta-t}--\eqref{boot-pssTheta-t} and \eqref{eq:jesus:just:left:chicago:a}--\eqref{eq:jesus:just:left:chicago:c} upon differentiating~\eqref{eq:xt:xs:chain:rule-Thd} one more time, and appealing to~\eqref{p2-tstar}, \eqref{p22-tstar}, \eqref{Qd-lower-upper-H}, and~\eqref{p2-Qd-H}. Since these bounds will not be crucially used in the subsequent analysis, we choose not to re-state these bounds.

\subsubsection{The function $\thd$ is strictly decreasing in $\s$.}
We recall the definition of $\thd(x_2,\s)$ from~\eqref{little-theta}. Differentiating the first identity in \eqref{little-theta} with respect to $\s$, and recalling that $\bar \Thd(x) = \Thd(x,0)$, we deduce that 
\begin{equation*}
\p_1 \Thd(\thd(x_2,\s),x_2,0) \cdot \p_\s \thd(x_2,\s) = 1 \,.
\end{equation*}
This identity, combined with the bound \eqref{p1-Theta-sign-t}, \eqref{eq:US:kappa:0:cond:0}, \eqref{eq:xt:xs:chain:rule-Thd:c}, and~\eqref{Qd-lower-upper-H}, yields
\begin{align} 
- \tfrac{4\alpha \kappa_0}{\Jg(\thd(x_2,\s),x_2,\s)} \leq \p_\s\thd(x_2,\s) \leq - \tfrac{10^5(1+\alpha)}{2} < 0\,. 
\label{ps-theta-sign}
\end{align} 
Similarly, for the function~$\thsin$ defined in~\eqref{thsin-def}, and appealing to the bound~\eqref{JJ-le-Jg:c}, we also have 
\begin{align} 
- 5 \alpha \kappa_0 \leq - \tfrac{4\alpha \kappa_0}{\Jg(\thsin(x_2,\s),x_2,\s)} \leq \p_\s\thsin(x_2,\s) = \tfrac{1}{\p_1 \Thd(\thsin(x_2,\s),x_2,\sin)} \leq - \tfrac{10^5(1+\alpha)}{2} < 0\,. \label{ps-thsin-sign}
\end{align} 
In view of the definition of $\nbs^*$ in~\eqref{eq:adjoints-H}, it is also convenient to record the estimate
\begin{equation}
\sabs{\p_2 \thd(x_2,\s)} \leq \tfrac{24 \cdot 10^3 (1+\alpha)^2 \alpha \kappa_0 }{\Jg(\thd(x_2,\s),x_2,\s)} \eps
\,,
 \label{p2-ths-bound}
\end{equation}
which follows upon differentiating the first identity in~\eqref{little-theta} with respect to $x_2$, and appealing to~\eqref{boot-p2Theta}.

\subsubsection{Lower bounds for $\JJ$}
We recall that the function $\JJ(x,t)$ satisfies the pointwise lower bounds~\eqref{JJ-and-t} and~\eqref{JJ-and-t-new} in $\tHdm$, and respectively $\tHdmp$. Via the definition of $\JJ(x,\s)$ in \eqref{eq:JJ:x:s:def}, appealing to the definitions~\eqref{qinv-H} and \eqref{eq:bar:tilde:Theta:delta:def}, and to the lower bound for $\Qd$ in~\eqref{Qd-lower-upper-H}, these pointwise lower bounds in~\eqref{JJ-and-t} become
\begin{equation}
\JJ(x,\s) \geq 
\begin{cases}
\bigl( \bar\Thd(x) - \s\bigr) \tfrac{1+\alpha}{2\eps}  \cdot \tfrac{22}{25}\,, &(x,\s) \in \Hdmp\,,\\
1\,, &(x,\s) \in \Hdmm\,.
\end{cases}
\label{JJ-and-s-new}
\end{equation}
By additionally appealing to~\eqref{JJ-def-plus}, \eqref{t-to-s-transform-H}, \eqref{eq:QQQ:a-H}, and~\eqref{Qd-lower-upper-H}, the lower bound~\eqref{JJ-and-t-new} becomes
\begin{align} 
\JJ(x,\Thd(x,\s))
\geq 
\tfrac{895}{1000} \Qd(x_2)^{-1} \cdot\tfrac{\s}{\sin}
\geq \tfrac{89}{100} \cdot\tfrac{\s}{\sin} \,,
 \label{JJ-and-s}
\end{align}  
for all $(x,\s) \in \Omega_{\mathsf{US},+}$, or equivalently, for all $(x,\Thd(x,\s)) \in \Hdmp$.

\subsubsection{Comparison of $\JJ$ and $\Jg$}
We only record the fact that the results of Lemma~\ref{lem:JJ-le-Jg} carry over directly to $(x,\s)$ variables. In particular, \eqref{JJ-le-Jg} implies that
\begin{equation}
\JJ(x,\s) \leq \tfrac{101}{100} \Jg(x,\s) {\bf 1}_{|x_1|\leq 13\pi\eps} + \tfrac{21}{10} \Jg(x,\s) {\bf 1}_{|x_1|>13 \pi \eps} 
\label{JJ-le-Jg-s}
\end{equation}
for all $(x,\s) \in \Hdm$. 

\subsubsection{Damping properties of $\JJ$}
Using~\eqref{eq:JJ:x:s:def}, \eqref{eq:xt:xs:chain:rule-H}, \eqref{transport-s-H}, and \eqref{eq:Q:all:bbq-H}, it is direct to transfer any derivative information for $\JJ$ from $(x,t)$ variables, into $(x,\s)$ variables. 
In particular, from \eqref{eq:waitin:for:the:bus} we directly deduce
\begin{subequations}
\label{eq:waitin:for:the:bus-s}
\begin{equation}
- (\Q \p_\s + V \p_2) \JJ(x,\s) \geq \tfrac{1+\alpha}{2 \eps} \cdot \tfrac{899}{1000} \geq \tfrac{11(1+\alpha)}{25 \eps}
 \label{qps-JJ-bound}
\end{equation}
for all $(x,\s) \in \Hdm$, and 
\begin{align}
-\tfrac{412(1+\alpha)}{\eps} \leq  \p_\s \JJ(x,\s)  &\leq - \tfrac{9(1+\alpha)}{25\eps}\,,
\label{mod-ps-JJ}
\\
\sabs{\nbs_2 \JJ(x,\s)} &\leq  \Cn \,,
\label{fat-marmot2b}
\\
\sabs{\nbs_1 \JJ(x,\s)} &\leq 10^{-3}\,,
\label{fat-marmot2}
\\
\sabs{(\Q \p_\s + V \p_2)^2 \JJ(x,\s)}  
&\leq \tfrac{202^2(1+\alpha)^2}{\eps^2}\,,
\label{qps2JJ-bound}
\end{align}
\end{subequations}
for all $(x,\s) \in \Hdmp \cup \Hdmm$.
Lastly, from Lemma~\ref{lem:damping:anti:damping-H-t} we obtain that 
\begin{equation}
-\tfrac 32 \JJh \Jg  (\Q \p_\s +V\p_2) \JJ+   \JJss(\Q\p_\s + V\p_2)   \Jg 
\ge \tfrac{1+\alpha}{6 \eps} \JJh\Jg     \,,
\label{eq:fakeJg:LB-H}
\end{equation}
pointwise for all $(x,\s) \in \Hdm$.

\subsection{Bounds for the geometry, sound speed, and ALE velocity for upstream development}
Just as we did in Section \ref{sec:downstreammaxdev}, we again
record all necessary upstream  modifications to the bounds obtained earlier in Section~\ref{sec:geometry:sound:ALE}.
We specifically highlight all the modification caused by the change in the weight function $\mathcal{J} \mapsto \JJ$.   
Bounds without the presence of such a weight function
remain identical to those in Section~\ref{sec:geometry:sound:ALE} (and we continue to make reference to equation numbers from Section~\ref{sec:geometry:sound:ALE}), and  bounds with such a weight function are modified with  $\mathcal{J}$ replaced with $\JJ$. 
For instance, the bounds in Proposition~\ref{prop:geometry} now become the bounds given in Proposition~\ref{prop:geometry-H} below. The corollaries and remarks which follow this proposition (in particular, the closure of the \eqref{bootstraps-Dnorm:h2} and~\eqref{bootstraps-Dnorm:Jg} bootstraps in Corollary~\ref{cor:Bj:Bh}) remain the same as in Section~\ref{sec:geometry:sound:ALE}, 
and to avoid redundancy we do not repeat those arguments here.
\begin{proposition}[Bounds for the geometry, sound speed, and ALE velocity]
\label{prop:geometry-H}
Assume that  the bootstrap assumptions \eqref{bootstraps-H} hold, and that $\eps$ is taken to be sufficiently small to ensure $\eps^{\frac 12}  ( \brak{\mathsf{B_J}} +  \brak{\mathsf{B_h}} +     \brak{\mathsf{B_6}} )  \leq 1$.  
Then, assuming $\eps$ is sufficiently small with respect to $\alpha,\kappa_0,$ and $\Cdata$, we have that 
\begin{subequations}
\label{geom-H-original}
\begin{align} 
 (\Jg, \nbs_1 h, \nbs_2 h, \Sigma, V) \text{ satisfy the bounds } \eqref{eq:signed:Jg},  \eqref{D5JgEnergy}, \eqref{D6h1Energy:new},  
\eqref{D5h2Energy}, \eqref{eq:Sigma:H6:new}\!-\!\eqref{eq:Sigma:H6:new:bdd}, & \eqref{eq:V:H6:new}\!-\!\eqref{eq:V:H6:new:bdd}, 
\\
  \snorm{ \JJof \nbs^6 \Jg(\cdot,\s) }^2_{L^\infty_\s L^2_{x}} 
+ \tfrac{1}{\eps}  \snorm{  \JJmof \nbs^6 \Jg }_{L^2_{x,\s}}^2 
&\les  \eps  \brak{\mathsf{B}_6}^2    \,,   \label{D6JgEnergy:new-H}
\\
 \| \JJof \nbs^{6} \nbs_2h(\cdot, \s) \|^2_{L^\infty_\s L^2_{x}} 
+ \tfrac{1}{\eps}  \|\JJmof \nbs^{6}  \nbs_2h \|_{L^2_{x,\s}}^2
&\les 
\mathsf{K}^2 \eps^3 \brak{\mathsf{B}_6}^2 \,,
\label{D6h2Energy:new-H}
\\
\sum_{ |\gamma|=3}^6 
 \|\JJmof\bigl(\nbs^{|\gamma|}  \nn + g^{-1}   \tt \nbs^{|\gamma|} \nbs_2 h\bigr) \|_{L^2_{x,\s}}
+
\|\JJmof \bigl( \nbs^{|\gamma|}  \tt - g^{-1}  \nn \nbs^{|\gamma|} \nbs_2 h \bigr) \|_{L^2_{x,\s}}
&\les \mathsf{K} \eps^3 \brak{\mathsf{B}_6}  \,,
\label{D6n-bound:a:new-H}
\\
\| \JJmof  \nbs^6  \nn  \|_{L^2_{x,\s}}
+
\|\JJmof  \nbs^6 \tt \|_{L^2_{x,\s}}
&\les \mathsf{K} \eps^2 \brak{\mathsf{B}_6}   \,,
\label{D6n-bound:b:new-H}
\\
\big( \|\JJof \nbs^6  \nn  \|_{L^\infty_\s L^2_{x}} 
+
\|\JJof  \nbs^6 \tt \|_{L^\infty_\s L^2_{x}} \big)
&\les \mathsf{K} \eps^{\frac 32} \brak{\mathsf{B}_6}   \,,
\label{D6n-bound:b:new:bdd-H}
\end{align} 
where the implicit constants in all the above inequalities depend only on $\alpha$, $\kappa_0$,  and $\Cdata$. 
\end{subequations}
\end{proposition}
\begin{proof}[Proof of Proposition~\ref{prop:geometry-H}]
We explain the upstream modifications required for the proof of the inequality  \eqref{D6JgEnergy:new-H}. 
We compute the $L^2$-inner product of \eqref{D5-Jg-s} with $\JJh \nbs^6\Jg$ to obtain that for any $\s \in (\sin,\sfin]$, 
\begin{align*} 
&\tfrac{1}{2} \tints \JJh (\Q\p_\s +V\p_2) |\nbs^6\Jg|^2 
\notag\\
&\qquad
= \tfrac{1+\alpha}{2}\tints \JJh \nbs^6(\Jg \Wbn)\nbs^6\Jg
+ \tfrac{1-\alpha}{2}\tints \JJh   \nbs^6(\Jg\Zbn)\nbs^6\Jg 
+ \tints \JJh   \Rj \ \nbs^6\Jg \,.
\end{align*} 
Using \eqref{eq:adjoints-H}, we have that
\begin{align*}
&
\tfrac{1}{2} \dint \Q  \JJh \sabs{\nbs^6 \Jg }^2\Big|^\s_\sin 
 -\tfrac{1}{4} \tints \JJhi  (\Q\p_\s+V\p_2)\JJ \  \sabs{\nbs^6 \Jg}^2
- \tfrac{1}{2}  \tints (\Qc +V,_2) \Jgh \sabs{\nbs^6 \Jg }^2
\notag \\
& \qquad
- \tfrac{1}{2} \int_\sin^\s\!\! \int_\mathbb{T}  \Big( \big( \Q\p_\s\theta  +V\p_2\theta   \big) \JJh \sabs{\nbs^6 \Jg}^2\Big)\Big|_{\theta(x_2,\s)}  {\rm d}x_2{\rm d}\s
\notag\\
&= \tfrac{1+\alpha}{2}\tints \JJh \nbs^6(\Jg \Wbn) \nbs^6\Jg
+ \tfrac{1-\alpha}{2}\tints \JJh  \nbs^6(\Jg\Zbn)\nbs^6\Jg 
+ \tints \JJh   \Rj \ \nbs^6\Jg \,.
\end{align*} 
Notice that the second line in this equality vanishes due to the fact that $\JJ=0$ along the surface $x_1=\thsd$  from \eqref{JJ-vanish}.
It follows that
\begin{align*}
&
\tfrac{1}{2} \dint \Q  \JJh \sabs{\nbs^6 \Jg }^2\Big|^\s_\sin 
 - \tfrac{1}{4} \tints  \JJhi  (\Q\p_\s+V\p_2)\JJ \  \sabs{\nbs^6 \Jg}^2
- \tfrac{1}{2}  \tints (\Qc +V,_2) \JJh \sabs{\nbs^6 \Jg }^2
\notag \\
&= \tfrac{1+\alpha}{2}\tints \JJh \nbs^6(\Jg \Wbn) \nbs^6\Jg
+ \tfrac{1-\alpha}{2}\tints \JJh  \nbs^6(\Jg\Zbn)\nbs^6\Jg 
+ \tints \JJh   \Rj \ \nbs^6\Jg \,.
\end{align*}

Then, using the bootstraps~\eqref{bootstraps}, and the bounds \eqref{eq:Q:all:bbq-H} and  \eqref{qps-JJ-bound}, we obtain that
for $\eps>0$ sufficiently small, 
\begin{align} 
&\tfrac{1-\eps^{\frac{7}{4}} }{2}  \dint \JJh  |\nbs^6\Jg|^2 \Big|^\s_\eps  
+ \tfrac{(1+\alpha)}{10\eps}  \| \JJmof  \nbs^6\Jg  \|_{L^2_{x,\s}}^2 
\notag\\
& \qquad
\leq
\| \JJmof \nbs^6\Jg  \|_{L^2_{x,\s}} \Bigl(
\tfrac{1+\alpha}{2} \|\JJtf  \nbs^6(\Jg \Wbn)\|_{L^2_{x,\s}}
+ \tfrac{|1-\alpha|}{2} \|\JJtf \nbs^6(\Jg\Zbn)\|_{L^2_{x,\s}} 
+\| \JJtf   \Rj \|_{L^2_{x,\s}} \Bigr)
\notag\\
&\qquad
\leq
\tfrac{1+\alpha}{20 \eps} \| \JJmof  \nbs^6\Jg  \|_{L^2_{x,\s}}^2
+4(1+\alpha) \eps \|\JJtf  \nbs^6(\Jg \Wbn,\Jg \Zbn)\|_{L^2_{x,\s}}^2 
+ \tfrac{16\eps}{1+\alpha} \| \JJtf \Rj \|_{L^2_{x,\s}}^2
\,.
\label{D6JgEnergy:temp1-H}
\end{align} 
The bound for $ \| \Jgtf \Rj \|_{L^2_{x,\s}}^2$ is obtained identically as in \eqref{D6JgEnergy:temp2}, 
and hence, it follows that
\begin{align} 
& \dint \JJh  |\nbs^6\Jg( \cdot , \s)|^2
+ \tfrac{1}{\eps}  \| \JJmof  \nbs^6\Jg  \|_{L^2_{x,\s}}^2 
\le
20 \dint \JJh  |\nbs^6\Jg( \cdot , \sin)|^2
+ \eps C \brak{\mathsf{B}_6}^2
 \,, \label{D6JgEnergy:temp3-H}
\end{align}
where $C$ depends only on $\alpha$.  
The bound for $\dint \JJh  |\nbs^6\Jg( \cdot , \sin)|^2$, given in 
\eqref{D6JgEnergy},  has an implicit constant that depends on $\alpha$ and $\Cdata$.  This concludes the 
proof of \eqref{D6JgEnergy:new-H}.

The upstream modifications required for proof of the inequalities \eqref{D6h2Energy:new-H}--\eqref{D6n-bound:b:new:bdd-H}  are identical  and these details will be omitted. To avoid redundancy, we also omit the proofs of the unweighted bounds for $(\Jg, \nbs_1 h, \nbs_2 h, \Sigma, V)$, as these are established exactly as in the proof of Proposition~\ref{prop:geometry}.
\end{proof}

\subsection{$L^2$-norms of functions evaluated along $\Thd$}
 
We first define the set of points at the intersection of each surface $\Thd(x,\s)$ and the initial time-slice $\{\s=\sin\}$. Recall the definitions~\eqref{eq:Thd:x:s:def} and~\eqref{eq:X:plus:stopping:time}, for each $\s\in [\sin,0]$ and $x_2\in\TT$, the intersection of the initial time-slice $\{\s=\sin\}$ with the slow characteristic surface $\Thd(x,\s)$ occurs at $x_1 = \mathfrak{X}_1^+(x_2,\mathfrak{q}^{-1}(x_2,\s))$, by the very definition of  the stopping time $\mathfrak{X}_1^+(x_2,t)$ in~\eqref{eq:X:plus:stopping:time}. Similarly, the intersection of $\Thd(x,\s)$ with the future temporal boundary of $\Hdm$ occurs at $x_1 = \mathfrak{X}_1^-(x_2,\mathfrak{q}^{-1}(x_2,\s))$, which may be seen by tracing back the definitions~\eqref{eq:Thd:x:s:def}, \eqref{t-to-s-transform-H}, and  \eqref{eq:X:minus:stopping:time}.  That is, upon defining
\begin{equation}
\tilde{\mathfrak{X}}_1^\pm(x_2,\s) :=  \mathfrak{X}_1^\pm(x_2,\mathfrak{q}^{-1}(x_2,\s)) \,,
\label{eq:X:pm:stopping:time:s}
\end{equation}
we have the identities
\begin{equation}
\Thd( \tilde{\mathfrak{X}}_1^-(x_2,\s),x_2,\s) = \mathfrak{q}(x_2,\final) \leq \sfin\,,
\qquad
\Thd( \tilde{\mathfrak{X}}_1^+(x_2,\s),x_2,\s) = \mathfrak{q}(x_2,\initial) = \sin\,.
\label{eq:never:thought:we:d:need:this:crap}
\end{equation}

Next, we recall that the domain of the function $\Thd(x,\s)$ is the spacetime set $\Omega_{\mathsf{US},+}$, defined in \eqref{eq:tilde:Omega:DS:+}. Keeping this in mind, we define:

\begin{definition} [The $L^2$ norms of a composite function $F(x,\Thd(x,\s))$]\label{def::L2-norm}
For any fixed time $\s\in[\sin,0)$, we define the spatial $L^2_x$ norm by  
\begin{align} 
\snorm{F(\cdot,\Thd(\cdot,\s))}^2_{L^2_x} = \int_{x_2=-\pi}^\pi \int_{x_1=\tilde{\mathfrak{X}}_1^-(x_2,\s)}^{\tilde{\mathfrak{X}}_1^+(x_2,\s)}  \sabs{F(x_1,x_2,\Thd(x_1,x_2,\s))}^2 {\rm d}x_1 {\rm d}x_2 \,,
\label{L2x-composite}
\end{align} 
where we recall the definitions in~\eqref{eq:X:pm:stopping:time:s}.
The spacetime  $L^2$ norm of  $F(x,\Thd(x,\s))$ is written as
\begin{subequations} 
\label{L2sx-composite}
\begin{align} 
\snorm{F(\cdot ,\Thd(\cdot,\cdot))}^2_{L^2_{x,\s}(\Omega_{\mathsf{US},+})} 
&= \int_{\s=\sin}^0 \snorm{F(\cdot,\Thd(\cdot,\s))}^2_{L^2_x} {\rm d} \s
\notag\\
&= \int_{\s=\sin}^0 \int_{x_2=-\pi}^\pi \int_{x_1=\tilde{\mathfrak{X}}_1^-(x_2,\s)}^{\tilde{\mathfrak{X}}_1^+(x_2,\s)}  \sabs{F(x,\Thd(x,\s))}^2 {\rm d}x_1 {\rm d}x_2 {\rm d}\s \notag\\
&= \int_{\Omega_{\mathsf{US},+}} \sabs{F(x,\Thd(x,\s))}^2 {\rm d}x_1 {\rm d}x_2 {\rm d}\s \,,
\end{align} 
upon recalling the definition~\eqref{eq:tilde:Omega:DS:+}.
By a slight abuse of notation, we shall write  
\begin{equation} 
\snorm{F(\cdot ,\Thd(\cdot,\cdot))}^2_{L^2_{x,\s}} := \snorm{F(\cdot ,\Thd(\cdot,\cdot))}^2_{L^2_{x,\s}(\Omega_{\mathsf{US},+})} 
\end{equation} 
\end{subequations} 
as $F(\cdot ,\Thd(\cdot,\cdot))$ is only defined on $\Omega_{\mathsf{US},+}$ and hence there is no other possible domain for this spacetime norm.
\end{definition} 

\subsubsection{A spacetime change of variables formula}
It will be necessary for us to make the spacetime change of variables $\s \mapsto \Thd(x,\s)$ in order to relate the norm
$\|F(\cdot ,\Thd(\cdot,\cdot))\|_{L^2_{x,\s}}$ (as defined by~\eqref{L2sx-composite}), to the norm $\|F\|_{L^2_{x,\s,+}}$ (as defined by~\eqref{int-Hdmp} and~\eqref{norms-H+-}). This is possible because according to~\eqref{eq:spacetime-Theta:aa} the spacetime $\Hdmp$ is indeed foliated by the surfaces $(x,\Thd(x,\s^\prime))$, with $(x,\s^\prime) \in \Omega_{\mathsf{US},+}$. That is, for every $(x,\s) \in \Hdmp$, there exists a unique $\s^\prime \in [\sin,0)$ such that $\s = \Thd(x,\s^\prime)$, and the map $(x,\s^\prime) \mapsto (x,\Thd(x,\s^\prime)) = (x,\s)$ is a bijection $\Omega_{\mathsf{US},+} \to \Hdmp$. The Jacobian determinant of this map is $\p_{\s^\prime} \Thd(x,\s^\prime)$, which according to~\eqref{ps-Theta} is uniformly close to $1$. Applying the change-of-variables theorem to this transformation, we find that 
\begin{align} 
\snorm{F(\cdot ,\Thd(\cdot,\cdot))}^2_{L^2_{x,\s}} 
&=
 \int_{\Omega_{\mathsf{US},+}} \sabs{F(x,\Thd(x,\s^\prime))}^2 {\rm d}x_1 {\rm d}x_2 {\rm d}\s^\prime
\notag\\
&=
\int_{\Hdmp} \sabs{F(x,\s)}^2 \tfrac{1}{\p_{\s^\prime} \Thd(x,\s^\prime)}\Bigr|_{\s = \Thd(x,\s^\prime)} {\rm d}x_1 {\rm d}x_2 {\rm d}\s 
=  \snorm{(\p_\s\Thd)^{-\frac{1}{2}} \; F }^2_{L^2_{x,\s}(\Hdmp)} \,,
\label{L2sx-fubini}
\end{align} 
the last equality coming from the definition of the $L^2_{{x,\s}}(\Hdmp)$-norm given by~\eqref{int-Hdmp} and~\eqref{norms-H+-}.

\subsubsection{Useful inequalities for fifth-derivative bounds}
We next develop some inequalities that will be used for fifth-derivative bounds.  We follow the methodology of Section~\ref{sec:D5:bootstrap}. The main estimates established below are \eqref{lazy-ass3}, \eqref{lazy-ass5}, and~\eqref{lazy-ass8}.
 
We first seek to obtain a bound for the  $L^2_x$-norm of the composite function $F(\cdot,\Thd(\cdot,s))$, as defined in~\eqref{L2x-composite}. In order to bound from above $\tfrac{d}{d\s} \| F(\cdot,\Thd(\cdot,\s))\|_{L^2_x}^2$, it is clear from~\eqref{L2x-composite} that we need to be able to differentiate the stopping times $\tilde{\mathfrak{X}}_1^{\pm}(x_2,\s)$ with respect to $\s$. This is achieved by differentiating the two identities in~\eqref{eq:never:thought:we:d:need:this:crap} with respect to $\s$, and appealing to \eqref{p1-Theta-sign-t}, \eqref{boot-psTheta-t}, \eqref{eq:xt:xs:chain:rule-Thd:a}, \eqref{eq:xt:xs:chain:rule-Thd:c}, \eqref{Qd-lower-upper-H}, and to the fact that $\Jg(\cdot,\sin) =1$; this leads to
\begin{subequations}
\label{eq:chris:rea}
\begin{align}
\tfrac{d}{d\s} \bigl( \tilde{\mathfrak{X}}_1^{+}(x_2,\s) \bigr)
&= - \tfrac{(\p_\s \Thd)(\tilde{\mathfrak{X}}_1^{+}(x_2,\s),x_2,\s)}{(\p_1 \Thd)(\tilde{\mathfrak{X}}_1^{+}(x_2,\s),x_2,\s)} \in [\tfrac{\alpha \kappa_0}{5}, 5\alpha \kappa_0\bigr]
\,,
\label{eq:chris:rea:a}
\\
\tfrac{d}{d\s} \bigl( \tilde{\mathfrak{X}}_1^{-}(x_2,\s) \bigr)
&= - \tfrac{(\p_\s \Thd)(\tilde{\mathfrak{X}}_1^{-}(x_2,\s),x_2,\s)}{(\p_1 \Thd)(\tilde{\mathfrak{X}}_1^{-}(x_2,\s),x_2,\s)} \in [\tfrac{\alpha \kappa_0}{5}, \infty)
\,.
\label{eq:chris:rea:b}
\end{align} 
\end{subequations}
The rough upper bound in~\eqref{eq:chris:rea:b} is due to the fact that the $\Jg \circ \Thd$ factor implicitly present in $\p_1 \Thd$ is evaluated at time $\mathfrak{q}(x_2,\final)$, and so $\Jg$ here could in principle be as close as possible to $0$, though still strictly positive. With \eqref{eq:chris:rea}, we simply differentiate \eqref{L2x-composite} and deduce that 
\begin{align}
&\tfrac{d}{d\s}  \snorm{F(\cdot,\Thd(\cdot,\s))}^2_{L^2_x} 
\notag\\
&\leq   5\alpha\kappa_0 \int_{\TT} \sabs{F(\tilde{\mathfrak{X}}_1^{+}(x_2,\s),x_2,\sin)}^2 {\rm d}x_2
\notag\\
&\quad + \int_{x_2=-\pi}^\pi \int_{x_1=\tilde{\mathfrak{X}}_1^{-}(x_2,\s)}^{\tilde{\mathfrak{X}}_1^{+}(x_2,\s)}  2  F(x_1,x_2,\Thd(x_1,x_2,\s)) \p_\s F(x_1,x_2,\Thd(x_1,x_2,\s)) \p_\s \Thd(x_1,x_2,\s)  {\rm d}x_1 {\rm d}x_2
\notag\\
&\leq 
 5\alpha\kappa_0 \|F(\tilde{\mathfrak{X}}_1^+(\cdot,\s),\cdot,\sin)\|_{L^2_{x_2}}^2
+2 \cdot \tfrac{1001}{999}\snorm{F(\cdot,\Thd(\cdot,\s))}_{L^2_x} \snorm{\p_\s F(\cdot,\Thd(\cdot,\s))}_{L^2_x}  
\label{eq:chris:rea:2}
\,.
\end{align}
In the last inequality above we have appealed to~\eqref{ps-Theta} to bound $\p_\s\Thd$. The estimate \eqref{eq:chris:rea:2} then implies upon integration in time, and by appealing to~\eqref{eq:chris:rea:a} and the change of variable formula for the map $(x_2,\s^\prime) \mapsto(x_2,\tilde{\mathfrak{X}}_1^+(x_2,\s^\prime))$, whose determinant Jacobian is $\frac{d}{d\s} \tilde{\mathfrak{X}}_1^+$, that 
\begin{align}
 \snorm{F(\cdot,\Thd(\cdot,\s))}_{L^2_x}^2
 &\leq 
  \snorm{F(\cdot,\Thd(\cdot,\sin))}_{L^2_x}^2
  + 5\alpha\kappa_0 \int_{\sin}^{\s}\int_{\TT} \sabs{F(\tilde{\mathfrak{X}}_1^+(x_2,\s^\prime),x_2,\sin)}^2 {\rm d}x_2 {\rm d}\s^\prime
  \notag\\
 &\qquad 
  + 2 \cdot \tfrac{1001}{999} \int_{\sin}^{\s}
  \snorm{F(\cdot,\Thd(\cdot,\s^\prime))}_{L^2_x} \snorm{\p_\s F(\cdot,\Thd(\cdot,\s^\prime))}_{L^2_x}  {\rm d} \s^\prime
  \notag\\
 &\leq 
  \snorm{F(\cdot,\Thd(\cdot,\sin))}_{L^2_x}^2
  + 25  \snorm{F(\cdot,\sin)}_{L^2_x}^2
  \notag\\
 &\qquad 
  + 2 \cdot \tfrac{1001}{999} \int_{\sin}^{\s}
  \snorm{F(\cdot,\Thd(\cdot,\s^\prime))}_{L^2_x} \snorm{\p_\s F(\cdot,\Thd(\cdot,\s^\prime))}_{L^2_x}  {\rm d} \s^\prime  
  \label{eq:chris:rea:3}
\end{align}

We now consider a  function $F(x,\s)$ such that $\JJ^{r} \p_\s F \in L^2_{x,\s}(\Hdm)$, for some $0\leq r \leq \frac 34$.
Our goal is to establish the bound
\begin{align} 
&  \snorm{F}_{L^2_{x,\s}(\Hdm)}^2
 \le  224 \eps \snorm{ F(\cdot,\sin)}_{L^2_x}^2
+  664 \eps^2  \snorm{\JJ^r\, \p_\s F}^2_{L^2_{x,\s}(\Hdm)}  
 \,, \qquad 0 \leq r \le \tfrac{3}{4} 
\,.
\label{lazy-ass3}
\end{align} 
In order to prove \eqref{lazy-ass3}, we first decompose the left side of this inequality as integrals over $\Hdmp$ and $\Hdmm$. From~\eqref{norms-H:b}, \eqref{norms-H+-}, \eqref{ps-Theta}, \eqref{L2sx-composite},  and~\eqref{L2sx-fubini}, we deduce that 
\begin{equation}
\snorm{F}_{L^2_{x,\s}(\Hdm)}^2
\leq \snorm{F}_{L^2_{x,\s}(\Hdmp)}^2 + \snorm{F}_{L^2_{x,\s}(\Hdmm)}^2
\leq \tfrac{1001}{999} \snorm{ F (\cdot,\Thd(\cdot,\cdot))}_{L^2_{x,\s}}^2 + \snorm{F}_{L^2_{x,\s}(\Hdmm)}^2
\,.
\label{eq:lazy-ass:0}
\end{equation}

For the first term on the right side of~\eqref{eq:lazy-ass:0}, we use definition~\eqref{L2sx-composite} the lower bound~\eqref{JJ-and-s} and the previously obtained upper bound \eqref{eq:chris:rea:3}, to deduce
\begin{align*}
 \snorm{ F (\cdot,\Thd(\cdot,\cdot))}_{L^2_{x,\s}}^2
 &= \int_{\sin}^0  \snorm{ F (\cdot,\Thd(\cdot,\s))}_{L^2_{x}}^2 {\rm d}\s
\notag\\
&\leq |\sin|  \snorm{F(\cdot,\Thd(\cdot,\sin))}_{L^2_x}^2
  + 25 |\sin|  \snorm{F(\cdot,\sin)}_{L^2_x}^2
  \notag\\
 &\qquad 
  + 2 \cdot \tfrac{1001}{999} \int_{\sin}^0 \int_{\sin}^{\s}
  \snorm{F(\cdot,\Thd(\cdot,\s^\prime))}_{L^2_x} \snorm{\p_\s F(\cdot,\Thd(\cdot,\s^\prime))}_{L^2_x}  {\rm d}\s^\prime  {\rm d}\s
 \notag\\
 &\leq |\sin|  \snorm{F(\cdot,\Thd(\cdot,\sin))}_{L^2_x}^2
  + 25 |\sin|  \snorm{F(\cdot,\sin)}_{L^2_x}^2
  \notag\\
 &\qquad 
  + 2 \cdot \tfrac{1001}{999} \int_{\sin}^0 \int_{\sin}^{\s}
  \snorm{F(\cdot,\Thd(\cdot,\s^\prime))}_{L^2_x} \snorm{(\JJ^r \p_\s F)(\cdot,\Thd(\cdot,\s^\prime))}_{L^2_x} \bigl( \tfrac{89}{100} \cdot \tfrac{\s^\prime}{\sin}\bigr)^{-r} {\rm d}\s^\prime  {\rm d}\s
  \notag\\
   &\leq |\sin|  \snorm{F(\cdot,\Thd(\cdot,\sin))}_{L^2_x}^2
  + 25 |\sin|  \snorm{F(\cdot,\sin)}_{L^2_x}^2
  \notag\\
 &\qquad 
  + 9 |\sin| \snorm{F(\cdot,\Thd(\cdot,\cdot))}_{L^2_{x,\s}} 
  \snorm{(\JJ^r \p_\s F)(\cdot,\Thd(\cdot,\cdot))}_{L^2_{x,\s}}
  \,.
\end{align*}
In the last inequality we have used the fact that for $r\in[0,3/4]$ we have $\int_{\sin}^0    (  \tfrac{\s }{\sin} )^{-r}   {\rm d}\s \leq4 |\sin|$. Using an Cauchy-Young argument, and then by appealing to~\eqref{L2sx-fubini} and \eqref{ps-Theta}, we deduce from the above estimate that 
\begin{align}
\snorm{ F (\cdot,\Thd(\cdot,\cdot))}_{L^2_{x,\s}}^2
 &\leq 2 |\sin|  \snorm{F(\cdot,\Thd(\cdot,\sin))}_{L^2_x}^2
  + 50 |\sin|  \snorm{F(\cdot,\sin)}_{L^2_x}^2
  + 81  |\sin|^2 \snorm{(\JJ^r \p_\s F)(\cdot,\Thd(\cdot,\cdot))}_{L^2_{x,\s}}^2
  \notag\\
   &\leq 2 |\sin|  \snorm{F(\cdot,\Thd(\cdot,\sin))}_{L^2_x}^2
  + 50 |\sin|  \snorm{F(\cdot,\sin)}_{L^2_x}^2
  + 82  |\sin|^2 \snorm{ \JJ^r \p_\s F}_{L^2_{x,\s}(\Hdmp)}^2
  \,.
\label{eq:lazy-ass:00}
\end{align}
We notice that the first term on the right side of \eqref{eq:lazy-ass:00} needs a further upper bound, in order to prove \eqref{lazy-ass3}. For this purpose, pointwise in $x$ we apply the fundamental theorem of calculus in time to write
\begin{equation*} 
F^2(x,\Thd(x,\sin))  = F^2(x,\sin) + 2 \int_\sin^{ \Thd(x,\sin)} F(x,\s) \p_\s F(x,\s) {\rm d}\s \,.
\end{equation*} 
Hence, recalling the definition~\eqref{L2x-composite}, appealing to \eqref{fubini-thsin}, and using the lower bound $\JJ \geq 1$ on the curve $(x,\Thd(x,\sin))$ (cf.~\eqref{JJ-and-s-new}) we deduce
\begin{align} 
\snorm{ F(\cdot,\Thd(\cdot ,\sin)) }^2_{L^2_x} 
&\leq \snorm{ F(\cdot,\sin)}_{L^2_x}^2
+ 2 \int_{x_2=-\pi}^\pi \int_{x_1=\tilde{\mathfrak{X}}_1^-(x_2,\sin)}^{\tilde{\mathfrak{X}}_1^+(x_2,\sin)} \int_{\s =\sin}^{ \Thd(x_1,x_2,\sin)} F(x,\s) \p_\s F(x,\s) {\rm d}\s {\rm d}x_1 {\rm d}x_2
\notag\\
&\leq \snorm{ F(\cdot,\sin)}_{L^2_x}^2
+ 2 \int_{x_1=-\pi}^\pi \int_{x_2=-\pi}^\pi  \int_{\s =\sin}^{ \Thd(x_1,x_2,\sin)} \sabs{F(x,\s)} \sabs{ (\JJ^r \p_\s F)(x,\s)} {\rm d}\s {\rm d}x_2{\rm d}x_1 
\notag\\
&= \snorm{ F(\cdot,\sin)}_{L^2_x}^2
+ 2 \int_{\s=\sin}^{\sfin} \int_{x_2=-\pi}^\pi   \int_{x_1=-\pi}^{\thsin(x_2,\s)} \sabs{F(x,\s)} \sabs{ (\JJ^r \p_\s F)(x,\s)} {\rm d}x_1 {\rm d}x_2 {\rm d}\s 
\notag\\ 
&\leq \snorm{ F(\cdot,\sin)}_{L^2_x}^2 + 2 \snorm{F}_{L^2_{x,\s}(\Hdmm)} \snorm{\JJ^r \p_\s F}_{L^2_{x,\s}(\Hdmm)}
\label{lazy-ass1}
\,.
\end{align} 
 
Let us now  estimate the second term on the right side of~\eqref{eq:lazy-ass:0}, namely $\|F\|_{L^2_{x,\s}(\Hdmm)}$. For the $L^2$ norm on $\Hdmm$, which was earlier defined in~\eqref{int-Hdmm}, we have that
\begin{align*} 
\int_{\TT} \int_{-\pi}^{\thsin(x_2,\s)} F^2(x,\s){\rm d}x 
& 
\leq \snorm{F( \cdot , \sin)}^2_{L^2_x} 
+ \int_\sin^\s \tfrac{d}{d\s'} \int_{\TT} \int_{-\pi}^{\thsin(x_2,\s')} F^2(x,\s'){\rm d}x {\rm d}\s'
\notag \\
&= \snorm{F( \cdot , \sin)}^2_{L^2_x} 
+ 2\int_\sin^\s  \int_{\TT} \int_{-\pi}^{\thsin(x_2,\s')} F(x,\s')\p_{\s'}F(x,\s') {\rm d}x{\rm d}\s' \notag \\
& \qquad
+ \int_\sin^\s  \int_{\TT} \p_{\s'}\thsin(x_2,\s') F^2(\thsin(x_2,\s'),x_2 ,\s'){\rm d}x_2 {\rm d}\s' \notag \\
&\le  \snorm{F( \cdot , \sin)}^2_{L^2_x} 
+ 2 \snorm{F}_{L^2_{x,\s}(\Hdmm)} \snorm{\JJ^r \p_\s F}_{L^2_{x,\s}(\Hdmm)}\,,
\end{align*} 
the last inequality following from \eqref{ps-thsin-sign} and the fact that $\JJ \geq 1$ in $\Hdmm$ due to~\eqref{JJ-and-s-new}. Integrating the above inequality for $\s\in [\sin,\sfin]$ and then using Cauchy-Young, we obtain
\begin{equation} 
\snorm{F}_{L^2_{x,\s}(\Hdmm)}^2  
\le 2 (\sfin-\sin)  \snorm{F( \cdot , \sin)}^2_{L^2_x} 
 +  4 (\sfin-\sin)^2  \snorm{\JJ^r \p_\s F}_{L^2_{x,\s}(\Hdmm)}^2
 \,.
\label{eq:lazy-ass:000}
 \end{equation} 
 
Combining the bounds~\eqref{eq:lazy-ass:0}, \eqref{eq:lazy-ass:00}, \eqref{lazy-ass1}, \eqref{eq:lazy-ass:000}, using that $|\sfin-\sin|\leq \tfrac{26}{25} |\sin|$, and applying the Cauchy-Young inequality, we deduce
\begin{align*}
\snorm{F}_{L^2_{x,\s}(\Hdm)}^2 
&\leq    
112 |\sin|  \snorm{F(\cdot,\sin)}_{L^2_x}^2
+ 166  |\sin|^2 \snorm{ \JJ^r \p_\s F}_{L^2_{x,\s}(\Hdmp)}^2  
+ 35 |\sin|^2 \snorm{ \JJ^r \p_\s F}_{L^2_{x,\s}(\Hdmm)}^2
\,.
\end{align*}
This bound clearly implies~\eqref{lazy-ass3}, as claimed.

Our next goal is to establish the bound
\begin{equation}
\sup_{\s \in [\sin,0]} 
\snorm{(\Jgh F)(\cdot,\Thd(\cdot,\s))}_{L^2_x}^2
\leq  40 e^{15}  \snorm{F(\cdot,\sin)}_{L^2_x}^2 
+ 20 e^{30} \eps  \snorm{\JJof \Jgh \p_\s  F}_{L^2_{x,\s}(\Hdm)}^2
\label{lazy-ass5}
\,,
\end{equation}
in direct analogy to~\eqref{eq:r:Linfty:time:2}. For this purpose, we revisit the computation in~\eqref{eq:chris:rea:2}--\eqref{eq:chris:rea:3}, appeal to the fact that via~\eqref{geom-H-original} (more precisely, \eqref{eq:signed:Jg}) and the bootstrap assumptions in~\eqref{boots-HH} we have that $\p_\s \Jg \leq - \frac 45 \cdot \frac{1+\alpha}{2\eps} + 14 \cdot \frac{1+\alpha}{2\eps} \Jg \leq 14|\sin|^{-1} \Jg$, and thus
\begin{align}
\snorm{(\Jgh F)(\cdot,\Thd(\cdot,\s))}_{L^2_x}^2
&\leq  
\snorm{(\Jgh F)(\cdot,\Thd(\cdot,\sin))}_{L^2_x}^2
+ 25 \snorm{(\Jgh F)(\cdot,\sin)}_{L^2_x}^2 
\notag\\
&\quad
+ \tfrac{2002}{999} \int_{\sin}^{\s} \snorm{(\Jgh F)(\cdot,\Thd(\cdot,\s^\prime))}_{L^2_x} \snorm{(\Jgh \p_\s F)(\cdot,\Thd(\cdot,\s^\prime))}_{L^2_x} {\rm d}\s^\prime
\notag\\
&\quad
+ \tfrac{1001}{999} \cdot \tfrac{14}{|\sin|} \int_{\sin}^{\s} \snorm{(\Jgh F)(\cdot,\Thd(\cdot,\s^\prime))}_{L^2_x}^2 {\rm d}\s^\prime
\,.
\label{eq:ramblin:on:my:mind:1}
\end{align}
Using that $\Jg(\cdot,\sin) \equiv 1$, that $|\Jg -1|\leq 5 \cdot 10^{-4}$ and $1\leq \JJ \leq 3$ in the closure of $\Hdmm$ via~\eqref{JJ-le-Jg:c}, appealing to the upper bounds~\eqref{lazy-ass1},~\eqref{eq:lazy-ass:000}, to the lower bound~\eqref{JJ-and-s}, and using the norm equivalence stated in~\eqref{L2sx-fubini}, we deduce from~\eqref{eq:ramblin:on:my:mind:1} that
\begin{align}
\snorm{(\Jgh F)(\cdot,\Thd(\cdot,\s))}_{L^2_x}^2
&\leq   26 \snorm{F(\cdot,\sin)}_{L^2_x}^2 
+ 3 \snorm{F}_{L^2_{x,\s}(\Hdmm)} \snorm{\JJof \Jgh \p_\s F}_{L^2_{x,\s}(\Hdmm)}
\notag\\
&\quad
+ \tfrac{21}{10} \int_{\sin}^{\s} \snorm{(\Jgh F)(\cdot,\Thd(\cdot,\s^\prime))}_{L^2_x} \snorm{(\JJof \Jgh \p_\s F)(\cdot,\Thd(\cdot,\s^\prime))}_{L^2_x}  \bigl( \tfrac{\s^\prime}{\sin} \bigr)^{-\frac 14} {\rm d}\s^\prime
\notag\\
&\quad
+  \tfrac{15}{|\sin|} \int_{\sin}^{\s} \snorm{(\Jgh F)(\cdot,\Thd(\cdot,\s^\prime))}_{L^2_x}^2 {\rm d}\s^\prime
\notag\\
&\leq 29 \snorm{F(\cdot,\sin)}_{L^2_x}^2 
+ 10 |\sin| \snorm{\JJof \Jgh \p_\s  F}_{L^2_{x,\s}(\Hdmm)}^2
\notag\\
&\quad
+ 3 |\sin|^{\frac 12} \Bigl( \sup_{\s \in [\sin,0]} \snorm{(\Jgh F)(\cdot,\Thd(\cdot,\s))}_{L^2_x} \Bigr)
\snorm{\JJof \Jgh \p_\s F}_{L^2_{x,\s}(\Hdmp)}
\notag\\
&\quad
+  \tfrac{15}{|\sin|} \int_{\sin}^{\s} \snorm{(\Jgh F)(\cdot,\Thd(\cdot,\s^\prime))}_{L^2_x}^2 {\rm d}\s^\prime
\,.
\label{eq:ramblin:on:my:mind:2}
\end{align}
Using Gr\"onwall's inequality for $\s \in [\sin,0]$ and the Cauchy-Schwartz inequality, we deduce from \eqref{eq:ramblin:on:my:mind:2} that
\begin{align*}
\sup_{\s \in [\sin,0]} 
\snorm{(\Jgh F)(\cdot,\Thd(\cdot,\s))}_{L^2_x}^2
&\leq  
29 e^{15} \snorm{F(\cdot,\sin)}_{L^2_x}^2 
+ 14 e^{15} |\sin| \snorm{\JJof \Jgh \p_\s  F}_{L^2_{x,\s}(\Hdmm)}^2
\notag\\
&\quad
+ 9 e^{30} |\sin| \snorm{\JJof \Jgh \p_\s  F}_{L^2_{x,\s}(\Hdmp)}^2
+ \tfrac 14 \sup_{\s \in [\sin,0]} \snorm{(\Jgh F)(\cdot,\Thd(\cdot,\s))}_{L^2_x}^2 
\,.
\end{align*}
Absorbing the last term on the right side into the left side of the above bound, then establishes~\eqref{lazy-ass5}.

The remaining bound that is established in this section is that for the norm defined in~\eqref{norms-H:a}, we have
\begin{align} 
\sup_{\s \in [\sin,\sfin]} 
\|( \Jgh F)(\cdot,\s)\|_{L^2_x}^2
\le 14 e^{28}  
\snorm{F( \cdot , \sin)}_{L^2_x}^2 
+ 92 \eps e^{56} 
\snorm{\JJof \Jgh\p_\s F}^2_{L^2_{x,\s}(\Hdm)} \,.
\label{lazy-ass8}
\end{align} 
In order to prove~\eqref{lazy-ass8}, we split the  integral defined in~\eqref{norms-H:a} similarly to~\eqref{int-Hdmp-Hdmm} as
\begin{align}
\|( \Jgh F)(\cdot,\s)\|_{L^2_x}^2
&=  \int_{\TT} \! \int_{-\pi}^{\thsin(x_2,\s)}
(\Jg F^2)(x_1,x_2,\s) {\rm d}x_1{\rm d}x_2
+ \int_{\TT} \! \int_{\thsin(x_2,\s)}^{\thd(x_2,\s)}
(\Jg F^2)(x_1,x_2,\s) {\rm d}x_1{\rm d}x_2
\notag\\
&=: \mathcal{I}(\s,-) + \mathcal{I}(\s,+)
\,.
\label{eq:pompeii:1}
\end{align}
For the first integral, a fundamental theorem of calculus in time is applied between time $\s$ and time $\sin$. The monotonicity with respect to $x_1$ of $\Thd(x_1,x_2,\sin)$ ensured by \eqref{p1-Theta-sign} implies that for every $x_1 \leq \thsin(x_2,\s)$, and for all $\s^\prime \in [\sin,\s)$, we have $(x_1,x_2,\s^\prime) \in \Hdmm$. Then, with the bounds for $\JJ$ and $\Jg$ provided by \eqref{JJ-le-Jg:c} in $\Hdmm$, and with the estimate $\p_\s \Jg \leq 14 |\sin|^{-1} \Jg$, we deduce
\begin{align}
 &\mathcal{I}(\s,-)
 = \int_{\TT} \! \int_{-\pi}^{\thsin(x_2,\s)}
(\Jg F^2)(x_1,x_2,\sin) {\rm d}x_1{\rm d}x_2 
+ 
\int_{\sin}^{\s} \int_{\TT} \! \int_{-\pi}^{\thsin(x_2,\s)}
\p_\s (\Jg F^2)(x_1,x_2,\s^\prime) {\rm d}x_1{\rm d}x_2 {\rm d}\s^\prime
\notag\\
 &\leq \snorm{F(\cdot,\sin)}_{L^2_x}^2
+ 2
\int_{\sin}^{\s} \|( \Jgh F)(\cdot,\s^\prime)\|_{L^2_x} \|(\JJof \Jgh \p_\s F)(\cdot,\s^\prime)\|_{L^2_x} {\rm d}\s^\prime
+ \tfrac{14}{|\sin|} 
\int_{\sin}^{\s} \|( \Jgh F)(\cdot,\s^\prime)\|_{L^2_x}^2 {\rm d}\s^\prime
\,.
\label{eq:pompeii:2}
\end{align}

In order to bound the $\mathcal{I}(\s,+)$ integral appearing in~\eqref{eq:pompeii:1}, we let $\s^\prime \in (\sin,0)$ be arbitrary (we will eventually pass $\s^\prime \to 0$). Then, for each $x_2 \in \TT$ fixed, there exists a unique  
\begin{equation*}
\nu(x_2,\s,\s^\prime):= x_1 \in [ \tilde{\mathfrak{X}}_1^-(x_2,\s^\prime),\tilde{\mathfrak{X}}_1^+(x_2,\s^\prime)] \cap [\thsin(x_2,\s),\thd(x_2,\s))
 \quad \mbox{such that} \quad 
 \s = \Thd(x_1,x_2,\s^\prime)
 \,.
\end{equation*} 
Note that $\lim_{\s^\prime \to \sin} \nu(x_2,\s,\s^\prime) = \thsin(x_2,\s)$, and more importantly, $\lim_{\s^\prime \to 0} \nu(x_2,\s,\s^\prime) = \thd(x_2,\s)$. Because of this, by the monotone convergence theorem, we may write 
\begin{equation}
\mathcal{I}(\s,+) = \lim_{\s^\prime \to 0} \int_{\TT} \! \int_{\thsin(x_2,\s)}^{\nu(x_2,\s,\s^\prime)} (\Jg F^2)(x_1,x_2,\s) {\rm d}x_1 {\rm d}x_2\,.
\label{eq:pompeii:3} 
\end{equation}
Our goal is to obtain uniform in $\s^\prime$ bounds for the integral expression in \eqref{eq:pompeii:3}. For this purpose, we use the fundamental theorem of calculus in time between the time $\s$, going down to time $\max\{ \Thd(x,\sin) , \sin \}$. We deduce
\begin{align}
\int_{\TT} \! \int_{\thsin(x_2,\s)}^{\nu(x_2,\s,\s^\prime)} (\Jg F^2)(x_1,x_2,\s) {\rm d}x_1 {\rm d}x_2
&= \int_{\TT} \! \int_{\thsin(x_2,\s)}^{\nu(x_2,\s,\s^\prime)} (\Jg F^2)(x_1,x_2,\max\{ \Thd(x,\sin) , \sin \}) {\rm d}x_1 {\rm d}x_2
\notag\\
&  + 
\int_{\TT} \! \int_{\thsin(x_2,\s)}^{\nu(x_2,\s,\s^\prime)} 
\int_{\max\{ \Thd(x,\sin) , \sin \}}^{\s} \p_\s(\Jg F^2)(x_1,x_2,\s^{\prime\prime}){\rm d}\s^{\prime\prime} {\rm d}x_1 {\rm d}x_2 
\,.
\label{eq:pompeii:4} 
\end{align}
Appealing to the triangle inequality, to the bound~\eqref{lazy-ass1} and the information on $\Jg$ in $\Hdmm$ available from~\eqref{JJ-le-Jg:c}, and to the bound~\eqref{eq:lazy-ass:000}, we deduce that the first term on the right side of 
\eqref{eq:pompeii:4} is bounded from above by
\begin{equation}
\snorm{F(\cdot,\sin)}_{L^2_x}^2 + \snorm{F(\cdot,\Thd(\cdot,\sin))}_{L^2_x}^2
\leq 
6 \snorm{F(\cdot,\sin)}_{L^2_x}^2 + 10 |\sin|  \snorm{\JJof \Jgh \p_\s F}_{L^2_{x,\s}(\Hdmm)}^2 
\,.
\label{eq:pompeii:5} 
\end{equation}
The above estimate is clearly independent of $\s^\prime$.
The second term on the right side of 
\eqref{eq:pompeii:4} is bounded similarly, using that $\p_\s \Jg \leq 14 |\sin|^{-1} \Jg$,  Fubini,
the change of variables formula, and the fact that  $\JJ$ satisfies the lower bound~\eqref{JJ-and-s}. An upper bound is provided by
\begin{align}
&\tfrac{14}{|\sin|} \int_{\TT} \! \int_{\thsin(x_2,\s)}^{\nu(x_2,\s,\s^\prime)} 
\int_{\max\{ \Thd(x,\sin) , \sin \}}^{\s} (\Jg F^2)(x_1,x_2,\s^{\prime\prime}) {\rm d}\s^{\prime\prime}{\rm d}x_1 {\rm d}x_2 
\notag\\
&\qquad
+ 
\int_{\TT} \! \int_{\thsin(x_2,\s)}^{\nu(x_2,\s,\s^\prime)} 
\int_{\max\{ \Thd(x,\sin) , \sin \}}^{\s} \bigl| (\Jgh F) \cdot (\JJof \Jgh \p_\s F)\cdot \JJ^{-\frac 14}\bigr|(x_1,x_2,\s^{\prime\prime}) {\rm d}\s^{\prime\prime}{\rm d}x_1 {\rm d}x_2 
\notag\\
&\leq \tfrac{14}{|\sin|} \int_{\sin}^{\s} \snorm{\Jgh F(\cdot,\s^{\prime\prime})}_{L^2_x}^2 {\rm d} \s^{\prime\prime}
+ \int_{\sin}^{\s} \int_{\TT}  \int_{\thsin(x_2,\s^{\prime\prime})}^{\nu(x_2,\s,\s^\prime)} 
\bigl| (\Jgh F) \cdot (\JJof \Jgh \p_\s  F)\cdot \JJ^{-\frac 14}\bigr|(x_1,x_2,\s^{\prime\prime}) {\rm d}x_1 {\rm d}x_2 {\rm d}\s^{\prime\prime} \,.
\label{eq:pompeii:6} 
\end{align}
The second term in~\eqref{eq:pompeii:6} requires a further careful analysis, in order to suitably bound $\JJ^{-\frac 14}(x,\s^{\prime\prime})$. First, by appealing to \eqref{JJ-formula}, \eqref{qps-JJ-bound}, and~\eqref{fat-marmot2b} we see that $2\alpha \Sigma \p_1 \JJ = (1-\dl) \Jg (\Q\p_\s + V \p_2) \JJ + 2\alpha \Sigma \Jg g^{-\frac 12} \nbs_2 h \nbs_2 \JJ < - \Jg  \tfrac{1+\alpha}{5\eps}  < 0$, and therefore  $\JJ$ is strictly decreasing with respect to $x_1$. As such, with~\eqref{mod-ps-JJ}, the fundamental theorem of calculus in time, the definition of $\nu(x_2,\s,\s^\prime)$ above, and the fact that $\JJ>0$ in $\Hdmp$, we deduce
\begin{align}
\min_{x_1 \leq \nu(x_2,\s,\s^\prime)} \JJ (x_1,x_2,\s^{\prime\prime}) 
&= \JJ(\nu(x_2,\s,\s^\prime),x_2,\s^{\prime\prime})
\notag\\
&= \JJ(\nu(x_2,\s,\s^\prime),x_2,\s) - (\s-\s^{\prime\prime}) \p_\s\JJ(\nu(x_2,\s,\s^\prime),x_2,\s^{\prime\prime\prime})
\geq  (\s-\s^{\prime\prime}) \tfrac{9(1+\alpha)}{25 \eps}
\label{eq:pompeii:7}
\,.
\end{align}
Using, \eqref{eq:pompeii:7} the right side of \eqref{eq:pompeii:6} may be further bounded from above by
\begin{align}
&\tfrac{14}{|\sin|} \int_{\sin}^{\s} \snorm{\Jgh F(\cdot,\s^{\prime\prime})}_{L^2_x}^2 {\rm d} \s^{\prime\prime}
+ \int_{\sin}^{\s} \snorm{(\Jgh F)(\cdot,\s^{\prime\prime})}_{L^2_x} 
\snorm{(\JJof \Jgh  \p_\s F)(\cdot,\s^{\prime\prime})}_{L^2_x} \bigl((\s-\s^{\prime\prime}) \tfrac{9(1+\alpha)}{25 \eps}\bigr)^{-\frac 14}  {\rm d}\s^{\prime\prime} 
\notag\\
&\leq \tfrac{14}{|\sin|} \int_{\sin}^{\s} \snorm{\Jgh F(\cdot,\s^{\prime\prime})}_{L^2_x}^2 {\rm d} \s^{\prime\prime}
+ 2 |\sin|^{\frac 14} (\s-\sin)^{\frac 14} \Bigl( \sup_{\s^{\prime\prime}\in[\sin,\s]}  \snorm{(\Jgh F)(\cdot,\s^{\prime\prime})}_{L^2_x}\Bigr)
\snorm{\JJof \Jgh  \p_\s F}_{L^2_{x,\s}(\Hdmp)}
\,.
\label{eq:pompeii:8} 
\end{align}
Importantly, this bound is also independent of $\s^\prime$. 
By combining \eqref{eq:pompeii:3}, \eqref{eq:pompeii:4}, \eqref{eq:pompeii:5}, and \eqref{eq:pompeii:8}, we thus obtain
\begin{align}
\mathcal{I}(\s,+)
&\leq 
6 \snorm{F(\cdot,\sin)}_{L^2_x}^2 + 10 |\sin|  \snorm{\JJof \Jgh \p_\s F}_{L^2_{x,\s}(\Hdmm)}^2
+ \tfrac{14}{|\sin|} \int_{\sin}^{\s} \snorm{\Jgh F(\cdot,\s^{\prime\prime})}_{L^2_x}^2 {\rm d} \s^{\prime\prime}
\notag\\
&\qquad + 3 |\sin|^{\frac 12}  \Bigl( \sup_{\s^{\prime\prime}\in[\sin,\s]}  \snorm{(\Jgh F)(\cdot,\s^{\prime\prime})}_{L^2_x}\Bigr)
\snorm{\JJof \Jgh F}_{L^2_{x,\s}(\Hdmp)}\,.
\label{eq:pompeii:9} 
\end{align}
From~\eqref{eq:pompeii:1}, \eqref{eq:pompeii:2}, \eqref{eq:pompeii:9} we obtain 
\begin{align*}
\|( \Jgh F)(\cdot,\s)\|_{L^2_x}^2
&\leq 
7 \snorm{F(\cdot,\sin)}_{L^2_x}^2 + 10 |\sin|  \snorm{\JJof \Jgh \p_\s F}_{L^2_{x,\s}(\Hdm)}^2
+ \tfrac{28}{|\sin|} \int_{\sin}^{\s} \snorm{\Jgh F(\cdot,\s^{\prime\prime})}_{L^2_x}^2 {\rm d} \s^{\prime\prime}
\notag\\
&\qquad + 5 |\sin|^{\frac 12}  \Bigl( \sup_{\s^{\prime\prime}\in[\sin,\s]}  \snorm{(\Jgh F)(\cdot,\s^{\prime\prime})}_{L^2_x}\Bigr)
\snorm{\JJof \Jgh F}_{L^2_{x,\s}(\Hdm)}\,.
\end{align*}
An application of the Gr\"onwall inequality in $\s$ then directly gives the proof of~\eqref{lazy-ass8}.  

\subsubsection{Bounds for $ \mathcal{E}_5$ and $ \mathcal{D}_5$}

We first show that \eqref{bootstraps-Dnorm:5} follows from \eqref{bootstraps-Dnorm:6}, assuming that  $\mathsf{B}_5$ is sufficiently large with 
respect to $\mathsf{B}_6$.  We will employ the inequality  \eqref{lazy-ass3} with $r=\frac 34$.   Using \eqref{table:derivatives}, 
recalling that $\p_\s = \frac{1}{\eps \Qd} \nbs_\s$, and that $\Qd^{-1}$ is bounded according to \eqref{Qd-lower-upper-H}, we deduce
from  \eqref{lazy-ass3} with $r=\frac 34$ that
\begin{align*}
\widetilde{\mathcal{D}}^2_{5,\nnn} (\sfin)
&\le  224 \Cdata^2+  680 \widetilde{\mathcal{D}}^2_{6,\nnn}(\sfin)\,,
\\
\widetilde{\mathcal{D}}^2_{5,\ttt} (\sfin)
&\le  224\eps^2 \Cdata^2 + 680  \widetilde{\mathcal{D}}^2_{6,\ttt}(\sfin)
\,.
\end{align*}
From the above bound and the definition \eqref{eq:tilde:D5D6-H}, it follows that 
\begin{align}
\widetilde{\mathcal{D}}^2_5(\sfin)  
&\le 224 \Cdata^2+  680 \widetilde{\mathcal{D}}^2_{6,\nnn}(\sfin)
+ (\mathsf{K} \eps)^{-2} \bigl( 224 \eps^2\Cdata^2+  680 \widetilde{\mathcal{D}}^2_{6,\ttt}(\sfin)\bigr)
\notag\\
&\le 448 \Cdatatwo  + 680 \widetilde{\mathcal{D}}^2_6(\sfin) 
\,,
\label{eq:D5:D6:dominate-H}
\end{align}
since $\mathsf{K}\geq 1$.
Similarly, by appealing to \eqref{lazy-ass8}, using \eqref{table:derivatives} and the fact that $\Jg(\cdot,\sin) = 1$, we deduce that 
\begin{align*}
\sup_{\s\in [\sin,\sfin]} \widetilde{\mathcal{E}}_{5,\nnn}(\s)^2
&\le \tfrac{14}{\eps} e^{28} \Cdata^2   +  \tfrac{94}{\eps} e^{56}\widetilde{\mathcal{D}}^2_{6,\nnn}(\sfin)
\\
\sup_{\s\in [\sin,\sfin]} \widetilde{\mathcal{E}}_{5,\ttt}(\s)^2
&\le 14 e^{36} \Cdata^2 \eps   +  \tfrac{94}{\eps} e^{56} \widetilde{\mathcal{D}}^2_{6,\ttt}(\sfin)
\end{align*}
and therefore, with \eqref{eq:tilde:E5E6:N+T-H}, we obtain that
\begin{align}
\eps \sup_{\s\in [\sin,\sfin]} \widetilde{\mathcal{E}}_{5}^2(\s)
&\le  14 e^{28} \Cdatatwo +  94 e^{56} \widetilde{\mathcal{D}}_{6,\nnn}^2(\sfin)
+
(\mathsf{K}\eps)^{-2} \bigl( 14 \eps^2  e^{28} \Cdatatwo + 94 e^{56}\widetilde{\mathcal{D}}_{6,\ttt}^2(\sfin))\big)
\notag\\
&\le 28 e^{28} \Cdatatwo+ 94 e^{56}\widetilde{\mathcal{D}}^2_6(\sfin) 
\,.
\label{eq:E5:D6:dominate-H}
\end{align}
From \eqref{eq:D5:D6:dominate-H} and \eqref{eq:E5:D6:dominate-H} we thus obtain
\begin{equation*}
\eps^{\frac 12} \sup_{\s\in [\sin,\sfin]} \widetilde{\mathcal{E}}_{5}(\s) 
+ \widetilde{\mathcal{D}}_5(\sfin)  
\leq  6 e^{14} \Cdata +10 e^{28} \widetilde{\mathcal{D}}_6(\sfin)
\end{equation*}
and so the bootstrap \eqref{bootstraps-Dnorm:6} implies that \eqref{bootstraps-Dnorm:5} holds (with a strict inequality) as soon as 
\begin{equation}
 6 e^{14} \Cdata +10 e^{28} \mathsf{B}_6 =: \mathsf{B}_5
  \,.
  \label{eq:B5B6:relation-H}
\end{equation}

\subsection{Improved estimates}
\subsubsection{The $H^6$ vorticity energy estimate}
\begin{proposition}[$H^6$ estimates for the vorticity]
\label{prop:vort:H6-H}
Let $\Omega$ be the ALE vorticity, defined in \eqref{vort-id-good}, and $\Upomega$ be the ALE specific vorticity given by \eqref{svort-id-good}.
Assume that the bootstrap assumptions~\eqref{bootstraps-H} hold, and that $\eps$ is taken to be sufficiently small to ensure $\eps^{\frac 12} \brak{\mathsf{B_J}} + \eps^{\frac 12} \brak{\mathsf{B_h}} +   \eps^{\frac 12} \brak{\mathsf{B_6}} \leq 1$.  Then, assuming $\eps$ is sufficiently small with respect to $\alpha,\kappa_0$, $\Cdata$,  and $\dl$, we have the bound
\begin{equation}
\sup_{\s\in[\sin,\sfin]} \snorm{\Jgh  \nbs^6\Upomega (\cdot,\s)}_{L^2_x}^2 
 +
 \tfrac{1}{ \eps}  \int_\sin^{\sfin}\snorm{   \nbs^6\Upomega (\cdot,\s)}_{L^2_x}^2 {\rm d}\s
\leq \Cn \eps  \brak{\mathsf{B_6}}^2  \,,
\label{eq:svort:H6-H}
\end{equation}
where the implicit constant depends only on $\alpha,\kappa_0$,  $\Cdata$, and $\dl$.
Additionally, we have that 
\begin{equation}
\sup_{\s\in[\sin,\sfin]} \snorm{\JJh  \nbs^6\Omega (\cdot,\s)}_{L^2_x}^2 
 +
 \tfrac{1}{ \eps}  \int_\sin^{\sfin}\snorm{   \nbs^6\Omega (\cdot,\s)}_{L^2_x}^2 {\rm d}\s
\leq \Cn \eps  \brak{\mathsf{B_6}}^2 \,,
\label{eq:vort:H6-H}
\end{equation}
where the implicit constant depends only on $\alpha,\kappa_0$, and  $\Cdata$.
Moreover, we have that 
\begin{subequations}
\label{eq:vorticity:pointwise-H}
\begin{align}
\snorm{\Omega}_{L^\infty_{x,\s}} 
&\leq  2^{3+\frac{2}{\alpha}} e^{18}  \Cdata
\,,
\label{eq:vorticity:pointwise:a-H}
\\
\snorm{\nbs \Omega}_{L^\infty_{x,\s}} 
&\leq 2 (4 e^{18})^{\frac{20 \cdot 23 (1+\alpha)}{\alpha}} \Cdata
\,.
\label{eq:vorticity:pointwise:b-H}
\end{align}
\end{subequations} 
\end{proposition}
 
\begin{proof}[Proof of Proposition~\ref{prop:vort:H6-H}]
First, we note that the proof of  the inequalities \eqref{eq:vorticity:pointwise:a-H} and \eqref{eq:vorticity:pointwise:b-H} of  is identical to that given in the proof of Proposition~\ref{prop:vort:H6} so that these pointwise bounds hold, and will be used for our energy estimate via the Gagliardo-Nirenberg-type
inequality  Lemma~\ref{lem:time:interpolation}.

Second, we will use the weight 
\begin{equation*}
\JJ^{2r} \,, \qquad  0<r \le \tfrac{3}{4}  \,,
\end{equation*}
for our energy method, with the intent of ultimately passing to the limit as $r \to 0$.

Third, by using \eqref{lazy-ass3}  with $F=\nbs^5\nbs_\s \Upomega$, we obtain the bound for the unweighted fifth-derivatives of
specific vorticity:
\begin{equation} 
\snorm{\nbs^5\Upomega}_{L^2_{\s,x}}   
\le 15 \eps^{\frac{1}{2}}  \| \nbs^5\Upomega(\cdot,\sin)\|_{L^2_x}
+ 28  \snorm{\JJ^r\ \nbs^5\nbs_\s\Upomega}_{L^2_{\s,x}}   
\,.
\label{svort-D5-1}
\end{equation} 

Fourth, we can now turn to the $H^6$ energy estimate.
From \eqref{D6-vort-s}, we have that 
\begin{equation} 
\tfrac{\Jg}{\Sigma} (\Q \p_\s + V \p_2)  \nbs^6  \Upomega -  \alpha    \p_1 \nbs^6  \Upomega + \alpha \Jg   g^{-\frac 12} \nbs_2 h  \nbs_2 \nbs^6  \Upomega = \RR_\Upomega
\,, \label{D6-vort-s-H}
\end{equation} 
where the remainder term $ \RR_\Upomega$ is defined in \eqref{D4-vort-s-remainder}.
For $\beta>0$ to be chosen below, we compute the spacetime $L^2$ inner-product of \eqref{D6-vort-s-H} with 
$\Sigma^{-2\beta+1} \JJ^{2r} \nbs^6 \Upomega$, 
use the identities \eqref{eq:QQQ:c-H},  \eqref{adjoint-1-H}, \eqref{adjoint-2-H}, and \eqref{adjoint-3-H}, to obtain that
\begin{align}
& \dint  \Q  \Sigma^{-2\beta} \JJrr\Jg \sabs{\nbs^6\Upomega}^2 (\cdot,\s) 
-  \dint  \Q  \Sigma^{-2\beta} \JJrr\Jg \sabs{\nbs^6\Upomega}^2 (\cdot,\sin) 
-   \tints \Sigma^{-2\beta}  \JJrr \sabs{\nbs^6\Upomega}^2 (\Q \p_\s + V \p_2)\Jg 
\notag \\
& \qquad\qquad
- 2r  \tints \Sigma^{-2\beta} \JJ^{2r-1} \Jg \sabs{\nbs^6\Upomega}^2
 \Big( (\Q \p_\s + V \p_2)\JJ - {\tfrac{\alpha }{\eps}} \Sigma \nbs_1\JJ + \alpha \Jg g^{- {\frac{1}{2}} } \nbs_2 h \, \nbs_2\JJ\Big)
\notag \\
& \qquad\qquad
- \alpha (2\beta-1) \tint  \Sigma^{-2\beta} \JJrr \Sigma,_1 \sabs{\nbs^6\Upomega}^2 
-   \tints   \JJrr \Jg \sabs{\nbs^6\Upomega}^2 (\Q \p_\s + V \p_2)\Sigma^{-2\beta}
\notag \\
&=   \tints  \Sigma^{-2\beta} \JJrr \Jg  \nbs_2 V \sabs{\nbs^6\Upomega}^2  
+ \alpha  \tints \JJrr \nbs_2\big( \Sigma^{-2\beta}  \Jg  g^{-\frac 12} \nbs_2 h\big) \sabs{\nbs^6\Upomega}^2
\notag\\
&\qquad
+ \alpha  \dint \Sigma^{-2\beta+1} \Qb\JJrr \Jg  g^{-\frac 12} \nbs_2 h\sabs{\nbs^6\Upomega}^2\Big|^{\s}_\sin
+  2 \tints \Sigma^{-2\beta+1} \JJrr \RR_\Upomega \nbs^6 \Upomega 
\,.
\label{eq:vorticity:energy:0-H}
\end{align}
We use the identities  \eqref{Jg-evo-s-P-US}, \eqref{p1-Sigma-s-P-US}--\eqref{Sigma0i-ALE-s-P-US},  and \eqref{JJ-formula}, together with  the bounds \eqref{bootstraps-H}, \eqref{eq:Q:all:bbq-H}, \eqref{eq:waitin:for:the:bus-s}; from \eqref{eq:vorticity:energy:0-H},  with $\eps$ taken sufficiently small,
 we obtain the inequality
\begin{align}
&(1-\eps^ {\frac{3}{2}} ) \dint   \Sigma^{-2\beta} \JJrr\Jg \sabs{\nbs^6\Upomega}^2 (\cdot,\s)  
+ \tfrac{  2+4\alpha\beta}{5\eps}  \tints \Sigma^{-2\beta}  \JJrr \sabs{\nbs^6\Upomega}^2 
\notag \\
&\le
(1+ \eps^{\frac{3}{2}} )  \dint   \Sigma^{-2\beta} \JJrr\Jg \sabs{\nbs^6\Upomega}^2 (\cdot,\sin)
+{\tfrac{7 + 14 \alpha \beta}{\eps}}   \tints \Sigma^{-2\beta}  \JJrr \Jg\sabs{\nbs^6\Upomega}^2
+  2 \tints \Sigma^{-2\beta+1} \JJrr \RR_\Upomega \nbs^6 \Upomega 
\,.
\label{fat-dink1}
\end{align}
In particular, in obtaining \eqref{fat-dink1}, we have used the fact that  by \eqref{JJ-formula} and \eqref{qps-JJ-bound}, 
the fourth integral on the left side of \eqref{eq:vorticity:energy:0-H} is non-negative.

We turn to the bound for $\snorm{\tfrac{\JJh}{\Sigma^\beta} \RR_\Upomega}_{L^2_{x,\s}}$ and use the decomposition
\begin{equation*} 
\RR_\Upomega
=: \RR_\Upomega^{(1)} + \RR_\Upomega^{(2)} + \RR_\Upomega^{(3)} + \RR_\Upomega^{(4)}
\end{equation*} 
which was defined in  \eqref{D4-vort-s-remainder}.  We begin with the most difficult term containing 
$\RR_\Upomega^{(1)}=- \tfrac{1}{\eps}  \jump{\nbs^6,\tfrac{\Jg}{\Sigma}} \nbs_s  \Upomega $.
From H\"older's inequality and the bootstrap assumption~\eqref{bs-Sigma},  
\begin{equation*}
 \snorm{\tfrac{\JJr}{\Sigma^\beta} \RR_\Upomega^{(1)}}_{L^2_{x,\s}} 
\les
\tfrac{1}{\eps} \snorm{\nbs \tfrac{\Jg}{\Sigma}}_{L^\infty_{x,\s}} \snorm{\tfrac{\JJr}{\Sigma^\beta}  \nbs^6 \Upomega}_{L^2_{x,\s}}
+ \tfrac{1}{\eps}
(4 \kappa_0^{-1})^\beta \sum_{i=0}^{4} \snorm{\nbs^{4-i} \nbs^2\tfrac{\Jg}{\Sigma}}_{L^{\frac{8}{4-i}}_{x,\s}} \snorm{\nbs^{i}  \nbs_s \Upomega}_{L^{\frac{8}{i}}_{x,\s}} 
\,,
\end{equation*}
where the implicit constant depends only on $\mathsf{C_{supp}}$, and hence on $\alpha$ and $\kappa_0$.    To estimate
the sum on the right side of the above inequality, we will use the interpolation  Lemma~\ref{lem:time:interpolation}.  For specific 
vorticity we will interpolate ``down to''  the pointwise bounds \eqref{eq:vorticity:pointwise-H} and   ``up to'' $\snorm{\nbs^5\Upomega}_{L^2_{x,\s}}$ and make use of the inequality
\eqref{svort-D5-1} to bound the fifth-derivatives by weighted sixth-derivatives.   By additionally
using the initial data assumptions~\eqref{table:derivatives}, the bootstrap assumptions~\eqref{bootstraps-H}, and the bounds in  Lemma~\ref{lem:time:interpolation}, we arrive at
\begin{align}
 \snorm{\tfrac{\JJr}{\Sigma^\beta} \RR_\Upomega^{(1)}}_{L^2_{x,\s}} 
&\leq
\tfrac{\Cn}{\eps}  \snorm{\tfrac{\JJr}{\Sigma^\beta}  \nbs^6 \Upomega}_{L^2_{x,\s}} 
+ \tfrac{\Cn}{\eps} (4 \kappa_0^{-1})^\beta \sum_{i=0}^{4} 
\bigl(\snorm{\nbs^{6}  \tfrac{\Jg}{\Sigma}}_{L^{2}_{x,\s}}^{\frac{4-i}{4}}  + \eps^{\frac{4-i}{4}}    \bigr)
\bigl((\kappa_0^\beta \snorm{\tfrac{\JJr}{\Sigma^\beta} \nbs^{6}  \Upomega}_{L^{2}_{x,\s}})^{\frac{i}{4}}  +  \eps^{\frac{i}{4}}  \bigr)
\notag\\
&\leq\tfrac{\Cn}{\eps} \Bigl(
 \snorm{\tfrac{\JJr}{\Sigma^\beta}  \nbs^6 \Upomega}_{L^2_{x,\s}} 
+ (\tfrac{4^4}{\kappa_0})^\beta   \snorm{\nbs^{6} \tfrac{\Jg}{\Sigma}}_{L^{2}_{x,\s}}
+ \eps      (\tfrac{4^4}{\kappa_0})^\beta \Bigr)
 \,,
 \label{corpulent1}
\end{align}
The bounds for  $\snorm{\tfrac{\JJr}{\Sigma^\beta} \RR_\Upomega^{(i)}}_{L^2_{x,\s}} $,  $i=2,3,4$ are obtained in the same manner (note that
these norms are a power of $\eps$ smaller than $\snorm{\tfrac{\JJr}{\Sigma^\beta} \RR_\Upomega^{(1)}}_{L^2_{x,\s}} $), and the modifications
are similar to those detailed in \eqref{eq:vorticity:reminder:2}.  We obtain that
\begin{align}
&{\textstyle \sum}_{i=2}^{4} \snorm{\tfrac{\JJr}{\Sigma^\beta} \RR_\Upomega^{(i)}}_{L^2_{x,\s}}
\leq 
(4 \kappa_0^{-1})^\beta {\textstyle \sum}_{i=2}^{4} \snorm{ \JJr \RR_\Upomega^{(i)}}_{L^2_{x,\s}}
\notag\\
&\leq \Cn (4 \kappa_0^{-1})^\beta
\Bigl(
 \eps \kappa_0^\beta \snorm{\tfrac{\JJr}{\Sigma^\beta} \nbs^6 \Upomega}_{L^2_{x,\s}}
+ \eps  \snorm{\nbs^6 \tfrac{\Jg}{\Sigma} }_{L^2_{x,\s}}
+ \snorm{\nbs^6  V }_{L^2_{x,\s}}
+ \eps \snorm{\nbs^6 \Jg}_{L^2_{x,\s}}
+  \snorm{\nbs^6 \nbs_2 h}_{L^2_{x,\s}} + \eps\Bigr)
 \label{corpulent2}
\,.
\end{align}
It follows from \eqref{corpulent1} and \eqref{corpulent2} and the inequality \eqref{eq:vorticity:reminder:3}, that
\begin{equation}
\snorm{\tfrac{\JJh}{\Sigma^\beta}  \RR_\Upomega}_{L^2_{x,\s}}
\leq 
\tfrac{\Cn}{\eps} (1 + \eps^2 4^\beta) \snorm{\tfrac{\JJh}{\Sigma^{\beta}}  \nbs^6 \Upomega}_{L^2_{x,\s}}  
+ \Cn (4^4 \kappa_0^{-1})^\beta \brak{\mathsf{B_6}}
\,,
\label{eq:vorticity:remainder:4-H}
\end{equation}
where $\Cn = \Cn(\alpha,\kappa_0,\Cdata) \geq 1$ is a constant independent of $\beta$.  By the Cauchy-Schwarz inequality
and \eqref{eq:vorticity:remainder:4-H}, we have that
that
\begin{align} 
2 \tints \Sigma^{-2\beta+1} \JJrr \RR_\Upomega \nbs^6 \Upomega 
\le 
\tfrac{\Cn}{\eps} (1 + \eps^2 4^\beta) \snorm{\tfrac{\JJr}{\Sigma^{\beta}}  \nbs^6 \Upomega}_{L^2_{x,\s}}^2 
+ {\tfrac{1}{\eps}}  \snorm{\JJr \nbs^6 \Upomega}_{L^2_{x,\s}}^2+  \Cn\eps (4^4 \kappa_0^{-1})^{2\beta} \brak{\mathsf{B_6}}^2 \,.
\label{R-vort-bound-H}
\end{align} 
The idea is to  choose $\beta$ sufficiently large so as to absorb the first two norms on the right side of \eqref{R-vort-bound-H} by the 
second integral on the left side of \eqref{fat-dink1}.
 Since $\Cn_{\eqref{R-vort-bound-H}}$ is independent of $\beta$, we may first  choose $\beta$ to be sufficiently large (in terms of $\alpha,\kappa_0,\Cdata$), and then $\eps$ to be sufficiently small (in terms of $\alpha,\kappa_0,\Cdata$), to ensure that
\begin{equation} 
\tfrac{4 \alpha \beta}{5} \ge \Cn_{\eqref{R-vort-bound-H}}(1 + \eps^2 4^\beta)  +  1   \,. 
\label{fat-dink2}
\end{equation} 
With $\beta=\beta(\alpha,\kappa_0,\Cdata)$ chosen so that \eqref{fat-dink2} holds,  from \eqref{fat-dink1}, we have that
\begin{align}
& \dint   \Sigma^{-2\beta} \JJrr\Jg \sabs{\nbs^6\Upomega}^2 (\cdot,\s)  
+ \tfrac{  1}{3\eps}  \tints \Sigma^{-2\beta}  \JJrr \sabs{\nbs^6\Upomega}^2 
\notag \\
&
\qquad\qquad
\le
\tfrac{1+ \eps^{\frac{3}{2}} }{1- \eps^{\frac{3}{2}} }  \dint   \Sigma^{-2\beta} \JJrr \Jg \sabs{\nbs^6\Upomega}^2 (\cdot,\sin)
+{\tfrac{\Cn}{\eps}}   \tints \Sigma^{-2\beta}  \JJrr \Jg\sabs{\nbs^6\Upomega}^2 
\,.
\label{fat-dink3}
\end{align}
From \eqref{eq:svort:H6} and the fact that $\JJ(x, \sin) =1 $, we have that
\begin{equation*}
\tfrac{1+ \eps^{\frac{3}{2}} }{1- \eps^{\frac{3}{2}} }  \dint   \Sigma^{-2\beta} \JJrr\Jg \sabs{\nbs^6\Upomega}^2 (\cdot,\sin)
\les \eps  \brak{\mathsf{B_6}}^2 \,,
\end{equation*}
and so it follows from \eqref{fat-dink3} and Gr\"onwall's inequality that
\begin{equation*}
\sup_{\s\in[\sin,\sfin]} \snorm{\tfrac{\JJr\Jgh}{\Sigma^{\beta}}   \nbs^6\Upomega (\cdot,\s)}_{L^2_x}^2 
 +
 \tfrac{1}{ \eps}  \int_\sin^{\sfin}\snorm{\tfrac{\JJr}{\Sigma^{\beta}}   \nbs^6\Upomega (\cdot,\s)}_{L^2_x}^2 {\rm d}\s
\leq \Cn \eps  \brak{\mathsf{B_6}}^2
\,,
\end{equation*}
where $\Cn = \Cn(\alpha,\kappa_0,\Cdata,\dl)>0$.   Since the right side is independent of $r$,  we may pass to the (weak) limit
as $r \to 0$ and obtain that
\begin{equation*}
\sup_{\s\in[\sin,\sfin]} \snorm{\tfrac{\Jgh}{\Sigma^{\beta}}   \nbs^6\Upomega (\cdot,\s)}_{L^2_x}^2 
 +
 \tfrac{1}{ \eps}  \int_\sin^{\sfin}\snorm{\tfrac{1}{\Sigma^{\beta}}   \nbs^6\Upomega (\cdot,\s)}_{L^2_x}^2 {\rm d}\s
\leq \Cn \eps  \brak{\mathsf{B_6}}^2
\,,
\end{equation*}
The proof of \eqref{eq:svort:H6-H} is concluded upon multiplying the above estimate by $\kappa_0^{2\beta}$ and using that
$\tfrac{\kappa_0}{4} \leq \Sigma \leq \kappa_0$.

The vorticity inequality \eqref{eq:vort:H6-H} is an immediate consequence of specific vorticity bound \eqref{eq:svort:H6-H}, the relation \eqref{svort-id-good}, the bootstrap bounds \eqref{bootstraps-H}, the sound speed inequalities \eqref{eq:Sigma:H6:new} and \eqref{eq:Sigma:H6:new:bdd}, the product and chain rules.
\end{proof}

\subsubsection{Improved estimates for $\Abn$ in the upstream spacetime}
The improved estimate for $\Abn$ is a corollary of Proposition~\ref{prop:vort:H6-H}.
\begin{corollary}
\label{cor:Abn:improve-H}
Under the standing bootstrap assumptions, we have that 
\begin{subequations}
\label{eq:Jg:Abn:D5:improve-H}
\begin{align}
\snorm{\Jgh  \nbs^5 \Abn}_{L^\infty_\s L^2_{x}}
&\les \eps^{\frac 12} \mathsf{K} \brak{\mathsf{B_6}}
\label{eq:Jg:Abn:D5:improve:a-H}
\,, \\
\snorm{\JJtf \Jgh  \nbs^6 \Abn}_{L^\infty_\s L^2_{x}}
&\les \eps^{\frac 12} \mathsf{K} \brak{\mathsf{B_6}}
\label{eq:Jg:Abn:D6:improve:b-H}
\,, \\
\snorm{\nbs^5 \Abn}_{L^2_{x,\s}} 
&\les 
\eps \mathsf{K} \brak{\mathsf{B_6}}
\label{eq:Jg:Abn:D5:improve:c-H}
\,, \\
\snorm{\JJof  \JJh \nbs^6  \Abn}_{L^2_{x,\s}} 
+\snorm{\JJtf \nbs^6  \Abn}_{L^2_{x,\s}} 
&\les 
\eps \mathsf{K} \brak{\mathsf{B_6}}
\label{eq:Jg:Abn:D6:improve:d-H}
\,.
\end{align}
\end{subequations}
where the implicit constant depends only on $\alpha,\kappa_0$,  and $\Cdata$.
\end{corollary}
\begin{proof}[Proof of Corollary~\ref{cor:Abn:improve-H}]
We first prove \eqref{eq:Jg:Abn:D6:improve:b-H} and  \eqref{eq:Jg:Abn:D6:improve:d-H}. 
From \eqref{vort-id-good}, \eqref{eq:tilde:E6},   \eqref{bootstraps-H}, \eqref{eq:vort:H6-H} and the bound $\JJ\leq 1$, we deduce that  
\begin{align*}
\snorm{\JJtf \Jgh  \nbs^6 \Abn}_{L^\infty_\s L^2_{x}}
& \les
\snorm{\Jgh  \nbs^6 \Omega}_{L^\infty_\s L^2_{x}}
+
\snorm{ \JJtf \Jgh  \nbs^6 (  \Wbt,  \Zbt)}_{L^\infty_\s L^2_{x}}
\les \eps^\frac 12  \brak{\mathsf{B_6}} +  \mathsf{K} \eps \brak{\mathsf{B_6}} \eps^{-\frac 12} \,, \\
\snorm{ \JJof\JJh \nbs^6  \Abn }_{L^2_{x,\s}} 
& \les \snorm{\nbs^6     \Omega }_{L^2_{x,\s}}  + \snorm{ \JJof\JJh \nbs^6    (  \Wbt,  \Zbt)}_{L^2_{x,\s}} 
\les \eps   \brak{\mathsf{B_6}}  + \mathsf{K} \eps \brak{\mathsf{B_6}}\,, \\
\snorm{ \JJtf \nbs^6  \Abn }_{L^2_{x,\s}} 
& \les \snorm{\nbs^6     \Omega }_{L^2_{x,\s}}  + \snorm{ \JJtf \nbs^6    (  \Wbt,  \Zbt)}_{L^2_{x,\s}} 
\les \eps   \brak{\mathsf{B_6}}  + \mathsf{K} \eps \brak{\mathsf{B_6}}\,,
\end{align*}
thereby proving  \eqref{eq:Jg:Abn:D6:improve:b-H} and  \eqref{eq:Jg:Abn:D6:improve:d-H}.

Next, we note that  \eqref{eq:Jg:Abn:D5:improve:c-H} immediately follows  from \eqref{eq:Jg:Abn:D6:improve:d-H} and  \eqref{lazy-ass3} with $r=\frac 34$ and and $F = \nbs^5 \Abn$.
Since $\p_\s = \tfrac{1}{\eps} \nbs_\s$, with \eqref{table:derivatives} and  \eqref{eq:Jg:Abn:D6:improve:d-H}, we find that
\begin{equation*}
 \| \nbs^5 \Abn \|_{L^2_{x,\s}} 
 \les \eps^{\frac 12}  \| \nbs^5 \Abn(\cdot,\sin)\|_{L^2_x} 
 + \| \JJtf \nbs_\s \nbs^5 \Abn \|_{L^2_{\s,x}} 
 \les \eps \mathsf{K} \brak{\mathsf{B}_6}
 \,,
\end{equation*}
thereby establishing \eqref{eq:Jg:Abn:D5:improve:c-H}. 
To establish \eqref{eq:Jg:Abn:D5:improve:a-H}, 
we first note that the first inequality in \eqref{eq:r:Linfty:time:2} holds in $\Hdm$ giving us
that for an $\s$-time-slice that intersects $\Hdm$, 
\begin{equation} 
\| \Jgh F (\cdot,\s) \|_{L^2_x} 
\leq 
e^9 \| \Jgh F (\cdot,\sin) \|_{L^2_x}  
+ e^9 \int_0^\eps \|\Jgh \p_\s  F(\cdot,\s)\|_{L^2_x} {\rm d}\s \,.
\label{fat-ass70}
\end{equation} 
Letting $F=\nbs^5\Omega$ in \eqref{fat-ass70}, we then have that
\begin{equation} 
\| \Jgh \nbs^5\Omega (\cdot,\s) \|_{L^2_x}
\leq 
e^9 \| \Jgh \nbs^5\Omega (\cdot,\sin) \|_{L^2_x}  
+ e^9 \eps^{-\frac{1}{2}} \|\Jgh \nbs_\s  \nbs^5\Omega\|_{L^2_{x,\s}}  \,.
\label{fat-ass71}
\end{equation} 
We now once again employ the identity \eqref{vort-id-good}. We have that
by \eqref{table:derivatives},  \eqref{bootstraps-H}, and \eqref{fat-ass71}, 
\begin{equation*}
 \sup_{\s \in [\sin,\sfin]} \|\Jgh \nbs^5 \Abn (\cdot,\s)\|_{L^2_{x}} 
 \les    \sup_{\s \in [\sin,\sfin]} \Big( \|\Jgh \nbs^5 \Omega (\cdot,\s)\|_{L^2_{x}} 
 + \tfrac{1}{2}  \|\Jgh  \nbs^5 \Wbt \|_{L^2_{\s,x}}   + \tfrac{1}{2}  \|\Jgh  \nbs^5 \Zbt \|_{L^2_{\s,x}} \Big)
 \les \eps^{\frac 12} \mathsf{K} \brak{\mathsf{B}_6}
 \,,
\end{equation*}
concluding the proof of the lemma.
\end{proof}

\subsubsection{Improved estimates for $\Zbn$ in the upstream spacetime}
\begin{lemma}
\label{lem:Zbn:improve-H}
Under the  assumptions of Proposition~\ref{prop:vort:H6-H}, we have that  for any $\bar \beta >0$, and $\bar a \in [0,\frac 12]$, 
\begin{subequations}
\label{eq:Jg:Zbn:D5:improve-H}
\begin{align}
\snorm{ \Jgh \nbs^5 \Zbn }_{L^\infty _\s L^2_x}^2
+ \tfrac{1}{\eps} \int_0^\eps
\snorm{\nbs^5  \Zbn (\s)}_{L^2_{x}}^2 {\rm d}\s
&\leq \Cn \eps \brak{\mathsf{B}_6}^2 \,, 
\label{eq:Jg:Zbn:D5:improve:a-H}\\
\snorm{\nbs^4 \Zbn }_{L^\infty _\s L^2_x}
&\leq \Cn \eps^{\frac{1}{2}}   \brak{\mathsf{B}_6} \,, 
\label{eq:Jg:Zbn:D3:improve-H} \\
\snorm{\Sigma^{-\bar \beta} \JJ^{\frac 34 - \bar a} \Jg^{\! \bar a} \nbs_1 \nbs^5 \Zbn }_{L^2_{x,\s}} 
&\le \tfrac{1}{2\alpha} \| \Sigma^{-1-\bar \beta} \JJ^{\frac 34 - \bar a} \Jg^{\! \bar a} \nbs^5 \nbs_\s (\Jg \Zbn) \|_{L^2_{x,\s}}
+ \Cn \eps  (\tfrac{4}{\kappa_0})^{\bar \beta}    \mathsf{K} \brak{\mathsf{B_6}}\,,
\label{eq:Jg:Zbn:p1D5:improve-H}  \\
\snorm{\Sigma^{-\bar \beta} \JJ^{\frac 34 - \bar a} \Jg^{\! \bar a}  \nbs_2 \nbs^5 \Zbn }_{L^2_{x,\s}}  
&\le \tfrac{1}{\alpha  \eps} \|\Sigma^{-1 - \bar \beta } \JJ^{\frac 34 - \bar a} \Jg^{\! \bar a} \nbs^5 \nbs_\s \Zbt \|_{L^2_{x,\s}} + \Cn \eps (\tfrac{4}{\kappa_0})^{\bar \beta} \mathsf{K} \brak{\mathsf{B_6}}  
\label{eq:Jg:Zbn:p2D5:improve:new-H} \,,
\end{align}
where the implicit constant $\Cn$ depends only on $\alpha,\kappa_0$, and $\Cdata$. 
It is also convenient to record the estimates 
\begin{align}
\snorm{\JJtf \Jgh \nbs_1 \nbs^5 \Zbn }_{L^\infty_\s L^2_{x}} 
&\les \eps^{-\frac 12} \brak{\mathsf{B_6}}\,,
\label{eq:Jg:Zbn:p1D6:sup:improve:new-H}  \\
\snorm{\JJtf \Jgh \nbs_2 \nbs^5 \Zbn }_{L^\infty_\s L^2_{x}}  
&\les  \eps^{-\frac 12}   \brak{\mathsf{B_6}} \,,
\label{eq:Jg:Zbn:p2D6:sup:improve:new-H}
\end{align}
\end{subequations}
where the implicit constant and the constant depends only on $\alpha,\kappa_0$, and $\Cdata$. 
\end{lemma}
\begin{proof}[Proof of Lemma~\ref{lem:Zbn:improve-H}]
The proof of \eqref{eq:Jg:Zbn:D3:improve-H}--\eqref{eq:Jg:Zbn:p2D5:improve:new-H}
is identical to that given in the proof of  Lemma~\ref{lem:Zbn:improve}.
The proof of \eqref{eq:Jg:Zbn:D5:improve:a-H} requires the following modification: We test the equation
\eqref{eq:Zbn:improv:1} with $\Sigma^{-2\beta+1} \JJrr \nbs^5\Zbn$ for $0<r\le \tfrac{3}{4} $. The energy estimate is performed in the same
manner as our weighted energy estimate for specific vorticity in the proof of Proposition~\ref{prop:vort:H6-H}; namely, we obtain uniform-in-$r$
estimates for norms weighted by $\JJr$ and then pass to the limit as $r \to 0$.  The  inequality \eqref{eq:Jg:Zbn:D5:improve:a-H} is an 
immediate consequence.
\end{proof}

\subsubsection{Improved estimates for $\Wbn$ in the upstream spacetime}
\begin{lemma}
\label{lem:duck:a:fck-H}
Under the  assumptions of Proposition~\ref{prop:vort:H6-H},
\begin{subequations}
\label{eq:Jg:Wbn:improve:material-H}
\begin{align}
 \|(\Q\p_\s + V\p_2) ( \Jg \Wbn)\|_{L^\infty_{x,\s}} 
 &\les 1
 \label{eq:Jg:Wbn:improve:material:a-H}
 \\
\|  \nbs^4(\Q\p_\s + V\p_2) (\Jg \Wbn)\|_{L^2_{x,\s}}
&\les  \mathsf{K} \brak{\mathsf{B_6}}
\label{eq:Jg:Wbn:improve:material:b-H}
\\
\| \JJof \nbs^5 (\Q\p_\s + V\p_2)( \Jg \Wbn) \|_{L^2_{x,\s}}
&\les \mathsf{K} \brak{\mathsf{B_6}}
\label{eq:Jg:Wbn:improve:material:c-H} 
\end{align}
\end{subequations}
and we also have that
\begin{subequations}
\label{Wbn:improved-H}
\begin{align}
\snorm{\nbs^5 (\Jg \Wbn)}_{L^\infty_\s L^2_x}
&\leq 2 \Cdata \eps^{-\frac 12}
\,,
\label{eq:D5:JgWbn-H}
\\
\snorm{ \JJof \Jgh \nbs_\s \nbs^5 (\Jg\Wbn )}_{L^2_{x,\s}} 
&\leq \Cn \eps \mathsf{K} \brak{\mathsf{B}_6}
\,.   
\label{eq:Ds:D4:Jg:Wbn:new-H}
\end{align} 
\end{subequations}
\end{lemma}
\begin{proof}[Proof of Lemma~\ref{lem:duck:a:fck-H}]
The proof of the inequalities \eqref{eq:Jg:Wbn:improve:material-H} is identical to the proof of \eqref{eq:Jg:Wbn:improve:material} and will not be repeated. 

The proof of \eqref{eq:D5:JgWbn-H} requires an upstream modification to the establishment of the inequality  \eqref{last-pain-in-the-ass}.  
We apply $\nb^5$ to \eqref{eq:Jg:Wb:nn} and then transform to upstream $(x,\s)$ variables, to obtain that
\begin{align}
(\Q  \p_\s + V \p_2) \nbs^5 (\Jg \Wbn)
+\alpha \Sigma g^{-\frac 12} \Jg \nbs^5  \nbs_2 \Abn 
= \mathsf{F}_{5} \,,
\label{fat-ass1}
\end{align}
and following the bounds \eqref{eq:duck:a:fck:2a} and \eqref{eq:duck:a:fck:2b}, we have that
\begin{align} 
\snorm{ \mathsf{F}_{5}}_{L^2_{x,\s}} \le  \mathsf{K}  \brak{\mathsf{B}_6} \,.
\label{fat-ass6}
\end{align} 
We test \eqref{fat-ass1} with $\nbs^5(\Jg\Wbn)$ which results in 
\begin{align} 
&\dint \Q\sabs{\nbs^5(\Jg\Wbn)}^2(x,\s)   
- \int_\sin^\s \int_{x_2=-\pi}^\pi  (\Q\p_\s+V\p_2) \thd(x_2,\s') \sabs{\nbs^5(\Jg\Wbn) (\thd(x_2,\s'),x_2,\s')}^2  {\rm d}x_2 {\rm d}\s'
\notag \\
& \qquad \qquad
= \dint\Q\sabs{\nbs^5(\Jg\Wbn)}^2(x,\sin)   
 -2\alpha \tints \Sigma g^{- {\frac{1}{2}} } \Jg \nbs^5\nbs_2\Abn \ \nbs^5(\Jg\Wbn)
+ \tints  \mathsf{F}_{5}  \nbs^5(\Jg\Wbn) 
\notag \\
&\qquad\qquad \qquad\qquad
+  \tints (\Qc +V,_2)  \sabs{\nbs^5(\Jg\Wbn)}^2  \,.
\label{fat-ass0}
\end{align}

The most difficult term in~\eqref{fat-ass0} is the second term on the right side, which we estimate next. We split the integral (see the integral notation given in \eqref{int-Hdmp-Hdmm})
\begin{align*} 
\int_{\Hdm} \Sigma g^{- {\frac{1}{2}} } \Jgh \nbs^5\nbs_2\Abn \ \Jgh\nbs^5(\Jg\Wbn)
&=\left( \int_{\Hdmm} +  \int_{\Hdmp} \right)\Sigma g^{- {\frac{1}{2}} } \Jgh \nbs^5\nbs_2\Abn \ \Jgh\nbs^5(\Jg\Wbn)
\,.
\end{align*} 
Since $\JJ\ge 1$ in $\Hdmm$, the analysis of that integral is straightforward and follows the analysis of \eqref{last-pain-in-the-ass} (and we will provide the bound for this integral below); hence, we shall
focus on the integral over $\Hdmp$. Using the change of variables formula from~\eqref{L2sx-fubini}, the bootstrap bounds, \eqref{ps-Theta}, the definition~\eqref{L2sx-composite}, the Cauchy-Schwartz inequality,   the lower bound in~\eqref{JJ-and-s},  and the estimate \eqref{lazy-ass5} with $F= \nbs^5(\Jg \Wbn)$, we deduce 
\begin{align}
&\left| \int_{\Hdmp} \Sigma g^{- {\frac{1}{2}} } \Jgh \nbs^5\nbs_2\Abn \ \Jgh\nbs^5(\Jg\Wbn)  \right|
\notag\\
&\leq \int_{\Omega_{\mathsf{US},+}} \p_\s \Thd(x,\s) \sabs{\Sigma g^{- {\frac{1}{2}} } \Jgh \nbs^5\nbs_2\Abn \ \Jgh\nbs^5(\Jg\Wbn) }(x,\Thd(x,\s)) {\rm d}x {\rm d}\s
\notag\\
&\leq 2 \kappa_0 \sup_{\s \in [\sin,0]} \snorm{ \Jgh\nbs^5(\Jg\Wbn)(\cdot, \Thd(\cdot ,\s)) }_{L^2_x} 
\snorm{(\JJof \Jgh \nbs^5\nbs_2\Abn)(\cdot,\Thd(\cdot,\cdot))}_{L^2_{x,\s}}
\left(\int_{\sin}^{0} \bigl(\tfrac{89}{100} \cdot \tfrac{\s}{\sin}\bigr)^{-\frac 12} {\rm d}\s\right)^{\frac 12} 
\notag\\
&\leq \Cn \eps^{\frac 12} \Bigl(  \snorm{\Jgh \nbs^5(\Jg\Wbn( \cdot , \sin))}_{L^2_x} + \eps^ {\frac{1}{2}} \snorm{\JJof \Jgh \p_\s \nbs^5(\Jg\Wbn)}_{L^2_{\s,x}}  \Bigr) 
\snorm{\JJof \Jgh \nbs^5\nbs_2\Abn}_{L^2_{x,\s}}
\notag\\
&\leq \Cn \eps \mathsf{K} \brak{\mathsf{B}_6}^2
\,.
\label{fat-ass4}
\end{align}
In the last inequality we have used~\eqref{boots-H} and \eqref{eq:Jg:Abn:D6:improve:d-H}.

Using a similar argument, appealing to \eqref{eq:lazy-ass:000} and the fact that $\JJ \geq 1$ in $\Hdmm$ due to~\eqref{JJ-and-s-new}, we have that
\begin{align*} 
\left| \int_{\Hdmm} \Sigma g^{- {\frac{1}{2}} } \Jgh \nbs^5\nbs_2\Abn \ \Jgh\nbs^5(\Jg\Wbn)\right|
 \le \eps \Cn   \mathsf{K} \brak{\mathsf{B_6}}^2 \,,
\end{align*} 
and hence together with  \eqref{fat-ass4}, we have established that
\begin{align} 
\int_{\Hdm} \Sigma g^{- {\frac{1}{2}} } \Jgh \nbs^5\nbs_2\Abn \ \Jgh\nbs^5(\Jg\Wbn)
 \le \eps \Cn   \mathsf{K} \brak{\mathsf{B_6}}^2 \,.
 \label{fat-ass5}
 \end{align} 
 
 The inequality \eqref{fat-ass5} bounds the second integral on the right side of \eqref{fat-ass0}.   To bound the third and fourth integrals
 on  the right side of \eqref{fat-ass0}, we use the bootstrap bounds \eqref{bootstraps-H}, the bound \eqref{fat-ass6}, and the Cauchy-Young inequality.
 The bound for the first integral on  the right side of \eqref{fat-ass0} is bounded by \eqref{table:derivatives}.    Collecting these bounds and taking
 $\eps$ sufficiently small, we establish the inequality \eqref{eq:D5:JgWbn-H}.
 
 The proof of the bound \eqref{eq:D5:JgWbn-H}
 is identical to the proof of the bound \eqref{eq:Ds:D4:Jg:Wbn:new} and will not be repeated.
\end{proof}

\subsection{Closing the pointwise bootstraps in the upstream geometry}
 
The pointwise bounds which were previously established in Section~\ref{sec:pointwise:bootstraps}, only relied on the evolution equations in $(x,\s)$ coordinates (cf.~\eqref{Jg-evo-s}--\eqref{nn-tt-evo-s}), the bootstrap assumptions, the functional analytic setup in Appendix~\ref{app:functional}, and of the $L^\infty$ estimates from Appendix~\ref{sec:app:transport}. 
These arguments apply {\bf as is} in the geometry of the downstream development, except that we refer to the evolution identities~\eqref{fundamental-H}, the bootstraps~\eqref{bootstraps-H}, to  the functional analytic bounds from Section~\ref{rem:app:upstream:flat}, and to the $L^\infty$ estimates from Section~\ref{app:upstream:Lp}. We omit these redundant details concerning the closure of the pointwise bootstraps.

\subsection{Upstream energy estimates}
It remains for us to close the bootstrap~\eqref{boots-H} (see~\eqref{bootstraps-Dnorm:6}) for the sixth order energy $\widetilde{\mathcal{E}}_6$ and damping $\widetilde{\mathcal{D}}_6$ norms, defined earlier in~\eqref{eq:tilde:E5E6:N+T-H} and~\eqref{eq:tilde:D5D6:N+T-H}. 

Previously, this was achieved by separately establishing a bound for the tangential parts of the energy $\widetilde{\mathcal{E}}_{6,\ttt}$ and damping  $\widetilde{\mathcal{D}}_{6,\ttt}$ in~Section~\ref{sec:sixth:order:energy-tangential}, and the normal parts of the energy $\widetilde{\mathcal{E}}_{6,\nnn}$ and damping  $\widetilde{\mathcal{D}}_{6,\nnn}$ in~Section~\ref{sec:sixth:order:energy}. In turn, these estimates required that we established improved energy bounds for six ``pure time derivatives'' in Section~\ref{sec:pure:time}. 

Just as for the downstream  development, for the upstream development,  we follow the same exact methodology.
As before, the tangential bounds from~Section~\ref{sec:sixth:order:energy-tangential} and normal energy estimates from Section~\ref{sec:sixth:order:energy} run in parallel, the only difference being that the fundamental variables are un-weighted for the tangential part (i.e.~$(\Wbt,\Zbt,\Abt)$) and are $\Jg$-weighted for the normal part (i.e.~$(\Jg \Wbn,\Jg \Zbn, \Jg \Abn)$). The special estimates for six pure time derivatives from Section~\ref{sec:pure:time} are used in the same way, to treat the remainders $\mathcal{R} _\Zb^\tt$ and $\mathcal{R} _\Zb^\nn$, in the $\Zbt$, and respectively $\Jg \Zbn$ equations. 

Since the tangential and normal energy estimates run in parallel (similarities and differences may be seen by comparing 
Sections~\ref{sec:sixth:order:energy-tangential} and~\ref{sec:sixth:order:energy}), we do not repeat both of these two sets of energy estimates for the 
upstream geometry.   Just as in our downstream analysis, 
 the upstream modifications to the tangential component energy estimates are identical to the modifications made to the normal component energy estimates. For conciseness, we shall therefore only provide details for the upstream normal-component energy estimates  (see Section~\ref{eq:upstream:normal} below). 

\subsubsection{Sixth order tangential energy estimates}
For the tangential energy estimates, at this point we simply record that by repeating the arguments from Section~\ref{sec:sixth:order:energy-tangential}, with the modifications outlined in Section~\ref{eq:upstream:normal} below (see the argument leading to~\eqref{eq:normal:conclusion:1-H}--\eqref{eq:normal:conclusion:3-H}), 
similarly to \eqref{eq:hate:11:a}--\eqref{eq:hate:12}, there exists a constant 
\begin{equation*}
 \hat{\mathsf{c}}_{\alpha,\kappa_0,\dl} > 0 \,,
\end{equation*}
which depends only on $\alpha$, $\dl$,  and  $\kappa_0$, and may be computed explicitly, such that 
\begin{align}
& 
\sup_{\s \in [0,\eps]}
\snorm{ \JJ^{\frac 34}  \Jg^{\!\frac 12} \nbs^6(\Wbt,\Zbt,\Abt)(\cdot,\s)}_{L^2_x}^2 
 + \tfrac{1}{\eps}   \int_0^{\eps}  
\snorm{\JJ^{\frac 14} \Jg^{\!\frac 12} \nbs^6 (\Wbt,\Zbt,\Abt)(\cdot,\s)}_{L^2_x}^2   {\rm d} \s
 \notag\\
 &\qquad 
+ \tfrac{1}{\eps^2} \sup_{\s \in [0,\eps]} 
 \snorm{ \JJ^{\frac 14} \nn \cdot  \nbs^6 \tt (\cdot,\s)}_{L^2_{x}}^2
+ \tfrac{1}{\eps^3} \int_0^{\eps} \snorm{ \JJ^{-\frac 14} \nn \cdot  \nbs^6 \tt (\cdot,\s)}_{L^2_{x}}^2 {\rm d} \s
\notag\\
 &  
 \leq \hat{\mathsf{c}}_{\alpha,\kappa_0,\dl} \cdot \eps  \Big( \Cdatatwo + \mathsf{B}_6^2 + \Cn \eps^{\frac 12}  \mathsf{K}^2  \brak{\mathsf{B}_6}^2\Bigr)
 \,.
 \label{eq:hate:12-H}
\end{align}
Then, as in~\eqref{eq:hate:13:aa}--\eqref{eq:hate:13:a}, upon ensuring that 
\begin{align}
\mathsf{B}_6 &\geq \max\{1, \Cdata \} \,,
\label{eq:B6:choice:1-H}
\end{align}
and upon defining
\begin{align}
\mathsf{K} &:= 8 \max\{1, \hat{\mathsf{c}}_{\alpha,\kappa_0,\dl}^{\frac 12} \} \,,
\label{eq:K:choice:1-H}
\end{align}
by letting $\eps$ be sufficiently small in terms of $\alpha,\kappa_0$, $\Cdata$, and $\dl$, we deduce from~\eqref{eq:hate:12-P} that 
\begin{align}
\eps \sup_{\s \in [0,\eps]} \widetilde{\mathcal{E}}_{6,\ttt}^2(\s)
+\widetilde{\mathcal{D}}_{6,\ttt}^2(\eps) 
\leq \tfrac{1}{8} (\eps \mathsf{K})^2 \mathsf{B}_6^2 
 \,.
 \label{eq:hate:13-H}
\end{align}
This bound is the same as~\eqref{eq:hate:13:a}. It closes the ``tangential part'' of the remaining bootstrap~\eqref{boots-H} for $\widetilde{\mathcal{E}}_{6}$ and $\widetilde{\mathcal{D}}_{6}$.

\subsubsection{Sixth order pure-time energy estimates}
The energy estimates for the case of pure-time derivatives following  the arguments from Section~\ref{sec:pure:time}, with the modifications outlined in Section~\ref{eq:upstream:normal} below, the same bound as given in \eqref{eq:madman:2} holds and therefore,
\begin{equation*} 
\eps^{\frac 12} 
\snorm{\JJ^{\! \frac 34}\Jgh \nbs_{\s}^6  \Zbn }_{L^\infty_x L^2_{x}}
+ 
\snorm{\JJ^{\! \frac 34} \nbs_{\s}^6  \Zbn }_{L^2_{x,\s}} 
\leq
\eps^{\frac 12} 
\snorm{\JJ^{\! \frac 34}\Jgh \nbs_{\s}^6  \Zbn }_{L^\infty_x L^2_{x}}
+
\snorm{\JJ^{\! \frac 14} \Jgh \nbs_{\s}^6  \Zbn }_{L^2_{x,\s}} 
\les \eps \mathsf{K} \brak{\mathsf{B}_6} 
\,.
\end{equation*}


\subsection{Upstream energy estimates for normal components}\label{eq:upstream:normal}
We continue to use the equation set \eqref{energy-WZA-s}, with the operator $\nbs^6$ employed in the following manner:
The energy identity \eqref{D6-L2-N} is then replaced with the upstream energy identity
\begin{align} 
&
\tints  \JJss \Big( \underbrace{\eqref{energy-Wn-s} \ \Jg  \nbs^6(\Jg\Wbn)} _{ I^{\WW_n}}
+ \underbrace{\eqref{energy-Zn-s} \  \nbs^6(\Jg\Zbn)} _{ I^{\ZZ_n}}
+ \underbrace{2\eqref{energy-An-s} \ \nbs^6(  \Jg \Abn) } _{ I^{\AA_n}}
 \Big) \Sigma^{-2\beta+1}{\rm d}x {\rm d}\s'=0 \,, \label{D5-L2-H}
\end{align} 
where once again $\beta>0$ is a constant,  whose value will be made precise below, and 
$\jb = \Sigma^{-2\beta+1}$.   

\subsubsection{The integral  $I^{\WW_n}$} 
 We additively decompose the integral  $I^{\WW_n}$.  For $\s \in [\sin,\sfin]$, we have that
 \begin{subequations} 
 \label{Integral-Wbn-H}
\begin{align}
I^{\WW_n}&= I^{\WW_n}_1+I^{\WW_n}_3+I^{\WW_n}_4 +I^{\WW_n}_5+I^{\WW_n}_6
\notag \,, \\
 I^{\WW_n}_1 &=
\tints \tfrac{1}{\Sigma^{2\beta}} \JJss \Jg (\Q\p_\s +V\p_2)\nbs^6( \Jg\Wbn  ) \  \nbs^6(\Jg\Wbn)
\,, \label{I1-Wbn-H} \\
I^{\WW_n}_3 &=
 \alpha \tints \jb \JJss \Jg g^{- {\frac{1}{2}} }  \nbs_2\nbs^6(\Jg\Abn)  \  \nbs^6(\Jg\Wbn)
\,, \label{I3-Wbn-H} \\
 I^{\WW_n}_4 &=
-\alpha \tints \jb  \JJss \Jg g^{- {\frac{1}{2}} }  \Abn \nbs_2\nbs^6 \Jg \  \nbs^6(\Jg\Wbn)
  \,,   \label{I4-Wbn-H} \\
 I^{\WW_n}_5 &=
-\tfrac{\alpha }{2}  \tints \jb  \JJss  \Jg g^{- {\frac{1}{2}} } (\Jg\Wbn + \Jg\Zbn - 2\Jg\Abt)\nbs_2 \nbs^6 \tt\cdo\nn \ \nbs^6(\Jg\Wbn)
 \,, \label{I5-Wbn-H} \\
 I^{\WW_n}_6 &=
-\tfrac{\alpha }{2}  \tints \jb \JJss \Jg \big(  \nbs^6 \Fwn + \mathcal{R}^\nn_{\Wb} + \mathcal{C}^\nn_{\Wb}\big) \ \nbs^6(\Jg\Wbn)
 \,. \label{I6-Wbn-H}
\end{align} 
\end{subequations}

\subsubsection{The integral  $I^{\ZZ_n}$} 
 We additively decompose the integral  $I^{\ZZ_n}$ as 
 \begin{subequations} 
 \label{Integral-Zbn-H}
 \begin{align}
  I^{\ZZ_n} & = I^{\ZZ_n}_1 + I^{\ZZ_n}_2 + I^{\ZZ_n}_3+ I^{\ZZ_n}_4+ I^{\ZZ_n}_5+ I^{\ZZ_n}_6+ I^{\ZZ_n}_7+ I^{\ZZ_n}_8
  + I^{\ZZ_n}_9+ I^{\ZZ_n}_{10}
  \,,  \notag \\
 I^{\ZZ_n}_1& =
 \tints \tfrac{1}{\Sigma^{2\beta}} \JJss  \Jg(\Q\p_\s +V\p_2)\nbs^6( \Jg\Zbn) \  \nbs^6(\Jg\Zbn)
  \,, \label{I1-Zbn-H}\\
 I^{\ZZ_n}_2& =
- \tints \tfrac{1}{\Sigma^{2\beta}}  \bubu{  \alpha (\Jg\Wbn - \Jg\Zbn)} \JJss\sabs{\nbs^6(\Jg\Zbn)  }^2
  \,, \label{I2-Zbn-H} \\
I^{\ZZ_n}_3 &=
- \alpha \tints \jb  g^{- {\frac{1}{2}} }\JJss \Jg  \nbs_2\nbs^6(\Jg\Abn)  \  \nbs^6(\Jg\Zbn)
\,, \label{I3-Zbn-H} \\
 I^{\ZZ_n}_4 &=
\alpha \tints \jb\JJss \Jg g^{- {\frac{1}{2}} }  \Abn \nbs_2\nbs^6 \Jg \  \nbs^6(\Jg\Zbn)
  \,,   \label{I4-Zbn-H} \\
 I^{\ZZ_n}_5 &=
\tfrac{\alpha }{2}  \tints \jb \JJss \Jg g^{- {\frac{1}{2}} } (\Jg\Wbn + \Jg\Zbn - 2\Jg\Abt)\nbs_2 \nbs^6 \tt\cdo\nn \ \nbs^6(\Jg\Zbn)
 \,, \label{I5-Zbn-H} 
\\
 I^{\ZZ_n}_6& =
-{\tfrac{2 \alpha }{\eps}} \tints  \jb \JJss   \nbs_1 \nbs^6 (\Jg\Zbn) \ \nbs^6(\Jg\Zbn)
\,,  \label{I6-Zbn-H}\\
 I^{\ZZ_n}_7& =
-{\tfrac{2 \alpha }{\eps}} \tints  \jb  \JJss  \Jg(\Abn+\Zbt) \Big( \nbs_1\nbs^6\tt \cdo\nn - \eps \Jg g^{- {\frac{1}{2}} } \nbs_2h\, \nbs_2\nbs^6\tt \cdo\nn   \Big)   \ \nbs^6(\Jg\Zbn)
\,,  \label{I7-Zbn-H}\\
 I^{\ZZ_n}_8& =
{\tfrac{2 \alpha }{\eps}} \tints  \jb \JJss  \Zbn  \Big( \nbs_1\nbs^6\Jg -  \eps \Jg g^{- {\frac{1}{2}} } \nbs_2 h\,  \nbs_2\nbs^6\Jg \Big) \ \nbs^6(\Jg\Zbn)
\,,  \label{I8-Zbn-H}\\
 I^{\ZZ_n}_9& =
2 \alpha \tints  \jb  \JJss \Jg g^{- {\frac{1}{2}} } \nbs_2h  \nbs_2 \nbs^6 (\Jg\Zbn) \ \nbs^6(\Jg\Zbn)
\,,  \label{I9-Zbn-H}\\
 I^{\ZZ_n}_{10}& =
 - \tints  \jb \JJss \big(  \nbs^6 \Fzn + \mathcal{R}^\nn_{\Zb} + \mathcal{C}^\nn_{\Zb}\big) \  \nbs^6(\Jg\Zbn)
\,.  \label{I10-Zbn-H}
\end{align} 
\end{subequations} 

 \subsubsection{The integral  $I^{\AA_n}$}  
We additively decompose the integral  $I^{\AA_n}$ as
\begin{subequations} 
\label{Integral-Abn-H}
 \begin{align}
  I^{\AA_n} & = I^{\AA_n}_1 + I^{\AA_n}_2 + I^{\AA_n}_3+ I^{\AA_n}_4+ I^{\AA_n}_5+ I^{\AA_n}_6+ I^{\AA_n}_7+ I^{\AA_n}_8
  + I^{\AA_n}_9+ I^{\AA_n}_{10}
  \,,  \notag \\
 I^{\AA_n}_1& =
2 \tints \tfrac{1}{\Sigma^{2\beta}}  \JJss \Jg(\Q\p_\s +V\p_2)\nbs^6( \Jg\Abn) \  \nbs^6(\Jg\Abn)
  \,, \label{I1-Abn-H}\\
 I^{\AA_n}_2& =
- \tints \tfrac{1}{\Sigma^{2\beta}}  \bubu{  \alpha  \big(\Jg\Wbn - \Jg\Zbn\big) } \JJss\sabs{\nbs^6(\Jg\Abn)}^2
  \,, \label{I2-Abn-H} \\
I^{\AA_n}_3 &=
2 \alpha \tints \jb  g^{- {\frac{1}{2}} } \JJss \Jg  \nbs_2\nbs^6(\Jg\Sbn)  \  \nbs^6(\Jg\Abn)
\,, \label{I3-Abn-H} \\
 I^{\AA_n}_4 &=
-2\alpha \tints \jb  g^{- {\frac{1}{2}} } \JJss (\Jg \Sbn) \nbs_2\nbs^6 \Jg \  \nbs^6(\Jg\Abn)
  \,,   \label{I4-Abn-H} \\
 I^{\AA_n}_5 &=
2 \alpha \tints \jb \JJss \Jg g^{- {\frac{1}{2}} }  \Sbt \nbs_2 \nbs^6 \tt\cdo\nn \ \nbs^6(\Jg\Abn)
 \,, \label{I5-Abn-H} 
\\
 I^{\AA_n}_6& =
-{\tfrac{2 \alpha }{\eps}} \tints  \jb   \JJss \nbs_1 \nbs^6 (\Jg\Abn) \ \nbs^6(\Jg\Abn)
\,,  \label{I6-Abn-H}\\
 I^{\AA_n}_7& =
{\tfrac{\alpha }{\eps}} \tints  \jb \JJss   (\Jg\Wbn+\Jg\Zbn - 2\Jg\Abt) \Big( \nbs_1\nbs^6\tt \cdo\nn - \eps \Jg g^{- {\frac{1}{2}} } \nbs_2h\, \nbs_2\nbs^6\tt \cdo\nn   \Big)  \ \nbs^6(\Jg\Abn)
\,,  \label{I7-Abn-H}\\
 I^{\AA_n}_8& =
{\tfrac{2 \alpha }{\eps}} \tints  \jb  \JJss  \Abn \Big( \nbs_1\nbs^6\Jg -  \eps \Jg g^{- {\frac{1}{2}} } \nbs_2 h\,  \nbs_2\nbs^6\Jg \Big)  \ \nbs^6(\Jg\Abn)
\,,  \label{I8-Abn-H}\\
 I^{\AA_n}_9& =
2 \alpha \tints  \jb  \JJss \Jg g^{- {\frac{1}{2}} } \nbs_2h  \nbs_2 \nbs^6 (\Jg\Abn) \ \nbs^6(\Jg\Abn)
\,,  \label{I9-Abn-H}\\
 I^{\AA_n}_{10}& =
 - \tints  \jb \JJss \big(  \nbs^6 \Fzn + \mathcal{R}^\nn_{\Zb} + \mathcal{C}^\nn_{\Zb}\big) \  \nbs^6(\Jg\Abn)
\,.  \label{I10-Abn-H}
 \end{align} 
 \end{subequations}

\subsubsection{Upstream analysis of $I^{\WW_n}_1$, $I^{\ZZ_n}_1$, and  $I^{\AA_n}_1$}  
 From \eqref{I1-Wbn-H}, by applying \eqref{eq:adjoints-H} and \eqref{qps-JJ-H} and using the fact that from \eqref{JJ-vanish} $\JJ =0$ on the surface $x_1=\ths$, we have that
 \begin{align} 
\big(I^{\WW_n}_1, I^{\ZZ_n}_1, I^{\AA_n}_1\big) 
&= \big(I^{\WW_n}_{1,a}, I^{\ZZ_n}_{1,a}, I^{\AA_n}_{1,a}\big)
+\big(I^{\WW_n}_{1,b}, I^{\ZZ_n}_{1,b}, I^{\AA_n}_{1,b}\big)
+\big(I^{\WW_n}_{1,c}, I^{\ZZ_n}_{1,c}, I^{\AA_n}_{1,c}\big)
+ \big(I^{\WW_n}_{1,d}, I^{\ZZ_n}_{1,d}, I^{\AA_n}_{1,d}\big) \,, 
\notag \\
\big(I^{\WW_n}_{1,a}, I^{\ZZ_n}_{1,a}, I^{\AA_n}_{1,a}\big) 
&=\tfrac{1}{2} \dint \tfrac{\Q}{\Sigma^{2\beta}}   \JJss \Jg \big( \sabs{\nbs^6 (\Jg\Wbn) }^2 , \sabs{\nbs^6( \Jg\Zbn )}^2, 2\sabs{\nbs^6 (\Jg\Abn) }^2\big)\Big|^\s_\sin \,, 
\notag \\
\big(I^{\WW_n}_{1,b}, I^{\ZZ_n}_{1,b}, I^{\AA_n}_{1,b}\big) 
&= -\tfrac{1}{2} \tints \tfrac{1}{\Sigma^{2\beta}} \JJss 
(\Q\p_\s+V\p_2) \Jg \big( \sabs{\nbs^6( \Jg\Wbn) }^2 , \sabs{\nbs^6 (\Jg\Zbn )}^2, 2\sabs{\nbs^6 (\Jg\Abn) }^2\big)   \,,
\notag \\
\big(I^{\WW_n}_{1,c}, I^{\ZZ_n}_{1,c}, I^{\AA_n}_{1,c}\big) 
&= -\tfrac{3}{4} \tints \tfrac{1}{\Sigma^{2\beta}} \JJh \Jg
(\Q\p_\s+V\p_2) \JJ \big( \sabs{\nbs^6( \Jg\Wbn) }^2 , \sabs{\nbs^6 (\Jg\Zbn) }^2, 2\sabs{\nbs^6 (\Jg\Abn) }^2\big)   \,,
\label{I-Wbn-Zbn-Abn-1c}\\
\big(I^{\WW_n}_{1,d}, I^{\ZZ_n}_{1,d}, I^{\AA_n}_{1,d}\big) 
&=- \!\! \tints \!\! \tfrac{1}{\Sigma^{2\beta}} \big( \alpha \beta  (\Zbn +\Abt) + \tfrac{1}{2} (\Qc +V,_2) \big) \JJss\Jg
\big( \sabs{\nbs^6 (\Jg\Wbn) }^2 , \sabs{\nbs^6 (\Jg\Zbn) }^2,2 \sabs{\nbs^6 (\Jg\Abn )}^2\big)
  \,. \notag
\end{align}

 \subsubsection{Upstream analysis of $\big(I^{\ZZ_n}_1,I^{\AA_n}_1\big)+\big(I^{\ZZ_n}_2,I^{\AA_n}_2\big)
 +\big(I^{\ZZ_n}_6,I^{\AA_n}_6\big)+\big(I^{\ZZ_n}_9,I^{\AA_n}_9\big)$} 
We now employ integration-by-parts with respect to $\nbs_1$ and $\nbs_2$ using \eqref{eq:adjoints-H}. From \eqref{JJ-vanish}, 
 $\JJ =0$ on the surface $x_1=\ths$.   Hence, we obtain that
\begin{subequations} 
 \label{Zbn-Zbn-6and9}
\begin{align} 
\big(I^{\ZZ_n}_6,I^{\AA_n}_6\big)& 
=  \big(I^{\ZZ_n}_{6,a} ,I^{\AA_n}_{6,a}\big)+ \big(I^{\ZZ_n}_{6,b} ,I^{\AA_n}_{6,b}\big) \,, \qquad 
\big(I^{\ZZ_n}_9,I^{\AA_n}_9\big) 
=  \big(I^{\ZZ_n}_{9,a} ,I^{\AA_n}_{9,a}\big)+ \big(I^{\ZZ_n}_{9,b} ,I^{\AA_n}_{9,b}\big) + \big(I^{\ZZ_n}_{9,c} ,I^{\AA_n}_{9,c}\big)  \,, 
\end{align} 
where
\begin{align} 
\big(I^{\ZZ_n}_{6,a} ,I^{\AA_n}_{6,a}\big)
&= {\tfrac{ \alpha(-2\beta+1) }{2}} \tints   \tfrac{1}{\Sigma^{2\beta}}  \JJss 
 \big(\Jg\Wbn - \Jg\Zbn  +  \Jg \nbs_2 h (\Wbt -\Zbt) \big) \big( \sabs{\nbs^6 (\Jg\Zbn) }^2, \sabs{\nbs^6 (\Jg\Abn) }^2\big)\,, 
\\
\big(I^{\ZZ_n}_{6,b} ,I^{\AA_n}_{6,b}\big)
&= {\tfrac{ 3\alpha }{2\eps}} \tints  \jb \JJh \nbs_1\JJ  \big( \sabs{\nbs^6( \Jg\Zbn) }^2, \sabs{\nbs^6 (\Jg\Abn )}^2\big)  \,,
\label{I-Zbn-Abn-6b}\\
  \big(I^{\ZZ_n}_{9,a} ,I^{\AA_n}_{9,a}\big)
 &=  - \alpha \dint  \tfrac{\Qb}{\Sigma^{2\beta}} \Sigma \JJss\Jg g^{- {\frac{1}{2}} } \nbs_2h \big( \sabs{\nbs^6( \Jg\Zbn) }^2, \sabs{\nbs^6 (\Jg\Abn )}^2\big)\Big|^\s_\sin \,, 
 \\
 \big(I^{\ZZ_n}_{9,b} ,I^{\AA_n}_{9,b}\big)
 &=  -\tfrac{3 \alpha}{2} \tints  \jb  \JJh  \nbs_2\JJ \Jg g^{- {\frac{1}{2}} } \nbs_2h  \big( \sabs{\nbs^6 (\Jg\Zbn) }^2, \sabs{\nbs^6 (\Jg\Abn )}^2\big) \,, 
 \label{I-Zbn-Abn-9b}\\
 \big(I^{\ZZ_n}_{9,c} ,I^{\AA_n}_{9,d}\big)
 &=  - \alpha \tints   \JJss \nbs_2(\jb\Jg g^{- {\frac{1}{2}} } \nbs_2h )\big( \sabs{\nbs^6 (\Jg\Zbn )}^2, \sabs{\nbs^6 (\Jg\Abn) }^2\big) \,, 
\end{align} 
\end{subequations} 
where we have used the identity \eqref{p1-Sigma-s-P-US} in the integral $ \big(I^{\ZZ_n}_{6,a} ,I^{\AA_n}_{6,a}\big)$ above.

Recalling the integral notation given in \eqref{int-Hdmp-Hdmm}, 
we analyze the integrals in \eqref{Zbn-Zbn-6and9} using the additive integral decomposition
$$
 \tints  f = \int_{\Hdmp}\!\! \!\!\!\!f (x,\s') {\rm d}x{\rm d}\s'+   \int_{\Hdmm} \!\! \!\!\!\!f (x,\s') {\rm d}x{\rm d}\s'
$$
and we shall write $\big(I^{\ZZ_n}_1, I^{\ZZ_n}_6,I^{\AA_n}_6\big) = \big(I^{\ZZ_n}_{1,+},I^{\ZZ_n}_{6,+},I^{\AA_n}_{6,+}\big) 
+ \big(I^{\ZZ_n}_{1,-}, I^{\ZZ_n}_{6,-},I^{\AA_n}_{6,-}\big)$
to indicate this additive decomposition.
For integrals over $\Hdmp$, 
the important observation is that by summing the integrals  \eqref{I-Wbn-Zbn-Abn-1c}, \eqref{I-Zbn-Abn-6b}, and 
\eqref{I-Zbn-Abn-9b}, and applying the basic identity \eqref{JJ-formula}, we obtain that
\begin{subequations} 
\label{nicechickennice}
\begin{align} 
I^{\ZZ_n}_{1,c,+}+I^{\ZZ_n}_{6,b,+}+ I^{\ZZ_n}_{9,b,+} & =
-\tfrac{3}{4} \dl   \int_{\Hdmp} \tfrac{1}{\Sigma^{2\beta}} \JJh \Jg
(\Q\p_\s+V\p_2) \JJ
\sabs{\nbs^6 (\Jg\Zbn )}^2 \,, \\
I^{\AA_n}_{1,c,+}+I^{\AA_n}_{6,b,+}+ I^{\AA_n}_{9,b,+} & 
=
-\tfrac{3}{4}(1+ \dl ) \int_{\Hdmp} \tfrac{1}{\Sigma^{2\beta}} \JJh \Jg(\Q\p_\s+V\p_2) \JJ \sabs{\nbs^6 (\Jg\Abn )}^2
 \,.
\end{align} 
\end{subequations} 

For integrals over $\Hdmm$, we use the lower bound \eqref{qps-JJ-bound} together with the upper-bounds \eqref{fat-marmot2b}--\eqref{fat-marmot2}, and the
fact that $\JJ \ge 1$ and $\Jg \ge \tfrac{1}{9} $ in $\Hdmm$.   We find that
\begin{subequations} 
\label{nicechickennice-minus}
\begin{align} 
I^{\ZZ_n}_{1,c,-}+I^{\ZZ_n}_{6,b,-}+ I^{\ZZ_n}_{9,b,-} & \ge
\bubu{+}\tfrac{33(1+ \alpha )}{100 \eps}   \int_{\Hdmm} \tfrac{1}{\Sigma^{2\beta}} \JJh \Jg\sabs{\nbs^6 (\Jg\Zbn )}^2 
- \tfrac{\Cn}{\eps}  \snorm{\tfrac{\Q^{\frac{1}{2}} \JJtf \Jgh}{\Sigma^\beta} \nbs^6(\Jg \Zbn)(\cdot,\s)}_{L^2_{x,\s}}^2  , \\
I^{\AA_n}_{1,c,-}+I^{\AA_n}_{6,b,-}+ I^{\AA_n}_{9,b,-} &  \ge
\bubu{+}\tfrac{33(1+ \alpha )}{50 \eps}   \int_{\Hdmm} \tfrac{1}{\Sigma^{2\beta}} \JJh \Jg\sabs{\nbs^6 (\Jg\Abn )}^2 
- \tfrac{\Cn}{\eps}  \snorm{\tfrac{\Q^{\frac{1}{2}} \JJtf \Jgh}{\Sigma^\beta} \nbs^6(\Jg \Abn)(\cdot,\s)}_{L^2_{x,\s}}^2  .
\end{align} 
\end{subequations}

Therefore, using \eqref{nicechickennice} and \eqref{nicechickennice-minus}, we have that
\begin{align} 
& I^{\WW_n}_1 +I^{\ZZ_n}_1 + I^{\AA_n}_1 
\bubu{ +  I^{\WW_n}_2 + I^{\ZZ_n}_2 + I^{\AA_n}_2  }
+  I^{\ZZ_n}_6 + I^{\AA_n}_6  
+  I^{\ZZ_n}_9 + I^{\AA_n}_9  
\notag \\
& \qquad
\bubu{ \ge }
 \tfrac{1}{2}   \dint \tfrac{\Q}{\Sigma^{2\beta}}   \JJss\Jg
\big( \sabs{\nbs^6( \Jg\Wbn  ) }^2 + \sabs{\nbs^6( \Jg\Zbn  ) }^2 + 2\sabs{\nbs^6( \Jg\Abn  ) }^2\big)\Big|^\s_\sin
\notag \\
&\qquad
 -\tfrac{1}{2} \tints \tfrac{1}{\Sigma^{2\beta}} \JJss 
(\Q\p_\s+V\p_2) \Jg
\big(\sabs{\nbs^6( \Jg\Wbn) }^2 + \sabs{\nbs^6 (\Jg\Zbn )}^2 + 2\sabs{\nbs^6 (\Jg\Abn )}^2\big)   
 \notag \\
& \qquad\qquad
- \tfrac{3}{4} \tints \tfrac{1}{\Sigma^{2\beta}} \JJh \Jg
(\Q\p_\s+V\p_2) \JJ  \big( \sabs{\nbs^6 (\Jg\Wbn) }^2 +  \sabs{\nbs^6 (\Jg\Abn) }^2\big)
\notag \\
& \qquad\qquad
-\tfrac{3}{4} \dl\tints \tfrac{1}{\Sigma^{2\beta}} \JJh \Jg
(\Q\p_\s+V\p_2) \JJ
\big( \sabs{\nbs^6 (\Jg\Zbn )}^2 + \sabs{\nbs^6 (\Jg\Abn )}^2 \big)
\notag \\
& \qquad
+{\tfrac{ \alpha(-2\beta+1) }{2}} \tints   \tfrac{1}{\Sigma^{2\beta}}  \JJss 
 \big(\Jg\Wbn - \Jg\Zbn  +  \Jg \nbs_2 h (\Wbt -\Zbt) \big)   \big( \sabs{\nbs^6 (\Jg\Zbn) }^2 +  \sabs{\nbs^6( \Jg\Abn) }^2\big)
\notag \\
&\qquad 
- \tints \tfrac{1}{\Sigma^{2\beta}} \big( \alpha \beta  (\Zbn +\Abt) + \tfrac{1}{2} (\Qc +V,_2)  \big)\JJss\Jg 
\big( \sabs{\nbs^6( \Jg\Wbn) }^2 + \sabs{\nbs^6 (\Jg\Zbn) }^2 +  2\sabs{\nbs^6( \Jg\Abn) }^2\big)
\notag \\
& \qquad
\bubu{ 
 - \alpha \tints\tfrac{1}{\Sigma^{2\beta}}   \Big( \Jg\Wbn-\Jg\Zbn+ \Sigma^{2\beta}\nbs_2(\jb  \Jg g^{- {\frac{1}{2}} } \nbs_2h )\Big)   \JJss 
 \big( \sabs{\nbs^6( \Jg\Zbn )}^2 + \sabs{\nbs^6 (\Jg\Abn) }^2\big) }
 \notag \\
& \qquad
 - \tfrac{\Cn}{\eps} \tints \tfrac{\Q}{\Sigma^{2\beta}} \JJss \Jg
 \big( \sabs{\nbs^6 (\Jg\Zbn) }^2 +  \sabs{\nbs^6 (\Jg\Abn) }^2\big)
 \,.
 \label{fat-marmot10}
 \end{align} 
 
 \begin{remark} [The $\JJ$ evolution in \eqref{JJ-formula} is fundamental]\label{rem:JJ-formula}
 Let us, at this point of the proof, make an important remark concerning the identities \eqref{nicechickennice}.
With the parameter $\dl$  chosen in \eqref{def-dl}, the choice of the {\it approximate} slow acoustic characteristic surface $\Thd$ is  
 made in  \eqref{eq:Thd:PDE} so that the weight function $\JJ$ satisfies the evolution equation \eqref{JJ-formula}, and thus yields the damping norm 
 obtain in \eqref{nicechickennice}.   In particular, the formula \eqref{JJ-formula} replaces the $x_1$-independence of the weight $\mathcal{J} $ in 
Section \ref{sec:sixth:order:energy} and replaces the crucial fact that  $\nbs_1\Jg \ge 0$ in a large region of the downstream spacetime in Section 
\ref{sec:downstreammaxdev}.
 \end{remark} 
 
\subsubsection{Upstream analysis of $I^{\WW_n}_3 +I^{\ZZ_n}_3 +I^{\AA_n}_3$}
Integrating-by-parts with respect to $\nbs_2$ and using \eqref{adjoint-2-H}, we find  
\begin{align} 
I^{\WW_n}_3 +I^{\ZZ_n}_3 +I^{\AA_n}_3
 &=  2\alpha \tints \jb g^{- {\frac{1}{2}} }   \ \JJss \Jg \,  \nbs_2\big(\nbs^6(\Jg\Abn)  \  \nbs^6(\Jg\Sbn) \big)
  \notag \\
  &=: I^{\WW_n+\ZZ_n+\AA_n}_{3,a}+ I^{\WW_n+\ZZ_n+\AA_n}_{3,b}+ I^{\WW_n+\ZZ_n+\AA_n}_{3,c}+ I^{\WW_n+\ZZ_n+\AA_n}_{3,d} \,,
  \notag \\
 I^{\WW_n+\ZZ_n+\AA_n}_{3,a}& = \tfrac{3 \alpha}{2} \tints  \jb  g^{- {\frac{1}{2}} } \nbs_2\JJ  \,  \JJh\Jg \, \nbs^6(\Jg\Abn)  
 \nbs^6(\Jg\Zbn)  \,,
\notag \\
 I^{\WW_n+\ZZ_n+\AA_n}_{3,b}& = -\tfrac{3 \alpha}{2} \tints  \jb  g^{- {\frac{1}{2}} } \nbs_2\JJ  \,  \JJh\Jg \, \nbs^6(\Jg\Abn)  
 \nbs^6(\Jg\Wbn)   \,,
\notag \\
 I^{\WW_n+\ZZ_n+\AA_n}_{3,c}&   =
- \alpha \tints   \jb   g^{- {\frac{1}{2}} }   \nbs_2 \Jg \, \JJss\,   \nbs^6(\Jg\Abn)  
\big( \nbs^6(\Jg\Wbn)- \nbs^6(\Jg\Zbn)\big)  \,,
\notag \\
 I^{\WW_n+\ZZ_n+\AA_n}_{3,d}
& =
- 2\alpha \dint \Qb\,   \jb \JJss \Jg g^{- {\frac{1}{2}} }  \nbs^6(\Jg\Abn) \big( \nbs^6(\Jg\Wbn)- \nbs^6(\Jg\Zbn)\big)\Big|^\s_\sin
 \,.
  \label{fat-marmot11}
\end{align} 
The novelty in estimating  the integral $ I^{\WW_n+\ZZ_n+\AA_n}_{3,a}$ is that its bound involves a constant that is proportional to $ \dl ^{-1} $.   
In particular, by an 
application of the Cauchy-Young inequality, \eqref{bs-h}, \eqref{bs-Sigma}, we have that
\begin{align*} 
\sabs{I^{\WW_n+\ZZ_n+\AA_n}_{3,a}} 
&\le  3 \alpha \kappa_0 \snorm{  \nbs_2 \JJ}_{L^\infty_{x,\s}}
\snorm{  \tfrac{\JJof \Jgh}{\Sigma^{\beta}} \nbs^6(\Jg\Zbn) }_{L^2_{x,\s}} \snorm{  \tfrac{\JJof \Jgh}{\Sigma^{\beta}} \nbs^6(\Jg\Abn) }_{L^2_{x,\s}}
\notag \\
& \le \tfrac{4 \dl (1+ \alpha )}{25\eps}  \snorm{  \tfrac{\JJof \Jgh}{\Sigma^{\beta}} \nbs^6(\Jg\Zbn) }_{L^2_{x,\s}}^2 +
 \tfrac{25 \eps}{16 \dl (1+ \alpha )} \cdot 9 \alpha^2\kappa_0^2 \snorm{  \nbs_2 \JJ}_{L^\infty_{x,\s}}^2
 \snorm{  \tfrac{\JJof \Jgh}{\Sigma^{\beta}} \nbs^6(\Jg\Abn) }_{L^2_{x,\s}}^2
 \notag \\
& \le \tfrac{4 \dl (1+ \alpha )}{25\eps}  \snorm{  \tfrac{\JJof \Jgh}{\Sigma^{\beta}} \nbs^6(\Jg\Zbn) }_{L^2_{x,\s}}^2 +
 \tfrac{ \eps}{\dl} \Cn  (\tfrac{4}{\kappa_0})^{2\beta}  \brak{\mathsf{B}_6}^2 \,,
\end{align*} 
where we have additionally used \eqref{bootstraps-Dnorm:6}, \eqref{fat-marmot2b}, and~\eqref{fat-marmot2} for the last inequality.
By choosing 
\begin{equation} 
\eps \le \dl \,, \label{eps-le-dl}
\end{equation} 
we have that
\begin{align*} 
\sabs{I^{\WW_n+\ZZ_n+\AA_n}_{3,a}} 
& \le \tfrac{4 \dl (1+ \alpha )}{25\eps}  \snorm{  \tfrac{\JJof \Jgh}{\Sigma^{\beta}} \nbs^6(\Jg\Zbn) }_{L^2_{x,\s}}^2 +
 \Cn  (\tfrac{4}{\kappa_0})^{2\beta}  \brak{\mathsf{B}_6}^2 \,,
\end{align*} 

The bootstrap bounds \eqref{eq:Qb:bbq-H} together with the bound for $\Qb$ in \eqref{eq:Qb:bbq-H} then show that
\begin{align*} 
\sabs{I^{\WW_n+\ZZ_n+\AA_n}_{3,b}} +\sabs{I^{\WW_n+\ZZ_n+\AA_n}_{3,c}} + \sabs{I^{\WW_n+\ZZ_n+\AA_n}_{3,d}} 
\les  (\tfrac{4}{\kappa_0})^{2\beta}  \brak{\mathsf{B}_6}^2 \,.
\end{align*} 

\subsubsection{Upstream lower bound for
$ I^{\WW_n}_1 + I^{\ZZ_n}_1  + I^{\AA_n}_1 + I^{\ZZ_n}_2 + I^{\AA_n}_2
 + I^{\WW_n}_3 + I^{\ZZ_n}_3 + I^{\AA_n}_3
 + I^{\ZZ_n}_6  + I^{\AA_n}_6 
 + I^{\ZZ_n}_9  + I^{\AA_n}_9 $}
We now obtain a lower bound for the sum of \eqref{fat-marmot10} and \eqref{fat-marmot11}.  
Making use of   \eqref{bootstraps-H}, \eqref{eq:Q:all:bbq-H},  \eqref{eq:waitin:for:the:bus-s}, and taking $\eps$ sufficiently small, 
we have that
\begin{align}
& I^{\WW_n}_1 + I^{\ZZ_n}_1  + I^{\AA_n}_1 + I^{\ZZ_n}_2 + I^{\AA_n}_2
 + I^{\WW_n}_3 + I^{\ZZ_n}_3 + I^{\AA_n}_3
 + I^{\ZZ_n}_6  + I^{\AA_n}_6 
 + I^{\ZZ_n}_9  + I^{\AA_n}_9 
 \notag\\
&\qquad
\geq \bigl( \tfrac 12 - \eps^{\frac{7}{4}} \bigr) 
 \snorm{\tfrac{\Q^{\frac{1}{2}} \JJtf \Jgh}{\Sigma^\beta} \nbs^6(\Jg \Wbn,\Jg \Zbn,\Jg \Abn)(\cdot,\s)}_{L^2_x}^2 
-  \tfrac{1}{2} \snorm{ \tfrac{1}{\Sigma^\beta} \nbs^6 (\Jg \Wbn,\Jg \Zbn,\Jg \Abn)(\cdot,\sin)}_{L^2_x}^2 
\notag\\
&\qquad\qquad
 +\tints\tfrac{1}{\Sigma^{2\beta}}   \mathsf{G_{good}}
\sabs{\nbs^6(\Jg \Wbn)}^2    
+  \tfrac{17(1+\alpha)  }{100 \eps} \tints \tfrac{1}{\Sigma^{2\beta}} 
 \JJh \Jg \big( \dl \sabs{\nbs^6(\Jg \Zbn)}^2 + (1+\dl) \sabs{\nbs^6(\Jg \Abn)}^2  \big)
\notag \\
& \qquad\qquad
+\tints   \tfrac{1}{\Sigma^{2\beta}} 
 \big( ({\tfrac{ 2\alpha(2\beta-1) + \bubu{  4 \alpha } }{5\eps}}  - \tfrac{13(1+3\alpha-4\alpha\beta)}{4\eps}  \Jg \big) \JJss   \big( \sabs{\nbs^6 (\Jg\Zbn) }^2 +  \sabs{\nbs^6( \Jg\Abn) }^2\big)
\notag
 \\
 & \qquad\qquad
 - \tfrac{\Cn}{\eps} \tints \tfrac{\Q}{\Sigma^{2\beta}} \JJss \Jg
 \big( \sabs{\nbs^6 (\Jg\Zbn) }^2 +  \sabs{\nbs^6 (\Jg\Abn) }^2\big)  -  \Cn  (\tfrac{4}{\kappa_0})^{2\beta}  \brak{\mathsf{B}_6}^2 \,,
 \label{eq:I:n:12369-H}
\end{align}
where
\begin{equation} 
 \mathsf{G_{good}} = - \tfrac{1}{2} \JJss (\Q\p_\s+V\p_2) \Jg - \tfrac{3}{4} \JJh \Jg (\Q\p_\s+V\p_2) \JJ - \Cn \,. \label{Ggood-H}
\end{equation}

\subsubsection{Geometric identities in the upstream spacetime}
We will need the upstream variants of Lemmas \ref{lem:tau-Jg-D2} and \ref{lem:tau-Jg-D1} to analyze the terms with
over-differentiated geometry.
\begin{lemma} \label{lem:tau-Jg-D2-H}
We have that
\begin{align} 
\left| \tints f \, \nn \cdo \nbs^6\tt  \ \nbs_2 \nbs^6 \Jg \right| 
+
\left| \tints f \, \nn \cdo \nbs_2 \nbs^6\tt  \  \nbs^6 \Jg \right| 
&\les  \eps^{3} \mathsf{K}^2 \brak{\mathsf{B}_6}^2 
\Bigl(\snorm{\JJh  \nbs f}_{L^\infty_{x,\s}}  
+\eps \snorm{\JJmh f}_{L^\infty_{x,\s}}  
   \Bigr) 
   \label{eq:tau-Jg-D2-H}
\end{align} 
for all differentiable functions $f(x,\s)$ such that the right side of \eqref{eq:tau-Jg-D2-H}.  
\end{lemma} 
\begin{proof} [Proof of Lemma \ref{lem:tau-Jg-D2-H}]
We begin with the modifications to the first integral on the left side of \eqref{eq:tau-Jg-D2-H}.
Substituting  the identity \eqref{big-chicken1} into the integral $\int_\eps^\s \! \iint^\theta \!\! f \, \nn \cdo \nbs^6\tt  \ \nbs_2 \nbs^6 \Jg$, and 
integrating-by-parts using \eqref{eq:adjoints-H}, we find that
\begin{align}
&\tints f \, \nn \cdo \nbs^6\tt  \, \nbs_2 \nbs^6 \Jg
\notag\\
& \qquad
=
- \tfrac{1}{2\eps} \tints  \nbs_1 (f \, g^{\frac 12})   \, (\nn \cdo \nbs^6 \tt)^2 
+ \tfrac 12 \tints \nbs_2 (f \, \Jg \nbs_2 h) \, (\nn \cdo \nbs^6 \tt)^2
+ \tfrac 12 \dint \Qb f \, \Jg \nbs_2 h   (\nn \cdo \nbs^6 \tt)^2 \Bigr|^{\s}_\sin
\notag\\
&\qquad
- \tfrac{1}{\eps} \tints f \, \nn \cdo \nbs^6\tt  \Big( g^{\frac 12} \tt \cdo \nbs_1  \nn\; \tt \cdo \nbs^6 \tt
-\nn \cdo \nbs^6(g^{\frac 12} \nn)  \nn \cdo \nbs_1 \tt\Big)
+ \tfrac{1}{\eps} \tint f \, \nn \cdo \nbs^6\tt  \,\doublecom{\nbs^6, g^{\frac 12} \nn_i, \nbs_1 \tt_i}
\notag\\
&\qquad 
+ \tints f \, \nn \cdo \nbs^6\tt  \Big( \Jg \nbs_2 h \tt \cdo \nbs_2 \nn \; \tt \cdo \nbs^6 \tt
- \nn \cdo \nbs^6(\Jg \nbs_2 h \nn)  \nn \cdo \nbs_2 \tt -  \doublecom{\nbs^6, \Jg \nbs_2 h \nn_i, \nbs_2 \tt_i}\Big)
\,.
\label{big-chicken2-H}
\end{align}
Here, we have used the fact that $f$ must contain at least a power of $\JJh$ in order for the right side of 
\eqref{eq:tau-Jg-D2-H} to be bounded, which in turn implied that $f=0$ on the surface $x_1=\thsd$.
Upon comparison,  every integral on the right side of \eqref{big-chicken2-H} is directly analogous to every integral on   the right side of 
\eqref{big-chicken2} (with $\JJ$ replacing $\mathcal{J} $ and with the upstream modification of the domains of integration).  Hence, 
the identical bounds hold for \eqref{big-chicken2-H} as for \eqref{big-chicken2}.

To bound  the second integral on the left side of \eqref{eq:tau-Jg-D2-H}, we
use the adjoint formula for $\nbs_2^*$ given in  
 \eqref{adjoint-2-H}, and obtain that
\begin{align*}
 \tint f \, \nn \cdo \nbs_2 \nbs^6\tt  \  \nbs^6 \Jg
 &= - \tint f \, \nn \cdo  \nbs^6\tt  \ \nbs_2 \nbs^6 \Jg
 - \tint \bigl(\nbs_2 f \, \nn \cdo  \nbs^6\tt + f \tt \cdot \nbs_2 \nn \tt \cdo \nbs^6 \tt \bigr)  \nbs^6 \Jg
 \notag\\
 &\qquad 
 - \dint \Qb f \, \nn \cdo  \nbs^6\tt  \  \nbs^6 \Jg \Bigr|^{\s}_\sin
 \,.
\end{align*}
Using the bootstrap inequalities~\eqref{bootstraps}, the estimates~\eqref{eq:Q:all:bbq}, and the bounds for the geometry~\eqref{geometry-bounds-new}, we deduce
\begin{align*}
 \left|\tint f \, \nn \cdo \nbs_2 \nbs^6\tt  \  \nbs^6 \Jg\right|
 &\leq \left|\tint f \, \nn \cdo  \nbs^6\tt  \ \nbs_2 \nbs^6 \Jg\right|
 + \Cn \snorm{\JJh \nbs_2 f}_{L^\infty_{x,\s}}  \eps^3 \mathsf{K}   \brak{\mathsf{B}_6}^2  
 + \Cn \snorm{\JJmh  f}_{L^\infty_{x,\s}}  \eps^3 \mathsf{K}   \brak{\mathsf{B}_6}^2 
 \,.
\end{align*}
This concludes the proof of the lemma.
\end{proof} 

\begin{lemma} \label{lem:tau-Jg-D1-H}
We have that
\begin{equation} 
\left|\tint f \, \nn \cdo \nbs_1 \nbs^6\tt   \  \nbs^6 \Jg\right|
+
\left|\tint f \, \nn \cdo \nbs^6\tt  \ \nbs_1 \nbs^6 \Jg\right|
\les
\eps^3 \mathsf{K} \brak{\mathsf{B}_6}^2
\Bigl( \|\JJh \nbs f\|_{L^\infty_{x,\s}} +  \|\JJmh f\|_{L^\infty_{x,\s}} \Bigr) 
\,,
\label{eq:tau-Jg-D1-H}
\end{equation} 
for all differentiable functions $f(x,\s)$ such that the right side of \eqref{eq:tau-Jg-D1-H}.  
\end{lemma} 
The proof of  Lemma \ref{lem:tau-Jg-D1-H} is identical to the proof of Lemma  \ref{lem:tau-Jg-D1}.

\subsubsection{Upstream bounds for the forcing, remainder, and commutator functions}
\label{subsub:FRC-H}
Following the analysis in Section \ref{subsub:FRC}, we have the following upstream bounds for the forcing, remainder, and commutator functions:
\begin{subequations}
\label{eq:I:W:nn:6:all-H}
\begin{align} 
\snorm{\tfrac{\JJtf(\Jg \Q)^{\frac 12}}{\Sigma^\beta} \nbs^6\Fwn}_{L^2_{x,\s}}
& \les  (\tfrac{4}{\kappa_0})^\beta  \brak{\mathsf{B}_6}
\,,
\label{eq:I:W:nn:6:b-H}
\\
\snorm{\tfrac{\JJtf(\Jg \Q)^{\frac 12}}{\Sigma^\beta}\mathcal{R}_\Wb^\nn}_{L^2_{x,\s}}
&\les  (\tfrac{4}{\kappa_0})^\beta  \brak{\mathsf{B}_6}
\label{eq:I:W:nn:6:c-H}
\\
\snorm{\tfrac{\JJtf (\Jg \Q)^{\frac 12}}{\Sigma^\beta}\mathcal{C}_\Wb^\nn}_{L^2_{x,\s}}
&\les    (\tfrac{4}{\kappa_0})^\beta  \brak{\mathsf{B}_6}
\,.
\label{eq:I:W:nn:6:d-H}
\end{align}
\end{subequations}
and 
\begin{subequations}
\label{eq:I:Z:nn:10:all-H}
\begin{align}
\snorm{\tfrac{\JJtf }{\Sigma^{\beta-1}} \nbs^6\Fzn}_{L^2_{x,\s}}
&\leq \Cn  (\tfrac{4}{\kappa_0})^\beta  \mathsf{K} \brak{\mathsf{B}_6}
\,,  \label{eq:I:Z:nn:10:b-H} \\
\snorm{\tfrac{\JJtf}{\Sigma^{\beta-1}}\mathcal{R}_\Zb^\nn}_{L^2_{x,\s}}
&\leq \Cn    (\tfrac{4}{\kappa_0})^\beta  \mathsf{K} \brak{\mathsf{B}_6}
+   \tfrac{4 (1+ \alpha )}{\eps} \snorm{ \tfrac{1}{\Sigma^{\beta}} \JJtf \nbs^6 (\Jg \Zbn)}_{L^2_{x,\s}} \,,
\label{eq:I:Z:nn:10:c-H} \\
\snorm{\tfrac{\JJtf }{\Sigma^{\beta-1}}\mathcal{C}_\Zb^\nn}_{L^2_{x,\s}}
&\leq \Cn    (\tfrac{4}{\kappa_0})^\beta  \mathsf{K} \brak{\mathsf{B}_6}
\,,
\label{eq:I:Z:nn:10:d-H}
\end{align}
\end{subequations}
and 
\begin{subequations}
\label{eq:I:A:nn:10:all-H}
\begin{align}
\snorm{\tfrac{\JJtf }{\Sigma^{\beta-1}} \nbs^6\Fan}_{L^2_{x,\s}}
&\leq \Cn  (\tfrac{4}{\kappa_0})^\beta  \mathsf{K} \brak{\mathsf{B}_6}
\,,  \label{eq:I:A:nn:10:b-H} \\
\snorm{\tfrac{\JJtf}{\Sigma^{\beta-1}}\mathcal{R}_\Ab^\nn}_{L^2_{x,\s}}
&\leq \Cn    (\tfrac{4}{\kappa_0})^\beta  \mathsf{K} \brak{\mathsf{B}_6}
+   \tfrac{4 (1+ \alpha )}{\eps} \snorm{ \tfrac{1}{\Sigma^{\beta}} \JJtf \nbs^6 (\Jg \Abn)}_{L^2_{x,\s}} \,,
\label{eq:I:A:nn:10:c-H} \\
\snorm{\tfrac{\JJtf }{\Sigma^{\beta-1}}\mathcal{C}_\Ab^\nn}_{L^2_{x,\s}}
&\leq \Cn    (\tfrac{4}{\kappa_0})^\beta  \mathsf{K} \brak{\mathsf{B}_6}
\,.
\label{eq:I:A:nn:10:d-H}
\end{align}
\end{subequations}
The proof of these inequalities relies on 
the bootstrap bounds  \eqref{bootstraps-H}, the bounds on geometry \eqref{geom-H-original}, 
the improved estimates \eqref{eq:Jg:Abn:D5:improve-H},  \eqref{eq:Jg:Zbn:D5:improve-H}, \eqref{eq:D5:JgWbn-H},  
as well as   \eqref{D2-Jg-Linfty}, \eqref{se3:time},  and \eqref{eq:Lynch:1}; with these bounds in hand, the inequalities
\eqref{eq:I:W:nn:6:all-H}--\eqref{eq:I:A:nn:10:all-H} follow in the identical manner as proven in Section \ref{subsub:FRC}.

\subsubsection{The upstream analysis of $I^{\WW_n}_5+I^{\ZZ_n}_5+I^{\AA_n}_7$}\label{sec::Wbn5+Zbn5+Abn7}
We first note that since $\Sbn= \tfrac{1}{2} (\Wbn-\Zbn)$, we have that
\begin{align} 
I^{\WW_n}_5+I^{\ZZ_n}_5 &=
-\alpha  \tint \jb \JJss  \Jg g^{- {\frac{1}{2}} } (\Jg\Wbn + \Jg\Zbn - 2\Jg\Abt)\nbs_2 \nbs^6 \tt\cdo\nn \ \nbs^6(\Jg\Sbn) \,.
\label{IWn5+IZn5-H}
\end{align} 

For the integral $I^{\AA_n}_7$ in \eqref{I7-Abn-H}, we use the adjoint formulas \eqref{eq:adjoints-H} and  \eqref{calc-dsd2} to find that
\begin{align*} 
I^{\AA_n}_7& = I^{\AA_n}_{7,a} + I^{\AA_n}_{7,b} + I^{\AA_n}_{7,c} +  I^{\AA_n}_{7,d } +  I^{\AA_n}_{7,e}+  I^{\AA_n}_{7,f} \,, 
\notag \\
I^{\AA_n}_{7,a} &=
-{\tfrac{\alpha }{\eps}} \tint  \jb  \JJss (\Jg\Wbn+\Jg\Zbn - 2\Jg\Abt) \nbs^6\tt \cdo\nn  
\ \Big( \nbs_1\nbs^6(\Jg\Abn) -  \eps\Jg g^{- {\frac{1}{2}} } \nbs_2h\, \nbs_2\nbs^6(\Jg\Abn) \Big) \,, 
\notag \\
I^{\AA_n}_{7,b} &=
-{\tfrac{\alpha }{\eps}} \tint  \nbs_1  \Big( \jb  \JJss (\Jg\Wbn+\Jg\Zbn - 2\Jg\Abt) \nn_i \Big) \nbs^6\tt_i  \ \nbs^6(\Jg\Abn) \,, 
\notag \\
I^{\AA_n}_{7,c} &=
 \alpha \tint  \nbs_2 \Big( \jb   g^{- {\frac{1}{2}} } \nbs_2 h \,  \JJss \Jg  (\Jg\Wbn+\Jg\Zbn - 2\Jg\Abt) \nbs^6\tt \cdo\nn  \ \nbs^6(\Jg\Abn) \Big) \,, 
\notag \\
I^{\AA_n}_{7,d} &=
\alpha  \dint  \Qb \   g^{- {\frac{1}{2}} } \nbs_2 h \,  \jb\JJss\Jg  (\Jg\Wbn+\Jg\Zbn - 2\Jg\Abt) \nbs^6\tt \cdo\nn  \ \nbs^6(\Jg\Abn)\Big|^\s_\sin\,.
\end{align*} 
Comparing the upstream decomposition of the integral $I^{\AA_n}_7$ with the decomposition \eqref{eq:IAAn:7:decompose} of the
corresponding integral  in Section \ref{sec:sixth:order:energy}, we
see that the integrals are the same once the weight $\Fg$ is replaced by the weight $\JJ$. 
Given the bounds in Section \ref{subsub:FRC-H}, 
 the estimates required here are identical to those provided
in Section \ref{sec:sixth:order:energy} so we will not repeat those detailed computations.  Just as in \eqref{I7-An-bcd}, we find that
$\sabs{I^{\AA_n}_{7,b} } +\sabs{I^{\AA_n}_{7,c} } + \sabs{I^{\AA_n}_{7,d} } 
\les  (\tfrac{4}{\kappa_0})^{2\beta}  \mathsf{K} \Bsix^2$, while the integral $I^{\AA_n}_{7,a}$ is given the corresponding decomposition as in 
\eqref{IAbn-7a-i-viii}.   We again find that $I^{\WW_n}_5+I^{\ZZ_n}_5+I^{\AA_n}_{7,a,i} $ cancel the derivative loss and that 
\begin{equation} 
\sabs{I^{\WW_n}_5+I^{\ZZ_n}_5+I^{\AA_n}_{7} }\les (\tfrac{4}{\kappa_0})^{2\beta}  \mathsf{K}^2 \Bsix^2  \,.
\label{IWn5+IZn5+IAn7-H}
\end{equation}

\subsubsection{The upstream analysis of $I^{\WW_n}_4+I^{\ZZ_n}_4+I^{\AA_n}_8$}
The estimates for the sum of these three integrals follows identically the analysis of Section \ref{sec:IW4+IZ4+IA8} and we do not repeat the details
here.  We simply record that 
\begin{equation}
\sabs{I^{\WW_n}_4+I^{\ZZ_n}_4+I^{\AA_n}_8} 
\les  
(\tfrac{4}{\kappa_0})^{2\beta}  \brak{\mathsf{B}_6}^2
\,.
\label{eq:small:duck:0-H}
\end{equation}

\subsubsection{The upstream analysis of $I^{\AA_n}_4$}\label{sec::Abn4}
For the integral $ I^{\AA_n}_4$ defined in~\eqref{I4-Abn-H} we first integrate-by-parts  with respect to the $\nbs_2$ derivative.  Using \eqref{adjoint-2-H}, we have that
\begin{align*} 
 I^{\AA_n}_{4}&=  I^{\AA_n}_{4,a} +  I^{\AA_n}_{4,b} +  I^{\AA_n}_{4,c}  \,,
 \notag \\
I^{\AA_n}_{4,a} &= 
2\alpha \tint \jb  \JJss ( \Jg \Sbn) \nbs^6 \Jg \  g^{- {\frac{1}{2}} } \nbs^6\nbs_2(\Jg\Abn) \,,
\notag \\
I^{\AA_n}_{4,b} &= 
2\alpha \tint \nbs_2 ( \jb \JJss g^{- {\frac{1}{2}} }  \Jg \Sbn ) \nbs^6 \Jg \   \nbs^6(\Jg\Abn) \,,
\notag \\
I^{\AA_n}_{4,c} &= 
2\alpha   \int  \Qb  \jb  g^{- {\frac{1}{2}} }  \JJss \Jg \Sbn  \nbs^6 \Jg \   \nbs^6(\Jg\Abn) \Big|^\s_\sin \,.
\end{align*} 
We bound  $I^{\AA_n}_{4,b}$ and $I^{\AA_n}_{4,c}$ in a straightforward manner by using \eqref{bootstraps-H}, \eqref{eq:Q:all:bbq-H},   \eqref{geom-H-original}, to obtain
\begin{align}
\sabs{I^{\AA_n}_{4,b}} + \sabs{I^{\AA_n}_{4,c}}
\les 
(\tfrac{4}{\kappa_0})^{2\beta}  \brak{\mathsf{B}_6}^2
\,.
\label{eq:small:duck:7-H} 
\end{align}
For the integral $I^{\AA_n}_{4,a} $, we use equation \eqref{energy-Wn-s} to substitute for  $\alpha g^{- {\frac{1}{2}} } \nbs^6\nbs_2(\Jg\Abn)$.   This leads to an
additive decomposition for $I^{\AA_n}_{4,a} $ which we write as
\begin{subequations} 
\label{I4-Abn-a-H}
\begin{align} 
I^{\AA_n}_{4,a} & = I^{\AA_n}_{4,a,i} + I^{\AA_n}_{4,a,ii} + I^{\AA_n}_{4,a,iii} + I^{\AA_n}_{4,a,iv} \,,
\notag \\
 I^{\AA_n}_{4,a,i}&= - 2\tint \tfrac{1}{\Sigma^{2\beta}}  \JJss  ( \Jg \Sbn) \nbs^6 \Jg \ (\Q\p_\s +V\p_2) \nbs^6(\Jg\Wbn  ) \,,
\label{I4-Abn-a,i-H} \\
 I^{\AA_n}_{4,a,ii}&= 
 \alpha \tint \jb\JJss    (\Jg \Sbn)  \  g^{- {\frac{1}{2}} } \Abn \nbs_2 \Bigl( \bigl(\nbs^6 \Jg\bigr)^2 \Bigr) \,,
\label{I4-Abn-a,ii-H} \\
 I^{\AA_n}_{4,a,iii}&= 
\alpha \tint \jb  \JJss  (\Jg \Sbn) \  g^{- {\frac{1}{2}} } (\Jg\Wbn + \Jg\Zbn - 2\Jg\Abt)\nbs^6 \Jg  \ \nbs_2 \nbs^6 \tt\cdo\nn \,,
\label{I4-Abn-a,iii-H} \\
 I^{\AA_n}_{4,a,iv}&= 
2 \tint \jb  \JJss ( \Jg \Sbn) \nbs^6 \Jg \   ( \nbs^6\Fwn + \mathcal{R}_\Wb^\nn + \mathcal{C}_\Wb^\nn) \,.
\label{I4-Abn-a,iv-H}
\end{align} 
\end{subequations} 
The last three terms in \eqref{I4-Abn-a-H} may be estimated directly: from 
\eqref{bootstraps-H}, \eqref{eq:Q:all:bbq-H}, \eqref{geom-H-original}, \eqref{eq:I:W:nn:6:all-H}, we obtain
\begin{subequations}
\label{eq:small:duck:8-H}
\begin{equation}
\sabs{I^{\AA_n}_{4,a,iv}}
\les    (\tfrac{4}{\kappa_0})^{2\beta}  \brak{\mathsf{B}_6}^2
\,,
\end{equation}
by additionally integrating by parts the $\nbs_2$ term via \eqref{adjoint-2-H} we obtain that
\begin{equation}
\sabs{I^{\AA_n}_{4,a,ii}}
\les \eps   (\tfrac{4}{\kappa_0})^{2\beta}  \brak{\mathsf{B}_6}^2
\,,
\end{equation}
and by also appealing to Lemma~\ref{lem:tau-Jg-D2-H}  we have that
\begin{equation}
\sabs{I^{\AA_n}_{4,a,iii}}
\les \eps  \mathsf{K} (\tfrac{4}{\kappa_0})^{2\beta}  \brak{\mathsf{B}_6}^2
\,.
\end{equation}
\end{subequations}

We now focus on the integral $ I^{\AA_n}_{4,a,i}$ defined in \eqref{I4-Abn-a,i-H} which produces an anti-damping term.  By once again using integration-by-parts via \eqref{adjoint-3-H}, we see that, 
\begin{subequations}
\label{I4-Abn-a-i-H}
\begin{align} 
 I^{\AA_n}_{4,a,i} &= J^{\AA_n}_{1} + J^{\AA_n}_{2}+ J^{\AA_n}_{3}+ J^{\AA_n}_{4}+ J^{\AA_n}_{5}+ J^{\AA_n}_{6}+ J^{\AA_n}_{7}
 + J^{\AA_n}_{8} \,
\notag \\
J^{\AA_n}_{1} &  = 2 \tints \tfrac{1}{\Sigma^{2\beta}}  \JJss   (\Jg\Sbn) (\Q\p_\s+V\p_2) \nbs^6\Jg  \  \nbs^6(\Jg\Wbn)  \,,
  \\
J^{\AA_n}_{2} &  =2 \tints \tfrac{1}{\Sigma^{2\beta}}  \JJss  (\Q\p_\s+V\p_2) (\Jg\Sbn)  \nbs^6\Jg  \  \nbs^6(\Jg\Wbn)   \,,
  \\
J^{\AA_n}_{3} &  =  \tints \tfrac{1}{\Sigma^{2\beta}} (\Q\p_\s+V\p_2) \JJss  \  (\Jg\Wbn)  \nbs^6\Jg  \  \nbs^6(\Jg\Wbn)  \,,
\label{J-Abn-3}
  \\
J^{\AA_n}_{4} &  =  2\tints      (\Q\p_\s +V\p_2)\big(  \tfrac{1}{\Sigma^{2\beta}}\big) \ \JJss (\Jg \Sbn) \   \nbs^6 \Jg \  \nbs^6(\Jg\Wbn  ) 
  \\
J^{\AA_n}_{5} & = 2 \tints \tfrac{1}{\Sigma^{2\beta}}  \JJss   (\Jg \Sbn)(-V \Qr_2+  \nbs_2 V) \nbs^6 \Jg \ \nbs^6(\Jg\Wbn  ) \,,
  \\
J^{\AA_n}_{6} & =
-  \dint \tfrac{\Q}{\Sigma^{2\beta}} \JJss    \Jg \Wbn \nbs^6 \Jg \  \nbs^6(\Jg\Wbn  ) \Big|_\s
 \,,
  \\
J^{\AA_n}_{7} & =
\dint \tfrac{\Q}{\Sigma^{2\beta}}    \Jg \Zbn\,   \JJof \nbs^6 \Jg \  \JJff \nbs^6(\Jg\Wbn  ) \Big|_\s
  \\  
J^{\AA_n}_{8} & =
 2 \dint \tfrac{\Q}{\Sigma^{2\beta}} \JJss    \Jg \Sbn \nbs^6 \Jg \  \nbs^6(\Jg\Wbn  ) \Big|_\sin \,,
  \\    
J^{\AA_n}_{9} &  = - \tints \tfrac{1}{\Sigma^{2\beta}} (\Q\p_\s+V\p_2) \JJss  \  (\Jg\Zbn)  \nbs^6\Jg  \  \nbs^6(\Jg\Wbn)  \,.
\end{align} 
\end{subequations}
In the same way that we bounded the analogous terms in \eqref{I4-Abn-a-i}, 
most of the terms in \eqref{I4-Abn-a-i-H} are estimated using \eqref{Sigma0i-ALE-s}, 
\eqref{bootstraps-H}, \eqref{eq:Q:all:bbq-H}, \eqref{eq:waitin:for:the:bus-s},   \eqref{geom-H-original}, and \eqref{eq:Jg:Wbn:improve:material:a-H} as 
\begin{align}
&\sabs{J^{\AA_n}_{2}} +  \sabs{J^{\AA_n}_{4}}+\sabs{J^{\AA_n}_{5}}+\sabs{J^{\AA_n}_{9}}
\les 
(\tfrac{4}{\kappa_0})^{2\beta}  \brak{\mathsf{B}_6}^2
\,.
\label{eq:small:duck:8a-H} 
\end{align}
By additionally employing \eqref{JJ-le-Jg}, we have that
\begin{align}
\sabs{J^{\AA_n}_{7}} & \le 
\big( \tfrac{101}{100} \big)^{\!\!\frac{1}{2}} 
\dint \tfrac{\Q}{\Sigma^{2\beta}}    \sabs{\Jg \Zbn} \,   \sabs{ \JJof \nbs^6 \Jg} \ \sabs{ \JJtf \Jgh \nbs^6(\Jg\Wbn  ) }\Big|_\s
\les 
(\tfrac{4}{\kappa_0})^{2\beta}  \brak{\mathsf{B}_6}^2
\,.
\label{eq:small:duck:8aa-H} 
\end{align}
and by also appealing to \eqref{sigma0-bound}, \eqref{table:derivatives} and \eqref{Q-lower-upper-H}, we have that for
$\eps$ sufficiently small, 
\begin{equation}
\sabs{J^{\AA_n}_{8}}
\le 2 \left( \tfrac{1}{\eps} +\Cn\right) (1+ \eps)
(\tfrac{3}{\kappa_0})^{2\beta}  \Cdatatwo
\leq \tfrac{4}{\eps} (\tfrac{3}{\kappa_0})^{2\beta}  \Cdatatwo
\,.
\label{eq:small:duck:8b-H} 
\end{equation}

Once again, the remaining integrals $J^{\AA_n}_{1}$, $J^{\AA_n}_{3}$, and $J^{\AA_n}_{6}$ require some care in their estimation. 
The integral $J^{\AA_n}_{1}$ produces an anti-damping term that must be combined with the last integral on the right
side of \eqref{eq:heavy:fuel:1n}.   Using \eqref{good-comm}, \eqref{Jg-evo-s},  a short computation shows that
\begin{align*} 
J^{\AA_n}_{1} & = J^{\AA_n}_{1,a} + J^{\AA_n}_{1,b}+ J^{\AA_n}_{1,c}  \,, 
\\
J^{\AA_n}_{1,a} &  = \tints \tfrac{1}{\Sigma^{2\beta}}  \JJss   (\Q\p_\s+V\p_2) \Jg\    \sabs{  \nbs^6(\Jg\Wbn)  }^2
  \\
J^{\AA_n}_{1,b} &  = \tfrac{1-\alpha}{2}
 \tints \tfrac{1}{\Sigma^{2\beta}}  \JJss    (\Jg\Wbn) \nbs^6(\Jg\Zbn)  
\nbs^6(\Jg\Wbn) \,,
  \\
J^{\AA_n}_{1,c} &  = 
- \tints \tfrac{1}{\Sigma^{2\beta}}  \JJss   (\Jg \Zbn) \bigl( \nbs^6(\Jg\Wbn) 
\bigr)^2 \,,
  \\
J^{\AA_n}_{1,d} &  = 
- 2 \tints \tfrac{1}{\Sigma^{2\beta}}   \JJss  ( \Jg \Sbn)  \big(\nbs^6 V \nbs_2\Jg + \doublecom{\nbs^6,V, \nbs_2\Jg} \big)   \  \nbs^6(\Jg\Wbn  ) 
  \,.
\end{align*} 
Using \eqref{bootstraps}, \eqref{eq:Q:all:bbq}, \eqref{geometry-bounds-new}, and Lemma~\ref{lem:comm:tangent} we deduce
\begin{subequations}
\label{eq:small:duck:big:duck-H}
\begin{equation}
\sabs{J^{\AA_n}_{1,c}} + \sabs{J^{\AA_n}_{1,d}}
\les (\tfrac{4}{\kappa_0})^{2\beta}  \brak{\mathsf{B}_6}^2 
+ \mathsf{K} \eps (\tfrac{4}{\kappa_0})^{2\beta}  \brak{\mathsf{B}_6}^2
\les (\tfrac{4}{\kappa_0})^{2\beta}  \brak{\mathsf{B}_6}^2
\label{eq:small:duck:9-H} 
\,.
\end{equation}
Moreover, by the Cauchy-Schwartz and Cauchy-Young inequalities, \eqref{bs-JgnnWb}, and \eqref{JJ-le-Jg},  we obtain that
\begin{align} 
\sabs{J^{\AA_n}_{1,b}}
&\le \tfrac{1+\alpha}{\eps} \sqrt{ \tfrac{101}{100} }
\int_0^{\s} \snorm{\tfrac{\JJof \Jgh}{\Sigma^\beta} \nbs^6(\Jg\Wbn)(\cdot,\s')}_{L^2_x}
 \snorm{\tfrac{\JJtf}{\Sigma^\beta} \nbs^6(\Jg\Zbn) (\cdot,\s')}_{L^2_x} {\rm d} \s'
 \notag \\
&
\le  \tfrac{1+ \alpha }{14 \eps}  \int_0^{\s} \snorm{\tfrac{\JJof \Jgh}{\Sigma^\beta} \nbs^6(\Jg\Wbn)(\cdot,\s')}_{L^2_x}^2 {\rm d} \s'
+  \tfrac{4(1+ \alpha )}{\eps}   \int_0^{\s}  \snorm{\tfrac{\JJtf}{\Sigma^\beta} \nbs^6(\Jg\Zbn) (\cdot,\s')}_{L^2_x}^2 {\rm d} \s'
\,,
\label{eq:small:duck:10-H}  
\end{align} 
whereas 
\begin{equation}
J^{\AA_n}_{1,a}
= \tints \tfrac{1}{\Sigma^{2\beta}}      \mathsf{G_{bad}}  \sabs{  \nbs^6(\Jg\Wbn)  }^2 \,,
\label{eq:small:duck:11-H}  
\end{equation}
with $ \mathsf{G_{bad}}  = \JJss (\Q\p_\s+V\p_2) \JJ$, 
produces an anti-damping integral 
that must be combined with the associated damping integral  on the right side of \eqref{eq:I:n:12369-H}. 
For this purpose, with  $ \mathsf{G_{good}}$ defined in \eqref{Ggood-H}, we   see that for $\eps$ taken small enough, 
\begin{align} 
 \mathsf{G_{good}} + \mathsf{G_{bad}}
 & =  - \tfrac{3}{4} \JJh \Jg (\Q\p_\s+V\p_2) \JJ + \tfrac{1}{2} \JJss (\Q\p_\s+V\p_2) \Jg  - \Cn 
 \ge \tfrac{1+\alpha}{7 \eps} \JJh\Jg  
 \,,
\label{eq:small:duck:11:a-H}
\end{align} 
\end{subequations}
where we have used 
and \eqref{eq:fakeJg:LB-H} to obtain this lower-bound.

Having estimated the integrals  $J^{\AA_n}_{1} $,  $J^{\AA_n}_{2} $, $J^{\AA_n}_{4} $,  $J^{\AA_n}_{5} $,  $J^{\AA_n}_{7} $,  and $J^{\AA_n}_{8} $
in  \eqref{I4-Abn-a-i-H}, it remains for us to estimate the integrals  $J^{\AA_n}_{3} $ and $J^{\AA_n}_{6} $.   We will first treat the
integral $J^{\AA_n}_{3} $ which produces new energy and damping norms that we will crucially rely upon.
As before, the following identity is the key to the construction of these new energy and damping norms: from 
\eqref{Jg-evo-s} and \eqref{good-comm}, we find that
\begin{align}
\nbs^6 (\Jg \Wbn) 
&= \tfrac{2}{1+\alpha} \bigl(\nbs^6 (\Q \p_\s + V\p_2) \Jg - \tfrac{1-\alpha}{2} \nbs^6 (\Jg \Zbn) \bigr)
\notag
\\
&=  \tfrac{2}{1+\alpha}(\Q \p_\s + V\p_2)  \nbs^6  \Jg + \tfrac{2}{1+\alpha} \bigl(\nbs^6 V \nbs_2 \Jg + \doublecom{\nbs^6,V,\nbs_2 \Jg} \bigr)- \tfrac{1-\alpha}{1+\alpha} \nbs^6 (\Jg \Zbn)  
\,.
\label{JgWbn-qpsJG-identity}
\end{align}
We substitute \eqref{JgWbn-qpsJG-identity} into the integral $J^{\AA_n}_{3} $ in \eqref{J-Abn-3}, employ
\eqref{Sigma0i-ALE-s} and the adjoint formula
\eqref{adjoint-3-H}, use that by \eqref{JJ-vanish} $\JJ=0$ on $\thsd$,  and arrive at the additive decomposition
\begin{align} 
J^{\AA_n}_{3} 
&= \tfrac{3}{2} \tint \tfrac{1}{\Sigma^{2\beta}} (\Q\p_\s+V\p_2) \JJ\, \JJh\,   (\Jg\Wbn)  \nbs^6\Jg  \  \nbs^6(\Jg\Wbn) 
\notag\\
&= J^{\AA_n}_{3,a} + J^{\AA_n}_{3,b}+ J^{\AA_n}_{3,c} + J^{\AA_n}_{3,d} + J^{\AA_n}_{3,e}+ J^{\AA_n}_{3,f}+ J^{\AA_n}_{3,g}+ J^{\AA_n}_{3,h}+ J^{\AA_n}_{3,i}\,, 
\label{eq:J3:An:decompose-H}\\
J^{\AA_n}_{3,a}
&=  \tfrac{3}{2(1+ \alpha )}   \dint \tfrac{1}{\Sigma^{2\beta}} \Big(\!\!\!-\!\!(\Q\p_\s+V\p_2) \JJ\Big) \JJh (-\Jg\Wbn + \tfrac{13}{\eps} \Jg )\sabs{  \nbs^6\Jg}^2\Big|_\s \,,
\notag \\
J^{\AA_n}_{3,b}
&
=  \tfrac{3}{4(1+ \alpha )}    \tints \tfrac{1}{\Sigma^{2\beta}}\Big((\Q\p_\s+V\p_2) \JJ\Big)^2  \JJmh (-\Jg\Wbn + \tfrac{13}{\eps} \Jg) \sabs{\nbs^6\Jg}^2 \,,
\notag \\
J^{\AA_n}_{3,c}
&=  \tfrac{3}{2(1+ \alpha )}  \tints \tfrac{1}{\Sigma^{2\beta}} (\Q\p_\s+V\p_2) \JJ \, \JJh (\Q\p_\s+V\p_2) (-\Jg\Wbn + \tfrac{13}{\eps} \Jg) 
\ \sabs{ \nbs^6\Jg}^2 \,,
\notag \\
J^{\AA_n}_{3,d}
&=  \tfrac{3}{2(1+ \alpha )}  \tints \tfrac{1}{\Sigma^{2\beta}} (\Q\p_\s+V\p_2)^2\JJ \,   \JJh (-\Jg\Wbn +   \tfrac{13}{\eps} \Jg) 
\ \sabs{ \nbs^6\Jg}^2 \,,
\notag \\
J^{\AA_n}_{3,e}
&=  \tfrac{3}{2(1+ \alpha )}  \tints \tfrac{1}{\Sigma^{2\beta}} (\Q\p_\s+V\p_2)\JJ \,  \big(\Qc+V,_2 +2 \alpha \beta(\Zbn+\Abt)\big) \,   \JJh (-\Jg\Wbn +   \tfrac{13}{\eps} \Jg) 
\ \sabs{ \nbs^6\Jg}^2 \,,
\notag \\
 J^{\AA_n}_{3,f} 
&=  \tfrac{3}{2(1+ \alpha )}   \dint \tfrac{1}{\Sigma^{2\beta}} (\Q\p_\s+V\p_2) \JJ\,  \JJh (-\Jg\Wbn + \tfrac{13}{\eps} \Jg )\sabs{  \nbs^6\Jg}^2\Big|_\sin \,,
\,,
\notag\\
J^{\AA_n}_{3,g}
&=  \tfrac{39}{2\eps} \tints \tfrac{1}{\Sigma^{2\beta}} (\Q\p_\s+V\p_2) \JJ\,  \JJh  \Jg  \nbs^6\Jg  \  \nbs^6(\Jg\Wbn) 
\,,
\notag\\
J^{\AA_n}_{3,h}
&=  \tfrac{3}{(1+\alpha) } \tints \tfrac{1}{\Sigma^{2\beta}} (\Q\p_\s+V\p_2) \JJ\, \JJh  (\Jg\Wbn- {\tfrac{13}{\eps}} \Jg)  \nbs^6\Jg  \  \bigl(\nbs^6 V \nbs_2 \Jg + \doublecom{\nbs^6,V,\nbs_2 \Jg} \bigr) \,,
\notag \\
J^{\AA_n}_{3,i}
&= - \tfrac{3(1-\alpha)}{2(1+\alpha)} \tints \tfrac{1}{\Sigma^{2\beta}} (\Q\p_\s+V\p_2) \JJ\, \JJh   (\Jg\Wbn- {\tfrac{13}{\eps}} \Jg)  \nbs^6\Jg  \  \nbs^6(\Jg\Zbn)
\notag
\,.
\end{align} 
At this stage, it is convenient to define the function
\begin{align} 
\QQ =- \eps (\Q\p_\s+V\p_2) \JJ \ge  \tfrac{11(1 + \alpha )}{25}  \,,
\label{QQ-def}
\end{align} 
the inequality following from \eqref{qps-JJ-bound}.

We shall first obtain lower-bounds for the integrals $J^{\AA_n}_{3,a}$ and $J^{\AA_n}_{3,b}$ using 
\eqref{eq:signed:Jg};  we find that
\begin{subequations}
\label{eq:fuck:yeah:0-H}
\begin{align}
J^{\AA_n}_{3,a}
&\geq 
\tfrac{27 }{20(1+ \alpha )} \tfrac{1}{\eps^2} \snorm{\tfrac{ \QQ^{\frac{1}{2}} \JJof}{\Sigma^\beta} \nbs^6 \Jg(\cdot,\s)}_{L^2_x}^2 
\,,
\label{eq:fuck:yeah:1-H}
\\
J^{\AA_n}_{3,b}
&\geq 
\tfrac{27 }{40(1+ \alpha )} \tfrac{1}{\eps^3}  \int_0^{\s} \snorm{\tfrac{ \QQ \JJmof}{\Sigma^\beta} \nbs^6 \Jg(\cdot,\s')}_{L^2_x}^2 {\rm d} \s'
\,.
\label{eq:fuck:yeah:2-H}
\end{align}

For the remaining terms in the additive decomposition of $J^{\AA_n}_{3}$, using the bounds \eqref{table:derivatives}, \eqref{eq:signed:Jg}, \eqref{bootstraps-H}, \eqref{eq:Q:all:bbq-H},   \eqref{eq:waitin:for:the:bus-s},  \eqref{geom-H-original},  \eqref{eq:Jg:Wbn:improve:material:a-H}, 
 together with the identitiy \eqref{Jg-evo-s} and choosing $\eps$ sufficiently small, we  obtain the bounds 
\begin{align}
J^{\AA_n}_{3,c}
&\ge \tfrac{3}{2(1+ \alpha )}  {\tfrac{1}{\eps^3}} \tints \tfrac{\QQ}{\Sigma^{2\beta}}  \JJh (\Q\p_\s+V\p_2) (\Jg\Wbn)
\ \sabs{ \nbs^6\Jg}^2
- \tfrac{39(1-\alpha)}{4(1+ \alpha )}  {\tfrac{1}{\eps^3}} \tints \tfrac{\QQ}{\Sigma^{2\beta}} \JJh ( \Jg \Zbn)  \sabs{ \nbs^6\Jg}^2
\notag\\
&\qquad
+ \tfrac{39 \cdot 9 }{4}  {\tfrac{1}{\eps^3}} \tints \tfrac{\QQ}{\Sigma^{2\beta}} \JJh    \sabs{ \nbs^6\Jg}^2
- \tfrac{39 \cdot 13 }{4}  {\tfrac{1}{\eps^3}} \tints \tfrac{\QQ}{\Sigma^{2\beta}} \JJh  \Jg \sabs{ \nbs^6\Jg}^2
\notag\\
&\ge
- \big( \tfrac{5733 }{ 40} +\Cn \eps \big)  {\tfrac{1}{\eps^3}} \int_\sin^{\s}   \snorm{\tfrac{\QQ^{\frac{1}{2}}  \JJof}{\Sigma^\beta} \nbs^6 \Jg(\cdot,\s')}_{L^2_x}^2 {\rm d} \s'
\ge 
- {\tfrac{144}{\eps^3}} \int_\sin^{\s}   \snorm{\tfrac{\QQ^{\frac{1}{2}}  \JJof}{\Sigma^\beta} \nbs^6 \Jg(\cdot,\s')}_{L^2_x}^2 {\rm d} \s'
\label{eq:fuck:yeah:3-H} \,,
\end{align}
and by additionally using \eqref{qps2JJ-bound},  
\begin{align} 
\sabs{J^{\AA_n}_{3,d}}
+
\sabs{J^{\AA_n}_{3,e}}
&\leq  \tfrac{1600^2}{\eps^3}  \tint \tfrac{\QQ}{\Sigma^{2\beta}} \JJh    \sabs{ \nbs^6\Jg}^2
\label{eq:fuck:yeah:4-H}
\,.
\end{align}
Continuing to rely on the above referenced bounds, we find that 
\begin{align} 
\sabs{J^{\AA_n}_{3,f}}
&\leq \tfrac{3}{2(1+\alpha)}  {\tfrac{402(1+ \alpha ) }{\eps}} \big( \tfrac{11}{10} \big)^{\frac{1}{2}}   \tfrac{14}{\eps} (\tfrac{3}{\kappa_0})^{2\beta} \eps \Cdatatwo
\le
\tfrac{9000} { \eps }  (\tfrac{3}{\kappa_0})^{2\beta}   \Cdatatwo
\label{eq:fuck:yeah:5-H}
\,,
\\
\sabs{J^{\AA_n}_{3,g}}
&\leq \tfrac{39}{ \eps^2} \big( \tfrac{51}{50} \cdot \tfrac{6}{5} \big)^{\frac{1}{2}} \int_\sin^{\s}
\snorm{\tfrac{ \QQ \JJmof }{\Sigma^\beta} \nbs^6 \Jg(\cdot,\s')}_{L^2_x} 
\snorm{\tfrac{\Q^{\frac{1}{2}} \JJtf \Jgh }{\Sigma^{\beta}} \nbs^6 (\Jg \Wbn)(\cdot,\s')}_{L^2_x} {\rm d}\s' 
\notag\\
&\leq 
\tfrac{1}{4(1+ \alpha )} \tfrac{1}{\eps^3} \int_\sin^{\s} \snorm{\tfrac{\QQ \JJmof}{\Sigma^\beta} \nbs^6 \Jg(\cdot,\s')}_{L^2_x}^2 {\rm d} \s'
+ 
\tfrac{44(1+\alpha )}{\eps} 
\int_\sin^{\s} 
\snorm{\tfrac{\Q^{\frac{1}{2}} \JJtf \Jgh}{\Sigma^{\beta}} \nbs^6 (\Jg \Wbn)(\cdot,\s')}_{L^2_x}^2 {\rm d}\s' 
\label{eq:fuck:yeah:6-H}
\,,
\\
\sabs{J^{\AA_n}_{3,h}}
&\les \eps \mathsf{K} (\tfrac{4}{\kappa_0})^{2\beta} \brak{\mathsf{B}_6}^2
\label{eq:fuck:yeah:7-H}
\,,
\\
\sabs{J^{\AA_n}_{3,i}}
&\leq  \tfrac{24}{\eps^2} \int_\sin^{\s}  
\snorm{\tfrac{\QQ \JJmof}{\Sigma^\beta} \nbs^6 \Jg (\cdot,\s')}_{L^2_x}
\snorm{\tfrac{\JJtf}{\Sigma^\beta} \nbs^6(\Jg\Zbn) (\cdot,\s')}_{L^2_x}
{\rm d} \s'
\notag\\
&\leq \tfrac{1}{4(1+\alpha)} \tfrac{1}{\eps^3} 
\int_\sin^{\s} 
\snorm{\tfrac{\QQ\JJmof}{\Sigma^\beta} \nbs^6 \Jg(\cdot,\s')}_{L^2_x}^2 
{\rm d} \s'
+ 
\tfrac{24^2(1+\alpha )}{\eps} 
\int_\sin^{\s} 
\snorm{\tfrac{\JJtf}{\Sigma^\beta} \nbs^6(\Jg\Zbn) (\cdot,\s')}_{L^2_x}^2 
{\rm d} \s'
\,.
\label{eq:fuck:yeah:8-H}
\end{align}
\end{subequations}
The bounds in~\eqref{eq:fuck:yeah:0-H}  complete our estimates for the nine terms in the decomposition of $J^{\AA_n}_{3}$.

To estimate $J^{\AA_n}_{6} $, we decompose this integral as
\begin{align*} 
J^{\AA_n}_{6} & =J^{\AA_n}_{6,\le}  + J^{\AA_n}_{6,>}  \,,
\notag \\
J^{\AA_n}_{6,\le} & = -  \dint {\bf 1}_{|x_1|\le 13\pi\eps} \tfrac{\Q}{\Sigma^{2\beta}} \JJss    \Jg \Wbn \nbs^6 \Jg \  \nbs^6(\Jg\Wbn  ) \Big|_\s
\notag \\
J^{\AA_n}_{6,>} & = -  \dint {\bf 1}_{|x_1|> 13\pi\eps} \tfrac{\Q}{\Sigma^{2\beta}} \JJss    \Jg \Wbn \nbs^6 \Jg \  \nbs^6(\Jg\Wbn  ) \Big|_\s
\end{align*} 
The integral $J^{\AA_n}_{6,\le}$ requires an application of the $\eps$-Young inequality.  More precisely, using \eqref{bs-JgnnWb}, 
\eqref{Qd-lower-upper-H}, \eqref{Q-lower-upper-H}, 
\eqref{JJ-le-Jg:a}, \eqref{QQ-def},  and choosing $\eps$ sufficiently small,  we obtain that
\begin{align} 
\sabs{J^{\AA_n}_{6,\le}}
& \le
\big( \tfrac{101}{100} \big)^{\!\!\frac{1}{2}} \big( \tfrac{1002}{1000} \big)^{\!\!\frac{1}{2}}\big( \tfrac{25}{11(1+ \alpha )} \big)^{\!\!\frac{1}{2}}
 {\tfrac{1+\eps}{\eps}}  \dint\Big( \tfrac{1}{\Sigma^{2\beta}}  \sabs{ \QQ^{\frac{1}{2}} \JJof \nbs^6 \Jg} \  \sabs{\Q^{\frac{1}{2}} \JJtf \Jgh  \nbs^6(\Jg\Wbn  )}\Big)(x,\s) {\rm d}x
\notag \\
&\le 
\tfrac{25}{52}   \snorm{\tfrac{\Q^{\frac{1}{2}} \JJtf \Jgh }{\Sigma^\beta} \nbs^6(\Jg\Wbn)(\cdot,\s)}_{L^2_x} ^2
+ \tfrac{13}{10(1+ \alpha )\eps^2}  \snorm{\tfrac{ \QQ^{\frac{1}{2}} \JJof}{\Sigma^\beta} \nbs^6 \Jg(\cdot,\s)}_{L^2_x}^2
\,.
\label{eq:fuck:yeah:9-H}
\end{align}
The specific pre-factors $\tfrac{25}{52}$ and $ \tfrac{13}{10(1+ \alpha )\eps^2} $ allow us to absorb the two terms on the right
side of \eqref{eq:fuck:yeah:9-H} by the associated energy norms.

Next, we use \eqref{JJ-le-Jg:b} in place of \eqref{JJ-le-Jg:a} and find that
\begin{align} 
\sabs{J^{\AA_n}_{6,>}}
& \le
10^{\frac{1}{2}}  \big( \tfrac{1002}{1000} \big)^{\!\!\frac{1}{2}}\big( \tfrac{25}{11(1+ \alpha )} \big)^{\!\!\frac{1}{2}}
\cdot \Cn \eps \dint\Big( \tfrac{1}{\Sigma^{2\beta}}  \sabs{ \QQ^{\frac{1}{2}} \JJof \nbs^6 \Jg} \  \sabs{\Q^{\frac{1}{2}} \JJtf \Jgh  \nbs^6(\Jg\Wbn  )}\Big)(x,\s) {\rm d}x
\notag \\
&\les
 (\tfrac{4}{\kappa_0})^{2\beta}  \brak{\mathsf{B}_6}^2 
\,.
\label{eq:fuck:yeah:9-H-<}
\end{align}

To summarize, with ~\eqref{eq:small:duck:7-H}, \eqref{eq:small:duck:8-H}, \eqref{eq:small:duck:8a-H}, \eqref{eq:small:duck:8b-H}, \eqref{eq:small:duck:big:duck-H}, \eqref{eq:fuck:yeah:0-H}, \eqref{eq:fuck:yeah:9-H} and~\eqref{eq:fuck:yeah:9-H-<}, the inequality  $\JJ\le \tfrac{101}{100}  \Jg$  for
$(x,\s) \in \Hdm$ such that $|x_1|\le 13  \pi \eps$, and 
a specific Caucy-Young inequality, we have found that
\begin{align}
 I^{\AA_n}_{4}
&\geq  
- \Cn (\tfrac{4}{\kappa_0})^{2\beta}  \brak{\mathsf{B}_6}^2
-  {\tfrac{4}{\eps}}  (\tfrac{3}{\kappa_0})^{2\beta}   \Cdatatwo
\notag\\
&\qquad 
+ \tint \tfrac{1}{\Sigma^{2\beta}} \bigl( \tfrac{1+ \alpha }{7 \eps}\JJh \Jg - \mathsf{G_{good}} \bigr)   \sabs{  \nbs^6(\Jg\Wbn)  }^2
- \tfrac{(24^2 +4)(1+\alpha )}{\eps} 
\int_0^{\s} \snorm{\tfrac{\JJtf }{\Sigma^\beta} \nbs^6(\Jg\Zbn) (\cdot,\s')}_{L^2_x}^2 {\rm d} \s'
\notag\\
&\qquad
+ \tfrac{1 }{20(1+ \alpha )} \tfrac{1}{\eps^2}  \snorm{\tfrac{\QQ^{\frac{1}{2}}  \JJof }{\Sigma^\beta} \nbs^6 \Jg(\cdot,\s)}_{L^2_x}^2 
-  {\tfrac{144 + 1600^2 }{\eps^3}} 
\int_\sin^{\s}  
 \snorm{\tfrac{\QQ^{\frac{1}{2}} \JJof }{\Sigma^\beta} \nbs^6 \Jg(\cdot,\s')}_{L^2_x}^2 {\rm d} \s'
\notag\\
&\qquad
+ \tfrac{7}{40(1+ \alpha )} \tfrac{1}{\eps^3}\int_\sin^{\s} \snorm{\tfrac{\QQ \JJmof }{\Sigma^\beta} \nbs^6 \Jg(\cdot,\s')}_{L^2_x}^2 {\rm d} \s'
- \tfrac{25}{52}
\snorm{\tfrac{\Q^{\frac{1}{2}} \JJtf \Jgh}{\Sigma^\beta} \nbs^6(\Jg\Wbn)(\cdot,\s)}_{L^2_x}^2
\notag\\
&\qquad
- 
\tfrac{44(1+\alpha )}{\eps} 
\int_\sin^{\s} 
\snorm{\tfrac{\Q^{\frac{1}{2}} \JJtf \Jgh }{\Sigma^{\beta}} \nbs^6 (\Jg \Wbn)(\cdot,\s')}_{L^2_x}^2 {\rm d}\s' 
\,. \label{eq:fuck:yeah:end-H}
\end{align}

\subsubsection{The integrals $I^{\ZZ_n}_7$, $I^{\AA_n}_5$, and $I^{\ZZ_n}_8$} 
Following the same analysis as in Sections \ref{sec:I:Zn:7:integral}--\ref{sec:IZ8}, we obtain
\begin{equation}
\sabs{I^{\ZZ_n}_{7}} + \sabs{ I^{\AA_n}_{5}}  + \sabs{I^{\ZZ_n}_8} 
\les (1+   \Cn\mathsf{K} \eps)   (\tfrac{4}{\kappa_0})^{2\beta} \brak{\mathsf{B}_6}^2
\,.
\label{I-Zbn7-Zbn8-Abn5-H}
\end{equation}

\subsection{The forcing and commutator terms}\label{sec:14:forcing:comm}
It remains for us to bound  the  integrals $ I^{\WW_n}_6$, $ I^{\ZZ_n}_{10}$, $ I^{\AA_n}_{10}$ in 
\eqref{I6-Wbn-H}, \eqref{I10-Zbn-H}, and \eqref{I10-Abn-H}, respectively.
In order to estimate these integrals, we use  the definitions  for the forcing functions in \eqref{forcing-nt} together with
the definitions of the so-called remainder and commutator functions in \eqref{Cw-Rw-comm}, \eqref{Cz-Rz-comm}, and  \eqref{Ca-Ra-comm},
whose bounds were given in~\eqref{eq:I:W:nn:6:all-H},~\eqref{eq:I:Z:nn:10:all-H}, and~\eqref{eq:I:A:nn:10:all-H}. 
Using the Cauchy-Schwartz inequality, the bound $\JJ \le \tfrac{21}{10}  \Jg$ from~\eqref{JJ-le-Jg-s}, and~\eqref{bootstraps-Dnorm:6}, we deduce from the aforementioned bounds that
\begin{subequations}
\begin{align}
\sabs{I^{\WW_n}_6}
&\leq \Cn (\tfrac{4}{\kappa_0})^{2\beta}  \brak{\mathsf{B}_6}^2 
\label{eq:I:W:nn:6:final-H}
\,.
\end{align}
Similarly, using \eqref{eq:I:Z:nn:10:all-H} we also have that
\begin{align}
\sabs{I^{\ZZ_n}_{10}}
&\le  \int_\sin^\s  \snorm{\tfrac{ \JJtf }{\Sigma^{\beta-1}} \nbs^6(\Jg \Zbn)(\cdot,\s')}_{L^2_x}  
 \Bigl( \snorm{\tfrac{\JJtf }{\Sigma^{\beta-1}} \nbs^6\Fzn (\cdot,\s')}_{L^2_x}  
+ \snorm{\tfrac{\JJtf }{\Sigma^{\beta-1}} \mathcal{R}_\Zb^\nn (\cdot,\s')}_{L^2_x}  
+ \snorm{\tfrac{\JJtf }{\Sigma^{\beta-1}} \mathcal{C}_\Zb^\nn(\cdot,\s')}_{L^2_x}  \Bigr){\rm d} \s'
\notag\\
&\le \tfrac{4 (1+ \alpha )}{\eps} \int_\sin^\s \snorm{ \tfrac{\JJtf }{\Sigma^{\beta}}   \nbs^6 (\Jg \Zbn)(\cdot,\s')}_{L^2_{x}}^2 {\rm d}\s'
+ \Cn   (\tfrac{4}{\kappa_0})^{2\beta}    \mathsf{K} \brak{\mathsf{B_6}}^2
   \,.  
\label{I10-Zn-H}
\end{align} 
In the same way, using \eqref{eq:I:A:nn:10:all-H} we  obtain the bound
\begin{equation} 
\sabs{ I^{\AA_n}_{10} } \le 
\tfrac{4 (1+ \alpha )}{\eps} \int_\sin^\s \snorm{ \tfrac{\JJtf}{\Sigma^{\beta}}   \nbs^6 (\Jg \Abn)(\cdot,\s')}_{L^2_{x}}^2 {\rm d}\s'
+ \Cn   (\tfrac{4}{\kappa_0})^{2\beta}    \mathsf{K} \brak{\mathsf{B_6}}^2
   \,.  \label{I10-An-H}
\end{equation} 
\end{subequations}
Collecting the bounds \eqref{eq:I:W:nn:6:final-H}, \eqref{I10-Zn-H}, and \eqref{I10-An-H},   we have that
\begin{equation} 
\sabs{ I^{\WW_n}_6}+ \sabs{ I^{\ZZ_n}_{10}} + \sabs{ I^{\AA_n}_{10}}
  \le
\tfrac{4 (1+ \alpha )}{\eps} \int_\sin^\s \snorm{ \tfrac{\JJtf}{\Sigma^{\beta}}   \nbs^6 (\Jg\Zbn,\Jg \Abn)(\cdot,\s')}_{L^2_{x}}^2 {\rm d}\s'
+ \Cn   (\tfrac{4}{\kappa_0})^{2\beta}  \mathsf{K} \brak{\mathsf{B_6}}^2\,.
\label{eq:I:forcing:normal-H}
\end{equation}

\subsection{Conclusion of the six derivative normal energy bounds}
\label{sec:D6:n:final-H}
We return to the energy identity~\eqref{D5-L2-H}, with the decompositions~\eqref{Integral-Wbn-H}, \eqref{Integral-Zbn-H}, 
and~\eqref{Integral-Abn-H}. 
We collect the lower bounds   \eqref{eq:I:n:12369-H}, \eqref{eq:fuck:yeah:end-H},    
the upper bounds~\eqref{IWn5+IZn5+IAn7-H}, 
\eqref{eq:small:duck:0-H}, \eqref{I-Zbn7-Zbn8-Abn5-H}, \eqref{eq:I:forcing:normal-H}, 
 and the initial data assumption~\eqref{table:derivatives}, to arrive at
\begin{align}
0
&\geq 
\bigl( \tfrac{1}{52} -  \eps^ {\frac{7}{4}} \bigr) 
\snorm{\tfrac{\Q^{\frac{1}{2}} \JJtf\Jgh }{\Sigma^\beta} \nbs^6(\Jg \Wbn,\Jg \Zbn,\Jg \Abn)(\cdot,\s)}_{L^2_x}^2
- \Cn (\tfrac{4}{\kappa_0})^{2\beta}  \brak{\mathsf{B}_6}^2
-  {\tfrac{10000}{\eps}}  (\tfrac{3}{\kappa_0})^{2\beta}   \Cdatatwo
\notag\\
&\qquad
+ \tfrac{1+\alpha}{7 \eps} 
\int_\sin^{\s} \snorm{\tfrac{\JJof \Jgh}{\Sigma^\beta}   \nbs^6(\Jg\Wbn)(\cdot,\s')}_{L^2_x}^2 {\rm d}\s'
+  \tfrac{17(1+\alpha)  }{100 \eps} \tints \tfrac{1}{\Sigma^{2\beta}} 
 \JJh \Jg \big( \dl \sabs{\nbs^6(\Jg \Zbn)}^2 + (1+\dl) \sabs{\nbs^6(\Jg \Abn)}^2  \big)
\notag\\
&\qquad 
+\tints   \tfrac{1}{\Sigma^{2\beta}} 
 \big({\tfrac{ 2\alpha(2\beta-1) + \bubu{  4 \alpha } }{5}} - (24^2+4)(1+\alpha ) \big) {\tfrac{1}{\eps}}    \JJss   \big( \sabs{\nbs^6 (\Jg\Zbn) }^2 +  \sabs{\nbs^6( \Jg\Abn) }^2\big)
\notag\\
&\qquad 
-   \tfrac{ 7(1+3\alpha-4\alpha\beta)   +  88(1+\alpha ) }{ 2\eps} 
\int_\sin^{\s} \snorm{\tfrac{\Q^{\frac{1}{2}} \JJtf \Jgh }{\Sigma^\beta} \nbs^6(\Jg\Wbn,\Jg \Zbn,\Jg\Abn)(\cdot,\s')}_{L^2_x}^2
{\rm d} \s'
+   \tfrac{1}{20\eps^2}  \snorm{\tfrac{ \QQ^{\frac{1}{2}} \JJof }{\Sigma^\beta} \nbs^6 \Jg(\cdot,\s)}_{L^2_x}^2
\notag\\
&\qquad  
-    \tfrac{144+ 1600^2}{\eps^3} 
             \int_\sin^{\s}   \snorm{\tfrac{\QQ^{\frac{1}{2}}  \JJof }{\Sigma^\beta} \nbs^6 \Jg(\cdot,\s')}_{L^2_x}^2 {\rm d} \s'
+ \tfrac{7}{40(1+ \alpha ) \eps^3} \int_\sin^{\s} \snorm{\tfrac{\QQ \JJmof }{\Sigma^\beta} \nbs^6 \Jg(\cdot,\s')}_{L^2_x}^2 {\rm d} \s'
\,,
\label{eq:normal:conclusion:1-H}
\end{align}
where
 $\Cn = \Cn(\alpha,\kappa_0,\Cdata)$ is independent of $\beta$ (and $\eps$, as always).

We choose $\beta = \beta(\alpha)$ to be sufficiently large to ensure that the damping for $\nbs^6(\Jg\Zbn,\Jg\Abn)$ is strong enough, i.e.,  
\begin{equation*}
{\tfrac{ 2\alpha(2\beta-1) + \bubu{ 4 \alpha } }{5}} - 580(1+\alpha ) \ge  1 \,.
\end{equation*}
More precisely, we choose $\beta$ to ensure equality in the above inequality; namely, we have that 
\begin{equation}
\bubu{ \beta_{\alpha}   :=  \tfrac{ 7(415+414 \alpha) } {4 \alpha } }  \,.
\label{eq:normal:bounds:beta-H}
\end{equation}
With this choice of $\beta = \beta_\alpha$, we return to \eqref{eq:normal:conclusion:1-H}, and choose $\eps$ to be sufficiently small in terms of $\alpha,\kappa_0,\Cdata$. After re-arranging,  we deduce that
\begin{align}
& \tfrac{1}{53}  \snorm{\tfrac{\Q^{\frac{1}{2}} \JJtf \Jgh }{\Sigma^{\beta_\alpha}} \nbs^6(\Jg \Wbn,\Jg \Zbn,\Jg \Abn)(\cdot,\s)}_{L^2_x}^2
 +   \tfrac{1}{20\eps^2}  \snorm{\tfrac{\QQ^{\frac{1}{2}}  \JJof }{\Sigma^\beta} \nbs^6 \Jg(\cdot,\s)}_{L^2_x}^2  
\notag \\
&\qquad
 + \tfrac{1+\alpha}{7 \eps} 
\int_\sin^{\s} \snorm{\tfrac{\JJof \Jgh}{\Sigma^\beta}   \nbs^6(\Jg\Wbn,\Jg\Abn)(\cdot,\s')}_{L^2_x}^2 {\rm d}\s'
+ \tfrac{17\dl(1+\alpha)}{100 \eps} 
\int_\sin^{\s} \snorm{\tfrac{\JJof \Jgh}{\Sigma^\beta}   \nbs^6(\Jg\Zbn)(\cdot,\s')}_{L^2_x}^2 {\rm d}\s'
 \notag\\
&\qquad
+ \tfrac{1}{ \eps} 
\int_\sin^{\s} \snorm{\tfrac{\JJtf }{\Sigma^{\beta_\alpha}}   \nbs^6(\Jg\Zbn,\Jg\Abn)(\cdot,\s')}_{L^2_x}^2 {\rm d}\s'
+ \tfrac{7}{40(1+ \alpha )\eps^3}\int_\sin^{\s} \snorm{\tfrac{ \QQ \JJmof }{\Sigma^\beta} \nbs^6 \Jg(\cdot,\s')}_{L^2_x}^2 {\rm d} \s'
\notag\\
&\le
{\tfrac{10000}{\eps}}  (\tfrac{3}{\kappa_0})^{2\beta}   \Cdatatwo
+\Cn (\tfrac{4}{\kappa_0})^{2\beta}  \brak{\mathsf{B}_6}^2
\notag\\
 &\qquad
+   \tfrac{C_ \alpha }{\eps} 
\int_\sin^{\s} \snorm{\tfrac{\Q^{\frac{1}{2}} \JJtf \Jgh }{\Sigma^\beta} \nbs^6(\Jg\Wbn,\Jg \Zbn,\Jg\Abn)(\cdot,\s')}_{L^2_x}^2
{\rm d} \s'
+ {\tfrac{C}{\eps^3}} 
             \int_\sin^{\s}   \snorm{\tfrac{\QQ^{\frac{1}{2}}  \JJof }{\Sigma^\beta} \nbs^6 \Jg(\cdot,\s')}_{L^2_x}^2 {\rm d} \s'
\,,
\label{eq:normal:conclusion:2-H}
\end{align}
where $C_ \alpha $ is a universal constant which is independent of $\kappa_0,\Cdata$, 
 $C $ is a universal constant which is independent of $\alpha,\kappa_0,\Cdata$,  and $\Cn$ is as defined above.

An application of Gr\"onwall's inequality to  \eqref{eq:normal:conclusion:2-H} shows that
there exists a constant 
\begin{equation*}
 \check{\mathsf{c}}_{\alpha,\kappa_0,\dl} > 0 
\end{equation*}
which only depends on $\alpha$, $\kappa_0$, and $\dl$, such that 
\begin{align} 
&\sup_{\s \in [\sin,\sfin]}\snorm{\tfrac{\JJtf \Jgh }{\Sigma^{\beta_\alpha}} \nbs^6(\Jg \Wbn,\Jg \Zbn,\Jg \Abn)(\cdot,\s)}_{L^2_x}^2
+   \tfrac{1}{\eps^2} \sup_{\s \in [\sin,\sfin]} \snorm{\tfrac{ \JJof }{\Sigma^\beta} \nbs^6 \Jg(\cdot,\s)}_{L^2_x}^2  
\notag \\
&\qquad
 + \tfrac{1}{ \eps} \int_\sin^{\sfin} \snorm{\tfrac{\JJof \Jgh}{\Sigma^\beta}   \nbs^6(\Jg\Wbn,\Jg\Zbn, \Jg\Abn)(\cdot,\s')}_{L^2_x}^2 {\rm d}\s'
+ \tfrac{1}{ \eps} 
\int_\sin^{\sfin} \snorm{\tfrac{\JJtf }{\Sigma^{\beta_\alpha}}   \nbs^6(\Jg\Zbn,\Jg\Abn)(\cdot,\s')}_{L^2_x}^2 {\rm d}\s'
\notag \\
&\qquad
+  \tfrac{1}{\eps^3}\int_\sin^{\sfin} \snorm{\tfrac{ \JJmof }{\Sigma^\beta} \nbs^6 \Jg(\cdot,\s')}_{L^2_x}^2 {\rm d} \s'
 \notag\\
 &  
 \leq  \check{\mathsf{c}}_{\alpha, \kappa_0,\dl}  \tfrac{1}{\eps} (\tfrac{4}{\kappa_0})^{2\beta_{\alpha,\kappa_0}  } 
\big( \Cdatatwo  + \eps \Cn  \mathsf{K}^2  \brak{\mathsf{B}_6}^2\big)
 \,.
 \label{eq:normal:conclusion:3-H}
\end{align}

Just as in the conclusion of our previously established energy estimates,  at this stage,  we multiply the above estimate by 
$\kappa_0^{2 \beta_\alpha}$, appeal to \eqref{bs-Sigma}, drop the energy and damping terms for $ \nbs^6 \Jg$ (since these were bounded already in Proposition~\ref{prop:geometry}),  and recall the definitions of $\widetilde{\mathcal{E}}_{6,\nnn}^2(\s)$ and 
$\widetilde{\mathcal{D}}_{6,\nnn}^2(\s) $ to deduce that 
\begin{align}
\eps \sup_{\s \in [\sin,\sfin]} \widetilde{\mathcal{E}}_{6,\nnn}^2(\s)
+\widetilde{\mathcal{D}}_{6,\nnn}^2(\sfin) 
&  \le  \check{\mathsf{c}}_{\alpha, \kappa_0,\dl}  4^{2\beta_\alpha} 
 \Big( \Cdatatwo  + \Cn \eps \mathsf{K}^2  \brak{\mathsf{B}_6}^2\Bigr)
 \notag\\
 &  \leq  \mathsf{B}_6^2  \check{\mathsf{c}}_{\alpha, \kappa_0,\dl}  4^{2\beta_\alpha} 
 \Big(\tfrac{\Cdatatwo}{ \mathsf{B}_6^2}  + \Cn \eps \mathsf{K}^2  \tfrac{\brak{\mathsf{B}_6}^2} {\mathsf{B}_6^2}\Bigr)
 \,.
 \label{eq:normal:conclusion:4-H}
\end{align}
Since $ \mathsf{B}_6 $ and $\mathsf{K}$ defined by \eqref{eq:B6:choice:1-H} and \eqref{eq:K:choice:1-H}, respectively, 
upon ensuring that 
\begin{align}
\mathsf{B}_6 &\geq 
4  \check{\mathsf{c}}_{\alpha, \kappa_0,\dl} ^{\frac 12} 
4^{\beta_\alpha}
\Cdata\,,
\label{eq:B6:choice:2-H}
\end{align}
where $\beta_\alpha$ is as defined in~\eqref{eq:normal:bounds:beta-H}, 
and taking $\eps$ sufficiently small in terms of $\alpha,\kappa_0,\Cdata$, and $\dl$, we deduce from \eqref{eq:normal:conclusion:4-H} that 
\begin{equation}
\eps \sup_{\s \in [\sin,\sfin]} \widetilde{\mathcal{E}}_{6,\nnn}^2(\s)
+\widetilde{\mathcal{D}}_{6,\nnn}^2(\sfin) 
\leq \tfrac{1}{8}   \mathsf{B}_6^2 
 \,,
\label{eq:normal:conclusion:5-H}
\end{equation}
which closes the ``normal part'' of the remaining bootstrap \eqref{bootstraps-Dnorm:6}.

\subsection{Closing the bootstrap for the sixth order energy}
\label{sec:bootstrap:closed:US}
Combining~\eqref{eq:normal:conclusion:5-H}  with \eqref{eq:hate:13-H}
we arrive at the same inequality as obtained in
\eqref{eq:normfinal1}
\begin{align}
\eps \sup_{\s \in [\sin,\sfin]} \widetilde{\mathcal{E}}_{6}(\s)
+\widetilde{\mathcal{D}}_{6}(\sfin) 
\leq \tfrac{1}{2}   \mathsf{B}_6 
\,,  \label{eq:normfinal2-H}
\end{align}
which closes the bootstrap \eqref{boots-H}  (cf.~\eqref{bootstraps-Dnorm:6}) in the upstream coordinate system  \eqref{t-to-s-transform-H}.

\def\Qb{ {\bar{\mathsf{Q}} }}

\section{Optimal regularity for velocity, sound speed, and ALE  map}
\label{sec:optimal:reg}

We discuss the optimal regularity of the fields $U = u\circ \psi$ and $\Sigma = \sigma \circ \psi$ (as defined in~\eqref{vsL}), and of the ALE map $\psi$ (as defined in~\eqref{psi-def}). That is, we show how the bootstrap bounds~\eqref{bootstraps} obtained for the differentiated Riemann variables $(\Wb,\Zb,\Ab)$ and the geometric quantities $(\Jg,h,_2)$ at the sixth derivative level, imply corresponding estimates for the un-differentiated unknowns $(U,\Sigma,h)$ at the seventh derivative level. Our main result is:

\begin{proposition}[\bf Optimal regularity for $(U,\Sigma,h)$]
\label{prop:D7:U:Sigma:h}
Let $\varphi$ be a weight function, defined as follows:
\begin{itemize}[leftmargin=16pt] 
\item For the ``shock formation'' in Sections~\ref{sec:formation:setup}--\ref{sec:sixth:order:energy}, let $\varphi = \mathcal{J}$ as defined in \eqref{eq:fake:Jg:def}. 
\item For the ``downstream development'' in Section~\ref{sec:downstreammaxdev}, let $\varphi = \JJ$, as defined in~\eqref{eq:fake:Jg:def:ds}. 
\item For the ``upstream development'' of Section~\ref{sec:upstreammaxdev}, let $\varphi = \JJ$, as defined by~\eqref{JJ-def-plus-t} and~\eqref{JJ-def-minus-t}. 
\end{itemize}
Then, for each of these weights $\varphi$ individually, we have the bounds
\begin{subequations}
\label{eq:D7:U:Sigma:h}
\begin{align}
\eps^{\frac 12} \sup_{\s\in[0,\eps]} 
\snorm{\varphi^{\frac 34} \Jgh \nbs^7 \Sigma (\cdot,\s)}_{L^2_x} 
+ 
\snorm{\varphi^{\frac 14} \Jgh \nbs^7 \Sigma}_{L^2_\s L^2_x} 
&\les 
\eps \mathsf{K} \brak{\mathsf{B}_6}
\\
\eps^{\frac 12} 
\snorm{\varphi^{\frac 34} \Jgh  \nbs^7 U }_{L^\infty_{\s} L^2_x} 
+ 
\snorm{\varphi^{\frac 14} \Jgh \nbs^7 U}_{L^2_\s L^2_x} 
&\les 
\eps \mathsf{K} \brak{\mathsf{B}_6}
\\
\eps^{\frac 12} 
\snorm{\varphi^{\frac 14} \nbs^7  h}_{L^\infty_{\s} L^2_{x}} 
+
\snorm{\nbs^7  h}_{L^2_{\s}L^2_x} 
&\les 
\eps^2 \mathsf{K} \brak{\mathsf{B}_6}
\,.
\end{align}
\end{subequations}
In~\eqref{eq:D7:U:Sigma:h} it is understood that the $(x,\s)$ coordinates are given by~\eqref{t-to-s-transform}  
for Sections~\ref{sec:formation:setup}--\ref{sec:sixth:order:energy}, by~\eqref{t-to-s-transform-P} for Section~\ref{sec:downstreammaxdev}, and by~\eqref{t-to-s-transform-H} for Section~\ref{sec:upstreammaxdev}.
\end{proposition}

\begin{remark}[\bf Bounds in terms of the $(x,t)$ variables]
The estimates provided by~\eqref{eq:D7:U:Sigma:h} only concern $(x,\s)$ variables, and we have dropped the tildes, as described in Remark~\ref{rem:no:tilde}, Remark~\ref{rem:no:tilde-P}, and Remark~\ref{rem:US:drop:tilde}. It is clear however that the $L^2_{x,\s}$ estimates in~\eqref{eq:D7:U:Sigma:h} may be converted into $L^2_{x,t}$ coordinates by using the change of variables formula, and that the $L^\infty_\s L^2_x$ imply a suitable bound in $(x,t)$ coordinates. For example, for the analysis in Sections~\ref{sec:formation:setup}--\ref{sec:sixth:order:energy}, the Jacobian of the map $(x,t) \mapsto (x,\s)$ present in~\eqref{eq:f:tilde:f} is easily seen to equal $|\p_t \mathfrak{q}| = \Qd$, and~\eqref{Qd-lower-upper} gives global upper and lower bounds for $\Qd$ (which are strictly positive and depend only on $\alpha$). This matter was previously addressed in Remark~\ref{rem:L2:norms:x:s:x:t:A}. As such, with the spacetime $\mathcal{P}$ defined in~\eqref{eq:spacetime:smooth}, we deduce from the $L^2_\s$ bounds in~\eqref{eq:D7:U:Sigma:h} that 
\begin{subequations}
\label{eq:nbs:7:shock:dev:x:t}
\begin{equation}
\label{eq:nbs:7:shock:dev:x:t:2}
\|\mathcal{J}^{\frac 14} \Jgh \nb^7 \Sigma\|_{L^2_{x,t}(\mathcal{P})}
+
\|\mathcal{J}^{\frac 14} \Jgh \nb^7 U\|_{L^2_{x,t}(\mathcal{P})}
+
\tfrac 1\eps
\| \Jgh \nb^7 h\|_{L^2_{x,t}(\mathcal{P})}
\les 
\eps \mathsf{K} \brak{\mathsf{B}_6}
\,.
\end{equation}
On the other hand, the uniform-in-$\s$ bounds~\eqref{eq:D7:U:Sigma:h} imply uniform bounds along the foliation of $\mathcal{P}$ with the (cylindrical) level sets $\{(x_1,x_2,t) \colon \mathcal{J}(x_2,t) = 1 -\tfrac{\s}{\eps} \} = \{(x_1,x_2,\mathfrak{q}^{-1}(x_2,\s))\}$ for $\s\in[0,\eps]$. That is, we have
\begin{align}
\label{eq:nbs:7:shock:dev:x:t:sup}
&\sup_{\s\in[0,\eps]}
(1-\tfrac{\s}{\eps})^{\frac 14} \|\Jgh \nb^7 \Sigma(x_1,x_2,\mathfrak{q}^{-1}(x_2,\s))\|_{L^2_{x}}
+
\sup_{\s\in[0,\eps]}
(1-\tfrac{\s}{\eps})^{\frac 14} \|\Jgh \nb^7 U(x_1,x_2,\mathfrak{q}^{-1}(x_2,\s))\|_{L^2_{x}}
\notag\\
&\qquad 
+
\tfrac 1\eps
\sup_{\s\in[0,\eps]}
\|\Jgh \nb^7 h(x_1,x_2,\mathfrak{q}^{-1}(x_2,\s))\|_{L^2_{x}}
\les 
\eps^{\frac 12} \mathsf{K} \brak{\mathsf{B}_6}
\,,
\end{align}
\end{subequations}
where as usual $\nb = (\eps \p_t , \eps x_1, x_2)$.

Bounds similar to \eqref{eq:nbs:7:shock:dev:x:t} hold also for the downstream development considered in Section~\ref{sec:downstreammaxdev}. The $L^2_{x,t}$ bounds analogous to \eqref{eq:nbs:7:shock:dev:x:t:2} (with $\mathcal{J}$ replaced by $\JJ$) hold over the spacetime $\mPds$ defined in \eqref{eq:spacetime:P}, while the bound analogous to \eqref{eq:nbs:7:shock:dev:x:t:sup} holds with $\mathfrak{q}^{-1}$ being replaced by $\qds^{-1}$.  Similarly, $L^2_{x,t}$ bounds for the upstream development considered in Section~\ref{sec:upstreammaxdev} hold over the spacetime $\Hdm$ defined in \eqref{eq:spacetime-Theta}.
\end{remark}

\begin{proof}[Proof of Proposition~\ref{prop:D7:U:Sigma:h}]
For simplicity, we only state here the bounds involving the weight function $\varphi = \mathcal{J}$, as defined in \eqref{eq:fake:Jg:def}, and which is used for shock formation, in Sections~\ref{sec:formation:setup}--\ref{sec:sixth:order:energy}. These bounds are stated in~\eqref{eq:nbs:7:Sigma},~\eqref{eq:nbs:7:U}, and~\eqref{eq:nbs:7:h} below. The same bounds hold for the downstream development discussed in Section~\ref{sec:downstreammaxdev}, upon replacing the weight $\mathcal{J}$ appearing in \eqref{eq:nbs:7:Sigma},~\eqref{eq:nbs:7:U}, and~\eqref{eq:nbs:7:h}, with the weight $\JJ$ defined in~\eqref{eq:fake:Jg:def:ds}, and  the same bounds hold for the upstream  development discussed in Section~\ref{sec:upstreammaxdev}, upon replacing the weight $\mathcal{J}$ appearing in \eqref{eq:nbs:7:Sigma},~\eqref{eq:nbs:7:U}, and~\eqref{eq:nbs:7:h}, with the weight $\JJ$, defined in~\eqref{JJ-def-plus-t} and~\eqref{JJ-def-minus-t}.

\subsection{Bounds for $\nbs^7 \Sigma$}
We note that the bounds~\eqref{eq:Sigma:H6:new} and~\eqref{eq:Sigma:H6:new:bdd} below, provide bounds for $\nbs^6 \Sigma$:
\begin{equation}
\eps^{\frac 12} 
\snorm{\Jgh \nbs^6 \Sigma}_{L^\infty_{\s} L^2_{x}} 
+
\snorm{\nbs^6 \Sigma}_{L^2_{\s}L^2_x} 
\les 
\eps \brak{\mathsf{B}_6} \,.
\label{eq:nbs:6:Sigma}
\end{equation}
In order to estimate $\nbs^7 \Sigma$, we recall from~\eqref{Sigma0-ALE} and~\eqref{grad-Sigma} that
\begin{align*} 
\tfrac{1}{\eps} \nbs_1 \Sigma &= \tfrac{1}{2} (\Jg\Wbn - \Jg\Zbn)  + \tfrac{1}{2} \Jg \nbs_2 h (\Wbt -\Zbt) 
\,,  
\\
\nbs_2 \Sigma &= \tfrac{1}{2} g^{\frac{1}{2}} (\Wbt-\Zbt) \\
\tfrac{1}{\eps} \nbs_\s \Sigma &= - \alpha \Sigma (\Zbn  + \Abt) -  \tfrac{1}{2}  V  g^{\frac{1}{2}} (\Wbt-\Zbt)
\,. 
\end{align*}
The above identities may be combined with the definitions~\eqref{eq:tilde:E6} and~\eqref{eq:tilde:D6}, the bootstraps~\eqref{bootstraps-Dnorm:6}, the bounds~\eqref{geometry-bounds-new} for the geometry, the improved estimates in~\eqref{eq:Jg:Zbn:D5:improve} and~\eqref{eq:madman:2:all}, the product bounds in~\eqref{eq:Lynch:1} and~\eqref{eq:Lynch:1:bdd} and the $L^\infty_{x,\s}$ bounds in~\eqref{bootstraps}, to yield
\begin{equation}
\label{eq:nbs:7:Sigma}
\eps^{\frac 12}   
\snorm{\mathcal{J}^{\frac 34} \Jgh \nbs^7 \Sigma}_{L^\infty_{\s} L^2_x} 
+ 
\snorm{\mathcal{J}^{\frac 14} \Jgh \nbs^7 \Sigma}_{L^2_\s L^2_x} 
\les 
\eps \mathsf{K} \brak{\mathsf{B}_6}
\,.
\end{equation}
When compared to the $\nbs^6 \Sigma(x,\s)$ bounds in~\eqref{eq:nbs:6:Sigma}, we note that the $\nbs^7 \Sigma(x,\s)$ estimates obtained in \eqref{eq:nbs:7:Sigma} come at the penalty of additional weights of $\mathcal{J}^{\frac 34}$, and respectively $\Jgh \mathcal{J}^{\frac 14}$, as is natural given our bootstrap assumptions.

\subsection{Bounds for $\nbs^7 U$}
First we note that bounds for $\nbs^6 U$ are available through estimates for $\nbs^6 (W,Z,A,\Jg, \nbs_2 h)$. More precisely, from~\eqref{component-identities} and \eqref{p1-n-tau} we deduce
\begin{align*}
\Jg \Wbn 
&= \tfrac{1}{\eps} \nbs_1 W - \nbs_2 h g^{-\frac 12}  \Jg \nbs_2 W + g^{-\frac 12} A \nbs_2  \Jg  
\,,
\\
\Jg \Zbn 
&= \tfrac{1}{\eps} \nbs_1 Z - \nbs_2 h g^{-\frac 12}  \Jg \nbs_2 Z + g^{-\frac 12} A \nbs_2 \Jg
\,,
\\
\Jg \Abn 
&=  \tfrac{1}{\eps} \nbs_1 A - \nbs_2 h g^{-\frac 12} \Jg \nbs_2 A - \tfrac 12 g^{-\frac 12} (W+Z) \nbs_2 \Jg
\,,
\\
g^{\frac 12} \Wbt 
&= \nbs_2 W + A g^{-1} \nbs_2^2 h
\,,
\\
g^{\frac 12} \Zbt 
&= \nbs_2 Z + A g^{-1} \nbs_2^2 h
\,,
\\
g^{\frac 12} \Abt 
&= \nbs_2  A - \tfrac 12 (W+Z) g^{-1} \nbs_2^2 h
\,.
\end{align*}
In the six identities above, we add the equations for $\Wb$ and $\Zb$, resulting in the four identities
\begin{align*}
\nn \cdo \p_1 U - g^{-\frac 12} \nbs_2 h \Jg \nn \cdo \nbs_2 U 
&= \tfrac 12 \bigl( \Jg \Wbn + \Jg \Zbn) 
\,\\
\tt \cdo \p_1 U - g^{-\frac 12} \nbs_2 h \Jg \tt \cdo \nbs_2 U 
&= \Jg \Abn
\,,\\
\nn \cdo \nbs_2 U 
&= \tfrac 12 g^{\frac 12} \bigl( \Wbt + \Zbt \bigr) 
\,,\\
\tt \cdot \nbs_2 U
&= g^{\frac 12} \Abt
\,.
\end{align*}
By substitution of the last two equations into the first two, and by applying $\nbs^k$ to both sides for $0 \leq k \leq 6$, we obtain 
\begin{subequations}
\label{eq:U:elliptic:0}
\begin{align}
\tfrac{1}{\eps} \nn \cdo \nbs_1 \nbs^k U 
&= \tfrac 12 \nbs^k \bigl( \Jg \Wbn + \Jg \Zbn) 
+ \tfrac 12  \nbs^k\bigl(  \nbs_2 h \Jg ( \Wbt + \Zbt) \bigr) 
- \tfrac{1}{\eps} \jump{\nbs^k,\nn } \cdot \nbs_1 U
\,
\label{eq:nn:cdot:nbs:1:U}\\
\tfrac{1}{\eps} \tt \cdo \nbs_1 \nbs^k  U
&= \nbs^k (\Jg \Abn)
+  \nbs^k  \bigl( \nbs_2 h \Jg \Abt\bigr)
- \tfrac{1}{\eps} \jump{\nbs^k,\tt} \cdot \nbs_1 U
\,,
\label{eq:tt:cdot:nbs:1:U}\\
\nn \cdo \nbs_2 \nbs^k U 
&= \tfrac 12 \nbs^k \bigl( g^{\frac 12} ( \Wbt + \Zbt) \bigr) 
\label{eq:nn:cdot:nbs:2:U}
- \jump{\nbs^k,\nn} \cdot \nbs_2 U
\,,\\
\tt \cdot \nbs_2 \nbs^k U
&=\nbs^k \bigl( g^{\frac 12} \Abt \bigr) 
- \jump{\nbs^k ,\tt} \cdot \nbs_2 U
\label{eq:tt:cdot:nbs:2:U}
\,.
\end{align}
The above four identities need to be supplemented for an identity for computing $\nn \cdo \nbs_\s \nbs^k U$ and $\tt \cdo\nbs_\s \nbs^k U$. We obtain the later by first noting that $\nbs_\s = \eps \bigl( (\Q\p_\s+V \p_2) - V\nbs_2\bigr)$, and then recalling that \eqref{U0-ALE} gives 
\begin{align*}
\nn \cdo (\Q \p_\s + V \p_2) U
&= 
\alpha\Sigma \Zbn
\,,\\
\tt \cdo (\Q \p_\s + V \p_2) U
&= 
\alpha \Sigma   \Abn - \tfrac{\alpha}{2} \Sigma (\Wbt - \Zbt)
\,.
\end{align*}
The two identities above may be differentiated by applying $\nbs^k$, and use \eqref{good-comm} to yield
\begin{align}
\nn \cdo (\Q \p_\s + V \p_2) \nbs^k U
&= 
\alpha \nbs^k \bigl( \Sigma \Zbn\bigr)
- \nn \cdot \jump{\nbs^k,V} \nbs_2 U
- \jump{\nbs^k,\nn} \cdot \bigl(\tfrac{1}{\eps} \nbs_\s U + V \nbs_2 U \bigr)
\,,
\label{eq:nn:cdot:nbs:s:U}\\
\tt \cdo (\Q \p_\s + V \p_2) \nbs^k U
&= 
\tfrac{\alpha}{2} \nbs^k \bigl(\Sigma (2 \Abn - \Wbt + \Zbt) \bigr)
- \tt \cdot \jump{\nbs^k,V} \nbs_2 U
- \jump{\nbs^k,\tt} \cdot \bigl(\tfrac{1}{\eps} \nbs_\s U + V \nbs_2 U \bigr)
\label{eq:tt:cdot:nbs:s:U}
\,.
\end{align}
\end{subequations}
Considering the six identities in \eqref{eq:U:elliptic:0} with $k=5$ together with the formula $\nbs_\s = \eps \bigl( (\Q\p_\s+V \p_2) - V\nbs_2\bigr)$, and using the definitions~\eqref{eq:tilde:E6} and~\eqref{eq:tilde:D6}, the bootstraps~\eqref{bootstraps-Dnorm:5}, the bounds~\eqref{geometry-bounds-new} for the geometry, the improved estimates in~\eqref{eq:Jg:Abn:D5:improve}, \eqref{eq:Jg:Zbn:D5:improve} and~\eqref{eq:madman:2:all}, the product bounds in~\eqref{eq:Lynch:1} and~\eqref{eq:Lynch:1:bdd} and the $L^\infty_{x,\s}$ bounds in~\eqref{bootstraps}, the commutator bounds in~\eqref{eq:Lynch:2}, \eqref{eq:Lynch:3}, and~\eqref{eq:Lynch:3:bdd}, and with Remark~\ref{rem:B5:B6}, we deduce  that $\nbs^6 U$ satisfies the estimates
\begin{equation}
\eps^{\frac 12} 
\snorm{\Jgh \nbs^6 U }_{L^\infty_{\s} L^2_{x}} 
+
\snorm{\nbs^6 U}_{L^2_{\s}L^2_x} 
\les 
\eps \mathsf{K} \brak{\mathsf{B}_6} \,.
\label{eq:nbs:6:U}
\end{equation}
The above estimate is in direct analogy with~\eqref{eq:nbs:6:Sigma}.

Next, we aim to bound $\nbs^7 U$ in a suitably weighted space. To achieve this, we consider the six identities in \eqref{eq:U:elliptic:0} with $k=6$. Here, it is convenient to replace $\nbs_1$ by $\eps \p_1$, and $\nbs_2$ by $\p_2$ plus a material derivative term, via the identity 
\begin{equation}
\p_2  
= \nbs_2 + \tfrac{1}{\eps} \Qb_2 \Qd^{-1} \nbs_\s 
= \tfrac{1}{1+\Qb_2 V \Q^{-1}}
\bigl(\nbs_2 + \Qb_2 \Q^{-1} (\Q \p_\s + V \p_2) \bigr)
\label{eq:p:2:nbs:2}
\,.
\end{equation}
Using \eqref{eq:U:elliptic:0} with $k=6$ and \eqref{eq:p:2:nbs:2}, we may thus derive
\begin{subequations}
\label{eq:U:elliptic:1}
\begin{align}
\nn \cdo \p_1 \nbs^6 U 
&= \mathcal{M}_{\nnn,1}
\\
\tt \cdo \p_1 \nbs^6  U
&= \mathcal{M}_{\ttt,1}
\\
\nn \cdo \p_2 \nbs^6 U 
&=
\tfrac{1}{1+\Qb_2 V \Q^{-1}}
\bigl( \mathcal{M}_{\nnn,2} + \Qb_2 \Q^{-1} \mathcal{M}_{\nnn,\s} \bigr)
\\
\tt \cdot \p_2 \nbs^6 U
&= 
\tfrac{1}{1+\Qb_2 V \Q^{-1}}
\bigl( \mathcal{M}_{\ttt,2} + \Qb_2 \Q^{-1} \mathcal{M}_{\ttt,\s} \bigr)
\,,
\end{align}
\end{subequations}
where we have denoted the emerging six forcing terms by
\begin{subequations}
\label{eq:U:elliptic:1a}
\begin{align}
\mathcal{M}_{\nnn,1}
&:= \tfrac 12 \nbs^6 \bigl( \Jg \Wbn + \Jg \Zbn) 
+ \tfrac 12  \nbs^6 \bigl(  \nbs_2 h \Jg ( \Wbt + \Zbt) \bigr) 
- \tfrac{1}{\eps} \jump{\nbs^6,\nn } \cdot \nbs_1 U
\\
\mathcal{M}_{\ttt,1}
&:=
\nbs^6 (\Jg \Abn)
+  \nbs^6  \bigl( \nbs_2 h \Jg \Abt\bigr)
- \tfrac{1}{\eps} \jump{\nbs^6,\tt} \cdot \nbs_1 U
\\
\mathcal{M}_{\nnn,2}
&:=
\tfrac 12 \nbs^6 \bigl( g^{\frac 12} ( \Wbt + \Zbt) \bigr) 
- \jump{\nbs^6,\nn} \cdot \nbs_2 U
\\
\mathcal{M}_{\ttt,2}
&:= 
\nbs^6 \bigl( g^{\frac 12} \Abt \bigr) 
- \jump{\nbs^6 ,\tt} \cdot \nbs_2 U
\\
\mathcal{M}_{\nnn,\s}
&:=
\alpha \nbs^6 \bigl( \Sigma \Zbn\bigr)
- \nn \cdot \jump{\nbs^6,V} \nbs_2 U
- \jump{\nbs^6,\nn} \cdot \bigl(\tfrac{1}{\eps} \nbs_\s U + V \nbs_2 U \bigr)
\\
\mathcal{M}_{\ttt,\s}
&:= 
\tfrac{\alpha}{2} \nbs^6 \bigl(\Sigma (2 \Abn - \Wbt + \Zbt) \bigr)
- \tt \cdot \jump{\nbs^6,V} \nbs_2 U
- \jump{\nbs^6,\tt} \cdot \bigl(\tfrac{1}{\eps} \nbs_\s U + V \nbs_2 U \bigr)
\,.
\end{align}
\end{subequations}
Next, using that $e_1 = g^{-\frac 12} (\nn + \nbs_2 h  \tt)$ and $e_2 = g^{-\frac 12} (\tt - \nbs_2 h \nn)$, we may perform linear combinations of the terms in \eqref{eq:U:elliptic:1}, finally arriving at
\begin{subequations}
\label{eq:U:elliptic:2}
\begin{align}
\p_1 \nbs^6  U^1
&= g^{-\frac 12}  \bigl( \mathcal{M}_{\nnn,1} +   \nbs_2 h \mathcal{M}_{\ttt,1}\bigr)
\,,
\\
\p_1  \nbs^6 U^2
&=g^{-\frac 12} \bigl( \mathcal{M}_{\ttt,1} - \nbs_2 h \mathcal{M}_{\nnn,1}\bigr)
\,,
\\
\p_2  \nbs^6  U^1
&=
\tfrac{g^{-\frac 12}}{1+\Qb_2 V \Q^{-1}}
\bigl( 
(\mathcal{M}_{\nnn,2} + \nbs_2 h \mathcal{M}_{\ttt,2})
+
\Qb_2 \Q^{-1} 
(\mathcal{M}_{\nnn,\s} + \nbs_2 h \mathcal{M}_{\ttt,\s})
\bigr)
\,,
\\
\p_2  \nbs^6 U^2
&=
\tfrac{g^{-\frac 12}}{1+\Qb_2 V \Q^{-1}}
\bigl( 
(\mathcal{M}_{\ttt,2} - \nbs_2 h \mathcal{M}_{\nnn,2})
+
\Qb_2 \Q^{-1} 
(\mathcal{M}_{\ttt,\s} - \nbs_2 h \mathcal{M}_{\nnn,\s})
\bigr)
\, .
\end{align}
\end{subequations}
In order to use \eqref{eq:U:elliptic:2}, we thus must bound the right side in suitably-weighted $L^2_x$ spaces, both in $L^\infty_\s$ (with weight $\mathcal{J}^{\frac 34} \Jgh$) and in $L^2_\s$ (with weight $\mathcal{J}^{\frac 14} \Jgh$). For this purpose, we record the bounds 
\begin{subequations}
\label{eq:U:elliptic:2a}
\begin{align}
\eps^{\frac 12} \snorm{\mathcal{J}^{\frac 34} \Jgh \mathcal{M}_{\nnn,1}}_{L^\infty_\s L^2_x}
+ \snorm{\mathcal{J}^{\frac 14} \Jgh \mathcal{M}_{\nnn,1}}_{L^2_\s L^2_x}
&\les \brak{\mathsf{B}_6}
\,,\\
\eps^{\frac 12} \snorm{\mathcal{J}^{\frac 34} \Jgh \mathcal{M}_{\ttt,1}}_{L^\infty_\s L^2_x}
+ \snorm{\mathcal{J}^{\frac 14} \Jgh \mathcal{M}_{\ttt,1}}_{L^2_\s L^2_x}
&\les \eps \mathsf{K} \brak{\mathsf{B}_6}
\,,\\
\eps^{\frac 12} \snorm{\mathcal{J}^{\frac 34} \Jgh \mathcal{M}_{\nnn,2}}_{L^\infty_\s L^2_x}
+ \snorm{\mathcal{J}^{\frac 14} \Jgh \mathcal{M}_{\nnn,2}}_{L^2_\s L^2_x}
&\les \eps \mathsf{K} \brak{\mathsf{B}_6}
\,,\\
\eps^{\frac 12} \snorm{\mathcal{J}^{\frac 34} \Jgh \mathcal{M}_{\ttt,2}}_{L^\infty_\s L^2_x}
+ \snorm{\mathcal{J}^{\frac 14} \Jgh  \mathcal{M}_{\ttt,2}}_{L^2_\s L^2_x}
&\les \eps \mathsf{K} \brak{\mathsf{B}_6}
\,,\\
\eps^{\frac 12} \snorm{\mathcal{J}^{\frac 34} \Jgh  \mathcal{M}_{\nnn,\s}}_{L^\infty_\s L^2_x}
+ \snorm{\mathcal{J}^{\frac 14} \Jgh \mathcal{M}_{\nnn,\s}}_{L^2_\s L^2_x}
&\les \mathsf{K} \brak{\mathsf{B}_6}
\,,\\
\eps^{\frac 12} \snorm{\mathcal{J}^{\frac 34} \Jgh \mathcal{M}_{\ttt,\s}}_{L^\infty_\s L^2_x}
+ \snorm{\mathcal{J}^{\frac 14} \Jgh \mathcal{M}_{\ttt,\s}}_{L^2_\s L^2_x}
&\les \eps \mathsf{K} \brak{\mathsf{B}_6}
\,,
\end{align}
\end{subequations}
which follow from the definitions~\eqref{eq:tilde:E6} and~\eqref{eq:tilde:D6}, the bootstrap~\eqref{bootstraps-Dnorm:5}, the bounds~\eqref{geometry-bounds-new} for the geometry, the improved estimates in~\eqref{eq:Jg:Abn:D5:improve}, \eqref{eq:Jg:Zbn:D5:improve}, and~\eqref{eq:madman:2:all}, the product bounds in~\eqref{eq:Lynch:1} and~\eqref{eq:Lynch:1:bdd},  the $L^\infty_{x,\s}$ bounds in~\eqref{bootstraps}, the commutator bounds in~\eqref{eq:Lynch:2}, \eqref{eq:Lynch:3}, and~\eqref{eq:Lynch:3:bdd}, and from the identities in~\eqref{eq:U:elliptic:0} to compute commutator terms.

In turn, the bootstraps~\eqref{bs-h}, \eqref{bs-V}, the bounds~\eqref{eq:Qb:bbq}, \eqref{eq:Q:bbq}, and estimates~\eqref{eq:U:elliptic:2a}, imply via \eqref{eq:U:elliptic:2} that 
\begin{equation}
\label{eq:nabla:D6:U:a}
\eps^{\frac 12} \snorm{\mathcal{J}^{\frac 34} \Jgh (\eps \p_1 , \p_2) \nbs^6 U}_{L^\infty_\s L^2_x}
+ \snorm{\mathcal{J}^{\frac 14} \Jgh (\eps \p_1 , \p_2) \nbs^6  U}_{L^2_\s L^2_x}
\les \eps \mathsf{K} \brak{\mathsf{B}_6} 
\,.
\end{equation}
It thus remains to convert the $\nabla \nbs^6$ bound in \eqref{eq:nabla:D6:U:a} to a full $\nbs^7 U$ estimate. 
This is achieved by appealing to \eqref{eq:nn:cdot:nbs:s:U}, which may be rewritten as $\nn \cdo (\Q \p_\s + V \p_2) \nbs^k U = \mathcal{M}_{\nnn,\s}$, and \eqref{eq:tt:cdot:nbs:s:U}, which may be rewritten as $\tt \cdo (\Q \p_\s + V \p_2) \nbs^k U = \mathcal{M}_{\ttt,\s}$. Together with~\eqref{eq:U:elliptic:2a}, these identities yield
\begin{equation}
\label{eq:nabla:D6:U:b}
\eps^{\frac 12} \snorm{\mathcal{J}^{\frac 34} \Jgh (\Q \p_\s + V \p_2) \nbs^6 U}_{L^\infty_\s L^2_x}
+ \snorm{\mathcal{J}^{\frac 14} \Jgh (\Q \p_\s + V \p_2) \nbs^6  U}_{L^2_\s L^2_x}
\les  \mathsf{K} \brak{\mathsf{B}_6} 
\,.
\end{equation}
The last ingredient is~\eqref{nb-s}, which may be rewritten as
\begin{equation*}
\nbs_1 = \eps \p_1\,,
\quad
\nbs_{2} = \bigl(1+\tfrac{V \Qb_2}{\Q}\bigr) \p_2 - \tfrac{\Qb_2}{\Q} (\Q\p_\s+V\p_2) \,,
\quad 
\nbs_{\s} = \eps \tfrac{\Qd}{\Q}  (\Q \p_\s + V\p_2) - \eps \tfrac{V  \Qd}{\Q}   \p_2 \,.
\end{equation*}
The above three identities, together with \eqref{eq:nabla:D6:U:a}, \eqref{eq:nabla:D6:U:b}, and Lemma~\ref{lem:Q:bnds}, give
\begin{equation}
\label{eq:nbs:7:U}
\eps^{\frac 12} 
\snorm{\mathcal{J}^{\frac 34} \Jgh  \nbs^7 U }_{L^\infty_{\s} L^2_x} 
+ 
\snorm{\mathcal{J}^{\frac 14} \Jgh \nbs^7 U}_{L^2_\s L^2_x} 
\les 
\eps \mathsf{K} \brak{\mathsf{B}_6}
\,.
\end{equation}
The above estimate is in direct analogy with~\eqref{eq:nbs:7:Sigma}.

\subsection{Bounds for $\nbs^7 h$}
We recall from \eqref{psi-def} that the ALE map $\psi$ is given by $h(x_1,x_2,t) e_1 + x_2 e_2$. As such, the regularity of $\psi$ is equivalent to the regularity of the map $h$. 
So far, the best available estimates for $h$ are \eqref{D6h2Energy:new} for $\nbs^6 \nbs_2 h$, and \eqref{D6h1Energy:new} for $\nbs^6 \nbs_1 h$. To obtain a full control of $\nbs^7 h$, it thus remains to estimate $\nbs_\s^7 h$. For this purpose, we recall from \eqref{eq:psi:evo:def}, written in $(x,\s)$ variables, that 
$\tfrac{1}{\eps} \p_\s h = g^{\frac 12}(U \cdot \nn + \Sigma)$. Upon applying $\nbs_\s^6$ to this expression, we determine that 
\begin{equation*}
\nbs_\s^7  h
= \eps \nbs_\s^6 \bigl( g^{\frac 12}(U \cdot \nn + \Sigma) \bigr)
\,.
\end{equation*}
Using the bootstraps~\eqref{bootstraps}, the inequality $1\leq \mathcal{J}^{-\frac 14}$, the geometric bounds in~\eqref{geometry-bounds-new}, the $\nbs^6 \Sigma$ bounds in~\eqref{eq:nbs:6:Sigma}, the $\nbs^6 U$ bounds in ~\eqref{eq:nbs:6:U}, and the Moser-type bounds in Lemmas~\ref{lem:Moser:tangent}, we deduce from the above identity that
\begin{subequations}
\begin{equation}
\snorm{\nbs_\s^7  h}_{L^2_{\s}L^2_x} 
\les 
\eps^2 \mathsf{K} \brak{\mathsf{B}_6} \,.
\label{eq:nbs:7:h:a}
\end{equation}
By additionally using the inequalities $\mathcal{J}^{\frac 14} \leq \Jg^{\! \frac 14} \leq \Jgh \les 1$, \eqref{table:derivatives}, \eqref{eq:sup:in:time:L2}, and \eqref{se3:time}, 
we have that 
\begin{equation}
\snorm{\mathcal{J}^{\frac 14} \nbs_\s^7  h}_{L^\infty_{\s} L^2_{x}} 
\les \eps^{\frac 32} \mathsf{K} \brak{\mathsf{B}_6}
+ 
\eps^{\frac 12} \displaystyle{\sum}_{k=1}^{6} 
\| \nbs_\s^{7-k} (g^{\frac 12} \nn) \cdo \nbs_\s^k U\|_{L^2_{x,\s}}
+
\|\nbs_\s^{7-k} g^{\frac 12}  \cdo \nbs_\s^k \Sigma\|_{L^2_{x,\s}}
\les 
\eps^{\frac 32} \mathsf{K} \brak{\mathsf{B}_6} 
\,.
\label{eq:nbs:7:h:b}
\end{equation}
\end{subequations}

To summarize, we combine \eqref{D6h2Energy:new}, \eqref{D6h1Energy:new}, \eqref{eq:nbs:7:h:a}, and \eqref{eq:nbs:7:h:b}, and deduce that 
\begin{equation}
\eps^{\frac 12} 
\snorm{\mathcal{J}^{\frac 14} \nbs^7  h}_{L^\infty_{\s} L^2_{x}} 
+
\snorm{\nbs^7  h}_{L^2_{\s}L^2_x} 
\les 
\eps^2 \mathsf{K} \brak{\mathsf{B}_6}
\,,
\label{eq:nbs:7:h}
\end{equation}
which gives the optimal regularity bounds for $h$, and thus also for the ALE map $\psi$.
\end{proof}



\appendix

\section{Maximal globally hyperbolic development in a box for the 1D Euler equations}\label{sec:usersguide}

The purpose of this appendix is to showcase the process of {\em shock formation from smooth initial data}, the  spacetime of \MGHDB, and to discuss the physical {\em shock development problem} in the context of the 1D Euler equations (cf.~\eqref{euler0} for $d=1$). While these results are surely known in one space dimension, we were not able to find a suitable exposition of these concepts in the literature.\footnote{After completing this work, we were made aware of the paper \cite{AbSp2023} by Abbrescia \& Speck, which 
constructs the \MGHDB\ for 1D relativistic Euler.}
Additionally, the discussion in this appendix highlights some of the main concepts discussed in the multi-dimensional setting of this paper, but without having to deal with a number of geometric difficulties. Lastly, it is easier to draw accurate images in $1+1$ dimensions (space \& time), so that the reader can build part of the necessary intuition required for the multi-dimensional setting.  

\subsection{The setup: equations, variables, initial data}
In one space dimension, the system~\eqref{euler0} becomes
\begin{subequations}
\label{eq:Euler:1D:conservation}
\begin{align}
 \partial_t (\rho u) + \partial_y( \rho u^2) + \partial_y p &=0\,,
 \label{eq:Euler:1D:conservation:a}\\
 \partial_t \rho + \partial_y (\rho u) &=0 \,,
  \label{eq:Euler:1D:conservation:b}\\
 \partial_t E + \partial_y (u (E+p)) &=0
  \label{eq:Euler:1D:conservation:c}\,,
\end{align}
\end{subequations}
where $E = \tfrac{1}{\gamma-1} p+ \tfrac 12 \rho u^2$ is the energy, and $\gamma>1$ is the adiabatic exponent. Here the space variable is $y\in \TT$, and \eqref{eq:Euler:1D:conservation} is supplemented with smooth Cauchy data $(\rho_0,u_0,E_0)$ at the initial time $t= \initial$.  

While the formulation \eqref{eq:Euler:1D:conservation}, as a system of conservation laws, is necessary in order to read off the correct Rankine-Hugoniot jump conditions in the presence of a shock singularity, when the Euler system~\eqref{eq:Euler:1D:conservation} is initiated with smooth initial data, an equivalent, more symmetric formulation is useful. For this purpose we remark that at least until the emergence of shocks,  the  evolution of the energy $E$ in~\eqref{eq:Euler:1D:conservation:c} may be replaced by the transport equation 
\begin{equation}
\partial_t S + u \partial_y S = 0\,, 
\label{eq:1D:entropy:transport}
\end{equation}
where $S$ is the specific entropy, defined as $ S = \log(\tfrac{\gamma p}{\rho^\gamma})$. As such, if the initial datum $S_0$ is a constant (which for simplicity we may take as zero), then $S$ remains a constant under the above transport equation,  as long as $u \in L^1_t W^{1,\infty}_x$ (prior to the formation of singularities). Since the \MGHDB\ spacetime does not contain any singular behavior in its interior, it is convenient for the sake of presentation to confine the discussion to the {\em isentropic} case, in which $S_0 \equiv 0$, and thus $S \equiv 0$ prior to the formation of shocks.\footnote{See the analysis in~\cite{NeRiShVi2023} for a parallel discussion in the case that $S_0 \not \equiv 0$.}

As discussed earlier in~\eqref{euler1}, in this isentropic case the pressure law becomes $p = \tfrac{1}{\gamma} \rho^\gamma$, and upon letting $\alpha = \tfrac{\gamma-1}{2}$, and denoting the rescaling sound speed by $\sigma= \tfrac{1}{\alpha} c = \tfrac{1}{\alpha} \rho^\alpha$, the conservation laws in~\eqref{eq:Euler:1D:conservation:a}--\eqref{eq:Euler:1D:conservation:b} may be written in a more symmetric way as
\begin{subequations}
\label{eq:Euler:1D:u:sigma}
\begin{align}
\partial_t u + u\p_y u + \alpha \sigma \p_y \sigma &=0 
\label{eq:Euler:1D:u:sigma:a}\,,\\
\partial_t \sigma + u \p_y \sigma + \alpha \sigma \p_y u &=0 
\label{eq:Euler:1D:u:sigma:b}\,.
\end{align}
\end{subequations}
The isentropic dynamics in~\eqref{eq:Euler:1D:u:sigma}  consists in fact of two coupled transport equations for the Riemann variables
\begin{equation*}
w = u + \sigma\,, \qquad z = u - \sigma \,.
\end{equation*}
The dominant Riemann variable $w$ is transported along the fast characteristic speed $\lambda_3$, while the  subdominant Riemann variable is transported along the slow characteristic speed $\lambda_1$. That is, 
\begin{subequations}
\label{eq:1D:Euler:Riemann}
\begin{align}
(\partial_t + \lambda_3 \partial_y) w &= 0
\label{eq:1D:Euler:Riemann:a}\,,
\\
(\partial_t + \lambda_1 \partial_y) z &= 0
\label{eq:1D:Euler:Riemann:b}\,,
\end{align}
where, rewriting~\eqref{wave-speeds} in one space dimension, we have denoted for $i\in\{1,2,3\}$ the wave speeds
\begin{equation}
\lambda_i = u + (i-2) \alpha \sigma = u + (i-2) c = \tfrac{1+(i-2) \alpha}{2} w + \tfrac{1-(i-2) \alpha}{2} z \,.
\end{equation}
\end{subequations}
The system~\eqref{eq:1D:Euler:Riemann} is equivalent to~\eqref{eq:Euler:1D:u:sigma}.

For consistency with the rest of the paper, we supplement \eqref{eq:1D:Euler:Riemann} with initial data $(w_0,z_0) = (u_0 + \sigma_0,u_0-\sigma_0)$ that has many of the same properties as the multi-dimensional data discussed in Section~\ref{cauchydata} (compressive and generic), while at the same aiming for simplicity of the presentation.  For definiteness, we consider initial data $w_0$ which is $\OO(1)$ and has an $\OO(\frac{1}{\eps})$ most negative slope, and $z_0$ which is $\OO(\eps)$ and has an $\OO(1)$ slope. For simplicity, let 
\begin{equation}
w_0(x) = 1 - \tfrac 12 \operatorname{sin}(\tfrac{2 x}{\eps})\,,
\qquad
z_0(x) = \tfrac{\eps}{2} \cos(\tfrac{2x}{\eps}) \,,
\label{eq:1D:Euler:IC}
\end{equation} 
where $\eps\ll 1$ is a small parameter. In particular, $\tfrac 18 \leq \sigma_0(x) \leq 1$, $w_0^\prime$ attains a global minimum of $-\tfrac 1\eps$ at $x=0$, and this minimum is non-degenerate, as $w_0^{\prime\prime\prime}(0) = \tfrac{4}{\eps^3}>0$.

For consistency with~\eqref{tin}, let us assume that the initial time is given by $\initial = - \tfrac{2\eps}{1+\alpha}$. We also define $\final = \tfrac{\eps}{1+\alpha}$. Our goal is to study the system~\eqref{eq:1D:Euler:Riemann}, with initial datum~\eqref{eq:1D:Euler:IC}, for $t\in [\initial,\final]$. To avoid any discussion of localization at ``large'' values of $|y|$,  let us only discuss the problem for $|y| \leq \frac{3\pi \eps}{4}$, appealing to the finite speed of propagation.

\subsection{The classical perspective on shock formation}

Shock formation is classically described as the process through which an initial negatively sloped disturbance in the graph of the density progressively steepens, and eventually $\p_y \sigma$ diverges towards $-\infty$. For initial data such as that in~\eqref{eq:1D:Euler:IC}, at the same time the slope of the dominant Riemann variable $w$ also diverges towards $-\infty$. A classical way to see this, following Lax~\cite{Lax1964}, is to differentiate~\eqref{eq:1D:Euler:Riemann:a} and to use that~\eqref{eq:1D:Euler:Riemann:b} may be rewritten as $(\p_t + \lambda_3 \p_y) \sigma = (\lambda_3 - \lambda_1) \p_y z = 2\alpha \sigma \p_y z$, leading to 
\begin{align*}
(\p_t + \lambda_3 \p_y) (\p_y w) 
= - \p_y \lambda_3 \, \p_y w
&= - \tfrac{1+\alpha}{2} (\p_y w)^2 - \tfrac{1-\alpha}{2} (\p_y w) (\p_y z) \notag\\
&=  - \tfrac{1+\alpha}{2} (\p_y w)^2 - \tfrac{1-\alpha}{4\alpha} (\p_y w) (\p_t + \lambda_3 \p_y) \log \sigma
\,. 
\end{align*}
Since vacuum cannot form on $[\initial,\final]$\footnote{To see this, note that for initial data as in~\eqref{eq:1D:Euler:IC}, the maximum principle applied to~\eqref{eq:1D:Euler:Riemann:a} yields $|w(y,t) - 1| \leq \tfrac 12$ and the maximum principle applied to~\eqref{eq:1D:Euler:Riemann:b} yields $|z(y,t)|\leq \tfrac{\eps}{2}$. Hence $\sigma = \tfrac 12 (w+z) \in [\tfrac 18, \tfrac 74]$ when $\eps \leq \frac 12 $.}
on may integrate the above PDE along the characteristics $\eta = \eta(x,t)$ associated to the fast wave-speed $\lambda_3$, i.e., the solution of
\begin{equation}
\label{eq:eta:flow:1D}
\partial_t \eta(x,t) = \lambda_3(\eta(x,t),t)\,, 
\qquad 
\eta(x,\initial) = x\,,
\end{equation}
leading to 
\begin{equation}
\p_t (\p_y w \cir \eta) = - \tfrac{1+\alpha}{2}  (\p_y w \cir \eta)^2 - \tfrac{1-\alpha}{4\alpha} (\p_y w \cir \eta) \p_t (\log \sigma \cir \eta)
\,.
\label{eq:good:ol:Lax}
\end{equation}
With $\log \sigma = \OO(1)$ given, the above ODE may now be directly integrated, leading to a proof that $\p_y w \circ \eta$ blows up towards $-\infty$ in finite time, for the label at which $w_0^\prime$ was the most negative. Since by assumption this label is $x=0$, we deduce that there exists a minimal time $t_* > \initial$ such that $\p_y w(\eta(0,t),t) \to -\infty$ as $t\to t_*$. In order to estimate $t_*$, we may compare the above ODE to the Ricatti ODE $\dot a = - \tfrac{1+\alpha}{2} a^2$, with datum $a(\initial) = -\tfrac 1 \eps$, and deduce that $t_* = \initial + \tfrac{2\eps}{1+\alpha} \pm \OO(\eps^2) = 0 \pm \OO(\eps^2) \in (\initial,\final)$. This is Lax's proof of shock formation for the system~\eqref{eq:1D:Euler:Riemann}, which may be found in~\cite{Lax1964}.  

A more geometric perspective on  shock formation is provided by the method of characteristics (see Figure~\ref{fig:1D:Eulerian:shock} below). Since in 1D we can rule out the formation of vacuum, it may be shown (see~\cite[Proposition 5.4]{NeRiShVi2023}) that the Eulerian blow-up criterion $\int_{\initial}^{t_*} \|\p_y w(\cdot,t)\|_{L^\infty} + \|\p_y \sigma(\cdot,t)\|_{L^\infty} dt = +\infty$, is equivalent to the Lagrangian criterion $\lim_{t\to t_*} \min_x \p_x \eta (x,t) = 0$. That is, shock formation occurs at the first time $t_*$ at which the map $x\mapsto \eta(x,t_*)$ loses strict monotonicity, because the fast characteristic curves $(\eta(x,t),t)$ impinge on each other (for $t>t_*$ the map $x\mapsto \eta(x,t)$ would fail to be injective). The time $t_*$, and the label $x_*$ at which this strict monotonicity is lost may be naturally computed by solving  
\begin{equation}
 \p_x \eta(x_*,t_*) =0\,,
 \qquad 
 \p_{xx} \eta(x_*,t_*) =0\,.
 \label{eta:shock:characterization}
\end{equation}
We note that for the initial data~\eqref{eq:1D:Euler:IC}, at $(x_*,t_*)$ the field $w$ develops a H\"older-$\frac 13$ cusp singularity,\footnote{\label{footnote:16}To see this, let $y = \eta(x,t_*)$. Due to~\eqref{eq:1D:Euler:Riemann:a} we have $w(y,t_*) = w_0(x)$. Thus we need to compose $w_0$ with the inverse flow $x = \eta^{-1}(y,t_*)$. Due to~\eqref{eta:shock:characterization} we have the power series expansion $\eta(x,t_*) = \eta(x_*,t_*) + \tfrac{1}{6} (x-x_*)^3 \p_x^3 \eta(x_*,t_*) + \OO((x-x_*)^4)$. The fact that $\p_x^3 w_0(x_*) > 0$, which is the genericity condition for the data, implies that $\p_x^3 \eta(x_*,t_*) \approx (t_* - \initial) \tfrac{1+\alpha}{2} \p_x^3 w_0(x_*) \approx 4 \eps^{-2} > 0$. Thus when  we are solving for $x$ the equation $y = \eta(x,t_*)$, we are in essence solving the cubic equation $y \approx y_* + \tfrac 23\eps^{-2} (x-x_*)^3$. This results in $x-x_* \approx \eps^{\frac 23} (y-y_*)^{\frac 13}$, and thus $w(y,t_*) \approx w_0(\eps^{\frac 23} (y-y_*)^{\frac 13})$. Since $w_0$ is smooth, this  explainins the H\"older-$\frac 13$ cusp.} not an actual shock-discontinuity, which is why the point $(x_*,t_*)$ is sometimes referred to as a {\em pre-shock}.

Also, we note that for initial datum of the type~\eqref{eq:1D:Euler:IC} the characteristics $\phi= \phi(x,t)$ associated to the slow wave-speed $\lambda_1$, i.e., the solution of 
\begin{equation}
\label{eq:phi:flow:1D}
\p_t \phi(x,t) = \lambda_1(\phi(x,t),t)\,,
\qquad 
\phi(x,\initial) = x\,, 
\end{equation}
remain injective and the map $x\mapsto \phi(x,t)$ remains strictly monotone uniformly for $(x,t)$ with $t\in [\initial,t_*]$. Equivalently, $\|\p_y z (\cdot,t)\|_{L^\infty}$ remains uniformly bounded as $t\to t^*$.  

\begin{figure}[htb!]
\centering
\includegraphics[width=.6\textwidth]{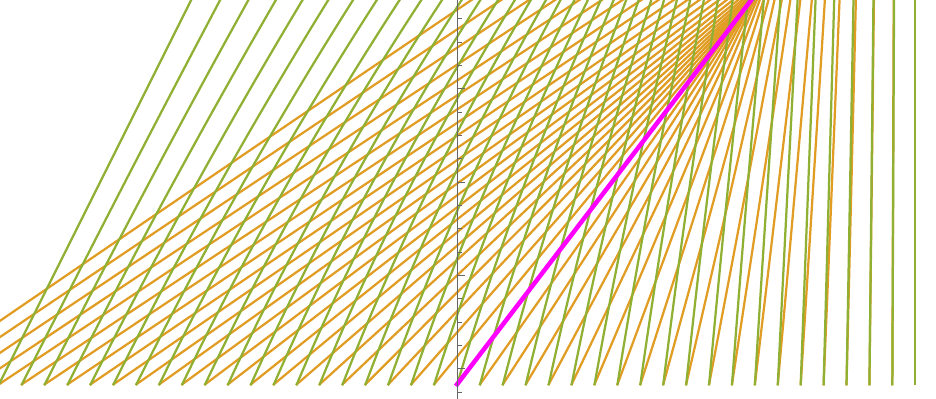};
\vspace{-.1in}
\caption{{\bf Shock formation for 1D Euler via classical characteristics.} Consider the 1D Euler system~\eqref{eq:1D:Euler:Riemann} with initial data~\eqref{eq:1D:Euler:IC} for $\alpha=\frac 12$ and $\eps = \frac{1}{16}$. The bounding box represents $(x,t) \in [-\frac{3\pi \eps}{4},\frac{3\pi\eps}{4}] \times [\initial,t_*]$. 
In orange we represent the ``fast'' characteristics $\eta(x,t)$  which solve~\eqref{eq:eta:flow:1D}, and in olive-green we represent the ``slow'' characteristics $\phi(x,t)$ that solve~\eqref{eq:phi:flow:1D}.
The first blow-up occurs at the Eulerian spacetime location $(y_*,t_*)$, where $y_*= \eta(t_*,x_*)$, and $(x_*,t_*)$ solve~\eqref{eta:shock:characterization}. In this figure $x_*=0$. In magenta we have represented the distinguished $\lambda_3$-characteristic curve leading to the very first singularity, $(\eta(x_*,t),t)_{t\in[\initial,t_*]}$.  }
\label{fig:1D:Eulerian:shock}
\end{figure}

\subsection{A characteristic description of the \MGHDB}
The aforementioned analysis only describes the solution up the time slice of the very first singularity, $\{ t=t_* \}$.  However, due to finite speed of propagation, even for times $t>t_*$ there exist regions of spacetime where a smooth extension of the solution may be uniquely computed. In hyperbolic PDEs it is a natural question to describe the largest such spacetime region, which is free of singularities. The \MGHDB\ spacetime $\mathcal{M}$ for~\eqref{eq:1D:Euler:Riemann} consists of all points $(y,t)$ for which the fields $w(y,t)$ and $z(y,t)$ can be computed in a {\em unique} and {\em smooth} way in terms of the initial data $(w_0,z_0)$, prescribed on the time slice $\{t = \initial\}$. 

In the particular case of~\eqref{eq:1D:Euler:Riemann}, which consists of two coupled transport equations, we have that 
\begin{equation*}
w(\eta(x,t),t) = w_0(x)\,,
\qquad
z(\phi(x,t),t) = z_0(x)\,,
\end{equation*}
assuming that $\eta$ and $\phi$ are well-defined and sufficiently regular. Therefore, the requirement that $w(y,t)$ and $z(y,t)$ may be computed in a unique and smooth fashion in terms of the Cauchy data is in fact the requirement that the backwards flows\footnote{\label{footnote:20}These may be defined either as inverse maps via $\eta(\eta^{-1}(y,t),t)=y$ and $\phi(\phi^{-1}(y,t),t) = y$, or as solutions of the transport PDEs $(\p_t + \lambda_3(y,t) \p_y) \eta = 0$ and $(\p_t + \lambda_1(y,t) \p_y) \phi = 0$, with initial data that is the identity map at time $\initial$.} $\eta^{-1}(y,t)$ and $\phi^{-1}(y,t)$ are smooth, and that 
\begin{equation}
\eta(\eta^{-1}(y,t),t^\prime), \phi(\phi^{-1}(y,t),t^\prime) \in \mathcal{M}\,,
\qquad \mbox{for all} \qquad
(y,t) \in \mathcal{M}\,, t^\prime \in [\initial,t]\,.
\label{eq:flows:remain:in:max:dev}
\end{equation}
Since $w(y,t) = w_0(\eta^{-1}(y,t))$ and $z(y,t) = z_0(\phi^{-1}(y,t))$, the condition in~\eqref{eq:flows:remain:in:max:dev} is indeed a full characterization of the \MGHDB\ spacetime $\mathcal{M}$, for 1D isentropic Euler~\eqref{eq:1D:Euler:Riemann}.

For the compressive and generic initial condition considered here (cf.~\eqref{eq:1D:Euler:IC}), one may give a precise definition of the spacetime $\mathcal{M}$ characterized by \eqref{eq:flows:remain:in:max:dev}. We refer to Figure~\ref{fig:1D:Eulerian} below, where the part of $\mathcal{M}$ which lies in the  bounding box $(x,t) \in [-\frac{3\pi \eps}{4},\frac{3\pi\eps}{4}] \times [\initial,\final]$ is represented:
\begin{itemize}[leftmargin=16pt]

\item It it is clear that that all time slices with $t < t_*$ lie in $\mathcal{M}$ since prior to the very first singularity all functions remain smooth; as such, Figure~\ref{fig:1D:Eulerian} extends the image in Figure~\ref{fig:1D:Eulerian:shock}.

\item Consider labels $x$ such that $x>x_*$. In our setup, where the waves move to the right, we refer to these labels as being {\em downstream} of the pre-shock label. For every $x>x_*$, we may smoothly and uniquely extend the fast acoustic characteristic $\eta$ emanating from $x$ for all times $t\in[\initial, \mathsf{T}_*(x)]$, where the time $\mathsf{T}_*(x)$ is characterized by the condition 
\begin{equation}
\p_x \eta(x,\mathsf{T}_*(x)) = 0 \,.
\label{eta:shock:characterization:2}
\end{equation} 
This is the same condition that we used in~\eqref{eta:shock:characterization} to characterize the very first blow-up time, and indeed $t_* = \mathsf{T}_*(x_*)$. For $t> \mathsf{T}_*(x)$ injectivity is lost and hence we cannot uniquely extend $\eta(x,\cdot)$ after this time. The  curve $\partial_{\sf top}^+\mathcal{M} := (\eta(x,\mathsf{T}_*(x)),  \mathsf{T}_*(x))_{x\geq x_*}$ which emanates from the very first singularity at $(x_*,t_*)$ (represented in red in Figure~\ref{fig:1D:Eulerian}) is smooth, as it is parametrizes the zero level set of the smooth function $\p_x \eta$, and it serves as the future temporal boundary of the ``downstream part'' of the spacetime $\mathcal{M}$. Indeed, it is clear from Figure~\ref{fig:1D:Eulerian} that for every $(y,t)$ which lies either ``below'' or ``to the right'' of the curve $\partial_{\sf top}^+\mathcal{M}$ the backwards trajectories corresponding to the labels $\eta^{-1}(y,t)$ and $\phi^{-1}(y,t)$ remain in $\mathcal{M}$, and therefore they satisfy~\eqref{eq:flows:remain:in:max:dev}. Lastly, we mention that for any point $(\bar{y}, \bar{t}) := (\eta(\bar{x},\mathsf{T}_*(\bar{x})),  \mathsf{T}_*(\bar{x})) \in \partial_{\sf top}^+\mathcal{M}$ with $\bar{x}\geq x_*$,  we have that $\lim_{\mathcal{M} \ni (y,t) \to (\bar{y},\bar{t})} \p_y w(y,t) \to -\infty$. In particular, the curve $\partial_{\sf top}^+\mathcal{M}$ parametrizes a family of successive gradient catastrophes in the dominant Riemann variable, upstream of the very first singularity. To see this, simply note that $\p_x w_0(x) = (\p_y w)(\eta(x,t),t) \p_x \eta(x,t)$, and hence $\p_y w(y,t) =w_0^\prime(\eta^{-1}(y,t)) (\p_x \eta(\eta^{-1}(y,t),t))^{-1}$. Since $(\eta^{-1}(y,t),t) \to (\bar{x}, \mathsf{T}_*(\bar{x}))$ as $(y,t)\to(\bar{y},\bar{t})$, we have that $\p_x \eta(\eta^{-1}(y,t),t)  \to\p_x \eta(\bar{x}, \mathsf{T}_*(\bar{x})) = 0$, and since for $x$ in the vicinity of $x_*$ we have $w_0^\prime(x) \approx - \frac{1}{\eps} < 0$, the claimed divergence towards $-\infty$ follows. We note that the gradient catastrophes lurking on $\partial_{\sf top}^+\mathcal{M}$ are not jump discontinuities, instead, the field $w$ has H\"older cusps on this curve.  Using the same argument as in Footnote~\ref{footnote:16}, one may verify that for $x>x_*$ the regularity of these cusps is H\"older $\frac 12$.  

\item Consider labels $x$ such that $x<x_*$. In our setup, where the waves move to the right, we refer to these labels as being {\em upstream} of the pre-shock label. Here, the main observation is that there exists a distinguished slow acoustic characteristic which emanates from the very first singularity $(y_*,t_*)$, the green curve in Figure~\ref{fig:1D:Eulerian}. Letting $x^* = \phi^{-1}(y_*,t_*)$, the slow acoustic characteristic emanating from the very first singularity is the curve $\partial_{\sf top}^-\mathcal{M} := (\phi(x^*,t),t)_{t\geq t_*}$. The notation we have chosen indicates that this curve serves as the future temporal boundary of the ``upstream part'' of the spacetime $\mathcal{M}$. Indeed, for any $(y,t)$ which lies either ``below'' or ``to the left'' of $\partial_{\sf top}^-\mathcal{M}$, the backwards trajectories corresponding to the labels $\eta^{-1}(y,t)$ and $\phi^{-1}(y,t)$ remain in $\mathcal{M}$, and therefore they satisfy~\eqref{eq:flows:remain:in:max:dev}. On the other hand, any point $(y,t)$ which lies ``to the right'' of $\partial_{\sf top}^-\mathcal{M}$ and ``to the left'' of $\partial_{\sf top}^+\mathcal{M}$, the white region in Figure~\ref{fig:1D:Eulerian}, is inaccessible from the perspective of the Cauchy data at time $t=\initial$. This is because following a $\phi$ characteristic backwards in time from any such point $(y,t)$, we are bound to intersect $\partial_{\sf top}^+\mathcal{M}$ (the red curve in Figure~\ref{fig:1D:Eulerian}), and we have just discussed that there the field $w$ experiences a gradient singularity. As such, as smooth continuation backwards in time ``through'' $\partial_{\sf top}^+\mathcal{M}$ and all the way back to the initial data is not possible. In the physical shock development problem, this issue is overcome by the introduction of a shock curve, see the discussion in Section~\ref{sec:1D:shock:dev} below.
\end{itemize}
In summary, the \MGHDB\ for the 1D isentropic Euler system is the spacetime $\mathcal{M}$ consisting of points $(y,t)$ which either lie to the left of the curve $(\eta(x_*,t),t)_{t\in[\initial,t_*]} \cup  \partial_{\sf top}^-\mathcal{M}$ (the upstream part), or they lie to the right of the curve $(\eta(x_*,t),t)_{t\in[\initial,t_*]} \cup  \partial_{\sf top}^+\mathcal{M}$ (the downstream part).

\begin{figure}[htb!]
\centering
\includegraphics[width=.6\textwidth]{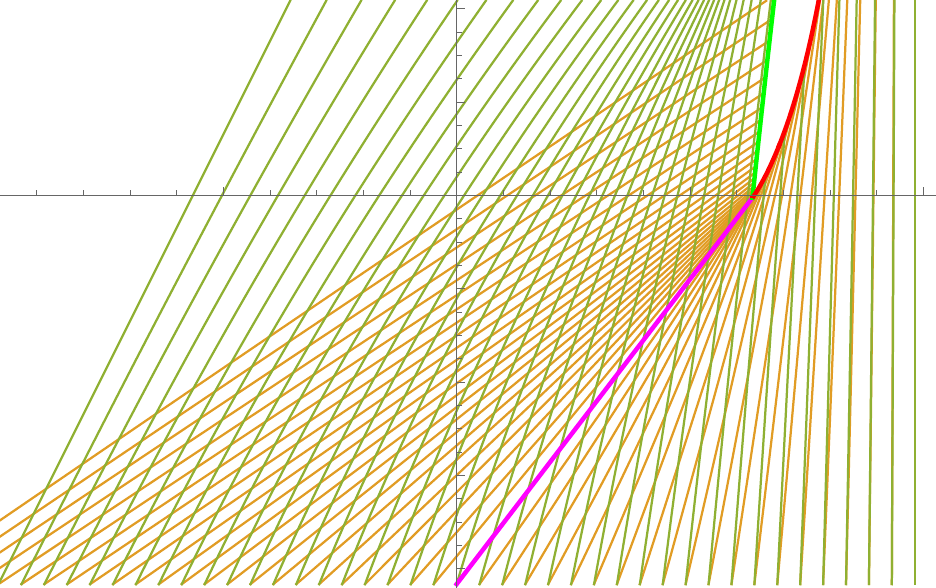}
\vspace{-.1in}
\caption{{\bf \MGHDB\ for 1D Euler via characteristics.} Consider the 1D Euler system~\eqref{eq:1D:Euler:Riemann} with initial data~\eqref{eq:1D:Euler:IC} for $\alpha=\frac 12$ and $\eps = \frac{1}{16}$. The bounding box represents $(x,t) \in [-\frac{3\pi \eps}{4},\frac{3\pi\eps}{4}] \times [\initial, \final]$. As in Figure~\ref{fig:1D:Eulerian:shock}, we represent in orange the ``fast'' characteristics $\eta(x,t)$ and in olive-green we represent the ``slow'' characteristics $\phi(x,t)$. The very first blow-up occurs at the spacetime location $(y_*,t_*) = ( \eta(x_*,t_*), t_*)$. Here $t_* = 0 \pm \OO(\eps^3)$ and $x_*=0$. The curve $\partial_{\sf top}^+\mathcal{M}$, characterizing the sequence of gradient catastrophes emerging from the very first singularity in the downstream region is represented in red. The curve $\partial_{\sf top}^-\mathcal{M}$, representing the slow acoustic characteristic emerging from the very first singularity in the upstream region is represented in green. The \MGHDB\ spacetime $\mathcal{M}$  is the complement of the white region.  In the white region of spacetime,  the acoustic characteristics cannot enter, so that the fields $(w,z)$ cannot be computed in a  smooth and unique way from $(w_0,z_0)$.}
\label{fig:1D:Eulerian}
\end{figure}

\begin{remark}[\bf Sound evolution, not the physical shock development]
\label{rem:sound:no:shock}
Its important to note that Figure~\ref{fig:1D:Eulerian} shows the dynamics of characteristics propagating as sound waves (and not shocks) and that the ``real'' Euler solution transitions from sound evolution to shock evolution (with associated the Rankine-Hugoniot jump conditions) instantaneously after the first singular time $t_*$. This issue is discussed in Section~\ref{sec:1D:shock:dev} below.
\end{remark}

\begin{remark}[\bf The boundary of $\mathcal{M}$ is not smooth]
We emphasize that in this characteristics language, the future temporal boundary of $\mathcal{M}$, i.e., the set $\partial_{\sf top}^-\mathcal{M} \cup \partial_{\sf top}^+\mathcal{M}$ is merely a Lipschitz surface. This lack of smoothness makes the analysis in this classical coordinate system rather cumbersome, as we shall describe next.  
\end{remark}

\subsection{A fast-acoustic characteristic perspective on shock formation}

The new perspective taken in this paper is to consider the entire analysis of shock formation and of \MGHDB\ in the Lagrangian coordinates\footnote{In multiple space dimensions, we in fact need to consider the Arbitrary-Eulerian-Lagrangian coordinates discussed in Section~\ref{sec:acoustic:geometry}.} $(\eta(x,t),t)$ of the fast acoustic characteristics. 

That this is a natural perspective may be seen from at least two points of view: first, the quantity which develops a singularity on the downstream part of the future temporal boundary of the spacetime $\mathcal{M}$ is the gradient of the dominant Riemann variable, namely $\p_y w$, and this quantity is naturally evolving along the flow $\eta$ (see~\eqref{eq:good:ol:Lax}); and second, the very first singularity and the subsequent gradient catastrophes in the downstream region are characterized in terms of properties of $\eta$ and its derivatives (see~\eqref{eta:shock:characterization} and \eqref{eta:shock:characterization:2}).

In order to build better intuition for the multi-D problem,\footnote{As explained in Section~\ref{sec:new:Euler:variables}, in multi-D, working with differentiated Riemann-type variables (which are then composed with the ALE map) is necessary in order to avoid derivative loss.} even in the 1D setting of this Appendix it is convenient to work with the spatially differentiated version of the isentropic Euler system~\eqref{eq:1D:Euler:Riemann}. That is, we define
\begin{equation*}
\mathring{W}(x,t):= (\p_y w)(\eta(x,t),t)\,,
\qquad
\mathring{Z}(x,t):= (\p_y z)(\eta(x,t),t)\,,
\qquad
\Sigma(x,t) := \sigma(\eta(x,t),t)\,.
\end{equation*}
With the above notation, the chain rule, and using the identity $\p_x \Sigma = \tfrac 12 w_0' - \tfrac 12 \p_x \eta \, \mathring{Z}$, the system~\eqref{eq:1D:Euler:Riemann} may be equivalently rewritten as\footnote{See also~\cite{NeRiShVi2023} for the similar variables in the non-isentropic setting.}
\begin{subequations}
\label{eq:Lagrangian:1D}
\begin{align}
\p_x \eta \, \mathring{W} 
&= w_0' \,,
\label{eq:Lagrangian:1D:a}\\
\p_t (\p_x \eta) &= \tfrac{1+\alpha}{2} w_0' + \tfrac{1-\alpha}{2} \p_x \eta \, \mathring{Z}\,,
\label{eq:Lagrangian:1D:b}\\ 
\p_x \eta \p_t \mathring{Z} - 2\alpha \Sigma \p_x \mathring{Z} 
&= - \mathring{Z} (\tfrac{1-\alpha}{2} w_0' + \tfrac{1+\alpha}{2} \p_x \eta \, \mathring{Z} ) \,,
\label{eq:Lagrangian:1D:c}\\
\p_t \Sigma 
&= - \alpha \Sigma \, \mathring{Z} \,.
\label{eq:Lagrangian:1D:d}
\end{align}
\end{subequations}
The remarkable features of the system~\eqref{eq:Lagrangian:1D} are as follows: \eqref{eq:Lagrangian:1D:a} shows that the quantity $\p_x \eta \, \mathring{W}$ is frozen into the flow, \eqref{eq:Lagrangian:1D:b} and \eqref{eq:Lagrangian:1D:d} are ODEs, while the only PDE left \eqref{eq:Lagrangian:1D:c} is the one for $\mathring{Z}$, which is a benign transport equation with a polynomial nonlinearity in the unknowns. It is thus natural that the analysis of this system for $(\mathring{W}, \mathring{Z}, \p_x\eta, \Sigma)$ should encounter no derivative loss (this turns out to be true even in multi-D). 

By working directly in the Lagrangian coordinate system associated to the fast acoustic characteristics $\eta$, the orange $\lambda_3$-characteristics from Figure~\ref{fig:1D:Eulerian:shock} and Figure~\ref{fig:1D:Eulerian} are now vertical straight lines (this holds ``by definition'', it is a tautology). What is more interesting is that the Lagrangian coordinates $\eta$, the characteristics associated to the slow wave speed $\lambda_1$ (which were the integral curves of the operator $\p_t + \lambda_1(y,t) \p_y$ in original coordinates, in olive-green in Figure~\ref{fig:1D:Eulerian:shock} and Figure~\ref{fig:1D:Eulerian}), are now the integral curves of the transport operator $\p_x \eta(x,t) \p_t - 2\alpha \Sigma(x,t) \p_x$. 

Finally, let us discuss the \MGHDB\ spacetime $\mathcal{M}$ in Lagrangian coordinates. This is represented in Figure~\ref{fig:1D:Lagrangian} below. Just as before, the location of the very first singularity at $(x_*,t_*)$ is determined by solving the system~\eqref{eta:shock:characterization}.
In the downstream region, i.e., for $x>x_*$, one still solves \eqref{eta:shock:characterization:2} in order to determine the time $\mathsf{T}_*(x)$ at which the curve $\{(x,t)\}_{t\geq \initial}$ intersects the downstream part of the future temporal boundary of $\mathcal{M}$. In Lagrangian coordinates, this boundary is characterized as $\partial_{\mathsf{top}}^+\mathcal{M} = \{ (x,\mathsf{T}_*(x))\}_{x>x_*}$. In light of \eqref{eq:Lagrangian:1D:a} it is immediately clear that the vanishing of $\p_x \eta(x,t)$ as $(x,t) \to (\bar{x},\mathsf{T}_*(\bar{x}))$, for any $\bar{x}\geq x_*$, implies the blow-up (towards $-\infty$) of $\mathring{W}(x,t)$; which is to say that $\partial_{\mathsf{top}}^+\mathcal{M}$ parametrizes a succession of gradient catastrophes past the time of the first singularity. In the upstream region, i.e., for $x<x_*$, we still solve for the distinguished slow acoustic characteristic (which is now an integral curve of the transport operator $\p_x \eta(x,t) \p_t - 2\alpha \Sigma(x,t) \p_x$), which passes through $(x_*,t_*)$. We denote the upstream part of this curve by $\partial_{\mathsf{top}}^-\mathcal{M}$, and this is the future temporal boundary of the upstream part of $\mathcal{M}$. 

The main observation is that in these fast acoustic Lagrangian coordinates the \MGHDB\ spacetime  becomes globally $W^{2,\infty}$ smooth (in essence, a cubic to the left, joining a parabola to the right, with matching first derivatives). This fact allows a PDE analysis in this spacetime which avoids derivative loss. 
The reader may compare Figure~\ref{fig:1D:Lagrangian}  (in Lagrangian coordinates)  and its direct analogue, Figure~\ref{fig:1D:Eulerian} (in Eulerian coordinates).

\begin{figure}[htb!]
\centering
\includegraphics[width=.6\textwidth]{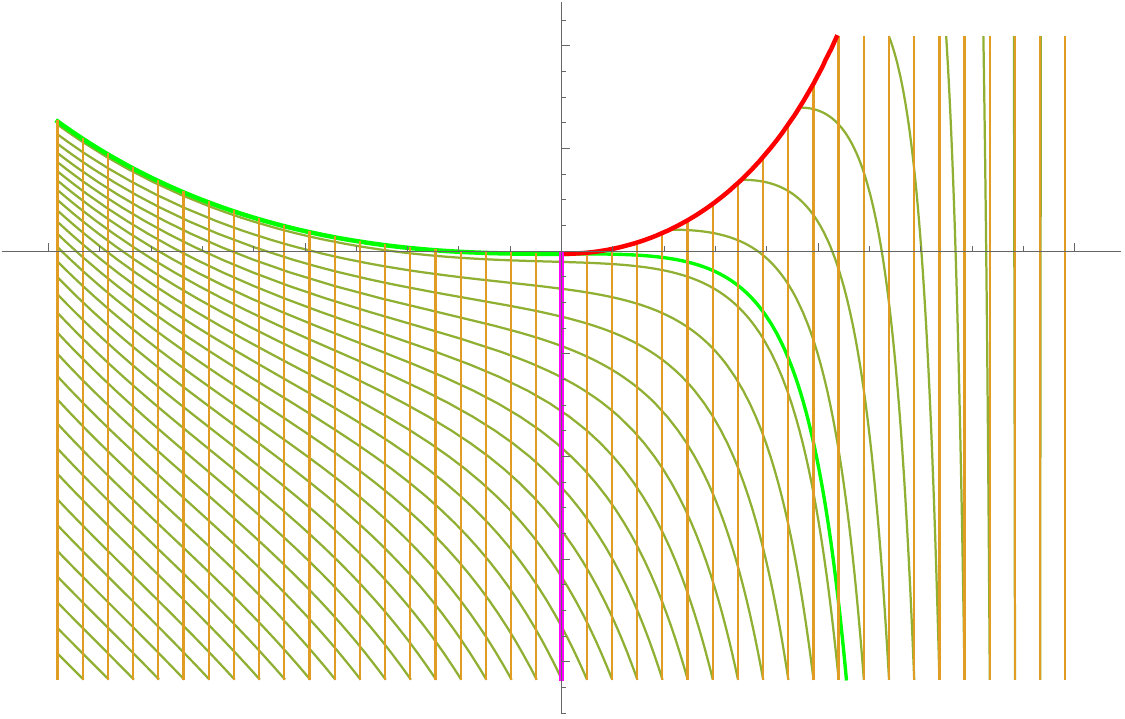};
\vspace{-.1in}
\caption{{\bf \MGHDB\ for 1D Euler -- the Lagrangian perspective.} Consider the 1D Euler system~\eqref{eq:1D:Euler:Riemann} with initial data as in \eqref{eq:1D:Euler:IC}. The bounding box represents $(x,t) \in [-\frac{3\pi \eps}{4},\frac{3\pi\eps}{4}] \times [\initial,\final]$, for the specific values $\alpha=\frac 12$ and $\eps = \frac{1}{16}$. We revisit Figure~\ref{fig:1D:Eulerian}, except that we work in the Lagrangian coordinates of the fast acoustic characteristic corresponding to $\lambda_3$. That is, a point $(y,t)$ in Figure~\ref{fig:1D:Eulerian} is replaced by the point $(x,t) = (\eta^{-1}(y,t),t)$ in the above figure. Then, the $\lambda_3$-characteristics become straight lines in orange, while the $\lambda_1$-characteristics become the curves in olive-green. The red curve represents the curve $\partial_{\mathsf{top}}^+\mathcal{M}$, while the green curve represents the curve $\partial_{\mathsf{top}}^-\mathcal{M}$ (which is extended backwards in time all the way up to $t=\initial$). In this Figure, the \MGHDB\ spacetime $\mathcal{M}$ is the region which lies ``below'' (in the causal past) of the union of  $\partial_{\mathsf{top}}^+\mathcal{M}$ and $\partial_{\mathsf{top}}^-\mathcal{M}$. This spacetime has a $W^{2,\infty}$-regular boundary.} 
\label{fig:1D:Lagrangian}
\end{figure}

\subsection{The physical shock development problem}
\label{sec:1D:shock:dev}

As already alluded to in Remark~\ref{rem:sound:no:shock}, when constructing the development of the Cauchy data past the time-slice of the very first singularity $\{t =t_*\}$,  the physical Euler evolution is replaced by the local hyperbolic propagation of sound waves via~\eqref{eq:1D:Euler:Riemann} (in the non-isentropic model, one would add the local hyperbolic propagation of entropy waves). While the extension of the solution to \eqref{eq:1D:Euler:Riemann} past the time-slice $\{t=t_*\}$ (see~Figure~\ref{fig:1D:Eulerian}) agrees with the physical, entropy-producing, shock solution\footnote{With Rankine-Hugoniot jump conditions to ensure that we are working with a weak solution of the full Euler system~\eqref{eq:Euler:1D:conservation}.} in certain regions of spacetime, it certainly {\em does not} agree with the ``real'' Euler solution globally, not even globally in $\mathcal{M}$!

In order to describe the physical, entropy-producing, shock solution to~\eqref{eq:Euler:1D:conservation} past the time slice $\{t=t_*\}$, we need to consider the shock curve, which may be parametrized as $(\mathfrak{s}(t),t)_{t>t_*}$, with $\mathfrak{s}(t_*)= y_* = \eta(x_*,t_*)$. While at the point of the very first singularity, $(\mathfrak{s}(t_*),t_*)$, the fields $(w,z,S)$ (and hence $(u,\rho,E)$) have in the worst case a H\"older $\frac 13$ cusp, for all $t>t_*$, the fields $(u,\rho,E)$ (and hence $(w,z,S)$) have a jump discontinuity across $(\mathfrak{s}(t),t)$. 

For convenience of notation, for $t>t_*$ let us denote the jump of a function $f(y,t)$ as $y$ crosses the shock location at $\mathfrak{s}(t)$ by $[[f(t)]] = f(\mathfrak{s}(t)^-,t) - f(\mathfrak{s}(t)^+,t)$. Here we have chosen the normal vector to the shock curve to point in the same direction as the propagation of the shock itself. With this convention, the physical entropy condition is that $[[S(t)]]>0$ for all $t>t_*$. 

\begin{figure}[htb!]
\centering
\includegraphics[width=.6\textwidth]{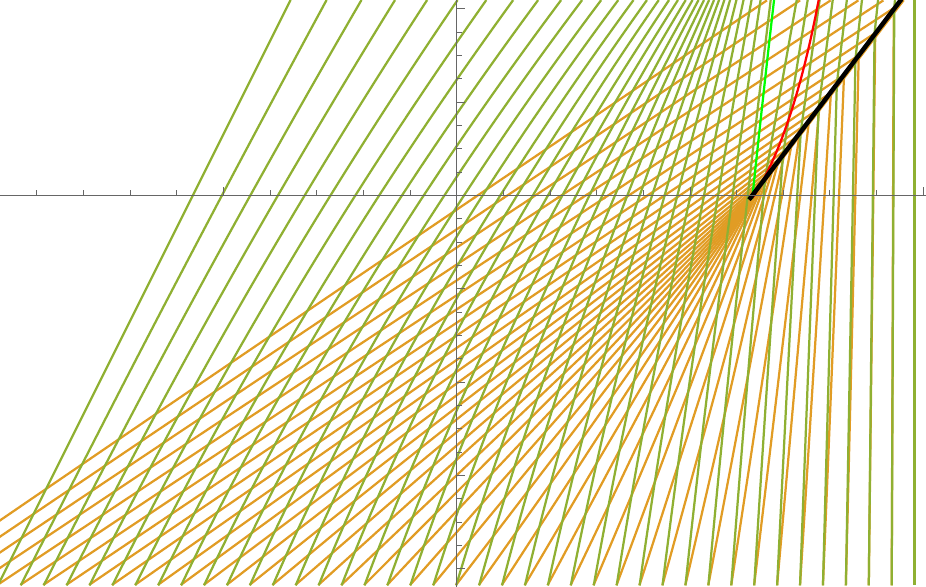};
\vspace{-.1in}
\caption{{\bf Shock development for 1D Euler.} Consider the 1D Euler system~\eqref{eq:1D:Euler:Riemann} with initial data as in \eqref{eq:1D:Euler:IC}. The bounding box represents $(x,t) \in [-\frac{3\pi \eps}{4},\frac{3\pi\eps}{4}] \times [\initial,\final]$, for the specific values $\alpha=\frac 12$ and $\eps = \frac{1}{16}$. We revisit Figure~\ref{fig:1D:Eulerian}, but instead of considering the \MGHDB, we consider the physical shock solution which emerges from the first singularity $(y_*,t_*)=(\eta(x_*,t_*),t_*)$. This physical shock curve, represented here in black, is parametrized as $(\mathfrak{s}(t),t)_{t>t_*}$. In order to not overcrowd the image, we have only represented the $\lambda_3$-characteristics (in orange), and the $\lambda_1$-characteristics (in olive-green), but have not depicted the $\lambda_2$-characteristics. From Figure~\ref{fig:1D:Eulerian} we have kept the green curve representing $\partial_{\sf top}^-\mathcal{M}$, and the red curve representing $\partial_{\sf top}^+\mathcal{M}$.}
\label{fig:1D:Shock}
\end{figure}

The shock speed (which uniquely determines the location of the shock) is denoted by $\dot{\mathfrak{s}}(t)$, and is computed via the Rankine-Hugoniot conditions as follows. The weak formulation of the conservation-law form of the 1D Euler equations~\eqref{eq:Euler:1D:conservation} yields a system of three equations in seven unknowns: $u(\mathfrak{s}(t)^\pm,t)$, $\rho(\mathfrak{s}(t)^\pm,t)$, $E(\mathfrak{s}(t)^\pm,t)$, and $\dot{\mathfrak{s}}(t)$. Equivalently, we may use the seven unknowns $w(\mathfrak{s}(t)^\pm,t)$, $z(\mathfrak{s}(t)^\pm,t)$, $S(\mathfrak{s}(t)^\pm,t)$, and $\dot{\mathfrak{s}}(t)$. The first observation is that the values of the fields downstream of the shock, i.e., $w(\mathfrak{s}(t)^+,t)$, $z(\mathfrak{s}(t)^+,t)$, $S(\mathfrak{s}(t)^+,t)$ may be computed uniquely and smoothly in terms of the initial data. With the isentropic data considered in this Appendix, this leads to $S(\mathfrak{s}(t)^+,t)=0$, $w(\mathfrak{s}(t)^+,t) = w_0(\eta^{-1}(\mathfrak{s}(t)^+,t),t)$ and $z(\mathfrak{s}(t)^+,t) = z_0(\phi^{-1}(\mathfrak{s}(t)^+,t),t)$. This is because upstream of the shock all characteristic families (corresponding to $\lambda_3$ for $w$, $\lambda_2$ for $S$, and $\lambda_1$ for $z$) are subsonic relative to the shock itself. The second observation (which together with the previous one are referred to as the {\em Lax geometric entropy conditions}) is that we may in fact also compute $w(\mathfrak{s}(t)^-,t)$ uniquely and smoothly in terms of the initial data as $w_0(\eta^{-1}(\mathfrak{s}(t)^-,t),t)$, because the fast acoustic characteristic upstream of the shock is supersonic relative to the shock itself.\footnote{This statement is a bit more subtle than it seems because the very introduction of the shock curve implies that the backwards map $\eta^{-1}(\cdot,t)$ is not continuous across $y=\mathfrak{s}(t)$, resulting in a nontrivial value for $[[w(t)]]$. The reader may for instance refer to Figure~\ref{fig:1D:Shock} and trace backwards-in-time the orange characteristic curves emanating from a point which is just to the left of $(\mathfrak{s}(t),t)$, and a point which is just to right of $(\mathfrak{s}(t),t)$. These backwards-in-time characteristics intersect $\{t=\initial\}$ at very different locations.} The third observation is that the Rankine-Hugoniot system reduces to a system in three equation and three unknowns, namely $z(\mathfrak{s}(t)^-,t)$, $S(\mathfrak{s}(t)^-,t)$, and $\dot{\mathfrak{s}}(t)$. For compressive and generic initial data $(w_0,z_0)$ as in~\eqref{eq:1D:Euler:IC} and with $S_0=0$, this 3$\times$3 system is uniquely solvable under the requirement that $S(\mathfrak{s}(t)^-,t) = [[S(t)]] > 0$, thus leading to a well-defined shock front evolution. The fourth and last observation is that the values of the fields $(w,z,S)$ on the left side (upstream) of the shock curve (that is, at $(\mathfrak{s}(t)^-,t)$) may now be used as Cauchy data for the system of equations~\eqref{eq:1D:entropy:transport} and~\eqref{eq:1D:Euler:Riemann}, which allows one to compute physical shock solution to 1D Euler globally in spacetime, even past the \MGHDB.

The shock development problem, as painted by the above paragraph, is illustrated in Figure~\ref{fig:1D:Shock} below. Mathematically, this picture has been analyzed in full-detail in~\cite{BuDrShVi2022}, in the context of the 2D Euler equations with azimuthal symmetry.\footnote{In~\cite{BuDrShVi2022} we have additionally described the weak characteristic singularities (weak-contact and weak-rarefaction) which emerge from the very first singularity at $(y_*,t_*)$, simultaneously with the shock curve in the upstream region.} See also~\cite{Yin2004} and~\cite{ChLi2016} in the context of 3D Euler with radial symmetry. 

\begin{remark}[\bf The shock development problem and the \MGHDB]
We close this appendix by discussing the location of the shock curve relative to the \MGHDB\ spacetime $\mathcal{M}$. Figure~\ref{fig:1D:Shock} is our point of reference. 

\begin{itemize}[leftmargin=16pt]

\item We emphasize that the shock curve is located in the interior of the downstream part of $\mathcal{M}$, i.e., the  curve $(\mathfrak{s}(t),t)_{t>t_*}$ (in black in Figure~\ref{fig:1D:Shock}) lies strictly below the curve $\partial_{\sf top}^+\mathcal{M}$ (in red in Figure~\ref{fig:1D:Shock}), and they intersect only at $(y_*,t_*)$. This is also why we may compute the values of all fields $(u,\rho,E)$ downstream of the shock: because this spacetime is fully contained in $\mathcal{M}$, and the interior of this spacetime contains no singularities whatsoever. 

\item This same fact, namely that the curve $(\mathfrak{s}(t),t)_{t>t_*}$ (in black in Figure~\ref{fig:1D:Shock}) lies strictly below the curve $\partial_{\sf top}^+\mathcal{M}$ (in red in Figure~\ref{fig:1D:Shock}) also shows that the \MGHDB\ computes a solution to Euler in a spacetime which is not physical (in the spacetime which lies between the red and black curves in Figure~\ref{fig:1D:Shock}). To see this, note merely that the physical shock solution is entropy producing, and thus $S\neq 0$ at points upstream of the shock. In contrast, the \MGHDB\ is computed via pure propagation of sound waves and thus the solution remains isentropic ($S=0$).

\item We note that in the upstream portion of $\mathcal{M}$, namely ``to the left'' or ``below''  the curve $\partial_{\sf top}^-\mathcal{M}$ (in green in Figure~\ref{fig:1D:Shock}), the physical shock solution in fact precisely matches the solution obtained as part of the \MGHDB. Again, the reason is that in this portion of spacetime there are no singularities, and each point therein is directly accessible from the initial data via $\lambda_i$-characteristics, for all $i\in \{1,2,3\}$. 

\item The last observation is the obvious one: solving the shock development problem allows one to compute the correct solution of 1D Euler (in particular, a weak solution, with shocks) in the portion of spacetime that lies beyond the \MGHDB\ (see the white region in Figure~\ref{fig:1D:Eulerian}). This is because the shock curve itself acts as a Cauchy hypersurface, on which the data is prescribed via a combination of the Lax geometric entropy conditions and the Rankine-Hugoniot jump conditions. The 1D Euler propagation away from the shock curve is via three transport equations: \eqref{eq:1D:entropy:transport}, \eqref{eq:1D:Euler:Riemann:a}, and~\eqref{eq:1D:Euler:Riemann:b}. Two of the corresponding characteristics (the slowest and the fastest) are represented in Figure~\ref{fig:1D:Shock}.
\end{itemize}
In summary, beyond the time slice of the very first singularity, $\{t = t_*\}$, the 1D Euler solution computed as part of the \MGHDB\ is useful (meaning, it describes the correct solution) in the spacetime which is either downstream of the shock (in black in Figure~\ref{fig:1D:Shock}), or which is upstream of the slow acoustic characteristic emanating from the very first singularity (in green in Figure~\ref{fig:1D:Shock}). Note however that since the shock curve is a-priori unknown (it must be nonlinearly solved for using the Rankine-Hugoniot jump conditions) a description of the Euler solution in the entire downstream part of $\mathcal{M}$ could prove to be useful in an iteration scheme  aimed at constructing the shock curve. 
\end{remark}

\section{Functional analysis in the remapped domain} 
\label{app:functional}
In this appendix we record a number of functional analytic bounds which are used throughout the paper. These bounds concern functions with arguments in the remapped spacetime $\TT^2 \times [0,\eps)$ (cf.~\eqref{eq:f:tilde:f}), which in turn is defined via the flattening map $\mathfrak{q} = \mathfrak{q}(x_2,t)$ from~\eqref{t-to-s-transform}. The associated differential operators given by the vector $\nbs$ (cf.~\eqref{nb-s}). This analytical framework is used in Sections~\ref{sec:formation:setup}--\ref{sec:sixth:order:energy}. In  Sections~\ref{rem:app:downstream:flat} and~\ref{rem:app:upstream:flat} below, we show that the same analytical framework developed here applies {\em as is} to the remapped domain from Section~\ref{sec:downstreammaxdev} and respectively Section~\ref{sec:upstreammaxdev}.

In view of definition~\eqref{eq:ic:supp} and the bootstrap~\eqref{bs-supp}, all the functions we consider in this Appendix are defined on $\mathcal{X}_{\rm fin} \times [0,\eps)$, where $\mathcal{X}_{\rm fin}  \subset \TT^2$ has Lebesgue measure bounded by $4\pi (13 \pi + \mathsf{C_{supp}}) \eps$, where the constant $\mathsf{C_{supp}} = \mathsf{C_{supp}}(\alpha,\kappa_0)>0$ was defined in \eqref{eq:Csupp:def}.  As such, throughout this section we will consider functions $f \colon \TT^2 \times [0,\eps) \to \mathbb{R}$ with the property that the support of $f$ in the $x_1$ direction has diameter $\les \eps$; more precisely, 
\begin{equation}
\label{eq:f:supp:ass}
\Bigg| \bigcup_{(x_2,\s) \in \TT \times [0,\eps)} \supp (f (\cdot,x_2,\s)) \Bigg|  \leq 4\pi (13 \pi + \mathsf{C_{supp}})   \eps
\,.
\end{equation} 
Additionally, throughout this section the  implicit constants in the $\les$ symbols depend on the constant $\mathsf{C_{supp}}$, and hence they depend only on $\alpha$ and $\kappa_0$.

\subsection{Sobolev and Poincar\'e}
Instead of the ``usual'' Sobolev embedding $H^s \subset L^\infty$, for $s>1$, we shall use the below bound, which avoids placing two copies of $\p_1 = \eps^{-1} \nbs_1$ derivatives on any term. 
\begin{lemma}
\label{lem:anisotropic:sobolev}
For sufficiently smooth $f \colon \TT^2 \times [0,\eps) \to \mathbb{R}$ satisfying~\eqref{eq:f:supp:ass}, we have
\begin{subequations}
\begin{align}
\snorm{f}_{L^2_{x,\s}} 
&\les \snorm{\nbs_1 f}_{L^2_{x,\s}}
\label{eq:x1:Poincare}
\\
\snorm{\p_2 f}_{L^2_{x,\s}} 
&\leq \snorm{\nbs_2 f}_{L^2_{x,\s}} + 8 \snorm{\nbs_\s f}_{L^2_{x,\s}} 
\label{eq:x2:Poincare}
\\
\snorm{f}_{L^ \infty_\s L^2_x}
&\les 
\snorm{f(\cdot,0)}_{L^2_x}
+ 
\eps^{-\frac 12}  \snorm{\nbs_\s f}_{L^2_{x,\s}} 
\label{eq:sup:in:time:L2} 
\\
 \snorm{f}_{L^\infty_{x,\s}} 
&\les 
\eps^{-\frac 12} \snorm{\nbs \nbs_1 f(\cdot,0)}_{L^2_x} 
+
\eps^{-1} \snorm{\nbs \nbs_\s \nbs_1 f}_{L^2_{x,\s}} \,,
\label{eq:Sobolev}
\\
\snorm{f}_{L^2_\s L^\infty_x}
&\les 
\eps^{-\frac 12} \snorm{\nbs \nbs_1 f}_{L^2_{x,\s}} 
\,,
\label{eq:Steve:punch}
\end{align}
\end{subequations}
where the implicit constant depends only on $\mathsf{C_{supp}}$, and hence on $\alpha$ and $\kappa_0$.
\end{lemma}
\begin{proof}[Proof of Lemma~\ref{lem:anisotropic:sobolev}]
The usual Poincar\'e inequality and the support assumption for $f(\cdot,\s)$ in \eqref{eq:f:supp:ass}  implies that 
\begin{equation} 
\snorm{f(\cdot,\s)}_{L^2_{x}} 
\les
\eps \snorm{\p_1 f(\cdot,\s)}_{L^2_{x}}
\,,
\label{eq:Poncare:temp:1}
\end{equation}
which yields \eqref{eq:x1:Poincare} since $\nbs_1 = \eps \p_1$.
The bound~\eqref{eq:x2:Poincare} follows since 
\begin{equation}
\p_2 = \nbs_2 + \Qb_2 \p_\s = \nbs_2 +  (\tfrac{1}{\eps} \Qb_2 \Qd^{-1})  \nbs_\s
\,,
\label{eq:Poncare:temp:2}
\end{equation}  
and  $|\frac{1}{\eps} \Qb_2 \Qd^{-1}| \leq \frac{3 \cdot 40}{17} \leq 8$.
Next, for any function $f$, the fundamental theorem of calculus in time implies that 
\begin{equation*}
\sup_{\s\in[0,\eps]} \norm{f(\cdot,\s)}_{L^2_x}^2 
\leq \norm{f(\cdot,0)}_{L^2_x}^2 
+ 2 \norm{f}_{L^2_{x,\s}} 
\norm{\p_\s f(\cdot,s)}_{L^2_{x,\s}}
.
\end{equation*}
Combined with the fact that $\norm{\p_\s f(\cdot,\s)}_{L^2_{x}} \les \eps^{-1} \snorm{\nbs_\s f(\cdot,\s)}_{L^2_{x}}$, which in turn is a consequence of $\snorm{\Qd^{-1}}_{L^\infty_x} \leq \frac{40}{17}$ (see~\eqref{Qd-lower-upper}), it follows that 
\begin{align}
\sup_{\s\in[0,\eps]} \snorm{f(\cdot,\s)}_{L^2_x}
&\les 
\snorm{f(\cdot,0)}_{L^2_x}
+ 
\eps^{-\frac 12} \snorm{f}_{L^2_{x,\s}}^{\frac 12}\snorm{\nbs_\s f}_{L^2_{x,\s}}^{\frac 12}
\notag\\
&\les 
\snorm{f(\cdot,0)}_{L^2_x}
+ 
\eps^{-\frac 14} \Bigl(\; \sup_{\s\in[0,\eps]} \snorm{f(\cdot,\s)}_{L^2_x}\Bigr)^{\frac 12}  \snorm{\nbs_\s f}_{L^2_{x,\s}}^{\frac 12}
\,,
\label{eq:L2:FTC:time}
\end{align}
and hence \eqref{eq:sup:in:time:L2} follows.

Lastly, for a fixed $\s\in [0,\eps]$, we may apply the fundamental theorem of calculus to $f(\cdot,\s)$, separately in $x_1$ and $x_2$, and then using \eqref{eq:Poncare:temp:1}--\eqref{eq:Poncare:temp:2}, we obtain 
\begin{align}
\norm{f(\cdot,\s)}_{L^\infty_x} 
&\les 
\bigl( \norm{f(\cdot,\s)}_{L^2_x} + \norm{f(\cdot,\s)}_{L^2_x}^{\frac 12} \norm{\p_2 f(\cdot,\s)}_{L^2_x}^{\frac 12} \bigr)^{\frac 12}
\bigl( \norm{\p_1 f(\cdot,\s)}_{L^2_x} + \norm{\p_1 f(\cdot,\s)}_{L^2_x}^{\frac 12} \norm{\p_{12} f(\cdot,\s)}_{L^2_x}^{\frac 12} \bigr)^{\frac 12}
\notag\\
&\les \eps^{-\frac 12}
\bigl( \snorm{\nbs_1^2 f(\cdot,\s)}_{L^2_x} + \snorm{\nbs_1^2 f(\cdot,\s)}_{L^2_x}^{\frac 12} \snorm{\nbs_1 \p_2 f(\cdot,\s)}_{L^2_x}^{\frac 12} \bigr) 
\notag\\
&\les  \eps^{-\frac 12}\snorm{\nbs_1 \nbs f(\cdot,\s)}_{L^2_x}
\,.
\label{eq:L2:FTC:space}
\end{align}  
Note that taking the $L^2_\s$ norm of both sides of the above estimate trivially yields \eqref{eq:Steve:punch}.
Then we apply \eqref{eq:sup:in:time:L2} with $f$ replaced by $\nbs_2 \nbs_1 f$ and respectively $\nbs_1 \nbs  f$ to deduce
\begin{equation}
\norm{f}_{L^\infty_{x,\s}} 
\les  
\eps^{-\frac 12} \snorm{\nbs  \nbs_1 f(\cdot,0)}_{L^2_x} 
+
\eps^{-1} \snorm{\nbs_\s \nbs \nbs_1 f}_{L^2_{x,\s}} 
\,,
\end{equation}
thereby establishing \eqref{eq:Sobolev}. 
\end{proof}

\begin{remark}
\label{rem:app:useful?}
Sometimes it is convenient to use a slightly modified variant of \eqref{eq:sup:in:time:L2}, which conforms to the transport operator $(\Qd \p_\s + V   \nbs_2)$, rather than the directional derivative $  \nbs_\s$. We claim that
\begin{equation}
\snorm{f}_{L^ \infty_\s L^2_x}
\les 
\snorm{f(\cdot,0)}_{L^2(\Omega)}
+ 
\eps^{\frac 12}  \snorm{(\Q \p_\s + V   \p_2) f}_{L^2_{x,\s}} 
\,,
\label{eq:sup:in:time:L2:alt}  
\end{equation}
where the implicit constant is universal.
In order to prove \eqref{eq:sup:in:time:L2:alt}, 
we use \eqref{adjoint-2} and then \eqref{the-time-der} to write 
\begin{align*}
\sup_{\s\in[0,\eps]} \norm{f(\cdot,\s)}_{L^2_x}^2 
&= \norm{f(\cdot,0)}_{L^2_x}^2 
+ 2 \int_0^{\s}\!\!\! \int f \p_\s f 
\notag\\
&= \norm{f(\cdot,0)}_{L^2_x}^2 
+ 2 \int_0^{\s}\!\!\! \int \tfrac{1}{\Qd} f \bigl(\Qd \p_\s + V \nbs_2 ) f 
-  \int_0^{\s}\!\!\! \int \tfrac{V}{\Qd} \nbs_2 (f^2)
\notag\\
&= \norm{f(\cdot,0)}_{L^2(\Omega)}^2 
+ 2 \int_0^{\s}\!\!\! \int \tfrac{1}{\Qd} f \bigl(\Q \p_\s + V \p_2 ) f 
+ \int_0^{\s}\!\!\! \int f^2 \bigl( -\Qr_2 + \nbs_2 \bigr) (  \tfrac{V}{\Qd} )
+   \! \int \bigl( f^2 \tfrac{V \Qb_2}{\Qd}\bigr)(x,\s) {\rm d}x
\,.
\end{align*}
Then, using  \eqref{bs-V}  and \eqref{eq:Q:all:bbq}  we arrive at
\begin{equation*}
\sup_{\s\in[0,\eps]} \norm{f(\cdot,\s)}_{L^2_x}^2 
\leq \norm{f(\cdot,0)}_{L^2_x}^2 
+ 5\eps^{\frac 12}   \snorm{\bigl(\Q \p_\s + V \p_2 ) f}_{L^2_{x,\s}} \sup_{\s\in[0,\eps]} \norm{f(\cdot,\s)}_{L^2_x}
+ \Cn \eps^2  \sup_{\s\in[0,\eps]} \norm{f(\cdot,\s)}_{L^2_x}^2.
\end{equation*}
The bound \eqref{eq:sup:in:time:L2:alt} now follows by absorbing terms involving $\sup_{\s\in[0,\eps]} \norm{f(\cdot,\s)}_{L^2_x}$ into the left side, upon taking $\eps$ to be sufficiently small.
\end{remark}

\subsection{Gagliardo-Nirenberg}

\begin{lemma}\label{lem:time:interpolation}
Let $f \colon \TT^2 \times [0,\eps) \to \mathbb{R}$ be a smooth function satisfying~\eqref{eq:f:supp:ass}. 
Then, for every $0 \leq i \leq m$ we have 
\begin{equation}
\snorm{\nbs^i f}_{L^{\frac{2m}{i}}_{x,\s}} 
\les 
\snorm{\nbs^m f}_{L^2_{x,\s}}^{\frac{i}{m}} 
\Bigl( \norm{f}_{L^\infty_{x,\s}} + \max_{1\leq j \leq i}  \snorm{\nbs^j f(\cdot,0)}_{L^\infty_x}    \Bigr)^{1- \frac{i}{m}} 
+  {\bf 1}_{0<i<m}  \eps^{\frac{i}{m}} 
\Bigl( \norm{f}_{L^\infty_{x,\s}} + \max_{1\leq j \leq i}  \snorm{\nbs^j f(\cdot,0)}_{L^\infty_x}    \Bigr) \,,
\label{se3:time}
\end{equation}
where the implicit constant depends only on $i$, $m$, and $\mathsf{C_{supp}}$ (and hence on $\alpha$ and $\kappa_0$).
\end{lemma} 
\begin{proof}[Proof of Lemma~\ref{lem:time:interpolation}]
The bound \eqref{se3:time} is trivial when $i=0$ or $i=m$, so we restrict the proof to $1\leq i \leq m-1$. The first nontrivial case is $i=1$ and $m=2$.
For this purpose we note that space integration by parts, the definition of the adjoint $\nbs^*$ in \eqref{eq:adjoints}, and the bounds \eqref{eq:Q:all:bbq}, yields
\begin{subequations}
\label{eq:waiting:for:the:worms:all}
\begin{align}
\snorm{\nbs_1 f}_{L^4_{x,\s}}^4 
&\les \norm{f}_{L^\infty_{x,\s}} \snorm{\nbs_1^2 f}_{L^2_{x,\s}} \snorm{\nbs_1 f}_{L^4_{x,\s}}^2
\label{eq:waiting:for:the:worms:a}
\\
\snorm{\nbs_2 f}_{L^4_{x,\s}}^4 
&\les \norm{f}_{L^\infty_{x,\s}} \snorm{\nbs_2^2 f}_{L^2_{x,\s}} \snorm{\nbs_2 f}_{L^4_{x,\s}}^2 
+ \eps^{\frac 12} \norm{f}_{L^\infty_{x,\s}} \snorm{\nbs_2 f}_{L^4_{x,\s}}^3
+ \eps \norm{f}_{L^\infty_{x,\s}} \sup_{\s\in[0,\eps]} \snorm{\nbs_2 f}_{L^3_x}^3\,,
\label{eq:waiting:for:the:worms:b}
\\
\snorm{\nbs_\s f}_{L^4_{x,\s}}^4 
&\les \norm{f}_{L^\infty_{x,\s}} \snorm{\nbs_\s^2 f}_{L^2_{x,\s}} \snorm{\nbs_\s f}_{L^4_{x,\s}}^2 
+ \eps \norm{f}_{L^\infty_{x,\s}} \sup_{\s\in[0,\eps]} \snorm{\nbs_\s f}_{L^3_x}^3\,,
\label{eq:waiting:for:the:worms:c}
\end{align}
where the last terms on the second and third lines arise from the temporal boundary term in the formula for $\nbs_2^*$ and $\nbs_\s^*$. For these terms, similarly to \eqref{eq:L2:FTC:time} we may prove
\begin{equation}
 \sup_{\s\in[0,\eps]} \snorm{\nbs  f}_{L^3_x}^3
 \les \eps \snorm{\nbs f(\cdot,0)}_{L^\infty_x}^3
 + \eps^{-1} \snorm{\nbs  f}_{L^4_{x,\s}}^2
 \snorm{\nbs_\s \nbs f}_{L^2_{x,\s}}
 \label{eq:waiting:for:the:worms:d}
\end{equation}
\end{subequations}
When combined, the above three displays yield
\begin{equation}
\label{eq:waiting:for:the:worms}
 \snorm{\nbs f}_{L^4_{x,\s}}
 \les \snorm{f}_{L^\infty_{x,\s}}^{\frac 12} \snorm{\nbs^2 f}_{L^2_{x,\s}}^{\frac 12} 
 + \eps^{\frac 12} \snorm{f}_{L^\infty_{x,\s}}
 + \eps^{\frac 12} \snorm{\nbs f(\cdot,0)}_{L^\infty_{x}}
\end{equation}
thereby probing \eqref{se3:time} for $i=1$ and $m=2$.

The general case $m\geq 3$ and $1\leq i \leq m-1$ in \eqref{se3:time} may be proven in a similar manner, by induction. We omit these redundant details.
\end{proof}

\subsection{Product and commutator estimates}

\begin{lemma}[Moser inequality]
\label{lem:Moser:tangent}
For $m\ge 1$, assume that  $f,g \colon \TT^2 \times [0,\eps] \to \mathbb{R}$ be smooth functions satisfying the standing support assumption~\eqref{eq:f:supp:ass}. 
Inspired by~\eqref{se3:time}, define the quantities
\begin{subequations}
\label{eq:pink:shirt}
\begin{align}
\mathfrak{B}_{f}
&:= 
\snorm{f}_{L^\infty_{x,\s}} + 
\max_{1\leq j \leq m-1}  \snorm{\nbs^j f(\cdot,0)}_{L^\infty_x}  \,,
\\
\mathfrak{B}_{g}
&:= 
\snorm{g}_{L^\infty_{x,\s}} + 
\max_{1\leq j \leq m-1}  \snorm{\nbs^j g(\cdot,0)}_{L^\infty_x}
\,.
\end{align}
\end{subequations}
Then, we have that
\begin{equation}
\snorm{\nbs^m (f g)}_{L^2_{x,\s}} 
\les 
\snorm{\nbs^m  f}_{L^2_{x,\s}}  \mathfrak{B}_{g}
+
\snorm{\nbs^m  g}_{L^2_{x,\s}}  \mathfrak{B}_{f}
+
\eps \; \mathfrak{B}_{f}    \mathfrak{B}_{g}  
\,,
\label{eq:Lynch:1}
\end{equation}
where the implicit constant only depends only on $m$ and $\mathsf{C_{supp}}$ (hence on $\alpha$ and $\kappa_0$).
\end{lemma}
\begin{proof}[Proof of Lemma~\ref{lem:Moser:tangent}]
Let $\gamma \in \mathbb{N}_0^3$ be such that $|\gamma|=m$.
From the Leibniz rule and the H\"older inequality we have
\begin{equation}
\snorm{\nbs^\gamma (f g)}_{L^2_{x,\s}}
\leq \sum_{\beta \leq \gamma} {\gamma \choose \beta} 
\snorm{\nbs^\beta f }_{L^{\frac{2m}{|\beta|}}_{x,\s}}
\snorm{\nbs^{\gamma-\beta} g}_{L^{\frac{2m}{|\gamma-\beta|}}_{x,\s}},
\end{equation}
where we recall that $\nbs^\gamma = \nbs_\s^{\gamma_0}  \nbs_1^{\gamma_1} \nbs_2^{\gamma_2}$. 
For each of the above terms we apply \eqref{se3:time} to deduce
\begin{align*}
\snorm{\nbs^m (f g)}_{L^2_{x,\s}}
&\les \sum_{i=0}^{m}
\bigl( \snorm{\nbs^m f}_{L^2_{x,\s}}^{\frac{i}{m}} \mathfrak{B}_{f}^{1- \frac{i}{m}} 
+  {\bf 1}_{0<i<m}  \eps^{\frac{i}{m}} \mathfrak{B}_{f} \bigr)
\bigl( \snorm{\nbs^m g}_{L^2_{x,\s}}^{1-\frac{i}{m}} \mathfrak{B}_{g}^{\frac{i}{m}} 
+  {\bf 1}_{0<i<m}  \eps^{1-\frac{i}{m}} \mathfrak{B}_{g} \bigr)
\notag\\
&\les \eps  \mathfrak{B}_{f}  \mathfrak{B}_{g}
+ \sum_{i=0}^m
(\snorm{\nbs^m f}_{L^2_{x,\s}}\mathfrak{B}_{g})^{\frac{i}{m}} 
(\snorm{\nbs^m g}_{L^2_{x,\s}} \mathfrak{B}_{f})^{1- \frac{i}{m}}
\notag\\
&\qquad
+ \sum_{i=1}^{m-1}
(\snorm{\nbs^m f}_{L^2_{x,\s}}\mathfrak{B}_{g})^{\frac{i}{m}} 
(\eps \mathfrak{B}_{f} \mathfrak{B}_{g})^{1- \frac{i}{m}}
+ \sum_{i=1}^{m-1}
(\eps \mathfrak{B}_{f} \mathfrak{B}_{g})^{\frac{i}{m}} 
(\snorm{\nbs^m g}_{L^2_{x,\s}} \mathfrak{B}_{f})^{1- \frac{i}{m}}
\,.
\end{align*}
The proof concludes by appealing to Young's inequality.
\end{proof}

The proof of Lemma~\ref{lem:Moser:tangent} also yields bounds for commutators, which are defined as in Remark~\ref{rem:notation:nb:doublecomm}. For $\gamma = (\gamma_0,\gamma_1, \gamma_2) \in {\mathbb N}_0^3$ with $|\gamma|=m$ and any scalar function $f$, we shall denote
\begin{align*}
 \jump{\nbs^\gamma,f} g 
 &= \nbs^\gamma (f g) - f \nbs^\gamma g = \sum_{\beta \leq \gamma, |\beta|\leq |\gamma|-1} {\gamma \choose \beta} \nbs^{\beta} f \; \nbs^{\gamma-\beta} g \,
 \,,
\end{align*}
and we shall use the notation
\begin{equation*}
\doublecom{\nbs^\gamma, f , g}
= \nbs^{\gamma}(f \, g) - f \nbs^\gamma g - g \nbs^\gamma f
= \sum_{\beta \leq \gamma, 1\leq |\beta|\leq |\gamma|-1} {\gamma \choose \beta} \nbs^{\beta} f \; \nbs^{\gamma-\beta} g \,,
\end{equation*}
so that 
\begin{equation*}
  \jump{\nbs^\gamma,f} g = \doublecom{\nbs^\gamma, f , g} + g \nbs^\gamma f \,.
\end{equation*}
The above identity implies that  bounds for $  \jump{\nbs^\gamma,f} g$ are direct consequences of bounds for $\doublecom{\nbs^\gamma, f , g}$, which are more symmetric in nature.

\begin{lemma}[Double commutator estimate]
\label{lem:comm:tangent}
For $m\ge 2$, assume that  $f,g \colon \TT^2 \times [0,\eps] \to \mathbb{R}$ are smooth functions satisfying the standing support assumption~\eqref{eq:f:supp:ass}. In analogy to \eqref{eq:pink:shirt}, define
\begin{subequations}
\label{eq:green:shirt}
\begin{align}
\mathfrak{B}_{\nbs f}
&:= 
\snorm{\nbs f}_{L^\infty_{x,\s}} + 
\max_{1\leq j \leq m-2}  \snorm{\nbs^j \nbs f(\cdot,0)}_{L^\infty_x}  \,,
\\
\mathfrak{B}_{\nbs g}
&:= 
\snorm{\nbs g}_{L^\infty_{x,\s}} + 
\max_{1\leq j \leq m-2}  \snorm{\nbs^j \nbs g(\cdot,0)}_{L^\infty_x}  
\,.
\end{align}
\end{subequations}
Then, we have that
\begin{equation}
\snorm{\doublecom{\nbs^m , f , g}}_{L^2_{x,\s}} 
\les 
\snorm{\nbs^{m-1}  f}_{L^2_{x,\s}}  \mathfrak{B}_{\nbs g}
+
\snorm{\nbs^{m-1}  g}_{L^2_{x,\s}}  \mathfrak{B}_{\nbs f}
+
\eps \; \mathfrak{B}_{\nbs f}    \mathfrak{B}_{\nbs g}  
\,,
\label{eq:Lynch:2}
\end{equation}
where the implicit constant depends only on $m$ and $\mathsf{C_{supp}}$ (hence on $\alpha$ and $\kappa_0$). 
In particular, the usual commutator $\jump{\nbs^m , f} g$ is bounded in $L^2_{x,\s}$ via \eqref{eq:Lynch:2} and the triangle inequality 
\begin{equation} 
\snorm{\jump{\nbs^m , f} g}_{L^2_{x,\s}} 
\leq \snorm{\doublecom{\nbs^m , f , g}}_{L^2_{x,\s}}  + \snorm{g}_{L^\infty_{t,x}}  \snorm{\nbs^m f}_{L^2_{x,\s}}
\,. \label{eq:Lynch:3}
\end{equation} 
\end{lemma}

\begin{proof}[Proof of Lemma~\ref{lem:comm:tangent}]
From the Leibniz rule and H\"older inequalities, we have that 
\begin{equation}
\snorm{\doublecom{\nbs^m, f , g}}_{L^2_{x,\s}}
\les
\sum_{i=1}^{m-1} 
\snorm{\nbs^{i-1} \nbs f}_{L^{\frac{2(m-2)}{i-1}}_{x,\s}} 
\snorm{\nbs^{m-i-1} \nbs g}_{L^{\frac{2(m-2)}{m-i-1}}_{x,\s}}
\,.
\end{equation}
Applying \eqref{se3:time}, we obtain
\begin{align}
\snorm{\doublecom{\nbs^m, f , g}}_{L^2_{x,\s}}
&\les
\sum_{i=1}^{m-1} 
\bigl( \snorm{\nbs^{m-1} f}_{L^2_{x,\s}}^{\frac{i-1}{m-2}} \mathfrak{B}_{\nbs f}^{1- \frac{i-1}{m-2}} 
+   \eps^{\frac{i-1}{m-2}} \mathfrak{B}_{\nbs f} \bigr)
\bigl( \snorm{\nbs^{m-1} g}_{L^2_{x,\s}}^{1-\frac{i-1}{m-2}} \mathfrak{B}_{\nbs g}^{\frac{i-1}{m-2}} 
+   \eps^{1-\frac{i-1}{m-2}} \mathfrak{B}_{\nbs g} \bigr)
\notag \\
&\les \eps \mathfrak{B}_{\nbs f}\mathfrak{B}_{\nbs g}
+ \sum_{i=1}^{m-1} 
 \bigl(\snorm{\nbs^{m-1} f}_{L^2_{x,\s}} \mathfrak{B}_{\nbs g}\bigr)^{\frac{i-1}{m-2}} 
 \bigl(\snorm{\nbs^{m-1} g}_{L^2_{x,\s}} 
 \mathfrak{B}_{\nbs f}\bigr)^{1- \frac{i-1}{m-2}} \notag\\
&\qquad 
+ \sum_{i=1}^{m-1} 
\bigl(\eps \mathfrak{B}_{\nbs f}\mathfrak{B}_{\nbs g}\bigr)^{\frac{i-1}{m-2}} 
\bigl(\mathfrak{B}_{\nbs f} \snorm{\nbs^{m-1} g}_{L^2_{x,\s}}\bigr)^{1-\frac{i-1}{m-2}} \notag\\
&\qquad 
+ \sum_{i=1}^{m-1} 
\bigl( \snorm{\nbs^{m-1} f}_{L^2_{x,\s}} \mathfrak{B}_{\nbs g} \bigr)^{\frac{i-1}{m-2}} 
\bigl(\eps \mathfrak{B}_{\nbs f} \mathfrak{B}_{\nbs g}  \bigr)^{1- \frac{i-1}{m-2}} 
\,.
\label{eq:poop:stinks:4}
\end{align}
The bound \eqref{eq:Lynch:1} now follows from Young's inequality.
\end{proof}

Next, we establish a version of the double-commutator bound from Lemma~\ref{lem:comm:tangent}, which is however uniform in $\s$.

\begin{lemma}[Double commutator estimate, uniform in time]
\label{lem:comm:bdd}
For $m\ge 2$, assume that  $f,g \colon \TT^2 \times [0,\eps] \to \mathbb{R}$ are smooth functions satisfying the support assumption~\eqref{eq:f:supp:ass}. Recall the definitions of $\mathfrak{B}_f$, $\mathfrak{B}_g$ in~\eqref{eq:pink:shirt}, and the definitions of $\mathfrak{B}_{\nbs f}$, $\mathfrak{B}_{\nbs g}$ in \eqref{eq:green:shirt}. Then,  
\begin{equation}
\snorm{ \doublecom{\nbs^m , f , g}}_{L^\infty_\s L^2_{x}} 
\les 
\eps^{-\frac 12} \Bigl( \snorm{\nbs^{m}  f}_{L^2_{x,\s}}  \mathfrak{B}_{\nbs g}
+
  \snorm{\nbs^{m}  g}_{L^2_{x,\s}}  \mathfrak{B}_{\nbs f}
+
\eps \; \mathfrak{B}_{\nbs f}    \mathfrak{B}_{\nbs g}  
\Bigr)
\,,
\label{eq:Lynch:2:bdd}
\end{equation}
where the implicit constant depends only on $m$ and $\mathsf{C_{supp}}$. 
For any function $\varphi \colon \TT^2 \times [0,\eps] \to [0,\infty)$ which is bounded from above, we have the commutator bound 
\begin{equation} 
\snorm{\varphi \jump{\nbs^m , f} g}_{L^\infty_\s L^2_{x}}  
\les 
\snorm{\varphi \nbs^m f}_{L^\infty_\s L^2_{x}} \mathfrak{B}_{g}
+
\eps^{-\frac 12} \Bigl( \snorm{\nbs^{m}  f}_{L^2_{x,\s}}  \mathfrak{B}_{\nbs g}
+
 \snorm{\nbs^{m}  g}_{L^2_{x,\s}}  \mathfrak{B}_{\nbs f}
+
\eps \; \mathfrak{B}_{\nbs f}    \mathfrak{B}_{\nbs g}  \Bigr)
\,,
 \label{eq:Lynch:3:bdd}
\end{equation} 
and the Moser-type bound 
\begin{align} 
\snorm{\varphi  \nbs^m (f \, g)}_{L^\infty_\s L^2_{x}} 
&\les 
 \snorm{\varphi \nbs^m f}_{L^\infty_\s L^2_{x}} \mathfrak{B}_{ g}
+ 
\snorm{\varphi \nbs^m g}_{L^\infty_\s L^2_{x}} \mathfrak{B}_{ f}
\notag\\
&\qquad 
+
\eps^{-\frac 12} \|\varphi\|_{L^\infty_{x,\s}}
\Bigl( \snorm{\nbs^{m}  f}_{L^2_{x,\s}}  \mathfrak{B}_{\nbs g}
+
 \snorm{\nbs^{m}  g}_{L^2_{x,\s}}  \mathfrak{B}_{\nbs f}
+
\eps \; \mathfrak{B}_{\nbs f}    \mathfrak{B}_{\nbs g} \Bigr)
\,, \label{eq:Lynch:1:bdd}
\end{align} 
where the implicit constant depends only on $m$,  and on $\mathsf{C_{supp}}$ (hence on $\alpha$ and $\kappa_0$).
\end{lemma}
\begin{remark}
\label{rem:shit:is:similar}
The important fact to notice about  \eqref{eq:Lynch:1:bdd} is that the terms on the second line are given by $\eps^{-\frac 12} \|\varphi\|_{L^\infty_{x,\s}}$ times the upper bound of $\snorm{\nbs^m (f \, g)}_{L^2_{x,\s}}$ given in \eqref{eq:Lynch:1}, except that $\mathfrak{B}_{f}$ and $\mathfrak{B}_{g}$ become $\mathfrak{B}_{\nbs f}$ and $\mathfrak{B}_{\nbs g}$. In all applications of this \eqref{eq:Lynch:1:bdd}, in view of \eqref{table:derivatives}  we have that $\mathfrak{B}_{f}$  obeys the same bounds as $\mathfrak{B}_{\nbs f}$ (and the same for $g$), and therefore the upper bound for $\snorm{\varphi  \nbs^m (f \, g)}_{L^\infty_\s L^2_{x}} $ will always be given by $\eps^{-\frac 12}\|\varphi\|_{L^\infty_{x,\s}}$ times the upper bound for $\snorm{   \nbs^m (f \, g)}_{L^2_{x,\s}}$, plus the first two terms on the right side of \eqref{eq:Lynch:1:bdd}, which retain the weight function $\varphi$ inside the norm.
\end{remark}

\begin{proof}[Proof of Lemma~\ref{lem:comm:bdd}]
From the triangle inequality we have that 
\begin{equation*}
\snorm{   \doublecom{\nbs^m , f , g}}_{L^\infty_\s L^2_{x}} 
\leq 2^m  {\textstyle\sum}_{j=1}^{m-1} \|  \nbs^j f \nbs^{m-j} g\|_{L^\infty_\s L^2_x} \,.  
\end{equation*}
In turn, for each fixed $1\leq j \leq m-1$ using \eqref{eq:sup:in:time:L2}, we have that 
 \begin{equation*}
 \| F \|_{L^\infty_\s L^2_x} 
 \les \|  F (0)\|_{L^2_x} 
 + \eps^{-\frac 12} \| \nbs_\s F\|_{L^2_{x,\s}} 
 \,,
\end{equation*}
where the implicit constant is universal. By appealing to the above bound with $F = \nbs^j f \nbs^{m-j} g$, we deduce that 
\begin{align*}
& {\textstyle\sum}_{j=1}^{m-1} \|  \nbs^j f \nbs^{m-j} g\|_{L^\infty_\s L^2_x} 
\notag\\
&\les  {\textstyle\sum}_{j=1}^{m-1} \Bigl( \eps^{\frac 12} \|  \nbs^j f(0)\|_{L^\infty_x} \| \nbs^{m-j} g(0)\|_{L^\infty_x}
+
\eps^{-\frac 12}  \|  \nbs_\s \nbs^j f \nbs^{m-j} g\|_{L^2_{x,\s}} 
+
\eps^{-\frac 12} \|   \nbs^j f \nbs_\s \nbs^{m-j} g\|_{L^2_{x,\s}} \Bigr)
\notag\\
&\les \eps^{\frac 12} \mathfrak{B}_{\nbs f} \mathfrak{B}_{\nbs g} 
+ \eps^{-\frac 12} \| \nbs^m f\|_{L^2_{x,\s}} \|\nbs g\|_{L^\infty_{x,\s}} 
+ \eps^{-\frac 12} \| \nbs^m g\|_{L^2_{x,\s}} \|\nbs f\|_{L^\infty_{x,\s}} 
+ \eps^{-\frac 12} {\textstyle\sum}_{j=2}^{m-1} \|  \nbs^{j-1} \nbs f \nbs^{m-j} \nbs g\|_{L^2_{x,\s}} 
\notag\\
&\les \eps^{\frac 12} \mathfrak{B}_{\nbs f} \mathfrak{B}_{\nbs g} 
+ \eps^{-\frac 12} \| \nbs^m f\|_{L^2_{x,\s}}  \mathfrak{B}_{\nbs g} 
+ \eps^{-\frac 12} \| \nbs^m g\|_{L^2_{x,\s}} \mathfrak{B}_{\nbs f} 
+ \eps^{-\frac 12} {\textstyle\sum}_{j=2}^{m-1} \|  \nbs^{j-1} \nbs f\|_{L^{\frac{2(m-1)}{j-1}}_{x,\s}} \|\nbs^{m-j} \nbs g\|_{L^{\frac{2(m-1)}{m-j}}_{x,\s}}
\end{align*}
where the implicit constant only depends on $m$. 
Using \eqref{se3:time}, the above   then furthermore implies
\begin{align*}
\snorm{  \doublecom{\nbs^m , f , g}}_{L^\infty_\s L^2_{x}} 
&\les \eps^{\frac 12} \mathfrak{B}_{\nbs f} \mathfrak{B}_{\nbs g} 
+ \eps^{-\frac 12} \| \nbs^m f\|_{L^2_{x,\s}}  \mathfrak{B}_{\nbs g} 
+ \eps^{-\frac 12} \| \nbs^m g\|_{L^2_{x,\s}} \mathfrak{B}_{\nbs f} 
\notag\\
&\qquad 
+ \eps^{-\frac 12} {\textstyle\sum}_{j=2}^{m-1} 
\Bigl(\|\nbs^m f\|_{L^2_{x,\s}}^{\frac{j-1}{m-1}} \mathfrak{B}_{\nbs f}^{\frac{m-j}{m-1}} + \eps^{\frac{j-1}{m-1}} \mathfrak{B}_{\nbs f}\Bigr)  
\Bigl(\|\nbs^m g\|_{L^2_{x,\s}}^{\frac{m-j}{m-1}} \mathfrak{B}_{\nbs g}^{\frac{j-1}{m-1}} + \eps^{\frac{m-j}{m-1}} \mathfrak{B}_{\nbs g}\Bigr)  
\notag\\
&\les \eps^{\frac 12} \mathfrak{B}_{\nbs f} \mathfrak{B}_{\nbs g} 
+ \eps^{-\frac 12} \| \nbs^m f\|_{L^2_{x,\s}}  \mathfrak{B}_{\nbs g} 
+ \eps^{-\frac 12} \| \nbs^m g\|_{L^2_{x,\s}} \mathfrak{B}_{\nbs f} 
\end{align*}
where the implicit constant only depends on $m$. This proves \eqref{eq:Lynch:2:bdd}. 
The bounds \eqref{eq:Lynch:3:bdd} and \eqref{eq:Lynch:1:bdd} now follow from the definitions of $ \doublecom{\nbs^m , f , g}$, $\jump{\nb^m,f} g$,   the triangle inequality, and the bound $\mathcal{J} \leq \Jgb \leq \Jg \leq \frac 32$.
\end{proof}

\subsection{The flattening map from Section~\ref{sec:downstreammaxdev}}
\label{rem:app:downstream:flat}
The map $\qds = \qds(x_1,x_2,t)$ defined in~\eqref{t-to-s-transform-P} is used in Section~\ref{sec:downstreammaxdev} to remap the spacetime which is downstream of the pre-shock onto $\TT^2 \times [0,\eps)$, via~\eqref{eq:f:tilde:f-P}. The associated differential operator $\nbs$ is given in \eqref{nb-s-P}. Save for the $x_1$-dependence of  $\qds$, the mappings $\mathfrak{q}$ and $\qds$ have nearly identical properties and bounds (see~Lemmas~\ref{lem:Q:bnds} and~\ref{lem:Q:bnds-DS}), which is why we have chosen to use nearly the same notation for these maps. 

The only difference between the analysis of the differential operators $\nbs$ from \eqref{nb-s} and the $\nbs$ from \eqref{nb-s-P} lies in the definition of $\nbs_1$. This difference results in the following modification to the upper bounds given in this section: every bound which contains a term of the type $\| \nbs_1 f\|$, with $\nbs_1$ as defined by \eqref{nb-s}, needs to be replaced by an upper of the type $\| \nbs_1 f\| + \|\nbs_\s f\|$, with $\nbs_1$ and $\nbs_\s$ as defined by \eqref{nb-s-P}. This change only affects the bounds \eqref{eq:x1:Poincare}, \eqref{eq:Sobolev}, and \eqref{eq:Steve:punch}; note however that when these estimates are used in the bulk of the paper, we in fact bound $\|\nbs_1 f\|$ with the full norm $\|\nbs f\|$, and therefore no single bound in the bulk of the paper is altered by this change. The estimates provided by the Remarks~\ref{rem:app:useful?},~\ref{rem:shit:is:similar}, and Lemmas~\ref{lem:time:interpolation}, \ref{lem:Moser:tangent}, \ref{lem:comm:tangent},~\ref{lem:comm:bdd} remain unchanged. In view of this fact, in Section~\ref{sec:downstreammaxdev} we still refer to the bounds in this Appendix, the above mentioned slight modification being implicit.
 
\subsection{The flattening map from Section~\ref{sec:upstreammaxdev}}
\label{rem:app:upstream:flat}
In Section~\ref{sec:upstreammaxdev} we work in  the spacetime $\tHdm$ defined in~\eqref{eq:spacetime-Theta-t}. While pointwise estimates are performed in $(x,t)\in\tHdm$ variables, the energy estimates at the center of our analysis are performed in the $(x,\s)$ variables that correspond to the remapped domain $\Hdm \subset \mathcal{X}_{\rm fin} \times [\sin,\sfin)$, defined in~\eqref{eq:spacetime-Theta}. This remapping (cf.~\eqref{eq:f:tilde:f-H}) is achieved by setting $\s = \mathfrak{q}(x_2,t)$, with the map $\mathfrak{q}$ being defined in~\eqref{t-to-s-transform-H}. In particular, the map $\mathfrak{q}$ is independent of $x_1$, and all functions  $f \colon \Hdm \to \mathbb{R}$ satisfy~\eqref{eq:f:supp:ass}. 

The differences between the $(x,\s)$ coordinates in Section~\ref{sec:upstreammaxdev} and the $(x,\s)$ coordinates in used in Sections~\ref{sec:formation:setup}--\ref{sec:sixth:order:energy} are as follows. First, the domain of the new time variable $\s$ is not $[0,\eps)$, but instead $[\sin,\sfin)$, as defined by~\eqref{sin} and~\eqref{sfin}. Since $\sfin-\sin \in [(\tfrac{51}{50} - \Cn \eps) \tfrac{2\eps}{1+\alpha}, (\tfrac{51}{50} + \Cn \eps) \tfrac{2\eps}{1+\alpha}]$ the length of the time interval is still $\OO(\eps)$, this time with an implicit constant depending only on $\alpha$ (on which the bounds in this Appendix are allowed to depend anyway). Second, the $L^2_x$ and $L^2_{x,\s}$ norms are not anymore defined on all of $\TT^2$, respectively $\TT^2 \times[0,\eps)$ (as was done in Sections~\ref{sec:formation:setup}--\ref{sec:sixth:order:energy}); instead, for $f \colon \Hdm \to \mathbb{R}$ we use the $L^2_{x,\s} = L^2_{x,\s}(\Hdm)$ norm defined in \eqref{norms-H:b}, while the $L^\infty_{\s} L^2_x = L^\infty_{\s} L^2_x(\Hdm)$ norm given by  $\sup_{\s\in[\sin,\sfin]} \|f(\cdot,\s)\|_{L^2_x}$, with $\|f(\cdot,\s)\|_{L^2_x}$ as defined by \eqref{norms-H:a}. Third, we now have that $\s$-time-slices of $\Hdm$ have a ``right'' boundary, at $x_1 = \thd(x_2,\s)$. Nonetheless, due to~\eqref{eq:f:supp:ass} it still holds that the norm in~\eqref{norms-H:a} satisfies a Poincar\'e inequality of the type $\|f(\cdot,\s)\|_{L^2_x} \les \eps \|\p_1 f(\cdot,\s)\|_{L^2_x} = \|\nbs_1 f(\cdot,\s)\|_{L^2_x}$, as claimed in~\eqref{eq:x1:Poincare}.  The new spacetime differential operator $\nbs$ is defined in~\eqref{nb-s-H}, with the adjoint $\nbs^*$ being given by~\eqref{eq:adjoints-H}. Since $\mathfrak{q}$ is independent of $x_1$, the estimates the $\Q$ coefficients appearing in the definition of the differential operator $\nbs$ are nearly identical to those in Sections~\ref{sec:formation:setup}--\ref{sec:sixth:order:energy}. Indeed, comparing the estimates in Lemma~\ref{lem:Q:bnds-H} to those in Lemma~\ref{lem:Q:bnds} we see that the only differences involve $\alpha$-dependent factors. 

Another potential difference concerns the  usage of the fundamental theorem of calculus in time (which is used to prove~\eqref{eq:sup:in:time:L2}, \eqref{eq:sup:in:time:L2:alt}), since the location of the right boundary of each $\s$-time-slice depends on $\s$. What saves the day is that the $x_1$ boundary term at $\thd(x_2,\s)$ (see~\eqref{calc-ds}) always has a sign since $\p_\s \thd \leq 0$ (see~\eqref{ps-theta-sign}). For example, 
\begin{align} 
\snorm{F(\cdot, \s)}_{L^p_x}^p 
&= \snorm{F(\cdot, \sin)}_{L^p_x}^p  
+ \int_\sin^\s \tfrac{d}{d\s^\prime}  \snorm{F(\cdot, \s^\prime)}_{L^2_x}^2 {\rm d}\s^\prime
\notag \\
&= \snorm{F(\cdot, \sin)}_{L^p_x}^p  \
+ \int_\sin^\s  \Big(\dint p \sabs{F(x,\s^\prime)}^{p-2} F(x,\s^\prime)  \p_{\s^\prime}F(x,\s^\prime){\rm d}x 
\notag \\
& \qquad \qquad  \qquad \qquad \qquad
+ \int_{x_2=-\pi}^\pi \underbrace{\p_{\s'}\thd(x_2,\s')}_{\leq 0} \sabs{F(\thd(x_2,\s^\prime) ,x_2,\s^\prime)}^p {\rm d}x_2\Big) {\rm d}\s^\prime 
\notag \\
&\le \snorm{F(\cdot, \sin)}_{L^p_x}^p + p  \int_\sin^\s  \snorm{F(\cdot, \s^\prime)}_{L^{\frac{q(p-1)}{q-1}}_x}^{p-1} \snorm{\p_{\s^\prime}F(\cdot,\s')}_{L^q_x} {\rm d}\s^\prime \,, 
\label{eq:got:my:mojo:workin}
\end{align} 
where in the last inequality we have used H\"older and \eqref{ps-theta-sign}. As such, no actual change arises when applying the fundamental theorem of calculus in $\s$, we are simply not using a helpful signed term in our upper bounds.

The last potential difference concerns the  integration by parts in $\s$ that is used to prove
~\eqref{se3:time}. To see that no difference arises in the final estimate, let us for instance consider the proof of~\eqref{eq:waiting:for:the:worms}, which consists of establishing the bounds~\eqref{eq:waiting:for:the:worms:a}--\eqref{eq:waiting:for:the:worms:d}. First, we note that~\eqref{eq:waiting:for:the:worms:d} remains unchanged, by using an argument similar to that  in~\eqref{eq:got:my:mojo:workin}, for $p=3$ and $q=2$. Second, using the definition of $\nbs_1^*$ in~\eqref{adjoint-1-H}, we see that the bound~\eqref{eq:waiting:for:the:worms:a} becomes
\begin{equation*}
\snorm{\nbs_1 f}_{L^4_{x,\s}}^4 
\leq 3 \norm{f}_{L^\infty_{x,\s}} \snorm{\nbs_1^2 f}_{L^2_{x,\s}} \snorm{\nbs_1 f}_{L^4_{x,\s}}^2
+ \eps \int_{\sin}^{\sfin}\int_{\TT} \bigl| f (\nbs_1 f)^3 \bigr|(\thd(x_2,\s),x_2,\s) {\rm d}x_2 {\rm d}\s
\,.
\end{equation*}
For the second term in the above estimate, using the support property~\eqref{eq:f:supp:ass}, the fundamental theorem of calculus in the $x_1$ variable, and the fact that $\nbs_1 = \eps \p_1$, we deduce that 
\begin{align*}
\eps \int_{\sin}^{\sfin}\int_{\TT} \bigl| f (\nbs_1 f)^3 \bigr|(\thd(x_2,\s),x_2,\s) {\rm d}x_2 {\rm d}\s
&\leq \eps \|f\|_{L^\infty_{x,\s}} \int_{\sin}^{\sfin}\int_{\TT} |\nbs_1 f|^3(\thd(x_2,\s),x_2,\s) {\rm d}x_2 {\rm d}\s
\notag\\
&\leq  \|f\|_{L^\infty_{x,\s}} \int_{\sin}^{\sfin}\int_{\TT} \! \int_{-\pi}^{\thd(x_2,\s)} \bigl| \nbs_1^2 f \cdot (\nbs_1 f)^2  \bigr|(x_1,x_2,\s) {\rm d}x_1 {\rm d}x_2 {\rm d}\s
\notag\\
&\leq   \|f\|_{L^\infty_{x,\s}} \snorm{\nbs_1 f}_{L^4_{x,\s}}^2 \snorm{\nbs_1^2 f}_{L^2_{x,\s}} \,.
\end{align*}
Combining the two estimates above shows that~\eqref{eq:waiting:for:the:worms:a} remains unchanged, i.e. 
\begin{equation*}
\snorm{\nbs_1 f}_{L^4_{x,\s}}^4 
\leq 4 \norm{f}_{L^\infty_{x,\s}} \snorm{\nbs_1^2 f}_{L^2_{x,\s}} \snorm{\nbs_1 f}_{L^4_{x,\s}}^2
\,.
\end{equation*}

A similar argument, using the definition of~$\nbs_2^*$ in~\eqref{adjoint-2-H}, and the bounds for $\Qr_2$ and $\Qb$ from~\eqref{eq:Q:all:bbq-H}, shows that
\begin{align*}
\snorm{\nbs_2 f}_{L^4_{x,\s}}^4  
&\leq 3 \snorm{f}_{L^\infty_{x,\s}} \snorm{\nbs_2 f}_{L^4_{x,\s}}^2   \snorm{\nbs_2^2 f}_{L^2_{x,\s}} 
+ 6(1+\alpha)  \snorm{f}_{L^\infty_{x,\s}} \snorm{\nbs_2 f}_{L^3_{x,\s}}^3 
+ 11 \eps \snorm{f}_{L^\infty_{x,\s}} \sup_{\s \in [\sin,\sfin]} \snorm{\nbs_2 f(\cdot,\s)}_{L^3_{x}}^3 
\notag\\
&\qquad 
+ \snorm{f}_{L^\infty_{x,\s}} \int_{\sin}^{\sfin} \int_{\TT} \sabs{\nbs_2 \thd(x_2,\s)} \sabs{\nbs_2 f(\thd(x_2,\s),x_2,\s)}^3 {\rm d}x_2 {\rm d}\s
\,.
\end{align*}
The first three terms on the right side of the above estimate already were contained in~\eqref{eq:waiting:for:the:worms:b}, and they are treated in the same way, by appealing to~\eqref{eq:got:my:mojo:workin} with $p=3$ and $q=2$. In order to estimate the last term in the above expression, we make two observations. First, from~\eqref{little-theta} and~\eqref{nb-s-H} we deduce
\begin{equation*}
\nbs_2 \thd(x_2,\s) = \tfrac{\Qb_2(x_2,\s) + \p_2 \bar \Thd(\thd(x_2,\s),x_2)}{-\p_1 \bar \Thd(\thd(x_2,\s),x_2)} 
\,.
\end{equation*}
Second, upon making the change of variables $\s\mapsto x_1$ via $\s = \bar\Thd(x_1,x_2)$, we have from~\eqref{little-theta} that 
\begin{equation*}
{\rm d}\s = \p_1 \bar \Thd(x_1,x_2) {\rm d}x_1 = \p_1 \bar \Thd(\thd(x_2,\s),x_2) {\rm d}x_1.
\end{equation*}
Hence, changing variables, recalling the definition~\eqref{L2x-composite}, and appealing to the bounds~\eqref{eq:Qb:bbq-H} and~\eqref{boot-p2Theta}, we deduce that 
\begin{align*}
 &\snorm{f}_{L^\infty_{x,\s}} \int_{\sin}^{\sfin} \int_{\TT} \sabs{\nbs_2 \thd(x_2,\s)} \sabs{\nbs_2 f(\thd(x_2,\s),x_2,\s)}^3 {\rm d}x_2 {\rm d}\s
 \notag\\
&\qquad = 
 \snorm{f}_{L^\infty_{x,\s}}  \int_{\TT} \int_{\tilde{\mathfrak{X}}_1^-(x_2,0)}^{\tilde{\mathfrak{X}}_1^+(x_2,0)} \sabs{\Qb_2(x_2,\bar\Thd(x_1,x_2)) + \p_2 \bar \Thd(x_1,x_2)} \sabs{\nbs_2 f(x_1,x_2,\bar\Thd(x_1,x_2))}^3 {\rm d}x_2 {\rm d}x_1
 \notag\\
 &\qquad \les \eps \snorm{f}_{L^\infty_{x,\s}}  \snorm{\nbs_2 f (\cdot,\Thd(\cdot,0))}^3_{L^3_x}
 \,.
\end{align*}
Then, re-doing the computation leading to~\eqref{eq:chris:rea:3} (in $L^3$ instead of $L^2$), using~\eqref{L2sx-fubini} (both with $L^2$ and with $L^4$, instead of just $L^2$), and re-doing the computation leading up to~\eqref{lazy-ass1} for $r=0$ (in $L^3$ instead of $L^2$), we deducethat 
\begin{align*}
\snorm{\nbs_2 f (\cdot,\Thd(\cdot,0))}_{L^3_x}^3
&\les \snorm{\nbs_2 f (\cdot,\Thd(\cdot,\sin))}_{L^3_x}^3  
+ \snorm{\nbs_2 f (\cdot,\sin)}_{L^3_x}^3 
+ \tfrac{1}{\eps}  \snorm{\nbs_2 f (\cdot,\Thd(\cdot,\cdot))}_{L^4_{x,\s}}^2  \snorm{\nbs_\s \nbs_2 f (\cdot,\Thd(\cdot,))}_{L^2_{x,\s}} 
\notag\\
&\les \snorm{\nbs_2 f (\cdot,\sin)}_{L^3_x}^3 
+ \tfrac{1}{\eps}  \snorm{\nbs_2 f}_{L^4_{x,\s}}^2  \snorm{\nbs_\s \nbs_2 f}_{L^2_{x,\s}} 
\,.
\end{align*}
Combining the bounds obtained in the five estimates above leads to 
\begin{align*}
\snorm{\nbs_2 f}_{L^4_{x,\s}}^4  
&\les  \snorm{f}_{L^\infty_{x,\s}} \snorm{\nbs_2 f}_{L^4_{x,\s}}^2   \snorm{\nbs_2^2 f}_{L^2_{x,\s}} 
+ \snorm{f}_{L^\infty_{x,\s}} \snorm{\nbs_2 f}_{L^4_{x,\s}}^2   \snorm{\nbs_\s \nbs_2 f}_{L^2_{x,\s}}
+   \eps^2 \snorm{f}_{L^\infty_{x,\s}}  \snorm{\nbs_2 f(\cdot,\sin)}_{L^\infty_{x}}^3 
\,,
\end{align*}
which precisely matches~\eqref{eq:waiting:for:the:worms}.
  
It remains to consider the modifications required to the bound~\eqref{eq:waiting:for:the:worms:c}. 
Recalling the definition of~$\nbs_\s^*$ in~\eqref{adjoint-s-H}, it is clear that only one new term emerges, on the right side of~\eqref{eq:waiting:for:the:worms}, which is
\begin{equation*}
\snorm{f}_{L^\infty_{x,\s}} \int_{\sin}^{\sfin} \int_{\TT} \sabs{\nbs_\s \thd(x_2,\s)} \sabs{\nbs_\s f(\thd(x_2,\s),x_2,\s)}^3 {\rm d}x_2 {\rm d}\s 
\,.
\end{equation*}
In a similar fashion to the $\nbs_2$ analysis above, we use that from~\eqref{little-theta} and~\eqref{nb-s-H} we deduce
\begin{equation*}
\nbs_\s \thd(x_2,\s) = \tfrac{\eps \Qd(x_2)}{-\p_1 \bar \Thd(\thd(x_2,\s),x_2)} 
\,.
\end{equation*}
The same change of variables as the one described above, followed up with an application of~\eqref{L2x-composite},  \eqref{eq:chris:rea:3}, \eqref{L2sx-fubini}, \eqref{lazy-ass1}, all adapted to the $L^3$ instead of the $L^2$ setting, shows that $\snorm{\nbs_\s f}_{L^4_{x,\s}}^4$ obeys a bound which is the same as~\eqref{eq:waiting:for:the:worms:c}.

Thus, we have shown that for the upstream geometry of Section~\ref{sec:upstreammaxdev} the bound~\eqref{eq:waiting:for:the:worms} still holds {\em as is}, albeit with a slightly more elaborate proof. This results in an unchanged statement of Lemma~\ref{lem:time:interpolation}.

In summary, {\em all the main bounds} obtained in this Appendix, namely Lemmas~\ref{lem:anisotropic:sobolev}, \ref{lem:time:interpolation}, \ref{lem:Moser:tangent}, \ref{lem:comm:tangent}, \ref{lem:comm:bdd}, and Remarks~\ref{rem:app:useful?},~\ref{rem:shit:is:similar}, hold {\em as is} for the $(x,\s)$ coordinates from Section~\ref{sec:upstreammaxdev}.

\section{Transport bounds}
\label{sec:app:transport}

The purpose of this appendix is to provide space-time $L^\infty$ bounds for objects which are transported and stretched along the $\lambda_1$ and $\lambda_2$ characteristics. As with Appendix~\ref{app:functional}, the estimates in this Appendix are written for the flattening map $\mathfrak{q}$ defined in Section~\ref{sec:remapping}, and utilized in Sections~\ref{sec:formation:setup}--\ref{sec:sixth:order:energy}.  In Sections~\ref{app:downstream:Lp} and~\ref{app:upstream:Lp} below, we show that the same analytical framework developed here applies {\em as is} to the remapped domains from Sections~\ref{sec:downstreammaxdev} and~\ref{sec:upstreammaxdev}.

Recall from \eqref{eq:Ab:nn:alt}, \eqref{eq:Ab:tt:alt}, and \eqref{vort-t} that in the ALE coordinates corresponding to the fast acoustic spped $\lambda_3$, the transport operator associated with the wave-speed $\lambda_2$ takes the form $ (\p_t +V \p_2) +  \alpha g^{-\frac 12} h,_2 \p_2 - \alpha \Sigma \Jg^{-1} \p_1$. Also, from~\eqref{eq:Zb:nn:alt} and \eqref{eq:Zb:tt:alt} we see that the $\lambda_1$ transport operator takes the form $(\p_t +V \p_2) +  2 \alpha g^{-\frac 12} h,_2 \p_2 - 2 \alpha \Sigma \Jg^{-1} \p_1$. 
For simplicity of the presentation, we only present the details of the $\lambda_2$ analysis (the change of an $\alpha$ to a $2\alpha$ does not affect our final estimate). 

We  recall from Section~\ref{sec:remapping} that the remapping of the space-time $\mathcal{P} \to \TT^2 \times [0,\eps]$, and the associated change of unknown $f(x,t) = \tilde f(x,\s)$, gives  $\|f\|_{L^\infty_{x,t}} = \|\tilde f\|_{L^\infty_{x,\s}}$. Under this change of coordinates, if a function $f$ that solves 
\begin{equation}
\bigl(
\Jg (\p_t + V\p_2) +  \alpha \Sigma g^{-\frac 12} \Jg h,_2 \p_2 - \alpha \Sigma \p_1 
\bigr) f
= 
q
\,,
\label{eq:abstract:0}
\end{equation}
then according to \eqref{eq:xt:xs:chain:rule} and \eqref{the-time-der} we have that the function $\tilde f$ solves 
\begin{equation}
  \bigl(
 \Jg (\Q \p_\s + V \p_2) +     \alpha \Sigma g^{-\frac 12} \nbs_2 h \Jg  \nbs_2  - \alpha   \Sigma  \p_1 
\bigr) 
\tilde f
= 
  \tilde q
\,.
\label{eq:abstract}
\end{equation}
As discussed in Remark~\ref{rem:no:tilde}, we drop the tildes on all the remapped functions of the unknown $(x,\s)$ present in \eqref{eq:abstract}. Our goal is to prove the $L^\infty$ bounds directly in $(x,\s)$ coordinates, working with \eqref{eq:abstract}. Clearly, the bounds established will also hold if $\alpha$ is replaced by $2\alpha$, i.e. the $\lambda_2$ transport operator is replaced by the $\lambda_1$ transport operator. The main result is Proposition~\ref{prop:transverse:bounds} below, whose proof has an equally useful consequence, see Corollary~\ref{cor:transverse:bounds}.
\begin{proposition}
\label{prop:transverse:bounds}
Assume that $f$ is a smooth solution \eqref{eq:abstract}, and that $q$ is bounded. Then, assuming the bootstrap inequalities~\eqref{bootstraps} hold, and that $\eps$ is sufficiently small in terms of $\alpha,\kappa_0$, and $\Cdata$, we have
\begin{equation}
\snorm{f}_{L^\infty_{x,\s}} 
\leq 4 e^{18} \snorm{f (\cdot,0)}_{L^\infty_x} 
+ \tfrac{20 \eps}{\alpha}  e^{18} \snorm{q}_{L^\infty_{x,\s}} 
\,.
\label{eq:abstract:transport}
\end{equation}
\end{proposition}
 
\begin{remark}
\label{rem:transverse:bounds}
In many applications of estimate~\eqref{eq:abstract:transport} it may be  convenient to appeal to  the fundamental theorem of calculus in time and \eqref{eq:Steve:punch} to bound the $L^\infty$ norm of $q$, resulting in the estimate 
\begin{equation}
\snorm{f}_{L^\infty_{x,\s}} 
\les \snorm{f (\cdot,0)}_{L^\infty_x} 
+ 
\tfrac{\eps}{\alpha}   \snorm{q(\cdot,0)}_{L^\infty_x}
+
\tfrac{1}{\alpha} \snorm{\nbs^2\nbs_1 q}_{L^2_{x,\s}} 
\,.
\label{eq:abstract:transport:1a}
\end{equation}
\end{remark}

\begin{proof}[Proof of Proposition~\ref{prop:transverse:bounds}]
First, we note that  \eqref{eq:abstract} (in which as usual we drop the tildes) implies 
\begin{equation*}
\tfrac{\Jg}{\Sigma^p} (\Q \p_\s + V \p_2)( f^p )
+ \tfrac{\alpha g^{-\frac 12} \Jg \nbs_2 h}{\Sigma^{p-1}} \nbs_2 (f^p)
- \tfrac{\alpha}{\Sigma^{p-1}} \p_1 (f^p)
=
\tfrac{p}{\Sigma^p}  f^{p-1} q
\,.
\end{equation*}
Integrating the above expression in space-time, and appealing to \eqref{eq:adjoints} we deduce  
\begin{equation}
\snorm{(\Q \Jg)^{\!\frac 1p} \Sigma^{-1} f(\cdot,\s)}_{L^p_x}^p
- \snorm{(\Q \Jg)^{\!\frac 1p} \Sigma^{-1} f(\cdot,0)}_{L^p_x}^p
=
\int_0^{\s}\!\!\! \int \tfrac{1}{\Sigma^p} f^p \mathsf{G}
+
 \int \Qb_2 \tfrac{\alpha g^{-\frac 12} \Jg \nbs_2 h}{\Sigma^{p-1}}  f^p \Bigr|_{\s}
+ 
 \int_0^{\s}\!\!\! \int \tfrac{p}{\Sigma^p} f^{p-1}  q ,
\label{eq:GSTQ:1}
\end{equation}
where, we have denoted
\begin{equation}
\mathsf{G} 
:=  
\Sigma^p \Bigl( \Qr_\s + \nbs_2 V - V \Qr_2 + (\Q \p_\s + V \p_2) \Bigr) \tfrac{\Jg}{\Sigma^p}
- 
\Sigma^p (\Qr_2 - \nbs_2) \Bigl(\tfrac{\alpha g^{-\frac 12} \Jg \nbs_2 h}{\Sigma^{p-1}}\Bigr)
- 
\Sigma^p \p_1 \Bigl( \tfrac{\alpha}{\Sigma^{p-1}} \Bigr)
\,.
\label{eq:GSTQ:2}
\end{equation}
Using \eqref{Jg-evo-s}, \eqref{Sigma0-ALE-s},   and  \eqref{grad-Sigma},  we may rewrite
\begin{align*}
\mathsf{G} 
&= \Jg \bigl( \Qr_\s + \nbs_2 V - V \Qr_2 \bigr)   
+ \tfrac{1+\alpha}{2} \Jg \Wbn + \tfrac{1-\alpha}{2} \Jg \Zbn 
+ p \alpha  \Jg (\Zbn + \Abt)
\notag\\
&\qquad 
- \alpha  \Qr_2  \Sigma  g^{-\frac 12} \Jg \nbs_2 h 
+
\alpha \Sigma   \nbs_2  \bigl( g^{-\frac 12} \Jg \nbs_2 h \bigr)
- (p-1) 
\alpha g^{-\frac 12} \Jg \nbs_2 h  \nbs_2 \Sigma  
\notag\\
&\qquad 
+ (p-1) \alpha \bigl(\tfrac{1}{2} (\Jg\Wbn -\Jg\Zbn)  + \tfrac{1}{2} \Jg \nbs_2 h (\Wbt -\Zbt) \bigr)
\notag\\
&=  
\bigl( \tfrac{1}{2} +  p \tfrac{\alpha}{2} \bigr) \Jg \Wbn
+  \Jg \bigl( \Qr_\s + \nbs_2 V - V \Qr_2 \bigr)   
+\bigl( \tfrac{1}{2}  + p \tfrac{\alpha}{2} \bigr) \Jg \Zbn 
\notag\\
&\qquad 
- \alpha  \Qr_2  \Sigma  g^{-\frac 12} \Jg \nbs_2 h 
+
\alpha \Sigma   \nbs_2  \bigl( g^{-\frac 12} \Jg \nbs_2 h \bigr)
+ p \alpha  \Jg  \Abt 
\,.
\end{align*}
Using \eqref{eq:Q:all:bbq}, \eqref{eq:signed:Jg}, and the bootstrap inequalities \eqref{bootstraps}, we deduce for $p\geq 1$ the pointwise bound
\begin{equation}
\mathsf{G}
\leq - p  \tfrac{9  \alpha  }{20 \eps }   +  \tfrac{17 p + 2 \cdot 250^2}{\eps} \Q \Jg + \Cn p
\label{eq:bingo:de:bongo:G}
\end{equation}
where $\Cn = \Cn(\alpha,\kappa_0,\Cdata)$ is independent of $\eps$ and $p$.
Moreover, by appealing to the same bounds we have that 
\begin{equation}
\sabs{\Qb_2 \tfrac{\alpha g^{-\frac 12} \Jg \nbs_2 h}{\Sigma^{p-1}} }
\leq \Cn \eps^2 \tfrac{\Jg \Q}{\Sigma^p}
\label{eq:bingo:de:bongo:Q}
\,.
\end{equation}
Therefore, upon taking $\eps$ to be sufficiently small, solely in terms of $\alpha,\kappa_0$, and $\Cdata$, we deduce from \eqref{eq:GSTQ:1}, \eqref{eq:bingo:de:bongo:G}, and \eqref{eq:bingo:de:bongo:Q} that
\begin{align*}
&\tfrac{17}{18} \snorm{(\Q\Jg)^{\!\frac 1p} \Sigma^{-1} f(\cdot,\s)}_{L^p_x}^p
+  \tfrac{2 \alpha p}{5 \eps}  \int_0^{\s} \snorm{\Sigma^{-1} f(\cdot,\s')}_{L^p_{x}}^p {\rm d} \s'
\notag\\
&\leq \snorm{(\Q\Jg)^{\!\frac 1p} \Sigma^{-1} f(\cdot,0)}_{L^p_x}^p
+ \tfrac{17 p + 2 \cdot 250^2}{\eps} \! \int_0^{\s} \snorm{(\Q \Jg)^{\!\frac 1p} \Sigma^{-1} f(\cdot,\s')}_{L^p_{x}}^p {\rm d} \s'
+ p \! \int_0^{\s} \snorm{\Sigma^{-1} f(\cdot,\s')}_{L^p_{x}}^{p-1} \snorm{\Sigma^{-1} q(\cdot,\s')}_{L^p_{x}}{\rm d}\s'
.
\end{align*}
Next, we apply the  $\eps$-Young  inequality  in the form
\begin{equation*}
 p \snorm{\Sigma^{-1} f}_{L^p_x}^{p-1} \snorm{\Sigma^{-1} q}_{L^p_x}
\leq 
 \tfrac{ p \alpha}{5\eps}     \snorm{\Sigma^{-1} f}_{L^p_x}^p
 + (\tfrac{5\eps}{ \alpha})^{p-1} \snorm{\Sigma^{-1} q}_{L^p_x}^p
\end{equation*}
to deduce
\begin{align}
&\snorm{(\Q\Jg)^{\!\frac 1p} \Sigma^{-1} f(\cdot,\s)}_{L^p_x}^p
+ \tfrac{\alpha p}{5 \eps}  \int_0^{\s} \snorm{\Sigma^{-1} f(\cdot,\s')}_{L^p_{x}}^p {\rm d} \s'
\notag\\
&\leq 2\snorm{(\Q\Jg)^{\!\frac 1p} \Sigma^{-1} f(\cdot,0)}_{L^p_x}^p
+ \tfrac{18 p + 2 \cdot 250^2}{\eps}  \int_0^{\s} \snorm{(\Q \Jg)^{\!\frac 1p} \Sigma^{-1} f(\cdot,\s')}_{L^p_{x}}^p {\rm d} \s'
+2 (\tfrac{5\eps}{ \alpha})^{p-1}   \int_0^{\s}   \snorm{\Sigma^{-1} q}_{L^p_x}^p {\rm d}\s'
\,. 
\label{eq:stupid:intermediary}
\end{align}
Next, we ignore the damping term on the left side of \eqref{eq:stupid:intermediary} and apply Gr\"onwall on the time interval $[0,\eps]$ to obtain
\begin{equation*}
\sup_{\s\in[0,\eps]}
\snorm{(\Q\Jg)^{\!\frac 1p} \Sigma^{-1} f(\cdot,\s)}_{L^p_x}^p
\leq  
2 e^{18p+ 2 \cdot 250^2} \snorm{(\Q\Jg)^{\!\frac 1p} \Sigma^{-1} f(\cdot,0)}_{L^p_x}^p
+2 e^{18p+ 2 \cdot 250^2} (\tfrac{5\eps}{ \alpha})^{p-1}    \snorm{\Sigma^{-1} q}_{L^p_{x,\s}}^p
\,. 
\end{equation*}
The second to last step is to take $p^{th}$ roots, use that $(a^p+b^p)^{\frac 1p} \leq a + b$, and  deduce
\begin{equation}
\sup_{\s\in[0,\eps]}
\snorm{(\Q \Jg)^{\!\frac 1p} \Sigma^{-1} f(\cdot,\s)}_{L^p_x}
\leq  
2^{\frac 1p} e^{18 + \frac{2 \cdot 250^2}{p}} \snorm{(\Q \Jg)^{\!\frac 1p} \Sigma^{-1} f(\cdot,0)}_{L^p_x}
+2^{\frac 1p} e^{18+ \frac{2 \cdot 250^2}{p}} (\tfrac{5\eps}{ \alpha})^{\frac{p-1}{p}}    \snorm{\Sigma^{-1} q}_{L^p_{x,\s}}
\,. 
\label{eq:very:stupid:intermediary}
\end{equation}
The last step is to pass $p\to \infty$. Via the dominated convergence theorem   
\begin{equation*}
 \snorm{\Sigma^{-1} f(\cdot,\s)}_{L^\infty_{x,\s}}
\leq  
 e^{18} \snorm{\Sigma^{-1} f (\cdot,0)}_{L^\infty_x} 
+ \tfrac{5\eps}{ \alpha}  e^{18}
\snorm{\Sigma^{-1} q}_{L^\infty_{x,\s}}
\,. 
\end{equation*}
The proof of \eqref{eq:abstract:transport} now follows from the above estimate and \eqref{bs-Sigma}.
\end{proof}

In a number on instances in the paper we need to bound the gradient of the function $f$ which solves~\eqref{eq:abstract}. In these cases the new forcing term is not merely the gradient of the function $q$, but there is a factor which is containing the gradient of $f$ itself, but without the necessary power of $\Jg$ in front of it. That is, we need to consider equations like
\begin{equation}
\bigl(
\Jg (\p_t + V \p_2) + \alpha  \Sigma g^{-\frac 12} \Jg h,_2  \p_2 - \alpha \Sigma \p_1 
\bigr) f
= m f + q
\,,
\label{eq:abstract:2} 
\end{equation}
where as before the function $q$ is bounded on $\mathcal{P}$, and we also assume that the function $m$ is bounded on $\mathcal{P}$. As before, in $(x,\s)$ variables the above equation transforms to 
\begin{equation}
  \bigl(
 \Jg (\Q \p_\s + V \p_2) +     \alpha \Sigma g^{-\frac 12} \nbs_2 h \Jg  \nbs_2  - \alpha   \Sigma  \p_1 
\bigr)
\tilde f
= 
\tilde m \tilde f +  \tilde q
\,.
\label{eq:abstract:2:xs}
\end{equation}
As usual, we drop the tilde from the unknowns, so that $(\tilde f,\tilde m,\tilde q)$ become $(f,m,q)$.
By a slight modification of the proof of Proposition~\ref{prop:transverse:bounds}, we obtain:
\begin{corollary}
\label{cor:transverse:bounds} 
Assume that $f$ is a smooth solution \eqref{eq:abstract:2:xs}, and that $q$ and $m$ are bounded. Then, assuming the bootstrap inequalities~\eqref{bootstraps} hold, and that $\eps$ is sufficiently small in terms of $\alpha,\kappa_0$, and $\Cdata$, we have
\begin{equation}
\snorm{f}_{L^\infty_{x,\s}} 
\leq (4 e^{18})^\beta \snorm{f (\cdot,0)}_{L^\infty_x} 
+ \eps \tfrac{20}{\alpha \beta} (4 e^{18})^\beta  \snorm{q}_{L^\infty_{x,\s}} 
\,.
\label{eq:abstract:transport:2}
\end{equation}
where the parameter $\beta$ is given cf.~\eqref{eq:bingo:de:bongo:beta} as $\beta = \max \{ \tfrac{20}{\alpha}\eps \|m^+\|_{L^\infty_{x,\s}}  ,  1\}$, and $m^+ = \max\{ m, 0\}$.
\end{corollary}
\begin{proof}[Proof of Corollary~\ref{cor:transverse:bounds}]
The main difference is that in the energy formed on the left side of \eqref{eq:GSTQ:1}, we need to replace $\Sigma^{-1}$ by $\Sigma^{-\beta}$ for a $\beta>0$ that is to be chosen in terms of $\norm{m^+}_{L^\infty_{x,\s}}$. With this change, \eqref{eq:GSTQ:1} becomes 
\begin{align}
& \snorm{(\Q\Jg)^{\!\frac 1p} \Sigma^{-\beta} f(\cdot,\s)}_{L^p_x}^p
- \snorm{(\Q\Jg)^{\!\frac 1p} \Sigma^{-\beta} f(\cdot,0)}_{L^p_x}^p
\notag\\
&\qquad =
\int_0^{\s}\!\!\! \int \tfrac{1}{\Sigma^{\beta  p}} f^p   \mathsf{G}  
+ \int \Qb_2 \tfrac{\alpha g^{-\frac 12} \Jg \nbs_2 h}{\Sigma^{\beta  p-1}}   f^p \Bigr|_{\s}
+ \int_0^{\s}\!\!\! \int \tfrac{p}{\Sigma^{\beta p}}  f^{p-1}   q     
\,,
\label{eq:bingo:de:bongo:energy}
\end{align}
where 
\begin{align*}
\mathsf{G} 
&=p m + \Sigma^{\beta p} \Bigl( \Qr_\s + \nbs_2 V - V \Qr_2 + (\Q \p_\s + V \p_2) \Bigr) \tfrac{\Jg}{\Sigma^{\beta p}}
- 
\Sigma^{\beta p} (\Qr_2 - \nbs_2) \Bigl(\tfrac{\alpha g^{-\frac 12} \Jg \nbs_2 h}{\Sigma^{\beta p-1}}\Bigr)
- 
\Sigma^{\beta p} \p_1 \Bigl( \tfrac{\alpha}{\Sigma^{\beta p-1}} \Bigr)
\notag\\
&=  pm +
\bigl( \tfrac{1}{2} +  \beta p \tfrac{\alpha}{2} \bigr) \Jg \Wbn
+  \Jg \bigl( \Qr_\s + \nbs_2 V - V \Qr_2 \bigr)   
+\bigl( \tfrac{1}{2}  + \beta p \tfrac{\alpha}{2} \bigr) \Jg \Zbn 
\notag\\
&\qquad 
- \alpha  \Qr_2  \Sigma  g^{-\frac 12} \Jg \nbs_2 h 
+
\alpha \Sigma   \nbs_2  \bigl( g^{-\frac 12} \Jg \nbs_2 h \bigr)
+ p \alpha \beta \Jg  \Abt 
\,.
\end{align*}
Similarly to \eqref{eq:bingo:de:bongo:G} we have the upper bound
\begin{equation*}
\mathsf{G} 
\leq  p \|m^+\|_{L^\infty_{x,\s}}
 -  \tfrac{9 \alpha }{20\eps } \beta p 
 + \bigl( 17  \beta p + 2 \cdot 250^2\bigr) \tfrac{1}{\eps} \Q \Jg + \Cn p \brak{\beta}
\end{equation*}
where $\Cn = \Cn(\alpha,\kappa_0,\Cdata)$ is independent of $\eps$, $p$, and $\beta$.
Choosing 
\begin{equation}
\beta = \max\bigl\{ \tfrac{20}{\alpha}\eps \|m^+\|_{L^\infty_{x,\s}}  ,  1 \bigr\}\,,
\label{eq:bingo:de:bongo:beta}
\end{equation}
and choosing $\eps$ to be sufficiently small, we arrive at 
\begin{equation}
\mathsf{G} 
\leq   
 -  \tfrac{2 \alpha\beta }{5 \eps }  p 
 +   \bigl( 17  \beta p + 2 \cdot 250^2\bigr) \tfrac{1}{\eps} \Q \Jg
 \label{eq:bingo:de:bongo:new:G}
 \,.
\end{equation}
From \eqref{eq:bingo:de:bongo:Q}, \eqref{eq:bingo:de:bongo:energy}, and \eqref{eq:bingo:de:bongo:new:G}, we deduce
\begin{align}
&\snorm{(\Q \Jg)^{\!\frac 1p} \Sigma^{-\beta} f(\cdot,\s)}_{L^p_x}^p
+ \tfrac{\alpha \beta p}{5 \eps}  \int_0^{\s} \snorm{\Sigma^{-\beta} f(\cdot,\s')}_{L^p_{x}}^p {\rm d} \s'
\notag\\
&\leq 2\snorm{(\Q\Jg)^{\!\frac 1p} \Sigma^{-\beta} f(\cdot,0)}_{L^p_x}^p
+ \tfrac{18 \beta p + 2 \cdot 250^2}{\eps}  \int_0^{\s} \snorm{(\Q\Jg)^{\!\frac 1p} \Sigma^{-\beta} f(\cdot,\s')}_{L^p_{x}}^p {\rm d} \s'
+2 (\tfrac{5\eps}{ \alpha \beta})^{p-1}   \int_0^{\s}   \snorm{\Sigma^{-\beta} q}_{L^p_x}^p {\rm d}\s'
\,. 
\label{eq:stupid:intermediary:2}
\end{align}
We then use Gr\"onwall's inequality, take $p^{th}$ roots and pass $p\to \infty$ to deduce
\begin{equation*}
 \snorm{\Sigma^{-\beta} f(\cdot,\s)}_{L^\infty_{x,\s}}
\leq  
 e^{18 \beta} \snorm{\Sigma^{-\beta} f (\cdot,0)}_{L^\infty_x} 
+ \tfrac{5\eps}{ \alpha\beta}  e^{18 \beta}
\snorm{\Sigma^{-\beta} q}_{L^\infty_{x,\s}}
\,. 
\end{equation*}
By further appealing to \eqref{bs-Sigma}, the proof is completed.
\end{proof}

\subsection{The flattening map from Section~\ref{sec:downstreammaxdev}}
\label{app:downstream:Lp}
In the flattened geometry utilized in Section~\ref{sec:downstreammaxdev}, the statements of Proposition~\ref{prop:transverse:bounds}, Remark~\ref{rem:transverse:bounds}, and Corollary~\ref{cor:transverse:bounds} remain unchanged. The changes to the proofs of these statements are entirely due to the fact that the $-\alpha \Sigma \p_1 \tilde f$ term in \eqref{eq:abstract} and \eqref{eq:abstract:2:xs} needs to be rewritten as $-\frac{\alpha}{\eps} \Sigma \nbs_1 \tilde f$, and instead of using that $\p_1^* = -\p_1$, we need to use the adjoint formula  $\frac{1}{\eps} \nbs_1^* = - \frac{1}{\eps} \nbs_1^* + \Qr_1 - \Qb_1 \delta_{\s}$ (cf.~\eqref{adjoint-1-P}). 

Let us discuss in detail the changes this causes to the proof of Proposition~\ref{prop:transverse:bounds}. Identities~\eqref{eq:GSTQ:1}--\eqref{eq:GSTQ:2} become 
\begin{align*}
&\snorm{(\Q \Jg)^{\!\frac 1p} \Sigma^{-1} f(\cdot,\s)}_{L^p_x}^p
- \snorm{(\Q \Jg)^{\!\frac 1p} \Sigma^{-1} f(\cdot,0)}_{L^p_x}^p
+   \int \Qb_1 \tfrac{\alpha}{\Sigma^{p-1}}  f^p \Bigr|_{\s}
\notag\\
&\qquad 
=
\int_0^{\s}\!\!\! \int \tfrac{1}{\Sigma^p} f^p \mathsf{G}
+
 \int \Qb_2 \tfrac{\alpha g^{-\frac 12} \Jg \nbs_2 h}{\Sigma^{p-1}}  f^p \Bigr|_{\s}
+ 
 \int_0^{\s}\!\!\! \int \tfrac{p}{\Sigma^p} f^{p-1}  q ,
\end{align*}
with 
\begin{equation*}
\mathsf{G} 
:=  
\Sigma^p \Bigl( \Qr_\s + \nbs_2 V - V \Qr_2 + (\Q \p_\s + V \p_2) \Bigr) \tfrac{\Jg}{\Sigma^p}
- 
\Sigma^p (\Qr_2 - \nbs_2) \Bigl(\tfrac{\alpha g^{-\frac 12} \Jg \nbs_2 h}{\Sigma^{p-1}}\Bigr)
- 
\tfrac{1}{\eps} \Sigma^p \nbs_1 \Bigl( \tfrac{\alpha}{\Sigma^{p-1}} \Bigr)
+
\alpha \Qr_1 \Sigma
\,.
\end{equation*}
The new term arising due to $\Qb_1$ is a helpful term that because \eqref{eq:Qrs1:bbq-DS} gives $\Qb_1 \geq 0$; we however ignore this helpful term. The only new contribution to $\mathsf{G}$ arising from the $\Qr_1$ term may be bounded from above using \eqref{eq:Qr1:bbq-DS} and \eqref{boots-PP} by $5 \alpha \kappa_0 \eps^{-1}$; a bound which is $p$-independent. As such, \eqref{eq:bingo:de:bongo:G} is replaced by 
\begin{equation}
\mathsf{G}
\leq - p  \tfrac{9  \alpha  }{20 \eps } +   \tfrac{5\alpha \kappa_0}{\eps} +  \tfrac{17 p + 2 \cdot 250^2}{\eps} \Q \Jg + \Cn p 
\leq - p   \tfrac{2  \alpha  }{5 \eps }   +  \tfrac{17 p + 2 \cdot 250^2}{\eps} \Q \Jg   
\,,
\label{eq:fall:thunder:storm}
\end{equation}
whenever $p \geq 200 \alpha \kappa_0$, and $\eps$ is taken to be sufficiently small in terms of $\alpha,\kappa_0$, and $\Cdata$. Since we pass $p\to \infty$, this restriction on a lower bound for $p$ that depends on $\alpha$ and $\kappa_0$ is irrelevant and may be assumed without loss of generality. Since none of the subsequent bounds are altered, we proceed as in the proof of Proposition~\ref{prop:transverse:bounds}, and upon passing $p\to \infty$ deduce that the bound  \eqref{eq:abstract:transport} remains unchanged. The changes to the proof of Corollary~\ref{cor:transverse:bounds} are   identical, and we omit these details.
 
\subsection{The flattening map from Section~\ref{sec:upstreammaxdev}}
\label{app:upstream:Lp}
In the flattened geometry utilized in Section~\ref{sec:upstreammaxdev}, the statements of 
Proposition~\ref{prop:transverse:bounds}, Remark~\ref{rem:transverse:bounds}, and Corollary~\ref{cor:transverse:bounds} remain unchanged. The changes to the proofs are due to the fact that every $\s$-time-slice of the spacetime $\Hdm$ has a right boundary at $x_1 = \thd(x_2,\s)$. These modifications were already discussed in the proof of Proposition~\ref{prop:vort:H6-H}, but for convenience we repeat here the main parts of the argument. 

For simplicity, we focus on the changes required by the proof of Proposition~\ref{prop:transverse:bounds}. Moreover, since the weight function $\JJ$ solves a PDE (see~\eqref{JJ-formula}) which is a $\dl$-modification of the $\lambda_1$ transport operator, we focus on the (potentially) more difficult case when the function $f$ solves a forced  $\lambda_1$ transport, namely
\begin{equation}
\bigl(
\Jg (\Q \p_\s + V \p_2) +   2\alpha \Sigma g^{-\frac 12} \nbs_2 h \Jg  \nbs_2  - 2\alpha   \Sigma  \p_1 
\bigr) 
 f
= 
 q
\,.
\label{eq:abstract:lambda:1}
\end{equation}
Here we recall that $\nbs$ is defined by~\eqref{nb-s-H}. 
The above PDE is considered instead of \eqref{eq:abstract}. When $2\alpha$ is replaced by $\alpha$ in~\eqref{eq:abstract:lambda:1} the same bounds will hold.
 
Let $r>0$ be arbitrary. We may multiply \eqref{eq:abstract:lambda:1} by $p \Sigma^{-p} \JJ^{r} f^{p-1}$ and  deduce that
\begin{equation*}
\tfrac{\Jg \JJ^{r}}{\Sigma^p} (\Q \p_\s + V \p_2)( f^p )
+ \tfrac{2\alpha g^{-\frac 12} \Jg \JJ^{r}\nbs_2 h}{\Sigma^{p-1}} \nbs_2 (f^p)
- \tfrac{2\alpha \JJ^{r}}{\Sigma^{p-1}} \p_1 (f^p)
=
\tfrac{p\JJ^{r}}{\Sigma^p}  f^{p-1} q
\,.
\end{equation*}
In analogy to~\eqref{eq:GSTQ:1}--\eqref{eq:GSTQ:2},  by integrating the above identity over the spacetime $\{ (x,\s^\prime) \colon \s^\prime\in[\sin,\s], x_2\in \TT, x_1 \leq \thd(x_2,\s^\prime)\} \subset \Hdm$, using the adjoint formulas in~\eqref{eq:adjoints-H} and the notation for norms from~\eqref{norms-H}, it follows that
\begin{subequations}
\label{eq:Lp:bounds:Appendix:B:15}
\begin{align}
&\snorm{(\Q \Jg \JJ^r )^{\!\frac 1p} \Sigma^{-1} f(\cdot,\s)}_{L^p_x}^p
- \snorm{(\Q \Jg\JJ^r )^{\!\frac 1p} \Sigma^{-1} f(\cdot,\sin)}_{L^p_x}^p
\notag\\
&\qquad 
=
\int_0^{\s}\!\!\! \int \tfrac{\JJ^r}{\Sigma^p} f^p \mathsf{G}
+
 \int \Qb_2 \tfrac{2 \alpha g^{-\frac 12} \Jg\JJ^r \nbs_2 h}{\Sigma^{p-1}}  f^p \Bigr|_{\s}
+ 
 \int_0^{\s}\!\!\! \int \tfrac{p \JJ^r}{\Sigma^p} f^{p-1}  q ,
 \label{eq:Lp:bounds:Appendix:B:15:a}
\end{align}
with 
\begin{equation}
\mathsf{G} 
:=  
\tfrac{\Sigma^p}{\JJ^r} \Bigl(\Qr_\s + \nbs_2 V - V \Qr_2 + (\Q \p_\s + V \p_2) \Bigr) \tfrac{\Jg \JJ^r}{\Sigma^p}
- 
\tfrac{\Sigma^p}{\JJ^r} (\Qr_2 - \nbs_2) \Bigl(\tfrac{2 \alpha g^{-\frac 12} \Jg \JJ^r \nbs_2 h}{\Sigma^{p-1}}\Bigr)
- 
\tfrac{\Sigma^p}{\JJ^r} \p_1 \Bigl( \tfrac{2\alpha \JJ^r}{\Sigma^{p-1}} \Bigr)
\label{eq:Lp:bounds:Appendix:B:15:b}
\,.
\end{equation}
\end{subequations}
In deriving~\eqref{eq:Lp:bounds:Appendix:B:15} we have crucially use the following fact. By~\eqref{JJ-vanish} we have that $\JJ(\thd(x_2,\s^\prime),x_2,\s^\prime) = 0$, and hence, since $r>0$ the boundary terms at $x_1 = \thd(x_2,\s^\prime)$ that should arise when integrating by parts via~\eqref{eq:adjoints-H}, are in fact all equal to $0$. We note that~\eqref{eq:Lp:bounds:Appendix:B:15:a} is identical to~\eqref{eq:GSTQ:1}, while the term $\mathsf{G}$ defined in \eqref{eq:Lp:bounds:Appendix:B:15:b} is nearly identical to the term $\mathsf{G}$ defined in~\eqref{eq:GSTQ:2}; the only difference arises from the differential operators landing on $\JJ^r$. That is, 
\begin{equation*}
\mathsf{G}_{\eqref{eq:Lp:bounds:Appendix:B:15:b}}
=
\mathsf{G}_{\eqref{eq:GSTQ:2}}
+
\tfrac{r}{\JJ} \Bigl( \Jg (\Q\p_\s+V\p_2) \JJ + 2\alpha \Sigma \Jg g^{-\frac 12} \nbs_2 h \nbs_2 \JJ - 2\alpha \Sigma \p_1 \JJ \Bigr)
\,.
\end{equation*}
The key observation is that by combining the PDE solved by $\JJ$ in $\Hdmp$ (recall~\eqref{JJ-formula}), with the derivative bounds satisfied by $\JJ$ in $\Hdm$ (recall~\eqref{qps-JJ-bound}), the pointwise bounds for $\JJ$ in $\Hdmm$ from~\eqref{JJ-and-s-new} and~\eqref{JJ-le-Jg-s}, and the previously obtained estimate~\eqref{eq:bingo:de:bongo:G} (which still holds because we work under teh same pointwise bootstraps), we obtain 
\begin{align*}
\mathsf{G}_{\eqref{eq:Lp:bounds:Appendix:B:15:b}}
&\leq 
\mathsf{G}_{\eqref{eq:GSTQ:2}}
+
\tfrac{r}{\JJ} \Bigl( - \dl \Jg \tfrac{11(1+\alpha)}{25\eps} {\bf 1}_{\Hdmp} 
+ \bigl( - \Jg\tfrac{10(1+\alpha)}{25\eps} + \tfrac{2\alpha \kappa_0}{10^{3}\eps} \bigr)  {\bf 1}_{\Hdmm} \Bigr)
\notag\\
&\leq 
\Bigl(- p  \tfrac{9  \alpha  }{20 \eps }   +  \tfrac{17 p + 2 \cdot 250^2}{\eps} \Q \Jg + \Cn p\Bigr)
-
\tfrac{2r \dl \Jg}{5 \JJ \eps} 
+ \tfrac{3 r \alpha \kappa_0}{10^{3}\eps} \Q \Jg
\notag\\
&\leq 
- p \tfrac{2\alpha}{5\eps} + \tfrac{17 p + 2 \cdot 250^2 +  r \alpha \kappa_0}{\eps} \Q \Jg
\,.
\end{align*}
We see that the modification to the upper bound for $\Jg$ is nearly identical to the one obtained in the previous Section in~\eqref{eq:fall:thunder:storm}. 
In fact, as $r\to 0^+$ (a limit which we will take in a moment), the above obtained upper bound is identical to the one from~\eqref{eq:fall:thunder:storm}, and to the one from the proof of Proposition~\ref{prop:transverse:bounds}.

The bounds for the last two terms on the right side of~\eqref{eq:Lp:bounds:Appendix:B:15:a} remain identical to those obtained in Proposition~\ref{prop:transverse:bounds}. Proceeding as  in the proof of Proposition~\ref{prop:transverse:bounds}, we arrive at the bound corresponding to~\eqref{eq:very:stupid:intermediary}, namely
\begin{align*}
\sup_{\s\in[\sin,\sfin]}
\snorm{(\Q \Jg \JJ^r)^{\!\frac 1p} \Sigma^{-1} f(\cdot,\s)}_{L^p_x}
&\leq  
2^{\frac 1p} e^{18 + \frac{2 \cdot 250^2 +  r \alpha\kappa_0}{p}} 
\snorm{(\Q \Jg \JJ^r)^{\!\frac 1p} \Sigma^{-1} f(\cdot,0)}_{L^p_x}
\notag\\
&\qquad 
+2^{\frac 1p} e^{18+ \frac{2 \cdot 250^2 + r \alpha \kappa_0}{p}} (\tfrac{5\eps}{ \alpha})^{\frac{p-1}{p}}    \snorm{\JJ^{\frac rp} \Sigma^{-1} q}_{L^p_{x,\s}}
\,. 
\end{align*}
Upon passing $r\to 0^+$ (which we can do since $0 \leq \JJ \leq 3$ uniformly in $\Hdm$) and letting $p\to \infty$, we deduce that the bound  \eqref{eq:abstract:transport} remains unchanged.  The changes to the proof of Corollary~\ref{cor:transverse:bounds} are   identical, and we omit these details.

\section*{Acknowledgments} 

The work of S.S.~was in part supported by the NSF grant DMS-2007606 and the Collaborative NSF grant DMS-2307680. The work of V.V.~was in part supported by the Collaborative NSF grant DMS-2307681 and a Simons Investigator Award.


\end{document}